\documentclass[a4paper,12pt]{book} 
\usepackage{subfigure}
\usepackage{mathrsfs} 
\usepackage{tikz}
\usepackage{dirtytalk} 
\usetikzlibrary{patterns}
\usetikzlibrary{positioning}
\usetikzlibrary{calc}
\usepackage[all]{xy}
\usepackage{bookmark} 
\usepackage{amssymb, xr, multicol,textcomp}
\usepackage{imakeidx,leftidx}
\usepackage{amsthm,amsmath,url}
\usepackage{enumerate}
\definecolor{oceanboatblue}{rgb}{0.0, 0.47, 0.75}
\definecolor{persianblue}{rgb}{0.11, 0.22, 0.73}
\definecolor{napiergreen}{rgb}{0.16, 0.5, 0.0}
\usepackage{hyperref}
\hypersetup{colorlinks=true,
    linkcolor=persianblue,
    filecolor=black,      
    urlcolor=persianblue,
    citecolor=persianblue}

\usepackage[toc,page]{appendix}
\makeindex
\makeindex[name=N,title=Index of notation,columns=2,intoc]
\makeindex[name=T,title=Terminology,columns=2,intoc]

\newtheorem{theorem}{Theorem}[chapter]
\newtheorem{lemma}[theorem]{Lemma}
\newtheorem{sublemma}[theorem]{Sublemma}
\newtheorem{proposition}[theorem]{Proposition}
\newtheorem{corollary}[theorem]{Corollary}
\theoremstyle{definition}
\newtheorem{definition}[theorem]{Definition}
\newtheorem{definitions}[theorem]{Definitions}
\newtheorem{remark}[theorem]{Remark}
\newtheorem{note}[theorem]{Note}
\newtheorem{remarks}[theorem]{Remarks}
\newtheorem{notation}[theorem]{Notation}

\newtheorem{example}[theorem]{Example}
\newtheorem{examples}[theorem]{Examples}
\newtheorem{exo}[theorem]{Exercise}
\newtheorem{conventions}[theorem]{Conventions}
\newtheorem{convention}[theorem]{Convention}

\numberwithin{equation}{chapter}

\newcommand{\bsp}{\noindent {\textsc{Sketch of proof:} }} 
\newcommand{\bp}{\noindent {\textsc{Proof:} }} 
\newcommand{\eop}{\nopagebreak[4] \hspace*{\fill}{$\square$} \medskip}
\newcommand{\eopwobp}{\nopagebreak[4] \hspace*{\fill}{$\square$} \medskip}
\newcommand{\bpo}[1]{\noindent {\textsc{#1:} }} 
\newcommand{\bfig}{\begin{figure}[!htbp]}                     
                    
\let\origdoublepage\cleardoublepage
\newcommand{\clearemptydoublepage}{\clearpage{\pagestyle{empty}\origdoublepage}}
\let\cleardoublepage\clearemptydoublepage

\newcommand{\hfl}[1]{\smash{\mathop{\hbox to 10mm{\rightarrowfill}}\limits^{\textstyle
#1}}}
\newcommand{\hflv}[1]{\smash{\mathop{\hbox to 8mm{\rightarrowfill}}\limits^{\scriptstyle
#1}}}
\newcommand{\hookfl}[1]{\smash{\mathop{\hookrightarrow}\limits^{\textstyle
#1}}}
\newcommand{\norm}[1]{{\lVert #1 \rVert}}
\newcommand{\lefnorm}[1]{\left\lVert #1 \right\rVert}
\newcommand{\bignorm}[1]{\bigl\lVert #1 \bigr\rVert}
\newcommand{\Bignorm}[1]{\Bigl\lVert #1 \Bigr\rVert}
\newcommand{\suchthat}{\, : \,} 
\newcommand{\indexT}[1]{\emph{#1}\index[T]{#1}}
\newcommand{\bigvec}[1]{\overrightarrow{#1}}

\newcommand{\cardlef}[1]{\left|#1\right|}
\newcommand{\card}[1]{|#1|}
\newcommand{\cardBig}[1]{\Bigl|#1\Bigr|}
\newcommand{\cardbig}[1]{\bigl|#1\bigr|}
\newcommand{\id}{\textbf{1}}
\newcommand\Aut{\operatorname{Aut}}
\newcommand{\CA}{\mathcal{A}} 
\newcommand{\elta}{a}
\newcommand{\funca}{a}
\newcommand{\veca}{a}
\newcommand{\curvalpha}{\alpha}
\newcommand{\crita}{a}
\newcommand{\aaug}{\mathcal{A}^{\mbox{\scriptsize aug}}}
\newcommand{\aaugc}{\mathcal{A}^{\mbox{\scriptsize aug,c}}}
\newcommand{\finseta}{A}
\newcommand{\Aman}{A}
\newcommand{\amba}{A} 
\newcommand{\Asc}{\mathcal{A}} 
\newcommand{\Aavis}{\mathcal{A}}
\newcommand{\Aavisc}{\mathcal{A}^c} 
\newcommand{\Assis}{\check{\mathcal{A}}} 

\newcommand{\Asimp}{\overline{\mathcal{A}}}
\newcommand{\eltb}{b}

\newcommand{\curvbeta}{\beta}
\newcommand{\funcb}{b}
\newcommand{\critb}{b}
\newcommand{\CB}{\mathcal{B}} 
\newcommand{\ballb}{B}
\newcommand{\submb}{B} 
\newcommand{\finsetb}{B}
\newcommand{\vertexb}{{\textbf{b}}}
\newcommand{\baseb}{B}
\newcommand{\bottom}{\operatorname{bot}}
\newcommand{\ansothree}{\beta}
\newcommand{\coefgambet}{\zeta}
\newcommand{\coefgambetred}{\zeta} 
\newcommand{\blowup}[2]{\mbox{B\!$\ell$}(#1,#2)}
\newcommand{\bigblowup}[2]{\mbox{B\!$\ell$}\bigl(#1,#2\bigr)}
\newcommand{\lefblowup}[2]{\mbox{B\!$\ell$}\left(#1,#2\right)}
\newcommand{\smallblowup}[2]{ \mbox{\scriptsize B\!$\ell$}(#1,#2)} 
\newcommand{\Bdesc}{\mathcal{B}} 

\newcommand{\hcylc}{\mathcal{C}}
\newcommand{\colorc}{\textbf{c}}
\newcommand{\cl}{\mbox{c$\ell$}}
\newcommand{\confc}{c}
\newcommand{\codimc}{c}
\newcommand{\subsubc}{C}
\newcommand{\finsetc}{C}
\newcommand{\CC}{\mathbb{C}}

\newcommand{\gamtwoleg}{\gamma}
\newcommand{\plainGamma}{\Gamma}
 
\newcommand{\CD}{\mathcal{D}} 
\newcommand{\drid}{d}
\newcommand{\deltb}{d} 
\newcommand{\finsetd}{D}

\newcommand{\diskstwo}{D}
\newcommand{\Davis}{\mathcal{D}}
\newcommand{\Davisred}{\mathcal{D}} 
\newcommand{\drad}[1]{D_{#1}} 
\newcommand{\dorad}[1]{\mathring{D}_{#1}} 
\newcommand{\dimd}{d} 
\newcommand{\dimdel}{\delta}
\newcommand{\diag}{\Delta} 

\newcommand\diagon{\operatorname{diag}}
\newcommand\diagonal{\operatorname{diagonal}}
\newcommand\dimension{\operatorname{dimension}}
\newcommand{\simplexr}{\Delta^{(r)}}
\newcommand{\simplexrm}{\Delta^{(r-1)}}
\newcommand\eo{\operatorname{eo}}
\newcommand{\CE}{\mathcal{E}} 
\newcommand{\ccride}{e}
\newcommand{\edgee}{e}
\newcommand{\mape}{{\textbf{e}}}
\newcommand{\finsete}{E}
\newcommand{\eastvec}{\vec{E}}
\newcommand{\CF}{\mathcal{F}} 
\newcommand{\funcf}{f}
\newcommand{\edgef}{f}
\newcommand{\fMorse}{f}

\newcommand{\ffuncspec}{f}

\newcommand{\facee}{F}
\newcommand{\finsetf}{F}
\newcommand{\chainproppseudo}[1]{G_2(#1)}
\newcommand{\funcg}{g}
\newcommand{\metrigo}{\mathfrak{g}}
\newcommand{\edgeh}{h}
\newcommand{\handlebody}{H}
\newcommand{\handleboda}{H_a}
\newcommand{\handlebodb}{H_b}
\newcommand{\aborb}{H}
\newcommand{\aborbp}{H^{\prime}}
\newcommand{\CH}{\mathcal{H}}
\newcommand{\HH}{\mathbb{H}} 
\newcommand\holl{\operatorname{hol}}
\newcommand{\hol}[2]{\holl_{#1}(#2)}
\newcommand{\lol}[2]{\widetilde{\holl}_{#1}(#2)}
\newcommand{\Biglol}[2]{\widetilde{\holl}_{#1}\Bigl(#2\Bigr)}
\newcommand{\biglol}[2]{\widetilde{\holl}_{#1}\bigl(#2\bigr)}
\newcommand{\lolsp}[2]{\widetilde{\holl}_{#1}#2}
\newcommand{\delol}[2]{\widetilde{\delta}_{#1}(#2)}
\newcommand{\bigdelol}[2]{\widetilde{\delta}_{#1}\bigl(#2\bigr)}

\newcommand\Hom{\operatorname{Hom}}
\newcommand{\KK}{\mathbb{K}} 
\newcommand{\dimk}{k}

\newcommand{\kids}{K}
\newcommand{\CI}{\mathcal{I}} 
\newcommand\Image{\operatorname{Im}}
\newcommand\Int{\operatorname{Int}}
\newcommand{\CJ}{\mathcal{J}}
\newcommand\Ker{\operatorname{Ker}}
\newcommand{\CK}{\mathcal{K}} 
\newcommand{\CL}{\mathcal{L}} 
\newcommand{\Link}{L}
\newcommand{\leafl}{\ell}
\newcommand{\source}{\mathcal{L}}
\newcommand{\sourcetl}{\mathcal{L}}
\newcommand{\tanghcyll}{L}
\newcommand{\linhp}{\mathcal{L}_+}
\newcommand{\linhn}{\mathcal{L}_-}
\newcommand{\linh}{\mathcal{L}}
\newcommand{\CM}{\mathcal{M}} 
\newcommand{\largem}{M}
\newcommand{\manifm}{M}
\newcommand{\cmanifm}{\check{\manifm}}
\newcommand{\merid}{m}
\newcommand{\mother}{m}
\newcommand{\NN}{\mathbb{N}} 
\newcommand{\largen}{N}

\newcommand{\neighn}{N}
\newcommand{\node}{n}
\newcommand{\nodemaj}{N}
\newcommand{\nnormalized}{n}
\newcommand{\upvec}{\vec{N}}
\newcommand{\invtrivcar}{\nu}

\newcommand{\CO}{\mathcal{O}} 
\newcommand{\origino}{{\textbf{o}}}
\newcommand{\CP}{\mathcal{P}} %
\newcommand{\parentp}{\mathcal{P}} 
\newcommand{\parntppx}{{\parentp}^{\prime}_{x}}
\newcommand{\parntplambda}{{\parentp}_{\overline{X}}}
\newcommand{\parntpplambda}{{\parentp}^{\prime}_{\overline{X}}}
\newcommand{\parntpwidehat}{{\parentp}_{\overline{X}}}
\newcommand{\projassis}{\check{p}}
\newcommand{\eltp}{p}
\newcommand{\projp}{p}

\newcommand{\pbl}{p_b}
\newcommand{\propP}{P}
\newcommand{\doublepointtang}{p}
\newcommand{\doublepointset}{P}
\newcommand{\preprop}{P}
\newcommand{\anomq}{Q}
\newcommand{\QQ}{\mathbb{Q}} 
\newcommand{\eltq}{q}
\newcommand\rev{\operatorname{rev}}
\newcommand{\CR}{\mathcal{R}}
\newcommand{\RR}{\mathbb{R}} 
\newcommand{\rats}{R} 
\newcommand{\vertexr}{{\textbf{r}}}
\newcommand{\rroot}{r}
\newcommand{\rhcylr}{r}
\newcommand{\crats}{\check{\rats}} 
\newcommand{\rhorot}{\rho} 
\newcommand{\rhomap}{\rho} 
\newcommand\signature{\operatorname{signature}}
\newcommand\Supp{\operatorname{Supp}}

\newcommand{\revmaps}{s}
\newcommand{\syms}{s}
\newcommand{\vertexs}{{\textbf{s}}}
\newcommand{\CS}{\mathcal{S}}
\newcommand{\sph}{S} 
\newcommand{\cinjuptd}[2]{\check{\mathcal{S}}_{#1}(#2)}
\newcommand{\ccompuptd}[2]{\mathcal{S}_{#1}(#2)}
\newcommand{\ccompuptdk}[3]{\mathcal{S}_{#1,#3}(#2)}
\newcommand{\ccompuptdempty}[1]{\mathcal{S}_{#1}}
\newcommand{\cuptd}[2]{\overline{\mathcal{S}}_{#1}(#2)}

\newcommand{\sunits}{\mathbb{S}}
\newcommand{\sinjupdtcs}{\check{\mathcal{S}}}

\newcommand{\ccompuptdanvec}{\mathcal{S}} 
\newcommand{\cinjuptdanvec}{\check{\mathcal{S}}} 

\newcommand{\setparent}{\mathcal{T}}
\newcommand{\twistriv}{\mathcal{T}} 
\newcommand{\CT}{\mathcal{T}}
\newcommand{\vecspt}{T}
\newcommand{\tang}{T}
\newcommand{\topt}{T}
\newcommand{\tanghcylt}{T}

\newcommand{\taust}{\tau_s}
\newcommand{\toptang}{\operatorname{top}}
\newcommand{\partau}{\tau} 
\newcommand{\ST}{U} 

\newcommand{\prodedge}{U}
\newcommand{\finsetu}{U}
\newcommand{\univu}{U}
\newcommand{\eltu}{u} 
\newcommand{\varedge}{u}
\newcommand{\CV}{\mathcal{V}} 
\newcommand{\eltv}{v}
\newcommand{\vertsetv}{V}
\newcommand{\finsetv}{V}
\newcommand{\confw}{w}
\newcommand{\qvarM}{W}
\newcommand{\vecw}{W}
\newcommand{\finsetw}{W}
\newcommand{\pointx}{x}
\newcommand{\confx}{x}
\newcommand{\eltx}{x}
\newcommand{\cutx}{X}
\newcommand{\setparamx}{X}
\newcommand{\vecx}{X}
\newcommand{\cvarM}{X}
\newcommand{\fvarM}{X}
\newcommand{\setx}{X}

\newcommand{\einftyxi}{\Xi}
\newcommand{\Zgen}{Y}
\newcommand{\vecy}{Y}
\newcommand{\svarM}{Y}
\newcommand{\finsety}{Y}
\newcommand{\vecZorth}{Z}
\newcommand{\finsetz}{Z}
\newcommand{\confy}{y} 
\newcommand{\edgez}{z}
\newcommand{\tvarM}{Z}
\newcommand{\ZZ}{\mathbb{Z}}
\newcommand{\Zinv}{Z}
\newcommand{\zinv}{z}
\newcommand{\Zinvlink}{\check{Z}}
\newcommand{\zinvlink}{\check{z}}
\newcommand{\Zinvuf}{\mathcal{Z}}
\newcommand{\Zinvufmodis}{\overline{\mathcal{Z}}} 
\newcommand{\zinvuf}{\mathfrak{z}}
\newcommand{\Zinvufrf}{\mathcal{Z}^f}

\newcommand{\Zinvufrfneg}{\mathcal{Z}^f\!}
\newcommand{\Zinvlinkufrfneg}{\check{\mathcal{Z}}^f\!}
\newcommand{\zinvufrfneg}{\mathfrak{z}^f\!}
\newcommand{\Zinvlinkuf}{\check{\mathcal{Z}}}
\newcommand{\Zinvlinkufmodis}{\overline{\check{\mathcal{Z}}}}

\newcommand{\Zinvlinkufov}{\overline{\overline{\mathcal{Z}}}} 

\newcommand{\bzaug}{Z^{\mbox{\scriptsize aug}}}
\newcommand{\zaug}{z^{\mbox{\scriptsize aug}}}

\newcommand{\cltetra}{\begin{tikzpicture}\useasboundingbox (-.1,0) rectangle (.7,.5);
\draw (0,0) -- (.6,0) -- (.3,.5) -- (0,0);
\draw (.6,0) -- (.3,.2) -- (.3,.5) (0,0) -- (.3,.2);
\fill (0,0) circle (1.5pt) (.6,0) circle (1.5pt) (.3,.5) circle (1.5pt) (.3,.2) circle (1.5pt);
\end{tikzpicture}}
\newcommand{\tetrasmall}{\begin{tikzpicture}[scale=.5]
\draw (0,0) -- (.6,0) -- (.3,.5) -- (0,0);
\draw (.6,0) -- (.3,.2) -- (.3,.5) (0,0) -- (.3,.2);
\fill (0,0) circle (1.5pt) (.6,0) circle (1.5pt) (.3,.5) circle (1.5pt) (.3,.2) circle (1.5pt);
\end{tikzpicture}}
\newcommand{\smalltripod}[3]{\begin{tikzpicture} \useasboundingbox (-.2,-.4) rectangle (.9,0);
\begin{scope}[yshift=-.3cm]
\draw [thick] (.4,0) -- (0,0) -- (.4,.3) (0,0) -- (.4,-.3);
\draw (.4,.3) node[right]{\footnotesize #3} (.4,0) node[right]{\footnotesize #2} (.4,-.3) node[right]{\footnotesize #1};
\fill (0,0) circle (1.5pt);
\end{scope}
\end{tikzpicture}}
\newcommand{\bigtripod}[3]{\begin{tikzpicture} [baseline=-.1cm] \useasboundingbox (-.1,-.7) rectangle (1.1,.9);
\draw [thick] (.4,0) -- (0,0) -- (.4,.55) (0,0) -- (.4,-.55);
\draw (.35,.65) node[right]{\footnotesize #3} (.4,.1) node[right]{\footnotesize #2} (.35,-.45) node[right]{\footnotesize #1};
\fill (0,0) circle (1.5pt);
\end{tikzpicture}}
\newcommand{\bigtripodrev}[3]{\begin{tikzpicture} [baseline=-.1cm] \useasboundingbox (-1.1,-.7) rectangle (.1,.9);
\draw [thick] (-.4,0) -- (0,0) -- (-.4,.55) (0,0) -- (-.4,-.55);
\draw (-.35,.65) node[left]{\footnotesize #3} (-.4,.1) node[left]{\footnotesize #2} (-.35,-.45) node[left]{\footnotesize #1};
\fill (0,0) circle (1.5pt);
\end{tikzpicture}}
\newcommand{\nccirc}{\begin{tikzpicture}
\useasboundingbox (-.4,.1) rectangle (.4,.55);
\draw [->] (-.2,0) -- (.2,.4);
\draw [->,draw=white,double=black,very thick]  (.2,0) -- (-.2,.4);
\draw [->]  (.2,0) -- (-.2,.4);
\draw [dash pattern=on 2pt off 2pt] (0,.2) circle (.283);
\end{tikzpicture}}
\newcommand{\pccirc}{\begin{tikzpicture} \useasboundingbox (-.4,.1) rectangle (.4,.55);
\draw [->] (.2,0) -- (-.2,.4);
\draw [->,draw=white,double=black,very thick]  (-.2,0) -- (.2,.4);
\draw [->]  (-.2,0) -- (.2,.4);
\draw [dash pattern=on 2pt off 2pt] (0,.2) circle (.283);
\end{tikzpicture}}
\newcommand{\doublepcirc}{\begin{tikzpicture} \useasboundingbox (-.4,.1) rectangle (.4,.55);
\draw [->] (.2,0) -- (-.2,.4);
\draw [->]  (-.2,0) -- (.2,.4);
\draw [dash pattern=on 2pt off 2pt] (0,.2) circle (.283);
\fill (0,.2) circle (1.5pt);
\end{tikzpicture}}

\newcommand{\doublecupnotwistsource}{\begin{tikzpicture} [baseline=-.05cm] \useasboundingbox (-.7,-.5) rectangle (.4,.9);
\draw [<-,dashed]  (0,.8) to (0,0);
\draw [<-,dashed] (.3,.8) to (.3,0);
\draw [dashed] (-.6,.8) -- (-.6,0) .. controls (-.6,-.6) and (.3,-.6) .. (.3,0);
\draw [dashed] (-.3,.8) -- (-.3,0) .. controls (-.3,-.3) and (0,-.3) .. (0,0);
\draw (.3,.4) .. controls (-.1,.1) and (-.3,.4) .. (0,.4);
\fill (0,.4) circle (1.5pt) (.3,.4) circle (1.5pt);
\end{tikzpicture}}
\newcommand{\doublecuptwistsing}{\begin{tikzpicture} [baseline=-.05cm] \useasboundingbox (-.7,-.6) rectangle (.4,.8);
\draw [out=90,in=-90] (.3,0) to (0,.4);
\draw [out=90,in=-90] (0,0) to (.3,.4);
\draw [->, out=90,in=-90] (.3,.4) to (0,.8);
\draw [out=90,in=-90,draw=white,double=black,very thick] (0,.4) to (.3,.8);
\draw [->, out=90,in=-90] (0,.4) to (.3,.8);
\draw (-.6,.8) -- (-.6,0) .. controls (-.6,-.6) and (.3,-.6) .. (.3,0);
\draw (-.3,.8) -- (-.3,0) .. controls (-.3,-.3) and (0,-.3) .. (0,0);
\fill (.15,.2) circle (1.5pt);
\end{tikzpicture}}
\newcommand{\twostraightstrands}{\begin{tikzpicture} [baseline=.2cm] \useasboundingbox (-.1,-.1) rectangle (.4,.9);
\draw [->] (0,0) to (0,.8);
\draw [->] (.3,0) to (.3,.8);
\end{tikzpicture}}
\newcommand{\twostraightstrandspetit}{\begin{tikzpicture} [baseline=0cm] \useasboundingbox (-.1,0) rectangle (.25,.3);
\draw [->] (0,0) to (0,.3);
\draw [->] (.15,0) to (.15,.3);
\end{tikzpicture}}

\newcommand{\twostraightstrandsnonorlong}{\begin{tikzpicture} [baseline=-.05cm] \useasboundingbox (-.1,-.7) rectangle (.3,.8);
\draw [-] (0,-.7) to (0,.8);
\draw [-] (.2,-.7) to (.2,.8);
\end{tikzpicture}}

\newcommand{\fulltwist}{\begin{tikzpicture} [baseline=.2cm] \useasboundingbox (-.1,-.1) rectangle (.4,.9);
\draw [out=90,in=-90] (.3,0) to (0,.4);
\draw [out=90,in=-90,draw=white,double=black,very thick] (0,0) to (.3,.4);
\draw [out=90,in=-90,->] (.3,.4) to (0,.8);
\draw [out=90,in=-90,->,draw=white,double=black,very thick] (0,.4) to (.3,.8);
\draw [out=90,in=-90,->] (0,.4) to (.3,.8);
\end{tikzpicture}}
\newcommand{\fulltwistfake}{\begin{tikzpicture} [baseline=.2cm] \useasboundingbox (-.1,-.1) rectangle (.4,.9);
\draw [out=90,in=-90] (0,0) to (.3,.4);
\draw [out=90,in=-90,draw=white,double=black,very thick] (.3,0) to (0,.4);
\draw [out=90,in=-90,->] (.3,.4) to (0,.8);
\draw [out=90,in=-90,->,draw=white,double=black,very thick] (0,.4) to (.3,.8);
\draw [out=90,in=-90,->] (0,.4) to (.3,.8);
\end{tikzpicture}}
\newcommand{\fulltwistsing}{\begin{tikzpicture} [baseline=.2cm] \useasboundingbox (-.1,-.1) rectangle (.4,.9);
\draw [out=90,in=-90] (.3,0) to (0,.4);
\draw [out=90,in=-90] (0,0) to (.3,.4);
\draw [out=90,in=-90,->] (.3,.4) to (0,.8);
\draw [out=90,in=-90,->,draw=white,double=black,very thick] (0,.4) to (.3,.8);
\draw [out=90,in=-90,->] (0,.4) to (.3,.8);
\fill (.15,.2) circle (1.5pt);
\end{tikzpicture}}
\newcommand{\fulltwistonerev}{\begin{tikzpicture} [baseline=.2cm] \useasboundingbox (-.1,-.1) rectangle (.4,.9);
\draw [out=90,in=-90,<-] (.3,0) to (0,.4);
\draw [out=90,in=-90,draw=white,double=black,very thick] (0,0) to (.3,.4);
\draw [out=90,in=-90,->] (.3,.4) to (0,.8);
\draw [out=90,in=-90,->,draw=white,double=black,very thick] (0,.4) to (.3,.8);
\draw [out=90,in=-90] (0,.4) to (.3,.8);
\end{tikzpicture}}
\newcommand{\onechordtwosor}{\begin{tikzpicture} [baseline=.2cm] \useasboundingbox (-.1,-.1) rectangle (.4,.9);
\draw [->, dashed] (0,0) to (0,.8);
\draw [->, dashed] (.3,0) to (.3,.8); 
\draw (.3,.4) .. controls (-.1,.1) and (-.3,.4) .. (0,.4);
\fill (0,.4) circle (1.5pt) (.3,.4) circle (1.5pt);
\end{tikzpicture}}
\newcommand{\twochordtwosor}{\begin{tikzpicture} [baseline=.2cm] \useasboundingbox (-.1,-.1) rectangle (.4,.9);
\draw [->, dashed] (0,0) to (0,.8);
\draw [->, dashed] (.3,0) to (.3,.8); 
\begin{scope}[yshift=-.15cm]
\draw (.3,.4) .. controls (-.1,.1) and (-.3,.4) .. (0,.4);
\fill (0,.4) circle (1.5pt) (.3,.4) circle (1.5pt);
\end{scope}
\begin{scope}[yshift=.15cm]
\draw (.3,.4) .. controls (-.1,.1) and (-.3,.4) .. (0,.4);
\fill (0,.4) circle (1.5pt) (.3,.4) circle (1.5pt);
\end{scope}
\end{tikzpicture}}
\newcommand{\threechordtwosor}{\begin{tikzpicture} [baseline=.2cm] \useasboundingbox (-.1,-.1) rectangle (.4,.9);
\draw [->, dashed] (0,0) to (0,.8);
\draw [->, dashed] (.3,0) to (.3,.8); 
\begin{scope}[yshift=.25cm]
\draw (.3,.4) .. controls (-.1,.1) and (-.3,.4) .. (0,.4);
\fill (0,.4) circle (1.5pt) (.3,.4) circle (1.5pt);
\end{scope}
\begin{scope}[yshift=.05cm]
\draw (.3,.4) .. controls (-.1,.1) and (-.3,.4) .. (0,.4);
\fill (0,.4) circle (1.5pt) (.3,.4) circle (1.5pt);
\end{scope}
\begin{scope}[yshift=-.15cm]
\draw (.3,.4) .. controls (-.1,.1) and (-.3,.4) .. (0,.4);
\fill (0,.4) circle (1.5pt) (.3,.4) circle (1.5pt);
\end{scope}
\end{tikzpicture}}
\newcommand{\onechordtwosoror}{\begin{tikzpicture} [baseline=.2cm] \useasboundingbox (-.1,-.1) rectangle (.4,.9);
\draw [->, dashed] (0,0) to (0,.8);
\draw [->, dashed] (.3,0) to (.3,.8); 
\draw [->] (0,.4) .. controls (-.3,.4) and (-.1,.2) .. (.1,.2);
\draw [out=0, in= -150]  (.1,.2) to (.3,.4);
\fill (0,.4) circle (1.5pt) (.3,.4) circle (1.5pt);
\end{tikzpicture}}
\newcommand{\onechordtwosororrev}{\begin{tikzpicture} [baseline=.2cm] \useasboundingbox (-.1,-.1) rectangle (.4,.9);
\draw [->, dashed] (0,0) to (0,.8);
\draw [-<] (0,.4) .. controls (-.3,.4) and (-.1,.2) .. (.1,.2);
\draw [out=0, in= -150]  (.1,.2) to (.3,.4);
\draw [->, dashed] (.3,0) to (.3,.8); 
\fill (0,.4) circle (1.5pt) (.3,.4) circle (1.5pt);
\end{tikzpicture}}
\newcommand{\onechordonverticalinterval}{
\begin{tikzpicture} \useasboundingbox (-.2,-.2) rectangle (.2,.2);
\draw (0,-.15) .. controls (-.25,-.15) and (-.25,.15) .. (0,.15);
\draw [dash pattern=on 2pt off 2pt]  (0,-.3) -- (0,.3);
\fill (0,-.15) circle (1.5pt) (0,.15) circle (1.5pt);
\end{tikzpicture}}
\newcommand{\onechordwithboxonverticalinterval}{
\begin{tikzpicture} [baseline=-.1cm] \useasboundingbox (-.6,-.2) rectangle (.2,.6);
\draw (0,-.2) .. controls (-.25,-.2) .. (-.25,-.1);
\draw (0,.35) .. controls (-.25,.35) .. (-.25,.25);
\draw (-.5,-.1) rectangle (-.05,.25) (-.25,.05) node{\scriptsize $\alpha$};
\draw [->, dashed]  (0,-.3) -- (0,.5);
\fill (0,-.2) circle (1.5pt) (0,.35) circle (1.5pt);
\end{tikzpicture}}
\newcommand{\onechordwithboxtwosor}{\begin{tikzpicture} [baseline=.2cm] \useasboundingbox (-.3,-.1) rectangle (.75,.9);
\draw [->, dashed] (0,0) to (0,.8);
\draw (.1,.4) .. controls (-.2,.4) and (-.2,.6) .. (0,.6);
\draw (.1,.2) rectangle (.55,.6) (.55,.4) -- (.65,.4) (.35,.4) node{\scriptsize $\alpha$};
\draw [->, dashed] (.65,0) to (.65,.8);
\fill (0,.6) circle (1.5pt) (.65,.4) circle (1.5pt);
\end{tikzpicture}}
\newcommand{\dasheddownw}{\begin{tikzpicture} \useasboundingbox (-.1,-.1) rectangle (.1,.2);
\draw [->,dash pattern=on 2pt off 2pt]  (0,.2) -- (0,-.2);
\end{tikzpicture}}
\newcommand{\dashedupw}{\begin{tikzpicture} \useasboundingbox (-.1,-.1) rectangle (.1,.2);
\draw [<-,dash pattern=on 2pt off 2pt]  (0,.2) -- (0,-.2);
\end{tikzpicture}}
\newcommand{\dashedcircleorrrevseul}{\begin{tikzpicture} \useasboundingbox (-.3,-.1) rectangle (.3,.2);
\draw [-<,dash pattern=on 2pt off 2pt]  (0,.2) arc (-270:90:.2);
\end{tikzpicture}}
\newcommand{\diagcross}{\begin{tikzpicture} \useasboundingbox (-.15,-.15) rectangle (.15,.15);
\draw [dash pattern=on 2pt off 2pt] (.15,0) arc (0:360:.15);
\draw (45:.15) -- (-135:.15);
\draw (135:.15) -- (-45:.15);
\fill (45:.15) circle (1pt) (135:.15) circle (1pt) (-135:.15) circle (1pt) (-45:.15) circle (1pt);
\end{tikzpicture}}
\newcommand{\diagcrossmini}{\begin{tikzpicture} [scale=.7] \useasboundingbox (-.25,-.14) rectangle (.2,.15);
\draw [dash pattern=on 2pt off 2pt] (.15,0) arc (0:360:.15);
\draw (45:.15) -- (-135:.15);
\draw (135:.15) -- (-45:.15);
\fill (45:.15) circle (1pt) (135:.15) circle (1pt) (-135:.15) circle (1pt) (-45:.15) circle (1pt);
\end{tikzpicture}}
\newcommand{\diagthetthet}{\begin{tikzpicture} \useasboundingbox (-.15,-.15) rectangle (.15,.15);
\draw [dash pattern=on 2pt off 2pt] (.15,0) arc (0:360:.15);
\draw (60:.15) -- (-60:.15);
\draw (120:.15) -- (-120:.15);
\fill (60:.15) circle (1pt) (120:.15) circle (1pt) (-120:.15) circle (1pt) (-60:.15) circle (1pt);
\end{tikzpicture}}
\newcommand{\diagtripod}{\begin{tikzpicture} \useasboundingbox (-.15,-.15) rectangle (.15,.15);
\draw [dash pattern=on 2pt off 2pt] (.15,0) arc (0:360:.15);
\draw (60:.15) -- (0,0);
\draw (-60:.15) -- (0,0);
\draw (180:.15) -- (0,0);
\fill (60:.15) circle (1pt) (-60:.15) circle (1pt) (180:.15) circle (1pt) (0,0) circle (1pt);
\end{tikzpicture}}
\newcommand{\diagtripodmini}{\begin{tikzpicture}[scale=.7] \useasboundingbox (-.25,-.14) rectangle (.2,.15);
\draw [dash pattern=on 2pt off 2pt] (.15,0) arc (0:360:.15);
\draw (60:.15) -- (0,0);
\draw (-60:.15) -- (0,0);
\draw (180:.15) -- (0,0);
\fill (60:.15) circle (1pt) (-60:.15) circle (1pt) (180:.15) circle (1pt) (0,0) circle (1pt);
\end{tikzpicture}}
\newcommand{\diagwheeltwo}{\begin{tikzpicture} \useasboundingbox (-.23,-.15) rectangle (.23,.15);
\begin{scope}[xscale=1.5]\draw [dash pattern=on 2pt off 2pt] (.15,0) arc (0:360:.15);
\draw (.15,0) -- (.075,0);
\draw (-.15,0) -- (-.075,0);
\draw (-.075,0) .. controls (0,-.1) .. (.075,0) (-.075,0) .. controls (0,.1) .. (.075,0);
\fill (.15,0) circle (1pt) (-.15,0) circle (1pt) (.07,0) circle (1pt) (-.07,0) circle (1pt);
\end{scope}
\end{tikzpicture}}
\newcommand{\diagv}[7]{
\begin{tikzpicture} \useasboundingbox (-1,-.5) rectangle (1,.5);
\draw [->,dash pattern=on 2pt off 2pt] (.5,0) arc (0:360:.5);
\draw (60:.5) -- (-60:.5) node[right]{\scriptsize $#7$};
\draw (60:.5) node[right]{\scriptsize $#5$} (245:.3) node[above]{\tiny $#4$} (0:.3)  node[left]{\tiny $#2$};
\draw (-60:.5) -- (180:.5) node[left]{\scriptsize $#6$};
\draw [-#1] (60:.5) -- (0:.25);
\draw [-#3] (180:.5) -- (240:.25);
\fill (60:.5) circle (1.6pt) (180:.5) circle (1.6pt) (-60:.5) circle (1.6pt);
\end{tikzpicture}}
\newcommand{\diagvhaut}[7]{
\begin{tikzpicture} [baseline=-.3cm] 
\draw [->,dash pattern=on 2pt off 2pt] (.5,0) arc (0:360:.5);
\draw (60:.5) -- (-60:.5) node[right]{\scriptsize $#7$};
\draw (60:.5) node[right]{\scriptsize $#5$} (-.2,0) node{\tiny $#4$} (.1,.15)  node{\tiny $#2$};
\draw (-60:.5) -- (180:.5) node[left]{\scriptsize $#6$};
\draw [-#1] (60:.5) -- (0:.25);
\draw [-#3] (180:.5) -- (240:.25);
\fill (60:.5) circle (1.6pt) (180:.5) circle (1.6pt) (-60:.5) circle (1.6pt);
\end{tikzpicture}}
\newcommand{\diagvv}[4]{
\begin{tikzpicture} \useasboundingbox (-1,-.5) rectangle (1,.5);
\draw [->,dash pattern=on 2pt off 2pt] (.5,0) arc (0:360:.5);
\draw (60:.5) -- (-60:.5);
\draw (-.42,-.05) node[right]{\tiny $#4$} (.32,.27)  node[left]{\tiny $#2$};
\draw (-60:.5) -- (180:.5);
\draw [-#1] (60:.5) -- (.25,.27);
\draw [-#3] (180:.5) -- (240:.25);
\fill (60:.5) circle (1.5pt) (180:.5) circle (1.5pt) (-60:.5) circle (1.5pt);
\end{tikzpicture}}
\newcommand{\diagtripoded}[4]{
\begin{tikzpicture} \useasboundingbox (-.7,-.4) rectangle (.7,.4);
\draw [->,dash pattern=on 2pt off 2pt] (.5,0) arc (0:360:.5);
\draw (0,0) -- (-60:.5);
\draw (0,0) -- (60:.5);
\draw (-30:.32) node{\tiny $1$}  (90:.32) node{\tiny $#2$} (210:.32) node{\tiny $#4$} (0,0) -- (180:.5);
\draw [->] (-60:.5) -- (-60:.25);
\draw [-#1] (60:.5) -- (60:.25);
\draw [-#3] (180:.5) -- (180:.25);
\fill (60:.5) circle (1.6pt) (180:.5) circle (1.6pt) (-60:.5) circle (1.6pt) (0,0) circle (1.6pt);
\end{tikzpicture}}
\newcommand{\nochordonesorshort}{\begin{tikzpicture} [baseline=0cm] \useasboundingbox (-.1,.05) rectangle (.1,.35);
\draw [->, dashed,dash pattern=on 2pt off 2pt] (0,0) to (0,.3);
\end{tikzpicture}}
\newcommand{\nochordonesor}{\begin{tikzpicture} [baseline=.2cm] \useasboundingbox (-.1,-.1) rectangle (.1,.9);
\draw [->, dashed,dash pattern=on 2pt off 2pt] (0,0) to (0,.8);
\end{tikzpicture}}
\newcommand{\nochordtwosor}{\begin{tikzpicture} [baseline=.2cm] \useasboundingbox (-.1,-.1) rectangle (.4,.9);
\draw [->, dashed] (0,0) to (0,.8);
\draw [->, dashed] (.3,0) to (.3,.8);
\end{tikzpicture}}
\newcommand{\nochordtwosorpetit}{\begin{tikzpicture} [baseline=0cm] 
\useasboundingbox (-.1,0) rectangle (.25,.3);
\draw [->, dashed] (0,0) to (0,.3);
\draw [->, dashed] (.15,0) to (.15,.3);
\end{tikzpicture}}
\newcommand{\fulltwistnonor}{\begin{tikzpicture} [baseline=.2cm] \useasboundingbox (-.1,-.1) rectangle (.4,.9);
\draw [out=90,in=-90] (.3,0) to (0,.4);
\draw [out=90,in=-90,draw=white,double=black,very thick] (0,0) to (.3,.4);
\draw [out=90,in=-90] (.3,.4) to (0,.8);
\draw [out=90,in=-90,draw=white,double=black,very thick] (0,.4) to (.3,.8);
\draw [out=90,in=-90] (0,.4) to (.3,.8);
\end{tikzpicture}}
\newcommand{\onechordtwos}{\begin{tikzpicture} [baseline=.2cm] \useasboundingbox (-.3,-.1) rectangle (.4,.9);
\draw [dashed] (0,0) to (0,.8);
\draw [dashed] (.3,0) to (.3,.8); 
\draw (.3,.4) .. controls (-.1,.1) and (-.3,.4) .. (0,.4);
\fill (0,.4) circle (1.5pt) (.3,.4) circle (1.5pt);
\end{tikzpicture}}
\newcommand{\symvogone}{\begin{tikzpicture} [baseline=.5cm] \useasboundingbox (-.9,.1) rectangle (.9,1.3);
\draw[very thick,gray] (0,1.2) .. controls (-.8,1.2) .. (-.8,1) .. controls (-.8,.8) .. (0,.8)
 .. controls (.8,.8) .. (.8,1) .. controls (.8,1.2) .. (0,1.2);
\draw (-.4,.9) -- (-.4,.2) (.4,.9) -- (.4,.2) (0,.9) .. controls (0,.5) .. (.4,.5) (-.55,.2) node{\scriptsize 1} (.55,.2) node{\scriptsize 2};
\fill (.4,.5) circle (1.5pt) (-.4,.2) circle (1.5pt) (.4,.2) circle (1.5pt);
\end{tikzpicture}}
\newcommand{\symvogtwo}{\begin{tikzpicture} [baseline=.5cm] \useasboundingbox (-.9,.1) rectangle (.9,1.3);
\draw[very thick,gray] (0,1.2) .. controls (-.8,1.2) .. (-.8,1) .. controls (-.8,.8) .. (0,.8)
 .. controls (.8,.8) .. (.8,1) .. controls (.8,1.2) .. (0,1.2);
\draw (-.4,.9) -- (-.4,.2) (.4,.9) -- (.4,.2) (0,.9) .. controls (0,.5) .. (-.4,.5) (-.55,.2) node{\scriptsize 1} (.55,.2) node{\scriptsize 2};
\fill (-.4,.5) circle (1.5pt) (-.4,.2) circle (1.5pt) (.4,.2) circle (1.5pt);
\end{tikzpicture}}
\newcommand{\symvogonesym}{\begin{tikzpicture} [baseline=.5cm] \useasboundingbox (-.9,.1) rectangle (.9,1.3);
\draw[very thick,gray] (0,1.2) .. controls (-.8,1.2) .. (-.8,1) .. controls (-.8,.8) .. (0,.8)
 .. controls (.8,.8) .. (.8,1) .. controls (.8,1.2) .. (0,1.2);
\draw (-.4,.9) -- (-.4,.2) (.4,.9) -- (.4,.2) (0,.9) .. controls (0,.5) .. (.4,.5) (-.55,.2) node{\scriptsize 2} (.55,.2) node{\scriptsize 1};
\fill (.4,.5) circle (1.5pt) (-.4,.2) circle (1.5pt) (.4,.2) circle (1.5pt);
\end{tikzpicture}}
\newcommand{\bulletpetit}{\begin{tikzpicture} [baseline=-.1cm]
\useasboundingbox (-.05,-.05) rectangle (.05,.05);
\fill (0,0) circle (1pt);
\end{tikzpicture}}
\newcommand{\bulletmoyen}{\begin{tikzpicture} [baseline=-.1cm]
\useasboundingbox (-.05,-.05) rectangle (.05,.05);
\fill (0,0) circle (1.5pt);
\end{tikzpicture}}
\newcommand{\ncpetit}{\begin{tikzpicture}
\useasboundingbox (-.25,.05) rectangle (.25,.35);
\draw [->] (-.15,0) -- (.15,.3);
\fill[white] (0,.15) circle (.07);
\draw [->]  (.15,0) -- (-.15,.3);
\end{tikzpicture}}
\newcommand{\pcpetit}{\begin{tikzpicture} \useasboundingbox (-.25,.05) rectangle (.25,.35);
\draw [->] (.15,0) -- (-.15,.3);
\fill[white] (0,.15) circle (.07);
\draw [->]  (-.15,0) -- (.15,.3);
\end{tikzpicture}}
\newcommand{\ncpetittexte}{\begin{tikzpicture}
\useasboundingbox (-.18,.05) rectangle (.33,.35);
\draw [->] (-.15,0) -- (.15,.3);
\fill[white] (0,.15) circle (.07);
\draw [->]  (.15,0) -- (-.15,.3);
\end{tikzpicture}}
\newcommand{\pcpetittexte}{\begin{tikzpicture} \useasboundingbox (-.18,0.05) rectangle (.33,.35);
\draw [->] (.15,0) -- (-.15,.3);
\fill[white] (0,.15) circle (.07);
\draw [->]  (-.15,0) -- (.15,.3);
\end{tikzpicture}}
\newcommand{\huitncpetit}{\begin{tikzpicture}
\useasboundingbox (-.35,0) rectangle (.35,.3);
\draw [->] (-.15,0) -- (.15,.3);
\fill[white] (0,.15) circle (.07);
\draw [->]  (.15,0) -- (-.15,.3);
\draw (.15,.3) .. controls (.25,.4) and (.3,.25) .. (.3,.15) .. controls (.3,.05) and (.25,-.1) .. (.15,0);
\draw (-.15,.3) .. controls (-.25,.4) and (-.3,.25) .. (-.3,.15) .. controls (-.3,.05) and (-.25,-.1) .. (-.15,0);
\end{tikzpicture}}
\newcommand{\huitpcpetit}{\begin{tikzpicture} \useasboundingbox (-.35,0) rectangle (.35,.3);
\draw [->] (.15,0) -- (-.15,.3);
\fill[white] (0,.15) circle (.07);
\draw [->]  (-.15,0) -- (.15,.3);
\draw (.15,.3) .. controls (.25,.4) and (.3,.25) .. (.3,.15) .. controls (.3,.05) and (.25,-.1) .. (.15,0);
\draw (-.15,.3) .. controls (-.25,.4) and (-.3,.25) .. (-.3,.15) .. controls (-.3,.05) and (-.25,-.1) .. (-.15,0);
\end{tikzpicture}}
\newcommand{\nc}{\begin{tikzpicture}
\useasboundingbox (-.3,.1) rectangle (.35,.55);
\draw [->] (-.2,0) -- (.2,.4);
\draw [->,draw=white,double=black,very thick]  (.2,0) -- (-.2,.4);
\draw [->]  (.2,0) -- (-.2,.4);
\end{tikzpicture}}
\newcommand{\pc}{\begin{tikzpicture} \useasboundingbox (-.3,.1) rectangle (.35,.55);
\draw [->] (.2,0) -- (-.2,.4);
\draw [->,draw=white,double=black,very thick]  (-.2,0) -- (.2,.4);
\draw [->]  (-.2,0) -- (.2,.4);
\end{tikzpicture}}
\newcommand{\pcbg}{\begin{tikzpicture}
\useasboundingbox (-.6,.1) rectangle (.6,.55);
\draw [->] (.2,0) -- (-.2,.4) node[left]{\tiny $J$};
\draw [->,draw=white,double=gray,very thick]  (-.2,0) -- (.2,.4) node[right]{\tiny $K$};
\draw [->,draw=gray]  (-.2,0) -- (.2,.4);
\end{tikzpicture}} 
\newcommand{\ncbg}{\begin{tikzpicture}
\useasboundingbox (-.6,.1) rectangle (.6,.55);
\draw [->] (-.2,0) -- (.2,.4) node[right]{\tiny $J$};
\draw [->,draw=white,double=gray,very thick]  (.2,0) -- (-.2,.4) node[left]{\tiny $K$};
\draw [->,draw=gray]  (.2,0) -- (-.2,.4);
\end{tikzpicture}}

\newcommand{\pcortrig}{\begin{tikzpicture} \useasboundingbox (-.5,.1) rectangle (.5,.55);
\draw [->] (.2,0) -- (-.2,.4);
\draw [->,draw=white,double=black,very thick]  (-.2,0) -- (.2,.4);
\draw [->]  (-.2,0) -- (.2,.4);
\draw [draw=gray,->] (.1,.45) .. controls (0,.5) .. (-.1,.45);
\end{tikzpicture}}
\newcommand{\pcbgkj}{\begin{tikzpicture}
\useasboundingbox (-.6,.1) rectangle (.6,.55);
\draw [->,draw=gray] (.2,0) -- (-.2,.4) node[left]{\tiny $K$};
\draw [->,draw=white,double=black,very thick]  (-.2,0) -- (.2,.4) node[right]{\tiny $J$};
\draw [->]  (-.2,0) -- (.2,.4);
\end{tikzpicture}}
\newcommand{\ncbgkj}{\begin{tikzpicture}
\useasboundingbox (-.6,.1) rectangle (.6,.55);
\draw [->,draw=gray] (-.2,0) -- (.2,.4) node[right]{\tiny $K$};
\draw [->,draw=white,double=black,very thick]  (.2,0) -- (-.2,.4) node[left]{\tiny $J$};
\draw [->]  (.2,0) -- (-.2,.4);
\end{tikzpicture}}
\newcommand{\singtrefoil}{\begin{tikzpicture}[scale=.5] \useasboundingbox (-1,-.3) rectangle (1,.2);
\draw [out=120,in=60] (30:.7) to (150:.3);
\draw [out=240,in=180] (150:.3) to (-90:.7);
\draw [out=0,in=-60] (-90:.7) to (30:.3);
\draw [out=120,in=60] (30:.3) to (150:.7);
\draw [out=240,in=180] (150:.7) to (-90:.3);
\draw [out=0,in=-60,->] (-90:.3) to (30:.7);
\fill (-30:.39) circle (2.5pt) (90:.39) circle (2.5pt) (-150:.39) circle (2.5pt);
\end{tikzpicture}}
\newcommand{\doublep}{\begin{tikzpicture} \useasboundingbox (-.3,.1) rectangle (.3,.55);
\draw [->] (.2,0) -- (-.2,.4);
\draw [->]  (-.2,0) -- (.2,.4);
\fill  (0,.2) circle (1.5pt);
\end{tikzpicture}}
\newcommand{\doubleppetit}{\begin{tikzpicture} \useasboundingbox (-.25,0.05) rectangle (.25,.35);
\draw [->] (.15,0) -- (-.15,.3);
\draw [->]  (-.15,0) -- (.15,.3);
\fill  (0,.15) circle (1.5pt);
\end{tikzpicture}}
\newcommand{\doubleppetittexte}{\begin{tikzpicture} \useasboundingbox (-.18,0.05) rectangle (.33,.35);
\draw [->] (.15,0) -- (-.15,.3);
\draw [->]  (-.15,0) -- (.15,.3);
\fill  (0,.15) circle (1.5pt);
\end{tikzpicture}}

\newcommand{\doublepdashchord}{\begin{tikzpicture} \useasboundingbox (-.5,.1) rectangle (.5,.55);
\draw [->, dash pattern=on 2pt off 2pt] (.2,-.1) .. controls (.2,0) and (-.2,.1) .. (-.2,.2) -- (-.2,.5);
\draw [->, dash pattern=on 2pt off 2pt]  (-.2,-.1) .. controls (-.2,0) and (.2,.1) .. (.2,.2) -- (.2,.5);
\draw (.2,.35) .. controls (0,.1) and (-.45,.35) .. (-.2,.35);
\fill (-.2,.35) circle (1.5pt) (.2,.35) circle (1.5pt);
\end{tikzpicture}}

\newcommand{\righttrefoil}{\begin{tikzpicture}
[scale=.5] \useasboundingbox (-1,-.3) rectangle (1,.2);
\draw [out=120,in=60] (30:.7) to (150:.3);
\draw [out=0,in=-60] (-90:.7) to (30:.3);
\draw [out=240,in=180] (150:.7) to (-90:.3);
\draw [out=240,in=180,draw=white,double=black,very thick] (150:.3) to (-90:.7);
\draw [out=120,in=60,draw=white,double=black,very thick] (30:.3) to (150:.7);
\draw [out=0,in=-60,draw=white,double=black,very thick] (-90:.3) to (30:.7);
\draw [out=0,in=-60,->] (-90:.3) to (30:.7);
\end{tikzpicture}}
\newcommand{\threexchord}{\begin{tikzpicture} \useasboundingbox (-.4,-.2) rectangle (.4,.5);
\draw [->,dash pattern=on 2pt off 2pt] (.3,0) arc (0:360:.3);
\draw (-150:.3) .. controls (180:.1) ..  (150:.3);
\draw (30:.3) -- (-90:.3);
\draw  (90:.3) -- (-30:.3);
\fill (-30:.3) circle (1.5pt) (30:.3) circle (1.5pt) (-90:.3) circle (1.5pt) (90:.3) circle (1.5pt) (-150:.3) circle (1.5pt) (150:.3) circle (1.5pt);
\end{tikzpicture}}
\newcommand{\onechord}{\begin{tikzpicture} \useasboundingbox (-.4,-.2) rectangle (.4,.3);
\draw [->,dash pattern=on 2pt off 2pt] (.3,0) arc (0:360:.3);
\draw (0,-.3) -- (0,.3);
\fill (0,-.3) circle (1.5pt) (0,.3) circle (1.5pt);
\end{tikzpicture}}
\newcommand{\onechordbas}{\begin{tikzpicture} \useasboundingbox (-.4,-.13) rectangle (.4,.2);
\draw [->,dash pattern=on 2pt off 2pt] (.3,0) arc (0:360:.3);
\draw (0,-.3) -- (0,.3);
\fill (0,-.3) circle (1.5pt) (0,.3) circle (1.5pt);
\end{tikzpicture}}
\newcommand{\twoischord}{\begin{tikzpicture} \useasboundingbox (-.4,-.2) rectangle (.4,.3);
\draw [->,dash pattern=on 2pt off 2pt] (.3,0) arc (0:360:.3);
\draw (-45:.3) .. controls (-90:.1) ..  (-135:.3);
\draw (45:.3) .. controls (90:.1) ..  (135:.3);
\fill (-45:.3) circle (1.5pt) (45:.3) circle (1.5pt) (-135:.3) circle (1.5pt) (135:.3) circle (1.5pt);
\end{tikzpicture}}
\newcommand{\twoischordbas}{\begin{tikzpicture} \useasboundingbox (-.4,-.13) rectangle (.4,.2);
\draw [->,dash pattern=on 2pt off 2pt] (.3,0) arc (0:360:.3);
\draw (-45:.3) .. controls (-90:.1) ..  (-135:.3);
\draw (45:.3) .. controls (90:.1) ..  (135:.3);
\fill (-45:.3) circle (1.5pt) (45:.3) circle (1.5pt) (-135:.3) circle (1.5pt) (135:.3) circle (1.5pt);
\end{tikzpicture}}
\newcommand{\twoxchord}{\begin{tikzpicture} \useasboundingbox (-.4,-.2) rectangle (.4,.3);
\draw [->,dash pattern=on 2pt off 2pt] (.3,0) arc (0:360:.3);
\draw (-45:.3) -- (135:.3);
\draw (-135:.3) -- (45:.3);
\fill (-45:.3) circle (1.5pt) (45:.3) circle (1.5pt) (-135:.3) circle (1.5pt) (135:.3) circle (1.5pt);
\end{tikzpicture}}
\newcommand{\twoxchordbas}{\begin{tikzpicture} \useasboundingbox (-.4,-.13) rectangle (.4,.2);
\draw [->,dash pattern=on 2pt off 2pt] (.3,0) arc (0:360:.3);
\draw (-45:.3) -- (135:.3);
\draw (-135:.3) -- (45:.3);
\fill (-45:.3) circle (1.5pt) (45:.3) circle (1.5pt) (-135:.3) circle (1.5pt) (135:.3) circle (1.5pt);
\end{tikzpicture}}
\newcommand{\twochordfreeleg}{\begin{tikzpicture} \useasboundingbox (-.4,-.2) rectangle (.4,.3);
\draw [->,dash pattern=on 2pt off 2pt] (.3,0) arc (0:360:.3);
\draw (-45:.3) -- (-45:.15);
\draw (135:.3) -- (135:.15);
\fill (-45:.3) circle (1.5pt) (135:.3) circle (1.5pt);
\end{tikzpicture}}
\newcommand{\twochordbubble}{\begin{tikzpicture} \useasboundingbox (-.4,-.2) rectangle (.4,.3);
\draw [->,dash pattern=on 2pt off 2pt] (.3,0) arc (360:60:.3) (.3,0) arc (0:60:.3);
\draw (.15,0) -- (.3,0);
\draw (-.15,0) -- (-.3,0);
\draw (-.15,0) .. controls (-.15,.1) and (.15,.1) .. (.15,0);
\draw (-.15,0) .. controls (-.15,-.1) and (.15,-.1) .. (.15,0);
\fill (.15,0) circle (1.5pt) (.3,0) circle (1.5pt) (-.15,0) circle (1.5pt) (-.3,0) circle (1.5pt);
\end{tikzpicture}}
\newcommand{\twochordbubblu}{\begin{tikzpicture} [baseline=-.1cm] \useasboundingbox (-.4,-.2) rectangle (.4,.3);
\draw [->,dash pattern=on 2pt off 2pt]  (.3,0) arc (360:60:.3) (.3,0) arc (0:60:.3);
\draw (.15,.46) node{\scriptsize $K_u$} (.15,0) -- (.3,0);
\draw (-.15,0) -- (-.3,0);
\draw (-.15,0) .. controls (-.15,.1) and (.15,.1) .. (.15,0);
\draw (-.15,0) .. controls (-.15,-.1) and (.15,-.1) .. (.15,0);
\fill (.15,0) circle (1.5pt) (.3,0) circle (1.5pt) (-.15,0) circle (1.5pt) (-.3,0) circle (1.5pt);
\end{tikzpicture}}
\newcommand{\vertbubblu}{\begin{tikzpicture} [baseline=-.1cm] \useasboundingbox (-.15,-.4) rectangle (.55,.45);
\draw [->,dash pattern=on 2pt off 2pt] (.44,.59) node[left]{\scriptsize $K_u$} (.3,-.4) -- (.3,.45);
\draw (0,.15) .. controls (0,.3) .. (.3,.3);
\draw (0,-.15) .. controls (0,-.3) .. (.3,-.3);
\draw (0,-.15) .. controls (.1,-.15) and (.1,.15) .. (0,.15);
\draw (0,-.15) .. controls (-.1,-.15) and (-.1,.15) .. (0,.15);
\fill (0,.15) circle (1.5pt) (.3,.3) circle (1.5pt) (0,-.15) circle (1.5pt) (.3,-.3) circle (1.5pt);
\end{tikzpicture}}
\newcommand{\threeischord}{\begin{tikzpicture} \useasboundingbox (-.4,-.2) rectangle (.4,.5);
\draw [->,dash pattern=on 2pt off 2pt] (.3,0) arc (0:360:.3);
\draw (-30:.3) .. controls (-60:.1) ..  (-90:.3);
\draw (30:.3) .. controls (60:.1) ..  (90:.3);
\draw  (-150:.3) .. controls (180:.1) ..  (150:.3);
\fill (-30:.3) circle (1.5pt) (30:.3) circle (1.5pt) (-90:.3) circle (1.5pt) (90:.3) circle (1.5pt) (-150:.3) circle (1.5pt) (150:.3) circle (1.5pt);
\end{tikzpicture}}
\newcommand{\threeparchord}{\begin{tikzpicture} \useasboundingbox (-.4,-.2) rectangle (.4,.5);
\draw [->,dash pattern=on 2pt off 2pt] (.3,0) arc (0:360:.3);
\draw (-30:.3) -- (150:.3);
\draw (30:.3) .. controls (60:.1) ..  (90:.3);
\draw  (-150:.3) .. controls (-120:.1) ..  (-90:.3);
\fill (-30:.3) circle (1.5pt) (30:.3) circle (1.5pt) (-90:.3) circle (1.5pt) (90:.3) circle (1.5pt) (-150:.3) circle (1.5pt) (150:.3) circle (1.5pt);
\end{tikzpicture}}
\newcommand{\threexparchord}{\begin{tikzpicture} \useasboundingbox (-.4,-.2) rectangle (.4,.5);
\draw [->,dash pattern=on 2pt off 2pt] (.3,0) arc (0:360:.3);
\draw (-150:.3) -- (-30:.3);
\draw (90:.3) -- (-90:.3);
\draw  (150:.3) -- (30:.3);
\fill (-30:.3) circle (1.5pt) (30:.3) circle (1.5pt) (-90:.3) circle (1.5pt) (90:.3) circle (1.5pt) (-150:.3) circle (1.5pt) (150:.3) circle (1.5pt);
\end{tikzpicture}}
\newcommand{\threexparchordbas}{\begin{tikzpicture} \useasboundingbox (-.4,-.13) rectangle (.4,.2);
\draw [->,dash pattern=on 2pt off 2pt] (.3,0) arc (0:360:.3);
\draw (-150:.3) -- (-30:.3);
\draw (90:.3) -- (-90:.3);
\draw  (150:.3) -- (30:.3);
\fill (-30:.3) circle (1.5pt) (30:.3) circle (1.5pt) (-90:.3) circle (1.5pt) (90:.3) circle (1.5pt) (-150:.3) circle (1.5pt) (150:.3) circle (1.5pt);
\end{tikzpicture}}
\newcommand{\threexxchord}{\begin{tikzpicture} \useasboundingbox (-.4,-.2) rectangle (.4,.5);
\draw [->,dash pattern=on 2pt off 2pt] (.3,0) arc (0:360:.3);
\draw (-150:.3) -- (30:.3);
\draw (-30:.3) -- (150:.3);
\draw  (90:.3) -- (-90:.3);
\fill (-30:.3) circle (1.5pt) (30:.3) circle (1.5pt) (-90:.3) circle (1.5pt) (90:.3) circle (1.5pt) (-150:.3) circle (1.5pt) (150:.3) circle (1.5pt);
\end{tikzpicture}}
\newcommand{\threexxchordbas}{\begin{tikzpicture} \useasboundingbox (-.4,-.13) rectangle (.4,.2);
\draw [->,dash pattern=on 2pt off 2pt] (.3,0) arc (0:360:.3);
\draw (-150:.3) -- (30:.3);
\draw (-30:.3) -- (150:.3);
\draw  (90:.3) -- (-90:.3);
\fill (-30:.3) circle (1.5pt) (30:.3) circle (1.5pt) (-90:.3) circle (1.5pt) (90:.3) circle (1.5pt) (-150:.3) circle (1.5pt) (150:.3) circle (1.5pt);
\end{tikzpicture}}
\newcommand{\zerochord}{\begin{tikzpicture} \useasboundingbox (-.4,-.2) rectangle (.4,.3);
\draw [->,dash pattern=on 2pt off 2pt] (.3,0) arc (0:360:.3);
\end{tikzpicture}}
\newcommand{\zerochordbas}{\begin{tikzpicture} \useasboundingbox (-.4,-.13) rectangle (.4,.2);
\draw [->,dash pattern=on 2pt off 2pt] (.3,0) arc (0:360:.3);
\end{tikzpicture}}
\newcommand{\twoischordor}{\begin{tikzpicture} \useasboundingbox (-.4,-.2) rectangle (.4,.3);
\draw [->,dash pattern=on 2pt off 2pt] (.3,0) arc (0:360:.3);
\draw [->] (-30:.3) -- (0,-.15);
\draw (0,-.15) -- (-150:.3);
\draw [->] (30:.3) -- (0,.15);
\draw [->] (0,.15) -- (150:.3);
\fill (-30:.3) circle (1pt) (30:.3) circle (1pt) (-150:.3) circle (1pt) (150:.3) circle (1pt);
\end{tikzpicture}}
\newcommand{\twoxchordor}{\begin{tikzpicture} \useasboundingbox (-.4,-.2) rectangle (.4,.3);
\draw [->,dash pattern=on 2pt off 2pt] (.3,0) arc (0:360:.3);
\draw [->] (-45:.3) -- (135:.15);
\draw (135:.15) -- (135:.3);
\draw [->] (-135:.3) -- (45:.15);
\draw (45:.15) -- (45:.3);
\fill (-45:.3) circle (1pt) (45:.3) circle (1pt) (-135:.3) circle (1pt) (135:.3) circle (1pt);
\end{tikzpicture}}
\newcommand{\twobullechordbis}{\begin{tikzpicture} \useasboundingbox (-.4,-.2) rectangle (.4,.3);
\draw [->,dash pattern=on 2pt off 2pt] (.3,0) arc (0:360:.3);
\draw [-] (-90:.3) -- (-90:.12);
\draw (-90:.12) arc (-90:270:.18);
\draw [-] (-.18,.06) -- (.18,.06);
\fill (-90:.3) circle (1pt) (-.18,.06) circle (1pt) (-90:.12) circle (1pt) (.18,.06) circle (1pt);
\end{tikzpicture}}
\newcommand{\twobullechord}{\begin{tikzpicture} \useasboundingbox (-.4,-.2) rectangle (.4,.3);
\draw [->,dash pattern=on 2pt off 2pt] (.3,0) arc (0:360:.3);
\draw (90:.3) -- (0,.12);
\draw (0,-.12) -- (-90:.3);
\draw [-] (-90:.3) -- (-90:.18);
\draw [-] (0,.12) ..  controls (.06,.12) and (.12,.06) .. (.12,0);
\draw (0,-.12) ..  controls (.06,-.12) and (.12,-.06) .. (.12,0);
\draw [-] (0,.12) ..  controls (-.06,.12) and (-.12,.06) .. (-.12,0);
\draw (0,-.12) ..  controls (-.06,-.12) and (-.12,-.06) .. (-.12,0);
\fill (-90:.3) circle (1pt) (90:.3) circle (1pt) (-90:.12) circle (1pt) (90:.12) circle (1pt);
\end{tikzpicture}}
\newcommand{\twobullechordor}{\begin{tikzpicture} \useasboundingbox (-.4,-.2) rectangle (.4,.3);
\draw [->,dash pattern=on 2pt off 2pt] (.3,0) arc (0:360:.3);
\draw (90:.3) -- (0,.12);
\draw (0,-.12) -- (-90:.3);
\draw [->] (90:.3) -- (90:.18);
\draw [->] (-90:.3) -- (-90:.18);
\draw [->] (0,.12) ..  controls (.06,.12) and (.12,.06) .. (.12,0);
\draw (0,-.12) ..  controls (.06,-.12) and (.12,-.06) .. (.12,0);
\draw [->] (0,.12) ..  controls (-.06,.12) and (-.12,.06) .. (-.12,0);
\draw (0,-.12) ..  controls (-.06,-.12) and (-.12,-.06) .. (-.12,0);
\fill (-90:.3) circle (1pt) (90:.3) circle (1pt) (-90:.12) circle (1pt) (90:.12) circle (1pt);
\end{tikzpicture}}
\newcommand{\tripodor}{\begin{tikzpicture} \useasboundingbox (-.4,-.2) rectangle (.4,.3);
\draw [->,dash pattern=on 2pt off 2pt] (.3,0) arc (0:360:.3);
\draw (-60:.3) -- (0,0);
\draw (60:.3) -- (0,0);
\draw (180:.3) -- (0,0);
\draw [->] (-60:.3) -- (-60:.15);
\draw [->] (60:.3) -- (60:.15);
\draw [->] (180:.3) -- (180:.15);
\fill (60:.3) circle (1.5pt) (180:.3) circle (1.5pt) (-60:.3) circle (1.5pt) (0,0) circle (1.5pt);
\end{tikzpicture}}
\newcommand{\tripodnonor}{\begin{tikzpicture} \useasboundingbox (-.4,-.2) rectangle (.4,.3);
\draw [->,dash pattern=on 2pt off 2pt] (.3,0) arc (0:360:.3);
\draw (-60:.3) -- (0,0);
\draw (60:.3) -- (0,0);
\draw (180:.3) -- (0,0);
\fill (60:.3) circle (1.5pt) (180:.3) circle (1.5pt) (-60:.3) circle (1.5pt) (0,0) circle (1.5pt);
\end{tikzpicture}}
\newcommand{\dquatTun}{\begin{tikzpicture} \useasboundingbox (-.8,-.4) rectangle (.8,.6);
\draw [->,dash pattern=on 2pt off 2pt] (70:.5) arc (70:110:.5);
\draw [->,dash pattern=on 2pt off 2pt] (170:.5) arc (170:240:.5);
\draw [->,dash pattern=on 2pt off 2pt] (-60:.5) arc (-60:10:.5);
\draw (90:.5) --(-40:.5) (-10:.5) -- (205:.5);
\fill (90:.5) circle (1.5pt) (205:.5) circle (1.5pt) (-40:.5) circle (1.5pt) (-10:.5) circle (1.5pt);
\end{tikzpicture}}
\newcommand{\dquatTunnum}{\begin{tikzpicture}
[baseline=-.1cm] 
\useasboundingbox (-.8,-.7) rectangle (.8,.7);
\draw [->,dash pattern=on 2pt off 2pt] (-.35,.5) node{\scriptsize 3} (70:.5) arc (70:110:.5);
\draw [->,dash pattern=on 2pt off 2pt] (-.5,-.5) node{\scriptsize 1} (170:.5) arc (170:240:.5);
\draw [->,dash pattern=on 2pt off 2pt] (.7,.1) node{\scriptsize 2} (-60:.5) arc (-60:10:.5);
\draw (90:.5) --(-40:.5) (-10:.5) -- (205:.5);
\fill (90:.5) circle (1.5pt) (205:.5) circle (1.5pt) (-40:.5) circle (1.5pt) (-10:.5) circle (1.5pt);
\end{tikzpicture}}
\newcommand{\dquatTdeux}{\begin{tikzpicture} \useasboundingbox (-.8,-.4) rectangle (.8,.6);
\draw [->,dash pattern=on 2pt off 2pt] (70:.5) arc (70:110:.5);
\draw [->,dash pattern=on 2pt off 2pt] (170:.5) arc (170:240:.5);
\draw [->,dash pattern=on 2pt off 2pt] (-60:.5) arc (-60:10:.5);
\draw (90:.5) -- (220:.5) (190:.5) -- (-25:.5);
\fill (90:.5) circle (1.5pt) (220:.5) circle (1.5pt) (-25:.5) circle (1.5pt) (190:.5) circle (1.5pt);
\end{tikzpicture}}
\newcommand{\dquatTtrois}{\begin{tikzpicture} \useasboundingbox (-.8,-.4) rectangle (.8,.6);
\draw [->,dash pattern=on 2pt off 2pt] (70:.5) arc (70:110:.5);
\draw [->,dash pattern=on 2pt off 2pt] (170:.5) arc (170:240:.5);
\draw [->,dash pattern=on 2pt off 2pt] (-60:.5) arc (-60:10:.5);
\draw (90:.5) -- (-10:.5) (-40:.5) -- (205:.5);
\fill (90:.5) circle (1.5pt) (205:.5) circle (1.5pt) (-40:.5) circle (1.5pt) (-10:.5) circle (1.5pt);
\end{tikzpicture}}
\newcommand{\dquatTqua}{\begin{tikzpicture} \useasboundingbox (-.8,-.4) rectangle (.8,.6);
\draw [->,dash pattern=on 2pt off 2pt] (70:.5) arc (70:110:.5);
\draw [->,dash pattern=on 2pt off 2pt] (170:.5) arc (170:240:.5);
\draw [->,dash pattern=on 2pt off 2pt] (-60:.5) arc (-60:10:.5);
\draw (90:.5) --(190:.5) (220:.5) -- (-25:.5);
\fill (90:.5) circle (1.5pt) (220:.5) circle (1.5pt) (-25:.5) circle (1.5pt) (190:.5) circle (1.5pt);
\end{tikzpicture}}
\newcommand{\ischord}{\begin{tikzpicture} \useasboundingbox (-.2,-.2) rectangle (.4,.3);
\draw [->,dash pattern=on 2pt off 2pt] (-90:.3) arc (-90:90:.3);
\draw (-45:.3) -- (45:.3);
\fill (-45:.3) circle (1.5pt) (45:.3) circle (1.5pt);
\end{tikzpicture}}
\newcommand{\pkink}{\begin{tikzpicture} \useasboundingbox (-.1,-.1) rectangle (.5,.2);
\draw [out=-90,in=-90,->] (.4,0) to (0,.2);
\draw [out=90,in=90,draw=white,double=black,very thick] (0,-.2) to (.4,0);
\end{tikzpicture}}
\newcommand{\nkink}{\begin{tikzpicture} \useasboundingbox (-.1,-.1) rectangle (.5,.2);
\draw [out=90,in=90] (0,-.2) to (.4,0);
\draw [out=-90,in=-90,draw=white,double=black,very thick] (.4,0) to (0,.2);
\draw [out=-90,in=-90,->] (.4,0) to (0,.2);
\end{tikzpicture}}
\newcommand{\nkinks}{\begin{tikzpicture} \useasboundingbox (-.5,-.1) rectangle (.1,.2);
\draw [out=-90,in=-90,->] (-.4,0) to (0,.2);
\draw [out=90,in=90,draw=white,double=black,very thick] (0,-.2) to (-.4,0);
\end{tikzpicture}}
\newcommand{\pkinks}{\begin{tikzpicture} \useasboundingbox (-.5,-.1) rectangle (.1,.2);
\draw [out=90,in=90] (0,-.2) to (-.4,0);
\draw [out=-90,in=-90,draw=white,double=black,very thick] (-.4,0) to (0,.2);
\draw [out=-90,in=-90,->] (-.4,0) to (0,.2);
\end{tikzpicture}}
\newcommand{\exjactwosonebis}{\begin{tikzpicture} \useasboundingbox (-1.4,-.3) rectangle (1.4,.3);
\draw [->,dash pattern=on 2pt off 2pt] (0,-.3) arc (-90:270:.3);
\draw [->,dash pattern=on 2pt off 2pt] (-1,-.3) arc (-90:270:.3);
\draw (-1.3,0) -- (30:.3);
\draw (0,0) -- (-140:.3) (0,0) -- (-40:.3) (0,0) -- (0,.3);
\fill (-40:.3) circle (1.5pt) (-140:.3) circle (1.5pt) (90:.3) circle (1.5pt) (30:.3) circle (1.5pt) (-1.3,0) circle (1.5pt) (0,0) circle (1.5pt) (1.2,0) circle (1.5pt) (.8,0) circle (1.5pt);
\draw (1,0) circle (.2) (.8,0) -- (1.2,0) (-.85,-.5) node{\scriptsize $C_1$} (.15,-.5) node{\scriptsize $C_2$};
\end{tikzpicture}}
\newcommand{\lolli}{\begin{tikzpicture} \useasboundingbox (-.4,-.2) rectangle (.4,.3);
\draw [->,dash pattern=on 2pt off 2pt] (.3,0) arc (0:360:.3);
\draw (0,-.3) -- (0,-.1) (0,.05) circle (.15);
\fill (0,-.3) circle (1.5pt) (0,-.1) circle (1.5pt);
\end{tikzpicture}}
\newcommand{\tatasone}{\begin{tikzpicture} \useasboundingbox (-.4,-.2) rectangle (.4,.3);
\draw [->,dash pattern=on 2pt off 2pt] (.3,0) arc (0:360:.3);
\draw (-.2,0) -- (.2,0) (0,0) circle (.2);
\fill (-.2,0) circle (1.5pt) (.2,0) circle (1.5pt);
\end{tikzpicture}}
\newcommand{\tatasonetw}{\begin{tikzpicture} \useasboundingbox (-.4,-.2) rectangle (.4,.3);
\draw [->,dash pattern=on 2pt off 2pt] (.3,0) arc (0:360:.3);
\draw (.2,0) arc (0:180:.2) (-.2,0) -- (-.05,0) .. controls (.05,0) and (.2,-.3) .. (.2,0) 
(-.2,0) arc (-180:-90:.2) (0,-.2) .. controls (.15,-.2) and (0,.1) .. (.2,0);
\fill (-.2,0) circle (1.5pt) (.2,0) circle (1.5pt);
\end{tikzpicture}}
\newcommand{\haltsone}{\begin{tikzpicture} \useasboundingbox (-.4,-.2) rectangle (.4,.3);
\draw [->,dash pattern=on 2pt off 2pt] (.35,0) arc (0:360:.35);
\draw (0,-.07) -- (0,.07) (0,.19) circle (.12) (0,-.19) circle (.12);
\fill (0,-.07) circle (1.5pt) (0,.07) circle (1.5pt);
\end{tikzpicture}}
\newcommand{\trivAS}{\begin{tikzpicture} \useasboundingbox (-.4,-.3) rectangle (.4,.4);
\draw (0,0) -- (0,-.4) (0,0) -- (50:.4) (0,0) -- (130:.4);
\fill (0,0) circle (1.5pt);
\end{tikzpicture}}
\newcommand{\trivASop}{\begin{tikzpicture} \useasboundingbox (-.4,-.2) rectangle (.4,.3);
\draw (0,0) -- (0,-.2);
\draw (0,0) .. controls (50:.3) ..  (130:.5);
\draw [draw=white,double=black,very thick] (0,0) .. controls (130:.3) ..  (50:.5);
\fill (0,0) circle (1.5pt);
\end{tikzpicture}}
\newcommand{\ihxone}{\begin{tikzpicture} \useasboundingbox (-.4,-.3) rectangle (.4,.4);
\draw (-90:.4) -- (0,0) -- (30:.4) (0,0) -- (150:.4);
\draw (-.3,-.25) -- (0,-.25);
\fill (0,-.25) circle (1.5pt) (0,0) circle (1.5pt);
\end{tikzpicture}}
\newcommand{\ihxtwo}{\begin{tikzpicture} \useasboundingbox (-.4,-.3) rectangle (.4,.4);
\draw (-90:.4) -- (0,0) -- (30:.4) (0,0) -- (150:.4);
\draw [out=0,in=-60] (-.3,-.25) to (30:.25);
\fill (30:.25) circle (1.5pt) (0,0) circle (1.5pt);
\end{tikzpicture}}
\newcommand{\ihxthree}{\begin{tikzpicture} \useasboundingbox (-.4,-.3) rectangle (.4,.4);
\draw (-90:.4) -- (0,0) -- (30:.4) (0,0) -- (150:.4);
\draw (-.3,-.25) .. controls (-.2,-.25) and  (.25,-.2)  .. (.25,0)  .. controls (.25,.2) and (120:.35) .. (150:.25);
\fill (150:.25) circle (1.5pt) (0,0) circle (1.5pt);
\end{tikzpicture}}
\newcommand{\ihxi}{\begin{tikzpicture}[scale=.7] \useasboundingbox (-.4,-.3) rectangle (.4,.4);
\draw (-.3,-.3) -- (.3,-.3) (-.3,.3) -- (.3,.3) (0,-.3) -- (0,.3);
\fill (0,-.3) circle (2pt) (0,.3) circle (2pt);
\end{tikzpicture}}
\newcommand{\ihxh}{\begin{tikzpicture}[scale=.7] \useasboundingbox (-.4,-.3) rectangle (.4,.4);
\draw (-.3,-.3) -- (-.3,.3) (.3,-.3) -- (.3,.3) (-.3,0) -- (.3,0);
\fill (-.3,0) circle (2pt) (.3,0) circle (2pt);
\end{tikzpicture}}
\newcommand{\ihxx}{\begin{tikzpicture}[scale=.7] \useasboundingbox (-.4,-.3) rectangle (.4,.4);
\draw (-.3,-.3) -- (.3,.3) (.3,-.3) -- (-.3,.3) (-.15,-.15) -- (.15,-.15);
\fill (-.15,-.15) circle (2pt) (.15,-.15) circle (2pt);
\end{tikzpicture}}
\newcommand{\stuy}{\begin{tikzpicture} \useasboundingbox (-.1,-.25) rectangle (.7,.4);
\draw [->,dash pattern=on 2pt off 2pt] (0,-.2) -- (.6,-.2);
\draw (.3,-.2) -- (.3,0) -- (.15,.2) (.3,0) -- (.45,.2);
\fill (.3,-.2) circle (1.5pt) (.3,0) circle (1.5pt);
\end{tikzpicture}}
\newcommand{\stuyunor}{\begin{tikzpicture} \useasboundingbox (-.1,-.25) rectangle (.7,.4);
\draw [dash pattern=on 2pt off 2pt] (0,-.2) -- (.6,-.2);
\draw (.3,-.2) -- (.3,0) -- (.15,.2) (.3,0) -- (.45,.2);
\fill (.3,-.2) circle (1.5pt) (.3,0) circle (1.5pt);
\end{tikzpicture}}
\newcommand{\stuyab}{\begin{tikzpicture} \useasboundingbox (-.3,-.25) rectangle (.8,.4);
\draw [->,dash pattern=on 2pt off 2pt] (0,-.2) -- (.6,-.2);
\draw (.3,-.2) -- (.3,0) -- (.15,.2) (.3,0) -- (.45,.2) (.4,.2) node[right]{\scriptsize $a$} (.2,.2) node[left]{\scriptsize $b$};
\fill (.3,-.2) circle (1.5pt) (.3,0) circle (1.5pt);
\end{tikzpicture}}
\newcommand{\stui}{\begin{tikzpicture} \useasboundingbox (-.1,-.25) rectangle (.8,.4);
\draw [->,dash pattern=on 2pt off 2pt] (0,-.2) -- (.7,-.2);
\draw (.15,-.2) -- (.15,.2)  (.45,.2) -- (.45,-.2);
\fill (.15,-.2) circle (1.5pt) (.45,-.2) circle (1.5pt);
\end{tikzpicture}}
\newcommand{\stuiunor}{\begin{tikzpicture} \useasboundingbox (-.1,-.25) rectangle (.7,.4);
\draw [dash pattern=on 2pt off 2pt] (0,-.2) -- (.6,-.2);
\draw (.15,-.2) -- (.15,.2)  (.45,.2) -- (.45,-.2);
\fill (.15,-.2) circle (1.5pt) (.45,-.2) circle (1.5pt);
\end{tikzpicture}}
\newcommand{\stuiab}{\begin{tikzpicture} \useasboundingbox (-.15,-.15) rectangle (.85,.2);
\draw [->,dash pattern=on 2pt off 2pt] (0,-.2) -- (.6,-.2);
\draw (.15,-.2) -- (.15,.2)  (.45,.2) -- (.45,-.2) (.4,.2) node[right]{\scriptsize $a$} (.2,.2) node[left]{\scriptsize $b$};
\fill (.15,-.2) circle (1.5pt) (.45,-.2) circle (1.5pt);
\end{tikzpicture}}
\newcommand{\stux}{\begin{tikzpicture} \useasboundingbox (-.1,-.25) rectangle (.8,.4);
\draw [->,dash pattern=on 2pt off 2pt] (0,-.2) -- (.7,-.2);
\draw (.15,.2) -- (.45,-.2);
\draw [draw=white,double=black,very thick] (.15,-.2) -- (.45,.2);
\fill (.15,-.2) circle (1.5pt) (.45,-.2) circle (1.5pt);
\end{tikzpicture}}
\newcommand{\stuxunor}{\begin{tikzpicture} \useasboundingbox (-.1,-.25) rectangle (.7,.4);
\draw [dash pattern=on 2pt off 2pt] (0,-.2) -- (.6,-.2);
\draw (.15,.2) -- (.45,-.2);
\draw [draw=white,double=black,very thick] (.15,-.2) -- (.45,.2);
\fill (.15,-.2) circle (1.5pt) (.45,-.2) circle (1.5pt);
\end{tikzpicture}}
\newcommand{\stuxdis}{\begin{tikzpicture} \useasboundingbox (-.1,-.2) rectangle (.8,.2);
\draw [->,dash pattern=on 2pt off 2pt] (0,-.2) -- (.7,-.2);
\draw (.15,.2) -- (.45,-.2);
\draw [draw=white,double=black,very thick] (.15,-.2) -- (.45,.2);
\fill (.15,-.2) circle (1.5pt) (.45,-.2) circle (1.5pt);
\end{tikzpicture}}
\newcommand{\stuidis}{\begin{tikzpicture} \useasboundingbox (-.1,-.2) rectangle (.8,.2);
\draw [->,dash pattern=on 2pt off 2pt] (0,-.2) -- (.7,-.2);
\draw (.15,-.2) -- (.15,.2)  (.45,.2) -- (.45,-.2);
\fill (.15,-.2) circle (1.5pt) (.45,-.2) circle (1.5pt);
\end{tikzpicture}}
\newcommand{\stuydis}{\begin{tikzpicture} \useasboundingbox (-.1,-.2) rectangle (.7,.2);
\draw [->,dash pattern=on 2pt off 2pt] (0,-.2) -- (.6,-.2);
\draw (.3,-.2) -- (.3,0) -- (.15,.2) (.3,0) -- (.45,.2);
\fill (.3,-.2) circle (1.5pt) (.3,0) circle (1.5pt);
\end{tikzpicture}}
\newcommand{\stuxab}{\begin{tikzpicture} \useasboundingbox (-.15,-.15) rectangle (.85,.2);
\draw [->,dash pattern=on 2pt off 2pt] (0,-.2) -- (.6,-.2);
\draw (.15,.2) -- (.45,-.2) (.4,.2) node[right]{\scriptsize $a$} (.2,.2) node[left]{\scriptsize $b$};
\draw [draw=white,double=black,very thick] (.15,-.2) -- (.45,.2);
\fill (.15,-.2) circle (1.5pt) (.45,-.2) circle (1.5pt);
\end{tikzpicture}}
\newcommand{\dquatSTU}{\begin{tikzpicture} \useasboundingbox (-.8,-.4) rectangle (.8,.6);
\draw [->,dash pattern=on 2pt off 2pt] (70:.5) arc (70:110:.5);
\draw [->,dash pattern=on 2pt off 2pt] (170:.5) arc (170:240:.5);
\draw [->,dash pattern=on 2pt off 2pt] (-60:.5) arc (-60:10:.5);
\draw (0,0) -- (-25:.5) (0,0) -- (90:.5) (0,0) -- (205:.5);
\fill (90:.5) circle (1.5pt) (-25:.5) circle (1.5pt) (205:.5) circle (1.5pt) (0,0) circle (1.5pt);
\end{tikzpicture}}
\newcommand{\asunor}{\begin{tikzpicture} \useasboundingbox (-.1,-.25) rectangle (.7,.4);
\draw [dash pattern=on 2pt off 2pt] (0,-.2) -- (.6,-.2);
\draw (.3,.2) -- (.3,-.2);
\fill (.3,-.2) circle (1.5pt);
\end{tikzpicture}}
\newcommand{\asunorc}{\begin{tikzpicture} \useasboundingbox (-.1,-.25) rectangle (.7,.4);
\draw [dash pattern=on 2pt off 2pt] (0,-.2) -- (.6,-.2);
\draw (.3,.2) .. controls (.1,0) and  (.05,-.35)  .. (.15,-.35)  .. controls (.25,-.35) .. (.3,-.2);
\fill (.3,-.2) circle (1.5pt);
\end{tikzpicture}}
\newcommand{\stubone}{\begin{tikzpicture} \useasboundingbox (-.1,-.25) rectangle (1.1,.4);
\draw [dash pattern=on 2pt off 2pt] (0,-.2) -- (1,-.2);
\draw (.3,.3) -- (.3,-.2) (.7,.3) -- (.7,-.2);
\fill (.3,-.2) circle (1.5pt) (.7,-.2) circle (1.5pt);
\end{tikzpicture}}
\newcommand{\stubtwo}{\begin{tikzpicture} \useasboundingbox (-.1,-.25) rectangle (1.1,.4);
\draw [dash pattern=on 2pt off 2pt] (0,-.2) -- (1,-.2);
\draw (.4,.3)  .. controls (.1,-.1) and  (.3,-.35)  .. (.4,-.35)  .. controls (.5,-.35) .. (.7,-.2) (.7,.3) -- (.4,-.2);
\fill (.4,-.2) circle (1.5pt) (.7,-.2) circle (1.5pt);
\end{tikzpicture}}
\newcommand{\stubthree}{\begin{tikzpicture} \useasboundingbox (-.1,-.25) rectangle (1.1,.4);
\draw [dash pattern=on 2pt off 2pt] (0,-.2) -- (1,-.2);
\draw (.2,.3)  .. controls (.2,-.2) and  (.7,-.1)  .. (.7,0)  .. controls (.7,.1) .. (.5,.15) (.5,.3) -- (.5,-.2);
\fill (.5,-.2) circle (1.5pt) (.5,.15) circle (1.5pt);
\end{tikzpicture}}
\newcommand{\graphkev}{\begin{tikzpicture} \useasboundingbox (-.45,-.2) rectangle (.45,.4);
\draw (-150:.18) -- (-30:.18) -- (90:.18) -- (-150:.18) -- (-150:.4) -- (-30:.4) -- (90:.4) -- (-150:.4) (90:.4) -- (90:.18) (-30:.4) -- (-30:.18);
\fill (-30:.18) circle (1.5pt) (90:.18) circle (1.5pt) (-150:.18) circle (1.5pt) (-30:.4) circle (1.5pt) (90:.4) circle (1.5pt) (-150:.4) circle (1.5pt);
\end{tikzpicture}}
\newcommand{\annulone}{\begin{tikzpicture} \useasboundingbox (-1.3,-.8) rectangle (1.3,1.2);
\draw [fill=gray!20,draw=white] (0,0) circle (.9);
\draw [fill=white,draw=white] (0,0) circle (.5) (.6,-.2) rectangle (.8,.2);
\draw [dash pattern=on 2pt off 2pt] (0,-1) -- (0,1);
\draw (.7,0) node{\tiny $A$} (0,-.2) -- (-.3,0) -- (0,0) (-.3,0) -- (0,.2) (-.5,-1) -- (0,-.7) (.1,-.6) node{\tiny $v$} (.2,.75) node{\tiny $\alpha_2$} (.2,-.8) node{\tiny $\alpha_1$} (.25,0) node{\tiny $D$} (-.3,-.7) node{\tiny $\beta$} (.8,-.8) node{$\Gamma_1$};
\fill (0,-.7) circle (1.5pt) (0,-.2) circle (1.5pt) (0,0) circle (1.5pt) (0,.2) circle (1.5pt) (-.3,0) circle (1.5pt);
\end{tikzpicture}}
\newcommand{\annultwo}{\begin{tikzpicture} \useasboundingbox (-1.3,-.8) rectangle (1.3,1.2);
\draw [fill=gray!20,draw=white] (0,0) circle (.9);
\draw [fill=white,draw=white] (0,0) circle (.5) (.6,-.2) rectangle (.8,.2);
\draw [dash pattern=on 2pt off 2pt] (0,-1) -- (0,1);
\draw (.7,0) node{\tiny $A$} (0,-.2) -- (-.3,0) -- (0,0) (-.3,0) -- (0,.2) (-.5,-1) .. controls (.1,-.7) and (.55,-.4) .. (.55,0) .. controls  (.55,.4) and  (.2,.7) ..  (0,.7) (-.1,.8) node{\tiny $v$} (.2,.75) node{\tiny $\alpha_2$} (.25,0) node{\tiny $D$} (.2,-.8) node{\tiny $\alpha_1$} (-.3,-.7) node{\tiny $\beta$} (.8,-.8) node{$\Gamma_2$};
\fill (0,.7) circle (1.5pt) (0,-.2) circle (1.5pt) (0,0) circle (1.5pt) (0,.2) circle (1.5pt) (-.3,0) circle (1.5pt);
\end{tikzpicture}}
\newcommand{\annulthree}{\begin{tikzpicture} \useasboundingbox (-1.3,-.8) rectangle (1.3,1.2);
\draw [fill=gray!20,draw=white] (0,0) circle (.9);
\draw [fill=white,draw=white] (0,0) circle (.5) (.6,-.2) rectangle (.8,.2);
\draw [dash pattern=on 2pt off 2pt] (0,-1) -- (0,1);
\draw (.7,0) node{\tiny $A$} (-1,0) -- (0,0) (-.5,-1) -- (0,-.7) (.1,-.6) node{\tiny $v$} (.25,0) node{\tiny $D$} (.2,.75) node{\tiny $\alpha_2$} (.2,-.8) node{\tiny $\alpha_1$} (-.7,.15) node{\tiny $\alpha_3$} (-.3,-.7) node{\tiny $\beta$} (.8,-.8) node{$\Gamma_1$};
\fill (0,-.7) circle (1.5pt) (0,0) circle (1.5pt);
\end{tikzpicture}}
\newcommand{\annulfour}{\begin{tikzpicture} \useasboundingbox (-1.3,-.8) rectangle (1.3,1.2);
\draw [fill=gray!20,draw=white] (0,0) circle (.9);
\draw [fill=white,draw=white] (0,0) circle (.5) (.6,-.2) rectangle (.8,.2);
\draw [dash pattern=on 2pt off 2pt] (0,-1) -- (0,1);
\draw  (-.5,-1) .. controls (.1,-.7) and (.55,-.4) .. (.55,0) .. controls  (.55,.4) and  (.2,.7) ..  (0,.7) (.7,0) node{\tiny $A$} (-1,0) -- (0,0) (-.1,.8) node{\tiny $v$} (.25,0) node{\tiny $D$} (.2,.75) node{\tiny $\alpha_2$} (.2,-.8) node{\tiny $\alpha_1$} (-.7,.15) node{\tiny $\alpha_3$} (-.3,-.7) node{\tiny $\beta$} (.8,-.8) node{$\Gamma_2$};
\fill (0,.7) circle (1.5pt) (0,0) circle (1.5pt);
\end{tikzpicture}}
\newcommand{\annulfive}{\begin{tikzpicture} \useasboundingbox (-1.3,-.8) rectangle (1.3,1.2);
\draw [fill=gray!20,draw=white] (0,0) circle (.9);
\draw [fill=white,draw=white] (0,0) circle (.5) (.6,-.2) rectangle (.8,.2);
\draw [dash pattern=on 2pt off 2pt] (0,-1) -- (0,1);
\draw (-.5,-1) .. controls (.1,-.7) and (.55,-.4) .. (.55,0) .. controls  (.55,.4) and  (.2,.7) ..  (0,.7) .. controls  (-.25,.7) and  (-.7,.35) ..  (-.7,0)  (.7,0) node{\tiny $A$} (-1,0) -- (0,0) (-.8,.1) node{\tiny $v$} (.25,0) node{\tiny $D$} (.2,.75) node{\tiny $\alpha_2$} (.2,-.8) node{\tiny $\alpha_1$} (-.7,-.15) node{\tiny $\alpha_3$} (-.3,-.7) node{\tiny $\beta$} (.8,-.8) node{$\Gamma_3$};
\fill (-.7,0) circle (1.5pt) (0,0) circle (1.5pt);
\end{tikzpicture}}
\newcommand{\dmer}{\begin{tikzpicture} [baseline=-.1cm] \useasboundingbox (-.5,-.8) rectangle (.2,1);
\draw [->,dash pattern=on 2pt off 2pt] (0,-.8) -- (0,.8);
\draw (0,-.2) -- (-.3,0) -- (0,0) (-.3,0) -- (0,.2);
\fill (0,-.2) circle (1.5pt) (0,0) circle (1.5pt) (0,.2) circle (1.5pt) (-.3,0) circle (1.5pt);
\end{tikzpicture}}
\newcommand{\cirt}{\begin{tikzpicture} [baseline=-.1cm] \useasboundingbox (-1,-.8) rectangle (1,1);
\draw [->,dash pattern=on 2pt off 2pt]  (-.8,0) arc (-180:180:.8);
\draw (-135:.8) -- (0,-.15) -- (0,.15) -- (135:.8)  (0,.15) -- (45:.8) (0,-.15) --  (-45:.8);
\fill (-135:.8) circle (1.5pt) (135:.8) circle (1.5pt) (45:.8) circle (1.5pt) (-45:.8) circle (1.5pt) (0,-.15) circle (1.5pt) (0,.15) circle (1.5pt);
\end{tikzpicture}}
\newcommand{\cirtd}{\begin{tikzpicture} [baseline=-.1cm] \useasboundingbox (-1,-.8) rectangle (1,1);
\draw [->,dash pattern=on 2pt off 2pt]  (-.8,0) arc (-180:180:.8);
\draw (-135:.8) -- (0,-.15) -- (0,.15) -- (135:.8)  (0,.15) -- (60:.8) (0,-.15) --  (-60:.8);
\draw (-30:.8) -- (.4,0) -- (30:.8) (.4,0) -- (.8,0);
\fill (-135:.8) circle (1.5pt) (135:.8) circle (1.5pt) (60:.8) circle (1.5pt) (-60:.8) circle (1.5pt) (30:.8) circle (1.5pt) (0:.8) circle (1.5pt) (-30:.8) circle (1.5pt) (.4,0) circle (1.5pt) (0,-.15) circle (1.5pt) (0,.15) circle (1.5pt);
\end{tikzpicture}}
\newcommand{\cirtdp}{\begin{tikzpicture}[baseline=-.1cm]  \useasboundingbox (-1,-.8) rectangle (1,1);
\draw [->,dash pattern=on 2pt off 2pt]  (-.8,0) arc (-180:180:.8);
\draw (-135:.8) -- (0,-.15) -- (0,.15) -- (150:.8)  (0,.15) -- (30:.8) (0,-.15) --  (-45:.8);
\draw (60:.8) -- (0,.45) -- (120:.8) (0,.45) -- (0,.8);
\fill (-135:.8) circle (1.5pt) (150:.8) circle (1.5pt) (30:.8) circle (1.5pt) (-45:.8) circle (1.5pt) (60:.8) circle (1.5pt) (90:.8) circle (1.5pt) (120:.8) circle (1.5pt) (0,.45) circle (1.5pt) (0,-.15) circle (1.5pt) (0,.15) circle (1.5pt);
\end{tikzpicture}}
\newcommand{\loopedge}{\begin{tikzpicture} [baseline=-.15cm] \useasboundingbox (-.1,-.12) rectangle (.65,.12);
\draw  (.5,0) arc (0:360:.15);
\draw (-.1,0) -- (.2,0);
\fill (.2,0) circle (1.5pt);
\end{tikzpicture}}

\newcommand{\tata}{ \begin{tikzpicture} \useasboundingbox (-.25,-0.02) rectangle (.25,.15);
\draw (-.15,.1) -- (.15,.1) (0,.1) circle (.15);
\fill (-.15,.1) circle (1pt) (.15,.1) circle (1pt);
\end{tikzpicture}} 
\newcommand{\tatak}{ \begin{tikzpicture}\useasboundingbox (0,0) rectangle (.3,.3);
\draw [dotted] (0.15,0.15) circle (.15);
\draw (0,0.15) -- (.3,0.15);
\fill (0,0.15) circle (1pt) (.3,0.15) circle (1pt);
\end{tikzpicture}}
\newcommand{\tataksmall}{ \begin{tikzpicture}\useasboundingbox (-.04,0) rectangle (.2,.16);
\draw [dotted] (0.08,0.08) circle (.08);
\draw (0,0.06) -- (.16,0.08);
\fill (0,0.08) circle (1pt) (.16,0.08) circle (1pt);
\end{tikzpicture}}
\newcommand{\tatanumone}{ \begin{tikzpicture}[scale=1.3] \useasboundingbox (-.5,-.1) rectangle (.5,.4);
\draw (0,0) -- (.3,0) (0,0) circle (.3) (-.35,.15) node{\scriptsize $5$} (-.35,-.15) node{\scriptsize $1$} (.35,.15) node{\scriptsize $6$} (.35,-.15) node{\scriptsize $2$} (-.15,.1) node{\scriptsize $3$} (.15,.1) node{\scriptsize $4$};
\draw [->] (-.3,0) -- (0,0);
\draw [->] (-.3,0) arc (-180:-90:.3);
\draw [->] (-.3,0) arc (180:90:.3);
\fill (-.3,0) circle (1.5pt) (.3,0) circle (1.5pt);
\end{tikzpicture}}
\newcommand{\tatanumtwo}{ \begin{tikzpicture}[scale=1.3] \useasboundingbox (-.5,-.1) rectangle (.5,.4);
\draw (0,0) -- (.3,0) (0,0) circle (.3) (-.35,.15) node{\scriptsize $3$} (-.35,-.15) node{\scriptsize $1$}  (.35,.15) node{\scriptsize $6$} (.35,-.15) node{\scriptsize $5$} (-.15,.1) node{\scriptsize $2$} (.15,.1) node{\scriptsize $4$};
\draw [->] (-.3,0) -- (0,0);
\draw [->] (-.3,0) arc (-180:-90:.3);
\draw [->] (-.3,0) arc (180:90:.3);
\fill (-.3,0) circle (1.5pt) (.3,0) circle (1.5pt);
\end{tikzpicture}}
\newcommand{\haltereor}{\begin{tikzpicture} \useasboundingbox (-1,-.1) rectangle (1,.2);
\draw [-<,dash pattern=on 2pt off 2pt]  (-.5,0) arc (-180:180:.15);
\draw [-<,dash pattern=on 2pt off 2pt]  (.5,0) arc (0:360:.15);
\draw [->] (-.2,0) -- (0,0);
\draw (0,0) -- (.2,0);
\fill (-.2,0) circle (1.5pt) (.2,0) circle (1.5pt) (.8,0) node{\scriptsize $S^1_j$} (-.8,0) node{\scriptsize $S^1_i$};
\end{tikzpicture}}
\newcommand{\haltereorKKp}{\begin{tikzpicture} \useasboundingbox (-1,-.1) rectangle (1.6,.2);
\draw [-<,dash pattern=on 2pt off 2pt]  (-.64,0) arc (-180:180:.22);
\draw [-<,dash pattern=on 2pt off 2pt]  (.54,0) arc (0:360:.22);
\draw [->] (-.2,0) -- (0,0);
\draw (0,0) -- (.1,0);
\fill (-.2,0) circle (1.5pt) (.1,0) circle (1.5pt) (1.1,0) node{\scriptsize $K_{\parallel,\tau,Y}$} (-.85,0) node{\scriptsize $K$};
\end{tikzpicture}}
\newcommand{\tataor}{ \begin{tikzpicture} \useasboundingbox (-.3,-.1) rectangle (.3,.2);
\draw [->] (-.2,0) -- (0,0);
\draw (0,0) -- (.2,0) (0,0) circle (.2);
\draw [->] (-.2,0) arc (-180:-90:.2);
\draw [->] (-.2,0) arc (180:90:.2);
\fill (-.2,0) circle (1.5pt) (.2,0) circle (1.5pt);
\end{tikzpicture}}
\newcommand{\tatanumor}[2]{ \begin{tikzpicture} \useasboundingbox (-.4,-.2) rectangle (.4,.2);
\begin{scope}[yshift=-.1cm]
\draw (0,0) -- (.3,0) (0,0) circle (.3) (.05,.38) node{\tiny $1$} (.05,.08) node{\tiny $2$} (.05,-.18) node{\tiny $3$};
\draw [-#1] (-.3,0) -- (0,0);
\draw [-#2] (-.3,0) arc (-180:-90:.3);
\draw [->] (-.3,0) arc (180:90:.3);
\fill (-.3,0) circle (1.5pt) (.3,0) circle (1.5pt);
\end{scope}
\end{tikzpicture}}
\newcommand{\onechordsegment}[2]{\begin{tikzpicture} \useasboundingbox (-.5,-.1) rectangle (.9,.2);
\draw (0,0) -- (.4,0) node[right]{\scriptsize $#2$};
\draw (0,0) node[left] {\scriptsize $#1$};
\draw [->] (0,0) -- (.2,0);
\fill (0,0) circle (1pt) (.4,0) circle (1pt);
\end{tikzpicture}}
\newcommand{\onechordnum}[2]{\begin{tikzpicture} \useasboundingbox (-.4,-.2) rectangle (.73,.2);
\begin{scope}[yshift=-.1cm]
\draw [->,dash pattern=on 2pt off 2pt]  (0,-.27) arc (-90:270:.27);
\draw [->] (-.27,0) -- (0,0);
\draw (0,0) -- (.27,0);
\draw  (.53,0) node{\scriptsize #1} (.04,.13) node{\tiny #2};
\fill (-.27,0) circle (1.5pt) (.27,0) circle (1.5pt);
\end{scope}
\end{tikzpicture}}
\newcommand{\onechordtwocircles}[2]{\begin{tikzpicture}[scale=1.5] \useasboundingbox (-.85,-.1) rectangle (.77,.2);
\draw [-<,dash pattern=on 2pt off 2pt]  (-.5,0) arc (-180:180:.15);
\draw [-<,dash pattern=on 2pt off 2pt]  (.4,0) arc (0:360:.15);
\draw [-#2] (-.2,0) -- (0,0);
\draw  (-.02,.17) node{\scriptsize $#1$};
\draw (0,0) -- (.1,0);
\fill (-.2,0) circle (1pt) (.1,0) circle (1pt) (.6,0) node{\scriptsize $S^1_2$} (-.7,0) node{\scriptsize $S^1_1$};
\end{tikzpicture}}
\newcommand{\onechordtwocirclesrightmoyen}[2]{\begin{tikzpicture} \useasboundingbox (-1,-.1) rectangle (1,.2);
\draw [-<,dash pattern=on 2pt off 2pt]  (-.64,0) arc (-180:180:.22);
\draw [-<,dash pattern=on 2pt off 2pt]  (.54,0) arc (0:360:.22);
\draw [->] (-.2,0) -- (0,0);
\draw (0,0) -- (.1,0);
\fill (-.2,0) circle (1.5pt) (.1,0) circle (1.5pt) (.8,0) node{\scriptsize $#2$} (-.85,0) node{\scriptsize $#1$};
\end{tikzpicture}}
\newcommand{\onechordsmalljnum}{\onechordnum{$S^1_j$}{$k$}} 
\newcommand{\onechordsmallnumt}{\begin{tikzpicture} \useasboundingbox (-.4,-.15) rectangle (.5,.15);
\begin{scope}[yshift=-.05cm]
 \draw [->,dash pattern=on 2pt off 2pt]  (0,-.2) -- (0,.25);
\draw (-.05,.15) node[right]{\scriptsize K};
\draw (-.3,0) node{\scriptsize k};
\draw [->]  (0,-.12) arc (270:180:.12);
\draw (0,.12) arc (90:180:.12);
\fill (0,-.12) circle (1pt) (0,.12) circle (1pt);\end{scope}
\end{tikzpicture}}

\newcommand{\onechordsmallnonnum}{\begin{tikzpicture} \useasboundingbox (-.25,-.15) rectangle (.4,.15);
\begin{scope}[yshift=-.05cm]
 \draw [->,dash pattern=on 2pt off 2pt]  (0,-.2) -- (0,.25);
\draw (-.05,.15) node[right]{\scriptsize K};
\draw [->]  (0,-.12) arc (270:180:.12);
\draw (0,.12) arc (90:180:.12);
\fill (0,-.12) circle (1pt) (0,.12) circle (1pt);
\end{scope}
\end{tikzpicture}}
\newcommand{\onechordsmallnonnumorrev}{\begin{tikzpicture} \useasboundingbox (-.25,-.15) rectangle (.4,.15);

\begin{scope}[yshift=-.1cm]
\draw [->,dash pattern=on 2pt off 2pt]  (0,-.2) -- (0,.3);
\draw (0,.2) node[right]{\scriptsize K};
\draw [-<]  (0,-.15) arc (270:180:.15);
\draw (0,.15) arc (90:180:.15);
\fill (0,-.15) circle (1pt) (0,.15) circle (1pt);
\end{scope}
\end{tikzpicture}}
\newcommand{\onechordsmall}[1]{\begin{tikzpicture} \useasboundingbox (-.3,-.1) rectangle (.7,.2);
\draw [->,dash pattern=on 2pt off 2pt]  (0,.2) arc (-270:90:.2);
\draw (-.2,0) -- (.2,0) (.5,0) node{\scriptsize #1};
\fill (-.2,0) circle (1.5pt) (.2,0) circle (1.5pt);
\end{tikzpicture}}
\newcommand{\onechordsmallh}[1]{\begin{tikzpicture} \useasboundingbox (-.3,-.1) rectangle (.6,.2);
\draw [->,dash pattern=on 2pt off 2pt]  (0,.2) arc (-270:90:.2);
\draw (-.2,0) -- (.2,0) (.35,.2) node{\scriptsize #1};
\fill (-.2,0) circle (1.5pt) (.2,0) circle (1.5pt);
\end{tikzpicture}}

\newcommand{\onechordsmallj}{\onechordsmall{$S^1_j$}}
\newcommand{\onechordsmallK}{\onechordsmall{$K$}}
\newcommand{\onechordsmalltseul}{\begin{tikzpicture} \useasboundingbox (-.25,-.15) rectangle (.2,.15);\begin{scope}[yshift=-.05cm]
\draw [->,dash pattern=on 2pt off 2pt]  (0,-.2) -- (0,.25);
\draw (0,-.12) arc (270:180:.12);
\draw (0,.12) arc (90:180:.12);
\fill (0,-.12) circle (1pt) (0,.12) circle (1pt);
\end{scope}
\end{tikzpicture}}
\newcommand{\onechordsmalltseullong}{\begin{tikzpicture} \useasboundingbox (-.25,-.2) rectangle (.2,.2);\begin{scope}[yshift=-.05cm]
\draw [->,dash pattern=on 2pt off 2pt]  (0,-.25) -- (0,.3);
\draw (0,-.12) arc (270:180:.12);
\draw (0,.12) arc (90:180:.12);
\fill (0,-.12) circle (1pt) (0,.12) circle (1pt);
\end{scope}
\end{tikzpicture}}

\newcommand{\onechordsmalltK}{\begin{tikzpicture} \useasboundingbox (-.2,-.15) rectangle (.4,.15);\begin{scope}[yshift=-.05cm]
\draw [->,dash pattern=on 2pt off 2pt]  (0,-.2) -- (0,.25);
\draw (-.05,.15) node[right]{\scriptsize $K$};
\draw (0,-.12) arc (270:180:.12);
\draw (0,.12) arc (90:180:.12);
\fill (0,-.12) circle (1pt) (0,.12) circle (1pt);
\end{scope}
\end{tikzpicture}}

\newcommand{\locnonor}{\begin{tikzpicture} \useasboundingbox (-.3,-.2) rectangle (.3,.2);
\draw (-.15,.15)  node{\scriptsize u} (+.2,-.25)  node{\scriptsize C};
\draw [dash pattern=on 2pt off 2pt]  (0,-.3) -- (0,.3);
\fill (0,0) circle (1.5pt);
\end{tikzpicture}}
\newcommand{\locleftup}{\begin{tikzpicture} \useasboundingbox (-.5,-.2) rectangle (.2,.2);
\draw (-.3,0) -- (0,0);
\draw [dash pattern=on 2pt off 2pt,->]  (0,-.3) -- (0,.3);
\fill (0,0) circle (1.5pt);
\end{tikzpicture}}
\newcommand{\locupor}{\begin{tikzpicture} \useasboundingbox (-.5,-.2) rectangle (.2,.2);
\draw (-.3,0) -- (0,0);
\draw [->] (-75:.15) arc (-75:75:.15);
\draw [->] (105:.15) arc (105:165:.15);
\draw [->] (-165:.15) arc (-165:-105:.15);
\draw [dash pattern=on 2pt off 2pt]  (0,-.3) -- (0,.3);
\fill (0,0) circle (1.5pt);
\end{tikzpicture}}
\newcommand{\locleftdown}{\begin{tikzpicture} \useasboundingbox (-.5,-.2) rectangle (.5,.2);
\draw (-.3,0) .. controls (-.2,0) and (-.1,-.15) .. (0,-.15)  .. controls (.1,-.15) and (.3,0) .. (0,0);
\draw [dash pattern=on 2pt off 2pt,<-]  (0,-.3) -- (0,.3);
\fill (0,0) circle (1.5pt);
\end{tikzpicture}}
\newcommand{\locdownor}{\begin{tikzpicture} \useasboundingbox (-.5,-.2) rectangle (.5,.2);
\draw (-.3,0) .. controls (-.2,0) and (-.1,-.15) .. (0,-.15)  .. controls (.1,-.15) and (.3,0) .. (0,0);
\draw [->] (15:.15) arc (15:75:.15);
\draw [->] (105:.15) arc (105:180:.15);
\draw [dash pattern=on 2pt off 2pt]  (0,-.3) -- (0,.3);
\fill (0,0) circle (1.5pt);
\end{tikzpicture}}

\newcommand{\univv}{\begin{tikzpicture} \useasboundingbox (-.5,.4) rectangle (.5,-.2);
\draw (-.15,.35) -- (0,.2) -- (0,-.1)  (.15,.35) -- (0,.2);
\draw [dashed, ->] (-.4,-.1) -- (.4,-.1);
\fill (0,-.1) circle (1.5pt) (0,.2) circle (1.5pt);
\end{tikzpicture}}
\newcommand{\Nunivv}{\begin{tikzpicture} \useasboundingbox (-.5,.4) rectangle (.5,-.2);
\draw [draw=gray!50,line width=3pt] (-.25,-.1) -- (.25,-.1) (-.15,.35) -- (0,.2) -- (0,-.1)  (.15,.35) -- (0,.2);
\draw (-.15,.35) -- (0,.2) -- (0,-.1) (.15,.35) -- (0,.2);
\draw [dash pattern=on 2pt off 2pt, ->] (-.4,-.1) -- (.4,-.1);
\fill (0,-.1) circle (1.5pt) (0,.2) circle (1.5pt);
\end{tikzpicture}} 
\newcommand{\Nunivvhopf}{\begin{tikzpicture} \useasboundingbox (-.5,.4) rectangle (.5,-.2);
\draw [draw=gray!50,line width=3pt] (-.15,.45) -- (0,.3) -- (0,.1) (.15,.45) -- (0,.3);
\draw (-.4,-.1) -- (-.2,-.1);
\draw [draw=gray!50,line width=3pt] (0,.1) arc (90:350:.2) (0,.1) arc (90:10:.2);
\fill[white] (.1,-.2) rectangle (.3,0);
\draw [->] (-.1,-.1) -- (.4,-.1);
\draw (-.15,.45) -- (0,.3) -- (0,.1) (.15,.45) -- (0,.3);
\draw [->] (0,.1) arc (90:20:.2) (0,.1) arc (90:340:.2);
\fill (0,.3) circle (1.5pt);
\end{tikzpicture}}
\newcommand{\univvhopf}{\begin{tikzpicture} \useasboundingbox (-.5,.4) rectangle (.5,-.2);
\draw (-.4,-.1) -- (-.3,-.1);
\draw [->] (-.1,-.1) -- (.4,-.1);
\draw (-.15,.45) -- (0,.3) -- (0,.1) (.15,.45) -- (0,.3);
\draw [->] (0,.1) arc (90:20:.2) (0,.1) arc (90:340:.2);
\fill (0,.3) circle (1.5pt);
\end{tikzpicture}}
\newcommand{\Nchord}{\begin{tikzpicture} [baseline=0cm] \useasboundingbox (-.4,.5) rectangle (.4,-.3);
\fill[gray!50] (-.3,-.2) rectangle (.3,.4);
\draw (-.2,.1) -- (.2,.1);
\draw [dashed, ->] (.2,-.3) -- (.2,.5);
\draw [dashed, ->] (-.2,.5) -- (-.2,-.3);
\fill (-.2,.1) circle (1.5pt)  (.2,.1) circle (1.5pt);
\fill[white] (-.1,.2) rectangle (.1,.5) (-.1,0) rectangle (.1,-.3);
\end{tikzpicture}} 
\newcommand{\Nchordrotate}{\begin{tikzpicture}[rotate=90]
\fill[gray!50] (-.3,-.2) rectangle (.3,.4);
\draw (-.2,.1) -- (.2,.1);
\draw [dashed, ->] (.2,-.3) -- (.2,.5);
\draw [dashed, ->] (-.2,.5) -- (-.2,-.3);
\fill (-.2,.1) circle (1.5pt)  (.2,.1) circle (1.5pt);
\fill[white] (-.1,.2) rectangle (.1,.5) (-.1,0) rectangle (.1,-.3);
\end{tikzpicture}}

\newcommand{\Nchordnegwithoutsurf}{\begin{tikzpicture} [baseline=0cm] \useasboundingbox (-.4,.5) rectangle (.4,-.3);
\draw [->] (.2,-.3) -- (.2,.5);
\draw [->] (-.2,.5) -- (-.2,-.3);
\end{tikzpicture}} 
\newcommand{\Nchordpos}{\begin{tikzpicture} [baseline=0cm] \useasboundingbox (-.3,.5) rectangle (.5,-.3);
\draw [->] (.2,-.3) -- (.2,.2) (.2,.4) -- (.2,.5);
\draw [->] (-.2,.5) .. controls (-.2,.3) .. (.2,.3) .. controls (.4,.3) .. (.4,.1) .. controls (.4,-.1) .. (.3,-.1) (.1,-.1) .. controls (-.2,-.1) .. (-.2,-.3);
\end{tikzpicture}}
\newcommand{\Nchordsing}{\begin{tikzpicture} [baseline=0cm] \useasboundingbox (-.3,.5) rectangle (.5,-.3);
\draw [->] (.2,-.3) -- (.2,.5);
\draw [->] (-.2,.5) .. controls (-.2,.3) .. (.2,.3) .. controls (.4,.3) .. (.4,.1) .. controls (.4,-.1) .. (.3,-.1) (.1,-.1) .. controls (-.2,-.1) .. (-.2,-.3);
\fill (.2,.3) circle (1.5pt);
\end{tikzpicture}}
\newcommand{\edget}{\begin{tikzpicture} \useasboundingbox (-1.15,-.2) rectangle (.35,.4);
\draw (-1,.1) -- (-.8,-.1) -- (-1,-.3) (-.8,-.1) -- (0,-.1) -- (.2,.1) (0,-.1) -- (.2,-.3);
\fill (-.8,-.1) circle (1.5pt) (0,-.1) circle (1.5pt);
\end{tikzpicture}}
\newcommand{\Nedget}{\begin{tikzpicture} \useasboundingbox (-1.2,-.2) rectangle (.4,.4);
\fill[gray!50] (-1,-.3) rectangle (.2,.1);
\draw (-1,.1) -- (-.8,-.1) -- (-1,-.3) (-.8,-.1) -- (0,-.1) -- (.2,.1) (0,-.1) -- (.2,-.3);
\fill (-.8,-.1) circle (1.5pt) (0,-.1) circle (1.5pt);
\fill[white] (-1,0) -- (-1,-.2) -- (-.9,-.1) -- (-1,0) (-.8,-.3) -- (-.7,-.2) -- (-.1,-.2) -- (0,-.3) -- (-.8,-.3) (.2,0) -- (.2,-.2) -- (.1,-.1) -- (.2,0) (-.8,.1) -- (-.7,0) -- (-.1,0) -- (0,.1) -- (-.8,.1);
\end{tikzpicture}}
\newcommand{\edgehopf}{\begin{tikzpicture} \useasboundingbox (-1.2,-.2) rectangle (.4,.4);
\draw [->] (-.3,-.1) arc (0:90:.2);
\draw [->] (-.1,-.1) arc (0:90:.2);
\draw (-.1,-.1) arc (0:-200:.2) (-.5,.1) arc (90:270:.2) (-.3,-.1) arc (0:-20:.2);
\draw (-1.1,.1) -- (-.9,-.1) -- (-1.1,-.3) (-.9,-.1) -- (-.7,-.1) (-.1,-.1) -- (.1,-.1) -- (.3,.1) (.1,-.1) -- (.3,-.3);
\fill (-.9,-.1) circle (1.5pt) (.1,-.1) circle (1.5pt);
\end{tikzpicture}}
\newcommand{\Nedgehopf}{\begin{tikzpicture} \useasboundingbox (-1.2,-.2) rectangle (.4,.4);
\fill[gray!50] (-1.1,-.3) rectangle (-.7,.1);
\fill[gray!50] (-.1,-.3) rectangle (.3,.1);
\fill[white] (-1.1,0) -- (-1.1,-.2) -- (-1,-.1) -- (-1.1,0) (-.9,-.3) -- (-.8,-.2) -- (0,-.2) -- (.1,-.3) -- (-.9,-.3) (.3,0) -- (.3,-.2) -- (.2,-.1) -- (.3,0) (-.9,.1) -- (-.8,0) -- (0,0) -- (.1,.1) -- (-.9,.1);
\draw [draw=gray!50,line width=3pt] (-.3,-.1) arc (0:90:.2) (-.1,-.1) arc (0:90:.2) (-.1,-.1) arc (0:-200:.2);
\draw [draw=gray!50,line width=3pt] (-.5,.1) arc (90:270:.2);
\draw [draw=gray!50,line width=3pt]  (-.3,-.1) arc (0:-20:.2);
\draw [->] (-.3,-.1) arc (0:90:.2);
\draw [->] (-.1,-.1) arc (0:90:.2);
\draw (-.1,-.1) arc (0:-200:.2) (-.5,.1) arc (90:270:.2) (-.3,-.1) arc (0:-20:.2);
\draw (-1.1,.1) -- (-.9,-.1) -- (-1.1,-.3) (-.9,-.1) -- (-.7,-.1) (-.1,-.1) -- (.1,-.1) -- (.3,.1) (.1,-.1) -- (.3,-.3);
\fill (-.9,-.1) circle (1.5pt) (.1,-.1) circle (1.5pt);
\end{tikzpicture}}
\newcommand{\smally}{\begin{tikzpicture} \useasboundingbox (-.4,-.2) rectangle (.4,.2);
\draw (150:.15) -- (0,0) -- (-90:.15) (0,0) -- (30:.15) (30:.25) circle (.1) (150:.25) circle (.1) (-90:.25) circle (.1);
\fill (0,0) circle (1pt);
\end{tikzpicture}}

\newcommand{\normyep}{\begin{tikzpicture} \useasboundingbox (-.6,-.3) rectangle (.6,.5);
\draw [draw=gray!50,line width=3pt] (150:.2) -- (0,0) -- (-90:.2) (0,0) -- (30:.2) (30:.35) circle (.15) (150:.35) circle (.15) (-90:.35) circle (.15);
\draw (150:.2) -- (0,0) -- (-90:.2) (0,0) -- (30:.2) (30:.35) circle (.15) (150:.35) circle (.15) (-90:.35) circle (.15);
\fill (0,0) circle (1pt);
\end{tikzpicture}}
\newcommand{\tancuphigh}{\begin{tikzpicture} [baseline=-.1 cm]
\useasboundingbox (-.25,-.3) rectangle (.25,.4);
\draw [->] (-.12,.15) .. controls (-.12,-.1) and (.12,-.1) .. (.12,.15);
\draw (-.2,.25) node{\scriptsize \textrm{(}} (.2,.25) node{\scriptsize \textrm{)}};
\fill (-.12,.25) circle (1pt) (.12,.25) circle (1pt);
\end{tikzpicture}}
\newcommand{\tancuphighwep}{\begin{tikzpicture}
\useasboundingbox (-.25,-.1) rectangle (.25,.2);
\draw [->] (-.12,.15) .. controls (-.12,-.1) and (.12,-.1) .. (.12,.15);
\end{tikzpicture}}

\newcommand{\tancuplow}{\begin{tikzpicture} [baseline=0cm]
\useasboundingbox (-.3,-.1) rectangle (.3,.35);
\draw [->] (-.12,.15) .. controls (-.12,-.1) and (.12,-.1) .. (.12,.15);
\draw (-.2,.25) node{\scriptsize \textrm{(}} (.2,.25) node{\scriptsize \textrm{)}};
\fill (-.12,.25) circle (1pt) (.12,.25) circle (1pt);
\end{tikzpicture}}
\newcommand{\tancapred}{\begin{tikzpicture}
\useasboundingbox (-.3,-.25) rectangle (.3,.1);
\draw [->] (-.12,-.15) .. controls (-.12,.1) and (.12,.1) .. (.12,-.15);
\draw (-.2,-.25) node{\scriptsize \textrm{(}} (.2,-.25) node{\scriptsize \textrm{)}};
\fill (-.12,-.25) circle (1pt) (.12,-.25) circle (1pt);
\end{tikzpicture}}
\newcommand{\tancap}{\begin{tikzpicture}  [baseline=-.1 cm]
\useasboundingbox (-.25,-.3) rectangle (.25,.4);
\draw [->] (-.12,-.15) .. controls (-.12,.1) and (.12,.1) .. (.12,-.15);
\draw (-.2,-.25) node{\scriptsize \textrm{(}} (.2,-.25) node{\scriptsize \textrm{)}};
\fill (-.12,-.25) circle (1pt) (.12,-.25) circle (1pt);
\end{tikzpicture}}
\newcommand{\tancapwep}{\begin{tikzpicture}
\useasboundingbox (-.25,-.1) rectangle (.25,.2);
\draw [->] (-.12,-.15) .. controls (-.12,.1) and (.12,.1) .. (.12,-.15);
\end{tikzpicture}}
\newcommand{\smaxo}{\begin{tikzpicture} [baseline=-.1 cm]
\useasboundingbox (-.2,-.2) rectangle (.2,.2);
\draw [->] (-.15,-.15) .. controls (-.15,.3) and (.15,.3) .. (.15,-.15);
\end{tikzpicture}}
\newcommand{\smaxno}{\begin{tikzpicture} [baseline=-.1 cm]
\useasboundingbox (-.2,-.2) rectangle (.2,.4);
\draw [->] (-.15,-.15) -- (-.07,-.07)  (.07,.07) .. controls (.15,.15) .. (.15,.2) .. controls (.15,.35) and (-.15,.35) .. (-.15,.2) .. controls (-.15,.15) .. (.15,-.15);
\end{tikzpicture}}
\newcommand{\assocnonnum}{\begin{tikzpicture} [baseline=-.1 cm]
\useasboundingbox (-.4,-.3) rectangle (.4,.4);
\draw [->] (-.2,-.15) -- (-.2,.15);
\draw [->] (.2,-.15) -- (.2,.15);
\draw [->] (-.1,-.15) -- (.1,.15);
\draw (-.33,-.25) node{\scriptsize \textrm{((}} (0,-.25) node{\scriptsize  \textrm{)}} (.27,-.25) node{\scriptsize \textrm{)}};
\fill (-.2,-.25) circle (1pt) (.2,-.25) circle (1pt)  (-.1,-.25) circle (1pt);
\draw (-.27,.25) node{\scriptsize \textrm{(}}  (0,.25) node{\scriptsize \textrm{(}} (.33,.25) node{\scriptsize  \textrm{))}};
\fill (-.2,.25) circle (1pt) (.2,.25) circle (1pt) (.1,.25) circle (1pt);
\end{tikzpicture}}
\newcommand{\tanposprod}{\begin{tikzpicture}  [baseline=-.1 cm]
\useasboundingbox (-.3,-.3) rectangle (.3,.4);
\draw [->] (.15,-.15) -- (-.15,.15);
\draw [draw=white,double=black,very thick]  (-.15,-.15) -- (.15,.15);
\draw [->]  (-.15,-.15) -- (.15,.15);
\draw (-.23,-.25) node{\scriptsize \textrm{(}} (.23,-.25) node{\scriptsize \textrm{)}};
\fill (-.15,-.25) circle (1pt) (.15,-.25) circle (1pt);
\draw (-.23,.25) node{\scriptsize \textrm{(}} (.23,.25) node{\scriptsize \textrm{)}};
\fill (-.15,.25) circle (1pt) (.15,.25) circle (1pt);
\end{tikzpicture}}
\newcommand{\assocnonnumwep}{\begin{tikzpicture}
\useasboundingbox (-.4,-.1) rectangle (.4,.2);
\draw [->] (-.2,-.15) -- (-.2,.15);
\draw [->] (.2,-.15) -- (.2,.15);
\draw [->] (-.1,-.15) -- (.1,.15);
\draw (-.33,-.25) node{\scriptsize \textrm{(}} (0,-.25) node{\scriptsize  \textrm{)}};
\fill (-.2,-.25) circle (1pt) (.2,-.25) circle (1pt)  (-.1,-.25) circle (1pt);
\draw (0,.25) node{\scriptsize \textrm{(}} (.33,.25) node{\scriptsize  \textrm{)}};
\fill (-.2,.25) circle (1pt) (.2,.25) circle (1pt) (.1,.25) circle (1pt);
\end{tikzpicture}}
\newcommand{\assocnum}{\begin{tikzpicture} [baseline=-.15 cm]
\useasboundingbox (-.55,-.4) rectangle (.55,.4);
\draw [->] (-.4,-.1) node{\scriptsize $1$} (-.3,-.25) -- (-.3,.25);
\draw [->] (.4,-.1) node{\scriptsize $3$} (.3,-.25) -- (.3,.25);
\draw [<-] (.1,.1) node{\scriptsize $2$} (.2,-.25) -- (-.2,.25);
\draw (-.44,.4) node{\scriptsize \textrm{((}} (-.12,.4) node{\scriptsize \textrm{)}} (.38,.4) node{\scriptsize \textrm{)}};
\fill (-.3,.4) circle (1pt) (.3,.4) circle (1pt)  (-.2,.4) circle (1pt);
\draw (-.38,-.35) node{\scriptsize \textrm{(}}  (.12,-.35) node{\scriptsize \textrm{(}} (.44,-.35) node{\scriptsize \textrm{))}};
\fill (-.3,-.35) circle (1pt) (.3,-.35) circle (1pt) (.2,-.35) circle (1pt);
\end{tikzpicture}}
\newcommand{\sourceass}{\begin{tikzpicture}
\useasboundingbox (-.4,-.15) rectangle (.5,.5);
\draw [->] (-.2,0) node{\scriptsize $1$} (-.3,-.25) -- (-.3,.25);
\draw [->] (.4,0) node{\scriptsize $3$} (.3,-.25) -- (.3,.25);
\draw [->] (.1,0) node{\scriptsize $2$} (0,-.25) -- (0,.25);
\end{tikzpicture}}
\newcommand{\sourceassone}[1]{\begin{tikzpicture}
\useasboundingbox (-.3,-.15) rectangle (0,.5);
\draw [->] (-.12,0) node{\tiny $#1$} (-.2,-.15) -- (-.2,.15);
\end{tikzpicture}}
\newcommand{\sourceassr}{\begin{tikzpicture} [baseline=0cm]
\useasboundingbox (-1,-.15) rectangle (1.4,.2);
\draw [->] (-.5,.2) node{\scriptsize $1$} (-.7,0) -- (-.3,.0);
\draw [->] (.2,.2) node{\scriptsize $2$} (0,0) -- (.4,.0);
\draw [->] (.9,.2) node{\scriptsize $3$} (.7,0) -- (1.1,.0);
\draw [dotted,->] (-1,0) -- (-.7,0) (-.3,.0) -- (0,0) (.4,.0) -- (.7,0) (1.1,.0) -- (1.4,0);
\end{tikzpicture}}
\newcommand{\snake}{\begin{tikzpicture} [baseline=-.2cm]
\useasboundingbox (-.35,-.5) rectangle (.35,.35);
\draw [->]  (-.2,-.4) -- (-.2,.15) .. controls (-.2,.3) and (-.2,.3) .. (-.1,.15) -- (.1,-.15)
.. controls (.2,-.3) and (.2,-.3) .. (.2,-.15) -- (.2,.4);
\end{tikzpicture}}
\newcommand{\snakedashed}{\begin{tikzpicture} [baseline=-.2cm]
\useasboundingbox (-.4,-.5) rectangle (.4,.35);
\draw [->]  (-.2,-.4) -- (-.2,.15) .. controls (-.2,.3) and (-.2,.3) .. (-.1,.15) -- (.1,-.15)
.. controls (.2,-.3) and (.2,-.3) .. (.2,-.15) -- (.2,.4);
\draw[dashed] (-.4,-.12) -- (.4,-.12)  (-.4,.12) -- (.4,.12);
\end{tikzpicture}}
\newcommand{\snakestr}{\begin{tikzpicture} [baseline=-.2cm]
\useasboundingbox (-.2,-.5) rectangle (.2,.35);
\draw [->]  (0,-.4) -- (0,.4);
\end{tikzpicture}}
\newcommand{\botsnake}{\begin{tikzpicture} [baseline=-.3cm]
\useasboundingbox (-.4,-.33) rectangle (.55,0);
\draw [->]  (-.2,-.4) -- (-.2,-.15);
\draw [->] (.1,-.15) .. controls (.1,-.35) and (.3,-.35) .. (.3,-.15);
\draw (-.27,-.05) node{\tiny \textrm{(}}  (0,-.05) node{\tiny \textrm{(}} (.43,-.05) node{\tiny \textrm{))}};
\fill (-.2,-.05) circle (1pt) (.3,-.05) circle (1pt) (.1,-.05) circle (1pt);
\end{tikzpicture}}
\newcommand{\topsnake}{\begin{tikzpicture} [baseline=.05cm]
\useasboundingbox (-.6,0) rectangle (.4,.35);
\draw [->]  (-.3,.15) .. controls (-.3,.35) and (-.1,.35) .. (-.1,.15);
\draw [->] (.2,.15) -- (.2,.4);
\draw (-.43,.05) node{\tiny \textrm{((}} (0,.05) node{\tiny \textrm{)}} (.27,.05) node{\tiny \textrm{)}};
\fill (-.3,.05) circle (1pt) (.2,.05) circle (1pt)  (-.1,.05) circle (1pt);
\end{tikzpicture}}
\newcommand{\guv}{\begin{tikzpicture} [baseline=-.1cm] \useasboundingbox (-.5,-.2) rectangle (.6,.2);
\draw (-.3,-.15) -- (-.3,.15) (.3,-.15) -- (.3,.15) (-.05,-.15) node{\scriptsize u} (.05,.15) node{\scriptsize v};
\draw [dashed, ->] (-.45,-.15) -- (-.15,-.15);
\draw [dashed, <-] (.15,.15) -- (.45,.15);
\fill (-.3,-.15) circle (1.5pt) (.3,.15) circle (1.5pt);
\end{tikzpicture}}
\newcommand{\guu}{\begin{tikzpicture} [baseline=-.15cm] \useasboundingbox (-.4,-.2) rectangle (.4,.2);
\draw (-.3,-.15) -- (-.3,.15) (0,-.15) -- (0,.15) (.25,-.15) node{\scriptsize u};
\draw [dashed, ->] (-.45,-.15) -- (.15,-.15);
\fill (-.3,-.15) circle (1.5pt) (0,-.15) circle (1.5pt);
\end{tikzpicture}}
\newcommand{\ddparallel}{\begin{tikzpicture} [baseline=-.1cm] \useasboundingbox (-.35,-.6) rectangle (.35,.6);
\draw [->,dash pattern=on 2pt off 2pt] (.25,-.5) -- (.25,.5);
\draw [<-,dash pattern=on 2pt off 2pt] (-.25,-.5) -- (-.25,.5);
\draw (-.25,-.15) -- (.25,-.15)  (-.25,.15) -- (.25,.15);
\fill (-.25,-.15) circle (1.5pt) (-.25,.15) circle (1.5pt) (.25,-.15) circle (1.5pt) (.25,.15) circle (1.5pt);
\end{tikzpicture}}
\newcommand{\ddcross}{\begin{tikzpicture} [baseline=-.1cm] \useasboundingbox (-.35,-.6) rectangle (.35,.6);
\draw [->,dash pattern=on 2pt off 2pt] (.25,-.5) -- (.25,.5);
\draw [<-,dash pattern=on 2pt off 2pt] (-.25,-.5) -- (-.25,.5);
\draw (-.25,-.15) -- (.25,.15)  (-.25,.15) -- (.25,-.15);
\fill (-.25,-.15) circle (1.5pt) (-.25,.15) circle (1.5pt) (.25,-.15) circle (1.5pt) (.25,.15) circle (1.5pt);
\end{tikzpicture}}

\title{Invariants of links and $3$-manifolds from graph configurations\\
\begin{tikzpicture} \useasboundingbox (-4,-4.8) rectangle (4,3.6);
\draw (-2,3) .. controls (-2,2) .. (-2,0) .. controls (-2,-2) and (-1.8,-3.8) .. (-1.2,-4.2) .. controls (-.6,-4.6) and (.6,-4.6) .. (1.2,-4.2) .. controls (1.8,-3.8) and (2,-2) .. (2,0) .. controls (2,2) .. (2,3);
\draw [very thick] (-2,0) -- (0,-1) -- (2,0) (1.2,-4.2) -- (0,-2.2) -- (-1.2,-4.2) 
(0,-2.2) -- (0,-1) -- (0,-.4) 
(60:.4)  .. controls (90:.2) .. (120:.4) (0,0) circle (.4);
\fill (0,-2.2) circle (.1) (0,-1) circle (.1) (0,-.4) circle (.1) (1.2,-4.2) circle (.1) (-1.2,-4.2) circle (.1) (2,0) circle (.1) (-2,0) circle (.1) (60:.4) circle (.1) (120:.4) circle (.1);
\end{tikzpicture}}

\author{Christine Lescop \thanks{Institut Fourier, CNRS, Universit\'e Grenoble Alpes, 38000 Grenoble, France}}
\begin{document}
\maketitle

\tableofcontents

\bigskip


\bigskip

\newpage
\thispagestyle{empty}

\vspace*{5cm}

\begin{flushright}
\`A Lucien Guillou,\\
\`A nos 26 ans de bonheur ensemble,\\
\`A nos deux enfants, Ga\"elle et Ronan.
\end{flushright}

\chapter*{Preface}

In this book, we explain how to count graph configurations to obtain invariants
for $3$-manifolds and knots in these $3$-manifolds, and we investigate the properties of the obtained invariants.

The simplest of these invariants is the linking number of two disjoint knots in the ambient space $\RR^3$. Gauss defined it in 1833 \cite{gauss}. As we review in Section~\ref{secslowbeg}, this linking number counts configurations

\begin{center}
\begin{tikzpicture}
\draw [blue] (.4,0) -- (.4,.2) .. controls (.4,.45) and (.6,.6).. (.8,.6) .. controls (1,.6) and (1.4,.2) .. (1.4,-.1);
\draw [blue] (0,0) -- (0,.2) .. controls (0,.6) and (.5,.9) .. (.8,.9) .. controls (1.1,.9) and (1.7,.5) .. (1.7,-.1);
\draw [->,draw=gray] (.3,.2) -- (.1,.2) (-.1,.2) .. controls (-.2,.2) and (-.35,.15) .. (-.35,.05) .. controls (-.35,-.05) and (0,-.1) .. (.2,-.1) .. controls (.4,-.1) and (.75,-.05) .. (.75,.05) .. controls (.75,.15) and (.6,.2) .. (.5,.2);
\draw[blue,->] (.2,-.3) .. controls (-.1,-.3) and (.2,-1.1) .. (.8,-1.1) .. controls (1.1,-1.1) and (1.7,-.7) .. (1.7,-.1);
\draw [draw=white,double=blue,very thick] (0,-.2) .. controls (0,-.3) and (.1,-.5) .. (.2,-.5) .. controls (.3,-.5) and (.4,-.3) ..  (.4,-.2);
\draw [draw=white,double=blue,very thick]  (.2,-.3) .. controls (.3,-.3) and (.5,-.8) .. (.8,-.8) .. controls (1.1,-.8) and (1.4,-.4) .. (1.4,-.1);
\draw (.65,.35) node{\scriptsize $J$} (1.7,-.1) node[right]{\scriptsize $K$};
\fill (.75,.05) circle (1.5pt) (1.4,-.1) circle (1.5pt) ;
\draw (.75,.05) -- (1.4,-.1) (2.85,0) node{of} (6.7,0) node{,};
\fill (4.5,.0) circle (1.5pt) (5.5,.0) circle (1.5pt);
\draw (4.5,.0) -- (5.5,.0); 
\draw [blue, dashed] (6,0) circle (.5);
\draw [gray, dashed] (4,0) circle (.5);
\draw (4,-.5)  node[below]{\scriptsize $J$} (6,-.5)  node[below]{\scriptsize $K$};
\end{tikzpicture}
\end{center}
as the degree of an associated Gauss map.

Many mysterious knot invariants called \say{quantum invariants} were introduced in the mid-80s, starting with the Jones polynomial. Witten explained how to obtain many of them from the perturbative expansion of the Chern-Simons theory in a seminal article \cite{witten}. This physicist viewpoint led Guadagnini, Martellini, Mintchev \cite{gmm} and Bar-Natan \cite{barnatanper} to show in what sense a coefficient $w_2$ of the Jones polynomial counts configurations of the graphs 
\begin{center}
 \begin{tikzpicture} \draw [gray,->,dash pattern=on 2pt off 2pt] (.5,0) arc (0:360:.5);
\draw (0,0) -- (-60:.5) node[right]{\scriptsize $v_1$};
\draw (0,0) -- (60:.5) node[right]{\scriptsize $v_2$};
\draw (.2,0) node{\scriptsize $w$}  (0,0) -- (180:.5) node[left]{\scriptsize $v_3$};
\draw [->] (-60:.5) -- (-60:.25);
\draw [->] (60:.5) -- (60:.25);
\draw [->] (180:.5) -- (180:.25);
\fill (60:.5) circle (1pt) (180:.5) circle (1pt) (-60:.5) circle (1pt) (0,0) circle (1pt);

\begin{scope}[xshift=2.8cm]
\draw [gray,->,dash pattern=on 2pt off 2pt] (.5,0) arc (0:360:.5);
\draw (45:.5) -- (-135:.5) node[left]{\scriptsize $b_2$};
\draw (135:.5) -- (-45:.5) node[right]{\scriptsize $b_3$};
\draw (135:.5) node[left]{\scriptsize $a_3$}  (45:.5)  node[right]{\scriptsize $a_2$};
\draw [->] (45:.5) -- (-135:.25);
\draw [->] (135:.5) -- (-45:.25);
\fill (45:.5) circle (1pt) (135:.5) circle (1pt) (-135:.5) circle (1pt) (-45:.5) circle (1pt);

\end{scope}
\draw (1.25,0) node{and} (4.5,0) node{such as} (8.3,0) node{and} (11.6,0) node{.};

\begin{scope}[xshift=8.4cm,yshift =-2.1cm]
\draw [>-,draw=gray, thick] (.8,1.5) .. controls (.8,1) and (1.7,1) .. (1.7,2);
\draw [draw=white,double=gray,very thick] (1.2,3.3) .. controls (1,3.3) and (.7,3.2) .. (.7,3) .. controls (.7,2.7) and (2,1.3) .. (2.2,1.5) .. controls (2.4,1.5) and (2.6,2) .. (2.6,2.2) .. controls (2.6,2.4) and (2.4,3.3) .. (2.2,3.3);
\draw [draw=white,double=gray,very thick] (2.2,3.3)  .. controls (1.5,3.3) and (.8,2) .. (.8,1.5);
\draw [draw=white,double=gray,very thick] (1.7,2)  .. controls (1.7,2.7) and (1.5,3.3) .. (1.2,3.3);
\draw  (.6,1.5) node{\scriptsize $K$}; 
\draw [black] (1.2,3.3) -- (2.2,3.3) ;
\draw [black] (2.6,2.2) -- (1.7,2);
\draw (1.2,3.3) node[above]{\scriptsize $b_2$} (2.6,2.2) node[right]{\scriptsize $b_3$} (2.2,3.3) node[above]{\scriptsize $a_2$} (1.95,2.1) node[below]{\scriptsize $a_3$};
\draw [->] (2.2,3.3) -- (1.7,3.3);
\draw [->] (1.7,2) -- (2.15,2.1);
\fill (1.7,2) circle (1.5pt) (2.6,2.2) circle (1.5pt) (2.2,3.3) circle (1.5pt) (1.2,3.3) circle (1.5pt);
\end{scope}

\begin{scope}[xshift=4.9cm,yshift =-2.1cm]
\draw [>-,draw=gray, thick] (.8,1.5) .. controls (.8,1) and (1.7,1) .. (1.7,2);
\draw [draw=white,double=gray,very thick] (1.2,3.3) .. controls (1,3.3) and (.7,3.2) .. (.7,3) .. controls (.7,2.7) and (2,1.3) .. (2.2,1.5) .. controls (2.4,1.5) and (2.6,2) .. (2.6,2.2) .. controls (2.6,2.4) and (2.4,3.3) .. (2.2,3.3);
\draw [draw=white,double=gray,very thick] (2.2,3.3)  .. controls (1.5,3.3) and (.8,2) .. (.8,1.5);
\draw [draw=white,double=gray,very thick] (1.7,2)  .. controls (1.7,2.7) and (1.5,3.3) .. (1.2,3.3);
\draw [draw=gray, thick, ->] (1.7,2)  .. controls (1.7,2.7) and (1.5,3.3) .. (1.2,3.3);
\draw  (.6,1.5) node{\scriptsize $K$}; 
\draw [draw=white,double=black,very thick] (3,3.2) -- (1.7,2);
\draw [black] (3,3.2) -- (2.2,3.3) ;
\draw [black] (3,3.2) -- (1.7,2) ;
\draw (3.2,3) node{\scriptsize $w$}  (3,3.2) -- (2.6,2.2) node[left]{\scriptsize $v_3$} (2.2,3.3) node[above]{\scriptsize $v_1$} (1.7,2) node[right]{\scriptsize $v_2$};
\draw [->] (2.2,3.3) -- (2.6,3.25);
\draw [->] (1.7,2) -- (2.35,2.6);
\draw [->] (2.6,2.2) -- (2.8,2.7);
\fill (1.7,2) circle (1.5pt) (2.6,2.2) circle (1.5pt) (2.2,3.3) circle (1.5pt) (3,3.2) circle (1.5pt);

\end{scope}
\end{tikzpicture}
 
\end{center}

The theory of Vassiliev invariants reviewed in Chapter~\ref{chapfintype} associates a degree in $(\NN \cup \{\infty\})$ to a numerical knot invariant. The only knot invariants of degree $0$ are the constant functions. The knot invariants of degree $2$ are linear combinations of $w_2$ and the constant function that maps every knot to $1$. The Jones polynomial can be renormalized into a series whose coefficients are finite-degree knot invariants. Altsch\"uler and Freidel showed that every degree $n$ real-valued knot invariant may be obtained by \say{counting} configurations
of graphs with at most $2n$ vertices as explained in this book \cite{af}.
The knot invariants counting graph configurations mentioned above are assembled in a \emph{universal Vassiliev invariant} $\Zinvuf(S^3,.)$ valued in a product of vector spaces generated by some unitrivalent graphs called \emph{Jacobi diagrams}. Kontsevich had constructed another universal Vassiliev invariant with similar properties called \emph{the Kontsevich integral} \cite{Kon,barnatan}. The Kontsevich integral $Z_K$ may be defined combinatorially from planar knot diagrams. It has been extensively studied. 
To my knowledge, the coincidence of the spatial invariant $\Zinvuf(S^3,.)$ with $Z_K$ is an open problem.

Developing the Witten approach further, Kontsevich outlined a way to count trivalent graphs in more general $3$-manifolds and define a topological graded invariant $\Zinvuf$ for them \cite{ko}.
These more general manifolds are the $3$-dimensional $\QQ$-spheres, simply called \emph{$\QQ$-spheres} in this book.
They are the closed $3$-manifolds with the same rational homology as the standard $3$-dimensional sphere $S^3$.
They include the $3$-manifolds with the same $\ZZ$-homology as $S^3$ called \emph{$\ZZ$-spheres}. A $\ZZ$-sphere is a closed orientable three-manifold in which knots bound orientable compact surfaces.
The degree one part $\Zinvuf_1$ of the Kontsevich invariant $\Zinvuf$ of $\QQ$-spheres is determined by a real-valued invariant $\Theta$ of $\QQ$-spheres, which counts configurations of the graph $\tata$ in the manifold.

Ohtsuki, Habiro, Goussarov, and others developed theories of finite type invariants of $\ZZ$-spheres analogous to the Vassiliev theory for knots in $\RR^3$ \cite{ohtkno,ggp}. Kuperberg and Thurston showed why the Kontsevich invariant $\Zinvuf$ of $\ZZ$-spheres obtained by counting graphs in these manifolds is also \emph{universal} with respect to the above theories \cite{kt}. This universality result implies that any real-valued invariant of $\ZZ$-spheres of degree $2n$  (with respect to one of the equivalent developed theories) is a combination of invariants counting configurations
of graphs with at most $2n$ vertices.

In 1985, Casson had defined an invariant of \emph{$\ZZ$-spheres}. The Casson invariant \say{counts} conjugacy classes of irreducible $SU(2)$-representations of the $\ZZ$-spheres fundamental groups.
Their universality result allowed Kuperberg and Thurston to show that $\Theta$ and the Casson invariant are proportional. In particular, the Casson invariant also \say{counts} configurations of the graph $\tata$.
For a $\ZZ$-sphere $\rats$ equipped with a basepoint $\infty$, the configurations are counted in a suitable compactification $C_2(\rats)$ of the space of ordered pairs of distinct points in the punctured manifold $(\crats=\rats \setminus \{\infty\})$. The set of counted configurations is the intersection in $C_2(\rats)$ of three transverse codimension-$2$ submanifolds called \emph{propagating chains}, and $\Theta(\rats)$ is their algebraic intersection number. Dually, the invariant $\Theta(\rats)$ is the integral over $C_2(\rats)$ of the cube of a \emph{propagating} closed $2$-form. Propagating chains and propagating forms both represent the linking form on $\rats$. We call them \emph{propagators}. They are the main ingredient used to count graph configurations in this book. They are associated with the graph edges. They are precisely defined in Chapter~\ref{chapprop}. When $\rats$ is $\ZZ$-sphere, results of Pontrjagin and Rohlin in the 1950s \cite{Roh} imply that the punctured $\crats$ can be equipped with a preferred homotopy class of parallelizations. For a general $\QQ$-sphere, the invariant $\Theta$ is first introduced as an invariant of a pair $(\rats, \tau)$, where $\tau$ is a parallelization of $\crats$. It is next corrected with the help of a relative first Pontrjagin class to become an invariant of $\rats$. Chapter~\ref{chapTheta} contains the complete construction of $\Theta$, and Chapter~\ref{chapfram} establishes the needed properties of Pontrjagin classes.

Kuperberg and I associated explicit propagating chains to Morse functions and associated Morse flows. These propagators reviewed in Section~\ref{subMorse} allowed me to express the Theta invariant in terms of Heegaard diagrams \cite{lesHC}. With this type of propagator, the \say{counted} graph configurations either map an edge of the graph into a flow line, or map the edge ends into descending manifolds or ascending manifolds of critical points of the Morse function. Fukaya proposed such a way of counting graphs \cite{FukayaMorse}. Many authors, including Watanabe and Shimizu, further studied it.

In the book's second part, we define and study an invariant $\Zinvuf(\rats,\Link)$ for a link $\Link$ in a \emph{$\QQ$-sphere} $\rats$. This invariant generalizes both $\Zinvuf(S^3,\Link)$ and $\Zinvuf(\rats)=\Zinvuf(\rats,\emptyset)$. Our definitions are more flexible than the original ones.
We prove generalizations of the mentioned universality results in the book's fourth part.

To get more properties of $\Zinvuf$, we cut our pairs $(\rats,\Link)$ of links $\Link$ in $\QQ$-spheres $\rats$ into pieces called \emph{tangles in $\QQ$-cylinders}.
These pieces can be composed in various ways. 
In the book's third part, we generalize $\Zinvuf$ to a functorial invariant of framed tangles in $\QQ$-cylinders, and we prove that it behaves well under the various allowed compositions.

Our first chapter is a more complete and much longer preface to this
book. It contains several introductions. 
Section \ref{secabs} is a short summary for experts. 
Other readers can start with Section~\ref{secslowbeg}, a slow informal introduction based on examples from which a broad audience
can get the flavor of the studied topics and hopefully become interested.
Section \ref{secmathov} is an independent, more formal, mathematical overview of the contents. It is accessible to beginners in topology after the warm-up of Section~\ref{secslowbeg}.
Section~\ref{secorganiz} describes the book organization. Section~\ref{secgenesis} outlines why I wrote this book and what I consider original and new. 

Apart from this preface and the first chapter, which has some parts written for experts and is sometimes imprecise, the rest of the book is precise, detailed, and mostly self-contained. The only prerequisites are basic notions of algebraic topology and de Rham cohomology, surveyed in the appendices.

In 2018, Watanabe disproved a long-standing conjecture called the $4$-dimensional Smale conjecture by constructing a topologically trivial $S^4$-bundle over $S^2$, which is not smoothly equivalent to the trivial bundle $S^4 \times S^2$ \cite{watanabe2018exotic}.
He distinguished his exotic $S^4 \times S^2$ from the standard $S^4 \times S^2$ using characteristic classes introduced by Kontsevich \cite{ko}. The involved Kontsevich--Watanabe characteristic class of a smooth topologically trivial $S^4$-bundle over $S^2$ counts configurations of the complete graph \tetrasmall\; with four vertices in the total space of the bundle. The ideas and techniques used by Watanabe are similar to those presented in this book. Even though we only count graph configurations in dimension $3$, this book can also serve as an introduction
to the work of Watanabe and other articles about invariants counting graph configurations in higher dimensions.

I thank the referees for their careful reading and their helpful comments.

\part{Introduction}

\chapter[Introductions]{Introductions}

In this chapter, we propose several introductions to this book:
\begin{itemize}
 \item a short abstract for experts in Section~\ref{secabs}, 
 \item a slow informal introduction based on examples in Section~\ref{secslowbeg}, from which a broad audience
can get the flavor of the studied topics and hopefully become interested,
 \item an independent, more formal, mathematical overview of the contents, to which the experts can go directly, in Section~\ref{secmathov}, and
 \item a section on the genesis of this book in Section~\ref{secgenesis}.
\end{itemize}

I apologize for some repetitions due to this structure.
Section~\ref{secorganiz} describes the book organization. The chapter ends with a list of open problems in Section~\ref{secopq}

Unlike this first chapter, which has some parts written for experts and is sometimes imprecise, the rest of the book is precise, detailed, and mostly self-contained. It relies only on the basic notions of algebraic topology and the basic notions of de Rham cohomology, surveyed in the appendices.

\section{An abstract for experts}
\label{secabs}

\paragraph{Very first conventions.}
Unless otherwise mentioned, manifolds are smooth, but they may have boundary and corners.
Let $\KK$ be $\ZZ$ or $\QQ$.
In this book, a \emph{$\KK$-sphere}\index[T]{Qsphere@$\QQ$-sphere}\index[T]{Zsphere@$\ZZ$-sphere} is a compact oriented $3$-dimensional manifold with the same homology with coefficients in $\KK$ as the standard unit sphere $S^3$ of $\RR^4$.
The unit disk of $\CC$ is denoted by $\drad{1}$.
A \emph{$\KK$-ball} (resp.  a \emph{$\KK$-cylinder}\index[T]{Qcylinder@$\QQ$-cylinder}, a \emph{genus $g$ $\KK$-handlebody}) is a compact oriented $3$-dimensional manifold $A$ with the same homology with coefficients in $\KK$ as the standard unit ball $B^3$ of $\RR^3$, (resp. the cylinder $\drad{1} \times \left[0,1\right]$, the standard solid handlebody $H_g$ of Figure~\ref{fighandlebodyg}), such that a neighborhood of the boundary of $A$ (which is necessarily homeomorphic to the boundary of its model---ball, cylinder or handlebody) is identified with a neighborhood of the boundary of its model by a smooth diffeomorphism.
\bfig
\centering
 \begin{tikzpicture} \useasboundingbox (.6,-1) rectangle (10.4,1); 
\draw [->] (1.8,.45) node[left]{\tiny $a_1$} (2,.9) .. controls (1.9,.9) and (1.75,.7) .. (1.75,.5);  
\draw (1.75,.5) .. controls (1.75,.3) and (1.9,.1) .. (2,.1);
\draw [dashed] (2,.9) .. controls (2.1,.9) and (2.25,.7) .. (2.25,.5) .. controls (2.25,.3) and (2.1,.1) .. (2,.1);
\draw plot[smooth] coordinates{(1.4,.1) (1.6,0) (2,-.1) (2.4,0) (2.6,.1)};
\draw plot[smooth] coordinates{(1.6,0) (2,.1) (2.4,0)};
\draw [dotted] (6.5,-.4) -- (7.5,-.4)  (6.5,.4) -- (7.5,.4);
\draw (7.5,.4) .. controls (8,.4) and (8.4,.9) .. (9,.9) .. controls (9.6,.9) and (10.4,.45) .. (10.4,0) .. controls (10.4,-.45) and (9.6,-.9) .. (9,-.9) .. controls (8.4,-.9) and (8,-.4) .. (7.5,-.4);
\draw (6.5,.4) .. controls (6,.4) and (5.6,.9) .. (5,.9) .. controls (4.4,.9) and (4,.4) .. (3.5,.4) .. controls (3,.4) and (2.6,.9) .. (2,.9) .. controls (1.4,.9) and (.6,.45) .. (.6,0) .. controls (.6,-.45) and (1.4,-.9) .. (2,-.9) .. controls (2.6,-.9) and (3,-.4) .. (3.5,-.4) .. controls (4,-.4) and (4.4,-.9) .. (5,-.9) .. controls (5.6,-.9) and (6,-.4) .. (6.5,-.4);
\begin{scope}[xshift=3cm]
 \draw [->] (1.8,.45) node[left]{\tiny $a_2$} (2,.9) .. controls (1.9,.9) and (1.75,.7) .. (1.75,.5);  
\draw (1.75,.5) .. controls (1.75,.3) and (1.9,.1) .. (2,.1);
\draw [dashed] (2,.9) .. controls (2.1,.9) and (2.25,.7) .. (2.25,.5) .. controls (2.25,.3) and (2.1,.1) .. (2,.1);
\draw plot[smooth] coordinates{(1.4,.1) (1.6,0) (2,-.1) (2.4,0) (2.6,.1)};
\draw plot[smooth] coordinates{(1.6,0) (2,.1) (2.4,0)};
\end{scope}
\begin{scope}[xshift=7cm]
 \draw [->] (1.8,.45) node[left]{\tiny $a_g$} (2,.9) .. controls (1.9,.9) and (1.75,.7) .. (1.75,.5);  
\draw (1.75,.5) .. controls (1.75,.3) and (1.9,.1) .. (2,.1);
\draw [dashed] (2,.9) .. controls (2.1,.9) and (2.25,.7) .. (2.25,.5) .. controls (2.25,.3) and (2.1,.1) .. (2,.1);
\draw plot[smooth] coordinates{(1.4,.1) (1.6,0) (2,-.1) (2.4,0) (2.6,.1)};
\draw plot[smooth] coordinates{(1.6,0) (2,.1) (2.4,0)};
\end{scope}
\end{tikzpicture}
\caption{The genus $g$ handlebody $H_g$}
\label{fighandlebodyg}
\end{figure}

In \say{$\QQ$-spheres}, $\QQ$ is a shortcut for \emph{rational homology}. So $\QQ$-spheres are also called \emph{rational homology spheres}\index[T]{Qsphere@$\QQ$-sphere} or \emph{rational homology $3$-spheres}\index[T]{rational homology!sphere}, while $\ZZ$-spheres are also called \emph{integer homology $3$-spheres.}\index[T]{integer homology sphere}

\paragraph{Abstract.}
In this book, following Edward Witten \cite{witten}, Maxim Kontsevich \cite{ko}, Greg Kuperberg and Dylan Thurston \cite{kt}, we define an invariant $\Zinvuf$ of $n$-component links $\Link$ in rational homology $3$-spheres $\rats$, and we study its properties.
The invariant $\Zinvuf$ is often called \say{the perturbative expansion of the Chern--Simons theory}. It is
valued in a graded space $\Aavis(\sqcup_{j=1}^k S^1)$ generated by Jacobi diagrams $\Gamma$ on $\sqcup_{j=1}^k S^1$. These diagrams are a special kind of Feynman diagrams. They are uni-trivalent.
The invariant $\Zinvuf(\Link)$ is a combination $\Zinvuf(\Link)=\sum_{\Gamma}\Zinvuf_{\Gamma}(\Link)\left[\Gamma\right]$ for coefficients $\Zinvuf_{\Gamma}(\Link)$ that \say{count} embeddings of $\Gamma$ in $\rats$ mapping the univalent vertices of $\Gamma$ to $\Link$, in a sense explained in the book. The coefficients $\Zinvuf_{\Gamma}(\Link)$ may be defined as integrals over configuration spaces or, in a dual way, as algebraic intersections in the same configuration spaces.

When $\rats=S^3$, the invariant $\Zinvuf$ is a universal Vassiliev link invariant studied by many authors, including Enore Guadagnini, Maurizio Martellini, and Mihail Mintchev \cite{gmm}, Dror Bar-Natan \cite{barnatanper}, Raoul Bott and Clifford Taubes \cite{botttaubes}, Daniel Altsch\"uler and Laurent Freidel \cite{af}, Dylan Thurston \cite{thurstonconf}, and Sylvain Poirier \cite{poirier}\dots This book contains a more flexible definition of this invariant.

\emph{Rational LP-surgeries} in $\QQ$-spheres are replacements of rational homology handlebodies by other rational homology handlebodies, in a way that does not change the linking number of curves outside the replaced handlebodies.
We prove that the restriction of $\Zinvuf$ to $\QQ$-spheres (equipped with empty links) is a universal finite type invariant with respect to rational LP-surgeries.
 Together with recent results of Gw\'ena\"el Massuyeau \cite{massuyeausplit} and Delphine Moussard \cite{moussardAGT}, this implies that the restriction of $\Zinvuf$ to $\QQ$-spheres contains the same information as the L\^e--Murakami--Ohtsuki LMO invariant \cite{lmo} for these manifolds.
This also implies that the degree one part of $\Zinvuf$ is the Casson--Walker invariant.

We extend $\Zinvuf$ to a functorial invariant of framed tangles in rational homology cylinders. We describe the behavior
of this functor under various operations including some cabling operations. We also compute iterated derivatives 
of this extended invariant with respect to discrete derivatives associated to the main theories of finite type invariants.

\section{A slow informal introduction for beginners}
\label{secslowbeg}

In this introduction, we describe the contents of this book to a broad audience, including graduate students.
We start with examples to give a flavor of the topics studied in this book and to introduce ideas, conventions, and methods used later.

This book is about invariants of links and 3-manifolds that count graph configurations.
The first example of such an invariant goes back to Carl Friedrich Gauss in 1833 \cite{gauss}. It is the linking number of two knots. We discuss it in Subsection~\ref{sublkGauss}.

At the end of the eighties, Edward Witten's insight into the perturbative expansion of the Chern--Simons theory \cite{witten} gave birth to trickier examples.
Following Witten's ideas, Enore Guadagnini, Maurizio Martellini, and Mihail Mintchev \cite{gmm}, and Dror Bar-Natan \cite{barnatanper} defined and studied another knot invariant $w_2$, which counts configurations of some uni-trivalent graphs with $4$ vertices. Maxim Kontsevich \cite{ko}, Raoul Bott and Clifford Taubes \cite{botttaubes} \cite{bottturkish}, Daniel Altsch\"uler and Laurent Freidel \cite{af}, and others revisited and generalized the definition of $w_2$.
In Subsection~\ref{subsecwtwo}, we detail the example of $w_2$.

\subsection{The linking number as a degree}
\label{sublkGauss}
Let $S^1$ denote the unit circle of the complex plane $\CC$. We use \say{$\suchthat$} as a symbol for \say{such that}. So, the circle 
$S^1$ is the set $\{z \suchthat z \in \CC,|z|=1\}$.
Consider a $C^{\infty}$ embedding 
 \begin{equation*}J\sqcup K \colon S^1\sqcup S^1\hookrightarrow \RR^3\end{equation*} 
of the disjoint union $S^1\sqcup S^1$ of two circles in the ambient space $\RR^3$, as the example pictured in Figure~\ref{figlinkdiagram}. Such an embedding represents a \emph{$2$-component link.}
\bfig
 \centering
\begin{tikzpicture}
\useasboundingbox (.5,.5) rectangle (4,3.5);
\draw [>-,draw=white,double=black,very thick] (.8,1.5) .. controls (.8,1) and (1.7,1) .. (1.7,2);
\draw [>-,draw=gray] (4,2) .. controls (4,3) and (2.5,2.7) .. (2.5,2.5) (3.5,2) .. controls (3.5,1.8) and (2.5,1.7) .. (2.5,1.5);
\draw [draw=white,double=black,very thick] (1.2,3.3) .. controls (1,3.3) and (.7,3.2) .. (.7,3) .. controls (.7,2.8) and (1.5,.7) .. (2.2,.7) .. controls (2.4,.7) and (3,.9) .. (3,1.1) -- (3,2.9) .. controls (3,3.1) and (2.4,3.3) .. (2.2,3.3);
\draw [>-,draw=white,double=gray,very thick] (2.5,2.5) .. controls (2.5,2.3) and (3.5,2.2) .. (3.5,2) (2.5,1.5) .. controls (2.5,1) and (4,1) .. (4,2);
\draw [draw=white,double=black,very thick] (2.2,3.3)  .. controls (1.5,3.3) and (.8,2) .. (.8,1.5);
\draw [draw=white,double=black,very thick] (1.7,2)  .. controls (1.7,2.7) and (1.5,3.3) .. (1.2,3.3);
\draw [>-] (.8,1.6)--(.8,1.3);
\draw [>-,draw=gray] (2.5,1.6)--(2.5,1.3);
\draw  (.6,1.5) node{\tiny $K$} (3.7,2) node{\tiny $J$};
\end{tikzpicture}

\caption{A $2$-component link in $\RR^3$}
\label{figlinkdiagram}
\end{figure}

\noindent It induces 
the {\em Gauss map\/}
\begin{center}
 \begin{tikzpicture} 
\begin{scope}[xshift=-5.1cm]
\draw [thick] (.4,0) .. controls (.4,-.5) and (1.4,-.9) .. (2,-.9) .. controls (2.6,-.9) and (3.6,-.5) .. (3.6,0) .. controls (3.6,.5) and (2.6,.9) .. (2,.9) .. controls (1.4,.9) and (.4,.5) .. (.4,0); 
\draw plot[smooth] coordinates{(1.4,.1) (1.6,0) (2,-.1) (2.4,0) (2.6,.1)};
\draw [out=25,in=155] (1.6,0) to (2.4,0);
\draw [->] (2.8,-.15) -- (3.1,-.15) node[right]{\tiny $1$};
\draw [->] (2.8,-.15) -- (2.8,.15) node[above]{\tiny $2$};
\end{scope}  
\draw (0,0) node {$\hfl{p_{JK}}$};
\begin{scope}[xshift=3.5cm] \useasboundingbox (-1,-1) rectangle (1,1);
\draw (0,0) circle (.9);
\draw (-.9,0)  .. controls (-.9,-.15) and (-.3,-.25) .. (0,-.25) .. controls (.3,-.25) and (.9,-.15) ..  (.9,0);
\draw [dashed] (-.9,0)  .. controls (-.9,.15) and (-.3,.25) .. (0,.25) .. controls (.3,.25) and (.9,.15) ..  (.9,0);
\draw [->] (-.15,-.7) -- (.15,-.7) node[right]{\tiny $1$};
\draw [->] (-.15,-.7) -- (-.15,-.4) node[left]{\tiny $2$};
\end{scope}
\draw (-5.8,-1.5) node {$p_{JK} \colon $} (-3.1,-1.5) node {$S^1 \times S^1$} (0,-1.5) node{$\rightarrow$} (3.5,-1.5) node{$S^2$}
(-3.1,-2.1) node {$(w,z)$} (0,-2.1) node{$\mapsto$} 
(3.5,-2.1) node{$\frac{1}{\lVert K(z)-J(w) \rVert}(K(z)-J(w)).$};
\end{tikzpicture}\end{center}

\begin{definition}
\label{defGausslk}
The {\em Gauss linking number\/} $lk_G(J,K)$ of the disjoint {\em knot embeddings\/} $J$ and $K$ is the degree of the Gauss map $p_{JK}$.
\end{definition}

Below, we give our favorite definition of the degree for this book. To do so, we first agree on conventions used throughout the book.

\subsection{On orientations of manifolds and degrees of maps}
\label{subdegree}

We work with smooth (i.e., $C^{\infty}$) manifolds locally diffeomorphic to open subspaces of $\left[0,1\right]^n$. These manifolds are precisely described in Sections~\ref{submanifolds} and \ref{secbackgroundmfdsbry}. They may have boundaries and ridges (or corners). The cube $\left[0,1\right]^3$ is an example of such a manifold. Its edges and its vertices are \emph{ridges}\index[T]{ridges}.

\begin{conventions}
Let $M$ be such a manifold. The \emph{interior} of $M$ consists of the points with a neighborhood diffeomorphic to an open subspace of $\RR^n$. The \emph{boundary} of $M$ is the complement of its interior in $M$.

An \emph{orientation}\index[T]{orientation!of a vector space} of a real vector space $V$ of positive dimension is a basis of $V$ up to a change of basis with positive determinant. When $V=\{0\}$, an orientation of $V$ is an element of $\{-1,1\}$. An  \emph{orientation}\index[T]{orientation!of a manifold} of a smooth $n$-manifold is an orientation of its tangent space at each point of its interior, defined in a continuous way.
(A local diffeomorphism $h$ of $\RR^n$ is \emph{orientation-preserving} at $x$ if and only if the Jacobian determinant of its derivative $T_xh$ is positive.
If the transition maps $\phi_j \circ \phi_i^{-1}$ of an \emph{atlas} $\left(\phi_i\right)_{i \in I}$  of a manifold $M$ (as in Subsection~\ref{submanifolds}) are orientation-preserving (at every point) for $\{i,j\} \subseteq I$, then the manifold $M$ is {\em oriented} by this atlas.) 
Unless otherwise mentioned, manifolds are smooth, oriented, compact, and considered up to orientation-preserving diffeomorphisms, in this book.
Products are oriented by the order of the factors. More generally, unless otherwise mentioned, the order of appearance of coordinates or parameters orients manifolds. When $M$ is an oriented manifold, $(-M)$ denotes the same manifold, equipped with the opposite orientation.

The boundary $\partial M$ of an oriented manifold $M$ is oriented by the \emph{outward normal first} convention: If $x \in \partial M$ is not in a ridge, then the outward normal to $M$ at $x$ followed by an oriented basis of $T_x \partial M$ induce the orientation of $M$. 
\end{conventions}

For example, the standard orientation of the disk in the plane induces the traditional counterclockwise orientation of the circle, as the following picture shows.
\begin{center}
\begin{tikzpicture}
\useasboundingbox (-.8,-.8) rectangle (.8,.8);
\draw [fill=gray!20] (0,0) circle (.7);
\draw [->] (-.15,-.3) -- (.15,-.3) node[right]{\tiny $1$};
\draw [->] (-.15,-.3) -- (-.15,0) node[left]{\tiny $2$};
\draw [->] (.7,0) -- (1,0) node[right]{\tiny $1$};
\draw [->] (.7,0) -- (.7,.3) node[left]{\tiny $2$};
\draw [thick,->] (-.7,0) arc (-180:180:.7);
\end{tikzpicture}
\end{center}
As another example, the sphere $S^2$ is oriented as the boundary of the unit ball $B^3$ of $\RR^3$, which has the standard orientation
induced by (thumb, index finger (2), middle finger (3)) of the right hand.

\begin{center}
 \begin{tikzpicture} \useasboundingbox (-1,-.7) rectangle (1,1);
\draw (0,0) circle (.9);
\draw (-.9,0)  .. controls (-.9,-.15) and (-.3,-.25) .. (0,-.25) .. controls (.3,-.25) and (.9,-.15) ..  (.9,0);
\draw [dashed] (-.9,0)  .. controls (-.9,.15) and (-.3,.25) .. (0,.25) .. controls (.3,.25) and (.9,.15) ..  (.9,0);
\draw [->] (-.15,-.7) -- (.15,-.7) node[right]{\tiny $2$};
\draw [->] (-.15,-.7) -- (-.15,-.4) node[left]{\tiny $3$};
\end{tikzpicture}
\end{center}

\begin{definitions}
Let $M$ and $N$ be smooth manifolds, and let $p \colon M \rightarrow N$ be a smooth map from $M$ to $N$.
Assume that the boundary of $M$ is empty.
A point $x$ of $M$ is called a \emph{regular point}\index[T]{regular point} of $p$ if the tangent map 
\begin{equation*}\tang_xp \colon \tang_x M \to \tang_y N\end{equation*} at $x$ is surjective.
A point $y$ of $N$ is called a \emph{regular value}\index[T]{regular value!of a smooth map} of $p$ if $p^{-1}(y)$ contains only regular points.
\end{definitions}

Our $d$-manifolds are covered by countably many open sets diffeomorphic to open subsets of $\left[0,1\right]^d$. The following \emph{Morse--Sard theorem} is proved in the book \cite{hirsch} by Morris W. Hirsch. See \cite[Chapter 3, Theorem 1.3, p. 69]{hirsch}.

\begin{theorem}[Morse--Sard theorem]
\label{thmMorseSard}
 The set of regular values of a smooth map from a manifold to another such is dense.\footnote{It is even \emph{residual}, i.e., it contains the intersection of a countable family of dense open sets.}
\end{theorem}
 
If the source $M$ of a smooth map $p \colon M \rightarrow N$ is compact, then the set of its regular values is furthermore open.
In general, the source may have boundary and ridges. Then a point $y$ is a \emph{regular value} \index[T]{regular value!for a general domain} of $p$, if $y$ is a regular value of the restrictions of $p$ to the interior of $M$ and to all the open, smooth faces (or strata) of $M$ (of any codimension). 
In particular, if the dimensions of $M$ and $N$ coincide, then a regular value is not in the image of the boundary $\partial M$.

\begin{definition}
\label{defdegr}
Assume that $M$ and $N$ are oriented, $M$ is compact, and the dimension of $M$ coincides with the dimension of $N$. Then the \emph{(differential) degree}\index[T]{degree!at a point} of $p$ at a regular value $y$ of $N$ is the (finite) sum running over the $x \in p^{-1}(y)$ of the signs of the determinants of $T_xp$ (in oriented charts).
In this case, the differential degree of $p$ extends as a continuous function $\deg(p)$ from $N \setminus p(\partial M)$ to $\ZZ$, as we prove in detail in Lemma~\ref{lemdeggen}.
In particular, if the boundary of $M$ is empty and if $N$ is connected, then $\deg(p)$ extends as a constant map from $N$ to $\ZZ$, whose value is called the \emph{degree}\index[T]{degree!of a map} of $p$. See \cite[Chapter 5]{Mil}.
\end{definition}

Figure~\ref{figdegvertproj} shows the values of $\deg(p)$ for the pictured vertical projection $p$ from the interval $\left[0,1\right]$ to $\RR$.

\bfig 
\centering
 \begin{tikzpicture}
\useasboundingbox (0,-.7) rectangle (6,1.5);
\draw [dashed] (1,0) -- (1,1.8);
\draw plot[smooth] coordinates{(4,1.5) (5,1.4) (3.6,1.2) (5.5,.95) (3,.8) (1,.7) (.5,.6) (1,.5) (2,.4)};
\draw [->] (3.15,.9) -- (3,.8) -- (3.15,.7) (0,0) -- (6,0);
\fill (2,.4) circle (1pt)  (4,1.5) circle (1pt) (1,.5)  circle (1pt)  (1,.7) circle (1pt);
\draw (1.2,.85) node{\tiny $-$} (1.2,.35) node{\tiny $+$}  (-.2,-.6) node[left]{\scriptsize values of $\deg(p)$:} (1,-.6) node{\scriptsize $0$} (3,-.6) node{\scriptsize $-1$} (5,-.6) node{\scriptsize $0$};
\draw [thin] (0,-.25) -- (1.8, -.25) -- (1.9,-.15) (2.1,-.15) -- (2.2,-.25) -- (3.8, -.25) -- (3.9,-.15)
(4.1,-.15) -- (4.2,-.25) -- (6, -.25);
\end{tikzpicture}
\caption{Degree of the vertical projection}
\label{figdegvertproj}
\end{figure}

Another easy example in higher dimensions is the case of an orientation-preserving embedding $p$. In this case, $\deg(p)$ is $1$ on the image of the interior of $M$ and $0$ outside the image of $M$.

\subsection{Back to the linking number}
\label{sublkGausstwo}

The Gauss linking number $lk_G(J,K)$ can be computed from a link diagram as in Figure~\ref{figlinkdiagram} as follows.
It is the differential degree of $p_{JK}$ at the vector $Y$ that points towards us. The set $p_{JK}^{-1}(Y)$ is the set of pairs $(w,z)$ of points for which the projections of $J(w)$ and $K(z)$ coincide, and $J(w)$ is under $K(z)$. They correspond to the {\em crossings\/} \pcbg and \ncbg of the diagram.

In a diagram, a crossing is \emph{positive} if we turn counterclockwise from the arrow at the end of the upper strand towards the arrow at the end of the lower strand like \pcortrig.
Otherwise, it is \emph{negative} like \nc.

Consider a positive crossing \pcbg. Moving $J(w)$ along $J$ following the orientation of $J$ moves $p_{JK}(w,z)$ towards the southeast direction $Tpdw$. In the meantime, moving $K(z)$ along $K$ following the orientation of $K$ moves $p_{JK}(w,z)$ towards the northeast direction $Tpdz$. So the local orientation \begin{tikzpicture}
\useasboundingbox (-.1,.1) rectangle (1.2,.55);
\draw [->] (0,.2) -- (.2,0) node[right]{\tiny $Tpdw$};
\draw [->,]  (0,.2) -- (.2,.4) node[right]{\tiny $Tpdz$};
\end{tikzpicture} induced by the image of $p_{JK}$ around $Y \in S^2$ is \begin{tikzpicture}
\useasboundingbox (-.1,.1) rectangle (.5,.55);
\draw [->] (0,.2) -- (.2,0) node[right]{\tiny $1$};
\draw [->,]  (0,.2) -- (.2,.4) node[right]{\tiny $2$};
\end{tikzpicture}.
Therefore, the contribution of a positive crossing to the degree is $1$. It is easy to deduce that the contribution of a negative crossing is $(-1)$.
We denote the cardinality of a set $A$ by $\cardlef{A}$. 

In particular, $\cardlef{ \pcbg }$ is the number of occurrences of $\pcbg$ in the diagram.
We have thus proved that \begin{equation*}\deg_Y(p_{JK})= \cardlef{ \pcbg } - \cardlef{ \ncbg}.\end{equation*}
 So we have
\begin{equation*}lk_G(J,K)= \cardlef{ \pcbg } - \cardlef{ \ncbg } .\end{equation*}
We similarly obtain $\deg_{-Y}(p_{JK})= \cardlef{ \pcbgkj } - \cardlef{  \ncbgkj}$.
This implies
\begin{equation*}\begin{array}{ll}lk_G(J,K)&= \cardlef{ \pcbgkj } - \cardlef{  \ncbgkj }\\
   & = \frac12\left(\cardlef{ \pcbg } +\cardlef{ \pcbgkj } - \cardlef{  \ncbg }- \cardlef{  \ncbgkj}\right)
  \end{array}
\end{equation*}
and $lk_G(J,K)=lk_G(K,J)$.

In the example of Figure~\ref{figlinkdiagram}, the Gauss linking number $lk_G(J,K)$ is $2$. For the \emph{positive Hopf link} of Figure~\ref{figHopfWhitehead}, we have $lk_G(J,K)=1$. The Gauss linking number of the components of the \emph{negative Hopf link} is equal to $-1$. It is zero for the \emph{Whitehead link}.

\bfig
 \centering
\begin{tikzpicture}
\begin{scope}[xshift=-3.3cm]
\begin{scope}[scale=1.5]
\draw[->] (0,0) -- (0,.2) .. controls (0,.4) and (-.15,.5).. (-.3,.5) .. controls(-.45,.5) and (-.6,.2) .. (-.6,0) node[left]{\scriptsize $K$};
\draw[->] (-.6,0) .. controls (-.6,-.2) and (-.45,-.5) .. (-.3,-.5) .. controls (-.15,-.5) and (0,-.35) .. (0,-.2);
\draw[->,draw=gray] (-.1,.2) .. controls (-.2,.2) and (-.25,.15) .. (-.25,.05) .. controls (-.25,-.05) and (-.2,-.1) .. (0,-.1) .. controls (.2,-.1) and (.25,-.05) .. (.25,.05) .. controls (.25,.15) and (.2,.2) .. (.1,.2);
\draw (.25,.35) node{\scriptsize $J$};\end{scope}
\end{scope}
\draw (-3.3,-1.5) node{\scriptsize The positive Hopf link};
\begin{scope}[scale=1.5]
\draw[->] (0,0) -- (0,.2) .. controls (0,.4) and (-.15,.5).. (-.3,.5) .. controls(-.45,.5) and (-.6,.2) .. (-.6,0) node[left]{\scriptsize $K$};
\draw[->] (-.6,0) .. controls (-.6,-.2) and (-.45,-.5) .. (-.3,-.5) .. controls (-.15,-.5) and (0,-.35) .. (0,-.2);
\draw[<-,draw=gray] (-.1,.2) .. controls (-.2,.2) and (-.25,.15) .. (-.25,.05) .. controls (-.25,-.05) and (-.2,-.1) .. (0,-.1) .. controls (.2,-.1) and (.25,-.05) .. (.25,.05) .. controls (.25,.15) and (.2,.2) .. (.1,.2);
\draw (-.25,.3)node{\scriptsize $J$};
\end{scope}
\draw  (0,-1.5) node{\scriptsize The negative Hopf link};
\begin{scope}[xshift=2.8cm]
\draw (.4,0) -- (.4,.2) .. controls (.4,.45) and (.6,.6).. (.8,.6) .. controls (1,.6) and (1.4,.2) .. (1.4,-.1);
\draw (0,0) -- (0,.2) .. controls (0,.6) and (.5,.9) .. (.8,.9) .. controls (1.1,.9) and (1.7,.5) .. (1.7,-.1);
\draw [->,draw=gray] (.3,.2) -- (.1,.2) (-.1,.2) .. controls (-.2,.2) and (-.35,.15) .. (-.35,.05) .. controls (-.35,-.05) and (0,-.1) .. (.2,-.1) .. controls (.4,-.1) and (.75,-.05) .. (.75,.05) .. controls (.75,.15) and (.6,.2) .. (.5,.2);
\draw[->] (.2,-.3) .. controls (-.1,-.3) and (.2,-1.1) .. (.8,-1.1) .. controls (1.1,-1.1) and (1.7,-.7) .. (1.7,-.1) node[right]{\scriptsize $K$};
\draw [draw=white,double=black,very thick] (0,-.2) .. controls (0,-.3) and (.1,-.5) .. (.2,-.5) .. controls (.3,-.5) and (.4,-.3) ..  (.4,-.2);
\draw [draw=white,double=black,very thick]  (.2,-.3) .. controls (.3,-.3) and (.5,-.8) .. (.8,-.8) .. controls (1.1,-.8) and (1.4,-.4) .. (1.4,-.1);
\draw (.65,.35) node{\scriptsize $J$} (1,-1.5) node{\scriptsize The Whitehead link};
\end{scope}
\end{tikzpicture}

\caption{The Hopf links and the Whitehead link}
\label{figHopfWhitehead}
\end{figure}
 
Let $\omega_S$ be a $2$-form on $S^2$. The integral of the pull-back $p_{JK}^{\ast}(\omega_S)$ of $\omega_S$ over $S^1 \times S^1$ is the integral over $S^2$ of $\deg(p_{JK})\omega_S$. 
Since the differential degree $\deg(p_{JK})$ of the Gauss map $p_{JK}$ is constant and equal to $lk_G(J,K)$ on the set of regular values of $p_{JK}$, we have
\begin{equation*}\int_{S^1 \times S^1}p_{JK}^{\ast}(\omega_S)=\int_{S^2} \deg(p_{JK})\omega_S=lk_G(J,K)\int_{S^2}\omega_S.\end{equation*}
Denote the standard area form of $S^2$ by $4\pi \omega_{S^2}$. So $\omega_{S^2}$ is the homogeneous volume form of $S^2$ such that
$\int_{S^2}\omega_{S^2}=1$.
In 1833, Carl Friedrich Gauss defined the linking number of $J$ and $K$ as
an integral \cite{gauss}.
 In modern notation, his definition may be written as
\begin{equation*}lk_G(J,K)=\int_{S^1 \times S^1}p_{JK}^{\ast}(\omega_{S^2}).\end{equation*}

\subsection{A first noninvariant count of graph configurations}
\label{subnewonechord}
The above Gauss linking number counts the \emph{configurations} of the graph \onechordsegment{v_J}{v_K}, for which $v_J$ is on the knot $J$ (or more precisely on the image $J(S^1)$ of the knot embedding $J$), the vertex $v_K$ is on the knot $K$, and the edge from $v_J$ to $v_K$ is a straight segment with an arbitrary generic fixed direction $X$ in $S^2$.
Here, \emph{generic} means that $X$ is a regular value of $p_{JK}$.\footnote{This is a generic condition thanks to the recalled Morse--Sard theorem \ref{thmMorseSard}.} 
A \indexT{configuration} of \onechordsegment{v_J}{v_K} is an injection $\confc \colon \{v_J,v_K\} \hookrightarrow (\RR^3)^2$ such that $\confc(v_J)=J(z_J)$ for some $z_J \in S^1$ and $\confc(v_K)=K(z_K)$ for some $z_K \in S^1$.
The corresponding \emph{configuration space}\index[T]{configuration!space} is parametrized by $(z_J,z_K) \in S^1 \times S^1$. It is diffeomorphic to $S^1 \times S^1$.
The configurations such that the edge from $v_J$ to $v_K$ is a straight segment with direction $X$ in $S^2$ are in one-to-one correspondence with $p_{JK}^{-1}(X)$. The local degree of $p_{JK}$ equips each of these configurations with a sign.
The \say{count} of configurations with direction $X$ is the sum of these signs, which is nothing but the degree of $p_{JK}$ at $X$. So any choice of a generic $X$ will give the same integral result, which will not be changed by a continuous deformation of our embedding among embeddings.
The graph \onechordsegment{v_J}{v_K} will also be denoted by \onechordtwocirclesrightmoyen{J}{K}. In this diagram, the dashed circles show to which component the vertices must map.

Let $K \colon S^1 \hookrightarrow \RR^3$ be a smooth embedding of the circle into $\RR^3$.
Such an embedding is called a \emph{knot embedding}. An \emph{isotopy} between two knot embeddings $K$ and $K_1$ is a smooth map $\psi \colon \left[0,1\right] \times S^1 \to \RR^3$ such that the restriction $\psi(t,.)$ of $\psi$ to  $\{t\} \times S^1$ is a knot embedding for any $t \in \left[0,1\right]$, 
$\psi(0,.)=K$, and $\psi(1,.)=K_1$. When there exists such an isotopy, the knot embeddings
$K$ and $K_1$ are said to be \emph{isotopic} or in the same \emph{isotopy class}. A \emph{knot} is an isotopy class of knot embeddings. 

Let us now try to count the configurations \begin{equation*}\confc \colon \{v_1,v_2\} \hookrightarrow (\RR^3)^2\end{equation*} of the graph \onechordsegment{v_1}{v_2},
for which $\confc(v_1)$ and $\confc(v_2)$ are two distinct points on the (image of) the knot embedding $K$, and the edge from $v_1$ to $v_2$ is a straight segment with an arbitrary direction $X$ in $S^2$.
The graph \onechordsegment{v_1}{v_2} is also denoted by $\tatak$. The associated \emph{configuration space} is
\begin{equation*}\check{C}(K;\tatak)=\Bigl\{\Bigl(K\bigl(z\bigr),K\bigl(z\exp(2i\pi t)\bigr)\Bigr) \suchthat z \in S^1, t \in \left]0,1\right[\Bigr\}.\end{equation*}
(In this book, open and half-open intervals are denoted with square brackets. For example, $\left]0,1\right[$ denotes the open interval $\{x \suchthat 0<x<1\}$. Similarly, $\left]0,1\right]$ denotes the half-open interval $\{x \suchthat 0<x\leq 1\}$.)
The configuration space $\check{C}(K;\tatak)$ is naturally identified
 with the open annulus $S^1 \times \left]0,1\right[$. We have a \emph{Gauss direction map}
 \begin{equation*}\begin{array}{llll}G_K \colon &\check{C}(K;\tatak) &\to &S^2\\
&c&\mapsto& \frac1{\norm{K(z\exp(2i\pi t))-K(z)}}\Bigl(K\bigl(z\exp(2i\pi t)\bigr)-K\bigl(z\bigr)\Bigr)
  \end{array}
\end{equation*}
and the degree of $G_K$ at an element $X$ of $S^2$ makes sense as soon as $X$ is a regular value of $G_K$ whose preimage is finite.

\bfig
 \centering
\begin{tikzpicture}
\draw [>-,draw=gray, thick] (.8,1.5) .. controls (.8,1) and (1.7,1) .. (1.7,2);
\draw [draw=white,double=gray,very thick] (1.2,3.3) .. controls (1,3.3) and (.7,3.2) .. (.7,3) .. controls (.7,2.7) and (2,1.3) .. (2.2,1.5) .. controls (2.4,1.5) and (2.6,2) .. (2.6,2.2) .. controls (2.6,2.4) and (2.4,3.3) .. (2.2,3.3);
\draw [draw=white,double=gray,very thick] (2.2,3.3)  .. controls (1.5,3.3) and (.8,2) .. (.8,1.5);
\draw [draw=white,double=gray,very thick] (1.7,2)  .. controls (1.7,2.7) and (1.5,3.3) .. (1.2,3.3);
\draw [draw=gray, thick, ->] (1.7,2)  .. controls (1.7,2.7) and (1.5,3.3) .. (1.2,3.3);
\draw  (.6,1.5) node{\scriptsize $K$}; 
\draw (2.2,3.3) node[above]{\scriptsize $v_1$};
\draw (2.2,3.3) -- (1.7,2) node[right]{\scriptsize $v_2$};
\draw [->] (2.2,3.3) -- (1.95,2.65);
\fill (1.7,2) circle (1.5pt) (2.2,3.3) circle (1.5pt);
\end{tikzpicture}

\caption{A configuration of a segment 
on $K$}
\label{figconfchord}
\end{figure}

The annulus $\check{C}(K;\tatak)$ can be compactified to the closed annulus ${C}(K;\tatak)=S^1 \times \left[0,1\right]$, to which $G_K$ extends smoothly. The extended $G_K$, still denoted by $G_K$, maps $(z,0) \in S^1 \times \{0\}$  to the direction of the tangent vector to $K$ at $z$. It maps $(z,1) \in S^1 \times \{1\}$ to the opposite direction.
The closed annulus ${C}(K;\tatak)$ is an example of a smooth manifold whose boundary is 
\begin{equation*}\partial{C}(K;\tatak) = S^1 \times \{0\} \cup (- S^1 \times \{1\}).\end{equation*}
The degree of $G_K$ is a continuous map from $S^2\setminus G_K(\partial {C}(K;\tatak))$ to $\ZZ$.
Let us compute it for the following embeddings of the trivial knot.

Let $O$ be an embedding of the circle to the horizontal plane.\footnote{Here, \say{the} horizontal plane is the plane $\CC \times \{0\}$ of $\RR^3$ viewed as $\CC \times \RR$.} The image under $G_O$ of the whole annulus is in the horizontal great circle of $S^2$. The set of regular values of $G_O$ is the complement of this circle, and hence the degree of $G_O$ is zero on all this set.

Let $K_1$ and $K_{-1}$ be embeddings of $S^1$ such that
\begin{itemize}
 \item the images of $K_1$ and $K_{-1}$ project to the horizontal plane as in Figure~\ref{figeightdiagram},
 \item they lie in the horizontal plane everywhere except in a small ball around where they cross over, and
 \item they lie in the union of two orthogonal planes.
\end{itemize}

\bfig
\centering
\begin{tikzpicture}
\begin{scope}[xshift=-3cm]
\draw [->] (.3,0) -- (-.3,0);
\draw [->] (-.3,0) .. controls (-.5,0) and (-.7,-.2) .. (-.7,-.4) .. controls  (-.7,-.6) and (0,-.4).. (0,-.15);
\draw [<-] (.3,0) .. controls (.5,0) and (.7,.2) .. (.7,.4) .. controls (.7,.6) and (0,.4) .. (0,.15);
\draw  (.35,-.3) node{\scriptsize $K_{-1}$};
\end{scope}

\draw [->] (.85,0) node{\scriptsize $O$} (.6,0) arc (0:360:.6) ;

\begin{scope}[xshift=3cm]
\draw [->] (0,-.3) -- (0,.3);
\draw [->] (0,.3) .. controls (0,.5) and (.2,.7) .. (.4,.7) .. controls (.6,.7) and (.4,0) .. (.15,0);
\draw [<-] (0,-.3) .. controls (0,-.5) and (-.2,-.7)  .. (-.4,-.7) .. controls  (-.6,-.7) and (-.4,0) .. (-.15,0);
\draw  (-.25,.35) node{\scriptsize $K_1$};
\end{scope}

\end{tikzpicture}

\caption{Diagrams of the trivial knot}
\label{figeightdiagram}
\end{figure}

The image of the boundary of ${C}(K_{\pm 1};\tatak)=S^1 \times \left[0,1\right]$ in $S^2$ lies in the union of the great circles of the two planes. More precisely, it lies in the union of the horizontal plane and two vertical arcs, as in the following figure.

\begin{center}
\begin{tikzpicture} \useasboundingbox (-1,-.7) rectangle (1,1);
\draw (0,0) circle (.9);
\draw [very thick] (-.9,0)  .. controls (-.9,-.15) and (-.3,-.25) .. (0,-.25) .. controls (.3,-.25) and (.9,-.15) ..  (.9,0);
\draw [dashed,very thick] (-.9,0)  .. controls (-.9,.15) and (-.3,.25) .. (0,.25) .. controls (.3,.25) and (.9,.15) ..  (.9,0);
\draw [very thick] (-.1,-.45) -- (-.1,-.07);
\draw [dashed,very thick] (.1,.07) -- (.1,.45);
\end{tikzpicture}
\end{center}

Therefore, when $K=K_{\varepsilon}$, for $\varepsilon=\pm 1$, the degree of $G_{K_{\varepsilon}}$ (extends as a map, which) is constant on each side of our horizontal equator.
Computing it at the North Pole $\upvec$ as in Subsection~\ref{sublkGausstwo}, we find that the degree of $G_{K_{\varepsilon}}$ is $\varepsilon$ on the Northern Hemisphere.
We similarly compute the degree of $G_{K_{\varepsilon}}$ on the Southern Hemisphere. It is also $\varepsilon$.

For any embedding  $K \colon S^1 \hookrightarrow \RR^3$, define 
\begin{equation*}I_{\theta}(K)=\int_{\check{C}(K;\tataksmall)} G_K^{\ast}(\omega_{S^2}).\end{equation*}
The real number $I_{\theta}(K)$ is the integral of $\deg(G_K)\omega_{S^2}$ over $S^2$. It can be seen as the \emph{algebraic area} of $G_K({C}(K;\tatak))$.
The above degree evaluation allows us to compute the integrals
$I_{\theta}(O)=0$, $I_{\theta}(K_1)=1$, and $I_{\theta}(K_{-1})=-1$.

More generally, say that
a knot embedding $K$ that lies in the union of the horizontal plane and a finite union of vertical planes so that the unit tangent vector to $K$ is never vertical is \emph{almost-horizontal}. The \emph{writhe} of a generic almost-horizontal knot embedding is the number of positive crossings minus the number of negative crossings of its orthogonal projection onto the horizontal plane.\footnote{The genericity of $K$ implies that the orthogonal projection onto the horizontal plane of $K$ is an immersion whose only multiple points are transverse double points.} An almost-horizontal embedding $K$ has a natural parallel $K_{\parallel}$ (up to isotopy) obtained from $K$ by (slightly) pushing it down.\footnote{A \emph{parallel}\index[T]{parallel!of a knot} of a knot embedding $K$ in a $3$-manifold $M$ (such as $\RR^3$) is a knot embedding $K_{\parallel} \colon S^1 \hookrightarrow M$ such that there exists an embedding $f$ from $\left[0,1\right] \times S^1$ into $M$ that restricts to $\{0\} \times S^1$ as $K$ and to $\{1\} \times S^1$ as $K_{\parallel}$.}
For any almost-horizontal knot embedding $K$, the degree of $G_K$ extends to a constant function of $S^2$. More precisely, we have the following lemma.

\begin{lemma}
\label{lemalmosthor}
 For any almost-horizontal knot embedding $K$, the degree of $G_K$ at any regular value of $G_K$ is the writhe of $K$, and we have
 \begin{equation*}I_{\theta}(K)=\int_{\check{C}(K;\tataksmall)} G_K^{\ast}(\omega_{S^2})=lk(K,K_{\parallel}).\end{equation*}
\end{lemma}
\bp
As in the previous examples,
for such a knot embedding $K$, the degree of $G_K$ extends to a constant function
on each hemisphere of $S^2$. It maps the regular values of these hemispheres to the writhe of $K$. In the complement of $K$, the parallel $K_{\parallel}$ is isotopic to the parallel $K_{\parallel,\ell}$ on the left-hand side of $K$. The formulas of Subsection~\ref{sublkGauss} show that $lk(K,K_{\parallel,\ell})$ is the writhe of $K$.
\eop

A \emph{knot invariant} is a function of embeddings that
takes the same value on isotopic knots.
Unlike the Gauss linking number, the integral $I_{\theta}()$ is \emph{not invariant under isotopy} since it takes distinct values on the isotopic knot embeddings $K_{-1}$ and $K_1$.

In every isotopy class of embeddings of $S^1$ into $\RR^3$,
we can construct an embedding $K$ such that the degree of the direction map $G_K$ extends as the constant map of $S^2$ with value $0$ (or any arbitrary integer) as follows: add kinks such as $\pkink$ or $\nkink$ to a generic horizontal projection, and take a corresponding almost-horizontal embedding.

Since $I_{\theta}$ varies continuously under an isotopy, it maps any isotopy class of embeddings of $S^1$ into $\RR^3$ onto $\RR$. 
In particular, there are embeddings $K$ for which $I_{\theta}(K)$ is not an integer. For such an embedding, the degree of the direction map $G_K$ cannot be extended to a constant map on $S^2$.

\subsection{A first knot invariant which counts graph configurations}
\label{subsecwtwo}
Let us now try to count configurations $\confc$ of the following graph:

\begin{center}
 \begin{tikzpicture} \useasboundingbox (-.7,-.4) rectangle (.7,.4);
\draw [->,dash pattern=on 2pt off 2pt] (.5,0) arc (0:360:.5);
\draw (45:.5) -- (-135:.5) node[left]{\scriptsize $b_2$};
\draw (135:.5) -- (-45:.5) node[right]{\scriptsize $b_3$};
\draw (135:.5) node[left]{\scriptsize $a_3$}  (45:.5)  node[right]{\scriptsize $a_2$};
\draw [->] (45:.5) -- (-135:.25);
\draw [->] (135:.5) -- (-45:.25);
\fill (45:.5) circle (1pt) (135:.5) circle (1pt) (-135:.5) circle (1pt) (-45:.5) circle (1pt);
\end{tikzpicture}
\end{center}
A \indexT{configuration} $\confc$ of that graph is an injection from the set $\{b_3,a_2,a_3,b_2\}$ of vertices of $\diagcross$ to $\RR^3$ such that the images $\confc(b_3)$, $\confc(a_2)$, $\confc(a_3)$, and $\confc(b_2)$ of the vertices are on the knot $K$, and we successively meet $\confc(b_3)=K(z)$, $\confc(a_2)=K(z\exp(2i\pi \alpha_2))$, $\confc(a_3)=K(z\exp(2i\pi \alpha_3))$, and $\confc(b_2)=K(z\exp(2i\pi \beta_2))$ along $K$ following the orientation of $K$. The dashed circle shows the cyclic order of the four vertices.
The associated \emph{configuration space}\index[T]{configuration!space} $\check{C}(K;\diagcross)$ is
\begin{equation*}\left\{\begin{array}{r}\Bigl(K\bigl(z\bigr),K\bigl(z\exp(2i\pi \alpha_2)\bigr),K\bigl(z\exp(2i\pi \alpha_3)\bigr),K\bigl(z\exp(2i\pi \beta_2)\bigr)\Bigr) \suchthat \\  z \in S^1, (\alpha_2,\alpha_3,\beta_2) \in \left]0,1\right[^3, \alpha_2<\alpha_3<\beta_2\end{array}\right\}.\end{equation*}
For $i \in \{2,3\}$, set $e_i=(a_i,b_i)$, and let $G_{e_i}(\confc) = \frac{\confc(b_i)-\confc(a_i)}{\norm{\confc(b_i)-\confc(a_i)}}$ denote the direction in $S^2$ of the image under the configuration $\confc$ of the edge $e_i$.
The open configuration space $\check{C}(K;\diagcross)$ has the natural compactification \begin{equation*}{C}(K;\diagcross)=S^1 \times \{(\alpha_2,\alpha_3,\beta_2) \in \left[0,1\right]^3 \suchthat \alpha_2\leq\alpha_3\leq\beta_2\}.\end{equation*} The maps $G_{e_2}$ and $G_{e_3}$ smoothly extend to ${C}(K;\diagcross)$ as before, and 
\begin{equation*}G_{\diagcross} =(G_{e_2},G_{e_3}) \colon {C}(K;\diagcross) \to (S^2)^2\end{equation*} is a smooth map between two compact $4$-manifolds.

The \emph{codimension-one faces} ${C}(K;\diagcross)$ are the four faces $(\alpha_2=0)$, $(\alpha_2=\alpha_3)$, $(\alpha_3=\beta_2)$, and $(\beta_2=1)$, on which $\confc$ maps (at least) two consecutive vertices to the same point on $K$.\footnote{In this introduction, codimension-one faces are closed. Later, they will be open parts of the boundaries.}  When $G_{\diagcross}$ is locally an embedding near such a face, the degree of $G_{\diagcross}$ changes by $\pm 1$ when we cross the image of that face.\footnote{See Lemma~\ref{lemdeggen} for a precise statement.} Thus, it suffices to determine the images of the interiors of these codimension-one faces and the local degree at one regular value to determine the degree of $G_{\diagcross}$, as a map from $(S^2)^2 \setminus G_{\diagcross}\left(\partial {C}(K;\diagcross)\right) $ to $\ZZ$. We associate
the following figure to the face $(\beta_2=1)$ of ${C}(K;\diagcross)$.

\begin{center}
\diagv{>}{2}{>}{3}{a_2}{a_3}{\,b_2=b_3}
\end{center}

Let us now try to count configurations $\confc$ of the following tripod $\diagtripod$.

\begin{center}
\begin{tikzpicture} \useasboundingbox (-.7,-.4) rectangle (.7,.4);
\draw [->,dash pattern=on 2pt off 2pt] (.5,0) arc (0:360:.5);
\draw (0,0) -- (-60:.5) node[right]{\scriptsize $v_1$};
\draw (0,0) -- (60:.5) node[right]{\scriptsize $v_2$};
\draw (.2,0) node{\scriptsize $w$}  (0,0) -- (180:.5) node[left]{\scriptsize $v_3$};
\draw [->] (-60:.5) -- (-60:.25);
\draw [->] (60:.5) -- (60:.25);
\draw [->] (180:.5) -- (180:.25);
\fill (60:.5) circle (1pt) (180:.5) circle (1pt) (-60:.5) circle (1pt) (0,0) circle (1pt);
\end{tikzpicture} 
\end{center}
A \emph{configuration} $\confc$ of this tripod $\diagtripod$ is an injection from the set $\{w,v_1,v_2,v_3\}$ of its vertices into $\RR^3$, where the images $\confc(v_1)$, $\confc(v_2)$, and $\confc(v_3)$ of the vertices $v_1$, $v_2$, and $v_3$ are on the knot $K$, and we successively meet $\confc(v_1)$, $\confc(v_2)$, and $\confc(v_3)$ along $K$ following the orientation of $K$.

\bfig
 \centering
\begin{tikzpicture}
\draw [>-,draw=gray, thick] (.8,1.5) .. controls (.8,1) and (1.7,1) .. (1.7,2);
\draw [draw=white,double=gray,very thick] (1.2,3.3) .. controls (1,3.3) and (.7,3.2) .. (.7,3) .. controls (.7,2.7) and (2,1.3) .. (2.2,1.5) .. controls (2.4,1.5) and (2.6,2) .. (2.6,2.2) .. controls (2.6,2.4) and (2.4,3.3) .. (2.2,3.3);
\draw [draw=white,double=gray,very thick] (2.2,3.3)  .. controls (1.5,3.3) and (.8,2) .. (.8,1.5);
\draw [draw=white,double=gray,very thick] (1.7,2)  .. controls (1.7,2.7) and (1.5,3.3) .. (1.2,3.3);
\draw [draw=gray, thick, ->] (1.7,2)  .. controls (1.7,2.7) and (1.5,3.3) .. (1.2,3.3);
\draw  (.6,1.5) node{\scriptsize $K$}; 
\draw [draw=white,double=black,very thick] (3,3.2) -- (1.7,2);
\draw [black] (3,3.2) -- (2.2,3.3) ;
\draw [black] (3,3.2) -- (1.7,2) ;
\draw (3.2,3) node{\scriptsize $w$}  (3,3.2) -- (2.6,2.2) node[left]{\scriptsize $v_3$} (2.2,3.3) node[above]{\scriptsize $v_1$} (1.7,2) node[right]{\scriptsize $v_2$};
\draw [->] (2.2,3.3) -- (2.6,3.25);
\draw [->] (1.7,2) -- (2.35,2.6);
\draw [->] (2.6,2.2) -- (2.8,2.7);
\fill (1.7,2) circle (1.5pt) (2.6,2.2) circle (1.5pt) (2.2,3.3) circle (1.5pt) (3,3.2) circle (1.5pt);
\end{tikzpicture}
\caption{A configuration of the tripod on $K$}
\label{figconftripod}
\end{figure}

Such a \emph{configuration} $\confc$ maps $w$ to $\confc(w) \in \RR^3$, $v_1$ to $\confc(v_1)=K(z)$ for some $z\in  S^1$, $v_2$ to $\confc(v_2)=K(z\exp(2i\pi t_2))$, and $v_3$ to $\confc(v_3)=K(z\exp(2i\pi t_3))$. The set of these configurations is the \emph{configuration space} $\check{C}(K;\diagtripod)$. It is an open $6$-manifold parametrized by an open subspace of $\RR^3 \times S^1\times \{(t_2,t_3) \in \left]0,1\right[^2 \suchthat t_2 <t_3\}$.
For $i \in \underline{3}=\{1,2,3\}$, set $e_i= (v_i,w)$, and let $G_{e_i}(\confc) = \frac{\confc(w)-\confc(v_i)}{\norm{\confc(w)-\confc(v_i)}}$ denote the direction of the image under $\confc$ of the edge $e_i$ in $S^2$.
These edge directions together provide a map
\begin{equation*}\begin{array}{llll}\check{G}_{\diagtripod} \colon&\check{C}(K;\diagtripod)&\to &(S^2)^3\\
&\confc &\mapsto &\left(G_{e_1}(\confc),G_{e_2}(\confc),G_{e_3}(\confc) \right).\end{array}\end{equation*}
 from our open $6$-manifold $\check{C}(K;\diagtripod)$ to the $6$-manifold $(S^2)^3$.
 
For a regular value $(X_1,X_2,X_3)$ of $\check{G}_{\diagtripod}$ whose preimage is finite, we can again count the configurations of the tripod such that the direction of the edge $e_i$ is $X_i$, as the degree of $\check{G}_{\diagtripod}$ at $(X_1,X_2,X_3)$. 

In Chapter~\ref{chapcompconf}, we construct a \emph{compactification} ${C}(K;\diagtripod)$ of $\check{C}(K;\diagtripod)$ (and of many similar configuration spaces), using blow-up techniques, as William  Fulton, Robert MacPherson \cite{FultonMcP}, Scott Axelrod, Isadore Singer \cite{axelsingI}, and Maxim Kontsevich \cite{ko} did. This compactification ${C}(K;\diagtripod)$ is a smooth compact $6$-manifold with boundary and ridges. Its interior is $\check{C}(K;\diagtripod)$.
The map $\check{G}_{\diagtripod}$ extends to a smooth map ${G}_{\diagtripod}$ over ${C}(K;\diagtripod)$.
\emph{Regular values} of ${G}_{\diagtripod}$ are regular values of $\check{G}_{\diagtripod}$ in the complement of $G_{\diagtripod}\bigl(\partial C(K;\diagtripod) ={C}(K;\diagtripod) \setminus \check{C}(K;\diagtripod)\bigr)$. They form an open dense
subset $\CO$ of $(S^2)^3$, for which the differential degree of $\check{G}_{\diagtripod}$ makes sense.
As mentioned above and proved in Lemma~\ref{lemdeggen}, this local integral degree extends to a continuous map on $(S^2)^3 \setminus {G}_{\diagtripod}(\partial C(K;\diagtripod))$. It restricts to a constant map on every connected component of $(S^2)^3 \setminus {G}_{\diagtripod}(\partial C(K;\diagtripod))$. In Lemma~\ref{lemItripod}, we explicitly compute this local degree  when $K$ is the round circle $O$ in a plane, as an example.
Our computation shows that this degree cannot be extended to a constant map in this case.

Let $p_i \colon (S^2)^3 \to S^2$ denote the projection to the $i^{th}$ factor.
Define the algebraic volume \begin{equation*} I(K;\diagtripod)= \int_{(S^2)^3}\deg(G_{\diagtripod}) \wedge_{i=1}^3 p_i^{\ast}(\omega_{S^2})=\int_{\check{C}(K;\diagtripodmini)} G_{\diagtripod}^{\ast}\left(\wedge_{i=1}^3 p_i^{\ast}(\omega_{S^2})\right) \end{equation*}
of the image of $ C(K;\diagtripod)$ under $G_{\diagtripod}$.
This algebraic volume $I(K;\diagtripod)$ has no reason to be a knot isotopy invariant, and it is not.

To compute the map \begin{equation*}\deg(G_{\diagtripod}) \colon (S^2)^3 \setminus {G}_{\diagtripod}(\partial C(K;\diagtripod)) \to \ZZ,\end{equation*}
we look at  ${G}_{\diagtripod}(\partial C(K;\diagtripod))$.
Let us describe the compactification $C(K;\diagtripod)$ on loci where $\confc(w)$ approaches $\confc(v_1)$
and $\confc(v_1)$ is far from $\confc(v_2)$ and $\confc(v_3)$. With our coordinates, this amounts to assuming 
$\confc(w)=\confc(v_1) + \eta x $ for some $\eta \in \left]0,\varepsilon\right]$ for a small $\varepsilon >0$ and some $x \in S^2$, and $(t_2,t_3) \in \left[\alpha,1-\alpha\right]^2$ for some $\alpha \in\left]0,1/2\right[$.
This part of $\check{C}(K;\diagtripod)$ is diffeomorphic to
$\left]0,\varepsilon\right] \times S^2 \times S^1 \times \{(t_2,t_3) \in \left[\alpha,1-\alpha\right]^2 \suchthat t_2 <t_3\}$. Its closure in the compactification $C(K;\diagtripod)$ is naturally diffeomorphic to $\left[0,\varepsilon\right] \times S^2 \times S^1 \times \{(t_2,t_3) \in \left[\alpha,1-\alpha\right]^2 \suchthat t_2 \leq t_3\}$. (We close $\left]0,\varepsilon\right[$ at $0$---we also relax the inequality $t_2 <t_3$ but this is not important for us now.)
In the compactification ${C}(K;\diagtripod)$, the image $\confc(w)$ may coincide with $\confc(v_1)$ (when $\eta=0$). The direction from $\confc(v_1)$ to $\confc(w)$ is still defined in this case. It is contained in the $S^2$ factor.
In particular, $G_{\diagtripod}$ extends to this part of ${C}(K;\diagtripod)$. The compactification creates the local boundary $\{0\} \times(- S^2) \times S^1 \times \{(t_2,t_3) \in \left]\alpha,1-\alpha\right[^2 \suchthat t_2 \leq t_3\}$. This local boundary is a 5-dimensional manifold. Its image under $G_{\diagtripod}$ is the product by $S^2$ of the image in $(S^2)^2$ of the configuration space associated to the following graph by the natural Gauss map \say{direction of the edges}.
\begin{center}
\diagv{>}{ }{>}{ }{v_2}{v_3}{w=v_1}
\end{center}

Again, the dashed circle represents the cyclic order of $\confc(v_1)$, $\confc(v_2)$, and $\confc(v_3)$ along $K$. The image of this local boundary creates a \say{wall} in $(S^2)^3$ across which the local degree changes by $\pm 1$. (See Lemma~\ref{lemdeggen} for a precise statement.)

We recognize the picture associated to the face $(\beta_2=1)$ denoted by $F_{\beta_2=1}$ of ${C}(K;\diagcross)$. We observe that the image of the corresponding face under 
$G_{\diagtripod}$ coincides with the image of the face $S^2 \times F_{\beta_2=1}$ of $S^2 \times C(K;\diagcross)$ under 
\begin{equation*}G^{\prime}_{\diagcross} =1_{S^2} \times G_{\diagcross} \colon S^2 \times C(K;\diagcross) \to (S^2)^3.\end{equation*}
So the combination $\left(\deg(1_{S^2} \times G_{\diagcross}) - \deg(G_{\diagtripod})\right)$ does not vary across
the images of the corresponding faces (or at least not because of them).
(There are some sign and orientation issues to check here, but we will carefully treat them in a broader generality in Section~\ref{secdefconfspace} and Lemma~\ref{lemstu}. Let the reader trust me that the signs are correct here.)
We glued the images of $G^{\prime}_{\diagcross}$ and of $G_{\diagtripod}$ along $G^{\prime}_{\diagcross}(S^2 \times F_{\beta_2=1})$ to make the union of these images behave as the image of a manifold without boundary, locally. 
Unfortunately, the three other faces of $S^2 \times C(K;\diagcross)$ created other walls in $(S^2)^3$ associated to the following figures: 
\begin{center}
 \diagv{>}{3}{<}{2}{a_3}{b_2}{\,b_3=a_2} \diagv{<}{2}{<}{3}{b_2}{b_3}{a_2=a_3} \diagv{<}{3}{>}{2}{b_3}{a_2}{a_3=b_2}
\end{center}
To cancel these walls with the same type of faces of ${C}(K;\diagtripod)$ as before, we use Gauss maps associated to the following diagrams:
\begin{center}
\diagtripoded{>}{3}{<}{2} \diagtripoded{<}{2}{<}{3} \diagtripoded{<}{3}{>}{2}
\end{center}
Let us describe these Gauss maps more precisely.
For any subset $I$ of $ \underline{3}$, let $\iota_I$ denote the diffeomorphism of $(S^2)^3$ which maps $(X_1,X_2,X_3)$ to $(\varepsilon_1 X_1,\varepsilon_2 X_2, \varepsilon_3 X_3)$, where $\varepsilon_i=-1$ when $i \in I$, and $\varepsilon_i=1$ when $ i \notin I$.
For two elements $i,j$ of $\underline{3}$, simply write $\iota_i=\iota_{\{i\}}$ and $\iota_{ij}= \iota_{\{i,j\}}$. For a permutation $\sigma$ of $\underline{3}$, let $\sigma_{\ast}$ denote the diffeomorphism of $(S^2)^3$ that maps $(Y_1,Y_2,Y_3)$ to $(Y_{\sigma^{-1}(1)},Y_{\sigma^{-1}(2)},Y_{\sigma^{-1}(3)})$. (We have $\sigma_{\ast}\bigl((X_{\sigma(1)},X_{\sigma(2)},X_{\sigma(3)}) \bigr)=(X_1,X_2,X_3)$.)
The Gauss maps associated to the above diagrams are $\iota_2 \circ (23)_{\ast} \circ G_{\diagtripod}$, 
$\iota_{23} \circ G_{\diagtripod}$, and $\iota_3 \circ  (23)_{\ast} \circ G_{\diagtripod}$, respectively. The combination
\begin{multline*}\deg(1_{S^2} \times G_{\diagcross}) - \deg(G_{\diagtripod}) - \deg(\iota_{23} \circ G_{\diagtripod}) \\  + \deg(\iota_2 \circ (23)_{\ast} \circ G_{\diagtripod}) + \deg(\iota_3 \circ  (23)_{\ast} \circ G_{\diagtripod})\end{multline*} does not vary across the 
 boundary of the image of $1_{S^2} \times G_{\diagcross}$ as Figure~\ref{figvardeg} suggests.
 
 \bfig
\centering
\begin{tikzpicture}
\fill[gray!20] (7,.5) rectangle (11,2);
\fill[gray!40] (3,.5) rectangle (7,2);
\fill[gray!40] (11,.5) rectangle (15,2); 
\fill[gray!40] (7,2) rectangle (11,4);
\fill[gray!40] (7,-1.5) rectangle (11,.5);
\draw (7,-1)-- (7,3.5) (11,-1)-- (11,3.5) (9,1) node{\scriptsize $S^2 \times G_{\diagcross}({C}(K;\diagcross))$} (4,.5) -- (14,.5) (4,2) -- (14,2) 
(9,-.5) node{\scriptsize $G_{\diagtripod} \left(-{C}(K;\diagtripod)\right)$}
(9,3) node{\scriptsize $\iota_{23} G_{\diagtripod} \left(-{C}(K;\diagtripod)\right)$}
(5,1) node{\scriptsize $\iota_3 (23)_{\ast} G_{\diagtripod} \left({C}(K;\diagtripod)\right)$}
(13,1) node{\scriptsize $\iota_2 (23)_{\ast}G_{\diagtripod} \left({C}(K;\diagtripod)\right)$};
\draw[dashed] (4,.5) -- (3,.5) (15,.5) -- (14,.5) (4,2) -- (3,2) (15,2) -- (14,2)
(7,-1)-- (7,-1.5) (11,-1)-- (11,-1.5) (7,4)-- (7,3.5) (11,4)-- (11,3.5);
\draw (12,2.8) node{\diagv{<}{2}{<}{3}{b_2}{b_3}{a_2=a_3}};
\draw [->,very thin] (11.45,2.6) -- (10,2.08);
\draw (6,2.8) node{\diagv{<}{3}{>}{2}{b_3}{a_2}{ }};
\draw [->,very thin] (6.18,2.28) -- (6.92,1.5);
\draw (6,-.3) node{\diagv{>}{2}{>}{3}{a_2}{a_3}{ }};
\draw [->,very thin] (6.55,-.1) -- (8,.42);
\draw (12,-.3) node{\diagv{>}{3}{<}{2}{a_3}{b_2}{\,b_3=a_2}};
\draw [->,very thin] (11.65,.2) -- (11.08,1);
\end{tikzpicture}
\caption{Variation of the image of $\id_{S^2} \times G$ across the boundary in $(S^2)^3$ (pictured as 2-dimensional)}
\label{figvardeg}
\end{figure}

Unfortunately, this process introduces other walls, and we have not yet cancelled the walls due to the faces $\confc(v_2)=\confc(w)$ and $\confc(v_3)=\confc(w)$. However, we have the following proposition, which is a corollary of Theorem~\ref{thmConwaycircZSthree}, as it will be seen right after Theorem~\ref{thmConwaycircZSthree}.

\begin{proposition}
\label{propwtwo}
The map $\tilde{w}_2$
\begin{multline*}\frac1{24}\sum_{\sigma \in \mathfrak{S}_3}\deg\Bigl(\sigma_{\ast}\bigl(1_{S^2} \times G_{\diagcross}\bigr)\Bigr)\\ - \frac1{48}\sum_{I \subseteq \underline{3}} (-1)^{\cardlef{ I}}\Bigl( \deg\bigl(\iota_I \circ G_{\diagtripod}\bigr) + \deg\bigl(\iota_I \circ (23)_{\ast} \circ G_{\diagtripod}\bigr) \Bigr),\end{multline*}
which is well-defined on an open dense subset of $(S^2)^3$, extends as a constant function of $(S^2)^3$ whose value $w_2(K)$ is in 
$\frac1{48}\ZZ$. 
\end{proposition}
\bsp
Let us show that the boundary \begin{equation*}\sigma_{\ast}\Bigl(S^2 \times G_{\diagcross}\bigl(\partial {C}(K;\diagcross)\bigr)\Bigr)\end{equation*} can be glued to the images of the faces $\confc(v_i)=\confc(w)$ of ${C}(K;\diagtripod)$ under the $\iota_I \circ G_{\diagtripod}$ or the $\iota_I \circ (23)_{\ast} \circ G_{\diagtripod}$, up to sign. For every permutation $\sigma$ of $\underline{3}$, the boundary $\sigma_{\ast}\left(S^2 \times G_{\diagcross}(\partial {C}(K;\diagcross))\right)$ consists of the four faces
\begin{center}
\diagvv{>}{\sigma(2)}{>}{\sigma(3)}
\diagvv{>}{\sigma(3)}{<}{\sigma(2)}
\diagvv{<}{\sigma(2)}{<}{\sigma(3)}
\diagvv{<}{\sigma(3)}{>}{\sigma(2)}
\end{center} where the first one and the third one come with a coefficient $- \frac1{24}$, and the second one and the fourth one come with a coefficient $\frac1{24}$.
The images of the (open) faces $(\confc(v_1)=\confc(w))$, $(\confc(v_2)=\confc(w))$, and $(\confc(v_3)=\confc(w))$ under
$G_{\diagtripod}$ are 
\begin{center}
\diagvhaut{>}{2}{>}{3}{ }{ }{ },
\diagvhaut{>}{3}{>}{1}{ }{ }{ }, and \diagvhaut{>}{1}{>}{2}{ }{ }{ }, respectively.
\end{center}
The images of the (open) faces $(\confc(v_1)=\confc(w))$, $(\confc(v_2)=\confc(w))$, and $(\confc(v_3)=\confc(w))$ under
$(23)_{\ast}  \circ G_{\diagtripod}$ are
\begin{center}
\diagvhaut{>}{3}{>}{2}{ }{ }{ },
\diagvhaut{>}{2}{>}{1}{ }{ }{ }, and \diagvhaut{>}{1}{>}{3}{ }{ }{ }, respectively.
\end{center}
In order to obtain the images under the compositions of these maps by some $\iota_I$, we reverse the edge $i$ when $i \in I$.
Each of these faces appears twice (once for each orientation of the collapsed edge), with the same sign (since the antipodal map of $S^2$ reverses the orientation) with a coefficient $(-1)^{n(F)}\frac1{48}$, where $n(F)$ is the number of edges towards the bivalent vertex.

We leave the general discussion of signs to the reader, and we do not discuss all the faces of ${C}(K;\diagtripod)$
in this sketch.
Let us just mention that the image of the faces $\confc(v_1)=\confc(v_2)$ is contained in the codimension-two 
subspace of $(S^2)^3$, for which at least two $S^2$-coordinates are equal or opposite. So these faces do not create walls and may be forgotten.
\eopwobp

The $\sigma_{\ast}$ preserve the volume of $(S^2)^3$, and the $\iota_I$ multiply the volume by $(-1)^{\cardlef{I}}$. Therefore, the combination $w_2(K)$ in Proposition~\ref{propwtwo} may also be written as
\begin{equation*}w_2(K)=\frac14 \int_{{C}(K;\diagcrossmini)}G_{\diagcross}^{\ast}\left(\wedge_{i=1}^2 p_i^{\ast}(\omega_{S^2})\right) -\frac13 \int_{{C}(K;\diagtripodmini)}G_{\diagtripod}^{\ast}\left(\wedge_{i=1}^3 p_i^{\ast}(\omega_{S^2})\right).\end{equation*}
Since $w_2(.)$ is valued in $\frac1{48}\ZZ$, and since $w_2(K)$ varies continuously under an isotopy of $K$, $w_2(.)$ is an isotopy invariant.
Enore Guadagnini, Maurizio Martellini, and Mihail Mintchev \cite{gmm}, and Dror Bar-Natan \cite{barnatanper} independently studied the invariant $w_2(.)$ under this integral form in 1990.
To prove its isotopy invariance, one can alternatively use Stokes' theorem to evaluate the variations of the integrals under a knot isotopy $(t,z) \mapsto K_t(z)$. In \cite{botttaubes}, Raoul Bott and Clifford Taubes used compactifications of the one-parameter configuration spaces $\cup_{t \in \left[0,1\right]} C(K_t;\Gamma)$, like those discussed above, to check invariance. 

The above formulation of Proposition~\ref{propwtwo} presents the invariant $w_2(K)$ as a discrete count of configurations of 
$\diagtripod$ and $\diagcross$, as Dylan Thurston \cite{thurstonconf} and Sylvain Poirier \cite{poirier} first did independently.\footnote{Michael Polyak and Oleg Viro obtained a similar result \cite[Theorem 3.A, Section 3.5]{pv} in the setting of long knots of $\RR^3$, with fewer involved gluings, and with induced combinatorial formulae in terms of knot diagrams.} 
For a generic triple $(a,b,c)$ of $(S^2)^3$, we count the configurations of $\diagcross$ for which the edge directions are a pair of noncolinear vectors in $\{a,b,c,-a,-b,-c\}$,
and the configurations of $\diagtripod$ for which the edge directions are a triple of pairwise noncolinear vectors in $\{a,b,c,-a,-b,-c\}$, with some coefficients and some signs determined by the corresponding above local degrees.

Any knot of $\RR^3$ is obtained from the trivial knot by (isotopies and) a finite number of crossing changes  \pcpetittexte $\rightarrow$ \ncpetittexte. So, in order to determine a real-valued knot invariant $w$, it suffices to know its value on the trivial knot and its variation
$w\left(\pcpetit \right) - w\left(\ncpetit \right)$, denoted by $w\left(  \doubleppetit \right)$,
under crossing change. 
(Here, \pcpetittexte and \ncpetittexte represent knot diagrams that coincide outside a disk that they intersect as in the figure, and \doubleppetittexte represents a diagram that coincides with the former diagrams outside this disk.)

The variation $w\left(  \doubleppetit \right)=w\left(\pcpetit \right) - w\left(\ncpetit \right)$ can be thought of as a discrete derivative of the invariant. The variation of this variation 
\begin{equation*}\begin{array}{lll}w\left(  \doubleppetit \doubleppetit \right)&= w\left(\pcpetit  \doubleppetit \right)&- w\left(\ncpetit  \doubleppetit \right)\\
&= w\left(\pcpetit  \pcpetit \right)  - w\left(\pcpetit  \ncpetit \right)& - w\left(\ncpetit  \pcpetit \right) +w\left(\ncpetit  \ncpetit \right)\end{array}\end{equation*}
under a disjoint crossing change is thought of as a discrete second derivative of this invariant.
In $\RR^3$, a real-valued knot invariant $w$ is actually determined by its value on the trivial knot and the discrete second derivative. (If the second derivative is zero, then the variation under a crossing change is independent of the knot, and it is the same as $(w(\huitpcpetit) -w(\huitncpetit)=0)$.) In \cite{gmm}, Enore Guadagnini, Maurizio Martellini, and Mihail Mintchev computed $w_2(O)=-\frac1{24}$.
This is reproved in Examples~\ref{examplecomconfintthree} and Lemma~\ref{lemItripod}.
In \cite[\S 6.3]{barnatanper}, Dror Bar-Natan computed the above \say{discrete second derivative} for the invariant $w_2$,
and he found:
\begin{equation*} w_2\left(\begin{tikzpicture} \useasboundingbox (-.7,.02) rectangle (1.6,.38);
\draw [->] (.2,0) -- (-.2,.4);
\draw [->]  (-.2,0) -- (.2,.4);
\fill  (0,.2) circle (1.5pt);
\draw [<-] (1.1,0) -- (.7,.4);
\draw [<-]  (.7,0) -- (1.1,.4);
\fill  (.9,.2) circle (1.5pt);
\draw [dashed] (-.2,.4) .. controls (-.4,.6) and  (-.5,.35) .. (-.5,.2)  .. controls (-.5,.05) and  (-.4,-.2) .. (-.2,0)
 (1.1,.4) .. controls (1.3,.6) and  (1.4,.35) .. (1.4,.2)  .. controls (1.4,.05) and  (1.3,-.2) .. (1.1,0)
 (.2,.4) .. controls (.35,.55) and (.55,.55) .. (.7,.4) 
  (.2,0) .. controls (.35,-.15) and (.55,-.15) .. (.7,0);
\end{tikzpicture} \right) =0 \;\;\;\; \mbox{and} \;\;\;\; 
 w_2\left( \begin{tikzpicture} \useasboundingbox (-.7,.02) rectangle (1.6,.38);
\draw [->] (.2,0) -- (-.2,.4);
\draw [->]  (-.2,0) -- (.2,.4);
\fill  (0,.2) circle (1.5pt);
\draw [<-] (1.1,0) -- (.7,.4);
\draw [<-]  (.7,0) -- (1.1,.4);
\fill  (.9,.2) circle (1.5pt);
\draw [dashed] (.2,.4) .. controls (.4,.6) and (1.3,.6) .. (1.1,.4)
(-.2,.4) .. controls (-.4,.6) and (.5,.6) .. (.7,.4)
(-.2,0) .. controls (-.4,-.2) and  (1.3,-.2) .. (1.1,0)
(.2,0)  .. controls (.3,-.1) and (.6,-.1) .. (.7,0);
\end{tikzpicture} \right) =1, \end{equation*}
where the dashed lines indicate the connections inside $K$ of the crossing strands.
This allowed him to identify $w_2$ with $a_2-\frac1{24}$, where $a_2(K)=\frac12\Delta^{\prime\prime}(K)(1)$ is half the second derivative of the Alexander polynomial of $K$ at one, since $a_2$ has the same discrete second derivative as $w_2$ and $a_2(O)=0$.

\subsection{On other similar invariants of knots in \texorpdfstring{$\RR^3$}{the ambient space}}
We can associate similar integrals over configuration space to every uni-trivalent graph $\Gamma$ whose univalent vertices are ordered cyclically, as in the following figure:

\begin{center}
\begin{tikzpicture} \useasboundingbox (-4,-1.5) rectangle (4,1.5);
\draw [->,dash pattern=on 2pt off 2pt] (-1.5,0) arc (-180:180:1.5);
\begin{scope}[yshift=.7cm]
\draw [thick] 
(0,0) circle (.5)
(70:.5) -- (-.1,.1)
(155:.5)  .. controls (-.2,.1) .. (-.1,.1) .. controls (0,.1) .. (35:.5) ;
\fill (-.1,.1) circle (1.5pt) (155:.5) circle (1.5pt) (70:.5) circle (1.5pt) (35:.5) circle (1.5pt) (-90:.5) circle (1.5pt);
\end{scope}
\draw [thick] (-125:1.5) -- (-.15,-.4) -- (-55:1.5) (-.15,-.4) -- (-.05,0) .. controls (.1,0) and (1.3,-.2) .. (1.5,0)(0,.2) -- (-.05,0) (0,1.2) -- (110:1.5);
\fill (-55:1.5) circle (1.5pt) (-125:1.5) circle (1.5pt) (0:1.5) circle (1.5pt) (-.05,0) circle (1.5pt) (-.15,-.4)  circle (1.5pt) (0,1.2) circle (1.5pt) (110:1.5) circle (1.5pt);
\end{tikzpicture}
\end{center}
We can also exhibit other similar combinations that provide isotopy invariants of knots in $\RR^3$, as several authors, including Maxim Kontsevich \cite{ko}, Daniel Altsch\"uler and Laurent Freidel \cite{af}, Dylan Thurston \cite{thurstonconf}, did.
These invariants are all finite type invariants with respect to the following definition.
An invariant valued in an abelian group is of \emph{degree less than $n$} if all its discrete derivatives of order $n$, which generalize the previously studied discrete second derivative, vanish.
A knot invariant is a \emph{Vassiliev invariant} or \emph{a finite type invariant} if it is of degree less than some integer $n$. (A precise definition is given in Section~\ref{secdeffintype}.)
In \cite{Kon,barnatan}, Dror Bar-Natan and Maxim Kontsevich proved the \say{fundamental theorem of Vassiliev invariants}, which determines the space of real-valued finite type invariants as the dual of an algebra $\CA$ generated by uni-trivalent graphs.
Daniel Altsch\"uler and Laurent Freidel \cite{af}---and Dylan Thurston \cite{thurstonconf} independently---constructed a \say{universal Vassiliev invariant} $\Zinvuf$ of knots in $\RR^3$, valued in $\CA$ such that, for any knot $K$ of $\RR^3$, $\Zinvuf(K)$ is a combination of classes of uni-trivalent graphs
whose coefficients are integrals over corresponding configuration spaces.\footnote{The original proof of the fundamental theorem of Vassiliev invariants relies on the construction of another---possibly equal---universal Vassiliev invariant called the \emph{Kontsevich integral}.} Any real-valued finite type invariant may be expressed as $\psi \circ \Zinvuf$ for some linear form $\psi \colon \CA \to \RR$.
Therefore, every real-valued finite type knot invariant is a combination of integrals over configuration spaces associated to uni-trivalent graphs.
It is still unknown whether finite type knot invariants distinguish all knots of $\RR^3$,
but many known polynomial invariants of knots of $\RR^3$, such as the Alexander polynomial, the Jones polynomial, and the HOMFLYPT polynomial, factor through $\Zinvuf$.

\subsection{On similar invariants of knots in other \texorpdfstring{$3$}{3}-manifolds}
\label{subsecsiminv}
Recall that $\drad{1} \times \left[0,1\right]$ denotes the standard cylinder in $\RR^3=\CC \times \RR$, and that a \emph{rational homology cylinder} $\hcylc$ is a compact oriented $3$-manifold with the same boundary and the same rational homology as $\drad{1} \times \left[0,1\right]$ (i.e., the rational homology of a point).\footnote{This homological condition can be rephrased as \say{The (compact oriented) $3$-manifold $\hcylc$ is connected, and every knot embedding $K$ in $\hcylc$ bounds a rational chain in $\hcylc$.}. Let $N(K)$ be a compact tubular neighborhood of $K$ in $\hcylc$, and let $\mathring{N}(K)$ denote its interior.
The latter condition is equivalent to the existence of a compact oriented surface in $\hcylc \setminus \mathring{N}(K)$ whose (oriented!) boundary is a disjoint union of curves in $\partial N(K)$ that does not bound a compact oriented surface in $N(K)$.}

In this book, following ideas of Maxim Kontsevich \cite{ko}, Greg Kuperberg, and Dylan Thurston \cite{kt}, we generalize $\Zinvuf$ to links
in more general $3$-manifolds denoted by $\crats=\crats(\hcylc)$. These manifolds are constructed from $\RR^3$ by replacing $\drad{1} \times \left[0,1\right]$ by a rational homology cylinder $\hcylc$. Such a more general $3$-manifold $\crats$ is called a \emph{rational homology $\RR^3$}. It looks like $\RR^3$ near $\infty$. The linking number can be defined as follows for a two-component link embedding $J \sqcup K \colon S^1 \sqcup S^1 \to \crats$.

Two submanifolds $\Aman$ and $B$ in a manifold $M$ are \emph{transverse}\index[T]{transverse!intersection} if at each intersection point $x \in \Aman \cap B$, we have $T_xM=T_x\Aman+T_xB$.
If two transverse oriented submanifolds $\Aman$ and $B$ in an oriented manifold $M$ are of \emph{complementary dimensions}
(i.e., if the sum of their dimensions is the dimension of $M$), then the \emph{sign of an intersection point} is $+1$ if  $T_xM=T_x\Aman \oplus T_xB$ as oriented vector spaces. Otherwise, the sign is $-1$. If $\Aman$ and $B$ are compact, and if $\Aman$ and $B$ are of complementary dimensions in $M$, then their \indexT{algebraic intersection} is the sum of the signs of the intersection points, it is denoted by
$\langle \Aman, B \rangle_{\!M\,}$.

When $K$ bounds a compact oriented embedded surface $\Sigma_K$ in $\crats$ transverse to $J$, the \indexT{linking number} of $J$ and $K$ is the algebraic intersection number $\langle J, \Sigma_{K} \rangle_{\!\crats\,}$ of $J$ and $\Sigma_K$ in $\crats$.
In general, there is an oriented surface $\Sigma_{nK}$ immersed in $\crats$ whose (oriented!) boundary is a positive multiple $nK$ of $K$, and $lk(J,K)=\frac1{n}\langle J, \Sigma_{nK} \rangle_{\!\crats\,}$.

For two-component links in $\RR^3$, this definition coincides with Definition~\ref{defGausslk} of the Gauss linking number. See Proposition~\ref{propdeflkeq}.

In the more general setting of a rational homology $\RR^3$, instead of counting uni-trivalent graphs whose edge directions belong to a finite set of directions, we use the notion of \emph{propagator} defined in Chapter~\ref{chapprop}.
A propagator is a rational combination of oriented compact $4$-manifolds in a suitable compactification $C_2(\rats)$ (defined in Section~\ref{secCtwo}) of the configuration space \begin{equation*}\check{C}_2(\rats)=\{(x,y) \in \crats^2 \suchthat x\neq y\},\end{equation*}
which shares many properties with our \emph{model propagator} \begin{equation*}p_{S^2}^{-1}(X)=\overline{\{(x,x+tX) \suchthat x \in \RR^3, t\in \left]0,+\infty\right[\}},\end{equation*} for $X \in S^2$ (and $\crats=\RR^3$).

With these model propagators $p_{S^2}^{-1}(X)$, the direction of a configured edge $( \confc(x),\confc(y) )$ is $X$ if and only if $( \confc(x),\confc(y) )$ belongs to the propagator $p_{S^2}^{-1}(X)$.
More general propagators produce similar codimension-two constraints on configurations. They allow us to count 
uni-trivalent graphs whose configured edges belong to a finite set of propagators, with signs, as above.

In general, the boundary of a propagator of $C_2(\rats)$ is in the boundary of $C_2(\rats)$, and, for any two-component link embedding $J \sqcup K \colon S^1 \sqcup S^1 \to \crats$, the algebraic intersection of $J\times K \subset \check{C}_2(\rats)$ with a propagator in $C_2(\rats)$ is the linking number of $J$ and $K$. (Note that the linking number of two knots in $\RR^3$ is indeed this algebraic intersection with a model propagator.)

\subsection{Morse propagators}
\label{subMorse}
Let us show examples of propagators in a general rational homology $\RR^3$.

View a rational homology $\RR^3$ as the union of two genus $g$ handlebodies as in Figure~\ref{figHa}, where the two pieces are glued to each other by an a priori nontrivial diffeomorphism of $\partial \handleboda$.\footnote{Unlike the handlebodies in the rest of this book, which are as in Section~\ref{secabs}, the handlebodies of this section are not compact.}

\bfig
\centering
\begin{tikzpicture}
\begin{scope}[xscale=.8, yscale=.8]
\draw (12.8,-3.5) -- (2.8,-3.5) -- (-.2,-1) -- (-.2,0) -- (9.8,0) -- (9.8,-1) -- (12.8,-3.5) -- (12.8,-2.5) -- (9.8,0) (9.8,-1) -- (-.2,-1);
\draw (12.3,-2.7) node{\scriptsize $\handlebodb$};
\begin{scope}[xshift=-5cm]
\draw (8,-1.95) ellipse (.15 and .25);
\draw [blue,thick] (8,-3.2) arc (-90:0:.5);
\draw [blue,thick,->] (8,-2.2) arc (90:0:.5);
\draw [blue,thick,dashed] (8,-2.2) arc (90:270:.5);
\draw [blue] (8.5,-2.7) node[right]{\scriptsize $\curvbeta_1$};
\draw (7,-1.1) .. controls (7.2,-1.35) and (7.2,-3.2) .. (8,-3.2) .. controls (8.5,-3.2) and (9.4,-3.2) .. (9.5,-3.45);
\end{scope}
\draw (5.5,-2.8) node{$\dots$};
\draw (8,-1.95) ellipse (.15 and .25);
\draw [blue,thick] (8,-3.2) arc (-90:0:.5);
\draw [blue,thick,->] (8,-2.2) arc (90:0:.5);
\draw [blue] (8.5,-2.7) node[right=-.01]{\scriptsize$\curvbeta_g$};
\draw [blue,thick,dashed] (8,-2.2) arc (90:270:.5);
\draw (7,-1.1) .. controls (7.2,-1.35) and (7.2,-3.2) .. (8,-3.2) .. controls (8.5,-3.2) and (9.4,-3.2) .. (9.5,-3.45);
\begin{scope}[yshift=-8.2cm]
\draw (1.3,1.5) .. controls (2,1.5) and (2,3.2) .. (3,3.2);
\draw [->] (4.7,1.5) .. controls (4,1.5) and (4,3.2) .. (3,3.2) ;
\draw (12.8,3.5) -- (2.8,3.5) -- (-.2,1) -- (-.2,0) -- (9.8,0) -- (9.8,1) -- (12.8,3.5) -- (12.8,2.5) -- (9.8,0) (9.8,1) -- (-.2,1);
\draw (12.3,2.6) node{\scriptsize $\handleboda$};
\draw (3,2) circle (.5);
\draw (5.5,1.5) node{$\dots$};
\draw (6.3,1.5) .. controls (7,1.5) and (7,3.2) .. (8,3.2) .. controls (9,3.2) and (9,1.5) .. (9.7,1.5);
\draw (8,2) circle (.5);
\draw [red,thick,->] (3,3.2) .. controls (2.9,3.2) and (2.7,3) .. (2.7,2.85);
\draw [red,thick] (3,2.5) .. controls (2.9,2.5) and (2.7,2.7) .. (2.65,2.85) node[right]{\scriptsize $\curvalpha_1$};
\draw [red,thick,dashed] (3,3.2) .. controls (3.1,3.2) and (3.3,3) .. (3.3,2.85) .. controls (3.3,2.7) and (3.1,2.5) .. (3,2.5);
\draw [red,thick,->] (8,3.2) .. controls (7.9,3.2) and (7.7,3) .. (7.7,2.85);
\draw [red,thick] (8,2.5) .. controls (7.9,2.5) and (7.7,2.7) .. (7.65,2.85) node[right]{\scriptsize $\curvalpha_g$};
\draw [red,thick,dashed] (8,3.2) .. controls (8.1,3.2) and (8.3,3) .. (8.3,2.85) .. controls (8.3,2.7) and (8.1,2.5) .. (8,2.5);
\end{scope}
\end{scope}
\end{tikzpicture}
\caption{$\handleboda$ and $\handlebodb$}
\label{figHa}
\end{figure}

The handlebody $\handleboda$ has $g$ arbitrarily oriented meridian disks $D(\curvalpha_i)$ centered at $\crita_i$, for $i$ in  $\underline{g}=\{1,2,\dots,g\}$. The handlebody $\handlebodb$ has $g$ arbitrarily oriented meridian disks $D(\curvbeta_j)$ centered at $\critb_j$, for $j \in \underline{g}$.
The topology of $\crats$ is determined by the curves $\curvalpha_i =\partial D(\curvalpha_i)$ and $\curvbeta_j = \partial D (\curvbeta_j)$ in the surface  $\partial \handleboda$. The data $(\partial \handleboda, (\curvalpha_i)_{i \in \underline{g}}, (\curvbeta_j)_{j \in \underline{g}})$ is called a \emph{Heegaard diagram}. From such data, we construct
\begin{itemize}
\item a \emph{Morse function} $\fMorse \colon \crats \to \RR$ such that $\partial \handleboda=\fMorse^{-1}(1/2)$ 
and the \emph{critical points of $\fMorse$} (the points at which the derivative of $\fMorse$ vanishes) are the $\crita_i$, which have index one and are in $\fMorse^{-1}(1/3)$, and the $\critb_j$, which have index two and are in $\fMorse^{-1}(2/3)$,
\item a gradient vector field $\nabla \colon \crats \to T\crats$ associated to $\fMorse$ and some metric $\metrigo$ such that \begin{equation*}T_xf(y \in T_x \crats)=\langle \nabla(x), y\rangle_{\!\metrigo\,},\end{equation*} and 
\item an associated \emph{gradient flow} $\phi \colon \RR \times \crats \to \crats$ such that 
$\phi(0,.)$ is the Identity map and $\frac{\partial}{\partial t} \phi(t,x)_{(u,x)}=\nabla(\phi(u,x))$.\footnote{The definition of this flow is justified in \cite[Chapter5]{SpivakI}, for example.}
\end{itemize}
Let $\CS$ be the set of critical points of $\fMorse$.
The \emph{flow lines} of the flow $\phi$ are the $\phi(\RR \times \{y\})$, for the $y$ in the complement $\crats\setminus \CS$ in $\crats$ of $\CS$.
For the standard height function $\fMorse_0$ and the standard metric of $\RR^3$, the associated flow is \begin{equation*}\left((\phi_t=\phi(t,.)) \colon x \mapsto x +t\upvec\right),\end{equation*} where $\upvec=(0,0,1)$.
Our rational homology $\RR^3$ (equipped with the metric $\metrigo$) is assumed to coincide with $\RR^3$ (equipped with its standard metric) outside  $\drad{1} \times \left[0,1\right]$, and our Morse function (can be and) is assumed to coincide with $\fMorse_0$ outside  $\drad{1} \times \left[0,1\right]$.

As a set, the complement $\crats\setminus \CS$ of the set of critical points is the disjoint union of the flow lines diffeomorphic to $\RR$, which behave as follows.
There are two flow lines $\linhp(\crita_i)$ and $\linhn(\crita_i)$ starting as vertical lines and ending at $\crita_i$ that approach $\crita_i$ at $+\infty$, as in Figure~\ref{figda}. (They start as vertical lines $\{.\} \times \left]-\infty,0\right[$.) The closure in $\crats$ of their union is a line $\linh(\crita_i)$, called the \emph{descending manifold} of $\crita_i$. It is oriented so that its algebraic intersection with $D(\curvalpha_i)$ is $1$. (This is not consistent with the bottom--top orientation of one of the flow lines. If the flow line with the orientation of the positive normal to $D(\curvalpha_i)$ is  $\linhp(\crita_i)$, then $\linh(\crita_i)=\overline{\linhp(\crita_i)}\cup (-\linhn(\crita_i))$.)
\bfig
\centering
\begin{tikzpicture}
\draw [thick,->] (1,1) -- (1,2) node[left]{\scriptsize $\linhp(\crita_i)$};
\draw [thick] (1,2) .. controls (1,2.7) and (1.3,3) .. (2,3);
\draw [fill=gray!20,draw=white] (2,3) ellipse (.6 and 1);
\draw [red]  (2,3) ellipse (.6 and 1);
\draw [thick,->] (3,1) -- (3,2) node[right]{\scriptsize $\linhn(\crita_i)$};
\draw [draw=white,double=black,very thick] (3,2) .. controls (3,2.7) and (2.7,3) .. (2,3);
\draw [->] (2,3) node[left] {\small $\crita_i$} -- (2,4.3);
\fill (2,3) circle (0.1);
\draw [->] (2,3) -- (2,1.7);
\draw [->] (2,3) -- (2.3,1.85);
\draw (2,2.5) node[rectangle,rotate=20,fill=white] {\tiny $D(\curvalpha_i)$};
\draw [->] (2,3) -- (2.3,4.15);
\draw [->] (2,3) -- (2.65,3.7);
\draw [->] (2,3) -- (2.8,3.2);
\draw [->] (2,3) -- (1.7,4.15);
\draw [red,->] (2,4) arc (90:180:.6 and 1);
\draw [red] (1.4,3) node[left]{\scriptsize $\curvalpha_i$};
\begin{scope}[xshift=-1.5cm]
\draw [thick,-<] (9,4.4) -- (9,3) node[right]{\scriptsize $\linhp(\critb_j)$};
\draw [thick] (8,2) .. controls (9,2.4) .. (9,3);
\draw [fill=gray!20,draw=white] (8,2) circle (1); \draw [blue]  (8,2) circle (1);
\draw [->] (6.7,2) -- (7.5,2);
\draw (7.5,2) -- (8,2);
\draw [->] (9.3,2) -- (8.5,2);
\draw (8.5,2) -- node[above]{\small $\critb_j$} (8,2);
\draw [->] (8,3.3) -- (8,2.5);
\draw (8,2.5) -- (8,2);
\draw [->] (8,.7) -- (8,1.2);
\draw (8,1.2) -- (8,2);
\draw [thick,-<] (7,3.6) -- (7,3) node[left]{\scriptsize $\linhn(\critb_j)$};
\draw [draw=white,double=black,very thick] (8,2) .. controls (7,1.6) .. (7,3);
\fill (8,2) circle (0.1);
\draw (8,1.5) node[rectangle,fill=white] {\tiny $D(\curvbeta_j)$};
\draw [blue] (7.4,3.1) node{\scriptsize $\curvbeta_j$};
\draw [blue,-<] (9,2) arc (0:130:1);
\end{scope}
\end{tikzpicture}
\caption{$\linhp(\crita_i)$, $\linhn(\crita_i)$, $\linhp(\critb_j)$, $\linhn(\critb_j)$}
\label{figda}

\end{figure}
There are two flow lines $\linhp(\critb_j)$ and $\linhn(\critb_j)$ that approach $\critb_j$ at $-\infty$. The closure of their union is a line $\linh(\critb_j)$, which is called the \emph{ascending manifold} of $\critb_j$. It is oriented so that its algebraic intersection with $D(\curvbeta_j)$ is $1$.

The closure of the union of the flow lines that approach $\crita_i$ at $-\infty$ is called the \emph{ascending manifold} of $\crita_i$. It is denoted by $\Asc_i$.  Its intersection with $\handleboda$ is $D(\curvalpha_i)$, and it is oriented like $D(\curvalpha_i)$.
The closure of the union of the flow lines that approach $\critb_j$ at $+\infty$ is called the \emph{descending manifold} of  $\critb_j$. It is denoted by $\Bdesc_j$. Its intersection with $\handlebodb$ is $D(\curvbeta_j)$, and it is oriented like $D(\curvbeta_j)$.
The ascending manifold $\Asc_i$ is an immersion of Figure~\ref{figinA}, which restricts to its interior as an embedding, where the flow lines $\gamma(c_k)$ are flow lines through crossings $c_k$ of $\curvalpha_i \cap \curvbeta_j$, which approach $\crita_i$ near $-\infty$, and $\critb_j$ near $+\infty$. A figure for $\Bdesc_j$ is obtained by reversing the arrows and changing $\crita_i$ to $\critb_j$.

Except for the flow lines of the descending manifolds of the $\crita_i$ and the flow lines of the ascending manifolds of the $\critb_j$, each flow line intersects $\partial \handleboda$ once, transversally, with a positive sign.

\bfig
\centering
\begin{tikzpicture}[scale=.6]
\useasboundingbox (-4.2,-4.1) rectangle (4.2,4.2);
\draw [fill=gray!20, draw=white] (0,0) circle (4);
\draw (4,0) arc (0:220:4);
\draw (4,0) arc (0:-40:4);
\draw [dotted] (0,-4) arc (-90:-40:4);
\draw [dotted] (0,-4) arc (-90:-140:4);
\draw [fill=white, draw=white] (4,0) circle (1.2);
\draw [fill=white, draw=white] (-1.2,4) arc (-180:0:1.2) -- (0,4.2) -- (-1.2,4);
\draw [fill=white, draw=white] (-4,0) circle (1.2);
\draw (-2.8,0) arc (0:83:1.2);
\draw (-2.8,0) arc (0:-83:1.2);
\draw [->] (-2.8,0) arc (0:60:1.2);
\draw [->] (-2.8,0) arc (0:-60:1.2);
\draw (2.8,0) arc (180:97:1.2);
\draw (2.8,0) arc (-180:-97:1.2);
\draw [->] (2.8,0) arc (180:120:1.2);
\draw [->] (2.8,0) arc (-180:-120:1.2);
\draw (3.9,-.7) node{\scriptsize $\linhn(\critb_j)$} (3.9,.67) node{\scriptsize $\linhp(\critb_j)$};
\draw (0,2.8) arc (-90:-173:1.2);
\draw (0,2.8) arc (-90:-7:1.2);
\draw [->] (0,2.8) arc (-90:-150:1.2);
\draw [->] (0,2.8) arc (-90:-30:1.2);
\draw (-1.15,3.56) node[right]{\scriptsize $\linhn(\critb_k)$} (.8,3.1) node[right]{\scriptsize $\linhp(\critb_k)$};
\draw (0,0) circle (1);
\draw [->] (-1,0) arc (180:270:1);
\draw (0,-1.5) node{\scriptsize $\curvalpha_i$};
\draw [->] (0,0) -- (-1.9,0);
\draw [->] (0,0) -- (1.9,0);
\draw [->] (0,0) -- (0,1.9);
\draw (0,0) -- (45:4);
\draw [->] (0,0) -- (45:1.9);
\draw (0,0) -- (35:4);
\draw [->] (0,0) -- (35:1.9);
\draw [->] (0,0) -- (15:1.9);
\draw (15:1.9) .. controls (15:2.5) and (25:3) .. (25:4);
\draw (-1.9,-.5) node{\scriptsize $\gamma(c_3)$} (1.9,-.5) node{\scriptsize $\gamma(c_1)$} (0,1.9) node[left]{\scriptsize $\gamma(c_2)$} (7,1.3);
\draw (-1.9,0) -- (-2.8,0);
\draw (1.9,0) -- (2.8,0);
\draw (0,1.9) -- (0,2.8);
\fill (0,0) node[below]{\scriptsize $\crita_i$} circle (0.1);
\end{tikzpicture}
\caption{The interior of $\Asc_i$.}
\label{figinA}

\end{figure}

The flow lines not in the above descending or ascending manifolds begin as vertical half-lines $x \times \left]-\infty ,0\right[$. They end as vertical half-lines $y \times \left]1,\infty\right[$ for some $x$, $y$ in $\RR^2$.
Except for the critical points, every point has a neighborhood diffeomorphic to a cube $\left]0,1\right[^3$ such that, with the induced identification, the flow maps $(t,x)$ to $x+t\upvec$, for any $(t,x)$ in $\RR \times \left]0,1\right[^3$ such that $x+t\upvec \in \left]0,1\right[^3$.

 In \cite[Theorem 4.2]{lesHC}, Greg Kuperberg and I constructed a propagator from the gradient flow $(\phi_t=\phi(t,.))$ of a Morse function $\fMorse$ without minima and maxima as above, 
as follows.
Let $P_{\phi}$ denote the closure
 in $C_2(\rats)$ of $\{(x,\phi_t(x)) \suchthat x \in \crats \setminus \CS, t \in \left]0,+\infty\right[\}$.
Note that when $(\phi_t \colon x \mapsto x +t\upvec)$ is the flow associated to the standard height function $\fMorse_0$ of $\RR^3$,
$P_{\phi}=p_{S^2}^{-1}(\upvec)$ is one of our model propagators.
Let \begin{equation*}\left[\CJ_{ji}\right]_{(j,i) \in \{1,\dots,g\}^2}=\Bigl[\langle \curvalpha_i,\curvbeta_j \rangle_{\!\partial \handleboda\,}\Bigr]^{-1}\end{equation*}
be the inverse matrix of the matrix of the algebraic intersection numbers $\langle \curvalpha_i,\curvbeta_j \rangle_{\!\partial \handleboda\,}$. (This matrix is invertible because $\crats$ is a rational homology $\RR^3$.)
 
Let $\left((\Bdesc_j \times \Asc_i)\cap C_2(\rats)\right)$ denote the closure of $\left((\Bdesc_j \times \Asc_i)\cap (\crats^2 \setminus \diagonal) \right)$ in $C_2(\rats)$,
then
\begin{equation*}\propP(\funcf,\metrigo)=\preprop_{\phi} + \sum_{(i,j) \in \{1,\dots,g\}^2} \CJ_{ji} \Bigl((\Bdesc_j \times \Asc_i)\cap C_2(\rats)\Bigr)\end{equation*}
is an example of a propagator.

Pick four small generic perturbations $\propP_1$, $\propP_2$, $\propP_3$, and $\propP_4$ of such a propagator
(or of more general propagators as precisely defined in Section~\ref{secprop}). The invariant $w_2$ can be extended to knots $K$ in $\crats$ as follows.
Let $\Gamma$ be one of the graphs  $\diagcross$, $\diagtripod$, or $\diagwheeltwo$.
Orient the edges of $\Gamma$. Number them by the data of an injection $j_E$ from the set $E(\Gamma)$ of the (plain) edges of $\Gamma$ to $\underline{4}$.
Thus, we can count the configurations of $\Gamma$ such that the configured oriented edge numbered by $i$ (viewed as the ordered pair of its ends) is in $\propP_i$, with signs precisely defined.
Denote the average over the choices of such edge-orientations and numberings by
$I_a(K,\Gamma, (\propP_1,\propP_2,\propP_3,\propP_4))$. (When averaging, we divide by the number $2^{\cardlef{E(\Gamma)}}$ of edge orientations and by the number $\frac{4!}{ (4-\cardlef{E(\Gamma)})!}$ of numberings.) Then 
\begin{multline*}w_2(K)= I_a\Bigl(K,\diagcross, \left(\propP_1,\propP_2,\propP_3,\propP_4\right)\Bigr) - I_a\Bigl(K,\diagtripod, \left(\propP_1,\propP_2,\propP_3,\propP_4\right)\Bigr)\\ -2 I_a\Bigl(K,\diagwheeltwo, \left(\propP_1,\propP_2,\propP_3,\propP_4\right)\Bigr)\end{multline*}
does not depend on the chosen propagators. Furthermore, if $K$ is null-homologous, then $w_2(K)$ is again $\frac12\Delta^{\prime\prime}(K)(1)-\frac1{24}$, in this more general setting, as we prove in 
Theorem~\ref{thmcompztwo}.\footnote{In \cite{Leturcq3}, David Leturcq completely expresses
the Alexander polynomial of long null-homologous knots in $\QQ$-spheres in terms of similar counts of configurations.}
 Note that when we compute $w_2(K\subset \RR^3)$ with model propagators, $I_a\left(K,\diagwheeltwo, \left(\propP_1,\propP_2,\propP_3,\propP_4\right)\right)$ vanishes because the double edge gives contradictory constraints.

\subsection{More about the contents of the book}

The cited Altsch\"uler--Freidel universal Vassiliev invariant of knots in $\RR^3$ also extends to knots in a rational homology $\RR^3$. This book describes this extended invariant $\Zinvuf$ using more general (and precisely defined) propagators. Kenji Fukaya proposed a way of counting configurations in a rational homology $\RR^3$ in \cite{FukayaMorse}. Other authors, including Tadayuki Watanabe \cite{watanabeMorse} further studied his approach and made it rigorous. In this book, we can view the resulting way of counting as a particular way of counting with the above \say{Morse propagators} associated to Heegaard diagrams. These Morse propagators are examples of the general \say{propagating chains} described in Chapter~\ref{chapprop}.
Counts of configurations as above also yield 
invariants of the ambient rational homology $\RR^3$ and of its one-point compactification. This compactification is a \emph{rational homology $3$-sphere}, i.e., a connected oriented closed $3$-manifold, where knots have a nontrivial multiple that bounds an immersed oriented compact surface.\footnote{A manifold is said to be \emph{closed} \index[T]{closed!manifold} if it is compact and connected, and if its boundary is empty.}
In \cite{moussardAGT}, Delphine Moussard developed a theory of finite type invariants for rational homology $3$-spheres. Her theory involves surgery operations called \emph{rational Lagrangian-preserving surgeries} instead of crossing changes. We describe these rational Lagrangian-preserving surgeries, which replace a piece of a manifold by another such, in Subsection~\ref{submmcontext} below. 
The above counts of configurations yield a universal finite type invariant $\Zinvuf$ of rational homology $3$-sphere, with respect to Moussard's theory. The universality of $\Zinvuf$ follows from surgery formulae proved in \cite{lessumgen}.

\emph{Integer homology $3$-spheres} are connected oriented closed $3$-manifolds where knots bound an embedded oriented compact surface, as in the standard $3$-sphere $S^3$.
The invariant $\Zinvuf$ also restricts to a universal finite type invariant for integer homology $3$-spheres, with respect to the Ohtsuki theory of finite type invariants \cite{ohtkno} and other equivalent theories described in \cite{ggp}, as first shown by Greg Kuperberg and Dylan Thurston in \cite{kt}.

Let us say a little more about the contents of this book, which
is mostly self-contained. We describe the background material in the appendices.
Unlike in this informal introduction, we will give details and precise statements and carefully check all the assertions.
Our general invariant $\Zinv$ is an infinite series of independent nontrivial invariants. To define it, we must equip any rational homology $\RR^3$ with a parallelization $\tau \colon \crats \times \RR^3 \to T\crats$. We describe the space of suitable parallelizations up to homotopy and the associated Pontrjagin numbers in Chapter~\ref{chapfram}.

This book has four parts.
In its present first part, we first state conventions and known facts about $3$-manifolds.
Then we precisely define propagating chains associated with parallelizations and discuss the invariant $\Theta$ of parallelized rational homology $3$-spheres.
This invariant $\Theta$ is associated to the graph $\tata$. It is the algebraic intersection of three transverse propagating chains.
We explain how to get rid of the dependence of the parallelization of the rational homology 3-spheres with the help of relative Pontrjagin numbers to get an invariant of (unparallelized) rational homology $3$-spheres in 
Section~\ref{secdefThetam}. This invariant is (six times) the Casson--Walker invariant of rational homology spheres. Andrew Casson and Kevin Walker originally defined it as a \say{count} of conjugacy classes of $SU(2)$-representations of the fundamental group of these manifolds \cite{akmc,gm,mar,wal}.

Our general invariant $\Zinvufrf$ is an invariant of parallelized links in rational homology 3-spheres.\footnote{\emph{Framed links} or \emph{parallelized links} are links equipped with a parallel (up to isotopy).} We define it in the second part of the book. It is valued in vector spaces generated by uni-trivalent graphs like $\tata$, $\diagtripod$, and $\diagcross$.
We describe these spaces of diagrams and their rich structures, which help formulate the properties of the invariants, in Chapter~\ref{chapfintype}.
We present the general definition of $\Zinvufrf$ for links in Chapter~\ref{chapdefzinv}. We first give it in terms of integrals rather than in terms of discrete counts because the results are easier to write and prove in the world of differential forms, where no genericity hypotheses are required.
We prove the consistency of our definition and the first properties of $\Zinvufrf$ in Chapters~\ref{chapindepform}, \ref{chapanom}, and \ref{chaprat},
after the needed study of the compactifications of the involved configuration spaces in Chapter~\ref{chapcompconf}.

To compute and use an invariant of links or manifolds, it is interesting to cut links or manifolds into elementary pieces and understand how one can recover the invariant from the invariants of the pieces.
In the third part of this book, we achieve this task with elementary pieces that are \emph{tangles} in rational homology cylinders as in Figure~\ref{figLTRbis}. These tangles are cobordisms between planar configurations of points. We define them precisely in Section~\ref{secintqtangle}.
\bfig
\centering
\begin{tikzpicture}
\useasboundingbox (3,-.8) rectangle (11,4.6);
\draw (4.4,2)node[left]{$\Link(\source) =$};
\draw [thin,->]  (9.2,1.2) node[right]{\scriptsize $\drad{1} \times \left[0,1\right]$} (9.2,1.2) -- (8.5,1.2);
\draw [thin,dashed] (5,0)  .. controls (5,.25) and (6.3,.5) .. (7,.5) .. controls (7.7,.5) and (9,.25) ..  (9,0);
\draw [thin,dashed] (5,3.5)  .. controls (5,3.75) and (6.3,4) .. (7,4) .. controls (7.7,4) and (9,3.75) ..  (9,3.5);
\draw [red,thick] (6.2,.95)  .. controls (6.2,1.05) and (6.8,1.2) .. (7.1,1.2) .. controls (7.4,1.2) and (8,1.05) ..  (8,.95);
\draw [blue, thick,->] (8.1,4.2) -- (8.1,4.6) (8.1,2.7) -- (8.1,4.2);
\draw [blue, thick] (8.6,2.1) .. controls (8.6,1.8) and (8.1,1.6) .. (7.8,1.6) .. controls (7.5,1.6) and (6.7,1.7) .. (6.7,1.9);
\draw [draw=white,double=blue,very thick] (7.6,-.6) -- (7.6,1.4) .. controls (7.6,1.6) and (7.1,2.2) .. (7.1,2.4) .. controls (7.1,2.6) and (7.8,2.9) .. (8,2.9)  .. controls (8.2,2.9) and  (8.6,2.4) .. (8.6,2.1);
\draw [draw=white,double=blue,very thick] (6.7,1.9)  .. controls (6.7,2.2) and  (8.1,2.4) .. (8.1,2.7);
\draw [blue,thick] (6.7,1.9)  .. controls (6.7,2.2) and  (8.1,2.4) .. (8.1,2.7);
\draw [draw=white,double=yellow!80!black,very thick] (5.6,-.8) -- (5.6,1.4) .. controls (5.6,1.6) and (5.9,1.9) .. (6.1,1.9) .. controls (6.3,1.9) and (6.6,1.6) ..  (6.6,1.4) -- (6.6,-.4);
\draw [yellow!80!black, thick,->] (6.6,-.4) -- (6.6,-.7);
\draw [yellow!80!black, thick] (6.6,-.6) -- (6.6,-.8);
\draw [blue, thick,->] (7.6,-.8) -- (7.6,-.6);
\draw [draw=white,double=red,very thick] (6.2,.95)  .. controls (6.2,.8) and (6.8,.7) .. (7.1,.7) .. controls (7.4,.7) and (8,.8) ..  (8,.95);
\draw [draw=white,double=black] (5,0)  .. controls (5,-.25) and (6.3,-.5) .. (7,-.5) .. controls (7.7,-.5) and (9,-.25) ..  (9,0);
\draw [draw=white,double=black] (5,3.5)  .. controls (5,3.25) and (6.3,3) .. (7,3) .. controls (7.7,3) and (9,3.25) ..  (9,3.5);
\draw [thin] (5,0) -- (5,3.5) (9,0) -- (9,3.5);
\end{tikzpicture}
\caption{A tangle representative  in $\drad{1} \times \left[0,1\right]$}
\label{figLTRbis}

\end{figure}
They may be composed in many ways, horizontally, vertically, and by insertions in tubular neighborhoods of other tangle representatives. Such insertions are called \emph{cablings}. We generalize $\Zinvufrf$ to tangles and describe the properties of our generalized $\Zinvufrf$ under the mentioned compositions in the third part of this book.

As already mentioned, a fundamental property of $\Zinvufrf$ is its universality among finite type invariants.
For the restriction of $\Zinvufrf$ to links in $\RR^3$, universality refers to the Vassiliev theory of finite type invariants based on crossing changes. For the restriction of $\Zinvufrf$ to rational homology spheres, universality refers to the Moussard theory based on rational Lagrangian-preserving surgeries. The proofs of universality involve computations of iterated discrete derivatives of $\Zinvufrf$ in the same spirit as the Bar-Natan result recalled in the end of Subsection~\ref{subsecwtwo}. The book's fourth part presents these computations and some consequences for the general invariant $\Zinvufrf$ of framed tangles in rational homology cylinders.

\section{A quicker introduction}
\label{secmathov}

This quicker introduction, for experienced topologists, is independent of the first one. Beginners can read it after the warm-up of the slower one.
Here, we describe the invariant $\Zinvuf$ of $n$-component links $\Link$ in rational homology $3$-spheres $\rats$ studied in this book, more precisely.
We also specify some notions vaguely introduced in the slow introduction of Section~\ref{secslowbeg}, and
we say more on the mathematical landscape around $\Zinvuf$.

\subsection{On the construction of \texorpdfstring{$\Zinvuf$}{Z}}

The invariant $\Zinvuf(\Link)=\bigl(\Zinvuf_k(\Link)\bigr)_{k\in \NN}$ of an $n$-component link $\Link$ in a rational homology $3$-sphere $\rats$ is valued in a graded space generated by uni-trivalent graphs. Its degree $k$ part is a sum \begin{equation*}\Zinvuf_k(\Link)=\sum_{\Gamma}\Zinvuf_{\Gamma}(\Link)\left[\Gamma\right],\end{equation*} running over such graphs $\Gamma$ with $2k$ vertices. The coefficient $\Zinvuf_{\Gamma}(\Link)$ \say{counts} embeddings of $\Gamma$ in $\rats$ mapping the univalent vertices of $\Gamma$ to $\Link$, in a sense explained in the book. 
Let us slightly specify that sense.

For technical reasons, we
remove a point $\infty$ from our rational homology $3$-spheres $\rats$ to transform them into open manifolds $\crats$\index[N]{Rmanif@$3$-manifolds!Rcheck@$\crats$ punctured $3$-manifold}. When $\rats$ is the standard sphere $S^3$, the punctured $\QQ$-sphere $\crats$ is $\RR^3$.

Let $\diag(\crats^2)$ denote the diagonal of $\crats^2$.
Following William Fulton, Robert MacPherson \cite{FultonMcP}, Maxim Kontsevich \cite{ko}, Scott Axelrod, Isadore Singer \cite[Section 5]{axelsingII}, and others,
we will introduce a suitable smooth compactification $C_2(\rats)$ (with boundary and ridges) of $\crats^2 \setminus \diag(\crats^2)$ such that the map \index[N]{Projections!pStwo@$p_{S^2}$ to $S^2$}
\begin{equation*}\begin{array}{llll}p_{S^2}\colon &\left(\RR^3\right)^2 \setminus \diag\bigl(\left(\RR^3\right)^2\bigr)&\rightarrow & S^2\\&(x,y) &\mapsto &\frac{1}{\lVert y-x \rVert}(y-x) \end{array}\end{equation*}
extends to $C_2(S^3)$.
We will introduce a notion of \emph{propagating chain} and the dual notion of  \emph{propagating form} for $\rats$.  When $\rats=S^3$, for any $X \in S^2$, the submanifold $p_{S^2}^{-1}(X)$ of $C_2(S^3)$ is an example of a propagating chain. For any $2$-form $\omega_{S}$ on $S^2$ such that $\int_{S^2}\omega_{S}=1$, the form $p_{S^2}^{\ast}(\omega_{S})$ is an example of a propagating form. The propagating chain $p_{S^2}^{-1}(X)$ is a \emph{model propagating chain}. The propagating form $p_{S^2}^{\ast}(\omega_{S})$ is a \emph{model propagating form}.
In general, a propagating chain is a $4$-dimensional rational chain (i.e., a finite rational combination of oriented compact $4$-manifolds with possible ridges) of $C_2(\rats)$, while a propagating form is a closed $2$-form on $C_2(\rats)$.
Both have to satisfy some conditions on the boundary of $C_2(\rats)$, which make them share sufficiently many properties with our model propagating chains or forms.

In particular, for any propagating form $\omega$, for any propagating chain $\propP$, and for any two-component link $(J,K) \colon S^1 \sqcup S^1 \rightarrow \crats$, we have 
\begin{equation*}\int_{J \times K \subset C_2(\rats)}\omega= \langle J \times K, \propP\rangle_{\!C_2(\rats)\,}= lk(J,K),\end{equation*}
where $\langle ., .\rangle_{\!C_2(\rats)\,}$ stands for the algebraic intersection in $C_2(\rats)$, and $lk$ is the linking number in $\rats$. 

Propagating forms and propagating chains are both called \emph{propagators} when their nature is clear from the context.
The above equalities tell us in which way \say{propagators represent the linking form}.

A \emph{Jacobi diagram}\index[T]{Jacobi diagram} $\Gamma$\index[N]{Gamma@$\Gamma$ Jacobi diagram} on $\sqcup_{i=1}^n S^1$ is a uni-trivalent graph $\Gamma$ equipped with an isotopy class of injections from its set $U(\Gamma)$\index[N]{U@$U(\Gamma)=\{\mbox{univalent vert. of }\Gamma\}$} of univalent vertices into the domain $\sqcup_{i=1}^n S^1$ of a link $\Link$. In the figures, we represent the domain of $\Link$ by dashed circles and put the univalent vertices of $\Gamma$ on their images under an injection of the given isotopy class as in Figure~\ref{figfirstJacdiag}.
\bfig
\centering
\cirtd
\caption{A (plain) Jacobi diagram on (the dashed) $S^1$}
\label{figfirstJacdiag}

\end{figure}
Let $V(\Gamma)$\index[N]{V@$V(\Gamma)=\{\mbox{vertices of }\Gamma\}$}, $T(\Gamma)$\index[N]{T@$T(\Gamma)=\{\mbox{trivalent vert. of }\Gamma\}$}, and $E(\Gamma)$\index[N]{E@$E(\Gamma)=\{\mbox{edges of }\Gamma\}$} respectively denote the set of vertices, trivalent vertices, and edges of $\Gamma$.
The \emph{configuration space}\index[T]{configuration!space} $\check{C}(\rats,\Link;\Gamma)$ is the set of injections from $V(\Gamma)$ to $\crats$ that map the set $U(\Gamma)$ of univalent vertices of $\Gamma$ to $\Link$ and induce the given isotopy class of injections. It is an open submanifold of $\crats^{T(\Gamma)}\times \Link^{U(\Gamma)}$. Assume that $\Gamma$ has no \emph{looped edge} like \loopedge and that its edges are oriented.
Then each edge $e$ of $\Gamma$ provides a natural restriction map $p(\Gamma,e)$ from $\check{C}(\rats,\Link;\Gamma)$ to $\crats^2\setminus \diag(\crats^2)$. When propagating forms $\omega(e)$ are associated to the edges, this allows one to define a real number \begin{equation*}I\Bigl(\rats,\Link,\Gamma,\bigl(\omega(e)\bigr)_e\Bigr)=\int_{\check{C}(\rats,\Link;\Gamma)} \bigwedge_{e \in E(\Gamma)} p(\Gamma,e)^{\ast}\bigl(\omega(e)\bigr).\end{equation*}
Similarly and dually, when propagating chains $\propP(e)$ in general position are associated to the edges, one can define a rational number 
$I\left(\rats,\Link,\Gamma,(\propP(e))_e\right)$ as the algebraic intersection of the codimension-two chains $p(\Gamma,e)^{-1}(\propP(e))$ in $\check{C}(\rats,\Link;\Gamma)$.

For example, the Jacobi diagram $\Gamma$ pictured as \onechordtwocirclesrightmoyen{i}{j},  is an oriented edge $e$ from a univalent vertex that must go to the component $K_i$ of $\Link$ to another univalent vertex, which must go to another component $K_j$ of $\Link$.
The associated configuration space $\check{C}(\rats,\Link;\Gamma)$ is $K_i \times K_j$,
According to the given property of our propagators, we have 
\begin{equation*}I\Bigl(\rats,\Link,\onechordtwocirclesrightmoyen{i}{j},\omega(e)\Bigr)=\int_{K_i \times K_j \subset C_2(\rats)}\omega(e)= lk(K_i,K_j)\end{equation*} 
and
\begin{equation*}I\Bigl(\rats,\Link,\onechordtwocirclesrightmoyen{i}{j},\propP(e)\Bigr)= \bigl\langle K_i \times K_j, \propP(e) \bigr\rangle_{\!C_2(\rats)\,}=lk(K_i,K_j)\end{equation*}
for any propagating form $\omega(e)$ and for any propagating chain $\propP(e)$.

As another example, consider the Jacobi diagram $\tata$ with two trivalent vertices and three edges $e_1$, $e_2$, and $e_3$ from one vertex to the other. When $\rats$ is a $\ZZ$-sphere, we will show how one can choose propagators $\propP(e_i)$ and $\omega(e_i)$ so that we have
\begin{equation*}I\Bigl(\rats,\emptyset,\tata,\bigl(\omega(e)\bigr)_{e \in \{e_1,e_2,e_3\}}\Bigr)=I\Bigl(\rats,\emptyset,\tata,\bigl(\propP(e)\bigr)_{e \in \{e_1,e_2,e_3\}}\Bigr)=6 \lambda_{CW}(\rats),\end{equation*} 
where $\lambda_{CW}$ \index[N]{lambdaCW@$\lambda_{CW}$ Casson--Walker invariant} is the Casson invariant normalized as in \cite{akmc,gm,mar}.
We can choose the above propagating forms so that $\omega(e_1)=\omega(e_2)= \omega(e_3)$. In particular, the Casson invariant, which may be written as \begin{equation*}\lambda_{CW}(\rats)=\frac16\int_{C_2(\rats)}\omega(e_1)^3,\end{equation*} may be viewed as a \say{cube of the linking number}.

In general, when $\rats$ is a $\ZZ$-sphere, one can choose a propagating form $\omega$ for $\rats$ and set $\omega(e)=\omega$ for all edges $e$ of Jacobi diagrams so that
the real coefficient $\Zinvuf_{\Gamma}(\Link)$ in \say{$\Zinvuf(\Link)=\sum_{\Gamma}\Zinvuf_{\Gamma}(\Link)\left[\Gamma\right]$} is the product of $I(\rats,\Link,\Gamma,(\omega)_e)$ and a constant which depends only on $\Gamma$. The coefficient $\Zinvuf_{\Gamma}(\Link)$ can alternatively be obtained by averaging some $I(\rats,\Link,\Gamma,(\propP(e))_e)$ over ways of equipping edges of $\Gamma$ by propagating chains in a fixed set of generic propagating chains, and over ways of orienting the edges of $\Gamma$. The coefficients $\Zinvuf_{\Gamma}(\Link)$ depend on propagator choices, but relations among Jacobi diagrams in the target space $\Aavis(\sqcup_{i=1}^n S^1)$ of $\Zinvuf$ ---generated by Jacobi diagrams---ensure that $\Zinvuf(\Link)$ is an isotopy invariant.
The invariant $\Zinvuf$ may be thought of as a series of higher-order linking invariants.

The definition of $\Zinvuf$ that is presented here is a generalization of the definition that I explained in detail for $\QQ$-spheres in the unsubmitted preprint \cite{lesconst}, which was inspired by \cite{kt} and discussions with Dylan Thurston in Kyoto in 2001. 

The present definition also includes links in $\QQ$-spheres. Most of the additional arguments involved in the construction for links already appear in many places. We repeat them to make the book as self-contained as possible. We present many variants of the definitions and make them as flexible as possible because the flexibility has proved useful in many generalizations and applications of these constructions,
such as equivariant constructions in \cite{lesbonn,lesuniveq},
or the recent explicit computations of integrals over configuration spaces by David Leturcq \cite{Leturcq2,Leturcq3}, in addition to the applications presented in this book. David Leturcq obtained an expression of the Alexander polynomial of knots in $\QQ$-spheres in terms of such integrals. His expression identifies some combinations of integrals over configuration spaces with coefficients of the Alexander polynomial.

With our flexible definition of a propagating chain,
there is a natural propagating chain associated to a generic Morse function on a punctured $\QQ$-sphere, and to a generic metric, as described in Subsection~\ref{subMorse}. The main part of such a propagator is the space of pairs of points on a gradient line such that the second point is after the first one.
Greg Kuperberg and I constructed such a Morse propagator in \cite{lesHC} for Morse functions without minima or maxima.
Independent work of Tadayuki Watanabe \cite{watanabeMorse} allows one to generalize these propagators to any Morse function. Thus, up to some corrections, $\Zinvuf$ counts embeddings of graphs whose edges embed in gradient lines of Morse functions as in a Fukaya article \cite{FukayaMorse}. Tadayuki Watanabe has used similar constructions in his recent construction of exotic elements in the fundamental group of the group of diffeomorphisms of $S^4$ \cite{watanabe2018exotic}.
This book contains a framework to study these questions precisely.\footnote{We mainly develop the framework for tangles in dimension $3$, but many arguments can be easily adapted in higher dimensions. See \cite{watanabe2018exotic, Leturcq1, Leturcq2}, for example.}

We will also show how the construction of $\Zinvuf$ extends to tangles in rational homology cylinders so that $\Zinvuf$ extends to a functor $\Zinvufrf$ on a category of framed tangles with many important properties.
These properties of the functor $\Zinvufrf$ are stated in Theorem~\ref{thmmainfunc}, one of this book's main original theorems. They provide tools to reduce the computation of $\Zinvufrf$ to its evaluation at elementary pieces.

\subsection{More mathematical context}
\label{submmcontext}

\paragraph{Finite type invariants.}

The \emph{finite type invariant} concept for knots was introduced in the 90's in order to classify knot invariants, with the work of Victor Vassiliev, Mikhail Goussarov, and Dror Bar-Natan, shortly after the birth of numerous quantum knot invariants, described by Vladimir Turaev in \cite{Turaevbookquantum}. Tomotada Ohtsuki extended this very useful concept to $3$-manifold invariants \cite{ohtkno}. See also \cite{ohtsukibookquantum}. Theories of finite type invariants in dimension $3$ are defined from a set $\CO$ of operations on links or $3$-manifolds. 
In the case of links in $\RR^3$, $\CO$ is the set $\CO_{V}$ of crossing changes \pcpetit $\leftrightarrow$ \ncpetit.
The variation of an invariant $\lambda$ under an operation of $\CO$ may be thought of as a discrete derivative. When $k$ independent operations $o_1,\dots, o_k$ on a pair $(\rats,\Link)$ consisting of a link $\Link$ in a $\QQ$-sphere $\rats$ are given, for a part $I$ of $\{1,\dots,k\}$ with cardinality $\cardlef{I}$, the pair $(\rats,\Link)((o_i)_{i \in I})$ is the pair obtained from $(\rats,\Link)$ by applying the operations $o_i$ for $i \in I$. Then the alternate sum \begin{equation*}\sum_{I \subseteq \{1,\dots,k\}} (-1)^{\cardlef{I}} \lambda\Bigl((\rats,\Link)\bigl((o_i)_{i \in I}\bigr)\Bigr)\end{equation*}
may be thought of as the $k^{th}$ derivative of $\lambda$ with respect to $\{o_1,\dots, o_k\}$ at $(\rats,\Link)$.
An \emph{invariant of degree at most $k$} with respect to $\CO$ is an invariant whose degree $k+1$ derivatives vanish. A \indexT{finite type invariant} with respect to $\CO$ is an invariant that is of degree at most $k$ for some positive integer $k$.
Finite type invariants of links in $\RR^3$ with respect to the set of crossing changes are called \emph{Vassiliev invariants}\index[T]{Vassiliev invariant}.

In this case of links in $\RR^3$, Daniel Altsch\"uler and Laurent Freidel \cite{af} proved that the invariant $\Zinvuf$ described in this book is a \emph{universal Vassiliev invariant}, meaning that all real-valued Vassiliev invariants of links in $\RR^3$ factor through $\Zinvuf$. Since all the quantum invariants of \cite{Turaevbookquantum} can be viewed as sequences of finite type invariants, $\Zinvuf$ also contains all these invariants such as the Jones polynomial, its colored versions, the HOMFLY polynomial...
Dylan Thurston proved similar universality results in \cite{thurstonconf} independently. He also showed that $\Zinvuf$ is rational.
Further substantial work of Sylvain Poirier in \cite{poirierv2} allowed me to identify the invariant $\Zinvuf$ with the famous Kontsevich integral of links in $\RR^3$---described in  \cite{barnatan}, in \cite[Chapter 6]{ohtsukibookquantum}, and in \cite{chmutovmostov} by Sergei Chmutov, Sergei Duzhin, and Jacob Mostovoy---up to a change of variables described in \cite{lesunikon} in terms of an \say{anomaly}, which is sometimes called the \emph{Bott--Taubes anomaly}.

Let us now decribe operations on $3$-manifolds.
The boundary $\partial A$ of a genus $g$ $\QQ$-handlebody is a closed oriented genus $g$ surface.
The \indexT{Lagrangian} $\CL_A$ \index[N]{Lagrangians!LA@$\CL_A$} of a compact $3$-manifold $A$ is the kernel of the map induced by the inclusion from $H_1(\partial A;\QQ)$ to $H_1(A;\QQ)$.
(In Figure~\ref{fighandlebodyg} of $H_g$, the Lagrangian of $H_g$ is freely generated by the classes of the curves $a_i$.)

An {\em integral (resp. rational) Lagrangian-Preserving (or LP) surgery\/} \index[T]{Lagrangian-Preserving surgery}\index[T]{LP surgery} $(A^{\prime}/A)$ is the replacement of an integral (resp. rational) homology handlebody $A$ embedded in the interior of a $3$-manifold $M$ by another such $A^{\prime}$ whose boundary $\partial A^{\prime}$ is identified with $\partial A$
by an orientation-preserving diffeomorphism that sends $\CL_{A^{\prime}}$ to $\CL_A$.

Theories of finite type invariants of integer (resp. rational) homology $3$-spheres $\rats$ can be defined from the set $\CO^{\ZZ}_{\CL}$  \index[N]{Operation sets!OLZ@$\CO^{\ZZ}_{\CL}$} (resp. $\CO^{\QQ}_{\CL}$ \index[N]{Operation sets!OLQ@$\CO^{\QQ}_{\CL}$}) of integral (resp. rational) LP-surgeries. For $\ZZ$-spheres, results of Kazuo Habiro \cite{habiro}, Stavros Garoufalidis, Mikhail Goussarov, and Michael Polyak \cite{ggp}, and Emmanuel Auclair and me \cite{aucles} imply that
the theory of real-valued finite type invariants with respect to $\CO^{\ZZ}_{\CL}$ is equivalent to the original theory defined by Tomotada Ohtsuki in \cite{ohtkno} using surgeries on algebraically split links. We will call this theory the \emph{Ohtsuki--Goussarov--Habiro theory.}

Greg Kuperberg and Dylan Thurston first showed that the restriction of $\Zinvuf$ to integer homology $3$-spheres (equipped with empty links) is a \emph{universal finite type invariant} of $\ZZ$-spheres with respect to $\CO^{\ZZ}_{\CL}$ in \cite{kt}.

As in the case of links in $\RR^3$, their proof of universality rests on a computation of the 
$k^{th}$ derivatives of the degree $k$ part $\Zinvuf_k$ of the invariant $\Zinvuf=(\Zinvuf_k)_{k \in \ZZ}$, which proves that $\Zinvuf_k$ is a degree $k$ invariant whose $k^{th}$ derivatives are universal in the following sense. All the $k^{th}$ derivatives of degree $k$ real-valued invariants factor through them.

The \say{Universality part} of this book will be devoted to a general computation of the $k^{th}$ derivatives of the extension of $\Zinvuf_k$ to tangles with respect both to $\CO_V$ and $\CO^{\QQ}_{\CL}$ (which contains $\CO^{\ZZ}_{\CL}$). The resulting formulae stated in Theorems~\ref{thmunivsingtang} and \ref{thmmainunivlag} are crucial properties of $\Zinvuf$. Theorem~\ref{thmmainunivlag} is one of the main original results of this book.

I first proved the splitting formulae, which compute the $k^{th}$ derivatives of the degree $k$ part $\Zinvuf_k$ with respect to $\CO^{\QQ}_{\CL}$, in \cite{lessumgen} for the restriction of $\Zinvuf_k$ to $\QQ$-spheres.
They allowed Delphine Moussard to classify finite type invariants of $\QQ$-spheres with respect to $\CO^{\QQ}_{\CL}$ in \cite{moussardAGT}. In particular, she proved that, when associated with the $p$-valuations of the cardinality $\cardlef{H_1}$ of the torsion first homology group, $\Zinvuf$ is a \emph{universal finite type invariant} of $\QQ$-spheres with respect to $\CO^{\QQ}_{\CL}$. Together with results of Gw\'ena\"el Massuyeau \cite{massuyeausplit} who proved that the LMO invariant $Z_{LMO}$ of Thang L\^e, Jun Murakami, and Tomotada Ohtsuki \cite{lmo} satisfies the same formulae, the Moussard classification implies that $\Zinvuf$ and $Z_{LMO}$ are equivalent in the sense that they distinguish the same $\QQ$-spheres with identical $\cardlef{H_1}$.

Thus, the invariant $\Zinvuf$ is as powerful as the famous LMO invariant for $\QQ$-spheres, and as the famous Kontsevich integral for links. 
The Kontsevich integral $Z^K$ and its relations with the theory of quantum groups developed by Drinfeld and Jimbo have been extensively studied. See the books \cite{Kasselbook,ohtsukibookquantum,chmutovmostov}, for instance.
Thang L\^e, Jun Murakami \cite[Theorem 10]{LeMur}, and Christian Kassel \cite[Theorem XX.8.3]{Kasselbook} independently showed how the Turaev quantum link invariants \cite{TuraevYB88}, which include the HOMFLYPT and Kauffman polynomials, can be recovered from  $Z^K$. See also \cite[Theorem 6.14]{ohtsukibookquantum}. 
An explicit way of recovering the Alexander polynomial from $Z^K$ for knots in $\RR^3$ based on \cite{bngaroufMMR} can be found in \cite[\S 11.2.4]{chmutovmostov}.
Tomotada Ohtsuki computed the two-loop part of the Kontsevich integral---which coincides with the two-loop part of $\Zinvuf$ according to Corollary~\ref{corlesunikontwoleg} and Note~\ref{rkalphafive}---for the genus one knots \cite{ohttwoloop}.
Dror Bar-Natan, Thang L\^e, and Dylan Thurston computed the Kontsevich integral of the trivial knot $O$ \cite{bltwheels}.\footnote{The computation for the torus knots follows as in \cite[Proposition 3.16]{lesintrokonsummer} or \cite[Section 9.3]{chmutovmostov}. See also the computations of Julien March\'e in \cite{MarcheTorusAGT}.}
 Since $\Zinvuf$ is obtained from the Kontsevich integral by a change of variables determined by the anomaly as in Corollary~\ref{corlesunikontwoleg}, these results have direct corollaries for $\Zinvuf$.

The LMO invariant and its generalizations for links in $\QQ$-spheres \cite{lmo, ohtsukibookquantum, bnaarhus1, bnaarhus2,bnaarhus3} are defined from the Kontsevich integral of links in $\RR^3$, in a combinatorial way. Any compact oriented $3$-manifold can be presented by a \emph{framed}\index[T]{framed!link} link of $\RR^3$, which is a link equipped with a favorite parallel, according to a theorem proved by Raymond Lickorish and Andrew Wallace, independently and nicely reproved by Colin Rourke in \cite{rourke}. The \emph{Kirby moves} are specific modifications of framed links that do not change the presented manifold. According to a theorem of Robion Kirby, two framed links present the same manifold if and only if they are related by a finite sequence of Kirby moves.
The LMO invariant of a $3$-manifold is defined from the Kontsevich integral of a framed link that presents such a manifold.
The proof of its invariance relies on the cited Kirby theorem. The combinatorial nature of $Z_{LMO}$ allowed Dror Bar-Natan and Ruth Lawrence to compute $Z_{LMO}$ for lens spaces and Seifert fibered spaces \cite{barnatanlawr}.

When restricted to braids, the Kontsevich integral has a natural geometric meaning. It measures how the strands turn around each other (see \cite{chmutovmostov} or \cite[Section 1]{lesintrokonsummer}), and it defines morphisms from braid groups to algebras of horizontal chord diagrams. Dror Bar-Natan extended the Kontsevich integral to links \cite{barnatan}. Thang L\^e and Jun Murakami extended the Kontsevich integral to a functor from framed tangles to a category of Jacobi diagrams \cite{LeMur}. 
Next, Thang L\^e, Jun Murakami, and Tomotada Ohtsuki defined the LMO invariant from the Le--Murakami--Kontsevich invariant of surgery presentations of the $3$-manifolds using tricky algebraic manipulations of Jacobi diagrams with the help of Kirby calculus \cite{lmo,ohtsukibookquantum}. Though some of the physical meaning of the LMO invariant can be recovered from its universality properties, much of it gets lost in the manipulations. 

The presented construction of $\Zinvuf$ is much more physical, geometric, and natural---at least to me.
It does not rely on the Kirby theorem and provides information about graph embeddings in $\QQ$-spheres. On the other hand, some explicit computations or properties available for the LMO invariant have yet to be performed or proved for $\Zinvuf$.

\section{Book organization}
\label{secorganiz}

Chapter~\ref{chaplk} completes our slow introduction to the linking number of Sections~\ref{sublkGauss} and \ref{sublkGausstwo}. It contains more conventions and arguments used throughout this book.  
Chapter~\ref{chapprop} introduces our definitions of propagators.
These propagators are the basic ingredients of all our constructions. Most of the time, they are associated with a parallelization of the $3$-manifold.
The Theta invariant is the simplest $3$-manifold invariant that can be derived from the techniques described in this book. We present it in detail in Chapter~\ref{chapTheta}. We first describe $\Theta$ as an invariant of a parallelized punctured $\QQ$-sphere $(\crats,\tau)$. The invariant $\Theta(\rats,\tau)$ is the intersection of three propagating chains associated with the given parallelization $\tau$ of $\crats$ in the two-point configuration space $C_2(\rats)$. Equivalently, the invariant $\Theta(\rats,\tau)$ is the integral over $C_2(\rats)$ of the cube of a propagating form associated with $\tau$. Next, we transform $\Theta$ to an invariant of $\QQ$-spheres using relative Pontrjagin classes, also called Hirzebruch defects, as Greg Kuperberg and Dylan Thurston did in \cite{kt}.
Like $\Theta$, the series $\Zinvuf$ comes from a more natural invariant $\Zinvufrf$ of parallelized links $(\Link,\Link_{\parallel})$ in parallelized punctured $\QQ$-spheres, constructed with associated propagators. We transform $\Zinvufrf\bigl(\rats,\Link,\Link_{\parallel},\tau\bigr)$ to an invariant $\Zinvuf$ of links $\Link$ in $\QQ$-spheres using a function of linking numbers associated to the link parallelizations, Pontrjagin numbers associated to the manifold parallelizations, and constants $\alpha$ and $\ansothree$ called \emph{anomalies}.
Chapter~\ref{chapfram} presents parallelizations of oriented $3$-manifolds with boundaries and associated Pontrjagin numbers in detail. It closes this introductory part.

Thus, the book's first part describes the degree one part of the graded invariant $\Zinvuf$ for links in $\QQ$-spheres, which is determined by the linking numbers of the components and the $\Theta$-invariant of the ambient manifold. The book's second part is devoted to the general presentation of $\Zinvuf$ for links in $\QQ$-spheres. 
In this part, we first review various theories of finite type invariants for which specializations of $\Zinvuf$ will be universal finite type invariants. This allows us to introduce 
the spaces of Jacobi diagrams in which $\Zinvuf$ takes its values, in a natural way, in Chapter~\ref{chapfintype}.
The complete definitions of $\Zinvuf$ for links in $\QQ$-spheres are given in Chapter~\ref{chapdefzinv} without proofs of consistency.
We show that these definitions make sense and do not depend on the involved
choices of propagating forms in Chapters~\ref{chapindepform} and \ref{chapanom}. Our proofs rely on the study of suitable compactifications of configuration spaces presented in Chapter~\ref{chapcompconf} and on some standard arguments of the subject already appearing in many places starting with \cite{ko}, \cite{botttaubes}\dots
This second part of the book ends with discrete equivalent definitions of $\Zinvuf$ in terms of propagating chains and algebraic intersections rather than propagating forms and integrals in Chapter~\ref{chaprat}. These definitions make clear that the invariant $\Zinvuf$ is rational.
The other main properties of $\Zinvuf$ are precisely described in Sections~\ref{secinvql}, \ref{secfstpropzinv}, \ref{secexppty}, and Chapter~\ref{chapuniv}. Some of them involve the extension of $\Zinvuf$ to tangles, which can be found in Theorem~\ref{thmfstconsttang}.

The book's third part is devoted to this extension of $\Zinvuf$ to tangles, still denoted by $\Zinvuf$. We introduce the framed version $\Zinvufrf$ of $\Zinvuf$ for framed tangles in Definition~\ref{deffstconsttangframed}. The spirit of the definition is the same. However, its justification is more difficult because the involved compactified configuration spaces are more complicated.
They are no longer smooth manifolds with ridges and have additional types of faces. We first present the definition and properties of the extension without proofs in Chapters~\ref{chapfirstintfunc} and \ref{chapzinvtang}.
Next, we justify them in Chapters~\ref{chapconszinvf} and \ref{chappropzinvffunc}, respectively.
In particular, Chapter~\ref{chappropzinvffunc} contains the proofs of many properties of the link invariant $\Zinvuf$. These proofs involve easy-to-discretize variants of the functor $\Zinvufrf$, which are interesting on their own and presented in Section~\ref{secvarzinvf}. The third part of the book ends with the computation of the iterated derivatives of the generalized $\Zinvuf$
with respect to crossing changes, in Section~\ref{secproofthmbn}. This computation proves that the restriction of $\Zinvuf$ to links in $S^3$ is a universal Vassiliev invariant.

The fourth part focuses on the computation of the iterated derivatives of the generalized $\Zinvuf$
with respect to rational LP-surgeries. It begins with Chapter~\ref{chapuniv}, stating the main results and their corollaries, and reducing their proofs to the proofs of two key propositions, presented in Chapter~\ref{chapsimnormprop}.
The proofs of these propositions involve the introduction of a more flexible definition of $\Zinvuf$ because the restriction of a parallelization of a $\QQ$-sphere to the exterior of a $\QQ$-handlebody does not necessarily extend to a $\QQ$-handlebody that replaces the former one during a rational LP-surgery. Chapter~\ref{chappseudopar} contains an extension of the notion of parallelization to a more flexible notion of \emph{pseudo-parallelization} and a corresponding more flexible definition of $\Zinvuf$. Pseudo-parallelizations also have associated propagators and Pontrjagin numbers. They easily extend to arbitrary $\QQ$-handlebodies.
More flexible variants of the definition of $\Zinvuf$ based on pseudo-parallelizations can be found in Chapter~\ref{chappseudoparmuchmore}.

The book ends with two appendices. Appendix~\ref{chaphomology} lists the basic results and techniques of algebraic topology used in the book. Appendix~\ref{chapDeRham} reviews the used properties of differential forms and de Rham cohomology.

Most chapters have their own detailed introduction. Many cross-references help the reader choose what she/he wants to read.

\section{Book genesis}
\label{secgenesis}
At first, this book aimed at presenting the results of two preprints \cite{lesconst,lessumgen} and lecture notes \cite{lesmek}. It contains generalizations of the results of these preprints to wider settings. I have never submitted the preprints \cite{lesconst,lessumgen} for publication.
Edward Witten's insight into the perturbative expansion of the Chern--Simons theory \cite{witten} inspired the mathematical guidelines for constructing the invariant $\Zinvuf$. Maxim Kontsevich gave them in \cite[Section~2]{ko}. Greg Kuperberg and Dylan Thurston developed these guidelines in \cite{kt}. They defined $\Zinvuf$ for $\QQ$-spheres and sketched a proof that the restriction of $\Zinvuf$ to $\ZZ$-spheres is a universal finite type invariant of $\ZZ$-spheres in the Ohtsuki--Goussarov--Habiro sense. This allowed them to identify the degree one part of $\Zinvuf$ with the Casson invariant for $\ZZ$-spheres. I thank Dylan Thurston for explaining to me his joint work with Greg Kuperberg in Kyoto in 2001. 

In \cite{lessumgen}, I proved splitting formulae for $\Zinvuf$. These formulae compute derivatives of $\Zinvuf$ with respect to rational LP-surgeries. They generalize similar Kuperberg--Thurston implicit formulae about Torelli surgeries. These formulae allowed me to identify the degree one part of $\Zinvuf$ with the Walker generalization of the Casson invariant for $\QQ$-spheres, in \cite[Section 6]{lessumgen}. They also allowed Delphine Moussard to classify finite type invariants with respect to these rational LP-surgeries and prove that all such real-valued finite type invariants factor through some \say{augmentation} of $\Zinvuf$ by invariants derived from the order of the $H_1(.;\ZZ)$, in \cite{moussardAGT}. In \cite{massuyeausplit},  Gw\'ena\"el Massuyeau proved that the LMO invariant of Thang L\^e, Jun Murakami, and Tomotada Ohtsuki \cite{lmo} satisfies the same splitting formulae as $\Zinvuf$. Thus, the Moussard classification implies that $\Zinvuf$ and $Z_{LMO}$ are equivalent in the sense that they distinguish the same $\QQ$-spheres with identical $\cardlef{H_1(.;\ZZ)}
$.
In order to write the proof of my splitting formulae, I needed to specify the definition of $\Zinvuf$
and I described the Kontsevich--Kuperberg--Thurston construction in detail in \cite{lesconst}.

In \cite{lesmek}, mixing known constructions in the case of links in $\RR^3$ with the construction of $\Zinvuf$ allowed me to define a natural extension of $\Zinvuf$ as an invariant of links in $\QQ$-spheres. This extension also generalizes invariants of links in $\RR^3$ defined by Enore Guadagnini, Maurizio Martellini, and Mihail Mintchev \cite{gmm}, Dror Bar-Natan \cite{barnatanper}, and by Raoul Bott and Clifford Taubes \cite{botttaubes}, which emerged after the Witten work \cite{witten}.\footnote{The relation between the perturbative expansion of the Chern--Simons theory of the Witten article and the configuration space integral viewpoint is explained by Michael Polyak in \cite{PolyakFeyn} and by Justin Sawon in \cite{Sawon}. See also \cite[Appendix F]{ohtsukibookquantum}}  I also gave more flexible definitions of $\Zinvuf$.

In addition to the revisited contents of the preprints \cite{lesconst,lessumgen} and of the notes \cite{lesmek}, this book contains an extension of $\Zinvuf$ as a functorial invariant of tangles in rational homology cylinders. It also contains the proofs of many properties of this extension, which imply simpler properties for $\Zinvuf$, such as the multiplicativity of $\Zinvuf$ under connected sum, for example. This functorial extension, which generalizes the Poirier extension in \cite{poirierv2}, and its properties are new. They appear only in this book (to my knowledge).

Most of the properties of $\Zinvuf$ are very intuitive and rather easy to accept after some hand-waving. Writing complete proofs is often more complicated than one would expect. I hope I have succeeded in this task, which was sometimes much more difficult than I expected.

\section{Some open questions}
\label{secopq}
\begin{enumerate}
\item A Vassiliev invariant is {\em odd\/} if it distinguishes some knot from the same knot with the opposite orientation. Are there odd Vassiliev invariants?
\item More generally, do Vassiliev invariants distinguish knots in $S^3$?
In \cite{Kupinvertib}, Greg Kuperberg proved that if they distinguish unoriented knots in $S^3$, then there exist odd Vassiliev invariants.
\item According to a theorem of Dror Bar-Natan and Ruth Lawrence \cite{barnatanlawr}, the LMO invariant fails to distinguish rational homology spheres with isomorphic $H_1$. So, according to a Moussard theorem \cite{moussardAGT}, rational finite type invariants fail to distinguish $\QQ$-spheres. Do finite type invariants distinguish $\ZZ$-spheres?
\item Compute the anomalies $\alpha$ and $\ansothree$ of Sections~\ref{secanomalpha} and \ref{secansothree}. For links in $\RR^3$, I expressed the invariant $\Zinvuf$  as a function of the Kontsevich integral $Z^K$, described in \cite{barnatan,ohtsukibookquantum,chmutovmostov}, and of the Bott--Taubes anomaly $\alpha$, in \cite{lesunikon}. The computation of $\alpha$ would finish clarifying the relationship between $\Zinvuf$ and $Z^K$ for links in $\RR^3$.
See Note~\ref{rkalphafive}.
\item Find surgery formulae for $\Zinvuf$. Do the surgery formulae that define $Z_{LMO}$
from $Z^K$ define $\Zinvuf$ from its restriction to links in $\RR^3$?
\item Compare $\Zinvuf$ with the LMO invariant $Z_{LMO}$ of Thang L\^e, Jun Murakami, and Tomotada Ohtsuki described in \cite[Chapter 10]{ohtsukibookquantum}.
\item Find relationships between $\Zinvuf$ or other finite type invariants and Heegaard Floer homologies. Recall the propagators associated to Heegaard diagrams of \cite{lesHC} from Subsection~\ref{subMorse}.
\item 
Andrew Kricker defined a lift $\tilde{Z}^K$ of the Kontsevich integral $Z^K$ (or the LMO invariant) for null-homologous knots in $\QQ$-spheres \cite{kricker,garkri}.
The Kricker lift is valued in a space $\tilde{A}$ of trivalent diagrams whose edges are decorated by rational functions whose denominators divide the Alexander polynomial.
Compare the Kricker lift $\tilde{Z}^K$ with the equivariant configuration space invariant $\tilde{Z}^c$ of \cite{lesbonn, lesuniveq} valued in the same diagram space  $\tilde{A}$. 
\item Does one obtain $\Zinvuf$ from  $\tilde{Z}^c$ in the same way as one obtains  $Z^K$ from $\tilde{Z}^K$? 
\item Study extensions of $\Zinvuf$ to manifolds with boundary. Dorin Cheptea, Kazuo Habiro, and  Gw\'{e}na\"{e}l Massuyeau introduced a functorial extension of the LMO invariant in \cite{ChepteaHM}. See also \cite{HabiroM21Quantum}.
\end{enumerate}

\chapter{More on manifolds and on the linking number}
\label{chaplk}

The first section of this chapter specifies some basic notions of differential topology, quickly and sometimes vaguely introduced in Subsection~\ref{subdegree}.
It also contains some additional notation and conventions. 
The second section completes our discussion about the linking number of Sections~\ref{sublkGauss} and \ref{sublkGausstwo}.

\section{More background material on manifolds}

\subsection{Manifolds without boundary}
\label{submanifolds}

This section presents a quick review of the notions of manifold and tangent bundle. We refer the reader to \cite[Chapter 1]{hirsch}  by Morris Hirsch for a clean and complete introduction.

A {\em topological $n$-dimensional manifold $M$ without boundary\/} is a Hausdorff topological space that is a union of open subsets $U_i$ labeled in a countable set $I$ 
($i \in I$), where every $U_i$ is identified with an open subset $V_i$ of $\RR^n$ by a homeomorphism $\phi_i\colon U_i \rightarrow V_i$, called a {\em chart.\/} 
Such a collection $(\phi_i\colon U_i \rightarrow V_i)_{i\in I}$ of charts, for which $\cup_{i\in I}U_i=M$, is called an
\emph{atlas} of $M$.
We consider manifolds up to homeomorphism. So homeomorphic manifolds are considered identical.

For $r=0,\dots,\infty$, the topological manifold $M$ {\em  has a $C^r$-structure\/} (induced by the atlas $(\phi_i)_{i\in I}$) or {\em is a $C^r$-manifold\/}, if, for each
pair $\{i,j\} \subset I$, the \emph{transition map} $\phi_j \circ \phi_i^{-1}$ defined on $\phi_i(U_i \cap U_j)$ is a $C^r$-diffeomorphism onto its image. The notion
of $C^s$-maps, $s \leq r$, from such a manifold to another one can be induced naturally from the known case for which the manifolds are open subsets of some $\RR^n$, thanks to the local identifications provided by the charts. Manifolds of class $C^r$ are considered up to $C^r$-diffeomorphism. They are called \emph{$C^r$-manifolds}. \emph{Smooth manifolds} are $C^{\infty}$-manifolds. 

A $C^r$ \emph{embedding} from a $C^r$ manifold $A$ into a $C^r$ manifold $M$ is an injective $C^r$ map $j\colon A \hookrightarrow M$ such that, for any point $a$ of $A$, there exist \begin{itemize}
                                                                                                                            \item a $C^r$ diffeomorphism $\phi$ from an open neighborhood $U$ of $j(a)$ in $M$ to an open subset $V$
of $\RR^{n}$ and
                                                                                                                                                        \item  an open neighborhood $U_A$ of $a$ in $A$ such that the restriction $j\vert_{U_A}$ of $j$ to $U_A$ is a $C^r$-diffeomorphism onto its image, which may be written as $j(U_A)=j(A)\cap U=\phi^{-1}\left(V\cap (\RR^d\times \{(0,\dots,0)\})\right)$. 
                                                                                                                                                                            \end{itemize}

A \emph{submanifold} of a manifold $M$ is the image of such an embedding into $M$. 

The \emph{tangent space} $T_xA$ to a $C^r$ submanifold $A$ of $\RR^n$ at a point $x$ of $A$, for $r\geq 1$ is the vector space of all tangent vectors to a curve (i.e., a $1$-dimensional submanifold) of $A$ at $x$. A well-known theorem \cite[Theorem 3.4, Chapter 1]{hirsch} asserts that any compact $C^r$-manifold, for $r\geq 1$ may be embedded in some $\RR^d$, and thus viewed as a submanifold of $\RR^d$.
The \emph{tangent bundle} $TA$ to $A$ is the union over the elements $x$ of $A$ of the $T_xA$. Its \emph{bundle projection} $p \colon TA \rightarrow A$ maps an element $v$ of $T_xA$ to $x$.
The tangent bundle to $\RR^n$ is canonically isomorphic to $\RR^n \times \RR^n$. A $C^r$ diffeomorphism between two open sets of $\RR^n$, together with its (first order) derivatives induces a canonical $C^{r-1}$-diffeomorphism between their tangent bundles. The notion of tangent bundle of any $C^r$ $n$-manifold, for $r\geq 1$, is naturally induced from the local identifications provided by the charts.
A $C^r$ map $f$ from a $C^r$-manifold $M$ to another one $N$ has a well-defined \emph{tangent map}, which is a map $Tf \colon TM \rightarrow TN$ restricting as a linear map $T_xf \colon T_xM \rightarrow T_{f(x)}N$ for any $x$ of $M$.

\subsection{More on low-dimensional manifolds}
\label{submorelowdif}

We now review classical results, which ensure that for $n=1$, $2$ or $3$, any topological $n$-manifold may be equipped with a unique smooth (i.e., $C^{\infty}$) structure (up to diffeomorphism).

A topological manifold $M$ as in the previous section {\em  has a piecewise linear (or PL) structure\/}  (induced by the atlas $(\phi_i)_{i\in I}$) or {\em is a $PL$-manifold\/}, if, for each
pair $\{i,j\} \subset I$, the transition map $\phi_j \circ \phi_i^{-1}$ is a piecewise linear homeomorphism onto its image.
$PL$-manifolds are considered up to $PL$-homeomorphism.

\emph{An $n$-dimensional simplex} is the convex hull of $(n+1)$ points that are not contained in an affine subspace of dimension $(n-1)$ in some $\RR^k$, with $k\geq n$. For example, a $1$-dimensional simplex is a closed interval, a $2$-dimensional simplex is a solid triangle, and a $3$-dimensional simplex is a solid tetrahedron.
A topological space $X$ has a {\em triangulation\/}, if it is a locally finite union of $k$-simplices (closed in $X$), which are the simplices of the triangulation, such that
\begin{itemize}
 \item the simplices are embedded in $X$,
 \item every face of a simplex of the triangulation is a simplex of the triangulation,
 \item when two simplices of the triangulation are not disjoint, their intersection is a simplex of the triangulation.\footnote{A collection of sets in $X$ is \emph{locally finite} if each point of $X$ has a neighborhood in $X$ that intersects finitely many sets of the collection.}
\end{itemize}

PL manifolds always have such triangulations.

When $n \leq 3$, the above notion of PL-manifold coincides with the notions of smooth and topological manifold, according to the following theorem. This is no longer true when $n>3$. See 
\cite{kui} by Nicolaas Kuiper.

\begin{theorem} 
\label{thmstructhree}
Let $n$ be a natural integer such that $n \leq 3$.
Any topological $n$-manifold has a unique PL structure (up to PL homeomorphism).
For any $r\in (\NN \setminus \{0\}) \cup \{\infty\}$, any topological $n$-manifold has a unique $C^{r}$-structure (up to $C^r$ diffeomorphism).
\end{theorem}

This statement contains several theorems (see \cite{kui}): 
The fact that any $C^1$-manifold in any dimension has a unique $C^{r}$-structure (up to $C^r$-diffeomorphism), for any $r\in (\NN \setminus \{0,1\}) \cup \{\infty\}$, follows from work of Hassler Whitney in 1936 \cite{whi}.
In 1934, Stewart Cairns \cite{cai1} equipped any $C^{1}$-manifold with a PL-structure. He proved that any PL manifold of dimension $3$ arises as the image of a $C^{1}$-manifold under such a process \cite[Theorem III]{cai2} in 1940. (This was already known in dimension less than $3$.) Edwin Moise \cite{moise} proved that any topological $3$-manifold has a unique PL structure in 1952. (This was also already known in dimension less than $3$.)
In dimension $3$, James Munkres \cite[Theorem 6.3]{munk} and Henry Whitehead \cite{white} completed this scheme in 1960 by their independent proofs of the uniqueness of a $C^1$-structure for any topological $3$-manifold.

\subsection{Connected sum}

Let $M_1$ and $M_2$ be two smooth closed manifolds of dimension $n$. The \indexT{connected sum} $M_1 \# M_2$ of $M_1$ and $M_2$ is defined as follows. For $i \in \{1,2\}$, let $\phi_i \colon 2\mathring{B}^n \hookrightarrow M_i$ be a smooth embedding of
the open ball $2\mathring{B}^n$ of radius $2$ (centered at the origin) of the Euclidean vector space $\RR^n$ into $M_i$, such that $\phi_1$ is orientation-preserving and $\phi_2$ is orientation-reversing. The elements of $(2\mathring{B}^n\setminus \{0\})$ may be written as $\lambda x$ for a unique pair $(\lambda,x) \in\left]0,2\right[ \times S^{n-1}$, where $S^{n-1}$ is the unit sphere of $\RR^n$. 
Let $h\colon \phi_1(2\mathring{B}^n\setminus \{0\})\to \phi_2(2\mathring{B}^n\setminus \{0\})$ be the diffeomorphism
such that $h\left(\phi_1(  \lambda x)\right)= \phi_2\left((2-\lambda)x\right)  $ 
for any $(\lambda,x) \in\left]0,2\right[ \times S^{n-1}$.

Then \begin{equation*}M_1 \# M_2=\bigl(M_1 \setminus \{\phi_1(0)\}\bigr) \cup_{h} \bigl(M_2 \setminus\{ \phi_2(0)\}\bigr)\end{equation*}
is the quotient space of $\left(M_1 \setminus \{\phi_1(0)\}\right) \sqcup \left(M_2 \setminus \{\phi_2(0)\}\right)$, in which an element of $\phi_1(2\mathring{B}^n\setminus \{0\})$ is identified with its image under $h$.
As a topological manifold, $M_1 \# M_2$ can be written as
\begin{equation*}M_1 \# M_2=\left(M_1 \setminus \phi_1(\mathring{B}^n)\right) \cup_{\phi_1(S^{n-1})\stackrel{h}{\sim} \phi_2(S^{n-1})} \left(M_2 \setminus \phi_2( \mathring{B}^n)\right).\end{equation*}

\subsection{Manifolds with boundary and ridges}
\label{secbackgroundmfdsbry}

A {\em topological $n$-dimensional manifold $M$ with possible boundary\/} is a Hausdorff topological space that is a union of open subsets $U_i$ labeled in a set $I$, 
($i \in I$), where every $U_i$ is identified with an open subset $V_i$ of $\left]-\infty,0\right]^k \times \RR^{n-k}$ by a chart $\phi_i\colon U_i \rightarrow V_i$. The {\em boundary\/} of $\left]-\infty,0\right]^k \times \RR^{n-k}$ consists of the points $(x_1,\dots,x_n)$ of $\left]-\infty,0\right]^k \times \RR^{n-k}$ for which there exists $i \leq k$ such that $x_i=0$.
 The {\em boundary\/} of $M$ consists of the points mapped to the boundary of $\left]-\infty,0\right]^k \times \RR^{n-k}$ by a chart.

A map from an open subset $O$ of $\left]-\infty,0\right]^k \times \RR^{n-k}$ to an open subset of $\left]-\infty,0\right]^{k^{\prime}} \times \RR^{n^{\prime}-k^{\prime}}$ is smooth
(resp. $C^r$) at a point $x \in O$ if it extends as a smooth (resp. $C^r$) map from an open neighborhood of $x$ in $\RR^{n}$ to $\RR^{n^{\prime}}$.

The topological manifold $M$ is a \emph{smooth manifold with ridges}, if, for each
pair $\{i,j\} \subset I$, the map $\phi_j \circ \phi_i^{-1}$ defined on $\phi_i(U_i \cap U_j)$ is a smooth diffeomorphism onto its image.\footnote{Manifolds with ridges are often called manifolds with corners.} When $k\leq 1$ for any $i$, we simply say that $M$ is
a \emph{smooth manifold with boundary}.
The \indexT{codimension $j$ boundary} of such a manifold $M$, which is denoted by $\partial_j(M)$\index[N]{Dzzel@$\partial_j$ codimension $j$ boundary},
consists of the points that are mapped to points $(x_1,\dots,x_n)$ of $\left]-\infty,0\right]^k \times \RR^{n-k}$ for which there are at least $j$ indices $i \leq k$ such that $x_i=0$. It is a closed subset of $M$. We have $\partial M= \partial_1(M)$.

The \emph{codimension $j$ faces} of such a smooth manifold $M$ with ridges are the connected components of $\partial_{j}(M) \setminus \partial_{j+1}(M)$. They are smooth manifolds of dimension $(n-j)$.
The \emph{interior} of $M$ is $M \setminus \partial M$.

For $j\geq 2$, the codimension $j$ faces are called \emph{ridges} of $M$.

\subsection{Algebraic intersections}
\label{subalgint}

We start this subsection with more orientation conventions and notation. Again, unless otherwise mentioned, manifolds are smooth, compact, and oriented.
The {\em normal bundle\/} to a submanifold $\Aman$ in a manifold $M$, at a point $x$, is the quotient $T_xM/T_x\Aman$ of tangent bundles at $x$. We denote it by $N_x\Aman$ or $N_x(\Aman)$. We orient $N_x\Aman$ so that (a lift of an oriented basis of) $N_x\Aman$ followed by (an oriented basis of) $T_x\Aman$ induce the orientation of $T_xM$. The orientation of $N_x(\Aman)$ is a \indexT{coorientation} of $\Aman$ at $x$. The regular preimage of a submanifold under a map $f$ is oriented so that $f$ preserves the coorientations.

Two submanifolds $\Aman$ and $B$ in a manifold $M$ are \emph{transverse} if we have $T_xM=T_x\Aman+T_xB$ at each $x \in \Aman \cap B$. As proved in \cite[Chapter 3 (Theorem 2.4 in particular)]{hirsch}, transversality is a generic condition.
The intersection $\Aman\cap B$ of two transverse submanifolds $\Aman$ and $B$ in a manifold $M$ is a manifold. We orient $\Aman\cap B$ so that the normal bundle to $\Aman\cap B$ is $(N(\Aman) \oplus N(B))$, fiberwise. 
In order to give a meaning to the sum $(N_x(\Aman) \oplus N_x(B))$ at $x \in A \cap B$, pick a Riemannian metric on $M$. Such a metric identifies $N_x(\Aman)$ with $T_x(\Aman)^{\perp}$, $N_x(B)$ with $T_x(B)^{\perp}$, and $N_x(\Aman\cap B)$ with $T_x(\Aman\cap B)^{\perp}=T_x(\Aman)^{\perp} \oplus T_x(B)^{\perp}$. Since the space of Riemannian metrics on $M$ is convex, and therefore connected, the induced orientation of $T_x(\Aman\cap B)$ does not depend on the choice of a Riemannian metric.

Let  $A$, $B$, $C$ be three pairwise transverse submanifolds in a manifold $M$ such that $A\cap B$ is transverse to $C$. The oriented intersection $(A\cap B )\cap C$ is a well-defined manifold.
Our assumptions imply that at any $x \in A\cap B \cap C$,
the sum $(T_xA)^{\perp} + (T_xB)^{\perp} +(T_xC)^{\perp}$ is a direct sum $(T_xA)^{\perp} \oplus (T_xB)^{\perp} \oplus (T_xC)^{\perp}$ for any Riemannian metric on $M$. So $A$ is also transverse to $B\cap C$, and $(A\cap B )\cap C=A\cap (B \cap C)$.
Thus, the intersection of transverse, oriented submanifolds is a well-defined associative operation, where \emph{transverse submanifolds} are manifolds such that the elementary pairwise intermediate possible intersections are well-defined, as above. This intersection is also commutative when the codimensions of the submanifolds are even.

Recall from Subsection~\ref{subsecsiminv} that,
for two transverse submanifolds $\Aman$ and $B$ of complementary dimensions in a manifold $M$,
the sign $\pm 1$ of a point $x \in \Aman \cap B$ is $+1$ if and only if $T_xM=T_x\Aman \oplus T_xB$ as oriented vector spaces. This is equivalent to the condition that the orientation of the normal bundle to $x \in \Aman \cap B$ coincides with the orientation of the ambient space $M$, that is that $T_xM=N_x\Aman \oplus N_xB$ (as oriented vector spaces again---exercise). 
If $\Aman$ and $B$ are of complementary dimensions in $M$, and if $\Aman \cap B$ if finite, then the \emph{algebraic intersection} $\langle \Aman, B \rangle_{\!M\,}$ of $\Aman$ and $B$ is the sum of the signs of the intersection points.

In a manifold $M$, a $k$-dimensional \emph{chain (resp. a $k$-dimensional rational chain)}\index[T]{chain} is a finite combination with coefficients in $\ZZ$ (resp. in $\QQ$) of (smooth, compact, oriented) $k$-dimensional submanifolds $C$ of $M$ with boundary and ridges, up to the identification of $(-1)C$ with $(-C)$ and other natural identifications (e.g., a $k$-manifold $A\cup B$ such that $A \cap B$ is a $(k-1)$-submanifold of $A \cup B$ is identified with the chain $A$ + $B$). The boundary $\partial$ of chains is the linear map that maps a submanifold to its oriented boundary (with respect to the usual outward normal first convention). This boundary is the sum of the closures of the codimension-one faces when there are ridges.  The canonical orientation of a point is the sign $+1$. So $\partial \left[0,1\right]= \{1\}-\{0\}$. A $k$-dimensional chain whose boundary vanishes is a $k$-dimensional \indexT{cycle}, or a \emph{$k$-cycle} for short.

If $A_1, \dots, A_k$ are $k$ transverse compact submanifolds of $M$ whose codimension sum is the dimension of $M$, then their \indexT{algebraic intersection} is defined to be
$\langle A_1, \dots, A_k \rangle_{\!M\,}=\langle \cap_{i=1}^{k-1}A_i, A_k \rangle_{\!M\,}$.
If $M$ is a connected manifold, which contains a point $x$, then the class of a $0$-cycle in $H_0(M;\QQ)=\QQ\left[x\right]=\QQ$ is a well-defined number. The algebraic intersection $\langle A_1, \dots, A_k \rangle_{\!M\,}$ can be equivalently defined to be the homology class of the (oriented) intersection $\cap_{i=1}^{k}A_i$. This algebraic intersection extends multilinearly to rational chains.

We will use the following lemma in Subsection~\ref{subgendeflinking}.

\begin{lemma}
\label{lemsignint}
 Let $\Aman$ and $B$ be two transverse submanifolds of a $d$-dimensional manifold $M$ with disjoint boundaries. Let $\beta$ denote the dimension of $B$. 
Then \begin{equation*}\partial (\Aman\cap B)=(-1)^{d-\beta}\partial \Aman \cap B + \Aman \cap \partial B.\end{equation*}
\end{lemma}
\bp Note that $\partial (\Aman\cap B) \subset \partial \Aman \cup \partial B$. At a point $a \in \partial \Aman$, $T_aM$ is oriented by $(N_a\Aman,o,T_a\partial \Aman)$, where $o$ is the outward normal to $\Aman$. If $a \in \partial \Aman \cap B$, then $o$ is also an outward normal for $\Aman\cap B$, and $\partial (\Aman \cap B)$ is cooriented by $(N_a\Aman,N_aB,o)$, while $\partial \Aman \cap B$ is cooriented by $(N_a\Aman,o,N_aB)$.
At a point $b \in \Aman \cap \partial B$, both $\partial (\Aman \cap B)$ and $\Aman \cap \partial  B$ are cooriented by $(N_b\Aman,N_bB,o)$.
\eop

\subsection{More on the degree}
\label{submoredeg}

Here, we make the notion of \emph{walls} in Subsection~\ref{subsecwtwo} more precise, by stating a lemma on the general behavior of the degree. We prove it with the Morse--Sard theorem~\ref{thmMorseSard}. The formula of the lemma could also be justified with Stokes' theorem. 

\begin{lemma}
\label{lemdeggen}
Let $n \in \NN$. Let $M$ be a compact (oriented) $n$-manifold with possible boundary. Let $N$ be a connected (oriented) $n$-manifold.
Let $f\colon M \to N$ be a smooth map. Let  $a$ and $b$
be two distinct regular values of $f$ in the interior $\Int(N)$ of $N$.

Then there exists an embedding $\gamma \colon \left[0,1\right] \to \Int(N)$ such that $\gamma(0)=a$, $\gamma(1)=b$, and, for any $x \in \partial M \cap f^{-1}\left(\gamma\left(\left[0,1\right]\right)\right)$,
\begin{itemize}
\item $x$ is in an open face of $\partial M$ (of codimension one in $M$) and
\item $T_{f(x)} N$ equals $\eta(x) T_{f(x)} \gamma \oplus T_xf(T_x \partial M)$ as an oriented vector space for some $\eta(x)=\pm 1$.
\end{itemize}

For any such embedding $\gamma$, we have
\begin{equation*}\deg_af-\deg_bf = \langle \gamma,f(\partial M) \rangle_{\!N\,} =\sum_{x \in f^{-1}(\gamma) \cap \partial M}\eta(x).\end{equation*}
\end{lemma}
\bfig
\centering
\begin{tikzpicture} \useasboundingbox (-5,-1) rectangle (5,1);
\draw (-4,-1) rectangle (4,1) ;
\fill [lightgray] (0,0) ellipse (2 and .8); 
\draw (0,0) ellipse (2 and .8); 
\fill [white, rounded corners] (-.8,.2) rectangle (0,.6);
\fill (.2,0) circle (.05) (3.4,0) circle (.05) (2,0) circle (.05);
\fill [white, rounded corners] (.05,-.05) rectangle (.35,-.35);
\draw [->] (2.7,0) -- (3.4,0) (.2,0) -- (2.7,0);
\draw (2.75,-.05) node[above]{\scriptsize $\gamma$} (3.5,0) node[below]{\scriptsize $b$} (-.4,.4) node{\scriptsize $f(M)$} (-3.7,-.8) node{\scriptsize $N$} (.2,-.2) node{\scriptsize $a$} (2.25,.05) node[below]{\scriptsize $f(x)$};
\end{tikzpicture}

\caption{Lemma~\ref{lemdeggen} for an embedding $f$ from $M$ to a rectangle $N$}\label{figlemdeggen}
\end{figure}

\bp Let $B^{n-1}$ be the unit ball of $\RR^{n-1}$.
Since $N$ is connected, there exists an (orientation-preserving) embedding 
$\Psi \colon B^{n-1} \times \left[-1,2\right] \to \Int(N)$ such that $\Psi(0,0) = a$, $\Psi(0,1) = b$, and $\Psi(\varepsilon B^{n-1} \times \left]-\varepsilon,\varepsilon\right[)$ and $\Psi(\varepsilon B^{n-1} \times \left]1-\varepsilon,1+\varepsilon\right[)$ consist of regular values of $f$ for some $\varepsilon\in\left]0,\frac1{100}\right[$.

If all the elements of $\Psi(\{0\} \times \left[0,1\right])$ are regular values of $f$, then the preimage $f^{-1}\left(\Psi( \{0\} \times\left[0,1\right])\right)$ is a disjoint union of intervals between a point of $f^{-1}(a)$ and a point of $f^{-1}(b)$, along which the sign of $\det(\tang_xf)$ is constant, and $\deg_af$ is equal to $\deg_bf$. This proves that the degree is constant on 
any ball of regular values. Let us return to the general case and try to make a similar argument work.
Set $N_{\varepsilon}=\varepsilon \mathring{B}^{n-1} \times \left[0,1\right]$, $M_{\varepsilon} =f^{-1}(\Psi(N_{\varepsilon}))$, and $f_{\varepsilon}= \Psi^{-1}\circ f\vert_{M_{\varepsilon}}$. For any $c \in \Psi(\mathring{N}_{\varepsilon})$, $\deg_cf=\deg_{\Psi^{-1}(c)}(f_{\varepsilon})$.

Let $p_B$ and $p_I$ respectively denote the natural projections of $\varepsilon \mathring{B}^{n-1} \times \left[0,1\right]$ to its factors $\varepsilon \mathring{B}^{n-1}$and $\left[0,1\right]$ .
The Morse--Sard theorem~\ref{thmMorseSard} guarantees the existence of a regular value $y$ of $p_B \circ f_{\varepsilon}$ in $\varepsilon \mathring{B}^{n-1}$. Let $\gamma_y\colon \left[0,1\right] \to N_{\varepsilon}$ map $t$ to $(y,t)$. Then $\gamma_y$ is cooriented by $p_B$ (via  $T_w p_B \colon T_w N_{\varepsilon} \to \RR^{n-1}$).
So its preimage
\begin{equation*}f_{\varepsilon}^{-1}(\gamma_y)=f_{\varepsilon}^{-1}\bigl(\gamma_y([0,1])\bigr)= (p_B \circ f_{\varepsilon})^{-1}(y)\end{equation*}
is a submanifold of dimension $1$ of $M$, cooriented by $p_B \circ f_{\varepsilon}$.
Let us show that the \emph{oriented} boundary of $f_{\varepsilon}^{-1}(\gamma_y)$ is 
\begin{equation*}\partial f_{\varepsilon}^{-1}(\gamma_y)=-f_{\varepsilon}^{-1}\bigl((0,y)\bigr)+f_{\varepsilon}^{-1}\bigl((1,y)\bigr)  +
\sum_{x \in f_{\varepsilon}^{-1}(\gamma_y) \cap \partial M}\eta(x)x.\end{equation*}
For a \emph{regular point} $x$ of $f_{\varepsilon}$ (i.e., a point such that $f_{\varepsilon}$ is a local diffeomorphism near $x$), let $\delta(x)=\pm 1$ denote the sign of the Jacobian determinant of $f_{\varepsilon}$ at $x$.
Near such a regular point $x$ in $(p_B \circ f_{\varepsilon})^{-1}(y)$, 
 the curve $\delta(x)f_{\varepsilon}^{-1}(\gamma_y)$ is oriented by $p_I \circ f_{\varepsilon}$. This shows that the above signs before $f_{\varepsilon}^{-1}((0,y))$ and $f_{\varepsilon}^{-1}((1,y))$ are correct.
 
Let us now check the sign associated to 
$x \in f_{\varepsilon}^{-1}(\gamma_y) \cap \partial M$. Let $N_{x,\partial M}$ denote the outward normal to $M$ at $x$. So we have $T_xM=\RR N_{x,\partial M} \oplus T_x\partial M$ and
\begin{equation*}T_{f_{\varepsilon}(x)}N_{\varepsilon}=\eta(x)T_{f_{\varepsilon}(x)} \gamma_y \oplus T_xf_{\varepsilon}(T_x\partial M)=\delta(x)\bigl(T_xf_{\varepsilon} (\RR N_{x,\partial M}) \oplus T_xf_{\varepsilon}(T_x\partial M)\bigr).\end{equation*} 
Therefore, the outward normal $N_{x,\partial M}$ to $\partial M$ is oriented by $\eta(x) \delta(x)p_I \circ f_{\varepsilon}$, and the sign of the scalar product of $T_x\left(f_{\varepsilon}^{-1}(\gamma_y)\right)$ and $N_{x,\partial M}$ is $\eta(x)$.

Writing that the homology class of $\partial f_{\varepsilon}^{-1}(\gamma_y)$ is zero in  $H_0(M)$ shows
\begin{equation*}\deg_{(0,y)}f_{\varepsilon}-\deg_{(1,y)}f_{\varepsilon} = \sum_{x \in  f_{\varepsilon}^{-1}(\gamma_y) \cap \partial M}\eta(x).\end{equation*} 

This proves that the lemma holds for $\Psi(a_y=\gamma_y(0)) $, $\Psi(b_y=\gamma_y(1))$, and $\Psi\circ \gamma_y$.
Since $\Psi(a_y)$ and $a$ are in a common open ball $B_a$ of regular values of $\funcf$  (in $\Int(N)\setminus f(\partial M)$), and since $\Psi(b_y)$ and $b$ are in another common ball of regular values of $\funcf$,
we may modify $\Psi \circ \gamma_y$ to a smooth path $\gamma$ from $a$ to $b$ with all the properties of the statement. 

Let us finish by proving that for any smooth embedded path $\gamma$ from $a$ to $b$ transverse to $f(\partial M)$, we have $\deg_a\funcf-\deg_b\funcf = \langle \gamma,\funcf(\partial M) \rangle_{\!N\,}.$ For any such path $\gamma$, we can construct a neighborhood embedding $\Psi$ of $\gamma$ as above, such that $\Psi (0,t)=\gamma(t)$ for all $t \in \left[0,1\right]$, $\Psi(y\times\left]-\varepsilon,1+\varepsilon\right[)$ is transverse to $f(\partial M)$ for any $y \in \varepsilon \mathring{B}^{n-1}$, and $\langle \gamma_y,f(\partial M)\rangle$ does not depend on $y \in \varepsilon \mathring{B}^{n-1}$. So the previous study allows us to conclude.
\eop

\section{On the linking number, again}
\label{seclinkingnumber}

\subsection{A general definition of the linking number}
\label{subgendeflinking}

\begin{lemma}
\label{lemdefgenlk}
Let $J$ and $K$ be two rationally null-homologous disjoint cycles of respective dimensions $j$ and $k$ in a $d$-manifold $M$, where $d=j+k+1$.
There exists a rational $(j+1)$-chain $\Sigma_J$ bounded by $J$ transverse to $K$ and a rational $(k+1)$-chain $\Sigma_K$ bounded by $K$ transverse to $J$. For any two such rational chains $\Sigma_J$ and $\Sigma_K$, we have $\langle J,\Sigma_K\rangle_{\!M\,}=(-1)^{j+1}\langle \Sigma_J,K\rangle_{\!M\,} $.
In particular, $\langle J,\Sigma_K\rangle_{\!M\,}$ is a topological invariant of $(J,K)$. It is denoted by $lk(J,K)$ and called the \indexT{linking number} of $J$ and $K$. We have
\begin{equation*}lk(J,K)=(-1)^{(j+1)(k+1)}lk(K,J).\end{equation*}
\end{lemma}
\bp
Since $K$ is rationally null-homologous, $K$ bounds a rational $(k+1)$-chain $\Sigma_K$. Without loss of generality, $\Sigma_K$ is assumed to be transverse to $\Sigma_J$, so $\Sigma_J \cap \Sigma_K$ is a rational $1$-chain (which is a rational combination of circles and intervals). According to Lemma~\ref{lemsignint}, we have
\begin{equation*}\partial (\Sigma_J \cap \Sigma_K) = (-1)^{d+k+1}J \cap \Sigma_K +\Sigma_J \cap K.\end{equation*}
Furthermore, the sum of the coefficients of the points in the left-hand side must be zero since this sum vanishes for the boundary of an interval. This proves $\langle J,\Sigma_K\rangle_{\!M\,}=(-1)^{d+k} \langle \Sigma_J,K\rangle_{\!M\,}$. Therefore, this rational number is independent of the chosen $\Sigma_J$ and $\Sigma_K$. Since we have
\begin{equation*}(-1)^{d+k}\langle \Sigma_J,K\rangle_{\!M\,}=(-1)^{j+1}(-1)^{k(j+1)}\langle K, \Sigma_J\rangle_{\!M\,},\end{equation*} we get $lk(J,K)=(-1)^{(j+1)(k+1)}lk(K,J).$
\eop

\begin{remark}
\label{rksignlinking}
 Our sign convention for the linking number differs from that in \cite[Section 77, page 288]{SeifertThrel}, where the linking number of cycles $J$ and $K$ as in the lemma is defined as $\langle \Sigma_J,K\rangle_{\!M\,}$, instead. The reason for our sign convention is justified in Remark~\ref{rksignlinkingjust}.
\end{remark}

In particular, the {\em linking number\/} of two rationally null-homologous disjoint links $J$ and $K$ in a $3$-manifold $M$ is the algebraic intersection of a rational chain bounded by one of the links and the other.

For $\KK=\ZZ$ or $\QQ$, any knot is rationally null-homologous
in a $\KK$-sphere or in a $\KK$-ball (defined in Section~\ref{secabs}). So the linking number of two disjoint knots always makes sense in such a $3$-manifold.

A \emph{meridian\index[T]{meridian of a knot}} $m_K$ of a knot $K$ is the (oriented) boundary of a disk that intersects $K$ once with a positive sign. See Figure~\ref{figmeridian}. Note $lk(K,m_K)=1$.

\bfig
\centering
\begin{tikzpicture}
\useasboundingbox (-1,-.3) rectangle (1,1.3);
\draw[->] (.5,0) .. controls (.8,0) and (.9,.3) .. (.9,.6);
\draw [draw=gray] (1.05,.8)  .. controls (1.05,.85) and (1,.9) .. (.9,.9) .. controls (.8,.9) and (.75,.85) .. (.75,.8); 
\draw [draw=white,double=black,very thick] (.9,.6) .. controls (.9,.9) and (.7,1.2) .. (.5,1.2) .. controls (.3,1.2) and (-.2,1.15) .. (-.2,1);
\draw [draw=white,double=gray,very thick] (.75,.8)  .. controls (.75,.75) and (.8,.7) .. (.9,.7) .. controls (1,.7) and (1.05,.7) .. (1.05,.8);
\draw [->,draw=gray] (.75,.8)  .. controls (.75,.75) and (.8,.7) .. (.9,.7) .. controls (1,.7) and (1.05,.7) .. (1.05,.8) node[right]{\tiny $m_K$}; 
\draw[->] (.4,.5) .. controls (.4,.3) and (.2,0) .. (0,0) node[below]{\tiny $K$};
\draw [draw=white,double=black,very thick] (.2,1) .. controls (.2,1.15) and (-.3,1.2) .. (-.5,1.2) .. controls (-.8,1.2) and (-.9,.9) .. (-.9,.6) .. controls (-.9,.3) and (-.8,0) .. (-.5,0) .. controls (-.3,0) and (-.2,.4) .. (0,.4) .. controls (.2,.4) and (.3,0) .. (.5,0);
\draw[draw=white,double=black,very thick] (0,0) .. controls (-.2,0) and (-.4,.3) .. (-.4,.5) .. controls (-.4,.7) and (.2,.85) .. (.2,1);
\draw[draw=white,double=black,very thick] (-.2,1) .. controls (-.2,.85) and (.4,.7) .. (.4,.5);
\end{tikzpicture}

\caption{A meridian $m_K$ of a knot $K$}
\label{figmeridian}
\end{figure}

\begin{lemma}
\label{lemdeflkhom} Let $\rats$ be a $\QQ$-sphere or a $\QQ$-ball.
 Let $K$ be a knot in $\rats$. Then $H_1(\rats\setminus K;\QQ)$ is equal to $\QQ\left[m_K\right]$. Let $J$ be a knot of $\rats$ disjoint from $K$. Then we have $\left[J\right]=lk(J,K)\left[m_K\right]$ in $H_1(\rats\setminus K;\QQ)$. This equality provides an alternative definition for the linking number.
\end{lemma}
\bp Exercise. Note that a chain bounded by $J$ transverse to $K$ in $\rats$ provides a rational cobordism between $J$ and a combination of meridians of $K$. \eop

The reader is also invited to check that the Gauss linking number $lk_G$ of Subsection~\ref{sublkGauss} coincides with the above linking number $lk$ for two-component links of $S^3$, as an exercise. This is proved in the following subsection, see Proposition~\ref{propdeflkeq}.

\subsection{Generalizing the Gauss definition of the linking number and identifying the definitions}
\label{subgenGaussdef}

Let $X$ and $Y$ be two topological spaces.
Recall that a \indexT{homotopy} from a continuous map $f$ from $X$ to $Y$ to another such $g$ is a continuous map $H \colon \left[0,1\right] \times X \to Y$ such that for any $x \in X$, $H(0,x)=f(x)$ and $H(1,x)=g(x)$. Two continuous maps $f$ and $g$ from $X$ to $Y$ are said to be \emph{homotopic} if there exists a homotopy from $f$ to $g$.
A continuous map $f$ from $X$ to $Y$ is a \indexT{homotopy equivalence} if there exists a continuous map $g$ from $Y$ to $X$ such that $g \circ f$ is homotopic to the identity map of $X$ and $f \circ g$ is homotopic to the identity map of $Y$.
The topological spaces $X$ and $Y$ are said to be \emph{homotopy equivalent}, or \emph{of the same homotopy type} if there exists a homotopy equivalence from $X$ to $Y$. Appendix~\ref{chaphomology} describes the homology $H_{\ast}$ of topological spaces and continuous maps.

Let $\diag\left((\RR^3)^2\right)$ \index[N]{Dzelta@Diagonals!Dzeltaa@$\diag\left((\RR^3)^2\right)$} denote the diagonal of $(\RR^3)^2$.
\begin{lemma}
\label{lempstt}
The map \begin{equation*}\begin{array}{llll}p_{S^2}\colon &(\RR^3)^2 \setminus \diag\left((\RR^3)^2\right)&\rightarrow & S^2\\&(x,y) &\mapsto &\frac{1}{\lVert y-x \rVert}(y-x) \end{array}\end{equation*}\index[N]{Projections!pStwo@$p_{S^2}$ to $S^2$}
is a homotopy equivalence.
In particular \begin{equation*}H_{i}(p_{S^2}) \colon H_i\Bigl((\RR^3)^2 \setminus \diag\bigl((\RR^3)^2\bigr);\ZZ\Bigr) \rightarrow H_i(S^2;\ZZ)\end{equation*} is an isomorphism for all integer $i$, the space
$(\RR^3)^2 \setminus \diag\left((\RR^3)^2\right)$ is a homology $S^2$, and $\left[S\right]=H_{2}\left(p_{S^2}\right)^{-1}\left[S^2\right]$ is a canonical generator of \begin{equation*}H_2\Bigl((\RR^3)^2 \setminus \diag\bigl((\RR^3)^2\bigr);\ZZ\Bigr)=\ZZ\left[S\right].\end{equation*}
\end{lemma}
\bp The configuration space $\bigl((\RR^3)^2 \setminus \diag\left((\RR^3)^2\right)\bigr)$ is homeomorphic to $\RR^3 \times \left]0,\infty\right[ \times S^2$ via the map \begin{equation*}(x,y) \mapsto \bigl(x,\lVert y-x \rVert, p_{S^2}(x,y)\bigr).\end{equation*}
\eop

As in Subsection~\ref{sublkGauss}, consider a two-component link
$J \sqcup  K\colon S^1 \sqcup S^1 \hookrightarrow \RR^3$. 
This embedding induces an embedding
\begin{equation*}\begin{array}{llll}J \times  K\colon & S^1 \times S^1 &\hookrightarrow &(\RR^3)^2 \setminus \diag\left((\RR^3)^2\right)\\
&(w,z) &\mapsto & \bigl(J(w),K(z)\bigr).\end{array}\end{equation*}
The map $p_{JK}$ of Subsection~\ref{sublkGauss} is the composition $p_{S^2} \circ (J \times  K)$. We have 
\begin{equation*}H_2(p_{JK})\left[S^1 \times S^1\right]=\deg(p_{JK})\left[S^2\right]=lk_G(J,K)\left[S^2\right]\end{equation*} in $H_2(S^2;\ZZ)=\ZZ\left[S^2\right]$. Thus we get
\begin{equation*}\left[(J\times K)(S^1\times S^1)\right]=H_2(J \times  K)\left[S^1 \times S^1\right]=lk_G(J,K)\left[S\right]\end{equation*}
in $H_2((\RR^3)^2 \setminus \diag\left((\RR^3)^2\right);\ZZ)=\ZZ\left[S\right].$
We will see that this definition of $lk_G$ generalizes to links in rational homology spheres.
Then we will prove that our generalized definition coincides with the general definition of linking numbers in this case.

For a manifold $M$, the normal bundle to the diagonal of $M^2$ in $M^2$ is identified with the tangent bundle to $M$, fiberwise, by the map
\begin{equation*}\left[(u,v)\right] \in \frac{(T_xM)^2}{\diag\left((T_xM)^2\right)} \mapsto (v-u) \in T_xM.\end{equation*}

A \emph{parallelization}\index[T]{parallelization!of a $3$-manifold} $\tau$ of an oriented $3$-manifold $M$ is a (smooth) bundle isomorphism $\tau \colon M \times \RR^3 \longrightarrow TM$ that restricts to $ x\times \RR^3$ as an orientation-preserving linear isomorphism from $ x\times \RR^3$ to $T_xM$, for any $x \in M$.
It has long been known that any oriented $3$-manifold is parallelizable (i.e., admits a parallelization). It is proved in Section~\ref{secproofpar}. 
Therefore, a tubular neighborhood of the diagonal $\diag(M^2)$ in $M^2$ is diffeomorphic to $M \times \RR^3$.

\begin{lemma}
\label{lemhomstwo}
Let $\rats$ be a rational homology sphere and let $\infty$ be a point of $\rats$. Set $\crats=\rats \setminus \{\infty\}$. Then $\crats^2 \setminus \diag\left(\crats^2\right)$ has the same rational homology as $S^2$. Let $B$ be a ball in $\crats$ and let $x$ be a point inside $B$, then
the class $\left[S\right]$ of $x \times \partial B$ is a canonical generator of $H_2(\crats^2 \setminus \diag\left(\crats^2\right);\QQ)=\QQ\left[S\right]$.
\end{lemma}
\bp
In this proof, the homology coefficients are in $\QQ$. We refer the reader to Section~\ref{sechomology}.
Since $\crats$ has the homology of a point, the K\"unneth Formula (Theorem~\ref{thmKun}) implies that
 $\crats^2$ has the homology of a point. 
The excision axiom yields
\begin{equation*}\begin{array}{lll}H_{\ast}(\crats^2,\crats^2 \setminus \diag\left(\crats^2\right)) &\cong& H_{\ast}(\crats \times \RR^3,\crats\times (\RR^3 \setminus 0))\\
                  &
                 
 \cong & H_{\ast}( \RR^3, S^2) \cong \left\{\begin{array}{ll} \QQ \;\;&\;\mbox{if} \;\ast =3,\\
0\;&\;\mbox{otherwise.} \end{array} \right. \end{array}\end{equation*}
Using the long exact sequence of the pair $(\crats^2,\crats^2 \setminus \diag\left(\crats^2\right))$, we get
$H_{\ast}(\crats^2 \setminus \diag\left(\crats^2\right))\cong H_{\ast}(S^2)$.
\eop

Define the {\em Gauss linking number\/} of two disjoint links $J$ and $K$ in $\crats$ so that \begin{equation*}\left[(J\times K)(S^1\times S^1)\right]=lk_G(J,K)\left[S\right]\end{equation*} in $H_2(\crats^2 \setminus \diag\left(\crats^2\right);\QQ)$. Note that the two definitions of $lk_G$ coincide when $\crats=\RR^3$.
\begin{proposition}
\label{propdeflkeq} For two disjoint links  $J$ and $K$ in $\crats$, we have
 \begin{equation*}lk_G(J,K)=lk(J,K)\end{equation*}
\end{proposition}
\bp
First recall that $lk(J,K)$ is the algebraic intersection $\langle J,\Sigma_K\rangle_{\!\rats\,}$ of $J$ and a rational chain $\Sigma_K$ bounded by $K$.
Note that the definitions of $lk(J,K)$ and $lk_G(J,K)$ make sense when $J$ and $K$ are disjoint links. If $J$ has several components $J_i$, for $i = 1, \dots, n$, then
$lk_G( \sqcup_{i=1}^nJ_i,K) = \sum_{i=1}^nlk_G(J_i,K)$ and  $lk( \sqcup_{i=1}^nJ_i,K) = \sum_{i=1}^nlk(J_i,K)$. There is no loss of generality in assuming that $J$ is a knot for the proof, which we do.

The chain $\Sigma_K$ provides a rational cobordism $C$ in $\crats \setminus J$ between $K$ and a combination of meridians of $J$. Thus, it provides the rational cobordism $C\times J$ in $\crats^2 \setminus \diag\left(\crats^2\right)$, which allows us to see that
$\left[J\times K\right]=lk(J,K) \left[J \times m_J\right]$ in $H_2(\crats^2 \setminus \diag\left(\crats^2\right);\QQ)$.
Similarly, $\Sigma_J$ provides a rational cobordism between $J$ and a meridian $m_{m_J}$ of $m_J$. So $\left[J \times m_J\right]=\left[m_{m_J} \times m_J \right]$ in $H_2(\crats^2 \setminus \diag\left(\crats^2\right);\QQ)$. This shows
\begin{equation*}lk_G(J,K)=lk(J,K) lk_G(m_{m_J}, m_J).\end{equation*}
Thus, it remains to prove that $lk_G(m_{m_J}, m_J)=1$ for a positive Hopf link $(m_{m_J}, m_J)$, as in Figure~\ref{figHopfWhitehead}, in a standard ball embedded in $\crats$.
Now, there is no loss of generality in assuming that our link is a Hopf link in $\RR^3$, so the equality follows from that for the positive Hopf link in $\RR^3$.
\eop

\begin{remark}
\label{rksignlinkingjust}
Under the assumptions of Lemma~\ref{lemdefgenlk}, the reader can prove as an exercise that if $M$ is connected and if $B$ is a compact ball of $M$ containing a point $x$ in its interior, then
$J\times K$ is homologous to $lk(J,K)(x \times \partial B)$ in $M^2 \setminus \diag\left(M^2\right)$.
In particular, Proposition~\ref{propdeflkeq} generalizes to all pairs $(J,K)$ as in Lemma~\ref{lemdefgenlk} naturally. This justifies our sign convention in Lemma~\ref{lemdefgenlk}.
\end{remark}

\chapter{Propagators}
\label{chapprop}

For a two-component link $(J,K)$ in $\RR^3$, the definition of the linking number $lk(J,K)$ can be rewritten as
\begin{equation*}lk(J,K)=\int_{J \times K}p_{S^2}^{\ast}(\omega)=\bigl\langle J \times K, p_{S^2}^{-1}(Y)\bigr\rangle_{\!\!(\RR^3)^2 \setminus \diag\left((\RR^3)^2\right)\;}\end{equation*} for any $2$-form $\omega$ of $S^2$ such that $\int_{S^2}\omega=1$, and for any regular value $Y$ of $p_{JK}$.
Thus, $lk(J,K)$ is the integral of a $2$-form $p_{S^2}^{\ast}(\omega)$ of $(\RR^3)^2 \setminus \diag\left((\RR^3)^2\right)$ along the $2$-cycle $\left[J \times K\right]$, or it is the intersection of the $2$-cycle $\left[J \times K\right]$ with the $4$-manifold $p_{S^2}^{-1}(Y)$.
In order to adapt these definitions of the linking number to punctured rational homology $3$-spheres $\crats =\rats \setminus \{\infty\}$ and to build other invariants of links and rational homology spheres $\rats$, we compactify $(\crats)^2 \setminus \diag\left((\crats)^2\right)$ to a compact $6$-manifold $C_2(\rats)$, in Section~\ref{secCtwo}, using differential blow-ups described in Section~\ref{secblowup}. The above form $p_{S^2}^{\ast}(\omega)$ extends to $C_2(S^3=\RR^3 \cup \{\infty\})$ as a model \emph{propagating form}. The closure of $p_{S^2}^{-1}(Y)$ in $C_2(S^3)$ is a model \emph{propagating chain}. We define general propagating forms and propagating chains in $C_2(\rats)$ in Section~\ref{secprop} by their behaviors on the created boundary of $C_2(\rats)$. The linking number in a rational homology sphere $\rats$ is expressed in terms of these \emph{propagators} as in the above equation, in Lemma~\ref{lemlkprop}. These propagators are the main ingredient in the definitions of the invariant $\Zinvuf$ studied in this book.

\section{Blowing up in real differential topology}
\label{secblowup}

For a vector space $\vecspt$, $\sph(\vecspt)$ denotes the quotient $\sph(\vecspt)={(\vecspt \setminus \{0\})}/{\RR^{+\ast}}$,\index[N]{S@$\sph(\vecspt)={(\vecspt \setminus \{0\})}/{\RR^{+\ast}}$} where $\RR^{+\ast}$ acts by scalar multiplication. 
Recall that the \indexT{unit normal bundle} of a submanifold $\subsubc$ in a smooth manifold is the fiber bundle whose fiber over $x \in \subsubc$ is $\sph N_x(\subsubc)=\sph (N_x(\subsubc))$.

In this book, \emph{blowing up} a submanifold $\subsubc$ in a smooth manifold $\amba$ is a canonical process, which transforms $\amba$ into a smooth manifold $\blowup{\amba}{\subsubc}$ by replacing $\subsubc$ with the total space of its unit normal bundle.
Unlike blow-ups in algebraic geometry, this blow-up of differential topology, which amounts to removing an open tubular neighborhood (thought of as infinitely small) of $\subsubc$, topologically, creates boundaries.
Let us define it formally.

A smooth \emph{submanifold transverse to the ridges} of a smooth manifold $\amba$ is a subset $\subsubc$ of $\amba$ such that for any point $x \in \subsubc$ there exists a smooth open embedding $\phi$ from  $\RR^{\codimc}  \times \RR^{\ccride} \times \left[0,1\right[^{\drid}$ into $\amba$ such that $\phi(0)=x$ and the image of $\phi$  intersects $\subsubc$ exactly along
$\phi(0  \times \RR^{\ccride} \times \left[0,1\right[^{\drid})$. Here $\codimc$ is the \emph{codimension} of $\subsubc$, $\drid$ and $\ccride$ are integers, which depend on $x$.

\begin{definition}
\label{defblodifun}
Let $\subsubc$ be a smooth submanifold transverse to the ridges of a smooth manifold $\amba$.
The \indexT{blow-up} $\blowup{\amba}{\subsubc}$ \index[N]{Blowup@$\blowup{\amba}{\subsubc}$ blow-up} is the unique smooth manifold $\blowup{\amba}{\subsubc}$ (with possible ridges)
equipped with a canonical smooth projection 
\begin{equation*}\pbl \colon \blowup{\amba}{\subsubc} \rightarrow \amba\end{equation*}
called the \indexT{blowdown map} such that 
\begin{enumerate}
\item the restriction of $\pbl$ to $\pbl^{-1}(\amba \setminus \subsubc)$ is a canonical diffeomorphism onto $\amba \setminus \subsubc$, which identifies $\pbl^{-1}(\amba \setminus \subsubc)$ with  $\amba \setminus \subsubc$ (we will simply regard $\amba \setminus \subsubc$ as a subset of $\blowup{\amba}{\subsubc}$ via this identification),
\item there is a canonical identification of $\pbl^{-1}(\subsubc)$ with the total space $\sph N(\subsubc)$ of the unit normal bundle to $\subsubc$ in $\amba$,
\item the restriction of $\pbl$ to $\pbl^{-1}(\subsubc)=\sph N(\subsubc)$ is the bundle projection from $\sph N(\subsubc)$ to $\subsubc$,
\item any smooth diffeomorphism $\phi$ from $\RR^{\codimc}  \times \RR^{\ccride} \times \left[0,1\right[^{\drid}$ onto an open subset $\phi(\RR^{\codimc}  \times \RR^{\ccride} \times \left[0,1\right[^{\drid})$ of $\amba$ whose image intersects $\subsubc$ exactly along
$\phi(0  \times \RR^{\ccride} \times \left[0,1\right[^{\drid})$, for natural integers $\codimc,\ccride,\drid$, provides a smooth embedding  
\begin{equation*}\begin{array}{lll}\left[0, \infty\right[ \times S^{\codimc-1} \times  (\RR^{\ccride} \times \left[0,1\right[^{\drid})  &\hfl{\tilde{\phi}} &\blowup{\amba}{\subsubc} \\ (\lambda \in \left]0, \infty\right[,v,x) &\mapsto & \phi(\lambda v,x) \\ (0,v,x) &\mapsto & T\phi(0,x)(v) \in \sph N(\subsubc) \end{array}\end{equation*} 
with open image in $\blowup{\amba}{\subsubc}$.
\end{enumerate}
\end{definition}
\bpo{Proof that the definition is consistent}
Use local diffeomorphisms of the form $\tilde{\phi}$ and charts
on $\amba \setminus \subsubc$ to construct an atlas for $\blowup{\amba}{\subsubc}$.
These charts are obviously compatible over $\amba \setminus \subsubc$. We 
check compatibility for charts 
 $\tilde{\phi}$ and $\tilde{\psi}$ induced by embeddings $\phi$ and $\psi$ as in the statement.
For those, transition maps 
may be written as
\begin{equation*}(\lambda,u,x) \mapsto \left(\tilde{\lambda}=\lVert p_1 \circ \psi^{-1} \circ \phi(\lambda u ,x )\rVert, \tilde{u},\tilde{x}=p_2 \circ \psi^{-1} \circ \phi(\lambda u,x)\right),\end{equation*}
where $p_1$ and $p_2$ respectively denote the projections on the first and the second factors of $\RR^{\codimc}  \times \left( \RR^{\ccride} \times \left[0,1\right[^{\drid}\right)$, and
\begin{equation*} \tilde{u} = \left\{ \begin{array}{ll}\frac{p_1 \circ \psi^{-1} \circ \phi(\lambda u ,x )}{\tilde{\lambda}} &\mbox{if} \; \lambda \neq 0 \\
\frac{T\left(p_1 \circ \psi^{-1} \circ \phi\right)(0,x)(u)}{\lVert T\left(p_1 \circ \psi^{-1} \circ \phi\right)(0,x)(u)\rVert} &\mbox{if} \; \lambda = 0\end{array} \right. \end{equation*}
In order to check that this is smooth, write
\begin{equation*}p_1 \circ \psi^{-1} \circ \phi(\lambda u ,x ) = \lambda \int_0^1 T\left(p_1 \circ \psi^{-1} \circ \phi\right)(t \lambda u,x )(u)dt,\end{equation*}
and check that the integral does not vanish when $\lambda$ is small enough.
The restriction to $S^{\codimc-1}$ of $T\left(p_1 \circ \psi^{-1} \circ \phi\right)(0,x)$ is an injection. So for any $u_0 \in S^{\codimc-1}$, there exists a neighborhood of $(0,u_0)$ in $\RR \times S^{\codimc-1}$ such that for any $(\lambda,u)$ in this neighborhood, we have the following condition about the scalar product \begin{equation*}\Bigl\langle T\left(p_1 \circ \psi^{-1} \circ \phi\right)(\lambda u ,x )(u),T\left(p_1 \circ \psi^{-1} \circ \phi\right)(0,x)(u) \Bigr\rangle \; > 0.\end{equation*}
Therefore, there exists $\varepsilon  > 0$ such that for any $\lambda \in \left]-\varepsilon,\varepsilon\right[$ and for any $u \in S^{\codimc-1}$, we have
\begin{equation*}\Bigl\langle T\left(p_1 \circ \psi^{-1} \circ \phi\right)(\lambda u ,x )(u),T\left(p_1 \circ \psi^{-1} \circ \phi\right)(0,x)(u) \Bigr\rangle \; > 0.\end{equation*}
Then 
\begin{equation*} \tilde{\lambda} = \lambda \lefnorm{\int_0^1 T\left(p_1 \circ \psi^{-1} \circ \phi\right)(t \lambda u,x )(u)dt}\end{equation*}
is a smooth function (defined even when $\lambda \leq 0$)
and 
\begin{equation*} \tilde{u} =
\frac{\int_0^1 T\left(p_1 \circ \psi^{-1} \circ \phi\right)(t \lambda u,x )(u)dt}{\lVert\int_0^1 T\left(p_1 \circ \psi^{-1} \circ \phi\right)(t \lambda u,x )(u)dt\rVert}\end{equation*} is smooth, too. Thus, our atlas is compatible. This defines $\blowup{\amba}{\subsubc}$ together with its smooth structure.
The projection $\pbl$ maps $\tilde{\phi}(\lambda, v,x)$ to ${\phi}(\lambda v,x)$ in a chart as above. Thus, it is obviously smooth, and it has the desired properties.
\eop

Note the following immediate proposition.
\begin{proposition}
\label{propblodifone}
 The blown-up manifold $\blowup{\amba}{\subsubc}$ is homeomorphic to the complement in $\amba$ of an open tubular neighborhood of $\subsubc$. In particular, $\blowup{\amba}{\subsubc}$ is homotopy equivalent to $\amba \setminus \subsubc$. If $\subsubc$ and $\amba$ are compact, then $\blowup{\amba}{\subsubc}$ is compact and it is a smooth compactification of $\amba \setminus \subsubc$.
\end{proposition}
\eopwobp

Figure~\ref{figcompblowup} first shows the result of blowing up $(0,0)$ in $\RR^2$. The closures in $\blowup{\RR^2}{(0,0)}$ of $\{0\} \times \RR^{\ast}$, $\RR^{\ast} \times \{0\}$, and the diagonal of $\left(\RR^{\ast}\right)^2$ are next blown up in $\blowup{\RR^2}{(0,0)}$.

\bfig
\centering
\begin{tikzpicture}
\begin{scope}[xshift=-1cm]
\draw [fill=gray!20, draw=white] (-1,-1) rectangle (1,1);
\draw [->] (-1,0) -- (1,0) node[right]{\scriptsize $\RR \times 0$};
\draw [->] (0,-1) -- (0,1) node[right]{\scriptsize $0 \times \RR$};
\draw [->] (-1,-1) -- (1,1) node[right]{\scriptsize $\diag$};
\end{scope}
\draw (2.3,0) node{\scriptsize Blow up $(0,0)$};
\draw [->] (1.3,-.2) -- (3.3,-.2);
\begin{scope}[xshift=4.7cm]
\draw [fill=gray!20, draw=white] (-1,-1) rectangle (1,1);
\draw [->] (-1,0) -- (1,0);
\draw [->] (0,-1) -- (0,1);
\draw [->] (-1,-1) -- (1,1);
\draw [fill=white] (0,0) circle (.3);
\draw [->,very thin] (.5,-.5) -- (.3,-.3);
\draw (.2,-.6) node[right]{\scriptsize unit normal bundle to $(0,0)$};

\end{scope}
\begin{scope}[xshift=7.5cm,yshift=-3cm] 
\draw (-2.9,0) node{\scriptsize Blow up the lines};
\draw [->] (-4.2,-.2) -- (-1.6,-.2);
\draw [fill=gray!20, draw=white] (-1,-1) rectangle (1,1);
\draw [very thick] (-.93,-1) -- (1,.93) (-1,-.93) -- (.93,1);
\draw [fill=black]  (-.12,-1) rectangle (.13,1) (-1,-.12) rectangle (1,.13);
\draw [fill=white] (0,0) circle (.3);
\draw [fill=white,draw=white] (-1,-1) -- (-.93,-1) -- (1,.93) -- (1,1) -- (.93,1) -- (-1,-.93) -- (-1,-1);
\draw [fill=white,draw=white]  (-.1,-1.1) rectangle (.1,1.1) (-1.1,-.1) rectangle (1.1,.1);
\end{scope}
\end{tikzpicture}
\caption{A composition of blow-ups}\label{figcompblowup}
\end{figure}

\begin{proposition}
\label{propblodifdeux}
Let $\submb$ and $\subsubc$ be two smooth submanifolds transverse to the ridges of a $C^{\infty}$ manifold $\amba$. Assume that $\subsubc$ is a smooth submanifold of $\submb$ transverse to the ridges of $\submb$.
\begin{enumerate}
\item The closure $\overline{\submb \setminus \subsubc}$ of $(\submb \setminus \subsubc)$ in $\blowup{\amba}{\subsubc}$ is a submanifold of $\blowup{\amba}{\subsubc}$. It intersects \begin{equation*} \sph N(\subsubc) \subseteq \partial \blowup{\amba}{\subsubc}\end{equation*} as the unit normal bundle $\sph N_{\submb}(\subsubc)$ to $\subsubc$ in $\submb$. It is canonically diffeomorphic to $\blowup{\submb}{\subsubc}$.
\item The blow-up $\blowup{\blowup{\amba}{\subsubc}}{\overline{\submb \setminus \subsubc}}$ of $\blowup{\amba}{\subsubc}$ along $\overline{\submb \setminus \subsubc}$ has a canonical differential structure of a manifold with ridges. The preimage of \/ $\overline{\submb \setminus \subsubc} \subset \blowup{\amba}{\subsubc}$ in $\blowup{\blowup{\amba}{\subsubc}}{\overline{\submb \setminus \subsubc}}$ under the canonical projection
\begin{equation*}\blowup{\blowup{\amba}{\subsubc}}{\overline{\submb \setminus \subsubc}} \longrightarrow \blowup{\amba}{\subsubc} \end{equation*}
is the pull-back via the blowdown projection from $\overline{\submb \setminus \subsubc}$ to $\submb$ of the unit normal bundle
to $\submb$ in $\amba$.
\end{enumerate}
\end{proposition}
\bp

\begin{enumerate}
\item Let $x\in\subsubc$. It is always possible to choose an embedding $\phi$ into $\amba$ as in Definition~\ref{defblodifun} such that $\phi(0)=x$ and
$\phi\bigl(\RR^{\codimc}  \times \RR^{\ccride} \times \left[0,1\right[^{\drid}\bigr)$ intersects $\subsubc$ exactly along
$\phi(0  \times \RR^{\ccride} \times \left[0,1\right[^{\drid})$
and $\submb$ exactly along
$\phi(0 \times \RR^{\dimk} \times \left[0,1\right[^{\drid} )$, with $\dimk > \ccride$. (First choose an embedding suitable for $\submb$, and then modify it to suit $\subsubc$.)
Look at the induced chart $\tilde{\phi}$ of $\blowup{\amba}{\subsubc}$ near a point of $\partial \blowup{\amba}{\subsubc}$. \\
The intersection of $(\submb \setminus \subsubc)$ with the image of $\tilde{\phi}$ is \begin{equation*}\tilde{\phi}\left(\left]0,\infty\right[  \times (0 \times S^{\dimk - \ccride-1}\subset S^{\codimc-1}) \times \RR^{\ccride} \times \left[0,1\right[^{\drid} \right).\end{equation*} Thus, the closure of $(\submb \setminus \subsubc)$ intersects the image of $\tilde{\phi}$ as \begin{equation*}\tilde{\phi}\left(\left[0, \infty\right[  \times (0 \times S^{\dimk - \ccride-1}\subset S^{\codimc-1})  \times \RR^{\ccride} \times \left[0,1\right[^{\drid} \right).\end{equation*}
\item
Together with the above given charts of $\overline{\submb \setminus \subsubc}$, the smooth  injective map \begin{equation*}\RR^{\codimc -\dimk + \ccride} \times S^{\dimk - \ccride-1}  \longrightarrow S^{\codimc-1}\end{equation*} \begin{equation*}(u,y) \mapsto \frac{(u,y)}{\lVert (u,y) \rVert}\end{equation*} 
identifies $\RR^{\codimc -\dimk + \ccride}$
with the fibers of the normal bundle to $\overline{\submb \setminus \subsubc}$ in $\blowup{\amba}{\subsubc}$. The blow-up process
will therefore replace $\overline{\submb \setminus \subsubc}$ with the quotient of the corresponding $(\RR^{\codimc -\dimk + \ccride} \setminus \{0\})$-bundle
by $\left]0,\infty\right[$, which is of course the pull-back under the blowdown projection $(\overline{\submb \setminus \subsubc} \longrightarrow \submb)$ of the unit normal bundle to $\submb$ in $\amba$.
\end{enumerate}
\eop

The fiber $\sph N_c(\subsubc)$ is oriented as the boundary of a unit ball of $N_c(\subsubc)$.

\section{The configuration space \texorpdfstring{$C_2(\rats)$}{C2(R)}}
\label{secCtwo}

Regard $S^3$ as $\RR^3 \cup \{\infty\}$ or as two copies of $\RR^3$ identified along $\RR^3 \setminus \{0\}$ by the (exceptionally orientation-reversing) diffeomorphism $x \mapsto x/\norm{x}^2$.

Let $(-S^2_{\infty})$ denote the unit normal bundle to $\infty$ in $S^3$, so $\blowup{S^3}{\infty}=\RR^3 \cup S^2_{\infty}$ and 
$\partial \blowup{S^3}{\infty} = S^2_{\infty}$.
There is a canonical orientation-preserving diffeomorphism $p_{\infty} \colon  S^2_{\infty} \rightarrow S^2$\index[N]{Projections!ptaainfini@$p_{\infty} \colon  S^2_{\infty} \rightarrow S^2$}, such that $x \in S^2_{\infty}$ is the limit in $\blowup{S^3}{\infty}$ of a sequence of points of $\RR^3$ tending to $\infty$ along any line of $\RR^3$ directed by $p_{\infty}(x) \in S^2$, in the direction of the line. 

Let $\mathring{B}_{1,\infty}$ \index[N]{Bminfinity@$\mathring{B}_{1,\infty}$ ball around $\infty$} (resp. $B_{1,\infty}$) denote the complement of the closed (resp. open) ball of radius one of $\RR^3$ in $S^3$. 
Let $\mathring{B}_{2,\infty}$ be the complement in $S^3$ of the closed ball ${B}(2)$ of radius $2$ in $\RR^3$.

Fix a rational homology sphere $\rats$ and a point $\infty$ of $\rats$. Set $\crats=\rats \setminus \{\infty\}$. 
Identify a neighborhood of $\infty$ in $\rats$ with $\mathring{B}_{1,\infty}$.  The ball $\mathring{B}_{2,\infty}$ is a smaller neighborhood of $\infty$ in $\rats$ via the understood identification.
Then $\ballb_{\rats}=\rats \setminus \mathring{B}_{2,\infty}$ is a compact rational homology ball diffeomorphic to $\blowup{\rats}{\infty}$.

Define the {\em configuration space\/} $C_2(\rats)$\index[N]{Configuration spaces!Ca@$C_2(\rats)$} to be the compact $6$-manifold with boundary and ridges obtained from $\rats^2$ by blowing up first $(\infty,\infty)$ in $\rats^2$, and next the closures of $\{\infty\} \times \crats$, $\crats \times \{\infty\}$, and the diagonal of $\crats^2$ in $\blowup{\rats^2}{(\infty,\infty)}$, as in Figure~\ref{figcompblowup}. 
Then $\partial C_2(\rats)$ contains the unit normal bundle $\bigl(({T\crats^2}/{\diag(T\crats^2)}) \setminus \{0\}\bigr)/\RR^{+\ast}$ to the diagonal of $\crats^2$. This bundle is identified with the unit tangent bundle $\ST \crats$ to $\crats$ by the map $\bigl(\left[(x,y)\right] \mapsto \left[y-x\right]\bigr)$.

\begin{lemma}
\label{lemstrucctwo}
Let $\check{C}_2(\rats)=\crats^2 \setminus \diag\left( \crats^2 \right)$.
The open manifold $C_2(\rats) \setminus \partial C_2(\rats)$ is $\check{C}_2(\rats)$, and the inclusion $\check{C}_2(\rats) \hookrightarrow C_2(\rats)$ is a homotopy equivalence. In particular, $C_2(\rats)$ is a compactification of $\check{C}_2(\rats)$ homotopy equivalent to $\check{C}_2(\rats)$. It has the same rational homology as the sphere $S^2$.
The manifold $C_2(\rats)$ is a smooth compact $6$-dimensional manifold with boundary and ridges. Its boundary is
\begin{equation*}\partial C_2(\rats)=p_{\rats^2}^{-1}(\infty,\infty) \cup  ( S^2_{\infty}\times \crats) \cup (-\crats \times S^2_{\infty}) \cup \ST \crats.\end{equation*}
Furthermore, there is a canonical smooth projection $p_{\rats^2}\colon C_2(\rats) \rightarrow \rats^2$.
\end{lemma}
\bp This lemma is a corollary of Propositions~\ref{propblodifone} and \ref{propblodifdeux}, and Lemma~\ref{lemhomstwo}. We give a few additional arguments to check that we can perform the three blow-ups in 
$\blowup{\rats^2}{(\infty,\infty)}$ simultaneously, and we take a closer look at the structure of $p_{\rats^2}^{-1}(\infty,\infty)$ below.

Blowing up $(\infty,\infty)$ in $B_{1,\infty}^2$ transforms a neighborhood of $(\infty,\infty)$
into the product $\left[0,1\right[ \times S^5$. Explicitly, there is a map
\begin{equation*}\begin{array}{llll}\psi \colon& \left[0,1\right[ \times S^5 & \rightarrow & \lefblowup{B_{1,\infty}^2}{(\infty,\infty)} \\
  & \left(\lambda \in \left]0,1\right[, (x \neq 0,y \neq 0) \in S^5 \subset (\RR^3)^2\right)  & \mapsto & \left(\frac1{\lambda\norm{x}^2}x,\frac1{\lambda\norm{y}^2}y\right)\\
& \left(\lambda \in \left]0,1\right[, (0,y \neq 0) \in S^5 \subset (\RR^3)^2\right)  & \mapsto & \left(\infty,\frac1{\lambda\norm{y}^2}y\right)\\
& \left(\lambda \in \left]0,1\right[, (x \neq 0,0) \in S^5 \subset (\RR^3)^2\right)  & \mapsto & \left(\frac1{\lambda\norm{x}^2}x,\infty\right),
  \end{array}\end{equation*}
which is a diffeomorphism onto its open image.

Here, the explicit image of $(\lambda \in \left]0,1\right[, (x \neq 0,y \neq 0) \in S^5 \subset (\RR^3)^2)$ is written in $\bigl(\mathring{B}_{1,\infty} \setminus \{\infty\}\bigr)^2 \subset \bigblowup{\mathring{B}_{1,\infty}^2}{(\infty,\infty)}$, where $\mathring{B}_{1,\infty} \setminus \{\infty\} \subset \RR^3$.
The image of $\psi$ is a neighborhood of the preimage of $(\infty,\infty)$ under the blowdown map 
\begin{equation*}\bigblowup{\rats^2}{(\infty,\infty)} \hflv{p_1} \rats^2.\end{equation*}
This neighborhood respectively intersects $\infty \times \crats$, $\crats \times \infty$, and $\diag\left(\crats^2\right)$ as $\psi(\left]0,1\right[ \times 0 \times S^2)$, $\psi(\left]0,1\right[ \times S^2\times 0)$, and $\psi\bigl(\left]0,1\right[ \times \left(S^5 \cap \diag\left((\RR^3)^2\right)\right)\bigr)$. 
In particular, the closures of $\infty \times \crats$, $\crats \times \infty$, and $\diag\left(\crats^2\right)$ in $\blowup{\rats^2}{(\infty,\infty)}$ intersect the boundary $\psi(0 \times S^5)$ of $\blowup{\rats^2}{(\infty,\infty)}$ as three disjoint spheres in $S^5$. They are respectively isomorphic to
 $\infty \times \blowup{\rats}{\infty}$, 
$\blowup{\rats}{\infty} \times \infty$, and $\diag\left(\blowup{\rats}{\infty}^2\right)$.
Thus, the next steps will be three blow-ups along these three disjoint smooth manifolds.

These blow-ups will preserve the product structure $\psi(\left[0,1\right[ \times .)$.
Thus, $C_2(\rats)$ is a smooth compact $6$-dimensional manifold with boundary, with three {\em ridges\/} $S^2 \times S^2$ in $p_{\rats^2}^{-1}(\infty,\infty)$. A neighborhood of these ridges in $C_2(\rats)$ is diffeomorphic to $\left[0,1\right[^2 \times S^2 \times S^2$.
\eop

Let $\iota_{S^2}$
denote the \emph{antipodal map} of $S^2$, which sends $x$ to $\iota_{S^2}(x)=-x$.

\begin{lemma}
\label{lemextproj}
The map $p_{S^2}$ of Lemma~\ref{lempstt} extends smoothly to $C_2(S^3)$. Its extension $p_{S^2}$
satisfies
\begin{equation*}p_{S^2}=\left\{\begin{array}{ll} \iota_{S^2} \circ p_{\infty} \circ p_1 \;\;& \mbox{on} \;S^2_{\infty} \times \RR^3\\ p_{\infty} \circ p_2 \;\;& \mbox{on} \;  \RR^3 \times S^2_{\infty}\\ p_2 \;\;& \mbox{on} \;\ST \RR^3 
{=} \RR^3 \times S^2,\end{array}\right.\end{equation*}
where $p_1$ and $p_2$ denote the projections on the first and second factor with respect to the above expressions.
\end{lemma}
\bp 
Near the diagonal of $\RR^3$, we have a chart of $C_2(S^3)$
\begin{equation*}\psi_d\colon \RR^3 \times \left[0, \infty\right[  \times S^2 \longrightarrow C_2(S^3),\end{equation*}
which maps $( x \in \RR^3, \lambda \in \left]0,\infty\right[ ,y \in S^2)$ to $(x, x + \lambda y) \in (\RR^3)^2$. Here, $p_{S^2}$ extends as the projection onto the $S^2$ factor.

Consider the orientation-reversing embedding $\phi_{\infty}$
\begin{equation*}\begin{array}{llll}\phi_{\infty}\colon &\RR^3 &\longrightarrow &\bigl(S^3=\RR^3 \cup \{\infty\}\bigr)\\
& \mu (x \in S^2) & \mapsto & \left\{\begin{array}{ll} \infty \;&\;\mbox{if}\; \mu=0\\
\frac{1}{\mu}x \;&\;\mbox{otherwise.} \end{array}\right.\end{array}\end{equation*}
Note that this chart induces the identification $p_{\infty}$ of the unit normal bundle $S^2_\infty$ to $\{\infty\}$ in $S^3$ with $S^2$. When $\mu \neq 0$,
\begin{equation*}p_{S^2}(\phi_{\infty}(\mu x), y \in \RR^3)= \frac{\mu y-x}{\norm{\mu y- x}}.\end{equation*}
Then $p_{S^2}$ can be extended smoothly on $S^2_\infty \times \RR^3$ (where $\mu=0$) by
\begin{equation*}p_{S^2}(x \in S^2_\infty, y \in \RR^3) = -x.\end{equation*}
Similarly, set
$p_{S^2}( x \in \RR^3, y \in S^2_\infty) = y.$
Now, with the map $\psi$ of the proof of Lemma~\ref{lemstrucctwo}, if $x$ and $y$ are not equal to zero and if they are distinct, then we have
\begin{equation*}p_{S^2}\circ \psi\Bigl(\bigl(\lambda ,(x,y)\bigr)\Bigr)=\frac{\frac{y}{\norm{y}^2}-\frac{x}{\norm{x}^2}}
{\norm{\frac{y}{\norm{y}^2}-\frac{x}{\norm{x}^2}}}
=\frac{\norm{x}^2y-\norm{y}^2x}
{\norm{\norm{x}^2y-\norm{y}^2x}}\end{equation*}
when $\lambda \neq 0$. This map extends to $\blowup{\rats^2}{(\infty,\infty)}$ outside the boundaries of $\infty \times \blowup{\rats}{\infty}$, 
$\blowup{\rats}{\infty} \times \infty$ and $\diag\left(\blowup{\rats}{\infty}^2\right)$, naturally, by keeping the same formula when $\lambda = 0$.

Let us check that the map $p_{S^2}$ extends smoothly over the boundary of $\diag\left(\blowup{\rats}{\infty}^2\right)$.
There is a chart of $C_2(\rats)$ near the preimage of this boundary
in $C_2(\rats)$
\begin{equation*}\psi_2 \colon \left[0, \infty\right[ \times \left[0, \infty\right[ \times S^2 \times S^2 \longrightarrow C_2(S^3),\end{equation*}
which maps $(\lambda \in \left]0,\infty\right[ , \mu \in  \left]0,\infty\right[,  x \in S^2, y \in S^2)$
to $(\phi_{\infty}(\lambda x), \phi_{\infty}(\lambda(x + \mu y)))$. With these coordinates, the map $p_{S^2}$ may be written as
\begin{equation*}(\lambda,\mu,x,y) \mapsto \frac{y- 2\langle x,y \rangle x -\mu x}
{\norm{y- 2\langle x,y \rangle x -\mu x}}.\end{equation*} Therefore, it extends smoothly along $\mu=0$. We check that $p_{S^2}$ extends smoothly over the boundaries of $\infty \times \blowup{\rats}{\infty}$ and 
$\blowup{\rats}{\infty} \times \infty$ similarly.
\eop

\begin{definition}
\label{defparasyst}
Let $\taust$ \index[N]{tzaus@$\taust$ standard parallelization of $\RR^3$} denote the standard parallelization of $\RR^3$. Say that a parallelization \begin{equation*}\tau \colon \crats \times \RR^3 \rightarrow T \crats\end{equation*} of $\crats$ that coincides with $\taust$ 
on $\mathring{B}_{2,\infty} \setminus \{\infty\}$ is {\em asymptotically standard.\/} \index[T]{asymptotically standard parallelization} According to Proposition~\ref{propexistparasymptriv}, such a parallelization exists.
Such a parallelization identifies $\ST \crats$ with $\crats \times S^2$.
\end{definition}

\begin{proposition}
\label{propprojbord}
For any asymptotically standard parallelization $\tau$ of $\crats$,
there exists a smooth map $p_{\tau} \colon \partial C_2(\rats) \rightarrow S^2$\index[N]{Projections!ptau@$p_{\tau} \colon  \partial C_2(\rats) \rightarrow S^2$} such that
\begin{equation*}p_{\tau}= 
\left\{\begin{array}{ll} 
p_{S^2} \;\;& \mbox{on} \;  p_{\rats^2}^{-1}(\infty,\infty)\\
\iota_{S^2} \circ p_{\infty} \circ p_1 \;\;& \mbox{on} \;S^2_{\infty} \times \crats\\
 p_{\infty} \circ p_2 \;\;& \mbox{on} \;  \crats \times S^2_{\infty}\\ 
p_2 \;\;& \mbox{on} \;\ST \crats \stackrel{\tau} {=} \crats \times S^2,\end{array}\right.\end{equation*}
where $p_1$ and $p_2$ denote the projections on the first and second factors with respect to the above expressions.
\end{proposition}
\bp This is a consequence of Lemma~\ref{lemextproj}.\eop

\begin{definition}
\label{defasyrathommRthree}
Define an \indexT{asymptotic rational homology $\RR^3$} to be a pair $(\crats,\tau)$, in which $\crats$ is the union over $\left]1,2\right] \times S^2$ of a rational homology ball $\ballb_{\rats}$\index[N]{BR@$\ballb_{\rats}$ rational homology ball} and the complement $\mathring{B}_{1,\infty} \setminus \{\infty\}$ of the unit ball in $\RR^3$, and $\tau$ is an asymptotically standard parallelization of $\crats$.
Such a pair $(\crats,\tau)$ defines the rational homology sphere $\rats=\crats \cup \{\infty\}$ canonically. So \say{Let $(\crats,\tau)$ be an asymptotic rational homology $\RR^3$} is a shortcut for
\say{Let $\rats$ be a rational homology sphere, equipped with an embedding  $\phi_{\rats} \colon \mathring{B}_{1,\infty} \hookrightarrow \rats$ and an asymptotically standard parallelization $\tau$ of the complement $\crats = \rats \setminus \{\phi_{\rats}(\infty)\}$, with respect to $\phi_{\rats}$}. (The embedding $\phi_{\rats}$ is understood. We will not mention it anymore.
We will just view it as an inclusion  $\mathring{B}_{1,\infty} \hookrightarrow \rats$ and denote $\phi_{\rats}(\infty)$ by $\infty$.)
\end{definition}

\section{On propagators}
\label{secprop}

A \indexT{volume-one form of $S^2$} is a $2$-form $\omega_S$ of $S^2$ such that $\int_{S^2}\omega_S=1$. (See Appendix~\ref{chapDeRham} for a short survey of differential forms and de Rham cohomology.)

Let $(\crats,\tau)$ be an asymptotic rational homology $\RR^3$. Recall the map \begin{equation*}p_{\tau} \colon \partial C_2(\rats) \rightarrow S^2\end{equation*} of Proposition~\ref{propprojbord}. 

\begin{definition}\label{defpropagatorone}
A \emph{propagating chain}\index[T]{propagating!chain} of $(C_2(\rats),\tau)$ is a $4$-chain $\propP$ of $C_2(\rats)$
such that $\partial \propP =p_{\tau}^{-1}(\veca)$ for some $\veca \in S^2$.
A \emph{propagating form}\index[T]{propagating!form} of $(C_2(\rats),\tau)$ is a closed $2$-form $\omega$ on $C_2(\rats)$ whose restriction to $\partial C_2(\rats)$ may be expressed as $p_{\tau}^{\ast}(\omega_S)$ for some volume-one form $\omega_S$ of $S^2$.
We will call propagating chains and propagating forms \emph{propagators}\index[T]{propagator} when their nature is clear from the context.
\end{definition}

\begin{examples}
\label{exapropsthree}
Recall the map $p_{S^2} \colon C_2(S^3) \rightarrow S^2$ of Lemma~\ref{lemextproj}. For any $a \in S^2$, $p_{S^2}^{-1}(\veca)$ is a propagating chain of $(C_2(S^3),\taust)$. For any $2$-form $\omega_S$ of $S^2$ such that $\int_{S^2}\omega_S=1$, $p_{S^2}^{\ast}(\omega_S)$ is a propagating form of $(C_2(S^3),\taust)$.
\end{examples}

For our general $\QQ$-sphere $\rats$, propagating chains of $(C_2(\rats),\tau)$ exist because the $3$-cycle $p_{\tau}^{-1}(\veca)$ of $\partial C_2(\rats)$ bounds in $C_2(\rats)$ since $H_3(C_2(\rats);\QQ)=0$. Dually, propagating forms of $(C_2(\rats),\tau)$ exist because the restriction induces a surjective map \begin{equation*}H^2(C_2(\rats);\RR) \rightarrow H^2(\partial C_2(\rats);\RR)\end{equation*} since $H^3(C_2(\rats),\partial C_2(\rats);\RR)$ is trivial.

When $\rats$ is a $\ZZ$-sphere, there exist propagating chains of $(C_2(\rats),\tau)$ that are smooth $4$-manifolds properly embedded in $C_2(\rats)$. See Corollary~\ref{corpropmanifint}. 

\begin{definition}
\label{defpropagatortwo}
A \emph{propagating form}\index[T]{propagating!form}  of $C_2(\rats)$ is a closed $2$-form $\omega_p$ on $C_2(\rats)$ whose restriction to $\partial C_2(\rats)\setminus \ST \ballb_{\rats}$ is equal to $\projp_{\partau}^{\ast}(\omega)$ for some volume-one form $\omega_S$ of $S^2$ and some asymptotically standard parallelization $\partau$. Similarly, a \emph{propagating chain}\index[T]{propagating!chain} of $C_2(\rats)$ is a $4$-chain $\propP$ of $C_2(\rats)$ such that $\partial P \subset \partial C_2(\rats)$ and  $\partial P \cap \left(\partial C_2(\rats)\setminus \ST \ballb_{\rats}\right)$ is equal to $\projp_{\partau}^{-1}(\veca)$ for some $\veca \in S^2$. (These definitions do not depend on $\partau$.)
\end{definition}
Explicit propagating chains associated to Heegaard splittings, constructed with Greg Kuperberg, are described in Subsection~\ref{subMorse}. They are integral chains multiplied by $\frac{1}{\cardlef{H_1(\rats;\ZZ)}}$, where $\cardlef{H_1(\rats;\ZZ)}$ is the cardinality of $H_1(\rats;\ZZ)$.

Since $C_2(\rats)$ is homotopy equivalent to $(\crats^2 \setminus \diag\left(\crats^2\right))$, Lemma~\ref{lemhomstwo} ensures that $H_2(C_2(\rats);\QQ)=\QQ\left[S\right]$, where the canonical generator $\left[S\right]$
is the homology class of a fiber of $\ST \crats \subset \partial C_2(\rats)$.
For a two-component link $(J,K)$ of $\crats$, the homology class $\left[J \times K\right]$ of
$J \times K$ in $H_2(C_2(\rats);\QQ)$ is $lk(J,K)\left[S\right]$, according to Proposition~\ref{propdeflkeq}.

\begin{lemma}
\label{lemlkprop}
Let $(\crats,\tau)$ be an asymptotic rational homology $\RR^3$. Let $C$ be a two-cycle of $C_2(\rats)$. For any propagating chain $\propP$ of $C_2(\rats)$ transverse to $C$ and for any propagating form $\omega$ of $C_2(\rats)$, we have
\begin{equation*}\left[C\right]=\left(\int_{C}\omega\right)\left[S\right]=\langle C, \propP \rangle_{\!C_2(\rats)\,}\left[S\right]\end{equation*} in $H_2(C_2(\rats);\QQ)=\QQ\left[S\right]$.
In particular, we have
\begin{equation*}lk(J,K)=\int_{J \times K}\omega=\langle J \times K, \propP \rangle_{\!C_2(\rats)\,}.\end{equation*} for any two-component link $(J,K)$ of $\crats$.
\end{lemma}
\bp
Fix a propagating chain $\propP$, the algebraic intersection $\langle C, \propP \rangle_{\!C_2(\rats)\,}$ depends only on the homology class $\left[C\right]$ of $C$ in $C_2(\rats)$. Similarly, since $\omega$ is closed, $\int_{C}\omega$ depends only on $\left[C\right]$.
(Indeed, if $C$ and $C^{\prime}$ cobound a chain $D$ transverse to $\propP$, then $C \cap \propP$ and $C^{\prime} \cap \propP$ cobound $\pm (D\cap \propP)$, and $\int_{\partial D=C^{\prime}-C}\omega=\int_Dd\omega$ according to Stokes' theorem.)
Furthermore, the dependence on $\left[C\right]$ is linear. Therefore, it suffices to check the lemma for a chain that represents the canonical generator $\left[S\right]$ of $H_2(C_2(\rats);\QQ)$. Any fiber of $\ST \crats$ is such a chain.
\eop

\begin{definition}
\label{defprophomogen}
A propagating form $\omega$ of $C_2(\rats)$ is \emph{homogeneous} if its restriction to $\partial C_2(\rats)\setminus \ST \ballb_{\rats}$ is $\projp_{\tau}^{\ast}(\omega_{S^2})$ for the homogeneous volume form $\omega_{S^2}$ of $S^2$ of total volume $1$. \index[T]{homogeneous!propagating form}
\end{definition}

\begin{definition}
\label{defantisympropform}
Let $\iota$ \index[N]{iota@$\iota$ involution of $C_2(\rats)$} be the involution of $C_2(\rats)$ that exchanges the two coordinates in $\crats^2 \setminus \diag\left(\crats^2\right)$.
An \emph{antisymmetric} propagating form is a propagating form such that $\iota_{\ast}(\omega)=-\omega$.\index[T]{antisymmetric!propagating form}
\end{definition}

\begin{example}
 The propagating form $p_{S^2}^{\ast}(\omega_{S^2})$ of $(C_2(S^3),\partau)$ is antisymmetric.
\end{example}

The following easy lemma ensures the existence of antisymmetric propagating forms of $\bigl(C_2(\rats),\partau\bigr)$.

\begin{lemma}
\label{lemiotaform}
Let $\omega_0$ be a propagating form of $(C_2(\rats),\tau)$.
Then $(-\iota^{\ast}(\omega_0))$
and $\omega=\frac12(\omega_{0}-\iota^{\ast}(\omega_0))$ are propagating forms of $(C_2(\rats),\tau)$, and we have
 $\iota^{\ast}(\omega)=-\omega$.
If $\omega_0$ is homogeneous, then $(-\iota^{\ast}(\omega_0))$
and $\omega=\frac12(\omega_{0}-\iota^{\ast}(\omega_0))$ are also homogeneous.
\end{lemma}
\bp Let $\omega_S$ be the volume-one form of $S^2$ such that $\omega_{0}\vert_{\partial C_2(\rats)}=p_{\tau}^{\ast}(\omega_S)$. Recall the antipodal map  $\iota_{S^2}$ of $S^2$. The form $(-\iota_{S^2}^{\ast}(\omega_S))$ is a volume-one form of $S^2$, and
we have $(-\iota^{\ast}(\omega_0))\vert_{\partial C_2(\rats)}=
p_{\tau}^{\ast}(-\iota_{S^2}^{\ast}(\omega_S))$.
\eop

Also note the following lemma.

\begin{lemma}
\label{lemetactwo}
Let $\omega_{\veca}$ and $\omega^{\prime}_{\veca}$ be two propagating forms of $(C_2(\rats),\tau)$, respectively restricting to $\partial C_2(\rats)$ as $p_{\tau}^{\ast}(\omega_A)$ and $p_{\tau}^{\ast}(\omega^{\prime}_A)$, for two volume-one forms $\omega_A$ and $\omega^{\prime}_A$ of $S^2$. 
There exists a one-form $\eta_A$ on $S^2$ such that $\omega^{\prime}_A=\omega_A +d\eta_A$.
For any such $\eta_A$, there exists a one-form $\eta$ on $C_2(\rats)$ such that $\omega^{\prime}_{\veca}-\omega_{\veca}=d\eta$ and the restriction of $\eta$ to $\partial C_2(\rats)$ is $p_{\tau}^{\ast}(\eta_A)$.
\end{lemma}
\bp
Since $\omega_{\veca}$ and  $\omega^{\prime}_{\veca}$ are cohomologous, there exists a one-form $\eta$ on $C_2(\rats)$ such that $\omega^{\prime}_{\veca}=\omega_{\veca} +d\eta$. Similarly, since $\int_{S^2}\omega^{\prime}_A=\int_{S^2}\omega_A$, there exists a one-form $\eta_A$ on $S^2$ such that $\omega^{\prime}_A=\omega_A +d\eta_A$.
We have $d(\eta-p_{\tau}^{\ast}(\eta_A))=0$ on $\partial C_2(\rats)$.
The exact sequence with real coefficients
\begin{equation*} 0=H^1(C_2(\rats)) \longrightarrow H^1(\partial C_2(\rats)) \longrightarrow H^2(C_2(\rats), \partial C_2(\rats)) \cong H_4(C_2(\rats))=0, \end{equation*}
implies $H^1(\partial C_2(\rats);\RR)=0$. 
Therefore, there exists a function $f$ from $\partial C_2(\rats)$ to $\RR$ such that \begin{equation*}df =\eta -p_{\tau}^{\ast}(\eta_A)\end{equation*} on $ \partial C_2(\rats)$.
Extend $f$ to a $C^{\infty}$ map on $C_2(\rats)$ and replace $\eta$ with $(\eta -df)$.
\eop

\chapter{The Theta invariant}
\label{chapTheta}

Let $(\crats,\tau)$ be an asymptotic rational homology $\RR^3$.
We are ready to define the topological invariant $\Theta$ for $(\crats,\tau)$ to be the algebraic triple intersection of three propagating chains of $(C_2(\rats),\tau)$, in Section~\ref{secdefTheta}. In Section~\ref{secdefThetam}, we use relative Pontrjagin classes introduced in Section~\ref{secannpont} to turn $\Theta$ into a topological invariant of rational homology $3$-spheres.

\section{The \texorpdfstring{$\Theta$}{Theta}-invariant of \texorpdfstring{$(\rats,\tau)$}{(R,tau)}}
\label{secdefTheta}

We defined the algebraic intersection $\langle A, B, C \rangle_{\!D\,}$ of three transverse compact submanifolds $A$, $B$, $C$ in a manifold $D$ such that the sum of the codimensions of $A$, $B$, and $C$ is the dimension of $D$
in Subsection~\ref{subalgint}. It is the sum over the intersection points $a$ of  $A\cap B \cap C$ of the associated signs, where the sign of $a$ is positive if and only if the orientation of $D$ is induced by the orientation of $N_aA \oplus N_aB \oplus N_aC$, where $N_aA$, $N_aB$, and $N_aC$ are respectively identified with $(T_aA)^{\perp}$, $(T_aB)^{\perp}$, and $(T_aC)^{\perp}$ with the help of a Riemannian metric.

\begin{theorem}
\label{thmdefTheta}
Let $(\crats,\tau)$ be an asymptotic rational homology $\RR^3$.
Let $\propP_a$, $\propP_b$, and $\propP_c$ be three transverse propagators of $(C_2(\rats),\tau)$ with respective boundaries $p_{\tau}^{-1}(\veca)$, $p_{\tau}^{-1}(b)$, and $p_{\tau}^{-1}(c)$ for three distinct points $a$, $b$, and $c$ of $S^2$.
Then 
\begin{equation*}\Theta(\rats,\tau)=\bigl\langle \propP_{\veca}, \propP_b, \propP_c \bigr\rangle_{\!C_2(\rats)\,}\end{equation*}
does not depend on the chosen propagators $\propP_{\veca}$, $\propP_b$, and $\propP_c$. It is a topological invariant of $(\rats,\tau)$. 
For any three propagating chains $\omega_{\veca}$, $\omega_b$, and $\omega_c$ of $(C_2(\rats),\tau)$, we have
\begin{equation*}\Theta(\rats,\tau)=\int_{C_2(\rats)} \omega_{\veca} \wedge \omega_b \wedge  \omega_c.\end{equation*}
\end{theorem}
\bp
Since $H_4(C_2(\rats);\QQ)=0$, if the propagator $\propP_{\veca}$ is replaced by a propagator $\propP^{\prime}_{\veca}$ with the same boundary, the $4$-dimensional cycle $(\propP^{\prime}_{\veca}-\propP_{\veca})$ bounds a $5$-dimensional rational chain $W$ transverse to $\propP_b \cap \propP_c$. The $1$-dimensional chain $W \cap \propP_b \cap \propP_c$ does not meet $\partial C_2(\rats)$ since $\propP_b \cap \propP_c$ does not meet $\partial C_2(\rats)$. Therefore, up to a well-determined sign, the boundary of $W \cap \propP_b \cap \propP_c$ is $\propP^{\prime}_{\veca} \cap \propP_b \cap \propP_c -\propP_{\veca} \cap \propP_b \cap \propP_c$, as in Figure~\ref{figthmdefTheta}. 
\bfig
\centering
\begin{tikzpicture} \useasboundingbox (-5,-2) rectangle (5,1);
\fill [lightgray]  (0,-.5) ellipse (3.2 and 1.5);
\draw (0,-.5) ellipse (3.2 and 1.5);
\fill [white] (0,.3) ellipse (.8 and .3);
\draw (0,.3) ellipse (.8 and .3);

\fill [white, rounded corners] (1.25,-.35) rectangle (2.45,.05);
\draw (1.85,-.15) node{\scriptsize $\propP_b \cap \propP_c$};
\fill [white, rounded corners] (-3,-.7) rectangle (-2,-.3);
\draw (-2.5,-.5) node{\scriptsize $C_2(\rats)$};
\fill [gray] (0,0) .. controls (0,-.3) and (1,-.7) .. (1,-1) .. controls (1,-1.3) and (0,-1.7) .. (0,-2)
 .. controls (0,-1.7) and (-1,-1.3) ..  (-1,-1)  .. controls (-1,-.7) and (0,-.3) ..  (0,0);
\draw [-<] (0,0) .. controls (0,-.3) and (1,-.7) .. (1,-1);
\draw [->] (1,-1) .. controls (1,-1.3) and (0,-1.7) .. (0,-2)
 .. controls (0,-1.7) and (-1,-1.3) ..  (-1,-1);
\draw (-1,-1) .. controls (-1,-.7) and (0,-.3) ..  (0,0);
\draw [->] (1.2,.1) arc (0:360:1.2 and .7);
\fill [white, rounded corners] (-.25,-1.65) rectangle (.25,-1.35);
\draw (0,-1.5) node{\scriptsize $W$};
\fill [white, rounded corners] (1.05,-.9) rectangle (1.45,-1.3);
\draw (1.25,-1.1) node{\scriptsize $\propP^{\prime}_{\veca}$};
\fill [white, rounded corners] (-1.05,-.9) rectangle (-1.45,-1.3);
\draw (-1.25,-1.1) node{\scriptsize $\propP_{\veca}$};
\end{tikzpicture}

\caption{Proof of Theorem~\ref{thmdefTheta}}\label{figthmdefTheta}
\end{figure}
This proves that $\langle \propP_{\veca}, \propP_b, \propP_c \rangle_{\!C_2(\rats)\,}$ is independent of $\propP_{\veca}$ when $a$ is fixed. Similarly, it is independent of $\propP_b$ and $\propP_c$ when $b$ and $c$ are fixed. Thus, $\langle \propP_{\veca}, \propP_b, \propP_c \rangle_{\!C_2(\rats)\,}$ is a rational function on the connected set of triples $(a,b,c)$ of distinct point of $S^2$. It is easy to see that this function is continuous. Therefore, it is constant.

Let us similarly prove that $\int_{C_2(\rats)} \omega_{\veca} \wedge \omega_b \wedge  \omega_c$ is independent of the propagating forms $\omega_{\veca}$, $\omega_b$, and $\omega_c$. 
Let $\omega^{\prime}_{\veca}$ be a propagating form of $(C_2(\rats),\tau)$.
Lemma~\ref{lemetactwo} implies the existence of forms $\eta$ and $\eta_A$ such that $\omega^{\prime}_{\veca}-\omega_{\veca}=d\eta$ and the restriction of $\eta$ to $\partial C_2(\rats)$ is $p_{\tau}^{\ast}(\eta_A)$. So we have
\begin{equation*}\begin{array}{ll}\int_{C_2(\rats)} \omega^{\prime}_{\veca} \wedge \omega_b \wedge  \omega_c - \int_{C_2(\rats)} \omega_{\veca} \wedge \omega_b \wedge  \omega_c & =\int_{C_2(\rats)}d(\eta \wedge \omega_b \wedge  \omega_c)\\&=\int_{\partial C_2(\rats)}\eta \wedge \omega_b \wedge  \omega_c\\
        &=\int_{\partial C_2(\rats)}p_{\tau}^{\ast}(\eta_A \wedge \omega_B \wedge  \omega_C)=0
       \end{array}\end{equation*}
since any $5$-form on $S^2$ vanishes.
Thus, $\int_{C_2(\rats)} \omega_{\veca} \wedge \omega_b \wedge  \omega_c$ is independent of the propagating forms $\omega_{\veca}$, $\omega_b$, and $\omega_c$. Now, we can choose the propagating forms $\omega_{\veca}$, $\omega_b$, and $\omega_c$ supported in tiny neighborhoods of $\propP_{\veca}$,
$\propP_b$, and $\propP_c$ and respectively Poincar\'e dual to $\propP_{\veca}$,
$\propP_b$, and $\propP_c$. So the intersection of the three supports is a very small neighborhood of $\propP_{\veca} \cap \propP_b \cap \propP_c$, and we can easily see that $\int_{C_2(\rats)} \omega_{\veca} \wedge \omega_b \wedge  \omega_c=\langle \propP_{\veca}, \propP_b, \propP_c \rangle_{\!C_2(\rats)\,}$. See Section~\ref{secexistformtransv}, Section~\ref{secderhamcohom}, and  Lemma~\ref{lemconstrformalphadual} in particular, for more details.
\eop

In particular, $\Theta(\rats,\tau)$ is equal to $\int_{C_2(\rats)} \omega^3$ for any propagating chain $\omega$ of $(C_2(\rats),\tau)$. Since such a propagating chain represents the linking number, $\Theta(\rats,\tau)$ can be thought of as the {\em cube of the linking number with respect to $\tau$.\/}

When $\tau$ varies continuously, $\Theta(\rats,\tau)$ varies continuously in $\QQ$. So $\Theta(\rats,\tau)$ is an invariant of the homotopy class of $\tau$.

\begin{remark} Define a \emph{combing} of $\crats$ to be a section of $\ST \crats$ that coincides with $\taust(v)$ outside $\ballb_{\rats}$, for some unit vector $v$ of $\RR^3$. Let $X$ be a combing of $\crats$.
Define a \emph{propagating chain of $(C_2(\rats),X)$} to be a propagating chain of $C_2(\rats)$ that intersects $\ST \crats$ along the image of $X$. Define $\tilde{\Theta}(\rats,X)$ to be the algebraic intersection of a propagating chain of $(C_2(\rats),X)$, a propagating chain of $(C_2(\rats),-X)$, and any other propagating chain of $C_2(\rats)$. It is easily proved in \cite[Theorem 2.1]{lesHC} that $\tilde{\Theta}(\rats,X)$ depends only on $\rats$ and on the homotopy class of $X$ among combings. In particular, $\Theta(\rats,\tau)=\tilde{\Theta}(\rats,\tau(v))$ depends only on the homotopy class of the combing $\tau(v)$ of $\ST \crats$, for some unit vector $v$ of $\RR^3$.
Further properties of the invariant $\tilde{\Theta}(\rats,.)$ of combings are studied in \cite{lescomb}.
I found an explicit formula \cite[Theorem 3.8]{lesHC} for the invariant $\tilde{\Theta}(\rats,.)$ from a Heegaard diagram of $\rats$. I computed it directly using the above definition of $\tilde{\Theta}(\rats,.)$ together with 
propagators associated to Morse functions constructed with Greg Kuperberg in \cite{lesHC}. 
\end{remark}

\begin{example}
\label{exathetasthreetaust}
Using (disjoint!) propagators $p_{S^2}^{-1}(\veca)$, $p_{S^2}^{-1}(b)$, $p_{S^2}^{-1}(c)$ associated to three distinct points $a$, $b$, and $c$ of $\RR^3$, as in Example~\ref{exapropsthree}, it is clear that
\begin{equation*}\Theta(S^3,\taust)=\bigl\langle p_{S^2}^{-1}(\veca),p_{S^2}^{-1}(b),p_{S^2}^{-1}(c) \bigr\rangle_{\!C_2(S^3)\,}\end{equation*} is zero.
\end{example}

\section{Parallelizations of \texorpdfstring{$3$}{3}-manifolds and Pontrjagin classes}
\label{secannpont}

In this subsection, $M$ denotes a smooth, compact oriented $3$-manifold with possible boundary $\partial M$.
Recall that a well-known theorem reproved in Section~\ref{secproofpar} ensures that such a $3$-manifold is parallelizable.

Let $GL^+(\RR^3)$ denote the group of orientation-preserving linear isomorphisms of $\RR^3$.
Let $\left[(M,\partial M),(GL^+(\RR^3),1)\right]_m$ denote the set of (continuous) maps \begin{equation*}g\colon (M, \partial M) \longrightarrow (GL^+(\RR^3),1)\end{equation*}
from $M$ to $GL^+(\RR^3)$ that send $\partial M$ to the identity element $1$ of $GL^+(\RR^3)$.\\
Let $\left[(M,\partial M),(GL^+(\RR^3),1)\right]$ denote the group of homotopy classes of such maps, with the group structure induced by the multiplication of maps, using the multiplication in $GL^+(\RR^3)$.
For a map $g$ in $\left[(M,\partial M),(GL^+(\RR^3),1)\right]_m$, define
\begin{equation*}\begin{array}{llll} 
\psi_{\RR}(g)\colon &M \times \RR^3 &\longrightarrow  &M \times \RR^3\\
&(x,y) & \mapsto &(x,g(x)(y)).\end{array}\end{equation*}
Let $\tau_M \colon M \times \RR^3 \rightarrow TM$ be a parallelization of $M$.
Then any parallelization $\tau$ of $M$ that coincides with $\tau_M$ on $\partial M$ may be written as
\begin{equation*}\tau = \tau_M \circ \psi_{\RR}(g)\end{equation*} for some $g \in \left[(M,\partial M),(GL^+(\RR^3),1)\right]_m$.

Thus, fixing $\tau_M$ identifies the set of homotopy classes of parallelizations of $M$ fixed on $\partial M$ with 
the group
$\left[(M,\partial M),(GL^+(\RR^3),1)\right]$.
Since $GL^+(\RR^3)$ deformation retracts onto the group $SO(3)$ of orientation-preserving linear isometries of $\RR^3$, the group $\left[(M,\partial M),(GL^+(\RR^3),1)\right]$ is isomorphic to $\left[(M,\partial M),(SO(3),1)\right]$.

\begin{definition}
\label{defrho}
We regard $S^3$ as $B^3/\partial B^3$, where $B^3$ is the standard unit ball of $\RR^3$ viewed as the quotient of $\left[0,1\right] \times S^2$ where all points of $\{0\} \times S^2$ are identified with each other.
Let $\chi_{\pi}\colon \left[0,1\right] \to \left[0,2\pi\right]$ be an increasing smooth bijection whose derivatives
vanish at $0$ and $1$ such that $\chi_{\pi}(1 -\theta)=2\pi -\chi_{\pi}(\theta)$ for any $\theta \in \left[0,1\right]$.
Let $\rhomap \colon B^3 \rightarrow SO(3)$ be the map that sends $(\theta \in \left[0,1\right], v \in S^2)$ to the rotation $\rhorot(\chi_{\pi}(\theta);v)$ with axis directed by $v$ and with angle $\chi_{\pi}(\theta)$\index[N]{rzho@$\rhomap \colon B^3 \rightarrow SO(3)$}. 
\end{definition}
This map $\rhomap$ induces the double covering map $\tilde{\rhomap} \colon S^3 \rightarrow SO(3)$,\index[N]{rzhotilde@$\tilde{\rhomap} \colon S^3 \rightarrow SO(3)$} which orients $SO(3)$ and which allows one to deduce the first three homotopy groups of $SO(3)$ from those of $S^3$.  The first three homotopy groups of $SO(3)$ are $\pi_1(SO(3))=\ZZ/2\ZZ$, $\pi_2(SO(3))=0$, and $\pi_3(SO(3))=\ZZ\left[\tilde{\rhomap}\right]$. For $v\in S^2$, $\pi_1(SO(3))$ is generated by the class of the loop that maps $\exp(i \theta) \in S^1$ to the rotation $\rhorot(\theta;v)$. See Section~\ref{sechomotopygroups} and Theorem~\ref{thmlongseqhomotopy} in particular.

Note that a map $g$ from $(M,\partial M)$ to $(SO(3),1)$ has a degree $\deg(g)$. The \emph{degree} $\deg(g)$ is the differential degree at a regular value (different from $1$) of $g$. It can also be defined homologically, by $H_3(g)\left[M,\partial M\right] =\deg(g)\left[SO(3),1\right]$.

We prove the following theorem, for which we claim no originality, in Chapter~\ref{chapfram}, as a direct consequence of Definition~\ref{deffirstpontcrats}, Lemmas~\ref{lempreptrivun}, \ref{lemdeg}, \ref{lemdegtwo}, and Propositions~\ref{proppontdef}, \ref{proppontdeg}, and \ref{proppont}.

\begin{theorem}
\label{thmpone}
For any compact connected oriented $3$-manifold $M$, the group \begin{equation*}\left[(M,\partial M),(SO(3),1)\right]\end{equation*} is abelian and
the degree \begin{equation*}\deg \colon \left[(M,\partial M),(SO(3),1)\right] \longrightarrow \ZZ\end{equation*} is a group homomorphism,
which induces an isomorphism 
\begin{equation*}\deg \colon \left[(M,\partial M),(SO(3),1)\right] \otimes_{\ZZ}\QQ \longrightarrow \QQ.\end{equation*}
When $\partial M$ is equal to $\emptyset$ or $S^2$, there exists a canonical map $p_1$ from the set of homotopy classes of parallelizations of $M$---that coincide with $\taust$ near $S^2$ if $\partial M=S^2$---to $\ZZ$  such that
\begin{itemize}
 \item we have $p_1((\taust)\vert_{B^3})=0$ and
 \item for any map $g$ in $\left[(M,\partial M),(SO(3),1)\right]_m$ and for any trivialization $\tau$ of $TM$, we have
\begin{equation*}p_1( \tau \circ \psi_{\RR}(g))-p_1(\tau)=2\deg(g).\end{equation*}
\end{itemize}
\end{theorem}

Definition~\ref{deffirstpontcrats} of the map $p_1$ involves relative Pontrjagin classes. When $\partial M=\emptyset$, the map $p_1$ coincides with the map $h$ studied by Friedrich Hirzebruch in \cite[\S 3.1]{hirzebruchEM}, and by Robion Kirby and Paul Melvin in \cite{km} under the name of {\em Hirzebruch defect.\/}

Since $\left[(M,\partial M),(SO(3),1)\right]$ is abelian,
the set of parallelizations of $M$ that are fixed on $\partial M$ is an affine space with translation group $\left[(M,\partial M),(SO(3),1)\right]$.

Recall the map $\rhomap \colon B^3 \rightarrow SO(3)$ from Definition~\ref{defrho}. 
Let $M$ be an oriented connected $3$-manifold with possible boundary. For a ball $B^3$ embedded in $M$, let $\rhomap_M(B^3) \in \left[(M,\partial M),(SO(3),1)\right]_m$ be a smooth map that coincides with $\rhomap$ on $B^3$ and that maps the complement of $B^3$ to the identity element of $SO(3)$. The homotopy class of $\rhomap_M(B^3)$ is well-defined.
\begin{lemma}
\label{lemdegrho}
 We have $\deg(\rhomap_M(B^3))=2$.
\end{lemma}
\bp Exercise.
\eop

\section{Defining a \texorpdfstring{$\QQ$}{Q}-sphere invariant from \texorpdfstring{$\Theta$}{Theta}}
\label{secdefThetam}

Recall that an asymptotic rational homology $\RR^3$ is a pair $(\crats,\tau)$, where $\crats$ is decomposed as the union over $\left]1,2\right] \times S^2$ of a rational homology ball $\ballb_{\rats}$ and the complement $\mathring{B}_{1,\infty}\setminus \{\infty\}$ of the unit ball of $\RR^3$, and $\tau$ is  an asymptotically standard parallelization of $\crats$.
In this subsection, we prove the following proposition.

\begin{proposition}
\label{propThetap}
Let $(\crats,\tau)$ be an asymptotic rational homology $\RR^3$.
For any map $g$ in $[(\ballb_{\rats},\ballb_{\rats} \cap \mathring{B}_{1,\infty}),(SO(3),1)]_m$ trivially extended to $\crats$, we have
\begin{equation*}\Theta(\rats, \tau \circ \psi_{\RR}(g))-\Theta(\rats, \tau)=\frac12\deg(g).\end{equation*}
\end{proposition}

Theorem~\ref{thmpone} allows us to derive the following corollary from Proposition~\ref{propThetap}.
\begin{corollary}
\label{corThetap}
 Set $\Theta(\rats)=\Theta(\rats,\tau) - \frac1{4}{p_1(\tau)}$. Then $\Theta(\rats)$ is an invariant of $\QQ$-spheres.
\end{corollary}
\eopwobp

Since $p_1(\taust)=0$, Example~\ref{exathetasthreetaust} shows that $\Theta(S^3)=0$. We will prove $\Theta(-\rats)=-\Theta(\rats)$ for any $\QQ$-sphere $\rats$ in Proposition~\ref{propThetaorrev}.

More properties of $\Theta$ will appear later in this book. We will first view this invariant as the degree one part of a much more general invariant $\Zinvuf$ or $\zinvuf$ (introduced in Theorem~\ref{thmfstconst} and in Corollary~\ref{cordefzinvuf}, respectively) in Corollary~\ref{corthetazone}. The multiplicativity of $\Zinvuf$ under connected sum stated in Theorem~\ref{thmconnsum} implies that $\Theta$ is additive under connected sum. The invariant $\Theta$ will be identified with $6\lambda_{CW}$, where $\lambda_{CW}$ is the Walker generalization of the Casson invariant to $\QQ$-spheres, in Section~\ref{seccasson}. See Theorem~\ref{thmThetaeqlambda}. (The Casson--Walker invariant $\lambda_{CW}$ is normalized as $\frac{1}{2}\lambda_W$ for rational homology spheres, where $\lambda_W$ is the Walker normalisation in \cite{wal}.) The equality $\Theta=6\lambda_{CW}$ will be obtained as a consequence of a \emph{universality property} of $\Zinvuf$ with respect to a theory of finite type invariants. We will also present a direct proof of a surgery formula satisfied by $\Theta$ in Section~\ref{secproofsurcas}.

Let us begin the proof of Proposition~\ref{propThetap} with the following lemma.

\begin{lemma}
\label{lemtheprimhom}
 The variation $\Theta(\rats, \tau \circ \psi_{\RR}(g))-\Theta(\rats, \tau)$ in Proposition~\ref{propThetap} is independent of $\tau$. Set $\Theta^{\prime}(g)=\Theta(\rats, \tau \circ \psi_{\RR}(g))-\Theta(\rats, \tau)$. Then $\Theta^{\prime}$ is a homomorphism from $\left[(\ballb_{\rats},\ballb_{\rats} \cap \mathring{B}_{1,\infty}),(SO(3),1)\right]$ to $\QQ$.
\end{lemma}
\bp
For $d=a$, $b$ or $c$, the propagator $\propP_d$ of $(C_2(\rats),\tau)$ of Theorem~\ref{thmdefTheta} can be assumed to be a product $\left[-1,0\right] \times p_{\tau }\vert_{\ST \ballb_{\rats}}^{-1}(d)$ on a collar $\left[-1,0\right] \times \ST \ballb_{\rats}$ of $\ST \ballb_{\rats}$ in $C_2(\rats)$. Since $H_3(\left[-1,0\right] \times \ST \ballb_{\rats};\QQ)=0$, the cycle \begin{equation*} \Bigl(\partial \bigl(\left[-1,0\right] \times p_{\tau }\vert_{\ST \ballb_{\rats}}^{-1}(d)\bigr) \setminus \bigl(\{0\} \times p_{\tau }\vert_{\ST \ballb_{\rats}}^{-1}(d)\bigr) \Bigr) \cup \Bigl(\{0\} \times p_{\tau \circ \psi_{\RR}(g) }\vert_{\ST \ballb_{\rats}}^{-1}(d)\Bigr)\end{equation*} bounds a chain $G_d$. 

The chains $G_{\veca}$, $G_b$, and $G_c$ can be assumed to be transverse. Construct the propagator $\propP_d(g)$ of $(C_2(\rats),\tau \circ \psi_{\RR}(g))$ from $\propP_d$ by replacing $\left[-1,0\right] \times p_{\tau }\vert_{\ST \ballb_{\rats}}^{-1}(d)$ with $G_d$ on $\left[-1,0\right] \times \ST \ballb_{\rats}$. Then we have
\begin{equation*}\Theta(\rats, \tau \circ \psi_{\RR}(g))-\Theta(\rats, \tau)=\langle G_{\veca},G_b,G_c\rangle_{\!\left[-1,0\right] \times \ST \ballb_{\rats}\,}.\end{equation*}
Using $\tau$ to identify $\ST \ballb_{\rats}$ with $\ballb_{\rats} \times S^2$, the cycle $\partial G_d$ may be written as
\begin{equation*} \partial G_d = \Bigl(\left(\partial (\left[-1,0\right]  \times \ballb_{\rats}) \setminus (\{0\} \times \ballb_{\rats})\right) \times \{d\}\Bigr) \cup \Bigl(\{0\} \times\left( \cup_{m \in \ballb_{\rats}}(m,g(m)(d))\right)\Bigr).\end{equation*} Therefore, $\Theta(\rats, \tau \circ \psi_{\RR}(g))-\Theta(\rats, \tau)$ is independent of $\tau$. Then it is easy to see that $\Theta^{\prime}$ is a homomorphism from $\left[(\ballb_{\rats},\partial \ballb_{\rats}),(SO(3),1)\right]$ to $\QQ$.
\eop

So, Theorem~\ref{thmpone} and Lemma~\ref{lemdegrho} reduce the proof of Proposition~\ref{propThetap} to the proof that $\Theta^{\prime}(\rhomap_{\rats}(B^3))=1$. It is easy to see that $\Theta^{\prime}(\rhomap_{\rats}(B^3))=\Theta^{\prime}(\rhomap)$. Thus, we have reduced the proof of Proposition~\ref{propThetap} to the proof of the following lemma.

\begin{lemma}
\label{lemtheprimrho}
We have $\Theta^{\prime}(\rhomap)=1$.
\end{lemma}
To prove this lemma, we will compute the expression $\langle G_{\veca},G_{-\veca},G_c\rangle_{\!\left[-1,0\right] \times \ST B^3\,}$ from the above proof of Lemma~\ref{lemtheprimhom} with explicit $G_{\veca}$ and $G_{-\veca}$ in $\left[-1,0\right] \times \ST B^3$. We will construct the chain $G_{\veca}$ from a chain $G(\veca)$ of $B^3 \times S^2$
described in Lemma~\ref{lemchainanobeta}.

Again, we regard $B^3$ as the quotient of $\left[0,1\right] \times S^2$ where all points of $\{0\} \times S^2$ are identified with each other.
We first prove the following lemma:

\begin{lemma}
\label{lemL} Recall the map $\rhomap \colon B^3 \to SO(3)$ from Definition~\ref{defrho}.
Let $a \in S^2$. The point $(-\veca)$ is a regular value of the map
\begin{equation*}\begin{array}{llll}\rhomap_{\veca} \colon &B^3 &\rightarrow &S^2\\
   & m &\mapsto & \rhomap(m)(\veca)
  \end{array}
\end{equation*}
and its preimage (cooriented by $S^2$ via $\rhomap_{\veca}$) is the knot $L_{a}=-\{\frac12\} \times S_{\veca}$, where $S_{\veca}$ is the circle of $S^2$ of vectors orthogonal to $a$, oriented as the boundary of the hemisphere that contains $a$.
\end{lemma}
\bp We prove the lemma when $a$ is the North Pole $\upvec$.
It is easy to see that $\rhomap_{\upvec}^{-1}(-\upvec)=L_{\upvec}$ up to orientation.
\begin{center}
\begin{tikzpicture} \useasboundingbox (-1,-1.2) rectangle (1,1.2);
\draw (0,0) circle (1) (-1.1,-.05) node[above]{\tiny $x$} (0,1) node[above]{\tiny $\upvec$} (0,-.95) node[below]{\tiny $-\upvec$} (0,-.2) node[below]{\small $L_{\upvec}$};
\draw [-<] (-1,0) .. controls (-.95,-.1) and (-.3,-.25) .. (0,-.25);
\draw (0,-.25) .. controls (.3,-.25) and  (.95,-.1) .. (1,0);
\draw [out=29,in=151,dashed] (-1,0) to (1,0); 
\draw [->] (-1,0) -- (-1,-.3) node[below]{\tiny $1$};
\draw [->] (0,-1) -- (.6,-1) node[right]{\small $v_1$};
\draw [->] (-1,0) -- (-1.3,0) node[left]{\tiny $2$};
\fill (0,1) circle (1.5pt) (0,-1) circle (1.5pt);
\end{tikzpicture}
\end{center}
Let $(\frac12,x) \in L_{\upvec}$. Let $C$ denote the great circle of $S^2$ that contains ${\upvec}$ and $x$. Orient $C$ from $x$ towards $(-\upvec)$.
When $m$ moves along $\{\frac12\} \times C$ from $(\frac12,x) $ towards $(\frac12, -\upvec)$, $\rhomap(m)(\upvec)$ moves from $(-\upvec)$ along $C$ following the orientation of $C$, too.
Let $v_1$ denote the oriented unit tangent vector to $C$ at $(-\upvec)$.
 In our figure, $x$ is on the left, $C$ is oriented counterclockwise, and $v_1$ points to the right. So $S^2$ is oriented at $(-\upvec)$ by $v_1$ and by the tangent vector $v_2$ at $(-\upvec)$ towards us. In order to move $\rhorot(\theta;v)(\upvec)$ in the $v_2$ direction, one increases $\theta$,
so $L_{\upvec}$ is cooriented and oriented as in the figure.
\eop

\begin{lemma}
\label{lemchainanobeta}
Let $\veca \in S^2$. Let $m \in B^3$. Recall the notation from Lemma~\ref{lemL} above where $\rhomap(m)(\veca)=\rhomap_{\veca}(m)$.
When $m \notin L_{\veca}$, let $[\veca,\rhomap_{\veca}(m)]$ denote the unique geodesic arc of $S^2$ with length $\left(\ell \in \left[0, \pi\right[\right)$ from $a$ to $\rhomap_{\veca}(m)$. For $t \in \left[0,1\right]$, let $X_t(m) \in \left[\veca,\rhomap_{\veca}(m)\right]$ be such that
the length of $\left[\veca,X_t(m)\right]$ is $t\ell$. 
Let $G_h(\veca)$ be the closure of 
\begin{equation*}\left(\cup_{t \in \left[0,1\right], m\in (B^3 \setminus L_{\veca})} \left(m,{X_t}(m)\right)\right)\end{equation*} in $B^3 \times S^2$.
Let $D_{\veca}$ be the disk image of $\left[0,\frac12\right]\times (-S_{\veca})$ bounded by $L_{\veca}$ in $B^3$. Set \begin{equation*}G(\veca)=G_h(\veca) +D_{\veca} \times S^2.\end{equation*}
Then $G(\veca)$ is a chain of $B^3 \times S^2$ such that \begin{equation*}\partial G(\veca)=-(B^3\times \veca) +\cup_{m\in B^3}\bigl(m,\rhomap_{\veca}(m)\bigr).\end{equation*}
\end{lemma}
\bp The map $X_t$ is well-defined on $(B^3 \setminus L_{\veca})$. We have $X_0(m)=\veca$ and $X_1(m)=\rhomap_{\veca}(m)$. Let us show how the definition of $X_t$ extends smoothly to the manifold $\blowup{B^3}{L_{\veca}}$ obtained from $B^3$ by blowing up $L_{\veca}$. The map $\rhomap_{\veca}$ maps the normal bundle to $L_{\veca}$ to a disk of $S^2$ around $(-\veca)$
by an orientation-preserving diffeomorphism on every fiber (near the origin). In particular, $\rhomap_{\veca}$ induces a map from the unit normal bundle $\sph N(L_{\veca})$ to $L_{\veca}$ in $\rats$ to the unit normal bundle to $(-\veca)$ in $S^2$. This map preserves the orientation of the fibers.
Then for an element $y$ of $\sph N(L_{\veca})$, define $X_t(y)$ as before on the half great circle $\left[a,-\veca\right]_{\rhomap_{\veca}(-y)}$ from $a$ to $(-\veca)$ that is tangent to $\rhomap_{\veca}(-y)$ at $(-\veca)$. (So $\rhomap_{\veca}(-y)$ is an outward normal to $\left[a,-\veca\right]_{\rhomap_{\veca}(-y)}$ at $(-\veca)$.)
This continuously extends the definition of $X_t$. So we have
\begin{equation*}G_h(\veca)=\cup_{t \in \left[0,1\right], m\in \smallblowup{B^3}{L_{\veca}}}\bigl(p_{B^3}(m),{X_t}(m)\bigr),\end{equation*}
and 
\begin{multline*}\partial G_h(\veca) = -(B^3\times a) +\cup_{m\in B^3}(m,\rhomap_{\veca}(m))\\ +
\cup_{t \in \left[0,1\right], m\in -\partial \smallblowup{B^3}{L_{\veca}}}\bigl(p_{B^3}(m),{X_t}(m)\bigr),\end{multline*}
where $(-\partial \blowup{S^3}{L_{\veca}})$ is oriented as the boundary of a tubular neighborhood of $L_{\veca}$. For any $x \in L_{\veca}$, the sphere $S^2$ is covered with degree $(-1)$ by the image of $\left[0,1\right] \times \sph N_x(L_{\veca})$,
where the fiber $\sph N_x(L_{\veca})$ of $\sph N(L_{\veca})$ is oriented as the boundary of a disk in the fiber $N_x(L_{\veca})$ of the normal bundle.
We may therefore write the last summand as $(-L_{\veca} \times S^2)$.
\eop

\bpo{Proof of Lemma~\ref{lemtheprimrho}}
We use the notation of the proof of Lemma~\ref{lemtheprimhom} and construct an explicit $G_{\veca}$ in $\left[-1,0\right] \times \ST B^3 \stackrel{\taust} {=}\left[-1,0\right] \times B^3 \times S^2$.
Let $\iota$ be the endomorphism of $\ST B^3$ over $B^3$ that maps a unit vector to the opposite one. Recall the chain $G(\veca)$ of Lemma~\ref{lemchainanobeta}.
Set \begin{multline*}G_{\veca}=\Bigl(\left[-1,-2/3\right]\times B^3\times a\Bigr) + \Bigl(\{-2/3\}\times G(\veca)\Bigr) \\+ \Bigl(\left[-2/3,0\right]\times \cup_{m\in B^3}\bigl(m,\rhomap_{\veca}(m)\bigr)\Bigr)\end{multline*}
and
\begin{multline*}
G_{-\veca}=\Bigl(\left[-1,-1/3\right]\times B^3 \times (-\veca)\Bigr) + \Bigl(\{-1/3\}\times \iota\bigl(G(\veca)\bigr) \Bigr)\\+ \Bigl(\left[-1/3,0\right]\times \cup_{m\in B^3}\bigl(m,\rhomap(m)(-\veca)\bigr)\Bigr).
\end{multline*}
Then we have \begin{multline*}G_{\veca} \cap G_{-\veca}=\bigl(\left[-2/3,-1/3\right] \times L_{\veca} \times (-\veca)\bigr) + \bigl(\{-2/3\} \times D_{\veca} \times (-\veca) \bigr)\\-\bigl(\{-1/3\} \times \cup_{m\in D_{\veca}}(m,\rhomap_{\veca}(m))\bigr). \end{multline*}
Finally, according to the proof of Lemma~\ref{lemtheprimhom}, $\Theta^{\prime}(\rhomap)$ is the algebraic intersection of $G_{\veca} \cap G_{-\veca}$ with $\propP_c(\rhomap)$ in $C_2(\rats)$. This intersection coincides with the algebraic intersection of $G_{\veca} \cap G_{-\veca}$ with any propagator of $C_2(\rats)$, according to Lemma~\ref{lemlkprop}. 
Therefore, we have \begin{equation*} \Theta^{\prime}(\rhomap)=\langle \propP_{c},G_{\veca} \cap G_{-\veca} \rangle_{\!\left[-1,0\right] \times B^3\times S^2\,}=-\deg_{c}(\rhomap_{\veca} \colon D_{\veca} \rightarrow S^2)\end{equation*}
for any regular value $c\neq-\veca$ of $\rhomap_{\veca}\vert_{D_{\veca}}$.
Since the image of the quotient $D_{\veca}$ of $\left[0,\frac12\right]\times (-S_{\veca})$ under $\rhomap_{\veca}$ covers the sphere with degree $(-1)$, we get $\Theta^{\prime}(\rhomap)=1$.
\eop

So, Proposition~\ref{propThetap} is proved.

\chapter{Parallelizations of \texorpdfstring{$3$}{3}-manifolds and Pontrjagin classes}
\label{chapfram}

In this chapter, we fix a smooth oriented connected $3$-manifold $M$ with possible boundary, and we study its parallelizations.
In particular, we prove Theorem~\ref{thmpone} in Sections~\ref{secmsothree}, \ref{sechomdegpithree}, and \ref{secpfpone}.
We will use this theorem in our general constructions of link invariants in $3$-manifolds as in the definition of $\Theta$ in Section~\ref{secdefThetam}.
This chapter also describes other properties of Pontrjagin classes, which will be used in the fourth part of this book in some
universality proofs. Section~\ref{secmorestrucpar} describes the structure of the space of parallelizations of oriented $3$-manifolds, more precisely. It is not used in other parts of the book.

\section{\texorpdfstring{$\left[(M,\partial M),(SO(3),1)\right]$}{[M,SO(3)]} is an abelian group.}
\label{secmsothree}

Again, we regard $S^3$ as $B^3/\partial B^3$, and $B^3$ as the quotient of $\left[0,1\right] \times S^2$ where all points of $\{0\} \times S^2$ are identified with each other.
Recall the map $\rhomap \colon B^3 \rightarrow SO(3)$ of Definition~\ref{defrho}, which maps $(\theta \in\left[0,1\right], v \in S^2)$ to the rotation $\rhorot(\chi_{\pi}(\theta);v)$ with axis directed by $v$ and with angle $\chi_{\pi}(\theta)$. 
Also recall that the group structure of $\left[(M,\partial M),(SO(3),1)\right]$ is induced by the multiplication of maps, using the multiplication of $SO(3)$.

Any $g \in \left[(M,\partial M),(SO(3),1)\right]_m$
 induces a map \begin{equation*}H_1(g;\ZZ) \colon H_1(M,\partial M;\ZZ)\longrightarrow (H_1(SO(3),1)=\ZZ/2\ZZ).\end{equation*}
Since \begin{equation*}\begin{array}{ll}H_1(M,\partial M;\ZZ/2\ZZ)&=H_1(M,\partial M;\ZZ)/2H_1(M,\partial M;\ZZ)\\&=H_1(M,\partial M;\ZZ)\otimes_{\ZZ}\ZZ/2\ZZ,\end{array}\end{equation*}
we have 
\begin{equation*}\begin{array}{ll}\Hom\bigl(H_1(M,\partial M;\ZZ),\ZZ/2\ZZ\bigr)&=\Hom\bigl(H_1(M,\partial M;\ZZ/2\ZZ),\ZZ/2\ZZ\bigr)\\&=H^1(M,\partial M;\ZZ/2\ZZ),\end{array}\end{equation*} and the image of $H_1(g;\ZZ)$ in $H^1(M,\partial M;\ZZ/2\ZZ)$ under the above isomorphisms is denoted by
$H^1(g;\ZZ/2\ZZ)$. (Equivalently, $H^1(g;\ZZ/2\ZZ)$ denotes the image of the generator of $H^1(SO(3),1;\ZZ/2\ZZ)=\ZZ/2\ZZ$ under $H^1(g;\ZZ/2\ZZ)$ in $H^1(M,\partial M;\ZZ/2\ZZ)$.)

For $v\in S^2$, let $\gamma_v$ denote the loop that maps $\exp(i \theta) \in S^1$ to the rotation $\rhorot(\theta;v)$ with axis directed by $v$ and with angle $\theta$. Let $\RR P^2_S$ denote the projective plane embedded in $SO(3)$ consisting of the rotations of $SO(3)$ of angle $\pi$. Note that $\gamma_v(S^1)=\rhomap\left(\left[0,1\right] \times \{v\}\right)$ and $\RR P^2_S=\rhomap(\{1/2\} \times S^2)$ intersect once transversally. Thus, we have $H_1(SO(3);\ZZ)=\ZZ/2\ZZ\left[\gamma_v\right]$ and $H_2(SO(3);\ZZ/2\ZZ)=\ZZ/2\ZZ\left[\RR P^2_S\right]$. 

\begin{lemma}
\label{lemfirsthomom}
The map
\begin{equation*}\begin{array}{llll}H^1(.;\ZZ/2\ZZ)  \colon & \left[(M,\partial M),(SO(3),1)\right] &\to &H^1(M,\partial M;\ZZ/2\ZZ)\\
& \left[g\right] & \mapsto & H^1(g;\ZZ/2\ZZ)
\end{array}\end{equation*}
is a group homomorphism.
\end{lemma}
\bp Let $f$ and $g$ be two elements of $\left[(M,\partial M),(SO(3),1)\right]_m$. In order to prove
$H^1(fg;\ZZ/2\ZZ)=H^1(f;\ZZ/2\ZZ)+H^1(g;\ZZ/2\ZZ)$, it suffices to prove that, for any path $\gamma \colon \left[0,1\right] \to M$ whose image $\gamma\left(\partial \left[0,1\right]\right)$ is in $\partial M$, we have \begin{equation*}H_1(fg;\ZZ/2\ZZ)\left(\left[\gamma\right]\right)=H_1(f;\ZZ/2\ZZ)\left(\left[\gamma\right]\right)+H_1(g;\ZZ/2\ZZ)\left(\left[\gamma\right]\right),\end{equation*}
where $H_1(f;\ZZ/2\ZZ)\left(\left[\gamma\right]\right)=\langle f \circ \gamma, \RR P^2_S \rangle_{\!SO(3)\,}\left[\gamma_v\right]$ for the mod $2$ algebraic intersection $\langle f \circ \gamma, \RR P^2_S \rangle_{\!SO(3)\,}$.
Both sides of this equation depend only on the homotopy classes of $f \circ \gamma$ and $g\circ \gamma$. Therefore, we may assume that the supports (closures of the preimages of $SO(3) \setminus \{1\}$) of $f \circ \gamma$ and $g\circ \gamma$ are disjoint and easily conclude.
\eop

\begin{lemma}
\label{lempreptrivun}
Let $M$ be an oriented connected $3$-manifold with possible boundary. Recall that $\rhomap_M(B^3) \in \left[(M,\partial M),(SO(3),1)\right]_m$ is a map that coincides with the map $\rhomap$ of Definition~\ref{defrho} on a ball $B^3$ embedded in $M$ and that maps the complement of $B^3$ to the unit of $SO(3)$.
\begin{enumerate}
\item Any homotopy class of a map $g$ from $(M,\partial M)$ to $(SO(3),1)$, such that $H^1(g;\ZZ/2\ZZ)$ 
is trivial, belongs to the subgroup \begin{equation*}\Bigl<\left[\rhomap_M(B^3)\right]\Bigr>\end{equation*} of 
$\left[(M,\partial M),(SO(3),1)\right]$
generated by $\left[\rhomap_M(B^3)\right]$. 
\item For any $\left[g\right] \in \left[(M,\partial M),(SO(3),1)\right]$, we have
$\left[g\right]^2 \in \bigl<\left[\rhomap_M(B^3)\right]\bigr>.$
\item The group $\left[(M,\partial M),(SO(3),1)\right]$ is abelian.
\end{enumerate}
\end{lemma}
\bp Let $g \in \left[(M,\partial M),(SO(3),1)\right]_m$. Assume that $H^1(g;\ZZ/2\ZZ)$ is trivial.
Choose a cell decomposition of $M$ relative to its boundary,
with only one three-cell, no zero-cell if $\partial M \neq \emptyset$, one zero-cell if $\partial M = \emptyset$,  one-cells, and two-cells. See \cite[Chapter 6, Section 3]{hirsch}. Then after a homotopy relative to $\partial M$, we may assume that $g$ maps the one-skeleton of $M$ to $1$. 
Next, since $\pi_2(SO(3)) = 0$, we may assume that $g$ maps the two-skeleton of $M$ to $1$, and therefore that $g$ maps the exterior of some $3$-ball to $1$. Now $g$ becomes a map from $B^3/\partial B^3=S^3$ to $SO(3)$, and its homotopy class is $k\left[\tilde{\rhomap}\right]$ in $\pi_3(SO(3))=\ZZ\left[\tilde{\rhomap}\right]$.
Therefore, the map $g$ is homotopic to $\rhomap_M(B^3)^k$. This proves the first assertion.

Since $H^1(g^2;\ZZ/2\ZZ)=2H^1(g;\ZZ/2\ZZ)$ is trivial, the second assertion follows.

For the third assertion, first note that $\left[\rhomap_M(B^3)\right]$ belongs to the center of
$\left[(M,\partial M),(SO(3),1)\right]$ because it can be supported in a small ball disjoint 
from the support (preimage of $SO(3) \setminus \{1\}$) of a representative of any other element. Thus, according to the second assertion, any square is in the center. 
In particular, if $f$ and $g$ are elements of $\left[(M,\partial M),(SO(3),1)\right]$, we have
\begin{equation*}(gf)^2=(fg)^2=(f^{-1}f^2g^2f)(f^{-1}g^{-1}fg),\end{equation*} 
where the first factor is equal to $f^2g^2=g^2f^2$. Exchanging $f$ and $g$ yields
$f^{-1}g^{-1}fg=g^{-1}f^{-1}gf$. Therefore, the commutator, which is a power of $\left[\rhomap_M(B^3)\right]$, thanks to Lemma~\ref{lemfirsthomom} and to the first assertion, has a vanishing square, and thus a vanishing degree. So it is trivial.
\eop

\section{Any oriented \texorpdfstring{$3$}{3}-manifold is parallelizable.}
\label{secproofpar}

In this subsection, we prove the following standard theorem. The spirit of our proof is the same as the Kirby proof in \cite[p.46]{Kir}. But instead of assuming familiarity with the obstruction theory described by Norman Steenrod in \cite[Part III]{St}, we use this proof to introduce this theory.

\begin{theorem}[Stiefel]
\label{thmpar}
Any orientable $3$-manifold is parallelizable.
\end{theorem}

\begin{lemma} Let $M$ be an oriented $3$-manifold.
Let $S$ be a closed surface, orientable or not, immersed in $M$.
 Then the restriction to $S$ of the tangent bundle $TM$ to $M$ is trivializable.
\end{lemma}
\bp Let us first prove that this bundle is independent of the immersion. It is the direct sum of the tangent bundle to the surface and of its normal
one-dimensional bundle. This normal bundle is trivial when $S$ is orientable Otherwise, its unit bundle is the $2$-fold orientation cover of the surface. (The orientation cover of $S$ is its $2$-fold orientable cover that is trivial over annuli embedded in the surface). Since any surface $S$ can be immersed in $\RR^3$, the restriction $TM\vert_{S}$ is the pull-back of the trivial bundle of $\RR^3$ by such an immersion, and it is trivial.
\eop

Using Stiefel-Whitney classes, the proof of Theorem~\ref{thmpar} goes as follows.

\bpo{Quick proof of Theorem~\ref{thmpar}}
Let $M$ be an orientable smooth $3$-manifold, equipped with a smooth triangulation.\footnote{A theorem of Henry Whitehead proved in the Munkres book \cite{Mu} ensures the existence of such a triangulation.}
By definition, the {\em first Stiefel-Whitney class\/} $w_1(TM) \in H^1(M;\ZZ/2\ZZ=\pi_0(GL(\RR^3)))$ viewed as a map from
$\pi_1(M)$ to $\ZZ/2\ZZ$ maps the class of a loop $c$ embedded in $M$ to $0$ if $TM\vert_{c}$ is orientable, and to $1$ otherwise.
It is the \emph{obstruction} to the existence of a trivialization of $TM$ over the one-skeleton of $M$.
Since $M$ is orientable, the first Stiefel-Whitney class $w_1(TM)$ vanishes. So $TM$ can be trivialized over the one-skeleton of $M$.
The {\em second Stiefel-Whitney class\/} $w_2(TM) \in H^2(M;\ZZ/2\ZZ=\pi_1(GL^+(\RR^3)))$ viewed as a map from $H_2(M;\ZZ/2\ZZ)$ to $\ZZ/2\ZZ$ maps the class of a connected closed surface $S$ to $0$ if $TM\vert_{S}$ is trivializable, and to $1$ otherwise.
The second Stiefel-Whitney class $w_2(TM)$ is the obstruction to the existence of a trivialization of $TM$ over the two-skeleton of $M$ when $w_1(TM)=0$.
According to the above lemma, we have $w_2(TM)=0$, and $TM$ can be trivialized over the two-skeleton of $M$.
Then since $\pi_2(GL^+(\RR^3))=0$, any parallelization over the two-skeleton of $M$ can be extended as a parallelization of $M$.
\eop

We detail the involved arguments without mentioning Stiefel-Whitney classes, (in fact, by almost defining $w_2(TM)$) below. The elementary proof below can serve as an introduction to the obstruction theory used above.

\bpo{Elementary proof of Theorem~\ref{thmpar}}
Let $M$ be an oriented $3$-manifold. Choose a triangulation of $M$. For any cell $c$ of the triangulation, define an arbitrary trivialization $\tau_c \colon c \times \RR^3 \rightarrow TM\vert_{c}$ such that $\tau_c$ induces the orientation of $M$.
This defines a trivialization $\tau^{(0)}\colon M^{(0)} \times \RR^3 \rightarrow TM\vert_{M^{(0)}}$ of $M$ over the $0$-skeleton $M^{(0)}$ of $M$, which is the set of $0$-dimensional cells of the triangulation. Let $C_k(M)$ be the set of $k$-cells of the triangulation. Equip every cell with an arbitrary orientation.
Let  $e \in C_1(M)$ be an edge of the triangulation. On $\partial e$, the trivialization $\tau^{(0)}$ may be written as $\tau^{(0)}= \tau_e \circ \psi_{\RR}(g_e)$, for a map $g_e \colon \partial e \rightarrow GL^+(\RR^3)$. Since $GL^+(\RR^3)$ is connected, the map $g_e$ extends to $e$, and $\tau^{(1)}= \tau_e  \circ \psi_{\RR}(g_e)$ extends $\tau^{(0)}$ to $e$. Doing so for all the edges extends $\tau^{(0)}$ to a trivialization $\tau^{(1)}$ of the one-skeleton $M^{(1)}$ of $M$, which is the union of the edges of the triangulation.

Let $t$ be a triangle of the triangulation. There is a map $g_t \colon \partial t \rightarrow GL^+(\RR^3)$ such that $\tau^{(1)}=  \tau_t \circ \psi_{\RR}(g_t)$ on $\partial t$. Let $E(t,\tau^{(1)})$ be the homotopy class of $g_t$ in $(\pi_1(GL^+(\RR^3))=\pi_1(SO(3))=\ZZ/2\ZZ)$.
Then $E(.,\tau^{(1)}) \colon C_2(M) \rightarrow \ZZ/2\ZZ$ extends linearly to a cochain, which is independent of the $\tau_t$. When $E(.,\tau^{(1)})=0$, $\tau^{(1)}$ extends to
a trivialization $\tau^{(2)}$ over the two-skeleton of $M$, as before.
Since $\pi_2(GL^+(\RR^3))=0$, the trivialization $\tau^{(2)}$ can next be extended over the three-skeleton of $M$. So it extends over $M$.

Let us now study the obstruction cochain $E(.,\tau^{(1)})$ whose vanishing guarantees the existence of a parallelization of $M$.

Choose a map  $d(e) \colon (e,\partial e) \rightarrow (GL^+(\RR^3),1)$ for every edge $e$ of the triangulation.
Let $\tau^{(1)\prime}$ be the trivialization of $TM\vert_{M^{(1)}}$ obtained as above by changing the map $g_e$ associated to $e$ to $d(e)g_e$.
Define the cochain 
\begin{equation*}D(\tau^{(1)},\tau^{(1)\prime}) \colon (\ZZ/2\ZZ)^{C_1(M)} \rightarrow \ZZ/2\ZZ\end{equation*} that maps $e$ to the homotopy class of $d(e)$.
Then $(E(.,\tau^{(1)\prime})-E(.,\tau^{(1)}))$ is the coboundary of $D(\tau^{(1)},\tau^{(1)\prime})$. (See Section~\ref{sechomology}, before Theorem~\ref{thmpoincdual}.)

Let us prove that $E(.,\tau^{(1)})$ is a cocycle.
Consider a $3$-simplex $T$. The trivialization $\tau^{(0)}$ extends to $T$. Assume that $\tau_T$
is an extension of $\tau^{(0)}$ to $T$, that $\tau_t$ is the restriction of $\tau_T$ to $t$ for any face $t$ of $T$, and that the above $\tau^{(1)\prime}$ coincides with $\tau_T$ on the edges of $\partial T$. Then we have $E(.,\tau^{(1)\prime})(\partial T)=0$. Since any coboundary
 maps $\partial T$ to $0$, we also have $E(.,\tau^{(1)})(\partial T)=0$.

Now, it suffices to prove that the cohomology class of $E(.,\tau^{(1)})$ (which is equal to $w_2(TM)$) vanishes in order to prove the existence of an extension $\tau^{(1)\prime}$ of $\tau^{(0)}$ on $M^{(1)}$ that extends on $M$.

Since $H^2(M;\ZZ/2\ZZ)=\Hom(H_2(M;\ZZ/2\ZZ);\ZZ/2\ZZ)$, it suffices to prove that $E(.,\tau^{(1)})$ maps any $2$-dimensional $\ZZ/2\ZZ$-cycle $C$ to $0$. 

We represent the class of such a cycle $C$ by a closed surface $S$, orientable or not, as follows.

Let $N(M^{(0)})$ bs a disjoint union of balls around the vertices and let $N(M^{(1)})$ be a small regular neighborhood of $M^{(1)}$ in $M$ as in Figure~\ref{figneighconstsurf}.
We assume that
$N(M^{(1)}) \cap (M \setminus N(M^{(0)}))$ is a disjoint union, over the edges $e$, of solid cylinders $B_e$ identified with $\left]0,1\right[ \times D^2$. The core $\left]0,1\right[ \times \{0\}$ of $B_e= \left]0,1\right[ \times D^2$ is a connected part of the interior of the edge $e$. (As in Figure~\ref{figneighconstsurf}, the neighborhood $N(M^{(1)})$ is thinner than $N(M^{(0)})$.)

\bfig
\centering
\begin{tikzpicture}[scale=.5] \useasboundingbox (-5,-1.2) rectangle (5,2);
\begin{scope}[xshift=-2.5cm]
\fill [lightgray] (-1.8,-1) rectangle (1.8,1.9);
\draw [thin] (-30:1.4) -- (90:1.4) -- (-150:1.4) -- (-30:1.4) (-1.8,-.7) -- (1.8,-.7) 
(90:1.4) -- (90:1.9);
\draw [thin] (-150:1.4) -- (-30:1.4);
\fill [white, rounded corners] (-1.9,2) -- (-.3,2) -- (-.3,1.4) -- (-1.32,-.4) -- (-1.9,-.4) -- (-1.9,2);
\fill [white, rounded corners] (1.9,2) -- (.3,2) -- (.3,1.4) -- (1.32,-.4) -- (1.9,-.4) -- (1.9,2);
\fill [white, rounded corners] (-30:1.1) -- (90:1.1) -- (-150:1.1) -- (-30:1.1);
\fill [white, rounded corners] (90:1.1) -- (-150:1.1) -- (-30:1.1) -- (90:1.1);
\draw [very thin, <-] (.25,1.6) -- (.6,1.6) node[right]{\scriptsize $N(M^{(1)})$};
\draw [very thin, <-] (-1.65,-.65) -- (-1.65,-.3) node[above]{\scriptsize $M^{(1)}$};
\draw [very thin, <-] (-.97,-.1) -- (-1.45,-.1);
\end{scope}
\begin{scope}[xshift=2.5cm]
\fill [lightgray] (-1.8,-1) rectangle (1.8,1.9);
\fill [white, rounded corners] (-1.9,2) -- (-.3,2) -- (-.3,1.4) -- (-1.32,-.4) -- (-1.9,-.4) -- (-1.9,2);
\fill [white, rounded corners] (1.9,2) -- (.3,2) -- (.3,1.4) -- (1.32,-.4) -- (1.9,-.4) -- (1.9,2);
\fill [white, rounded corners] (-30:1.1) -- (90:1.1) -- (-150:1.1) -- (-30:1.1);
\fill [gray] (90:1.4) circle (.45);
\fill [gray] (-30:1.4) circle (.45);
\fill [gray] (-150:1.4) circle (.45);
\draw [thin] (-30:1.4) -- (90:1.4) -- (-150:1.4) -- (-30:1.4) (-1.8,-.7) -- (1.8,-.7) 
(90:1.4) -- (90:1.9);
\draw [very thin, <-] (.25,1.45) -- (.6,1.6) node[right]{\scriptsize $N(M^{(0)})$};
\draw [very thin, <-] (1.2,-.35) -- (1.2,1.3);
\end{scope}
\end{tikzpicture}

\caption{The neighborhoods $N(M^{(1)})$ and $N(M^{(0)})$}\label{figneighconstsurf}
\end{figure}

Construct $S$ in the complement of $N(M^{(0)}) \cup N(M^{(1)})$ as the intersection of the support of $C$ with this complement.
Then the closure of $S$ meets the part $\left[0,1\right] \times S^1$ of every $\overline{B_e}$ as an even number of parallel intervals from $\{0\} \times S^1$ to $\{1\} \times S^1$. Complete $S$ in $M  \setminus N(M^{(0)})$ by connecting the intervals pairwise in $\overline{B_e}$ by disjoint bands. After this operation, the boundary of the closure of $S$ is a disjoint union of circles in the boundary of $N(M^{(0)})$. Glue disjoint disks of $N(M^{(0)})$ along these circles to finish the construction of $S$.

Extend $\tau^{(0)}$ to $N(M^{(0)})$, assume that $\tau^{(1)}$ coincides with this extension over $M^{(1)}\cap N(M^{(0)})$, and extend $\tau^{(1)}$ to $N(M^{(1)})$.
The bundle $TM\vert_{S}$ is trivial, and we may choose a trivialization $\tau_S$ of $TM$ over $S$ that coincides with our extension of $\tau^{(0)}$ over $N(M^{(0)})$, over $S \cap N(M^{(0)})$. We have a cell decomposition of $(S,S \cap N(M^{(0)}))$ with only 
$1$-cells and $2$-cells, for which the $2$-cells of $S$ are in one-to-one canonical correspondence with the $2$-cells of $C$, and one-cells correspond bijectively to bands connecting two-cells in the cylinders $B_e$. These one-cells are equipped with the trivialization of $TM$ induced by $\tau^{(1)}$.
Then we can define $2$-dimensional cochains $E_S(.,\tau^{(1)})$ and $E_S(.,\tau_S)$ as before, with respect to this cellular decomposition of $S$. The cochain $(E_S(.,\tau^{(1)})-E_S(.,\tau_S))$
is again a coboundary, and $E_S(.,\tau_S)=0$. So $E_S(C,\tau^{(1)})=0$, and we also have $E(C,\tau^{(1)})=0$ since $E(C,\tau^{(1)})=E_S(C,\tau^{(1)})$.
\eop

Theorem~\ref{thmpar} has the following immediate corollary.
\begin{proposition}
\label{propexistparasymptriv}
Any punctured oriented $3$-manifold $\crats$ as in Definition~\ref{defparasyst} can be equipped with an asymptotically standard parallelization.
 \end{proposition}
\bp The oriented manifold $\rats$ admits a parallelization $\tau_0 \colon \rats \times \RR^3 \to T\rats$.
Over $\mathring{B}_{1,\infty} \setminus \{\infty\}$, $\taust =\tau_0 \circ \psi_{\RR}(g)$ for a map $g \colon \mathring{B}_{1,\infty} \setminus \{\infty\} \rightarrow GL^+(\RR^3)$.
For $r \in \left[1,2\right]$,
let $\mathring{B}_{r,\infty}$ (resp. ${B}_{r,\infty}$) be the complement in $S^3$ of the closed (resp. open) ball ${B}(r)$ of radius $r$ in $\RR^3$. Since $\pi_2(GL^+(\RR^3))=\{0\}$, the restriction of $g$ to ${B}_{7/4,\infty} \setminus \mathring{B}_{2,\infty}$
extends to a map of $\mathring{B}_{1,\infty} \setminus \mathring{B}_{2,\infty}$ that maps $\mathring{B}_{1,\infty} \setminus \mathring{B}_{5/3,\infty}$ to $1$. After smoothing, we get a smooth map $\tilde{g}\colon \mathring{B}_{1,\infty} \setminus \{\infty\} \rightarrow GL^+(\RR^3)$ that coincides with $g$ on $\mathring{B}_{2,\infty}$ and that maps  $\mathring{B}_{1,\infty} \setminus \mathring{B}_{3/2,\infty}$ to $1$. Extend $\tilde{g}$ to $\crats$ so that it maps $\crats \setminus \mathring{B}_{3/2,\infty} $ to $1$. So $\tau_0 \circ \psi_{\RR}(\tilde{g})$ is an asymptotically standard parallelization as wanted.
\eop

\section{The homomorphism induced by the degree}
\label{sechomdegpithree}

In this section, $M$ is a compact connected oriented $3$-manifold, with or without boundary.
Let $S$ be a closed surface, orientable or not, embedded in the interior of our manifold $M$. Let $\tau$ be a parallelization of our $3$-manifold $M$.
We define a twist $g(S,\tau) \in \left[(M,\partial M),(SO(3),1)\right]_m$ below.

The surface $S$ has a tubular neighborhood $N(S)$, which is a $[-1,1]$-bundle over $S$. This bundle admits (orientation-preserving) bundle charts $\phi \colon [-1,1]\times D \to N(S)$ for disks $D$ of $S$, such that the changes of coordinates restrict to the fibers as $\pm \id$, where $\id$ denotes the identity map. Let $\phi(t,s)$ be in the image of such a chart.
Let $N_s=T_{\phi(0,s)}\phi\left([-1,1] \times s\right)$ be the tangent vector to the fiber $\phi\left([-1,1] \times s\right)$ at $\phi(0,s)$.
Define $g(S,\tau)\bigl(\phi(t,s)\bigr)$ to be the rotation with angle $\pi(t+1)$ and with axis
\begin{equation*}p_2\Bigl(\tau^{-1}(N_s)=\bigl(\phi(0,s),p_2(\tau^{-1}(N_s))\bigr)\Bigr).\end{equation*} Since this rotation coincides with the rotation with opposite axis and with opposite angle $\pi(1-t)$,
this provides a consistent definition of $g(S,\tau)\vert_{N(S)}$. Extend $g(S,\tau)\vert_{N(S)}$ to $M$ so that 
\begin{equation*}g(S,\tau)\colon (M, \partial M) \longrightarrow (SO(3),1)\end{equation*} maps $M\setminus N(S)$ to $1$.

The homotopy class of $g(S,\tau)$ depends only on the homotopy class of $\tau$ and on the isotopy class of $S$. When $M=B^3$, when $\tau$ is the standard parallelization of $\RR^3$, and when $\frac12 S^2$ denotes the sphere $\frac12 \partial B^3$ inside $B^3$, the homotopy class of $g(\frac12 S^2,\tau)$ coincides with the homotopy class of $\rhomap$.

We will see later (Proposition~\ref{propgStau}) that the homotopy class of $g(S,\tau)$ depends only on the Euler characteristic $\chi(S)$ of $S$ and on the class of $S$ in $H_2(M;\ZZ/2\ZZ)$. Thus, we will simply denote $g(S,\tau)$ by $g(S)$. We will also see (Corollary~\ref{corgStautwo}) that any element of $\left[(M,\partial M),(SO(3),1)\right]$ can be represented by some $g(S)$.

\begin{lemma}
\label{lemgStau}
The morphism $H^1(g(S,\tau);\ZZ/2\ZZ)$ maps the generator of \begin{equation*}H^1(SO(3);\ZZ/2\ZZ)\end{equation*} to the mod $2$ intersection with $S$ in \begin{equation*}\Hom(H_1(M,\partial M;\ZZ/2\ZZ),\ZZ/2\ZZ)=H^1(M,\partial M;\ZZ/2\ZZ).\end{equation*}
Thus, the morphism $H^1(.;\ZZ/2\ZZ)\colon \left[(M,\partial M),(SO(3),1)\right] \rightarrow H^1(M,\partial M;\ZZ/2\ZZ)$ is onto.
\end{lemma}
\bp The first assertion is obvious. 

The second one follows since $H^1(M,\partial M;\ZZ/2\ZZ)$ is the Poincar\'e dual of $H_2(M;\ZZ/2\ZZ)$ and since any element of $H_2(M;\ZZ/2\ZZ)$ is the class of a closed surface.
\eop

\begin{lemma}
\label{lemdeg}
The degree is a group homomorphism \begin{equation*}\deg \colon \left[(M,\partial M),(SO(3),1)\right] \longrightarrow \ZZ\end{equation*}
and $\deg(\rhomap_M(B^3)^k)=2k$.
\end{lemma}
\bp
It is easy to see that $\deg(fg)=\deg(f)+\deg(g)$ when $f$ or $g$
is a power of $\left[\rhomap_M(B^3)\right]$.

Let us prove that $\deg(f^2)=2\deg(f)$ for any $f$.
According to Lemma~\ref{lemgStau}, there is an unoriented embedded surface $S_f$ of the interior of $C$ such that $H^1\left(f;\ZZ/2\ZZ\right)=H^1\left(g(S_f,\tau);\ZZ/2\ZZ\right)$ for some trivialization $\tau$ of $TM$. According to Lemmas~\ref{lemfirsthomom} and \ref{lempreptrivun}, $fg(S_f,\tau)^{-1}$ is homotopic to some power of $\rhomap_M(B^3)$. So, it suffices to prove that the degree of $g^2$ is $2\deg(g)$
for $g=g(S_f,\tau)$. This can be done easily, by noticing that $g^2$ is homotopic to $g(S^{(2)}_f,\tau)$, where $S^{(2)}_f$ is the boundary of the tubular neighborhood of $S_f$.
In general, we have
\begin{equation*}\deg(fg)=\frac{1}2\deg((fg)^2)=\frac{1}2\deg(f^2g^2)=\frac{1}2\left( \deg(f^2)+\deg(g^2)\right),\end{equation*} and the lemma is proved.
\eop

Lemmas~\ref{lempreptrivun} and \ref{lemdeg} imply the following lemma.

\begin{lemma}
\label{lemdegtwo}
The degree induces an isomorphism 
\begin{equation*}\deg \colon \left[(M,\partial M),(SO(3),1)\right] \otimes_{\ZZ}\QQ \longrightarrow \QQ.\end{equation*}
Any group homomorphism $\psi \colon \left[(M,\partial M),(SO(3),1)\right] \longrightarrow \QQ$ may be expressed as
\begin{equation*}\frac12\psi(\rhomap_M(B^3))\deg.\end{equation*}
\end{lemma}
\eopwobp

Recall that $M$ is a compact connected oriented $3$-manifold in this section.
So $H^3(M,\partial M;\pi_3(SO(3)))$ is canonically isomorphic to $\ZZ$.

\begin{proposition}
\label{propstructwo} Let $P \colon H^1(M,\partial M;\ZZ/2\ZZ) \to H_2(M;\ZZ/2\ZZ)$ be the Poincar\'e duality isomorphism.

Let
$i \colon H^3(M,\partial M;\pi_3(SO(3))) \to \left[(M,\partial M),(SO(3),1)\right]$ be the group morphism that sends the canonical generator $\left[M,\partial M\right]$ of $H^3(M,\partial M;\pi_3(SO(3)))$ to $\left[\rhomap_M(B^3)\right]$.
Then the sequence 
\begin{multline*}
0 \rightarrow H^3(M,\partial M;\pi_3(SO(3))) \hflv{i} \left[(M,\partial M),(SO(3),1)\right] \\ \vspace{5pt}
\\
\hflv{P\circ H^1(\left[.\right];\ZZ/2\ZZ)} H_2(M;\ZZ/2\ZZ) \rightarrow 0,\end{multline*}
is a canonical exact sequence.
\end{proposition}
\bp The proposition
is a consequence of Lemmas~\ref{lempreptrivun}, \ref{lemgStau}, and \ref{lemdeg}.
\eop

\section{On the groups \texorpdfstring{$SU(n)$}{SU(n)}}
\label{secsun}

Let $\KK= \RR$ or $\CC$. Let $n \in \NN$. The stabilization maps induced
by the inclusions
\begin{equation*}\begin{array}{llll}i\colon & GL(\KK^n) & \longrightarrow & GL(\KK \oplus \KK^n)\\
& g & \mapsto & i(g)\colon (x,y) \mapsto (x,g(y))\end{array}\end{equation*}
are denoted by $i$. We represent elements of $GL(\KK^n)$ by matrices whose columns contain the coordinates of the images of the basis elements with respect to the standard basis of $\KK^n$.
View $S^3$ as the unit sphere of $\CC^2$. So, its elements are the pairs $(z_1,z_2)$ of complex numbers such that $|z_1|^2 + |z_2|^2=1$.
The group $SU(2)$ may be identified with $S^3$ by the homeomorphisms
\begin{equation*}\begin{array}{lllllllll} m^{\CC}_r \colon & S^3 &\rightarrow & SU(2)&\;\;\mbox{and} \;\;\overline{m}^{\CC}_r \colon & S^3 &\rightarrow & SU(2)\\ 
 &(z_1,z_2) &\mapsto & \left[\begin{array}{cc} z_1&-\overline{z}_2\\z_2& \overline{z}_1\end{array} \right]& &(z_1,z_2) &\mapsto & \left[\begin{array}{cc} \overline{z}_1&\overline{z}_2\\ -z_2& z_1\end{array} \right]. \end{array}\end{equation*}
So the first nontrivial homotopy group of $SU(2)$ is
$\pi_3(SU(2))=\ZZ\left[\overline{m}^{\CC}_r\right]$, where $\left[\overline{m}^{\CC}_r\right]=-[{m}^{\CC}_r]$ and $\overline{m}^{\CC}_r$ is a group homomorphism (it induces the group structure of $S^3$).
The long exact sequence associated to the fibration \begin{equation*} SU(n-1) \hookfl{i} SU(n) \rightarrow S^{2n-1},\end{equation*}
described in Theorem~\ref{thmlongseqhomotopy},
shows that $i^n_{\ast}\colon \pi_j(SU(2)) \longrightarrow \pi_j(SU(n+2))$ is an isomorphism for $j \leq 4$ and $n \geq 0$.
In particular, the group $\pi_j(SU(4))$ is trivial for $j\leq 2$ and we have 
\begin{equation*}\pi_3(SU(4))=\ZZ\left[i^2_{\ast}(\overline{m}^{\CC}_r)\right],\end{equation*}
where $i^2_{\ast}(\overline{m}^{\CC}_r)$ 
is the following map
\begin{equation*}\begin{array}{llll} i^2_{\ast}(\overline{m}^{\CC}_r)\colon &(S^3 \subset \CC^2) & \longrightarrow & SU(4)\\
& (z_1,z_2) & \mapsto & \left[\begin{array}{cccc} 1&0&0&0\\0&1&0&0\\0&0&\overline{z}_1&\overline{z}_2\\0&0&-z_2& z_1\end{array} \right] \end{array}.\end{equation*}

\section{Definition of relative Pontrjagin numbers}
\label{secdefpont}
Let ${\manifm}_0$ and ${\manifm}_1$ be two compact connected oriented $3$-manifolds whose boundaries have collars identified by a diffeomorphism. Let $\tau_0\colon {\manifm}_0 \times \RR^3 \rightarrow T{\manifm}_0$ and $\tau_1\colon {\manifm}_1 \times \RR^3 \rightarrow T{\manifm}_1$ be two parallelizations (which respect the orientations) that agree on the collar neighborhoods of $\partial {\manifm}_0=\partial {\manifm}_1$. 
Then the {\em relative Pontrjagin number $p_1(\tau_0,\tau_1)$\/} is the Pontrjagin obstruction to extending the trivialization of $TW \otimes \CC$ induced by $\tau_0$ and $\tau_1$ across the interior of a signature $0$ cobordism $W$ from ${\manifm}_0$ to ${\manifm}_1$. Details follow.

Let ${\manifm}$ be a compact connected oriented $3$-manifold. A \emph{special complex trivialization}\index[T]{special!complex trivialization} of $T{\manifm}$ is a trivialization of $T{\manifm} \otimes \CC$ that is obtained from a trivialization $\tau_{\manifm} \colon {\manifm} \times \RR^3 \rightarrow T{\manifm}$ by composing $(\tau^{\CC}_{\manifm} = \tau_{\manifm} \otimes_{\RR} \CC) \colon {\manifm} \times \CC^3 \rightarrow T{\manifm} \otimes \CC$ by \begin{equation*}\begin{array}{llll} 
\psi(G)\colon &{\manifm} \times \CC^3 &\longrightarrow  &{\manifm} \times \CC^3\\
&(x,y) & \mapsto &\bigl(x,G(x)(y)\bigr)\end{array}\end{equation*}
for a map $G\colon ({\manifm},\partial \manifm) \rightarrow (SL(3,\CC),1)$, which is a map $G\colon {\manifm} \rightarrow SL(3,\CC)$ that maps $\partial \manifm$ to $1$. The definition and properties of relative Pontrjagin numbers $p_1(\tau_0,\tau_1)$, given with more details below, are valid for pairs $(\tau_0,\tau_1)$ of special complex trivializations.

The {\em signature\/} of a $4$-manifold is the signature of the intersection form on its $H_2(.;\RR)$ (i.e., the number of positive entries minus the number of negative entries in a diagonalized version of this form). It is well known that any closed oriented three-manifold bounds a compact oriented $4$-dimensional manifold. See \cite{rourke} for an elegant elementary proof. The signature of such a bounded $4$-manifold may be  changed arbitrarily by connected
sums with copies of $\CC P^2$ or $-\CC P^2$.
A {\em cobordism from ${\manifm}_0$ to ${\manifm}_1$\/} is a compact oriented $4$-dimensional manifold $W$ with ridges such that
\begin{equation*}\partial W=-{\manifm}_0\cup_{\partial {\manifm}_0 \sim 0 \times \partial {\manifm}_0}  (-\left[0,1\right] \times \partial {\manifm}_0)  \cup_{\partial {\manifm}_1 \sim 1 \times \partial {\manifm}_0} {\manifm}_1,\end{equation*} 
where $W$ is identified with an open subspace of one of the products $\left[0,1\right[ \times {\manifm}_0$ or $\left]0,1\right] \times {\manifm}_1$ near $\partial W$, as the following picture suggests.

\begin{center}
\begin{tikzpicture}[scale=.8]
\useasboundingbox (-3,-.7) rectangle (5.5,2.5);
\draw (0,0) -- (4,0) (4,2) -- (0,2) (2,1) node{\scriptsize $W^4$} (0,1) node[left]{\scriptsize $\{0\} \times {\manifm}_0={\manifm}_0 $} (4,1) node[right]{\scriptsize $\{1\} \times {\manifm}_1={\manifm}_1$};
\draw [very thick] (0,0) -- (0,2) (4,0) -- (4,2);
\draw [dashed] (1.2,0) -- (.8,-.4) -- (1.2,2);
\draw (.9,-.5) node[left]{\scriptsize $\left[0,1\right] \times (-\partial {\manifm}_0)$};
\draw (4,.1) node{\scriptsize $ \rightarrow $} (3.2,.1) node{\scriptsize $ \rightarrow $} (2.4,.1) node{\scriptsize $ \rightarrow $} (3.2,.1) node[above]{\scriptsize $ \vec{N}$};
\end{tikzpicture}
\end{center}

Let $W=W^4$ be such a cobordism from ${\manifm}_0$ to ${\manifm}_1$, with signature $0$.
Consider the complex $4$-bundle $TW \otimes \CC$ over $W$.
Let $\vec{N}$ be the tangent vector to $\left[0,1\right] \times \{\mbox{pt}\}$ over $\partial W$ (under the above identifications). Let 
$\tau(\tau_0,\tau_1)$ denote the trivialization of $TW \otimes \CC$ over $\partial W$ obtained by stabilizing either $\tau_0$ or $\tau_1$ into $\vec{N} \oplus \tau_0$ or $\vec{N} \oplus \tau_1$. Then the obstruction to extending this trivialization to $W$ is the relative first \emph{Pontrjagin class}\index[T]{Pontrjagin class!first}\begin{equation*}p_1\bigl(W;\tau(\tau_0,\tau_1)\bigr)\left[W,\partial W\right] \in \Bigl(H^4\bigl(W,\partial W;\ZZ=\pi_3(SU(4))\bigr)=\ZZ\left[W,\partial W\right]\Bigr)\end{equation*} of the trivialization. 

Now, we specify our sign conventions for this Pontrjagin class. They are the same as in \cite{milnorsta}. 
In particular, $p_1$ is the opposite of the complexified tangent bundle's second Chern class $c_2$. See \cite[p. 174]{milnorsta}. Let us describe these conventions. The determinant bundle of $TW$
is trivial because $W$ is oriented, and $\det(TW \otimes \CC)$ is also trivial.
Our parallelization $\tau(\tau_0,\tau_1)$ over $\partial W$ is special with respect to
the trivialization of $\det(TW \otimes \CC)$.
Equip ${\manifm}_0$ and ${\manifm}_1$ with Riemannian metrics that coincide near $\partial {\manifm}_0$. Equip $W$ with a Riemannian metric that coincides with the orthogonal product metric of one of the products $\left[0,1\right] \times {\manifm}_0$ or $\left[0,1\right] \times {\manifm}_1$ near $\partial W$. 
Equip $TW \otimes \CC$ with the associated Hermitian structure.
Up to homotopy, assume that  $\tau(\tau_0,\tau_1)$
is unitary with respect to the Hermitian structure of $TW \otimes \CC$ and the standard Hermitian form of $\CC^4$.
Since $\pi_i(SU(4))=\{0\}$ when $i<3$, the trivialization $\tau(\tau_0,\tau_1)$
extends to a special unitary trivialization $\tau$ outside the interior of a $4$-ball $B^4$
and defines 
\begin{equation*}\tau\colon S^3 \times \CC^4 \longrightarrow (TW \otimes \CC)\vert_{S^3}\end{equation*}
over the boundary $S^3=\partial B^4$ of this $4$-ball $B^4$.
Over this $4$-ball $B^4$, the bundle $TW \otimes \CC$ admits a special unitary trivialization
\begin{equation*}\tau_B\colon B^4 \times \CC^4 \longrightarrow (TW \otimes \CC)\vert_{B^4}.\end{equation*}
Then $\tau_B^{-1} \circ \tau(v \in S^3, w \in \CC^4)=(v, \phi(v)(w))$ for a map $\phi\colon S^3 \longrightarrow  SU(4)$ whose homotopy class may be written as
\begin{equation*}\left[\phi\right]=p_1\bigl(W;\tau(\tau_0,\tau_1)\bigr)\left[i^2_{\ast}(\overline{m}^{\CC}_r)\right] \in \pi_3(SU(4)),\end{equation*}
where $i^2_{\ast}(\overline{m}^{\CC}_r)$ was defined at the end of Section~\ref{secsun}.

Define $p_1(\tau_0,\tau_1)= p_1\bigl(W;\tau(\tau_0,\tau_1)\bigr)$.

\begin{proposition}
\label{proppontdef}
Let ${\manifm}_0$ and ${\manifm}_1$ be two compact connected oriented $3$-manifolds whose boundaries have collars identified by a diffeomorphism. Let $\tau_0 \colon {\manifm}_0 \times \CC^3 \rightarrow T{\manifm}_0 \otimes \CC  $ and $\tau_1\colon {\manifm}_1 \times \CC^3 \rightarrow T{\manifm}_1 \otimes \CC$ be two special complex trivializations (which
 respect the orientations) that agree on the collar neighborhoods of $\partial {\manifm}_0=\partial {\manifm}_1$. 
The (first) \indexT{Pontrjagin number} $p_1(\tau_0,\tau_1)$ is well defined by the above 
conditions.
\end{proposition}
\bp 
According to the Nokivov additivity theorem, if a closed $4$-manifold $Y$ is decomposed as $Y= Y^+ \cup_{X} Y^-$,
where $Y^+$ and $Y^-$ are two $4$-manifolds with boundary, embedded in $Y$, which intersect along a closed $3$-manifold $X$ (their common boundary, up to orientation), then we have  \begin{equation*}\signature(Y)=\signature(Y^+)+\signature(Y^-).\end{equation*} According to a theorem of Vladimir Rohlin (see \cite{Roh} or \cite[p. 18]{gm2}),
when $Y$ is a compact oriented $4$-manifold without boundary, we have 
\begin{equation*}p_1(Y)=3\,\signature(Y).\end{equation*}
We only need to prove that $p_1(\tau_0,\tau_1)$ is independent of the signature $0$ cobordism $W$. Let $W_E$
be a $4$-manifold of signature $0$ bounded by $(-\partial W)$. Then $W \cup_{\partial W} W_E$ is a $4$-dimensional manifold without boundary whose signature is \begin{equation*}\signature(W_E) + \signature(W)=0\end{equation*} by the Novikov additivity theorem.
According to the Rohlin theorem, the first Pontrjagin class of $W \cup_{\partial W} W_E$ is also zero. On the other hand,
this first Pontrjagin class is the sum of the relative first Pontrjagin classes of $W$ and $W_E$ with respect to $\tau(\tau_0,\tau_1)$. These two relative Pontrjagin classes are opposite. Therefore, the relative first Pontrjagin class of $W$ with respect to $\tau(\tau_0,\tau_1)$ does not depend on $W$.
\eop

Similarly, it is easy to prove the following proposition.
\begin{proposition}
\label{proppontdefsignon}
Under the above assumptions, except for the assumption on the signature of the cobordism $W$, we have
\begin{equation*}p_1(\tau_0,\tau_1)=p_1(W;\tau(\tau_0,\tau_1))-3\,\textrm{signature}(W).\end{equation*}
\end{proposition}
\eopwobp

\begin{remark} 
 When $\partial {\manifm}_1=\emptyset$ and ${\manifm}_0=\emptyset$, the map $p_1(=p_1(\tau(\emptyset),.))$ coincides with the map $h$ studied by Friedrich Hirzebruch in \cite[\S 3.1]{hirzebruchEM}, and by Robion Kirby and Paul Melvin in \cite{km} under the name of {\em Hirzebruch defect.\/} 
\end{remark}

\begin{definition}
\label{deffirstpontcrats}
 When $(\crats,\tau)$ is an asymptotic rational homology $\RR^3$, set 
\begin{equation*}p_1(\tau)=p_1\bigl((\taust)\vert_{B^3},\tau\vert_{\ballb_{\rats}}\bigr),\end{equation*}
with the notation of Proposition~\ref{proppontdef}.
\end{definition}

\begin{lemma} 
\label{lemrevor}
Let $(\crats,\partau)$ be an asymptotic rational homology $\RR^3$ as in Definition~\ref{defasyrathommRthree}. 
The parallelization $\partau \colon \crats \times \RR^3 \to T\crats$ induces the parallelization
$\overline{\partau} \colon (-\crats) \times \RR^3 \to T(-\crats)$ such that $\overline{\partau}(x,v)=-\partau(x,v)$.

Compose the orientation-preserving identification of a neighborhood of $\infty$ in $\rats$ with $\mathring{\ballb}_{1,\infty}$ by the (restriction of the) multiplication by $(-1)$ in $\RR^3 \cup \{\infty\}$ in order to get an orientation-preserving identification of a neighborhood of $\infty$ in $(-\rats)$ with $\mathring{\ballb}_{1,\infty}$.
Then \begin{equation*}\left(-\crats=\check{\overbrace{(-\rats)}},\overline{\partau}\right)\end{equation*} is an asymptotic rational homology $\RR^3$, and we have
$p_1(\overline{\partau})=-p_1(\partau)$.
\end{lemma}
\bp Use a signature $0$ cobordism $W$ from $\{0\} \times B^3$ to $\{1\} \times {\ballb}_{\rats}$ to compute $p_1(\partau)$. Extend the trivialization of $TW \otimes \CC$ on $\partial W$, which may be expressed as $\vec{N} \oplus \taust$ on $\{0\} \times B^3 \cup (-\left[0,1\right] \times \partial B^3)$, and $\vec{N} \oplus \partau$ on $\{1\} \times {\ballb}_{\rats}$, to a special trivialization on the complement of an open ball $\mathring{\ballb}^4$ in $W$.
Let $\overline{W}$ be the cobordism obtained from $W$ by reversing the orientation of $W$. Equip $\overline{W \setminus \mathring{\ballb}^4}$ with the trivialization obtained from the above trivialization by a composition by $\id_{\RR} \times (-\id)_{\RR^3}$. Then the changes of trivializations $\phi \colon \partial B^4 \to SU(4)$ and $\overline{\phi} \colon \partial \overline{B^4} \to SU(4)$ are obtained from one another by the orientation-preserving conjugation by $\id_{\RR} \times (-\id)_{\RR^3}$. Since $\partial \overline{B^4}$ and $\partial B^4$ have opposite orientations, we get the result.
\eop

Back to the invariant $\Theta$ defined in Corollary~\ref{corThetap}, we can now prove the following proposition.
\begin{proposition}
\label{propThetaorrev} For any $\QQ$-sphere $\rats$, we have $\Theta(-\rats)=-\Theta(\rats)$.
\end{proposition}
\bp Note that $C_2(-\rats)$ is naturally identified with $C_2(\rats)$, with the same orientation. Recall the orientation-reversing diffeomorphism $\iota$ of $C_2(\rats)$ that exchanges the two coordinates in $\crats^2 \setminus \diag\left(\crats^2\right)$.
If $\omega$ is a propagating form of $\bigl(C_2(\rats),\partau\bigr)$, then
$\iota^{\ast}(\omega)$ is a propagating form of $(C_2(-\rats),\overline{\partau})$. So we have
$\Theta(-\rats,\overline{\partau})=\int_{C_2(-\rats)} \iota^{\ast}(\omega^3)$.  This proves $\Theta(-\rats,\overline{\partau})=-\Theta(\rats,{\partau})$. Corollary~\ref{corThetap} and Lemma~\ref{lemrevor} yield the conclusion.
\eop

\section{On the groups \texorpdfstring{$SO(3)$}{SO(3)} and \texorpdfstring{$SO(4)$}{SO(4)}}

Let $\HH$ denote the vector space $\CC \oplus \CC j$. Set $k=ij$. The \emph{conjugate} of an element $(z_1+z_2j)$ of $\HH$ is
\begin{equation*}\overline{z_1+z_2j}=\overline{z_1}-z_2j.\end{equation*}

\begin{lemma}
The bilinear map that maps $(z_1+z_2j,z^{\prime}_1+z^{\prime}_2j)$ to $(z_1z^{\prime}_1-z_2\overline{z^{\prime}_2}) +(z_2\overline{z^{\prime}_1} +z_1z^{\prime}_2)j$
maps $(i,j)$ to $k$, $(j,k)$ to $i$, $(k,i)$ to $j$, $(j,i)$ to $(-k)$, $(k,j)$ to $(-i)$, $(i,k)$ to $(-j)$, $(i,i)$, $(j,j)$, and $(k,k)$ to $(-1)$, and $(z_1+z_2j,\overline{z_1+z_2j})$ to $|z_1|^2+|z_2|^2$.
It defines an associative product on $\HH$ such that this product and the addition make $\HH$ a field.
\end{lemma}
\bp Exercise.\eop

The noncommutative field $\HH$, which contains $\CC$ as a subfield, is the \emph{field of quaternions}. It is equipped with the scalar product $\langle.,.\rangle$ that makes $(1,i,j,k)$ an orthonormal basis of $\HH$.  The associated norm maps $(z_1+z_2j)$ to 
\begin{equation*}\sqrt{(z_1+z_2j)\overline{z_1+z_2j}}.\end{equation*} It is multiplicative. The unit sphere of $\HH$ is the sphere $S^3$. It is equipped with the group structure induced by the product of $\HH$. The elements of $\HH$ are the \indexT{quaternions}. The \emph{real part} of a quaternion $(z_1+z_2j)$ is the real part of $z_1$. The \emph{pure quaternions} are the quaternions with zero real part.

For $\KK= \RR$ or $\CC$ and $n \in \NN$, the $\KK$ (Euclidean or Hermitian) oriented vector space with the direct orthonormal basis $(v_1, \dots, v_n)$ is denoted by $\KK\langle v_1, \dots, v_n \rangle $.
The \emph{cross product} or \emph{vector product} of two elements $v$ and $w$ of $\RR^3=\RR\langle i , j , k \rangle$ is the element $v \times w$ of $\RR^3$ such that for any $x \in \RR^3$, $x \wedge v \wedge w =\langle x,v \times w \rangle i \wedge j \wedge k$ in $\bigwedge^3 \RR^3 = \RR$.

\begin{lemma}
\label{lemrotquat}
The product of two pure quaternions $v$ and $w$ is 
\begin{equation*}vw = -\langle v, w \rangle + v \times w.\end{equation*}
Every element of $S^3$ may be expressed as $\cos(\theta) + \sin(\theta) v$ for a unique $\theta \in\left[0,\pi\right]$ and a pure quaternion $v$ of norm $1$, which is unique when $\theta \notin \{0,\pi\}$. 
For such an element, the restriction to $\RR\langle i , j , k \rangle$ of the conjugation \begin{equation*}R(\theta,v) \colon w \mapsto \bigl(\cos(\theta) + \sin(\theta) v\bigr)w\overline{\bigl(\cos(\theta) + \sin(\theta) v\bigr)}\end{equation*}
is the rotation with axis directed by $v$ and with angle $2 \theta$.
\end{lemma}
\bp It is easy to check the first assertion. The conjugation $R(\theta,v)$ preserves the scalar product of $\HH$ and fixes $\RR \oplus \RR v$, pointwise. Therefore, it restricts to $\RR\langle i , j , k \rangle$ as an orthonormal transformation that fixes $v$.
Let $w$ be a pure quaternion orthogonal to $v$.
\begin{equation*}R(\theta,v)(w)=\bigl(\cos(\theta) + \sin(\theta) v\bigr)w\bigl(\cos(\theta) - \sin(\theta) v\bigr)\end{equation*}
is equal to
\begin{equation*}\begin{array}{ll}R(\theta,v)(w)&=\cos^2(\theta)w -\sin^2(\theta)vwv + \cos(\theta)\sin(\theta)(vw-wv)\\
&=\cos(2\theta)w + \sin(2\theta) v \times w.\end{array}\end{equation*}
\eop

\begin{lemma} The group morphism \begin{equation*}\begin{array}{llll} \tilde{\rhomap} \colon & S^3 &\rightarrow & SO(\RR\langle i , j , k \rangle)=SO(3)\\ 
& x &\mapsto &\left(w \mapsto x.w.\overline{x} \right)\end{array}\end{equation*}
is surjective, and its kernel is $\{-1,+1\}$. The morphism $\tilde{\rhomap}$ is a two-fold covering map, and
this definition of $\tilde{\rhomap}$ coincides with the previous one (after Definition~\ref{defrho}), up to homotopy.
\end{lemma}
\bp 
According to Lemma~\ref{lemrotquat}, the map
$\tilde{\rhomap}$ is surjective. Its kernel is the center of the group of unit quaternions, which is $\{-1,+1\}$. Thus, the map $\tilde{\rhomap}$ is a two-fold covering map.

This two-fold covering map clearly coincides with the previous one, up to homotopy and orientation, since both classes generate $\pi_3(SO(3))=\ZZ$. We take care of the orientation using the outward normal first convention to orient boundaries, as usual.
Consider the diffeomorphism 
\begin{equation*}\begin{array}{llll}\psi \colon &\left]0,\pi\right[\times S^2 &\to &S^3\setminus \{-1,1\}\\
&(\theta,v) & \mapsto & \cos(\theta) + \sin(\theta) v.\end{array}\end{equation*}
We study $\psi$ at $(\pi/2,i)$.
At $\bigl(\psi(\pi/2,i)=i\bigr)$, the space $\RR\langle i , j , k \rangle$ is oriented by the outward normal $o_i$ to $S^2$ followed by the orientation of $S^2$. The field
$\HH$ is oriented as $\RR \oplus \RR\langle i , j , k \rangle$. Since $o_i$ coincides with the outward normal to $S^3$ in $\HH$, the sphere $S^3$ is oriented by the opposite of the real part followed by the orientation of $S^2$. Since $\cos$ is an orientation-reversing diffeomorphism at $\pi/2$, $\psi$ preserves the orientation near $(\pi/2,i)$. So the diffeomorphism $\psi$ preserves the orientation everywhere, and
the two maps $\tilde{\rhomap}$ are homotopic.
\eop

The following two group morphisms from $S^3$ to $SO(4)$ induced by the multiplication in $\HH$
\begin{equation*}\begin{array}{llll} m_{\ell} \colon & S^3 &\rightarrow & (SO(\HH)=SO(4))\\ 
& x &\mapsto &(m_{\ell}(x) \colon v \mapsto x.v)\end{array}\end{equation*}

\begin{equation*}\begin{array}{llll} \overline{m_r} \colon & S^3 &\rightarrow & SO(\HH)\\ 
& y &\mapsto & (\overline{m_r}(y)\colon v \mapsto v.\overline{y})\end{array}\end{equation*}
together induce the surjective group morphism
\begin{equation*}\begin{array}{lll} 
S^3 \times S^3 &\rightarrow & SO(4)\\ 
(x,y) &\mapsto & (v \mapsto x.v.\overline{y}).\end{array}\end{equation*}
The kernel of this group morphism is $\{(-1,-1),(1,1)\}$. So this morphism is a two-fold covering map.
In particular, we have $\pi_3(SO(4))=\ZZ\left[m_{\ell}\right] \oplus \ZZ\left[\overline{m_r}\right].$
Define
\begin{equation*}\begin{array}{llll} {m_r} \colon & S^3 &\rightarrow & (SO(\HH)=SO(4))\\ 

& y &\mapsto & ({m_r}(y)\colon v \mapsto v.{y}).\end{array}\end{equation*}

\begin{lemma}
\label{lempitroissoquatrestr} In
$\pi_3(SO(4))=\ZZ\left[m_{\ell}\right] \oplus \ZZ\left[\overline{m_r}\right]$, we have
\begin{equation*}i_{\ast}\left(\left[\tilde{\rhomap}\right]\right)=\left[m_{\ell}\right] + \left[\overline{m}_r\right]=
\left[m_{\ell}\right] - \left[{m}_r\right].\end{equation*}
\end{lemma}
\bp
The $\pi_3$-product in $\pi_3(SO(4))$ coincides with the product induced by the group structure of $SO(4)$.
\eop

\begin{lemma} 
\label{lem237}
Recall that $m_r$
denotes the map from the unit sphere $S^3$ of $\HH$ to $SO(\HH)$ induced by the
right-multiplication. 
Denote the inclusions $SO(n) \subset SU(n)$ by $c$. Then we have 
\begin{equation*}c_{\ast}\bigl(\left[m_r\right]\bigr)=2\left[i^2_{\ast}(m^{\CC}_r)\right]\end{equation*} in $\pi_3(SU(4))$.
\end{lemma}
\bp
Let $\HH + I \HH$ denote the complexification of $\RR^4= \HH= \RR \langle 1,i,j,k \rangle$.
Here, we have $\CC=\RR \oplus I\RR$ and $I^2=-1$.
When $x \in \HH$ and $v \in S^3$, we have $c(m_r(v))(Ix)=Ixv$.
Let $\varepsilon=\pm 1$. Define
\begin{equation*} \CC^2(\varepsilon)=\CC\biggl\langle \frac{\sqrt{2}}{2}(1+ \varepsilon Ii),\frac{\sqrt{2}}{2}(j+ \varepsilon Ik)\biggr\rangle.\end{equation*}
Consider the quotient $\CC^4/\CC^2(\varepsilon)$.
In this quotient, we have $Ii=-\varepsilon 1$, $Ik =-\varepsilon j$, and, since $I^2=-1$,
$I1=\varepsilon i$ and $Ij=\varepsilon k$. Therefore, this quotient is isomorphic to $\HH$ as a real vector space with its complex structure $I= \varepsilon i$.
Then it is easy to see that $c(m_r(v))$ maps $\CC^2(\varepsilon)$ to $0$ in this quotient, for any $v \in S^3$. We get $c(m_r(v))(\CC^2(\varepsilon)) = \CC^2(\varepsilon)$.
Now, observe that $\HH + I \HH$ is the orthogonal sum of $\CC^2(-1)$ and $\CC^2(1)$.
In particular, $\CC^2(\varepsilon)$ is isomorphic to the quotient $\CC^4/\CC^2(-\varepsilon)$, which is isomorphic to $(\HH;I= -\varepsilon i)$, and 
$c(m_r)$ acts on it by the right multiplication. Therefore,
with respect to the orthonormal basis $\frac{\sqrt{2}}{2}(1-Ii, j-Ik, 1+Ii, j+Ik )$, $c(m_r(z_1+z_2j=x_1+y_1i+x_2j+y_2k))$ may be written as
\begin{equation*}c\bigl(m_r(x_1+y_1i+x_2j+y_2k)\bigr) =\left[\begin{array}{cccc} 
x_1+y_1I & -x_2+y_2I & 0 & 0\\ 
x_2+y_2I & x_1-y_1I & 0 & 0\\
0 & 0 & x_1-y_1I & -x_2-y_2I\\
0 & 0 & x_2-y_2I & x_1+y_1I\\
\end{array} \right].\end{equation*}
Therefore, the homotopy class of $c(m_r)$ 
is the sum of the homotopy classes of
\begin{equation*}(z_1+z_2j) \mapsto \left[\begin{array}{cc} 
m^{\CC}_r(z_1,z_2) & 0 \\ 
0 & 1

\end{array} \right] \;\; \mbox{and}\;\; (z_1+z_2j) \mapsto \left[\begin{array}{cc} 
1 & 0 \\ 
0 & m^{\CC}_r \circ \iota(z_1,z_2)
\end{array} \right], \end{equation*}
where $\iota (z_1,z_2)= (\overline{z_1},\overline{z_2})$.
Since the first map is conjugate to $i^2_{\ast}(m^{\CC}_r)$ by a fixed element of $SU(4)$, it is homotopic to $i^2_{\ast}(m^{\CC}_r)$. Since $\iota$ induces the 
identity on $\pi_3(S^3)$, the second map is homotopic to $i^2_{\ast}(m^{\CC}_r)$, too.
\eop

The following lemma finishes to determine the maps \begin{equation*}c_{\ast}\colon \pi_3(SO(4)) \longrightarrow \pi_3(SU(4))\end{equation*}
and $c_{\ast}i_{\ast}\colon \pi_3(SO(3)) \longrightarrow \pi_3(SU(4)).$
\begin{lemma}
\label{lempitroissoquatre} We have
\begin{equation*}c_{\ast}\bigl(\left[\overline{m_r}\right]\bigr)=c_{\ast}\bigl(\left[m_{\ell}\right]\bigr)=-2\left[i^2_{\ast}(m^{\CC}_r)\right]=2\left[i^2_{\ast}(\overline{m}^{\CC}_r)\right],\end{equation*} and
\begin{equation*}c_{\ast}\Bigl(i_{\ast}\bigl(\left[\tilde{\rhomap}\right]\bigr)\Bigr)=4\left[i^2_{\ast}(\overline{m}^{\CC}_r)\right].\end{equation*}
\end{lemma}
\bp
According to Lemma~\ref{lempitroissoquatrestr}, $i_{\ast}\left(\left[\tilde{\rhomap}\right]\right)=\left[m_{\ell}\right] + \left[\overline{m}_r\right].$
Using the conjugacy of quaternions, we have $m_{\ell}(x)(v)=x.v=\overline{\overline{v}.\overline{x}}=\overline{\overline{m}_{r}(x)(\overline{v})}$.
Therefore, $m_{\ell}$ is conjugated to $\overline{m}_{r}$ via 
the conjugacy of quaternions, which lies in $(O(4)\subset U(4))$.

Since $U(4)$ is connected, the conjugacy by an element of $U(4)$ induces the identity on $\pi_3(SU(4))$. 
Thus, we get \begin{equation*}c_{\ast}\bigl(\left[m_{\ell}\right]\bigr)=c_{\ast}\bigl(\left[\overline{m}_{r}\right]\bigr)=-c_{\ast}\bigl(\left[{m}_{r}\right]\bigr)=-2\left[i^2_{\ast}(m^{\CC}_r)\right]=2\left[i^2_{\ast}(\overline{m}^{\CC}_r)\right],\end{equation*}
and 
$c_{\ast}\bigl(i_{\ast}\left(\left[\tilde{\rhomap}\right]\right)\bigr)=c_{\ast}\bigl(\left[m_{\ell}\right]\bigr)+c_{\ast}\bigl(\left[\overline{m}_{r}\right]\bigr)=4\left[i^2_{\ast}(\overline{m}^{\CC}_r)\right]$.
\eop

\section{Relating the Pontrjagin number to the degree}
\label{secpfpone}

We finish proving Theorem~\ref{thmpone} by proving the following proposition. See Lemmas~\ref{lempreptrivun}, \ref{lemdeg}, and \ref{lemdegtwo}.

\begin{proposition}
\label{proppontdeg}
Let $M_0$ and $M$ be two compact connected oriented $3$-manifolds whose boundaries have collars identified by a diffeomorphism. Let $\tau_0\colon M_0 \times \CC^3 \rightarrow TM_0 \otimes \CC$ and $\tau\colon M \times \CC^3 \rightarrow TM \otimes \CC$ be two special complex trivializations (which respect the orientations) that coincide on the collar neighborhoods of $\partial M_0=\partial M$. Let $\left[(M, \partial M),(SU(3),1)\right]$ denote the group of homotopy classes of maps from $M$ to $SU(3)$ that map $\partial M$ to $1$. For any 
\begin{equation*}g\colon (M, \partial M) \longrightarrow (SU(3),1),\end{equation*}
define
\begin{equation*}\begin{array}{llll} 
\psi(g)\colon &M \times \CC^3 &\longrightarrow  &M \times \CC^3\\
&(x,y) & \mapsto &\bigl(x,g(x)(y)\bigr).\end{array}\end{equation*}
Then 
$\left(p_1(\tau_0,\tau \circ \psi(g))-p_1(\tau_0,\tau)\right)$ is independent of $\tau_0$ and $\tau$. Set \begin{equation*}p^{\prime}_1(g)=p_1\bigl(\tau_0,\tau \circ \psi(g)\bigr)-p_1(\tau_0,\tau).\end{equation*}
The map $p^{\prime}_1$ induces an isomorphism from the group $\left[(M, \partial M),(SU(3),1)\right]$ to $\ZZ$,
and, if $g$ is valued in $SO(3)$, then we have \begin{equation*}p^{\prime}_1(g)=2\deg(g).\end{equation*}
\end{proposition}
To prove Proposition~\ref{proppontdeg}, we first prove the following lemma.
\begin{lemma}
\label{lempunind}
Under the hypotheses of Proposition~\ref{proppontdeg}, \begin{equation*}p_1\bigl(\tau_0,\tau \circ \psi(g)\bigr)-p_1\bigl(\tau_0,\tau\bigr)=p_1\bigl(\tau,\tau \circ \psi(g)\bigr)=-p_1\bigl(\tau \circ \psi(g),\tau\bigr)\end{equation*} is independent of $\tau_0$ and $\tau$.
\end{lemma}
\bp 
The equalities of the statement are easy to observe. Let us prove that $p_1(\tau,\tau \circ \psi(g))$ is independent of $\tau$. Let $\tau_W$ be the trivialization of the restriction to $\partial(\left[0,1\right] \times M)$ of the complexified tangent bundle $T\left[0,1\right] \times M \otimes_{\RR} \CC$ to
$\left[0,1\right] \times M$ given by
$T\left[0,1\right] \oplus \tau$ on $(\{0\} \times M) \cup (\left[0,1\right] \times \partial M)$, and $T\left[0,1\right] \oplus \tau \circ \psi(g)$ on $\{1\} \times M$. 
Let $\tilde{g}\colon \partial(\left[0,1\right] \times M)\to SU(4)$ map $(\{0\} \times M) \cup \left(\left[0,1\right] \times \partial M\right)$ to $1$ and coincide with $i \circ g$ on $\{1\} \times M$.
Then the obstruction $p_1(\tau,\tau \circ \psi(g))$ to extending $\tau_W$ to $\left[0,1\right] \times M$ is the obstruction
to extending the map $\tilde{g}$ to $\left[0,1\right] \times M$. It lies in 
$\pi_3(SU(4))$ since $\pi_i(SU(4))=0$ for $i<3$. It is independent of $\tau$.
\eop

Lemma~\ref{lempunind} guarantees that $p^{\prime}_1$ defines two group homomorphisms to $\ZZ$ from $\left[(M, \partial M),(SU(3),1)\right]$ and from $\left[(M, \partial M),(SO(3),1)\right]$. 
\begin{lemma} \label{lempunindbis}
Under the hypotheses of Proposition~\ref{proppontdeg}
\begin{equation*}p^{\prime}_1\colon \bigl[(M, \partial M),(SO(3),1)\bigr] \to \ZZ\end{equation*} is a group homomorphism and
 \begin{equation*}p^{\prime}_1 \colon\bigl[(M, \partial M),(SU(3),1)\bigr] \to \ZZ\end{equation*}
 is a group isomorphism.
\end{lemma} \bp
Since $\pi_i(SU(3))$ is trivial for $i<3$ and since $\pi_3(SU(3))=\ZZ$, the group of homotopy classes  $\left[(M, \partial M) , (SU(3),1)\right]$ is generated by the class of a map that maps the complement of a $3$-ball $B$ to $1$ and that factors through a map whose homotopy class generates $\pi_3(SU(3))$ on $B$. By definition of the Pontrjagin classes, $p^{\prime}_1$ sends such a generator to $\pm 1$, and it induces an isomorphism from the group $\left[(M, \partial M),(SU(3),1)\right]$ to $\ZZ$. 
\eop

\begin{lemma}
\label{lemvarpun} We have
\begin{equation*}p^{\prime}_1\bigl(\rhomap_M(B^3)\bigr)=4.\end{equation*}
\end{lemma}
\bp Let $g=\rhomap_M(B^3)$. Extend the map $\tilde{g}$ of the proof of Lemma~\ref{lempunind} by the constant map with value 1 outside
$\left[\varepsilon, 1\right] \times B^3 \cong B^4$, for some $\varepsilon \in \left]0,1\right[$. We have
\begin{equation*}\left[\tilde{g}\vert_{\partial B^4}\right]=p_1\bigl(\tau,\tau \circ \psi(g)\bigr)\left[i^2_{\ast}(\overline{m}^{\CC}_r)\right]\end{equation*} in $\pi_3(SU(4))$.
Since $\tilde{g}\vert_{\partial B^4}$ is homotopic to $c \circ i \circ \tilde{\rhomap}$, Lemma~\ref{lempitroissoquatre} allows us to conclude.
\eop 

\bpo{Proof of Proposition~\ref{proppontdeg}}
According to Lemmas~\ref{lemdegtwo} and \ref{lempunindbis}, the restriction of $p^{\prime}_1$ to $\left[(M, \partial M),(SO(3),1)\right]$ is equal to
$p^{\prime}_1(\rhomap_M(B^3)) \frac{\deg}{2}$. Conclude with Lemma~\ref{lemvarpun}.
\eop

\section{Properties of Pontrjagin numbers}
\label{secproppontnumb}

\begin{proposition}
\label{proppont}
Let ${\manifm}_0$ and ${\manifm}_1$ be two compact connected oriented $3$-manifolds whose boundaries  have collars identified by a diffeomorphism. Let $\tau_0 \colon {\manifm}_0 \times \CC^3 \rightarrow T{\manifm}_0 \otimes \CC  $ and $\tau_1 \colon {\manifm}_1 \times \CC^3 \rightarrow T{\manifm}_1 \otimes \CC  $ be two special complex trivializations (which respect the orientations) that agree on the collar neighborhoods of $\partial {\manifm}_0=\partial {\manifm}_1$. 

Then the first Pontrjagin number $p_1(\tau_0,\tau_1)$  satisfies the following properties.
\begin{enumerate}
\item Let ${\manifm}_2$ be a compact $3$-manifold whose boundary has a collar neighborhood identified with a collar neighborhood of $\partial {\manifm}_0$.
Let $\tau_2$ be a special complex trivialization of $T{\manifm}_2$ that agrees with $\tau_0$ near $\partial {\manifm}_2$.
If two of the Lagrangians of ${\manifm}_0$, ${\manifm}_1$, and ${\manifm}_2$ coincide in $H_1(\partial {\manifm}_0;\QQ)$, then we have \begin{equation*}p_1(\tau_0,\tau_2)=p_1(\tau_0,\tau_1)+p_1(\tau_1,\tau_2).\end{equation*}
In particular, we also have $p_1(\tau_1,\tau_0)=-p_1(\tau_0,\tau_1)$ since $p_1(\tau_0,\tau_0)=0$.
 \item 
Let $D$ be a connected compact $3$-manifold that contains ${\manifm}_0$ in its interior. Let $\tau_D$ be a special complex trivialization of $TD$ that restricts as the special complex trivialization $\tau_0$ on $T{\manifm}_0$. Let $D_1$ be obtained from $D$ by replacing ${\manifm}_0$ by ${\manifm}_1$. Let $\tau_{D_1}$ be the trivialization
of $TD_1$ that agrees with $\tau_1$ on $T{\manifm}_1$, and with $\tau_D$ on $T(D\setminus {\manifm}_0)$.
\bfig
\centering
\begin{tikzpicture} 
\begin{scope}[xshift=-2.8cm]
\fill [gray!20] (0,0) ellipse (2.6 and 1.2); 
\fill [lightgray] (0,-.25) ellipse (1 and .5); 
\draw (0,0) ellipse (2.6 and 1.2); 
\draw (0,-.25) ellipse (1 and .5); 
\fill [white, rounded corners] (-.5,.55) rectangle (.5,1.05);
\draw (0,.8) node{\scriptsize $(D,\tau_D)$};
\fill [white, rounded corners] (-.55,-.45) rectangle (.55,.05);
\draw (0,-.2) node{\scriptsize $({\manifm}_0,\tau_0)$};
\end{scope}
\begin{scope}[xshift=2.8cm]
 \fill [gray!20] (0,0) ellipse (2.6 and 1.2); 
\fill [lightgray] (0,-.25) ellipse (1 and .5); 
\draw (0,0) ellipse (2.6 and 1.2); 
\draw (0,-.25) ellipse (1 and .5); 
\fill [white, rounded corners] (-.65,.55) rectangle (.65,1.05);
\draw (0,.8) node{\scriptsize $(D_1,\tau_{D_1})$};
\fill [white, rounded corners] (-.55,-.45) rectangle (.55,.05);
\draw (0,-.2) node{\scriptsize $({\manifm}_1,\tau_1)$};
\fill [white, rounded corners] (1.6,-.2) rectangle (2,.2);
\draw (1.8,0) node{\scriptsize $\tau_D$};
\end{scope}
\end{tikzpicture}
\caption{The manifolds $D$ and $D_1$ of Proposition~\ref{proppont}}
\end{figure}
If the Lagrangians of ${\manifm}_0$ and ${\manifm}_1$ coincide,
then we have \begin{equation*}p_1(\tau_D,\tau_{D_1})=p_1(\tau_0,\tau_1).\end{equation*}
\end{enumerate}
\end{proposition}

The proof uses a weak form of the Wall Non-Additivity theorem. We quote the weak form we need.

\begin{theorem}[\cite{wall}]
\label{thmwall}
Let $Y$ be a compact oriented $4$-manifold (with possible boundary). Let $X$ be a $3$-manifold properly embedded in $Y$ that separates $Y$ and that induces the splitting  $Y = Y^+ \cup_{X} Y^-$, for two $4$-manifolds $Y^+$ and $Y^-$ in $Y$, as in the following figure of $Y$:
\begin{center}
\begin{tikzpicture}
\useasboundingbox (-3,-.4) rectangle (3,.4);
\draw (0,-.4) -- (0,.4);
\draw (0,-.4) .. controls (1.2,-.4)  .. (1.2,0)  .. controls (1.2,.4) .. (0,.4);
\draw (0,-.4) .. controls (-1.2,-.4)  .. (-1.2,0)  .. controls (-1.2,.4) .. (0,.4);
\draw (-1.1,-.1) node[left]{\footnotesize $X^-$} (1.1,-.1) node[right]{\footnotesize $X^+$} (-.1,-.1) node[right]{\footnotesize $X$} (.7,0) node{\footnotesize $Y^+$} (-.5,0) node{\footnotesize $Y^-$};
\end{tikzpicture}
\end{center}

The manifold $X$ is the intersection  $Y^+\cap Y^-$. It is oriented as a part of the boundary of $Y^-$.
Set
\begin{equation*}X^+ =\overline {\partial Y^+ \setminus (-X)}\;\;\;\mbox{and}\;\;\;X^- =-\overline {\partial Y^- \setminus X}.\end{equation*} 
Let $\CL$, $\CL^-$, and $\CL^+$ denote the Lagrangians of $X$, $X^-$, and $X^+$, respectively. They are Lagrangian subspaces of $H_1(\partial X, \QQ)$.
Then \begin{equation*}\left(\signature(Y)-\signature(Y^+)-\signature(Y^-)\right)\end{equation*} is the signature of an explicit quadratic form on 
\begin{equation*}\frac{\CL \cap \left(\CL^- + \CL^+ \right)}{(\CL \cap \CL^-) + (\CL \cap \CL^+)}.\end{equation*}
Furthermore, this space is isomorphic to $\frac{\CL^+ \cap \left(\CL + \CL^- \right)}{(\CL^+ \cap \CL) + (\CL^+ \cap \CL^-)}$ and $\frac{\CL^- \cap \left(\CL + \CL^+ \right)}{(\CL^- \cap \CL) + (\CL^- \cap \CL^+)}$.
\end{theorem}
We do not describe the involved quadratic form because
we use this theorem only when the above space is trivial.

\bpo{Proof of Proposition~\ref{proppont}} 
Let us prove the first property. Let $Y^-=W$ be a signature $0$ cobordism from $X^-={\manifm}_0$ to $X={\manifm}_1$. Let $Y^+$ be a signature $0$ cobordism from ${\manifm}_1$ to $X^+={\manifm}_2$. Then it is enough to prove that the signature of $Y=Y^+\cup_X Y^-$ is zero. With the notation of Theorem~\ref{thmwall}, under our assumptions, the space $\frac{\CL \cap (\CL^- + \CL^+ )}{(\CL \cap \CL^-) + (\CL \cap \CL^+)}$
is trivial. Therefore, according to the Wall theorem, the signature of $Y$ is zero. The first property follows.

We now prove that under the assumptions of the second property, we have \begin{equation*}p_1(\tau_D,\tau_{D_1})=p_1(\tau_0,\tau_1).\end{equation*}
Set $Y^+=\left(\left[0,1\right] \times (D\setminus \mathring{{\manifm}}_0)\right)$ and $X=-\left[0,1\right]\times \partial {\manifm}_0$. Let  $Y^-=W$ be a signature $0$ cobordism from ${\manifm}_0$ to ${\manifm}_1$. Note that the signature of $Y^+$ is zero.
In order to prove the desired equality, it is enough to prove that the signature of $Y=Y^+\cup_X Y^-$ is zero.
Here, we have $H_1(\partial X;\QQ)=H_1(\partial {\manifm}_0) \oplus H_1(\partial {\manifm}_0)$.
Let $j \colon H_1(\partial {\manifm}_0) \rightarrow H_1(D\setminus \mathring{{\manifm}}_0)$ and $j_{\partial D} \colon H_1(\partial D) \rightarrow H_1(D\setminus \mathring{{\manifm}}_0)$ be the maps induced by inclusions. With the notation of Theorem~\ref{thmwall}, we have
\begin{equation*}\begin{array}{rl}
\partial X&=-\bigl(\partial\left[0,1\right]\bigr) \times \partial {\manifm}_0,\\
X^- &=-\{1\} \times {\manifm}_1\cup \bigl(\{0\} \times {\manifm}_0 \bigr),\\
X^+&= -\left[0,1\right] \times \partial D\cup \Bigl(\bigl(\partial\left[0,1\right]\bigr) \times (D\setminus \mathring{\manifm}_0)\Bigr),\\
\CL&=\bigl\{(x,-x) \suchthat x \in H_1(\partial {\manifm}_0)\bigr\},\\ 
\CL^-&=\bigl\{(x,y) \suchthat x \in \CL_{{\manifm}_0}, y \in \CL_{{\manifm}_1}\bigr\},\\ 
\CL^+&= \bigl\{(x,y) \suchthat (j(x),j(y))=(j_{\partial D}(z \in H_1(\partial D)),-j_{\partial D}(z))\bigr\}\\
&=\bigl\{(y,-y) \suchthat j(y) \in \Image(j_{\partial D})\bigr\} \oplus \bigl\{(x,0) \suchthat j(x) =0\bigr\},\\
\CL \cap \CL^-  &= \bigl\{(x,-x) \suchthat x \in (\CL_{{\manifm}_0} \cap \CL_{{\manifm}_1}=\CL_{{\manifm}_0})\bigr\},\\
\CL \cap \CL^+  &= \bigl\{(x,-x) \suchthat j(x) \in \Image(j_{\partial D})\bigr\}.
\end{array}\end{equation*}
Let us prove $\CL \cap (\CL^- + \CL^+)=(\CL \cap \CL^-) + (\CL \cap \CL^+)$.
For a subspace $K$ of $H_1(\partial {\manifm}_0;\QQ)$, let $j_{MV}(K)$ be the subspace $\{(x,-x) \suchthat x \in K\}$ of $H_1(\partial X;\QQ)$. Then we have $\CL=j_{MV}(H_1(\partial {\manifm}_0))$, $\CL \cap \CL^+=j_{MV}\left(j^{-1}\bigl(\Image(j_{\partial D})\right)\bigr)$, and 
\begin{equation*}\CL \cap (\CL^- + \CL^+) = \CL \cap \CL^+  + j_{MV}\Bigl(\CL_{{\manifm}_1} \cap \bigl(\CL_{{\manifm}_0} + \Ker(j)\bigr)\Bigr).\end{equation*}
Since $\CL_{{\manifm}_0}=\CL_{{\manifm}_1}$, we have $\left(\CL_{{\manifm}_1} \cap (\CL_{{\manifm}_0} + \Ker(j))\right)= \CL_{{\manifm}_0}$. So $\CL \cap (\CL^- + \CL^+)=(\CL \cap \CL^+) +j_{MV}(\CL_{{\manifm}_0}).$
Then the second property is proved because Wall's theorem guarantees the additivity of the signature in this case.
\eop

\paragraph{The parallelizations of \texorpdfstring{$S^3$}{the sphere of dimension 3}.}

As a Lie group, $S^3$ has two natural homotopy classes of parallelizations $\tau_{\ell}$ and $\tau_{r}$, which we describe below. Identify the tangent space $T_1S^3$ to $S^3$ at $1$ with $\RR^3$ (arbitrarily, with respect to the orientation).
For $g \in S^3$, the multiplication induces two diffeomorphisms $m_{\ell}(g)$ and $m_{r}(g)$ of $S^3$, $m_{\ell}(g)(h)=gh$ and $m_{r}(g)(h)
=hg$. Let $T(m_{\ell}(g))$ and $T(m_{r}(g))$ denote their respective tangent maps at $1$. 
Then we have 
\begin{equation*}\tau_{\ell}\Bigl(h \in S^3, v \in \bigl(\RR^3=T_1S^3\bigr)\Bigr)=\Bigl(h,T_1\bigl(m_{\ell}(h)\bigr)(v) \Bigr)\end{equation*}
and $\tau_{r}(h, v)=\bigl(h,T_1(m_{r}(h))(v)\bigr)$.

\begin{proposition}
\label{proppontstrois} We have
$p_1(\tau_{\ell})=2$ and $p_1(\tau_{r})=-2$.
\end{proposition}
\bp
Regard $S^3$ as the unit sphere of $\HH$. So $T_1S^3=\RR\langle i,j,k \rangle$.
The unit ball $B(\HH)$ of $\HH$ has the standard parallelization of a real vector space equipped with a basis. The trivialization $\tau(\tau_{\ell})$ induced by $\tau_{\ell}$ on 
$\partial B(\HH)$ is such that $\tau_{\ell}(h \in S^3, v \in \HH)=(h,hv) \in (S^3 \times \HH=T\HH\vert_{S^3})$. So, we have
$c_{\ast}\left(\left[m_{\ell}\right]\right)=p_1(\tau_{\ell})\left[i^2_{\ast}(\overline{m}^{\CC}_r)\right]$ in $\pi_3(SU(4))$ by Definition~\ref{deffirstpontcrats} of $p_1$.
According to Lemma~\ref{lempitroissoquatre}, we get $p_1(\tau_{\ell})=2$.
We similarly get $p_1(\tau_{r})=-2$.
\eop

\paragraph{On the image of \texorpdfstring{$p_1$}{the Pontrjagin number}.}
For $n\geq 3$, a {\em spin structure\/} of a smooth $n$-manifold is a homotopy class of parallelizations over a $2$-skeleton of $M$ (or, equivalently, over the complement of a point, if $n=3$ and if $M$ is connected).

The class of the covering map $\tilde{\rhomap}$ described after Definition~\ref{defrho} is the 
standard generator of $\pi_3(SO(3))=\ZZ\left[\tilde{\rhomap}\right]$.
Recall the map $\rhomap_M(B^3)$ of Lemma~\ref{lemdegrho}.
Set $\left[\tilde{\rhomap}\right]\left[\tau\right]=\left[\tau\psi_{\RR}(\rhomap_M(B^3))\right]$. The homotopy classes of parallelizations that induce a given spin structure constitute an affine space with translation group $\pi_3(SO(3))$. According to Theorem~\ref{thmpone} and Lemma~\ref{lemdegrho}, $p_1\left(\left[\tilde{\rhomap}\right]\left[\tau\right]\right)=p_1(\tau)+4$.

\begin{definition}
\label{defmuRohlin}
 The {\em Rohlin invariant\/} $\mu(M,\sigma)$ of a smooth closed $3$-manifold $M$, equipped with a spin structure $\sigma$, is the mod $16$ signature of a compact $4$-manifold $W$, bounded by $M$, equipped with a spin structure that restricts to $M$ as a stabilization of $\sigma$.
\end{definition}
The \emph{first Betti number} of $M$ is the dimension of $H_1(M;\QQ)$. It is denoted by $\beta_1(M)$. 
Robion Kirby and Paul Melvin proved the following theorem \cite[Theorem 2.6]{km}.

\begin{theorem}
\label{thmim} 
For any closed oriented $3$-manifold $M$, for any parallelization $\tau$ of $M$, we have
\begin{equation*}\Bigl(p_1(\tau)-\dimension\bigl(H_1(M;\ZZ/2\ZZ)\bigr)-\beta_1(M)\Bigr) \in 2 \ZZ.\end{equation*}
Let $M$ be a closed $3$-manifold equipped with a given spin structure $\sigma$. Then $p_1$ is a bijection from the set of homotopy classes of parallelizations of $M$ that induce $\sigma$ to \begin{equation*}2\Bigl(\dimension\bigl(H_1(M;\ZZ/2\ZZ)\bigr)+1\Bigr) +\mu(M,\sigma)+4 \ZZ.\end{equation*}
When $M$ is a $\ZZ$-sphere, $p_1$ is a bijection from the set of homotopy classes of parallelizations of $M$ to $(2+4 \ZZ)$.
\end{theorem}

Thanks to Proposition~\ref{proppont}(2), Theorem~\ref{thmim} implies Proposition~\ref{propims} below.

\begin{proposition}
\label{propims}
Let $M_0$ be the unit ball of $\RR^3$. Let $\taust$ be the standard parallelization of $\RR^3$,
\begin{itemize}
\item for any given $\ZZ$-ball $M$, $p_1(.)=p_1((\taust)\vert_{B^3},.)$ defines a bijection from the set of homotopy classes of parallelizations of $M$ that are standard near $\partial M=S^2$ to $4 \ZZ$.
\item For any $\QQ$-ball $M$, for any trivialization $\tau$ of $M$ that is standard near $\partial M=S^2$, we have
\begin{equation*}\Bigl(p_1(\tau)-\dimension\bigl(H_1(M;\ZZ/2\ZZ)\bigr)\Bigr) \in 2 \ZZ.\end{equation*}
\end{itemize}
\end{proposition}

\section{More on \texorpdfstring{$\left[(M,\partial M),(SO(3),1)\right]$}{[M,SO(3)]}}
\label{secmorestrucpar}

This section is a complement to the study of $\left[(M,\partial M),(SO(3),1)\right]$ started in Sections~\ref{secmsothree} and \ref{sechomdegpithree}. It is not used later in this book.
We show how to describe all the elements of $\left[(M,\partial M),(SO(3),1)\right]$ as twists across surfaces, and 
we describe the structure of $\left[(M,\partial M),(SO(3),1)\right]$ precisely, by proving the following theorem.

\begin{theorem}
\label{mainthmtriv}
Let $M$ be a compact oriented $3$-manifold.

If all the closed surfaces embedded in $M$ have an even Euler characteristic,
then $\left[(M,\partial M),(SO(3),1)\right]$ is canonically isomorphic to $H^3(M, \partial M ;\ZZ) \oplus H_2(M;\ZZ/2\ZZ),$ and
the degree maps $\left[(M,\partial M),(SO(3),1)\right]$ onto $2\ZZ$.

If $H_1(M;\ZZ)$ has no $2$-torsion, then all the closed surfaces embedded in $M$ have an even Euler characteristic.

If $M$ is connected, and if there exists a closed surface $S$ of $M$ with odd Euler characteristic, then
the degree maps $\left[(M,\partial M),(SO(3),1)\right]$ onto $\ZZ$, and
$\left[(M,\partial M),(SO(3),1)\right]$ is isomorphic to $\ZZ \oplus \Ker(e\partial_B)$,
where \begin{equation*}e\partial_B \colon H_2(M;\ZZ/2\ZZ) \rightarrow \ZZ/2\ZZ\end{equation*} maps the class of a surface to its Euler characteristic modulo $2$,
and the kernel of $e\partial_B$ has a canonical image in $\left[(M,\partial M),(SO(3),1)\right]$.
\end{theorem}

\paragraph{Representing the elements of \texorpdfstring{$\left[(M,\partial M),(SO(3),1)\right]$}{\left[M,SO(3)\right]} by surfaces.}

\begin{proposition}
\label{propgStau} Let $S$ be a surface, orientable or not, embedded in a $3$-manifold $M$ equipped with a parallelization $\tau$. Recall the map $g(S,\tau)$ from the beginning of Section~\ref{sechomdegpithree}.
If $M$ is connected, then $g(S,\tau)^2$ is homotopic to $\rhomap_M(B^3)^{\chi(S)}$. In particular, the homotopy class of $g(S,\tau)$ depends only on $\chi(S)$ and on the class of $S$ in $H_2(M;\ZZ/2\ZZ)$.
\end{proposition}
\bp Assume that $S$ is connected and oriented. 
Perform a homotopy of $\tau$ so that $p_2 \circ \tau^{-1}$ maps the positive normal  $N^+(S)=T_{(u,s)}([-1,1] \times s)$ to $u \times S$ to a fixed vector $v$ of $S^2$, on $[-1,1] \times (S\setminus D)$, for a disk $D$ of $S$. Then there is a homotopy from $\left[0,1\right] \times [-1,1] \times (S\setminus D)$ to $SO(3)$,
\begin{itemize}
\item which factors through the projection onto  $\left[0,1\right] \times [-1,1]$,
\item which maps $(1,u,s)$ to the rotation $g(S,\tau)^2(u,s)$ with axis $v$ and angle $2{\pi}(u+1)$, for any $(u,s) \in [-1,1] \times (S\setminus D)$, and
\item which maps $\bigl(\partial \left(\left[0,1\right] \times [-1,1]\right) \setminus \left(\{1\}  \times [-1,1]\right)\bigr) \times (S\setminus D) $ to $1$.
\end{itemize}
This homotopy extends to a homotopy $h \colon \left[0,1\right] \times [-1,1] \times S \rightarrow SO(3)$
from $h_0=\rhomap_{[-1,1] \times S}(B^3)^k$ for some $k \in \ZZ$, to $h_1=g(S,\tau)^2\vert_{[-1,1] \times S}$, such that $h(\left[0,1\right] \times \{-1,1\}  \times D)=1$.
Thus, the map $g(S,\tau)^2$ is homotopic to $\rhomap_M(B^3)^{k}$.

Let us show that $k$ depends only on the degree of the Gauss map $G$ from $(D, \partial D)$ to $(S^2,v)$ that maps $s \in (D=\{0\} \times D)$ to $p_2(\tau^{-1}(N^+(s)))$.  Endow $\RR^3$ with an orthonormal basis beginning with $v$. Identify the tangent space to $D$ with $(D \times)\RR^2$ and $T\bigl([-1,1] \times D\bigr)\vert_{D}$ with $([-1,1] \times D \times)\left( \RR N^+(D) \oplus\RR^2\right)$. These identifications allow us to regard $p_2 \circ \tau^{-1}$ as a map $\Psi$ from $[-1,1] \times D$ to $SO(3)$.
The map $\Psi$ determines $\tau$. It sends
 $[-1,1] \times \partial D$ to $\bigl(i(SO(2))=SO(2)\bigr)$.
Since $\pi_2(SO(3))$ is trivial, the homotopy class of $\Psi\vert_{D}$ is determined by its restriction $\Psi\vert_{\partial D}$ to $\partial D$.
The map from $\pi_2(S^2)$ to $\pi_1(SO(2)=S^1)$ in the long exact sequence associated to the fibration \begin{equation*} SO(2) \hookfl{i} SO(3) \rightarrow S^{2},\end{equation*}
described in Theorem~\ref{thmlongseqhomotopy} sends the class of $G$ to the
class of  $\Psi\vert_{\partial D}$. So the degree of $G$ determines the homotopy class of the restriction of $\tau$ to $D$ relative to $\partial D$.
This homotopy class determines the homotopy class relative to $[-1,1] \times \partial D$ of the restriction of $\tau$ to $[-1,1] \times D$, which determines the homotopy class relative to $ \partial \bigl([-1,1] \times D\bigr)$ of the restriction of $g(S,\tau)^2$ to $[-1,1] \times D$.
Therefore, the degree $d$ of $G$ determines $k$. Cutting $D$ into smaller disks shows that $k$ depends linearly on $d$. Note that $d$ is the degree of the Gauss map from $S$ to $S^2$ before the homotopy of $\tau$. For a standard sphere $S^2$, we have $d=1$ and $g(S^2,\tau)=\rhomap_M(B^3)$. So, we get
$k=2d$. It remains to see that the degree $d$ of the Gauss map is $\frac{\chi(S)}{2}$. This is easy to observe for a standard embedding of $S$ into $\RR^3$ equipped with its standard trivialization. Up to homotopy, the trivializations of $TM\vert_{ S}$ are obtained from the standard one by 
compositions by rotations with fixed axis $v$ supported in neighborhoods of curves outside the preimage of $v$. So the degree (at $v$) is independent of the trivialization. Thus, the proposition is proved when $S$ is orientable and connected.
When $S$ is not orientable, the map $g(S,\tau)^2$ is homotopic to $\rhomap_M(B^3)^{k}$, for some $k$, according to Lemma~\ref{lempreptrivun}.
Furthermore, it is homotopic to $g(S^{(2)},\tau)$, where $S^{(2)}$ is the orientable boundary of the tubular neighborhood of $S$. The Euler characteristic of $S^{(2)}$ is $2\chi(S)$. So $g(S,\tau)^4$ is homotopic to $\rhomap_M(B^3)^{2k}$ and to $\rhomap_M(B^3)^{2\chi(S)}$. Since the arguments are local, they extend to the disconnected case and prove that $g(S,\tau)^2$ is homotopic to $\rhomap_M(B^3)^{\chi(S)}$ for any $S$. Then Proposition~\ref{propstructwo} and Lemma~\ref{lemgStau} allow us to conclude that the homotopy class of $g(S,\tau)$ depends only on $\chi(S)$ and on the class of $S$ in $H_2(M;\ZZ/2\ZZ)$.
\eop

Hence, $g(S,\tau)$ will be denoted by $g(S)$.
Lemma~\ref{lemgStau} and Propositions~\ref{propstructwo} and \ref{propgStau} easily imply the following corollary.
\begin{corollary}
\label{corgStautwo}
 All elements of $\left[(M,\partial M),(SO(3),1)\right]$ can be represented by $g(S)$ for some embedded disjoint union of closed surfaces $S$ of $M$.
\end{corollary}
\eopwobp

\paragraph{Structure of \texorpdfstring{$\left[(M,\partial M),(SO(3),1)\right]$}{\left[M,SO(3)\right]}.}
Tensoring a free chain complex $C_{\ast}(M;\ZZ)$ whose homology is $H_{\ast}(M;\ZZ)$ by the short exact sequence \begin{equation*}0 \rightarrow \ZZ \hfl{\times 2} \ZZ  \rightarrow \ZZ/2\ZZ  \rightarrow 0\end{equation*} yields the associated long exact homology sequence 
\begin{equation*}\dots  \rightarrow H_{\ast}(M;\ZZ) \hfl{\times 2} H_{\ast}(M;\ZZ)\rightarrow H_{\ast}(M;\ZZ/2\ZZ) \hfl{\partial_B} H_{\ast-1}(M;\ZZ)  \rightarrow\dots\; ,\end{equation*}
where $\partial_B$ is the {\em Bockstein morphism\/}.

\begin{definition}
\label{defparselflk}
The \emph{self-linking number} of a torsion element $x$ of $H_1(M;\ZZ)$ is the linking number $lk(c,c^{\prime})$ of a curve $c$ that represents $x$ and a parallel $c^{\prime}$ of $c$, modulo $\ZZ$.\footnote{Equivalently, it is the linking number of two disjoint representatives of $x$, modulo $\ZZ$.} It belongs to $\QQ/\ZZ$.
\end{definition}

\begin{proposition}
\label{propdefbock}
There is a canonical group homomorphism \begin{equation*}e\partial_B \colon H_2(M;\ZZ/2\ZZ) \longrightarrow \ZZ/2\ZZ\end{equation*} that admits the following two equivalent definitions:
\begin{enumerate}
 \item For any embedded surface $S$, $e\partial_B$ maps the class of $S$ to the Euler characteristic of $S$ modulo $2$.
\item The map $e\partial_B$ is the composition of the Bockstein morphism
\begin{equation*}\partial_B \colon H_2(M;\ZZ/2\ZZ) \longrightarrow \mbox{$2$-torsion of}\; H_1(M;\ZZ)\end{equation*}
and the map 
\begin{equation*}e \colon \mbox{$2$-torsion of}\; H_1(M;\ZZ) \longrightarrow \ZZ/2\ZZ\end{equation*}
that maps the class of a curve $x$ to $1$ if the self-linking number of $x$ is $\frac{1}{2}$ modulo $\ZZ$, and to $0$ otherwise.
\end{enumerate}
\end{proposition}
\bp The map $e\partial_B$ is well defined by the second definition. It is a group homomorphism. Let $S$ be a connected closed surface.
If $S$ is orientable, then the Euler characteristic of $S$ is even, and the long exact sequence shows that $\partial_B\left(\left[S\right]\right)=0$. 
Otherwise, there is a curve $x$ (Poincar\'e dual to $w_1(S)$) such that $S \setminus x$ is orientable, and the boundary of the closure of the domain of the embedding $S \setminus x$ maps to $(\pm 2x)$, so $\partial_B\left(\left[S\right]\right)=\left[x\right]$, by definition. The characteristic curve $x$ may be assumed to be connected. Then the tubular neighborhood of $x$ in $S$ is either a M\"obius band or an annulus. In the first case, we have $e\left(\left[x\right]\right)=1$ and $\chi(S)$ is odd. Otherwise, we have $e\left(\left[x\right]\right)=0$ and $\chi(S)$ is even.
\eop

\begin{proposition}
\label{propstructhree}
Let $M$ be an oriented connected $3$-manifold.
Then we have \begin{equation*}\left[(M,\partial M),(SO(3),1)\right] \cong \ZZ \oplus \Ker(e\partial_B),\end{equation*} and the degree maps
$\left[(M,\partial M),(SO(3),1)\right]$ onto $2\ZZ$ when $e\partial_B=0$, and onto $\ZZ$ otherwise.
\end{proposition}
\bpo{Proof of Proposition~\ref{propstructhree} and Theorem~\ref{mainthmtriv}}
The class of a surface with even Euler characteristic in $H_2(M;\ZZ/2\ZZ)$ can be represented by a surface $S$ with null Euler characteristic  (a disjoint union with trivial bounding surfaces). According to Proposition~\ref{propgStau}, for such an $S$, the class of $g(S)$ is a $2$-torsion element of $\left[(M,\partial M),(SO(3),1)\right]$ called $\sigma\left(\left[S\right]\right)$.
This defines a canonical partial section 
\begin{equation*}\sigma \colon \bigl(\Ker(e\partial_B)\subset H_2(M;\ZZ/2\ZZ)\bigr)
\rightarrow \Ker\bigl(\deg \colon \left[(M,\partial M),(SO(3),1)\right] \rightarrow \ZZ\bigr)\end{equation*}
of the sequence of Proposition~\ref{propstructwo}.
Therefore, if $e\partial_B = 0$, we have \begin{equation*}\Bigl[(M,\partial M),(SO(3),1)\Bigr]=\ZZ\left[\rhomap_M(B^3)\right] \oplus \sigma\Bigl(\Ker(e\partial_B)=H_2(M;\ZZ/2\ZZ)\Bigr).\end{equation*}
If $e\partial_B \neq 0$, there exists a closed surface $S_1$ with $\chi(S_1)=1$ in $M$. Since the degree is a group homomorphism from $\bigl[(M,\partial M),(SO(3),1)\bigr]$ to $\ZZ$, Proposition~\ref{propgStau} implies $\deg(g(S_1))=1$ for such an $S_1$. Thus we get
\begin{equation*}\bigl[(M,\partial M),(SO(3),1)\bigr]=\ZZ\bigl[g(S_1)\bigr]\oplus \sigma\bigl(\Ker(e\partial_B)\bigr).\end{equation*} 
\eop

\part{The general invariants}

\chapter{Introduction to finite type invariants and Jacobi diagrams}
\label{chapfintype}

This chapter introduces the target space of the invariant $\Zinvuf$ studied in this book. It is a space generated by uni-trivalent graphs called Jacobi diagrams. Dror Bar-Natan has studied this space in his fundational article \cite{barnatan}, where most of the results of this chapter come from.
In this chapter, the field $\KK$ is $\QQ$ or $\RR$.

\section{Definition of finite type invariants}
\label{secdeffintype}

A $\KK$-valued {\em invariant\/} of oriented $3$-manifolds is a function from the set of $3$-manifolds, considered up to orientation-preserving diffeomorphism, to $\KK$. Let $\sqcup_{i=1}^n S_i^1$ denote a disjoint union of $n$ circles, where each $S^1_i$ is a copy of $S^1$.
Here, an {\em $n$-component link\/} in a $3$-manifold $\rats$ is an equivalence class of smooth embeddings
$L \colon \sqcup_{i=1}^n S_i^1 \hookrightarrow \rats$ under the equivalence relation that identifies two embeddings $L$ and $L^{\prime}$ if and only if there is an orientation-preserving diffeomorphism $h$ of $\rats$ such that $h(L)=L^{\prime}$.\footnote{This relation is equivalent to the usual equivalence relation defined by isotopies when $\rats$ is $\RR^3$ or $S^3$. In general $3$-manifolds, two equivalent links are not necessarily isotopic, but the link invariants described in this book are invariant under the above equivalence relation.}
A {\em knot\/} is a one-component link. 
A {\em link invariant\/} (resp. a {\em knot invariant\/}) is a function of links (resp. knots).
For example, $\Theta$ is an invariant of $\QQ$-spheres, and the linking number is a rational invariant of two-component links in rational homology spheres.

In order to study a function, it is common to study its derivative, and the derivatives of its derivative.
The derivative of a function is defined from its variations. 
For a function $f$ from $\ZZ^d = \oplus_{i=1}^d\ZZ e_i$ to $\KK$, one can define its first-order derivatives $\frac {\partial f}{\partial e_i}\colon \ZZ^d \rightarrow \KK$ 
by \begin{equation*}\frac {\partial f}{\partial e_i}(z)=f(z+e_i)-f(z),\end{equation*}
and check that all the first-order derivatives of $f$ vanish if and only if $f$ is constant.
Inductively define an $n$-order derivative to be a first-order derivative of an $(n-1)$-order derivative for a positive integer $n$.
Then it can be checked that all the $(n+1)$-order derivatives of a function $f \colon \ZZ^d \to \KK$ vanish if and only if $f$ is a polynomial of degree not greater than $n$ in the coordinates.
In order to study topological invariants, we can similarly study their variations under {\em simple operations.\/}

Below, $X$ denotes one of the following sets
\begin{itemize}
\item $\ZZ^d$,
\item  the set $\CK$ of knots in $\RR^3$, the set $\CK_k$ of $k$-component links in $\RR^3$,\footnote{Recall that a knot is an isotopy class of knot embeddings.}
\item the set $\CM$ of $\ZZ$-spheres, the set $\CM_{\QQ}$ of $\QQ$-spheres (up to orientation-preserving diffeomorphism).
\end{itemize}
and $\CO(X)$ denotes a set of {\em simple operations\/} acting on some elements of $X$.

For $X=\ZZ^d$, $\CO(X)$ consists of the operations $(z \rightarrow z \pm e_i)$.

For knots or links in $\RR^3$, the {\em simple operations\/} are {\em crossing changes\/}. 
A  {\em crossing change ball\/} of a link $L$ is a ball $B$ of the ambient space, where $L\cap B$ is a disjoint union of two arcs $\alpha_1$ and $\alpha_2$ properly embedded in $B$, and there exist two disjoint topological disks $D_1$ and $D_2$ embedded in $B$, such that, for $i\in \{1,2\}$, the topological circle $\partial D_i$ is the union of $\alpha_i$ and an arc of $\partial B$ as in the following picture:

\begin{center}
{\begin{tikzpicture}
\begin{scope}[yshift=.2cm]
\fill [lightgray] (.15,0) .. controls (.15,-.1) .. (-65:.64) arc (-65:65:.64) (65:.64) .. controls (.15,.1) .. (.15,0);
\draw [-<] (.15,0) .. controls (.15,.1) .. (65:.64) arc (65:-65:.64) (295:.64) .. controls (.15,-.1) .. (.15,0);
\fill [lightgray] (-.15,0) .. controls (-.15,-.1) .. (245:.64) arc (245:115:.64) (115:.64) .. controls (-.15,.1) .. (-.15,0);
\draw [-<] (-.15,0) .. controls (-.15,-.1) .. (245:.64) arc (245:115:.64)  (115:.64) .. controls (-.15,.1) .. (-.15,0);
\draw [dash pattern=on 2pt off 2pt] (0,0) circle (.64);
\draw (-.37,-.2) node{\scriptsize $D_1$} (.37,.2) node{\scriptsize $D_2$} (-.32,.15) node{\scriptsize $\alpha_1$} (.37,-.15) node{\scriptsize $\alpha_2$}; 
\end{scope}
\end{tikzpicture}}
\end{center}
After an isotopy, a projection of $(B,\alpha_1,\alpha_2)$ looks like $\nccirc$ or $\pccirc$. A {\em crossing change\/} is a change that does not change $L$ outside $B$ and that modifies it inside $B$ by a local move \begin{equation*}\left(\nccirc \rightarrow \pccirc\right)\mbox{ or }\left(\pccirc \rightarrow \nccirc\right).\end{equation*} For the left move, the crossing change is {\em positive\/}. It is {\em negative\/} for the move of the right-hand side.

For integer (resp. rational) homology $3$-spheres, the simple operations are integral (resp. rational) {\em LP-surgeries\/}, defined in Subsection~\ref{submmcontext}, and $\CO(\CM)$ (resp. $\CO(\CM_{\QQ})$) is denoted by  $\CO^{\ZZ}_{\CL}$  (resp. $\CO^{\QQ}_{\CL}$).

Say that crossing changes are {\em disjoint\/} if they sit inside disjoint $3$-balls. Say that LP-surgeries $(A^{\prime}/A)$ and $(B^{\prime}/B)$ in a manifold $\rats$ are {\em disjoint\/} if $A$ and $B$ are disjoint in $\rats$.
Two operations on $\ZZ^d$ are always {\em disjoint\/} (even if they look identical).
In particular, disjoint operations \emph{commute} (i.e., their result does not depend on which one is performed first).
Set $\underline{n}=\{1,2,\dots,n\}$. Consider the vector space $\CF_0(X)=\CF_0(X;\KK)$ \index[N]{FnX@$\CF_n(X;\KK)$ filtration} freely generated by $X$ over $\KK$.
For an element $x$ of $X$ and $n$ pairwise disjoint operations $o_1,\dots,o_n$ acting on $x$, let $x((o_i)_{i \in I})$ denote the element of $X$ obtained by performing the operations $o_i$ on $x$ for $i \in I$.
Define \begin{equation*}\left[x;o_1,\dots,o_n\right]=\sum_{I \subseteq \underline{n}} (-1)^{\cardlef{I}}x\bigl((o_i)_{i \in I}\bigr) \in \CF_0(X).\end{equation*}
Then define $\CF_n(X)=\CF_n(X;\KK)$ as the $\KK$-subspace of $\CF_0(X)$ generated by the $\left[x;o_1,\dots,o_n\right]$, for all $x \in X$ equipped with $n$ pairwise disjoint simple operations  $o_1,\dots,o_n$ acting on $x$.
Since we have \begin{equation*}\left[x;o_1,\dots,o_n,o_{n+1}\right]= \left[x;o_1,\dots,o_n\right] - \left[x(o_{n+1});o_1,\dots,o_n\right],\end{equation*}
we get $\CF_{n+1}(X) \subseteq \CF_{n}(X)$ for all $n \in \NN$.

\begin{definition}
A $\KK$-valued function $f$ on $X$, extends uniquely as a $\KK$-linear map on 
\begin{equation*}\CF_0(X)^{\ast}=\Hom(\CF_0(X);\KK),\end{equation*} which is still denoted by $f$.
For an integer $n \in \NN$, the invariant (or function) $f$ is of {\em degree $\leq n$\/} if and only if $f(\CF_{n+1}(X))=0$. The \emph{degree}\index[T]{degree!of an invariant} of such an invariant is the smallest integer $n \in \NN$ such that $f(\CF_{n+1}(X))=0$. 
An invariant is of {\em finite type\/} if it is of degree $n$ for some $n \in \NN$. This definition depends on the chosen set of operations $\CO(X)$. We fixed our choices for our sets $X$, but other choices could lead to different notions. See \cite{ggp} for $X=\CM$.
\end{definition}

Let $\CI_n(X)=(\CF_0(X)/\CF_{n+1}(X))^{\ast}$ be the space of invariants of degree at most $n$.
Of course, we have $\CI_n(X) \subseteq \CI_{n+1}(X)$ for all $n \in \NN$.

\begin{exo}
Prove that $\CI_n(\ZZ^d)$ is the space of polynomials of degree at most $n$ on $\ZZ^d$.
\end{exo}

\begin{lemma}
\label{lemmultinvfin}
 If $\funcf\in\mathcal{I}_m(X)$ and $\funcg\in\mathcal{I}_n(X)$, then $\funcf\funcg\in\mathcal{I}_{m+n}(X)$.
\end{lemma}
\bp
Let $\left[x;(o_i)_{i\in \underline{m+n+1}}\right] \in \CF_{m+n+1}(X)$. The lemma is a direct consequence of the equality
 \begin{equation*}\funcf\funcg\Bigl(\bigl[x;(o_i)_{i\in \underline{m+n+1}}\bigr]\Bigr)=\sum_{J\subseteq \underline{m+n+1}}
 \funcf\Bigl(\bigl[x;(o_j)_{j\in J}\bigr]\Bigr)\funcg\Bigl(\bigl[x((o_j)_{j\in J});(o_i)_{i\in \underline{m+n+1} \setminus J}\bigr]\Bigr),\end{equation*}
which is proved as follows. The right-hand side is equal to
\begin{multline*}\sum_{J\subseteq \underline{m+n+1}} (-1)^{|J|} \left( \sum_{K \suchthat  K\subseteq J} (-1)^{|K|} \funcf\Bigl(x\bigl((o_i)_{i\in K}\bigr)\Bigr)\right) 
 \left( \sum_{L \suchthat  J \subseteq L} (-1)^{|L|} \funcg\Bigl(x\bigl((o_i)_{i\in L}\bigr)\Bigr)\right)\\
=\sum_{(K,L) \suchthat  K \subseteq L \subseteq \underline{m+n+1}} (-1)^{|K|+|L|} \funcf\Bigl(x\bigl((o_i)_{i\in K}\bigr)\Bigr) \funcg\Bigl(x\bigl((o_i)_{i\in L}\bigr)\Bigr)
 \left( \sum_{J \suchthat  K\subseteq J\subseteq L} (-1)^{|J|} \right),
\end{multline*}
where $\sum_{J \suchthat K\subseteq J\subseteq L} (-1)^{|J|}=\left\lbrace 
\begin{array}{l l} 0 & if\ K\subsetneq L \\ (-1)^{|K|} & if\ K=L . \end{array}\right.$
\eop

\begin{lemma}
 Any $n$-component link in $\RR^3$ can be transformed to the trivial $n$-component link below by a finite number of disjoint crossing changes.

\begin{center}
 \begin{tikzpicture}
\useasboundingbox (-.5,-.35) rectangle (5.5,.35);
\draw [thick,->] (-.35,0) arc (-180:180:.35);
\draw  (-.35,0) node[left]{\scriptsize $U_1$};
\draw [thick,->] (2,0) arc (-180:180:.35);
\draw  (2,0) node[left]{\scriptsize $U_2$};
\draw  (3.5,0) node{\dots};
\draw [thick,->] (5,0) arc (-180:180:.35);
\draw  (5,0) node[left]{\scriptsize $U_n$};
\end{tikzpicture}\end{center}
\end{lemma}
\bp
Let $L$ be an (embedding representing an) $n$-component link in $\RR^3$.
Since $\RR^3$ is simply connected, there is a homotopy that carries $L$ to the trivial link. Such a homotopy $h \colon \left[0,1\right]\times \sqcup_{i=1}^n S^1_i\rightarrow \RR^3$ can be chosen to be smooth and such that $h(t,.)$ is an embedding, except for finitely many times $t_i$, $0<t_1 < \dots <t_i < t_{i+1} < \dots <1$, at which $h(t_i,.)$ is an immersion with one double point and no other multiple points, and the link $h(t,.)$ changes exactly by a crossing change when $t$ crosses a $t_i$. (For an alternative elementary proof of this fact, see \cite[Subsection~7.1]{lescol} before Definition~7.5, for example).
\eop

In particular, a degree $0$ invariant of $n$-component links of $\RR^3$ must be constant since it does not vary under a crossing change.

\begin{exo}
1. Check that  $\CI_1(\CK)=\KK c_0$, where $c_0$ is the constant map that maps any knot to $1$.\\
2. Check that the linking number is a degree $1$ invariant of $2$-component links of $\RR^3$.\\
3. Check that $\CI_1(\CK_2)=\KK c_0 \oplus \KK lk$, where $c_0$ is the constant map that maps any two-component link to $1$.
\end{exo}

\section{Introduction to chord diagrams}
\label{secchordd}

A {\em singular knot with $n$ double points\/} in $\RR^3$ is an immersion of a circle with $n$ transverse double points in $\RR^3$.
Such a double point \doubleppetittexte can be \emph{desingularized} in two ways, the \emph{positive} one \pcpetittexte and the \emph{negative} one \ncpetittexte.
For example, desingularizing the double points of the singular knot \singtrefoil in the positive way produces the knot \righttrefoil. Note that the sign of the desingularization is defined from the orientations of the knot and of the ambient space.

Let $o_i$ be disjoint negative crossing changes  $\pcpetit \to \ncpetit$  to be performed on a knot $K$.
We represent $\left[K;o_1,\dots,o_n\right]$ as a singular knot with $n$ double points such that $K\bigl((o_i)_{i \in I}\bigr)$ is obtained from the singular knot by desingularizing the crossings of $I$ in the negative way, and the others in the positive way, for $I \subseteq \underline{n}$.
Thus, singular knots represent elements of $\CF_0(\CK)$. Three singular knots that coincide outside a ball, inside which they look as in the following \emph{skein relation} 
\begin{equation*}\doublep =\pc -\nc, \end{equation*}
satisfy this relation in $\CF_0(\CK)$.

Define the {\em chord diagram\/} $\Gamma_C\left(\left[K;o_1,\dots,o_n\right]\right)$ associated to $\left[K;o_1,\dots,o_n\right]$ as follows. Draw the preimage of the associated singular knot with $n$ double points as an oriented dashed circle equipped with the $2n$ preimages of the double points, and join the pairs of preimages of a double point by a plain segment called a {\em chord\/}. For example, we have
\begin{equation*}\Gamma_C\Bigl(\singtrefoil\Bigr)=\threexxchordbas.\end{equation*}
Formally, a \indexT{chord diagram} with $n$ chords (on a circle) is a cyclic order of the $2n$ ends of the $n$ chords, up to a permutation of the chords
and up to exchanging the two ends of a chord.

\begin{lemma}
\label{leminvtyf} If $f$ is an invariant of knots in $\RR^3$ of degree at most $n$,
then $f\left(\left[K;o_1,\dots,o_n\right]\right)$ depends only on $\Gamma_C\left(\left[K;o_1,\dots,o_n\right]\right)$.
\end{lemma}
\bp
Since $f$ is of degree $n$, $f\left(\left[K;o_1,\dots,o_n\right]\right)$ is invariant under a crossing change outside the balls of the $o_i$, that is outside the double points of the associated singular knot. Therefore, $f\left(\left[K;o_1,\dots,o_n\right]\right)$ depends only on the cyclic order of the $2n$ arcs involved in the $o_i$ on $K$. (A more detailed proof can be found in \cite[Subsection 7.3]{lescol}.)
\eop

Let $\CD_n$ be the $\KK$-vector space freely generated by the $n$-chord diagrams (i.e., the diagrams with $n$ chords)  on $S^1$.
For example, we have
\begin{multline*}
\CD_0=\KK \zerochord \mbox{, } \CD_1=\KK \onechord\mbox{, } \CD_2=\KK \twoischord \oplus \KK \twoxchord\mbox{, and }   \\ \CD_3=\KK \threeischord \oplus \KK \threeparchord \oplus \KK \threexchord \oplus \KK \threexparchord \oplus \KK \threexxchord.\end{multline*}

\begin{lemma}
\label{lemphinsurj}
The map $\phi_n$ from $\CD_n$ to ${\CF_n(\CK)}/{\CF_{n+1}(\CK)}$ that maps an $n$-chord diagram $\Gamma$ to some $\left[K;o_1,\dots,o_n\right]$ whose diagram is $\Gamma$ is well-defined and surjective.
\end{lemma}
\bp Use the arguments of the proof of Lemma~\ref{leminvtyf}.
\eop

For example, we have \begin{equation*}\phi_3\Bigl(\threexxchordbas\Bigr)=\Bigl[\singtrefoil\Bigr].\end{equation*}
Lemma~\ref{lemphinsurj} implies that \begin{equation*}\phi_n^{\ast} \colon \left(\frac{\CF_n(\CK)}{\CF_{n+1}(\CK)}\right)^{\ast} \to \CD_n^{\ast}\end{equation*}
is injective.
The kernel of the restriction below
 \begin{equation*}\left(\CI_n(\CK)=\left(\frac{\CF_0(\CK)}{\CF_{n+1}(\CK)}\right)^{\ast}\right)
\rightarrow \left(\frac{\CF_n(\CK)}{\CF_{n+1}(\CK)}\right)^{\ast} \end{equation*}
is $\CI_{n-1}(\CK)$. Thus, ${\CI_n(\CK)}/{\CI_{n-1}(\CK)}$ injects into $\CD_n^{\ast}$, and the dimension of $\CI_n(\CK)$ is finite for all $n$.
In particular, ${\CF_0(\CK)}/{\CF_{n+1}(\CK)}$ is finite-dimensional, and the above restriction is surjective.
Therefore, we have \begin{equation*}\frac{\CI_n(\CK)}{\CI_{n-1}(\CK)}=\Hom\left(\frac{\CF_n(\CK)}{\CF_{n+1}(\CK)};\KK\right).\end{equation*}

An {\em isolated chord\/} in a chord diagram is a chord between two points of $S^1$ that are consecutive on the circle $S^1$.

\begin{lemma}
\label{lem14T}
Let $D$ be a diagram on $S^1$ that contains an isolated chord. Then
$\phi_n(D)=0.$
Let $D^1$, $D^2$, $D^3$, and $D^4$ be four n-chord diagrams identical outside three portions of circles, inside which they look like
\begin{equation*}D^1=\dquatTun \mbox{, } D^2= \dquatTdeux \mbox{, } D^3= \dquatTtrois\mbox{, and } D^4=\dquatTqua.\end{equation*}
Then we have
\begin{equation*}\phi_n(-D^1+D^2+D^3-D^4)=0.\end{equation*}
\end{lemma}
\bp For the first assertion, observe
\begin{equation*}\phi_n(\ischord)=\left[\pkink\right]-\left[\nkink\right].\end{equation*}
Let us prove the second one.
We may represent \begin{equation*}D^1=\dquatTunnum \end{equation*} by a singular knot $K^1$ with n double points, which intersects
a ball as 
\begin{equation*}K^1=\begin{tikzpicture} [baseline=.2cm] \useasboundingbox (0,-.1) rectangle (1,1.1);
\draw [->]  (1.05,.8) node{\scriptsize 1} (.1,0) -- (.9,.8);
\draw [->,draw=white,double=black,very thick] (.5,.95) node{\scriptsize 3} (.7,.2) .. controls (.9,.4) .. (.7,.6) -- (.5,.8);
\draw [->] (.7,.6) -- (.5,.8);
\draw [->] (-.05,.8) node{\scriptsize 2} (.9,0) -- (.1,.8);
\draw (.5,0) -- (.7,.2);
\fill (.5,.4) circle (1.5pt) (.7,.2) circle (1.5pt);
\end{tikzpicture}.\end{equation*}

Let $K^2$, $K^3$, $K^4$ be the singular knots with $n$ double points that coincide with $K^1$ outside this ball, and that intersect this ball as shown in the pictures:
\begin{equation*}K^2= \begin{tikzpicture} [baseline=.2cm] \useasboundingbox (0,-.1) rectangle (1,1.1);
\draw [->] (.5,.95) node{\scriptsize 3} (.7,.2) .. controls (.9,.4) .. (.7,.6) -- (.5,.8);
\draw (.5,0) -- (.7,.2);
\draw [->] (.7,.6) -- (.5,.8);
\draw [->,draw=white,double=black,very thick] (-.05,.8) node{\scriptsize 2} (.9,0) -- (.1,.8);
\draw [->]  (1.05,.8) node{\scriptsize 1} (.1,0) -- (.9,.8);
\draw [->] (.9,0) -- (.1,.8);

\fill (.5,.4) circle (1.5pt) (.7,.6) circle (1.5pt);
\end{tikzpicture}
\mbox{, }
K^3= \begin{tikzpicture} [baseline=.2cm] \useasboundingbox (0,-.1) rectangle (1,1.1);
\draw [->] (.5,.95) node{\scriptsize 3} (.3,.2) .. controls (.1,.4) .. (.3,.6) -- (.5,.8);
\draw (.5,0) -- (.3,.2);
\draw [->,draw=white,double=black,very thick]  (1.05,.8) node{\scriptsize 1} (.1,0) -- (.9,.8);
\draw [->] (.1,0) -- (.9,.8);
\draw [->] (-.05,.8) node{\scriptsize 2} (.9,0) -- (.1,.8);
\draw [->] (.3,.6) -- (.5,.8);

\fill (.5,.4) circle (1.5pt) (.3,.6) circle (1.5pt);
\end{tikzpicture}
\mbox{, }
K^4= \begin{tikzpicture} [baseline=.2cm] \useasboundingbox (0,-.1) rectangle (1,1.1);
\draw [->] (-.05,.8) node{\scriptsize 2} (.9,0) -- (.1,.8);
\draw [->,draw=white,double=black,very thick] (.5,.95) node{\scriptsize 3} (.3,.2) .. controls (.1,.4) .. (.3,.6) -- (.5,.8);
\draw [->] (.3,.6) -- (.5,.8);
\draw [->]  (1.05,.8) node{\scriptsize 1} (.1,0) -- (.9,.8);
\draw (.5,0) -- (.3,.2);
\fill (.5,.4) circle (1.5pt) (.3,.2) circle (1.5pt);
\end{tikzpicture}. 
\end{equation*}

Then the chord diagram $D(K^2)$ associated to $K^2$ is $D^2$. Similarly, we have $D(K^3)=D^3$ and $D(K^4)=D^4$. 
Therefore, we have $\phi_n(-D^1+D^2+D^3-D^4)=-\left[K^1\right]+\left[K^2\right]+\left[K^3\right]-\left[K^4\right].$

Thus, it is enough to prove that we have
\begin{equation*}-\left[K^1\right]+\left[K^2\right]+\left[K^3\right]-\left[K^4\right]=0\end{equation*}
in $\CF_n(\CK)$. Let us prove it.

Let $K_0$ be the singular knot with $(n-1)$ double points that intersects our ball
as \begin{equation*}K_0= \begin{tikzpicture} [baseline=.2cm] \useasboundingbox (0,-.1) rectangle (2,1.1);
\draw [->]  (1.05,.8) node{\scriptsize 1} (.1,0) -- (.9,.8);
\draw [->] (-.05,.8) node{\scriptsize 2} (.9,0) -- (.1,.8);
\draw [->,draw=white,double=black,very thick] (.5,.95) node{\scriptsize 3} (.5,0) -- (.7,.2) .. controls (.9,.4) .. (.7,.6) -- (.5,.8);
\draw [->] (.7,.6) -- (.5,.8);
\fill (.5,.4) circle (1.5pt);
\end{tikzpicture}
\end{equation*}
and that coincides with $K^1$ outside this ball.

The strands $1$ and $2$ in the pictured double point are in the
horizontal plane. They orient it. The strand $3$ is vertical. It intersects the horizontal plane in a positive way between the tails of $1$ and $2$.
Now, make $3$ turn around the double point counterclockwise so that it successively becomes the knots with $(n-1)$ double points:
\begin{equation*}K_1= \begin{tikzpicture} [baseline=.2cm] \useasboundingbox (0,-.1) rectangle (1,1.1);
\draw [->]  (1.05,.8) node{\scriptsize 2} (.1,0) -- (.9,.8);
\draw [->,draw=white,double=black,very thick] (.5,.95) node{\scriptsize 3} (.5,0) -- (.8,.3) .. controls (.9,.4) .. (.6,.7) -- (.5,.8);
\draw [->] (.6,.7) -- (.5,.8);
\draw [->,draw=white,double=black,very thick] (.9,0) -- (.5,.4);
\draw [->] (-.05,.8) node{\scriptsize 1} (.9,0) -- (.1,.8);
\fill (.5,.4) circle (1.5pt);
\end{tikzpicture}
\mbox{, }
K_2= \begin{tikzpicture} [baseline=.2cm] \useasboundingbox (0,-.1) rectangle (1,1.1);
\draw [->] (.5,.95) node{\scriptsize 3} (.5,0) -- (.8,.3) .. controls (.9,.4) .. (.6,.7) -- (.5,.8);
\draw [->,draw=white,double=black,very thick]  (1.05,.8) node{\scriptsize 2} (.1,0) -- (.9,.8);
\draw [->,draw=white,double=black,very thick] (-.05,.8) node{\scriptsize 1} (.9,0) -- (.1,.8);
\draw [->] (.9,0) -- (.1,.8);
\draw [->] (.1,0) -- (.9,.8);
\fill (.5,.4) circle (1.5pt);
\end{tikzpicture}
\mbox{, and }
K_3=\begin{tikzpicture} [baseline=.2cm] \useasboundingbox (0,-.1) rectangle (1,1.1);
\draw [->] (-.05,.8) node{\scriptsize 1} (.9,0) -- (.1,.8) ;
\draw [->,draw=white,double=black,very thick] (.5,.95) node{\scriptsize 3} (.5,0) -- (.2,.3) .. controls (.1,.4) .. (.4,.7) -- (.5,.8);
\draw [->] (.4,.7) -- (.5,.8);
\draw [->,draw=white,double=black,very thick] (.1,0) -- (.4,.3);
\draw [->] (1.05,.8) node{\scriptsize 2} (.1,0) -- (.9,.8);
\fill (.5,.4) circle (1.5pt);
\end{tikzpicture}.\end{equation*}
On its way, it goes successively
through our four knots $K^1$, $K^2$, $K^3$, and $K^4$ with $n$ double points, which appear inside matching parentheses, in the following obvious identity
in $\CF_{n-1}(\CK)$
\begin{equation*}\bigl(\left[K_1\right]-\left[K_0\right]\bigr)+\bigl(\left[K_2\right]-\left[K_1\right]\bigr)+\bigl(\left[K_3\right]-\left[K_2\right]\bigr)+\bigl(\left[K_0\right]-\left[K_3\right]\bigr)=0.\end{equation*}
Now, we have $\left[K^i\right]=\pm \left(\left[K_i\right]-\left[K_{i-1}\right]\right)$, where the sign $\pm$ is $+$ when the vertical strand goes through an arrow
from $K_{i-1}$ to $K_i$, and minus when it goes through a tail.
Therefore, the above equality can be written as
\begin{equation*}-\left[K^1\right]+\left[K^2\right]+\left[K^3\right]-\left[K^4\right]=0.\end{equation*}
\eop

Let $\Asimp_n$ \index[N]{Apb@$\Asimp_n$ chord diagrams mod 4T} denote the quotient of $\CD_n$ by the {\em four--term relation $(4T)$.\/} It is the quotient of $\CD_n$
by the vector space generated by the $(-D^1+D^2+D^3-D^4)$ for all the $4$-tuples $(D^1,D^2,D^3,D^4)$ as in Lemma~\ref{lem14T}.
Call $(1T)$ the relation that identifies a diagram with an isolated chord with $0$. So $\Asimp_n/(1T)$ is the quotient of $\Asimp_n$ by the vector space generated by diagrams with an isolated chord.

According to Lemma~\ref{lem14T} above, the map $\phi_n$ induces a map 
\begin{equation*}\overline{\phi}_n \colon\Asimp_n/(1T) \longrightarrow \frac{\CF_n(\CK)}{\CF_{n+1}(\CK)}.\end{equation*}

The fundamental theorem of {\em Vassiliev invariants\/} (which are finite type knot invariants) can now be stated.

\begin{theorem}[Bar-Natan, Kontsevich]
\label{thmbn}
There exists a family \begin{equation*}\left( \Zinvlinkufmodis_n \colon \CF_0(\CK) \rightarrow \Asimp_n/(1T) \right)_{n \in \NN}\end{equation*} of linear maps satisfying:
\begin{itemize}
\item $\Zinvlinkufmodis_n(\CF_{n+1}(\CK))=0$.
\item Let $\Zinvlinkufov_n$ be the map induced by $\Zinvlinkufmodis_n$ from ${\CF_n(\CK)}/{\CF_{n+1}(\CK)}$ to $\Asimp_n/(1T)$. Then $\Zinvlinkufov_n \circ \overline{\phi}_n$ is the identity map of $\Asimp_n/(1T)$.
\end{itemize}
In particular, the inverse isomorphisms $\Zinvlinkufov_n$ and $\overline{\phi}_n$ identify ${\CF_n(\CK)}/{\CF_{n+1}(\CK)}$ with $\Asimp_n/(1T)$, and we have \begin{equation*}\frac{\CI_n(\CK)}{\CI_{n-1}(\CK)}\cong (\Asimp_n/(1T))^{\ast}.\end{equation*}
\end{theorem}
Maxim Kontsevich and Dror Bar-Natan proved this theorem by using the {\em Kontsevich integral\/} $Z^K=(Z^K_n)_{n\in \NN}$ \cite{barnatan}, for $\KK=\RR$. It is also true when $\KK=\QQ$. It is reproved in Section~\ref{secproofthmbn} using the invariant $\Zinvuf$ studied in this book. An invariant $\Zinvlinkufmodis$ as in the above statement has the following universality property.

\begin{corollary}
\label{coruniv}
For any real-valued finite type invariant $f$ of knots in $\RR^3$ of degree at most $n$, there exist linear forms $\psi_i \colon \Asimp_i/(1T) \rightarrow \RR$, for $i=0, \dots, n$, such that
\begin{equation*}f= \sum_{i=0}^n \psi_i \circ \Zinvlinkufmodis_i.\end{equation*}
\end{corollary}
\bp 
Let $\psi_n = f\vert_{\CF_{n}} \circ \overline{\phi}_n$, then $\bigl(f-\psi_n \circ \Zinvlinkufmodis_n\bigr)$ is an invariant of degree at most $n-1$. Conclude by induction.
\eop

By projection (or up to (1T)), the invariant $\Zinvlinkufmodis$ defines a universal Vassiliev knot invariant with respect to the following definition.

An invariant $\Zgen \colon \CF_0(\CK) \rightarrow \prod_{n \in \NN} \Asimp_n/(1T)$ such that
\begin{itemize}
\item $\Zgen_n(\CF_{n+1}(\CK))=0$ and
\item $\Zgen_n$ induces a left inverse to $\overline{\phi}_n$ from ${\CF_n(\CK;\RR)}/{\CF_{n+1}(\CK;\RR)}$ to $\Asimp_n/(1T)$
\end{itemize}
is called a \emph{universal Vassiliev knot invariant}\index[T]{universal!Vassiliev knot invariant}. 

The terminology is justified because such an invariant contains all the real-valued Vassiliev knot invariants as in Corollary~\ref{coruniv}.

As Daniel Altsch\"uler and Laurent Freidel proved in \cite{af}, the restriction of the invariant $\Zinvuf=(\Zinvuf_n)_{n\in \NN}$ to knots of $\RR^3$ also satisfies the properties of Theorem~\ref{thmbn}. So it is also a \emph{universal Vassiliev knot invariant.}
We will give alternative proofs of generalizations of this Altsch\"uler--Freidel theorem in Section~\ref{secproofthmbn}.

Similar characterizations of the spaces of finite type invariants of links in $\RR^3$, integer homology $3$-spheres, and rational homology $3$-spheres will be presented in Section~\ref{secproofthmbn} and Chapter~\ref{chapuniv}, respectively. For integer homology $3$-spheres and rational homology $3$-spheres, the most difficult parts of the proofs are consequences of the splitting formulae satisfied by $\Zinvuf$, stated in Theorem~\ref{thmmainunivlag}.

We end this subsection with an example of a nontrivial linear form on $\Asimp_n/(1T)$.
\begin{example}
\label{exaweightConway}
We first define a function $\check{w}_C$ of chord diagrams. Let $\Gamma$ be a
chord diagram. Immerse $\Gamma$ in the unit disk $\drad{1}$ of the plane so that the chords of $\Gamma$ are embedded and attached to the left-hand side of the boundary $S^1$ of $\drad{1}$, as in our former pictures like
\begin{equation*}\twoxchord\mbox{ or }\threexparchord.\end{equation*}

Attach disjoint oriented bands $\left[0,1\right]^2$, one \say{around each chord}, so that $\left(\partial \left[0,1\right]\right) \times \left[0,1\right]$ is a well-oriented neighborhood of the ends of the chord in the dashed circle $S^1$, as in the following figure. Perform the \emph{surgery} on the dashed circle $S^1$ that replaces $\left(\partial \left[0,1\right]\right) \times \left[0,1\right]$ with $\left[0,1\right] \times \partial \left[0,1\right]$ as in the following figure.

\begin{center}
\begin{tikzpicture} \useasboundingbox (-.4,-.2) rectangle (.4,.3);
\draw [->,dash pattern=on 2pt off 1pt] (-130:.3) arc (-130:-50:.3);
\draw [->,dash pattern=on 2pt off 1pt] (50:.3) arc (50:130:.3);
\draw (0,-.3) -- (0,.3);
\fill (0,-.3) circle (1.5pt) (0,.3) circle (1.5pt);
\end{tikzpicture}
$\longrightarrow$
\begin{tikzpicture} \useasboundingbox (-.4,-.2) rectangle (.4,.3);
\fill[gray!50] (75:.3) arc (75:105:.3) -- (-105:.3) arc (-105:-75:.3) -- (75:.3);
\draw [->,dash pattern=on 2pt off 1pt] (-130:.3) arc (-130:-50:.3);
\draw [->,dash pattern=on 2pt off 1pt] (50:.3) arc (50:130:.3);
\end{tikzpicture}
$\longrightarrow$
\begin{tikzpicture} \useasboundingbox (-.4,-.2) rectangle (.4,.3);
\draw [dash pattern=on 2pt off 1pt] (-130:.3) arc (-130:-105:.3);
\draw [dash pattern=on 2pt off 1pt] (-105:.3) -- (105:.3);
\draw [->,dash pattern=on 2pt off 1pt] (105:.3) arc (105:130:.3);
\draw [dash pattern=on 2pt off 1pt] (50:.3) arc (50:75:.3);
\draw [dash pattern=on 2pt off 1pt] (75:.3) -- (-75:.3);
\draw [->,dash pattern=on 2pt off 1pt] (-75:.3) arc (-75:-50:.3);
\end{tikzpicture}
\end{center}
If the resulting naturally oriented one-manifold is connected, then $\check{w}_C(\Gamma)=1$. Otherwise, $\check{w}_C(\Gamma)=0$. For example, we have
\begin{equation*}\check{w}_C\Bigl(\zerochordbas\Bigr)=1, \;\check{w}_C\Bigl(\onechordbas\Bigr)=0, \;\mbox{and}\;\check{w}_C\Bigl(\twoischordbas\Bigr)=0.\end{equation*} Since 
\begin{tikzpicture} \useasboundingbox (-.4,-.12) rectangle (.4,.2);
\draw (120:.3) -- (-30:.3) (150:.3) -- (-60:.3);
\draw [draw=white,double=black,very thick] (-120:.3) -- (30:.3) (-150:.3) -- (60:.3);
\draw (-120:.3) -- (30:.3) (60:.3) -- (-150:.3);
\draw (-30:.3) arc (-30:30:.3) (60:.3) arc (60:120:.3) (150:.3) arc (150:210:.3) (-120:.3)
arc (-120:-60:.3);
\end{tikzpicture}
is connected, we have $\check{w}_C\Bigl(\twoxchordbas\Bigr)=1$. We also have
\begin{equation*}\check{w}_C\Bigl(\threexparchordbas\Bigr)= \check{w}_C\Bigl(\threexxchordbas\Bigr)=0.\end{equation*}
More generally, the reader can check that one of our surgeries changes the mod $2$ congruence class of the number of connected components of the surgered manifold. So $\check{w}_C(\Gamma)= 0$ for any
chord diagram $\Gamma$ with an odd number of chords. Extend $\check{w}_C$ linearly over $\CD_n$. For $4$-tuples $(D^1,D^2,D^3,D^4)$ as in Lemma~\ref{lem14T}, the reader can check that the extended $\check{w}_C$ maps $(-D^1+D^2+D^3-D^4)$ to zero. It also maps diagrams with isolated chords to zero. Therefore, $\check{w}_C$ induces a linear map $w_C$ on $\Asimp_n$ (and on $\Asimp_n/(1T)$) for any $n$. For any chord diagram $\Gamma$, we have $w_C\left(\left[\Gamma\right]\right)=\check{w}_C(\Gamma)$. Dror Bar-Natan and Stavros Garoufalidis introduced this linear map in \cite{bngaroufMMR}. They called it the \emph{Conway weight system}.
It is zero when $n$ is odd. It is not zero when $n$ is even since we have \begin{equation*}w_C\left(\left[
\begin{tikzpicture}  [baseline=.2cm] \useasboundingbox (-.4,0) rectangle (4,.6);
\draw [->,dash pattern=on 2pt off 1pt] (3.6,.6) -- (0,.6) arc (90:270:.3) (2.6,0) -- (3.6,0) arc (-90:90:.3) ;
\draw [->,dash pattern=on 2pt off 1pt] (0,0) -- (2,0);
\draw (2.3,0) node{\dots};
\draw (.6,0)  arc (0:180:.2) (.8,0)  arc (0:180:.2);
\fill (.2,0) circle (1.5pt) (.4,0) circle (1.5pt) (.6,0) circle (1.5pt) (.8,0) circle (1.5pt);
\begin{scope}[xshift=1cm]
\draw (.6,0)  arc (0:180:.2) (.8,0)  arc (0:180:.2);
\fill (.2,0) circle (1.5pt) (.4,0) circle (1.5pt) (.6,0) circle (1.5pt) (.8,0) circle (1.5pt);
 \end{scope}
 \begin{scope}[xshift=2.6cm]
\draw (.6,0)  arc (0:180:.2) (.8,0)  arc (0:180:.2);
\fill (.2,0) circle (1.5pt) (.4,0) circle (1.5pt) (.6,0) circle (1.5pt) (.8,0) circle (1.5pt);
 \end{scope}
\end{tikzpicture}\right]\right)\neq 0.\end{equation*}
\end{example}

\section{More spaces of diagrams}

\begin{definition}
\label{defunitrivgra}
A {\em uni-trivalent graph\/} $\Gamma$ is a $6$-tuple \begin{equation*}(H(\Gamma),E(\Gamma),U(\Gamma),T(\Gamma),p_E,p_V),\end{equation*} where
\begin{itemize}
 \item 
$H(\Gamma)$, $E(\Gamma)$, $U(\Gamma)$, and $T(\Gamma)$ are finite sets, respectively called the \emph{set of half-edges} of $\Gamma$, the \emph{set of edges} of $\Gamma$, the \emph{set of univalent vertices} of $\Gamma$, and the set of trivalent vertices of $\Gamma$, 
\item $p_E\colon H(\Gamma) \rightarrow E(\Gamma)$ is a two-to-one map (every element of $E(\Gamma)$ has two preimages under $p_E$) and
\item  $p_V\colon  H(\Gamma) \rightarrow U(\Gamma) \sqcup T(\Gamma)$ is a map such that every element of $U(\Gamma)$ has one preimage under $p_V$ and every element of $T(\Gamma)$ has three preimages under $p_V$, \end{itemize} up to isomorphism.
In other words, $\Gamma$ is a set $H(\Gamma)$ equipped with two partitions, a partition into pairs (induced by $p_E$), and a partition into singletons and triples (induced by $p_V$), up to the bijections that preserve the partitions. These bijections are the {\em automorphisms \/} of $\Gamma$.
\end{definition}

\begin{definition}
\label{defdia}
Let ${\source}$ be a one-manifold, oriented or not.
A 
{\em Jacobi diagram $\Gamma$ with support ${\source}$,\/} also called \emph{Jacobi diagram on ${\source}$}\index[T]{Jacobi diagram}, is a finite uni-trivalent graph $\Gamma$ 
equipped with an isotopy class $\left[i_{\Gamma}\right]$ of injections $i_{\Gamma}$ from the set $U(\Gamma)$ of univalent vertices of $\Gamma$ into the interior of ${\source}$. For such a $\Gamma$, a \emph{$\Gamma$-compatible injection} is an injection in the class $\left[i_{\Gamma}\right]$.
An \emph{orientation}\index[T]{orientation!of a trivalent vertex} of a trivalent vertex of $\Gamma$ is a cyclic order
on the set of the three half-edges that meet at this vertex.
An \emph{orientation}\index[T]{orientation!of a univalent vertex} of a univalent vertex $\eltu$ of $\Gamma$ is an orientation of the connected component ${\source}(\eltu)$ of $i_{\Gamma}(\eltu)$ in ${\source}$, for a choice of $\Gamma$-compatible $i_{\Gamma}$. This orientation associated to $\eltu$ is also called (and thought of as) a \emph{local orientation of ${\source}$ at $\eltu$}.\footnote{A \emph{local orientation} of ${\source}$ at $\eltu$ is simply an orientation of ${\source}(u)$. However, since different vertices are allowed to induce different orientations, we think of these orientations as being \emph{local}, i.e., defined in a neighborhood of $i_{\Gamma}(\eltu)$ for a choice of $\Gamma$-compatible $i_{\Gamma}$.} When ${\source}$ is oriented, the orientation of ${\source}$ orients the univalent vertices of $\Gamma$ naturally.

A \emph{vertex-orientation}\index[T]{vertex-orientation!of a Jacobi diagram} of a Jacobi diagram $\Gamma$ is an \emph{orientation} of every vertex of $\Gamma$.
A Jacobi diagram is \emph{oriented}\index[T]{oriented!Jacobi diagram} if it is equipped with a vertex-orientation.\footnote{When ${\source}$ is oriented, it suffices to specify the orientations of the trivalent vertices since ${\source}$ orients the univalent vertices.}
\end{definition}

Unless otherwise mentioned, the supports of Jacobi diagrams are oriented, and we use the induced orientations of univalent vertices without mentioning them. Nevertheless, the above notion of local orientations will prove useful to state some properties of the invariant $\Zinv$ studied in this book, such as the behavior under cablings in Theorem~\ref{thmmainfunc}.

We represent an oriented Jacobi diagram $\Gamma$ by a planar immersion of $\Gamma \cup {\source}=\Gamma \cup_{U(\Gamma)} {\source}$, with the following conventions.
We represent the (oriented) one-manifold ${\source}$ by dashed lines, the edges of the diagram $\Gamma$ by plain segments, and the vertices by big dots.
The univalent vertices of $U(\Gamma)$ are located at their images under 
a $\Gamma$-compatible injection $i_{\Gamma}$. We represent the orientation of a trivalent vertex by the counterclockwise order of the three half-edges that meet at the vertex.
Here is an example of a picture of a Jacobi diagram $\Gamma$ on the disjoint 
union ${\source}=S^1 \sqcup S^1$ of two (oriented) circles:

\begin{center}
 \exjactwosonebis
\end{center}

The \emph{degree}\index[T]{degree!of a Jacobi diagram} of such a diagram is 
half the number of all its vertices. 
Note that a chord diagram of $\CD_n$ is a degree $n$ Jacobi diagram on $S^1$ without trivalent vertices.
For an (oriented) one-manifold ${\source}$, $\Davis_n({\source})$ denotes the $\KK$-vector space freely generated by the degree $n$ oriented Jacobi diagrams
on ${\source}$. For the (oriented) circle $S^1$, we have
\begin{equation*}\Davis_1(S^1)=\KK \onechord \oplus \KK \lolli \oplus \KK \tatasone \oplus \KK \tatasonetw \oplus \KK \haltsone.\end{equation*}
For an (oriented) one-manifold ${\source}$, $\Aavis_n({\source})$ \index[N]{Am@$\Aavis_n(.)$ target quotiented diagram space of $\Zinvuf_n$} denotes the quotient of $\Davis_n({\source})$ by the following relations AS, Jacobi, and STU:

\begin{center}
AS (or antisymmetry): \trivAS + \trivASop $=0$\\
Jacobi: \ihxone + \ihxtwo + \ihxthree $=0$\\
STU: \stuy $=$ \stui -\stux
\end{center}
Each of these relations relates oriented Jacobi diagrams which are identical outside the pictures.
The quotient $\Aavis_n({\source})$ is the largest quotient of $\Davis_n({\source})$ in which these relations hold. It is obtained by quotienting $\Davis_n({\source})$ by the $\KK$-vector space generated by elements of $\Davis_n({\source})$ of the form \begin{equation*}\biggl(\trivAS + \trivASop\biggr) \mbox{, }  \biggl(\ihxone + \ihxtwo + \ihxthree \biggr)\mbox{, or }
\biggl(\stuy - \stui +\stux\biggr).\end{equation*}

\begin{example} We have
 \begin{equation*}\Aavis_1(S^1)=\KK \onechord \oplus \KK \tatasone.\end{equation*}
\end{example}

\begin{note}
When $\partial {\source} =\emptyset$, any finite-dimensional Lie algebra
equipped with a finite-dimensional representation and a nondegenerate
bilinear symmetric invariant form provides a nontrivial linear map from $\Aavis_n({\source})$ to $\KK$. Such a map is called a \indexT{weight system}. See \cite{barnatan}, \cite[Chapter 6]{chmutovmostov} or \cite[Section 6]{lescol}, for example. In the weight system constructions, the Jacobi relation for the Lie bracket ensures that the maps defined for oriented Jacobi diagrams factor through the Jacobi relation. In \cite{vogel}, Pierre Vogel proved that the maps associated to Lie (super)algebras are sufficient to detect any nontrivial element of $\Aavis_n({\source})$ until degree $15$. He exhibited a nontrivial element of $\Aavis_{16}(\emptyset)$ that such maps cannot detect.
Dror Bar-Natan originally called the Jacobi relation IHX \index[N]{IHX} in \cite{barnatan} because we can write it as 
\begin{equation*}\ihxi=\ihxh - \ihxx,\end{equation*} up to AS.
Note that the four entries in this IHX relation play the same role up to AS.
\end{note}

\begin{definition}
\label{defrkoruniv}
In figures, the orientation of a univalent vertex $\eltu$ of a Jacobi diagram on a nonoriented one-manifold ${\source}$ 
is again represented by the counterclockwise cyclic order of the three half-edges that meet at $\eltu$ in a planar immersion of $\Gamma \cup_{U(\Gamma)} {\source}$, with respect to the following convention. 
The half-edge of $\eltu$ in $\Gamma$ is attached to the left-hand side of ${\source}$, with respect to the local orientation of ${\source}$ at $\eltu$, as in the following picture:
\begin{equation*}\locleftup  \leftrightarrow  \locupor\end{equation*}
In other words, to represent the upward local orientation of $\locnonor$ at $u$, we attach the half-edge of $\eltu$ in $\Gamma$ as in the above figure, and to represent its downward orientation, we attach the half-edge of $\eltu$ as follows: $\locleftdown \leftrightarrow  \locdownor $.

For a nonoriented one-manifold ${\source}$, the space $\Davis_n({\source})$ is the $\KK$-vector space generated by the degree $n$ oriented Jacobi diagrams on ${\source}$---where there are additional orientation choices for univalent vertices, and $\Aavis_n({\source})$ \index[N]{Am@$\Aavis_n(.)$ target quotiented diagram space of $\Zinvuf_n$} is the quotient of $\Davis_n({\source})$ by the previous relations AS, Jacobi, and STU together with the additional antisymmetry relation
 \begin{equation*}\asunor +\asunorc=0,\end{equation*}
where the (unoriented) STU relation may be written as
\begin{center}
STU: \stuyunor $=$ \stuiunor -\stuxunor
\end{center}
\end{definition}

\begin{remark}
We can draw the unoriented STU relation above like the Jacobi relation up to AS:
\begin{equation*}\stubone + \stubtwo + \stubthree=0.\end{equation*}
\end{remark}

\begin{lemma}
\label{lemrkoruniv}
 Let ${\source}$ be an oriented one-manifold. Let ${\source}^f$ denote the nonoriented manifold obtained from ${\source}$ by forgetting the orientation. Let $\Gamma^f$ be an oriented Jacobi diagram on ${\source}^f$.
A univalent vertex $\eltu$ of $\Gamma^f$ is ${\source}$-oriented if it induces the orientation of ${\source}$. Otherwise, it is $(-{\source})$-oriented. Let $\Gamma(\Gamma^f,{\source})$ be the oriented Jacobi diagram on ${\source}$ obtained from $\Gamma^f$ by reversing the local orientation of the $(-{\source})$-oriented univalent vertices. Let $k(\Gamma^f,{\source})$ be the number of $(-{\source})$-oriented univalent vertices of $\Gamma^f$. The linear map from $\Davis_n({\source}^f)$ into $\Davis_n({\source})$ that maps any oriented Jacobi diagram $\Gamma^f$ to $(-1)^{k(\Gamma^f,{\source})}\Gamma(\Gamma^f,{\source})$ and the linear canonical injection from $\Davis_n({\source})$ into $\Davis_n({\source}^f)$ induce canonical isomorphisms between $\Aavis_n({\source})$ and $\Aavis_n({\source}^f)$, which are inverse to each other, for any integer $n \in \NN$.
\end{lemma}
\bp Exercise. \eop

The above lemma justifies using the same notation $\Aavis_n(.)$ for oriented and unoriented supports. We draw Jacobi diagrams on oriented supports by attaching the half-edges of univalent vertices to the left-hand side of the support, to avoid confusion and get rid of the support orientation more easily.

\begin{notation}
\label{notAcheck}
When ${\source} \neq \emptyset$, let $\Assis_n({\source})=\Assis_n({\source};\KK)$ \index[N]{Anotriv@$\Assis_n(.)$ no trivalent component} denote the quotient of $\Aavis_n({\source})=\Aavis_n({\source};\KK)$ by the vector space generated by the diagrams that have at least one connected component without univalent vertices.
\end{notation}
So $\Assis_n({\source})$ is generated by the degree $n$ oriented Jacobi diagrams whose (plain) connected components contain at least one univalent vertex.

\begin{lemma} The space
$\Assis_n(S^1)$ is the quotient by the relations AS and STU of the vector space generated by the degree $n$ oriented Jacobi diagrams whose connected components contain at least one univalent vertex. In other words, the Jacobi relation is a consequence of the relations AS and STU in this vector space.
\end{lemma}
\bp
We want to prove that the Jacobi relation holds in the quotient by AS and STU of the space of uni-trivalent diagrams on $S^1$ with at least one univalent vertex in each connected component. Consider three diagrams represented
by three immersions which coincide outside a disk $D$, inside which they are as in the pictures involved in the Jacobi relation.
Use STU as much as possible to
remove all trivalent vertices that can be removed without changing the two vertices in $D$, on the three
diagrams simultaneously.
This transforms the Jacobi relation to be proved to a sum of similar relations,
where one of the four entries of the disk is directly connected to $S^1$. Since the four entries play the same role in the Jacobi relation, we may assume that the Jacobi relation to be proved is

\begin{equation*} \begin{tikzpicture} \useasboundingbox (-1,-.6) rectangle (.7,.6);
\draw (-.7,.4) node{\scriptsize 1} (-.1,.4) node{\scriptsize 2} (.5,.4) node{\scriptsize 3} (0,-.5) -- (0,0) -- (.3,.4) (0,0) -- (-.3,.4);
\draw (-.9,.4) -- (0,-.25);
\fill (0,-.25) circle (1.5pt) (0,0) circle (1.5pt) (0,-.5) circle (1.5pt);
\draw [->,dash pattern=on 2pt off 2pt] (-.9,-.5) -- (.7,-.5);
\end{tikzpicture}
+ \begin{tikzpicture} \useasboundingbox (-1,-.6) rectangle (.7,.6);
\draw (-.7,.4) node{\scriptsize 1} (-.1,.4) node{\scriptsize 2} (.5,.4) node{\scriptsize 3} (0,-.5) -- (0,0) -- (.3,.4) (0,0) -- (-.3,.4);
\draw (-.9,.4) .. controls (-.2,-.25) .. (0,-.25) .. controls (.2,-.25) and (.3,.05) .. (.15,.2);
\fill (.15,.2) circle (1.5pt) (0,0) circle (1.5pt) (0,-.5) circle (1.5pt);
\draw [->,dash pattern=on 2pt off 2pt] (-.9,-.5) -- (.7,-.5);
\end{tikzpicture}
+ \begin{tikzpicture} \useasboundingbox (-1,-.6) rectangle (.7,.6);
\draw (-.7,.4) node{\scriptsize 1} (-.1,.4) node{\scriptsize 2} (.5,.4) node{\scriptsize 3} (0,-.5) -- (0,0) -- (.3,.4) (0,0) -- (-.3,.4);
\draw (-.9,.4) .. controls (-.2,-.25) .. (0,-.25) .. controls (.2,-.25) and (.3,.05) .. (.15,.2) .. controls (0.05,.3) .. (-.15,.2);
\fill (-.15,.2) circle (1.5pt) (0,0) circle (1.5pt) (0,-.5) circle (1.5pt);
\draw [->,dash pattern=on 2pt off 2pt] (-.9,-.5) -- (.7,-.5);
\end{tikzpicture}
=0.\end{equation*}
Using STU twice and AS transforms the summands of the left-hand side to diagrams that can be represented by three straight lines from the entries $1$, $2$, $3$ to three fixed points of the horizontal line numbered from left to right. When the entry $i \in \{1,2,3\}$ is connected to the point $\sigma(i)$ of the horizontal dashed line, where $\sigma$ is a permutation of $\{1,2,3\}$, the corresponding diagram is denoted by $(\sigma(1)\sigma(2)\sigma(3))$.
Thus, the expansion of the left-hand side of the above equation is
\begin{equation*}\begin{array}{r}((123)-(132)-(231)+(321))\\-((213)-(231)-(132)+(312))\\-((123)-(213)-(312)+(321))\end{array},\end{equation*}
which vanishes, and the lemma is proved.
\eop

\begin{proposition}
\label{propchordtriv}
The natural map from $\CD_n$ to $\Assis_n(S^1)$ induces an isomorphism from the space
$\Asimp_n$ of chord diagrams to $\Assis_n(S^1)$.
\end{proposition}
\bpo{First part of the proof}
According to STU, we have
\begin{equation*}\dquatTtrois - \dquatTun = \dquatSTU= \dquatTqua -\dquatTdeux\end{equation*}
in $\Assis_n(S^1)$. So the natural map from $\CD_n$ to $\Assis_n(S^1)$ factors though $4T$.
Since STU allows us to inductively write any oriented Jacobi diagram whose connected components contain at least one univalent vertex as a combination of chord diagrams, the induced map from $\Asimp_n$ to $\Assis_n(S^1)$ is surjective. We will prove injectivity in Section~\ref{secprod} by constructing an inverse map. \eop

\section{Multiplying diagrams}
\label{secprod}

Set $\Aavis({\source})=\prod_{n \in \NN}\Aavis_n({\source})$, $\Assis({\source})=\prod_{n \in \NN}\Assis_n({\source})$, and $\Asimp=\prod_{n \in \NN}\Asimp_n$.

Assume that a one-manifold ${\source}$ is decomposed as a union of two one-manifolds ${\source} = {\source}_1 \cup {\source}_2$ whose interiors in ${\source}$ do not intersect. Let $\Gamma_1$ be  a Jacobi diagram with support ${\source}_1$ and let $\Gamma_2$ be a Jacobi diagram with support ${\source}_2$. Define the Jacobi diagram $\Gamma_1 \sqcup \Gamma_2$ on ${\source}$ to be the disjoint union of $\Gamma_1$ and $\Gamma_2$. (Formally, a $\Gamma_1 \sqcup \Gamma_2$-compatible injection restricts to $U(\Gamma_i)$ as a $\Gamma_i$-compatible injection, for $i \in \{1,2\}$.)
Define the
{\em product associated to the decomposition ${\source} = {\source}_1 \cup {\source}_2$\/}:
\begin{equation*}\Aavis({\source}_1) \times \Aavis({\source}_2) \longrightarrow \Aavis({\source})\end{equation*}
to be the continuous bilinear map that maps $\left(\left[\Gamma_1\right],\left[\Gamma_2\right]\right)$ to $\left[\Gamma_1 \sqcup \Gamma_2\right]$, for two diagrams $\Gamma_1$ and $\Gamma_2$ as above.
In particular, the disjoint union of diagrams turns $\Aavis(\emptyset)$ into a commutative algebra graded by the degree, and it turns $\Aavis({\source})$ into an $\Aavis(\emptyset)$-module, for any $1$-dimensional manifold ${\source}$.

An orientation-preserving diffeomorphism from a manifold ${\source}$ to another one ${\source}^{\prime}$ induces
a natural isomorphism from $\Aavis_n({\source})$ to $\Aavis_n({\source}^{\prime})$ (and from $\Assis_n({\source})$ to $\Assis_n({\source}^{\prime})$), for all $n$.
Let $I=\left[0,1\right]$ be the compact oriented interval. If $I={\source}$, and if we identify $I$ with ${\source}_1=\left[0,1/2\right]$
and with ${\source}_2=\left[1/2,1\right]$ with respect to the orientation, then
the above process turns $\Aavis(I)$ and $\Assis(I)$ into algebras. The elements of $\Aavis(I)$ with
nonzero degree zero part admit an inverse.

\begin{proposition}
\label{propdiagcom}
The algebra $\Aavis\left(\left[0,1\right]\right)$ is commutative.
 The projection from $\left[0,1\right]$ to $S^1=\left[0,1\right]/(0\sim 1)$ induces an isomorphism from 
$\Aavis_n\left(\left[0,1\right]\right)$ to $\Aavis_n(S^1)$ for all $n$. So $\Aavis(S^1)$ inherits a commutative algebra structure from this isomorphism.
The choice of an oriented connected component ${\source}_j$ of ${\source}$ equips $\Aavis({\source})$ with an \emph{$\Aavis\left(\left[0,1\right]\right)$-module structure $\#_j$}, induced by the orientation-preserving inclusion from $\left[0,1\right]$ to a small part of ${\source}_j$ outside the vertices, and the insertion of diagrams with support $\left[0,1\right]$ there.
\end{proposition}

In order to prove this proposition, we present a useful trick in diagram spaces.

\begin{lemma}
\label{lemcom} Let ${\source}$ be a nonoriented one-manifold.
Let $\Gamma_1$ be an oriented Jacobi diagram (resp. a chord diagram) with support ${\source}$ as in Definitions~\ref{defdia} and \ref{defrkoruniv}. Assume that $\Gamma_1 \cup {\source}$ is immersed in the plane so that  $\Gamma_1 \cup {\source}$ meets an open annulus $A$ embedded in the plane
exactly along $n+1$ embedded arcs $\alpha_1$, $\alpha_2$, \dots, $\alpha_n$, and $\beta$,
and one vertex $v$ on $\alpha_1$, as in the examples below, so that
\begin{enumerate}
\item the $\alpha_i$ are disjoint, they may be dashed or plain (they are dashed in the case of chord diagrams), they go from a boundary component of $A$ to the other,
\item $\beta$ is a plain arc going from the boundary of $A$ to
$v$ on $\alpha_1$,
\item the bounded component $D$ of the complement of $A$ does not contain a boundary point of ${\source}$,
\item the vertex-orientations are induced by the planar immersion by the local counterclockwise orders in the neighborhoods of vertices (as usual).
\end{enumerate}
Let $\Gamma_i$ be the diagram obtained from $\Gamma_1$ by attaching the 
endpoint $v$ of $\beta$ to $\alpha_i$ instead of $\alpha_1$ on the same side,
where
the side of an arc is its side when going from the outside boundary component
of $A$ to the inside one $\partial D$, as in the examples below.
Then we have \begin{equation*}\sum_{i=1}^n\Gamma_i=0\end{equation*} in $\Aavis({\source})$ (resp. in the space $\Asimp$ of Section~\ref{secchordd}).
\end{lemma}
\begin{examples} We have
\begin{equation*}\annulone + \annultwo =0\end{equation*} and
\begin{equation*}\annulthree + \annulfour + \annulfive =0.\end{equation*}
\end{examples}

\begin{remarks}
The second example shows that the STU relation is equivalent to the relation of the statement when the bounded component $D$ of $\RR^2 \setminus A$ intersects $\Gamma_1$ in the neighborhood of a univalent vertex on ${\source}$. Similarly, the Jacobi relation is given by the statement's relation when $D$ intersects $\Gamma_1$ in the neighborhood of a trivalent vertex. 
Also note that AS corresponds to the case in which $D$ intersects $\Gamma_1$ along a
dashed or plain arc.
\end{remarks}
\bpo{Proof of Lemma~\ref{lemcom}}
Let us give the Bar-Natan \cite[Lemma 3.1]{barnatan} proof. See also \cite[Lemma 3.3]{vogel}. Without loss of generality, assume that $v$ is always attached on the left-hand side of 
the $\alpha$'s. 

We first treat the case of trivalent diagrams. Add to the sum the contribution of the sum of the diagrams
obtained from $\Gamma_1$ by attaching $v$ to each of the three (dashed or plain) half-edges of each vertex $w$
of $\Gamma_1 \cup {\source}$ in $D$ on the right-hand side when the half-edges are oriented towards $w$ (i.e., by attaching $v$ to the hooks in \begin{tikzpicture} \useasboundingbox (0,-.15) rectangle (.7,.15);
\draw (.52,0) node{\scriptsize w} (.3,-.25) -- (.3,0) -- (.1,.2) (.3,0) -- (.5,.2) (.4,-.15) -- (.3,-.15)  (.13,.03) -- (.2,.1) (.33,.17) -- (.4,.1);
\fill (.3,0) circle (1.5pt);
\end{tikzpicture}). The Jacobi and STU relations ensure that this contribution is trivial. Now, group the terms of the obtained sum according to the edges of $\Gamma_1 \cup {\source}$,
where $v$ is attached. Observe that the sum is zero, edge by edge, by AS.

For chord diagrams, similarly add to the sum the contribution of the sum of the diagrams
obtained from $\Gamma_1$ by attaching $v$ to each of the four (dashed) half-edges adjacent to each chord $W$ of $\Gamma_1 \cup {\source}$ in $D$, on the right-hand side when the half-edges are oriented towards $W$ (i.e., by attaching $v$ to the hooks in \begin{tikzpicture} \useasboundingbox (-.1,-.25) rectangle (.75,.25);
\draw [dash pattern=on 1.5pt off 1.5pt] (0,-.25) -- (0,.25) (.6,-.25) -- (.6,.25);
\draw (.3,.15) node{\scriptsize W} (0,0) -- (.6,0) (0,-.1) -- (.1,-.1) (.6,-.1) -- (.7,-.1) (0,.1) -- (-.1,.1) (.6,.1) -- (.5,.1);
\fill (0,0) circle (1.5pt) (.6,0) circle (1.5pt);
\end{tikzpicture}). Thanks to the 4T relation, this contribution is trivial. Again, group the terms of the obtained sum according to the dashed edges of $\Gamma_1 \cup {\source}$
where $v$ is attached. Again observe that the sum is zero, edge by edge, by AS.
\eop

\bpo{End of proof of Proposition~\ref{propchordtriv}}

As promised, we construct a map $f$ from $\Assis_n(S^1)$ to the space $\Asimp_n$ of chord diagrams up to 4T and AS, and we prove that it is an inverse of the natural surjective map $g$ from $\Asimp_n$ to $\Assis_n(S^1)$.
Let $\CD_{n,k}$ denote the vector space generated by the oriented uni-trivalent degree $n$ diagrams on $S^1$ that have at most $k$ trivalent vertices and at least one univalent vertex per connected component. Then $\CD_{n,2n}$ is the vector space generated by the oriented uni-trivalent degree $n$ diagrams on $S^1$ with least one univalent vertex per connected component.

We will define linear maps $\lambda_k$ from $\CD_{n,k}$ to $\Asimp_n$ by induction on $k$ so that 
\begin{enumerate}
\item $\lambda_0$ maps a chord diagram to its class in $\Asimp_n$,
\item the restriction of $\lambda_k$ to $\CD_{n,k-1}$ is $\lambda_{k-1}$, and
\item $\lambda_k$ maps all the relations AS and STU that involve only
elements of $\CD_{n,k}$ to zero.
\end{enumerate}
When we have succeeded in such a task, the linear map that maps a diagram $d$ of $\CD_{n,2n}$ with $k$ trivalent vertices
to $\lambda_k(d)$ will factor through STU and AS. Then the induced map
$\overline{\lambda}$ will provide the desired inverse map and allow us to conclude the proof. Let us define
our maps $\lambda_k$ with the announced properties.

Let $k \geq 1$, assume that $\lambda_{k-1}$ is defined on $\CD_{n,k-1}$ and that $\lambda_{k-1}$ maps all the relations AS and STU that involve only
elements of $\CD_{n,k-1}$ to zero. We want to extend $\lambda_{k-1}$
on $\CD_{n,k}$ to a linear map $\lambda_k$ that maps all the relations AS and STU that involve only
elements of $\CD_{n,k}$ to zero.

Let $d$ be a diagram with $k$ trivalent vertices, and let $e$ be an edge of $d$
that contains one univalent vertex and one trivalent vertex.
Set \begin{equation*}\lambda\left((d,e)=\begin{tikzpicture} \useasboundingbox (-.1,-.25) rectangle (.7,.4);
\draw [->,dash pattern=on 2pt off 2pt] (0,-.2) -- (.6,-.2);
\draw (.15,-.1) node{\scriptsize e} (.3,-.2) -- (.3,0) -- (.15,.2) (.3,0) -- (.45,.2);
\fill (.3,-.2) circle (1.5pt) (.3,0) circle (1.5pt);
\end{tikzpicture}\right)
=\lambda_{k-1}\left(
\stui
-
\stux 
\right).\end{equation*}
It suffices to prove that $\lambda(d,e)$ is independent of our chosen edge $e$
to conclude the proof by defining the linear map $\lambda_{k}$, which will obviously satisfy the desired properties, by
\begin{equation*}\lambda_{k}(d)=\lambda(d,e).\end{equation*}

Assume that two different edges $e$ and $f$ of $d$ connect a trivalent vertex to a univalent vertex. We prove that $\lambda(d,e)=\lambda(d,f)$.
If $e$ and $f$ are disjoint, then the fact that $\lambda_{k-1}$ satisfies STU
allows us to express both $\lambda(d,e)$ and $\lambda(d,f)$ as the same combination
of four diagrams with $(k-2)$ vertices and conclude. Thus, we 
assume that $e$ and $f$ are two different edges that share a trivalent vertex $t$.
If there exists another trivalent vertex that is connected to $S^1$ by an edge $g$, then $\lambda(d,e)=\lambda(d,g)=\lambda(d,f)$. Thus, we furthermore assume that $t$ is the unique trivalent vertex connected
to $S^1$ by an edge. If $t$ is the unique trivalent vertex, then its component is necessarily like 
\dquatSTU
and the fact that $\lambda(d,e)=\lambda(d,f)$ is a consequence of 
(4T). Otherwise, the component of $t$ is of the form
\begin{tikzpicture} \useasboundingbox (-.1,-.25) rectangle (1.1,.65);
\draw [->,dash pattern=on 2pt off 1pt] (0,-.2) -- (.4,-.2);
\draw [->,dash pattern=on 2pt off 1pt] (.6,-.2) -- (1,-.2);
\draw (.4,.2) node {\scriptsize t} (.1,-.05) node {\scriptsize e} (.9,0) node {\scriptsize f} (.2,-.2) -- (.5,.1) -- (.8,-.2) (.5,.1) -- (.5,.3);
\draw [dotted] (.5,.45) circle (.15);
\fill (.2,-.2) circle (1.5pt) (.5,.1) circle (1.5pt) (.8,-.2) circle (1.5pt);
\end{tikzpicture} where the dotted circle represents a dashed diagram with only one pictured entry.
So we have \begin{equation*}\lambda(d,e)=\lambda_{k-1}\left(
\begin{tikzpicture} \useasboundingbox (-.1,-.25) rectangle (1.1,.55);

\draw [->,dash pattern=on 2pt off 1pt] (0,-.2) -- (.5,-.2);
\draw [->,dash pattern=on 2pt off 1pt] (.7,-.2) -- (1,-.2);
\draw (.15,-.2) -- (.15,.2);
\draw [out=90,in=90] (.35,-.2) to (.85,-.2);
\draw [dotted] (.15,.35) circle (.15);
\fill (.15,-.2) circle (1.5pt) (.35,-.2) circle (1.5pt) (.85,-.2) circle (1.5pt);
\end{tikzpicture}
-
\begin{tikzpicture} \useasboundingbox (-.1,-.25) rectangle (1.1,.55);
\draw [->,dash pattern=on 2pt off 1pt] (0,-.2) -- (.5,-.2);
\draw [->,dash pattern=on 2pt off 1pt] (.7,-.2) -- (1,-.2);
\draw (.35,-.2) -- (.35,.2);
\draw [out=90,in=90] (.15,-.2) to (.85,-.2);
\draw [dotted] (.35,.35) circle (.15);
\fill (.15,-.2) circle (1.5pt) (.35,-.2) circle (1.5pt) (.85,-.2) circle (1.5pt);
\end{tikzpicture}\right).\end{equation*}

Then $\lambda(d,e)$ is zero because the expansion of \begin{tikzpicture} \useasboundingbox (-.1,-.25) rectangle (.6,.45);
\draw [->,dash pattern=on 2pt off 1pt] (0,-.2) -- (.5,-.2);
\draw (.25,-.2) -- (.25,.2);
\draw [dotted] (.25,.35) circle (.15);
\fill (.25,-.2) circle (1.5pt);
\end{tikzpicture}
as a sum of chord diagrams commutes with any vertex in $\Asimp_n$, according to Lemma~\ref{lemcom}.
We similarly have $\lambda(d,f)=0$. So we get $\lambda(d,e)=\lambda(d,f)$ in this last case.
\eop

\begin{lemma}
\label{lemoneleg}
For any one-manifold ${\source}$, the class of a Jacobi diagram with one univalent vertex vanishes in $\Aavis_n({\source})$.
\end{lemma}
\bp Exercise. (Use Lemma~\ref{lemcom}.)
\eop

\bpo{Proof of Proposition~\ref{propdiagcom}}
Let $\Gamma^{\prime}$ be a diagram with support ${\source}$ and let $\Gamma$ is a diagram with support $I$. Define $\left[\Gamma\right]\#_j\left[\Gamma^{\prime}\right]$ to be the class in $\Aavis({\source})$ of the diagram obtained by inserting $\Gamma$ along ${\source}_j$ outside the vertices of $\Gamma$, according to the given orientation.
For example, we have \begin{equation*}\left[\dmer\right] \#_j\left[\cirt\right] =\left[ \cirtd\right] =\left[\cirtdp\right] .\end{equation*}
As shown in the first example that illustrates  Lemma \ref{lemcom}, the independence of the choice of the insertion locus is a consequence of Lemma~\ref{lemcom}, where $\Gamma_1$ is the disjoint union $\Gamma \sqcup \Gamma^{\prime}$, and $\Gamma_1$ intersects $D$ along $\Gamma \cup I$.
The continuous bilinear map \begin{equation*} \begin{array}{lll} \Aavis(I) \times \Aavis({\source}) & \longrightarrow &\Aavis({\source}) \\
                               
   \left(  \left[\Gamma\right],\left[\Gamma^{\prime}\right]   \right)              & \mapsto &         \left[\Gamma\right]\#_j\left[\Gamma^{\prime}\right]\end{array}\end{equation*}
endows $\Aavis({\source})$ with the \emph{$\Aavis(I)$-module structure} $\#_j$.

Lemma \ref{lemcom} allows us to similarly prove that $\Aavis(I)$ is a commutative algebra.
Now, it suffices to prove that the morphism from $\Aavis(I)$ to $\Aavis(S^1)$ induced by the identification of the two endpoints of $I$ is an isomorphism. This is proved in the more general proposition below.
\eop

\begin{proposition}
\label{propdiagrs}
Let $n \in \NN$. Let $\source$ be a disjoint union of circles.
The projection from $\left[0,1\right]$ to $S^1=\left[0,1\right]/(0\sim 1)$ induces an isomorphism from 
$\Aavis_n\left(\left[0,1\right] \sqcup \source \right)$ to $\Aavis_n(S^1 \sqcup \source)$.
\end{proposition}
\bp The morphism from $\Aavis\left(\left[0,1\right]\sqcup \source\right)$ to $\Aavis(S^1\sqcup \source)$ induced by the identification of the two endpoints of $\left[0,1\right]$ amounts to mod out $\Aavis\left(\left[0,1\right]\sqcup \source\right)$ by the relation that identifies two diagrams obtained from one another by moving the nearest univalent vertex to an endpoint of $\left[0,1\right]$ near the other endpoint. Applying Lemma~\ref{lemcom} (with $\beta$ coming from the inside boundary of the annulus) shows that this relation is a consequence of the relations in $\Aavis\left(\left[0,1\right]\sqcup \source\right)$. So this morphism is an isomorphism from $\Aavis\left(\left[0,1\right]\sqcup \source\right)$ to $\Aavis\left(S^1\sqcup \source\right)$.
\eop

As the following exercise shows, Proposition~\ref{propdiagrs} would not be satisfied if $\source$ were replaced by an interval.
\begin{exo}
\label{exocounterexcut}
1. Prove 
\begin{equation*}\left[\ddparallel \right] \neq \left[\ddcross \right].\end{equation*}
in $\Aavis_2(\dasheddownw \dashedupw)$.\\
\noindent 2. Prove that $\Aavis_2\left(\dasheddownw \dashedupw \right)\ncong \Aavis_2\left(\dashedcircleorrrevseul \dashedupw \right)$.
\end{exo}

\begin{lemma}
\label{lemlifdiag}
Let $\pi:{\source}^{\prime} \longrightarrow {\source}$ be a smooth map between
two unoriented compact one-manifolds ${\source}$ and ${\source}^{\prime}$. Assume $\pi(\partial {\source}^{\prime}) \subset \partial {\source}$.
Let $\Gamma$ be an oriented Jacobi diagram on ${\source}$ (as in Definitions~\ref{defdia} and \ref{defrkoruniv}) equipped with a $\Gamma$-compatible injection whose image avoids the critical values of $\pi$.
Define $\pi^{\ast}(\Gamma)$ to be the class in $\Aavis({\source}^{\prime})$ of the sum of all diagrams on ${\source}^{\prime}$ obtained from $\Gamma$ by lifting
each univalent vertex to one of its preimages under $\pi$. (These diagrams have the
same vertices and edges as $\Gamma$, and the local orientations at univalent vertices are naturally induced by the local orientations of the corresponding univalent vertices of $\Gamma$.) 
Then $\pi^{\ast}(\Gamma)$ depends only on the class of $\Gamma$ in $\Aavis({\source})$ and on the homotopy class of $\pi$.
So one can unambiguously define the linear degree-preserving map \begin{equation*}\pi^{\ast}\colon\Aavis({\source}) \longrightarrow  \Aavis({\source}^{\prime})\end{equation*}
that maps the class of a diagram $\Gamma$ as above to $\pi^{\ast}(\Gamma)$.
Furthermore, the map $\pi^{\ast}$ depends only on the homotopy class of $\pi$ relative to the boundary of ${\source}^{\prime}$.
\end{lemma}
\bp It suffices to see that this operation is compatible with STU and that nothing bad happens when a univalent vertex of ${\source}$ moves across a critical value of $\pi$.
\eop

\begin{remark}
 Propositions~\ref{propdiagcom} and \ref{propdiagrs} and Lemma~\ref{lemlifdiag} remain valid if $\Aavis$ is replaced by $\Assis$.
\end{remark}

\begin{notation}
\label{notationduplication}
Let ${\source}$ be a one-manifold. Let ${\source}_0$ be a connected component of ${\source}$.
Let \begin{equation*}{\source}(r \times {\source}_0) = ({\source} \setminus {\source}_0) \sqcup \left(\sqcup_{i=1}^r {\source}_0^{(i)}\right)\end{equation*} be the manifold
obtained from ${\source}$ by \emph{duplicating} ${\source}_0$ (r-1) times, that is by replacing ${\source}_0$ by r copies of ${\source}_0$. Let $\pi(r \times {\source}_0) \colon {\source}(r \times {\source}_0) \longrightarrow {\source}$ be the associated map, which is the identity on $({\source} \setminus {\source}_0)$, and
the trivial r-fold covering from $\sqcup_{i=1}^r {\source}_0^{(i)}$ to ${\source}_0$.
The associated map is the \emph{duplication map:}
\begin{equation*}\pi(r \times {\source}_0)^{\ast}\colon \CA({\source}) \longrightarrow \CA({\source}(r \times {\source}_0)).\end{equation*}
\end{notation}

\begin{example} We have
\begin{equation*}\pi(2 \times I)^{\ast} \left(\onechordonverticalinterval \right) = 
\begin{tikzpicture} \useasboundingbox (-.2,-.2) rectangle (.4,.2);
\draw (0,-.15) .. controls (-.25,-.15) and (-.25,.15) .. (0,.15);
\draw [dash pattern=on 2pt off 2pt]  (0,-.3) -- (0,.3);
\draw [dash pattern=on 2pt off 2pt]  (.3,-.3) -- (.3,.3);
\fill (0,-.15) circle (1.5pt) (0,.15) circle (1.5pt);
\end{tikzpicture}
+
\begin{tikzpicture} \useasboundingbox (-.2,-.2) rectangle (.4,.2);
\draw (.3,-.15) .. controls (.2,-.15) and (-.25,.15) .. (0,.15);
\draw [dash pattern=on 2pt off 2pt]  (0,-.3) -- (0,.3);
\draw [dash pattern=on 2pt off 2pt]  (.3,-.3) -- (.3,.3);
\fill (.3,-.15) circle (1.5pt) (0,.15) circle (1.5pt);
\end{tikzpicture}
+
\begin{tikzpicture} \useasboundingbox (-.2,-.2) rectangle (.4,.2);
\draw (0,-.15) .. controls (-.25,-.15) and (.2,.15) .. (.3,.15);
\draw [dash pattern=on 2pt off 2pt]  (0,-.3) -- (0,.3);
\draw [dash pattern=on 2pt off 2pt]  (.3,-.3) -- (.3,.3);
\fill (0,-.15) circle (1.5pt) (.3,.15) circle (1.5pt);
\end{tikzpicture}
+
\begin{tikzpicture} \useasboundingbox (-.2,-.2) rectangle (.4,.2);
\draw (.3,-.15) .. controls (.05,-.15) and (.05,.15) .. (.3,.15);
\draw [dash pattern=on 2pt off 2pt]  (0,-.3) -- (0,.3);
\draw [dash pattern=on 2pt off 2pt]  (.3,-.3) -- (.3,.3);
\fill (.3,-.15) circle (1.5pt) (.3,.15) circle (1.5pt);
\end{tikzpicture} .
\end{equation*}
\end{example}

Note the following lemma.

\begin{lemma}
\label{lemdupcom}
 When $\source$ is a disjoint union of $r$ intervals, an $r$-duplicated vertex commutes with an element of $\CA(\source)$.
This sentence is explained by the pictures below.
In the first picture
\begin{center}
\begin{tikzpicture}[scale=.6] 
\draw (0,1) -- (4,1) -- (4,2) -- (0,2) -- (0,1) (.8,.6) -- (1,.6)  (1.3,.6) -- (1.5,.6) (1.8,.6) -- (2,.6) (2.3,.6) -- (2.5,.6) (2.8,.6) -- (3,.6) (-1,1.5) -- (-.5,1.5) circle (3pt);
\draw[dashed,dash pattern=on 3pt off 1pt] (1,.2) -- (1,1) (1.5,.2) -- (1.5,1) (2,.2) -- (2,1) (2.5,.2) -- (2.5,1) (3,.2) -- (3,1)  (1,2) -- (1,2.8) (1.5,2) -- (1.5,2.8) (2,2) -- (2,2.8) (2.5,2) -- (2.5,2.8) (3,2) -- (3,2.8) (4.5,1.5) node{$=$};
\begin{scope}[xshift=6.3cm]
\draw (0,1) -- (4,1) -- (4,2) -- (0,2) -- (0,1) (.8,2.4) -- (1,2.4)  (1.3,2.4) -- (1.5,2.4) (1.8,2.4) -- (2,2.4) (2.3,2.4) -- (2.5,2.4) (2.8,2.4) -- (3,2.4) (-1,1.5) -- (-.5,1.5) circle (3pt);
\draw[dashed,dash pattern=on 3pt off 1pt] (1,.2) -- (1,1) (1.5,.2) -- (1.5,1) (2,.2) -- (2,1) (2.5,.2) -- (2.5,1) (3,.2) -- (3,1)  (1,2) -- (1,2.8) (1.5,2) -- (1.5,2.8) (2,2) -- (2,2.8) (2.5,2) -- (2.5,2.8) (3,2) -- (3,2.8);
\end{scope}
\end{tikzpicture}\end{center}
there is a Jacobi diagram on $\source$ inside the rectangle, and the picture represents the sum of the diagrams obtained by attaching the free end of an edge (the end with the empty circle) of some other part of a Jacobi diagram to each of the hooks attached to the vertical strands. The second picture
\begin{center}
\begin{tikzpicture}[scale=.6]
\begin{scope}
\draw (0,1) -- (4,1) -- (4,2) -- (0,2) -- (0,1) (.8,.6) -- (1,.6)  (1.3,.6) -- (1.5,.6) (1.8,.6) -- (2,.6) (2.3,.6) -- (2.5,.6) (2.8,.6) -- (3,.6) (0,1.5) -- (-.5,1.5) circle (3pt) (4.5,1.5) node{$=$};
\draw[dashed,dash pattern=on 3pt off 1pt]  (1,.2) -- (1,1) (1.5,.2) -- (1.5,1) (2,.2) -- (2,1) (2.5,.2) -- (2.5,1) (3,.2) -- (3,1)  (1,2) -- (1,2.8) (1.5,2) -- (1.5,2.8) (2,2) -- (2,2.8) (2.5,2) -- (2.5,2.8) (3,2) -- (3,2.8);
\end{scope}
\begin{scope}[xshift=6cm]
\draw (0,1) -- (4,1) -- (4,2) -- (0,2) -- (0,1) (.8,2.4) -- (1,2.4)  (1.3,2.4) -- (1.5,2.4) (1.8,2.4) -- (2,2.4) (2.3,2.4) -- (2.5,2.4) (2.8,2.4) -- (3,2.4) (0,1.5) -- (-.5,1.5) circle (3pt);
\draw[dashed,dash pattern=on 3pt off 1pt]  (1,.2) -- (1,1) (1.5,.2) -- (1.5,1) (2,.2) -- (2,1) (2.5,.2) -- (2.5,1) (3,.2) -- (3,1)  (1,2) -- (1,2.8) (1.5,2) -- (1.5,2.8) (2,2) -- (2,2.8) (2.5,2) -- (2.5,2.8) (3,2) -- (3,2.8);
\end{scope}
\end{tikzpicture}
\end{center}
 is similar, except that the edge with the free end is a part of a Jacobi diagram that is inside the box apart from this half-edge.
\end{lemma}
\bp This is a direct consequence of Lemma~\ref{lemcom} as the pictures show.
\eop

\section{Coproduct on \texorpdfstring{$\Aavis({\source})$}{diagram spaces}}
\label{seccoprod}

Recall that $\KK$ is $\QQ$ or $\RR$ in this chapter.
All the results are valid over a commutative field $\KK$ of characteristic zero.
Below, all tensor products are over the ground field $\KK$. The canonical identification of $V \otimes \KK$ with $V$ for a finite-dimensional vector space over $\KK$ will always be implicit. In this section, ${\source}$ denotes a one-manifold.

For $n \in \NN$, set \begin{equation*}(\Aavis({\source}) \otimes  \Aavis({\source}))_n=\oplus_{i=0}^n \Aavis_i({\source}) \otimes \Aavis_{n-i}({\source})\end{equation*}
and \begin{equation*}\Aavis({\source}) \hat{\otimes}  \Aavis({\source}) =\prod_{n \in \NN}(\Aavis({\source}) \otimes  \Aavis({\source}))_n.\end{equation*}
The topological product $\Aavis({\source})$ is equipped with the following collection
of linear maps
\begin{equation*}\Delta_n \colon \Aavis_n({\source}) \to (\Aavis({\source}) \otimes  \Aavis({\source}))_n.\end{equation*}
The image of the class of a Jacobi diagram $\Gamma=\sqcup_{i \in I}\Gamma_i$ with $\cardlef{I}$ nonempty connected components $\Gamma_i$, numbered arbitrarily in a set $I$, is
\begin{equation*}\Delta_n\left(\left[\Gamma\right]\right)=\sum_{J \subseteq I} \left[\sqcup_{i \in J}\Gamma_i\right] \otimes \left[\sqcup_{i \in (I\setminus J)}\Gamma_i\right].\end{equation*}
Here, the connected components are not the connected components of $\Gamma \cup_{U(\Gamma)} {\source}$. They are the connected components of $\Gamma$, i.e., the connected components of the solid (i.e., nondashed) part of the figures.
It is easy to check that $\Delta_n$ is well-defined.
The family $\Delta=(\Delta_n)_{n\in \NN}$ defines a degree-preserving map from $\Aavis({\source})$ to 
$\Aavis({\source}) \hat{\otimes}  \Aavis({\source})$.

There is a well-defined continuous linear map $\varepsilon \colon \Aavis({\source}) \to \KK$ that
maps the class of the empty diagram to $1$ and $\Aavis_i({\source})$ to $0$ for any $i>0$. The ground field $\KK$ is considered as a degree $0$ vector space. So the map $\varepsilon$ is a degree-preserving homomorphism. Recall that $\id$ denotes the identity map.

In the following statement, we omit degrees. However, the following  identities, which express the fact that $\Delta$ is a graded \indexT{coproduct} with associated \indexT{counit} $\varepsilon$, are collections of identities between collections of degree-preserving linear maps between finite-dimensional vector spaces. For example, the \emph{coassociativity} identity
\begin{equation*}(\Delta  \otimes \id) \circ \Delta = (\id  \otimes \Delta)\circ \Delta\end{equation*} means that
\begin{equation*}(\Delta  \otimes \id)_n \circ \Delta_n = (\id  \otimes \Delta)_n\circ \Delta_n\end{equation*}
for any $n \in \NN$,
where both maps are valued in \begin{equation*}(\Aavis({\source}) \otimes  \Aavis({\source}) \otimes  \Aavis({\source}))_n = \oplus_{i,j,k \suchthat (i,j,k)\in \NN^3, i+j+k=n} \Aavis_i({\source}) \otimes \Aavis_{j}({\source}) \otimes \Aavis_{k}({\source}).\end{equation*}
\begin{lemma} We have
 $(\varepsilon \otimes \id) \circ \Delta= (\id \otimes \varepsilon) \circ \Delta=\id$ and
\begin{equation*}(\Delta  \otimes \id) \circ \Delta = (\id  \otimes \Delta)\circ \Delta.\end{equation*}
\end{lemma}
\bp Exercise. \eop

Let \begin{equation*}\begin{array}{llll}\tau_n \colon & (\Aavis({\source}) \otimes  \Aavis({\source}))_n&\to &(\Aavis({\source}) \otimes  \Aavis({\source}))_n\\
      &x\otimes y &\mapsto & y \otimes x.
     \end{array}\end{equation*}
Then we also immediately have the identity
\begin{equation*}\tau \circ \Delta =\Delta,\end{equation*}
which expresses the \emph{cocommutativity} of $\Delta$.

\section{Bialgebra structures}
\label{secHopfalg}

\begin{definition}
A \emph{connected, finite type, commutative, cocommutative, graded bialgebra} over a field $\KK$ is the topological product $\CH=\prod_{n\in \NN}\CH_n$ of finite-dimensional vector spaces $\CH_n$ over $\KK$
equipped with
\begin{itemize}
 \item a \emph{multiplication} $m=(m_n \colon (\CH \otimes  \CH)_n \to \CH_n)_{n \in \NN}$, where $(\CH \otimes  \CH)_n = \oplus_{i=0}^n\CH_i \otimes \CH_{n-i}$,
\item a \emph{coproduct} $\Delta=(\Delta_n \colon \CH_n \to (\CH \otimes  \CH)_n)_{n \in \NN}$,
\item a \emph{unit} $\upsilon \colon \KK \to \CH$, which maps $\KK$ to $\CH_0$, and which is an isomorphism from $\KK$ to $\CH_0$ (connectedness),
\item a \emph{counit} $\varepsilon \colon \CH \to \KK$, where $\KK$ is again assumed to be of degree $0$,
\end{itemize}
where $(m,\Delta,\upsilon,\varepsilon)$ are families of degree-preserving linear maps 
that satisfy 
\begin{itemize}
\item the following identities, which express that $(m,\upsilon)$ is an associative and commutative product with unit $\upsilon(1)$:
\begin{equation*} \begin{array}{l}
m \circ (m \otimes \id)= m \circ (\id \otimes m)\\
m \circ (\upsilon \otimes \id)=m \circ ( \id \otimes \upsilon)=\id\\
m \circ \tau=m, 
   \end{array}
 \end{equation*}
 where $\tau_n\colon (\CH \otimes  \CH)_n \to (\CH \otimes  \CH)_n$ maps $x \otimes y$ to $y \otimes x$.
\item the following identities, which express that $(\Delta,\varepsilon)$ is a coassociative and cocommutative coproduct with counit $\varepsilon$:
\begin{equation*} \begin{array}{l}  
(\Delta  \otimes \id) \circ \Delta = (\id  \otimes \Delta)\circ \Delta\\
(\varepsilon \otimes \id) \circ \Delta= (\id \otimes \varepsilon) \circ \Delta=\id\\
\tau \circ \Delta =\Delta.
\end{array}
 \end{equation*}
\item the following \emph{compatibility identity,} which expresses the fact that $\Delta$ is an algebra morphism and that $m$ is a coalgebra morphism
\begin{equation*}\Delta \circ m =(m \otimes m) \circ (\id \otimes \tau \otimes \id) \circ (\Delta \otimes \Delta),\end{equation*}
where the product on $\CH \hat{\otimes} \CH=\prod_{n\in \NN}(\CH \otimes \CH)_n$ is defined from $m$ 
so that it maps $(a \otimes b) \otimes (a^{\prime} \otimes b^{\prime})$ to $m(a \otimes a^{\prime}) \otimes m(b \otimes b^{\prime})$.
\end{itemize}
\end{definition}

\begin{lemma}\label{lemcoprodone}
 In a connected, finite type, commutative, cocommutative, graded bialgebra, we have $\varepsilon \circ \upsilon =\id$ and 
$\Delta(\upsilon(1))=\upsilon(1) \otimes \upsilon(1)$. The element $\upsilon(1)$ is denoted by $1$.
\end{lemma}
\bp Since $\Delta$ is degree-preserving, there exists $k \in \KK$ such that $\Delta(\upsilon(1))= k\upsilon(1) \otimes \upsilon(1)$.
Since $(\varepsilon \otimes \id) \circ \Delta=\id$,
we have $k\varepsilon(\upsilon(1))=1$.
Applying the compatibility identity to 
$\upsilon(1) \otimes x$ yields $\Delta(x)=k\Delta(x)$. So we get $k=1=\varepsilon(\upsilon(1))$.
\eop

In a connected, finite type, commutative, cocommutative, graded bialgebra, a \indexT{primitive} element is an element such that
$\Delta(x)=1\otimes x+ x \otimes 1$,
and a \indexT{group-like} element is an element such that
$\Delta(x)=x\otimes x$ and $\varepsilon(x) \neq 0$.

The proof of the following lemma is straightforward and left to the reader.
\begin{lemma}
\label{lemAHopf}
Equipped with the product of Section~\ref{secprod}, with the coproduct of Section~\ref{seccoprod}, and with the counit that maps the class of the empty diagram to $1$, the graded vector spaces $\Assis(S^1)$, $\Aavis(S^1)$, and $\Aavis(\emptyset)$ are connected, finite type, commutative, cocommutative, graded bialgebras.\footnote{These spaces can be furthermore equipped with the linear \emph{antipode} that maps a product $\Pi$ of $p$ primitive elements to $(-1)^p \Pi$. So they are \emph{Hopf algebras}.} The unit $\upsilon(1)$ of these algebras is the class of the empty diagram. Furthermore, connected Jacobi diagrams are primitive elements in these algebras.
\end{lemma}
\eopwobp

Note the elementary lemma.

\begin{lemma}
\label{lemexpprimgrouplike}
If $y$ is a primitive element of a connected, finite type, commutative, cocommutative, graded bialgebra, then $\exp(y)$ is group-like.
\end{lemma}
\bp
It suffices to prove that $\Delta(y^n)=\sum_{k=0}^n \frac{n!}{k!(n-k)!}y^k \otimes y^{n-k}$.
When $n=0$, this is Lemma~\ref{lemcoprodone}.
The compatibility identity implies $\Delta(y^n)=\Delta(y^{n-1})(y\otimes \upsilon(1) +\upsilon(1) \otimes y)$.
\eop

We can  now state a version of the Milnor-Moore theorem.
\begin{theorem}
\label{thmmilmoo} 
Let $(\CH;m,\Delta,\upsilon,\varepsilon)$ be a connected, finite type, commutative, cocommutative, graded bialgebra over a field $\KK$ of characteristic zero.
 Let $\CP_n$ denote the set of primitive elements of $\CH_n$. It is a finite-dimensional vector space.  Pick a basis $b_n$ of each $\CP_n$, for each $n$.
For all $n \in \NN$, the space $\CH_n$ is the vector space freely generated by the degree $n$ monic monomials in the elements of $b_{\leq n}=\cup_{k \in \NN \suchthat k\leq n}{b_k}$.\footnote{By \emph{monic monomials}, we mean monomials with coefficient one. So degree $n$ monic monomials in the elements of $b_{\leq n}$ are of the form $\prod_{i \in I}p_i^{r(i)}$, for elements $p_i$ of $b_{\leq n}$ of degree $d(i)$, and positive integers $r(i)$ such that $\sum_{i \in I}d(i){r(i)}=n$.}
\end{theorem}
\bp
Since $\CP_0=\{0\}$ and $\CH_0$ is freely generated by the monomial $\upsilon(1)$, the theorem holds for $n=0$.
Let $n \geq 1$, let $d_n$ be the set of degree $n$ monic monomials in the elements of $b_{\leq n-1}$.
We want to prove that $\CH_n$ is freely generated by $d_n \sqcup b_n$, by induction on $n$. 

For $x \in \CH_n$,
set $\Delta^{\prime}(x)= \Delta(x)-x\otimes \upsilon(1) - \upsilon(1)\otimes x$.
According to Lemma~\ref{lemcoprodone}, since $(\varepsilon \otimes \id) \circ \Delta= (\id \otimes \varepsilon) \circ \Delta=\id$, we have
\begin{equation*}\Delta^{\prime}(x) \in \oplus_{i=1}^{n-1} \CH_i \otimes \CH_{n-i}\end{equation*}
and $\CP_n$ is the kernel of $\Delta^{\prime}$.

By induction, $\CH_i$ (resp. $\CH_{n-i}$) has a basis consisting of degree $i$ (resp. $(n-i)$) monic monomials in the elements of $b_{\leq n-1}$.
Thus, $\CH_i \otimes \CH_{n-i}$ has a basis consisting of tensor products of these monomials.
Multiplying two such monomials yields an element \begin{equation*}\prod_{i \in I}p_i^{r(i)}\end{equation*} of $d_n$, where the $p_i$ are distinct elements of $b_{\leq n-1}$, and the $r(i)$ are positive integers. We have
\begin{multline*}
 \Delta\biggl(\prod_{i \in I}p_i^{r(i)}\biggr)=\\
 \sum_{k \colon I \to \NN \suchthat  0\leq k(i) \leq r(i), \forall i}\biggl(\prod_{i \in I} \frac{r(i)!} {k(i)!(r(i)-k(i))!}\biggr)\biggl(\prod_{i \in I} p_i^{k(i)} \biggr)\otimes \biggl(\prod_{i \in I}p_i^{r(i)-k(i)} \biggr).
\end{multline*}

This formula proves that $\Delta^{\prime}$ injects the vector space freely generated by the degree $n$ monic monomials in the elements of $b_{\leq n-1}$ into $\oplus_{i=1}^{n-1} \CH_i \otimes \CH_{n-i}$.
Thus, the degree $n$ monic monomials in the elements of $b_{\leq n}$ form a free system of $\CH_n$, and it suffices to prove that they generate $\CH_n$.
To do so, we only need to check that for every $x \in \CH_n$, for every element $d \in d_n$, there exists a constant $a(x,d)$ such that \begin{equation*}\Delta^{\prime}(x)=\sum_{d \in d_n}a(x,d)\Delta^{\prime}(d).\end{equation*}
Indeed, in this case, $\left(x -\sum_{d \in d_n}a(x,d)d\right)$ is primitive. Fix $x \in \CH_n$.

Let $d=\prod_{i \in I}p_i^{r_d(i)} \in \CH_n$, where $r_d(i)>0$, for any  $i \in I$. Let $E(d)$ be the set of maps $k\colon I \to \NN$ such that $0< \sum_{i \in I} k(i) $, $k(i)\leq r_d(i)$, and $d(k)\stackrel{\mbox{\scriptsize def}}{=}\prod_{i \in I} p_i^{k(i)} \in  \CH_{j}$ for some $j$ such that $j<n$. Then we have
\begin{equation*}\Delta^{\prime}(d)=\sum_{k\in E(d)}c_{k,r_d}d(k) \otimes d(r_d-k)\end{equation*} with
\begin{equation*}c_{k,r_d}=\prod_{i \in I} \frac{r_d(i)!} {k(i)!(r_d(i)-k(i))!}\neq 0.\end{equation*}
We also have \begin{equation*}\Delta^{\prime}(x)=\sum_{d \in d_n}\Bigl(\sum_{k\in E(d)}c_k(x,d) d(k) \otimes d(r_d-k)\Bigr).\end{equation*}
Now, it suffices to prove that, for any $d \in d_n$, there exists $a(x,d)$ such that $c_k(x,d)=a(x,d)c_{k,r_d}$ for any $k \in E(d)$. Fix $d \in d_n$ and set $r=r_d$.

Thanks to the coassociativity of $\Delta$, we have $(\Delta  \otimes \id) \circ \Delta(x) = (\id  \otimes \Delta)\circ \Delta(x)$.
Therefore, if $h(i) \leq k(i)$ for all $i\in I$ and $\sum_{i \in I}h(i)>0$, then the coefficient 
\begin{equation*}c_{h,k}c_k(x,d)=c_h(x,d)c_{k-h,r-h}\end{equation*}
of $d(h)\otimes d(k-h)\otimes d(r-k)$ in $(\Delta  \otimes \id) \circ \Delta(x)$ determines both the coefficient $c_h(x,d)$ of $d(h)\otimes d(r-h)$ and the coefficient $c_k(x,d)$ in $\Delta^{\prime}(x)$. So, we have
\begin{equation*}c_h(x,d)=\frac{c_{h,k}}{c_{k-h,r-h}}c_k(x,d).\end{equation*}
The coassociativity applied to $d$ similarly implies \begin{equation*}c_h(d,d)=\frac{c_{h,k}}{c_{k-h,r-h}}c_k(d,d),\end{equation*} where $c_h(d,d)=c_{h,r}$. So, we get
\begin{equation*}c_h(x,d)=\frac{c_{h,r}}{c_{k,r}}c_k(x,d).\end{equation*}

Choose $j \in I$. Let $\delta_{j} \in E(d)$ be such that $\delta_{j}(j)=1$ and $\delta_{j}(i)=0$ for all $i \in I \setminus \{j\}$.
Set $a(x,d)={c_{\delta_{j}}(x,d)}/{c_{\delta_{j},r}}$. Then for all $k \in E(d)$ such that $k(j)\neq 0$, we have $c_k(x,d)=a(x,d)c_{k,r}$.
If $\sum_{i \in I} r(i)>2$, then for any $i \in I \setminus \{j\}$, the map $\delta_i +\delta_{j}$ is in $E(d)$. So we have $c_{\delta_i}(x,d)=a(x,d)c_{\delta_i,r}$ for any $i \in I$, and therefore $c_k(x,d)=a(x,d)c_{k,r}$ for any $k \in E(d)$.
The only untreated case is $I=\{i,j\}$ with $r(i)=r(j)=1$. In this case, the cocommutativity of $\Delta$ leads to the result.
\eop

\begin{corollary}
\label{corprojprim}
Under the assumptions of Theorem~\ref{thmmilmoo}, there is a well-defined unique linear projection $p^c$ from $\CH$ to $\CP=\prod_{n\in \NN}\CP_n$ that maps the products of two homogeneous elements of positive degree to $0$, and that maps $\CH_0$ to $0$.\footnote{An element of $\CH$ is \emph{homogeneous} if it belongs to some $\CH_j$ for some $j$.} 
\end{corollary}
\eopwobp

\begin{theorem}\label{thmgrouplikexp} Let $\CH$ be a connected, finite type, commutative, cocommutative, graded bialgebra.
Let $\CP$ be the space of its primitive elements, and let $p^c \colon \CH \to \CP$ be the projection of Corollary~\ref{corprojprim}.
Any group-like element $x$ of $\CH$ is the exponential of a unique primitive element of $\CP$.
This element is $p^c(x)$.
\end{theorem}
\bp First note that $p^c(\exp(y))=y$ for any primitive element $y$ of $\CH$. Therefore, if $\exp(y)=\exp(y^{\prime})$ for two primitive elements $y$ and $y^{\prime}$ of $\CH$, then $y=y^{\prime}$.

Set $y_n=p^c(x)_n$, and let us prove $x =\exp(y)$.
Since $\varepsilon(x) \neq 0$, we have $x_0=k \upsilon(1)$ with $k\neq 0$.
Since $\Delta(x)=x\otimes x$, we get $k=1$. So the equality $x =\exp(y)$ holds in degree $0$.
Assume that it holds until degree $(n-1)$. We have
$\Delta_n(x_n)=\sum_{i=0}^n x_i \otimes x_{n-i}$.
According to  Lemma~\ref{lemexpprimgrouplike}, we also have
$\Delta_n(\exp(y)_n)=\sum_{i=0}^n \exp(y)_i \otimes \exp(y)_{n-i}$.
Our induction hypothesis ensures \begin{equation*}\Delta_n(x_n - \exp(y)_n)=1 \otimes (x_n - \exp(y)_n) + (x_n - \exp(y)_n) \otimes 1.\end{equation*}
So $(x_n - \exp(y)_n)$ is primitive, and we get $(x_n - \exp(y)_n)=p^c(x_n - \exp(y)_n)=0$.
\eop

\chapter{First definitions of \texorpdfstring{$\Zinv$}{Z}}
\label{chapdefzinv}

In this chapter, we introduce the invariant $\Zinvuf$ of links in $\QQ$-spheres, which is the main
object of this book. We illustrate the required definitions with many examples, skippable by a reader who only wants the definition.

\section{Configuration spaces of Jacobi diagrams in \texorpdfstring{$3$}{3}-manifolds}
\label{secdefconfspace}

Let $(\crats,\tau)$ be an asymptotic rational homology $\RR^3$, as in Definition~\ref{defasyrathommRthree}.
Let $\source$ be a disjoint union of $k$ circles $S^1_i$, $i \in \underline{k}$. Let 
\begin{equation*}\Link \colon \source \longrightarrow \crats\end{equation*} denote a $C^{\infty}$ embedding from $\source$ to $\crats$. The embedding $\Link$ is a \emph{link embedding}. Let $\Gamma$ be a Jacobi diagram with support $\source$ as in Definition~\ref{defdia}. Let $U=U(\Gamma)$ denote the set of univalent vertices of $\Gamma$, and let $T=T(\Gamma)$ denote the set of trivalent vertices of $\Gamma$. A \indexT{configuration} of $\Gamma$ (with respect to $\Link$) is an embedding 
\begin{equation*}c\colon U \cup T \hookrightarrow \crats\end{equation*}
whose restriction $c\vert_{U}$ to $U$ may be written as $\Link \circ j$ for some $\Gamma$-compatible  injection 
\begin{equation*}j\colon U \hookrightarrow \source.\end{equation*} Denote the set of these configurations by
$\check{C}(\rats,\Link;\Gamma)$, \index[N]{Configuration spaces!CcheckLGamma@$\check{C}(\rats,\Link;\Gamma)$ open} we have
\begin{equation*}\check{C}(\rats,\Link;\Gamma)=\Bigl\{c \colon U \cup T \hookrightarrow \crats \suchthat
\bigl(\exists j \in \left[i_{\Gamma}\right] \suchthat c\vert_{U}=\Link \circ j\bigr)\Bigr\}.\end{equation*}
In $\check{C}(\rats,\Link;\Gamma)$, the univalent vertices move along $\Link(\source)$, while the trivalent vertices move in the ambient space $\crats$.
The configuration space $\check{C}(\rats,\Link;\Gamma)$ is naturally an open submanifold of 
$\source^U \times \crats^T$.

An \emph{orientation}\index[T]{orientation!of a set} of a set of cardinality at least $2$ is a total order of its elements up to an even permutation.
When $\Link$ is oriented, such an orientation of the set $V(\Gamma)$ of vertices of $\Gamma$ orients $\check{C}(\rats,\Link;\Gamma)$ naturally since it orders the oriented odd-dimensional factors of $\source^U \times \crats^T$. Below, we associate an orientation of 
$\check{C}(\rats,\Link;\Gamma)$ to a vertex-orientation of $\Gamma$ and an orientation of the set $H(\Gamma)$ of half-edges of $\Gamma$.

Cut each edge of $\Gamma$ into two half-edges. When an edge is oriented, define its \emph{first} half-edge and its \emph{second} one so that we meet  the first half-edge first when following the orientation of the edge.
 When the edges of $\Gamma$ are oriented, the orientations of the edges of $\Gamma$ induce the following orientation of the set $H(\Gamma)$. Order the set $E(\Gamma)$ of edges of $\Gamma$ arbitrarily, and order $H(\Gamma)$ as 
(First half-edge of the first edge, second half-edge of the first edge, \dots, second half-edge of the last edge). The induced orientation of $H(\Gamma)$ is called the \indexT{edge-orientation} of $H(\Gamma)$. Note that it does not depend on the order of $E(\Gamma)$.

\begin{lemma}
\label{lemorc}
When $\Gamma$ is equipped with a vertex-orientation, orientations of the manifold $\check{C}(\rats,\Link;\Gamma)$ are in canonical
one-to-one correspondence with orientations of the set $H(\Gamma)$.
\end{lemma}
\bp
Since $\check{C}(\rats,\Link;\Gamma)$ is naturally an open submanifold of 
$\source^U \times \crats^T$, it inherits $\RR^{\cardlef{U} + 3\cardlef{T}}$-valued charts from $\RR$-valued orientation-preserving charts of $\source$,
with respect to the (possibly local as in Definition~\ref{defdia}) orientation(s) of $\source$, and $\RR^3$-valued orientation-preserving charts of $\crats$. To define the orientation of $\RR^{\cardlef{U} + 3\cardlef{T}}$, it suffices to identify its factors and order them (up to even permutation).
Each factor is labeled by an element of $H(\Gamma)$ as follows. The $\RR$-valued local coordinate of an element of $\source$ corresponding to the image under $j$ of an element $u$ of $U$ sits in the factor labeled by the half-edge that contains $u$.
The three $\RR$-valued coordinates of the image under a configuration $c$ of an element $t$ of $T$, with respect to an arbitrary oriented local chart, belong to the factors labeled by the three half-edges containing $t$ so that
the cyclic order of the three half-edges induced by the vertex-orientation of $\Gamma$ matches the order of the three factors.
\eop

We will use Lemma~\ref{lemorc} to orient $\check{C}(\rats,\Link;\Gamma)$ as summarized in the following immediate corollary.

\begin{corollary} \label{cororc}
If $\Gamma$ is equipped with a vertex-orientation $o(\Gamma)$ and if the edges of\, $\Gamma$ are oriented, then the induced edge-orientation of $H(\Gamma)$ orients $\check{C}(\rats,\Link;\Gamma)$, via the canonical correspondence described in the proof of Lemma~\ref{lemorc}.
\end{corollary}

\begin{examples}
\label{exaorconfknot}
In Subsection~\ref{subsecwtwo}, the orientations of the configuration spaces $\check{C}(K;\diagcross)$ and $\check{C}(K;\diagtripod)$ were induced by the order of the given coordinates. We can check that these orientations are also induced by the edge-orientations and by the orientation of the vertex $w$ in the following figures
\begin{center}
 \begin{tikzpicture}[scale=.8] \useasboundingbox (-.7,-.5) rectangle (.7,.5);
\begin{scope}[xshift=-2cm]
\draw [->,dash pattern=on 2pt off 2pt] (.6,0) arc (0:360:.6);
\draw (45:.6) -- (-135:.6) node[left]{\scriptsize $b_2$};
\draw (135:.6) -- (-45:.6) node[right]{\scriptsize $b_3$};
\draw (135:.6) node[left]{\scriptsize $a_3$}  (45:.6)  node[right]{\scriptsize $a_2$};
\draw [->] (45:.6) -- (-135:.25);
\draw [->] (135:.6) -- (-45:.25);
\fill (45:.6) circle (2pt) (135:.6) circle (2pt) (-135:.6) circle (2pt) (-45:.6) circle (2pt);
\end{scope}
\draw (0,0) node{and};
\begin{scope}[xshift=2cm]
\draw [->,dash pattern=on 2pt off 2pt] (.6,0) arc (0:360:.6);
\draw (0,0) -- (-60:.6) node[right]{\scriptsize $v_1$};
\draw (0,0) -- (60:.6) node[right]{\scriptsize $v_2$};
\draw (.25,0) node{\scriptsize $w$}  (0,0) -- (180:.6) node[left]{\scriptsize $v_3$};
\draw [->] (-60:.6) -- (-60:.25);
\draw [->] (60:.6) -- (60:.25);
\draw [->] (180:.6) -- (180:.25);
\fill (60:.6) circle (2pt) (180:.6) circle (2pt) (-60:.6) circle (2pt) (0,0) circle (2pt);
\end{scope}
\draw (3,0) node{.};
\end{tikzpicture}
\end{center}
Recall that $\check{C}(K;\diagtripod)$, which may be described as
\begin{equation*}\left\{\confc \colon \{w,v_1,v_2,v_3\} \hookrightarrow \RR^3 \suchthat \begin{array}{l} \confc(v_i)=K(z_i \in S^1),
 z_2=\exp(2i\pi t_2)z_1,\\ z_3=\exp(2i\pi t_3)z_1, 0<t_2<t_3<1 \end{array}\right\},\end{equation*}
was regarded as an open submanifold of 
$\RR^3 \times (S^1)^3=\{(X_1,X_2,X_3,z_1,z_2,z_3)\}$,
where $(X_1,X_2,X_3)=\confc(w)$.

Let us check that the above coordinates orient $\check{C}(K;\diagtripod)$ as the orientation of edges does. For $i \in \underline{3}$, let $e(i)$ denote the edge from $v_i$ to $w$. This orders the three factors of $(S^2)^{E(\Gamma)}$. Distribute the coordinates $X_1$, $X_2$, $X_3$ so that $X_i$ is on the second half-edge of $e(i)$. Then the vertex-orientation and the edge-orientation of $\Gamma$ induce the orientation of $\check{C}(K;\diagtripod)$ represented by $(z_1,X_1,z_2,X_2,z_3,X_3)$, which is the same as the orientation represented by $(X_1,X_2,X_3,z_1,z_2,z_3)$. The case of $\check{C}(K;\diagcross)$ is left to the reader.
\end{examples}

\begin{example}
\label{exaortheta}
Equip the diagram $\tata$ with its vertex-orientation induced by the picture.
Orient its three edges so that they start from the same vertex. Then the orientation of $\check{C}(\rats,\Link;\tata)$ induced by this edge-orientation of $\tata$ matches the orientation of $(\crats \times \crats) \setminus \diag$ induced by the order of the two factors, where the first factor corresponds to the position of the vertex from which the three edges start, as shown in the following picture. 
\begin{equation*}\tatanumone \cong \tatanumtwo \end{equation*} 
\end{example}

\begin{remark}
\label{rkedgeorvert}
For a Jacobi diagram $\Gamma$ equipped with a vertex-orientation $o(\Gamma)$, an orientation of $V(\Gamma)$ induces the following orientation of $H(\Gamma)$. Fix a total order of $V(\Gamma)$ that induces its given orientation.
Then the corresponding orientation of $H(\Gamma)$ is induced by a total order which starts with the half-edges adjacent to the first vertex, ordered with respect to $o(\Gamma)$ if the vertex is trivalent, and continues with the half-edges adjacent to the second vertex, to the third one, \dots
This orientation is called the \emph{vertex-orientation}\index[T]{vertex-orientation!of $H(\Gamma)$} of $H(\Gamma)$ associated to $o(\Gamma)$ and to the orientation of $V(\Gamma)$.
In particular, an orientation of $H(\Gamma)$ (such as the edge-orientation of $H(\Gamma)$ when the edges of $\Gamma$ are oriented) and a vertex-orientation $o(\Gamma)$ together induce an orientation of $V(\Gamma)$, namely the orientation of $V(\Gamma)$ such that the induced vertex-orientation of $H(\Gamma)$ matches the given orientation of $H(\Gamma)$.
\end{remark}

The dimension of $\check{C}(\rats,\Link;\Gamma)$ is
\begin{equation*}\cardlef{U(\Gamma)} + 3\cardlef{T(\Gamma)}= 2 \cardlef{E(\Gamma)}.\end{equation*}
Since the degree of $\Gamma$ is $n=n(\Gamma)=\frac12(\cardlef{U(\Gamma)} +\cardlef{T(\Gamma)})$, we have \begin{equation*}\cardlef{E(\Gamma)}=3n-\cardlef{U(\Gamma)}.\end{equation*}

\section{Configuration space integrals}
\label{secconfint}

\begin{definition}
\label{defnumdia}
Let $\finseta$ be a finite set.
An \emph{$\finseta$-numbered Jacobi diagram}\index[T]{Anumbered@$\finseta$-numbered Jacobi diagram} is a Jacobi diagram $\Gamma$ whose edges are oriented, equipped with an injection $j_E \colon E(\Gamma) \hookrightarrow \finseta$. \index[N]{jE@$j_E \colon E(\Gamma) \hookrightarrow \finseta$}
Such an injection numbers the edges when $\finseta \subset \NN$.
Let $\Davis^e_n(\source)$\index[N]{Diag@Diagram sets!Denn@$\Davis^e_n$ $\underline{3n}$-numbered} denote the set of $\underline{3n}$-numbered degree $n$ Jacobi diagrams with support $\source$ without \emph{looped edges} like \loopedge.
\end{definition}
Note that the injection $j_E$ is a bijection for any diagram of $\Davis^e_n(\source)$ without univalent vertices.

\begin{examples}
\label{exanumdia} We have
\begin{multline*}
 \begin{array}{l}\Davis^e_1(\emptyset)=\left\{\tatanumor{>}{>},\tatanumor{>}{<},\tatanumor{<}{>},\tatanumor{<}{<}\right\},\\
 \Davis^e_1(S^1)=\Davis^e_1(\emptyset) \sqcup \left\{ \onechordnum{$S^1$}{$1$},\onechordnum{$S^1$}{$2$},\onechordnum{$S^1$}{$3$} \right\}, \mbox{ and} \\
\Davis^e_1(S^1_1 \sqcup S^1_2)=\Davis^e_1(\emptyset) \sqcup \left(\Davis^e_1(S^1_1) \setminus \Davis^e_1(\emptyset)\right)
 \sqcup \left(\Davis^e_1(S^1_2) \setminus \Davis^e_1(\emptyset)\right) \end{array}\\
\sqcup \Bigl\{\onechordtwocircles{1}{>},\onechordtwocircles{2}{>},\onechordtwocircles{3}{>},\\\onechordtwocircles{1}{<},\onechordtwocircles{2}{<},\onechordtwocircles{3}{<}\Bigr\}.
\end{multline*}
\end{examples}

Let $\Gamma$ be a numbered degree $n$ Jacobi diagram with support $\source$.
An edge $e$ oriented from a vertex $v_1$ to a vertex $v_2$ of $\Gamma$ induces the map \begin{equation*}\begin{array}{llll}p_e \colon &\check{C}(\rats,\Link;\Gamma) &\rightarrow &C_2(\rats)\\
           & c & \mapsto & (c(v_1),c(v_2)).\end{array}\end{equation*}

Let $o(\Gamma)$ be a  vertex-orientation of $\Gamma$.
For any $i \in \underline{3n}$, let $\omega(i)$ be a propagating form of $(C_2(\rats),\tau)$. Let $(\check{C}(\rats,\Link;\Gamma),o(\Gamma))$ denote the manifold $\check{C}(\rats,\Link;\Gamma)$, equipped with the orientation induced by $o(\Gamma)$ and by the edge-orientation of $\Gamma$, as in Corollary~\ref{cororc}.
Define \begin{equation*}I\left(\rats,\Link,\Gamma,o(\Gamma),(\omega(i))_{i \in \underline{3n}}\right)=\int_{(\check{C}(\rats,\Link;\Gamma),o(\Gamma))} \bigwedge_{e \in E(\Gamma)}p_e^{\ast}\Bigl(\omega\bigl(j_E(e)\bigr)\Bigr) \index[N]{Integrals over configuration spaces!IRLGammaa@$I(\rats,\Link,\Gamma,o(\Gamma),(\omega(i))_{i \in \underline{3n}})$}.\end{equation*} 
The convergence of this integral is a consequence of the following proposition, which will be proved in Chapter~\ref{chapcompconf}. (See the end of Section~\ref{secfirstprescomp}.)

\begin{proposition}
\label{propcompext}
 There exists a smooth compactification of $\check{C}(\rats,\Link;\Gamma)$, which will be denoted by $C(\rats,\Link;\Gamma)$, to which the maps $p_e$ extend smoothly.
\end{proposition}

According to this proposition, $\bigwedge_{e \in E(\Gamma)}p_e^{\ast}\left(\omega\bigl(j_E(e)\bigr)\right)$ extends smoothly to $C(\rats,\Link;\Gamma)$, and we have
\begin{equation*}\int_{(\check{C}(\rats,\Link;\Gamma),o(\Gamma))} \bigwedge_{e \in E(\Gamma)}p_e^{\ast}\Bigl(\omega\bigl(j_E(e)\bigr)\Bigr)=
\int_{(C(\rats,\Link;\Gamma),o(\Gamma))} \bigwedge_{e \in E(\Gamma)}p_e^{\ast}\Bigl(\omega\bigl(j_E(e)\bigr)\Bigr).\end{equation*}
\begin{examples}
For any three propagating forms $\omega(1)$, $\omega(2)$, and $\omega(3)$ of $(C_2(\rats),\tau)$, we have
 \begin{equation*}I\Bigl(\rats, K_i \sqcup K_j \colon S^1_i \sqcup S^1_j \hookrightarrow \crats,\haltereor,\bigl(\omega(i)\bigr)_{i \in \underline{3}}\Bigr)=lk(K_i,K_j)\end{equation*}
and \begin{equation*}I\Bigl(\rats,\emptyset,\tataor,\bigl(\omega(i)\bigr)_{i \in \underline{3}}\Bigr)=\Theta(\rats,\tau)\end{equation*}
for any numbering of the (nondashed) diagrams (exercise).
\end{examples}

\begin{examples}
\label{examplecomconfinttwo}
For any oriented trivalent numbered degree $n$ Jacobi diagram $\bigl(\Gamma,o(\Gamma)\bigr)$, we have \begin{equation*}I\bigl(\Gamma,o(\Gamma)\bigr)=I\left(S^3,\emptyset,\Gamma,o(\Gamma),(p_{S^2}^{\ast}(\omega_{S^2}))_{i \in \underline{3n}}\right)=0.\end{equation*}
Indeed, $I\bigl(\Gamma,o(\Gamma)\bigr)$ is equal to
\begin{equation*}\int_{(\check{C}(S^3,\emptyset;\Gamma),o(\Gamma))}\Bigl(\prod_{e \in E(\Gamma)}p_{S^2} \circ p_e\Bigr)^{\ast}\Bigl( \bigwedge_{e \in E(\Gamma)}\omega_{S^2}\Bigr),\end{equation*} where 
\begin{itemize}
\item $\bigwedge_{e \in E(\Gamma)}\omega_{S^2}$ is a product volume form of $\left(S^2\right)^{E(\Gamma)}$ with total volume one,
 \item $\check{C}(S^3,\emptyset;\Gamma)$ is the space 
of injections of $\underline{3n}$ into $\RR^3$,
\item the degree of $\wedge_{e \in E(\Gamma)}\omega_{S^2}$ is equal to the dimension of $\check{C}(S^3,\emptyset;\Gamma)$, and
\item the map $\prod_{e \in E(\Gamma)}p_{S^2} \circ p_e$ is never a local diffeomorphism since it is invariant under the action of global translations on $\check{C}(S^3,\emptyset;\Gamma)$.
\end{itemize}
\end{examples}

\begin{examples}
\label{examplecomconfintthree}
Let $O$ denote the representative of the unknot of $S^3$, that is the image of the embedding of the unit circle $S^1$ of $\CC$ regarded as $\CC\times\{0\}$, into $\RR^3$ regarded as $\CC\times\RR$.
Let us compute 
\begin{equation*}I\left(S^3,O,\Gamma,o(\Gamma),(p_{S^2}^{\ast}(\omega_{S^2}))_{i \in \underline{3n}}\right),\end{equation*}
 for the oriented Jacobi diagrams
\begin{equation*}\Gamma_1=\twoischordor\mbox{, }\Gamma_2=\twoxchordor\mbox{, }\Gamma_3=\twobullechordor\mbox{, and }\Gamma_4=\tripodor.\end{equation*}
Since all edges are equipped with the same standard propagating form $p_{S^2}^{\ast}\left(\omega_{S^2}\right)$, we do not number the edges. For $i\in\underline{4}$, set
\begin{equation*}I(\Gamma_i)=I\left(S^3,O,\Gamma_i,o(\Gamma_i),(p_{S^2}^{\ast}(\omega_{S^2}))\right).\end{equation*}
We are about to prove
$I(\Gamma_1)=I(\Gamma_2)=I(\Gamma_3)=0$ and $I(\Gamma_4)=\frac18$.
For $i\in\underline{4}$, set $\Gamma=\Gamma_i$. Recall that $I(\Gamma)$ is equal to
\begin{equation*}\int_{(\check{C}(S^3,O;\Gamma),o(\Gamma))}\biggl(\prod_{e \in E(\Gamma)}p_{S^2} \circ p_e\biggr)^{\ast}\biggl( \bigwedge_{e \in E(\Gamma)}\omega_{S^2}\biggr).\end{equation*}
When $i \in\underline{2}$, the image of $\prod_{e \in E(\Gamma)}p_{S^2} \circ p_e$ lies in the subset of $(S^2)^2$ consisting of the pair of horizontal vectors. Since the interior of this subset is empty, the form $\bigl(\prod_{e \in E(\Gamma)}p_{S^2} \circ p_e\bigr)^{\ast}\bigl(\omega_{S^2}^{E(\Gamma)}\bigr)$ vanishes identically. So we get $I(\Gamma_i)=0$.
When $i=3$, the two edges with the same endpoints must have the same direction. So the image of $\prod_{e \in E(\Gamma)}p_{S^2} \circ p_e$ lies in the subset of $(S^2)^{E(\Gamma)}$, where two $S^2$-coordinates are identical (namely those in the $S^2$-factors corresponding to the above pair of edges), and we obtain $I(\Gamma_3)=0$ as before.
\end{examples}
\begin{lemma} \label{lemItripod} Let $\Gamma=\tripodor$. Then we have
 \begin{equation*}I\left(S^3,O,\Gamma,o(\Gamma),(p_{S^2}^{\ast}(\omega_{S^2}))\right)=\frac18.\end{equation*}
\end{lemma}
\bp 
Let $(S^1)^3_+$ be the subset of $(S^1)^3$ consisting of triples $(z_1,z_2,z_3)$ of pairwise distinct elements of $S^1$ such that the orientation of $S^1$ induces the cyclic order $(z_1,z_2,z_3)$.
Then we write \begin{equation*}\check{C}=\check{C}(S^3,O;\Gamma)=\Bigl\{(X,z_1,z_2,z_3) \suchthat   X \in \RR^3 \setminus \{z_1,z_2,z_3\},(z_1,z_2,z_3)\in (S^1)^3_+\Bigr\},\end{equation*}
where $X=(X_1,X_2,X_3)$.
As explained in Example~\ref{exaorconfknot}, these coordinates orient $\check{C}$ as does the orientation of edges.
Let $\check{C}^+=\{(X,z_1,z_2,z_3) \in \check{C} \suchthat  X_3>0\}$.
The reflection $\sigma_h$ with respect to the horizontal plane acts on $\check{C}$ by an orientation-reversing diffeomorphism, which changes $X_3$ to $(-X_3)$ and leaves the other coordinates unchanged. It also acts on $S^2$ by an orientation-reversing diffeomorphism, which preserves the volume up to sign.
Therefore, we have \begin{equation*}I(\Gamma)=2\int_{\check{C}^+}\Bigl(\prod_{e \in E(\Gamma)}p_{S^2} \circ p_e\Bigr)^{\ast}\Bigl( \bigwedge_{e \in E(\Gamma)}\omega_{S^2}\Bigr).\end{equation*}
Let $S^2_+$ denote the set of elements of $S^2$ with positive \emph{height} (third coordinate), and let $(S^2)^3_+$ be the set of elements of $\left(S^2_+\right)^3$ which form a direct basis.
Let us check that the volume of $(S^2)^3_+$ is $\frac1{16}$ (with respect to $\bigwedge_{e \in E(\Gamma)}\omega_{S^2}$). The volume of $\left(S^2_+\right)^3$ is $\frac18$. This is also the volume of the subset of $\left(S^2_+\right)^3$ consisting of triples of noncoplanar vectors. The involution that exchanges the last two vectors in the latter set sends the direct bases to the indirect ones, and it preserves the volume.
Therefore, in order to prove $I(\Gamma_4)=\frac18$, it suffices to prove that 
\begin{equation*}\begin{array}{llll}\Psi \colon & \check{C}^+ & \to &(S^2)^3 \\
   & c & \mapsto & \left(\prod_{e \in E(\Gamma)}p_{S^2} \circ p_e\right)(c)
  \end{array}\end{equation*}
is an orientation-preserving diffeomorphism onto $(S^2)^3_+$.
Let $c=(X,z_1,z_2,z_3)$ be a point of $\check{C}^+$. Let us first check that $\Psi(c)\in (S^2)^3_+$. For $j \in \underline{3}$, the vector $\bigvec{z_jX}$ may be written as $\lambda_jV_j$ for some $\lambda_j \in \left]0,+\infty\right[$ and for some $V_j \in S^2_+$. Since $(V_1,\bigvec{z_1z_2}=\lambda_1V_1-\lambda_2V_2, \bigvec{z_1z_3}=\lambda_1V_1-\lambda_3V_3)$ is a direct basis of $\RR^3$, so 
is $(V_1,V_2,V_3)$.
Let us now compute the sign of the Jacobian of $\Psi$ at $c$. Let $T_c\Psi$ denote the tangent map to $\Psi$ at $c$.
For $j \in \underline{3}$, let $Z_j$ denote the unit tangent vector to $S^1$ at $z_j$, and let 
$p_j\colon (\RR^3)^3 \to \RR^3/\RR V_j$ be the projection onto the $j^{th}$ factor composed by the projection onto the tangent space $\RR^3/\RR V_j$ to $S^2$ at $V_j$.
 Then
$p_j(T_c\Psi(Z_j))=-Z_j$ in the tangent space to $S^2$ at $V_j$, which is generated by (the projections onto $\RR^3/\RR V_j$ of) $-Z_j$ and any vector $W_j$ such that $\det(V_j,-Z_j,W_j) =1$.
Reorder the oriented basis $(-Z_1,W_1,-Z_2,W_2,-Z_3,W_3)$ of $T(S^2)^3_+$ at $(V_1,V_2,V_3)$ by a positive permutation of the coordinates as $(W_1,W_2,W_3,-Z_1,-Z_2,-Z_3)$. 
Writing the matrix of $T_c\Psi$ with respect to this basis $(W_1,W_2,W_3,-Z_1,-Z_2,-Z_3)$ of the target space and the basis $(V_1,V_2,V_3)$ of $T_X\RR^3$ followed by $(Z_1,Z_2,Z_3)$ for the domain $T_c\check{C}^+$ produces a matrix whose last three columns contain $1$ on the diagonals as only nonzero entries. In the quotient \begin{equation*}{T(S^2)^3_+}/{\bigl(\RR (-Z_1,0,0) \oplus \RR (0,-Z_2,0)\oplus \RR (0,0,-Z_3)\bigr)},\end{equation*} $T_c\Psi(V_1)$ may be expressed as
\begin{equation*}T_c\Psi(V_1) = \det(V_2,-Z_2,V_1)W_2 + \det(V_3,-Z_3,V_1)W_3,\end{equation*}
where $\det(V_2,-Z_2,V_1)=\det(Z_2,V_2,V_1)$. Let us prove that $\det(Z_2,V_2,V_1)>0$. When $z_2=-1$, projecting $c$ on $\left(\RR=\frac{\CC}{i\RR}\right) \times \RR$ produces a picture as in Figure~\ref{figprojcomptriv}
\bfig
\centering
\begin{tikzpicture} \useasboundingbox (-1.5,-.7) rectangle (1.5,1.2);
\draw [-] (-1,0) -- (1,0) (-1,0) -- (-1.4,1.2) -- (.4,0) (-.3,.6) node[right]{\scriptsize $V_1$} (-1.2,.6) node[left]{\scriptsize $V_2$} (-1,0) node[left]{\scriptsize $z_2$} (.4,0)  node[above]{\scriptsize $z_1$} (-1.4,1.2) node[left]{\scriptsize $X$} (1,0) node[right]{\scriptsize $1$} (0,0)  node[below]{\scriptsize $S^1$};
\draw [->] (-1,0) -- (-1.2,.6);
\draw [->] (.4,0) -- (-.5,.6);
\fill (-1,0) circle (1.5 pt) (1,0) circle (1.5pt) (-1.4,1.2) circle (1pt) (.4,0) circle (1.5pt);
\end{tikzpicture}
\caption{Partial projection of $(X,z_1,z_2,z_3)$ on $\frac{\CC}{i\RR} \times \RR$ when $z_2=-1$}
\label{figprojcomptriv}
\end{figure}
 that makes this result clear. The general result follows easily.
Finally, the Jacobian of $\Psi$ at $(X,z_1,z_2,z_3)$ is the determinant of 
\begin{equation*}\left[\begin{array}{ccc} 0 &\det(Z_1,V_1,V_2) &\det(Z_1,V_1,V_3)\\
 \det(Z_2,V_2,V_1)& 0 & \det(Z_2,V_2,V_3)\\ 
\det(Z_3,V_3,V_1)& \det(Z_3,V_3,V_2)& 0
\end{array}\right],\end{equation*}
which is \begin{equation*}\begin{array}{ll}&\det(Z_1,V_1,V_2)\det(Z_2,V_2,V_3)\det(Z_3,V_3,V_1)\\+&\det(Z_1,V_1,V_3) \det(Z_2,V_2,V_1)\det(Z_3,V_3,V_2).\end{array}\end{equation*} It is positive since all the involved terms are.
Let us finally check that every element $(V_1,V_2,V_3)$ of $(S^2)^3_+$ has a unique element in its preimage. Construct the three lines of $\RR^3$ directed by $V_1$, $V_2$, and $V_3$ through the origin of $\RR^3$. The line directed by $V_i$ intersects the horizontal plane at height $(-1)$ at a point $w_i$.
There is a unique circle in this horizontal plane that contains $w_1$, $w_2$, and $w_3$. Let $\frac1{\lambda}$ be the radius of this circle, and let $w_0$ be its center. Then the unique element of $\Psi^{-1}(V_1,V_2,V_3)$ is $(-\lambda w_0, \lambda w_1 - \lambda w_0,\lambda w_2 - \lambda w_0,\lambda w_3 - \lambda w_0)$.
\eop

\section{Configuration space integrals associated to a chord}
\label{secthetaknot}

Let us now study the case of $I(\onechordsmalljnum,(\omega(i))_{i \in \underline{3}})$. We will see that this integral depends on the chosen propagating forms and on the diagram numbering.

A \emph{dilation}\index[T]{dilation} is a homothety with a positive ratio.
Let $\ST^+K_j$\index[N]{Uplus@$\ST^+K$ unit positive tangent vectors to $K$} denote the fiber space over $K_j$ consisting of the tangent vectors to the knot $K_j$ of $\crats$ that orient $K_j$ up to dilation. The fiber of $\ST^+K_j$ consists of one point. So the total space of this {\em unit positive tangent bundle to $K_j$\/} is canonically diffeomorphic to $K_j$. Set $\ST^-K_j = \ST^+(-K_j)$\index[N]{Uminus@$\ST^-K$ unit negative tangent vectors to $K$}. 

For a knot $K_j$ in $\crats$, we have \begin{equation*}\check{C}(K_j;\onechordsmallj)=\Bigl\{\bigl(K_j(z),K_j(z\exp(i\theta))\bigr) \suchthat  (z,\theta) \in S^1 \times \left]0,2\pi\right[\Bigr\} .\end{equation*}
Let $C_j={C}(K_j;\onechordsmallj)$ be the closure of $\check{C}(K_j;\onechordsmallj)$ in $C_2(\rats)$. This closure is
diffeomorphic to $S^1 \times \left[0,2\pi\right]$, where $S^1 \times \{0\}$ is identified with $\ST^+K_j$, $S^1 \times \{2\pi\}$ is identified with $\ST^-K_j$, and $\partial C(K_j;\onechordsmallj) =\ST^+K_j - \ST^-K_j$.

\begin{lemma}
\label{lemvaritheta}
For any $i\in \underline{3}$, let $\omega(i)$ and $\omega^{\prime}(i)$ 
be propagating forms of $(C_2(\rats),\tau)$, which restrict to $\partial C_2(\rats)$ as $\projp_{\tau}^{\ast}(\omega(i)_{S^2})$ and $\projp_{\tau}^{\ast}(\omega^{\prime}(i)_{S^2})$, respectively.
Then there exists a one-form $\eta(i)_{S^2}$ on $S^2$ such that $\omega^{\prime}(i)_{S^2}=\omega(i)_{S^2} + d\eta(i)_{S^2}$, and we have
\begin{equation*}\begin{array}{ll}I\Bigl(\onechordsmalljnum,\bigl(\omega^{\prime}(i)\bigr)_{i \in \underline{3}}\Bigr)-I\Bigl(\onechordsmalljnum,\bigl(\omega(i)\bigr)_{i \in \underline{3}}\Bigr)=& \int_{\ST^+K_j}\projp_{\partau}^{\ast}\bigl(\eta(k)_{S^2}\bigr) \\
& -\int_{\ST^-K_j}\projp_{\partau}^{\ast}\bigl(\eta(k)_{S^2}\bigr).\end{array} \end{equation*}

\end{lemma}
\bp Such a form $\eta(i)_{S^2}$ exists for any $i$. According to Lemma~\ref{lemetactwo}, there exists a one-form $\eta(i)$ on $C_2(\rats)$ such that $\omega^{\prime}(i)=\omega(i)+d\eta(i)$ and the restriction of $\eta(i)$ to $\partial C_2(\rats)$ is $\projp_{\partau}^{\ast}(\eta(i)_{S^2})$. Apply Stokes' theorem to $\int_{C_j}\bigl(\omega^{\prime}(k)-\omega(k)\bigr)=\int_{C_j}d\eta(k)$. \eop

\begin{exo}
 Find a knot $K_j$ of $\RR^3$ and a form $\eta(k)$ of $C_2(\RR^3)$ such that the right-hand side of Lemma~\ref{lemvaritheta} does not vanish. (Use Lemma~\ref{lemetactwo}. Hints can be found in Section~\ref{secstraight}.)
\end{exo}

Recall that a propagating form $\omega$ of $\bigl(C_2(\rats),\partau\bigr)$ is \emph{homogeneous} if its restriction to $\partial C_2(\rats)$ is $\projp_{\partau}^{\ast}(\omega_{S^2})$ for the homogeneous volume form $\omega_{S^2}$ of $S^2$ of total volume $1$.

\begin{lemma}
\label{lemdefItheta}
For any $i\in \underline{3}$, let $\omega(i)$ be a homogeneous propagating form of $\bigl(C_2(\rats),\partau\bigr)$.
Then \begin{equation*}I\Bigl(\onechordsmalljnum,\bigl(\omega(i)\bigr)_{i \in \underline{3}}\Bigr)\end{equation*} does not depend on the choices of the $\omega(i)$, it is denoted by $I_{\theta}(K_j,\partau)$.\index[N]{Integrals over configuration spaces!Ithetaknot@$I_{\theta}(K,\partau)$ associated to \tatak \;}
\end{lemma}
\bp It is a corollary of Lemma~\ref{lemvaritheta}.
\eop

\section{First definition of \texorpdfstring{$\Zinv$}{Z}}
\label{secfirstdef}

From now on, the coefficients of our spaces of Jacobi diagrams are in $\RR$. ($\KK=\RR$.)
Let $(\crats,\partau)$ be an asymptotic rational homology $\RR^3$. 
Let $\Link \colon \source \hookrightarrow \crats$ be a link embedding.

Let $\Gamma$ be a numbered Jacobi diagram $\Gamma$ in the space $\Davis^e_n(\source)$ of Definition~\ref{defnumdia}.
Let
$\left[\Gamma,o(\Gamma)\right]$ denote the class in $\Aavis_n(\source)$ of $\Gamma$ equipped with a vertex-orientation $o(\Gamma)$.
Then \begin{equation*}I\Bigl(\rats,\Link,\Gamma,o(\Gamma),\bigl(\omega(i)\bigr)_{i \in \underline{3n}}\Bigr)\left[\Gamma,o(\Gamma)\right] \in \Aavis_n(\source)\end{equation*} is independent of the vertex-orientation $o(\Gamma)$ of $\Gamma$. It is simply denoted by 
\begin{equation*}I\Bigl(\rats,\Link,\Gamma,\bigl(\omega(i)\bigr)_{i \in \underline{3n}}\Bigr)\left[\Gamma\right].\end{equation*}

\begin{notation}
\label{notationzZ}
For $\Gamma \in \Davis^e_n(\source)$, set\index[N]{zwetaGamma@$\coefgambet_{\Gamma}$ averaging coefficient} \begin{equation*}\coefgambet_{\Gamma}=\frac{(3n-\cardlef{E(\Gamma)})!}{(3n)!2^{\cardlef{E(\Gamma)}}}.\end{equation*}
Recall the definitions of propagating forms from Section~\ref{secprop}.
For any $i\in \underline{3n}$, let $\omega(i)$ be a propagating form of $C_2(\rats)$.
For $n \in \NN$, set \begin{equation*}\Zinv_n\Bigl(\crats,\Link,\bigl(\omega(i)\bigr)_{i \in \underline{3n}}\Bigr)=\sum_{\Gamma \in \Davis^e_n(\source)}\coefgambet_{\Gamma}I\Bigl(\rats,\Link,\Gamma,\bigl(\omega(i)\bigr)_{i \in \underline{3n}}\Bigr)\left[\Gamma\right] \in \Aavis_n(\source).\end{equation*} \index[N]{ZZ@$\Zinvuf$ and some variants
(see also the summary in the next pages)!ZiomR@$\Zinv_n(\crats,\Link,(\omega(i)))$}

This $\Zinv_n\bigl(\crats,\Link,(\omega(i))_{i \in \underline{3n}}\bigr)$ is the hero of this book. Let us describe some of its variants.

Let $\Aavis^c_n(\emptyset)$ \index[N]{Ap@$\Aavis^c_n(\emptyset)$ connected} denote the subspace of $\Aavis_n(\emptyset)$ generated by connected trivalent Jacobi diagrams. Set $\Aavis^c(\emptyset)=\prod_{n \in \NN}\Aavis^c_n(\emptyset)$, and let $p^c \colon \Aavis(\emptyset) \to \Aavis^c(\emptyset)$ be the projection that maps the empty diagram and diagrams with several connected components to $0$. \index[N]{Projections!pc@$p^c \colon \Aavis(\emptyset) \to \Aavis^c(\emptyset)$}
Let $\Davis^{c}_n$ denote the subset of $\Davis^{e}_n(\emptyset)$ that contains the connected diagrams of $\Davis^{e}_n(\emptyset)$. \index[N]{Diag@Diagram sets!Dcn@$\Davis^{c}_n$ connected}
For $n \in \NN$, set \begin{equation*}\zinv_n\Bigl(\crats,\bigl(\omega(i)\bigr)_{i \in \underline{3n}}\Bigr)=p^c\biggl(\Zinv_{n}\Bigl(\crats,\emptyset,\bigl(\omega(i)\bigr)\Bigr)\biggr).\end{equation*} \index[N]{ZZ@$\Zinvuf$ and some variants (see also the summary in the next pages)!zaRomega@$\zinv(\crats,(\omega(i)))$} \begin{equation*}\zinv_n\Bigl(\crats,\bigl(\omega(i)\bigr)\Bigr)=\sum_{\Gamma \in \Davis^{c}_n}\coefgambet_{\Gamma}I\Bigl(\rats,\emptyset,\Gamma,\bigl(\omega(i)\bigr)_{i \in \underline{3n}}\Bigr)\left[\Gamma\right] \in \Aavis^c_n(\emptyset).\end{equation*} 

When all the forms $\omega(i)$ are equal to $\omega(1)$, $\Zinv_n(\crats,\Link,\omega(1))$ and $\zinv_n(\crats,\omega(1))$ respectively denote
$\Zinv_n\bigl(\crats,\Link,(\omega(i))_{i \in \underline{3n}}\bigr)$ and $\zinv_n\bigl(\crats,(\omega(i))_{i \in \underline{3n}}\bigr)$.

We also use the projection $\projassis \colon \Aavis(\source) \to \Assis(\source)$\index[N]{Projections!pcheck@$\projassis\colon \Aavis(\source) \to \Assis(\source)$} that maps the diagrams with connected components without univalent vertices to zero and that maps the other diagrams to themselves.
Set $\Zinvlink_n = \projassis \circ \Zinv_n$. For example, we have
\begin{equation*}\Zinvlink_n(\crats,\Link,\omega(1))=\projassis\left(\Zinv_n(\crats,\Link,\omega(1))\right).\end{equation*}
We also remove the subscript $n$ to denote the collection (or the sum) of the $\Zinv_n$ for $n \in \NN$.
For example, we write \begin{equation*}\Zinvlink\left(\crats,\Link,\omega(1)\right)=\Bigl(\Zinvlink_n\bigl(\crats,\Link,\omega(1)\bigr)\Bigr)_{n\in \NN}=
\sum_{n\in \NN}\Zinvlink_n\left(\crats,\Link,\omega(1)\right) \in \Assis(\source).\end{equation*}

\end{notation}
As a first example, let us prove the following proposition.
\begin{proposition}
\label{propthetazoneone}
Let $(\crats,\partau)$ be an asymptotic rational homology $\RR^3$, then for any triple $(\omega(i))_{i \in \underline{3}}$ of propagating forms of $\bigl(C_2(\rats),\partau\bigr)$, we have
\begin{equation*}\Zinv_1\Bigl(\crats,\emptyset,\bigl(\omega(i)\bigr)_{i \in \underline{3}}\Bigr)=\zinv_{1}\Bigl(\crats,\bigl(\omega(i)\bigr)\Bigr)=
\frac{\Theta(\rats,\partau)}{12}\left[ \tata \right]\end{equation*}
in $\Aavis_1(\emptyset)=\Aavis_1(\emptyset;\RR)=\RR\left[\tata\right]$.
\end{proposition}
\bp The diagram $\tata$ is the only trivalent diagram with $2$ vertices without looped edges, and it is easy to check that $\Aavis_1(\emptyset)=\RR\left[\tata\right]$.
Each of the three edges goes from one vertex to the other.
There are $4$ elements in $\Davis^e_1(\emptyset)$, depending on whether the orientations of Edge $2$ and Edge $3$ coincide with the orientation of Edge $1$. See Example~\ref{exanumdia}.
When the three edges start from the same vertex, call the corresponding element $\theta_{++}$, and order the vertices so that the vertex from which the edges start is first. Recall from Example~\ref{exaortheta} that the vertex-orientation $o(\theta_{++})$ of $\tata$ induced by the picture $\tata$ and the edge-orientation of $H(\theta_{++})$ orient $C_2(\rats)$ as $\crats^2 \setminus \diag(\crats^2)$, as in Corollary~\ref{cororc}.
Then, according to Theorem~\ref{thmdefTheta},
 for any triple $(\omega(i))_{i \in \underline{3}}$ of propagating forms  of $\bigl(C_2(\rats),\partau\bigr)$, we have
\begin{equation*}I\Bigl(\theta_{++},o(\theta_{++}),\bigl(\omega(i)\bigr)_{i \in \underline{3}}\Bigr)=\Theta(\rats,\partau).\end{equation*}
Reversing the edge-orientation of Edge $2$ transforms $\theta_{++}$ to $\theta_{-+}$ and $\omega(2)$ to $\iota^{\ast}(\omega(2))$.
It changes the edge-orientation of $H(\theta_{\pm +})$. 
According to Lemma~\ref{lemiotaform}, the form $\bigl(-\iota^{\ast}(\omega(2))\bigr)$ is a propagating form of $\bigl(C_2(\rats),\partau\bigr)$.
Therefore, we have
\begin{equation*}I\Bigl(\theta_{++},\bigl(\omega(i)\bigr)_{i \in \underline{3}}\Bigr)\left[\theta_{++}\right]=I\Bigl(\theta_{-+},\bigl(\omega(i)\bigr)_{i \in \underline{3}}\Bigr)\left[\theta_{-+}\right].\end{equation*}
Similarly, the four graphs of $\Davis^e_1(\emptyset)$ contribute in the same way 
to \begin{equation*}\Zinv_1\Bigl(\crats,\emptyset,\bigl(\omega(i)\bigr)_{i \in \underline{3}}\Bigr)=\frac{4}{3!2^3}\Theta(\rats,\partau)\left[ \tata \right].\end{equation*}
\eop

\begin{examples}
 \label{examplecomconfinttwocont}
According to the computations in Example~\ref{examplecomconfinttwo}, we have \begin{equation*}\Zinv_n\bigl(\RR^3,\emptyset,p_{S^2}^{\ast}(\omega_{S^2})\bigr)=\zinv_n\bigl(\RR^3,p_{S^2}^{\ast}(\omega_{S^2})\bigr)=0\end{equation*} for $n>0$, and $\Zinv_0(\RR^3,\emptyset,p_{S^2}^{\ast}(\omega_{S^2}))=\left[\emptyset\right]$, while $\zinv_0(\RR^3,p_{S^2}^{\ast}(\omega_{S^2}))=0$.

For the embedding $O$ of the trivial knot in $\RR^3$ of Example~\ref{examplecomconfintthree}, we have
\begin{equation*}\Zinv_0\bigl(\RR^3,O,p_{S^2}^{\ast}(\omega_{S^2})\bigr)=1=\left[\emptyset\right].\end{equation*}
Since $I_{\theta}(O,\taust)=0$, we have $\Zinv_1(\RR^3,O,p_{S^2}^{\ast}(\omega_{S^2}))=0$.

Let us now prove
\begin{equation*}\Zinv_2\left(\RR^3,O,p_{S^2}^{\ast}(\omega_{S^2})\right)
=\frac1{24}\left[\tripodnonor\right].\end{equation*}
Note $\iota_{S^2}^{\ast}(\omega_{S^2})=-\omega_{S^2}$, where $\iota_{S^2}$
is the antipodal map of $S^2$.
So, reversing the orientation of an edge does not change
$I(S^3,O,\Gamma,p_{S^2}^{\ast}(\omega_{S^2}))\left[\Gamma\right]$, for a degree $2$ numbered Jacobi diagram $\Gamma$, since it changes both the orientation of $\check{C}(S^3,0;\Gamma)$ and the sign of the form to be integrated. Thus, $I(S^3,O,\Gamma,p_{S^2}^{\ast}(\omega_{S^2}))\left[\Gamma\right]$ depends only on the underlying Jacobi diagram. The degree $2$ Jacobi diagrams all components of which have univalent vertices and without looped edges are
\begin{equation*}
\twoischord, \twoxchord, \twobullechord, \twobullechordbis\mbox{, and } \tripodnonor.
\end{equation*}
As proved in Example~\ref{examplecomconfintthree}, the diagrams $\twoischord$, $\twoxchord$, $\twobullechord$ do not contribute to $\Zinv_2\left(\RR^3,O,p_{S^2}^{\ast}(\omega_{S^2})\right)$. 

Since Lemma~\ref{lemoneleg} implies $[\twobullechordbis]=0$, the diagram $\twobullechordbis$ does not contribute either. 
Lemma~\ref{lemItripod} implies
 \begin{equation*}I\Bigl(S^3,O,\tripodnonor,\bigl(p_{S^2}^{\ast}(\omega_{S^2})\bigr)\Bigr)\left[\tripodnonor\right]=\frac18\left[\tripodnonor\right].\end{equation*}
When $\Gamma=\tripodnonor$, we have $\coefgambet_{\Gamma}=\frac{3!}{6!2^3}$, and there are $\frac13 \frac{6!2^3}{3!}$ numbered graphs of $\Davis^e_2(S^1)$ that are isomorphic to $\Gamma$ as a Jacobi diagram. This concludes the 
computation of $\Zinv_2\bigl(\RR^3,O,p_{S^2}^{\ast}(\omega_{S^2})\bigr)$.

See also Proposition~\ref{propdefhomog} and Example~\ref{examplecomconfinttwocontcont}.
Enore Guadagnini, Maurizio Martellini, and Mihail Mintchev have performed alternative computations of similar quantities in \cite{gmm}.
\end{examples}

The following theorem is proved in Chapter~\ref{chapindepform}. See Section~\ref{secintroindepform}, and Corollary~\ref{corinvone} and Lemma~\ref{lemformprod}, in particular.
\begin{theorem}
\label{thmevenz}
Let $(\crats,\partau)$ be an asymptotic rational homology $\RR^3$. Let $n \in \NN$.
For any $i\in \underline{6n}$, let $\omega(i)$ be a propagating form of $C_2(\rats)$.
Then $\zinv_{2n}(\crats,(\omega(i))_{i \in \underline{6n}})$ is independent of the chosen
$\omega(i)$, it depends only on the diffeomorphism class of $\rats$. It is denoted by $\zinv_{2n}(\rats)$.
\end{theorem}

For odd $n$, $\zinv_{n}(\crats,(\omega(i))_{i \in \underline{3n}})$ depends on the chosen $\omega(i)$. Theorem~\ref{thmfstconst} explains how to deal with this dependence when the $\omega(i)$ are homogeneous propagating forms of $\bigl(C_2(\rats),\partau\bigr)$. Tadayuki Watanabe \cite{watanabeMorse} and Tatsuro Shimizu \cite{shimizu} have studied alternative compensations for this dependence.

We are going to prove the following theorem in the next chapters. The proof will be concluded in the end of Section~\ref{secdephom}.

\begin{theorem}
\label{thmfstconst}
Let $(\crats,\partau)$ be an asymptotic rational homology $\RR^3$. Let $\source=\sqcup_{j=1}^kS^1_j$ be a disjoint union of $k$ circles. 
Let $\Link \colon \source \hookrightarrow \crats$ be an embedding. Let $n \in \NN$.
For any $i\in \underline{3n}$, let $\omega(i)$ be a homogeneous propagating form of $\bigl(C_2(\rats),\partau\bigr)$.
 
Then $\Zinv_n\bigl(\crats,\Link,(\omega(i))_{i \in \underline{3n}}\bigr)$ is independent of the chosen
$\omega(i)$. It depends only on the diffeomorphism class of $(\rats,\Link)$, on $p_1(\partau)$, and on the $I_{\theta}(K_j,\partau)$, for the components $K_j$ of $\Link$. It is denoted by $\Zinv_n(\crats,\Link,\partau)$.\index[N]{ZZ@$\Zinvuf$ and some variants (see also the summary in the next pages)!ZitauR@$\Zinv(\crats,\Link,\partau)$}
More precisely, 
set \begin{equation*}\Zinv(\crats,\Link,\partau)=\bigl(\Zinv_n(\crats,\Link,\partau)\bigr)_{n\in \NN} \in \Aavis(\sqcup_{j=1}^kS^1_j).\end{equation*} 
There exist two constants $\alpha \in \Assis(S^1;\RR)$ \index[N]{azlpha@$\alpha$ anomaly} and $\ansothree \in \Aavis(\emptyset;\RR)$, called \emph{anomalies},\footnote{The anomaly $\alpha$ is defined in Section~\ref{secanomalpha}, and $\beta$ is defined in Definition~\ref{defxin}.} such that $\alpha_{2n}=0$ and $\ansothree_{2n}=0$ for any $n \in \NN$, and
\begin{equation*}\exp\Bigl(-\frac14 p_1(\partau)\ansothree\Bigr)\prod_{j=1}^k\Bigl(\exp\bigl(-I_{\theta}(K_j,\partau)\alpha\bigr)\#_j\Bigr) \Zinv(\crats,\Link,\partau)=\Zinvuf(\rats,\Link)\end{equation*}
depends only on the diffeomorphism class of $(\rats,\Link)$. \index[N]{ZZ@$\Zinvuf$ and some variants (see also the summary in the next pages)!ZuRL@$\Zinvuf(\rats,\Link)$}
Here $\exp(- I_{\theta}(K_j,\partau)\alpha)$ acts on $\Zinv(\crats,\Link,\partau)$, on the copy $S^1_j$ of $S^1$ as indicated by the subscript $j$, as in Proposition~\ref{propdiagcom}.
\end{theorem}
Recall Notation~\ref{notationzZ}.
If $\crats=\RR^3$, then the projection $\Zinvlinkuf(S^3,\
\Link)$ of $\Zinv(S^3,\Link)$ on $\Assis(\sqcup_{j=1}^kS^1_j)$ is a \emph{ universal finite type invariant of links in $\RR^3$}, i.e., the projection of $\Zinvlinkuf_n$ onto $\Asimp_n/(1T)$ satisfies the properties stated for $\Zinvlinkufmodis_n$ in Theorem~\ref{thmbn}. This result, due to Daniel Altsch\"uler and Laurent Freidel \cite{af}, is proved in Section~\ref{secproofthmbn}. See Theorem~\ref{thmunivsingtang}. 

\begin{notation} \label{notzZcheck}
Let $\Zinvlinkuf(\rats,\Link)$ \index[N]{ZZ@$\Zinvuf$ and some variants (see also the summary in the next pages)!ZuTT@$\Zinvlinkuf$}
denote the projection of $\Zinvuf(\rats,\Link)$ on $\Assis(\sqcup_{j=1}^kS^1_j)$, and set
\begin{equation*}\Zinvuf(\rats)= \Zinvuf(\rats,\emptyset).\end{equation*} \index[N]{ZZ@$\Zinvuf$ and some variants (see also the summary in the next pages)!Zurrats@$\Zinvuf(\rats)$}
\end{notation}

The invariant $\Zinvuf$ of rational homology $3$-spheres is the Kontsevich configuration space invariant studied by Greg Kuperberg and Dylan Thurston in \cite{kt}, and described in \cite{lesconst} as $Z_{KKT}$.

We will see in Chapter~\ref{chaprat} that the anomalies $\alpha$ and $\ansothree$ are \emph{rational}, i.e., that $\alpha \in \Assis(S^1;\QQ)$ (in Proposition~\ref{propratalpha}) and $\ansothree \in \Aavis(\emptyset;\QQ)$ (in Theorem~\ref{thmZinvrat}).

\begin{examples}
\label{examplecomconfinttwobis}
 According to Example~\ref{examplecomconfinttwocont}, we have $\Zinv(\RR^3,\emptyset,\taust)= 1 =\left[\emptyset\right]$. Since $p_1(\taust)=0$, we also have $\Zinvuf(S^3)=1$.

For the embedding $O$ of the trivial knot in $\RR^3$ of Example~\ref{examplecomconfintthree}, we have \begin{equation*}I_{\theta}(O,\taust)=0.\end{equation*} So, according to Example~\ref{examplecomconfinttwocont}, we also have $\Zinvuf_0(S^3,O)=1$,  $\Zinvuf_1(S^3,O)=0$, and
\begin{equation*}\Zinvuf_2(S^3,O)=\frac1{24}\left[\tripodnonor\right].\end{equation*}
\end{examples}

\begin{definition}
 Let $\Gamma$ be a Jacobi diagram on an oriented one-manifold $\source$ as in Definition~\ref{defdia}. Let $j \colon U(\Gamma) \hookrightarrow \source$ be a $\Gamma$-compatible injection from its set $U(\Gamma)$ of univalent vertices to its support $\source$.
An \emph{automorphism}\index[T]{automorphism!of a Jacobi diagram} of $\Gamma$ is a permutation of the set $H(\Gamma)$ of half-edges of $\Gamma$ that maps a pair of half-edges of an edge to another such and a triple of half-edges that contain a vertex to another such, and such that, 
for the induced bijection $b$ of $U(\Gamma)$, $j\circ b$ is isotopic to $j$. (So the automorphisms preserve the components of the univalent vertices. They also preserve their linear order on intervals and their cyclic order on circles.) Let $\Aut(\Gamma)$ denote the set of automorphisms of $\Gamma$.
\end{definition}

\begin{examples}
There are six automorphisms of $\tata$ that fix each vertex. They correspond to the permutations of the edges. The cardinality $\cardlef{ \Aut(\tata)}$ of $\Aut(\tata)$ is $12$, and we have
\begin{equation*} \cardBig{ \Aut(\onechordsmalltseul)}=1 \mbox{ and }\cardBig{ \Aut(\tripodnonor)}=3 \end{equation*}
\end{examples}

Recall from Notation~\ref{notationzZ} that
$\coefgambet_{\Gamma}=\frac{(3n-\cardlef{E(\Gamma)})!}{(3n)!2^{\cardlef{E(\Gamma)}}}.$
Also recall that
Lemma~\ref{lemiotaform} ensures the existence of antisymmetric homogeneous propagating forms $\omega$ of $\bigl(C_2(\rats),\partau\bigr)$.

\begin{proposition}
 \label{propdefhomog} Let $\Davis^u_n(\source)$ \index[N]{Diag@Diagram sets!Du@$\Davis^u_n$ unnumbered} denote the set of unnumbered, unoriented degree $n$ Jacobi diagrams on $\source$ without looped edges.
Under the assumptions of Theorem~\ref{thmfstconst}, let $\omega$ be an antisymmetric homogeneous propagating form of $C_2(\rats)$.
Then we have \begin{equation*}\Zinv_n(\crats,\Link,\omega)=\sum_{\Gamma \in \Davis^u_n(\source)} \frac1{\cardlef{\Aut(\Gamma)}}I(\rats,\Link,\Gamma,(\omega)_{i \in \underline{3n}})\left[\Gamma\right].\end{equation*}
\end{proposition}
\bp Set $\omega(i)=\omega$ for any $i$. For a \emph{numbered} graph $\Gamma$ (i.e., a graph equipped with the structure described in Definition~\ref{defnumdia}), there are $\frac{1}{\coefgambet_{\Gamma}}$ ways of \emph{renumbering} it (i.e., changing this structure), and $\cardlef{\Aut(\Gamma)}$ of them will produce the same numbered graph.  Therefore, we have \begin{equation*}\sum_{\Gamma \in \Davis^e_n(\source)}\coefgambet_{\Gamma}I\Bigl(\rats,\Link,\Gamma,\bigl(\omega(i)\bigr)_{i \in \underline{3n}}\Bigr)\left[\Gamma\right]=\sum_{\Gamma \in \Davis^u_n(\source)} \frac1{\cardlef{\Aut(\Gamma)}}I\bigl(\rats,\Link,\Gamma,(\omega)_{i \in \underline{3n}}\bigr)\left[\Gamma\right].\end{equation*} 
\eop

\begin{remark}
\label{rkhomog}
Let $\omega$ be an antisymmetric homogeneous propagating form of $\bigl(C_2(\rats),\partau\bigr)$.
 The homogeneous definition of $\Zinv_n(\crats,\Link,\partau)=\Zinv_n(\crats,\Link,\omega)$ above makes clear that $\Zinv_n(\crats,\Link,\partau)$ is a \say{measure} of graph configurations, where a graph configuration is an embedding of the set of vertices of a uni-trivalent graph into $\crats$, which maps univalent vertices to $\Link(\source)$. The embedded vertices are connected by a set of abstract solid edges, which represent the measuring form. The factor $\frac1{\cardlef{\Aut(\Gamma)}}$ ensures that every such configuration of an unnumbered, unoriented graph is \say{measured} once.
\end{remark}

\begin{lemma}
 \label{lemmultcheck} Recall Notation~\ref{notzZcheck}.
 Under the assumptions of Theorem~\ref{thmfstconst}, we have
 \begin{equation*}\Zinv\bigl(\crats,\Link,\omega(1)\bigr)=\Zinv\bigl(\crats,\emptyset,\omega(1)\bigr) \Zinvlink\bigl(\crats,\Link,\omega(1)\bigr).\end{equation*}
In particular, Theorem~\ref{thmfstconst} will imply
\begin{equation*}\Zinvuf(\rats,\Link)=\Zinvuf(\rats)\Zinvlinkuf(\rats,\Link).\end{equation*}
\end{lemma}
\bp Remark~\ref{rkhomog} shows that the coefficients are correct.
\eop

\begin{example}
 \label{examplecomconfinttwocontcont}
For the embedding $O$ of the trivial knot in $\RR^3$ of Examples~\ref{examplecomconfintthree}, according to the computations performed in the series of examples \ref{examplecomconfinttwocont}, and with the notation of Theorem~\ref{thmfstconst} and Notation~\ref{notationzZ}, as in Example~\ref{examplecomconfinttwobis}, we have
\begin{equation*}\begin{array}{ll}\Zinvuf_2\left(S^3,O\right)
&=\Zinv_2\left(\RR^3,O,p_{S^2}^{\ast}(\omega_{S^2})\right)
\\&=\frac13 I\left(S^3,O,\tripodor,p_{S^2}^{\ast}(\omega_{S^2})\right)\left[\tripodnonor\right]\\&=\frac1{24}\left[\tripodnonor\right].\end{array}\end{equation*}
\end{example}

We end this section by stating Theorems~\ref{thmConwaycircZ} and \ref{thmConwaycircZSthree} about the numerical invariants obtained from $\Zinvlinkuf_n$ by applying the Conway weight system $w_C$ of Example~\ref{exaweightConway}. For $n=2$, we get the invariant $w_2$ discussed in Subsection~\ref{subsecwtwo}.

Since $w_C$ is multiplicative, and since $w_C$ sends elements of odd degree to zero, $w_C$ sends
$\exp(-I_{\theta}(K_j,\partau)\alpha)$ to the unit of $\Assis(S^1)$. So, with the notation of Theorem~\ref{thmfstconst}, we have $w_C\left(\Zinvlinkuf_n(\rats,\Link)\right)=w_C\left(\Zinvlink_n(\crats,\Link,\partau)\right)$ for any $n \in \NN$, and we can forget the anomalies for $w_C \circ \Zinvlinkuf$. Theorem~\ref{thmConwaycircZ} will tell us that we can furthermore omit the homogeneity assumptions on the forms for $w_C \circ \Zinvlinkuf$ and reduce
our averaging process.
We will average only over some degree $n$ graphs whose edges are numbered in $\underline{3n-2}$, or even in $\underline{3}$ when $n=2$, in some cases.

(Since $\projassis \colon \Aavis(\source) \to \Assis(\source)$ maps all graphs with less than $2$ univalent vertices to zero (thanks to Lemma~\ref{lemoneleg}), the graphs that contribute to  $\Zinvlink_n(\crats,\Link,\omega(1))$ have at most $3n-2$ edges. So it is natural to average only over these graphs.)

\begin{notation}
\label{notlessav}
For a finite set $\finseta$, let $\Davisred^e_{n,\finseta}(\source)$\index[N]{Diag@Diagram sets!DeA@$\Davis^e_{n,\finseta}(\sourcetl)$ $\finseta$-numbered} denote the set of $\finseta$-numbered degree $n$ Jacobi diagrams with support $\source$ without looped edges, as in Definition~\ref{defnumdia}. (These diagrams have at most $\cardlef{\finseta}$ edges.) The coefficient $\coefgambetred_{\Gamma}$ associated to a diagram $\Gamma \in \Davisred^e_{n,\finseta}(\source)$ is \index[N]{zwetaGamma@$\coefgambet_{\Gamma}$ averaging coefficient} \begin{equation*}\coefgambet_{\Gamma}=\frac{\bigl(\cardlef{\finseta}-\cardlef{E(\Gamma)}\bigr)!}{\cardlef{\finseta}!2^{\cardlef{E(\Gamma)}}}.\end{equation*} 
For any $i\in \finseta$, let $\omega(i)$ be a propagating form of $C_2(\rats)$. Set
\begin{equation*}\Zinv_{n,\finseta}\Bigl(\crats,\Link,\bigl(\omega(i)\bigr)_{i \in \finseta}\Bigr)=\sum_{\Gamma \in \Davisred^e_{n,\finseta}(\source)}\coefgambetred_{\Gamma}I\Bigl(\rats,\Link,\Gamma,\bigl(\omega(i)\bigr)_{i \in \finseta}\Bigr)\left[\Gamma\right] \in \Aavis_n(\source).\end{equation*} \index[N]{ZZ@$\Zinvuf$ and some variants (see also the summary in the next pages)!ZiomRA@$\Zinv_{n,\finseta}(\crats,\Link,(\omega(i)))$}
For $m \in \NN$, set $\Davisred^e_{n,m}(\source)=\Davisred^e_{n,\underline{m}}(\source)$\index[N]{Diag@Diagram sets!DeAn@$\Davis^e_{n,m}(\sourcetl)$} and \begin{equation*}\Zinv_{n,m}\Bigl(\crats,\Link,\bigl(\omega(i)\bigr)_{i \in \underline{m}}\Bigr)= \Zinv_{n,\underline{m}}\Bigl(\crats,\Link,\bigl(\omega(i)\bigr)_{i \in \underline{m}}\Bigr).\end{equation*} Note that $\Davis^e_n=\Davisred^e_{n,3n}$ and that $\Zinv_{n,3n}=\Zinv_n$.
With the projection $\projassis \colon \Aavis_n(S^1) \to \Assis_n(S^1)$ of Notation~\ref{notationzZ}, set
\begin{equation*}\Zinvlink_{n,m}\Bigl(\crats,K,\bigl(\omega(i)\bigr)_{i \in \underline{m}}\Bigr)=\projassis \biggl( \Zinv_{n,m}\Bigl(\crats,K,\bigl(\omega(i)\bigr)_{i \in \underline{m}}\Bigr)\biggr).\end{equation*}
\end{notation}

\begin{theorem}
\label{thmConwaycircZ}
Let $(\crats,\partau)$ be an asymptotic rational homology $\RR^3$. Let $K \colon S^1 \hookrightarrow \crats$ be an embedding. Let $n \in \NN$. For any $i\in \underline{3n-2}$, let $\omega(i)$ be a propagating form of $C_2(\rats)$.

Let $p^c$ be the projection given by Corollary~\ref{corprojprim} from $\Assis(S^1)$ to the space $\Assis^c(S^1)$ of its primitive elements.
Recall the linear form $w_C \colon (\Assis_n(S^1)=\Asimp_n) \to \RR$ induced by the Conway weight system of Example~\ref{exaweightConway}.  

Then $w_C \left( \Zinvlink_{n,3n-2}\left(\crats,K,(\omega(i))_{i \in \underline{3n-2}}\right)\right)$
and \begin{equation*}w_C \circ p^c  \left( \Zinvlink_{n,3n-2}\left(\crats,K,(\omega(i))_{i \in \underline{3n-2}}\right)\right)\end{equation*}
are independent of the chosen
$\omega(i)$. They depend only on the diffeomorphism class of $(\crats,K)$. They are respectively denoted by $w_C \Zinvlinkuf_n(\rats,K)$ and \begin{equation*}w_C p^c\Zinvlinkuf_n(\rats,K).\end{equation*}
\end{theorem}
At the end of Subsection~\ref{secintroindepform}, the proof of Theorem~\ref{thmConwaycircZ} will be reduced to the proof of Proposition~\ref{propinvonebis}, which is proved in Subsection~\ref{seccancelfaces}.

\begin{remark}
\label{rknoconflict}
Assuming both Theorems~\ref{thmConwaycircZ} and \ref{thmfstconst}, there is no notation conflict between them.
Indeed, with the notation of Theorem~\ref{thmfstconst}, we have
\begin{equation*}
\begin{array}{ll}w_C\bigl(\Zinvlinkuf_n(\rats,\Link)\bigr)&=w_C\bigl(\Zinvlink_n(\crats,\Link,\partau)\bigr)=w_C\bigl(\Zinvlink_n(\crats,\Link,\omega)\bigr)\\&=\sum_{\Gamma \in \Davis^u_n(\source)} \frac1{\cardlef{\Aut(\Gamma)}}I\bigl(\rats,\Link,\Gamma,(\omega)_{i \in \underline{3n}}\bigr)w_C\bigl(\left[\Gamma\right]\bigr)\end{array}\end{equation*} 
for any homogeneous propagating form $\omega$ of $\bigl(C_2(\rats),\partau\bigr)$ and for any $n \in \NN$, thanks to Proposition~\ref{propdefhomog}. When $\omega(i)=\omega$ for all $i$, this is also the expression of \begin{equation*}w_C \Bigl( \Zinvlink_{n,3n-2}\bigl(\crats,K,(\omega(i))_{i \in \underline{3n-2}}\bigr)\Bigr). \end{equation*}
\end{remark}

We prove the following theorem at the end of Section~\ref{seccancelfaces}. It implies Proposition~\ref{propwtwo}.
\begin{theorem}
\label{thmConwaycircZSthree}
Let $K \colon S^1 \hookrightarrow \RR^3$ be a knot embedding. For $i\in \underline{3}$, let 
$\omega_{S^2}(i)$ be a volume-one form of $S^2$, then $w_C \bigl( \Zinvlink_{2,3}\bigl(\RR^3,K,(p_{S^2}^{\ast}(\omega_{S^2}(i)))_{i \in \underline{3}}\bigr)\bigr)$ is independent of the chosen
$\omega_{S^2}(i)$. It is an isotopy invariant of $K$, which coincides with $w_C \Zinvlinkuf_2(S^3,K)$ and $w_C p^c\Zinvlinkuf_2(S^3,K)$.
\end{theorem}

\bpo{Proof of Proposition~\ref{propwtwo} assuming Theorem~\ref{thmConwaycircZSthree}}
The only degree $2$ Jacobi diagrams on $S^1$ without looped edges, with at most three edges, and with no trivalent component
are $\diagcross$, $\diagthetthet$, and $\diagtripod$. We have $w_C\left(\left[\diagcross\right]\right)=1$, $w_C\left(\left[\diagthetthet\right]\right)=0$, $w_C\left(\left[\diagtripod\right]\right) =w_C\left(\left[\diagthetthet\right]\right) -w_C\left(\left[\diagcross\right]\right) =-1$.
There are six elements of $\Davisred^e_{2,3}(S^1)$ isomorphic to $\diagcross$, one for each permutation $\sigma$ of $\underline{3}$. 
They may be drawn as
\begin{center}
 $\Gamma(\sigma)=\begin{tikzpicture}[xscale=1.6] \useasboundingbox (-.7,-.2) rectangle (.7,.5);
\draw [->,dash pattern=on 2pt off 2pt] (.6,0) arc (0:360:.6);
\draw (45:.6) -- (-135:.6);
\draw (135:.6) -- (-45:.6);
\draw [->] (45:.6) -- (-135:.25);
\draw [->] (135:.6) -- (-45:.25) ;
\draw (-.1,-.1) node[left]{\scriptsize $\sigma(2)$} (.1,-.1) node[right]{\scriptsize $\sigma(3)$};
\fill (45:.6) circle (2pt) (135:.6) circle (2pt) (-135:.6) circle (2pt) (-45:.6) circle (2pt);
\end{tikzpicture}$
\end{center}
and we have $\coefgambet_{\Gamma(\sigma)}=\frac{1}{3!\times 4}$.
For such a graph $\Gamma(\sigma)$, we also have
\begin{multline*}I\biggl(S^3,K,\Gamma(\sigma),\Bigl(p_{S^2}^{\ast}\bigl(\omega_{S^2}(i)\bigr)\Bigr)_{i \in \underline{3}}\biggr)\\
=\int_{(S^2)^3} \deg\left(\sigma_{\ast}\left(1_{S^2} \times G_{\diagcross}\right)\right)
\wedge_{i=1}^3 p_i^{\ast}\bigl(\omega_{S^2}(i)\bigr).\end{multline*}
There are 16 elements of $\Davisred^e_{2,3}(S^1)$ isomorphic to $\diagtripod$. They are obtained from the two diagrams
\begin{equation*}
\diagtripoded{>}{2}{>}{3}\mbox{ and }\diagtripoded{>}{3}{>}{2}\end{equation*} by reversing the directions of some edges. For these diagrams $\coefgambet_{\Gamma}=\frac{1}{3!\times 8}$.
Theorem~\ref{thmConwaycircZSthree} implies that for any embedding $K \colon S^1 \hookrightarrow \RR^3$, we have
\begin{equation*}w_C \Zinvlinkuf_2(S^3,K)=\int_{(S^2)^3}\tilde{w}_2(K)\wedge_{i=1}^3 p_i^{\ast}\bigl(\omega_{S^2}(i)\bigr)\end{equation*}
for the locally constant degree map $\tilde{w}_2(K)$ of Proposition~\ref{propwtwo} defined on an open dense subset of $(S^2)^3$.
Any point $(X_1,X_2,X_3)$ of this open dense subset of regular values of $(S^2)^3$, has an open connected neighborhood $\prod_{i=1}^3W_i$ of regular values, and there are volume-one $2$-forms $\omega_{S^2}(i)$ supported on $W_i$.
For such forms, $w_C \Zinvlinkuf_2(S^3,K)=\tilde{w}_2(K)(X_1,X_2,X_3)$. So Theorem~\ref{thmConwaycircZSthree} implies that $\tilde{w}_2(K)$ is constant.
\eop

\section{Straight links}
\label{secstraight}

A one-cycle $c$ of $S^2$ is \emph{algebraically trivial}\index[T]{algebraically trivial!one-cycle} if, for any two points $x$ and $y$ outside its support, the algebraic intersection of an arc from $x$ to $y$ transverse to $c$ with $c$ is zero, or equivalently if the integral of any one-form of $S^2$ along $c$ is zero.

Let $(\crats,\partau)$ be an asymptotic rational homology $\RR^3$. Let $K \colon S^1 \hookrightarrow \crats$ be a knot embedding.
Recall the notation from Proposition~\ref{propprojbord} and Section~\ref{secthetaknot}.
Set $C_K=C(K;\onechordsmall{$S^1$})$. We have $C_K \subset C_2(\rats)$.

\begin{lemma}
\label{lemindepstraight} 
If $\projp_{\partau}(\partial C_K)$ is algebraically trivial, then for any propagating chain $F$ of $\bigl(C_2(\rats),\partau\bigr)$ transverse to $C_K$ and for any propagating form $\omega_p$ of $\bigl(C_2(\rats),\partau\bigr)$, we have
\begin{equation*}\int_{C_K}\omega_p=\langle C_K, F \rangle_{\!C_2(\rats)\,}=I_{\theta}(K,\partau),\end{equation*}
with respect to the definition of $I_{\theta}(K,\partau)$ in Lemma~\ref{lemdefItheta}.
In particular, we have $I_{\theta}(K,\partau) \in \QQ$, and $I_{\theta}(K,\partau)$ is an integer when $\rats$ is an integral homology sphere.
\end{lemma}
\bp According to Lemma~\ref{lemvaritheta}, changing the propagating form $\omega_p$ to $\omega^{\prime}_p$ adds some $\int_{\partial C_K}\projp_{\partau}^{\ast}(\eta_S)=\int_{\projp_{\partau}(\partial C_K)}\eta_S$ for some one-form on $S^2$.
Then by definition, $\int_{C_K}\omega_p$ is independent of the propagating form $\omega_p$ of $\bigl(C_2(\rats),\partau\bigr)$. 
For a propagating chain $F$ of $\bigl(C_2(\rats),\partau\bigr)$ transverse to $C_K$, one can choose a propagating form dual to $F$ and supported near $F$ such that $\int_{C_K}\omega_p=\langle C_K, F \rangle_{\!C_2(\rats)\,}$.
(See the end of Section~\ref{secderhamcohom} and Lemma~\ref{lemconstrformalphadual} in particular, for details.)
The rationality of $I_{\theta}(K,\partau)$ follows from the rationality of $F$. Since $F$ can be chosen to be an integral chain when $\rats$ is a $\ZZ$-sphere, we get $I_{\theta}(K,\partau) \in \ZZ$ in this case. \eop

\begin{remark}
One could have proved that $\langle C_K, F \rangle_{\!C_2(\rats)\,}$ is independent of the chosen propagating chain $F$ of $\bigl(C_2(\rats),\partau\bigr)$ transverse to $C_K$ directly, as follows.
Let $F$ and $F^{\prime}$ be two such propagating chains with respective boundaries $\projp_{\partau}^{-1}(a)$ and $\projp_{\partau}^{-1}(a^{\prime})$, where $a$ and $a^{\prime}$ are in $S^2$.
Let $\left[a,a^{\prime}\right]$ be a path from $a$ to $a^{\prime}$ in $S^2$ transverse to $\projp_{\partau}(\partial C_K)$.
There exists a rational chain $W$ of $C_2(\rats)$ transverse to $C_K$ and to $\partial C_2(\rats)$, whose boundary is $F^{\prime}-F-\projp_{\partau}^{-1}\left(\left[a,a^{\prime}\right]\right)$. Then $\partial C_K \cap W \subset C_K \cap \partial W$, and $\partial (C_K \cap W) = \pm C_K \cap \partial W$. So we have
\begin{equation*}\bigl\langle C_K, \partial W \bigr\rangle_{\!C_2(\rats)} =\bigl\langle C_K, F^{\prime}-F-\projp_{\partau}^{-1}\left(\left[a,a^{\prime}\right]\right) \bigr\rangle_{\!C_2(\rats)}=0.\end{equation*}
So we get
\begin{equation*}\begin{array}{ll}\langle C_K, F^{\prime} \rangle_{\!C_2(\rats)\,} - \langle C_K, F \rangle_{\!C_2(\rats)\,} &=  \langle C_K, \projp_{\partau}^{-1}\left(\left[a,a^{\prime}\right]\right) \rangle_{\!C_2(\rats)\,}\\&=  \langle \partial C_K, \projp_{\partau}^{-1}\left(\left[a,a^{\prime}\right]\right) \rangle_{\!\partial C_2(\rats)\,}
\\&=  \langle \projp_{\partau}(\partial C_K), \left[a,a^{\prime}\right] \rangle_{\!S^2\,}=0.\end{array}\end{equation*}
\end{remark}

\begin{lemma} \label{lemdegloctheta}
Assume that $\projp_{\partau}(\partial C_K)$ is algebraically trivial. 
 Let $Y \in S^2 \setminus \projp_{\partau}(\partial C_K)$.
 Let $\vecZorth \colon K \to S^2$ map $k \in K$ to the vector $\vecZorth(k)$ of $S^2$ orthogonal to $\projp_{\partau}(T_kK)$ in the half great circle of $S^2$ that contains $\projp_{\partau}(T_kK)$, $Y$, and $\projp_{\partau}(-T_kK)$.
Define the parallel $K_{\parallel,\partau,Y}$ by pushing a point $k$ of $K$ in the direction $\partau(\vecZorth(k))$.
Then we have
\begin{equation*}I_{\theta}(K,\partau)=lk(K,K_{\parallel,\partau,Y}).\end{equation*}
\end{lemma}
\bp Thanks to Lemma~\ref{lemlkprop}, $lk(K,K_{\parallel,\partau,Y})$ is the evaluation of any propagator of $C_2(\rats)$ on $K \times K_{\parallel,\partau,Y}$. Let $\omega_p$ be a propagating form of $C_2(\rats)$ (as in Definition~\ref{defpropagatortwo}), which may be expressed as $\projp_{\partau}^{\ast}(\omega_{-Y})$  on $\ST N(K)$, for a $2$-form $\omega_{-Y}$ of $S^2$ supported in a geometric disk in $S^2 \setminus \projp_{\partau}(\partial C_K)$ centered at $(-Y)$. Observe that the intersection of such a disk with the half great circle of $S^2$ that contains $\projp_{\partau}(T_kK)$, $Y$, and $\projp_{\partau}(-T_kK)$ is empty for any $k \in K$.

Compute $lk(K,K_{\parallel,\partau,Y})$ as the limit of $I(\haltereorKKp,\omega_p)$,
when $K_{\parallel,\partau,Y}$ tends to $K$. The configuration space $K \times K_{\parallel,\partau,Y}$ is a torus of $C_2(\rats)$. When $K_{\parallel,\partau,Y}$ tends to $K$, this torus tends to the union of the annulus $C(K;\onechordsmallK)$ and an annulus $K \times_{\partau,Y} J$ contained in $\ST \rats\vert_{K}$, which fibers over $K$ and whose fiber over $k \in K$ contains all the limit directions from $k$ to a close point on $K_{\parallel,\partau,Y}$. This fiber is the half great circle of $\ST \rats\vert_{k}$ that $\projp_{\partau}$ maps to the half great circle of $S^2$ that contains $\projp_{\partau}(T_kK)$, $Y$, and $\projp_{\partau}(-T_kK)$. Thus $K \times K_{\parallel,\partau,Y}$ is homologous to the torus 
\begin{equation*}T=C_K\cup \left(K \times_{\partau,Y} J\right).\end{equation*}
 The integral of $\omega_p$ on $K \times_{\partau,Y} J$ is the integral of $\omega_{-Y}$ along $\projp_{\partau}(K \times_{\partau,Y} J)$. It is zero since $\projp_{\partau}(K \times_{\partau,Y} J)$ does not meet the support of $\omega_{-Y}$. Therefore, we have $lk(K,K_{\parallel,\partau,Y})=\int_{C_K}\omega_p=I_{\theta}(K,\partau)$, thanks to Lemma~\ref{lemindepstraight}.
\eop

An isotopy class of parallels of a knot is called a \emph{framing}\index[T]{framing!of a knot} (or a \emph{parallelization}) of a knot.

\begin{corollary} \label{cordegloctheta}
A knot embedding $K$ such that $\projp_{\partau}\left(\partial C_K\right)$ is algebraically trivial, with respect to a parallelization $\partau$, has a canonical framing induced by $\partau$, which is the framing induced by a parallel $K_{\parallel,\partau,Y}$ for an arbitrary $Y \in S^2 \setminus \projp_{\partau}(\partial C_K)$.
For such a knot embedding $K$, for any propagating chain $F$ of $\bigl(C_2(\rats),\partau\bigr)$ transverse to $C_K$, and for any propagating form $\omega_p$ of $\bigl(C_2(\rats),\partau\bigr)$, we have
\begin{equation*}\int_{C_K}\omega_p=\langle C_K, F \rangle_{\!C_2(\rats)\,}=I_{\theta}(K,\partau)=lk(K,K_{\parallel,\partau,Y}).\end{equation*}
\end{corollary}
\bp Since the linking number $lk(K,K_{\parallel,\partau,Y})$ determines the framing, the corollary is a direct consequence of Lemmas~\ref{lemindepstraight} and \ref{lemdegloctheta}.
\eop

\begin{definition} \label{defstraightlink}
 A knot embedding $K \colon S^1 \hookrightarrow \crats$ is {\em straight\/} with respect to $\partau$ if the curve $\projp_{\partau}(\ST^+K)$ of $S^2$ is algebraically trivial (with the notation from Proposition~\ref{propprojbord} and Section~\ref{secthetaknot}).\index[T]{straight!link} A link embedding is {\em straight\/} with respect to $\partau$ if all its components are.
\end{definition}

Straight knot embeddings and almost-horizontal knot embeddings in $\RR^3$ (defined before Lemma~\ref{lemalmosthor}) are examples of knot embeddings $K$ such that $\projp_{\partau}(\partial C_K)$ is algebraically trivial.
Therefore, Lemma~\ref{lemalmosthor} is a particular case of the above corollary.
As a second corollary, we get the following lemma.
\begin{lemma}
\label{lemstrongstraight}
For any knot embedding $K$ in $\crats$, there exists an asymptotically standard parallelization $\tilde{\partau}$ homotopic to $\partau$ such that the first vector $\tilde{\partau}(.;(1,0,0))$ of $\tilde{\partau}$ is tangent to $K$ and orients $K$. In this case, let $K_{\tilde{\partau}}$ be the parallel of $K$ obtained by pushing $K$ in the direction of the second vector $\tilde{\partau}(.;(0,1,0))$ of $\tilde{\partau}$.
Then we have
\begin{equation*}I_{\theta}(K,\tilde{\partau})=lk(K,K_{\tilde{\partau}}).\end{equation*}
\end{lemma}
\bp
In order to obtain $\tilde{\partau}$, it suffices to perform a homotopy of $\partau$ around the image of $K$ so that the first vector of $\tilde{\partau}$ becomes tangent to $K$ along $K$. Thus, $K$ is straight with respect to $(\rats,\tilde{\partau})$. Apply Corollary~\ref{cordegloctheta}.
\eop

\begin{lemma}
\label{lemstraightlk}
Let $K_0\colon \{0\} \times S^1 \rightarrow \crats$ be a straight embedding with respect to $\partau$.
Let $K\colon \left[0,1\right] \times S^1 \rightarrow \crats$ be an embedding such that its restriction  $K_1\colon \{1\} \times S^1 \rightarrow \crats$ is straight with respect to an asymptotically standard parallelization $\partau_1$ homotopic to $\partau$. Then $\bigl(I_{\theta}(K_1,\partau_1)-I_{\theta}(K_0,\partau)\bigr)$ is an even integer.
Furthermore, for any real number $x$ congruent to $I_{\theta}(K_0,\partau)$ modulo $2$, there exists a straight embedding $K_1$ isotopic to $K_0$ such that $I_{\theta}(K_1,\partau)=x$.
\end{lemma}
\bp Let $H \colon t \mapsto \partau_t$ be a smooth homotopy from $\partau=\partau_0$ to $\partau_1$.
Let $p_H \colon \left[0,1\right] \times \partial C_2(\rats) \to S^2$ be the smooth map that restricts
to $\{t\} \times \partial C_2(\rats)$ as $\projp_{\partau_t}$.
There is a closed $2$-form $\omega$ on $\left[0,1\right] \times C_2(\rats)$ that restricts to 
$\left[0,1\right] \times  \partial C_2(\rats)$ as $p_{H}^{\ast}(\omega_{S^2})$.
(Such a form may be obtained by modifying $p_{C_2(\rats)}^{\ast}\projp_{\partau}^{\ast}(\omega_{S^2})$ in a collar neighborhood of $\ST \crats$ using the homotopy $H$.)
Then the integral of $\omega$ over
\begin{equation*}\partial \bigl( \cup_{t \in \left[0,1\right]} C(K;\onechordsmall{$K_t$})\bigr) = C(K;\onechordsmall{$K_1$}) - C(K;\onechordsmall{$K_0$}) - \cup_{t \in \left[0,1\right]}\partial C(K;\onechordsmall{$K_t$}) \end{equation*}
vanishes. So we have
\begin{equation*}I_{\theta}(K_1,\partau_1)-I_{\theta}(K_0,\partau)=\int_{\cup_{t \in \left[0,1\right]}\partial C(K;\onechordsmall{$K_t$})}\omega.\end{equation*} This is the area in $S
^2$ of the integral cycle
$\cup_{t \in \left[0,1\right]}\projp_{\partau_t}\bigl(\partial C(K;\onechordsmall{$K_t$})\bigr)$.
This cycle is the union of the two integral cycles \begin{equation*}\cup_{t \in \left[0,1\right]}\projp_{\partau_t}\bigl(\ST^+(K_t)\bigr)\mbox{ and } \cup_{t \in \left[0,1\right]}\projp_{\partau_t}\bigl(-\ST^-(K_t)\bigr)\end{equation*} which have the same integral area. So $\bigl(I_{\theta}(K_1,\partau_1)-I_{\theta}(K_0,\partau)\bigr)$ is an even integer.

Adding two small almost-horizontal kinks in a standard ball as in the end of Subsection~\ref{subnewonechord} 
which turn in opposite direction and contribute with the same crossing sign like ($\pkink$ and $\pkinks$) or ($\nkink$ and $\nkinks$) preserves straightness and adds $\pm 2$ to $I_{\theta}$.
\eop

\section{Second definition of \texorpdfstring{$\Zinv$}{Z}}
\label{secdefdeux}

Let us state another version of Theorem~\ref{thmfstconst} using straight links instead of homogeneous propagating forms. Recall $\coefgambet_{\Gamma}=\frac{(3n-\cardlef{E(\Gamma)})!}{(3n)!2^{\cardlef{E(\Gamma)}}}.$ 

\begin{theorem}
\label{thmconststraight}
Let $(\crats,\partau)$ be an asymptotic rational homology $\RR^3$. Let $\source=\sqcup_{j=1}^kS^1_j$ be a disjoint union of $k$ circles.
Let $\Link \colon \source \hookrightarrow \crats$ be a straight embedding with respect to $\partau$.
For any $i\in \underline{3n}$, let $\omega(i)$ be a propagating form of $\bigl(C_2(\rats),\partau\bigr)$.
Then the element $\Zinv_n\left(\crats,\Link,(\omega(i))_{i \in \underline{3n}}\right)$ of $\Aavis_n(\source)$ defined in Notation~\ref{notationzZ}  
is independent of the chosen propagating forms
$\omega(i)$ of $\bigl(C_2(\rats),\partau\bigr)$. It is denoted by $\Zinv^s_n(\crats,\Link,\partau)$.
In particular, with the notation of Theorem~\ref{thmfstconst}, we have \begin{equation*}\Zinv^s_n(\crats,\Link,\partau)=\Zinv_n(\crats,\Link,\partau).\end{equation*}
\end{theorem}

We will give the proof of this theorem in Section~\ref{secdepstraight}.
Theorem~\ref{thmtangconstcompar} shows how $\Zinv_n\left(\crats,\Link,(\omega(i))_{i \in \underline{3n}}\right)$ depends on the propagators without the straightness assumption.

Straight links $\Link$ with respect to $\partau$ are framed links $(\Link,\Link_{\parallel,\partau})$ according to Corollary~\ref{cordegloctheta} and Lemma~\ref{lemdegloctheta}. So we can keep the information from the link framing and define the
invariant of straight links
\begin{equation*}\Zinvufrfneg(\crats,\Link,\Link_{\parallel,\partau})=\exp\Bigl(-\frac14 p_1(\partau)\ansothree\Bigr)\Zinv(\crats,\Link,\partau),\end{equation*}
which depends only on $(\crats,\Link)$ and on the $lk(K_j,K_{j\parallel,\partau})$ for the components $K_j$ of $\Link$, according to Theorem~\ref{thmfstconst} and Corollary~\ref{cordegloctheta}. 

\begin{definition}
\label{defZinvufrflink}Recall the invariant $\Zinvuf$ of Theorem~\ref{thmfstconst}.
Define the framed link invariant $\Zinvufrf$ to satisfy \index[N]{ZZ@$\Zinvuf$ and some variants (see also the summary in the next pages)!ZzfL@$\Zinvufrfneg(\crats,\Link,\Link_{\parallel})$}
\begin{equation*}\Zinvufrfneg\Bigl(\crats,\sqcup_{j=1}^k K_j,\sqcup_{j=1}^k K_{j\parallel}\Bigr)=
\prod_{j=1}^k\Bigl(\exp\bigl(lk(K_j,K_{j\parallel})\alpha\bigr)\#_j\Bigr) \Zinvuf\bigl(\rats,\sqcup_{j=1}^k K_j\bigr),\end{equation*}
for a link $\sqcup_{j=1}^k K_j$ equipped with a parallel $K_{j\parallel}$ for each component $K_j$.
\end{definition}

Thanks to Theorem~\ref{thmfstconst}, Corollary~\ref{cordegloctheta}, and Theorem~\ref{thmconststraight}, both definitions coincide for straight framed links $(\Link,\Link_{\parallel,\partau})$. 

Again, we can reduce our averaging process when projecting to $\Assis(\source)$ and get the following theorem, also proved in Section~\ref{secdepstraight}.
\begin{theorem}
\label{thmconststraightred}
Under the assumptions of Theorem~\ref{thmconststraight}, we have
\begin{equation*}\projassis\Bigl(\Zinv^s_n(\crats,\Link,\partau)\Bigr)=\sum_{\Gamma \in \Davisred^e_{n,3n-2}(\source)}\coefgambetred_{\Gamma}I\Bigl(\rats,\Link,\Gamma,\bigl(\omega(i)\bigr)_{i \in \underline{3n-2}}\Bigr)\projassis\bigl(\left[\Gamma\right]\bigr).\end{equation*}
\end{theorem}

\chapter{Compactifications of configuration spaces}
\label{chapcompconf}

We use compactifications of our open configuration spaces $\check{C}(\rats,\Link;\Gamma)$ 
associated to links $\Link$ in a $\QQ$-sphere $\rats$ to study the behavior of our integrals and their dependence on the choice of propagating forms.
More specifically, we prove the convergence of the integrals involved in the definitions of $\Zinv$ by finding a smooth
compactification (with boundary and ridges), to which the integrated forms extend smoothly.
The variation of an integral under the addition of an exact form $d\eta$ is the integral of $\eta$ on the codimension-one faces of the boundary, which need to be identified precisely. 
Therefore, the proofs of Theorems~\ref{thmconststraight} and \ref{thmfstconst} 
(and their variants with reduced averages) require a deep knowledge of configuration spaces. We give all the useful statements in Sections~\ref{secfirstprescomp} to \ref{secfacecodimone}. We prove all of them in Sections~\ref{secprooffacanom} to \ref{secproofblowup}. 

Before giving all the required general statements, we present the main features of the involved compactifications 
in some examples. William Fulton and Robert MacPherson introduced similar compactifications in \cite{FultonMcP}. Scott Axelrod and Isadore Singer \cite[Section 5]{axelsingII}, Maxim Kontsevich~\cite{ko}, Raoul Bott and Clifford Taubes \cite{botttaubes}, and Dev Sinha \cite{Sinha} also investigated these compactifications.

\section{An informal introduction}
\label{secinformalcompactif}

Our first example of a compactification ${C}(\rats,\Link;\Gamma)$ of a configuration space $\check{C}(\rats,\Link;\Gamma)$ is the closed annulus ${C}(K;\tatak)={C}(S^3,K;\tatak)$ of Subsection~\ref{subnewonechord}. 
Our (more interesting) second example is the compactification $C_2(\rats)={C}(\rats,\emptyset;\tata)$ studied in Section~\ref{secCtwo}. Note that ${C}(K;\tatak)$ is the closure of $\check{C}(K;\tatak)$ in $C_2(\rats)$.

As in the example of $C_2(\rats)={C}(\rats,\emptyset;\tata)$, for a trivalent Jacobi diagram $\Gamma$, our general compactification ${C}(\rats,\emptyset;\Gamma)$  of the space $\check{C}(\rats,\emptyset;\Gamma)=\check{C}_{\finsetv(\Gamma)}({\rats})$ of injective maps from ${\finsetv(\Gamma)}$ to $\crats$ depends only on the finite set $\finsetv=\finsetv(\Gamma)$ of vertices of $\Gamma$. It is denoted by ${C}_{\finsetv}({\rats})$. As in the example of ${C}(S^3,K;\tatak)$, for a Jacobi diagram $\Gamma$ on the domain of a link $\Link$ of $\crats$, we define the compactification ${C}(\rats,\Link;\Gamma)$ of $\check{C}(\rats,\Link;\Gamma)$ to be the closure of 
$\check{C}(\rats,\Link;\Gamma)$ in ${C}_{\finsetv(\Gamma)}({\rats})$ in Proposition~\ref{propcompconfL}. This is why we first study ${C}_{\finsetv}({\rats})$ by generalizing the construction of $C_2(\rats)$ performed in Section~\ref{secCtwo}. In this general case, we start with $\rats^{\finsetv}$ and blow up all the diagonals and all the loci that involve $\infty$, in the sense of Section~\ref{secblowup}, following a process precisely described in Theorem~\ref{thmcompconf}.

In this informal introduction, we forget about $\infty$ and first discuss how we successively blow up the diagonals in the manifold $(\RR^3)^{\finsetv}$ of maps from $\finsetv$ to $\RR^3$. Thus, we get the preimage ${C}_{\finsetv}\!\left[\RR^3\right]$ of $(\RR^3)^{\finsetv}$ in ${C}_{\finsetv(\Gamma)}({S^3})$ under the composition valued in $(S^3)^{\finsetv}$ of the blowdown maps.

\subsection{On \texorpdfstring{the}{a} configuration space \texorpdfstring{${C}_{\finsetv}\!\left[\RR^3\right]$}{} of four points in \texorpdfstring{$\RR^3$}{the ambient space}}

The diagonal $\diag_{\finsetv}((\RR^3)^{\finsetv})$ is the set of constant maps from $\finsetv$ to $\RR^3$.
The fiber of the normal bundle to the vector space $\diag_{\finsetv}((\RR^3)^{\finsetv})$ in $(\RR^3)^{\finsetv}$ is the vector space $(\RR^3)^{\finsetv}/\diag_{\finsetv}((\RR^3)^{\finsetv})$. It consists of maps $\confc$ from $\finsetv$ to $\RR^3$ up to global translation.\footnote{A map $\confc$ is identified with the map $(v \mapsto \confc(v) + W)$ for any $W \in \RR^3$.} 
The fiber of the unit normal bundle to the diagonal $\diag_{\finsetv}((\RR^3)^{\finsetv})$ is the space $\cuptd{\finsetv}{\RR^3}$ of nonconstant maps from $\finsetv=\{v_1,v_2,\dots, v_{\cardlef{\finsetv}}\}$ to $\RR^3$ up to (global) translation and up to dilation.\footnote{In $\cuptd{\finsetv}{\RR^3}$, a map $\confc$ is furthermore identified with the map $(v \mapsto \lambda\confc(v))$ for any $\lambda \in \left]0,+\infty\right[$.} The space $\cuptd{\finsetv}{\RR^3}$, studied in Section~\ref{secpresfacanom}, can be identified with the space of maps $w \colon \finsetv \to \RR^3$ that map $v_1$ to $0$ and such that $\sum_{i=2}^{\cardlef{\finsetv}}\norm{w(v_i)}^2 =1$. It is diffeomorphic to a sphere of dimension $\left(3 \cardlef{\finsetv} -4\right)$. Let $B\!\ell_1=\blowup{(\RR^3)^{\finsetv}}{\diag_{\finsetv}((\RR^3)^{\finsetv})}$ be obtained from $(\RR^3)^{\finsetv}$ by blowing up 
$\diag_{\finsetv}((\RR^3)^{\finsetv})$.
We have a diffeomorphism
\begin{equation*}\psi \colon \RR^3 \times \left[0, \infty\right[ \times \cuptd{\finsetv}{\RR^3} \to B\!\ell_1,\end{equation*}
which maps $(u,\mu,w)$ to the map $\confc \colon \finsetv \to \RR^3$ such that
$\confc(v_i)=u+\mu w(v_i)$ for any $i \in \underline{\cardlef{\finsetv}}$, with respect to the identification above. This map $\confc$ is furthermore equipped with the data of the map $w$ when $\confc$ is constant, or equivalently when $\mu=0$.
In particular, this first blow-up equips each constant map $\confc_0$ in the manifold $B\!\ell_1$ with the additional data of a nonconstant map $w \colon \finsetv \to \RR^3$ up to translation and dilation.  Let $\confc_0$ be the constant map which maps $\finsetv$ to $u$, and let $(\confc_0,w)$ denote $\psi(u,0,w)$, then $(\confc_0,w)$ is the limit in $B\!\ell_1$ of $\psi(u,t,w)$, when $t>0$  tends to $0$. Therefore, we can think of the map $w$ as an infinitesimal configuration. The first blow-up provides a magnifying glass, allowing us to see this infinitesimal configuration $w$ of the vertices mapped to the same point $u$ in $\RR^3$. 

When $\cardlef{\finsetv}=2$, we are done and $B\!\ell_1$ is the preimage ${C}_{\finsetv}\!\left[\RR^3\right]$ of $(\RR^3)^{2}$ in ${C}_{2}({S^3})$. In general, we blow up the other \emph{diagonals} $\diag_{\finseta}\left((\RR^3)^{\finsetv}\right)$ for 
all subsets $\finseta$ of $\finsetv$, where $\diag_{\finseta}\left((\RR^3)^{\finsetv}\right)$
is the subspace of $(\RR^3)^{\finsetv}$ consisting of the maps $\confc$ that map $\finseta$ to a single element and such that $\confc(\finsetv \setminus \finseta) \subset \RR^3 \setminus \confc(\finseta)$, as in Section~\ref{secfirstprescomp}.

Let us describe the process when $\finsetv=\finsetv(\Gamma)=\{v_1,v_2,v_3,v_4\}$. 
From now on, we restrict to this case in this subsection.
The closure of the diagonal $\diag_{\finseta}\left((\RR^3)^{\finsetv}\right)$ in  $B\!\ell_1$ is $\psi \left( \RR^3 \times \left[0, \infty\right[ \times \diag_{\finseta}\left(\cuptd{\finsetv}{\RR^3}\right)\right)$,
where $\diag_{\finseta}\left(\cuptd{\finsetv}{\RR^3}\right)$ is the subspace of $\cuptd{\finsetv}{\RR^3}$ made of the classes of its maps $\confc$ constant on $\finseta$ such that $\confc(\finsetv \setminus \finseta) \subset \RR^3 \setminus \confc(\finseta)$. In particular, the closures of the diagonals $\diag_{\finseta}\left((\RR^3)^{\finsetv}\right)$ in  $B\!\ell_1$ for the subsets $\finseta$ of $\finsetv$ of cardinality $3$ are disjoint in $B\!\ell_1$.
We blow up these closures in an arbitrary order. Since the closures are disjoint, the order of the blow-ups does not affect the result.

The fiber of the unit normal bundle of $\diag_{\finseta}\left((\RR^3)^{\finsetv}\right)$ in $B\!\ell_1$ is the space $\cuptd{\finseta}{\RR^3}$ of nonconstant maps from $\finseta$ to $\RR^3$ up to (global) translation and dilation.
Let $\finsetb=\finseta_{123}= \{v_1,v_2,v_3\}$. View $\cuptd{123}{\RR^3}=\cuptd{\finsetb}{\RR^3}$ as the space of maps $w_{123} \colon \finsetv \to \RR^3$ that map $v_1$ and $v_4$ to $0$ and such that $\norm{w_{123}(v_2)}^2 +\norm{w_{123}(v_3)}^2 =1$.
We have a smooth embedding
\begin{equation*}\psi_2 \colon \RR^3 \times \left[0, \infty\right[^2 \times S^2 \times \cuptd{\finsetb}{\RR^3} \to \bigblowup{B\!\ell_1}{\overline{\diag_{\finsetb}((\RR^3)^{\finsetv})}}\end{equation*}
that maps $(u,\mu,\mu_{123},W_4,w_{123})$ to the map $\confc \colon \finsetv \to \RR^3$ such that
$\confc(v_4)=u+\frac{\mu}{\sqrt{1+\mu_{123}^2}}W_4$ and $\confc(v_i)=u+\frac{\mu}{\sqrt{1+\mu_{123}^2}} \mu_{123} w_{123}(v_i)$ for $i \in \underline{3}$. 
The preimage of $\overline{\diag_{\finsetb}((\RR^3)^{\finsetv})}$ under the blowdown map $\blowup{B\!\ell_1}{\overline{\diag_{\finsetb}((\RR^3)^{\finsetv})}} \to B\!\ell_1$ is $\psi_2 \left(\RR^3 \times \left[0, \infty\right[ \times \{0\} \times S^2 \times \cuptd{\finsetb}{\RR^3}\right)$, where $\mu_{123}=0$. 
The image of $\psi_2$ contains a neighborhood of this preimage in $\blowup{B\!\ell_1}{\overline{\diag_{\finsetb}((\RR^3)^{\finsetv})}}$.

Here, this blow-up equips a map $\confc$ of $\diag_{\finsetb}((\RR^3)^{\finsetv})$ with the additional data of the (infinitesimal, nonconstant) configuration $w_{123}\vert_{\finsetb} \colon \finsetb \to \RR^3$ up to translation and dilation. It equips a constant map $\confc_0$ with value $u$ in the closure of $\diag_{\finsetb}((\RR^3)^{\finsetv})$ in $B\!\ell_1$
with such a configuration $w_{123}\vert_{\finsetb}$ in addition to the former $w$, which maps $\finsetb$ to $0$ and $v_4$ to some $W_4 \in S^2$. In this case, we denote the obtained configuration $\psi_2(u,0,0,W_4,w_{123})$ by $(\confc,w,w_{123}\vert_{\finsetb})$, and we have three observation scales.
All the $\confc(v_i)$ coincide at the first initial scale.
At the second (infinitely smaller) scale $w$ in $B\!\ell_1$, the $w(v_i)$ coincide for $i \in \underline{3}$, but $\left(w(v_4)-w(v_1)\right)$ is not zero, and its direction is a vector $W_4 \in S^2$.
At the third scale (infinitely smaller than the second one) in $\blowup{B\!\ell_1}{\overline{\diag_{\finsetb}((\RR^3)^{\finsetv})}}$, the configuration $w_{123}\vert_{\finsetb}$ of $\{v_1,v_2,v_3\}$ is visible up to global translation and dilation. The first three pictures of Figure~\ref{figmagnifyingglasses}
show these three observation scales.

\bfig
\centering
\begin{tikzpicture} 

\draw (0,-1.8) node{In $(\RR^3)^{\finsetv}$} (0,0) node[below]{\scriptsize $\finsetv$} (0,-.25) node[below]{\scriptsize $\shortparallel$} (0,-.5) node[below]{\scriptsize $\{v_1,v_2,v_3,v_4\}$};
 \fill (0,0) circle (1.5 pt);
\begin{scope}[xshift=5cm]
\draw [draw=black,double=white,very thick] (0,0) circle (1.5 cm) (-128:2.4) -- (-128:1.5);
\draw (0,-1.8) node{In $B\!\ell_1$} (0,0) node[below]{ $\finseta_{123}$} (48:1.35) node[below]{\scriptsize $v_4$} (48:1.4) node[left]{\scriptsize $W_4$};
\fill (0,0) circle (1.5 pt) (48:1.35) circle (1.5 pt);
\end{scope}
\begin{scope}[xshift=2.5cm, yshift=-4.2cm]
\draw [draw=gray!50,double=white,very thick] (0,0) circle (1.5 cm) (-128:2.4) -- (-128:1.5);
\draw [draw=black,double=white,very thick] (0,0) circle (1.1 cm) (-138:2.2) -- (-138:1.1);
\draw (0,-1.8) node{In $B\!\ell_2$} (0,0) node[below]{\scriptsize $v_1$} (-.5,0) node[above]{\scriptsize $\finseta_{23}$}(48:1.42) node[left]{\scriptsize $v_4$} ;
 \fill (0,0) circle (1.5 pt) (48:1.35) circle (1.5 pt) (-.55,.1) circle (1.5 pt);
\end{scope}
\begin{scope}[xshift=7.5cm, yshift=-4.2cm]
\draw [draw=gray!20,double=white,very thick] (0,0) circle (1.5 cm) (-128:2.4) -- (-128:1.5);
\draw [draw=gray!50,double=white,very thick] (0,0) circle (1.1 cm) (-138:2.2) -- (-138:1.1);
\begin{scope}[xshift=-.55cm,yshift=.1]
\draw [draw=black,double=white,very thick] (0,0) circle (.4 cm) (-138:1.7) -- (-138:.4);
\draw (-.05,0.05) node[below]{\scriptsize $v_2$} (.3,.2) node[left]{\scriptsize $v_3$};
\fill (0,0) circle (1.5 pt) (.2,.1) circle (1.5 pt);
\end{scope}
\draw (0,-1.8) node{In ${C}_{\finsetv}\!\left[\RR^3\right]$} (0,0) node[right]{\scriptsize $v_1$} (48:1.42) node[left]{\scriptsize $v_4$};
 \fill (0,0) circle (1.5 pt) (48:1.35) circle (1.5 pt);
\end{scope}
\end{tikzpicture}
\caption{The magnifying glasses provided by the successive blow-ups from $(\RR^3)^{\finsetv}$ to ${C}_{\finsetv}\!\left[\RR^3\right]$, for a configuration $(\confc,w,w_{123},w_{23})$ }
\label{figmagnifyingglasses}
\end{figure}

Let $B\!\ell_2$ be the manifold obtained by blowing up the four closures of the diagonals $\diag_{\finseta}\left((\RR^3)^{\finsetv}\right)$ in  $B\!\ell_1$ for the subsets $\finseta$ of $\finsetv$ of cardinality $3$. We have local charts similar to $\psi_2$ for $B\!\ell_2$ in the neighborhoods of the disjoint blown-up loci.

Finally, we blow up (the preimages under the composition of the previous blowdown maps of) the closures of the diagonals $\diag_{\finseta}\left((\RR^3)^{\finsetv}\right)$ in $B\!\ell_2$ for the subsets $\finseta$ of $\finsetv$ of cardinality $2$, in an arbitrary order, to get the manifold ${C}_{\finsetv}\!\left[\RR^3\right]$ of Section~\ref{secblodiag}. Here the diagonals $\diag_{\{v_1,v_2\}}((\RR^3)^{\finsetv})$ and $\diag_{\{v_3,v_4\}}((\RR^3)^{\finsetv})$ are no longer disjoint.
Nevertheless, the blow-up operations associated to $\finseta_{12}=\{v_1,v_2\}$ and $\finseta_{34}=\{v_3,v_4\}$, which act on different coordinates, commute. In the neighborhood of the intersection of the corresponding blown-up loci in ${C}_{\finsetv}\!\left[\RR^3\right]$,
we have an embedding \begin{equation*}\psi_3 \colon \RR^3 \times \left[0, \infty\right[ \times \left[0,1/3\right[^2 \times (S^2)^3  \to {C}_{\finsetv}\!\left[\RR^3\right],\end{equation*}
which maps \begin{equation*}(u, \mu^{\prime}, \mu_{12},\mu_{34},W_3, W_{12}, W_{34})\end{equation*}
to the map $\confc \colon \finsetv \to \RR^3$ such that
$\confc(v_1)=u$, $\confc(v_2)=u+\mu^{\prime} \mu_{12} W_{12}$, $\confc(v_3)=u+\mu^{\prime} W_3$, $\confc(v_4)=u+\mu^{\prime} (W_3 +\mu_{34} W_{34})$ where $W_3\in S^2 = \cuptd{\{\finseta_{12},\finseta_{34}\}}{\RR^3}$, $W_{12}\in S^2 = \cuptd{\finseta_{12}}{\RR^3}$, and $W_{34}\in S^2 = \cuptd{\finseta_{34}}{\RR^3}$.
The configuration $\confc$ is
equipped with the map $w$ up to translation and dilation when $\confc$ is constant. (A representative $w$ of this map sends
the elements of $V$ to $w(v_1)=0$, $w(v_2)=\mu_{12} W_{12}$, $w(v_3)=W_3$, and  $w(v_4)=W_3 + \mu_{34} W_{34}$.) 
The configuration $\confc$ is equipped with $W_{12}$ when the restriction of $w$ (or $\confc$) to $\finseta_{12}$ is constant.
It is equipped with $W_{34}$, when the restriction of $w$ (or $\confc$) to $\finseta_{34}$ is constant. Figure~\ref{figmagnifyingglassestwo} shows the three magnifying glasses that have popped up for a configuration $(\confc, w, W_{12}, W_{34})=\psi_3(u,0,0,0,W_3, W_{12}, W_{34})$.

\bfig
\centering

\begin{tikzpicture} 
\begin{scope}[xshift=-4cm]
\draw (0,-1.8) node{In $(\RR^3)^{\finsetv}$} (0,0) node[below]{\scriptsize $\finsetv=\{v_1,v_2,v_3,v_4\}$};
 \fill (0,0) circle (1.5 pt);
\end{scope}
\draw [draw=black,double=white,very thick] (0,0) circle (1.5 cm) (-158:2.4) -- (-158:1.5);
\draw (0,-1.8) node{In $B\!\ell_1$ and in $B\!\ell_2$} (0,.05) node[below]{\scriptsize $\finseta_{12}$} (-1.54,.32) node[right]{\scriptsize $\finseta_{34}$} (-1.2,.1) node[below]{\scriptsize $W_3$};
\fill (0,0) circle (1.5 pt) (-1.3,.1) circle (1.5 pt);
\begin{scope}[xshift=4cm]
\draw [draw=gray!30,double=white,very thick] (0,0) circle (1.5 cm) (-158:2.4) -- (-158:1.5);

\begin{scope}[xshift=-.65cm,yshift=.05]
\draw [draw=black,double=white,very thick] (0,0) circle (.57 cm) (-170:1.7) -- (-170:.57);
\draw (-.05,0.05) node[below]{\scriptsize $v_3$} (-.2,.3) node[below]{\scriptsize $v_4$} (-.3,.3) node[right]{\scriptsize $W_{34}$};
\fill (0,0) circle (1.5 pt) (-.2,.25) circle (1.5 pt);
\end{scope}
\begin{scope}[xshift=.65cm,yshift=-.05]
\draw [draw=black,double=white,very thick] (0,0) circle (.57 cm) (-42:1.7) -- (-42:.57);
\draw (-.05,0.05) node[below]{\scriptsize $v_1$} (.5,.28) node[left]{\scriptsize $W_{12}$} (.1,-.05) node[right]{\scriptsize $v_2$};
\fill (0,0) circle (1.5 pt) (.2,.1) circle (1.5 pt);
\end{scope}
\draw (0,-1.8) node{In ${C}_{\finsetv}\!\left[\RR^3\right]$};
\end{scope}
\end{tikzpicture}
\caption{The magnifying glasses provided by the successive blow-ups from $(\RR^3)^{\finsetv}$ to ${C}_{\finsetv}\!\left[\RR^3\right]$, for a configuration $(\confc, w, W_{12}, W_{34})$}
\label{figmagnifyingglassestwo}
\end{figure}

For a constant map $\confc$, whose associated infinitesimal $w$ is constant on $\finseta_{123}$, and whose next associated $w_{123}$ is constant on $\finseta_{23}=\{v_2,v_3\}$, the third blow-up family provides a fourth smaller observation scale.
With this additional scale, we see a nonconstant map $w_{23} \colon \finseta_{23} \to \RR^3$ up to global translation and dilation, as in Figure~\ref{figmagnifyingglasses}. The map $w_{23}$ gives the direction from $w_{23}(v_2)$ to $w_{23}(v_3)$.

\subsection{More configuration spaces and their stratifications}
\label{submoreconfstrat}

In general, for a finite set $\finsetv$ and for an integer $d \in \NN \setminus \{0\}$, we transform $(\RR^d)^{\finsetv}$ to a manifold ${C}_{\finsetv}\!\left[\RR^d\right]$ by successively blowing up the closures of (the preimages under the composition of the previous blowdown maps of) the diagonals $\diag_{\finseta}((\RR^d)^{\finsetv})$ associated to the subsets $\finseta$ of $\finsetv$ of cardinality $k$, for $k= \cardlef{\finsetv}, \cardlef{\finsetv} -1, \dots, 2$ in this decreasing order, in Section~\ref{secblodiag}. It will follow from Theorem~\ref{thmcompdiaggen} and Proposition~\ref{proprest} that there are natural smooth (quotients of) restriction maps from ${C}_{\finsetv}\!\left[\RR^d\right]$ to the space $\cuptd{\finseta}{\RR^d}$ of nonconstant maps from $\finseta$ to $\RR^d$ up to translation and dilation, for every subset $\finseta$ of $\finsetv$ of cardinality at least $2$. These restriction maps are smooth extensions of the natural (quotients of) restriction maps from the space $\check{C}_{\finsetv}\!\left[\RR^d\right]$ of injective maps from $\finsetv$ to $\RR^d$ to the compact space $\cuptd{\finseta}{\RR^d}$. We could define the configuration space  ${C}_{\finsetv}\!\left[\RR^d\right]$ to be the closure of the image of $\check{C}_{\finsetv}\!\left[\RR^d\right]$ in the product $(\RR^d)^{\finsetv} \times \prod_{\finseta \subseteq \finsetv; \cardlef{\finseta} \geq 2}\cuptd{\finseta}{\RR^d}$, or in smaller spaces, as Dev Sinha did in \cite{Sinha}. This would define its topology. However, the differential structures of our configuration spaces are essential for our purposes. This is why we study them in detail in this chapter.

\begin{definition}
\label{defpartitionone}

The \emph{partition associated to a map} $\funcf$ from a finite set $\finsetv$ to some set $X$
is the following set $\kids(\finsetv;\funcf)$ of subsets of $\finsetv$.
\begin{equation*}\kids(\finsetv;\funcf)=\{\funcf^{-1}(x);x \in \funcf({\finsetv})\}.\end{equation*} 
\index[N]{K(\finsetv;\funcf)@ $\kids(\finsetv;\funcf)$ kid set}

In this book, a \emph{partition} of a finite set $\finsetv$ is a set of disjoint nonempty subsets of $\finsetv$ whose union is $\finsetv$.
The elements of a partition $\kids(\finsetv)$ are called the \emph{kids}\index[T]{kid} of ${\finsetv}$ (with respect to the partition). (We do not call them children because we use the initial of children in the notation of configuration spaces.)
The \emph{daughters}\index[T]{daughter} of $\finsetv$ with respect to such a partition are its kids with cardinality at least $2$, and its \emph{sons}\index[T]{son} are the singletons of $\kids(\finsetv)$. (The daughters might bear smaller kids in other partitions.) We denote the set of daughters of $\finsetv$ by $D(\finsetv,\kids(\finsetv))$, 
 \index[N]{Dau@Daughter sets!D(\finsetv,\kids(\finsetv))@ $D(\finsetv,\kids(\finsetv))$} 
or by $D(\finsetv)$ \index[N]{Dau@Daughter sets!D(\finsetv)@ $D(\finsetv)$}
when $\kids(\finsetv)$ is understood. We simply denote $D(\finsetv,\kids(\finsetv;\funcf))$ by $D(\finsetv;\funcf)$.
\end{definition}

\begin{definition}
\label{defparenthesizationone}

A \indexT{parenthesization} $\parentp$ of a finite set $\finsetv$ is a set $\parentp=\{\finseta_i;i \in I\}$ of
 subsets of $\finsetv$, each of cardinality greater than one, such that, for any two distinct elements $i$, $j$ of $I$, one of the following holds $\finseta_i \subset \finseta_j$,  $\finseta_j \subset \finseta_i$ or $\finseta_i \cap \finseta_j=\emptyset$.
 \end{definition}

Every element $\confx$ of the space ${C}_{\finsetv}\!\left[\RR^d\right]$ defines the following parenthesization $\parentp(\confx)$ of $\finsetv$.
The maximal elements (with respect to the inclusion) of $\parentp(\confx)$ are the daughters of $\finsetv$ with respect to $(\pbl(\confx) \in (\RR^d)^{\finsetv})$.
For any element $\finseta$ of $\parentp(\confx)$, the maximal strict subsets of $\finseta$ in $\parentp(\confx)$ are the daughters of $\finseta$ with respect to the restriction of $\confx$ to $\finseta$.
In our examples of Figures~\ref{figmagnifyingglasses} and \ref{figmagnifyingglassestwo}, the parenthesizations are $\{\finsetv,\finseta_{123},\finseta_{23}\}$ and $\{\finsetv,\finseta_{12},\finseta_{34}\}$, respectively. They are in one-to-one correspondences with the magnifying glasses provided by the iterated blow-ups, or, equivalently, with the blow-ups that affected $\confx$.
For a parenthesization $\parentp$ of $\finsetv$,
define the \emph{stratum} $C_{\finsetv,\parentp}\left[\RR^d\right]=\{ \confx\in C_{\finsetv}\!\left[\RR^d\right] \suchthat \parentp(\confx)=\parentp\}$.
As we can see in the above examples, and as we will state in a larger generality in Proposition~\ref{propdescconftautwo}, such a stratum is a smooth manifold of codimension $\cardlef{\parentp}$ in $C_{\finsetv}\!\left[\RR^d\right]$.
As a set, the configuration space ${C}_{\finsetv}\!\left[\RR^d\right]$ is the disjoint union of the strata $C_{\finsetv,\parentp}\left[\RR^d\right]$.
So these strata define a \emph{stratification} of ${C}_{\finsetv}\!\left[\RR^d\right]$. 
The open codimension-one faces of $C_{\finsetv}\!\left[\RR^d\right]$ are the strata of codimension one. They are in one-to-one correspondence with the subsets $\finseta$ of $\finsetv$ of cardinality greater than $1$. The open face associated to $\finseta$ consists of the pairs $(\confc \in \diag_{\finseta}((\RR^d)^{\finsetv}) ,w_{\finseta} \in \cuptd{\finseta}{\RR^d})$ such that the restriction of $\confc$ to $(\finsetv \setminus \finseta) \cup \{a \in \finseta\}$ and the map $w_{\finseta}$ are injective.

The space $\cuptd{\finsetv}{\RR^d}$ of nonconstant maps from $\finsetv$ to $\RR^d$ up to translation and dilation is the preimage of a constant map under the blowdown map from $\blowup{(\RR^d)^{\finsetv}}{\diag_{\finsetv}((\RR^d)^{\finsetv})}$ to $(\RR^d)^{\finsetv}$.  In  Theorem~\ref{thmcompfacanom}, we successively blow up
the diagonals $\diag_{\finseta}(\cuptd{\finsetv}{\RR^d})$ of this space, for strict subsets $\finseta$ of $\finsetv$ of cardinality
greater than $1$, as above, to transform $\cuptd{\finsetv}{\RR^d}$ to the preimage $\ccompuptd{\finsetv}{\RR^d}$ of a constant map under the composed blowdown map from  ${C}_{\finsetv}\!\left[\RR^d\right]$ to $(\RR^d)^{\finsetv}$. 
The space $\ccompuptd{\finsetv}{\RR^d}$ is presented in detail in Section~\ref{secpresfacanom}. It is a compactification of the space $\cinjuptd{\finsetv}{\RR^d}$ of injective maps from $\finsetv$ to $\RR^d$ up to translation and dilation.

The spaces $\cinjuptd{\underline{2}}{\RR^d}$, $\cuptd{\underline{2}}{\RR^d}$, $\ccompuptd{\underline{2}}{\RR^d}$ are identical, they consist of the classes of the maps $w_X \colon \underline{2} \to \RR^d$ such that $w_X(1)=0$ and $w_X(2)=X$ for $X \in S^{d-1}$. So they are diffeomorphic to the unit sphere $S^{d-1}$ of $\RR^d$.
In particular, the space $\cuptd{\underline{2}}{\RR}$ has two elements, which are the classes of $w^+$ and $w^-$, where $w^+$ and $w^-$ map $1$ to $0$, $w^+(2)=1$, and $w^-(2)=-1$. 

Let us discuss the spaces $\cinjuptd{\finsetv}{\RR}$, $\cuptd{\finsetv}{\RR}$, and $\ccompuptd{\finsetv}{\RR}$ in more detail.
\begin{example}
\label{exaStasheff}
In general, for an integer $k\geq 2$, $\cinjuptd{\underline{k}}{\RR}$ and its compactification $\ccompuptd{\underline{k}}{\RR}$ have $k!$ connected components, which correspond to the orders of the $\confc(i)$ in $\RR$, for configurations $\confc$ in $\cinjuptd{\underline{k}}{\RR}$.
Denote the connected component of $\cinjuptd{\underline{k}}{\RR}$ in which the configurations $\confc$ satisfy $\confc(1) <\confc(2) <\dots<\confc(k)$, by $\cinjuptd{<,\underline{k}}{\RR}$. Its respective closures in $\ccompuptd{\underline{k}}{\RR}$ and in $\cuptd{\underline{k}}{\RR}$ are denoted by $\ccompuptd{<,\underline{k}}{\RR}$ and $\cuptd{<,\underline{k}}{\RR}$.
Then we have $\cinjuptd{<,\underline{3}}{\RR}= \{(0,t,1) \suchthat t \in \left]0,1\right[\}$. The spaces $\ccompuptd{<,\underline{3}}{\RR}$ and $\cuptd{<,\underline{3}}{\RR}$  coincide. They are the natural compactification $\left[0,1\right]$ of $\cinjuptd{<,\underline{3}}{\RR}$.
In this space, an element $t$ of $\left]0,1\right[$ represents the injective configuration $(0,t,1)$.
In $\left[0,1\right]$, the extremity $0$ represents the limit configuration $((..).)=\lim_{t \to 0}(0,t,1)$ and $1$ represents the limit configuration $(.(..))=\lim_{t \to 0}(0,1-t,1)$.
The configuration space $\ccompuptd{<,\underline{4}}{\RR}$ is diffeomorphic to the following well-known pentagon:
 \begin{center}
 \begin{tikzpicture} \useasboundingbox (-3,0) rectangle (3,1);
\draw [-] (-1,0) node[left]{ \scriptsize ((.(..)).)} -- (1,0) node[right]{ \scriptsize (.((..).))} -- (1.2,.7) node[right]{ \scriptsize (.(.(..)))} -- (0,1) -- (-1.2,.7) node[left]{ \scriptsize (((..).).)} -- (-1,0) (.1,1.1)node[right]{ \scriptsize ((..)(..))};
\end{tikzpicture}
\end{center}
The edge from $((.(..)).)$ to $(((..).).)$ is the preimage of the diagonal $\diag_{\underline{3}}(\cuptd{<,\finsetv}{\RR})$ under the composed blowdown map from $\ccompuptd{<,\underline{4}}{\RR}$ to $\cuptd{<,\underline{4}}{\RR}$. It is naturally diffeomorphic to $\ccompuptd{<,\underline{3}}{\RR}$.
The edge from $(((..).).)$ to $((..)(..))$ is the closure of the preimage of the diagonal $\diag_{\underline{2}}(\cuptd{<,\finsetv}{\RR})$. Its interior is naturally diffeomorphic to $\cinjuptd{<,\{2,3,4\}}{\RR}$.

In general, for $k \geq 3$, the configuration space $\ccompuptd{<,\underline{k}}{\RR}$ is a \emph{Stasheff polyhedron} \cite{Stasheff63} of dimension $(k-2)$ whose \emph{corners} (i.e., ridges of dimension $0$) are labeled by \emph{nonassociative words} in the letter $\,\bulletmoyen\,$,
as in the above example. For any integer $k\geq 2$, a nonassociative word $w$ with $k$ letters represents a limit configuration $w=\lim_{t \to 0}w(t)$, where $w(t)=(w_1(t)=0,w_2(t), \dots, w_{k-1}(t), w_k(t)=1)$ is an injective configuration inductively defined as follows for $t \in \left]0,\frac12\right[$.
If $w$ is the product $uv$ of a nonassociative word $u$ of length $j\geq 1$ and a nonassociative word $v$ of length $(k-j)\geq 1$, then $w_i(t)=tu_i(t)$ when $1<i \leq j$, and $w_i(t)=1-t + tv_{i-j}(t)$ when $j<i <k$.
For example, we have $(((..).).)(t)=(0,t^2,t,1)$. 
In a limit configuration associated to such a nonassociative word, points inside matching parentheses 
are thought of as infinitely closer to each other than they are to points outside these matching parentheses.
The parenthesization associated as above to a nonassociative word is the set of strict subsets inside matching parentheses.\footnote{In Theorem~\ref{thmcompuptd}, we will rather associate the set of all subsets inside matching parentheses.}
\end{example}

\section{General presentation of \texorpdfstring{${C}_{\finsetv}({\rats})$}{the compactifications}} \label{secfirstprescomp}

Let $\finsetv$ denote a finite set. We use this notation for a generic finite set since our sets will end up being sets of vertices of Jacobi diagrams.
The space of maps from $\finsetv$ to $\setx$ is denoted by $\setx^{\finsetv}$ as usual. For a subset $\finseta$ of $\finsetv$ of cardinality at least $2$, recall that 
the subspace of $\setx^{\finsetv}$ consisting of maps $\confc$ that map $\finseta$ to a single element and such that $\confc(\finsetv \setminus \finseta) \subset \setx \setminus \confc(\finseta)$ is a diagonal denoted by $\diag_{\finseta}(\setx^{\finsetv})$. \index[N]{Dzelta@Diagonals!Dzeltab@$\diag_{\finseta}(\setx^{\finsetv})$} In particular, if $\cardlef{\finsetv} \geq 2$, the \indexT{small diagonal} consisting of constant maps is denoted by $\diag_{\finsetv}(\setx^{\finsetv})$.

As often in this book, we fix a rational homology sphere $\rats$ and a point $\infty$ of $\rats$.
Recall that $\check{C}_{\finsetv}({\rats})$ is the space of injective maps from ${\finsetv}$ to $\left(\crats=\rats \setminus \{\infty\}\right)$.

\begin{theorem}
 \label{thmcompconf}
Define a compactification ${C}_{\finsetv}({\rats})$ of $\check{C}_{\finsetv}({\rats})$ as follows.
For a nonempty ${\finseta} \subseteq {\finsetv}$, let $\einftyxi_{\finseta}$ be the set of maps from $\finsetv$ to $\rats$ that map ${\finseta}$ to $\infty$ and $\finsetv \setminus \finseta$ to $\rats \setminus \{\infty\}$. 
Start with $\rats^{\finsetv}$.
Blow up $\einftyxi_{\finsetv}$ (which is reduced to the point $m=\infty^{\finsetv}$ such that $m^{-1}(\infty)={\finsetv}$). Set \begin{equation*}{C}_{\finsetv,\card{\finsetv}+1}({\rats})= \bigblowup{\rats^{\finsetv}}{\einftyxi_{\finsetv}}.\end{equation*}
For $k=\cardlef{\finsetv},\cardlef{\finsetv}-1, \dots, 3, 2$, let 
${C}_{\finsetv,k}({\rats})$ be obtained from ${C}_{\finsetv,k+1}({\rats})$ by blowing up
the closures of (the preimages under the composition of the previous blowdown maps of) the $\diag_{{\finseta}}(\crats^{\finsetv})$ such that $\cardlef{\finseta}=k$ and the closures of (the preimages under the composition of the previous blowdown maps of) the $\einftyxi_J$ such that $\cardlef{J}=k-1$. At each step, the blown-up  manifolds are smooth and transverse to the ridges, and ${C}_{\finsetv,k}({\rats})$ is independent of the order of the blow-ups.
The obtained manifold ${C}_{\finsetv}({\rats})= {C}_{\finsetv,2}({\rats})$ is a smooth compact $(3\cardlef{\finsetv})$-manifold with ridges. 
The interior of ${C}_{\finsetv}({\rats})$ is $\check{C}_{\finsetv}({\rats})$, and the composition of the blowdown maps gives rise to a canonical smooth blowdown projection
$\pbl \colon {C}_{\finsetv}({\rats}) \to \rats^{\finsetv}$.
\end{theorem}
We prove the generalization Theorem~\ref{thmcompconfbis} of Theorem~\ref{thmcompconf}
in Section~\ref{secbloinf} using the results of Sections~\ref{secprooffacanom} and \ref{secblodiag}. This generalization includes an alternative definition of ${C}_{\finsetv}({\rats})$.

Set $C_n(\rats)=C_{\underline{n}}(\rats)$.

As already announced, with the above definition, we have $C_1(\rats)=\blowup{\rats}{\infty}$, and the configuration space $C_2(\rats)$ is the compactification studied in Section~\ref{secCtwo}. In particular, Theorem~\ref{thmcompconf} is true when $\cardlef{V} \leq 2$.

\begin{theorem}
\label{thmcompconftwo}
The configuration spaces $\check{C}_{\finsetv}({\rats})$ have the following properties.
\begin{enumerate}
 \item 
Under the assumptions of Theorem~\ref{thmcompconf},
for $\finseta \subset {\finsetv}$,
the restriction map \begin{equation*}p_{\finseta} \colon \check{C}_{\finsetv}({\rats}) \rightarrow \check{C}_{\finseta}({\rats})\end{equation*} extends to a smooth restriction map still denoted by $p_{\finseta}$ from ${C}_{\finsetv}({\rats})$
to ${C}_{\finseta}({\rats})$ such that the following square commutes:
  \begin{equation*}\xymatrix{
    {C}_{\finsetv}({\rats}) \ar[r]^{p_{\finseta}} \ar[d]_{\pbl} & {C}_{\finseta}({\rats}) \ar[d]^{\pbl} \\
    \rats^{\finsetv} \ar[r]^{p_{\finseta}} & \rats^{\finseta}
  }\end{equation*}
  \item
For an open subset $U$ of $\rats$, let ${C}_{\finsetv}(U)$ denote $\pbl^{-1}(U^{\finsetv})$.
If $\finsetv = \sqcup_{i \in I}\finseta_i$ and if $(U_i)_{i \in I}$ is a family of disjoint open sets of $\rats$, then the product
\begin{equation*}\pbl^{-1}\Bigl(\prod_{i \in I}U_i^{\finseta_i}\Bigr) \hfl{\prod_{i \in I}p_{\finseta_i} }\prod_{i \in I}{C}_{\finseta_i}(U_i)\end{equation*}
 of the above restriction maps is a diffeomorphism.\end{enumerate}
\end{theorem}

The first part of this theorem will be a direct consequence of Proposition~\ref{proprest} and Theorem~\ref{thmcompconfbis}. Its second part comes from the locality of the blow-up operations. 
The spaces ${C}_{\finsetv}({\crats})$, which involve only blow-ups along the diagonals, have been
studied by Scott Axelrod and Isadore Singer \cite[Section 5]{axelsingII}, and with more details by Dev Sinha \cite{Sinha}. Their properties that are useful in this book are proved in Sections~\ref{secprooffacanom} and \ref{secblodiag}. William Fulton and Robert MacPherson studied similar compactifications of configuration spaces in an algebraic geometry setting in \cite{FultonMcP}.

Recall the configuration space
\begin{equation*}\check{C}(\rats,\Link;\Gamma)=\Bigl\{c \colon U \cup T \hookrightarrow \crats \suchthat
\bigl(\exists j \in \left[i_{\Gamma}\right] \suchthat c\vert_{U}=\Link \circ j\bigr)\Bigr\}\end{equation*}
of Section~\ref{secdefconfspace}.

\begin{proposition}
\label{propcompconfL}
The closure of $\check{C}(\rats,\Link;\Gamma)$ in $C_{V(\Gamma)}(\rats)$ is a smooth compact submanifold of $C_{V(\Gamma)}(\rats)$ transverse to the ridges. It is denoted by ${C}(\rats,\Link;\Gamma)$.
\end{proposition}
We will prove Proposition~\ref{propcompconfL} in Section~\ref{secproofblowup}.
Theorems~\ref{thmcompconf} and~\ref{thmcompconftwo}, and Proposition~\ref{propcompconfL} imply Proposition~\ref{propcompext}.

\section{Configuration spaces associated to unit normal bundles to diagonals}
\label{secpresfacanom}
For a vector space $\vecspt$, recall from Section~\ref{secblowup} that $\sph(\vecspt)$ denotes the quotient $\sph(\vecspt)=(\vecspt \setminus \{0\})/{\RR^{+\ast}}$ of $\vecspt \setminus \{0\}$ by the dilations. If $\vecspt$ is equipped with a Euclidean norm, then $\sunits(\vecspt)$ \index[N]{Sbf@$\sunits$ unit sphere} denotes the unit sphere of $\vecspt$. In this case, $\sph(\vecspt)$ and $\sunits(\vecspt)$ are diffeomorphic.

\begin{definition}
\label{defcuptdfinsetvvecspt}
Let $\finsetv$ denote a finite set of cardinality at least $2$.
We use the notation $\vecspt$ for a generic vector space since $\vecspt$ will end up being a tangent space.
Let
\begin{equation*}\cuptd{\finsetv}{\vecspt}= \sph\bigl(\vecspt^{\finsetv}/\diag_{\finsetv}(\vecspt^{\finsetv})\bigr)=\frac{\left(\vecspt^{\finsetv}/\diag_{\finsetv}(\vecspt^{\finsetv})\right) \setminus \{0\}}{\RR^{+\ast}}\end{equation*}
be the space of nonconstant maps from ${\finsetv}$ to ${\vecspt}$ up to translation and dilation.\index[N]{Sebar@$\cuptd{\finsetv}{\vecspt}$ configuration space}
\end{definition} 

\begin{lemma}
\label{lemnormfibdiagzero}
Let $\finseta$ be a subset of $\finsetv$.
 The fiber of the unit normal bundle to $\diag_{\finseta}(\rats^{\finsetv})$ 
 in $\rats^{\finsetv}$ over a configuration $m$ is
$\cuptd{\finseta}{T_{m({\finseta})}\rats}$.
\end{lemma}
\bp Exercise. \eop

As expected for the fiber of a unit normal bundle, we have the following lemma.

\begin{lemma}
\label{lemnormfibdiag}
The topological space $\cuptd{\finsetv}{\vecspt}$ has a canonical smooth structure. When $\vecspt$ is a vector space of dimension $\dimdel$, the space $\cuptd{\finsetv}{\vecspt}$
is diffeomorphic to a sphere of dimension $\left((\cardlef{\finsetv}-1)\dimdel -1\right)$. 
\end{lemma}
\bp Choosing a basepoint $b(\finsetv)$ for $\finsetv$ and a basis for $\vecspt$ identifies  
$\cuptd{\finsetv}{\vecspt}$ with the set $\sunits\left(\vecspt^{\finsetv \setminus \{b(\finsetv)\}}\right)$ of maps $\confc \colon \finsetv \to \vecspt$ such that
\begin{itemize}
\item $\confc\bigl(b(\finsetv)\bigr)=0$ and
\item $\sum_{\eltv \in \finsetv} \norm{\confc(\eltv)}^2=1$,
\end{itemize}
with respect to the Euclidean norm for which our basis is orthonormal.
It is easy to see that changes of basepoints, or bases of $\vecspt$ give rise to diffeomorphisms of spheres.
\eop

\begin{convention}
 \label{convordefcuptdfinsetvvecspt}
In this chapter, we do not take care of orientations. Later, we will associate the following orientation of $\cuptd{\finsetv}{\vecspt}$ to an order of $\finsetv$ and an orientation of $\vecspt$.
Assume that $\finsetv=\{v_1,\dots,v_n\}$ and that $\vecspt$ is oriented. The order on $\finsetv$ orients $\vecspt^{\finsetv}$ and $\vecspt^{\finsetv \setminus\{v_1\}}$. Then the map from the boundary of the unit ball of $\vecspt^{\finsetv \setminus\{v_1\}}$ to $\cuptd{\finsetv}{\vecspt}$ that maps an element $(x_2,\dots,x_n)$ of $\sunits\left(\vecspt^{\finsetv \setminus\{v_1\}}\right)$ to the class of $(0,x_2,\dots,x_n)$ is an orientation-preserving diffeomorphism. \end{convention}

When $\finsetv=\underline{2}$, fixing the basepoint $b(\underline{2})$ of $\underline{2}$ to be $1$ identifies $\cuptd{\finsetv}{\vecspt}$ with the unit sphere $\sunits(\vecspt)$ of $\vecspt$, if $\vecspt$ is equipped with a Euclidean norm, as in the proof of Lemma~\ref{lemnormfibdiag}.

For a strict subset $\finseta$ of $\finsetv$ of cardinality at least $2$, define the diagonal
$\diag_{\finseta}\bigl(\cuptd{\finsetv}{\vecspt}\bigr)$ to be the subset of $\cuptd{\finsetv}{\vecspt}$ consisting of classes of maps $\confc$ from $\finsetv$ to $\vecspt$ that are constant on $\finseta$ and such that $\confc(\finsetv \setminus \finseta) \cap\confc(\finseta) =\emptyset$.
Let
$\cinjuptd{\finsetv}{\vecspt}$\index[N]{Scheck@$\cinjuptd{\finsetv}{\vecspt}$ open configuration space} denote the subspace of $\cuptd{\finsetv}{\vecspt}$ consisting of \emph{injective} maps from ${\finsetv}$ to ${\vecspt}$ up to translation and dilation.
The following theorem defines a bigger compactification of $\cinjuptd{\finsetv}{\vecspt}$, which is also used in our study of the variations of $\Zinv$, and especially in the definition of the anomalies.
We will prove the following two theorems in Section~\ref{secprooffacanom}. See Theorem~\ref{thmcompuptd} and Proposition~\ref{proprestanom}.

\begin{theorem}
\label{thmcompfacanom} Set $\ccompuptdk{\finsetv}{\vecspt}{\card{\finsetv}}=\cuptd{\finsetv}{\vecspt}$.
For $k=\cardlef{\finsetv}-1, \dots, 3, 2$, define $\ccompuptdk{\finsetv}{\vecspt}{k}$ from $\ccompuptdk{\finsetv}{\vecspt}{k+1}$ by blowing up 
the closures of (the preimages under the composition of the previous blowdown maps of) the $\diag_{{\finseta}}(\cuptd{\finsetv}{\vecspt})$ such that $\cardlef{\finseta}=k$. At each step, the  blown-up manifolds are smooth and transverse to the ridges, and
the resulting blown-up manifold is independent of the order choice.

Thus, this process gives rise to a canonical compact smooth manifold 
$\ccompuptd{\finsetv}{\vecspt}=\ccompuptdk{\finsetv}{\vecspt}{2}$ with ridges. The interior of $\ccompuptd{\finsetv}{\vecspt}$ is
$\cinjuptd{\finsetv}{\vecspt}$. 
\end{theorem}

When $\finsetv$ has two elements, we have $\ccompuptd{\finsetv}{\vecspt}=\cinjuptd{\finsetv}{\vecspt}=\cuptd{\finsetv}{\vecspt}$. In general,
the manifold $\ccompuptd{\finsetv}{\vecspt}$ satisfies the following properties.
\begin{theorem}
\label{thmcompfacanomtwo}
With the notation of Theorem~\ref{thmcompfacanom}, for any subset $\finseta$ of $\finsetv$, the restriction from $\cinjuptd{\finsetv}{\vecspt}$ to 
$\cinjuptd{\finseta}{\vecspt}$ extends to a smooth map from $\ccompuptd{\finsetv}{\vecspt}$ to 
$\ccompuptd{\finseta}{\vecspt}$.
\end{theorem}

The open codimension-one faces of $\ccompuptd{\finsetv}{\vecspt}$ will be the loci for which only one blow-up along some $\diag_{{\finseta}}\bigl(\cuptd{\finsetv}{\vecspt}\bigr)$ is involved, for a strict subset $\finseta$ of $\finsetv$ such that $\cardlef{\finseta} \geq 2$. The blowdown map sends such a face $f(\finseta)(\vecspt)$ into $\diag_{{\finseta}}(\cuptd{\finsetv}{\vecspt})$. 
As it will be seen in Lemma~\ref{lemnormfibdiagtwo}, the fiber of the unit normal bundle
of $\diag_{\finseta}\bigl(\cuptd{\finsetv}{\vecspt}\bigr)$ is $\cuptd{\finseta}{\vecspt}$.
Thus, the following proposition will be clear in the end of Section~\ref{secprooffacanom}, where a one-line proof is given.

\begin{proposition}
\label{propfacfacanom}
The codimension-one faces of $\ccompuptd{\finsetv}{\vecspt}$ are in one-to-one correspondence with the strict subsets $\finseta$ of $\finsetv$ with cardinality at least $2$. The (open) face $f(\finseta)(\vecspt)$ corresponding to such an $\finseta$ is
\begin{equation*}f(\finseta)(\vecspt)=\cinjuptd{\finseta}{\vecspt} \times \cinjuptd{\{\elta\} \cup (\finsetv \setminus {\finseta})}{\vecspt}\end{equation*} for an element $\elta$ of $\finseta$.
For a subset $e$ of cardinality $2$ of $\finsetv$, the restriction to $f(\finseta)(\vecspt)$ of the extended restriction
\begin{equation*}p_e \colon \ccompuptd{\finsetv}{\vecspt} \rightarrow \ccompuptd{e}{\vecspt}\end{equation*}
may be described as follows:
\begin{itemize}
\item If $e \subseteq {\finseta}$, then $p_e$ is the composition of the natural
projections
\begin{equation*}f({\finseta})(\vecspt) \longrightarrow \cinjuptd{\finseta}{\vecspt} \longrightarrow S_e(\vecspt).\end{equation*}
\item If $e \subseteq (\finsetv \setminus {\finseta}) \cup \{\elta\}$, then $p_e$ is the composition of the natural
projections
\begin{equation*}f({\finseta})(\vecspt) \longrightarrow \cinjuptd{\{\elta\} \cup(\finsetv \setminus {\finseta})}{\vecspt} \longrightarrow  S_e(\vecspt).\end{equation*}
\item If $e \cap {\finseta} =\{\elta^{\prime}\}$, let $\tilde{e}$ be obtained from $e$ by replacing  $\elta^{\prime}$ by $\elta$, then we have $p_e=p_{\tilde{e}}$.
\end{itemize}
\end{proposition}

The space $\ccompuptd{\finsetv}{T_x \crats}$ described in the above theorem is related to 
${C}_{\finsetv}(\rats)$ by the following proposition, which is a corollary of Theorem~\ref{thmcompdiaggen}.

\begin{proposition}
\label{propfiberc}
For $x \in \crats$,  for $\pbl \colon {C}_{\finsetv}(\rats) \to \rats^{\finsetv}$, we have
\begin{equation*}\pbl^{-1}(x^{\finsetv})=\ccompuptd{\finsetv}{T_x \crats}.\end{equation*}
\end{proposition}

\section{The codimension-one faces of \texorpdfstring{${C}(\rats,\Link;\Gamma)$}{C(R,L;Gamma)}}
\label{secfacecodimone}

Recall that our terminology in Section~\ref{secbackgroundmfdsbry} makes codimension-one faces open in the boundary of a manifold with ridges. Our codimension-one faces may be called \emph{open codimension-one faces} in other references.
We describe the codimension-one faces of the configuration space ${C}(\rats,\Link;\Gamma)$ for a link $\Link$ in a $\QQ$-sphere $\rats$ below with an outline of justification. Details are provided in Sections~\ref{secprooffacanom} to \ref{secproofblowup}.  Again, the codimension-one faces are the loci where only one blow-up has been performed.

In this book, the sign $\subset$ stands for \say{is a strict subset of} or \say{$\subseteq$ and $\neq$}.

\paragraph{Faces that involve blow-ups along diagonals \texorpdfstring{$\diag_{\finseta}$}{Delta}.}

Let $\finseta$ be a subset of the set $V(\Gamma)$ of vertices of a Jacobi diagram $\Gamma$ on the domain $\source$ of a link $\Link$. Assume $\cardlef{\finseta} \geq 2$. Let us describe the (open) face $\facee(\finseta,\Link,\Gamma)$, which comes from the blow-up along $\diag_{\finseta}(\crats^{V(\Gamma)})$. \index[N]{Faces!FALGamma@$\facee(\finseta,\Link,\Gamma)$}
Such a face contains
limit configurations that map $\finseta$ to a point of $\crats$.
An element of such a face is described by its image $m$ under the blowdown map
\begin{equation*}\pbl \colon {C}_{V(\Gamma)}(\rats) \to \rats^{V(\Gamma)},\end{equation*}
which maps $\facee(\finseta,\Link,\Gamma)$ to $\diag_{\finseta}(\crats^{V(\Gamma)})$,
together with an element of the fiber 
\begin{equation*}\cuptd{\finseta}{T_{m({\finseta})}\crats}.\end{equation*}
Let $\elta \in \finseta$.
Let $\check{\diag}_{\finseta}(\crats^{V(\Gamma)})$ denote the set of maps of $\diag_{\finseta}(\crats^{V(\Gamma)})$ whose restriction to $\{\elta\} \cup \bigl(V(\Gamma) \setminus \finseta\bigr)$
is injective.
Set
\begin{equation*}\baseb(\finseta,\Link,\Gamma) = \check{\diag}_{\finseta}\bigl(\crats^{V(\Gamma)}\bigr) \cap \pbl({C}(\rats,\Link;\Gamma)).\end{equation*}
The face $\facee(\finseta,\Link,\Gamma)$ fibers over the subspace $\baseb(\finseta,\Link,\Gamma)$.
When $\finseta$ contains no univalent vertex, the fiber over a point $m$ is $\cinjuptd{\finseta}{T_{m({\finseta})}\crats}$.

Let $i_{\Gamma}$ be a $\Gamma$-compatible injection. Let $\source_1$ be an oriented 
 connected component of $\source$. Let $U_1=i_{\Gamma}^{-1}(\source_1)$.
The restriction to $U_1$ of the injection $i_{\Gamma}$ into $\source_1$ induces a permutation $\sigma$ of $\finsetu_1$, such that, travelling along $\source_1$, we meet 
$i_{\Gamma}(\eltv)$, $i_{\Gamma}\bigl(\sigma(\eltv)\bigr)$, \dots, $i_{\Gamma}\bigl(\sigma^{\cardlef{\finsetu_1}}(\eltv)=\eltv\bigr)$, successively, for any element $\eltv$ of $\finsetu_1$.
A \emph{set of consecutive\index[T]{consecutive elements} elements} of $\finsetu_1$, with respect to $i_{\Gamma}$, is a subset $\finseta_U$ of $\finsetu_1$ that may be written as
$\{\eltv, \sigma(\eltv), \dots, \sigma^{k}(\eltv)\}$ for some element $\eltv \in \finsetu_1$ and for $k \leq \cardlef{\finsetu_1} -1$.
If $\finseta_U \neq \finsetu_1$, then the first element $\eltv$ in such an $\finseta_U$ is unique, and $\sigma$ induces the following
unique linear order
 \begin{equation*} \eltv < \sigma(\eltv)< \dots< \sigma^{k}(\eltv)\end{equation*}
on such a set $\finseta_U$ of consecutive elements $\finsetu_1$, which is said to be \emph{compatible} with the isotopy class $\left[i_{\Gamma}\right]$ of $i_{\Gamma}$.
If $\finseta_U = \finsetu_1$, every choice of an element $\eltv$ in $\finseta_U$ induces a \emph{linear} (i.e., total) order 
 \begin{equation*} \eltv < \sigma(\eltv)< \dots< \sigma^{k}(\eltv),\end{equation*}
which is said to be \emph{compatible} with $\left[i_{\Gamma}\right]$.

Let $\finseta_U =\finseta \cap U(\Gamma)$.
When $\finseta$ contains univalent vertices, if $\facee(\finseta,\Link,\Gamma)$ is nonempty, then $\finseta_U$ must be a set of consecutive vertices on a component $\source_1$ of $\source$ with respect to the given class $\left[i_{\Gamma}\right]$ of injections of $U(\Gamma)$ into $\source_1$.
The fiber over a point $m$ is the subset $\cinjuptd{\finseta}{T_{m({\finseta})}\crats,\Link,\Gamma}$ of $\cinjuptd{\finseta}{T_{m({\finseta})}\crats}$ consisting of injections that map $\finseta_U$ on a line directed by $T_{m(\finseta)}\Link$,
so that the order induced by the line on $\finseta_U$ coincides with 
\begin{itemize}
\item  the linear order induced by $\left[i_{\Gamma}\right]$, if $\finseta_U$ is not the whole $i_{\Gamma}^{-1}(\source_1)$,
\item one of the $\cardlef{i_{\Gamma}^{-1}(\source_1)}$ linear orders compatible with the cyclic order induced by $\left[i_{\Gamma}\right]$, if $\finseta_U=i_{\Gamma}^{-1}(\source_1)$.
\end{itemize}
In this latter case, neither the fiber nor $\facee(\finseta,\Link,\Gamma)$ is connected. Their connected components are in one-to-one correspondence with the compatible orders.

According to Theorem~\ref{thmcompconftwo}, for any pair $e$ of $V(\Gamma)$, there exists a smooth restriction map from ${C}_{V(\Gamma)}$ to ${C}_e(\rats)$. An order on $e$ identifies ${C}_e(\rats)$ with ${C}_2(\rats)$. We describe the natural restriction $p_e$ to the (open) face $\facee(\finseta,\Link,\Gamma)$
below for a pair $e$ of $V(\Gamma)$.

\begin{itemize}
\item If $\cardlef{e \cap \finseta} \leq 1$, then $p_e$ is the composition of the natural
projections
\begin{equation*}\facee(\finseta,\Link,\Gamma) \longrightarrow \check{\diag}_{\finseta}(\crats^{V(\Gamma)}) \longrightarrow  \check{C}_e(\rats).\end{equation*}
\item If $e \subseteq \finseta$, then $p_e$ maps an element of $\cinjuptd{\finseta}{T_{m({\finseta})}\crats,\Link,\Gamma}$ to its projection
in $ \cinjuptd{e}{T_{m({\finseta})}\crats} \subset C_e(\rats)$.
\end{itemize}

\paragraph{Faces \texorpdfstring{$\facee(\finsetv(\Gamma),\Link,\Gamma)$}{F(V(Gamma),L,Gamma)}.}

Let us study the following special case of the previous paragraph's faces, which will play a particular role.
Consider a face $\facee(\finsetv(\Gamma),\Link,\Gamma)$ where a connected graph $\Gamma$ with at least one univalent vertex collapses and $\Link$ is a knot embedding. 
Such a face has one connected component for each linear order of $U(\Gamma)$ compatible with the cyclic order of $U(\Gamma)$.
\begin{notation}
\label{notorderlift}
A Jacobi diagram $\check{\Gamma}$ on $\RR$ yields a diagram $\cl(\check{\Gamma})$ on $S^1$, where $S^1$ is viewed as $\RR \cup \{\infty\}$ by adding $\infty$ to $\RR$. The natural orientation of $\RR$ orders $U(\check{\Gamma})$.
For a Jacobi diagram ${\Gamma}$ on $S^1$, represent a linear order of $U(\Gamma)$ compatible with the cyclic order of $U(\Gamma)$ by a Jacobi diagram $\check{\Gamma}$ on $\RR$ such that $\cl(\check{\Gamma})=\Gamma$ and the order of $U(\Gamma)$ matches the order of $U(\check{\Gamma})$.
Denote the corresponding connected component of $\facee(\finsetv(\Gamma),\Link,\Gamma)$ by
$\facee(\finsetv(\Gamma),\Link,\check{\Gamma})$.
\end{notation}

Such a connected component fibers over the domain $S^1$ of the knot. We denote the fiber over $z\in S^1$ by $\cinjuptdanvec(T_{\Link(z)}\crats,\vec{t}_{\Link(z)};\check{\Gamma})$,
where $\vec{t}_{\Link(z)}$ is an oriented tangent vector to $\Link$ at $\Link(z)$.

The space $\cinjuptdanvec(T_{\Link(z)}\crats,\vec{t}_{\Link(z)};\check{\Gamma})$ is the space of injections of $V(\check{\Gamma})$ into the vector space $T_{\Link(z)}\crats$ that map  $U(\check{\Gamma})$ to the oriented line $\RR\vec{t}_{\Link(z)}$ directed by $\vec{t}_{\Link(z)}$ with respect to the linear order of $U(\check{\Gamma})$, up to dilation and translation with respect to a vector in $\RR\vec{t}_{\Link(z)}$.
It is naturally a subspace of $\cinjuptd{\finsetv(\Gamma)}{T_{\Link(z)}\crats}$.
Indeed, moding out by all translations is equivalent to considering only configurations that map a univalent vertex to $\RR\vec{t}_{\Link(z)}$ and moding them out by translations along $\RR\vec{t}_{\Link(z)}$.

\begin{lemma}
\label{lemfacfacone}
The closure $\overline{\facee}(\finsetv(\Gamma),\Link,\check{\Gamma})$ of $\facee(\finsetv(\Gamma),\Link,\check{\Gamma})$ in ${C}(\rats,\Link;\Gamma)$ is a manifold with ridges.
The closure $\ccompuptdanvec(T_{\Link(z)}\crats,\vec{t}_{\Link(z)};\check{\Gamma})$ of $\cinjuptdanvec(T_{\Link(z)}\crats,\vec{t}_{\Link(z)};\check{\Gamma})$ in ${C}(\rats,\Link;\Gamma)$ is canonically diffeomorphic to its closure in $\ccompuptd{\finsetv(\Gamma)}{T_{\Link(z)}\crats}$.
It is a smooth manifold with ridges.
\end{lemma}
\bp The first assertion comes from Proposition~\ref{propcompconfL} and from the fact that the closed faces of manifolds with ridges are manifolds with ridges (or from the proof of Lemma~\ref{lemfacfacan} at the end of Section~\ref{secprooffacanom}).
The space $\ccompuptdanvec(T_{\Link(z)}\crats,\vec{t}_{\Link(z)};\check{\Gamma})$ is the closure of
$\cinjuptdanvec(T_{\Link(z)}\crats,\vec{t}_{\Link(z)};\check{\Gamma})$ 
in the fiber over $\Link(z)^{\finsetv(\Gamma)}$ of $C_{\finsetv(\Gamma)}(\crats)$, which is $\ccompuptd{\finsetv(\Gamma)}{T_{\Link(z)}\crats}$ according to Proposition~\ref{propfiberc}.

Now $\overline{\facee}(\finsetv(\Gamma),\Link,\check{\Gamma})$ fibers over $S^1$, and the fiber over $S^1$ is also a manifold with ridges.
\eop

Let $\finseta$ be a strict subset of $\finsetv(\check{\Gamma})$ with cardinality at least $2$ whose univalent vertices are consecutive on $\RR$.
Let $\elta \in \finseta$. Let $\check{\diag}_{\finseta}(\ccompuptd{\finsetv(\Gamma)}{T_{\Link(z)}\crats})$ denote the set of (classes of) maps of $\ccompuptd{\finsetv(\Gamma)}{T_{\Link(z)}\crats}$ whose restriction to $\finseta$ is constant and whose restriction to $\{\elta\} \cup (V(\Gamma) \setminus \finseta)$
is injective.
Set \begin{equation*}\baseb(\finseta,\vec{t}_{\Link(z)};\check{\Gamma})= \ccompuptdanvec(T_{\Link(z)}\crats,\vec{t}_{\Link(z)};\check{\Gamma}) \cap \check{\diag}_{\finseta}\bigl(\ccompuptd{\finsetv(\Gamma)}{T_{\Link(z)}\crats}\bigr).\end{equation*}

Define the (open) face $f(\finseta,\vec{t}_{\Link(z)};\check{\Gamma})$ of $\ccompuptdanvec(T_{\Link(z)}\crats,\vec{t}_{\Link(z)};\check{\Gamma})$ to be the space that fibers over the subspace $\baseb(\finseta,\vec{t}_{\Link(z)};\check{\Gamma})$, whose fiber is the space of injections from $\finseta$ into $T_{\Link(z)}\crats$ that map the univalent vertices of $\finseta$ to the oriented line $\RR\vec{t}_{\Link(z)}$ with respect to the linear order of $U(\check{\Gamma}) \cap \finseta$, up to dilation and translation by a vector in $\RR\vec{t}_{\Link(z)}$.

We will prove the following lemma at the end of Section~\ref{secprooffacanom}.

\begin{lemma}
\label{lemfacfacan}
The codimension-one faces of $\ccompuptdanvec(T_{\Link(z)}\crats,\vec{t}_{\Link(z)};\check{\Gamma})$ are the faces $f(\finseta,\vec{t}_{\Link(z)};\check{\Gamma})$
for the strict subsets $\finseta$ of $\finsetv(\check{\Gamma})$ with cardinality at least $2$ whose univalent vertices are consecutive on $\RR$.
The faces $f(\finseta,\vec{t}_{\Link(z)};\check{\Gamma})$ are the intersections of $\ccompuptdanvec(T_{\Link(z)}\crats,\vec{t}_{\Link(z)};\check{\Gamma})$ with the codimension-one faces $f(\finseta)(T_{\Link(z)}\crats)$ of $\ccompuptd{\finsetv(\check{\Gamma})}{T_{\Link(z)}\crats}$ listed in Proposition~\ref{propfacfacanom}.
In particular, the restriction maps $p_e$ from $f(\finseta)(T_{\Link(z)}\crats)$ to $\ccompuptd{e}{T_{\Link(z)}\crats}$ of Proposition~\ref{propfacfacanom} restrict as restriction maps still denoted by $p_e$ from $f(\finseta,\vec{t}_{\Link(z)};\check{\Gamma})$ to $\ccompuptd{e}{T_{\Link(z)}\crats}$.
\end{lemma}

\paragraph{Faces that involve $\infty$.}

Let $\finseta$ be a nonempty subset of the set $V(\Gamma)$ of vertices of a Jacobi diagram $\Gamma$. Let $\elta \in \finseta$.
 Let us describe the (open) face $\facee_{\infty}(\finseta,\Link,\Gamma)$ \index[N]{Faces!Finfty@$\facee_{\infty}(\finseta,\Link,\Gamma)$} that comes from the blow-up along $\einftyxi_{\finseta}$. It contains
limit configurations, which map $\finseta$ to $\infty$. If it is nonempty, then $\finseta$ contains no univalent vertices.

Let $\check{\einftyxi}_{\finseta}$ denote the set of maps of $\einftyxi_{\finseta}$
that restrict to injective maps on $\{\elta\} \cup (V(\Gamma) \setminus \finseta)$,
and set
\begin{equation*}\baseb_{\infty}(\finseta,\Link,\Gamma) = \check{\einftyxi}_{\finseta} \cap \pbl\bigl({C}(\rats,\Link;\Gamma)\bigr).\end{equation*}
The face $\facee_{\infty}(\finseta,\Link,\Gamma)$ fibers over the subspace $\baseb_{\infty}(\finseta,\Link,\Gamma)$.
\begin{notation}
\label{notsinjupdtcs}
 Let $\sinjupdtcs(T_{\infty}\rats,{\finseta})$ 
\index[N]{SiTinfty@$\sinjupdtcs(T_{\infty}\rats,{\finseta})$ open configuration space}
denote the set of injective maps from $\finseta$ to $(T_{\infty}\rats \setminus 0)$ up to dilation.
\end{notation}

 Note that $\sinjupdtcs(T_{\infty}\rats,{\finseta})$ is an open submanifold of the unit normal bundle  of $\einftyxi_{\finseta}$, which is nothing but $\sph((T_{\infty}\rats)^{\finseta})$.
We have 
\begin{equation*}\facee_{\infty}(\finseta,\Link,\Gamma) = \baseb_{\infty}(\finseta,\Link,\Gamma) \times \sinjupdtcs(T_{\infty}\rats,{\finseta}).\end{equation*} 
For a pair $e$ of $V(\Gamma)$, the natural restriction to $\facee_{\infty}(\finseta,\Link,\Gamma)$ of \begin{equation*}p_e \colon  {C}_{V(\Gamma)}
\to C_e(\rats)\end{equation*} behaves in the following way.
\begin{itemize}
\item If $e \subseteq V(\Gamma) \setminus \finseta$, then $p_e$ is the composition of the natural
maps
\begin{equation*}\facee_{\infty}(\finseta,\Link,\Gamma) \rightarrow \baseb_{\infty}(\finseta,\Link,\Gamma) \rightarrow \breve{C}_{V(\Gamma) \setminus \finseta}(\rats) \rightarrow C_e(\rats).\end{equation*}
\item If $e \subseteq \finseta$, then $p_e$ is the composition of the natural
maps
\begin{equation*}\facee_{\infty}(\finseta,\Link,\Gamma) \longrightarrow \sinjupdtcs(T_{\infty}\rats,{\finseta}) \longrightarrow  \sinjupdtcs(T_{\infty}\rats,e)  \hookrightarrow C_e(\rats).\end{equation*}
\item If $e \cap \finseta =\{\elta^{\prime}\}$, then $p_e$ is the composition of the natural
maps
\begin{multline*}\facee_{\infty}(\finseta,\Link,\Gamma) \longrightarrow \breve{C}_{e \setminus \{\elta^{\prime}\}}(\rats) \times  \sinjupdtcs(T_{\infty}\rats,{\{\elta^{\prime}\}}) \longrightarrow \\ \longrightarrow  \crats^{e \setminus \{\elta^{\prime}\}} \times  S(T_{\infty}\rats^{\{\elta^{\prime}\}}) \hookrightarrow C_e(\rats).\end{multline*}
\end{itemize}

\paragraph{Summary.}

We have just outlined a proof of the following proposition.

\begin{proposition}
\label{propconffaceun}
The disjoint union of
\begin{itemize}
 \item the faces $\facee_{\infty}(\finseta,\Link,\Gamma)$ for nonempty subsets $\finseta$ of $T(\Gamma)$ and
\item the faces $\facee(\finseta,\Link,\Gamma)$ for subsets $\finseta$ of $V(\Gamma)$ with cardinality at least $2$ such that $\finseta \cap U(\Gamma)$ is a (possibly empty) set of consecutive vertices on a connected component of $\source$,
\end{itemize}
described above, embeds canonically in ${C}(\rats,\Link;\Gamma)$.
Its image is the open codimension-one boundary
$\partial_1{C}(\rats,\Link;\Gamma) \setminus \partial_2{C}(\rats,\Link;\Gamma)$
of ${C}(\rats,\Link;\Gamma)$.
Furthermore, for any ordered pair $e$ of $V(\Gamma)$, the restriction to this codimension-one boundary of the map
\begin{equation*}p_e \colon {C}(\rats,\Link;\Gamma) \to C_2(\rats)\end{equation*}
given by Theorem~\ref{thmcompconf}
is as described above.
\end{proposition}

Proposition~\ref{propconffaceun} is proved in Section~\ref{secproofblowup}.
It will be used to prove that $\Zinv^s$ and $\Zinv$ are independent of the used propagating forms of Theorems~\ref{thmfstconst} and \ref{thmconststraight}, in Chapters~\ref{chapindepform} and \ref{chapanom}.

All the results stated so far in this chapter are sufficient to understand the proofs of Theorems~\ref{thmfstconst} and \ref{thmconststraight}. We prove them in detail in
the rest of this chapter, which can be skipped by a reader already convinced by the outline above. 

\section{Detailed study of \texorpdfstring{$\ccompuptd{\finsetv}{\vecspt}$}{SV(T)}}
\label{secprooffacanom}

In this section, we study the configuration space $\ccompuptd{\finsetv}{\vecspt}$ presented in Section~\ref{secpresfacanom}.
We first prove Theorem~\ref{thmcompfacanom}.
Let us first describe the transformations operated by the first blow-ups, locally.

Fix $\vecspt$, equip $\vecspt$ and its powers with Euclidean norms.
Fix $\tilde{w}_0 \in \cuptd{\finsetv}{\vecspt}$ and $b(\finsetv) \in \finsetv$.
Identify $\cuptd{\finsetv}{\vecspt}$ with the unit sphere $\sunits\left(\vecspt^{\finsetv \setminus \{b(\finsetv)\}}\right)$ of $\vecspt^{\finsetv \setminus \{b(\finsetv)\}}$.
Then $\tilde{w}_0$ is viewed as a map from $\finsetv$ to $\vecspt$ such that
$\tilde{w}_0\bigl(b({\finsetv})\bigr)=0$ and $\sum_{\eltv \in \finsetv} \norm{\tilde{w}_0(\eltv)}^2=1$.
 Recall Definition~\ref{defpartitionone} for a partition associated to a map and the associated notation.
\begin{definition}
\label{defbasedpartition}
A \emph{based partition} of a finite set $\finsetv$ equipped with a basepoint $b(\finsetv) \in \finsetv$ is a partition $\kids(\finsetv)$ of $\finsetv$ into nonempty subsets $\finseta$ equipped with basepoints $b(\finseta)$ such that
\begin{itemize}
 \item $b(\finseta)\in \finseta$, and
\item if $b(\finsetv) \in \finseta$, then $b(\finseta)=b(\finsetv)$. 
\end{itemize}
\end{definition}
Fix $\tilde{w}_0$. Let $\kids(\finsetv)=\kids(\finsetv;\tilde{w}_0)$ be the associated fixed partition.
Fix associated basepoints so that $\kids(\finsetv)$ becomes a based partition.

In general for a based subset $\finseta$ of $\finsetv$ equipped with a based partition $(\kids({\finseta}),b)$, define the set $O({\finseta},\kids({\finseta}),b,\vecspt)$ \index[N]{Open configuration sets!OCSb@$O({\finseta},\kids({\finseta}),b,\vecspt)$} 
of maps $w:{\finsetv} \longrightarrow {\vecspt}$ such that
\begin{itemize}
\item $\sum_{{\finsetb} \in \kids({\finseta})} \lefnorm{w\bigl(b({\finsetb})\bigr)}^2=1$,
\item $w(b({\finseta}))=0$, $w({\finsetv} \setminus {\finseta})=\{0\}$, and
\item two elements of ${\finseta}$ that belong to different kids
of ${\finseta}$ are mapped to different points of ${\vecspt}$ by $w$.
\end{itemize}

There is a straightforward identification of $O({\finsetv},\kids({\finsetv}),b,\vecspt)$ with an open subset of $\cuptd{\finsetv}{{\vecspt}}$,
which contains $\tilde{w}_0$.
Let $w_0 \in O({\finsetv},\kids({\finsetv}),b,\vecspt)$ denote the element that corresponds to $\tilde{w}_0$ under this identification.
Set 
\begin{equation*}W_{\finsetv}=O\bigl({\finsetv},\kids({\finsetv}),b,\vecspt\bigr) \cap \left( \cap_{\finseta \in D(\finsetv)}\diag_{\finseta}\left(\cuptd{\finsetv}{{\vecspt}}\right)\right).\end{equation*}

Note that $\kids({\finsetv})$ is a set that is naturally based by the element $b\bigl(\kids({\finsetv})\bigr)$ of $\kids({\finsetv})$ that contains $b(\finsetv)$. Then $W_{\finsetv}$ is an open subset of 
$\sunits(T^{\kids({\finsetv}) \setminus \{b(\kids(\finsetv))\}})$. It is the image $\check{\sunits}(T^{\kids({\finsetv}) \setminus \{b(\kids(\finsetv))\}})$ of $\cinjuptd{\kids(\finsetv)}{T}$ under the canonical identification of $\cuptd{\kids(\finsetv)}{T}$ with $\sunits(T^{\kids({\finsetv}) \setminus \{b(\kids(\finsetv))\}})$.

For $\finseta \in \kids({\finsetv})$, view the elements of $\vecspt^{\finseta \setminus \{b({\finseta})\}}$ as the maps from ${\finsetv}$ to $\vecspt$ that map $({\finsetv} \setminus {\finseta}) \cup \{b({\finseta})\}$ to $0$, and let $\vecspt^{\finseta \setminus \{b({\finseta})\}}_{< \varepsilon}$ denote 
the ball of its elements of norm smaller than $\varepsilon$.
Note the easy lemma.
\begin{lemma} \label{lem819}
 Let $\kids({\finsetv})$ be a based partition of $\finsetv$.
Let $w_0 \in W_{\finsetv}$. Then there exist an open neighborhood $\neighn(w_0)$ of $w_0$ in $W_{\finsetv}$ and an $\varepsilon \in \left]0,\infty\right[$ such that
the map
\begin{equation*}\begin{array}{lll}\neighn(w_0) \times \prod_{\finseta \in D(\finsetv)} \vecspt^{\finseta \setminus \{b({\finseta})\}}_{< \varepsilon} &\to & \cuptd{\finsetv}{{\vecspt}}\\
   \bigl(w,(\mu_{\finseta}\tilde{w}_{\finseta})_{{\finseta} \in D({\finsetv})}\bigr) &\mapsto& w +\sum_{{\finseta} \in D({\finsetv})} \mu_{\finseta}\tilde{w}_{\finseta}\end{array}\end{equation*}
is an open embedding whose image does not meet diagonals that do not correspond to (nonstrict) subsets $\finsetb$ of daughters of $\finsetv$.
\end{lemma}
\eopwobp

In particular, the first blow-ups that will affect this neighborhood of $\tilde{w}_0$ in $\cuptd{\finsetv}{{\vecspt}}$ are blow-ups along diagonals corresponding to daughters of $\finsetv$.
For any daughter $\finseta$ of $\finsetv$, the above identification identifies the normal bundle to $\diag_{\finseta}\bigl(\cuptd{\finsetv}{{\vecspt}}\bigr)$ with $\vecspt^{\finseta \setminus \{b({\finseta})\}}$, and the corresponding blow-up replaces the factor $\vecspt^{\finseta \setminus \{b({\finseta})\}}_{< \varepsilon}$ with $\left[0,\varepsilon\right[ \times \sunits(\vecspt^{\finseta \setminus \{b({\finseta})\}})$. Thus, it is clear that the blow-ups corresponding to different daughters of $\finsetv$ commute.
Note that our argument also proves the following lemma.
\begin{lemma}
\label{lemnormfibdiagtwo}
Let $\finseta$ be a subset of $\finsetv$.
 The fiber of the unit normal bundle to $\diag_{\finseta}(\cuptd{\finsetv}{{\vecspt}})$ 
is
$\cuptd{\finseta}{\vecspt}$.
\end{lemma}
\eopwobp

When performing the blow-ups along the diagonals corresponding to the daughters of $\finsetv$,
we replace $\mu_{\finseta}\tilde{w}_{\finseta} \in \vecspt^{\finseta \setminus \{b({\finseta})\}}$, for $\mu_{\finseta} \in \left[0,\varepsilon\right[$ and $\tilde{w}_{\finseta} \in \sunits(\vecspt^{\finseta \setminus \{b({\finseta})\}})$, with $(\mu_{\finseta},\tilde{w}_{\finseta})$.
Thus, we replace $0 \in \prod_{\finseta \in D(\finsetv)} \vecspt^{\finseta \setminus \{b({\finseta})\}}_{< \varepsilon}$ with the set of normal vectors $\tilde{w}_{\finseta}$ that pop up during the blow-ups.

\begin{lemma}
\label{lemfirstchartsinj}
In particular, with the notation of Lemma~\ref{lem819}, we get a chart of 
\begin{equation*}\lefblowup{O\bigl({\finsetv},\kids({\finsetv}),b,\vecspt\bigr)}{\bigl(\diag_{\finseta}(\cuptd{\finsetv}{{\vecspt}})\cap O({\finsetv},\kids({\finsetv}),b,\vecspt)\bigr)_{{\finseta} \in D({\finsetv})}},\end{equation*}
which maps
\begin{equation*} \left(w,(\mu_{\finseta},\tilde{w}_{\finseta})_{{\finseta} \in D({\finsetv})}\right) \in \neighn(w_0) \times \prod_{\finseta \in D(\finsetv)} \left(\left[0,\varepsilon\right[ \times \sunits(\vecspt^{\finseta \setminus \{b({\finseta})\}})\right)\end{equation*}
to the element
\begin{equation*}\Bigl(w +\sum_{{\finseta} \in D({\finsetv})} \mu_{\finseta}\tilde{w}_{\finseta},(\tilde{w}_{\finseta})_{{\finseta} \in D({\finsetv}) \suchthat  \mu_{\finseta}=0}\Bigr)\end{equation*}
of \begin{equation*}\lefblowup{O\bigl({\finsetv},\kids({\finsetv}),b,\vecspt\bigr)}{\left(\diag_{\finseta}(\cuptd{\finsetv}{{\vecspt}})\left(\cap O({\finsetv},\kids({\finsetv}),b,\vecspt)\right)\right)_{{\finseta} \in D({\finsetv})}},\end{equation*}
where the $\tilde{w}_{\finseta}$ are the normal vectors viewed in $\cuptd{\finseta}{\vecspt}$ that popped up during the blow-ups. 

We can construct an atlas of $\bigblowup{O({\finsetv},\kids({\finsetv}),b,\vecspt)}{(\diag_{\finseta}(\cuptd{\finsetv}{{\vecspt}}))_{{\finseta} \in D({\finsetv})}}$ with charts of this form.
\end{lemma}
\eopwobp

In order to conclude and get charts of manifolds with ridges, we blow up the
\begin{equation*}\cuptd{\finseta}{{\vecspt}} \cong\sunits(\vecspt^{\finseta \setminus \{b(A)\}})\end{equation*}
for the daughters $A$ of $V$, and we iterate.
Such an iteration produces a parenthesization of $\finsetv$ as in Definition~\ref{defparenthesizationone}.

\begin{definition}
\label{defparenthesizationtwo}
A \emph{$\Delta$-parenthesization}\index[T]{parenthesization!$\Delta$-parenthesization} of $\finsetv$ is a parenthesization $\parentp$ of $\finsetv$
such that $\finsetv \in \parentp$. The \emph{daughters}\index[T]{daughter} of an element $\finseta$ of a parenthesization $\parentp$ (with respect to $\parentp$) are the maximal elements of $\parentp$ strictly included in $\finseta$. The \indexT{mother} of an element $\finseta$ in $\parentp$ is the smallest element
of $\parentp$ that strictly contains $\finseta$ (if there is one). A $\Delta$-parenthesization $\parentp$ is organized as a tree, in which the vertices are the elements of $\parentp$ and the edges are in one-to-one correspondence with the pairs (daughter, mother) of $\parentp^2$.
We orient its edges from the daughter to her mother, as in Figure~\ref{figtreedeltaparone}.
The \emph{sons} of an element $\finseta$ of $\parentp$ are the singletons consisting of elements of $\finseta$ that do not belong to a daughter of $\finseta$ (with respect to $\parentp$). Any element ${\finseta}$ of $\parentp$ is equipped with the set $\kids({\finseta},\parentp)$ ($=\kids({\finseta})$ when $\parentp$ is fixed) of the
\emph{kids}\index[T]{kid} of ${\finseta}$, which are its daughters and its sons.
\end{definition}

\begin{example} 
\label{exaafterdefparenthesizationtwo}
The trees associated to the parenthesizations $\{\finsetv,\finseta_{123},\finseta_{23}\}$ and $\{\finsetv,\finseta_{12},\finseta_{34}\}$, which correspond to Figures~\ref{figmagnifyingglasses} and \ref{figmagnifyingglassestwo} in Section~\ref{secinformalcompactif}, are pictured in Figure~\ref{figtreedeltaparone}.
With respect to the parenthesization $\{\finsetv,\finseta_{123},\finseta_{23}\}$, the daughter $\finseta_{123}$ of $\finsetv$ has two kids, its son $\{1\}$ and its daughter $\finseta_{23}$.

\bfig \centering
\begin{tikzpicture}
\useasboundingbox (-1,-.2) rectangle (7.8,1);
\draw  (0,0) -- (2.8,.4) (4,0) -- (5.4,.4) -- (6.8,0) 
(2.8,.4) node[above]{\scriptsize $\finsetv$} (5.4,.4) node[above]{\scriptsize $\finsetv$}
(1.4,.2) node[above]{\scriptsize $\finseta_{123}$} (0,0) node[above]{\scriptsize $\finseta_{23}$}
(4,0) node[left]{\scriptsize $\finseta_{12}$} (6.8,0) node[right]{\scriptsize $\finseta_{34}$};
\draw[->] (0,0) -- (.7,.1);
\draw[->] (1.4,.2) -- (2.1,.3);
\draw[->] (4,0) -- (4.7,.2);
\draw[->] (6.8,0) -- (6.1,.2);
\fill (0,0) circle (1.5pt) (1.4,.2) circle (1.5pt) (2.8,.4) circle (1.5pt) (4,0) circle (1.5pt) (5.4,.4) circle (1.5pt) (6.8,0) circle (1.5pt);
\end{tikzpicture}
\caption{Trees associated to $\Delta$-parenthesizations}
\label{figtreedeltaparone}
\end{figure} 
\end{example}

\begin{definition}
\label{defbasedparenthesization}
Let $\setparent(\finsetv)$ (resp. $\setparent_{\Delta}(\finsetv)$) denote the set of parenthesizations (resp. $\Delta$-parenthesizations) of $\finsetv$.
Fix $\parentp \in \setparent_{\Delta}(\finsetv)$. 
Choose a basepoint $b({\finseta})=b({\finseta};\parentp)$\index[N]{basepoint@$b({\finseta})$ basepoint} for any ${\finseta} \in \parentp$ 
so that if ${\finseta}$ and $\finsetb$ are in $\parentp$,
if ${\finsetb} \subset {\finseta}$, and if $b({\finseta}) \in {\finsetb}$, then $b({\finsetb})=b({\finseta})$.
When ${\finseta} \in \parentp$, $D({\finseta})$ denotes the set of daughters of ${\finseta}$.
The basepoint $b({\finsetb})$ of a son ${\finsetb}$ of ${\finseta} \in \parentp$ is its unique point.
A $\Delta$-parenthesization equipped with basepoints as above is called a \indexT{based $\Delta$-parenthesization}.
\end{definition}

Recall the canonical identification of the set $O\bigl({\finseta},\kids({\finseta}),b,\vecspt\bigr)$ of Definition~\ref{defbasedpartition} with an open subset of $\cuptd{\finseta}{{\vecspt}}$.
Set 
\begin{equation*}W_{\finseta}=O\bigl({\finseta},\kids({\finseta}),b,\vecspt\bigr) \cap \left(\cap_{B \in D(A)}\diag_{\finsetb}(\cuptd{\finseta}{{\vecspt}})\right).\end{equation*}
Note that $W_{\finseta}$ may be identified canonically with an open subset of the sphere $\cuptd{\kids(A)}{\vecspt}$.

For $((\mu_{\finseta})_{{\finseta} \in \parentp \setminus \{\finsetv\}}, (w_{\finseta})_{{\finseta} \in \parentp}) \in (\RR^+)^{\parentp \setminus \{\finsetv\}} \times \prod_{{\finseta} \in \parentp} W_{\finseta}$, and for ${\finsetb} \in \parentp$,
define $v_{\finsetb}((\mu_{\finseta})_{{\finseta} \in \parentp \setminus \{\finsetv\}}, (w_{\finseta})_{{\finseta} \in \parentp})$ as the following map from $\finsetb$ to ${\vecspt}$:
\begin{equation*}\begin{array}{ll}v_{\finsetb}\bigl((\mu_{\finseta}), (w_{\finseta})\bigr)&=\sum_{\finsetc \in \parentp \suchthat  \finsetc \subseteq \finsetb}\left(\prod_{\finsetd \in \parentp \suchthat  \finsetc \subseteq \finsetd \subset \finsetb} \mu_{\finsetd} \right) w_{\finsetc}\\
   &=w_{\finsetb}+ \sum_{\finsetc \in D(\finsetb)}\mu_{\finsetc}\left(w_{\finsetc} + \sum_{\finsetd \in D(\finsetc)}\mu_{\finsetd}\left(w_{\finsetd}+\dots\right)\right).
  \end{array}
\end{equation*}
The construction of $v_{\finsetb}$ is illustrated in Figure~\ref{figtreedeltapartwo}.
We can see it on the subtree whose vertices are the subsets of $\finsetb$. The map $v_{\finsetb}$ is the sum over the vertices $\finsetc$ of this subtree of the maps $w_{\finsetc}$ associated to its vertices, multiplied by the products of the coefficients $\mu_{\finsetd}$ associated to the edges of the path from $\finsetc$ to $\finsetb$.

\bfig \centering
\begin{tikzpicture}
\useasboundingbox (-1,-.2) rectangle (7.8,1);
\draw  (0,0) -- (2.8,.4) (4,0) -- (5.4,.4) -- (6.8,0) 
(2.8,.4) node[above]{\scriptsize $w_{\finsetv}$} (5.4,.4) node[above]{\scriptsize $w_{\finsetv}$}
(1.4,.2) node[above]{\scriptsize $w_{\finseta_{123}}$} (0,0) node[above]{\scriptsize $w_{\finseta_{23}}$}
(4,0) node[left]{\scriptsize $w_{\finseta_{12}}$} (6.8,0) node[right]{\scriptsize $w_{\finseta_{34}}$};
\draw[->] (.7,.1) node[below]{\scriptsize $\mu_{\finseta_{23}}$} (0,0) -- (.7,.1);
\draw[->] (2.1,.3) node[below]{\scriptsize $\mu_{\finseta_{123}}$}(1.4,.2) -- (2.1,.3);
\draw[->] (4.7,.2) node[below]{\scriptsize $\mu_{\finseta_{12}}$} (4,0) -- (4.7,.2);
\draw[->] (6.1,.2) node[below]{\scriptsize $\mu_{\finseta_{34}}$} (6.8,0) -- (6.1,.2);
\fill (0,0) circle (1.5pt) (1.4,.2) circle (1.5pt) (2.8,.4) circle (1.5pt) (4,0) circle (1.5pt) (5.4,.4) circle (1.5pt) (6.8,0) circle (1.5pt);
\end{tikzpicture}
\caption{About the construction of $v_{\finsetb}$}
\label{figtreedeltapartwo}
\end{figure}

Note that $v_{\finsetb}((\mu_{\finseta}), (w_{\finseta})) \in O({\finsetb},\kids({\finsetb}),b,\vecspt)$ when the $\mu_{\finseta}$ are small enough.
Also note the following easy lemma.
\begin{lemma}
\label{lemneighchartcomf}
For any $(w^0_{\finseta})_{{\finseta} \in \parentp} \in \prod_{{\finseta} \in \parentp} W_{\finseta}$, there exists a neighborhood $\neighn((w^0_{\finseta}))$ of 
$0 \times (w^0_{\finseta})$ in  $(\RR^+)^{\parentp \setminus \{\finsetv\}} \times \prod_{{\finseta} \in \parentp} W_{\finseta}$ such that,\\ for any $((\mu_{\finseta})_{{\finseta} \in \parentp \setminus \{\finsetv\}}, (w_{\finseta})_{{\finseta} \in \parentp}) \in \neighn((w^0_{\finseta}))$,
\begin{itemize}
\item if $v_{\finsetv}((\mu_{\finseta}), (w_{\finseta}))$ is constant on $\finsetb$ for a subset $\finsetb$ of $\finsetv$, then
$\finsetb$ is included in (or equal to) a daughter of $\finsetv$,
\item if $\mu_{\finseta} \neq 0$ for any ${\finseta}\in \parentp$, then $v_{\finsetv}\bigl((\mu_{\finseta}), (w_{\finseta})\bigr)$ is an injective map from $\finsetv$ to ${\vecspt}$.
\end{itemize}
\end{lemma}
\eopwobp

When the construction of Theorem~\ref{thmcompfacanom} makes sense, we denote a point of our configuration space $\ccompuptd{\finsetv}{\vecspt}$ as a tuple \begin{equation*}\biggl(v_{\finsetv}\bigl((\mu_{\finseta}), (w_{\finseta})\bigr) ,\Bigl(v_{\finsetb}\bigl((\mu_{\finseta}), (w_{\finseta})\bigr)\Bigr)_{\finsetb \suchthat \mu_{\finsetb}=0}\biggr),\end{equation*} which contains its blowdown projection $v_{\finsetv}\bigl((\mu_{\finseta}), (w_{\finseta})\bigr)$ in $\cuptd{\finsetv}{\vecspt}$ followed by the normal vectors $v_{\finsetb}\bigl((\mu_{\finseta}), (w_{\finseta})\bigr)$ that have popped up during the blow-ups.
Lemma~\ref{lemnormfibdiagtwo} ensures that the normal vectors are elements of $\cuptd{\finsetb}{\vecspt}$ for some $\finsetb$. These  normal vectors are nonconstant maps from $\finsetb$ to ${\vecspt}$ up to translation and dilation.

\begin{theorem} 
\label{thmcompuptd}
Theorem~\ref{thmcompfacanom} is correct and defines $\ccompuptd{\finsetv}{\vecspt}$.
Let $\parentp$ be a $\Delta$-parenthesization of $\finsetv$.
Let $\ccompuptd{\finsetv,\parentp}{\vecspt}$ be the space of the elements that have been transformed by the blow-ups along the closures of the $(\diag_{\finseta}(\vecspt))_{\finseta \in \parentp \setminus \{\finsetv\}}$ and that have not been transformed by other blow-ups.
Then $\ccompuptd{\finsetv,\parentp}{\vecspt}$ is canonically diffeomorphic to 
\begin{equation*}\prod_{\finseta \in \parentp}\cinjuptd{\kids(\finseta)}{\vecspt}.\end{equation*}
The composed blowdown map sends $((w_{\finseta})_{\finseta \in \parentp})\in \prod_{\finseta \in \parentp}\cinjuptd{\kids(\finseta)}{\vecspt}$ to the map that sends
an element of a kid $\finsetb$ of $\finsetv$ to $w_{\finsetv}(\finsetb)$. For $\finseta \in \parentp \setminus \{\finsetv\}$, the $w_{\finseta}$ are (similarly identified with) the unit normal vectors that have appeared during the blow-ups.

As a set, $\ccompuptd{\finsetv}{\vecspt}$ is the disjoint union over the $\Delta$-parenthesizations $\parentp$ of $\finsetv$ of the $\ccompuptd{\finsetv,\parentp}{\vecspt}$.

 Any based $\Delta$-parenthesization $\parentp$ of $\finsetv$ and any $(w^0_{\finseta})_{{\finseta} \in \parentp} \in \prod_{{\finseta} \in \parentp} W_{\finseta}$ provide the following smooth open embedding
$\psi\bigl(\parentp,(w^0_{\finseta})_{{\finseta} \in \parentp}\bigr)$ from a neighborhood $\neighn\bigl((w^0_{\finseta})\bigr)$ as in Lemma~\ref{lemneighchartcomf} to $\ccompuptd{\finsetv}{\vecspt}$:
\begin{equation*}\begin{array}{lll}\neighn\bigl((w^0_{\finseta})\bigr) &\hookrightarrow &\ccompuptd{\finsetv}{\vecspt}\\
   \bigl((\mu_{\finseta})_{{\finseta} \in \parentp \setminus \{\finsetv\}}, (w_{\finseta})_{{\finseta} \in \parentp}\bigr) & \mapsto & \left( \begin{array}{l}v_{\finsetv}\bigl((\mu_{\finseta}), (w_{\finseta})\bigr),\\
\left(v_{\finsetb}\bigl((\mu_{\finseta}), (w_{\finseta})\bigr)\right)_{\{{\finsetb} \in \parentp \setminus \{\finsetv\} \suchthat \mu_{\finsetb} =0\}}\end{array}\right).
  \end{array}
\end{equation*}
This embedding restricts to $\neighn\bigl((w^0_{\finseta})\bigr) \cap \left( (\RR^{+\ast})^{\parentp\setminus \{\finsetv\}} \times \prod_{{\finseta} \in \parentp} W_{\finseta}\right)$ as a diffeomorphism onto an open subset of $\cinjuptd{\finsetv}{\vecspt}$.
Furthermore, the open images of the embeddings $\psi\bigl(\parentp,(w^0_{\finseta})_{{\finseta} \in \parentp}\bigr)$ for different $\bigl(\parentp,(w^0_{\finseta
})_{{\finseta} \in \parentp}\bigr)$ cover $\ccompuptd{\finsetv}{\vecspt}$, and the codimension of $\ccompuptd{\finsetv,\parentp}{\vecspt}$ in $\ccompuptd{\finsetv}{\vecspt}$ equals the cardinality of $\parentp \setminus \{\finsetv\}$.
\end{theorem}
\bp Note that the images of the embeddings corresponding to $\parentp=\{\finsetv\}$  cover $\cinjuptd{\finsetv}{\vecspt}$ trivially. The theorem is obviously true when $\cardlef{\finsetv}=2$. We proceed by induction on $\cardlef{\finsetv}$. By induction, for each daughter $\finseta$ of $\finsetv$, the space $\ccompuptd{\finseta}{\vecspt}$ is covered by images of smooth open embeddings $\psi(\parentp_{\finseta},(w^0_{\finsetb})_{{\finsetb} \in \parentp_{\finseta}})$ from $\neighn\left((w^0_{\finsetb})_{\finsetb \in \parentp_{\finseta}}\right)$ to some open subset $U_{\finseta}$ of $\ccompuptd{\finseta}{\vecspt}$
associated to parenthesizations $\parentp_{\finseta}$ of $\finseta$,
as in the statement.
So Lemma~\ref{lemfirstchartsinj} leaves us with the proof that $\psi\bigl(\parentp,(w^0_{\finseta})_{{\finseta} \in \parentp}\bigr)$ is smooth and open for any $(w^0_{\finseta})_{{\finseta} \in \parentp} \in \prod_{{\finseta} \in \parentp} W_{\finseta}$, for a small neighborhood $\neighn\bigl((w^0_{\finseta})\bigr)$ satisfying the conditions of Lemma~\ref{lemneighchartcomf}. 
It suffices to check that the change of coordinates replacing 
the coordinates $(\tilde{\mu}_A, \tilde{w}_A)$ written as $({\mu}_A, \tilde{w}_A)$ in Lemma~\ref{lemfirstchartsinj} by the coordinates $(\mu_A,v_A)$ compatible with $\psi\bigl(\parentp,(w^0_{\finseta})_{{\finseta} \in \parentp}\bigr)$ is smooth and open for each daughter $\finseta$ of $\finsetv$.
We have $\tilde{\mu}_A\tilde{w}_A= \mu_Av_A$, $\sum_{\elta \in \finseta}\norm{\tilde{w}_A(\elta)}^2=1$, and $\sum_{\finsetb \in D(\finseta)}\norm{v_A(b(\finsetb))}^2=1$.
We get \begin{equation*}
        \tilde{w}_A=\left(\sum_{\elta \in \finseta}\norm{v_A(\elta)}^2 \right)^{-1/2}v_A \mbox{, }v_A = \left(\sum_{\finsetb \in D(\finseta)}\norm{\tilde{w}_A(b(\finsetb))}^2\right)^{-1/2}\tilde{w}_A,
       \end{equation*}
\begin{equation*}
  \tilde{\mu}_A=\left(\sum_{\elta \in \finseta}\norm{v_A(\elta)}^2 \right)^{1/2}\mu_A\mbox{, and } \mu_A= \left(\sum_{\finsetb \in D(\finseta)}\norm{\tilde{w}_A(b(\finsetb))}^2\right)^{1/2} 
\tilde{\mu}_A,
 \end{equation*}
where $\sum_{\elta \in \finseta}\norm{v_A(\elta)}^2 \geq 1$, and $\sum_{\finsetb \in D(\finseta)}\norm{\tilde{w}_A(b(\finsetb))}^2>0$. 
\eop

The above proof also proves the following two propositions, with the notation $\partial_{r}$ introduced in the beginning of Section~\ref{secbackgroundmfdsbry}.
\begin{proposition}
\label{propdescconftauone}
Let $\setparent_{r,\Delta}(\finsetv)$ be the subset of $\setparent_{\Delta}(\finsetv)$ consisting of the $\Delta$-parenthesizations of $\finsetv$ of cardinality $r$.
Then \begin{equation*}\partial_{r-1}\bigl(\ccompuptd{\finsetv}{\vecspt}\bigr)\setminus \partial_{r}\bigl(\ccompuptd{\finsetv}{\vecspt}\bigr)
=\sqcup_{\parentp \in \setparent_{r,\Delta}(\finsetv)}\ccompuptd{\finsetv,\parentp}{\vecspt}.\end{equation*}
For two $\Delta$-parenthesizations $\parentp$ and $\parentp^{\prime}$, if $\parentp \subset \parentp^{\prime}$, then 
$\ccompuptd{\finsetv,\parentp^{\prime}}{\vecspt} \subset \overline{\ccompuptd{\finsetv,\parentp}{\vecspt}}$.
\end{proposition}
\bp We can deduce the first assertion from the charts of Theorem~\ref{thmcompuptd}.
Let $c_0=(w_{\finseta})_{{\finseta} \in \parentp^{\prime}} \in \ccompuptd{\finsetv,\parentp^{\prime}}{\vecspt}$ 
and let $\Psi=\psi(\parentp^{\prime},(w_{\finseta})_{{\finseta} \in \parentp^{\prime}}) \colon \neighn((w_{\finseta})) \to  \ccompuptd{\finsetv}{\vecspt}$ be a smooth open embedding as in Theorem~\ref{thmcompuptd}.
Let $\varepsilon \in \left]0,\infty\right[$ be such that $\left[0,\varepsilon\right]^{\parentp^{\prime} \setminus \{\finsetv\}} \times \{(w_{\finseta})_{{\finseta} \in \parentp^{\prime}}\} \subset \neighn((w_{\finseta})_{{\finseta} \in \parentp^{\prime}})$.
For $t \in \left[0,1\right]$, set 
\begin{equation*}\mu_{\finseta}(t)=\left\{\begin{array}{ll} \varepsilon t & \mbox{if}\; \finseta \in \parentp^{\prime} \setminus \parentp \\
                       0 &  \mbox{if}\; \finseta \in \parentp \setminus \{\finsetv\}
                      \end{array}\right.\end{equation*}
and let $c(t)=\Psi\left((\mu_{\finseta}(t))_{\finseta \in \parentp^{\prime} \setminus \{\finsetv\}},(w_{\finseta})_{{\finseta} \in \parentp^{\prime}}\right) $.
Then $c(t) \in \ccompuptd{\finsetv,\parentp}{\vecspt}$ for any $t \in \left]0,1\right]$, and $\lim_{t\to 0}c(t)=c_0$.
\eop

\begin{proposition}
 Any injective linear map $\phi$ from a vector space $\vecspt$ to another such $\vecspt^{\prime}$ induces a canonical embedding $\phi_{\ast} \colon \ccompuptd{\finsetv}{\vecspt} \to \ccompuptd{\finsetv}{\vecspt^{\prime}}$.
This embedding maps an element $\bigl((w_{\finseta})_{\finseta \in \parentp}\bigr)$ of $\ccompuptd{\finsetv,\parentp}{\vecspt}$ to the element $\bigl((\phi \circ w_{\finseta})_{\finseta \in \parentp}\bigr)$ of $\ccompuptd{\finsetv,\parentp}{\vecspt^{\prime}}$.
If $\psi$ is another injective linear map from a vector space $\vecspt^{\prime}$ to a third vector space $\vecspt^{\prime \prime}$, then we have $(\psi \circ \phi)_{\ast}= \psi_{\ast} \circ \phi_{\ast}$.
\end{proposition}
\eopwobp

Finally, let us check the following proposition, which implies Theorem~\ref{thmcompfacanomtwo}.
\begin{proposition}
\label{proprestanom}
Let  $\finseta$ be a finite subset of cardinality at least $2$ of a finite set $\finsetv$.
For $\parentp \in \setparent_{\Delta}(\finsetv)$, define 
\begin{equation*}\parentp_{\finseta} =\bigl\{\finsetb \cap \finseta \suchthat  \finsetb \in \parentp, \cardlef{\finsetb \cap \finseta} \geq 2\bigr\}\end{equation*}
in $\setparent_{\Delta}(\finseta)$. For $\finsetc \in \parentp_{\finseta}$, let $\hat{\finsetc}$ be the smallest element of $\parentp$ that contains $\finsetc$ or is equal to $\finsetc$.
For
\begin{equation*}w=\bigl(w_{\finsetb} \in \cinjuptd{\kids(\finsetb)}{\vecspt}\bigr)_{\finsetb \in \parentp} \in \ccompuptd{\finsetv,\parentp}{\vecspt},\end{equation*}
and for $\finsetc \in \parentp_{\finseta}$, define $w^{\prime}_{\finsetc}$ to be the natural restriction of $w_{\hat{\finsetc}}$ to $\kids(\finsetc,\parentp_{\finseta})$. Then set
\begin{equation*}p_{\finseta}(w)=\bigl((w^{\prime}_{\finsetc})_{\finsetc \in \parentp_{\finseta}}\bigr)
\in \ccompuptd{\finseta,\parentp_{\finseta}}{\vecspt}.\end{equation*}
This consistently defines a smooth map \begin{equation*}p_{\finseta} \colon \ccompuptd{\finsetv}{\vecspt} \rightarrow \ccompuptd{\finseta}{\vecspt}.\end{equation*}
The map $p_{\finseta}$ is the unique continuous extension from $\ccompuptd{\finsetv}{\vecspt}$ to $\ccompuptd{\finseta}{\vecspt}$ of the restriction map from $\cinjuptd{\finsetv}{\vecspt}$ to $\cinjuptd{\finseta}{\vecspt}$.
\end{proposition}
\bp 
It is easy to see that this restriction map is well-defined. In order to prove that it is smooth, use the charts of Theorem~\ref{thmcompuptd}.
Fix a $\Delta$-parenthesization $\parentp$ of $\finsetv$. Let $\parentp_{\finseta}$ be the induced $\Delta$-parenthesization of $\finseta$. 
Fix basepoints $b_{\finseta}$ for $\parentp_{\finseta}$ according to the rule in Definition~\ref{defbasedparenthesization}. Fix basepoints for the elements $\finsetb$ of
$\parentp$ so that $b(\finsetb) =b_{\finseta}( \finseta \cap \finsetb)$ when $\finsetb \cap \finseta \in \parentp_{\finseta}$, and $b(\finsetb) \in \finsetb \cap \finseta$ when $\finsetb \cap \finseta \neq \emptyset$.
According to Theorem~\ref{thmcompuptd}, it suffices to prove smoothness in charts involving the based $\Delta$-parenthesizations $\parentp$ and $\parentp_{\finseta}$. So it suffices to prove that the projections on the factors that contain the restrictions maps $w^{\prime}_{\finsetc}$ are smooth for $\finsetc \in\parentp_{\finseta}$,
and that the projections on the factors that contain the dilation factors $\mu_{\finsetc}$ are smooth for $\finsetc \in\parentp_{\finseta} \setminus \{\finseta\}$.

With our charts and our conditions on the basepoints, for any $\finsetc \in \parentp_{\finseta}$, we have
\begin{equation*}w^{\prime}_{\finsetc}=\frac{1}{g(\finsetc,w_{\hat{\finsetc}})}w_{\hat{\finsetc}}\vert_{\kids(\finsetc,\parentp_{\finseta})},\end{equation*}
where \begin{equation*}g(\finsetc,w_{\hat{\finsetc}})=\sqrt{\sum_{\finsetd \in \kids(\finsetc,\parentp_{\finseta})}\norm{w_{\hat{\finsetc}}(\finsetd)}^2}\end{equation*}
is not zero since $w_{\hat{\finsetc}}$ is nonconstant on $\kids(\finsetc,\parentp_{\finseta})$.
So, the projection on the factor of $w^{\prime}_{\finsetc}$ is smooth.

Let $\finsetc \in \parentp_{\finseta}\setminus \{\finseta\}$, and let $\mother(\finsetc)$ denote the \emph{mother} of $\finsetc$ in $\parentp_{\finseta}$. 
Let $\finsete$ be the set of basepoints of the kids of $\mother(\finsetc)$ distinct from $\finsetc$ with respect to $\parentp_{\finseta}$.
Consider the elements $\finsetb_i \in \parentp$, for $i =1, 2, \dots, k(\finsetb)$, such that $\finsetc =\finsetb_i \cap \finseta$ for any $i\leq k(\finsetc)$, where $(\hat{\finsetc}=\finsetb_1) \subset \finsetb_2 \dots \subset \finsetb_{k(\finsetc)}$. We have $b(\finsetb_i) =b_{\finseta}( \finsetc)$  for any $i\in \underline {k(\finsetc)}$ and $b(\widehat{\mother(\finsetc)}) =b_{\finseta}(\mother(\finsetc))$.

Then the restrictions of the configurations $w^{\prime}_{\mother(\finsetc)} + \mu_{\finsetc}w^{\prime}_{\finsetc}$ and 
\begin{equation*}w_{\widehat{\mother(\finsetc)}} +\left( \prod_{i=1}^{k(\finsetc)} \mu_{\finsetb_i} \right) w_{\hat{\finsetc}}\end{equation*}
to $\finsetc \cup \finsete$ coincide up to dilation.
So $g(\mother(\finsetc),w_{\widehat{\mother(\finsetc)}})\mu_{\finsetc}w^{\prime}_{\finsetc}$
coincides with \begin{equation*}\left( \prod_{i=1}^{k(\finsetc)} \mu_{\finsetb_i} \right)g(\finsetc,w_{\hat{\finsetc}})w^{\prime}_{\finsetc}\end{equation*} on $\finsetc$, and we have
\begin{equation*}\mu_{\finsetc}=\left( \prod_{i=1}^{k(\finsetc)} \mu_{\finsetb_i} \right)\frac{g(\finsetc,w_{\hat{\finsetc}})}{g(\mother(\finsetc),w_{\widehat{\mother(\finsetc)}})}.\end{equation*}
Thus $\mu_{\finsetc}$ is smooth (it is defined even when the $\mu_{\finsetb_i}$ are negative).
\eop

Proposition~\ref{propfacfacanom} follows from Propositions~\ref{propdescconftauone} and \ref{proprestanom}. \eopwobp

\bpo{Proof of Lemma~\ref{lemfacfacan}} 
The structure of a smooth manifold with ridges of
$\ccompuptdanvec(T_{\Link(z)}\crats,\vec{t}_{\Link(z)};\check{\Gamma})$ in Lemma~\ref{lemfacfacone}
can be alternatively deduced from the charts of Theorem~\ref{thmcompuptd}.
These charts also show that $\ccompuptdanvec(T_{\Link(z)}\crats,\vec{t}_{\Link(z)};\check{\Gamma})$
is a submanifold transverse to the ridges of $\ccompuptd{\finsetv(\check{\Gamma})}{T_{\Link(z)}\crats}$. So its codimension-one faces are the intersections of $\ccompuptdanvec(T_{\Link(z)}\crats,\vec{t}_{\Link(z)};\check{\Gamma})$ with the codimension-one faces of $\ccompuptd{\finsetv(\check{\Gamma})}{T_{\Link(z)}\crats}$.
Then Lemma~\ref{lemfacfacan} follows from Proposition~\ref{propfacfacanom}.
\eop

\section{Blowing up diagonals}
\label{secblodiag}

In the rest of this chapter, $\manifm$ is a smooth manifold without boundary of dimension $\dimdel>0$.
It is not necessarily oriented. The set of injective maps from $\finsetv$ to $\manifm$ is denoted
by $\check{C}_{\finsetv}\!\left[\manifm\right]$\index[N]{Configuration spaces!Ccheckbrac@$\check{C}_{\finsetv}[\manifm]$ open} with brackets instead of parentheses (we have $\check{C}_{\finsetv}\!\left[\crats\right]=\check{C}_{\finsetv}(\crats)$, but $\check{C}_{\finsetv}\!\left[\rats\right] \neq \check{C}_{\finsetv}(\rats)$).

\begin{theorem}
\label{thmcompdiaggen}
Let $\manifm$ be a manifold. Let $\finsetv$ be a finite set.

Set $C_{\finsetv,\card{\finsetv}+1}\!\left[\manifm\right]=\manifm^{\finsetv}$.
For $k=\cardlef{\finsetv}, \dots, 3, 2$,  define $C_{\finsetv,k}\!\left[\manifm\right]$ from $C_{\finsetv,k+1}\!\left[\manifm\right]$
by blowing up the closures of
(the preimages under the composition of the previous blowdown maps of) the $\diag_{{\finseta}}(\manifm^{\finsetv})$ such that $\cardlef{\finseta}=k$. At each step, the blown-up manifolds are  smooth and transverse to the ridges, and the resulting blown-up manifold is independent of the order choice.
Thus, this process gives rise to a canonical smooth manifold $C_{\finsetv}\!\left[\manifm\right]=C_{\finsetv,2}\!\left[\manifm\right]$ with ridges equipped with its composed blowdown projection
\begin{equation*}\pbl \colon C_{\finsetv}\!\left[\manifm\right] \to \manifm^{\finsetv}.\end{equation*}
\begin{itemize}
\item 
Let $\funcf \in \manifm^{\finsetv}$ be a map from
$\finsetv$ to $\manifm$. Then
$\pbl^{-1}(\funcf)$ is canonically diffeomorphic to $\prod_{\finseta \in D(\finsetv;\funcf)} \ccompuptd{\finseta}{T_{\funcf(\finseta)}\manifm}$.
We will denote an element $x \in \pbl^{-1}(\funcf)$ of $C_{\finsetv}\!\left[\manifm\right]$ by $\bigl(\funcf \in \manifm^{\finsetv},\bigl(w_{\finseta}\in \ccompuptd{\finseta}{T_{\funcf(\finseta)}\manifm}\bigr)_{\finseta \in D(\finsetv;\funcf)}\bigr)$, with the notation of Definition~\ref{defpartitionone}.

\item For any open subset $U$ of $\manifm$, $C_{\finsetv}\!\left[U\right]=\pbl^{-1}(U^{\finsetv})$.

\item The space $\check{C}_{\finsetv}\!\left[\manifm\right]$ is dense in ${C}_{\finsetv}\!\left[\manifm\right]$.

\item If $\manifm$ is compact, then ${C}_{\finsetv}\!\left[\manifm\right]$ is compact, too.

\item Any choice of a basepoint $b=b(\finsetv)$ of $\finsetv$ and of an open embedding
\begin{equation*}\phi \colon \RR^{\dimdel} \to \manifm\end{equation*} 
induces the diffeomorphism
\begin{equation*}\psi(\phi,b) \colon \RR^{\dimdel} \times \RR^+  \times \ccompuptd{\finsetv}{\RR^{\dimdel}} \to C_{\finsetv}\!\left[\phi\left(\RR^{\dimdel}\right)\right]\end{equation*}
described below.\\
Let $u \in \RR^{\dimdel}$, $\mu \in \RR^{+}$, and $n \in \ccompuptd{\finsetv}{\RR^{\dimdel}}$.
Then $\psi(\phi,b)$ satisfies \begin{equation*}\pbl\bigl(\psi(\phi,b)(u,\mu,n)\bigr)=\phi \circ \bigl( u + \mu\pbl(n)\bigr),\end{equation*}
where $\pbl(n)$ is viewed as a map from $\finsetv$ to $\RR^{\dimdel}$ such that $\pbl(n)(b)=0$ and $\sum_{\eltv \in \finsetv}\norm{\pbl(n)(\eltv)}^2=1$, and
 $(u + \mu \pbl(n))$ denotes the map from
$\finsetv$ to $\RR^{\dimdel}$ obtained from $\pbl(n)$ by composition by the homothety with ratio $\mu$, followed by the translation of vector $u$.
For $u \in \RR^{\dimdel}$ and $\nnormalized \in \cinjuptd{\finsetv}{\RR^{\dimdel}}$, the map $\pbl(\psi(\phi,b)(u,0,\nnormalized))$ is constant with value $\phi(u)$ and we have
\begin{equation*}\psi(\phi,b)(u,0,\nnormalized)=\bigl(\phi(u)^{\finsetv}, (T_u\phi)_{\ast}(\nnormalized)\bigr).\end{equation*}
The restriction of $\psi(\phi,b)$ to $\RR^{\dimdel} \times \RR^{+\ast} \times \cinjuptd{\finsetv}{\RR^{\dimdel}}$ is a diffeomorphism 
onto $\check{C}_{\finsetv}\!\left[\phi\left(\RR^{\dimdel}\right)\right]$.
\item An embedding $\phi$ of a manifold $M_1$ into another such $M_2$ induces the canonical embedding $\phi_{\ast}$ such that
\begin{equation*}\phi_{\ast}\Bigl(\funcf \in M_1^{\finsetv},\bigl(w_{\finseta}\in \ccompuptd{\finseta}{T_{\funcf(\finseta)}M_1}\bigr)_{\finseta \in D(\finsetv;\funcf)}\Bigr)=\bigl(\phi \circ f, (T_{\funcf(\finseta)}\phi)_{\ast}(w_{\finseta})\bigr)\end{equation*}
 from $C_{\finsetv}\!\left[M_1\right]$ to $C_{\finsetv}\!\left[M_2\right]$. If $\psi$ is an embedding from $M_2$ to another manifold $M_3$, then we have $(\psi \circ \phi)_{\ast}=\psi_{\ast} \circ \phi_{\ast}$.
\end{itemize}
\end{theorem}
\bp Start with $\manifm =\RR^{\dimdel}$ equipped with its usual Euclidean norm. Fix $b(\finsetv)$.
Any map $\funcf$ from $\finsetv$ to $\RR^{\dimdel}$ may be written as $\funcf\bigl(b(\finsetv)\bigr) + y$, for a unique element $y \in (\RR^{\dimdel})^{\finsetv \setminus \{b(\finsetv)\}}$.
Then blowing up $(\RR^{\dimdel})^{\finsetv}$ along the diagonal $\diag_{\finsetv}((\RR^{\dimdel})^{\finsetv} )$, described by the equation $y=0$, replaces $(\RR^{\dimdel})^{\finsetv \setminus \{b(\finsetv)\}}$ with $\RR^+ \times \sunits((\RR^{\dimdel})^{\finsetv \setminus \{b(\finsetv)\}})$ and provides a diffeomorphism
from $\RR^{\dimdel}  \times \RR^+ \times \sunits\bigl((\RR^{\dimdel})^{\finsetv \setminus \{b(\finsetv)\}}\bigr)$
to $\bigblowup{(\RR^{\dimdel})^{\finsetv}}{\diag_{\finsetv}\bigl((\RR^{\dimdel})^{\finsetv} \bigr)}$.

The diagonals corresponding to strict subsets of $\finsetv$ are products by $\RR^{\dimdel}  \times \RR^+$ of the diagonals corresponding to the same subsets for the manifold
$\sunits\bigl((\RR^{\dimdel})^{\finsetv \setminus \{b(\finsetv)\}}\bigr)\cong \cuptd{\finsetv}{\RR^{\dimdel}}$, which was studied in the previous subsection.

Thus, ${C}_{\finsetv}\!\left[\RR^{\dimdel}\right]$ is well-defined, and our diffeomorphism lifts to a diffeomorphism from $\RR^{\dimdel}  \times \RR^+ \times \ccompuptd{\finsetv}{\RR^{\dimdel}}$
to ${C}_{\finsetv}\!\left[\RR^{\dimdel}\right]$. So the composition of this diffeomorphism with the product of the charts obtained in Theorem~\ref{thmcompuptd} by the identity map yields an atlas of ${C}_{\finsetv}\!\left[\RR^{\dimdel}\right]$.

For a diffeomorphism $\phi$ from $\RR^{\dimdel}$ to an open subspace $U$ of a manifold, the diffeomorphism \begin{equation*}\phi^{\finsetv} \colon (\RR^{\dimdel})^{\finsetv} \to U^{\finsetv}\end{equation*} preserves diagonals. So ${C}_{\finsetv}\!\left[U\right]$ is well-defined for any open subset $U$ diffeomorphic to $\RR^{\dimdel}$, in a manifold $\manifm$.
Furthermore, $\phi^{\finsetv}$ lifts as a natural diffeomorphism $\phi_{\ast} \colon {C}_{\finsetv}\!\left[\RR^{\dimdel}\right] \to {C}_{\finsetv}\!\left[U\right]$.

Note that the elements of ${C}_{\finsetv}\!\left[\RR^{\dimdel}\right]$ have the prescribed form.
Since the normal bundle to a diagonal $\diag_{{\finseta}}\bigl(\check{C}_{\finsetv}\!\left[U\right]\bigr)$ is $\bigl(({TU^{\finseta}}/{\diag_{\finseta}(TU)}) \setminus \{0\}\bigr)/\RR^{+\ast}$,
the diffeomorphism $\phi_{\ast}$ from
${C}_{\finsetv}\!\left[\RR^{\dimdel}\right]$ to ${C}_{\finsetv}\!\left[U\right]$ maps 
\begin{equation*}x=\Bigl(\funcf \in (\RR^{\dimdel})^{\finsetv},\bigl(w_{\finseta}\in \ccompuptd{\finseta}{\RR^{\dimdel}}\bigr)_{\finseta \in D(\finsetv;\funcf)}\Bigr)\end{equation*} to
\begin{equation*}\phi_{\ast}(x)=\Bigl(\phi \circ f, \bigl((T_{\funcf(\finseta)}\phi)_{\ast}(w_{\finseta})\bigr)_{\finseta \in D(\finsetv;\funcf)}\Bigr).\end{equation*}
Then the elements of ${C}_{\finsetv}\!\left[U\right]$ have the prescribed form, too.

In order to prove that ${C}_{\finsetv}\!\left[\manifm\right]$ is well-defined for a manifold $\manifm$, it suffices to see that it is well-defined over an open neighborhood of any point $\confc$ of $\manifm^{\finsetv}$. Such a map $\confc$ defines the partition $\kids(\finsetv;\confc)$.
There exist pairwise disjoint open neighborhoods $U_{\finseta}$ diffeomorphic to $\RR^{\dimdel}$ of the $\confc(\finseta)$, for $\finseta \in \kids(\finsetv)=\kids(\finsetv;\confc)$.
It is easy to see that ${C}_{\finsetv}\!\left[\manifm\right]$ is well-defined over $\prod_{\finseta \in \kids(\finsetv)} U_{\finseta}^{\finseta}$ and that it is canonically isomorphic to $\prod_{\finseta \in \kids(\finsetv)} {C}_{\finseta}\left[U_{\finseta}\right]$, there.
Thus, ${C}_{\finsetv}\!\left[\manifm\right]$ is well-defined, and its elements have the prescribed form.

If $\phi$ is a diffeomorphism from a manifold $M_1$ to another such $M_2$, then the diffeomorphism \begin{equation*}\phi^{\finsetv} \colon M_1^{\finsetv} \to M_2^{\finsetv}\end{equation*} preserves diagonals. So it lifts as a natural diffeomorphism $\phi_{\ast} \colon {C}_{\finsetv}\!\left[M_1\right] \to {C}_{\finsetv}\!\left[M_2\right]$, which behaves as stated.

Then the study of the map induced by an embedding from a manifold $M_1$ into another such $M_2$ can be easily reduced to the case of a linear embedding from $\RR^{k}$ to $\RR^{\dimdel}$. For such an embedding, use the identification
of $C_{\finsetv}\!\left[\RR^{\dimdel}\right]$ with $\RR^{\dimdel} \times \RR^+  \times \ccompuptd{\finsetv}{\RR^{\dimdel}}$.

The other statements are easy to check.
\eop

\begin{proposition}
\label{proprest} 
Let $\finsetb$ be a subset of $\finsetv$, then the restriction from  $\check{C}_{\finsetv}\!\left[\manifm\right]$ to $\check{C}_{\finsetb}\left[\manifm\right]$ extends uniquely to a smooth map $p_{\finsetb}$ from ${C}_{\finsetv}\!\left[\manifm\right]$ to ${C}_{\finsetb}\left[\manifm\right]$. 
Let $\funcf \in \manifm^{\finsetv}$.
The elements $\finseta$ of $D(\finsetb;\funcf\vert_{\finsetb})$ are of the form $\finsetb \cap \hat{\finseta}$ for a unique $\hat{\finseta}$ in $D(\finsetv;\funcf)$. Then we have
\begin{equation*}p_{\finsetb}\Bigl(\funcf,\bigl(w_{\finsetc}\in \ccompuptd{\finsetc}{T_{\funcf(\finsetc)}\manifm}\bigr)_{\finsetc \in D(\finsetv;\funcf)}\Bigr)=
\Bigl(\funcf\vert_{\finsetb},\bigl(w_{\hat{\finseta}}\vert_{\finseta}\bigr)_{\finseta \in D(\finsetb;\funcf\vert_{\finsetb})}\Bigr).\end{equation*}
\end{proposition}
\bp It is obvious when $\cardlef{\finsetb} =1$. Assume that $\cardlef{\finsetb} \geq 2$.
In order to check that this restriction map is smooth, it is enough to prove that it is smooth when $\manifm=\RR^{\dimdel}$. Assume that $b(\finsetv) \in \finsetb$.
Use the diffeomorphism of Theorem~\ref{thmcompdiaggen} to write an element of ${C}_{\finsetv}\!\left[\RR^{\dimdel}\right]$
as $(u \in \RR^{\dimdel},\mu \in \RR^+,n \in \ccompuptd{\finsetv}{\RR^{\dimdel}})$.
Then the restriction maps
$(u,\mu,n)$ to $\bigl(u,{\norm{\pbl(n)\vert_{\finsetb}}}{\mu},n\vert_{\finsetb}\bigr)$. It is smooth according to Proposition~\ref{proprestanom}.
\eop

An element $x=\bigl(\funcf \in \manifm^{\finsetv},(w_{\finseta}\in \ccompuptd{\finseta}{T_{\funcf(\finseta)}\manifm})_{\finseta \in D(\finsetv;\funcf)}\bigr)$ of 
${C}_{\finsetv}\!\left[\manifm\right]$ induces the parenthesization $\parentp(x)$ of $\finsetv$ 
that is the union over the elements $\finseta$ of $D(\finsetv;\funcf)$ of the $\parentp_{\finseta}(x)$ such that $w_{\finseta}\in \ccompuptd{\finseta,\parentp_{\finseta}(x)}{T_{\funcf(\finseta)}\manifm}$. (See Theorem~\ref{thmcompuptd}.)
Let $\parentp$ be a parenthesization of $\finsetv$.
Set $C_{\finsetv,\parentp}\left[\manifm\right]=\{ x\in C_{\finsetv}\left[\manifm\right] \suchthat \parentp(x)=\parentp\}$.

The following proposition is easy to observe. (See Proposition~\ref{propdescconftauone}.)
\begin{proposition}
\label{propdescconftautwo}
Let $\setparent_{r}(\finsetv)$ be the subset of $\setparent(\finsetv)$ consisting of the parenthesizations of $\finsetv$ of cardinality $r$.
Then we have \begin{equation*}\partial_{r}\left(C_{\finsetv}\!\left[\manifm\right]\right)\setminus \partial_{r+1}\left(C_{\finsetv}\!\left[\manifm\right]\right)
=\sqcup_{\parentp \in \setparent_{r}(\finsetv)}C_{\finsetv,\parentp}\!\left[\manifm\right].\end{equation*}
Furthermore, for any two parenthesizations $\parentp$ and $\parentp^{\prime}$ such that $\parentp \subset \parentp^{\prime}$, we have 
$C_{\finsetv,\parentp^{\prime}}\!\left[\manifm\right] \subset \overline{C_{\finsetv,\parentp}\!\left[\manifm\right]}$.
\end{proposition}
\eopwobp

\section{Blowing up \texorpdfstring{$\infty$}{the point at infinity}}
\label{secbloinf}
For a finite set $V$, let $V^+$ be obtained from $V$ by adding a special element $v_{\infty}$ to $V$. We have $V^+ = V \sqcup\{v_{\infty}\}$.

We state and prove the following generalization of Theorem~\ref{thmcompconf}, where ${C}_{\finsetv}(\rats)={C}_{\finsetv}\!\left[\rats,\infty\right]$. 

\begin{theorem}
 \label{thmcompconfbis} Let $\finsetv$ be a finite set.
Let $\manifm$ be a manifold without boundary of dimension $\dimdel$, and let $\infty \in \manifm$. Set $\cmanifm=\manifm \setminus \{\infty\}$.
Recall the manifold $C_{V^+}[\manifm]$ of Theorem~\ref{thmcompdiaggen}.
Let $\pbl\vert_{\{v_{\infty}\}} \colon C_{V^+}[\manifm] \to \manifm$ map a configuration to its value at $v_{\infty}$.
Define ${C}^+_{\finsetv}\!\left[\manifm,\infty\right]$ to be the preimage of $\infty$ under this map. 
Then ${C}^+_{\finsetv}\!\left[\manifm,\infty\right]$ is a smooth submanifold of $C_{V^+}[\manifm]$ transverse to the ridges equipped with charts induced by the local models of Theorems~\ref{thmcompdiaggen} 
and \ref{thmcompuptd}. 
It is the closure in $C_{V^+}[\manifm]$ of $\check{C}_{\finsetv}\!\left[\cmanifm\right] \times \{\infty\}^{\{v_{\infty}\}}$, where $\check{C}_{\finsetv}\!\left[\cmanifm\right] \subset \cmanifm^{\finsetv}$ is the space
 of injective maps from $\finsetv$ to $\cmanifm$.
 
Furthermore, the manifold ${C}^+_{\finsetv}\!\left[\manifm,\infty\right]$ is canonically diffeomorphic to the manifold ${C}_{\finsetv}\!\left[\manifm,\infty\right]$ obtained from $\manifm^{\finsetv}$ by the following process.

For a nonempty ${\finseta} \subseteq {\finsetv}$, let $\einftyxi_{\finseta}$ be the set of maps from $\finsetv$ to $\manifm$  that map ${\finseta}$ to $\infty$, and $\finsetv \setminus \finseta$ to $\cmanifm$. 

Start with $\manifm^{\finsetv}$.
Blow up $\einftyxi_{\finsetv}$, which is reduced to the point $m=\infty^{\finsetv}$ such that $m^{-1}(\infty)={\finsetv}$. Set 
$${C}_{\finsetv,\card{\finsetv}+1}\!\left[\manifm,\infty\right]=\bigblowup{\manifm^{\finsetv}}{\infty^{\finsetv}}.$$\\
Then for $k=\cardlef{\finsetv},\cardlef{\finsetv}-1, \dots, 3, 2$, define
${C}_{\finsetv,k}\!\left[\manifm,\infty\right]$ from ${C}_{\finsetv,k+1}\!\left[\manifm,\infty\right]$ by blowing up 
the closures of (the preimages under the composition of the previous blowdown maps of) the $\diag_{{\finseta}}(\cmanifm^{\finsetv})$ such that $\cardlef{\finseta}=k$ and the closures of (the preimages under the composition of the previous blowdown maps of) the $\einftyxi_J$ such that $\cardlef{J}=k-1$. 
At each step, the  blown-up manifolds are smooth and transverse to the ridges, and ${C}_{\finsetv,k}\!\left[\manifm,\infty\right]$ is independent of the possible order choice of the blow-ups.
The obtained manifold ${C}_{\finsetv}\!\left[\manifm,\infty\right]={C}_{\finsetv,2}\!\left[\manifm,\infty\right] $ is a smooth manifold of dimension $\dimdel\!\cardlef{\finsetv}$, with ridges. 
It is compact if $\manifm$ is compact.
The interior of ${C}_{\finsetv}\!\left[\manifm,\infty\right]$ is $\check{C}_{\finsetv}\!\left[\cmanifm\right]$, and the composition of the blowdown maps gives rise to a canonical smooth blowdown projection
$\pbl \colon {C}_{\finsetv}\!\left[\manifm,\infty\right] \to \manifm^{\finsetv}$.
\end{theorem}
\bp Let $U_{\infty}$ be a small open neighborhood of $\infty$ in $\manifm$. Then the local models given by a composition of the local models of Theorems~\ref{thmcompdiaggen}
and \ref{thmcompuptd}  
make clear that ${C}^+_{\finsetv}\!\left[U_{\infty},\infty\right]$ has a canonical smooth structure and that ${C}^+_{\finsetv}\!\left[U_{\infty},\infty\right]$ is the closure in $C_{V^+}[U_{\infty}]$ of $\check{C}_{\finsetv}\!\left[U_{\infty} \setminus \{\infty\}\right] \times \{\infty\}^{\{v_{\infty}\}}$.

Let $U$ be an open subset of $\manifm$ disjoint from $U_{\infty}$. Let $A$ be a subset of $V$.
Consider the map $\pbl\vert_{V}\colon {C}^+_{\finsetv}\!\left[\manifm,\infty\right] \to \manifm^{V}$.
Observe \begin{equation*}\bigl(\pbl\vert_{V}\bigr)^{-1}\bigl(U_{\infty}^{A^+} \times U^{V \setminus A}\bigr) ={C}^+_{\finseta}\!\left[U_{\infty},\infty\right] \times C_{V \setminus A}[U].\end{equation*}
The preimage of $\bigl(U_{\infty}^{A} \times U^{V \setminus A}\bigr)$ under $\pbl \colon {C}_{\finsetv}\!\left[\manifm,\infty\right] \to \manifm^{V}$ is similarly canonically diffeomorphic to ${C}_{A}\!\left[U_{\infty},\infty\right] \times  C_{V \setminus A}[U]$. 

So we are left with the identification of a small open neighborhood of $\bigl(\pbl\vert_{V}\bigr)^{-1}(\infty^V)$ in ${C}^+_{\finsetv}\!\left[U_{\infty},\infty\right]$ and a small open neighborhood of $\pbl^{-1}(\infty^V)$ in ${C}_{V}\!\left[U_{\infty},\infty\right]$ for an arbitrarily small neighborhood
$U_{\infty}$ of $\infty$ and a finite set $V$. Assume without loss that $U_{\infty}$ is identified with $\RR^{\dimdel}$ by a diffeomorphism $\varphi \colon \RR^{\dimdel} \to U_{\infty}$ such that $\varphi(0)=\infty$.
The naturality of the constructions leaves us with the case $(U_{\infty},\infty)=(\RR^{\dimdel},0)$.

Equip $\finsetv^+$ with the basepoint $v_{\infty}$.
The first blow-up along the small diagonal of $\finsetv^+$ transforms
$(\RR^{\dimdel})^{\finsetv^+}$ to 
\begin{equation*}C_{\finsetv^+,\card{\finsetv^+}}\!\left[\RR^{\dimdel}\right] \cong
(\RR^{\dimdel})^{\{v_{\infty}\}} \times \RR^+ \times \cuptd{\finsetv^+}{\RR^{\dimdel}},
\end{equation*}
where the space of $\cuptd{\finsetv^+}{\RR^{\dimdel}}$ of nonconstant maps from $V^+$ to 
$\RR^{\dimdel}$ up to translation and dilation is naturally identified with the space $\sph\bigl((\RR^{\dimdel})^{\finsetv}\bigr)$ of nonzero maps from $\finsetv$ to 
$\RR^{\dimdel}$ up to dilation. The first factor $(\RR^{\dimdel})^{\{v_{\infty}\}}$ contains the value of the configuration at $v_{\infty}$. It is zero on $\pbl\bigl({C}^+_{\finsetv}\!\left[\RR^{\dimdel},0\right]\bigr)$.
Since $\sph\bigl((\RR^{\dimdel})^{\finsetv}\bigr)$ is the unit normal bundle of $\{0\}^{\finsetv}$ in $\manifm^{\finsetv}$, we get a canonical diffeomorphism from
$C_{\finsetv^+,\card{\finsetv^+}}\!\left[\RR^{\dimdel}\right] \cap \bigl(\pbl\vert_{\{v_{\infty}\}}\bigr)^{-1}(0)$ to ${C}_{\finsetv,\card{\finsetv}+1}\!\left[\RR^{\dimdel},0\right]$.

Proceed by induction to get a canonical diffeomorphism from $C_{\finsetv^+,k}\!\left[\RR^{\dimdel}\right] \cap \bigl(\pbl\vert_{\{v_{\infty}\}}\bigr)^{-1}(0)$ to ${C}_{\finsetv,k}\!\left[\RR^{\dimdel},0\right]$ for all $k$.
\eop

Theorem~\ref{thmcompconftwo} is a direct consequence of Proposition~\ref{proprest} and Theorem~\ref{thmcompconfbis}. (Recall ${C}_{\finsetv}(\rats)={C}_{\finsetv}\!\left[\rats,\infty\right]$.) 
\eopwobp

An element $\confx$ of ${C}^+_{\finsetv}\!\left[\manifm,\infty\right]$ induces a parenthesization $\parentp^+(\confx)=\parentp(\confx \in C_{V^+}[\manifm])$ of $\finsetv^+$  as before Proposition~\ref{propdescconftautwo}.
By Theorems~\ref{thmcompuptd} 
and \ref{thmcompdiaggen} 
the element $\confx$ can be written as \begin{equation*}
\Bigl(\pbl(\confx)\vert_{V}, \bigl(w_{\finsetb} \in \cinjuptd{\kids(\finsetb)}{T_{\confx(\finsetb)}\manifm}\bigr)_{\finsetb \in \parentp^+(\confx)}\Bigr).
\end{equation*}
Let $\parentp^+_s(\confx)$ denote the set of elements of $\parentp^+(x)$ containing $v_{\infty}$.
This set is totally ordered by the inclusion. So is
\begin{equation*}
\parentp_s(\confx)= \{\finseta \setminus \{v_{\infty}\} \suchthat \finseta \in \parentp^+_s(\confx)\}=\{\finsetv(1), \finsetv(2), \dots, \finsetv(\sigma)\}\end{equation*}
with $\finsetv(i+1) \subset \finsetv(i)$.
We have $\confx(\finsetv(i))=\infty$ for all $i$. Let $\kids^s_d(\finsetv(i))$ denote the set of kids of $\finsetv(i)^+$ that do not contain $v_{\infty}$.

Recall from Notation~\ref{notsinjupdtcs} that $\sinjupdtcs(T_{\infty}\manifm,{\finseta})$ 
denotes the set of injective maps from $\finseta$ to $(T_{\infty}\manifm \setminus 0)$ up to dilation.
The natural basepoint choice of $v_{\infty}$ identifies the set $\cinjuptd{\kids(\finsetv(i)^+)}{T_{\infty}\manifm}$ of injective maps from 
$\kids(\finsetv(i)^+)$ to $T_{\infty}\manifm$ up to dilation and translation with $\sinjupdtcs\bigl(T_{\infty}\manifm,{\kids^s_d(\finsetv(i))}\bigr)$.

Set $\parentp_d(\confx)=\parentp^+(\confx) \setminus \parentp^+_s(\confx)$.
When $\confx$ is seen as an element of ${C}_{\finsetv}\!\left[\manifm,\infty\right]$, 
the set $\parentp_s(\confx)$ is the set of subsets $\finsetb$ of $\finsetv$ such that $\confx$ has been transformed by the blow-up along (the closure of the preimage of) $\einftyxi_{\finsetb}$
and $\parentp_d(\confx)$ is the set of subsets $\finsetb$ of $\finsetv$ such that $\confx$ has been transformed by the blow-up along (the closure of the preimage of) $\diag_{\finsetb}(\manifm^{\finsetv})$.

The proof of Theorem~\ref{thmcompconfbis} also proves the following proposition.
\begin{proposition}
\label{propstratinfty} Let $\parentp^+$ be a parenthesization of $\finsetv^+ = V \sqcup\{v_{\infty}\}$.
Set \begin{equation*}C_{\finsetv,\parentp^+}\left[\manifm,\infty\right]=\{x\in {C}^+_{\finsetv}\!\left[\manifm,\infty\right] \suchthat \parentp^+(x)=\parentp^+\}.\end{equation*}
Use the canonical identification of ${C}^+_{\finsetv}\!\left[\manifm,\infty\right]$ with ${C}_{\finsetv}\!\left[\manifm,\infty\right]$.
Then the stratum $C_{\finsetv,\parentp^+}\!\left[\manifm,\infty\right]$ is an open part of \begin{equation*}\partial_{\cardlef{\parentp^+}}\bigl(C_{\finsetv}\!\left[\manifm,\infty\right]\bigr)\setminus \partial_{\cardlef{\parentp^+} +1}\bigl(C_{\finsetv}\!\left[\manifm,\infty\right]\bigr).\end{equation*}

Let $\parentp^+_s$ be the set of elements of $\parentp^+$ containing $v_{\infty}$. 
Set $\parentp_d=\parentp^+ \setminus \parentp^+_s$ and
$\parentp_s= \{\finseta \setminus \{v_{\infty}\} \suchthat \finseta \in \parentp^+_s\}$.
For an element $\finsetb$ of $\parentp_s \cup \{\finsetv\}$, let $\kids^s_d(\finsetb)$ denote the set of kids of $\finsetb^+$ that do not contain $v_{\infty}$. 
The stratum $C_{\finsetv,\parentp^+}\left[\manifm,\infty\right]$ fibers 
 over the space $\check{C}_{\kids^s_d(\finsetv)}\!\left[\cmanifm\right]$ of injective maps $\confc$ from $\kids^s_d(\finsetv)$ to $\cmanifm$.
 Its fiber over such an injective map $\confc$ 
 is \begin{equation*}\Bigl(\prod_{\finsetb \in \parentp_s} \sinjupdtcs\bigl(T_{\infty}\manifm,\kids^s_d(\finsetb)\bigr)\Bigr) \times \Bigl(\prod_{\finseta \in \parentp_d} \check{S}_{\kids(\finseta)}\bigl(T_{\confc(\finseta)}\manifm \bigr)\Bigr) = \prod_{\finseta \in \parentp^+} \check{S}_{\kids(\finseta)}\bigl(T_{\confc(\finseta)}\manifm \bigr).\end{equation*}
\end{proposition}

\begin{corollary}
 \label{cordescconftaucomp}
An element of $C_{\finsetv}\!\left[\manifm,\infty\right]$ is a map $\confc$ from ${\finsetv}$ to $\manifm$, equipped with
\begin{itemize} 
\item a parenthesization $\parentp+$ of $V^+ = V \sqcup\{v_{\infty}\}$ and induced parenthesizations $\parentp_s$ and $\parentp_d$ of $V$, as in Proposition~\ref{propstratinfty},
\item an element $\ffuncspec_{\finsetb}$ of $\sinjupdtcs\left(T_{\infty}\manifm,{\kids^s_d\bigl(\finsetb\bigr)}\right)$ for any element $\finsetb$ of $\parentp_s$, 
\item an element $w_{\finseta}$ of $\check{S}_{\kids(\finseta)}\bigl(T_{\confc(\finseta)}\manifm \bigr)$ for each $\finseta \in \parentp_d$.
\end{itemize}
\end{corollary}

\section{Finishing the proofs of the statements of Sections~\ref{secfirstprescomp} and \ref{secfacecodimone}}
\label{secproofblowup}

\bpo{Proof of Proposition~\ref{propcompconfL}}
To study the closure of $\check{C}(\rats,\Link;\Gamma)$ in $C_{V(\Gamma)}(\rats)$,
we study its intersection with some $\pbl^{-1}(\prod_{i \in I} U_i^{{\finsetv}_i})$ for disjoint small compact $U_i$. We assume that at most one $U_i$ contains $\infty$ and that this $U_i$ does not meet the link. So the corresponding ${\finsetv}_i$ does not contain univalent vertices, and Theorem~\ref{thmcompconf} gives the structure of the corresponding factor.

Thus, it is enough to study $p_{\finseta}\bigl(\overline{\check{C}(\rats,\Link;\Gamma)} \bigr) \cap C_{\finseta}\bigl(\phi(\RR^3)\bigr)$
when 
\begin{itemize}
 \item 
 $\phi$ is an embedding from $\RR^3$ to $\crats$ that maps the vertical line $\RR\vec{v}$ through the origin oriented from bottom to top onto $\phi(\RR^3)\cap \Link$, so that $\phi$ identifies
$(\RR^3,\RR\vec{v})$ with $\bigl(\phi(\RR^3),\phi(\RR^3)\cap \Link\bigr)$,\footnote{Here, $\Link$ also denotes the image $\Link(\source)$ of $\Link$.}
\item 
the univalent vertices of ${\finseta}$ form a nonempty subset $\finseta_U=U(\Gamma) \cap {\finseta}$ 
of consecutive vertices on the component $\Link_i$ of $\Link$ such that $\phi(\RR^3)\cap \Link=\phi(\RR^3)\cap \Link_i$.
\end{itemize} 
Let $\CO(\finseta_U)$ denote the set of the linear orders $<$ on $\finseta_U$ compatible with $\left[i_{\Gamma}\right]$. (This is a singleton unless $\finseta_U$ contains all the univalent vertices of $\Link_i$.)
Via the natural maps induced by $\phi$, the space $\check{C}_{\finseta}(\phi(\RR^3))$ is identified with  the set of injections from
${\finseta}$ to $\RR^3$, and
$p_{\finseta}\bigl(\check{C}(\rats,\Link;\Gamma) \bigr) \cap \check{C}_{\finseta}\bigl(\phi(\RR^3)\bigr)$
is identified with the disjoint union over $\CO(\finseta_U)$ of the subsets $\check{C}_{\finseta}\bigl(\RR^3,\finseta_U,<\; \in \CO(\finseta_U)\bigr)$ of injections that map $\finseta_U$ to $\RR\vec{v}$ so that the order induced by $\RR\vec{v}$ coincides with $<$. Fix $< \; \in \CO(\finseta_U)$, and write $\finseta_U=\{v_1,\dots,v_k\}$ so that $v_1<v_2\dots<v_k$.

Fix $b(\finseta)=v_1$, and study the closure of $\check{C}_{\finseta}(\RR^3,\finseta_U,<)$ in $\RR^3 \times \RR^+ \times \ccompuptd{\finseta}{\RR^3}$, using the diffeomorphism $\psi(\phi,v_1)$ of Theorem~\ref{thmcompdiaggen}.
This closure is the closure of $\RR\vec{v} \times \RR^+ \times \cinjuptd{\finseta}{\RR^3,\finseta_U,<}$, where $\cinjuptd{\finseta}{\RR^3,\finseta_U,<}$ is the quotient of $\check{C}_{\finseta}(\RR^3,\finseta_U,<)$ by translations by vectors of $\RR\vec{v}$ and dilations.
Then the charts of Theorem~\ref{thmcompuptd} (used with basepoints that are as much as possible in $\finseta_U$) make clear that the closure of $\cinjuptd{\finseta}{\RR^3,\finseta_U,<}$ in $\ccompuptd{\finseta}{\RR^3}$ consists of the limit configurations $c$ such that $\left(c(v_j)-c(v_i)\right)$ is nonnegatively colinear with $\vec{v}$ at any scale (i.e., in any infinitesimal configuration $w$ that has popped up during the blow-ups) for any $i$ and $j$ in $\{1, \dots, k\}$ such that $i<j$, and that this closure is a smooth submanifold of $\ccompuptd{\finseta}{\RR^3}$ transverse to the ridges.
\eop

This proof also proves the following lemma.

\begin{lemma}
The codimension-one faces of ${C}(\rats,\Link;\Gamma)$ are the intersections of ${C}(\rats,\Link;\Gamma)$ with the codimension-one faces of ${C}_{V(\Gamma)}(\rats)$.
\end{lemma}
\eopwobp

Proposition~\ref{propconffaceun} then follows, with the help of Propositions~\ref{proprest} 
and \ref{propstratinfty}.
\eopwobp

\section{Alternative descriptions of configuration spaces}
\label{secfsttwoblow}

Apart from Lemma~\ref{lemdescclosccompuptd}, this section will not be used in this book. It mentions other presentations of the configuration spaces studied in Sections~\ref{secprooffacanom} to \ref{secbloinf} without proofs. Most of the proofs are left to the reader as exercises.

Let $\finsetv$ be a finite set of cardinality at least $2$. Let $\CP_{\geq 2}=\CP_{\geq 2}(\finsetv)$\index[N]{Pzzgreater2@$\CP_{\geq 2}(\finsetv)$} be the set of its (nonstrict) subsets of cardinality at least $2$.
The smooth blowdown projection from $\ccompuptd{\finseta}{\vecspt}$ to $\cuptd{\finseta}{\vecspt}$ for an $\finseta \in \CP_{\geq 2}(\finsetv)$ may be composed with the smooth restriction map from $\ccompuptd{\finsetv}{\vecspt}$ to $\ccompuptd{\finseta}{\vecspt}$ to produce a smooth map $\pi_{\finseta}$ from $\ccompuptd{\finsetv}{\vecspt}$ to $\cuptd{\finseta}{\vecspt}$.

\begin{lemma}
\label{lemdescclosccompuptd}
The product over the subsets $\finseta$ of $\finsetv$ with cardinality at least $2$ of the $\pi_{\finseta}$ is a smooth embedding 
\begin{equation*}\prod_{\finseta \in \CP_{\geq 2}}\pi_{\finseta} \colon \ccompuptd{\finsetv}{\vecspt} \hookrightarrow \prod_{\finseta \in \CP_{\geq 2}}\cuptd{\finseta}{\vecspt}\end{equation*}
and the image of $\ccompuptd{\finsetv}{\vecspt}$ is the closure of the image of the restriction of $\prod_{\finseta \in \CP_{\geq 2}}\pi_{\finseta}$ to $\cinjuptd{\finsetv}{\vecspt}$.
\end{lemma}
\bp The injectivity of $\prod_{\finseta \in \CP_{\geq 2}}\pi_{\finseta}$ can be seen from the description of $\ccompuptd{\finsetv}{\vecspt}$ as a set, which is given in Theorem~\ref{thmcompuptd}.
\eop

\begin{proposition}
 The image $\bigl(\prod_{\finseta \in \CP_{\geq 2}}\pi_{\finseta}\bigr)\bigl(\ccompuptd{\finsetv}{\vecspt}\bigr)$
is the subset of \begin{equation*}\prod_{\finseta \in \CP_{\geq 2}}\cuptd{\finseta}{\vecspt}\end{equation*} consisting of the elements $((c_{\finseta})_{\finseta \in \CP_{\geq 2}})$ such that for any two elements $\finseta$ and $\finsetb$ of $\CP_{\geq 2}$ such that $\finsetb \subset \finseta$, the restriction of $c_{\finseta}$ to $\finsetb$ coincides with $c_{\finsetb}$ if it is not constant.
\end{proposition}
\bp Exercise. \eop

Thus, $\ccompuptd{\finsetv}{\vecspt}$ can be defined as its image described in the above proposition.
Similar definitions involving only cardinality $2$ or $3$ subsets of $\finsetv$ can be found in \cite{Sinha}.

For a smooth manifold $\manifm$ without boundary, we have similar smooth maps $\pi_{\finseta}$ from $C_{\finsetv}\!\left[\manifm\right]$ to $\blowup{\manifm^{\finseta}}{\diag_{\finseta}(\manifm^{\finseta})}$.
They also define a smooth map
\begin{equation*}\prod_{\finseta \in \CP_{\geq 2}}\pi_{\finseta} \colon C_{\finsetv}\!\left[\manifm\right] \hookrightarrow \prod_{\finseta \in \CP_{\geq 2}}\lefblowup{\manifm^{\finseta}}{\diag_{\finseta}(\manifm^{\finseta})}.\end{equation*}

The elements of $\blowup{\manifm^{\finseta}}{\diag_{\finseta}(\manifm^{\finseta})}$ are maps $\confc$ from $\finseta$ to $\manifm$ that are equipped with an element $w \in \cuptd{\finseta}{T_{\confc(\finseta)}\manifm}$ when they are constant.

\begin{proposition}
 The map \begin{equation*}\prod_{\finseta \in \CP_{\geq 2}}\pi_{\finseta} \colon C_{\finsetv}\!\left[\manifm\right] \hookrightarrow \prod_{\finseta \in \CP_{\geq 2}}\lefblowup{\manifm^{\finseta}}{\diag_{\finseta}(\manifm^{\finseta})}\end{equation*} is an embedding.
Its image is the subset of $\prod_{\finseta \in \CP_{\geq 2}}\blowup{\manifm^{\finseta}}{\diag_{\finseta}(\manifm^{\finseta})}$ consisting of the elements $((c_{\finseta})_{\finseta \in \CP_{\geq 2}})$ such that
 for any two elements $\finseta$ and $\finsetb$ of $\CP_{\geq 2}$ such that $\finsetb \subset \finseta$,
\begin{itemize}
 \item the restriction to $\finsetb$ of the map $\pbl(c_{\finseta})\colon \finseta \to \manifm$ coincides with $\pbl(c_{\finsetb})$, and,
\item if $\pbl(c_{\finseta})$ is constant, then the restriction to $\finsetb$ of $\bigl(w_{\finseta}(c_{\finseta}) \in \cuptd{\finseta}{T_{\confc_{\finseta}(\finseta)}\manifm}\bigr)$ is $w_{\finsetb}(c_{\finsetb})$ if this restriction is not constant.
\end{itemize}
\end{proposition}
\bp Exercise. \eop

Again, $C_{\finsetv}\!\left[\manifm\right]$ can be defined as its image described in the above proposition, and
similar definitions involving only cardinality $2$ and $3$ subsets of $\finsetv$ can be found in \cite{Sinha}.

We may obtain similar statements for $C_{\finsetv}\!\left[\manifm,\infty\right]$, where $\infty \in \manifm$, from the definition of $C_{\finsetv}\!\left[\manifm,\infty\right]$ as ${C}^+_{\finsetv}\!\left[\manifm,\infty\right]$ in Theorem~\ref{thmcompconfbis}.

More information about the homotopy groups and the homology of the configuration spaces $\check{C}_{\finsetv}\!\left[\RR^d\right]$ and $\check{C}_{\finsetv}\!\left[S^d\right]$ can be found in the book \cite{FadellHusseinibook} by Edward Fadell and Sufian Husseini.
See also the book \cite{CamposLambrechts} by Ricardo Campos, Najib Idrissi, Pascal Lambrechts, and Thomas Willwacher.

\chapter{Dependence on the propagating forms}
\label{chapindepform}

This chapter shows how our combinations of integrals over configuration spaces depend on the chosen propagating forms.

\section{Introduction}
\label{secintroindepform}

In this section, we give a first general description of the variation of $\Zinv$ when propagating forms change in Proposition~\ref{propinvone}.
Then we show how this proposition and a preliminary lemma (\ref{lemformprod}) apply to prove Theorem~\ref{thmevenz} and two other lemmas (\ref{leminvfstconst} and \ref{leminvwoboun}), about independence of chosen propagating forms as in Definition~\ref{defpropagatortwo}. The typical proof of Proposition~\ref{propinvone} will occupy the next sections.

Again, any closed $2$-form on $\partial C_2(\rats)$ extends to $C_2(\rats)$ because the restriction induces a surjective map $H^2(C_2(\rats);\RR) \rightarrow H^2(\partial C_2(\rats);\RR)$ since we have 
\begin{equation*}H^3\bigl(C_2(\rats),\partial C_2(\rats);\RR\bigr)=0.\end{equation*}

\begin{lemma}
\label{lemformprod}
Let $(\crats,\partau_0)$ be an asymptotic rational homology $\RR^3$ as in Definition~\ref{defasyrathommRthree}. 
Let $\partau \colon \left[0,1\right] \times \crats  \times \RR^3 \to T \crats$ be a smooth map whose restriction 
to $\{t\} \times  \crats  \times \RR^3$ is an asymptotically standard parallelization $\partau_t$ of $\crats$ for any $t \in \left[0,1\right]$.
Define $\projp_{\partau} \colon \left[0,1\right] \times \partial C_2(\rats) \to \left[0,1\right] \times S^2$ by $\projp_{\partau}(t,x)=\bigl(t,p_{\partau_t}(x)\bigr)$.

Let $\omega_0$ and $\omega_1$ be two propagating forms of $C_2(\rats)$ that restrict to 
$\partial C_2(\rats) \setminus \ST B_{\rats}$ as $p_{\partau_0}^{\ast}(\omega_{0,S^2})$ and as $p_{\partau_1}^{\ast}(\omega_{1,S^2})$, respectively, for two volume-one forms $\omega_{0,S^2}$ and $\omega_{1,S^2}$ of $S^2$. 

Then there exist \begin{itemize}
\item a closed $2$-form $\tilde{\omega}_{S^2}$ on $\left[0,1\right] \times S^2$ whose restriction to $\{t\} \times S^2$ is $\omega_{t,S^2}$ for $t \in \{0,1\}$,
\item for any such $\tilde{\omega}_{S^2}$, a closed $2$-form $\omega^{\partial}$ on $\left[0,1\right] \times \partial C_2(\rats)$ whose restriction to $\{t\} \times \partial C_2(\rats)$ is $\omega_{t}\vert_{\partial C_2(\rats)}$ for $t \in \{0,1\}$,
and whose restriction to $\left[0,1\right] \times (\partial C_2(\rats) \setminus \ST(\ballb_{\rats}))$ is $\projp_{\partau}^{\ast}(\tilde{\omega}_{S^2})$, and,
\item for any such compatible $\tilde{\omega}_{S^2}$ and $\omega^{\partial}$, a closed $2$-form ${\omega}$ on $\left[0,1\right] \times C_2(\rats)$ whose restriction to $\{t\} \times C_2(\rats)$ is $\omega_t$ for $t \in \{0,1\}$,
and whose restriction to $\left[0,1\right] \times \partial C_2(\rats)$ is $\omega^{\partial}$.
\end{itemize}
If $\omega_0$ and $\omega_1$ are propagating forms of $\bigl(C_2(\rats),\partau_0\bigr)$ and $\bigl(C_2(\rats),\partau_1\bigr)$, then
we may choose ${\omega}^{\partial}=\projp_{\partau}^{\ast}(\tilde{\omega}_{S^2})$ on $\left[0,1\right] \times \partial C_2(\rats)$.
\end{lemma}
\bp As in Lemma~\ref{lemetactwo}, there exists a one-form $\eta_{S^2}$ on $S^2$ such that
$d\eta_{S^2}=\omega_{1,S^2}-\omega_{0,S^2}$. Define the closed $2$-form $\tilde{\omega}_{S^2}$ on $\left[0,1\right] \times S^2$ by 
\begin{equation*}\tilde{\omega}_{S^2} = p_{S^2}^{\ast}({\omega}_{0,S^2}) + d \bigl(t p_{S^2}^{\ast}(\eta_{S^2})\bigr),\end{equation*}
where $t$ is the coordinate on $\left[0,1\right]$.

Now, the form $\omega^{\partial}$ is defined on the boundary of $\left[0,1\right] \times \ST(\ballb_{\rats})$, and it extends as a closed $2$-form $\omega^{\partial}$ as desired there because the restriction induces a surjective map $H^2(\left[0,1\right] \times \ST(\ballb_{\rats});\RR) \rightarrow H^2\left(\partial \left(\left[0,1\right] \times \ST(\ballb_{\rats})\right);\RR\right)$ since we have \begin{equation*}H^3\Bigl(\left[0,1\right] \times \ST(\ballb_{\rats}),\partial \bigl(\left[0,1\right] \times \ST(\ballb_{\rats})\bigr);\RR\Bigr)=0.\end{equation*}

Finally, the desired form ${\omega}$ is defined on the boundary of $\left[0,1\right] \times C_2(\rats)$. It similarly extends as a closed $2$-form to $\left[0,1\right] \times C_2(\rats)$.
\eop

When ${\finseta}$ is a subset of the set of vertices ${\finsetv(\Gamma)}$ of a numbered Jacobi diagram $\Gamma$ with support a one-manifold $\source$,
the set of edges of $\Gamma$ between two elements of ${\finseta}$ is denoted by $E(\Gamma_{\finseta})$ (edges of $\Gamma$ are plain), and $\Gamma_{\finseta}$ is the subgraph of $\Gamma$ consisting of the vertices of ${\finseta}$ and the edges of $E(\Gamma_{\finseta})$ together with the natural restriction to $U(\Gamma) \cap \finseta$ of the isotopy class of injections from $U(\Gamma)$ to $\source$ associated to $\Gamma$.\index[N]{Gamma{\finseta}@$\Gamma_{\finseta}$ subgraph of $\Gamma$}

The following proposition, whose proof occupies most of this chapter, is crucial in the study of the variations of $\Zinv$.

\begin{proposition}
\label{propinvone}
Let $(\crats,\partau)$ be an asymptotic rational homology $\RR^3$. Let $\Link=\sqcup_{j=1}^kK_j$ be an embedding of $\source=\sqcup_{j=1}^kS^1_j$ into $\crats$.
Let $\partau \colon \left[0,1\right] \times \crats  \times \RR^3 \to T \crats$ be a smooth map whose restriction 
to $\{t\} \times  \crats  \times \RR^3$ is an asymptotically standard parallelization $\partau_t$ of $\crats$ for any $t \in \left[0,1\right]$.\footnote{This homotopy $\partau$ is not useful for this statement, but we will use this notation later.}
Define $\projp_{\partau} \colon \left[0,1\right] \times \partial C_2(\rats) \to \left[0,1\right] \times S^2$ by $\projp_{\partau}(t,x)=(t,p_{\partau_t}(x))$. 

Let $n$ and $m$ be positive integers.
For $i \in \underline{m}$, let $\tilde{\omega}(i)$ be a closed $2$-form on $\left[0,1\right] \times C_2(\rats)$ whose restriction to $\{t\} \times C_2(\rats)$ is denoted by $\tilde{\omega}(i,t)$, for any $t \in \left[0,1\right]$.

Assume that $\tilde{\omega}(i)$ restricts to $\left[0,1\right] \times\bigl(\partial C_2(\rats) \setminus \ST {\ballb}_{\rats}\bigr)$ as $\projp_{\partau}^{\ast}(\tilde{\omega}_{S^2}(i))$, for some closed two-form $\tilde{\omega}_{S^2}(i)$ on $\left[0,1\right] \times S^2$ such that $\int_{\{t\} \times S^2}\tilde{\omega}_{S^2}(i)=1$ for $t \in \left[0,1\right]$.
Recall Notation~\ref{notlessav} and set
\begin{equation*}\Zinv_{n,m}(t)=\sum_{\Gamma \in \Davisred^e_{n,m}(\source)}\coefgambetred_{\Gamma}I\Bigl(\rats,\Link,\Gamma,\bigl(\tilde{\omega}(i,t)\bigr)_{i \in \underline{m}}\Bigr)\left[\Gamma\right] \in \Aavis_n(\sqcup_{j=1}^kS^1_j)\end{equation*}
 and $\Zinv_n(t)=\Zinv_{n,3n}(t)$.
Let $\Davis^{e,F}_{n,m}(\source)$ denote the set of pairs $(\Gamma,{\finseta})$
such that 
\begin{itemize}
 \item $\Gamma \in \Davisred^e_{n,m}(\source)$,
 \item ${\finseta} \subseteq {\finsetv(\Gamma)}$, $\cardlef{\finseta} \geq 2$,
 \item $\Gamma_{\finseta}$ is a connected component of $\Gamma$,
 \item $\cardlef{\finseta} \equiv 2$ mod $4$ if $\finseta \cap U(\Gamma)=\emptyset$, and 
 \item $\finseta \cap U(\Gamma)$
is a set of consecutive vertices on a component $\source_{\finseta}$ of $\source$ if $\finseta \cap U(\Gamma)\neq \emptyset$.
\end{itemize}
For $(\Gamma,{\finseta}) \in \Davis^{e,F}_{n,m}(\source)$,
set \begin{equation*}I(\Gamma,\finseta)= \int_{\left[0,1\right]\times \facee(\finseta,\Link,\Gamma)}\bigwedge_{e \in E(\Gamma)}p_e^{\ast}\Bigl(\tilde{\omega}\bigl(j_E(e)\bigr)\Bigr)\left[\Gamma\right],\end{equation*}
\index[N]{Integrals over configuration spaces!IRLGammaA@$I(\Gamma,\finseta)$}
where $p_e \colon \left[0,1\right] \times {C}(\rats,\Link;\Gamma) \to \left[0,1\right] \times C_2({\rats})$ is the product by the identity map $\id_{\left[0,1\right]}$ of $\left[0,1\right]$ of the previous $p_e$, $j_E$ is the edge-numbering map of Definition~\ref{defnumdia}, and the face $\facee(\finseta,\Link,\Gamma)$ of ${C}(\rats,\Link;\Gamma)$ is described in Section~\ref{secfacecodimone}. Set $\Davis^{e,F}_n(\source)=\Davis^{e,F}_{n,3n}(\source)$.
Then we have \begin{equation*}\Zinv_n(1)-\Zinv_n(0)=\sum_{ (\Gamma,{\finseta}) \in \Davis^{e,F}_n(\source)}\coefgambet_{\Gamma} I(\Gamma,\finseta)\end{equation*}
and 
\begin{equation*}\Zinv_{n,3n-2}(1)-\Zinv_{n,3n-2}(0)=\sum_{ (\Gamma,{\finseta}) \in \Davis^{e,F}_{n,3n-2}(\source)}\coefgambetred_{\Gamma} I(\Gamma,\finseta).\end{equation*}

\end{proposition}

This statement simplifies as in Corollary~\ref{corinvone}
when $L=\emptyset$ using the projection 
$p^c \colon \Aavis(\emptyset) \to \Aavis^c(\emptyset)$, which maps diagrams with several connected components to $0$. Recall that $\Davis^{c}_n$ is the subset of $\Davis^{e}_n(\emptyset)$ whose elements are the numbered diagrams of $\Davis^{e}_n(\emptyset)$ with one connected component.

For an oriented connected trivalent diagram $\Gamma$, the face $\facee\bigl(\finsetv(\Gamma),\emptyset,\Gamma\bigr)$ fibers over $\crats$, and the fiber over $x \in \crats$ is the space $\cinjuptd{\finsetv(\Gamma)}{T_x\crats}$ of injections from $\finsetv(\Gamma)$ to $T_x\crats$, up to translation and dilation. See Section~\ref{secpresfacanom}.
We also denote this face by $\cinjuptd{\finsetv(\Gamma)}{T\crats}$. 

\begin{lemma}
\label{lemoriface}
Let $\Gamma \in \Davis^c_n$ be equipped with a vertex-orientation, which induces an orientation of ${C}(\rats,\emptyset;\Gamma)$ as in Corollary~\ref{cororc}.
These orientations induce the orientation of $\finsetv(\Gamma)$ described in Remark~\ref{rkedgeorvert}.
The orientation of $\facee\bigl(\finsetv(\Gamma),\emptyset,\Gamma\bigr)$ as part of the boundary of ${C}(\rats,\emptyset;\Gamma)$ can be alternatively described as follows.
The face $\facee\bigl(\finsetv(\Gamma),\emptyset,\Gamma\bigr)$ is oriented as the local product $\crats \times \mbox{fiber}$, where the fiber is oriented as in Convention~\ref{convordefcuptdfinsetvvecspt}, using the above orientation of $\finsetv(\Gamma)$.
\end{lemma}
\bp The dilation factor for the quotient $\cinjuptd{\finsetv(\Gamma)}{T_x\crats}$ plays the role of an inward normal for ${C}(\rats,\emptyset;\Gamma)$. The orientation of ${C}(\rats,\emptyset;\Gamma)$ near the boundary is given by the orientation of $\crats$, followed by this inward normal, followed by the fiber orientation.
\eop

For any pair $e$ of $\finsetv(\Gamma)$, we have a natural restriction map \begin{equation*}p_e \colon \cinjuptd{\finsetv(\Gamma)}{T\crats} \to \cinjuptd{e}{T\crats} \cong \ST \crats,\end{equation*}
which provides natural restriction maps 
\begin{equation*}p_e \colon \left[0,1\right] \times \cinjuptd{\finsetv(\Gamma)}{T\crats} \to \left[0,1\right] \times \cinjuptd{e}{T\crats}\end{equation*}
by multiplication by $\id_{\left[0,1\right]}$.

Proposition~\ref{propinvone} has the following corollary.

\begin{corollary}
\label{corinvone} Assume $\Link=\emptyset$. Recall  Notation~\ref{notationzZ}.
Under the assumptions of Proposition~\ref{propinvone}, set $\zinv_n(t)=p^c\left(\Zinv_n(t)\right)$ so that we have
\begin{equation*}\zinv_n(t)=\zinv_n\left(\rats,\bigl(\tilde{\omega}(i,t)\bigr)_{i \in \underline{3n}}\right).\end{equation*}

Then $\zinv_n(1)=\zinv_n(0)$ if $n$ is even, and
$(\zinv_n(1)-\zinv_n(0))$ depends only on $\left(\tilde{\omega}(i)\vert_{\left[0,1\right] \times \ST {\ballb}_{\rats}}\right)_{i \in \underline{3n}}$ if $n$ is odd.
Define
\begin{equation*}\zinv_n\left(\left[0,1\right] \times \ST {\ballb}_{\rats}; \bigl(\tilde{\omega}(i)\bigr)_{i \in \underline{3n}}\right)=\sum_{\Gamma \in \Davis^{c}_n}\coefgambet_{\Gamma} I\bigl(\Gamma,V(\Gamma)\bigr),
\end{equation*}
where
\begin{equation*}I\bigl(\Gamma,V(\Gamma)\bigr)=\int_{\left[0,1\right]\times \cinjuptd{\finsetv(\Gamma)}{T{\ballb}_{\rats}}}\bigwedge_{e \in E(\Gamma)}p_e^{\ast}\Bigl(\tilde{\omega}\bigl(j_E(e)\bigr)\Bigr)\left[\Gamma\right].\end{equation*}
Then we have
\begin{equation*}\zinv_n(1)-\zinv_n(0)=\zinv_n\left(\left[0,1\right] \times \ST {\ballb}_{\rats}; \bigl(\tilde{\omega}(i)\bigr)_{i \in \underline{3n}}\right)\end{equation*}
for any odd integer $n$.
\end{corollary}
\bp We have \begin{equation*}\int_{\left[0,1\right]\times \cinjuptd{\finsetv(\Gamma)}{T(\crats \setminus {\ballb}_{\rats})}}\bigwedge_{e \in E(\Gamma)}p_e^{\ast}\Bigl(\tilde{\omega}\bigl(j_E(e)\bigr)\Bigr)=0\end{equation*}
because the integrated form factors through $\left[0,1\right]\times \cinjuptd{\finsetv(\Gamma)}{\RR^3}$ via a map induced by $p_{\tau}$, which is fixed and independent of $\tau$, there.
In particular, the $I\bigl(\Gamma,V(\Gamma)\bigr)$ depend only on $\left(\tilde{\omega}(i)\vert_{\left[0,1\right] \times \ST {\ballb}_{\rats}}\right)_{i \in \underline{3n}}$.
\eop

\bpo{Proof of Theorem~\ref{thmevenz} assuming Proposition~\ref{propinvone}}
Changing propagating forms $\omega(i)_0$ of $C_2(\rats)$ to other ones $\omega(i)_1$ provides
forms $\tilde{\omega}(i)$ on $\left[0,1\right] \times C_2(\rats)$ as in Lemma~\ref{lemformprod}. Then Corollary~\ref{corinvone} guarantees that $\zinv_{2n}\bigl(\crats,\emptyset,(\omega(i))\bigr)$ does not depend on the used propagating forms (which are not normalized on $\ST(\ballb_{\rats})$ and hence do not depend on parallelizations).
\eop

\begin{lemma}
\label{leminvfstconst}
Let $(\crats,\partau)$ be an asymptotic rational homology $\RR^3$. 
Let $\Link \colon \source \hookrightarrow \crats$ be a link embedding.
For any $i\in \underline{3n}$, let $\omega(i)$ be a homogeneous propagating form of $\bigl(C_2(\rats),\partau\bigr)$.
Then, as stated in Theorem~\ref{thmfstconst}, $\Zinv_n\bigl(\crats,\Link,(\omega(i))\bigr)$ is independent of the chosen
$\omega(i)$.\footnote{The proof of Theorem~\ref{thmfstconst} will be concluded in the end of Section~\ref{secdephom}.} Denote it by $\Zinv_n(\crats,\Link,\partau)$.
\end{lemma}
\bpo{Proof assuming Proposition~\ref{propinvone}}
To prove this lemma, it suffices to prove that if some homogeneous $\omega(i)=\tilde{\omega}(i,0)$ is changed to another homogeneous propagating form
$\tilde{\omega}(i,1)$, then $\Zinv_n\bigl(\crats,\Link,(\omega(i))\bigr)$ is unchanged.
According to Lemma~\ref{lemetactwo}, under these assumptions, there exists a one-form $\eta$ on $C_2(\rats)$ such that 
\begin{itemize}
\item $\tilde{\omega}(i,1) = \tilde{\omega}(i,0) +d \eta$ and
\item $\eta\vert_{\partial C_2({\rats})}=0$.
\end{itemize}

Let $p_{C_2} \colon \left[0,1\right] \times C_2(\rats) \to C_2(\rats)$ denote the projection on the second factor.
Define closed $2$-forms $\tilde{\omega}(j)$ on $\left[0,1\right] \times C_2(\rats)$ by
\begin{itemize}
 \item $\tilde{\omega}(j) = p_{C_2}^{\ast}\bigl(\omega(j)\bigr)$ if $j \neq i$, and
\item $\tilde{\omega}(i) = p_{C_2}^{\ast}\bigl(\tilde{\omega}(i,0)\bigr) + d \bigl(t p_{C_2}^{\ast}(\eta)\bigr).$
\end{itemize}
Then the variation of $\Zinv_n\bigl(\crats,\Link,(\omega(i))\bigr)$ is $\Zinv_n(1)-\Zinv_n(0)$, with the notation of Proposition~\ref{propinvone}, where all the forms involved in some $I(\Gamma,\finseta)$, except $p_{e(i)}^{\ast}(\tilde{\omega}(i))$, for the possible edge $e(i)$ such that $j_E\bigl(e(i)\bigr)=i$, factor through $p_{C_2}^{\ast}$. Thus, if $i \notin \Image(j_E)$, then all the forms factor through $p_{C_2}^{\ast}$, and $I(\Gamma,\finseta)$ vanishes.

Locally, $\facee(\finseta,\Link,\Gamma)$ is diffeomorphic to the product of $\facee(\finseta,\Link,\Gamma_{\finseta})$ by ${C}(\rats,\Link;\Gamma \setminus \Gamma_{\finseta})$.
If $e(i)$ is not an edge of $\Gamma_{\finseta}$, then the form $\bigwedge_{e \in E(\Gamma_{\finseta})}p_e^{\ast}\bigl(\tilde{\omega}(j_E(e))\bigr)$ factors through $\facee(\finseta,\Link,\Gamma_{\finseta})$ whose dimension is $2 \cardlef{E(\Gamma_{\finseta})} -1$, and the form vanishes.
If $e(i)$ is an edge of $\Gamma_{\finseta}$, then the part $p_{e(i)}^{\ast}\bigl(d (t p_{C_2}^{\ast}(\eta))\bigr)$ vanishes since $\eta$ vanishes on $\partial C_2({\rats})$.
Thus $\bigwedge_{e \in E(\Gamma_{\finseta})}p_e^{\ast}\bigl(\omega(j_E(e))\bigr)$ still factors through $\facee(\finseta,\Link,\Gamma_{\finseta})$.
So, when $(\rats,\Link,\partau)$ is fixed, 
$\Zinv_n\bigl(\crats,\Link,(\omega(i))\bigr)$ is independent of the chosen homogeneous $\omega(i)$.
\eop

\begin{lemma}
\label{leminvwoboun}
Let $(\crats,\partau_0)$ be an asymptotic rational homology $\RR^3$. Let $n \in \NN$.
For any $i\in \underline{3n}$, let $\omega(i)$ be a propagating form of $\bigl(C_2(\rats),\partau_0\bigr)$.
 Then $\Zinv_n\bigl(\crats,\emptyset,(\omega(i))\bigr)=\Zinv_n(\crats,\emptyset,\partau_0)$ with the notation of Lemma~\ref{leminvfstconst}. Furthermore, $\Zinv_n(\crats,\emptyset,\partau_0)$ depends only on the homotopy class of $\partau_0$.
\end{lemma}
\bpo{Proof assuming Proposition~\ref{propinvone}}
Let $\partau \colon \left[0,1\right] \times \crats  \times \RR^3 \to T \crats$ be a smooth map whose restriction 
to $\{t\} \times  \crats  \times \RR^3$ is an asymptotically standard parallelization $\partau_t$ of $\crats$ for any $t \in \left[0,1\right]$.
Define $\projp_{\partau} \colon \left[0,1\right] \times \partial C_2(\rats) \to \left[0,1\right] \times S^2$ by $\projp_{\partau}(t,x)=(t,p_{\partau_t}(x))$.
For any $i\in \underline{3n}$, let $\omega(i)_0$ be a (non-necessarily homogeneous) propagating form of $\bigl(C_2(\rats),\partau_0\bigr)$, and let $\omega(i)_1$ be a propagating form of $\bigl(C_2(\rats),\partau_1\bigr)$.
It suffices to prove that 
\begin{equation*}\Zinv_n\Bigl(\crats,\emptyset,\bigl(\omega(i)_1\bigr)\Bigr)-\Zinv_n\Bigl(\crats,\emptyset,\bigl(\omega(i)_0\bigr)\Bigr)=0.\end{equation*}
Use forms ${\omega}(i)$ on $\left[0,1\right] \times C_2(\rats)$ provided by Lemma~\ref{lemformprod},
which restrict to $\left[0,1\right] \times \partial C_2(\rats)$ as $\projp_{\partau}^{\ast}(\tilde{\omega}_{S^2})$,
to express this variation as in Proposition~\ref{propinvone}.
Here, a face $\facee(\finseta,\emptyset,\Gamma)$ is an open dense subset of the product of $\facee(\finseta,\emptyset,\Gamma_{\finseta})$ by $\check{C}(\rats,\emptyset;\Gamma \setminus \Gamma_{\finseta})$,
and $\partau$ identifies $\left[0,1\right] \times \facee(\finseta,\emptyset,\Gamma_{\finseta})$ with 
$\left[0,1\right] \times \crats \times \cinjuptd{\finseta}{\RR^3}$.
The form $\bigwedge_{e \in E(\Gamma)}p_e^{\ast}\bigl(\omega(j_E(e))\bigr)$ pulls back through $\left[0,1\right] \times \cinjuptd{\finseta}{\RR^3} \times \check{C}(\rats,\emptyset;\Gamma \setminus \Gamma_{\finseta})$,
and it vanishes.
\eop

The following variant of Proposition~\ref{propinvone} implies Theorem~\ref{thmConwaycircZ}. We prove it in Section~\ref{seccancelfaces}.

\begin{proposition}
\label{propinvonebis}
Under the assumptions of Proposition~\ref{propinvone},
the following statement
is also true.
For $\Gamma \in \Davis^e_{n,m}(\source)$ and for a connected component $\Gamma_A$ of $\Gamma$ (with no univalent vertices or) whose univalent vertices are consecutive on one component $\source_{\finseta}$ of $\source$,
let $\Gamma^{\rev}(\finseta)$ denote the graph obtained from $\Gamma$ by reversing the order on the univalent vertices of $\Gamma_A$ induced by $i_{\Gamma}$.\footnote{When $\source_{\finseta}$ is oriented, this order is a linear order if $U(\Gamma)$ has vertices of $V(\Gamma) \setminus A$ on $\source_{\finseta}$; it is cyclic, otherwise. When $\source_{\finseta}$ is not oriented, the order is not defined, but reversing the order is well-defined in any case. When $\source_{\finseta}$ is oriented, we have $\left[\Gamma\right]=\left[\Gamma_{\finseta} \#_{\source_{\finseta}}\Gamma_{\finsetv(\Gamma) \setminus \finseta}\right]$ and $\left[\Gamma^{\rev}(\finseta)\right]=\left[\Gamma_{\finseta}^{\rev}(\finseta) \#_{\source_{\finseta}}\Gamma_{\finsetv \setminus \finseta}\right]$.} 
Set
\begin{equation*}I^{\prime}(\Gamma,\finseta)= \frac{1}{2}\int_{\left[0,1\right]\times \facee(\finseta,\Link,\Gamma)}\bigwedge_{e \in E(\Gamma)}p_e^{\ast}\Bigl(\tilde{\omega}\bigl(j_E(e)\bigr)\Bigr)\Bigl(\left[\Gamma\right]-(-1)^{\cardlef{E(\Gamma_{\finseta})}}\left[\Gamma^{\rev}(\finseta)\right]\Bigr).\end{equation*}
Then we have \begin{equation*}\Zinv_n(1)-\Zinv_n(0)=\sum_{ (\Gamma,{\finseta}) \in \Davis^{e,F}_n(\source)}\coefgambet_{\Gamma} I^{\prime}(\Gamma,\finseta)\end{equation*}
and 
\begin{equation*}\Zinv_{n,3n-2}(1)-\Zinv_{n,3n-2}(0)=\sum_{ (\Gamma,{\finseta}) \in \Davis^{e,F}_{n,3n-2}(\source)}\coefgambetred_{\Gamma} I^{\prime}(\Gamma,\finseta).\end{equation*}

\end{proposition}

\bpo{Proof of Theorem~\ref{thmConwaycircZ} assuming Proposition~\ref{propinvonebis}}
There is a map $\rev$ from $\Aavis(S^1)$ to itself that sends the class of a diagram $\Gamma$ to the class of the diagram obtained from $\Gamma$ by reversing the order of the univalent vertices on $S^1$ and by multiplying the class by $(-1)^{\cardlef{T(\Gamma)}}$. The composition $w_C \circ \rev$ equals $w_C$. Furthermore, $w_C$ sends odd-degree diagrams to zero, and $w_C$ is multiplicative with respect to the multiplication of $\Aavis(S^1)$. So for any $(\Gamma,{\finseta}) \in \Davis^{e,F}_{n,m}(S^1)$ as in the above statement such that $\Gamma$ has no component without univalent vertices, $\bigl(w_C\left(\left[\Gamma\right]\right)-(-1)^{\cardlef{E(\Gamma_{\finseta})}}w_C\bigl(\left[\Gamma^{\rev}(\finseta)\right]\bigr)\bigr)$ is equal to
\begin{equation*}w_C\Bigl(\bigl[\Gamma_{\finsetv(\Gamma) \setminus \finseta}\bigr]\Bigr) \Bigl(w_C\left(\left[\Gamma_{\finseta}\right]\right) -(-1)^{\cardlef{E(\Gamma_{\finseta})}+\cardlef{T(\Gamma_{\finseta})}}w_C\bigl(\rev\left(\left[\Gamma_{\finseta}\right]\right)\bigr)\Bigr),\end{equation*} which is zero
since $\cardlef{E(\Gamma_{\finseta})}+\cardlef{T(\Gamma_{\finseta})}$ is even when the degree of $\Gamma_{\finseta}$ is even.
So $w_C(I^{\prime}(\Gamma,\finseta))$ is always zero, and $w_C \bigl( \Zinvlink_{n,3n-2}\bigl(\crats,K,(\omega(i))_{i \in \underline{3n-2}}\bigr)\bigr)$ is independent of the chosen
$\omega(i)$. The variation of $p^c\bigl(\Zinvlink_{n,3n-2}\bigl(\crats,K,(\omega(i))_{i \in \underline{3n-2}}\bigr)\bigr)$
is \begin{equation*}\sum_{ (\Gamma,{\finseta}) \in \Davis^{e,F}_{n,3n-2}(S^1)}\coefgambetred_{\Gamma} p^c\projassis\bigl(I^{\prime}(\Gamma,\finseta)\bigr),\end{equation*} where $w_C$ sends $p^c\projassis\bigl(I^{\prime}(\Gamma,\finseta)\bigr)$ to zero, for the same reasons as above.
\eop

\section{Sketch of proof of Proposition~\ref{propinvone}}
\label{seccancelfacessk}

According to Stokes' theorem, for any $\Gamma \in \Davis^e_{n,m}(\source)$, where $m=3n$ or $m=3n-2$, we have
\begin{multline*}I\left(\rats,\Link,\Gamma,(\tilde{\omega}(i,1))_{i \in \underline{m}}\right)=I\left(\rats,\Link,\Gamma,(\tilde{\omega}(i,0))_{i \in \underline{m}}\right)\\+\sum_F\int_{\left[0,1\right]\times F} \bigwedge_{e \in E(\Gamma)}p_e^{\ast}\bigl(\tilde{\omega}(j_E(e))\bigr),\end{multline*}
where the sum runs over the codimension-one faces $F$ of $C(\rats,\Link;\Gamma)$, which are described in Proposition~\ref{propconffaceun}.
Let $\tilde{p}_e \colon C(\rats,\Link;\Gamma) \to C_e(\rats)$ be the natural restriction
and set \begin{equation*}\projp(\Gamma)=\id_{\left[0,1\right]} \times \prod_{e \in E(\Gamma)} \tilde{p}_e \colon \left[0,1\right] \times C(\rats,\Link;\Gamma) \rightarrow  \left[0,1\right] \times \prod_{e \in E(\Gamma)} C_e(\rats).\end{equation*}
For an edge $e_0$ of $E(\Gamma)$, let $p_{j_E(e_0)} \colon \left[0,1\right] \times \prod_{e \in E(\Gamma)} C_e(\rats) \to \left[0,1\right] \times C_2(\rats)$ be the composition of the natural projection onto $\left[0,1\right] \times C_{e_0}(\rats)$ and the natural identification of $\left[0,1\right] \times C_{e_0}(\rats)$ with $\left[0,1\right] \times C_2(\rats)$. Define the form 
\begin{equation*}\Omega_{E(\Gamma)}=\bigwedge_{e \in E(\Gamma)}p_{j_E(e)}^{\ast}\Bigl(\tilde{\omega}\bigl(j_E(e)\bigr)\Bigr)\end{equation*}
on $\left[0,1\right] \times \prod_{e \in E(\Gamma)} C_e(\rats)$.
Set 
\begin{equation*}I(\Gamma,\finseta)= \int_{\left[0,1\right]\times \facee(\finseta,\Link,\Gamma)}\projp(\Gamma)^{\ast}(\Omega_{E(\Gamma)})\left[\Gamma\right]\end{equation*}\index[N]{Integrals over configuration spaces!IRLGammaA@$I(\Gamma,\finseta)$}for any subset $\finseta$ of ${\finsetv(\Gamma)}$ of cardinality at least $2$, where $\facee(\finseta,\Link,\Gamma)$ is empty (and hence $I(\Gamma,\finseta)=0$) if $\finseta \cap U(\Gamma)$ is not a set of consecutive vertices on one component of $\source$. Set
\begin{equation*}I(\Gamma,\finseta,\infty)= \int_{\left[0,1\right]\times \facee_{\infty}(\finseta,\Link,\Gamma)}\projp(\Gamma)^{\ast}(\Omega_{E(\Gamma)})\left[\Gamma\right]\end{equation*} 
for any subset $\finseta$ of ${\finsetv(\Gamma)}$ of cardinality at least $1$,
where $\facee_{\infty}(\finseta,\Link,\Gamma)$ is empty (and hence $I(\Gamma,\finseta,\infty)=0$) if $\finseta \cap U(\Gamma)$ is not empty. We have
\begin{multline*}\Zinv_{n,m}(1)-\Zinv_{n,m}(0)=\\
\sum_{\Gamma \in \Davis^e_{n,m}(\source)} \left(\sum_{{\finseta} \in \CP_{\geq 2}(\finsetv(\Gamma))}\coefgambet_{\Gamma}I(\Gamma,\finseta)+\sum_{{\finseta} \in \CP_{\geq 1}(\finsetv(\Gamma))}\coefgambet_{\Gamma}I(\Gamma,\finseta,\infty) \right),\end{multline*}
where $\CP_{\geq 1}(\finsetv(\Gamma))$\index[N]{Pzzgreater1@$\CP_{\geq 1}(\finsetv)$ set of nonempty subsets} denotes the set of the (nonstrict) nonempty subsets of $\finsetv(\Gamma)$.
In order to prove Proposition~\ref{propinvone}, it suffices to prove that 
the codimension-one faces $F$ of the $C(\rats,\Link;\Gamma)$ that do not appear in the statement of Proposition~\ref{propinvone} do not contribute.

This is the consequence of Lemmas~\ref{lemfaceinf} to \ref{lemstu}, together with the analysis before Lemma~\ref{lemihx}.

\begin{lemma}
\label{lemfaceinf} For any $\Gamma \in \Davis^e_{n,m}(\source)$,
For any nonempty subset ${\finseta}$ of ${\finsetv(\Gamma)}$, we have
$I(\Gamma,\finseta,\infty)=0$.
\end{lemma}
\bp
Recall 
$\facee_{\infty}(\finseta,\Link,\Gamma) = \baseb_{\infty}(\finseta,\Link,\Gamma) \times \sinjupdtcs(T_{\infty}\rats,{\finseta})$ from Section~\ref{secfacecodimone}.
Let $E_{\finsetc}$ be the set of the edges of $\Gamma$ that contain an element of ${\finsetv(\Gamma)} \setminus {\finseta}$ and an element of ${\finseta}$.
Let $p_2$ denote the projection of $\facee_{\infty}(\finseta,\Link,\Gamma)$ onto $ \sinjupdtcs(T_{\infty}\rats,{\finseta})$.
For $e \in E_{\finseta} \cup E_{\finsetc}$, the map $P_e:(S^2)^{E_{\finseta} \cup E_{\finsetc}} \longrightarrow S^2$ is the projection onto the factor indexed by $e$.
We prove that there exists a smooth map
\begin{equation*}g \colon \sinjupdtcs(T_{\infty}\rats,{\finseta}) \longrightarrow (S^2)^{E_{\finseta} \cup E_{\finsetc}}\end{equation*}
such that 
$\bigwedge_{e \in E_{\finseta} \cup E_{\finsetc}}p_e^{\ast}\Bigl(\tilde{\omega}\bigl(j_E(e)\bigr)\Bigr)$ is equal to
\begin{equation*}\bigl(\id_{\left[0,1\right]} \times (g \circ p_2)\bigr)^{\ast}\left( \bigwedge_{e \in E_{\finseta} \cup E_{\finsetc}}(\id_{\left[0,1\right]} \times P_e)^{\ast}\Bigl(\tilde{\omega}_{S^2}\bigl(j_E(e)\bigr)\Bigr) \right).\end{equation*}
If $e \in E_{\finseta} \cup E_{\finsetc}$, then $p_e(\facee_{\infty}(\finseta,\Link,\Gamma)) \subset \partial C_2({\rats}) \setminus \ST(\crats)$, and we have \begin{equation*}\left(\id_{\left[0,1\right]} \times p_e\right)^{\ast}\Bigl(\tilde{\omega}\bigl(j_E(e)\bigr)\Bigr)=\left(\id_{\left[0,1\right]} \times (\projp_{\partau} \circ p_e)\right)^{\ast}\Bigl(\tilde{\omega}_{S^2}\bigl(j_E(e)\bigr)\Bigr).\end{equation*}
If $e \in E_{\finsetc}$, then $\projp_{\partau} \circ p_e$ depends only on the projection on $\sph(T_{\infty}{\rats})$ of the vertex at $\infty$ (of ${\finseta}$).
If $e \in E_{\finseta}$, then $\projp_{\partau} \circ p_e$ factors through $\sinjupdtcs(T_{\infty}\rats,e)$.
So in both cases, the map $\projp_{\partau} \circ p_e$ factors through $\sinjupdtcs(T_{\infty}\rats,{\finseta})$. Thus it may be expressed as
$((P_e \circ g) \circ p_2)$.
Therefore, if the degree of the form $\bigwedge_{e \in E_{\finseta} \cup E_{\finsetc}}p_e^{\ast}\bigl(\tilde{\omega}_{S^2}(j_E(e))\bigr)$ is bigger than the dimension $3 \cardlef{\finseta}$ of $\left[0,1\right] \times \sinjupdtcs(T_{\infty}\rats,{\finseta})$, this form vanishes on $\left[0,1\right] \times\facee_{\infty}(\finseta,\Link,\Gamma)$.
The degree of the form 
is $( 2\cardlef{E_{\finseta}} + 2\cardlef{E_{\finsetc}})$, and we have
\begin{equation*}3 \cardlef{\finseta}= 2\cardlef{E_{\finseta}} + \cardlef{E_{\finsetc}}.\end{equation*}
So the integral vanishes unless $E_{\finsetc}$ is empty.
In this case, all the $\projp_{\partau} \circ p_e$, for $e \in E_{\finseta}$
factor through the conjugates under the inversion $(x \mapsto x/\norm{x}^2)$ of the translations that make sense, and the form $ \bigwedge_{e \in E_{\finseta}}p_e^{\ast}\bigl(\tilde{\omega}_{S^2}(j_E(e))\bigr)$ factors through the product by $\left[0,1\right]$ of the quotient of $\sinjupdtcs(T_{\infty}\rats,{\finseta})$ by these translation conjugates. So it vanishes, too.
\eop

If there exists a smooth map from $\left[0,1\right] \times \facee(\finseta,\Link,\Gamma)$ to a manifold of strictly smaller
dimension that factorizes the restriction of \begin{equation*}\projp(\Gamma)=\Bigl( \id_{\left[0,1\right]} \times \prod_{e \in E(\Gamma)} \tilde{p}_e \Bigr) \colon \left[0,1\right] \times C(\rats,\Link;\Gamma) \rightarrow   \left[0,1\right] \times \prod_{e \in E(\Gamma)} C_e(\rats)\end{equation*} to $\left[0,1\right] \times \facee(\finseta,\Link,\Gamma)$, then we have $I(\Gamma,\finseta)=0$.
We use this principle to get rid of some faces.

\begin{lemma}
\label{lemdiscon} 
Let $\Gamma \in \Davis^e_{n,m}(\source)$. For any subset ${\finseta}$ of ${\finsetv(\Gamma)}$ such that the graph $\Gamma_{\finseta}$ defined before Proposition~\ref{propinvone} is not connected and $\Gamma_{\finseta}$ is not a pair of univalent vertices, we have
$I(\Gamma,\finseta)=0$.
\end{lemma}
\bp
In the fiber $\cinjuptd{\finseta}{T_{m({\finseta})}\crats,\Link,\Gamma}$ of $\facee(\finseta,\Link,\Gamma)$ we may translate
one connected component of $\Gamma_{\finseta}$ whose set of vertices is $\finsetc$, independently, without changing the restriction of $\projp(\Gamma)$ to $\facee(\finseta,\Link,\Gamma)$. The translation vector is in $T_{m({\finseta})}\Link$ when $\finsetc$ contains univalent vertices.
Unless $\finsetc$ and $\finseta \setminus \finsetc$ are reduced to a univalent vertex, the quotient of $\cinjuptd{\finseta}{T_{m({\finseta})}\crats,\Link,\Gamma}$ by these translations has a smaller dimension than $\cinjuptd{\finseta}{T_{m({\finseta})}\crats,\Link,\Gamma}$, and 
$\projp(\Gamma)$ factors through the corresponding quotient of $\left[0,1\right] \times \facee(\finseta,\Link,\Gamma)$.
\eop

\begin{lemma}
\label{lemedge} 
Let $\Gamma \in \Davis^e_{n,m}(\source)$. Let ${\finseta}$ be a subset of ${\finsetv(\Gamma)}$ such that $\cardlef{\finseta} \geq 3$. If some trivalent vertex of ${\finseta}$ belongs to exactly one edge of $\Gamma_{\finseta}$, then
$I(\Gamma,\finseta)=0$.
\end{lemma}
\bp Let $b$ be the mentioned trivalent vertex. Let $e$ be its edge in $\Gamma_{\finseta}$, and let $d \in {\finseta}$ be the other element of $e$. The group $\left]0,\infty\right[$ acts on the map $t$ from ${\finseta}$ to $T_{c(b)}{\rats}$ by moving $t(b)$ on the half-line from $t(d)$ through $t(b)$, by multiplying $(t(b)-t(d))$ by a scalar. When $\cardlef{\finseta} \geq 3$, this action defined on an open dense subset of $\cinjuptd{\finseta}{T_{m({\finseta})}\crats,\Link,\Gamma}$ is not trivial, and $\projp(\Gamma)$
factors through the corresponding quotient of an open dense subset of $\left[0,1\right] \times \facee(\finseta,\Link,\Gamma)$, which is of smaller dimension.
\eop

\section{Cancellations of nondegenerate faces}
\label{seccancelfaces}

From now on, we will study cancellations that are no longer individual, and
orientations must be seriously considered. Recall that the codimension-one faces are oriented as parts of the boundary of $C(\rats,\Link;\Gamma)$, with the outward normal first convention, where $C(\rats,\Link;\Gamma)$ is oriented 
by an orientation of $\Link$ and an order on $V(\Gamma)$.
The relations between an orientation of $V(\Gamma)$, which orients $C(\rats,\Link;\Gamma)$, a vertex-orientation of $\Gamma$, and an edge-orientation of the set $H(\Gamma)$ of half-edges of $\Gamma$ are explained in Lemma~\ref{lemorc}, Corollary~\ref{cororc}, and Remark~\ref{rkedgeorvert}. Fortunately, we do not have to fix everything to compare similar orientations.

\begin{lemma}
\label{lemsym} 
Let $\Gamma \in \Davis^e_{n,m}(\source)$. Let ${\finseta}$ be a subset of ${\finsetv(\Gamma)}$ such that at least one element of ${\finseta}$ belongs to exactly two edges of $\Gamma_{\finseta}$. Let $\CE(\Gamma,\finseta)$  denote the set of graphs of $\Davis^e_{n,m}(\source)$ that are isomorphic to $\Gamma$ by an isomorphism that is only allowed to change the labels and the orientations of the edges of $\Gamma_{\finseta}$. Such an isomorphism preserves $\finseta$, and we have
\begin{equation*}\sum_{\tilde{\Gamma} \suchthat  \tilde{\Gamma} \in \CE(\Gamma,\finseta)}\coefgambetred_{\tilde{\Gamma}}I(\tilde{\Gamma},\finseta)=0.\end{equation*} 
\end{lemma}
\bp Let us first check that the isomorphisms of the statement preserve $\finseta$.
The vertices of the elements of $\Davis^e_{n,m}(\source)$ are not numbered. A vertex is characterized by the half-edges
that contain it. Therefore, the isomorphisms of the statement preserve the vertices of $\finsetv(\Gamma) \setminus \finseta$. So they preserve $\finseta$ setwise. These isomorphisms also preserve the vertices that have adjacent edges outside $E(\Gamma_{\finseta})$ pointwise. The isomorphisms described below actually induce the identity map on $\finsetv(\Gamma)$.

Among the vertices of ${\finseta}$ that belong to exactly two edges of $\Gamma_{\finseta}$ and one edge $j_E^{-1}(k)$ of $\Gamma$ outside $\Gamma_{\finseta}$,
choose the vertices such that $k$ is minimal. If there is one such vertex, then call this vertex $v_m$. Otherwise, there are two choices, and $v_m$ is chosen to be the vertex that belongs to the first half-edge of $j_E^{-1}(k)$.

We first describe an orientation-reversing diffeomorphism of the complement of a codimension-three submanifold of $\facee(\finseta,\Link,\Gamma)$. Let $v_j$ and $v_k$ denote the (possibly 
equal) two other vertices of the two edges of $\Gamma_{\finseta}$ that contain $v_m$.
Consider the transformation $S$ of the space $\cuptd{\finseta}{T_{c(\finseta)}\rats}$ of nonconstant maps $f$ from ${\finseta}$ to $T_{c(\finseta)}{\rats}$
up to translations and dilations that maps $f$ to $S(f)$, where
\begin{center}
$S(f)(v_{\ell})=f(v_{\ell})$ if $v_\ell \neq v_m$, and\\
$S(f)(v_m)=f(v_j)+f(v_k)-f(v_m)$.\\
\end{center}
This is an orientation-reversing involution of $\cuptd{\finseta}{T_{c(\finseta)}\rats}$.
The set of elements of $\cinjuptd{\finseta}{T_{c(\finseta)}{\rats}}$
whose image under $S$ is not in $\cinjuptd{\finseta}{T_{c(\finseta)}{\rats}}$ is a codimension-three submanifold of $\cinjuptd{\finseta}{T_{c(\finseta)}{\rats}}$.
The fibered product of $S$ by the identity of the base $\baseb(\finseta,\Link,\Gamma)$ is an orientation-reversing smooth involution outside a codimension-three submanifold $F_S$ of $\facee(\finseta,\Link,\Gamma)$. It is still denoted by $S$, as its product by $\id_{\left[0,1\right]}$ is, too.

Now, let $\sigma({\finseta};\Gamma)(\tilde{\Gamma})$ be obtained from $(\tilde{\Gamma} \in \CE(\Gamma,\finseta))$ by exchanging 
the labels of the two edges of $\Gamma_{\finseta}$ that contain $v_m$ and by reversing their orientations if (and only if) they both start or end at $v_m$.
Then, as Figure~\ref{figparallelogram} shows, we have
\begin{equation*}\projp(\tilde{\Gamma}) \circ S=\projp\bigl(\sigma({\finseta};\Gamma)(\tilde{\Gamma})\bigr).\end{equation*}

\bfig
\centering
\begin{tikzpicture}[scale=.7] 
\draw [->] (4,.5) -- (5.4,1.2);
\draw [->] (1,.5) -- (2.4,1.2);
\draw [->] (1,.5) -- (2.5,.5);
\draw [->] (4,2) -- (5.5,2);
\draw (5.4,1.2) -- (7,2) -- (5.5,2) (2.4,1.2) -- (4,2) (2.5,.5) -- (4,.5)
(5.5,2) node[above]{$\scriptstyle a$}
(2.5,.4) node[below]{$\scriptstyle \sigma({\finseta};\Gamma)(a)$}
(5.6,1.1) node[right]{$\scriptstyle b$}
(2.5,1.5) node[left]{$\scriptstyle \sigma({\finseta};\Gamma)(b)$} (.9,.6) node[below]{$\scriptstyle f(v_m)$}
(3.9,2.1) node[left]{$\scriptstyle f(v_k)$}
(4.2,.4) node[right]{$\scriptstyle f(v_j)$}
(7,1.9) node[above]{$\scriptstyle S(f(v_m))$};
\fill (4,.5) circle (2pt) (1,.5) circle (2pt) (4,2) circle (2pt) (7,2) circle (2pt);
\begin{scope}[xshift=7.5cm]
\draw [->] (4,.5) -- (5.4,1.2);
\draw [->] (1,.5) -- (2.4,1.2);
\draw [-<] (1,.5) -- (2.5,.5);
\draw [-<] (4,2) -- (5.5,2);
\draw (5.4,1.2) -- (7,2) -- (5.5,2) (2.4,1.2) -- (4,2) (2.5,.5) -- (4,.5)
(5.5,2) node[above]{$\scriptstyle a$}
(2.5,.4) node[below]{$\scriptstyle \sigma({\finseta};\Gamma)(a)$}
(5.6,1.1) node[right]{$\scriptstyle b$}
(2.5,1.5) node[left]{$\scriptstyle \sigma({\finseta};\Gamma)(b)$} (.9,.6) node[below]{$\scriptstyle f(v_m)$}
(3.9,2.1) node[left]{$\scriptstyle f(v_k)$}
(4.2,.4) node[right]{$\scriptstyle f(v_j)$}
(7,1.9) node[above]{$\scriptstyle S(f(v_m))$};
\fill (4,.5) circle (2pt) (1,.5) circle (2pt) (4,2) circle (2pt) (7,2) circle (2pt);
\end{scope}

\end{tikzpicture}
\caption{The parallelogram transformation $S$}
\label{figparallelogram}
 
\end{figure}
Recall \begin{equation*}I(\Gamma,\finseta)= \int_{\left[0,1\right]\times \facee(\finseta,\Link,\Gamma)}\projp({\Gamma})^{\ast}(\Omega_{E(\Gamma)})\left[\Gamma\right],\end{equation*} 
with the map \begin{equation*}\projp(\Gamma)=\left(\id_{\left[0,1\right]} \times \prod_{e \in E(\Gamma)} \tilde{p}_e \right)\colon \left[0,1\right] \times C(\rats,\Link;\Gamma) \rightarrow  \left[0,1\right] \times \prod_{e \in E(\Gamma)} C_e(\rats)\end{equation*} and
$\Omega_{E(\Gamma)}=\bigwedge_{e \in E(\Gamma)}p_{j_E(e)}^{\ast}\bigl(\tilde{\omega}(j_E(e))\bigr)$. We have
\begin{equation*}\begin{array}{ll}
I(\tilde{\Gamma},\finseta) &=\int_{\left[0,1\right]\times (\facee(\finseta,\Link,\tilde{\Gamma}) \setminus F_S)}\projp(\tilde{\Gamma})^{\ast}({\Omega_{E(\Gamma)}})[\tilde{\Gamma}]\\
&=-\int_{\left[0,1\right]\times (\facee(\finseta,\Link,\tilde{\Gamma}) \setminus F_S)}S^{\ast}\left(\projp(\tilde{\Gamma})^{\ast}({\Omega_{E(\Gamma)}})\right)[\tilde{\Gamma}]\\
&=-\int_{\left[0,1\right]\times(\facee(\finseta,\Link,\tilde{\Gamma}) \setminus F_S)}(\projp(\tilde{\Gamma}) \circ S)^{\ast}({\Omega_{E(\Gamma)}})[\tilde{\Gamma}]\\ 
&=-\int_{\left[0,1\right]\times (\facee(\finseta,\Link,\tilde{\Gamma}) \setminus F_S)}\projp(\sigma({\finseta};\Gamma)(\tilde{\Gamma}))^{\ast}({\Omega_{E(\Gamma)}})[\tilde{\Gamma}]\\
&=-I\bigl(\sigma({\finseta};\Gamma)(\tilde{\Gamma}),\finseta\bigr)\end{array}\end{equation*}
since we have $\bigl[\tilde{\Gamma}\bigr]=\bigl[\sigma({\finseta};\Gamma)(\tilde{\Gamma})\bigr]$.
Now, $\sigma({\finseta};\Gamma)$ defines an involution of $\CE(\Gamma,\finseta)$, and we get
\begin{equation*}\sum_{\tilde{\Gamma} \in \CE(\Gamma,\finseta)}I(\tilde{\Gamma},\finseta)=\sum_{\tilde{\Gamma} \in \CE(\Gamma,\finseta)}I\bigl(\sigma({\finseta};\Gamma)(\tilde{\Gamma}),\finseta\bigr)=-\sum_{\tilde{\Gamma} \in \CE(\Gamma,\finseta)}I(\tilde{\Gamma},\finseta)=0.\end{equation*}
\eop

Maxim Kontsevich observed the symmetry of the above proof in \cite{ko}.

\begin{lemma}
\label{lemsymeven}
Let $\Gamma \in \Davis^e_n(\source)$. Let ${\finseta}$ be a subset of ${\finsetv(\Gamma)}$ such that $\Gamma_{\finseta}$ is a connected component of $\Gamma$. Let ${\Gamma}^{\eo}(\finseta)$ denote the graph obtained from $\Gamma$ by reversing all the orientations of the edges of $\Gamma_{\finseta}$.
Recall the notation of Proposition~\ref{propinvone}.
If $\Gamma_{\finseta}$ is a diagram of even degree without univalent vertices, then we have
\begin{equation*}I\bigl({\Gamma},\finseta\bigr)+I\bigl({\Gamma}^{\eo}(\finseta),\finseta\bigr) =0.\end{equation*}
If $(\Gamma,\finseta)$ is an element of the set $\Davis^{e,F}_{n,m}(\source)$ defined in Proposition~\ref{propinvone} and if $\Gamma_{\finseta}$ has univalent vertices,
let ${\Gamma}^{\eo , \rev}(\finseta)$ denote the graph obtained from ${\Gamma}^{\eo}(\finseta)$ by reversing the order of the univalent vertices of $\Gamma_{\finseta}$ on $\source_{\finseta}$.
Recall \begin{equation*}I^{\prime}(\Gamma,\finseta)=\frac12 \int_{\left[0,1\right]\times \facee(\finseta,\Link,\Gamma)}\bigwedge_{e \in E(\Gamma)}p_e^{\ast}\Bigl(\tilde{\omega}\bigl(j_E(e)\bigr)\Bigr)\Bigl(\left[\Gamma\right]-(-1)^{\cardlef{E(\Gamma_{\finseta})}}\left[\Gamma^{\rev}(\finseta)\right]\Bigr)\end{equation*} from Proposition~\ref{propinvonebis}.
Then we have \begin{equation*}I\bigl({\Gamma},\finseta\bigr)+I\bigl({\Gamma}^{\eo , \rev}(\finseta),\finseta\bigr) =I^{\prime}\bigl({\Gamma},\finseta\bigr)+I^{\prime}\bigl({\Gamma}^{\eo , \rev}(\finseta),\finseta\bigr).\end{equation*}
\end{lemma}
\bp 
Set $\overline{\Gamma}(A)=\Gamma^{\eo , \rev}(\finseta)$ in both cases.
The opposite of the identity map of $T_{c(\finseta)}\crats$ induces a diffeomorphism from the fiber 
of $\facee(\finseta,\Link,\Gamma)$ to the fiber of  $\facee(\finseta,\Link,\overline{\Gamma}(A))$, which induces a diffeomorphism from $\facee(\finseta,\Link,\Gamma)$ to
$\facee(\finseta,\Link,\overline{\Gamma}(A))$ over the identity map of the base. Denote by $\CS$ the product of this diffeomorphism by $\id_{\left[0,1\right]}$. Let us carefully discuss orientations to determine when this diffeomorphism preserves the orientation.

Order and orient the vertices of $\Gamma$ so that the corresponding orientation of $H(\Gamma)$, as in Remark~\ref{rkedgeorvert}, is induced by the edge-orientation of $\Gamma$.
There is a natural bijection from $V(\Gamma)$ to $V(\overline{\Gamma}(A))$. 
This bijection is the identity on the set $V(\Gamma) \setminus \finseta$ of vertices of $\Gamma \setminus \Gamma_{\finseta}$, unaffected by the modifications.
When $\finseta$ is not a pair of vertices in a $\theta$-component, 
a vertex of $\Gamma_{\finseta}$ is characterized by the labels of the edges that contain it. Such a vertex is sent to the vertex of $\overline{\Gamma}(A)$ with the same set of labels of adjacent edges.
When $\finseta$ is a pair of vertices in a $\theta$-component, the vertex at which an edge of $\Gamma_{\finseta}$ labeled by $i$ begins is sent to the vertex of $\overline{\Gamma}(A)$ at which the edge of 
$\overline{\Gamma}(A)$ labeled by $i$ ends.

Order the vertices of $\overline{\Gamma}(A)$ like the vertices of $\Gamma$ if $\cardlef{E(\Gamma_{\finseta})}$ is even, and permute two vertices if $\cardlef{E(\Gamma_{\finseta})}$ is odd.
Orient $C(\rats,\Link;\Gamma)$ and $C(\rats,\Link;\overline{\Gamma}(A))$ with respect to the above orders of $V(\Gamma)$ and $V(\overline{\Gamma}(A))$, using the orientations of $\rats$ and $\Link$. Then $\CS$ reverses the orientation if and only if $\cardlef{E(\Gamma_{\finseta})}$ is even since $\cardlef{V(\Gamma_{\finseta})}$ is even.

Orient the vertices of $\overline{\Gamma}(A)$ like the vertices of $\Gamma$. So the orientation of $H(\overline{\Gamma}(A))$ associated to that vertex-orientation and to the above order of vertices is induced by the edge-orientation of $\overline{\Gamma}(A)$. 
When orientations are fixed as above, set $I_0(\Gamma,\finseta)=\int_{\left[0,1\right]\times \facee(\finseta,\Link,{\Gamma})}\projp({\Gamma})^{\ast}({\Omega_{E(\Gamma)}})$ so that $I(\Gamma,\finseta)=I_0(\Gamma,\finseta)\left[\Gamma\right]$. Let $e$ be an edge of $\Gamma$. Up to its edge-orientation, the edge $e$ is also an edge of $\overline{\Gamma}(A)$, and $p_{j_E(e)} \circ \projp(\overline{\Gamma}(A))$ 
restricts to $\left[0,1\right] \times \facee(\finseta,\Link,\overline{\Gamma}(A))$ as $p_{j_E(e)} \circ \projp(\Gamma) \circ \CS^{-1}$. We have
\begin{equation*}\begin{array}{ll}I_0\bigl(\overline{\Gamma}(A),\finseta\bigr)&=\int_{\left[0,1\right]\times \facee(\finseta,\Link,\overline{\Gamma}(A))}\projp(\overline{\Gamma}(A))^{\ast}({\Omega_{E(\overline{\Gamma}(A))}})\\
&=\int_{\left[0,1\right]\times \facee(\finseta,\Link,\overline{\Gamma}(A))}\bigwedge_{e \in E(\Gamma)}\Bigl(p_{j_E(e)}\circ \projp(\overline{\Gamma}(A))\Bigr)^{\ast}\Bigl(\tilde{\omega}\bigl(j_E(e)\bigr)\Bigr)
\\
&=\int_{\left[0,1\right]\times \facee(\finseta,\Link,\overline{\Gamma}(A))}\left(\CS^{-1}\right)^{\ast}\Bigl(\bigwedge_{e \in E(\Gamma)}\left(p_{j_E(e)}\circ \projp({\Gamma})\right)^{\ast}\bigl(\tilde{\omega}\bigl(j_E(e)\bigr)\bigr)\Bigr)
\\
&=(-1)^{\cardlef{E(\Gamma_{\finseta})}+1} \int_{\left[0,1\right]\times \facee(\finseta,\Link,\Gamma)}\bigwedge_{e \in E(\Gamma)}\left(p_{j_E(e)}\circ \projp({\Gamma})\right)^{\ast}\Bigl(\tilde{\omega}\bigl(j_E(e)\bigr)\Bigr)
\\
&=(-1)^{\cardlef{E(\Gamma_{\finseta})}+1}I_0(\Gamma,\finseta).\end{array}\end{equation*}
In particular, if $\Gamma_{\finseta}$ is a connected diagram of even degree without univalent vertices, 
then $\cardlef{E(\Gamma_{\finseta})}$ is even, and $\left[\Gamma\right]$ is equal to $\left[\overline{\Gamma}(A)\right]$ and $\left[{\Gamma}^{\eo}(\finseta)\right]$. So we have $I({\Gamma},\finseta)=-I({\Gamma}^{\eo}(\finseta),\finseta)$.
Otherwise, we have \begin{equation*}\begin{array}{lll} I({\Gamma},\finseta)+I(\overline{\Gamma}(A),\finseta)& =&
\frac12\left(I_0(\Gamma,\finseta) -(-1)^{\cardlef{E(\Gamma_{\finseta})}} I_0(\overline{\Gamma}(A)),\finseta)\right)\left[\Gamma\right]\\&& + \frac12\left(I_0(\overline{\Gamma}(A)),\finseta) -(-1)^{\cardlef{E(\Gamma_{\finseta})}}I_0(\Gamma,\finseta) \right)\left[{\Gamma}^{\rev}(\finseta)\right]\\
&=&I^{\prime}({\Gamma},\finseta)+I^{\prime}(\overline{\Gamma}(A),\finseta).\end{array}\end{equation*}
\eop

Lemmas \ref{lemdiscon}, \ref{lemedge}, and \ref{lemsym}  allow us to get rid of the pairs $(\Gamma,{\finseta})$ with $\cardlef{\finseta} \geq 3$ such that 
\begin{itemize}
 \item at least one element of ${\finseta}$ does not have all its adjacent edges in $E(\Gamma_{\finseta})$, or
\item  $\Gamma_{\finseta}$ is disconnected.
\end{itemize}
Lemma~\ref{lemsymeven} rules out the pairs $(\Gamma,{\finseta})$ such that $\Gamma_{\finseta}$ is an even degree connected component of $\Gamma$, without univalent vertices (where $\cardlef{\finseta} \equiv 0$ mod $4$).
Therefore, according to Lemma~\ref{lemfaceinf}, we are left
with \begin{itemize}
 \item 
the pairs $(\Gamma,{\finseta})$ of the statement of Proposition~\ref{propinvone}, for which $\Gamma_{\finseta}$ is a connected component of $\Gamma$ (which may be an edge between two univalent vertices) and
\item the following pairs, for which $\cardlef{\finseta} =2$ and
\begin{itemize}
\item $\Gamma_{\finseta}$ is an edge between two trivalent vertices,
\item $\Gamma_{\finseta}$ is an edge between a trivalent vertex and a univalent one, or
\item $\Gamma_{\finseta}$ is a pair of isolated consecutive univalent vertices,
\end{itemize}
since Lemma~\ref{lemdiscon} rules out the disconnected $\Gamma_{\finseta}$ with a trivalent vertex, and Lemma~\ref{lemsym} rules out $\Gamma_{\finseta}= \begin{tikzpicture} \useasboundingbox (-.3,-.1) rectangle (.3,.2);
\draw (-.2,0) .. controls (0,.15) .. (.2,0);
\draw (-.2,0) .. controls (0,-.15) .. (.2,0);
\fill (-.2,0) circle (1.5pt) (.2,0) circle (1.5pt);
\end{tikzpicture}$.
\end{itemize}

The following lemma allows us to eliminate the pairs $(\Gamma,{\finseta})$ such that $\Gamma_{\finseta}$ is an edge between two trivalent vertices using the Jacobi relation.

\begin{lemma}
\label{lemihx} The contributions to $(\Zinv_{n,m}(1)-\Zinv_{n,m}(0))$ of the faces $\facee(\finseta,\Link,{\Gamma})$ for which $\Gamma_{\finseta}$ is an edge between two trivalent vertices cancel. More precisely,
let $\Gamma \in \Davis^e_{n,m}(\source)$. Let ${\finseta}$ be a subset of ${\finsetv(\Gamma)}$
such that $\Gamma_{\finseta}$ is an edge $e(\ell)$ with label $\ell$.
Let $\Gamma/\Gamma_{\finseta}$ be the labeled edge-oriented graph obtained from $\Gamma$
by contracting $\Gamma_{\finseta}$ to one point, as in Figure~\ref{figihx}. 
(The labels of the edges of $\;\Gamma/\Gamma_{\finseta}$ belong to $\underline{m} \setminus \{\ell\}$. The graph $\Gamma/\Gamma_{\finseta}$ has one four-valent vertex, and its other vertices are univalent or trivalent.) 
Let $\CE(\Gamma;{\finseta})$ be the subset of $\Davis^e_{n,m}(\source)$ that contains the graphs $\tilde{\Gamma}$ equipped with a pair $\finseta$ of vertices joined by an edge $e(\ell)$ with label $\ell$ such that $\tilde{\Gamma}/\tilde{\Gamma}_{\finseta}$ is equal to $\Gamma/\Gamma_{\finseta}$. Then we have
\begin{equation*}\sum_{\tilde{\Gamma} \suchthat  \tilde{\Gamma} \in \CE(\Gamma;{\finseta})}\coefgambet_{\tilde{\Gamma}}I(\tilde{\Gamma},\finseta)=0.\end{equation*} 
\end{lemma}

\bfig \centering
\begin{tikzpicture}
\draw [very thick, ->] (2,.5) -- (3,.5) (1,.5) -- (2,.5);
\draw (0,.8) -- (1,.5) -- (0,.2) (4,.8) -- (3,.5) -- (4,.2) node[right]{ \scriptsize $c$} (4,.8) node[right]{ \scriptsize $d$}
(0,.8) node[left]{ \scriptsize $a$} (0,.2) node[left]{ \scriptsize $b$}
(1.2,.5) node[below]{ \scriptsize $v(\ell,1)$}
(2.8,.5) node[below]{ \scriptsize $v(\ell,2)$}
(2,.5) node[above]{\scriptsize $e(\ell)$}
(2,.2) node[below]{\scriptsize $\Gamma$};
\fill (1,.5) circle (2pt) (3,.5) circle (2pt);
\begin{scope}[xshift=-2cm]
 \draw (8,.8) -- (9,.5) -- (8,.2) (10,.8) -- (9,.5) -- (10,.2) node[right]{ \scriptsize $c$} 
(10,.8) node[right]{ \scriptsize $d$}
(8,.8) node[left]{ \scriptsize $a$} (8,.2) node[left]{ \scriptsize $b$}
(9,.2) node[below]{\scriptsize $\Gamma/\Gamma_{\finseta}$};
\fill (9,.5) circle (1.5pt);
\end{scope}
\end{tikzpicture}
\caption{The graph $\Gamma$, its bold subgraph $\Gamma_{\finseta}$, and $\Gamma/\Gamma_{\finseta}$}
\label{figihx}
\end{figure} 
\bp
Let us prove that there are $6$ graphs in $\CE(\Gamma;{\finseta})$.
Let $a,b,c,d$ be the four half-edges of 
$\Gamma/\Gamma_{\finseta}$ that contain its four-valent vertex. 
In $\tilde{\Gamma}$, the edge $e(\ell)$ goes from a vertex $v(\ell,1)$ to a vertex $v(\ell,2)$.
The vertex $v(\ell,1)$ is adjacent to the first half-edge of $e(\ell)$ and to two half-edges of $\{a,b,c,d\}$.
The unordered pair of $\{a,b,c,d\}$ adjacent to $v(\ell,1)$
determines $\tilde{\Gamma}$ as an element of $\Davis^e_{n,m}(\source)$ 
and
there are $6$ graphs in $\CE(\Gamma;{\finseta})$ labeled by the pairs of elements
of $\{a,b,c,d\}$. They are $\Gamma=\Gamma_{ab}$, $\Gamma_{ac}$, $\Gamma_{ad}$, $\Gamma_{bc}$, $\Gamma_{bd}$, and $\Gamma_{cd}$.

The face $\facee(\finseta,\Link,\Gamma)$ is fibered over $\baseb(\finseta,\Link,\Gamma)$ with fiber $\cinjuptd{\finseta}{T_{c(v(\ell,1))}\rats}=_{\partau}S^2$, which contains the direction of the vector from $c(v(\ell,1))$ to $c(v(\ell,2))$. 
Consistently order the vertices of the $\Gamma_{..}$ starting with $v(\ell,1)$, $v(\ell,2)$ (the other vertices are in natural correspondences for different $\Gamma_{..}$). Use these orders to orient the
configuration spaces ${C}(\rats,\Link;\Gamma_{..})$.

The oriented face $\facee(\finseta,\Link,\Gamma_{..})$ and the map
\begin{equation*}\projp(\Gamma_{..})\colon \left[0,1\right] \times \bigl(\facee(\finseta,\Link,\Gamma_{..}) \subset {C}(\rats,\Link;\Gamma_{..})\bigr) \longrightarrow \left[0,1\right]\times \prod_{e \in E(\Gamma_{..})} C_2({\rats})^{e}\end{equation*}
are the same for all the elements $\Gamma_{..}$ of $\CE(\Gamma;{\finseta})$.
Therefore, the \begin{equation*}I_0(\Gamma_{..},\finseta)=\int_{\left[0,1\right]\times \facee(\finseta,\Link,\Gamma_{..})}\projp(\Gamma_{..})^{\ast}(\Omega_{E(\Gamma_{..})})\end{equation*} are the same for all the elements $\Gamma_{..}$ of $\CE(\Gamma;{\finseta})$ (for our consistent orders of the vertices), and the sum of the statement is 
\begin{equation*}\sum_{\tilde{\Gamma} \suchthat  \tilde{\Gamma} \in \CE(\Gamma;{\finseta})}\coefgambet_{\tilde{\Gamma}}I_0(\tilde{\Gamma},\finseta)[\tilde{\Gamma}]=
\coefgambet_{\Gamma}I_0(\Gamma,\finseta) \sum_{\tilde{\Gamma} \suchthat  \tilde{\Gamma} \in \CE(\Gamma;{\finseta})}[\tilde{\Gamma}].\end{equation*}

Let $e_1$ be the first half-edge of $e(\ell)$, and let $e_2$ be the other half-edge of $e(\ell)$.
Equip $\Gamma=\Gamma_{ab}$ with a vertex-orientation, represented by
$(a,b,e_1)$ at $v(\ell,1)$, and $(c,d,e_2)$ at $v(\ell,2)$, which is consistent
with its given edge-orientation (i.e., such that the edge-orientation of $H(\Gamma)$
is equivalent to its vertex-orientation, with respect to the above order of vertices).
The orientation of $H(\Gamma)$ is represented by
$(a,b,e_1,c,d,e_2, \dots )$.
It induces the edge-orientation of $H(\Gamma)$, which is the same for all the elements of $\CE(\Gamma;{\finseta})$.

Thus, permuting the letters $b,c,d$ cyclically gives rise to two other graphs ($\Gamma_{ac}$ and $\Gamma_{ad}$)
in $\CE(\Gamma;{\finseta})$ equipped with suitable vertex-orientations, respectively represented by 
\begin{equation*}
\begin{array}{l}
 (a,c,e_1)\mbox{ at }v(\ell,1)\mbox{, and }(d,b,e_2)\mbox{ at }v(\ell,2)\mbox{, or }\\
 (a,d,e_1)\mbox{ at }v(\ell,1)\mbox{, and }(b,c,e_2)\mbox{ at }v(\ell,2)\mbox{.}
\end{array}
\end{equation*}

The three other elements of $\CE(\Gamma;{\finseta})$ with their suitable vertex-orientation are obtained from the three previous ones by exchanging the ordered pair before $e_1$ with the ordered pair before $e_2$. This does not change the unlabeled vertex-oriented graph. The first three graphs can be represented by the following three graphs identical outside the pictured disk:
\begin{equation*}\begin{tikzpicture} \useasboundingbox (-.5,.2) rectangle (1.5,1.1);
\draw [->] (.25,0) -- (.45,.4);
\draw (.1,1) -- (.35,.2) (.75,0) -- (.5,.5) -- (.5,1) (.45,.4) -- (.5,.5)
(.7,-.05) node[right]{ \scriptsize $c$} (.45,1.1) node[right]{ \scriptsize $d$}
(.15,1) node[left]{ \scriptsize $a$} (.3,0) node[left]{ \scriptsize $b$}
(.45,.55) node[right]{ \scriptsize $v(\ell,2)$};
\fill (.35,.2) circle (1.5pt) (.5,.5) circle (1.5pt);
\end{tikzpicture}\mbox{, } 
\begin{tikzpicture} \useasboundingbox (-.5,.2) rectangle (1.5,1.1);
\draw (.2,0) -- (.5,.6) -- (.5,1) (.6,.4) -- (.5,.6)
(.15,1) node[left]{ \scriptsize $a$} (.3,0) node[left]{ \scriptsize $b$}
(.7,-.05) node[right]{ \scriptsize $c$} (.45,1.1) node[right]{ \scriptsize $d$}
(.45,.65) node[right]{ \scriptsize $v(\ell,2)$};
\draw [draw=white,double=black,very thick] (.1,1) .. controls (.3,.3) .. (.7,.2);
\draw (.1,1) .. controls (.3,.3) .. (.7,.2);
\draw [->] (.8,0) -- (.6,.4);
\fill (.5,.6) circle (1.5pt) (.7,.2) circle (1.5pt);
\end{tikzpicture}
 \mbox{, and} 
\begin{tikzpicture} \useasboundingbox (-.5,.2) rectangle (1.5,1.1);
\draw (.75,0) -- (.5,.35) -- (.25,0) (.5,.45) -- (.5,.35)
(.15,1) node[left]{ \scriptsize $a$} (.3,0) node[left]{ \scriptsize $b$}
(.7,-.05) node[right]{ \scriptsize $c$} (.45,1.15) node[right]{ \scriptsize $d$}
(.45,.4) node[right]{ \scriptsize $v(\ell,2)$};
\draw [->] (.5,1) -- (.5,.45);
\draw [draw=white,double=black,very thick] (.1,1) .. controls (.3,.5) and (.8,.6) .. (.5,.85);
\draw (.1,1) .. controls (.3,.5) and (.8,.6) .. (.5,.85);
\fill (.5,.35) circle (1.5pt) (.5,.85) circle (1.5pt);
\end{tikzpicture}\mbox{.} \end{equation*}
Then the sum $\sum_{\tilde{\Gamma} \suchthat  \tilde{\Gamma} \in \CE(\Gamma;{\finseta})}[\tilde{\Gamma}]$ is zero thanks to the Jacobi relation.
\eop

Now, we get rid of the remaining faces with the help of the STU relation.
\begin{lemma}
\label{lemstu}
The contributions to $(\Zinv_n(1)-\Zinv_n(0))$ or to $(\Zinv_{n,3n-2}(1)-\Zinv_{n,3n-2}(0))$ of the faces $\facee(\finseta,\Link,{\Gamma})$ such that 
\begin{itemize}
\item $\Gamma_{\finseta}$ is an edge between a trivalent vertex and a univalent vertex or
\item $\finseta$ is a pair of consecutive univalent vertices and $\Gamma_{\finseta}$ is not an edge of $\Gamma$                                                                                                                                       \end{itemize}
cancel. More precisely,
let $\Gamma \in \Davis^e_{n,m}(\source)$, let $\finseta$ be a pair of consecutive univalent vertices of $\Gamma$ on a component of $\source$, and assume that $\Gamma_{\finseta}$ is not an edge of $\Gamma$.
Let $\Gamma/\Gamma_{\finseta}$ be the labeled edge-oriented graph obtained from $\Gamma$
by contracting $\Gamma_{\finseta}$ to one point. 
(The labels of the edges of $\;\Gamma/\Gamma_{\finseta}$ belong to $\underline{m} $, $\Gamma/\Gamma_{\finseta}$ has one bivalent vertex injected on $\source$.) 
Let $\CE(\Gamma/\Gamma_{\finseta})$ be the subset of $\Davis^e_{n,m}(\source)$ that contains the graphs $\tilde{\Gamma}$ equipped with a pair $\finseta$ of vertices that are either 
\begin{itemize}
\item two consecutive univalent vertices or
\item a univalent vertex and a trivalent vertex connected by an edge, 
\end{itemize}
such that $\Gamma/\Gamma_{\finseta}$ is equal to $\tilde{\Gamma}/\tilde{\Gamma}_{\finseta}$. If $m=3n$ or if $m=3n-2$, then 
\begin{equation*}\sum_{\tilde{\Gamma} \suchthat  \tilde{\Gamma} \in \CE(\Gamma;{\finseta})}\coefgambet_{\tilde{\Gamma}}I(\tilde{\Gamma},\finseta)=0.\end{equation*}
\end{lemma}
\bp
Note that the face $\facee(\finseta,\Link,\Gamma)$ has two connected components if the only univalent vertices of $\Gamma$ on the component of $\source$ of the univalent vertices of $\finseta$ are the two vertices  of $\finseta$. The two connected components correspond to the two possible linear orders of $\finseta$ at the collapse.

Below, we consider these connected components as two different faces, and a face corresponds to a subset $\finseta$ equipped with a linear order compatible with $i_{\Gamma}$.
In particular, the graph and its face are determined by the labeled edge-oriented graph $\Gamma/\Gamma_{\finseta}$ obtained from $\Gamma$
by contracting ${\finseta}$ to one point, together with a linear order of the two half-edges of the bivalent vertex.
Let $k \in \underline{m} \setminus j_E(E(\Gamma))$. Define $\Gamma_k^+$ (resp. $\Gamma_k^-$) to be the graph in $\Davis^e_{n,m}(\source)$ with an edge $e(k)$ such that $j_E(e(k))=k$, which goes from a univalent vertex $u$ to a trivalent vertex $t$ (resp. from a trivalent vertex $t$ to a univalent vertex $u$) forming a pair $\finseta=\{u,t\}$ such that $\Gamma_k^+/(\Gamma^+_k)_{\finseta}$ (resp. $\Gamma_k^-/(\Gamma^-_k)_{\finseta}$) coincides with $\Gamma/\Gamma_{\finseta}$.

Order the sets of vertices of the $\Gamma_k^{\pm}$ by putting the vertices of $\finseta$ first, with respect to the order induced by the edge orientation (source first), and so that the order of the remaining vertices is the same for all $\Gamma_k^{\pm}$. Let $o(V(\Gamma) \setminus A)$ denote this order of $V(\Gamma) \setminus A$.
For $(\tilde{\Gamma},\finseta) \in \CE(\Gamma/\Gamma_{\finseta})$ such that $\finseta$ is an ordered pair of univalent vertices of $\tilde{\Gamma}$,
order $\finsetv(\tilde{\Gamma})$ by putting the vertices of $\finseta$ first with respect to the linear order induced by the collapse, and next the others with the same order  $o(V(\Gamma) \setminus A)$ as for the $\Gamma_k^{\pm}$.

Let $\phi \colon \RR^3 \to \crats$ be an orientation-preserving diffeomorphism onto a neighborhood of the image of $\finseta$ in a configuration of $\facee(\finseta,\Link,{\Gamma})$. There exists $\rho_{\phi} \colon  \RR^3 \to GL^+(\RR^3)$ such that $T_x\phi\left(\rho_{\phi}(x)\left(\vec{v}\right)\right)=\partau(\phi(x),\vec{v})$, for any $\vec{v}\in (\RR^3=T_x\RR^3)$.
Then for $\Gamma_k^+$, the configuration space is locally diffeomorphic to
$\source \times \crats \times \dots$, where $\source$ contains the position $c(u)=\phi(x)$ of the univalent vertex $u$ of $\finseta$, and $\crats$ contains the position $c(t)$ of the trivalent vertex $t$ of $\finseta$. This position $c(t)=\phi(x+\lambda \rho_{\phi}(x)(\vec{v}))$ is described by a small positive $\lambda$, which plays the role of an inward normal near the collapse, and a vector $\vec{v}$ of $S^2$, which is equal to $\projp_{\partau} \circ p_{e(k)}(c)$ when $\lambda$ reaches $0$. The face is diffeomorphic to
$S^2 \times \source \times \dots$, where the projection onto the factor $S^2$ is $\projp_{\partau} \circ p_{e(k)}$, and the dots contain the coordinates of the remaining vertices, which are the same for all the considered diagrams.

For $\Gamma_k^-$, the configuration space is locally diffeomorphic to
$\crats \times \source \times \dots$, where $\source$ contains the position $c(u)=\phi(x)$ of the univalent vertex $u$ of $\finseta$, and $\crats$ contains the position $c(t)$ of the trivalent vertex $t$ of $\finseta$. This position $c(t)=\phi(x-\lambda \rho_{\phi}(x)(\vec{v}))$ is described by a small positive $\lambda$, which still plays the role of an inward normal near the collapse, and a vector $\vec{v}$ of $S^2$, which is $\projp_{\partau} \circ p_{e(k)}(c)$ when $\lambda$ reaches $0$. The face is again diffeomorphic to
$S^2 \times \source \times \dots$, where the projection onto the factor $S^2$ is $\projp_{\partau} \circ p_{e(k)}$, and the dots contain the coordinates of the remaining vertices, which are the same as for the $\Gamma_k^+$.

For $(\tilde{\Gamma},\finseta) \in \CE(\Gamma/\Gamma_{\finseta})$ such that $\finseta$ is an ordered pair of univalent vertices of $\tilde{\Gamma}$, the configuration space is locally diffeomorphic to
$\source \times \source \times \dots$, where the first $\source$ contains the position $c(u_1)=\phi(x)$ of the first univalent vertex $u_1$ of $\finseta$, and the second $\source$ contains the position $c(u_2)$ of the vertex $u_2$ that follows $u_1$ along $\source$. The position $c(u_2)=\phi(x+\lambda\rho_{\phi}(x)(\vec{t}))$ is again described by a small positive $\lambda$, which plays the role of an inward normal near the collapse, where $\projp_{\partau}$ maps the oriented unit tangent vector to $\Link$ at $c(u_1)$ to $\vec{t} \in S^2$ when $\lambda$ reaches $0$. The face is diffeomorphic to
$\source \times \dots$ with the same notation as before. So the previous faces are the products by $S^2$ of this one,
and $\projp_{\partau} \circ p_{e(k)}$ is the projection to the factor $S^2$. Since the other $p_{e}$ do not depend on this factor $S^2$, $\int_{\left[0,1\right]\times \facee(\finseta,\Link,\Gamma_k^+)}\projp(\Gamma_k^+)^{\ast}(\Omega_{E(\Gamma_k^+)})$ is equal to
\begin{equation*}\int_{(t,\confc) \in\left[0,1\right]\times \facee(\finseta,\Link,\Gamma)}\left(\int_{\{t\} \times\cinjuptd{e(k)}{T_{\confc({\finseta})}\crats}}\tilde{\omega}(k)\right)\projp(\Gamma)^{\ast}(\Omega_{E(\Gamma)}),\end{equation*}
where $\{t\} \times\cinjuptd{e(k)}{T_{\confc({\finseta})}\crats}$ is the factor $S^2$ above.
Furthermore, we have
\begin{equation*}\int_{\{t\} \times\cinjuptd{e(k)}{T_{\confc({\finseta})}\crats}}\tilde{\omega}(k)=1\end{equation*}
since the integral of the closed form $\tilde{\omega}(k)$ over any representative of the homology class of the fiber of the unit tangent bundle of $\crats$ in $\left[0,1\right] \times \partial C_2(\rats) $ is $1$.

This argument, which also works for $\Gamma_k^-$, implies that all the integrals 
$I_0(\tilde{\Gamma},\finseta)=\int_{\left[0,1\right]\times \facee(\finseta,\Link,\tilde{\Gamma})}\projp(\tilde{\Gamma})^{\ast}(\Omega_{E(\tilde{\Gamma})})$ coincide for all the $(\tilde{\Gamma},\finseta) \in \CE(\Gamma/\Gamma_{\finseta})$ equipped with orders of their vertices as above. So it suffices to prove
\begin{equation*}\coefgambet_{\Gamma}\Bigl(\left[\Gamma\right] + \left[\Gamma^{\prime}\right]\Bigr) + \sum_{k \in \underline{m} \setminus j_E(E(\Gamma))}\coefgambet_{\Gamma_k^{+}}\Bigl(\left[\Gamma_k^{+}\right]+\left[\Gamma_k^{-}\right]\Bigr)=0,\end{equation*}
where $\Gamma^{\prime}$ is the graph obtained from $\Gamma$ by permuting the order of the two univalent vertices on $\source$, and all the graphs $\tilde{\Gamma}$ are vertex-oriented so that the vertex-orientation of $H(\tilde{\Gamma})$ induced by the fixed order of the vertices coincides with the edge-orientation of $H(\tilde{\Gamma})$ (as in Remark~\ref{rkedgeorvert}), for $m=3n-2$ and for $m=3n$.

Let $a$ and $b$ denote the half-edges of $\Gamma$ that contain the vertices of $\finseta$.
Without loss of generality, assume that the vertex of $b$ follows the vertex of $a$ on $\source$ for $\Gamma$ (near the connected face). Let $o_V(H(\Gamma) \setminus\{a,b\})$ be an order of $H(\Gamma)\setminus\{a,b\}$, such that the order $(a,b,o_V(H(\Gamma) \setminus\{a,b\}))$ (i.e., $(a,b)$ followed by the elements $H(\Gamma)\setminus\{a,b\}$ ordered by $o_V(H(\Gamma) \setminus\{a,b\})$) induces the edge-orientation of $H(\Gamma)$.

Orient the vertices of $V(\Gamma)\setminus A$ in $\Gamma$ so that $o_V(H(\Gamma)\setminus\{a,b\})$ is induced by the order $o(V(\Gamma) \setminus A)$ and this vertex-orientation.
Let $f$ (resp. $s$) denote the first (resp. second) half-edge of $e(k)$ in $\Gamma^{\pm}_{k}$. Then $(f,s,a,b,o_V(H(\Gamma) \setminus\{a,b\}))$ induces the edge-orientation of $H(\Gamma^{\pm}_{k})$.
Equip the trivalent vertex of $\finseta$ in $\Gamma^{\pm}_{k}$ with the vertex-orientation $( (f \;\mbox{or}\; s),a,b)$, which corresponds to the picture 
\begin{equation}\stuyab,\end{equation} and equip the other vertices of $\Gamma^{\pm}_{k}$ with the same vertex-orientation as their vertex-orientation in $\Gamma$. 
Then the vertex-orientation of $H(\Gamma_k^{+})$ is induced by $(f,s,a,b, o_V(H(\Gamma) \setminus\{a,b\}))$ and coincides with its edge-orientation.
Similarly, the vertex-orientation of $H(\Gamma_k^{-})$ is induced by $(f,a,b,s, o_V(H(\Gamma) \setminus\{a,b\}))$ and coincides with its edge-orientation.
Thus, we have
$\left[\Gamma_k^{+}\right]=\left[\Gamma_k^{-}\right]$ for any $k$, and $\left[\Gamma_k^{+}\right]$ is independent of $k$.

Note that $\left[\Gamma\right]$ looks like \stuxab locally. It coincides with $\left[\Gamma_k^{+}\right]$ outside the pictured part. But $\Gamma^{\prime}$ must be equipped with the opposite vertex-orientation, and $\bigl(-\left[\Gamma^{\prime}\right]\bigr)$ looks like \stuiab.

Thus, it suffices to prove
\begin{equation*}\coefgambet_{\Gamma}\left(\left[\stuxdis\right] -\left[\stuidis\right]\right) + 2(m-\cardlef{E(\Gamma)})\coefgambet_{\Gamma_k^{+}}\left[\stuydis\right]=0\end{equation*}
if $m\in \{3n-2,3n\}$.
With the expression of the $\coefgambet_{\Gamma}$ in Notation~\ref{notationzZ} and Notation~\ref{notlessav}, this equality is equivalent to the STU relation when $m> \cardlef{E(\Gamma)}$. In particular, it is equivalent to the STU relation when $m=3n$.
When $m=3n-2$, if $m = \cardlef{E(\Gamma)}$, then $\Gamma$ has exactly $2$ univalent vertices, so we have $\bigl([\Gamma] + [\Gamma^{\prime}]\bigr) =0$, and the equality is still true.
\eop

Proposition~\ref{propinvone} is now proved. Proposition~\ref{propinvonebis} follows from Proposition~\ref{propinvone} and Lemma~\ref{lemsymeven}. \eopwobp

\bpo{Proof of Theorem~\ref{thmConwaycircZSthree}}
We follow the face cancellations in the proof of Proposition~\ref{propinvone}, to study the effect of changing $\omega_{S^2}(i)=\tilde{\omega}_{S^2}\vert_{\{0\} \times S^2}$ to $\omega^{\prime}_{S^2}(i)=\tilde{\omega}_{S^2}\vert_{\{1\} \times S^2}$, for a closed two-form $\tilde{\omega}_{S^2}$ on $\left[0,1\right] \times S^2$ as in Lemma~\ref{lemformprod}, for some $i \in \underline{3}$. We use the form ${\omega}(i)=p_{\left[0,1\right] \times S^2}^{\ast}(\tilde{\omega}_{S^2})$ on $\left[0,1\right] \times C_2(S^3)$. 
Here, the involved graphs have no looped edges, $4$ vertices, at most $3$ edges, and hence at most one trivalent vertex. They are $\diagcross$, $\diagthetthet$,  and $\diagtripod$.
The only cancellation that requires an additional argument is the cancellation of the faces, for which $\Gamma$ is isomorphic to \diagtripod \; and $\finseta$ is a pair of univalent vertices of $\Gamma$. (The cancellation of Lemma~\ref{lemstu} would involve \diagwheeltwo \, .) In this case, the open face $\facee(\finseta,\Link,\Gamma)$ is the configuration space of $2$ vertices on the knot (one of them stands for the two vertices of $\finseta$) and a trivalent vertex in $\RR^3$. The integral $\int_{\left[0,1\right]\times \facee(\finseta,\Link,\Gamma)}\bigwedge_{e \in E(\Gamma)}p_e^{\ast}\bigl({\omega}(j_E(e))\bigr)$ is the pull-back of a $6$-form on $\left[0,1\right] \times (S^2)^3$, by a map 
whose image is in the codimension $2$ subspace of $\left[0,1\right] \times (S^2)^3$ in which two $S^2$-coordinates coincide. Indeed, the two edges that contain the vertices of $\finseta$ have the same direction. Therefore, we have $I(\Gamma,\finseta)=0$ for these faces, and $w_C \bigl( \Zinvlink_{2,3}\bigl(\RR^3,K,(p_{S^2}^{\ast}(\omega_{S^2}(i)))_{i \in \underline{3}}\bigr)\bigr)$ is independent of the chosen
$\omega_{S^2}(i)$. Conclude with the arguments of Remark~\ref{rknoconflict}.
\eop

\chapter{First properties of \texorpdfstring{$\Zinv$}{Z} and anomalies}
\label{chapanom}

\section{Some properties of \texorpdfstring{$\Zinv(\crats,\Link,\partau)$}{Z(R,L,tau)}}
\label{secinvql}

Lemma~\ref{leminvfstconst} allows us to set 
\begin{equation*}\Zinv_n(\crats,\Link,\partau)=\Zinv_n\Bigl(\crats,\Link,\bigl(\omega(i)\bigr)\Bigr)\end{equation*}
for any collection $(\omega(i))$ of homogeneous propagating forms of $\bigl(C_2(\rats),\partau\bigr)$, under the assumptions of Theorem~\ref{thmfstconst}.
We still have to study how $\Zinv_n(\crats,\Link,\partau)$ varies when $\partau$ varies inside its homotopy class when $\Link \neq \emptyset$, but the naturality of the construction of $\Zinv_n$ already implies the following proposition.

\begin{proposition}
\label{propinvdiffnat}
Let $\rats$ be the $\QQ$-sphere equipped with its neighborhood $\mathring{B}_{1,\infty}$ of $\infty$ of the beginning of Section~\ref{secCtwo}.
Let $\psi$ be an orientation-preserving diffeomorphism from $\rats$ to $\psi(\rats)$. Use the restriction of $\psi$ to the ball $\mathring{B}_{1,\infty}$ as an identification of $\mathring{B}_{1,\infty}$ with a neighborhood of $\psi(\infty)$ in $\psi(\rats)$. Define $\psi_{\ast}(\partau)=T\psi \circ \partau \circ (\psi^{-1} \times \id_{\RR^3})$.
Then we have \begin{equation*}\Zinv_n\bigl(\psi(\crats),\psi(\Link),\psi_{\ast}(\partau)\bigr)=\Zinv_n\bigl(\crats,\Link,\partau\bigr)\end{equation*} for all $n \in \NN$,
where  $p_1(\psi_{\ast}(\partau))=p_1(\partau)$.
\end{proposition}
\bp The diffeomorphism $\psi$ induces natural diffeomorphisms $\psi_{\ast}$ from $C_2(\rats)$ 
to $C_2(\psi(\rats))$, and from the $\check{C}(\rats,\Link;\Gamma)$ to the $\check{C}(\psi(\rats),\psi(\Link);\Gamma)$. If $\omega$ is a homogeneous propagating form of $\bigl(C_2(\rats),\partau\bigr)$, then
$(\psi_{\ast}^{-1})^{\ast}(\omega)$ is a homogeneous propagating form of $(C_2(\psi(\rats)),\psi_{\ast}(\partau))$ since the restriction of $(\psi_{\ast}^{-1})^{\ast}(\omega)$ to $\ST \psi(\crats)$ is $(T\psi ^{-1})^{\ast}(\projp_{\partau}^{\ast}(\omega_{S^2}))=(\projp_{\partau}\circ T\psi ^{-1})^{\ast}(\omega_{S^2})=\projp_{\psi_{\ast}(\partau)}^{\ast}(\omega_{S^2})$. 
For any Jacobi diagram $\Gamma$ on the domain $\sourcetl$ of $\Link$, we have
\begin{equation*}\begin{array}{ll}I\bigl(\psi(\rats),\psi(\Link),\Gamma,(\psi_{\ast}^{-1})^{\ast}(\omega)\bigr)&=\int_{\check{C}(\psi(\rats),\psi(\Link);\Gamma)} \bigwedge_{e \in E(\Gamma)}p_e^{\ast}\bigl((\psi_{\ast}^{-1})^{\ast}(\omega)\bigr)\\
  &=\int_{\check{C}(\psi(\rats),\psi(\Link);\Gamma)} (\psi_{\ast}^{-1})^{\ast}\left(\bigwedge_{e \in E(\Gamma)}p_e^{\ast}(\omega)\right)\\
&=I\left(\rats,\Link,\Gamma,\omega\right),
  \end{array}
 \end{equation*}
  where  $\Gamma$ is equipped with an implicit orientation $o(\Gamma)$.
Therefore, we have \begin{equation*}\Zinv_n\bigl(\psi(\crats),\psi(\Link),\psi_{\ast}(\partau)\bigr)=\Zinv_n(\crats,\Link,\partau)\end{equation*} for all $n \in \NN$.
\eop

We study some other properties of $\Zinv_n(\crats,\Link,\partau)$.

Let $(\crats,\partau)$ be an asymptotic rational homology $\RR^3$.
Thanks to Lemma~\ref{leminvwoboun},
$\Zinv_n(\crats,\emptyset,\partau)$
depends only on the homotopy class of $\partau$ for any integer $n$.
Set $\Zinv_n(\rats,\partau)=\Zinv_n(\crats,\emptyset,\partau)$ and $\Zinv(\rats,\partau)=\bigl(\Zinv_n(\rats,\partau)\bigr)_{n \in \NN}$.
Using Notation~\ref{notationzZ}, let
$\zinv_n(\rats,\partau)=p^c\bigl(\Zinv_n(\rats,\partau)\bigr)$ be the connected part of $\Zinv_n(\rats,\partau)$, and set $\zinv(\rats,\partau)=\bigl(\zinv_n(\rats,\partau)\bigr)_{n \in \NN}$.\index[N]{ZZ@$\Zinvuf$ and some variants (see also the summary in the next pages)!za@$\zinv$}

We give a direct elementary proof of the following proposition, which could also be proved in the same way as Corollary~\ref{corlogZinv} below.

\begin{proposition}
\label{propexpon}
For any propagating form $\omega$ of $C_2(\rats)$, we have
\begin{equation*}\Zinv(\rats,\omega)=\exp\bigl(\zinv(\rats,\omega)\bigr).\end{equation*}
In particular, for any asymptotically standard parallelization $\partau$ of $\rats$, we have
\begin{equation*}\Zinv(\rats,\partau)=\exp\bigl(\zinv(\rats,\partau)\bigr).\end{equation*}
\end{proposition}
\bp Let $\Gamma$ be a trivalent Jacobi diagram whose components are isomorphic to some $\Gamma_i$ for $i\in \underline{r}$ and such that
$\Gamma$ has $k_i$ connected components isomorphic to $\Gamma_i$.
Then we have \begin{equation*}I\bigl(\rats,\emptyset,\Gamma,(\omega)_{i \in \underline{3n}}\bigr)\left[\Gamma\right]=\prod_{i=1}^r\left(I(\rats,\emptyset,\Gamma_i,(\omega)_{i \in \underline{3\deg(\Gamma_i)}})^{k_i}\left[\Gamma_i\right]^{k_i}\right)\end{equation*}
and $\cardlef{\Aut(\Gamma)} = \prod_{i=1}^r\left(k_i!(\Aut(\Gamma_i))^{k_i}\right)$. Conclude with Proposition~\ref{propdefhomog}.
\eop

Recall the coproduct maps $\Delta_n$ defined in Section~\ref{seccoprod}.
\begin{proposition}
\label{propgrouplike} Under the assumptions of Theorem~\ref{thmfstconst}, we have
\begin{equation*}\Delta_n\bigl(\Zinv_n(\crats,\Link,\partau)\bigr)= \sum_{i=0}^n \Zinv_i(\crats,\Link,\partau) \otimes \Zinv_{n-i}(\crats,\Link,\partau).\end{equation*}
\end{proposition}
\bp Let $T_i=\Zinv_i(\rats,\Link,\partau) \otimes \Zinv_{n-i}(\rats,\Link,\partau)$. We have

\begin{equation*}T_i=\sum \frac1{\cardlef{\Aut(\Gamma^{\prime})}}I\bigl(\Gamma^{\prime},(\omega)_{j \in \underline{3i}}\bigr)\frac1{\cardlef{\Aut(\Gamma^{\prime \prime})}}I\bigl(\Gamma^{\prime \prime},(\omega)_{j \in \underline{3(n-i)}}\bigr)\left[\Gamma^{\prime}\right] \otimes [\Gamma^{\prime \prime}],\end{equation*}
where the sum runs over the pairs $(\Gamma^{\prime},\Gamma^{\prime \prime})$ in $\Davis^u_i\bigl(\source(\Link)\bigr) \times \Davis^u_{n-i}\bigl(\source(\Link)\bigr)$.
Use Remark~\ref{rkhomog} to view the summands as a measure of configurations of graphs $\Gamma^{\prime} \sqcup \Gamma^{\prime \prime}$ (which may correspond to several elements of $\Davis^u_n(\source(\Link))$) 
together with a choice of an embedded subgraph $\Gamma^{\prime}$.
\eop

\begin{corollary}
\label{corlogZinv} Recall the projection $p^c$ of Corollary~\ref{corprojprim} from $\Assis(S^1)$ to the space $\Assis^c(S^1)$ of its primitive elements. 
If $\Link$ has one component, 
set \begin{equation*}\zinvlink(\rats,\Link,\partau)=p^c\bigl(\Zinvlink(\crats,\Link,\partau)\bigr).\index[N]{ZZ@$\Zinvuf$ and some variants (see also the summary in the next pages)!zb@$\zinvlink$}\end{equation*}
Then we have \begin{equation*}\Zinvlink(\crats,\Link,\partau)=\exp\bigl(\zinvlink(\rats,\Link,\partau)\bigr).\end{equation*}
\end{corollary}
\bp This is a direct consequence of Lemma~\ref{lemAHopf}, Proposition~\ref{propgrouplike} and Theorem~\ref{thmgrouplikexp}.\eop

\section{On the anomaly \texorpdfstring{$\ansothree$}{beta}}
\label{secansothree}

We now study how $\zinv_n(\rats,\partau)$, which is defined before Proposition~\ref{propexpon}, depends on $\partau$.

\begin{definition}
\label{defxin}
 Let $\rhomap \colon (B^3,\partial B^3) \rightarrow (SO(3),1)$ be the map of Definition~\ref{defrho}. Extend it to $\RR^3$ by considering $B^3$ as the unit ball of $\RR^3$ and by letting $\rho$ map $(\RR^3 \setminus B^3)$ to $1$.
Consider the parallelization $\taust \circ \psi_{\RR}(\rho)$, where $\psi_{\RR}(\rho)(x,v)=\bigl(x,\rho(x)(v)\bigr)$. Set
\begin{equation*}\ansothree_n=\zinv_n\bigl(S^3,\taust \circ \psi_{\RR}(\rho)\bigr).\end{equation*}\index[N]{bzeta@$\ansothree$ anomaly}
\end{definition}

\begin{proposition}
\label{propdefoneanom}
Let $(\crats,\partau_0)$ be an asymptotic rational homology $\RR^3$, and let $\partau_1$ be a parallelization of $\crats$ that coincides with $\partau_0$ outside ${\ballb}_{\rats}$.
Then we have
\begin{equation*}\zinv_n(\rats,\partau_1) - \zinv_n(\rats,\partau_0) =\frac{p_1(\partau_1)-p_1(\partau_0)}{4}\ansothree_n\end{equation*} for any integer $n$.
\end{proposition}

Proposition~\ref{propdefoneanom} is an easy consequence of Proposition~\ref{propdeftwoanom} below. The latter looks more complicated,
but 
it is very useful since it offers more practical definitions of the \indexT{anomaly} \begin{equation*}\ansothree=(\ansothree_n)_{n\in \NN}\end{equation*}
when applied to $(\crats,\partau_0,\partau_1)=\bigl(\RR^3,\taust,\taust \circ \psi_{\RR}(\rho)\bigr)$ (and to the case in which $\tilde{\omega}_{S^2}(i)$ is the pull-back of $\omega_{0,S^2}(i)$ under the natural projection from $\left[0,1\right] \times S^2$ to $S^2$).

\begin{proposition}
\label{propdeftwoanom}
Let $(\crats,\partau_0)$ be an asymptotic rational homology $\RR^3$, and let $\partau_1$ be a parallelization of $\crats$ that coincides with $\partau_0$ outside ${\ballb}_{\rats}$.
For $i \in \underline{3n}$, let $\omega_{0,S^2}(i)$ and $\omega_{1,S^2}(i)$ be two volume-one forms on $S^2$.
Then there exists a closed two-form $\tilde{\omega}_{S^2}(i)$ on $\left[0,1\right] \times S^2$ such that the restriction 
of $\tilde{\omega}_{S^2}(i)$ to $\{t\} \times S^2$ is ${\omega}_{t,S^2}$ for $t \in \{0,1\}$.
For any such forms $\tilde{\omega}_{S^2}(i)$, there exist closed $2$-forms $\tilde{\omega}(i)$ on $\left[0,1\right] \times \ST \crats$
such that
\begin{itemize}
 \item the restriction  
of $\tilde{\omega}(i)$ to $\{t\} \times  \ST \crats$ is 
$p_{\partau_t}^{\ast}\bigl({\omega}_{t,S^2}(i)\bigr)$ for $t \in \{0,1\}$,
\item the restriction of $\tilde{\omega}(i)$ to $\left[0,1\right] \times  \bigl(\ST (\crats \setminus {\ballb}_{\rats}) \bigr)$ is $\bigl(\id_{\left[0,1\right]} \times p_{\partau_0}\bigr)^{\ast}\bigl(\tilde{\omega}_{S^2}(i)\bigr)$.
\end{itemize}
Then we have
\begin{multline*}\zinv_n(\rats,\partau_1) - \zinv_n(\rats,\partau_0)=\zinv_n\left(\left[0,1\right] \times \ST {\ballb}_{\rats}; \bigl(\tilde{\omega}(i)\bigr)_{i \in \underline{3n}}\right)\\
=\sum_{\Gamma \in \Davis^c_n}\coefgambet_{\Gamma}\int_{\left[0,1\right] \times \cinjuptd{\finsetv(\Gamma)}{T{\ballb}_{\rats}}}\bigwedge_{e \in E(\Gamma)}p_e^{\ast}\bigl(\tilde{\omega}\bigl(j_E(e)\bigr)\bigr)\left[\Gamma\right]\\
=\frac{p_1(\partau_1)-p_1(\partau_0)}{4}\ansothree_n,\end{multline*}
and $\ansothree_n=0$ if $n$ is even.
(Recall that the orientation of $\cinjuptd{\finsetv(\Gamma)}{T\crats}$ is defined in Lemma~\ref{lemoriface}.)
\end{proposition}
\bp 
The existence of $\tilde{\omega}_{S^2}(i)$ comes from Lemma~\ref{lemformprod}. In order to prove the existence of $\tilde{\omega}(i)$, which is defined on $\partial \left(\left[0,1\right]\times S^2 \times {\ballb}_{\rats}\right)$ by the conditions, we need to extend it to $\left[0,1\right]\times S^2 \times {\ballb}_{\rats}$. The obstruction belongs to 
\begin{equation*}H^3\left(\left[0,1\right]\times S^2 \times {\ballb}_{\rats},\partial (\left[0,1\right]\times S^2 \times {\ballb}_{\rats})\right) \cong H_3(\left[0,1\right]\times S^2 \times {\ballb}_{\rats}),\end{equation*}
which is trivial. So $\tilde{\omega}(i)$ extends as desired.
In order to prove that the first equality is a consequence of Corollary~\ref{corinvone}, extend the forms $\tilde{\omega}(i)$ of the statement to $\left[0,1\right] \times C_2(\rats)$ as forms that satisfy the conditions in Proposition~\ref{propinvone}. First extend the $\tilde{\omega}(i)$ to $\left[0,1\right]\times (\partial C_2(\rats) \setminus \ST{\ballb}_{\rats}) $ as $(\id_{\left[0,1\right]} \times p_{\partau_0})^{\ast}(\tilde{\omega}_{S^2}(i))$. Next extend the restriction of $\tilde{\omega}(i)$
to $\{0\} \times \partial C_2(\rats)$ (resp. to $\{1\} \times \partial C_2(\rats)$) on $\{0\} \times C_2(\rats)$ (resp. on $\{1\} \times C_2(\rats)$) as a propagating form of $\bigl(C_2(\rats),\partau_0\bigr)$ (resp. of $\bigl(C_2(\rats),\partau_1\bigr)$) as in Section~\ref{secprop}.
Thus, $\tilde{\omega}(i)$ is defined consistently on $\partial (\left[0,1\right]\times C_2(\rats))$, and it extends as a closed form that satisfies the assumptions in Proposition~\ref{propinvone} as in Lemma~\ref{lemformprod}. 
Corollary~\ref{corinvone}, Lemma~\ref{leminvwoboun}, and Lemma~\ref{lemsymeven} yield 
\begin{equation*}\zinv_n(\rats,\partau_1) - \zinv_n(\rats,\partau_0)
=\sum_{\Gamma \in \Davis^c_n}\frac{1}{(3n)!2^{3n}}\int_{\left[0,1\right] \times \cinjuptd{\finsetv(\Gamma)}{T{\ballb}_{\rats}}}\bigwedge_{e \in E(\Gamma)}p_e^{\ast}\Bigl(\tilde{\omega}\bigl(j_E(e)\bigr)\Bigr)\left[\Gamma\right],\end{equation*}
which is zero 
if $n$ is even. So everything is proved when $n$ is even. Assume that $n$ is odd.

There exists a map $g \colon (\crats, \crats \setminus {\ballb}_{\rats}) \to (SO(3),1)$ such that
$\partau_1=\partau_0 \circ \psi_{\RR}(g)$.
Using $\partau_0$ to identify $\cinjuptd{\finsetv(\Gamma)}{T\crats}$ with $\crats \times \cinjuptd{\finsetv(\Gamma)}{\RR^3}$ makes clear that
$\left(\zinv_n(\rats,\partau_0 \circ \psi_{\RR}(g)) - \zinv_n(\rats,\partau_0)\right)$ does not depend on $\partau_0$.
For any $g \colon (\crats, \crats \setminus {\ballb}_{\rats}) \to (SO(3),1)$, set $\zinv_n^{\prime}(g)=\zinv_n(\rats,\partau_0 \circ \psi_{\RR}(g)) - \zinv_n(\rats,\partau_0)$.
Then $\zinv_n^{\prime}$ is a homomorphism from $\left[\left({\ballb}_{\rats},\partial{\ballb}_{\rats}\right),(SO(3),1)\right]$ to the vector space $\Aavis^c_n(\emptyset)$ over $\RR$.
Theorem~\ref{thmpone} and Lemma~\ref{lemdegrho} yield \begin{equation*}\zinv_n^{\prime}(g)=\frac{\deg(g)}{2}\zinv_n^{\prime}\bigl(\rho_{{\ballb}_{\rats}}(B^3)\bigr).\end{equation*}
It is easy to see that $\zinv_n^{\prime}\bigl(\rho_{{\ballb}_{\rats}}(B^3)\bigr)$ is independent of $\crats$. Since Example~\ref{examplecomconfinttwocont} shows $ \zinv_n(S^3,\partau_s) =0$, we have $\zinv_n^{\prime}\bigl(\rho_{{\ballb}_{\rats}}(B^3)\bigr)=\ansothree_n$ by Definition~\ref{defxin}. Recall $p_1(\partau_0 \circ \psi_{\RR}(g))- p_1(\partau_0)=2\deg(g)$ from Theorem~\ref{thmpone}.
\eop

\begin{remark}
 The anomaly $\ansothree$ is the opposite of the constant $\xi$ defined in \cite[Section 1.6]{lesconst}.
\end{remark}

\begin{corollary}
\label{cordefzinvuf}
Let $(\crats,\partau)$ be an asymptotic rational homology $\RR^3$, then
\begin{equation*}\zinv_n(\rats,\partau) - \frac{p_1(\partau)}{4}\ansothree_n\end{equation*} is independent of $\partau$.
Set $\zinvuf_n(\rats)=\zinv_n(\rats,\partau) - \frac{p_1(\partau)}{4}\ansothree_n$, $\zinvuf(\rats)=(\zinvuf_n(\rats))_{n \in \NN}$, and 
\begin{equation*}\Zinvuf(\rats)=\exp(\zinvuf(\rats)).\end{equation*}
Then \begin{equation*}\Zinvuf(\rats)=\Zinv(\rats,\partau)\exp\left(- \frac{p_1(\partau)}{4}\ansothree\right)\end{equation*}
is the invariant $\Zinvuf(\rats, \emptyset)$ that was announced in Theorem~\ref{thmfstconst}.
\end{corollary}

\bp See Proposition~\ref{propexpon}.\eop

\begin{proposition}
\label{propanodegone} The degree-one part of the anomaly $\ansothree$ is
 $\ansothree_1=\frac{1}{12}\left[ \tata \right]$.
\end{proposition}
\bp According to Proposition~\ref{propthetazoneone}, we have $\zinv_1(\rats,\partau)=\frac{\Theta(\rats,\partau)}{12}\left[ \tata \right]$.
Proposition~\ref{propdefoneanom} implies $\zinv_1(\rats,\partau_1) - \zinv_1(\rats,\partau_0) =\frac{p_1(\partau_1)-p_1(\partau_0)}{4}\ansothree_1$, while Corollary~\ref{corThetap} implies
 \begin{equation*}\Theta(\rats,\partau_1)-\Theta(\rats,\partau_0) = \frac1{4}\bigl(p_1(\partau_1)-p_1(\partau_0)\bigr).\end{equation*}
\eop

\begin{corollary}
\label{corthetazone}
Let $(\crats,\partau)$ be an asymptotic rational homology $\RR^3$, then we have
$\Zinv_1(\rats,\partau)=\zinv_1(\rats,\partau)=\frac{\Theta(\rats,\partau)}{12}\left[ \tata \right]$ and
\begin{equation*}\Zinvuf_1(\rats)=\zinvuf_1(\rats)=\frac{\Theta(\rats)}{12}\left[ \tata \right]\end{equation*}
in $\Aavis_1(\emptyset)=\RR\left[\tata\right]$.
\end{corollary}
\bp The first equality is a direct consequence of  Proposition~\ref{propthetazoneone}. The second one follows from Corollary~\ref{corThetap}, Corollary~\ref{cordefzinvuf}, and Proposition~\ref{propanodegone}.
\eop

\begin{remark} According to Proposition~\ref{propdeftwoanom}, the even-degree part of the anomaly $\ansothree$ vanishes.
 The values of $\ansothree_{2n+1}$ are unknown when $n \geq 1$. We may hope them to be zero, but I do not know any conjecture about them.
\end{remark}

\section{On the anomaly \texorpdfstring{$\alpha$}{alpha}}
\label{secanomalpha}

We define the \indexT{anomaly} \begin{equation*}\alpha=(\alpha_n)_{n\in \NN},\end{equation*}
which is sometimes called the \emph{Bott--Taubes anomaly,} below.
Let $v \in S^2$. Let $D_v$ denote the linear map
\begin{equation*}\begin{array}{llll}D_v\colon &\RR &\longrightarrow& \RR^3\\
& 1&\mapsto&v.\end{array}\end{equation*} 
Let $\CD^c_n(\RR)$ denote the set of degree $n$, connected, $\underline{3n-2}$-numbered diagrams on $\RR$ with at least one univalent vertex, without looped edges. 
As in Definition~\ref{defnumdia}, a degree $n$ diagram $\check{\Gamma}$ is \emph{numbered} if the edges of $\check{\Gamma}$ are oriented and if $E(\check{\Gamma})$ is equipped with an injection $j_E \colon E(\check{\Gamma}) \hookrightarrow \underline{3n-2}$,
which numbers its edges.
Let $\check{\Gamma} \in \CD^c_n(\RR)$.
Define $\check{C}(D_v;\check{\Gamma})$ as in Section~\ref{secdefconfspace}, where the line $D_v$ of $\RR^3$ replaces the link $L$ of $\crats$, and $\RR$ replaces the domain $\source$.
Let $\check{\anomq}(v;\check{\Gamma})$ be the quotient of $\check{C}(D_v;\check{\Gamma})$ by the translations parallel to $D_v$ and by the dilations. Then $\check{\anomq}(v;\check{\Gamma})$ is the space denoted by $\cinjuptdanvec(\RR^3,v;\check{\Gamma})$ before Lemma~\ref{lemfacfacone}.
Let ${\anomq}(v;\check{\Gamma})$ \index[N]{Qv@${\anomq}(v;\check{\Gamma})$ configuration space} denote the closure of $\check{\anomq}(v;\check{\Gamma})$ in $\ccompuptd{\finsetv(\check{\Gamma})}{\RR^3}$.
According to Lemma~\ref{lemfacfacone}, the space ${\anomq}(v;\check{\Gamma})$, which coincides with $\ccompuptdanvec(\RR^3,v;\check{\Gamma})$, is a compact smooth manifold with ridges.

To each edge $e$ of $\check{\Gamma}$, associate a map $p_{e,S^2}$, which maps a configuration of $\check{\anomq}(v;\check{\Gamma})$ to the direction of the vector from the origin of $e$ to its end in $S^2$. This map
extends to ${\anomq}(v;\check{\Gamma})$ according to Theorem~\ref{thmcompfacanomtwo}.

 Now, define $\check{\anomq}(\check{\Gamma})$ (resp. ${\anomq}(\check{\Gamma})$) as the total space
of the fibration over $S^2$ whose fiber over $v$ is $\check{\anomq}(v;\check{\Gamma})$ (resp. ${\anomq}(v;\check{\Gamma})$). The configuration space $\check{\anomq}(\check{\Gamma})$ and its compactification ${\anomq}(\check{\Gamma})$ carry natural smooth structures.
The configuration space $\check{\anomq}(\check{\Gamma})$ is oriented as follows, when a vertex-orientation $o(\check{\Gamma})$ is given. Equip $\check{C}(D_v;\check{\Gamma})$ with its orientation induced by Corollary~\ref{cororc}, as before. Orient $\check{\anomq}(v;\check{\Gamma})$ so that $\check{C}(D_v;\check{\Gamma})$ is locally homeomorphic to the oriented product
(translation vector $z$ in $\RR v$, ratio of homothety $\lambda \in \left]0,\infty\right[$) $\times \check{\anomq}(v;\check{\Gamma})$. Next orient $\check{\anomq}(\check{\Gamma})$ with the $(\mbox{base} (=S^2) \oplus \mbox{fiber})$ convention. (We can summarize this by saying that the $S^2$-coordinates replace $(z,\lambda)$.)

\begin{proposition}
\label{propanom}
For $i \in \underline{3n-2}$, let $\omega(i,S^2)$ be a volume-one form of $S^2$. Define 
\begin{equation*}I\bigl(\check{\Gamma},o(\check{\Gamma}),\omega(i,S^2)\bigr)=\int_{\check{\anomq}(\check{\Gamma})}\bigwedge_{e \in E(\check{\Gamma})}p_{e,S^2}^{\ast}\bigl(\omega(j_E(e),S^2)\bigr).\end{equation*}
Set
\begin{equation*}\alpha_n=\frac12\sum_{\check{\Gamma} \in \CD^c_n(\RR)}\coefgambet_{\check{\Gamma}}I\bigl(\check{\Gamma},o(\check{\Gamma}),\omega(i,S^2)\bigr)\left[\check{\Gamma},o(\check{\Gamma})\right] \;\; \in \mathcal{A}(\RR),\end{equation*}
where $\coefgambet_{\check{\Gamma}}=\frac{(3n-2-\cardlef{E(\check{\Gamma})})!}{(3n-2)!2^{\cardlef{E(\check{\Gamma})}}}$.
Then $\alpha_n$ does not depend on the chosen $\omega(i,S^2)$, we have $\alpha_1= \frac12\left[ \onechordsmalltseul \right]$, and $\alpha_{2n}=0$ for all $n$.
\end{proposition}
\bp Let us first prove that $\alpha_n$ does not depend on the chosen $\omega(i,S^2)$,
by proving that its variation vanishes when $\omega(i,S^2)$ is changed to some $\tilde{\omega}(i,1,S^2)$.
According to Lemma~\ref{lemformprod}, there exists a closed $2$-form
$\tilde{\omega}(i,S^2)$ on $\left[0,1\right] \times S^2$ whose restriction to $\{0\}\times S^2$ is $\omega(i,S^2)=\tilde{\omega}(i,0,S^2)$ and whose restriction to $\{1\}\times S^2$ is $\tilde{\omega}(i,1,S^2)$.
According to Stokes' theorem, for any $\check{\Gamma} \in \CD^c_n(\RR)$, we have
\begin{multline*}I\bigl(\check{\Gamma},o(\check{\Gamma}),(\tilde{\omega}(i,1,S^2))_{i \in \underline{3n-2}}\bigr)-I\bigl(\check{\Gamma},o(\check{\Gamma}),(\tilde{\omega}(i,0,S^2))_{i \in \underline{3n-2}}\bigr)\\=\sum_F\int_{\left[0,1\right]\times F} \bigwedge_{e \in E(\check{\Gamma})}p_{e,S^2}^{\ast}\Bigl(\tilde{\omega}\bigl(j_E(e),S^2\bigr)\Bigr),\end{multline*}
where $p_{e,S^2}\colon \left[0,1\right]\times \check{\anomq}(\check{\Gamma}) \to \left[0,1\right]\times S^2$ denotes the product by $\id_{\left[0,1\right]}$ of $p_{e,S^2}$,
and the sum runs over the codimension-one faces $F$ of ${\anomq}(\check{\Gamma})$.
These faces fiber over $S^2$, and the fibers over $v \in S^2$ are the codimension-one faces
$f(\finseta,v;\check{\Gamma})$ of ${\anomq}(v,\check{\Gamma})=\ccompuptdanvec(\RR^3,v;\check{\Gamma})$
for the strict subsets $\finseta$ of $\finsetv(\check{\Gamma})$ with cardinality at least $2$ whose univalent vertices are consecutive on $\RR$, as in Lemma~\ref{lemfacfacan}. Let $\facee(\finseta,\check{\Gamma})$ denote the face with fiber $f(\finseta,v;\check{\Gamma})$.
Now, it suffices to prove that the contributions of all the $\facee(\finseta,\check{\Gamma})$ vanish.

When the product of all the $p_{e,S^2}$ factors through a quotient of $\left[0,1\right]\times \facee(\finseta,\check{\Gamma})$ of smaller dimension, the face $\facee(\finseta,\check{\Gamma})$ does not contribute. This allows us to get rid of 
\begin{itemize}
 \item the faces $\facee(\finseta,\check{\Gamma})$ for which $\check{\Gamma}_{\finseta}$ is not connected and ${\finseta}$ is not a pair of univalent vertices of $\check{\Gamma}$, as in Lemma~\ref{lemdiscon}, and
\item the faces $\facee(\finseta,\check{\Gamma})$ for which $\cardlef{\finseta} \geq 3$ and $\check{\Gamma}_{\finseta}$ has a univalent vertex that was trivalent in $\check{\Gamma}$, as in Lemma~\ref{lemedge}.
\end{itemize}
We also have faces that cancel each other for graphs that are identical outside their $\check{\Gamma}_{\finseta}$ part.
\begin{itemize}
 \item The faces $\facee(\finseta,\check{\Gamma})$ (which are not already listed) such that $\check{\Gamma}_{\finseta}$ has at least one bivalent vertex cancel (by pairs) by the parallelogram identification as in Lemma~\ref{lemsym}.
\item The faces $\facee(\finseta,\check{\Gamma})$ for which $\check{\Gamma}_{\finseta}$ is an edge between two trivalent vertices, cancel by triples, thanks to the Jacobi relation as in Lemma~\ref{lemihx}.
\item Similarly, two faces for which ${\finseta}$ is a pair of (necessarily consecutive in $\RR$) univalent vertices of $\check{\Gamma}$, cancel $(3n-2-\cardlef{E(\Gamma)})$ faces $F(\check{\Gamma}^{\prime},{\finseta}^{\prime})$ for which $\check{\Gamma}^{\prime}_{{\finseta}^{\prime}}$ is an edge between a univalent vertex of $\check{\Gamma}$ and a trivalent vertex of $\check{\Gamma}$, thanks to the STU relation (and to Lemma~\ref{lemoneleg}) as in Lemma~\ref{lemstu}.
\end{itemize}
Here, there are no faces left, and $\alpha_n$ does not depend on the chosen $\omega(i,S^2)$.

The computation of $\alpha_1$ is straightforward.

Let us prove that $\alpha_{n}=0$ for any even $n$.
Let $\check{\Gamma}$ be a numbered graph, and let ${\check{\Gamma}}^{eo}$ be obtained from $\check{\Gamma}$ by reversing the orientations of the $\cardlef{E}$ edges of $\check{\Gamma}$.
Consider the map $r$ from $\check{\anomq}({\check{\Gamma}}^{eo})$ to $\check{\anomq}({\check{\Gamma}})$ that composes a configuration by the multiplication by $(-1)$ in $\RR^3$. It sends a configuration over $v \in S^2$ to a configuration over $(-v)$. It is a fibered space map over the orientation-reversing antipode of $S^2$. Equip $\check{\Gamma}$ and ${\check{\Gamma}}^{eo}$ with the same vertex-orientation and with the same orders on their vertex sets. Then our map $r$ is orientation-preserving if and only if $\cardlef{T(\check{\Gamma})} +1$ is even. The vertex-orientations of $H({\check{\Gamma}})$ and $H({\check{\Gamma}}^{eo})$ can be consistent with both of the edge-orientations of $H({\check{\Gamma}})$ and $H({\check{\Gamma}}^{eo})$ if and only if $\cardlef{E(\check{\Gamma})}$ is even. We have $p_{e,S^2,{\check{\Gamma}}^{eo}}=p_{e,S^2,{\check{\Gamma}}} \circ r$ for all the edges $e$ of ${\check{\Gamma}}^{eo}$. Since $\cardlef{E(\check{\Gamma})} =n+\cardlef{T(\check{\Gamma})}$, we get \begin{equation*}I\bigl({\check{\Gamma}}^{eo},\omega(i,S^2)\bigr)\left[{\check{\Gamma}}^{eo}\right]=(-1)^{n+1}I\bigl({\check{\Gamma}},\omega(i,S^2)\bigr)\left[\check{\Gamma}\right].\end{equation*}
\eop

\begin{note} \label{rkalphafive}
It is known that $\alpha_3=0$ \cite[Proposition 1.4]{poirier}. Sylvain Poirier also found that $\alpha_5=0$ with the help of a Maple program.
Furthermore, according to \cite[Corollary 1.4]{lesunikon}, $\alpha_{2n+1}$ is a combination of diagrams with two univalent vertices, and $\Zinvuf(S^3,\Link)=\Zinvlinkuf(S^3,\Link)$ is obtained from the Kontsevich integral $Z^K$ by inserting $d$
times the plain part of $2\alpha$ on each degree $d$ connected component of a diagram. A similar statement is valid for the functorial extension of $\Zinvuf$ to tangles described in the book's third part. See Corollary~\ref{corlesunikontwoleg} and Lemma~\ref{lemalphatwoleg}, which provides an alternative definition of $\alpha$. It has been conjectured that $\alpha_n=0$ for $n>1$ and incorrect proofs have circulated.
\end{note}

\section{Dependence on the forms for straight links}
\label{secdepstraight}

In this section, we prove Theorems~\ref{thmconststraight} and \ref{thmconststraightred}.
In order to do it,
we will prove the following lemma.

\begin{lemma}
 \label{lemconststraightred}
 Under the assumptions of Theorem~\ref{thmconststraight}, the element \begin{equation*}\Zinv_{n,3n-2}\bigl(\crats,\Link,(\omega(i))_{i \in \underline{3n-2}}\bigr)=\sum_{\Gamma \in \Davisred^e_{n,3n-2}(\source)}\coefgambetred_{\Gamma}I\bigl(\rats,\Link,\Gamma,(\omega(i))_{i \in \underline{3n-2}}\bigr)\left[\Gamma\right]\end{equation*}
of $\Aavis_n(\source)$ is independent of the chosen propagating forms
$\omega(i)$ of $\bigl(C_2(\rats),\partau\bigr)$. 
\end{lemma}

\begin{lemma}
 \label{lemconststraightredimp}
Lemma~\ref{lemconststraightred} implies Theorem~\ref{thmconststraight} and Theorem~\ref{thmconststraightred}. 
\end{lemma}
\bp
Let us prove that Lemma~\ref{lemconststraightred} implies Theorem~\ref{thmconststraight}.
$\Aavis_n(\source)=\Aavis_n(\emptyset) \oplus \Aavis^{\prime}_n(\source)$, where $\Aavis^{\prime}_n(\source)$ is the subspace of $\Aavis_n(\source)$ generated by the Jacobi diagrams with at least one univalent vertex. Since we know from
Lemma~\ref{leminvwoboun}
that $\Zinv_{n}(\crats,\emptyset,(\omega(i))_{i \in \underline{3n}})$ is independent of the chosen propagating forms
$\omega(i)$ of $\bigl(C_2(\rats),\partau\bigr)$,
we focus on the projection $\Zinv^{\prime}_{n}\bigl(\crats,\Link,(\omega(i))_{i \in \underline{3n}}\bigr)$
of \begin{equation*}\Zinv_{n}\Bigl(\crats,\Link,\bigl(\omega(i)\bigr)_{i \in \underline{3n}}\Bigr)=\sum_{\Gamma \in \Davis^e_{n}(\source)}\coefgambet_{\Gamma}I\Bigl(\rats,\Link,\Gamma,\bigl(\omega(i)\bigr)_{i \in \underline{3n}}\Bigr)\left[\Gamma\right]
\in \Aavis_n(\source)\end{equation*} onto $\Aavis^{\prime}_n(\source)$. This sum is a sum over diagrams with at least two univalent vertices, according to Lemma~\ref{lemoneleg}. Recall Notation~\ref{notlessav}.  Lemma~\ref{lemconststraightred} guarantees that
\begin{equation*}\Zinv_{n,I}\Bigl(\crats,\Link,\bigl(\omega(i)\bigr)_{i \in I}\Bigr)=\sum_{\Gamma \in \Davisred^e_{n,I}(\source)}\coefgambetred_{\Gamma}I\Bigl(\rats,\Link,\Gamma,\bigl(\omega(i)\bigr)_{i \in I}\Bigr)\left[\Gamma\right]
\in \Aavis_n(\source)\end{equation*} is independent of the chosen
$\omega(i)$ for any subset $I$ of $\underline{3n}$ with cardinality $(3n-2)$. 
Observe \begin{equation*}\Zinv^{\prime}_{n}\Bigl(\crats,\Link,\bigl(\omega(i)\bigr)_{i \in \underline{3n}}\Bigr)=\frac{2\bigl((3n-2)!\bigr)}{(3n)!}\sum_{I \subset \underline{3n} \suchthat \cardlef{I} =3n-2}\Zinv_{n,I}\Bigl(\crats,\Link,\bigl(\omega(i)\bigr)_{i \in I}\Bigr)\end{equation*}
since for a numbered graph $\Gamma$ of $\Davis^e_{n}(\source)$, the coefficient of $I(\rats,\Link,\Gamma,(\omega(i))_{i \in \underline{3n}})$
 is $\frac{(3n-\cardlef{E(\Gamma)})!}{(3n)!2^{\cardlef{E(\Gamma)}}}\left[\Gamma\right]$ in the left-hand side, and 
\begin{multline*}\frac{2\bigl((3n-2)!\bigr)}{(3n)!}\sum_{I \subset \underline{3n} \suchthat \cardlef{I} 
=3n-2, j_{E}(E(\Gamma)) \subseteq I}\frac{(3n-2-\cardlef{E(\Gamma)})!}{(3n-2)!2^{\cardlef{E(\Gamma)}}}\left[\Gamma\right]\\=\frac{2\bigl((3n-2)!\bigr)}{(3n)!}  \frac{(3n-\cardlef{E(\Gamma)})!}{2((3n-2-\cardlef{E(\Gamma)})!)}   \frac{(3n-2-\cardlef{E(\Gamma)})!}{(3n-2)!2^{\cardlef{E(\Gamma)}}}\left[\Gamma\right]\end{multline*} in the right-hand-side.
Thus, Lemma~\ref{lemconststraightred} and Lemma~\ref{leminvfstconst} imply Theorem~\ref{thmconststraight}.

Lemma~\ref{lemconststraightred} also directly implies that
\begin{equation*}\sum_{\Gamma \in \Davisred^e_{n,3n-2}(\source)}\coefgambetred_{\Gamma}I\Bigl(\rats,\Link,\Gamma,\bigl(\omega(i)\bigr)_{i \in \underline{3n-2}}\Bigr)\projassis\bigl([\Gamma]\bigr)\end{equation*} is independent of the chosen
$\omega(i)$, so it implies Theorem~\ref{thmconststraightred} as above.
\eop

\bpo{Proof of Lemma~\ref{lemconststraightred}}
Let $\Link \colon \source \hookrightarrow \crats$ be a straight embedding  with respect to $\partau$. It suffices to prove that
$\Zinv_{n,3n-2}(\crats,\Link,(\omega(i))_{i \in \underline{3n-2}})$
does not change when some $\omega(i)$ is changed to $\omega(i)+d\eta$ for some one-form $\eta$ on $C_2(\rats)$, which restricts to $\partial C_2(\rats)$ as $\projp_{\partau}^{\ast}(\eta_{S^2})$ for some one-form $\eta_{S^2}$ on $S^2$, as in Lemma~\ref{lemetactwo}.
Assume that the forms ${\omega}(j)$ restrict to $\partial C_2(\rats)$ as $\projp_{\partau}^{\ast}({\omega}_{S^2}(j))$.
Set $\tilde{\omega}(i,0)=\omega(i)$.
Let $p_{C_2} \colon \left[0,1\right] \times C_2(\rats) \to C_2(\rats)$ and $p_{S^2} \colon \left[0,1\right] \times S^2 \to S^2$ denote the projections onto the second factor.
Define the closed $2$-form $\tilde{\omega}_{S^2}(i)$ on $\left[0,1\right] \times S^2$ by 
\begin{equation*}\tilde{\omega}_{S^2}(i) = p_{S^2}^{\ast}\bigl({\omega}_{S^2}(i)\bigr) + d \bigl(t p_{S^2}^{\ast}(\eta_{S^2})\bigr),\end{equation*}
where $t$ is the coordinate on $\left[0,1\right]$.
Define the closed $2$-form $\tilde{\omega}(i)$ on $\left[0,1\right] \times C_2(\rats)$ to be
\begin{equation*}\tilde{\omega}(i) = p_{C_2}^{\ast}\bigl(\omega(i)\bigr) + d \bigl(t p_{C_2}^{\ast}(\eta)\bigr).\end{equation*}
For $j \in \underline{3n}\setminus \{i\}$, define $\tilde{\omega}_{S^2}(j)=p_{S^2}^{\ast}({\omega}_{S^2}(j))$
and $\tilde{\omega}(j) = p_{C_2}^{\ast}(\omega(j))$. For any $j  \in \underline{3n},$ let
$\tilde{\omega}(j,t)$ denote the restriction of $\tilde{\omega}(j)$ to $\{t\} \times C_2(\rats)$.
Thus, it suffices to prove that $\Zinv_{n,3n-2}(1)=\Zinv_{n,3n-2}(0)$, with \begin{equation*}\Zinv_{n,3n-2}(t)=\sum_{\Gamma \in \Davis^e_{n,3n-2}(\source)}\coefgambet_{\Gamma}I\Bigl(\rats,\Link,\Gamma,\bigl(\tilde{\omega}(i,t)\bigr)_{i \in \underline{3n-2}}\Bigr)\left[\Gamma\right] \in \Aavis_n(\source).\end{equation*}

Proposition~\ref{propinvone} expresses $\bigl(\Zinv_{n,3n-2}(1)-\Zinv_{n,3n-2}(0)\bigr)$ as a sum over numbered graphs $\Gamma$ equipped with a connected component $\Gamma_{\finseta}$, with no univalent vertex, or whose univalent vertices form a nonempty set of consecutive vertices in $\Gamma$ on some component $S^1_j$ of $\source$.
The faces for which $\finseta$ has no univalent vertex do not contribute, as in the proof of Lemma~\ref{leminvwoboun}. So we focus on the remaining faces.

Such a $\facee(\finseta,\Link,\Gamma)$ may split according to the possible compatible linear orders of the univalent vertices of $\Gamma_{\finseta}$,  represented by lifts $\check{\Gamma}_{\finseta}$ of  $\Gamma_{\finseta}$ on $\RR$ as in Notation~\ref{notorderlift}. We view
$\Zinv_{n,3n-2}(1)-\Zinv_{n,3n-2}(0)$ as a sum over pairs $(\Gamma,\check{\Gamma}_{\finseta})$ of terms
\begin{equation*} \coefgambet_{\Gamma}I(\Gamma,\check{\Gamma}_{\finseta} )= \coefgambet_{\Gamma}\int_{\left[0,1\right]\times\facee(\check{\Gamma}_{\finseta},\Link,\Gamma)}\bigwedge_{e \in E(\Gamma)}p_e^{\ast}\Bigl(\tilde{\omega}\bigl(j_E(e)\bigr)\Bigr)\left[\Gamma\right]\end{equation*}
 associated to the corresponding face components denoted by $\facee(\check{\Gamma}_{\finseta},\Link,\Gamma)$.

For a fixed numbered graph $\check{\Gamma}_{\finseta}$ on $\RR$ as above, we study the sum of the contributions $\coefgambet_{\Gamma} I(\tilde{\Gamma},\check{\Gamma}_{\finseta})$ 
running over the graphs $\tilde{\Gamma}$ such that
\begin{itemize}
 \item  the graph $\check{\Gamma}_{\finseta}$ is a subgraph of $\tilde{\Gamma}$ when $\RR$ is identified with a part of $S^1_j$,
 \item the univalent vertices of $\check{\Gamma}_{\finseta}$ are consecutive on $S^1_j$ in $\tilde{\Gamma}$, with respect to their linear order, and
 \item $\tilde{\Gamma} \setminus \check{\Gamma}_{\finseta}$ is equal to a fixed ${\Gamma} \setminus \check{\Gamma}_{\finseta}$, as above
\end{itemize}
Recall from Proposition~\ref{propdiagcom} that $[\tilde{\Gamma}]=\left[\check{\Gamma}_{\finseta}\right]\#_j\left[{\Gamma} \setminus \check{\Gamma}_{\finseta}\right]$ in $\Assis(\source)$ for any such pair $(\tilde{\Gamma},\check{\Gamma}_{\finseta})$, provided that all the vertices of the graph $\tilde{\Gamma}$ inherit their orientations from a fixed vertex-orientation of ${\Gamma}$.

If $i \notin j_E(E(\check{\Gamma}_{\finseta}))$, then the  $ I(\tilde{\Gamma},\check{\Gamma}_{\finseta} )$ vanish 
because all the $p_e^{\ast}(\tilde{\omega}(j_E(e)))$ for $e \in E(\check{\Gamma}_{\finseta})$ factor through the projection onto $\facee(\check{\Gamma}_{\finseta},\Link,\Gamma)$ whose dimension is
$\bigl(2\card{E(\check{\Gamma}_{\finseta})} -1\bigr)$.

If $i \in j_E(E(\check{\Gamma}_{\finseta}))$, let $e(i)$ be the edge such that $j_E(e(i))=i$. The sum of the contributions $\coefgambet_{.} I(.,\check{\Gamma}_{\finseta} )$ involving $\check{\Gamma}_{\finseta}$ factors through
\begin{equation*}I=\int_{\left[0,1\right] \times \cup_{c(\eltv)\in K_j}\cinjuptdanvec(T_{c(\eltv)}\crats,\vec{t}_{c(\eltv)};\check{\Gamma}_{\finseta})}p_{e(i)}^{\ast}\Bigl(d\bigl(t p_{C_2}^{\ast}(\eta) \bigr)\Bigr) \bigwedge_{e \in E(\check{\Gamma}_{\finseta}) \setminus e(i)}p_e^{\ast}\Bigl(\tilde{\omega}\bigl(j_E(e)\bigr)\Bigr),\end{equation*}
where $\vec{t}_{c(\eltv)}$ denotes the unit tangent vector to $K_j$ at $c(\eltv)$.

Recall that $\check{\anomq}(\check{\Gamma}_{\finseta})$ was defined in Section~\ref{secanomalpha}, together with natural maps 
\begin{equation*}p_{e,S^2} \colon \check{\anomq}(\check{\Gamma}_{\finseta}) \to S^2.\end{equation*} Let $p_{e,S^2}$ also denote $\id_{\left[0,1\right]} \times p_{e,S^2} \colon \left[0,1\right] \times \check{\anomq}(\check{\Gamma}_{\finseta}) \to \left[0,1\right] \times S^2$. The form
\begin{equation*}p_{e(i)}^{\ast}\Bigl(d\bigl(t p_{C_2}^{\ast}(\eta) \bigr)\Bigr) \bigwedge_{e \in E(\check{\Gamma}_{\finseta}) \setminus e(i)}p_{e}^{\ast}\Bigl(\tilde{\omega}\bigl(j_E(e)\bigr)\Bigr)\end{equation*}
is the pull-back of the closed form \begin{equation*}
\Omega= p_{e(i),S^2}^{\ast}\bigl(d(t \eta_{S^2} )\bigr) \bigwedge_{e \in E(\check{\Gamma}_{\finseta}) \setminus e(i)}p_{e,S^2}^{\ast}\Bigl(\tilde{\omega}_{S^2}\bigl(j_E(e)\bigr)\Bigr)\end{equation*}
on $\left[0,1\right] \times \check{\anomq}(\check{\Gamma}_{\finseta})$ under the projection
\begin{equation*}\left[0,1\right] \times \cup_{c(\eltv)\in K_j}\cinjuptdanvec\bigl(T_{c(\eltv)}\crats,\vec{t}_{c(\eltv)};\check{\Gamma}_{\finseta}\bigr) \to \left[0,1\right] \times \check{\anomq}(\check{\Gamma}_{\finseta}).\end{equation*}

The image of this projection is the product by $\left[0,1\right]$ of
the restriction of the bundle $\check{\anomq}(\check{\Gamma}_{\finseta})$ over $\projp_{\partau}(\ST^+K_j)$, and $I$ is the integral of $\Omega$ along this image. Compute the integral by integrating first along the fibers of $\check{\anomq}(\check{\Gamma}_{\finseta})$, next along $\left[0,1\right]$. Afterwards, the integral $I$ becomes the integral of a one-form along $\projp_{\partau}(\ST^+K_j) \subset S^2$. So it vanishes because $K_j$ is straight.\eop

\section{The general variation for homogeneous propagating forms}
\label{secdephom}
Set $\CD^c(\RR)=\cup_{n \in \NN} \CD^c_n(\RR)$, where $\CD^c_n(\RR)$ is the set of degree $n$ connected $(3n-2)$-numbered Jacobi diagrams on $\RR$ introduced in the beginning of Section~\ref{secanomalpha}.
In this section, we write various sums over numbered diagrams, but all the edges of a diagram are equipped with the same propagating forms. So neither the set in which the edges are numbered nor its cardinality matters, provided that the cardinality is greater than the possible number of edges for a given degree. (See Proposition~\ref{propdefhomog} and Remark~\ref{rknoconflict}.)

\begin{proposition}
\label{propinvtwomieux}
Let $(\crats,\partau)$ be an asymptotic rational homology $\RR^3$. Let $\Link=\sqcup_{j=1}^kK_j$ be an embedding of $\source=\sqcup_{j=1}^kS^1_j$ into $\crats$.
Let $\tilde{\omega}(0)$ and $\tilde{\omega}(1)$ be two homogeneous propagating forms of $C_2(\rats)$. Let $\tilde{\omega}$ be a closed $2$-form on $\left[0,1\right] \times \partial C_2(\rats)$ whose restriction $\tilde{\omega}(t)$ to $\{t\} \times \left(\partial C_2(\rats) \setminus \ST  {\ballb}_{\rats}\right)$ is $\projp_{\partau}^{\ast}(\omega_{S^2})$ for any $t \in \left[0,1\right]$, and whose restriction $\tilde{\omega}(t)$ to $\{t\} \times \partial C_2(\rats) $ coincides with the restriction to $\partial C_2(\rats) $ of the given $\tilde{\omega}(t)$ for $t \in \{0,1\}$.
For any $j \in \underline{k}$, define
$I_j=\sum_{\Gamma_{\finsetb} \in \CD^c(\RR)}\coefgambet_{\Gamma_{\finsetb}}I(\Gamma_{\finsetb},K_j,\tilde{\omega}),$
where \begin{equation*}I(\Gamma_{\finsetb},K_j,\tilde{\omega})=\int_{u\in \left[0,1\right]}\int_{w\in K_j}\int_{\cinjuptdanvec(T_{w}\crats,\vec{t}_w;{\Gamma}_{\finsetb})} \bigwedge_{e \in E(\Gamma_{\finsetb})}p_e^{\ast}\bigl(\tilde{\omega}(u)\bigr)\left[\Gamma_{\finsetb}\right],\end{equation*} 
and $\vec{t}_w$ denotes the unit tangent vector to $K_j$ at $w$.
Define \begin{equation*}\zinv(\tilde{\omega})=\sum_{n\in \NN}\zinv_n\bigl([0,1] \times \ST {\ballb}_{\rats}; \tilde{\omega}\bigr)\end{equation*}
as in Corollary~\ref{corinvone}.
Then we have
\begin{equation*}\Zinv\bigl(\crats,\Link,\tilde{\omega}(1)\bigr)=\left(\prod_{j=1}^k \exp\left(I_j\right) \#_j\right)\Zinv\bigl(\crats,\Link,\tilde{\omega}(0)\bigr) \exp\bigl(\zinv(\tilde{\omega})\bigr),\end{equation*}
where $\#_j$ stands for the insertion on a diagram on $\RR$ on the component $S^1_j$ of $\source$.
\end{proposition}

We will actually prove the following lemma.
\begin{lemma}
\label{lempropinvtwomieux} Recall Notation~\ref{notationzZ}.
Under the assumptions of Proposition~\ref{propinvtwomieux}, we have
\begin{equation*}\Zinvlink\left(\crats,\Link,\tilde{\omega}(1)\right)=\left(\prod_{j=1}^k \exp\left(I_j\right) \#_j\right)\Zinvlink\left(\crats,\Link,\tilde{\omega}(0)\right).\end{equation*}
\end{lemma}

\begin{lemma}
 Lemma~\ref{lempropinvtwomieux} implies Proposition~\ref{propinvtwomieux}.
\end{lemma}
\bp When $\Link =\emptyset$, Proposition~\ref{propinvtwomieux} follows from Corollary~\ref{corinvone}, Proposition~\ref{propexpon}, and Lemma~\ref{lemformprod}, which ensures that there exists a closed $2$-form $\tilde{\omega}$ on $\left[0,1\right] \times C_2(\rats)$ that extends the $2$-form $\tilde{\omega}$ of the statement. Conclude with Lemma~\ref{lemmultcheck}.
\eop

Let us begin the proof of Lemma~\ref{lempropinvtwomieux} with the proof of the following corollary of Proposition~\ref{propinvone}.
\begin{lemma}
\label{leminvtwo}
Under the assumptions of Proposition~\ref{propinvtwomieux},
let $\tilde{\omega}$ be a closed $2$-form on $\left[0,1\right] \times C_2(\rats)$ which extends the $2$-form $\tilde{\omega}$ of Proposition~\ref{propinvtwomieux}. For any $t \in \left[0,1\right]$, let $\tilde{\omega}(t)$ denote the restriction to $\{t\} \times C_2(\rats)$ of $\tilde{\omega}$.
Set
\begin{equation*}\Zinvlink(t)=\Bigl(\Zinvlink_n\bigl(\rats,\Link,\tilde{\omega}(t)\bigr)\Bigr)_{n\in \NN}.\end{equation*}
For $\Gamma_{\finsetb} \in \CD^c(\RR)$, $u \in\left[0,1\right]$, and $j \in \underline{k}$,
set \begin{equation*}\eta(\rats,\Link,\Gamma_{\finsetb},K_j,\tilde{\omega})(u)=\int_{w\in K_j}\int_{\cinjuptdanvec(T_{w}\crats,\vec{t}_w;{\Gamma}_{\finsetb})} \bigwedge_{e \in E(\Gamma_{\finsetb})}p_e^{\ast}\bigl(\tilde{\omega}(u)\bigr)\left[\Gamma_{\finsetb}\right],\end{equation*} 
where $\vec{t}_w$ denotes the unit tangent vector to $K_j$ at $w$,
and set
\begin{equation*}\gamma_j(u)=\sum_{\Gamma_{\finsetb} \in \CD^c(\RR)}\coefgambet_{\Gamma_{\finsetb}}\eta(\rats,\Link,\Gamma_{\finsetb},K_j,\tilde{\omega})(u).\end{equation*}
Then $\Zinvlink(t)$ is differentiable, and we have
\begin{equation*}\Zinvlink^{\prime}(t)dt = \left(\sum_{j=1}^k \gamma_j(t) \#_j\right) \Zinvlink(t).\end{equation*}
\end{lemma}
\bpo{Proof of Lemma~\ref{leminvtwo}} 
The variations of $\Zinvlink_n(t)$ are given by Proposition~\ref{propinvone}, by sending the diagrams with components without univalent vertices to $0$.
They involve only faces $\facee(\finseta,\Link,\Gamma)$ for which $\Gamma_{\finseta}$ is a connected component of $\Gamma$, with univalent vertices on one component of $\source$. Again, such a face may split according to the possible compatible linear orders of the univalent vertices of $\Gamma_{\finseta}$, represented by lifts $\check{\Gamma}_{\finseta}$ of $\Gamma_{\finseta}$ on $\RR$, as in Notation~\ref{notorderlift}. The corresponding face component is denoted by $\facee(\check{\Gamma}_{\finseta},\Link,\Gamma)$, and the corresponding integral is
\begin{equation*}I(\Gamma,\check{\Gamma}_{\finseta})=\int_{\left[0,t\right]\times \facee(\check{\Gamma}_{\finseta},\Link,\Gamma)}\bigwedge_{e \in E(\Gamma)}p_e^{\ast}(\tilde{\omega})\left[\Gamma\right].\end{equation*}
Proposition~\ref{propinvone} implies
\begin{equation*}\Zinvlink_n(t)-\Zinvlink_n(0)=\sum_{\substack{\scriptstyle(\Gamma,{\finseta}) \suchthat \Gamma \in {\Davis}^e_n(\source), {\finseta} \subseteq {\finsetv(\Gamma)}, \cardlef{\finseta} \geq 2,
\\ \mbox{\scriptsize every component of $ \Gamma $ has univalent vertices,}\\
\mbox{\scriptsize $\Gamma_{\finseta}$ is a connected component of $\Gamma$,}\\
\mbox{\scriptsize the univalent vertices of $\Gamma_{\finseta}$ are consecutive}\\
\mbox{\scriptsize on one component of $\Gamma$, and}\\
\mbox{\scriptsize $\check{\Gamma}_{\finseta}$ is a compatible lift of  $\Gamma_{\finseta}$  on $\RR$.}
}}\coefgambet_{\Gamma}I\bigl(\Gamma,\check{\Gamma}_{\finseta}\bigr),\end{equation*}
where the set of univalent vertices of $\check{\Gamma}_{\finseta}$ is equipped with the unique linear order induced by $\Gamma$ if there are univalent vertices of $\Gamma \setminus \Gamma_{\finseta}$ on the component of $\Gamma_{\finseta}$, and with one of the linear orders compatible with $\Gamma$ otherwise.

This expression implies that $\Zinvlink_n$ (valued in a finite-dimensional vector space) is differentiable. (For any smooth compact $d$-dimensional manifold $C$ and for any smooth $(d+1)$-form $\omega$ on $\left[0,1\right]\times C$, the function $\bigl(t \mapsto \int_{\left[0,t\right]\times C}\omega \bigr)$ is differentiable.)
Assume that the vertices of $\Gamma_{\finseta}$ are on a component $K_j$ of $\Link(\sqcup_{j=1}^kS^1_j)$.
The forms associated to edges of $\Gamma_{\finseta}$ are integrated along $\left[0,1\right] \times (\cup_{c(\eltv)\in K_j}\cinjuptdanvec(T_{c(\eltv)}\crats,\vec{t}_{c(\eltv)};\check{\Gamma}_{\finseta}))$, where $\vec{t}_{c(\eltv)}$ denotes the unit tangent vector to $K_j$ at $c(\eltv)$. So they do not depend on the configuration of $\bigl({\finsetv(\Gamma)} \setminus {\finseta}\bigr)$.
The other forms are integrated along $\check{C}(\rats,\Link;\Gamma \setminus \Gamma_{\finseta})$ at $u \in \left[0,1\right]$. 

Group the contributions of the pairs $(\Gamma,\check{\Gamma}_{\finseta})$ with common
$(\Gamma \setminus \Gamma_{\finseta},\check{\Gamma}_{\finseta})$ to view the
 global variation $\bigl(\Zinvlink(t)-\Zinvlink(0)\bigr)$ as
\begin{equation*}\sum_{j=1}^k\int_0^t \Bigl(\sum_{\check{\Gamma}_{\finseta} \in \CD^c(\RR)}\coefgambet_{\check{\Gamma}_{\finseta}}\eta(\rats,\Link,\check{\Gamma}_{\finseta},K_j,\tilde{\omega})(u)\#_j\Bigr) \Zinvlink(u) .\end{equation*}
Use Proposition~\ref{propdefhomog} and Remark~\ref{rkhomog} to check that the coefficients are correct.
So we get
\begin{equation*}\Zinvlink(t)-\Zinvlink(0)=\int_0^t\Bigl(\sum_{j=1}^k \gamma_j(u)\#_j\Bigr) \Zinvlink(u).\end{equation*}
\eop

\bpo{Proof of Lemma~\ref{lempropinvtwomieux}}
Set
$I_j(t)=\int_{0}^t\gamma_j(u)$. We have $\Zinvlink_0(t)=1$.
The equation \begin{equation*}\Zinvlink^{\prime}(t)dt = \bigl(\sum_{j=1}^k \gamma_j(t) \#_j\bigr) \Zinvlink(t)\end{equation*} of Lemma~\ref{leminvtwo} determines $\Zinvlink(t)$ as a function of $\Zinvlink(0)$ by induction on the degree, and we have
$\Zinvlink(t)=\prod_{j=1}^k\exp\bigl(I_j(t)\bigr) \#_j \Zinvlink(0)$.
\eop

Let us now apply Lemma~\ref{leminvtwo} to study the variation of the quantity $\zinvlink(\rats,\Link,\partau)$ of Corollary~\ref{corlogZinv} when $\partau$ varies smoothly.

\begin{lemma}
\label{lemvartauanom}
Let $\bigl(\partau(t)\bigr)_{t \in \left[0,1\right]}$ define a smooth homotopy of asymptotically standard parallelizations of $\crats$. We have
\begin{equation*}\frac{\partial}{\partial t}\Zinvlink\bigl(\crats,\Link,\partau(t)\bigr)=\left(\sum_{j=1}^k\frac{\partial}{\partial t}\left( 2 \int_{\left[0,t\right]\times \ST^+K_j}p_{\partau(.)}^{\ast}(\omega_{S^2})\right)\alpha \#_j\right)\Zinvlink\bigl(\crats,\Link,\partau(t)\bigr).\end{equation*} 
\end{lemma}
\bp
Fix a homogeneous propagating form $\omega$ of $(C_2(\rats),\partau(0))$ and
a form $\tilde{\omega}$ on $\left[0,1\right]\times C_2(\rats)$ such that $\tilde{\omega}(t)$ is a homogeneous propagating form $\omega$ of $(C_2(\rats),\partau(t))$ for all $t \in \left[0,1\right]$ as in 
Lemma~\ref{lemformprod}.
Lemma~\ref{leminvtwo} ensures 
\begin{equation*}\frac{\partial}{\partial t}\Zinvlink\bigl(\crats,\Link,\partau(t)\bigr)dt=\left(\sum_{j=1}^k \gamma_j(t) \#_j\right) \Zinvlink\bigl(\crats,\Link,\partau(t)\bigr),\end{equation*}
with 
\begin{equation*}\gamma_j(u)=\sum_{\Gamma_{\finsetb} \in \CD^c(\RR)}\coefgambet_{\Gamma_{\finsetb}}\eta(\rats,\Link,\Gamma_{\finsetb},K_j,\tilde{\omega})(u),\end{equation*}
and
\begin{equation*}\eta(\rats,\Link,\Gamma_{\finsetb},K_j,\tilde{\omega})(u)=\int_{w\in K_j}\int_{\cinjuptdanvec(T_{w}\crats,\vec{t}_w;{\Gamma}_{\finsetb})} \bigwedge_{e \in E(\Gamma_{\finsetb})}p_e^{\ast}\bigl(p_{\partau(u)}^{\ast}(\omega_{S^2})\bigr)\left[\Gamma_{\finsetb}\right].\end{equation*} 
The restriction of $p_{\partau(.)}$ from $\left[0,1\right]\times \ST^+K_j$ to $S^2$ induces a map 
\begin{equation*}p_{a,\partau,{\Gamma}_{\finsetb}}\colon \left[0,1\right]  \times \cup_{w\in K_j}\cinjuptdanvec(T_{w}\crats,\vec{t}_{w};{\Gamma}_{\finsetb}) \rightarrow \check{\anomq}({\Gamma}_{\finsetb})\end{equation*} 
over $\left(p_{\partau(.)} \colon \left[0,1\right]\times \ST^+K_j \to S^2\right)$, which restricts to the fibers as the identity map, for any ${\Gamma}_{\finsetb} \in \CD^c(\RR)$. (Recall the definition of $\check{\anomq}({\Gamma}_{\finsetb})$ from the beginning of Section~\ref{secanomalpha}.) We have
\begin{equation*}\int_0^1\eta(\rats,\Link,\Gamma_{\finsetb},K_j,\tilde{\omega})(u)=\int_{\Image(p_{a,\partau,{\Gamma}_{\finsetb}})}\left(\bigwedge_{e \in E({\Gamma}_{\finsetb})}p_{e,S^2}^{\ast}(\omega_{S^2})\right)\left[\Gamma_{\finsetb}\right].\end{equation*}
Integrating $\bigl(\bigwedge_{e \in E({\Gamma}_{\finsetb})}p_{e,S^2}^{\ast}(\omega_{S^2})\bigr)\left[{\Gamma}_{\finsetb}\right]$ along the fiber in $\check{\anomq}({\Gamma}_{\finsetb})$ yields a two-form
on $S^2$. This two-form is homogeneous because everything is. Thus, it may be expressed as $2\alpha({\Gamma}_{\finsetb})\omega_{S^2}\left[{\Gamma}_{\finsetb}\right]$, where $\alpha({\Gamma}_{\finsetb}) \in \RR$ and \begin{equation*}\sum_{{\Gamma}_{\finsetb} \in \CD^c(\RR)}\coefgambet_{{\Gamma}_{\finsetb}}\alpha(\Gamma_{\finsetb})\left[{\Gamma}_{\finsetb}\right]=\alpha.\end{equation*}
So we get
\begin{equation*}\int_0^t\eta(\rats,\Link,\Gamma_{\finsetb},K_j,\tilde{\omega})(u)=2\alpha({\Gamma}_{\finsetb})\int_{\left[0,t\right]\times \ST^+K_j}p_{\partau(.)}^{\ast}(\omega_{S^2})\left[{\Gamma}_{\finsetb}\right],\end{equation*}
and $\gamma_j(t)=2\frac{\partial}{\partial t}\left(\int_{\left[0,t\right]\times \ST^+K_j}p_{\partau(.)}^{\ast}(\omega_{S^2})\right)\alpha dt$.
\eop

\begin{corollary}
\label{corvartauanom} The product
\begin{equation*}\prod_{j=1}^k\exp\Bigl(-I_{\theta}\bigl(K_j,\partau(t)\bigr)\alpha\Bigr)\#_j \Zinvlink\bigl(\crats,\Link,\partau(t)\bigr)\end{equation*} does not change when $\partau$ varies by a smooth homotopy.
\end{corollary}
\bp With the notation of Lemma~\ref{lemdefItheta}, we have \begin{equation*}\Zinvlink_1\bigl(\crats,\Link,\partau(t)\bigr)=\frac12\sum_{j=1}^k I_{\theta}\bigl(K_j,\partau(t)\bigr) \left[ \onechordsmallh{$S^1_j$}\right].\end{equation*}
Therefore, Lemma~\ref{lemvartauanom} and Proposition~\ref{propanom} imply 
\begin{equation*}\frac{\partial}{\partial t}I_{\theta}\bigl(K_j,\partau(t)\bigr)=2\frac{\partial}{\partial t}\int_{\left[0,t\right]\times \ST^+K_j}p_{\partau(.)}^{\ast}(\omega_{S^2}).\end{equation*}                                              (The reader can also check it directly as an exercise.)
So Lemma~\ref{lemvartauanom} implies
\begin{equation*}\frac{\partial}{\partial t}\Zinvlink\bigl(\crats,\Link,\partau(t)\bigr)=\sum_{j=1}^k\Bigl(\frac{\partial}{\partial t}I_{\theta}\bigl(K_j,\partau(t)\bigr)\alpha \#_j\Bigr)\Zinvlink\bigl(\crats,\Link,\partau(t)\bigr).\end{equation*}
Therefore, the derivative of \begin{equation*}\prod_{j=1}^k\exp\Bigl(-I_{\theta}\bigl(K_j,\partau(t)\bigr)\alpha\Bigr)\#_j \Zinvlink\bigl(\crats,\Link,\partau(t)\bigr)\end{equation*} with respect to $t$ vanishes.
\eop

\bpo{Proof of Theorem~\ref{thmfstconst}}
According to the naturality of Proposition~\ref{propinvdiffnat}, Lemma~\ref{leminvfstconst}, Proposition~\ref{propdeftwoanom}, Corollary~\ref{cordefzinvuf}, and Proposition~\ref{propanom}, it
suffices to prove that \begin{equation*}\prod_{j=1}^k\Bigl(\exp\bigl(-I_{\theta}(K_j,\partau)\alpha\bigr)\#_j\Bigr) \Zinvlink(\crats,\Link,\partau) \in \Assis(\source)\end{equation*}
is independent of the homotopy class of parallelization $\partau$.

When $\partau$ changes in a ball that does not meet the link, the forms can be changed only in the neighborhoods of the unit tangent bundle to this ball. Apply 
Proposition~\ref{propinvtwomieux}
again to $\Zinvlink$, where the $p_e^{\ast}(\tilde{\omega}(u))$ are independent of $u$ over $K_j$.
So we get that
\begin{equation*}\prod_{j=1}^k\Bigl(\exp\bigl(-I_{\theta}(K_j,\partau)\alpha\bigr)\#_j\Bigr) \Zinvlink(\crats,\Link,\partau) \in \Assis(\source)\end{equation*}
is invariant under the natural action of $\pi_3(SO(3))$ on the homotopy classes of parallelizations, thanks to Corollary~\ref{corvartauanom}.

Now examine the effect of the twist of the parallelization by a map $g \colon ({\ballb}_{\rats},1) \rightarrow (SO(3),1)$. Without loss of generality, assume that $\projp_{\partau}(\ST^+K_j)=v$ for some $v \in S^2$ and $g$ maps $K_j$ to rotations with axis $v$ for any $j \in \underline{k}$.
We want to compare $\Zinvlink(\crats,\Link,\partau \circ \psi_{\RR}(g))$ with $\Zinvlink(\crats,\Link,\partau)$.
There exists a closed form $\omega$ on
$\left[0,1\right]\times \ST {\ballb}_{\rats}$ equal to $\projp_{\partau}^{\ast}(\omega_{S^2})$ on $\partial (\left[0,1\right]\times \ST {\ballb}_{\rats}) \setminus (\{1\} \times \ST {\ballb}_{\rats})$ and equal to $p_{\partau \circ \psi_{\RR}(g)}^{\ast}(\omega_{S^2})$ on $ \{1\} \times \ST {\ballb}_{\rats}$.
Extend this form to a closed form $\Omega$ on $\left[0,1\right]\times C_2(\rats)$, which restricts to $\left[0,1\right] \times (\partial C_2(\rats) \setminus \ST {\ballb}_{\rats})$ as $\projp_{\partau}^{\ast}(\omega_{S^2})$, 
as in 
Lemma~\ref{lemformprod}.
Let $\Omega(t)$  denote the restriction of $\Omega$ to  $\{t\} \times C_2(\rats)$.
According to Proposition~\ref{propinvtwomieux}, we have
\begin{equation*}\Zinvlink\bigl(\crats,\Link,\partau \circ \psi_{\RR}(g)\bigr)=\prod_{j=1}^k\bigl(\exp(I_j)\#_j \bigr)\Zinvlink(\crats,\Link,\partau),\end{equation*}
where $I_j=\int_{0}^1\gamma_j(u)$, with
$\gamma_j(t)=\sum_{\Gamma_{\finsetb} \in \CD^c(\RR)}\coefgambet_{{\Gamma}_{\finsetb}}\eta(\rats,\Link,\Gamma_{\finsetb},K_j,\Omega)(t)$
and \begin{equation*}\eta(\rats,\Link,\Gamma_{\finsetb},K_j,\Omega)(t)=\int_{w\in K_j}\int_{\cinjuptdanvec(T_{w}\crats,\vec{t}_w;{\Gamma}_{\finsetb})} \left(\bigwedge_{e \in E(\Gamma_{\finsetb})}p_e^{\ast}\bigl(\Omega(t)\bigr)\right)\left[\Gamma_{\finsetb}\right].\end{equation*} It suffices to prove $I_j=\bigl(I_{\theta}(K_j,\partau \circ \psi_{\RR}(g))-I_{\theta}(K_j,\partau)\bigr)\alpha$.
Proposition~\ref{propinvtwomieux} implies that the degree one part $I_{1,j}$ of $I_j$ is
\begin{equation*}\begin{array}{ll}I_{1,j}&=\Zinvlink_1\bigl(\crats,K_j,\partau \circ \psi_{\RR}(g)\bigr) -\Zinvlink_1\bigl(\crats,K_j,\partau\bigr)\\&=\frac12\left(I_{\theta}\bigl(K_j,\partau \circ \psi_{\RR}(g)\bigr)-I_{\theta}\bigl(K_j,\partau\bigr)\right)\left[ \onechordsmalltseul\right].\end{array}\end{equation*} 

Let $^{\partau}\psi(g^{-1}) \colon \ST {\ballb}_{\rats} \to \ST {\ballb}_{\rats}$ denote the map induced by $\partau  \circ \psi_{\RR}(g^{-1}) \circ \partau^{-1}$. Recall $\projp_{\partau}=p_{S^2} \circ \partau^{-1}$. So we have \begin{equation*}\projp_{\partau \circ \psi_{\RR}(g)}=p_{S^2} \circ \psi_{\RR}(g^{-1}) \circ \partau^{-1} =\projp_{\partau} \circ \leftidx{^{\partau}}\psi(g^{-1}).\end{equation*}
Let $(.-1) \colon \left[1,2\right] \to \left[0,1\right]$ map $x$ to $x-1$. Set $^{\;\;\,\partau}_{-1}\psi(g^{-1})=\left((.-1)  \times \leftidx{^{\partau}}\psi(g^{-1})\right)$.
Extend $\Omega$ over $\left[0,2\right]\times C_2(\rats)$ so that $\Omega$ restricts to $\left[1,2\right]\times \ST {\ballb}_{\rats}$ as $^{\;\;\,\partau}_{-1}\psi(g^{-1})^{\ast}(\Omega)$.
For any ${\Gamma}_{\finsetb}$, the map $^{\;\;\,\partau}_{-1}\psi(g^{-1})$ induces an orientation-preserving diffeomorphism
\begin{equation*}^{\;\;\,\partau}_{-1}\psi(g^{-1})_{\ast}\colon \left[1,2\right] \times \cup_{w\in K_j}{\cinjuptdanvec(T_{w}\crats,\vec{t}_w;{\Gamma}_{\finsetb})} \to  \left[0,1\right] \times \cup_{w\in K_j}{\cinjuptdanvec(T_{w}\crats,\vec{t}_w;{\Gamma}_{\finsetb})}\end{equation*}
such that $p_e \circ ^{\;\;\,\partau}_{-1}\psi(g^{-1})_{\ast} = ^{\;\;\,\partau}_{-1}\psi(g^{-1}) \circ p_e$ for any edge $e$ of ${\Gamma}_{\finsetb}$.
Using these diffeomorphisms $^{\;\;\,\partau}_{-1}\psi(g^{-1})_{\ast}$ to pull back $\bigl(\bigwedge_{e \in E(\Gamma_{\finsetb})}p_e^{\ast}(\Omega(t))\bigr)$
proves $\gamma_j(t+1)=\gamma_j(t)$.
In particular, we get \begin{equation*}I_j(2)=\int_0^2\gamma_j(u)= 2I_j.\end{equation*}
Set $\Zinvlink(2)=\Zinvlink(\rats,\Link,\partau \circ \psi_{\RR}(g)^2)$.
We have
\begin{equation*}\Zinvlink(2)=\prod_{j=1}^k\exp\biggl(\Bigl(I_{\theta}(K_j,\partau \circ \psi_{\RR}(g)^2)-I_{\theta}\bigl(K_j,\partau\bigr)\Bigr)\alpha\biggr)\#_j \Zinvlink(\crats,\Link,\partau),\end{equation*} since $g^2$ is homotopic to the trivial map outside a ball (see Lemma~\ref{lempreptrivun},~2).
Proposition~\ref{propinvtwomieux} implies
\begin{equation*}\Zinvlink(2)=\prod_{j=1}^k\bigl(\exp(2I_j)\#_j \bigr)\Zinvlink(\crats,\Link,\partau).\end{equation*}
By induction on the degree, we get
\begin{equation*}2I_j=\bigl(I_{\theta}(K_j,\partau \circ \psi_{\RR}(g)^2)-I_{\theta}(K_j,\partau)\bigr)\alpha.\end{equation*}
The degree one part of this equality implies 
\begin{equation*}2\Bigl(I_{\theta}\bigl(K_j,\partau \circ \psi_{\RR}(g)\bigr)-I_{\theta}\bigl(K_j,\partau\bigr)\Bigr)=I_{\theta}\bigl(K_j,\partau \circ \psi_{\RR}(g)^2\bigr)-I_{\theta}\bigl(K_j,\partau\bigr).\end{equation*}
We get $I_j=\left(I_{\theta}(K_j,\partau \circ \psi_{\RR}(g))-I_{\theta}(K_j,\partau)\right)\alpha$, as desired.
\eop

\section{Some more properties of \texorpdfstring{$\Zinv$}{Z}}
\label{secfstpropzinv}

When $\crats=\RR^3$, then $\Zinvuf(S^3,\Link)=\Zinvlinkuf(S^3,\Link)$ is the configuration space invariant studied by Daniel Altsch\"uler, Laurent Freidel \cite{af}, Dylan Thurston \cite{thurstonconf}, Sylvain Poirier \cite{poirier}, and others, after the work of many people including Edward Witten \cite{witten}, Enore Guadagnini, Maurizio Martellini, Mihail Mintchev~\cite{gmm}, Maxim Kontsevich~\cite{ko,Kon}, Raoul Bott and Clifford Taubes \cite{botttaubes}, Dror Bar-Natan~\cite{barnatanper}\dots

\paragraph{Reversing a link component orientation.}

The following proposition is obvious from the definition of $\Zinv$.
\begin{proposition}
\label{proporcomp}
Let $\Link \colon \sqcup_{j=1}^kS^1_j \to \rats$ be a link in a $\QQ$-sphere $\rats$.
For a Jacobi diagram $\Gamma$ on $\sqcup_{j=1}^kS^1_j$, let $U_j(\Gamma)$ denote the set of univalent vertices of $\Gamma$ mapped to $S^1_j$. This set is cyclically ordered by $S^1_j$.
When the orientation of the component $\Link(S^1_j)$ is reversed, $\Zinvuf(\Link)$ is modified by reversing the circle $S^1_j$ (that is reversing the cyclic order of $U_j(\Gamma)$) in classes $\left[\Gamma\right]$ of diagrams $\Gamma$ on $\sqcup_{j=1}^kS^1_j$ and multiplying them by $(-1)^{\cardlef{U_j(\Gamma)}}$ in $\Aavis(\sqcup_{j=1}^kS^1_j)$.

In other words, we can forget the orientation of the link $\Link$ and view $\Zinvuf(\Link)$ as valued in $\Aavis(\sqcup_{j=1}^kS^1_j)$, where the $S^1_j$ are not oriented, as in Definitions~\ref{defdia} and \ref{defrkoruniv}.
\end{proposition}

\begin{remark}
\label{rkorunivtwo}
The orientation of a component $\Link(S^1_j)$ is used in two ways. It defines a cyclic order 
on $U_j(\Gamma)$, and it defines the orientation of the vertices of $U_j(\Gamma)$ as in Definition~\ref{defdia}.  The local orientation of $S^1_j$ near the image of a vertex orients the corresponding local factor of the configuration space. The cyclic order is encoded in the isotopy class of the injection of $U_j$ into the domain $S^1_j$. 
\end{remark}

\paragraph{Numbering of the link components.}
The following proposition is obvious from the definition of $\Zinvuf$.
\begin{proposition}
When the numbering of the components of $\Link$ is changed, $\Zinvuf(\Link)$ is modified by the corresponding change of numbering of the circles $S^1_j$ in diagram classes of $\Aavis(\sqcup_{j=1}^kS^1_j)$.
\end{proposition}

For a link $\Link \colon \source(\Link) \to \rats$ in a $\QQ$-sphere $\rats$,
$\Zinvuf(\rats,\Link)$ is valued in $\Aavis(\source(\Link))$. The one-manifold $\source(\Link)$ is a disjoint union of oriented circles, which have been numbered so far. However, the numbers may be changed to any decoration that marks the component.

This gives sense to the statement of the following theorem.

\paragraph{Connected sums.}

\begin{theorem}
\label{thmconnsum}
For any two links $\Link_1$ and $\Link_2$ in rational homology spheres $\rats_1$ and $\rats_2$, we have
\begin{equation*}\Zinvuf(\rats_1 \# \rats_2, \Link_1 \sqcup \Link_2) = \Zinvuf(\rats_1,\Link_1)\Zinvuf(\rats_2,\Link_2).\end{equation*}
\end{theorem}

We prove a generalization of Theorem~\ref{thmconnsum} in Section~\ref{secfunc}. See Theorem~\ref{thmfuncgen} in particular. See also Section~\ref{secexppty}. The proof given in Section~\ref{secfunc} can be read without reading the intermediate chapters.
Theorem~\ref{thmconnsum} and Corollary~\ref{corthetazone} yield the following corollary.

\begin{corollary}
 \label{corconnsumTheta}
For any two rational homology spheres $\rats_1$ and $\rats_2$, we have
\begin{equation*}\Theta(\rats_1 \# \rats_2)=\Theta(\rats_1) + \Theta(\rats_2).\end{equation*}
\end{corollary}

\paragraph{Reversing the orientation of the ambient space.}

\begin{lemma}
\label{lemorchangman} Under the assumptions of Lemma~\ref{lemrevor}, we have
\begin{equation*}\Zinv_n(-\crats,\Link,\overline{\partau})=(-1)^n\Zinv_n(\crats,\Link,\partau).\end{equation*}
\end{lemma}
\bp
If $\omega$ is a homogeneous propagating form of $\bigl(C_2(\rats),\partau\bigr)$, then
$(-\omega)$ is a homogeneous propagating form of $(C_2(-\rats),\overline{\partau})$.
Let $\Gamma$ be a degree $n$ numbered Jacobi diagram. When the orientation of $\rats$ is reversed, the orientation of $\check{C}(\rats,\Link;\Gamma)$ is reversed if and only if $\cardlef{T(\Gamma)}$ is odd.
Thus, the integrals will be multiplied by $(-1)^{\cardlef{E(\Gamma)} + \cardlef{T(\Gamma)}}$. We have
$2\cardlef{E(\Gamma)} = 3 \cardlef{T(\Gamma)} + \cardlef{U(\Gamma)}$, and hence $2n=2(\cardlef{E(\Gamma)} -  \cardlef{T(\Gamma)})$.
\eop

\begin{theorem}
\label{thmorchangman} For any link $\Link$ in a rational homology sphere $\rats$, we have
\begin{equation*}\Zinvuf_n(-\rats, \Link) = (-1)^n\Zinvuf_n(\rats, \Link).\end{equation*}
\end{theorem}
\bp Theorem~\ref{thmfstconst} implies
\begin{equation*}\Zinvuf(-\rats,\Link)=\exp\Bigl(-\frac14 p_1(\overline{\partau})\ansothree\Bigr)\prod_{j=1}^k\Bigl(\exp\bigl(-I_{\theta}(K_j,\overline{\partau})\alpha\bigr)\#_j\Bigr) \Zinv(-\crats,\Link,\overline{\partau}).\end{equation*}
Lemma~\ref{lemorchangman} implies $I_{\theta}(K_j,\overline{\partau})=-I_{\theta}(K_j,{\partau})$, and Lemma~\ref{lemrevor} yields $p_1(\overline{\partau})=-p_1(\partau)$. So Lemma~\ref{lemorchangman} implies the result since $\alpha$ and $\ansothree$ vanish in even degrees.
\eop

\chapter{Rationality}
\label{chaprat}

In this chapter, we give equivalent definitions of $\Zinvuf$ based on algebraic intersections of propagating chains, and we prove that $\Zinvuf$ and the anomalies $\alpha$ and $\ansothree$ are rational.

\section{From integrals to algebraic intersections}
\label{secfromintegtointer}

In order to warm up, we first prove the following rationality result, which is due to Sylvain Poirier \cite{poirier} and Dylan Thurston \cite{thurstonconf}, independently.

\begin{proposition}
\label{propratalpha}
The anomaly $\alpha$ of Section~\ref{secanomalpha} is rational, i.e., it belongs to
$\Assis(\RR;\QQ)$.
For any link $\Link \colon \source \to \RR^3$, the value $\Zinvuf(S^3,\Link)$ is also rational. It belongs to
$\Aavis(\source;\QQ)$.
\end{proposition}
\bp Let us fix $n$ and prove that $\alpha_n$ is in $\Assis_n(\RR;\QQ)$.
For any degree $n$ numbered Jacobi diagram $\check{\Gamma}$ on $\RR$, define the smooth map 
\begin{equation*}g(\check{\Gamma}) \colon {\anomq}(\check{\Gamma}) \times (S^2)^{\underline{3n-2} \setminus j_E(E(\check{\Gamma}))} \to (S^2)^{3n-2}\end{equation*}
as the product $\bigl(\prod_{e \in E(\check{\Gamma})} p_{e,S^2}\bigr) \times \id \bigl((S^2)^{\underline{3n-2} \setminus j_E(E(\check{\Gamma}))}\bigr)$. 
Note that a regular value of $g(\check{\Gamma})$ is not in the image of $\partial {\anomq}(\check{\Gamma}) \times  (S^2)^{\underline{3n-2} \setminus j_E(E(\check{\Gamma}))}$.
According to the Morse--Sard theorem~\ref{thmMorseSard}, the set of regular values of $g(\check{\Gamma})$ is dense. 
Since ${\anomq}(\check{\Gamma})$ is compact, so are $\partial {\anomq}(\check{\Gamma})$, the boundary of the domain of $g(\check{\Gamma})$, and the subset of the domain of $g(\check{\Gamma})$ consisting of the points at which the derivative of 
$g(\check{\Gamma})$ is not surjective. Therefore, the set of regular values of $g(\check{\Gamma})$ is open.
Thus, the finite intersection over all the $ \check{\Gamma} \in \CD^c_n(\RR)$ of the sets of regular values of the $g(\check{\Gamma})$ is also open and dense.
Let $\prod_{i=1}^{3n-2} B(x_i)$ be a product of open balls of $S^2$ that is in this intersection.
Then for any $ \check{\Gamma} \in \CD^c_n(\RR)$, the local degree of $g(\check{\Gamma})$ (which is an integer) is constant over $\prod_{i=1}^{3n-2} B(X_i)$.
In particular, if $\omega(i,S^2)$ is a volume-one form of $S^2$ that is supported on $B(X_i)$ for each $i \in \underline{3n-2}$, then $I(\check{\Gamma},o(\check{\Gamma}),\omega(i,S^2))$, which is nothing but this integral local degree, is an integer for any $\check{\Gamma}$ in $\CD^c_n(\RR)$. Thus $\alpha_n$, which is defined in Proposition~\ref{propanom},  is in $\Assis_n(\RR;\QQ)$.

For a fixed $n$ and a given $k$-component link $\Link$ of $S^3$, there exists a similar product $\prod_{i=1}^{3n} B_{\Link}(Y_i)$ of open balls of $S^2$ consisting of points of $(S^2)^{3n}$ that are regular values for all maps 
\begin{equation*}\Bigl(\prod_{e \in E(\Gamma)} p_{e,S^2}\Bigr) \times \id \Bigl((S^2)^{\underline{3 n} \setminus j_E(E({\Gamma}))}\Bigr) \colon {C}(S^3,\Link;\Gamma) \times (S^2)^{\underline{3n} \setminus j_E(E({\Gamma}))} \to (S^2)^{3n}\end{equation*}
associated to Jacobi diagrams $\Gamma$ of $\Davis^e_n(\source)$.
Then if $\omega(i,S^2)$ is a volume-one form of $S^2$ that is supported on $B_{\Link}(Y_i)$ for each $i \in \underline{3 n}$, \begin{equation*}I\biggl(S^3,\Link,\Gamma,\Bigl(p_{S^2}^{\ast}\bigl(\omega(i,S^2)\bigr)\Bigr)\biggr)\end{equation*} is an integer for every $\Gamma$ of $\Davis^e_n(\source)$.

If the link is straight, then Theorem~\ref{thmconststraight} implies that 
$\Zinv_n(\RR^3,\Link,\taust)$ is rational.
Thus, $\Zinv(\RR^3,\Link,\taust)$ is rational for any straight link $\Link$ of $\RR^3$.
In particular $I_{\theta}(K,\taust)$ is rational for any component $K$ of a straight link  $\Link$,
and Theorem~\ref{thmfstconst} together with the rationality of $\alpha$ implies
that $\Zinvuf(S^3,\Link)$ is rational.
\eop

With the notation of the above proof, the
$\propP(i)=p_{S^2}^{-1}(y_i) \subset C_2(S^3)$ for $y_i \in B_{\Link}(Y_i)$ are propagating chains 
such that, for any $\Gamma$ of $\Davis^e_n(\source)$, the intersection over $E(\Gamma)$ of the $p_e^{-1}(\propP(j_E(e)))$ in ${C}(S^3,\Link;\Gamma)$ is transverse.
The integral
$I(S^3,\Link,\Gamma,(p_{S^2}^{\ast}(\omega(i,S^2))))$
is nothing but their algebraic intersection. 

We will use Version~\ref{thmconststraight} of Theorem~\ref{thmfstconst} to replace the configuration space integrals with algebraic intersections in configuration spaces and thus prove the rationality of $\Zinv^s$ for straight links in any rational homology sphere as follows.
\begin{definition}
\label{deftransversealong}
A smooth map $f \colon \submb \to \amba$ is \emph{transverse} \index[T]{transverse!map} to a submanifold $\subsubc$ of $\amba$ along a subset $K$ of $\submb$ if \begin{equation*}\tang_{f(x)}\amba =\tang_xf(T_x\submb)+\tang_{f(x)}\subsubc\end{equation*} for any point $x$ of $K\cap f^{-1}(\subsubc)$.
When $\amba$ or $\submb$ have ridges, we furthermore require this equality to hold when $\amba$ or $\submb$ are replaced by all their open faces (of any dimension).

A smooth map $f \colon \submb \to \amba$ is \emph{transverse} to a submanifold $\subsubc$ of $\amba$ if it is transverse to $\subsubc$ along $\submb$.
Say that a smooth map $f \colon \submb \to \amba$ is \emph{transverse} to a rational chain $\subsubc$ of $\amba$, which is a multiple of a union of compact smooth embedded submanifolds with boundaries and ridges $\cup_{k \in J}\subsubc_k$, if $f$ is transverse to $\subsubc_k$ for any $k \in J$.
\end{definition}

A \emph{rational simplicial chain}, which is a rational combination of simplices in a triangulated smooth manifold, is an example of what we call a rational chain.
Rational multiples of compact immersion images provide other examples of chains. An immersion image will be represented as a union of embedded manifolds by decomposing the domain as a union of compact manifolds with boundaries and ridges glued along their boundaries. 
 
Recall from Notation~\ref{notlessav} that $\Davis^e_{k,\underline{3n}}(\source)$ denotes the set of $\underline{3n}$-numbered degree $k$ Jacobi diagrams with support $\source$ without looped edges. Let $\Davis^e_{\underline{3n}}(\source)=\cup_{k \in \NN}\Davis^e_{k,\underline{3n}}(\source)$\index[N]{Diag@Diagram sets!Dethreen@$\Davis^e_{\underline{3n}}$ $\underline{3n}$-numbered}.
Note that $\Davis^e_n(\source)=\Davis^e_{n,\underline{3n}}(\source)$ but $\underline{3n}$-numbered Jacobi diagrams may have a degree different from $n$.

\begin{definition}
\label{defgenthreenpos}
Say that a family $(\propP(i))_{i \in \underline{3n}}$ of propagating chains of $(C_2(\rats),\tau)$ is in \indexT{general $3n$ position} with respect to a link $\Link \colon \source \to \crats$ if for any $\Gamma \in \Davis^e_{\underline{3n}}(\source)$ and for any subset $E$ of $E(\Gamma)$, 
the map
\begin{equation*}\projp({\Gamma},E)=\prod_{e \in E}p_e \colon C(\rats,\Link;\Gamma) \to \bigl(C_2(\rats)\bigr)^{j_E(E)}\end{equation*}
is transverse to $\prod_{e \in E}\propP(j_E(e))$.
\end{definition}

For such a family $(\propP(i))_{i \in \underline{3n}}$ in general $3n$ position, 
the intersection \begin{equation*}\cap_{e \in E(\Gamma)}p_e^{-1}\Bigl(\propP\bigl(j_E(e)\bigr)\Bigr)\end{equation*} consists of a finite number of points
$x$, which sit in the interior of $C(\rats,\Link;\Gamma)$, and, for each such $x$,
the following conditions are satisfied.
\begin{itemize}
 \item For every edge $e \in E(\Gamma)$, 
 $p_e(x)$ meets the union of smooth embedded $4$-manifolds with boundaries that constitute $\propP(j_E(e))$ in the interior of finitely many of these manifolds. The family $(\Delta_{j_E(e),i})_{i \in J(e,x)}$ of met manifolds is indexed by a finite set $J(e,x)$.\footnote{The $\Delta_{j_E(e),i}$ are smooth embedded $4$-simplices when $\propP(j_E(e))$ is a simplicial chain.}
\item For every map $i \colon E(\Gamma) \to \cup_{e \in E(\Gamma)}J(e,x)$ such that $i(e) \in J(e,x)$, the local maps, from small open neighborhoods of $x$ in $C(\rats,\Link;\Gamma)$ to the product over $E(\Gamma)$ of the fibers of the locally trivialized normal bundles to the $\Delta_{j_E(e),i(e)}$, are local diffeomorphisms.
\end{itemize}
We will prove the following lemma in Section~\ref{secexistransv}.
\begin{lemma}
\label{lemexisttransv} Let $(\crats,\tau)$ be  an asymptotic rational homology $\RR^3$.
For any link $\Link \colon \source \hookrightarrow \crats$, for any integer $n$, there exists a family $(\propP(i))_{i \in \underline{3n}}$ of propagating chains of $(C_2(\rats),\tau)$ in general $3n$-position with respect to $\Link$.
\end{lemma}

Let $(\propP(i))_{i \in \underline{3n}}$ be such a family.
The codimension-two rational chains $p_e^{-1}\bigl(\propP(j_E(e))\bigr)$ are cooriented by the coorientation of $\propP(j_E(e))$ in $C_2(\rats)$.
Define $I\bigl(\Gamma,o(\Gamma),(\propP(i))_{i \in \underline{3n}}\bigr)$ to be the algebraic intersection in $\bigl(C(\rats,\Link;\Gamma),o(\Gamma)\bigr)$ of the chains $p_e^{-1}\bigl(\propP(j_E(e))\bigr)$ over the edges $e$ of $E(\Gamma)$.

For any finite set $\finsetv$, equip $C_{\finsetv}(\rats)$ with a Riemannian metric that is symmetric with respect to permutations of elements of $\finsetv$. Let $d$ denote the associated distance. Our choice of distance will not matter thanks to the following easy lemma.

\begin{lemma}
\label{lemequivdist}
All the distances associated to Riemannian metrics are equivalent on a compact smooth manifold.
\end{lemma}
\bp
Let $g_1$ and $g_2$ be two Riemannian metrics on the compact manifold $M$, let $\norm{.}_1$ and $\norm{.}_2$ be the two associated norms on tangent vectors, and let $d_1$ and $d_2$ be the two associated distances.
View the unit tangent bundle $\ST M$ of $M$ as the set of unit tangent vectors to $M$ with respect to $\norm{.}_1$. Then the image of $\ST M$ under the continuous map $\norm{.}_2$ is a compact interval $[a,b]$ with $a>0$, and we have
\begin{equation*}a \leq \frac{\norm{x}_2}{\norm{x}_1} \leq b\end{equation*}
for any nonzero tangent vector $x$ of $M$.
Let $p$ and $q$ be two distinct points of $M$. For any smooth path $\gamma \colon \left[0,1\right] \to M$ such that $\gamma(0)=p$ and $\gamma(1)=q$, we have
\begin{equation*}d_2(p,q)\leq\int_0^1\norm{\gamma^{\prime}(t)}_2 dt \leq b \int_0^1\norm{\gamma^{\prime}(t)}_1 dt.\end{equation*}
Therefore, we get $d_2(p,q)\leq b d_1(p,q)$.
Similarly, we have $d_1(p,q)\leq \frac{d_2(p,q)}{a}$.
\eop

\begin{definition}
\label{defformdual} For a subset $X$ of $C_{\finsetv}(\rats)$
 and for $\varepsilon >0$, set \begin{equation*}N_{\varepsilon}(X)=\{x \in C_{\finsetv}(\rats) \suchthat  d(x,X) < \varepsilon\}.\end{equation*}\index[N]{NepsilonX@$N_{\varepsilon}(X)$ $\varepsilon$-neighborhood}
 
For a small positive number $\eta$, a closed $2$-form $\omega(i)$ on $C_2(\rats)$ is said to be \emph{$\eta$-dual} to $\propP(i)$, if it is supported in $N_{\eta}(\propP(i))$ and if
$\int_D\omega(i)=\langle D,\propP(i)\rangle_{\!C_2(\rats)\,}$ for any $2$-dimensional disk $D$ embedded in $C_2(\rats)$ transverse to $\propP(i)$ whose boundary sits outside $N_{\eta}(\propP(i))$.
\end{definition}

We will prove the following lemma in Section~\ref{secexistformtransv}.
\begin{lemma}
\label{lemexistformtransv}
Assume Lemma~\ref{lemexisttransv}. Under its hypotheses,
for any $\eta >0$, there exist 
propagating forms $\omega(i)$ of $(C_2(\rats),\tau)$ $\eta$-dual to the $\propP(i)$ of Lemma~\ref{lemexisttransv}.
If $\eta$ is small enough, then we have
\begin{equation*}I\Bigl(\Gamma,o(\Gamma),\bigl(\propP(i)\bigr)_{i \in \underline{3n}}\Bigr)=I\Bigl(\Gamma,o(\Gamma),\bigl(\omega(i)\bigr)_{i \in \underline{3n}}\Bigr)\end{equation*} for any $\Gamma \in \Davis^e_{k,\underline{3n}}(\source)$, where $k \leq n$.
\end{lemma}

Thus, $I(\Gamma,o(\Gamma),(\omega(i))_{i \in \underline{3n}})$ is rational in this case, and we get the following theorem.

\begin{theorem}
\label{thmZinvrat} The anomaly $\ansothree$ is rational, i.e., it belongs to
$\Aavis(\emptyset;\QQ)$.
 Let $(\crats,\tau)$ be  an asymptotic rational homology $\RR^3$.
For any link $\Link \colon \source \to \crats$, the value $\Zinvuf(\rats,\Link)$ belongs to $\Aavis(\source;\QQ)$.
\end{theorem}
\bpo{Proof of the theorem assuming Lemmas~\ref{lemexisttransv} and \ref{lemexistformtransv}}
Theorems~\ref{thmconststraight} and \ref{thmfstconst}, Lemmas~\ref{lemexisttransv} and \ref{lemexistformtransv} imply that $\Zinv(S^3,\partau)=\Zinv(\RR^3,\emptyset,\tau)$ is rational for any $\partau$. So 
$\zinv(S^3,\partau)=p^c(\Zinv(S^3,\partau))$ is rational for any $\tau$, too, and,  by Definition~\ref{defxin}, the anomaly $\ansothree$ is rational.
Therefore, Theorem~\ref{thmconststraight}, Lemmas~\ref{lemexisttransv}, and \ref{lemexistformtransv} also imply that $\Zinv(\crats,\Link,\tau)$ is rational for any pair $((\crats,\tau),\Link)$ such that $(\crats,\tau)$ is an asymptotic rational homology $\RR^3$ and $\Link$ is a straight link with respect to $\tau$. In particular, the integral $I_{\theta}(K,\tau)$ is rational for any component $K$ of a straight link $\Link$. 
Since $p_1(\tau)$ and the anomalies $\alpha$ and $\ansothree$ are rational, 
Theorem~\ref{thmfstconst} now implies that $\Zinvuf(\rats,\Link)$ is rational.
\eop

\section{More on general propagating chains}
\label{secthom}

By a theorem of Ren\'e Thom \cite[Th\'eor\`eme II.27, p. 55]{Thom}, any integral codimension $2$ homology class in a manifold can be represented as the class of an embedded closed (oriented) submanifold. We prove a relative version of this result below, following Thom's original proof in this particular case of his theorem.

\begin{theorem}
\label{thmthom}
 Let $\amba$ be a compact smooth (oriented) manifold with boundary. Let $\subsubc$ be a smooth codimension $2$ closed (oriented) submanifold of $\partial \amba$ such that the homology class of $\subsubc$ vanishes in $H_{\dim(\amba)-3}(\amba;\ZZ)$. Then there exists a compact smooth codimension $2$ submanifold $\submb$
of $\amba$ transverse to $\partial \amba$ whose boundary is $\subsubc$.
\end{theorem}
\bp 
Let us first sketch Thom's proof with his notation.
The normal bundle to $\subsubc$ in $\partial \amba$ is an oriented disk bundle.
It is the pull-back of a universal disk bundle $A_{SO(2)}$ over a compact classifying space $B_{SO(2)}$ via a map $f_{\subsubc}$ from $\subsubc$ to $B_{SO(2)}$. See \cite[p. 145]{milnorsta}.
Like Ren\'e Thom \cite[p. 28, 29]{Thom}, define the Thom space $M(SO(2))$ of $SO(2)$ to be the space obtained from the total space $A_{SO(2)}$ by identifying its subspace $E_{SO(2)}$ consisting of the points in the boundaries of the fibers $D^2$ of $A_{SO(2)}$ with a single point $a$. Regard $B_{SO(2)}$ as the zero section of $A_{SO(2)}$. So $B_{SO(2)}$ sits inside $M(SO(2))$.

The map $f_{\subsubc}$ extends canonically to $\partial \amba$. Its extension $f_{\partial \amba}$ injects the fibers of an open tubular neighborhood of $\subsubc$ in $\partial \amba$ to fibers of $A_{SO(2)}$ and maps the complement of such a neighborhood to $a$.
Thus $\subsubc$ is the preimage of $B_{SO(2)}$.
In order to prove the theorem, it suffices to extend the map $f_{\partial \amba}$ to a map $f_{\amba}$ from $\amba$ to $M(SO(2))$ 
so that, in a neighborhood of any point 
of $f_{\amba}^{-1}(B_{SO(2)})$, the differential of a local projection to the fiber of the normal bundle to $B_{SO(2)}$ composed with 
$f_{\amba}$ 
is well-defined and surjective. (The fiber of the normal bundle to $B_{SO(2)}$ is isomorphic to the tangent space to a fiber of $A_{SO(2)}$.)
Indeed, the compact submanifold $\submb=f_{\amba}^{-1}(B_{SO(2)})$ of $\amba$, with respect to such an extension, has the desired properties.

The map $f_{\partial \amba}$ can be extended as a continuous map, using the fact that $M(SO(2))$ is a $K(\ZZ;2)$
\cite[ii), p. 50]{Thom}. In other words, the only nontrivial homotopy group of $M(SO(2))$ is its $\pi_2$, which is isomorphic to $\ZZ$. 

Let us now give some details about the above sketch and show how the pair $\bigl(B_{SO(2)}, M(SO(2))\bigr)$ can be replaced by
$(\CC P ^{N}, \CC P ^{N+1})$ for some large integer $N$, following  \cite[ii), p. 50]{Thom}.
View the fiber of a disk bundle as the unit disk of $\CC$. The corresponding complex line bundle over $\subsubc$ injects into a trivial complex bundle $\CC^{N+1} \times \subsubc$ as in \cite[Lemma 5.3, p. 61]{milnorsta} for some integer $N$, by some map $(f_{1,\subsubc},\id(\subsubc))$. Therefore, it is the pull-back of the tautological complex line bundle $\gamma^1_N$ over $B^{\prime}_{SO(2)}=\CC P^{N}$ by the map $f^{\prime}_{\subsubc}\colon \subsubc \to \CC P^{N}$ that sends a point $x$ of $\subsubc$ to the image of the fiber over $x$ under $f_{1,\subsubc}$.
The disk bundle $A^{\prime}_{SO(2)}$ associated to $\gamma^1_N$ is diffeomorphic to the normal bundle to 
$\CC P^{N}$ in $\CC P^{N+1}$ by the inverse of the following (orientation-reversing) map:
\begin{equation*}\begin{array}{lll}\CC P^{N+1} \setminus \bigl\{[(0,\dots,0,1)] \bigr\}&\to &\gamma^1_N \\
\left[z_1,\dots, z_{N+1},z\right] & \mapsto & \left(\frac{\overline{z}}{\sum_{i=1}^{N+1}|z_i|^2}(z_1,\dots,z_{N+1}),\left[z_1,\dots,z_{N+1}\right]\right).\end{array}\end{equation*}
This map also shows that the space $M^{\prime}(SO(2))$, obtained from $A^{\prime}_{SO(2)}$ by identifying $E^{\prime}_{SO(2)}=\partial A^{\prime}_{SO(2)}$ to a point, is homeomorphic to the whole $\CC P^{N+1}$.

The long exact sequence associated to the fibration $S^1 \hookrightarrow S^{2N+3} \to \CC P^{N+1}$ implies that $\pi_2(\CC P^{N+1})=\ZZ \left[\CC P^1\right]$ and that $\pi_i(\CC P^{N+1}))$ is trivial for any $i \in \underline{2N+2} \setminus \{2\}$. See Theorem~\ref{thmlongseqhomotopy}. Assume that $2N+2$ is bigger than the dimension of $\amba$ without loss of generality. 

It is not hard to see that $\pi_2(M^{\prime}(SO(2)))=H_2(M^{\prime}(SO(2));\ZZ)$ is freely generated by the class of the image $\left[F\right]$ of a fiber under the identification of $E^{\prime}_{SO(2)}$ with the point $a$.
(Indeed, since $B^{\prime}_{SO(2)}$ is connected, the homology class of $\left[F\right]$ is well-defined. Since $M^{\prime}(SO(2)) \setminus B^{\prime}_{SO(2)}$ is contractible and since any $2$-cycle is homologous to a $2$-cycle that is transverse to $B^{\prime}_{SO(2)}$, any homology class of degree $2$ is a multiple of $\left[F\right]$, which therefore generates $H_2(M^{\prime}(SO(2));\ZZ)$.) 

Extend  $f^{\prime}_{\subsubc}$ to a map $f_{\partial \amba}$ valued in $M^{\prime}(SO(2))$ as before so that $\subsubc = f_{\partial \amba}^{-1}(B^{\prime}_{SO(2)})$.
Recall that any smooth manifold is triangulable \cite{cai1}, \cite{WhiteheadCone}. Fix a triangulation for $(\amba,\partial \amba)$ 
transverse to $\subsubc$. In particular,
$\subsubc=f_{\partial \amba}^{-1}(B^{\prime}_{SO(2)})$ avoids the $1$-skeleton.
Extend $f_{\partial \amba}$ skeleton by skeleton starting with the zero and one-skeleta, for which there is no obstruction to extending $f_{\partial \amba}$ to a map valued in $M^{\prime}(SO(2)) \setminus B^{\prime}_{SO(2)}$, which is connected. 
There is no obstruction to extending $f_{\partial \amba}$, as a map valued in the contractible $M^{\prime}(SO(2)) \setminus B^{\prime}_{SO(2)}$, to the two-skeleton of $(\amba,\partial \amba)$, but such a map would not necessarily extend to the three-skeleton.
Let $f_{\amba}^{(2)}$ be an arbitrary generic extension valued in $M^{\prime}(SO(2))$ of $f_{\partial \amba}$ to the two-skeleton of $(\amba,\partial \amba)$. 
Define the $2$-cochain $c(f_{\amba}^{(2)})$ with $\ZZ$-coefficients such that $c(f_{\amba}^{(2)})(D)$ is the algebraic intersection of $B^{\prime}_{SO(2)}$ and $f_{\amba}^{(2)}(D)$ in $M^{\prime}(SO(2))$ for any $2$-cell $D$ of $\amba$.
Then $f_{\amba}^{(2)}$ extends to the $3$-skeleton if and only if this cochain (which is fixed on $\partial \amba$ and Poincar\'e dual to $\subsubc$ on $\partial \amba$) is a cocycle.
Thus, to prove that $f_{\partial \amba}$ extends to the $3$-skeleton, it suffices to prove that the class of $c(f_{\amba}^{(2)})\vert_{\partial \amba}$ in $H^2(\partial \amba;\ZZ)$ is in the natural image of $H^2(\amba;\ZZ)$, or, equivalently, that its image in $H^3(\amba,\partial \amba;\ZZ)$, by the boundary map of the long cohomology exact sequence of $(\amba,\partial \amba)$, vanishes.
This image is represented by a cochain that maps a $3$-cell $\CB$ of $(\amba,\partial \amba)$ to the algebraic intersection of $\partial \CB$ and $\subsubc$, which is, up to a fixed sign, the algebraic intersection of $\CB$ and $\subsubc$ (pushed inside $\amba$).
Therefore, the class in $H^3(\amba,\partial \amba;\ZZ)$ of this relative cocycle is Poincar\'e dual to the class of $\subsubc$ in $H_{\dim(\amba)-3}(\amba)$, which vanishes. So $f_{\partial \amba}$ can be extended to the $3$-skeleton. Since the next homotopy groups $\pi_i\bigl(M^{\prime}(SO(2))\bigr)$, for $3 \leq i < \dim (\amba)$, vanish, there is no obstruction to extending $f_{\partial \amba}$ to the manifold $\amba$.

Finally, make $f_{\amba}$ smooth, using an approximation theorem \cite[Chapter 2, Theorem 2.6]{hirsch} of continuous maps by smooth maps, and make $f_{\amba}$ transverse to $B^{\prime}_{SO(2)}$, with the help of a transversality theorem \cite[Chapter 3, Theorem 2.1]{hirsch}.
\eop

\begin{corollary}
\label{corpropmanifint}
 If $\rats$ is a $\ZZ$-sphere, for any asymptotically standard parallelization $\tau$ of $\crats$, for any $X \in S^2$, there exists a $4$-dimensional submanifold of $C_2(\rats)$ transverse to the ridges whose boundary is $\projp_{\tau}^{-1}(X)$.
\end{corollary}
\bp First extend $\projp_{\tau}$ as a regular map from a regular neighborhood $N(\partial C_2(\rats))$ of $\partial C_2(\rats)$, where $N(\partial C_2(\rats))$ is a smooth cobordism with ridges embedded in $C_2(\rats)$ from a smooth manifold $\partial C^{\prime}_2(\rats)$ without ridges to $\partial C_2(\rats)$, and $N(\partial C_2(\rats))$ is homeomorphic to the product $\left[0,1\right] \times \partial C_2(\rats)$. Then apply Theorem~\ref{thmthom} to
$C^{\prime}_2(\rats)=C_2(\rats) \setminus \Int(N(\partial C_2(\rats)))$ and to $p_{\tau}\vert_{\partial C^{\prime}_2(\rats)}^{-1}(X)$.
\eop

When $\rats$ is a $\QQ$-sphere, perform the same first step as in the above proof. Take a collar neighborhood of $\partial C^{\prime}_2(\rats)$ in $N(\partial C_2(\rats))$, which is (diffeomorphic to and) identified with $\left[0,8\right] \times \partial C^{\prime}_2(\rats)$ so that $\partial C^{\prime}_2(\rats)=\{0\} \times \partial C^{\prime}_2(\rats)$. Assume that $\projp_{\tau}$ factors through the projection to $\partial C^{\prime}_2(\rats)$ on $\left[0,8\right] \times \partial C^{\prime}_2(\rats)$. 
There exists a positive integer $k$ such that $k p_{\tau}\vert_{\partial C^{\prime}_2(\rats)}^{-1}(X)$ is null-homologous in $C^{\prime}_2(\rats)$. Let $p_k \colon S^2 \to S^2$ be a degree $k$ map that does not fix $X$ and such that $X$ is a regular value of $p_k$ with $k$ preimages. Then $(p_k \circ p_{\tau}\vert_{\partial C^{\prime}_2(\rats)})^{-1}(X)$ bounds a $4$-manifold $\propP^{\prime}$ properly embedded in $C^{\prime}_2(\rats)$ according to Theorem~\ref{thmthom}.
For $j\in \underline{k}$, let $\{\gamma_j \colon \left[0,4\right] \to S^2\}_{j\in \underline{k}}$ be a collection of smooth injective paths ending at $X=\gamma_j(4)$ whose images do not meet outside $X$ and such that $p_k^{-1}(X)=\{\gamma_j(0) \suchthat j\in \underline{k}\}$. 
Also assume that all the derivatives of $\gamma_j$ vanish at $0$ and $4$.
Consider \begin{equation*}\begin{array}{llll}p_{\left[0,8\right]} \times \projp_{\tau} \colon & \left[0,8\right] \times \partial C^{\prime}_2(\rats) & \to &\left[0,8\right] \times S^2\\ & (t,x) & \mapsto &(t,\projp_{\tau}(x)). \end{array}\end{equation*}
Then \begin{multline*}\propP= p^{-1}_{\tau}\vert_{N(\partial C_2(\rats)) \setminus  ([0,4[ \times \partial C^{\prime}_2(\rats)) }(X) +\frac{1}{k}\propP^{\prime}\\  + \frac{1}{k}(p_{\left[0,8\right]} \times \projp_{\tau})\vert_{\left[0,4\right] \times \partial C^{\prime}_2(\rats)}^{-1}\Bigl(\bigl\{\bigl(t,\gamma_j(t)\bigr) \suchthat j\in \underline{k},t\in\left[0,4\right]  \bigr\}\Bigr)\end{multline*}
is a propagating chain of $(\crats,\tau)$. See Figure~\ref{figproprat}.

\bfig
\centering
\begin{tikzpicture}[scale=.9] 
\draw (3,0) -- (3,5) (7,1) -- (7,4) (11,1.6) -- (11,4) (12.6,0) -- (13,.4) -- (13,4.4) -- (12.6,5)
(1.65,4.1) node{\scriptsize $C^{\prime}_2(\rats)$} (8,4.5) node{\scriptsize $N(\partial C_2(\rats))$} (2.25,3) node{\scriptsize $\propP^{\prime}$} (5,3.25) node{\scriptsize $\left[0,4\right] \times \partial C^{\prime}_2(\rats)$} (9,3.25) node{\scriptsize $\left[4,8\right] \times \partial C^{\prime}_2(\rats)$}
(7.6,-.4) node{$(p_{\left[0,8\right]} \times \projp_{\tau})\vert_{\left[0,4\right] \times \partial C^{\prime}_2(\rats)}^{-1}\left(\{(t,\gamma_j(t));j\in \underline{k},t\in\left[0,4\right]  \}\right)$} (10,1.1) node{$kp_{\tau}\vert_{N(\partial C_2(\rats)) \setminus  ([0,4[ \times \partial C^{\prime}_2(\rats)) }^{-1}(X)$} (3,4.75) node[left]{\scriptsize $\partial C^{\prime}_2(\rats)$} (12.8,4.75) node[left]{\scriptsize $\partial C_2(\rats)$} ;
\draw [very thick] (7,2.4) -- (13,2.4) (7,2.4) -- (1.5,2.4) (7,2.4) .. controls (6.6,2.4) and (4,2.8) ..  (3,2.8) -- (1.5,2.8) (7,2.4) .. controls (6.6,2.4) and (4,2) .. (3,2) -- (1.5,2) (7,2.4) .. controls (6.6,2.4) and (4,1.6) .. (3,1.6) -- (1.5,1.6);
\draw [very thick,dotted] (.6,1.6) -- (1.5,1.6) (.6,2) -- (1.5,2)  (.6,2.4) -- (1.5,2.4) (.6,2.8) -- (1.5,2.8);
\draw [very thin,dotted] (.6,3.85) -- (1.5,3.85);
\draw [very thin,->] (1.5,3.85) -- (2.9,3.85);
\draw [very thin,<->] (3.1,4.25) -- (12.9,4.25);
\draw [very thin,<->] (3.1,3) -- (6.9,3);
\draw [very thin,<->] (7.1,3) -- (10.9,3);
\draw [very thin,->] (7.6,.1) -- (5,1.85);
\draw [very thin,->] (9.8,1.4) -- (9,2.3);
\draw [very thin,->] (10.2,1.4) -- (11.8,2.3);
\end{tikzpicture}
\caption{A multiple of a propagating chain of $(\crats,\tau)$}
\label{figproprat}

\end{figure}

\section{Existence of transverse propagating chains} 
\label{secexistransv}

In this section, we prove Lemma~\ref{lemexisttransv}.

In order to warm up, we first prove a weak version of this lemma. The proof is a straightforward adaptation of a proof of Thom \cite[p. 23, 24, Lemma I.4]{Thom}.

Assume that $\rats$ is an integer homology $3$-sphere. Let $(\propP(i))_{i \in \underline{3n}}$ be a family of propagating chains of $(C_2(\rats),\tau)$ for an asymptotically standard parallelization $\tau$ of $\crats$.
Assume that these chains are submanifolds of $C_2(\rats)$ transverse to $\partial C_2(\rats)$ as in Corollary~\ref{corpropmanifint}.

Let $N(\propP(i))$ denote the normal bundle to $\propP(i)$ embedded in $C_2(\rats)$ as a tubular neighborhood whose fibers $N_x(\propP(i))$ over a point $x \in \propP(i)$ are disks embedded in $C_2(\rats)$.
Let $(K_{i,j})_{j \in J(i)}$ be a finite cover of $\propP(i)$ by compact subspaces $K_{i,j}$ embedded in open subspaces $\CO_{i,j}$ of $\propP(i)$ equipped with bundle charts $\bigl(\psi_{i,j} \colon N(\propP(i))\vert_{\CO_{i,j}} \rightarrow A_{i,j} \times D^2\bigr)_{j \in J}$.
We assume that the above $A_{i,j}$ are $\RR^4$, $\RR^+ \times \RR^3$
or $(\RR^+)^2 \times \RR^2$, and that the $\psi_{i,j}$ are bundle charts over 
diffeomorphisms
$\phi_{i,j} \colon \CO_{i,j} \rightarrow A_{i,j}$.
We also assume that for any $\{j,k\} \subset J(i)$, for any $x \in \CO_{i,j} \cap \CO_{i,k}$, the map \begin{equation*}v \mapsto p_{D^2} \circ \psi_{i,k} \circ \psi_{i,j}^{-1}\bigl(\phi_{i,j}(x),v\bigr)\end{equation*} is a linear map of $SO(2)$, where $p_{D^2}$ denotes the natural projection onto $D^2$.
Consider the space $\CH_i$ of smooth diffeomorphisms of $N(\propP(i))$ that are isotopic to the identity map, that fix a neighborhood of $\partial N(\propP(i))$ pointwise, and that map any fiber of $N(\propP(i))$ to itself.
Equip this space $\CH_i$ with the following distance $d$.\footnote{This distance induces the strong (or weak, which is the same since $N(\propP(i))$ is compact) $C^1$-topology. See \cite[Chapter 2, p.35]{hirsch}.} 
Each fiber is equipped with the distance $d_{\propP}$ induced by the norm of $\RR^2$.
This allows us to define a $C^0$ distance $d_0$ between two elements $h$ and $k$ of $\CH_i$ by
\begin{equation*}d_0(h,k)=\sup_{x \in N(\propP(i))}d_{\propP}\bigl(h(x),k(x)\bigr).\end{equation*}
Since $A_{i,j} \times D^2$ is a subset of $\RR^6$, the differential of a map $\psi_{i,j} \circ h \circ \psi_{i,j}^{-1}$ for $h \in \CH_i$ maps every element $x$ of $A_{i,j} \times D^2$ to a linear map of $\RR^6$.
The norm $\norm{L}$ of a linear map $L$ of $\RR^6$ is defined to be $\norm{L}=\sup_{x \in S^5}\norm{L(x)}$.
For $h$ and $k$ in $\CH_i$, set 
\begin{equation*}d_{(1)}(h,k)=\sup_{j \in J(i),x \in \phi_{i,j}(K_{i,j}) \times D^2}\left(\norm{\tang_x(\psi_{i,j} \circ h \circ \psi_{i,j}^{-1}) - \tang_x(\psi_{i,j} \circ k \circ \psi_{i,j}^{-1})}\right)\end{equation*}
and
\begin{equation*}d(h,k)=\sup\left(d_0(h,k),d_{(1)}(h,k)\right).\end{equation*}

\begin{lemma}
\label{lemexisttransvbaby}
Under the above hypotheses,
there is a dense open subset of $\prod_{i=1}^{3n}\CH_i$ such that for any $(h_{i})$ in this subset, the chains obtained from the $\propP(i)$ by replacing $\propP(i)$ with $h_i^{-1}(\propP(i))$ are in general $3n$ position with respect to $\Link$ in the sense of Definition~\ref{defgenthreenpos}.
\end{lemma}

\bp We will first list finitely many sufficient conditions on the $(h_i)$, which guarantee the conclusion \say{For any graph $\Gamma$ of $\Davis^e_{\underline{3n}}(\source)$, and for any subset $E$ of $E(\Gamma)$, the map 
\begin{equation*}\projp({\Gamma},E)=\prod_{e \in E}p_e \colon C(\rats,\Link;\Gamma) \to \bigl(C_2(\rats)\bigr)^{j_E(E)}\end{equation*}
is transverse to $\prod_{e \in E}h_{j_E(e)}^{-1}\propP(j_E(e))$.}.
Next we will prove that each of these conditions is realized in an open dense subset of $\prod_{i=1}^{3n}\CH_i$.

Extend the elements of $\CH_i$ to diffeomorphisms of $C_2(\rats)$, by the identity map of $C_2(\rats)\setminus \mathring{N}(\propP(i))$.
The propagating chains obtained from the $\propP(i)$ by replacing $\propP(i)$ with $h_i^{-1}(\propP(i))$ are in general $3n$ position with respect to $\Link$ if and only if the following condition
$(\ast)(\Gamma,E,\ell)$ holds for any triple $(\Gamma,E,\ell)$,
where $\Gamma \in \Davis^e_{\underline{3n}}(\source)$, $E \subseteq E(\Gamma)$,
and $\ell$ is a map $\ell \colon E \to \cup_{i \in \underline{3n}}J(i)$ such that $\ell(e) \in J(j_E(e))$.

$(\ast)(\Gamma,E,\ell) \colon$ The map $\projp(\Gamma,E)$ is transverse to $\prod_{e \in E}h_{j_E(e)}^{-1}\bigl(\propP(j_E(e))\bigr)$ along $\projp(\Gamma,E)^{-1}\bigl(\prod_{e \in E}h_{j_E(e)}^{-1}(K_{j_E(e),\ell(e)})\bigr)$ (as in Definition~\ref{deftransversealong}). 

In order to prove our lemma, it suffices to prove that, for any of the finitely many $(\Gamma,E,\ell)$ as above, the 
set $\CH(\Gamma,E,\ell)$ in which the condition $(\ast)(\Gamma,E,\ell)$ is realized is a dense open subset of $\prod_{i=1}^{3n}\CH_i$.
This condition is equivalent to

$(\ast)(\Gamma,E,\ell) :$
the map $\bigl(\prod_{e \in E}h_{j_E(e)}\bigr) \circ \projp(\Gamma,E)$ is transverse to $\prod_{e \in E}\propP\bigl(j_E(e)\bigr)$ along $\projp(\Gamma,E)^{-1}\bigl(\prod_{e \in E}h_{j_E(e)}^{-1}(K_{j_E(e),\ell(e)})\bigr)$.\\
\noindent Set 
\begin{equation*}C_{E,\ell}=C(\rats,\Link;\Gamma) \cap \projp(\Gamma,E)^{-1} \left(\prod_{e \in E}\psi_{j_E(e),\ell(e)}^{-1}\bigl(\phi_{j_E(e),\ell(e)}(K_{j_E(e),\ell(e)}) \times D^2\bigr) \right).\end{equation*}
The condition $(\ast)(\Gamma,E,\ell)$  can equivalently be written as \say{ $(0)_{e \in E}$ is a regular value of the map \begin{equation*}\prod_{e \in E}\bigl(p_{D^2} \circ \psi_{j_E(e),\ell(e)} \circ h_{j_E(e)} \circ p_e\bigr)\end{equation*}
on $C_{E,\ell}$.}

Note that the set of regular values of this map on the compact domain $C_{E,\ell}$ is open. Therefore, if $\left(h_i\right)_{i\in \underline{3n}} \in \CH(\Gamma,E,\ell)$ and if the $d_0(h_i,h^{\prime}_i)$ are small enough, the preimage
of $(0)_{e \in E}$ under the restriction of $\prod_{e \in E}p_{D^2} \circ \psi_{j_E(e),\ell(e)} \circ h^{\prime}_{j_E(e)} \circ \projp(\Gamma,E)$ to $C_{E,\ell}$ consists of regular points of $\prod_{e \in E} p_{D^2} \circ \psi_{j_E(e),\ell(e)} \circ h_{j_E(e)} \circ \projp(\Gamma,E)$. These points are regular for $\prod_{e \in E}p_{D^2} \circ \psi_{j_E(e),\ell(e)} \circ h^{\prime}_{j_E(e)} \circ \projp(\Gamma,E)$ provided that the $d(h_i,h^{\prime}_i)$ are small enough.
Therefore, the set $\CH(\Gamma,E,\ell)$ is open.

To prove density, we use explicit deformations of the $h_i \in \CH_i$, for a given $\left(h_i\right)_{i\in \underline{3n}} \in \prod_{i=1}^{3n}\CH_i$.
Fix a smooth map $\chi \colon D^2 \to \left[0,1\right]$, which maps the disk of radius $\frac12$ to $1$ and the complement of the disk of radius $\frac34$ to $0$.
For each compact $K_{i,j}$, such that $j \in J(i)$, fix a smooth map $\chi_{i,j} \colon A_{i,j} \to \left[0,1\right]$ that maps $\phi_{i,j}(K_{i,j})$ to $1$, and that vanishes outside a compact of $A_{i,j}$.
For $w \in D^2$, define
\begin{equation*}\begin{array}{llll}h_{i,j,w} \colon & A_{i,j} \times D^2 & \to & A_{i,j} \times D^2\\
&(x,v) &\mapsto & \bigl(x,v+\chi(v)\chi_{i,j}(x)w\bigr).
\end{array}\end{equation*}
Note that $h_{i,j,w}$ is a diffeomorphism as soon as $\norm{w}$ is smaller than a fixed positive number $\eta <\frac12$.
Extend $\psi_{i,j}^{-1} \circ h_{i,j,w} \circ \psi_{i,j}$ by the identity map outside $N(\propP(i))\vert_{\CO_{i,j}}$. Note that there exists a constant $C$ such that $d(\psi_{i,j}^{-1} \circ h_{i,j,w} \circ \psi_{i,j} \circ h_i, h_i) \leq C \norm{w}$.

Thus, it suffices to prove that, for any $\varepsilon$ such that $0<\varepsilon<\eta$, there exists $(w_e)_{e \in E}$ with $\norm{w_e} < \varepsilon$ such that the restriction of \begin{equation*}\prod_{e \in E}\bigl(\psi_{j_E(e),\ell(e)}^{-1} \circ h_{{j_E(e)},{\ell(e)},w_e} \circ \psi_{j_E(e),\ell(e)} \circ h_{j_E(e)}\bigr)\circ \projp(\Gamma,E)\end{equation*} to 
$(C_{E,\ell})$ is transverse to $\prod_{e \in E}\propP(j_E(e))$ along $C_{E,\ell}$.
Since this happens when $(-w_e)_{e \in E}$ is a regular value of the restriction of \begin{equation*}\prod_{e \in E}\left(p_{D^2}\circ\psi_{j_E(e),\ell(e)} \circ h_{j_E(e)}\right)\circ \projp(\Gamma,E)\end{equation*} to 
$C_{E,\ell}$, and since such regular values form a dense set according to the Morse--Sard theorem, the lemma is proved.
\eop

The magic in the Thom proof above is that it proves the density of manifolds in general $3n$ position without bothering to construct a single one.

Lemma~\ref{lemexisttransvbaby} does not quite prove Lemma~\ref{lemexisttransv} for two reasons. First, the $h_i$ do not fix the boundary of $\partial C_2(\rats)$ pointwise, so the perturbations $h_i^{-1}(\propP(i))$ are no longer propagating chains of $(C_2(\rats),\tau)$.
Second, we have to deal with immersed manifolds (multiplied by an element of $\QQ$) rather than embedded ones when $\rats$ is not an integer homology sphere.

To deal with this latter issue, we start with immersions $f_i$ of manifolds $\tilde{\propP}(i)$ to $C_2(\rats)$ whose images $f_i(\tilde{\propP}(i))$ represent chains $k\propP$ as in the end of Section~\ref{secthom},
and (extended) immersions $f_i$ of the pull-backs $N(\tilde{\propP}(i))$ of the normal bundles to their images.
Our immersions $f_i$ have the following properties. The restriction to
$f_i^{-1}\bigl(C^{\prime}_2(\rats)\bigr)$ of each immersion $f_i$ is an embedding. The preimage $f_i^{-1}\bigl(N(\partial C_2(\rats))\bigr)$ has $k$ connected components $C_{j,i}$ ($j \in \underline{k}$)
in $N(\tilde{\propP}(i))$, and $f_i$ embeds each of these $k$ connected components into $N(\partial C_2(\rats))$. We will think of the intersection with a preimage of $f_i(\tilde{\propP}(i)) \cap N(\partial C_2(\rats))$ as the sum of the intersections
with the preimages of the $f_i(C_{i,j})$, and argue with covers of $\tilde{\propP}(i)$ rather than covers of its image. So this latter issue is not a big one---if we do not require the boundary conditions. We keep this in mind, and we no longer discuss this issue.

The first issue is more serious. We want the boundaries of our propagating chains to be equal to $\projp_{\tau}\vert_{\partial C_2(\rats)}^{-1}(X_i)$ for some $X_i \in S^2$.
Recall that $\projp_{\tau}$ also denotes a regular extension of $\projp_{\tau}$ on $N(\partial C_2(\rats))$, that a collar $\left[0,8\right]\times \partial C^{\prime}_2(\rats)$ of $\partial C^{\prime}_2(\rats)$ in $N(\partial C_2(\rats))$ has been fixed, and that $\projp_{\tau}$ factors through the natural projection onto $\partial C^{\prime}_2(\rats)$ in this collar.
For an interval $I$ included in $\left[0,8\right]$, set 
\begin{equation*}N_{I}= I \times \partial C^{\prime}_2(\rats).\end{equation*}
For $a \in \left[1,8\right]$, set
\begin{equation*}N_{\left[a,9\right]}=N\bigl(\partial C_2(\rats)\bigr) \setminus \left[0,a\right[ \times \partial C^{\prime}_2(\rats).\end{equation*}

We will actually impose that our propagating chains intersect $N_{\left[7,9\right]}$ as $\projp_{\tau}\vert_{N_{\left[7,9\right]}}^{-1}(X_i)$, by modifying our immersions $f_i$ provided by the construction of the end of Section~\ref{secthom}, only on $N_{[4,7[}$.

We first describe appropriate choices for the $X_i$, to allow transversality near the boundaries.

Let $\Gamma$ of $\Davis^e_{\underline{3n}}(\source)$, 
let $E$ be a subset of $E(\Gamma)$.
A condition on $(X_i)_{i \in \underline{3n}}$ is that $(X_i)_{i \in j_E(E)}$ is a regular (for the restriction to any stratum of $C(\rats,\Link;\Gamma)$) value of the map 
\begin{equation*}\prod_{e \in E}\projp_{\tau}\circ p_e\end{equation*}
from \begin{equation*}C(\Gamma,E)= C(\rats,\Link;\Gamma)\cap \bigcap_{e \in E}p_e^{-1}\Bigl(N\bigl(\partial C_2(\rats)\bigr)\Bigr)\end{equation*} to $(S^2)^{j_E(E)}$.
According to the Morse--Sard theorem~\ref{thmMorseSard}, this condition holds when  $(X_i)_{i \in \underline{3n}}$ is in a dense subset of $(S^2)^{3n}$, which is furthermore open since $C(\Gamma,E)$ is compact.
Thus, this condition holds for any of the finitely many pairs $(\Gamma,E)$ as above, if $(X_i)_{i \in \underline{3n}}$ belongs to the intersection of the corresponding open dense subsets of $(S^2)^{3n}$, which is still open and dense.

Fix $(X_i)_{i \in \underline{3n}}$ in this open dense subset of $(S^2)^{3n}$.
Now, we refer to the proof of Lemma~\ref{lemexisttransvbaby} and adapt it to produce the desired family of propagating chains of Lemma~\ref{lemexisttransv}.
We fix a finite cover $(K_{i,j})_{j \in J(i)}$ of $\tilde{\propP}(i) \cap f_i^{-1}\left(C^{\prime}_2(\rats) \cup N_{\left[0,7\right]}\right)$. We assume that this cover contains a special element $K_{i,0}=\tilde{\propP}(i) \cap f_i^{-1}\left(N_{\left[5,7\right]}\right)$ and that $K_{i,j}$ is a compact subset of
$\tilde{\propP}(i) \cap f_i^{-1}\left(C^{\prime}_2(\rats) \cup N_{\left[0,5\right]}\right)$ for any $j \in J^{\prime}(i)=J(i) \setminus \{0\}$.

When $j\in J^{\prime}(i)$, $K_{i,j}$ is embedded in an open subspace $\CO_{i,j}$ of $\tilde{\propP}(i) \cap f_i^{-1}\left(C^{\prime}_2(\rats) \cup N_{\left[0,6\right]}\right)$. These $\CO_{i,j}$ are  diffeomorphic to $\RR^4$ via diffeomorphisms $\phi_{i,j} \colon \CO_{i,j} \rightarrow \RR^4$, and we have bundle charts $\bigl(\psi_{i,j} \colon N(\tilde{\propP}(i))\vert_{\CO_{i,j}} \rightarrow \RR^4 \times D^2\bigr)$, for ${j \in J^{\prime}(i)}$, as in the proof of Lemma~\ref{lemexisttransvbaby}.

The bundle $N(K_{i,0})$ is trivialized by $\projp_{\tau}$ in the following way. Fix a small neighborhood $D_i$ of $X_i$ in $S^2$ diffeomorphic to the standard disk $D^2$ and a diffeomorphism $\psi_{D,i}$ from $D_i$ to $D^2$. Without loss of generality, assume that $N(K_{i,0})=f_i^{-1}(\projp_{\tau}^{-1}(D_i) \cap N_{\left[5,7\right]})$, and identify $N(K_{i,0})$ with $K_{i,0} \times D^2$ so that the projection onto $D^2$ may be expressed as $p_{D^2}= \psi_{D,i} \circ \projp_{\tau} \circ f_i$.

The space $\CH_i$ is now the space of smooth diffeomorphisms of $N(\tilde{\propP}(i))$
that are isotopic to the identity map, that fix a neighborhood of $\partial N(\tilde{\propP}(i))$ and a neighborhood of $f_i^{-1}\bigl(f_i(N(\tilde{\propP}(i)))\cap (N_{\left[7,9\right]})\bigr)$ pointwise, and that map any fiber of $N(\tilde{\propP}(i))$ to itself.
The space $\CH_i$ is equipped with a distance similar to that described before Lemma~\ref{lemexisttransvbaby}.

We want to prove that the subspace of $\prod_{i \in \underline{3n}}\CH_i$ consisting of the $(h_i)_{i\in \underline{3n}}$ such that the
$f_i\bigl(h_i^{-1}(\tilde{\propP}(i))\bigr)$ are in general $3n$ position with respect to $\Link$, in the sense of Definition~\ref{defgenthreenpos}, is open and dense.

It is open as in the proof of Lemma~\ref{lemexisttransvbaby}.\footnote{If we ask only for transversality of the
\begin{equation*}\projp({\Gamma},E)=\prod_{e \in E}p_e \colon C(\rats,\Link;\Gamma) \to (C_2(\rats))^{j_E(E)}\end{equation*} to $\prod_{e \in E}f_{j_E(e)}\bigl(h_{j_E(e)}^{-1}(\tilde{\propP}(j_E(e)))\bigr)$
along 
$\prod_{e \in E} p_e^{-1}(C^{\prime}_2(\rats) \cup N_{\left[0,5\right]})$, then density could also be proved as in Lemma~\ref{lemexisttransvbaby}.}

Moreover, for any $\Gamma$ of $\Davis^e_{\underline{3n}}(\source)$, for any triple 
$(E_X,E_N,E_C)$ of pairwise disjoint subsets of $E(\Gamma)$, the subset $\CH(\Gamma,E_X,E_N,E_C)$ of $\prod_{i \in \underline{3n}}\CH_i$ such that
the restriction of $\projp({\Gamma},E_N\cup E_C)$ to $C(\rats,\Link;\Gamma) \cap \projp({\Gamma},E_X)^{-1}\prod_{e \in E_X}(N_{\left[5,9\right]} \cap \projp_{\tau}^{-1}(X_{j_E(e)}))$ is transverse to 
\begin{equation*}\prod_{e \in E_N \cup E_C}\Biggl(f_{j_E(e)}\biggl(h_{j_E(e)}^{-1}\Bigl(\tilde{\propP}\bigl(j_E(e)\bigr)\Bigr)\biggr)\Biggr)\end{equation*}
along \begin{multline*} \projp({\Gamma},E_X)^{-1}\prod_{e \in E_X}\bigl(N_{\left[5,9\right]} \cap \projp_{\tau}^{-1}(X_{j_E(e)})\bigr)
\cap   \projp({\Gamma},E_N)^{-1}\prod_{e \in E_N}\bigl(N_{\left[5,9\right]}\bigr)\\
\cap \projp({\Gamma},E_C)^{-1}\prod_{e \in E_C}\bigl(C^{\prime}_2(\rats) \cup N_{\left[0,5\right]}\bigr)
 \end{multline*}
is open.\footnote{Our hypotheses on $(X_i)_{i \in \underline{3n}}$ guarantee that $C(\rats,\Link;\Gamma) \cap \projp({\Gamma},E_X)^{-1}\prod_{e \in E_X}(N_{\left[5,9\right]} \cap \projp_{\tau}^{-1}(X_{j_E(e)}))$ is a manifold.}
The $\CH(\Gamma,E_X,\emptyset,E_C)$ are furthermore dense as in the proof of Lemma~\ref{lemexisttransvbaby}.

In order to prove Lemma~\ref{lemexisttransv}, it suffices to prove that for any $\Gamma$ of $\Davis^e_{\underline{3n}}(\source)$, for any pair $(E_N,E_C)$ of disjoint subsets of $E(\Gamma)$, the subset $\CH(\Gamma,\emptyset,E_N,E_C)$ of $\prod_{i \in \underline{3n}}\CH_i$ is dense.
To do that, we fix $(h_i)_{i \in \underline{3n}} \in \prod_{i \in \underline{3n}}\CH_i$ and $\varepsilon \in \left]0,1\right[$, and we prove that there exists $(h^{\prime}_i)_{i \in \underline{3n}} \in \CH(\Gamma,\emptyset,E_N,E_C)$ such that $\max_{i \in \underline{3n}}(d(h_i,h^{\prime}_i)) < \varepsilon$.
There exists $\eta \in\left]0,1\right[$ such that the restriction of $h_i$ to \begin{equation*}f_i^{-1}\Bigl(f_i\bigl(N(\tilde{\propP}(i))\bigr) \cap N_{\left[7-2\eta,9\right]}\Bigr)\end{equation*} is the identity map for any $i \in j_E(E_N)$.

For $i \in J_E(E_N)$, our $h^{\prime}_i$ will be constructed as some $h_{i,\eta,w} \circ h_i$. 
Let $\chi_{\eta}$ be a smooth map from $\left[4,9\right]$ to $\left[0,1\right]$ that maps $\left[5,7-2\eta\right]$ to $1$ and that maps the complement of $\left[5-\eta,7-\eta\right]$ to $0$.
Recall our smooth map $\chi \colon D^2 \to \left[0,1\right]$, which maps the disk of radius $\frac12$ to $1$ and the complement of the disk of radius $\frac34$ to $0$.
For $w \in D^2$ define
\begin{equation*}\begin{array}{llll}h_{\eta,w} \colon & \left[4,9\right] \times D^2 & \to &  D^2\\
&(t,v) &\mapsto & v+\chi(v)\chi_{\eta}(t)w.
\end{array}\end{equation*}
Define $h_{i,\eta,w} \in \CH_i$, for $w$ sufficiently small, to coincide with the identity map outside 
\begin{equation*}f_i^{-1}\Bigl(f_i\bigl(N(\tilde{\propP}(i))\bigr) \cap N_{\left[5-\eta,7-\eta\right]}\Bigr),\end{equation*} and with the map that sends
$(p,v) \in\left(\tilde{\propP}(i)\cap f_i^{-1}(N_{\{t\}})\right) \times D^2 $ \begin{equation*}  \biggl(\subset \Bigl( \bigl(\tilde{\propP}(i) \cap f_i^{-1}(N_{\left[5-\eta,7-\eta\right]})\bigr) \times D^2 =N\bigl(\tilde{\propP}(i)\bigr) \cap f_i^{-1}(N_{\left[5-\eta,7-\eta\right]}) \Bigr) \biggr)\end{equation*}
to $(p,h_{\eta,w}(t,v))$, for $t \in \left[5-\eta,7-\eta\right]$.\footnote{For any $t \in \left[0,9\right]$, we assume $f_i\left(\left(\tilde{\propP}(i)\cap f_i^{-1}(N_{\{t\}})\right) \times D^2\right) \subset N_{\{t\}}$, without loss of generality.}
There exists $u \in \left]0,1\right[$ such that, as soon as $\norm{w} <u$, $h_{i,\eta,w}$ is indeed a diffeomorphism and $d(h_{i,\eta,w}\circ h_i, h_i )< \varepsilon$.
Fix $(h^{\prime}_i)_{i \in J_E(E_C)}$ such that $d(h_i,h^{\prime}_i) < \varepsilon$ and $(h^{\prime}_i)_{i \in J_E(E_C)} \times (h_i)_{i \notin J_E(E_C)}$ is in the dense open set $\cap_{E_x \subseteq E_N}\CH(\Gamma,E_x,\emptyset,E_C)$ (this does not impose anything on  $(h_i)_{i \notin J_E(E_C)}$).

After reducing $u$, we may now assume that as soon as $\norm{w} <u$, for any $E_x \subseteq E_N$,
$\projp({\Gamma},E_C\cup E_x)$ is transverse to 
\begin{multline*}\prod_{e \in E_C}f_{j_E(e)}\biggl(\Bigl(h^{\prime}_{j_E(e)}\Bigr)^{-1}\Bigl(\tilde{\propP}\bigl(j_E(e)\bigr)\Bigr)\biggr)\\
\times \prod_{e \in E_x}f_{j_E(e)}\biggl(\Bigl(h_{j_E(e),\eta,w_{j_E(e)}} \circ h_{j_E(e)}\Bigr)^{-1}\Bigl(\tilde{\propP}\bigl(j_E(e)\bigr)\Bigr)\biggr)
\end{multline*}
along \begin{equation*}\projp({\Gamma},E_C)^{-1}\left(\prod_{e \in E_C}\left(C^{\prime}_2(\rats) \cup N_{\left[0,5\right]}\right)\right)
\cap \projp({\Gamma},E_x)^{-1}\left(\prod_{e \in E_x}\left(N_{\left[7-2\eta,9\right]}\right)\right)\end{equation*}
since $ h_{j_E(e)}^{-1}\bigl(\tilde{\propP}(j_E(e))\bigr)=\projp_{\tau}^{-1}(X_{j_E(e)}))$ on $N_{\left[7-2\eta,9\right]}$.
Furthermore, $\projp({\Gamma},E_C\cup E_N)$ is transverse to 
\begin{multline*} M\left((h^{\prime}_i)_{i \in J_E(E_C)}, (h_{i,\eta,w_i})_{i \in J_E(E_N)}\right)
=\prod_{e \in E_C}\Biggl(f_{j_E(e)}\biggl(\Bigl(h^{\prime}_{j_E(e)}\Bigr)^{-1}\Bigl(\tilde{\propP}\bigl(j_E(e)\bigr)\Bigr)\biggr)\Biggr)\\
\times \prod_{e \in E_N}f_{j_E(e)}\biggl(\Bigl(h_{j_E(e),\eta,w_{j_E(e)}} \circ h_{j_E(e)}\Bigr)^{-1}\Bigl(\tilde{\propP}\bigl(j_E(e)\bigr)\Bigr)\biggr)
\end{multline*}
along \begin{equation*}\projp({\Gamma},E_C)^{-1}\left(\prod_{e \in E_C}\left(C^{\prime}_2(\rats) \cup N_{\left[0,5\right]}\right)\right)
\cap \projp({\Gamma},E_N)^{-1}\left(\prod_{e \in E_N}N_{\left[5,9\right]}\right),\end{equation*}
if and only if, for any subset $E_x$ of $E_N$, the following condition $(\ast)(E_x)$ holds.\\
$(\ast)(E_x)$ :
$\projp({\Gamma},E_C\cup E_N)$ is transverse to $M\left((h^{\prime}_i)_{i \in J_E(E_C)}, (h_{i,\eta,w_i})_{i \in J_E(E_N)}\right)$ along \begin{multline*} \projp({\Gamma},E_C)^{-1}\left(\prod_{e \in E_C}\left(C^{\prime}_2(\rats) \cup N_{\left[0,5\right]}\right)\right)
\cap  \projp({\Gamma},E_x)^{-1}\left(\prod_{e \in E_x}N_{\left[7-2\eta,9\right]}\right)\\
\cap  \projp({\Gamma},E_N \setminus E_x)^{-1}\left(\prod_{e \in E_N \setminus E_x}N_{\left[5,7-2\eta\right]}\right).\end{multline*}

Let $\mathring{D}_u$ denote the open disk of $\RR^2$ centered at $0$ of radius $u$.
Our former hypotheses guarantee transversality of $\projp({\Gamma},E_C \cup E_x)$ as soon as the $\norm{w_i}$ are smaller than $u$ for $i \in E_x$. So
the condition $(\ast)(E_x)$ is realized as soon as
$(w_i)_{i\in j_E(E_N)}$ is in an open dense subset $\CD(E_x)$ of $\mathring{D}_u^{j_E(E_N)}$. Thus, we have the desired transversality 
when $(w_i)_{i\in j_E(E_N)}$ is in the intersection of the open dense subsets $\CD(E_x)$ over the subsets $E_x$ of $E_N$.
\eopwobp

\section{More on forms dual to transverse propagating chains}
\label{secexistformtransv}

Though Lemma~\ref{lemexistformtransv} is not surprising, we prove it and refine it in this section. We use its refinement in  Chapter~\ref{chappropzinvffunc}.
Recall the notation of Definition~\ref{defformdual}, and let $\drad{\varepsilon}$ \index[N]{Deps@$\drad{\varepsilon}$ disk of radius $\varepsilon$} ( resp. $\dorad{\varepsilon}$) denote the closed (resp. open) disk of $\CC$ centered at $0$ with radius $\varepsilon$.

\begin{lemma}
\label{lemepsneigprop} Recall that our configuration spaces are equipped with Riemannian metrics.
Let $(\crats,\tau)$ be an asymptotic rational homology $\RR^3$.
Let $\Link \colon \source \to \crats$ be a link in $\crats$.  Let $n \in \NN$, and let $(\propP(i))_{i \in \underline{3n}}$ be a family of propagating chains of $(C_2(\rats),\tau)$ in general $3n$ position with respect to $\Link$.
For any $\varepsilon >0$, there exists $\eta >0$ such that for any $i \in \underline{3n}$,
for any $\Gamma \in \Davis^e_{\underline{3n}}(\source)$, and for any $e \in E(\Gamma)$ with associated restriction map $p_e \colon C(\rats,\Link;\Gamma) \to C_{2}(\rats)$, we have
\begin{equation*}p_e^{-1}\left(N_{\eta}\Bigl(\propP\bigl(j_E(e)\bigr)\Bigr)\right) \subset N_{\varepsilon}\left(p_e^{-1}\Bigl(\propP\bigl(j_E(e)\bigr)\Bigr)\right).\end{equation*}
\end{lemma}
\bp Of course, it is enough to prove the lemma for a fixed $(\Gamma,e)$. 
Set $i=j_E(e)$. The compact $p_e\left( C(\rats,\Link;\Gamma) \setminus N_{\varepsilon}(p_e^{-1}(\propP(i))) \right)$ does not meet $\propP(i)$. So there exists $\eta>0$ such that this compact does not meet 
$N_{\eta}(\propP(i))$ either. This implies $p_e^{-1}(N_{\eta}(\propP(i))) \subset N_{\varepsilon}(p_e^{-1}(\propP(i)))$.
\eop

\begin{lemma}
\label{lemepsneigint}
Let $\Gamma \in \Davis^e_{\underline{3n}}(\source)$.
Assume that the hypotheses of Lemma~\ref{lemepsneigprop} are satisfied.
Then the intersection in $C(\rats,\Link;\Gamma)$ over the edges $e$ of $E(\Gamma)$ of the codimension $2$ rational chains $p_e^{-1}(\propP(j_E(e)))$ is a finite set $I_S(\Gamma,(\propP(i))_{i \in \underline{3n}})$.
Furthermore,
for any $\varepsilon >0$, there exists $\eta >0$ such that we have
\begin{equation*}\bigcap_{e \in E(\Gamma)}p_e^{-1}\left(N_{\eta}\Bigl(\propP\bigl(j_E(e)\bigr)\Bigr)\right) \subset N_{\varepsilon}\left(I_S\Bigl(\Gamma,\bigl(\propP(i)\bigr)_{i \in \underline{3n}}\Bigr)\right)\end{equation*} for any $\Gamma \in \Davis^e_{\underline{3n}}(\source)$.
So,
for any family $(\omega(i))_{i \in \underline{3n}}$ of propagating forms of $(C_2(\rats),\tau)$ $\eta$-dual to the $\propP(i)$, the form $\bigwedge_{e \in E(\Gamma)}p_e^{\ast}(\omega(j_E(e)))$ is supported in
$N_{\varepsilon}(I_S(\Gamma,(\propP(i))_{i \in \underline{3n}}))$.

Moreover, if $N_{\varepsilon}(I_S(\Gamma,(\propP(i))))$ is a disjoint union over the points $\pointx$ of $I_S(\Gamma,(\propP(i)))$ of the $N_{\varepsilon}(\pointx)$, then the integral over $N_{\varepsilon}(\pointx)$ of $\bigwedge_{e \in E(\Gamma)}p_e^{\ast}(\omega(j_E(e)))$ is the rational intersection number of the rational chains $p_e^{-1}(\propP(j_E(e)))$ at $\pointx$. When all the $\propP(j_E(e))$ are embedded manifolds with coefficient $1$ near $p_e(\pointx)$, this intersection number is the sign of $\pointx$ with respect to an orientation $o(\Gamma)$ of $\Gamma$.
\end{lemma}
\bp Again, it suffices to prove the lemma for a fixed $\Gamma \in \Davis^e_{\underline{3n}}(\source)$.
We refer to the description of the image under $\projp(\Gamma)=\prod_{e \in E(\Gamma)}p_e$ of an intersection point $\pointx$ after  Definition~\ref{defgenthreenpos}.

Fix such an $\pointx$. For each edge $e$, $p_e(\pointx)$ sits inside a nonsingular open $4$-dimensional smooth ball $\delta_e$ of a smooth piece $\Delta_{j_E(e),k}$ of $\propP(j_E(e))$.
Consider a tubular neighborhood $N_u(\delta_e)$ whose fibers are disks $\drad{\theta}$ orthogonal to $\delta_e$ of radius $\theta$. The bundle $N_u(\delta_e)$ is isomorphic to $\delta_e \times \drad{\theta}$, with respect to a trivialization of $N_u(\delta_e)$. Another trivialization would compose the diffeomorphism from $N_u(\delta_e)$  to $\delta_e \times \drad{\theta}$ by a map $(x,v) \mapsto (v,\phi(x)(v))$ for some $\phi \colon \delta_e \to SO(2)$.

The projection $p_e(\pointx)$ may sit simultaneously in different nonsingular $4$-dimensional smooth parts $\Delta_{j_E(e),k}$ of $\propP(j_E(e))$. Let $K(e,x)$ be the finite set of components $\Delta_{j_E(e),k}$ of $\propP(j_E(e))$
such that $p_e(\pointx)\in \Delta_{j_E(e),k}$. We first focus on one element of $K(e,x)$ 
for each $e$, and next take the sum over all the choices in $\prod_{e \in E(\Gamma)} K(e,x)$ multiplied by the products of the coefficients of the elements of $K(e,x)$ in the rational chains $\propP(j_E(e))$. Similarly, our forms $\eta$-dual to the $\propP(i)$ are thought of and constructed as linear combinations of forms $\eta$-dual to the elements of $K(e,x)$.

Without loss of generality, assume that $\varepsilon$ is small enough so that we have
\begin{equation*}p_e\bigl(N_{\varepsilon}(\pointx)\bigr) \subset \delta_e \times \drad{\theta}\end{equation*}
for any edge $e$ of $\Gamma$ and for any $\delta_e=\delta_{e,k}$ associated to an element $\Delta_{j_E(e),k}$ of $K(e,x)$,
and so that $p_e(N_{\varepsilon}(\pointx))$ does not meet the components $\Delta_{j_E(e),k}$ of $\propP(j_E(e))$ that are not in $K(e,x)$. Reduce $\varepsilon$ and choose $\eta<\theta$ small enough so that $p_e(N_{\varepsilon}(\pointx))$ does not meet the neighborhoods $N_{\eta}(\Delta_{j_E(e),k})$ of these components, either.

Let $p_{\drad{\theta}} \colon \delta_e \times \drad{\theta} \to \drad{\theta}$ denote the natural projection.
Let $\omega_{\eta}$ be a volume-one form supported on $\dorad{\eta}$. Forms $\eta$-dual to $\propP(j_E(e))$ can be constructed by patching forms $(p_{\drad{\theta}})^{\ast}(\omega_{\eta})$ (multiplied by the coefficients of the $\Delta_{j_E(e),k}$) together, as in Lemma~\ref{lemconstrformalphadual}.
Conversely, for any form $\omega(j_E(e))$ $\eta$-dual to a piece $\Delta_{j_E(e),k}$ of $\propP(j_E(e))$ that contains $\delta_e$, there exists a one-form $\alpha_e$ on $\delta_e \times \drad{\theta}$, such that 
$\omega(j_E(e))=p_{\drad{\theta}}^{\ast}(\omega_{\eta}) + d\alpha_e$
on $\delta_e \times \drad{\eta}$. Then we have 
\begin{equation*}\int_{\{x \in \delta_e\}\times \partial \drad{\eta}}\alpha_e =\int_{\{x \in \delta_e\}\times \drad{\eta}} \omega(j_E(e))-p_{\drad{\theta}}^{\ast}(\omega_{\eta})=0.\end{equation*} So $\alpha_e$ is exact on $\delta_e \times (\drad{\theta}\setminus \drad{\eta})$, and $\alpha_e$ can and will be assumed to be supported on $\delta_e \times \drad{\eta}$.

In the neighborhood $N_{\varepsilon}(\pointx)$ of $\pointx$, $\prod_{e \in E(\Gamma)}p_{\drad{\theta}} \circ p_e$ is a local diffeomorphism around $\pointx$. Without loss of generality, assume that $\eta$ and $\varepsilon$ are small enough so that \begin{equation*}\Pi_p=\prod_{e \in E(\Gamma)}p_{\drad{\theta}} \circ p_e \colon N_{\varepsilon}(\pointx) \to \drad{\theta}^{E(\Gamma)}\end{equation*} 
restricts to a diffeomorphism from $\Pi_p^{-1}\left(\drad{2\eta}^{E(\Gamma)}\right)$ 
 to $\drad{2\eta}^{E(\Gamma)}$, for each $\pointx$ (and for each choice in $\prod_{e \in E(\Gamma)} K(e,x)$).
If $\prod_{e \in E(\Gamma)} K(e,x)$ has one element, and if the coefficient of the element of $K(e,x)$ in $\propP(j_E(e))$ is $1$ for any edge $e$, then we have
\begin{equation*}\int_{N_{\varepsilon}(\pointx)}\bigwedge_{e \in E(\Gamma)}p_e^{\ast}\Bigl(\omega\bigl(j_E(e)\bigr)\Bigr)=\int_{N_{\varepsilon}(\pointx)}\bigwedge_{e \in E(\Gamma)}p_e^{\ast}(p_{\drad{\theta}}^{\ast}\bigl(\omega_{\eta})\bigr).\end{equation*}
Indeed, changing one $\omega(j_E(e))$ to $(p_{\drad{\theta}})^{\ast}(\omega_{\eta})$ amounts to add the integral obtained by replacing $\omega(j_E(e))$ by $d\alpha_e$. Since all the forms are closed, this latter integral is the integral over $\Pi_p^{-1}(\partial (\drad{2\eta}^{E(\Gamma)}))$ of the form obtained by replacing $d\alpha_e$ by $\alpha_e$, which is zero since the whole form is supported in $\Pi_p^{-1}(\drad{\eta}^{E(\Gamma)})$.
Therefore, the integral is the sign of the intersection point $\pointx$ with respect to the given orientation and coorientations.

The open neighborhoods $N_{\varepsilon}(\pointx)$ may be assumed to be disjoint from each other for distinct $\pointx$.
Consequently, since $C(\rats,\Link;\Gamma)$ is compact, the set of intersection points $\pointx$ is finite.
Consider the complement $C^c(\eta_0)$ in $C(\rats,\Link;\Gamma)$ of the union over the intersection points $\pointx$ of the $N_{\varepsilon}(\pointx)$. 
Since $p_{e_1}^{-1}\bigl(\propP(j_E(e_1))\bigr)$ does not meet $\bigcap_{e \in E(\Gamma) \setminus \{e_1\}}p_{e}^{-1}\bigl(\propP(j_E(e))\bigr)$ in $C^c(\eta_0)$, there is an $\varepsilon_1 >0$ such that $\overline{N_{\varepsilon_1}}\bigl(p_{e_1}^{-1}\bigl(\propP(j_E(e_1))\bigr)\bigr)$ does not meet $\bigcap_{e \in E(\Gamma) \setminus \{e_1\}}p_{e}^{-1}\bigl(\propP(j_E(e))\bigr)$ either in $C^c(\eta_0)$. Iterating, we find $\varepsilon_2 >0$ such that
\begin{equation*}C^c(\eta_0) \cap \bigcap_{e \in E(\Gamma)}N_{\varepsilon_2}\biggl(p_{e}^{-1}\Bigl(\propP\bigl(j_E(e)\bigr)\Bigr)\biggr) =\emptyset.\end{equation*}
According to Lemma~\ref{lemepsneigprop}, $\eta$ can be reduced so that 
$p_e^{-1}\bigl(N_{\eta}\bigl(\propP(i)\bigr)\bigr) \subset N_{\varepsilon_2}\bigl(p_e^{-1}\bigl(\propP(i)\bigr)\bigr)$ for any $i$.
Then $\bigwedge_{e \in E(\Gamma)}p_e^{\ast}(\omega(j_E(e)))$ is supported where we want it to be.
\eop

Lemma~\ref{lemexistformtransv} follows.
\eopwobp

Theorem~\ref{thmZinvrat} is now proved.\eopwobp

\section{A discrete definition of the anomaly \texorpdfstring{$\ansothree$}{beta}}

In this section, we give a discrete definition of the anomaly $\ansothree$ and mention a few recent results of K\'evin Corbineau on $\ansothree_3$.
\begin{lemma}
\label{lemansothreeadmis} Let $n \in \NN$. Let $\Gamma \in \Davis^c_n$. Recall the compactification $\ccompuptd{\finsetv(\Gamma)}{\RR^3}$ of $\cinjuptd{\finsetv(\Gamma)}{\RR^3}$ from Theorem~\ref{thmcompfacanom}.
For any edge $e=j_E^{-1}(i)$ of $\Gamma$, we have a canonical projection \begin{equation*}p_e\colon B^3 \times \ccompuptd{\finsetv(\Gamma)}{\RR^3} \to  B^3 \times S^2.\end{equation*}
 Let $i \in \underline{3n}$. When $\Gamma$ is fixed, set $p_i=p_{j_E^{-1}(i)}$.
For any $\veca_i \in S^2$, define the following cooriented chains of $B^3 \times \ccompuptd{\finsetv(\Gamma)}{\RR^3}$:
\begin{equation*}A(\Gamma,i,\veca_i)=p_i^{-1}\bigl(B^3 \times \{\veca_i\}\bigr),\end{equation*}
\begin{equation*}B(\Gamma,i,\veca_i)=p_i^{-1}\Bigl(\Bigl\{\cup_{m\in B^3}\bigl(m,\rho(m)(\veca_i)\bigr)\Bigr\}\Bigr),\end{equation*} and\\
\begin{equation*}H(\Gamma,i,\veca_i)=p_i^{-1}\bigl(G(\veca_i)\bigr),\end{equation*} 
where $\rho$ is introduced in Definition~\ref{defrho}, and the chain $G(\veca_i)$ of $B^3 \times S^2$ is introduced in Lemma~\ref{lemchainanobeta}.
The codimension of $A(\Gamma,i,\veca_i)$ and $B(\Gamma,i,\veca_i)$ is
$2$, while the codimension of
$H(\Gamma,i,\veca_i)$ is $1$.
An element $(a_1,\dots ,a_{3n})$ of $\left(S^2\right)^{3n}$ is \emph{$\ansothree_n$-admissible} if for $h \in \underline{3n}$ and for any $\Gamma \in \Davis^c_n$, the intersection of the $A(\Gamma,i,\veca_i)$ for $i \in \underline{h-1}$, the $B(\Gamma,i,\veca_i)$ for $i \in \underline{3n} \setminus \underline{h} $, and $H(\Gamma,h,\veca_h)$ is transverse.
Then the sets of elements of $\left(S^2\right)^{3n}$ that are $\ansothree_n$-admissible is an open dense subset of $\left(S^2\right)^{3n}$.
\end{lemma}
\bp The principle of the proof is the same as the proof of Proposition~\ref{propratalpha}. See also Section~\ref{secexistransv}. This lemma is proved in detail in \cite{Corbineau}.
\eop

\begin{proposition} Recall the orientation of $\cinjuptd{\finsetv(\Gamma)}{\RR^3}$ of Lemma~\ref{lemoriface}.
For any $\ansothree_n$-admissible element $(a_1,\dots ,a_{3n})$ of $\left(S^2\right)^{3n}$, we have
\begin{equation*}\ansothree_n=\sum_{h=1}^{3n}\sum_{\Gamma \in \Davis^c_n}\frac{1}{(3n)!2^{3n}}
I(\Gamma,h)\left[\Gamma\right],\end{equation*}
with \begin{equation*}I(\Gamma,h)\left[\Gamma\right]=\Bigl\langle \cap_{i=1}^{h-1}A(\Gamma,i,\veca_i),H(\Gamma,h,\veca_h),\cap_{i=h+1}^{3n}B(\Gamma,i,\veca_i)   \Bigr\rangle_{\!\!B^3 \times \ccompuptd{\finsetv(\Gamma)}{\RR^3}\,}\left[\Gamma\right].\end{equation*}
\end{proposition}
\bp For $\veca \in S^2$ and $t \in \left[0,1\right]$, define the following chain 
\begin{multline*}G(\veca,t)=\\ \left[0,t\right] \times B^3 \times \{\veca\} + \bigl(\{t\} \times G(\veca) \bigr) + \Bigl\{\bigl(u,m,\rho(m)(\veca)\bigr) \suchthat  u \in \left[t,1\right], m \in B^3\Bigr\}\end{multline*} of $\left[0,1\right] \times B^3 \times S^2$.
Let $(t_i)_{i \in \underline{3n}}$ be a strictly decreasing sequence of $\left]0,1\right[$.
Let $\Gamma \in \Davis^c_n$. For $i \in \underline{3n}$, let $p_i$ also denote the canonical projection associated to  $e=j_E^{-1}(i)$ from $\left[0,1\right] \times B^3 \times \ccompuptd{\finsetv(\Gamma)}{\RR^3}$
to $\left[0,1\right] \times B^3 \times S^2$.

If $(a_1,\dots ,a_{3n})$ is $\ansothree_n$-admissible, then for any $\Gamma \in \Davis^c_n$,
the intersection of the $p_i^{-1}(G(\veca_i,t_i))$ is transverse and equal to
\begin{equation*}\sqcup_{h=1}^{3n} \{t_h\} \times \biggl(\cap_{i=1}^{h-1}A(\Gamma,i,\veca_i) \cap H(\Gamma,h,\veca_h) \cap \Bigl(\cap_{i=h+1}^{3n}B(\Gamma,i,\veca_i)\Bigr)\biggr).\end{equation*}
Indeed, it is clear that the intersection may be expressed as above at the times $t \in \{t_h \suchthat  h\in \underline{3n}\}$. Since $(a_1,\dots ,a_{3n})$ is $\ansothree_n$-admissible, this intersection is transverse at these times. So it does not intersect the boundaries of the $H(\Gamma,h,\veca_h)$. Therefore, there is no intersection in $\left(\left[0,1\right] \setminus \{t_h \suchthat  h\in \underline{3n}\} \right) \times B^3 \times \ccompuptd{\finsetv(\Gamma)}{\RR^3}$.

Then for any $\alpha >0$, there exist closed $2$-forms $\tilde{\omega}(i)$ on $\left[0,1\right] \times \RR^3 \times S^2$, as in Proposition~\ref{propdeftwoanom}, applied when $\tau_0=\taust$ and $\tau_1=\tau_0 \circ \psi_{\RR}(\rho)$ on $\ST B^3$, such that $\tilde{\omega}(i)$ is
$\alpha$-dual to $G(\veca_i,t_i)$, for any $i$.

Theorem~\ref{thmpone} yields $p_1(\tau_0 \circ \psi_{\RR}(\rho))- p_1(\tau_0)=2\deg(\rho)=4$.
Therefore, Proposition~\ref{propdeftwoanom} implies
\begin{equation*}\ansothree_n=\sum_{\Gamma \in \Davis^c_n}\frac{1}{(3n)!2^{3n}}\int_{\left[0,1\right] \times \cinjuptd{\finsetv(\Gamma)}{T{\ballb}_{\rats}}}\bigwedge_{e \in E(\Gamma)}p_e^{\ast}\Bigl(\tilde{\omega}\bigl(j_E(e)\bigr)\Bigr)\left[\Gamma\right],\end{equation*}
where $\bigwedge_{e \in E(\Gamma)}p_e^{\ast}\bigl(\tilde{\omega}(j_E(e))\bigr)=\bigwedge_{i=1}^{3n} p_i^{\ast}\bigl(\tilde{\omega}(i)\bigr)$.

As in Section~\ref{secexistformtransv}, for $\alpha$ small enough, $\int_{\left[0,1\right] \times \cinjuptd{\finsetv(\Gamma)}{TB^3}}\bigwedge_{i=1}^{3n} p_i^{\ast}\bigl(\tilde{\omega}(i)\bigr)$ is the algebraic intersection of the $p_i^{-1}\bigl(G(\veca_i,t_i)\bigr)$. 

For the signs, note that the coorientation of $\{t_h\} \times G(\veca_h)$ in $\left[0,1\right] \times B^3 \times S^2$ is represented by the orientation of $\left[0,1\right]$, followed by the coorientation of $G(\veca_h)$ in $B^3 \times S^2$.
\eop

In his Ph. D. thesis \cite[Th\'eor\`eme 2.15]{Corbineau}, K\'evin Corbineau obtained the following simplified expression for $\ansothree_3$.

\begin{theorem}
For $j \in \underline{n}$, set
\begin{equation*}H_h(\Gamma,j,\veca_j)=p_j^{-1}\bigl(G_h(\veca_j)\bigr),\end{equation*} 
where the chain $G_h(\veca_j)$ of $B^3 \times S^2$ is introduced in Lemma~\ref{lemchainanobeta}.
Let $\Davis^c_3(T)$ be the set of numbered graphs in $\Davis^c_3$ isomorphic to \begin{equation*}\graphkev. \end{equation*}
For any element $\ansothree_3$-admissible $(a_1,\dots ,a_{9})$ of $\left(S^2\right)^{9}$, we have
\begin{equation*}\ansothree_3=\sum_{j=2}^{8}\sum_{\Gamma \in \Davis^c_3(T)}\frac{1}{(9)!2^{9}}
I_h(\Gamma,j)\left[\Gamma\right],\end{equation*}
with \begin{equation*}I_h(\Gamma,j)\left[\Gamma\right]=\Bigl\langle \cap_{i=1}^{j-1}A(\Gamma,i,\veca_i),H_h(\Gamma,j,\veca_j),\cap_{i=j+1}^{9}B(\Gamma,i,\veca_i)   \Bigr\rangle_{\!\!B^3 \times \ccompuptd{\finsetv(\Gamma)}{\RR^3}\;}\left[\Gamma\right].\end{equation*}

 \end{theorem}

The Ph. D. thesis of K\'evin Corbineau also contains an algorithm to compute $\ansothree_3$.

\part{Functoriality}

Recall that $\drad{1}$ denotes the closed disk of $\CC$ centered at $0$ with radius $1$.
In this book, a \emph{rational homology cylinder}\index[T]{rational homology!cylinder} (or \emph{$\QQ$-cylinder}) is a compact oriented $3$-manifold, whose boundary neighborhood is identified with a boundary neighborhood $N\bigl(\partial (\drad{1} \times [0,1])\bigr)$ of $\drad{1} \times \left[0,1\right]$, and which has the same rational homology as a point.

Roughly speaking, \emph{$q$-tangles} are parallelized cobordisms between limit planar configurations of points up to dilation and translation in rational homology cylinders.
We describe the category of $q$-tangles and its structures precisely in Section~\ref{secintqtangle}.
Framed links in rational homology spheres are particular $q$-tangles. They are cobordisms between empty configurations.

In this third part of the book, we define a functorial extension to $q$-tangles of the invariant $\Zinvufrf$ of framed links in $\QQ$-spheres defined in Section~\ref{secdefdeux}, and we prove that it has a lot of useful properties.
These properties are listed in Theorem~\ref{thmmainfunc}. They ensure that $\Zinvufrf$ is a functor, which behaves naturally with respect to other structures of the category of $q$-tangles, such as cabling or duplication. They allow one to reduce the computation of $\Zinvufrf$ for links to its computation for elementary pieces of the links.

Section~\ref{sectanginj} introduces particular $q$-tangles, for which the involved planar configurations are injective. Section~\ref{sectangwords} introduces other particular $q$-tangles, for which the involved planar configurations are corners of the Stasheff polyhedra of Example~\ref{exaStasheff}.
In Chapter~\ref{chapfirstintfunc}, we define $\Zinvufrf$ for these particular $q$-tangles without proofs. We also state a functoriality result, a monoidality result, and a duplication property, under simple hypotheses, to introduce the involved structures and motivate their introduction. These results are just particular cases of Theorem~\ref{thmmainfunc}.

In Chapter~\ref{chapzinvtang}, we state our general Theorem~\ref{thmmainfunc}, and we describe our strategy towards a consistent definition of $\Zinvufrf$ for general $q$-tangles in Section~\ref{secdefzinvf}.
Our proofs involve convergence results, which rely on intricate compactifications of configuration spaces described in Chapter~\ref{chapconszinvf}. In Chapter~\ref{chapzinvfbraid}, we study
$\Zinvufrf$ as a holonomy for the $q$-tangles that are paths in spaces of planar configurations. In Chapter~\ref{chapdiscrex},
we introduce discretizable versions of $\Zinvufrf$. We use these discretizable versions in the proofs of some important properties of $\Zinvufrf$ given in Chapter~\ref{chappropzinvffunc}.
The consistency of our strategy for the definition of $\Zinvufrf$ is shown in Chapters~\ref{chapconszinvf}, \ref{chapzinvfbraid}, and \ref{chapdiscrex}. The proof of Theorem~\ref{thmmainfunc} will be finished in Chapter~\ref{chappropzinvffunc}. 

This functoriality part contains a generalization of results of Sylvain Poirier \cite{poirierv2}, who constructed the functor $\Zinvufrf$ and proved Theorem~\ref{thmmainfunc} for combinatorial $q$-tangles of $\RR^3$. We recall his results in Section~\ref{secgoodmonfonc}.

\chapter{A first introduction to the functor \texorpdfstring{$\Zinvufrf$}{Z}}
\label{chapfirstintfunc}

In Section~\ref{sectanginj}, we extend the definition of the invariant $\Zinvuf$ of Theorem~\ref{thmfstconst} to long tangle representatives as in Figure~\ref{figLTRbis}. Then we define the framed version $\Zinvufrf$ of $\Zinvuf$ and state that it is multiplicative under the allowed vertical compositions in Section~\ref{secframtanginj}.

In Section~\ref{sectangwords}, we state that $\Zinvufrf$ reaches a limit with nice cabling properties when some vertical infinite strands of the long tangle representatives approach each other. We thus define the restriction of $\Zinvufrf$ to \emph{combinatorial $q$-tangles}, which are parallelized cobordisms between limit configurations on the real line in rational homology cylinders.
This definition is due to Sylvain Poirier \cite{poirierv2} when the involved rational homology cylinder is the standard one $\drad{1} \times\left[0,1\right]$.
In Section~\ref{secgoodmonfonc}, we list sufficiently many properties of the Poirier restriction of $\Zinvufrf$ to characterize the restriction of $\Zinvufrf$ to combinatorial $q$-tangles, in terms of the anomaly $\alpha$.

\section{Extension of \texorpdfstring{$\Zinvuf$}{Z} to long tangles}
\label{sectanginj}
View $\RR^3$ as $\CC \times \RR$, where $\CC$ is horizontal and $\RR$ is vertical, oriented upwards.
For a rational homology cylinder $\hcylc$, $\crats(\hcylc)$ \index[N]{Rmanif@$3$-manifolds!RH@$\crats(\hcylc)$} denotes the asymptotically standard $\QQ$-homology $\RR^3$ obtained by replacing the standard cylinder $\hcylc_0=\drad{1} \times \left[0,1\right]$ in $\RR^3$ by $\hcylc$.

\begin{definition}
\label{defLTR}
A \emph{long tangle representative}\index[T]{long tangle representative} (or LTR\index[T]{LTR} for short) in $\crats(\hcylc)$ is an embedding $\tanghcyll \colon \source \hookrightarrow \crats(\hcylc)$ of a one-manifold  $\source$, as in Figure~\ref{figLTRbis}, such that
\begin{itemize}
\item $\tanghcyll(\source)$ intersects the closure $\check{\hcylc}_0^c$ of the complement of $\hcylc_0$ in $\RR^3$ as
\begin{equation*}\tanghcyll(\source) \cap \check{\hcylc}_0^c= \bigl(c^-(\finsetb^-) \times \left]-\infty,0\right]\bigr) \cup \bigl(c^+(\finsetb^+) \times \left[1,\infty\right[\bigr)\end{equation*}
for two finite sets $\finsetb^-$ and $\finsetb^+$ and two injective maps $c^- \colon \finsetb^- \hookrightarrow \Int(\drad{1})$ and
$c^+ \colon \finsetb^+ \hookrightarrow \Int(\drad{1})$, which are called the \emph{bottom configuration} and the \emph{top configuration} of $\tanghcyll$, respectively, and 
\item $\tanghcyll(\source) \cap \hcylc$ is a compact one-manifold whose unoriented boundary is $\bigl(c^-(\finsetb^-) \times \{0\}\bigr) \cup\bigl (c^+(\finsetb^+) \times \{1\}\bigr)$.
\end{itemize}
\end{definition}

For a $\underline{3n}$-numbered degree $n$ Jacobi diagram with support $\sourcetl$ without looped edges, let $\check{C}\bigl(\crats(\hcylc),\tanghcyll;\Gamma\bigr)$ be its configuration space defined as in Section~\ref{secdefconfspace}.\footnote{The only differences are that $\source$ is not necessarily a disjoint union of circles and that $\check{C}\bigl(\crats(\hcylc),\tanghcyll;\Gamma\bigr)$ was denoted by $\check{C}\bigl(\rats(\hcylc),\tanghcyll;\Gamma\bigr)$.}  The univalent vertices on a \emph{strand}, which is the image under $\tanghcyll$ of an open connected component of $\sourcetl$ (diffeomorphic to $\RR$), move along this whole long component, as in Figure~\ref{figLTRJac}.

\bfig
\centering
\begin{tikzpicture}
\draw [white,very thick] (8.4,3.3) -- (7.3,1.7);
\draw (.4,2)node[left]{$\Gamma =$};
\draw (.6,1.5) node[right]{\cltetra}; 
\draw [yellow!80!black,thick,->] (2.2,.5) -- (2.2,3.5);
\draw [blue,thick,->] (2.6,.5) -- (2.6,3.5);
 \begin{scope}[yshift=2 cm]
\draw [red,thick,->] (1.2,1.6) arc (90:45:.6) (1.2,1.6) arc (90:405:.6);
\draw (1.2,1.6) -- (1.2,1.2) -- (1.8,1) (1.2,1.2) -- (.6,1) (1.2,.4)  .. controls (1.2,.7) and (1.6,.8) .. (2.6,.8);
\fill (1.2,1.6) circle (1.5pt) (1.2,1.2) circle (1.5pt) (1.8,1) circle (1.5pt) (.6,1) circle (1.5pt)
(1.2,.4) circle (1.5pt) (2.6,.8) circle (1.5pt);
\end{scope}
 \begin{scope}[xshift=1.2 cm]
 \draw (5.5,3.7) -- (6.2,2.3);
\draw (4.4,2)node[left]{$c(V(\Gamma)) =$};
\draw [thin,->]  (8.6,-.6) node[right]{\scriptsize $c^-(\finsetb^-) \times \left]-\infty,0\right]$} (8.6,-.6) -- (7.8,-.6);
\draw [thin,->]  (8.6,4.3) node[right]{\scriptsize $c^+(\finsetb^+) \times \left[1,\infty\right[$} (8.6,4.3) -- (8.2,4.3);
\draw [thin,->]  (9.2,1.2) node[right]{\scriptsize ${\hcylc}_0=\drad{1} \times \left[0,1\right]$} (9.2,1.2) -- (8.75,1.2);
\draw [thin,dashed] (5,0)  .. controls (5,.25) and (6.3,.5) .. (7,.5) .. controls (7.7,.5) and (9,.25) ..  (9,0);
\draw [thin,dashed] (5,3.5)  .. controls (5,3.75) and (6.3,4) .. (7,4) .. controls (7.7,4) and (9,3.75) ..  (9,3.5);
\draw [red,thick] (6.2,.95)  .. controls (6.2,1.05) and (6.8,1.2) .. (7.1,1.2) .. controls (7.4,1.2) and (8,1.05) ..  (8,.95);
\draw [blue, thick,->] (8.1,4.2) -- (8.1,4.8) (8.1,2.7) -- (8.1,4.2);
\draw [blue, thick] (8.6,2.1) .. controls (8.6,1.8) and (8.1,1.6) .. (7.8,1.6) .. controls (7.5,1.6) and (6.9,1.7) .. (6.9,1.9);
\draw [draw=white,double=blue,very thick] (7.6,-.6) -- (7.6,1.4) .. controls (7.6,1.6) and (7.1,2.2) .. (7.1,2.4) .. controls (7.1,2.6) and (7.8,2.9) .. (8,2.9)  .. controls (8.2,2.9) and  (8.6,2.4) .. (8.6,2.1);
\draw [draw=white,double=blue,very thick] (6.9,1.9)  .. controls (6.9,2.2) and  (8.1,2.4) .. (8.1,2.7);
\draw [blue,thick] (6.9,1.9)  .. controls (6.9,2.2) and  (8.1,2.4) .. (8.1,2.7);
\draw (8.1,4.6) -- (6.2,.95);
\draw [draw=white,double=yellow!80!black,very thick] (5.6,-.8) -- (5.6,1.4) .. controls (5.6,1.6) and (5.9,1.9) .. (6.1,1.9) .. controls (6.3,1.9) and (6.6,1.6) ..  (6.6,1.4) -- (6.6,-.4);
\draw [yellow!80!black, thick,->] (6.6,-.4) -- (6.6,-.7);
\draw [yellow!80!black, thick] (6.6,-.6) -- (6.6,-.8);
\draw [blue, thick,->] (7.6,-.8) -- (7.6,-.6);
\draw [white,very thick] (8,.95) -- (8.3,2.1) -- (7.1,1.2) (8.3,2.1) -- (7.1,.7);
\draw [draw=white,double=red,very thick] (6.2,.95)  .. controls (6.2,.8) and (6.8,.7) .. (7.1,.7) .. controls (7.4,.7) and (8,.8) ..  (8,.95);
\draw [draw=white,double=black,very thick] (5,0)  .. controls (5,-.25) and (6.3,-.5) .. (7,-.5) .. controls (7.7,-.5) and (9,-.25) ..  (9,0);
\draw (5.3,1.7) -- (6.2,2.3) -- (6.4,3.3) -- (5.5,3.7) -- (5.3,1.7) -- (6.4,3.3);
\draw [draw=white,double=black,very thick] (5,3.5)  .. controls (5,3.25) and (6.3,3) .. (7,3) .. controls (7.7,3) and (9,3.25) ..  (9,3.5);
\draw [thin] (5,0) -- (5,3.5) (9,0) -- (9,3.5);
\draw (8,.95) -- (8.3,2.1) -- (7.1,1.2) (8.3,2.1) -- (7.1,.7);
\fill (5.3,1.7)  circle (1.5pt) (6.2,2.3) circle (1.5pt) (6.4,3.3) circle (1.5pt) (5.5,3.7) circle (1.5pt)
(6.2,.95) circle (1.5pt) (8,.95) circle (1.5pt) (7.1,1.2) circle (1.5pt) 
(7.1,.7) circle (1.5pt) (8.3,2.1) circle (1.5pt)  (8.1,4.6) circle (1.5pt);
\end{scope}
\end{tikzpicture}
\caption{A (black) Jacobi diagram $\Gamma$ on the domain of an LTR $\tanghcyll$ and a configuration $\confc$ of $\check{C}(\crats(\hcylc),\tanghcyll;\Gamma)$}
\label{figLTRJac}

\end{figure}

For any $i \in \underline{3n}$, let $\omega(i)$ be a propagating form of $\bigl(C_2(\rats(\hcylc)),\tau\bigr)$. Let $o(\Gamma)$ be a vertex-orientation of $\Gamma$.
As in Section~\ref{secconfint}, define \begin{equation*}I\left(\hcylc,\tanghcyll,\Gamma,o(\Gamma),\bigl(\omega(i)\bigr)_{i \in \underline{3n}}\right)=\int_{(\check{C}(\crats(\hcylc),\tanghcyll;\Gamma),o(\Gamma))} \bigwedge_{e \in E(\Gamma)}p_e^{\ast}\left(\omega\bigl(j_E(e)\bigr)\right)\index[N]{Integrals over configuration spaces!IRLGammab@$I(\hcylc,\tanghcyll,\Gamma,o(\Gamma),(\omega(i)))$},\end{equation*} 
where $(\check{C}(\crats(\hcylc),\tanghcyll;\Gamma),o(\Gamma))$ denotes the manifold $\check{C}(\crats(\hcylc),\tanghcyll;\Gamma)$, equipped with the orientation induced by the vertex-orientation $o(\Gamma)$ and by the edge-orientation of $\Gamma$, as in Corollary~\ref{cororc}, and
\begin{equation*}I\left(\hcylc,\tanghcyll,\Gamma,\bigl(\omega(i)\bigr)_{i \in \underline{3n}}\right)\left[\Gamma\right]=I\left(\hcylc,\tanghcyll,\Gamma,o(\Gamma),\bigl(\omega(i)\bigr)_{i \in \underline{3n}}\right)\left[\Gamma,o(\Gamma)\right].\end{equation*}

\begin{theorem} \label{thmconvint}
The above integral is convergent.
\end{theorem}

We prove this theorem in Section~\ref{secstrucbyhand}. See Lemma~\ref{lemStokesconftang}. Again, its proof involves appropriate compactifications ${C}(\crats(\hcylc),\tanghcyll;\Gamma)$ of the configuration spaces $\check{C}(\crats(\hcylc),\tanghcyll;\Gamma)$. The compactifications are more complicated in this case. We study them in Chapter~\ref{chapconszinvf}.

As an example, let us compute $I\bigl(\hcylc,\tanghcyll,\Gamma,o(\Gamma),(\omega(i))_{i \in \underline{3n}}\bigr)$ when
\begin{itemize}
 \item $\hcylc=\hcylc_0=\drad{1} \times \left[0,1\right]$,
 \item $\tanghcyll$ is an LTR $\fulltwist$ whose bottom and top configurations coincide
and map $\finsetb^-=\finsetb^+$ to $\{-\frac12,\frac12\}$, 
\item $\tanghcyll(\source)\cap \hcylc_0$ projects to $\RR^2$ as $\fulltwist$,
\item $(\Gamma,o(\Gamma))$ is the vertex-oriented diagram $\,\onechordtwosoror$ whose chord is oriented and numbered, and
\item the propagating forms $\omega(i)$ pull back through $p_{S^2} \colon C_2(S^3) \to S^2$.
\end{itemize}

\begin{lemma} \label{lemfulltwistdegone} We have
\begin{equation*}I\left(\hcylc_0,\fulltwist,\onechordtwosoror,\left(p_{S^2}^{\ast}\left(\omega_{i,S}\right)\right)_{i \in \underline{3}}\right)=I\left(\hcylc_0,\fulltwist,\onechordtwosororrev,\left(p_{S^2}^{\ast}\left(\omega_{i,S}\right)\right)_{i \in \underline{3}}\right)=1\end{equation*} for any arbitrary numbering of the edge of the involved Jacobi diagram and any choice of volume-one forms $\omega_{i,S}$ of $S^2$. 
\end{lemma}
\bp Let us compute \begin{equation*}I\left(\fulltwist,\onechordtwosoror,\left(p_{S^2}^{\ast}\left(\omega_{i,S}\right)\right)_{i \in \underline{3}}\right).\end{equation*}
The configuration space \begin{equation*}\check{C}=\biggl(\check{C}\Bigl(\RR^3,\fulltwist;\Gamma\Bigr),o(\Gamma)\biggr)\end{equation*} is naturally diffeomorphic to $\left]-\infty,\infty\right[ \times \left]-\infty,\infty\right[$, where the first factor parametrizes the height of the vertex on the left strand oriented from bottom to top and the second one parametrizes the height of the vertex on the right strand. 
 
The map $p_{S^2}$ maps $\left]-\infty,0\right]^2$ and $\left[1,\infty\right[^2$ to the vertical circle through the horizontal real direction. Therefore, the integral of $p_{S^2}^{\ast}\left(\omega_{i,S}\right)$ vanishes there, and the integral of $p_{S^2}^{\ast}\left(\omega_{i,S}\right)$ over $\check{C}$ is the integral of $p_{S^2}^{\ast}\left(\omega_{i,S}\right)$ over $\left]-\infty,\infty\right[^2 \setminus \bigl(\left]-\infty,0\right[^2 \cup \left]1,\infty\right[^2\bigr)$ or over $\left[-\infty,\infty\right]^2 \setminus \bigl(\left[-\infty,0\right[^2 \cup \left]1,\infty\right]^2\bigr)$, to which $p_{S^2}^{\ast}\left(\omega_{i,S}\right)$ extends naturally. The boundary of this domain, which is  drawn in Figure~\ref{figTsetminusTCbis}, is mapped to the vertical half circle between the two vertical directions north $\upvec$ and south $(-\upvec)$ through the horizontal east direction $\eastvec$ towards the right. 

\bfig
\centering
\begin{tikzpicture}
\begin{scope}[xshift=9.5cm]
 \draw [out=90,in=-90, dashed, blue,very thick] (.6,1.8) to (0,2.4);
\draw [out=90,in=-90, dashed,blue,very thick] (0,1.8) to (.6,2.4);
\draw [dashed,blue,very thick] (.6,0) -- (.6,1.8);
\draw [dotted,blue,very thick] (0,0) -- (0,.4);
\draw [blue,very thick] (0,.4) -- (0,1.8);
\draw (.3,1.8) -- (.6,1.8) (0,.6) -- (.6,1.8);
\fill (0,1.8) circle (2pt) (.6,1.8) circle (2pt) (0,.6) circle (2pt);
\draw [->] (.25,1.6) node{\scriptsize $\eastvec$}  (.8,1.8) node{\scriptsize $0$} (0,1.8) -- (.3,1.8);
\draw [->] (0,.6) -- (.3,1.2);
\draw [->] (1.2,.6) node{\scriptsize $\upvec$} (1,.2) -- (1,.6);
\draw [very thin] (-.35,.85) node[left]{\scriptsize $\left[-\infty,0\right]$} (-.2,1.8) -- (-.35,1.65) -- (-.35,0);
\end{scope}
\fill[lightgray] (0,0) -- (0,.8) -- (.8,.8) -- (.8,0) -- (0,0) (1.6,1.6) -- (2.4,1.6) -- (2.4,2.4) -- (1.6,2.4) -- (1.6,1.6);
\draw [very thick,->] (0,1.6) -- (0,.8) -- (.8,.8) -- (.8,0) -- (2.4,0) -- (2.4,.8); 
\draw [very thick,->]  (2.4,.8) node[right]{\scriptsize $\partial \left(\left[-\infty,\infty\right]^2 \setminus \left(\left[-\infty,0\right[^2 \cup \left]1,\infty\right]^2\right)\right)$}  (2.4,.8)-- (2.4,1.6) -- (1.6,1.6) -- (1.6,2.4) -- (0,2.4) -- (0,1.6); 
\draw (.2,2.2) node{\scriptsize $\upvec$} (.2,1) node{\scriptsize $\upvec$} (1.4,2.2) node{\scriptsize $\upvec$};
\draw (2.05,.2) node{\scriptsize $-\upvec$} (2.05,1.4) node{\scriptsize $-\upvec$} (1.05,.2) node{\scriptsize $-\upvec$};
\draw (.95,.95) node{\scriptsize $\eastvec$} (1.45,1.45) node{\scriptsize $\eastvec$};
\draw [line width=.12cm] (0,.75) -- (.8,.75);
\end{tikzpicture}
\caption{Images of boundary points of $\left[-\infty,\infty\right]^2 \setminus \left(\left[-\infty,0\right[^2 \cup \left]1,\infty\right]^2\right)$ under $p_{S^2}$. The right-hand side shows the computation for the very thick part $\left[-\infty,0\right] \times \{0\}$.}
\label{figTsetminusTCbis}

\end{figure}

Thus, $p_{S^2}\bigl(\partial{C}(\RR^3,\fulltwist;\Gamma)\bigr)$ is algebraically trivial (as in the beginning of Section~\ref{secstraight}), and the differential degree of $p_{S^2}$ is constant on the set of regular values of $p_{S^2}$, according to Lemma~\ref{lemdeggen}. It can be computed as in Subsection~\ref{sublkGausstwo} at the vector that points towards the reader. It is equal to one. Thus we have
\begin{equation*}I\left(\fulltwist,\onechordtwosoror,\left(p_{S^2}^{\ast}\left(\omega_{i,S}\right)\right)_{i \in \underline{3}}\right)=1\end{equation*} for any arbitrary numbering of the edge of $\Gamma=\onechordtwosoror$ and for any choice of volume-one forms $\omega_{i,S}$ of $S^2$. Similarly, for the opposite orientation of the edge of $\Gamma$, we have \begin{equation*}I\left(\fulltwist,\onechordtwosororrev,\left(p_{S^2}^{\ast}\left(\omega_{i,S}\right)\right)_{i \in \underline{3}}\right)=1.\end{equation*}
\eop

\begin{definition} \label{defparacyl}
 A \emph{parallelization} of $\hcylc$ is a parallelization of $\crats(\hcylc)$ that agrees with the standard parallelization of $\RR^3$ outside $\hcylc$.
A \emph{parallelized rational homology cylinder} $(\hcylc,\tau)$ is a rational homology cylinder
equipped with such a parallelization.
\end{definition}

The following lemma shows other important examples of computations.

\begin{lemma}\label{lemvarithetatang}
Let $K \colon \RR \hookrightarrow \crats(\hcylc)$ be a component of $\tanghcyll$. Let $\tau$ be a parallelization of $\hcylc$ (which is standard near $\partial \hcylc$ by Definition~\ref{defparacyl}).
For any $i\in \underline{3}$, let $\omega(i)$ and $\omega^{\prime}(i)$ 
be propagating forms of $(C_2(\rats(\hcylc)),\tau)$, which restrict to $\partial C_2(\rats(\hcylc))$ as $p_{\tau}^{\ast}(\omega(i)_{S^2})$ and $p_{\tau}^{\ast}(\omega^{\prime}(i)_{S^2})$, respectively.
Let $\eta(i)_{S^2}$ be a one-form on $S^2$ such that $\omega^{\prime}(i)_{S^2}=\omega(i)_{S^2} + d\eta(i)_{S^2}$. Then when $K$ goes from bottom to top or from top to bottom, we have
\begin{equation*}\begin{array}{ll}I\left(\onechordsmallnumt,\bigl(\omega^{\prime}(i)\bigr)_{i \in \underline{3}}\right)-I\left(\onechordsmallnumt,\bigl(\omega(i)\bigr)_{i \in \underline{3}}\right)&=\int_{\ST^+K}p_{\tau}^{\ast}(\eta(k)_{S^2})\\&=\int_{p_{\tau}(\ST^+K)}\eta(k)_{S^2}.\end{array} \end{equation*}
When $K$ goes from bottom to bottom (resp. from top to top), let $S(K)$ be the half-circle from $-\upvec$ to $\upvec$ (resp. from $\upvec$ to $-\upvec$) through the horizontal direction from the initial vertical half-line of $K$ (the first encountered one) to the final one, then we have \begin{equation*}I\left(\onechordsmallnumt,\bigl(\omega^{\prime}(i)\bigr)_{i \in \underline{3}}\right)-I\left(\onechordsmallnumt,\bigl(\omega(i)\bigr)_{i \in \underline{3}}\right)=\int_{p_{\tau}(\ST^+K) \cup S(K)}\eta(k)_{S^2}.\end{equation*}
In particular, $I(\onechordsmallnumt,(\omega(i))_{i \in \underline{3}})$ depends only on $\omega(k)_{S^2}$. It is also denoted by $I\left(\onechordsmallnonnum,\omega(k)_{S^2}\right) $. 
\end{lemma}
\bp
In any case, the configuration space $\check{C}(\crats(\hcylc),\tanghcyll;\onechordsmallnonnum)$ is identified naturally with the interior of a triangle, as in the left part of Figure~\ref{figbasictriangle}.
When $K$ goes from bottom to top, $p_{\tau}$ extends smoothly to the triangle, as a map that sends the horizontal side and the vertical side to $\upvec$. Conclude as in Lemmas~\ref{lemvaritheta} and \ref{lemdefItheta}.
The case in which $K$ goes from top to bottom is similar.

\bfig
\centering
\begin{tikzpicture}
\fill[lightgray,draw=black] (0,0) -- (1,1) -- (0,1) -- (0,0);
\draw (1,1) node[right]{\scriptsize $(\infty,\infty)$};
\draw (0,1) node[left]{\scriptsize $(-\infty,\infty)$};
\draw (0,0) node[left]{\scriptsize $(-\infty,-\infty)$};
\draw (-.4,.5) node[left]{\scriptsize $\check{C}(\crats(\hcylc),\tanghcyll;\onechordsmallnonnum)$};
\draw [very thin, ->] (-.4,.5) -- (.25,.5);
\draw [->] (0,0) -- (.5,.5);
\draw (.5,.5) node[right]{\scriptsize $\ST^+K$};
\begin{scope}[xshift=4.5cm]
 \fill[lightgray,draw=black] (0,0) -- (1,1) -- (.3,1) arc (0:-90:.3) -- (0,0);
\draw (1,1) node[right]{\scriptsize $(1,1)$};
\draw (0,0) node[left]{\scriptsize $(0,0)$};
\draw [->] (0,0) -- (.5,.5);
\draw (.5,.5) node[right]{\scriptsize $\ST^+K$};
\fill[white] (.1,.4) rectangle (.4,.7);
\draw (.25,.55) node{\scriptsize $\overline{T}$};
\end{scope}
\end{tikzpicture}
\caption{Compactifications of configuration spaces for the proof of Lemma~\ref{lemvarithetatang}}
\label{figbasictriangle}

\end{figure}

Let us study the case in which $K$ goes from top to top.
Let $d_1=-\{z_1\} \times \left[1,\infty\right]$ and $d_2=\{z_2\} \times \left[1,\infty\right]$ denote the vertical half-lines of $K$ above $\hcylc$, where $d_1$ is before $d_2$. 
View $K$ as a path composition $d_1\bigl(K\cap \left(\drad{1}\times\left[0,1\right]\right)\bigr)d_2$ and
parametrize 
\begin{equation*}\begin{array}{llll}K=\begin{tikzpicture}
\useasboundingbox (-.6,0) rectangle (1,.5);
\draw (0,.5) -- (0,0) (.4,.5) -- (.4,0);
\draw [dashed] (0,0) arc (-180:0:.2);
\draw [->] (0,.5) -- (0,.25);
\draw [->] (.4,0) -- (.4,.25);
\draw (0,.25) node[left]{\scriptsize $d_1$};
\draw (.4,.25) node[right]{\scriptsize $d_2$};
\end{tikzpicture} \mbox{ by }m \colon &\left]0,1\right[ &\to &K\\
                           & t \in\left]0,1/3\right]& \mapsto &(z_1,1/(3t))\\
 & t \in\left[2/3,1\right[ & \mapsto &(z_2,1/(3(1-t))).
                           \end{array}\end{equation*}
Let $T_0=\{(t_1,t_2)\in \left]0,1\right[^2 \suchthat  t_1 < t_2\}$.
We study the integral of $\omega^{\prime}(k)-\omega(k)=\mbox{d}\eta(k)$ over $\check{C}(\crats(\hcylc),\tanghcyll;\onechordsmalltK)=m^2({T_0})$.
View $T_0$ as the set
\begin{equation*}T=\left\{(t_1,t_2,\alpha)\in \left]0,1\right[^2\times \left]-\frac{\pi}{2},\frac{\pi}{2}\right[  \suchthat  t_1 < t_2,\tan(\alpha)=\frac{1/(3(1-t_2))-1/(3t_1)}{|z_2-z_1|}\right\}\end{equation*}
denoted by $T$.
Note that when $(t_1,t_2) \in \left]0,1/3\right[ \times \left]2/3,1\right[$, and when $(t_1,t_2,\alpha)\in T$, we have
\begin{equation*}p_{\tau}\bigl(m(t_1),m(t_2)\bigr)=\cos \alpha \frac{z_2-z_1}{|z_2-z_1|} + \sin \alpha \upvec \in S(K).\end{equation*}
Let $\overline{T}$ be the closure of $T$ in $\left[0,1\right]^2\times \left[-\frac{\pi}{2},\frac{\pi}{2}\right]$. This closure, drawn in the right part of Figure~\ref{figbasictriangle}, is a smooth blow-up of $\{(t_1,t_2)\in \left[0,1\right]^2 \suchthat  t_1 \leq t_2\}$ at $(0,1)$ (with corners), where $(0,1)$ lifts as $(0,1)\times \left[-\frac{\pi}{2},\frac{\pi}{2}\right]$ in $\overline{T}$. The map $m^2$ extends 
as a smooth map valued in $C_2(\rats(\hcylc))$ on $\overline{T}$. Its composition with $p_{\tau}$ sends the vertical side of $\overline{T}$ to $-\upvec$, its horizontal side to $\upvec$, and the blown-up upper-left corner to $S(K)$.
The integral of $\mbox{d}\eta(k)$ over $\check{C}(\crats(\hcylc),\tanghcyll;\onechordsmalltK)$ is the integral of $\mbox{d}\eta(k)$ over
$m^2\left(\overline{T} \right)$. So it is the integral of $\eta(k)$ over $m^2(\partial \overline{T})$, where $\eta(k)$ can be assumed to be equal to
$p_{\tau}^{\ast}(\eta(k)_{S^2})$, as in Lemma~\ref{lemetactwo}. Furthermore, $p_{\tau}\circ m^2$ restricts to $\partial \overline{T}$ as a degree one map onto $p_{\tau}(\ST^+K)\cup S(K)$. So the stated conclusion follows.
The case in which $K$ goes from bottom to bottom can be treated similarly.
 \eop
 
 \begin{definition}
\label{defIthetalong} Recall the notation of Lemma~\ref{lemvarithetatang}.
For a long component $K \colon \RR \hookrightarrow \crats(\hcylc)$ of a tangle in a parallelized $\QQ$-cylinder $(\hcylc,\tau)$, define \index[N]{Integrals over configuration spaces!Ithetalong@$I_{\theta}(K,\partau)$ for a long knot}
\begin{equation*}I_{\theta}(K,\tau) =2 I\left(\onechordsmallnonnum,\omega_{S^2}\right).\end{equation*}
Recall the definition of $I_{\theta}(K,\partau)$ for a closed component $K$ of $\crats(\hcylc)$ from Lemma~\ref{lemdefItheta}.
 \end{definition}

The factor $2$ in the definition of $I_{\theta}$ for long components in 
Definition~\ref{defIthetalong} 
may seem unnatural. It allows to get homogeneous formulas
in 
Theorem~\ref{thmfstconsttang} below. Theorem~\ref{thmfstconsttang} generalizes Theorem~\ref{thmfstconst} to long tangle representatives.
It will be proved in Section~\ref{secvarzinf}.

\begin{theorem}
\label{thmfstconsttang} Let $\hcylc$ be a rational homology cylinder equipped with a parallelization $\tau$ (standard near $\partial \hcylc$).
Let $\tanghcyll \colon \sourcetl \hookrightarrow \crats(\hcylc)$ be a long tangle representative in $\crats(\hcylc)$.
Let $n \in \NN$.
For any $i\in \underline{3n}$, let $\omega(i)$ be a homogeneous propagating form of $(C_2(\rats(\hcylc)),\tau)$.
Set \begin{equation*}\Zinv_n\Bigl(\hcylc,\tanghcyll,\bigl(\omega(i)\bigr)\Bigr)=\sum_{\Gamma \in \Davis^e_n(\sourcetl)}\coefgambet_{\Gamma}I\Bigl(\hcylc,\tanghcyll,\Gamma,\bigl(\omega(i)\bigr)_{i \in \underline{3n}}\Bigr)\left[\Gamma\right] \in \Aavis_n(\sourcetl)\index[N]{ZZ@$\Zinvuf$ and some variants (see also the summary in the next pages)!ZiomT@$\Zinv_n(\hcylc,\tanghcyll,(\omega(i)))$},\end{equation*} 
where $\coefgambet_{\Gamma}=\frac{(3n-\cardlef{E(\Gamma)})!}{(3n)!2^{\cardlef{E(\Gamma)}}}$, with Definition~\ref{defnumdia} for $\Davis^e_n(\sourcetl)$.
Then $\Zinv_n(\hcylc,\tanghcyll,(\omega(i)))$ is independent of the chosen
$\omega(i)$. It depends only on $(\hcylc,\tanghcyll(\sourcetl) \cap \hcylc)$ up to diffeomorphisms that fix $\partial \hcylc$ (and $\tanghcyll(\sourcetl) \cap \partial \hcylc$), pointwise, on $p_1(\tau)$, and on the $I_{\theta}(K_j,\tau)$ for the components $K_j$, $j \in \underline{k}$, of $\tanghcyll$. We denote it by $\Zinv_n(\hcylc,\tanghcyll,\tau)$.
Set \begin{equation*}\Zinv(\hcylc,\tanghcyll,\tau)=\bigl(\Zinv_n(\hcylc,\tanghcyll,\tau)\bigr)_{n\in \NN} \in \Aavis(\sourcetl), \index[N]{ZZ@$\Zinvuf$ and some variants (see also the summary in the next pages)!ZitauT@$\Zinv(\hcylc,\tanghcyll,\tau)$}\end{equation*} and recall the anomalies $\alpha \in \Assis(S^1;\RR)$ and $\ansothree \in \Aavis(\emptyset;\RR)$ from Sections~\ref{secanomalpha} and \ref{secansothree}.
Then the expression \begin{equation*}\exp\Bigl(-\frac14 p_1(\tau)\ansothree\Bigr)\prod_{j=1}^k\Bigl(\exp\bigl(-I_{\theta}(K_j,\tau)\alpha\bigr)\#_j\Bigr) \Zinv(\hcylc,\tanghcyll,\tau)\end{equation*}
depends only on the boundary-preserving diffeomorphism class of $(\hcylc,\tanghcyll)$.\footnote{Again, the subscript $j$ of $\#_j$ indicates that $\exp(- I_{\theta}(K_j,\tau)\alpha)$ is inserted on the component of $K_j$ of the domain $\sourcetl$ of $\tanghcyll$.} We denote it by $\Zinvuf(\hcylc,\tanghcyll)$.\index[N]{ZZ@$\Zinvuf$ and some variants (see also the summary in the next pages)!ZuTL@$\Zinvuf(\hcylc,\tanghcyll)$}
\end{theorem}

\section{Definition of \texorpdfstring{$\Zinvufrf$}{Z} for framed tangles}
\label{secframtanginj}

Recall Definition~\ref{defparselflk} of knot parallels.

\begin{definition}
 \label{defparalleltang}
A \emph{parallel}\index[T]{parallel!of a long component} $K_{\parallel}$ of an embedding of a long component $K \colon \RR \hookrightarrow \crats(\hcylc)$
parametrized so that $K(\RR) \cap \hcylc =K\left(\left[0,1\right]\right)$ is the image of an embedding $K_{\parallel} \colon \RR \hookrightarrow \crats(\hcylc)$ such that there exists an embedding 
\begin{equation*}k \colon [-1,1] \times \RR \to  \rats(\hcylc) \setminus \bigl(\Link(\source) \setminus K(\RR)\bigr)\end{equation*}
such that $K=k\vert_{\{0\} \times \RR}$,$K_{\parallel} =k\vert_{\{1\} \times \RR}$
and $k(u,t)=K(t) + u\varepsilon(t)(1,0,0)$ for any $(u,t) \in [-1,1]\times \left(\RR \setminus \left]0,1\right[\right)$,
 for a small smooth function $\varepsilon \colon \RR \setminus \left]0,1\right[ \to \RR \setminus \{0\}$ such that $\varepsilon(0)\varepsilon(1)$ is positive for components going from bottom to top or from top to bottom, and negative 
for components going from bottom to bottom or from top to top. (We push in one of the two horizontal real directions in a way consistent with the orientation.)
Parallels are considered up to isotopies that stay within these parallels and up to the exchange of $k\vert_{\{1\} \times \RR}$ and $k\vert_{\{-1\} \times \RR}$.
A component $K$ of an LTR is \emph{framed}\index[T]{framed!tangle} if $K$ is equipped with such a class of parallels, called a \emph{framing} of $K$. An LTR is \emph{framed} if all its components are.
The \indexT{self-linking number} of a circle component $K$ in a framed LTR is the linking number $lk(K,K_{\parallel})$ of $K$ and its parallel $K_{\parallel}$.
\end{definition}

For a long component $K$ equipped with a parallel, we define its \emph{self-linking number} $lk(K,K_{\parallel})$ in Definitions~\ref{defselflkqtangone} and \ref{defselflkqtangtwo} below, and in Definition~\ref{defselflkqtanggen}, which covers the remaining cases (when $K$ goes from bottom to bottom or from top to top and when $K(1)-K(0)$ is not in the direction of the real line).

\begin{definition}
 \label{defselflkqtangone}
When $K$ goes from bottom to top, and when $\varepsilon(0)$ is positive, let $\left[K_{\parallel}(1),(1,1)\right]$ be the straight segment from $K_{\parallel}(1)$ to $(1,1) \in \drad{1} \times \{1\}$. Similarly define $\left[K(1),(-1,1)\right]$, $\left[(1,0),K_{\parallel}(0)\right]$ and $\left[(-1,0),K(0)\right]$, and note that they are pairwise disjoint. Define the topological circle embeddings \begin{equation*}\begin{array}{ll} \hat{K}&=K([0,1]) \cup \left[K(1),(-1,1)\right] \cup \bigl(-\{-1\} \times \left[0,1\right]\bigr) \cup \bigl[(-1,0),K(0)\bigr]\\
\hat{K}_{\parallel}&=K_{\parallel}([0,1]) \cup \left[K_{\parallel}(1),(1,1)\right] \cup \bigl(-\{1\} \times \left[0,1\right]\bigr) \cup \bigl[(1,0),K_{\parallel}(0)\bigr],\end{array}\end{equation*} as in the figure below,\footnote{As we often abusively do, we identify embeddings with their images.} and
set \begin{equation*}lk\bigl(K,K_{\parallel}\bigr)=lk\bigl(\hat{K},\hat{K}_{\parallel}\bigr).\end{equation*}
\begin{center}
\begin{tikzpicture}
\useasboundingbox (-3,-.8) rectangle (3,.8);
\draw [->,very thick] (0,.2) .. controls (.2,.4) and (.1,.5) .. (.1,.7) -- (-1.4,.7) -- (-1.4,-.7) -- (-.5,-.7) .. controls (-.5,-.2) and (-.2,0) .. (0,.2);
\draw [->] (.4,.2) .. controls (.6,.4) and (.5,.5) .. (.5,.7) -- (1.4,.7) -- (1.4,-.7) -- (-.1,-.7) .. controls (-.1,-.2) and (.2,0) .. (.4,.2);
\draw (-1.4,-.7) node[left]{\scriptsize $(-1,0)$} (-1.4,.7) node[left]{\scriptsize $(-1,1)$}
(1.4,-.7) node[right]{\scriptsize $(1,0)$} (1.4,.7) node[right]{\scriptsize $(1,1)$}
(-.85,-.2) node[below]{\scriptsize $K(0)$} (.4,-.2) node[below]{\scriptsize $K_{\parallel}(0)$}
(-.3,.5) node[below]{\scriptsize $\hat{K}$} (.6,.3) node[below]{\scriptsize $\hat{K}_{\parallel}$} (.95,.5) node{\scriptsize ${K}_{\parallel}(1)$};
\end{tikzpicture}
\end{center}
When $\varepsilon(0)$ is negative, the above definition determines $lk\left(K_{\parallel},K\right)$. Then set \begin{equation*}lk\left(K,K_{\parallel}\right)=lk\left(K_{\parallel},K\right).\end{equation*} So $lk\left(K,K_{\parallel}\right)$ is defined when $K$ goes from bottom to top.
When $K$ goes from top to bottom, $(-K)$ goes from bottom to top, and $-K_{\parallel}$ is a parallel of $(-K)$. Then set $lk\bigl(K,K_{\parallel}\bigr)=lk\bigl(-K,-K_{\parallel}\bigr)$.
\end{definition}

\begin{definition}
 \label{defselflkqtangtwo}
When $K$ goes from bottom to bottom or from top to top, and when $K(1)-K(0)$ is equal to $(v,0,0)$ for $v \in \RR \setminus \{0\}$, we define the self-linking number $lk(K,K_{\parallel})$ as follows.
Let us first assume that $v>0$ and that $\varepsilon(0)>0$, as in Figure~\ref{figgammaparthetazero}. In this case, define 
topological circle embeddings
$\hat{K}_{\parallel} = K_{\parallel}\left(\left[0,1\right]\right) \cup \left[K_{\parallel}(1),K_{\parallel}(0)\right]$, where $\left[K_{\parallel}(1),K_{\parallel}(0)\right]$ is the straight segment from $K_{\parallel}(1)$ to $K_{\parallel}(0)$ in $\drad{1} \times \{0\}$ or in $\drad{1} \times \{1\}$, and $\hat{K} = K\left(\left[0,1\right]\right) \cup \gamma\left(\left[0,1\right]\right)$ for an arbitrary path $\gamma$ from $\gamma(0)=K(1)$ to $\gamma(1)=K(0)$ such that 
$\gamma\left(\left]0,1\right[\right) \subset \crats(\hcylc) \setminus \hcylc$ as in Figure~\ref{figgammaparthetazero}, and set \begin{equation*}lk\bigl(K,K_{\parallel}\bigr)=lk\bigl(\hat{K},\hat{K}_{\parallel}\bigr).\end{equation*}
\bfig 
\centering
\begin{tikzpicture}
\useasboundingbox (-.1,-.4) rectangle (7.6,1.7);
\draw [<-, dash pattern=on 2pt off 1pt] (1.7,.25) arc (0:180:.7);
\draw [ very thick, dash pattern=on 2pt off 1pt] (2,.25) arc (0:180:1);
\draw [->] (1,0) -- (.3,0) -- (.3,.25) (1.7,.25) -- (1.7,0) -- (1,0);
\draw [->] (1,-.25) node[below]{\scriptsize $\gamma$} (1,-.1) node[above]{\scriptsize $\hat{K}_{\parallel}$} (0.05,0) node[left] {\scriptsize $K(0)$} (1.95,0) node[right] {\scriptsize $K(1)$} (1.95,.3) node[right] {\scriptsize $K$} (1,-.3) -- (0,-.3) -- (0,0) (2,0) -- (2,-.3) -- (1,-.3);
\draw [very thick] (0,0) -- (0,.25) (2,.25) -- (2,0);
\fill (0,0) circle (1.25pt) (2,0) circle (1.25pt);

\begin{scope}[xshift=5.5cm]
\draw [->, dash pattern=on 2pt off 1pt] (.3,.75) arc (-180:0:.7);
\draw [->, very thick, dash pattern=on 2pt off 1pt] (0,.75) arc (-180:0:1);
\draw [->] (1,1) -- (.3,1) -- (.3,.75) (1.7,.75) -- (1.7,1) -- (1,1);
\draw [->] (1,1.25) node[above]{\scriptsize $\gamma$} (1,1.1) node[below]{\scriptsize $\hat{K}_{\parallel}$} (0.05,1) node[left] {\scriptsize $K(0)$} (1.95,1) node[right] {\scriptsize $K(1)$} (1.95,.7) node[right] {\scriptsize $K$} (1,1.3) -- (0,1.3) -- (0,1) (2,1) -- (2,1.3) -- (1,1.3);
\draw [very thick] (0,1) -- (0,.75) (2,.75) -- (2,1);
\fill (0,1) circle (1.25pt) (2,1) circle (1.25pt);
\end{scope}

\end{tikzpicture}
\caption{Pictures of $\gamma$ and $\hat{K}_{\parallel}$}
\label{figgammaparthetazero}

\end{figure}

In the other cases, in which $K$ goes from bottom to bottom or from top to top, and $K(1)-K(0)$ is equal to $(v,0,0)$ for $v \in \RR \setminus \{0\}$, $lk(K,K_{\parallel})$ is defined so that we have again
$lk(K,K_{\parallel})=lk(K_{\parallel},K)=lk(-K,-K_{\parallel})$.
\end{definition}

\begin{example} Let $K$ be a framed component of an LTR in $\RR^3$, with a regular projection on $(\RR \subset \CC) \times \RR$, and such that $K(1)-K(0)$ is equal to $(v,0,0)$ for $v \in \RR \setminus \{0\}$ if $K$ goes from bottom to bottom or from top to top.
The \emph{self-linking number} of $K$ is its \emph{writhe}, which is the \emph{algebraic number} of its crossings, i.e., the number of positive crossings minus the number of negative crossings in its above regular projection. (Check it as an exercise.)
\end{example}

\begin{definition}
\label{deffstconsttangframed}
When a long tangle representative $\tanghcyll=(K_j)_{j \in \underline{k}}$ is framed by some $\tanghcyll_{\parallel}=(K_{j\parallel})_{j \in \underline{k}}$, set \index[N]{ZZ@$\Zinvuf$ and some variants (see also the summary in the next pages)!ZZfT@$\Zinvufrfneg(\hcylc,(\tanghcyll,\tanghcyll_{\parallel}))$}
\begin{equation*}\Zinvufrfneg\bigl(\hcylc,(\tanghcyll,\tanghcyll_{\parallel})\bigr)=\prod_{j=1}^k\Bigl(\exp\bigl(lk(K_j,K_{j\parallel})\alpha\bigr)\#_j\Bigr)\Zinvuf(\hcylc,\tanghcyll)\end{equation*} with the notation of Theorem~\ref{thmfstconsttang}.
\end{definition}

We will give general relations between $I_{\theta}$ and self-linking numbers in Section~\ref{secstraighttang}. See Proposition~\ref{proplktangles}.

A \emph{tangle representative}\index[T]{tangle!representative} is a pair $(\hcylc,\tanghcyll(\sourcetl) \cap \hcylc)$ for a rational homology cylinder $\hcylc$ and a long tangle representative $\tanghcyll \colon \sourcetl \hookrightarrow \crats(\hcylc)$ as in Definition~\ref{defLTR}. 
Tangle representatives and LTR are in natural one-to-one correspondence, and we also view $\Zinvuf$ as a function of tangle representatives.

A tangle representative $(\hcylc_1,\tanghcyll_1)$ is \emph{right-composable} by a tangle representative $(\hcylc_2,\tanghcyll_2)$ when the top configuration of $(\hcylc_1,\tanghcyll_1)$ coincides with the bottom configuration of $(\hcylc_2,\tanghcyll_2)$. In this case, the \emph{product} \begin{equation*}(\hcylc_1\hcylc_2,\tanghcyll_1\tanghcyll_2)=(\hcylc_1,\tanghcyll_1)(\hcylc_2,\tanghcyll_2)\end{equation*} is obtained by stacking $(\hcylc_2,\tanghcyll_2)$ above $(\hcylc_1,\tanghcyll_1)$, after affine reparametrizations of $\drad{1} \times \left[0,1\right]$, which becomes $\drad{1} \times \left[0,1/2\right]$ for $(\hcylc_1,\tanghcyll_1)$ and $\drad{1} \times \left[1/2,1\right]$ for $(\hcylc_2,\tanghcyll_2)$.

The product of two framed tangle representatives is naturally framed.
We will prove the following functoriality theorem for $\Zinvufrf$ in Section~\ref{secfunc}.

\begin{theorem}
\label{thmfunc}
 $\Zinvufrf$ is \emph{functorial:} For two framed tangle representatives $\tanghcyll_1=(\hcylc_1,\tanghcyll_1)$ and $\tanghcyll_2=(\hcylc_2,\tanghcyll_2)$ such that the top configuration of $\tanghcyll_1$ coincides with the bottom configuration of $\tanghcyll_2$, we have
\begin{equation*}\Zinvufrfneg(\tanghcyll_1\tanghcyll_2)=\Zinvufrfneg\left(\begin{tikzpicture} [baseline=-.1cm] \useasboundingbox (-.4,-.7) rectangle (.4,.7);
\draw (-.3,-.6) rectangle (.3,.6);
\draw (-.3,0) to (.3,0) (0,.3) node{\scriptsize $\tanghcyll_2$} (0,-.3) node{\scriptsize $\tanghcyll_1$};
\end{tikzpicture}\right)=\begin{tikzpicture} [baseline=-.1cm] \useasboundingbox (-.7,-.7) rectangle (.7,.7);
\draw (-.6,-.6) rectangle (.6,.6);
\draw (-.6,0) to (.6,0) (0,.3) node{\scriptsize $\Zinvufrfneg(\tanghcyll_2)$} (0,-.3) node{\scriptsize $\Zinvufrfneg(\tanghcyll_1)$};
\end{tikzpicture}=\Zinvufrfneg(\tanghcyll_1)\Zinvufrfneg(\tanghcyll_2).\end{equation*}
\end{theorem}

The product $\Zinvufrfneg(\tanghcyll_1)\Zinvufrfneg(\tanghcyll_2)$ is the natural product of Section~\ref{secprod}.
When applied to the case in which $\sourcetl_1$ and $\sourcetl_2$ are empty, the above theorem implies that the invariant $\Zinvuf$ of $\QQ$-spheres is multiplicative under connected sum.
\begin{remark}
\label{rkznonfunc}
 The multiplicativity of Theorem~\ref{thmfunc} does not hold for the unframed version $\Zinvuf$ of $\Zinvufrf$.
 Indeed, as an unframed tangle, the vertical product \begin{equation*}
                                      \fulltwistonerev \smaxo                 
                                                     \end{equation*}
equals $\smaxo$. However, the reader can prove 
\begin{equation*}\Zinvuf(\fulltwistonerev)\Zinvuf(\smaxo) \neq \Zinvuf(\smaxo)\end{equation*}
as an exercise, using Lemma~\ref{lemfulltwistdegone} and the behavior of $\Zinvufrf$ under component orientation reversals described in Theorem~\ref{thmmainfunc}. (Here, the image of any involved boundary planar two-point configuration is $\{-1/2,1/2\} \subset \CC$.) 
\end{remark}

Lemma~\ref{lemconnuptd} will imply 
that $\Zinvuf$ is invariant under any isotopy of a tangle representative in a rational homology cylinder during which the bottom and top configurations are constant up to translation and dilation. 

\begin{definition}
 \label{deftangleone}
 In this book, a \indexT{tangle} is an equivalence class of tangle representatives under the equivalence relation that identifies two representatives if and only if they can be obtained from one another by
a diffeomorphism $h$ from the pair $(\hcylc,\tanghcyll)$ to another such $(\hcylc,\tanghcyll^{\prime})$, 
\begin{itemize}
\item which fixes a neighborhood of $\partial \left(\drad{1}  \times\left[0,1\right]\right)$ setwise,
\item which fixes a neighborhood of $(\partial \drad{1})  \times\left[0,1\right]$ pointwise,
\item such that $h(c^-(\finsetb^-) \times \{0\}) \subset \drad{1}\times \{0\}$ coincides with $c^-(\finsetb^-) \times \{0\}$ up to translation and dilation (i.e., as a planar configuration of $\cinjuptd{c^-(\finsetb^-)}{\CC\times \{0\}}$),
and the classes of $h(c^+(\finsetb^+) \times \{1\}) \subset \drad{1}\times \{1\}$ and $c^+(\finsetb^+) \times \{1\}$ in $\cinjuptd{c^+(\finsetb^+)}{\CC\times \{1\}}$ coincide, too,
\end{itemize}
and which is isotopic to the identity map through such diffeomorphisms.
\end{definition}

\section{Defining \texorpdfstring{$\Zinvufrf$}{Z} for combinatorial \texorpdfstring{$q$-tangles}{q-tangles}}
\label{sectangwords}

The tangles of Definition~\ref{deftangleone} can be framed by parallels of their components as before to become \emph{framed tangles}\index[T]{framed!tangle}. These framed tangles are framed cobordisms (up to isotopy) in $\QQ$-cylinders from an injective configuration $c^-(\finsetb^-) \in \cinjuptd{c^-(\finsetb^-)}{\CC}$ of points in $\CC$, up to dilation and translation, to another planar configuration $c^+(\finsetb^+) \in \cinjuptd{c^+(\finsetb^+)}{\CC}$.
Then $\Zinvufrf$ is an invariant of these framed tangles. In Section~\ref{zinvfqbraid}, we extend $\Zinvufrf$ to framed cobordisms between limit configurations in the compactifications $\ccompuptd{c(\finsetb)}{\CC}$ of the spaces $\cinjuptd{c(\finsetb)}{\CC}$ introduced in Subsection~\ref{submoreconfstrat} and studied in Section~\ref{secpresfacanom}. The category of these limit cobordisms called \emph{$q$-tangles} is equipped with interesting cabling operations described in Section~\ref{secintqtangle}. As stated in Theorem~\ref{thmmainfunc}, $\Zinvufrf$ behaves nicely under these cabling operations. In this section, we define $\Zinvufrf$ for particular $q$-tangles, called \emph{combinatorial $q$-tangles}, defined below,
and we give two examples of
cabling properties for these combinatorial $q$-tangles, to motivate and introduce our more general presentation of $q$-tangles in Section~\ref{secintqtangle}.
Dror Bar-Natan, Thang L\^e, and Jun Murakami used these combinatorial $q$-tangles to define a functorial extension of the Kontsevich integral of framed oriented links in \cite{BarNatanNAT} and \cite{LeMur}. 

\begin{definition}
\label{defcombqtang}
A \emph{combinatorial $q$-tangle} is a triple $(\Link,w^-,w^+)$, where
\begin{itemize}
 \item $\Link$ is a framed tangle representative whose bottom and top configurations are on the real line, up to isotopies of $\hcylc$ that globally preserve the intersection of the bottom disk $\drad{1} \times \{0\}$ with $\RR \times \{0\}$ and the intersection of the top disk $\drad{1} \times \{1\}$ with $\RR \times \{1\}$,
 \item $w^-$ and $w^+$ are nonassociative words in the letter $\,\bulletmoyen\,$ (as in Example~\ref{exaStasheff}) such that the letters of $w^-$ (resp. $w^+$) are in canonical one-to-one correspondence with the elements of the domain $\finsetb^-$ (resp. $\finsetb^+$) of the bottom (resp. top) configuration of $\Link$.
\end{itemize} 
 The nonassociative words $w^-$ and $w^+$ are respectively called the \emph{bottom and top configurations} of the combinatorial $q$-tangle $(\Link,w^-,w^+)$. The combinatorial $q$-tangle $(\Link,w^-,w^+)$ is also called the combinatorial $q$-tangle $\Link$ from $w^-$ to $w^+$.
\end{definition}

\begin{examples} \label{excombqtang} 
Combinatorial $q$-tangles in the standard cylinder are unambiguously represented by one of their regular projections to $(\RR  \subset \CC) \times \left[0,1\right]$ such that the parallels of their components are parallel in the figures, together with their bottom and top configurations.
Examples of these combinatorial $q$-tangles include 
\begin{equation*}\tanposprod\mbox{, } \assocnonnum\mbox{, }
\tancuphigh
\mbox{, and }
\tancap. \end{equation*}
\end{examples}

Recall from Example~\ref{exaStasheff} that the involved nonassociative words are corners of $\ccompuptd{<,\underline{k}}{\RR} \subset\ccompuptd{\underline{k}}{\CC}$.
A combinatorial $q$-tangle $(\Link,w^-,w^+)$ from a bottom word $w^-$ to a top word $w^+$ is thought of as the limit,
when $t$ tends to $0$, 
of framed tangles $\Link(t)$  whose bottom and top configurations are the configurations $w^-(t)$ and $w^+(t)$ defined in Example~\ref{exaStasheff}, in the isotopy class of $\Link$ with respect to the isotopies of Definition~\ref{defcombqtang}.

In Theorem~\ref{thmpoirbraid} and Remark~\ref{deflimzqtang}, following Sylvain Poirier \cite{poirierv2}, we prove that $\lim_{t \to 0}\Zinvufrfneg(\Link(t))$ exists and that the formula
\begin{equation*}\Zinvufrfneg(\Link)=\lim_{t \to 0}\Zinvufrfneg\bigl(\Link(t)\bigr)\end{equation*}
defines an isotopy invariant $\Zinvufrf$ of these (framed) combinatorial $q$-tangles. The isotopy invariant $\Zinvufrf$ will thus behave naturally with respect to the deletion of components. If $\tanghcyll^{\prime}$ is a subtangle of a tangle $\tanghcyll$ with domain  $\sourcetl^{\prime}$, then $\Zinvufrfneg(\tanghcyll^{\prime})$ will be obtained from $\Zinvufrfneg(\tanghcyll)$ by forgetting $\sourcetl \setminus \sourcetl^{\prime}$ and all the diagrams with univalent vertices on $\sourcetl \setminus \sourcetl^{\prime}$.

As a first example, let us compute $\Zinvufrf$ for the combinatorial $q$-tangle $\twostraightstrandspetit$. 
\begin{lemma}
\label{lemtwostraightstrands}
We have \begin{equation*}\Zinvufrfneg\left(\twostraightstrandspetit\right)=1=\left[\nochordtwosorpetit\right].\end{equation*}
\end{lemma}
\bp By definition, the left-hand side is the limit, when $t$ tends to $0$, of the evaluation of $\Zinv$ of the LTR whose image is $\{0,t\}\times \RR$. 
There is an action of $\RR$ by vertical translation on the involved configuration spaces. The integrated forms factor through the quotients by this action of $\RR$ of the configuration spaces, whose dimensions are smaller (by one) than the degrees of the integrated forms. So the integrals vanish for all nonempty diagrams.
\eop

As a second example, we compute $\Zinvufrf_{\leq 1}$, which is the truncation of $\Zinvufrf$ in degrees lower than $2$, for the combinatorial $q$-tangle $\fulltwist$.

\begin{lemma} \label{lemfulltwist} We have
\begin{equation*}\Zinvufrf_{\leq 1}\left(\fulltwist\right)=\left[\nochordtwosor\right]+\left[\onechordtwosor\right].\end{equation*}
\end{lemma}
\bp Since $\Zinvuf$ is invariant under the isotopies that preserve the bottom and the top configurations,
the contributions of the nonempty diagrams whose univalent vertices are on one strand of the tangle
vanish. So we are left with the contribution of the numbered graphs with one vertex on each strand, treated in Lemma~\ref{lemfulltwistdegone}.
\eop

The obtained invariant $\Zinvufrf$ is still multiplicative under vertical composition as in Theorem~\ref{thmfunc}, and we can now define other interesting operations.

For two combinatorial $q$-tangles $\tanghcyll_1=(\hcylc_1,\tanghcyll_1)$ from $w_1^-$ to $w_1^+$, and $\tanghcyll_2=(\hcylc_2,\tanghcyll_2)$ from $w_2^-$ to $w_2^+$, define 
the product $\tanghcyll_1 \otimes \tanghcyll_2$, from the bottom configuration $w_1^-w_2^-$ to the top configuration $w_1^+w_2^+$, by shrinking $\hcylc_1$ and $\hcylc_2$ to make them respectively replace the products by $\left[0,1\right]$ of the horizontal disks with radius $\frac14$ and respective centers $-\frac12$ and $\frac12$.

\begin{examples} \label{excombqtangtensor} 
We have
\begin{equation*}\tanposprod \otimes \assocnonnum=
\begin{tikzpicture}  [baseline=-.1 cm]
\useasboundingbox (-.4,-.3) rectangle (1.4,.4);
\draw [->] (.15,-.15) -- (-.15,.15);
\draw [draw=white,double=black,very thick]  (-.15,-.15) -- (.15,.15);
\draw [->]  (-.15,-.15) -- (.15,.15);
\draw (-.23,-.25) node{\scriptsize \textrm{(}} (.23,-.25) node{\scriptsize \textrm{)}};
\fill (-.15,-.25) circle (1pt) (.15,-.25) circle (1pt);
\draw (-.23,.25) node{\scriptsize \textrm{(}} (.23,.25) node{\scriptsize \textrm{)}};
\fill (-.15,.25) circle (1pt) (.15,.25) circle (1pt);
\draw (-.3,.25) node{\scriptsize \textrm{(}} (1.25,.25) node{\scriptsize \textrm{)}};
\draw (-.3,-.25) node{\scriptsize \textrm{(}} (1.15,-.25) node{\scriptsize \textrm{)}};
\begin{scope}[xshift=+.75cm]
\draw [->] (-.2,-.15) -- (-.2,.15);
\draw [->] (.2,-.15) -- (.2,.15);
\draw [->] (-.1,-.15) -- (.1,.15);
\draw (-.33,-.25) node{\scriptsize \textrm{((}} (0,-.25) node{\scriptsize  \textrm{)}} (.27,-.25) node{\scriptsize \textrm{)}};
\fill (-.2,-.25) circle (1pt) (.2,-.25) circle (1pt)  (-.1,-.25) circle (1pt);
\draw (-.27,.25) node{\scriptsize \textrm{(}}  (0,.25) node{\scriptsize \textrm{(}} (.33,.25) node{\scriptsize  \textrm{))}};
\fill (-.2,.25) circle (1pt) (.2,.25) circle (1pt) (.1,.25) circle (1pt);\end{scope}
\end{tikzpicture}\mbox{, and } 
\tancuphigh \otimes \tancap = \tancuphigh \tancap.\end{equation*}
\end{examples}

We can easily deduce the following theorem from the cabling property and the functoriality property of Theorem~\ref{thmmainfunc}.
\begin{theorem}
 \label{thmmonoid}
 $\Zinvufrf$ is \emph{monoidal:} For two combinatorial $q$-tangles $\tanghcyll_1$ and $\tanghcyll_2$, we have
\begin{equation*}\Zinvufrfneg(\tanghcyll_1 \otimes \tanghcyll_2)=\Zinvufrfneg\biggl(\begin{tikzpicture} [baseline=-.1cm] \useasboundingbox (-.9,-.3) rectangle (.9,.5);
\draw (-.8,-.2) rectangle (.8,.4);
\draw (0,-.2) to (0,.4) (.4,.1) node{\scriptsize $\tanghcyll_2$} (-.4,.1) node{\scriptsize $\tanghcyll_1$};
\end{tikzpicture}\biggr)=\begin{tikzpicture} [baseline=-.1cm] \useasboundingbox (-1.3,-.3) rectangle (1.3,.5);
\draw (-1.2,-.2) rectangle (1.2,.4);
\draw (0,-.2) to (0,.4) (.6,.1) node{\scriptsize $\Zinvufrfneg(\tanghcyll_2)$} (-.6,.1) node{\scriptsize $\Zinvufrfneg(\tanghcyll_1)$};
\end{tikzpicture}=\Zinvufrfneg(\tanghcyll_1) \otimes \Zinvufrfneg(\tanghcyll_2),\end{equation*}
where $\Zinvufrfneg(\tanghcyll_1) \otimes \Zinvufrfneg(\tanghcyll_2)$ denotes the image of 
$\Zinvufrfneg(\tanghcyll_1) \otimes \Zinvufrfneg(\tanghcyll_2)$ under the natural product from $\Aavis(\sourcetl_1) \otimes_{\RR} \Aavis(\sourcetl_2)$ to $\Aavis(\sourcetl_1 \sqcup \sourcetl_2)$ induced by the disjoint union of diagrams.
\end{theorem}

We can also \emph{double} a component $K$ according to its parallelization in a combinatorial $q$-tangle $\Link$. This operation replaces a component by two parallel components, and if this component has boundary points, it replaces the corresponding letters in the nonassociative words with $(\bulletpetit\bulletpetit)$. The obtained combinatorial $q$-tangle is denoted by $\tanghcyll( 2 \times K)$. For example, we have
\newcommand{\tanposprodessai}{\begin{tikzpicture} [baseline=-.1 cm]
\useasboundingbox (-.3,-.3) rectangle (.3,.4);
\draw [->] (.15,-.15) -- (-.15,.15);
\draw [draw=white,double=black,very thick]  (-.15,-.15) -- (.15,.15);
\draw [->]  (-.15,-.15) -- (.15,.15);
\draw (-.23,-.25) node{\scriptsize \textrm{(}} (.23,-.25) node{\scriptsize \textrm{)}};
\fill (-.15,-.25) circle (1pt) (.15,-.25) circle (1pt);
\draw (-.23,.25) node{\scriptsize \textrm{(}} (.23,.25) node{\scriptsize \textrm{)}};
\fill (-.15,.25) circle (1pt) (.15,.25) circle (1pt);
\end{tikzpicture}}
\begin{equation*}\tanposprodessai\left(2 \times \begin{tikzpicture}
\useasboundingbox (-.2,-.1) rectangle (.2,.2);
\draw [->]  (-.15,-.15) -- (.15,.15);
\end{tikzpicture}\right) = 
\begin{tikzpicture} [baseline=-.1 cm] 
\useasboundingbox (-.4,-.3) rectangle (1.3,.4);
\draw [->] (.25,-.15) -- (-.15,.15);
\draw [draw=white,double=black,very thick]  (-.15,-.15) -- (.15,.15);
\draw [->]  (-.15,-.15) -- (.15,.15);
\draw [draw=white,double=black,very thick]  (-.05,-.15) -- (.25,.15);
\draw [->]  (-.05,-.15) -- (.25,.15);
\fill (-.15,-.25) circle (1pt) (-.05,-.25) circle (1pt) (.25,-.25) circle (1pt);
\draw (-.28,-.25) node{\scriptsize \textrm{((}} (.05,-.25) node{\scriptsize  \textrm{)}} (.35,-.25) node{\scriptsize \textrm{)}};
\fill (-.15,.25) circle (1pt) (.15,.25) circle (1pt) (.25,.25) circle (1pt);
\draw (-.25,.25) node{\scriptsize \textrm{(}}  (0.05,.25) node{\scriptsize \textrm{(}} (.36,.25) node{\scriptsize  \textrm{))}};
\end{tikzpicture} \end{equation*}

The following duplication property for $\Zinvufrf$ is a part of Theorem~\ref{thmmainfunc},
which is proved in Section~\ref{secduprop}. 

\begin{theorem}
 \label{thmdup}
 Let $K$ be a component of a combinatorial $q$-tangle $\tanghcyll$.\footnote{The theorem holds for components going from bottom to bottom or from top to top. It would not hold in this generality if $\Zinvufrf$ were replaced by certain functorial extensions of the Kontsevich integral.} Recall Notation~\ref{notationduplication}. Then we have
\begin{equation*}\Zinvufrfneg\bigl(\tanghcyll( 2 \times K)\bigr)=\pi(2 \times K)^{\ast}\Zinvufrfneg(\tanghcyll).\end{equation*}
\end{theorem}

These results are some particular cases of the properties of $\Zinvufrf$, which are listed in Theorem~\ref{thmmainfunc} and proved in Chapter~\ref{chappropzinvffunc}.

\section{Good monoidal functors for combinatorial \texorpdfstring{$q$-tangles}{q--tangles}}
\label{secgoodmonfonc}

Recall Notation~\ref{notAcheck}. In this section, all combinatorial $q$-tangles are combinatorial $q$-tangles in $\RR^3$ or in the standard cylinder, and coefficients for spaces of Jacobi diagrams are in $\KK=\CC$.
In \cite{poirierv2}, Sylvain Poirier extended the natural projection $\Zinvlinkufrfneg(\RR^3,.)$ of  $\Zinvufrfneg(\RR^3,.)$ in $\Assis(.;\RR) \subset \Assis(.;\CC)$, from framed links of $\RR^3$ to these combinatorial $q$-tangles, in an elegant way, and he proved that his extension $Z^l$ is a \emph{good monoidal functor} with respect to the definition below.
Good monoidal functors on the category of combinatorial $q$-tangles (in $\RR^3$) are characterized in \cite{lesunikon}.
This section reviews these results of \cite{poirierv2} and \cite{lesunikon}.
The quoted results of Sylvain Poirier will be reproved (with much more details) and generalized in this book. In contrast, the proofs of the results of \cite{lesunikon} will not be reproduced in this book since they do not involve analysis on configuration spaces.

\begin{definition}
\label{defgoodmon}
A {\em good monoidal functor\/} from the category of combinatorial $q$-tangles (in $\RR^3$) to the category
of spaces of Jacobi diagrams
is a map $\Zgen$ sending a combinatorial $q$-tangle $\tanghcyll$ with domain $\sourcetl$ to an element of the space $\Assis(\sourcetl;\CC)$ of Notation~\ref{notAcheck} with the following properties.
\begin{itemize}
\item For any combinatorial $q$-tangle $\tanghcyll$, the degree zero part $\Zgen_0(\tanghcyll)$ of $\Zgen(\tanghcyll)$ is 1, which is the class of the empty diagram.
\item $\Zgen$ is \emph{functorial:} For two combinatorial $q$-tangles $\tanghcyll_1$ and $\tanghcyll_2$ such that the top configuration $\toptang(\tanghcyll_1)$ of $\tanghcyll_1$ coincides with the bottom configuration $\bottom(\tanghcyll_2)$ of $\tanghcyll_2$, we have
\begin{equation*}\Zgen(\tanghcyll_1\tanghcyll_2)=\Zgen\left(\begin{tikzpicture} [baseline=-.1cm] \useasboundingbox (-.4,-.7) rectangle (.4,.7);
\draw (-.3,-.6) rectangle (.3,.6);
\draw (-.3,0) to (.3,0) (0,.3) node{\scriptsize $\tanghcyll_2$} (0,-.3) node{\scriptsize $\tanghcyll_1$};
\end{tikzpicture}\right)=\begin{tikzpicture} [baseline=-.1cm] \useasboundingbox (-.7,-.7) rectangle (.7,.7);
\draw (-.6,-.6) rectangle (.6,.6);
\draw (-.6,0) to (.6,0) (0,.3) node{\scriptsize $\Zgen(\tanghcyll_2)$} (0,-.3) node{\scriptsize $\Zgen(\tanghcyll_1)$};
\end{tikzpicture}=\Zgen(\tanghcyll_1)\Zgen(\tanghcyll_2).\end{equation*}
\item $\Zgen$ is \emph{monoidal:} For two combinatorial $q$-tangles $\tanghcyll_1$ and $\tanghcyll_2$, we have
\begin{equation*}\Zgen(\tanghcyll_1 \otimes \tanghcyll_2)=\Zgen\biggl(\begin{tikzpicture} [baseline=-.1cm] \useasboundingbox (-.9,-.3) rectangle (.9,.5);
\draw (-.8,-.2) rectangle (.8,.4);
\draw (0,-.2) to (0,.4) (.4,.1) node{\scriptsize $\tanghcyll_2$} (-.4,.1) node{\scriptsize $\tanghcyll_1$};
\end{tikzpicture}\biggr)=\begin{tikzpicture} [baseline=-.1cm] \useasboundingbox (-1.1,-.3) rectangle (1.1,.5);
\draw (-1,-.2) rectangle (1,.4);
\draw (0,-.2) to (0,.4) (.5,.1) node{\scriptsize $\Zgen(\tanghcyll_2)$} (-.5,.1) node{\scriptsize $\Zgen(\tanghcyll_1)$};
\end{tikzpicture}=\Zgen(\tanghcyll_1) \otimes \Zgen(\tanghcyll_2).\end{equation*}
\item $\Zgen$ is \emph{compatible with the deletion} of a component: If $\tanghcyll^{\prime}$ is a subtangle of $\tanghcyll$ with domain  $\sourcetl^{\prime}$, then $\Zgen(\tanghcyll^{\prime})$ is obtained from $\Zgen(\tanghcyll)$ by forgetting $\sourcetl \setminus \sourcetl^{\prime}$ and all the diagrams with univalent vertices on $\sourcetl \setminus \sourcetl^{\prime}$.
\item $\Zgen$ is \emph{compatible with the duplication} of a \emph{regular} component, which is a component that can be represented without horizontal tangent vectors: For such a component $K$ of a combinatorial $q$-tangle $\tanghcyll$, we have
\begin{equation*}\Zgen\bigl(\tanghcyll( 2 \times K)\bigr)=\pi(2 \times K)^{\ast}\Zgen(\tanghcyll)\end{equation*}
with respect to Notation~\ref{notationduplication}.
\item $\Zgen$ is \emph{invariant under the 180-degree rotation around a vertical axis} through the real line. 
\item Let $s_{\frac12}$ be the reflection with respect to the horizontal plane at height $\frac12$.
Let $\sigma_{\frac12}$ be the linear endomorphism of the topological vector space $\Assis(S^1)$ such that $\sigma_{\frac12}\left(\left[\Gamma\right]\right) = (-1)^d \left[\Gamma\right]$, for any element $\left[\Gamma\right]$ of $\Assis_d(S^1)$. For any framed knot $K=K(S^1)$, $\Zgen(K)$ belongs to $\Assis(S^1;\RR)$ and we have
\begin{equation*}\Zgen \circ s_{\frac12}(K) =\sigma_{\frac12} \circ \Zgen(K).\end{equation*} 
\item $\Zgen$ behaves as in Proposition~\ref{proporcomp} with respect to reversals of component orientations. It can be defined as an invariant of unoriented tangles valued in a space of diagrams whose support is the unoriented domain of the tangle, as in Definitions~\ref{defdia} and \ref{defrkoruniv}.
\item The degree one part $a^\Zgen_1$ of the element $a^\Zgen \in \Assis\left(\left[0,1\right]\right)$ such that $a^\Zgen_0=0$ and \begin{equation*}\Zgen\left(\smaxno \right) = \exp(a^\Zgen)\Zgen\left(\smaxo \right)\end{equation*} is \begin{equation*}a^\Zgen_1 = \frac12\left[ \onechordsmalltseullong \right].\end{equation*}
\end{itemize}
\end{definition}

As shown in \cite[Proposition 4.2, 10]{lesunikon}, a good monoidal functor $\Zgen$ is determined by its values 
\begin{equation*}\Zgen\left(\tanposprod \right)\mbox{,} \; \Zgen\left(\assocnonnum \right) \mbox{, and } \; \Zgen\left( \tancap\right).\end{equation*}
Moreover, $\Zgen$ can be computed combinatorially from these three ingredients. See also \cite{LeMur,BarNatanNAT}.
The value $\Zgen\left(\assocnonnum \right)$ is called an \indexT{associator}.

In \cite{LeMur}, Thang L\^e and Jun Murakami constructed the first example of such a good monoidal functor from the Kontsevich integral of links in $\RR^3$ described in \cite{chmutovduzhin}.\footnote{They also explain \cite[Section 5]{LeMur} how to construct such a functor combinatorially, from an element of $\Aavis(|||)$ satisfying some equations \cite[Section 3, A1--A4]{LeMur}, which allow it to be the associator of a good monoidal functor. See also \cite{BarNatanNAT}.} See also \cite{lesintrokonsummer}. We denote the L\^e--Murakami functor by $Z^K$ and call it the \emph{Kontsevich integral} of combinatorial $q$-tangles. This Kontsevich integral furthermore satisfies
\begin{equation*}Z^K\left(\fulltwist \right)= \exp\left(\onechordtwosor \right)=\nochordtwosor+\left[\onechordtwosor\right]+\frac12 \left[\twochordtwosor\right] + \frac16 \left[\threechordtwosor\right] +\dots \end{equation*}
This easily implies that the element $a^{Z^K}$ of the above definition is 
\begin{equation*}a^{Z^K} = \frac12\left[ \onechordsmalltseullong \right].\end{equation*}
So it vanishes in all degrees greater than one.

In \cite{poirierv2}, Sylvain Poirier extended $\Zinvlinkufrfneg(\RR^3,.)$ from framed links of $\RR^3$ to combinatorial $q$-tangles of $\RR^3$ and he proved that his extension $Z^l$ satisfies the above properties with
\begin{equation*}a^{Z^l}=\alpha,\end{equation*}
where $\alpha$ is the anomaly of Section~\ref{secanomalpha}.

\begin{remark}
The published version \cite{poirier} of \cite{poirierv2} does not contain the cited important results of \cite{poirierv2}, which will be generalized and proved with much more details in the present book.
\end{remark}

\begin{definition}
\label{deftwoleg}
Say that an element $\gamtwoleg=(\gamtwoleg_n)_{n \in {\NN}}$ in $\Assis\left(\left[0,1\right]\right)$ is a {\em two-leg element\/} if, for any $n \in {\NN}$, $\gamtwoleg_n$ is a combination of diagrams with two univalent 
vertices.
\end{definition}
 
Forgetting $\left[0,1\right]$ from such a two-leg element gives rise to a unique series $\gamtwoleg^s$ of diagrams
with two distinguished univalent vertices $v_1$ and $v_2$, such that $\gamtwoleg^s$ is symmetric with respect to the exchange of $v_1$ and $v_2$, according to the following lemma due to Pierre Vogel. See \cite[Corollary 4.2]{vogel}.

\begin{lemma}[Vogel] 
\label{lemsymvo}
Two-leg Jacobi diagrams are symmetric with respect to the exchange of their
two legs in a diagram space quotiented by the AS and Jacobi relations.
\end{lemma}
\bp Since a chord is obviously symmetric, we can 
restrict ourselves to 
a two-leg diagram with at least one trivalent vertex, and whose two univalent vertices are respectively numbered by $1$ and $2$. We draw it as 
\begin{equation*} \symvogone \end{equation*}
where we do not represent the trivalent part inside the disk bounded by the thick gray topological circle. Applying Lemma~\ref{lemcom}, when the annulus is a neighborhood of the thick topological circle that contains the pictured trivalent vertex, proves
\begin{equation*} \left[ \symvogone \right] = \left[\symvogtwo \right]. \end{equation*}

Similarly, we have
\begin{equation*}\left[ \begin{tikzpicture} [baseline=.5cm] \useasboundingbox (-.9,.1) rectangle (.9,1.3);
\draw[very thick,gray] (0,1.2) .. controls (-.8,1.2) .. (-.8,1) .. controls (-.8,.8) .. (0,.8)
 .. controls (.8,.8) .. (.8,1) .. controls (.8,1.2) .. (0,1.2);
\draw (0,.9) -- (-.4,.2) (.4,.9) -- (.4,.2) (-.4,.9) .. controls (-.4,.55) .. (-.2,.55) (-.55,.2) node{\scriptsize 1} (.55,.2) node{\scriptsize 2};
\fill (-.2,.55) circle (1.5pt) (-.4,.2) circle (1.5pt) (.4,.2) circle (1.5pt);
\end{tikzpicture}\right]  
= -\left[\begin{tikzpicture} [baseline=.5cm] \useasboundingbox (-.9,.1) rectangle (.9,1.3);
\draw[very thick,gray] (0,1.2) .. controls (-.8,1.2) .. (-.8,1) .. controls (-.8,.8) .. (0,.8)
 .. controls (.8,.8) .. (.8,1) .. controls (.8,1.2) .. (0,1.2);
\draw (0,.9) -- (-.4,.2) (.4,.9) -- (.4,.2) (-.4,.9) .. controls (-.4,.55) .. (.4,.5) (-.55,.2) node{\scriptsize 1} (.55,.2) node{\scriptsize 2};
\fill (.4,.5) circle (1.5pt) (-.4,.2) circle (1.5pt) (.4,.2) circle (1.5pt);
\end{tikzpicture} \right] 
\;\;\mbox{and}\;\;
\left[\begin{tikzpicture} [baseline=.5cm] \useasboundingbox (-.9,.1) rectangle (.9,1.3);
\draw[very thick,gray] (0,1.2) .. controls (-.8,1.2) .. (-.8,1) .. controls (-.8,.8) .. (0,.8)
 .. controls (.8,.8) .. (.8,1) .. controls (.8,1.2) .. (0,1.2);
\draw (0,.9) -- (-.4,.2) (-.4,.9) -- (.4,.2) (.4,.9) .. controls (.4,.55) .. (0,.55) (-.55,.2) node{\scriptsize 1} (.55,.2) node{\scriptsize 2};
\fill (0,.55) circle (1.5pt) (-.4,.2) circle (1.5pt) (.4,.2) circle (1.5pt);
\end{tikzpicture}\right] 
= -\left[\begin{tikzpicture} [baseline=.5cm] \useasboundingbox (-.9,.1) rectangle (.9,1.3);
\draw[very thick,gray] (0,1.2) .. controls (-.8,1.2) .. (-.8,1) .. controls (-.8,.8) .. (0,.8)
 .. controls (.8,.8) .. (.8,1) .. controls (.8,1.2) .. (0,1.2);
\draw (0,.9) -- (-.4,.2) (-.4,.9) -- (.4,.2) (.4,.9) .. controls (.4,.55) .. (-.2,.55) (-.55,.2) node{\scriptsize 1} (.55,.2) node{\scriptsize 2};
\fill (-.2,.55) circle (1.5pt) (-.4,.2) circle (1.5pt) (.4,.2) circle (1.5pt);
\end{tikzpicture}\right].\end{equation*}
So we get
\begin{equation*} \left[\symvogone \right]= \left[\symvogonesym\right].\end{equation*}
\eop

\begin{definition}
\label{defbetas}
Let  $\gamtwoleg$ be a two-leg element of $\Assis\left(\left[0,1\right]\right)$, recall that $\gamtwoleg^s$ is the series obtained from $\gamtwoleg$ by erasing $\left[0,1\right]$.
For a chord diagram $\Gamma$, define $\Psi(\gamtwoleg)(\Gamma)$ by replacing each chord
by $\gamtwoleg^s$. As it is proved in \cite[Lemmas~6.1 and 6.2]{lesunikon}, $\Psi(\gamtwoleg)$ is a well-defined morphism of topological vector spaces from $\Assis(C)$ to $\Assis(C)$ for any one-manifold $C$, and $\Psi(\gamtwoleg)$ is an isomorphism as soon as $\gamtwoleg_1 \neq 0$.
\end{definition}

The following theorem is Theorem~1.3 in \cite{lesunikon}.
\begin{theorem}
\label{thmlesunikon}
Let $\Zgen$ be a good monoidal functor as above. Then $a^\Zgen$ is a
 two-leg element of $\Assis\left(\left[0,1\right]\right)$ such that $a^\Zgen_{2i}=0$ for any integer $i$,
and we have
\begin{equation*}\Zgen(L) = \Psi(2a^\Zgen)\bigl(Z^K(L)\bigr)\end{equation*} for any framed link $L$, where $Z^K$ denotes the Kontsevich integral of framed links (denoted by $\hat{Z}_f$ in \cite{LeMur}, and by $Z$ in \cite{lesintrokonsummer}).
\end{theorem}

The following corollary is a particular case of \cite[Corollary~1.4]{lesunikon}.

\begin{corollary}
 \label{corlesunikontwoleg}
 The anomaly $\alpha$ is a two-leg element of $\Assis\left(\left[0,1\right]\right)$, and we have \begin{equation*}\Zinvlinkufrfneg(\RR^3,\Link)=\Psi(2\alpha)\bigl(Z^K(L)\bigr)\end{equation*} for any framed link $\Link$ of $\RR^3$.
\end{corollary}

\begin{note} \label{notelesunikonfunc}
Theorem~\ref{thmlesunikon} would not hold if framed links were replaced by arbitrary tangles. For pure combinatorial $q$-tangles from a nonassociative word to itself, whose components go from top to bottom or from bottom to top and connect a letter of the bottom word to the corresponding letter of the top word, it is true, up to conjugation by a \emph{twisting function} of nonassociative words, introduced by Thang L\^e and Jun Murakami in \cite[Section 7]{LeMur}. A generalization of Theorem~\ref{thmlesunikon} to all combinatorial $q$-tangles
whose components go from top to bottom or from bottom to top is given in \cite[Theorem 8.5]{lesunikon}. It involves the L\^e--Murakami twisting functions.
\end{note}

\chapter{More on the functor \texorpdfstring{$\Zinvufrf$}{Z}}
\label{chapzinvtang}

In this chapter, we state our general Theorem~\ref{thmmainfunc}, which ensures that $\Zinvufrf$ is a functor behaving naturally with respect to various structures of the category of $q$-tangles, such as cabling or duplication.
We first describe the category of $q$-tangles in Section~\ref{secintqtangle} 
before stating Theorem~\ref{thmmainfunc} in Section~\ref{secexppty}.

We describe the main steps of the generalization of the construction of $\Zinvufrf$ to $q$-tangles in Section~\ref{secdefzinvf}.
We will give the details of these steps in Chapters~\ref{chapconszinvf}, \ref{chapzinvfbraid}, and \ref{chapdiscrex}. We will finish the proof of Theorem~\ref{thmmainfunc} in Chapter~\ref{chappropzinvffunc}.

\section{Tangles and \texorpdfstring{$q$-tangles}{q--tangles}}
\label{secintqtangle}

Recall that a tangle representative is a pair $(\hcylc,\tanghcyll(\sourcetl) \cap \hcylc)$ for a rational homology cylinder $\hcylc$ and a long tangle representative $\tanghcyll \colon \sourcetl \hookrightarrow \crats(\hcylc)$ as in Definition~\ref{defLTR}. 

\begin{definition}
\label{deftanglelong}

In this book, a \indexT{braid representative} is a tangle representative $T(\tilde{\gamma})$ of the standard cylinder $\drad{1} \times \left[0,1\right]$ whose components called \emph{strands}
may be expressed as \begin{equation*}\Bigl\{\bigl(\tilde{\gamma}_{\eltb}(t),t\bigr) \suchthat t \in \left[0,1\right]\Bigr\},\end{equation*} for an element $\eltb$ of a finite set $\finsetb$, which labels the strands.
In the above expression $\tilde{\gamma}_{\eltb} \colon \left[0,1\right] \to \dorad{1}$ is a path, and, for any $t\in \left[0,1\right]$ and any pair $(\eltb,\eltb^{\prime})$ of distinct elements of $\finsetb$, we have $\tilde{\gamma}_{\eltb}(t)\neq\tilde{\gamma}_{\eltb^{\prime}}(t)$ as in Figure~\ref{figbraidthreestrands}.
Such a braid representative can be viewed naturally as a path $\tilde{\gamma} \colon \left[0,1\right] \to \check{C}_{\finsetb}[\dorad{1}]$, where $\check{C}_{\finsetb}[\dorad{1}]$ is the space of injections of $\finsetb$ into $\dorad{1}$, defined in the beginning of Section~\ref{secblodiag}.

\bfig
\centering
\begin{tikzpicture}
\useasboundingbox (0,-.3) rectangle (2,2.2);
\draw [thin,dashed] (0,0)  .. controls (0,.13) and (.85,.3) .. (1.2,.3) .. controls (1.5,.3) and (2,.13) ..  (2,0);
\draw [thin] (0,0)  .. controls (0,-.13) and (.5,-.3) .. (.8,-.3) .. controls (1.15,-.3) and (2,-.13) ..  (2,0);
\draw [thin,dashed] (0,1.9)  .. controls (0,2.03) and (.85,2.2) .. (1.2,2.2) .. controls (1.5,2.2) and (2,2.03) ..  (2,1.9);
\fill (1.5,0)  circle (1pt) (.85,.1)  circle (1pt) (.65,-.1)  circle (1pt)
 (1.5,1.9)  circle (1pt) (.85,2)  circle (1pt) (.65,1.8)  circle (1pt);
 \draw[->] (1.5,0)  .. controls (1.5,.5) and (.85,1.5) .. (.85,2);
 \draw[->] (.85,.1)   .. controls (.85,.6) and (.65,1.3) .. (.65,1.8);
 \draw [draw=white,double=black,very thick,->] (.65,-.1)  .. controls (.65,.4) and (1.5,1.4) .. (1.5,1.9);
 \draw [->] (.65,-.1)  .. controls (.65,.4) and (1.5,1.4) .. (1.5,1.9);
 \draw [draw=white,double=black,very thick] (0,1.9)  .. controls (0,1.77) and (.5,1.6) .. (.8,1.6) .. controls (1.15,1.6) and (2,1.77) ..  (2,1.9);
 \draw  [thin] (0,0) -- (0,1.9) (2,0) -- (2,1.9);
 \draw (2,.95) node[right]{$\gamma_{1,3}$};
\end{tikzpicture}
\caption{A braid representative with three strands}
\label{figbraidthreestrands}

\end{figure}

In this book, a \indexT{braid} (resp. a \indexT{q--braid}) is a homotopy class of paths $\gamma \colon \left[0,1\right] \to \cinjuptd{\finsetb}{\CC}$  (resp. of paths $\gamma \colon \left[0,1\right] \to \ccompuptd{\finsetb}{\CC}$) for some finite set $\finsetb$, where $\ccompuptd{\finsetb}{\CC}$ is the compactification of $\cinjuptd{\finsetb}{\CC}$ described in Theorem~\ref{thmcompfacanom}.
A braid $\gamma$ induces the tangle $T(\gamma)$, which is also called a braid, as above.
The path $\overline{\gamma}$ is the path such that $\overline{\gamma}(t)=\gamma(1-t)$.
A braid is naturally framed by the parallels obtained by pushing it in the direction of the real line of $\CC$.

Tangles, as in Definition~\ref{deftangleone}, can be multiplied if they have representatives that can be, i.e., if the top configuration of the first tangle agrees with the bottom configuration of the second one, up to dilation and translation. The product is associative.
Framed tangles multiply vertically to give rise to framed tangles.

A \emph{$q$-tangle} \index[T]{qtangle@$q$-tangle} is a framed tangle whose bottom and top configurations are 
allowed to be limit configurations in some $\ccompuptd{\finsetb^-}{\CC}$ and in some $\ccompuptd{\finsetb^+}{\CC}$. 
More precisely, a $q$-tangle is represented by a product \begin{equation*}T(\gamma^-)(\hcylc,\tanghcyll)T(\gamma^+),\end{equation*}
as in Figure~\ref{figqtangrep}, where $\gamma^-$ and $\gamma^+$ are $q$-braids, $(\hcylc,\tanghcyll)$ is a framed tangle whose bottom configuration is $\gamma^-(1)$ and whose top configuration is $\gamma^+(0)$, and the strands of $T(\gamma^-)$ and $T(\gamma^+)$ get their orientations from the orientation of $\tanghcyll$. For consistency, we allow braids with $0$ or $1$ strand, and we agree that $\ccompuptd{\emptyset}{\CC}$ and $\ccompuptd{\{b\}}{\CC}$ each have one element, which is the unique configuration of one point in $\CC$ up to translation in the latter case. Note that the restriction of a $q$-tangle to one of its components is a framed tangle since configurations of at most two points are always injective. The components of a $q$-tangle representative are framed since braids are.

\bfig
\centering
\begin{tikzpicture}
\useasboundingbox (0,0) rectangle (1.8,2);
\draw (0,.6) rectangle (1.8,1.4);
\draw (.3,.6) arc (180:90:.4);
\draw [draw=white,double=black,very thick] (.2,1.4)  .. controls (.2,.8) .. (.4,.8);
\draw (.6,1.4)  .. controls (.6,.8) .. (.4,.8);
\draw (.7,.6)  .. controls (.7,.9) and (1.5,1) .. (1.5,1.4);
\draw [draw=white,double=black,very thick] (.7,1) arc (90:0:.4);
\draw (.9,1.4) arc (-180:0:.15);
\draw (.3,0) -- (.3,.6) (1.1,0) -- (1.1,.6) (.7,.6)  .. controls (.7,.3) .. (.3,0);
\draw (.2,1.4)  .. controls (.2,1.7) .. (.6,2) (.6,1.4) -- (.6,2) (.9,1.4)  .. controls (.9,1.6) and (.6,1.7) .. (.6,2) (1.2,1.4) .. controls (1.2,1.6) .. (1.5,2) -- (1.5,1.4);
\draw (0,1) node[left]{\scriptsize $(\hcylc,\tanghcyll)$} (0,.3) node[left]{\scriptsize $T(\gamma^-)$} (0,1.7) node[left]{\scriptsize $T(\gamma^+)$};
\end{tikzpicture}
\caption{A $q$-tangle representative}
\label{figqtangrep}

\end{figure}

Now, $q$-tangles are classes of these representatives under the equivalence relation that identifies $T(\gamma^-)(\hcylc,\tanghcyll)T(\gamma^+)$ with  $T(\gamma^{-\prime})(\hcylc^{\prime},\tanghcyll^{\prime})T(\gamma^{+\prime})$
if and only if
$\gamma^-(0)=\gamma^{-\prime}(0)$, $\gamma^+(1)=\gamma^{+\prime}(1)$, and the framed tangles $(\hcylc,\tanghcyll)$ and $T(\alpha)(\hcylc^{\prime},\tanghcyll^{\prime})T(\beta)$ represent the same framed tangle
 for any braids $\alpha$ and $\beta$ such that 
\begin{itemize}
\item the composition $T(\alpha)(\hcylc^{\prime},\tanghcyll^{\prime})T(\beta)$ is well-defined,  
\item the path $\alpha$ of $\cinjuptd{\finsetb^-}{\CC}$ is homotopic to $\overline{\gamma^-}\gamma^{-\prime}$ in $\ccompuptd{\finsetb^-}{\CC}$,
and \item the path $\beta$ of $\cinjuptd{\finsetb^+}{\CC}$ is homotopic to $\gamma^{+\prime}\overline{\gamma^+}$ in $\ccompuptd{\finsetb^+}{\CC}$ (by homotopies that fix the boundary points).
\end{itemize}

The \emph{domain} of a $q$-tangle (represented by) $T(\gamma^-)(\hcylc,\tanghcyll)T(\gamma^+)$ is (identified with) the domain $\sourcetl$ of $\tanghcyll \colon \sourcetl \hookrightarrow \hcylc$, its \emph{bottom configuration} is $\gamma^-(0)$, and its \emph{top configuration} is $\gamma^+(1)$.

\end{definition}

\begin{examples} \label{exaqtang} A \emph{combinatorial configuration} is a nonassociative word in the letter $\,\bulletmoyen\,$.
The combinatorial $q$-tangles $(\Link,w^-,w^+)$ of Definition~\ref{defcombqtang} are examples of $q$-tangles. They are the $q$-tangles of the standard cylinder whose bottom and top configurations are combinatorial. 
Recall the notation from Example~\ref{exaStasheff}. A combinatorial $q$-tangle $(\tanghcyll,w^-,w^+)$ may be represented as
$T(\gamma^-)(\hcylc,\tanghcyll(1/4))T(\gamma^+)$ 
for a $q$-braid $\gamma^-$ from $\gamma^-(0)=w^-$ to  $\gamma^-(1)=w^-(1/4)$ such that $\gamma^-\left(\left]0,1\right]\right)$ is $w^-\left(\left]0,1/4\right]\right)$,
a $q$-braid $\gamma^+$ from $\gamma^+(0)=w^+(1/4)$ to 
$\gamma^+(1)=w^+$ such that the sets $\gamma^+\left(\left[0,1\right[\right)$ and $w^+\left(\left]0,1/4\right]\right)$ are equal,
and a representative $\tanghcyll(1/4)$ of $\tanghcyll$ from $w^-(1/4)$ to $w^+(1/4)$ in the isotopy class of $\tanghcyll$ of Definition~\ref{defcombqtang}.

We may omit the external pair of parentheses in a combinatorial configuration since it is always present. We may similarly omit the only possible two-point combinatorial configuration $(\bulletpetit \bulletpetit)$ from the notation in combinatorial $q$-tangles.
So we may represent the examples of Example~\ref{excombqtang} by
\begin{equation*}\pcpetit\mbox{, } \assocnonnumwep\mbox{, }
\tancuphighwep
 \mbox{, and }
\tancapwep. \end{equation*}
When the boundary of a $q$-tangle is empty, the $q$-tangle is a framed link in $\crats(\hcylc)$.
Conversely, any asymptotically standard $\QQ$-homology $\RR^3$ equipped with a framed link may be obtained in this way, up to diffeomorphism.
\end{examples}

In addition to the (vertical) product, which extends to $q$-tangles naturally, $q$-tangles support a \indexT{cabling} operation
of a component $K$ of a $q$-tangle $T_m$ (the main tangle) by a $q$-tangle $T_i$ (the inserted tangle), which produces a $q$-tangle $T_m(T_i/K)$. This operation roughly consists in replacing the strand $K$ in $T_m$ by a tangle $T_i$ with respect to the framing of $K$, as in Figure~\ref{figcabling}.

For example, if $T_i$ is the trivial braid $\mid_1\mid_2$ with two strands and if $T_m$ is a combinatorial $q$-tangle, then
$T_m(T_i/K)$ is the tangle $T_m( 2 \times K)$ described before Theorem~\ref{thmdup}. 
As another example, the product $\tanghcyll_1 \otimes \tanghcyll_2$ described before Example~\ref{excombqtangtensor} can be written as
$(\tanghcyll_1 \otimes \mid_2)(\tanghcyll_2/\mid_2)$, where $\tanghcyll_1 \otimes \mid_2=
(\mid_1\mid_2)(\tanghcyll_1/\mid_1)$.

Let us describe the cabling operation in general in (lengthy) detail.\footnote{Notation~\ref{notcabling} ends 4 lines before Definition~\ref{defselflkqtanggen}.}

\begin{notation} \label{notcabling} 

A \emph{semi-pure $q$-tangle} is a $q$-tangle with identical bottom and top configuration (up to dilation and translation). 
A \emph{pure $q$-tangle} is a semi-pure $q$-tangle whose interval components connect a bottom configuration point to the top configuration's corresponding point.
The cabling operation that produces $T_m(T_i/K)$ is defined for any pair $(T_m,T_i)$ of $q$-tangles equipped with an interval component $K$ of $T_m$ going from bottom to top. It is also defined for any pair $(T_m,T_i)$ of $q$-tangles equipped with a framed circle component $K$ of $T_m$ provided that $T_i$ is semi-pure. 

\bfig
\centering
\begin{tikzpicture}
\begin{scope}
 \draw [thin,dashed] (0,0)  .. controls (0,.13) and (.7,.3) .. (1,.3) .. controls (1.3,.3) and (2,.13) ..  (2,0);
\draw [thin] (0,0)  .. controls (0,-.13) and (.7,-.3) .. (1,-.3) .. controls (1.3,-.3) and (2,-.13) ..  (2,0);
\draw [thin,dashed] (0,1.9)  .. controls (0,2.03) and (.7,2.2) .. (1,2.2) .. controls (1.3,2.2) and (2,2.03) ..  (2,1.9);
\draw (1,.95) node{\scriptsize{$T_i$}} (.4,.7) rectangle (1.6,1.2) (.6,0) -- (.6,.7) (.6,1.2) -- (.6,1.9) (1.4,0) -- (1.4,.7) (1.4,1.2) -- (1.4,1.9);
 \draw [draw=white,double=black,very thick] (0,1.9)  .. controls (0,1.77) and (.7,1.6) .. (1,1.6) .. controls (1.3,1.6) and (2,1.77) ..  (2,1.9);
 \draw  [thin] (0,0) -- (0,1.9) (2,0) -- (2,1.9);
 \end{scope}
 
 \begin{scope}[xshift=4.5 cm]
 \draw [thin,dashed] (0,0)  .. controls (0,.13) and (.7,.3) .. (1,.3) .. controls (1.3,.3) and (2,.13) ..  (2,0);
\draw [thin] (0,0)  .. controls (0,-.13) and (.7,-.3) .. (1,-.3) .. controls (1.3,-.3) and (2,-.13) ..  (2,0);
\draw [thin,dashed] (0,1.9)  .. controls (0,2.03) and (.7,2.2) .. (1,2.2) .. controls (1.3,2.2) and (2,2.03) ..  (2,1.9);
\draw (.7,.5) arc (-90:80:.5 and .3);
\draw [->] (.7,.5) arc (-90:0:.5 and .3);
\draw (1.15,.8) node[right]{\scriptsize $K_1$};
\draw (.7,.5) arc (-90:-260:.5 and .3);
\draw [->] (.7,0) -- (.7,.43) (.7,.57) -- (.7,1.9) (.7,.57) -- (.7,1.3);
\draw (.65,1.3) node[right]{\scriptsize $K_2$};
\draw (1.2,0) arc (180:0:.2 and .4);
\draw [->] (1.2,0) arc (180:50:.2 and .4);
\draw (1.45,.35) node[right]{\scriptsize $K_3$};
\draw [draw=white,double=black,very thick] (0,1.9)  .. controls (0,1.77) and (.7,1.6) .. (1,1.6) .. controls (1.3,1.6) and (2,1.77) ..  (2,1.9);
 \draw  [thin] (0,0) -- (0,1.9) (2,0) -- (2,1.9);
 \draw (0,.95) node[left]{\scriptsize $T_m=$};
 \end{scope}
 
 \begin{scope}[xshift=2 cm,yshift=-3.2 cm]
 \draw [thin,dashed] (0,0)  .. controls (0,.13) and (.7,.3) .. (1,.3) .. controls (1.3,.3) and (2,.13) ..  (2,0);
\draw [thin] (0,0)  .. controls (0,-.13) and (.7,-.3) .. (1,-.3) .. controls (1.3,-.3) and (2,-.13) ..  (2,0);
\draw [thin,dashed] (0,1.9)  .. controls (0,2.03) and (.7,2.2) .. (1,2.2) .. controls (1.3,2.2) and (2,2.03) ..  (2,1.9);
\draw (.7,.5) arc (-90:80:.55 and .35);
\draw (.7,.5) arc (-90:-260:.55 and .35);
\draw (.7,.6) arc (-90:80:.45 and .25);
\draw (.7,.6) arc (-90:-260:.45 and .25);
\draw [->] (.7,0) -- (.7,.43) (.7,.54) -- (.7,.56) (.7,.67) -- (.7,1.9) (.7,.9) -- (.7,1.45);
\fill[white,draw=black] (.9,.7) rectangle (1.35,1.05) ;
\draw (1.15,.85) node{\scriptsize $T_i$};
\draw (.65,1.35) node[right]{\scriptsize $K_2$};
\draw (1.2,0) arc (180:0:.2 and .4);
\draw [->] (1.2,0) arc (180:50:.2 and .4);
\draw (1.45,.35) node[right]{\scriptsize $K_3$};
\draw [draw=white,double=black,very thick] (0,1.9)  .. controls (0,1.77) and (.7,1.6) .. (1,1.6) .. controls (1.3,1.6) and (2,1.77) ..  (2,1.9);
 \draw  [thin] (0,0) -- (0,1.9) (2,0) -- (2,1.9);
 \draw (0,.95) node[left]{\scriptsize $T_m\left(\frac{T_i}{K_1}\right)=$};
 \end{scope}
 
  \begin{scope}[xshift=7 cm,yshift=-3.2cm]
 \draw [thin,dashed] (0,0)  .. controls (0,.13) and (.7,.3) .. (1,.3) .. controls (1.3,.3) and (2,.13) ..  (2,0);
\draw [thin,dashed] (0,1.9)  .. controls (0,2.03) and (.7,2.2) .. (1,2.2) .. controls (1.3,2.2) and (2,2.03) ..  (2,1.9);
\draw (.7,.45) arc (-90:75:.5 and .3);
\draw [->] (.7,.45) arc (-90:0:.5 and .3);
\draw (1.15,.8) node[right]{\scriptsize $K_1$};
\draw (.7,.45) arc (-90:-255:.5 and .3);
\draw (.68,1.05) -- (.72,1.05);
\draw (.7,-.45) .. controls (.65, -.2) .. (.65,0) -- (.65,.4) (.65,.5) -- (.65,1.9) .. controls (.65,2) .. (.7,2.3);
\draw (.7,-.45) .. controls (.75, -.2) .. (.75,0) -- (.75,.4) (.75,.5) -- (.75,1.9) .. controls (.75,2) .. (.7,2.3);
\fill[white,draw=black] (.5,1.17) rectangle (.9,1.52);
\draw (.7,1.32) node{\scriptsize $T_i$};
\draw (1.6,-.45) -- (1.6,0) (1.2,-.45) -- (1.2,0) arc (180:0:.2 and .4);
\draw [->] (1.2,0) arc (180:50:.2 and .4);
\draw (1.45,.35) node[right]{\scriptsize $K_3$};
\draw [draw=white,double=black,very thick] (0,1.9)  .. controls (0,1.77) and (.7,1.6) .. (1,1.6) .. controls (1.3,1.6) and (2,1.77) ..  (2,1.9);
\draw [draw=white,double=black,very thick] (0,0)  .. controls (0,-.13) and (.7,-.3) .. (1,-.3) .. controls (1.3,-.3) and (2,-.13) ..  (2,0);
\draw  [thin] (0,0) -- (0,1.9) (2,0) -- (2,1.9);
\draw (0,.95) node[left]{\scriptsize $T_m\left(\frac{T_i}{K_2}\right)=$};
 \end{scope}
 
 \end{tikzpicture}
\caption{Examples of cablings}
\label{figcabling}

\end{figure}

We begin with the details when $K$ is a circle, because they are lighter in this case.
When $K$ is a closed component and $T_i$ is a framed tangle $(\hcylc_i, \tanghcyll_i)$ with identical injective bottom and top configuration, pick a tubular neighborhood $D^2 \times K$ of $K$, trivialized with respect to the parallelization of $K$, that does not meet the other components of $T_m$. Write this neighborhood as $D^2 \times \left[0,1\right]/(0 \sim 1)$, and replace it by $(\hcylc_i, \tanghcyll_i)$ using the identification of $N\left(\partial\left( \drad{1} \times \left[0,1\right]\right)\right)$ with a neighborhood of $\partial \left(\drad{1} \times \left[0,1\right]\right)$ in order to obtain $T_m(T_i/K)$. Note that when $\gamma$ is a braid such that $\gamma(1)$ is the bottom configuration of $T_i$, we have
\begin{equation*}T_m\left(\frac{T_i}{K}\right)=T_m\left(\frac{T(\gamma)T_iT(\overline{\gamma})}{K}\right).\end{equation*}
Any semi-pure $q$-tangle $T_q$ can be written as $T(\gamma)(\hcylc_i,\tanghcyll_i)T(\overline{\gamma})$, for some $q$-braid  $\gamma$ and some framed tangle $(\hcylc_i,\tanghcyll_i)$. For such a tangle, set $T_m(T_q/K)=T_m((\hcylc_i,\tanghcyll_i)/K).$
It is easy to check that this definition is consistent.

Let us now define \emph{cabling} or \emph{duplication} for configurations.
Let $\finsetb$ and $\finsetb_i$ be nonempty finite sets, let $\eltb \in \finsetb$, and let $\finsetb(\finsetb_i/\eltb) = (\finsetb \setminus \{\eltb \}) \cup \finsetb_i$.
Let $c_m$ be an element of $\ccompuptd{\finsetb}{\CC}$, and let $c_i$ be an element of $\ccompuptd{\finsetb_i}{\CC}$.
The configuration $c_m(c_i/\eltb)$ of $\ccompuptd{\finsetb(\finsetb_i/\eltb)}{\CC}$ is the configuration obtained by letting $c_i$ replace $\eltb$. 
Up to translation, there is only one configuration of a set of one element. So $c_m(c_i/\eltb)$ is this unique configuration
if $\cardlef{\finsetb(\finsetb_i/\eltb)}=1$.
If $\cardlef{\finsetb_i}=1$, then $c_m(c_i/\eltb)=c_m$. If $\cardlef{\finsetb}=1$, then $c_m(c_i/\eltb)=c_i$ (with natural identifications).
Assume that $\cardlef{\finsetb_i} \geq 2$ and $\cardlef{\finsetb} \geq 2$. When $c_m$ and $c_i$ are both combinatorial configurations, it makes natural sense to let $c_i$ replace $\eltb$ to produce $c_m(c_i/\eltb)$.
In general, recall Definition~\ref{defparenthesizationtwo} and define the $\Delta$-parenthesization $\tau(c_m)(\tau(c_i)/\eltb)$ of $\finsetb(\finsetb_i/\eltb)$ from the respective $\Delta$-parenthesizations $\tau(c_m)$ and $\tau(c_i)$ of $c_m$ and $c_i$ by the
following one-to-one correspondence \begin{equation*}\begin{array}{llll}\phi \colon &\tau(c_m) \sqcup \tau(c_i) &\to &\tau(c_m)(\tau(c_i)/\eltb)\\
                                      & A &\mapsto &\left\{\begin{array}{ll} A 
\;&\mbox{if}\; A \in \tau(c_i) \;\mbox{or}\; (A \in \tau(c_m) \;\mbox{and}\; \eltb \notin A)\\
A(\finsetb_i/\eltb) \;&\mbox{if}\; A \in \tau(c_m) \;\mbox{and}\; \eltb \in A.
\end{array}\right.
\end{array}\end{equation*} 
With the notation of Theorem~\ref{thmcompuptd}, the configuration $c_m(c_i/\eltb)$ is the configuration of $\ccompuptd{\finsetb(\finsetb_i/\eltb),\tau(c_m)(\tau(c_i)/\eltb)}{\CC}$ that restricts to $\finsetb_i$ as $c_i$, and to $\finsetb(\{\eltb^{\prime}\}/\eltb)$ as $c_m$ for any element $\eltb^{\prime}$ of $\finsetb_i$.

Assume that $K$ is an interval component going from bottom to top of a $q$-tangle \begin{equation*}T_m=T(\gamma_m^-)(\hcylc_m,\tanghcyll_m)T(\gamma_m^+),\end{equation*} and define
$T_m(T_i/K)$ for a $q$-tangle $T_i=T(\gamma_i^-)(\hcylc_i,\tanghcyll_i)T(\gamma_i^+)$. 
Let $\finsetb_i^-$, (resp. $\finsetb_i^+$, $\finsetb^-$, $\finsetb^+$) denote the set of strand indices of $\gamma_i^-$ (resp. $\gamma_i^+$, $\gamma_m^-$, $\gamma_m^+$). Let $b^-_K$ (resp. $b^+_K$) denote the strand index of $K$ in $\finsetb^-$ (resp. $\finsetb^+$). Let $c_m^-$ (resp. $c_i^-$) denote the bottom configuration of $T_m$ (resp. $T_i$). Let $c_m^+$ (resp. $c_i^+$) denote the top configuration of $T_m$ (resp. $T_i$). Assume that $\gamma_m^-\left(\left]0,1\right]\right) \subset \cinjuptd{\finsetb^-}{\CC}$, $\gamma_i^-\left(\left]0,1\right]\right) \subset \cinjuptd{\finsetb_i^-}{\CC}$, $\gamma_m^+\left(\left[0,1\right[\right) \subset \cinjuptd{\finsetb^+}{\CC}$, and $\gamma_i^+\left(\left[0,1\right[\right) \subset \cinjuptd{\finsetb_i^+}{\CC}$.
Let $I_K$ denote the intersection of $K$ with $\hcylc_m$. Identify $I_K$ with $\left[0,1\right]$. Let $D^{(i)}$ be a copy of the disk $D^2$. Let
$D^{(i)} \times \left[0,1\right]$ be a tubular neighborhood of $I_K$ in $\hcylc_m$ that does not meet the other components of $\tanghcyll_m$ and that meets
$\partial \hcylc_m$ along $D^{(i)} \times \partial \left[0,1\right]$ inside $\drad{1} \times \partial\left[0,1\right]$, such that $\left(\{\pm 1\} \times \left[0,1\right]\right) \subset \left(\partial D^{(i)} \times  \left[0,1\right]\right) $ is the given parallel of $I_K$. Replace $D^{(i)} \times \left[0,1\right]$ by $(\hcylc_i,\tanghcyll_i)$ in order to get a tangle $(\hcylc_m,\tanghcyll_m)\left((\hcylc_i,\tanghcyll_i)/I_K\right)$.

Let $\gamma_m^-(\gamma_i^-/K)$ be the path composition $\gamma_m^-(c_i^-/K)\left(\gamma_m^-(1)(\gamma_i^-/K)\right)$ of the paths $\gamma_m^-(c_i^-/K)$ and $\gamma_m^-(1)(\gamma_i^-/K)$ in $\ccompuptd{\finsetb^-(\finsetb_i^-/b^-_K)}{\CC}$,
where $\gamma_m^-(c_i^-/K)(t)=\gamma_m^-(t)(c_i^-/b^-_K)$ for any $t\in \left[0,1\right]$,
and the restriction to $\left]0,1\right]$ of $\gamma_m^-(1)(\gamma_i^-/K)$ is represented by a map from $\left]0,1\right]$ to  $\check{C}_{\finsetb^-(\finsetb_i^-/b^-_K)}[\dorad{1}]$, 
\begin{itemize}
 \item which maps $1$ to the bottom configuration of $(\hcylc_m,\tanghcyll_m)\bigl((\hcylc_i,\tanghcyll_i)/I_K\bigr)$,
\item whose restriction to $\finsetb^- \setminus \{b^-_K\}$ is constant,\footnote{It is thus located in $\check{C}_{\finsetb^- \setminus \{b^-_K\}}\left[\dorad{1} \setminus p_{\CC}(D^{(i)}\times \{0\})\right]$.} and,
\item whose restriction to $\finsetb_i^-$ is a lift of $\gamma^-_{i}\vert_{\left]0,1\right]}$ in $\check{C}_{\finsetb_i^-}[p_{\CC}(D^{(i)}\times \{0\})]$, such that 
$\gamma_m^-(1)(\gamma_i^-/K)$ is composable by $\gamma_m^-(c_i^-/K)$ on its left.
\end{itemize}

The tangle $T_m(T_i/K)$ is defined to be
\begin{equation*}T\bigl(\gamma_m^-(\gamma_i^-/K)\bigr) \Bigl((\hcylc_m,\tanghcyll_m)\bigl((\hcylc_i,\tanghcyll_i)/I_K\bigr)\Bigr)T\bigl(\gamma_m^+(\gamma_i^+/K)\bigr),\end{equation*}
where the definition of $T(\gamma_m^+(\gamma_i^+/K))$, similar to that of $T(\gamma_m^-(\gamma_i^-/K))$, follows (and can be skipped...).

The path $\gamma_m^+(\gamma_i^+/K)$ 
is the path composition $\left(\gamma_m^+(0)(\gamma_i^+/K)\right)\gamma_m^+(c_i^+/K)$ in $\ccompuptd{\finsetb^+(\finsetb_i^+/b^+_K)}{\CC}$,
where $\gamma_m^+(c_i^+/K)(t)=\gamma_m^+(t)(c_i^+/b^+_K)$ for any $t \in \left[0,1\right]$,
and the restriction to $\left[0,1\right[$ of $\gamma_m^+(0)(\gamma_i^+/K)$ is represented by a map from $\left[0,1\right[$ to  $\check{C}_{\finsetb^+(\finsetb_i^+/b^+_K)}[\dorad{1}]$, mapping $0$ to the top configuration of $(\hcylc_m,\tanghcyll_m)\left((\hcylc_i,\tanghcyll_i)/I_K\right)$,
whose restriction to $\finsetb^+ \setminus \{b^+_K\}$ is constant, and whose restriction to $\finsetb_i^+$ is a lift of $\gamma^+_{i}\vert_{\left[0,1\right[}$ in $\check{C}_{\finsetb_i^+}\bigl[p_{\CC}(D^{(i)}\times \{1\})\bigr]$, such that $\gamma_m^+(0)(\gamma_i^+/K)$ is composable by $\gamma_m^+(c_i^+/K)$ on its right.

A particular case of cablings is the case in which the inserted $q$-tangle $T_i$ is just the $q$-tangle $\confy \times \left[0,1\right]$ associated to the constant path of $\ccompuptd{\finsetb_i}{\CC}$ that maps $\left[0,1\right]$ to a configuration $\confy$ of $\ccompuptd{\finsetb_i}{\CC}$.
(Formally, this $q$-tangle is represented by $T(\gamma) \left(\gamma(1) \times \left[0,1\right]\right) T(\overline{\gamma})$ for some path $\gamma$ of $\ccompuptd{\finsetb_i}{\CC}$ such that $\gamma(0)=\confy$ and $\gamma(1) \in \cinjuptd{\finsetb_i}{\CC}$.)
Set 
\begin{equation*}T_m(\confy \times K)= T_m\bigl((\confy \times \left[0,1\right])/K\bigr).\end{equation*}
If $K$ is a closed component, then $T_m(\confy \times K)$ depends only on $\cardlef{\finsetb_i}$. It is denoted by $T_m(\cardlef{\finsetb_i}\times K)$.
If $K$ is an interval component and $\confy$ is the unique configuration of $\ccompuptd{\{1,2\}}{\RR}$,
then $T_m(\confy \times K)$ is again denoted by $T_m(2 \times K)$.
These special cablings introduced before Theorem~\ref{thmdup} are called \emph{duplications} or \emph{doublings}. Our functor $\Zinvufrf$ will behave well under all cablings. 
\end{notation}

We end this section by completing Definition~\ref{defselflkqtangtwo} of the self-linking number for a framed $q$-tangle component, which goes from bottom to bottom or from top to top.

\begin{definition}
 \label{defselflkqtanggen}
The \emph{self-linking number} $lk(K,K_{\parallel})$ of a framed component $K$ of a $q$-tangle going from bottom to bottom (resp. from top to top) is defined as follows.
The self-linking number depends only on the component. So there is no loss of generality in representing 
$K$ by a tangle representative with injective bottom (resp. top) configuration whose ends are at a distance bigger than $2\varepsilon$ for a small positive
$\varepsilon$. Also assume that ${K}_{\parallel}(0)= K(0) + (\varepsilon,0,0)$ and ${K}_{\parallel}(1)= K(1) - (\varepsilon,0,0)$. (There is no loss of generality in this assumption either. It suffices to choose a parallel that satisfies this assumption to define the self-linking number. Recall Definition~\ref{defparalleltang}.) Let $\hat{K} = K\left(\left[0,1\right]\right) \cup \gamma\left(\left[0,1\right]\right)$ for an arbitrary path $\gamma$ from $\gamma(0)=K(1)$ to $\gamma(1)=K(0)$ such that 
$\gamma\left(\left]0,1\right[\right) \subset \crats(\hcylc) \setminus \hcylc$.

Let $\left[K(1),K(0)\right]$ denote the straight segment from $K(1)$ to $K(0)$ in $\drad{1}$. 
Let $\alpha_1\colon \left[0,1\right] \to \drad{1}$ be an arc from $K_{\parallel}(1)$ to a point $a_1$ inside $\left[K(1),K(0)\right]$
such that $\alpha_1(t)=K(1)-\varepsilon\exp(2i\pi\theta_1 t)$ for some real number $\theta_1$. Similarly, let $\alpha_0 \colon \left[0,1\right] \to \drad{1}$ be an arc from $K_{\parallel}(0)$ to a point $a_0$ inside $\left[K(1),K(0)\right]$
such that $\alpha_0(t)=K(0)+\varepsilon\exp(2i\pi\theta_0 t)$ for some $\theta_0 \in \RR$.

\bfig
\centering
\begin{tikzpicture}
\useasboundingbox (-.1,-.1) rectangle (8.5,2);
\draw [->] (28:1.5) -- (28:3) node[right]{\scriptsize $K(1)$} (.8,0) -- (0,0) -- (28:1.5);
\draw [->] (14:.8)  arc (14:28:.8) (0:.8)  arc (0:14:.8); 
\draw [-<] (28:2.2) arc (208:180:.8) -- (28:3) (28:2.2) arc (208:194:.8);
\draw (35:2.05) node{\scriptsize $\alpha_1$} (10:1.05) node{\scriptsize $\alpha_0$} (2.2,.8) node[below]{\scriptsize $\left[K(0),K(1)\right]$}
(-.5,0) node{\scriptsize $K(0)$};
\fill (28:3) circle (1pt) (0,0) circle (1pt);
\begin{scope}[xshift=5.5cm]
\draw [->, very thick] (1,0) -- (2,0) .. controls (2,1) and (1.4,1.6) .. (1,1.6);
\draw [->, very thick] (1,1.6) .. controls (.6,1.6) and (0,1) .. (0,0) -- (1,0);
\draw [draw=white,double=black,very thick] (-.3,.4) .. controls (-.3,.3) and (-.15,.2) .. (0,.2) .. controls (.15,.2) and (.3,.3) .. (.3,.4) (1.7,.4) .. controls (1.7,.3) and (1.85,.2) .. (2,.2) .. controls (2.15,.2) and (2.3,.3) ..  (2.3,.4);
\draw [->] (-.3,.4) .. controls (-.3,.3) and (-.15,.2) .. (0,.2);
\draw [->] (2.3,.4) .. controls (2.3,.3) and (2.15,.2) .. (2,.2);
\draw [->] (1,.4) -- (1.7,.4) (2.3,.4) .. controls (2.3,1.3) and (1.55,1.9) .. (1,1.9);
\draw [->]  (1,1.9) .. controls (.45,1.9) and (-.3,1.3) .. (-.3,.4) (.3,.4) -- (1,.4);
\fill (2,0) circle (1pt) (0,0) circle (1pt);
\draw (1,1.6) node[below]{\scriptsize $\hat{K}$} (1,.3) node[above]{\scriptsize $\hat{K}_{\parallel,-\frac12,\frac12}$} (-.2,.05) node{\scriptsize $\alpha_1$} (2.25,.05) node{\scriptsize $\alpha_0$} (1,0) node[below]{\scriptsize $\gamma$}; 
\end{scope}

\end{tikzpicture}
\caption{A general picture of $\alpha_0$ and $\alpha_1$, and a picture of $\hat{K}_{\parallel}$ when $\theta_0=-\frac12$, $\theta_1=\frac12$, and $K$ goes from bottom to bottom}
\label{figgammapar}

\end{figure}
If $K$ goes from bottom to bottom, define
\begin{equation*}\hat{K}_{\parallel,\theta_0,\theta_1} =K_{\parallel}\bigl(\left[0,1\right]\bigr) \cup \bigl((\alpha_1 \cup \left[a_1,a_0\right] \cup \overline{\alpha}_0) \times \{0\}\bigr),\end{equation*}
and set \begin{equation*}lk(K,K_{\parallel})=lk(\hat{K},\hat{K}_{\parallel,\theta_0,\theta_1}) +\theta_1 +\theta_0.\end{equation*}
 If $K$ goes from top to top, 
define
\begin{equation*}\hat{K}_{\parallel,\theta_0,\theta_1} =K_{\parallel}\bigl(\left[0,1\right]\bigr) \cup \bigl((\alpha_1 \cup \left[a_1,a_0\right] \cup \overline{\alpha}_0) \times \{1\}\bigr)\end{equation*}
 and set \begin{equation*}lk(K,K_{\parallel})=lk(\hat{K},\hat{K}_{\parallel,\theta_0,\theta_1}) -(\theta_1 +\theta_0).\end{equation*}
\end{definition}

Note that these definitions do not depend on the chosen $\theta_1$ and $\theta_0$, which are well-determined modulo $\ZZ$. The angles $2\pi\theta_1$ and $2\pi\theta_0$ are both congruent to the angle from the oriented real line to $\bigvec{K(0)K(1)}$ modulo $2\pi$.
So $(lk(K,K_{\parallel})-2\theta_1)\in \QQ$ when $K$ goes from bottom to bottom, and $(lk(K,K_{\parallel})+2\theta_1)\in \QQ$ when $K$ goes from top to top.

When $\bigvec{K(0)K(1)}$ directs and orients the real line, we can choose $\theta_0=\theta_1=0$, and this definition coincides with Definition~\ref{defselflkqtangtwo}. When $\bigvec{K(1)K(0)}$ directs and orients the real line, we can choose $\theta_1=\frac12=-\theta_0$ so that $\hat{K}_{\parallel}$ is simply as in Figure~\ref{figgammapar}. So the present definition is again consistent with Definition~\ref{defselflkqtangtwo}.

\begin{lemma}
\label{lemlkindepor}
 The self-linking number does not depend on the orientations of the components.
\end{lemma}
\bp
This is easy to see for closed components, and this is part of the definition for components going from bottom to top or from top to bottom. When $K$ goes from top to top, let $K^{\prime}$ stand for $(-K)$, and let $K^{\prime}_{\parallel}$ be the parallel obtained from $K_{\parallel}$ by a rotation of angle $\pi$ around $K$ (and by reversing the orientation). Choose the corresponding angles $\theta^{\prime}_0$ and $\theta^{\prime}_1$ to be $\theta
^{\prime}_0=\theta_1-\frac12$ and $\theta^{\prime}_1=\theta_0+\frac12$. So
$\hat{K}^{\prime}_{\parallel,\theta^{\prime}_0,\theta^{\prime}_1}$ is isotopic to $(-\hat{K}_{\parallel,\theta_0,\theta_1})$ in the complement of $\hat{K}$.
See Figure~\ref{figgammaparttt}.
\bfig
\centering
\begin{tikzpicture}
\useasboundingbox (-1,-.9) rectangle (9,1.5);
\begin{scope}[yshift=-.4cm]
\draw [dashed, -<] (28:2.2) arc (-152:-180:.8) -- (28:3) (28:2.2) arc (-152:-166:.8);
\draw [-<] (28:1.5) -- (28:3) (-.8,0) -- (0,0) -- (28:1.5);
\draw [->] (-76:.8)  arc (-76:28:.8) (-180:.8)  arc (-180:-76:.8); 
\draw [-<] (28:2.2) arc (-152:0:.8) -- (28:3) (28:2.2) arc (-152:-76:.8);
\draw (34:2) node{\scriptsize $\theta_1$} (17:3) node{\scriptsize $\theta^{\prime}_0$} (-76:.6) node{\scriptsize $\theta^{\prime}_1$} (1.9,.6) node[below]{\scriptsize $\left[K^{\prime}(0),K^{\prime}(1)\right]$}
(-.3,.2) node{\scriptsize $K^{\prime}(1)$} (-.8,0) node[left]{\scriptsize $\hat{K}^{\prime}_{\parallel}(1)$} (3.44,1.41) node[right]{\scriptsize $\hat{K}^{\prime}_{\parallel}(0)$} (26:3) node[above]{\scriptsize $K^{\prime}(0)=K(1)$};
\fill (28:3) circle (1pt) (0,0) circle (1pt) (-.8,0) circle (1pt) (3.44,1.41) circle (1pt);
\end{scope}

\begin{scope}[xshift=6.5cm, yshift=-1cm]
\draw [->, very thick] (1,2) -- (2,2) .. controls (2,1) and (1.4,.4) .. (1,.4);
\draw [->, very thick] (1,.4) .. controls (.6,.4) and (0,1) .. (0,2) -- (1,2);
\draw [draw=white,double=black,very thick] (-.3,1.6) .. controls (-.3,1.5) and (-.15,1.4) .. (0,1.4) .. controls (.15,1.4) and (.3,1.5) .. (.3,1.6) (1.7,1.6) .. controls (1.7,1.5) and (1.85,1.4) .. (2,1.4) .. controls (2.15,1.4) and (2.3,1.5) ..  (2.3,1.6);
\draw [->] (-.3,1.6) .. controls (-.3,1.5) and (-.15,1.4) .. (0,1.4);
\draw [->] (2.3,1.6) .. controls (2.3,1.5) and (2.15,1.4) .. (2,1.4);
\draw [->] (1,1.6) -- (1.7,1.6) (2.3,1.6) .. controls (2.3,.7) and (1.55,.1) .. (1,.1);
\draw [->]  (1,.1) .. controls (.45,.1) and (-.3,.7) .. (-.3,1.6) (.3,1.6) -- (1,1.6);
\fill (2,2) circle (1pt) (0,2) circle (1pt);
\draw (1,.4) node[above]{\scriptsize $K^{\prime}$} (1.13,1.65) node[below]{\scriptsize $\hat{K}^{\prime}_{\parallel,-\frac12,\frac12}$} (-.2,1.75) node{\scriptsize $\theta^{\prime}_1$} (2.25,1.75) node{\scriptsize $\theta^{\prime}_0$} (1,2) node[above]{\scriptsize $\hat{K}^{\prime}$}; 
\end{scope}
\end{tikzpicture}
\caption{A general picture of $\theta^{\prime}_0$ and $\theta^{\prime}_1$, and a picture of $\hat{K}^{\prime}_{\parallel}$ when $\theta^{\prime}_0=-\frac12$, $\theta^{\prime}_1=\frac12$, and $K$ goes from top to top}
\label{figgammaparttt}

\end{figure}
\eop

Proposition~\ref{proplktangles} will show that the real-valued self-linking numbers $lk(K,K_{\parallel})$ coincide with $I_{\theta}(K,\tau)$ for the interval components $K$ of the straight tangles of Definition~\ref{defstraighttang}. Its proof relies only on the beginning of Section~\ref{secstraighttang}.

\section{Definition of \texorpdfstring{$\Zinvufrf$}{Zf} for all  \texorpdfstring{$q$-tangles}{q--tangles}}
\label{secdefzinvf}

Recall Definition~\ref{deffstconsttangframed} of the extension of $\Zinvufrf$ for $q$-tangles whose bottom and top configurations are injective. 
In this section, we extend the definition of $\Zinvufrf$ to all $q$-tangles.

In Chapter~\ref{chapzinvfbraid}, we will prove the following particular case
of the functoriality property stated in Theorem~\ref{thmfunc}. It is a direct corollary of Proposition~\ref{propfuncbraid}.

\begin{proposition} \label{propfuncbraidprel}
Let  $(\hcylc_1,\tanghcyll_1)$ and $(\hcylc_2,\tanghcyll_2)$ be two framed tangle representatives such that the bottom of $\tanghcyll_2$ coincides with the top of $\tanghcyll_1$. If one of them is a braid, then we have \begin{equation*}\Zinvufrfneg\bigl(\hcylc_1\hcylc_2,(\tanghcyll_1\tanghcyll_2)_{\parallel}\bigr) = \Zinvufrfneg\bigl(\hcylc_1,\tanghcyll_1,\tanghcyll_{1\parallel}\bigr) \Zinvufrfneg\bigl(\hcylc_2,\tanghcyll_2,\tanghcyll_{2\parallel}\bigr).\end{equation*}
\end{proposition}

The following lemma allows us to consider the tangles' bottom and top configurations up to translation and dilation. It will also be a direct corollary of Proposition~\ref{propfuncbraid}.

\begin{lemma}
\label{lemconnuptd}
Let $\gamma \colon \left[0,1\right] \to \check{C}_{\finsetb}[\dorad{1}]$. Let $p_{CS}\circ \gamma \colon \left[0,1\right] \to \cinjuptd{\finsetb}{\CC}$ be the composition of $\gamma$ by the natural projection $p_{CS} \colon \check{C}_{\finsetb}[\dorad{1}]  \to \cinjuptd{\finsetb}{\CC}$ (which mods out by translations and dilations). Then $\Zinvuf(\gamma)$ and $\Zinvufrfneg(\gamma)$ depend only on $p_{CS}\circ \gamma$. 
\end{lemma}

Under the assumptions of the lemma, we set $\Zinvuf(p_{CS}\circ \gamma)=\Zinvuf(\gamma)$ and $\Zinvufrfneg(p_{CS}\circ \gamma)=\Zinvufrfneg(\gamma)$. Recall that $\Zinvuf$ and $\Zinvufrf$ coincide for braids.
We extend the definition of $\Zinvuf$ to piecewise smooth paths of $\cinjuptd{\finsetb}{\CC}$ so that $\Zinvuf$ is multiplicative with respect to path composition of smooth paths.
The following theorem is essentially due to Sylvain Poirier \cite{poirierv2}. It allows us to extend $\Zinvufrf$ to $q$-braids, and to $q$-tangles in Definition~\ref{defzqtang}.

\begin{theorem}
\label{thmpoirbraid}
Let $p_{CS}\circ \gamma \colon \left[0,1\right] \to \ccompuptd{\finsetb}{\CC}$ be a path whose restriction to $\left]0,1\right[$
is the projection of some $\gamma \colon \left]0,1\right[ \to \check{C}_{\finsetb}[\dorad{1}]$ that can be described by a collection of piecewise polynomial continuous maps $(\gamma_{\eltb} \colon \left[0,1\right] \to \dorad{1})_{\eltb \in \finsetb}$.\footnote{Every $\gamma_{\eltb}$ is polynomial over a finite number of intervals covering $\left[0,1\right]$.}
Then $\lim_{\varepsilon \to 0} \Zinvuf(p_{CS}\circ \gamma\vert_{\left[\varepsilon, 1-\varepsilon\right]})$ makes sense, and it depends only on the homotopy class of $p_{CS}\circ \gamma$ relatively to its boundary. It is denoted by $\Zinvuf(p_{CS}\circ \gamma)$ or $\Zinvufrfneg(p_{CS}\circ \gamma)$.
\end{theorem}
Theorem~\ref{thmpoirbraid} and Theorem~\ref{thmconnecgen}, which generalizes its homotopy invariance part, will be proved in Section~\ref{zinvfqbraid}.

\begin{proposition}\label{propfuncqbraid}
Any $q$-braid $\gamma \colon \left[0,1\right] \to \ccompuptd{\finsetb}{\CC}$ is homotopic relatively to its boundary to a $q$-braid $p_{CS}\circ \tilde{\gamma}$ as in the statement of Theorem~\ref{thmpoirbraid}.
Setting $\Zinvuf(\gamma)=\Zinvuf(p_{CS}\circ \tilde{\gamma})$ consistently extends the definition of $\Zinvuf$ to all $q$-braids.
Furthermore, $\Zinvuf$ is multiplicative with respect to the $q$-braid composition: For two composable paths $\gamma_1$ and $\gamma_2$ of $\ccompuptd{\finsetb}{\CC}$, we have
\begin{equation*}\Zinvuf(\gamma_1\gamma_2)=\Zinvuf(\gamma_1)\Zinvuf(\gamma_2).\end{equation*}
\end{proposition}
\bp Let us first exhibit a $q$-braid $p_{CS}\circ \tilde{\gamma}$ homotopic to a given $q$-braid $\gamma \colon \left[0,1\right] \to \ccompuptd{\finsetb}{\CC}$, with the desired properties.
Define a path $\tilde{\gamma}_1 \colon \left[0,1/3\right] \to {C}_{\finsetb}[\dorad{1}]$, such that $p_{CS}\circ\tilde{\gamma}_1(0)=\gamma(0)$, $\tilde{\gamma}_1(\left]0,1/3\right]) \subset \check{C}_{\finsetb}[\dorad{1}]$, and $\tilde{\gamma}_1$ is a path obtained by replacing all the parameters $\mu_{\finseta}$ in the charts of Lemma~\ref{lemneighchartcomf} by $\varepsilon t$ for $t \in \left[0,1/3\right]$ for some small $\varepsilon > 0$. So $\tilde{\gamma}_1$ is described by a collection of polynomial maps $(\tilde{\gamma}_{1,\eltb} \colon \left[0,1/3\right] \to \dorad{1})_{\eltb \in \finsetb}$. 
Similarly define a polynomial path $\tilde{\gamma}_3 \colon \left[2/3,1\right] \to {C}_{\finsetb}[\dorad{1}]$ such that $p_{CS}\circ\tilde{\gamma}_3(1)=\gamma(1)$ and $\tilde{\gamma}_3(\left[2/3,1\right[) \subset \check{C}_{\finsetb}[\dorad{1}]$. 
Define a path $\tilde{\gamma}_2^{\prime} \colon \left[1/3,2/3\right] \to {C}_{\finsetb}[\dorad{1}]$, such that
 $\tilde{\gamma}_2^{\prime}(1/3)=\tilde{\gamma}_1(1/3)$, $\tilde{\gamma}_2^{\prime}(2/3)=\tilde{\gamma}_3(2/3)$, and $p_{CS} \circ \tilde{\gamma}_2^{\prime}$ is a path composition $(p_{CS} \circ \overline{\tilde{\gamma}_1})\gamma( p_{CS} \circ \overline{\tilde{\gamma}_3})$. 
This path $\tilde{\gamma}_2^{\prime}$ of ${C}_{\finsetb}[\dorad{1}]$ is homotopic to a path of the interior $\check{C}_{\finsetb}[\dorad{1}]$ of the manifold ${C}_{\finsetb}\left[\drad{1}\right]$ with ridges, and it is homotopic to a polynomial path $\tilde{\gamma}_2$ in $\check{C}_{\finsetb}[\dorad{1}]$. Now, the path composition $\tilde{\gamma}=\tilde{\gamma}_1\tilde{\gamma}_2\tilde{\gamma}_3$ satisfies the hypotheses of Theorem~\ref{thmpoirbraid}, and $p_{CS}\circ \tilde{\gamma}$ is homotopic to $\gamma$, relatively to its boundary. Furthermore, for any other path $\tilde{\gamma}^{\prime}$ satisfying these properties, 
 $p_{CS}\circ \tilde{\gamma}$ and $p_{CS}\circ \tilde{\gamma}^{\prime}$ are homotopic relatively to the boundary. So the definition of $\Zinvuf(\gamma)$ is consistent.
 
Let us prove the multiplicativity.
Pick a piecewise polynomial path
$\tilde{\gamma}_1 \colon [0,1] \to {C}_{\finsetb}[\dorad{1}]$, such that  $\gamma_1$ and $p_{CS} \circ \tilde{\gamma}_1$ are homotopic relatively to the boundary, and $\tilde{\gamma}_1(]0,1[) \subset \check{C}_{\finsetb}[\dorad{1}]$. Next, pick a piecewise polynomial path
$\tilde{\gamma}_2 \colon [0,1] \to {C}_{\finsetb}[\dorad{1}]$, such that, 
$\tilde{\gamma}_{2}(t)=\tilde{\gamma}_{1}(1-t)$ for any $t \in [0,1/2]$,
$\gamma_2$ and $p_{CS} \circ \tilde{\gamma}_2$ are homotopic relatively to the boundary, and $\tilde{\gamma}_2(]0,1[) \subset \check{C}_{\finsetb}[\dorad{1}]$.
Thus we have \begin{equation*}\Zinvuf(\gamma_1\gamma_2)=\lim_{\varepsilon \to 0} \Zinvuf\Bigl(p_{CS}\circ \bigl(\tilde{\gamma}_{1}\vert_{\left[\varepsilon, 1/2\right]}\tilde{\gamma}_{2}\vert_{\left[1/2, 1 -\varepsilon\right]}\bigr)\Bigr).\end{equation*}
We also have \begin{equation*}\begin{array}{lll}\Zinvuf(\gamma_1)\Zinvuf(\gamma_2)&=&\lim_{\varepsilon \to 0}\Zinvuf\bigl(p_{CS}\circ \tilde{\gamma}_{1}\vert_{\left[\varepsilon, 1-\varepsilon\right]}\bigr)\Zinvuf\bigl(p_{CS}\circ \tilde{\gamma}_{2}\vert_{\left[\varepsilon, 1 -\varepsilon\right]}\bigr)
\\
&=&\lim_{\varepsilon \to 0}\Zinvuf\bigl(p_{CS}\circ \tilde{\gamma}_{1}\vert_{\left[\varepsilon, 1-\varepsilon\right]}p_{CS}\circ \tilde{\gamma}_{2}\vert_{\left[\varepsilon, 1 -\varepsilon\right]}\bigr)\\
&=&\lim_{\varepsilon \to 0} \Zinvuf\bigl(p_{CS}\circ \bigl(\tilde{\gamma}_{1}\vert_{\left[\varepsilon, 1/2\right]}\tilde{\gamma}_{2}\vert_{\left[1/2, 1 -\varepsilon\right]}\bigr)\bigr),\end{array}\end{equation*}
where the second equality comes from Proposition~\ref{propfuncbraidprel}.
\eop

Theorem~\ref{thmfstconsttang}, Proposition~\ref{propfuncbraidprel}, Theorem~\ref{thmpoirbraid}, and Proposition~\ref{propfuncqbraid} allow us to extend $\Zinvufrf$ unambiguously to the $q$-tangles of Definition~\ref{deftanglelong} as follows.

\begin{definition}
\label{defzqtang} Let $\gamma^-$ and $\gamma^+$ be $q$-braids. Let $(\hcylc,\tanghcyll)$ be a framed tangle whose bottom and top configurations are $\gamma^-(1)$ and $\gamma^+(0)$, respectively.
Orient the strands of $T(\gamma^-)$ and $T(\gamma^+)$ so that their orientations are consistent with the orientation of $\tanghcyll$.
Define $\Zinvufrfneg\left(T(\gamma^-)(\hcylc,\tanghcyll)T(\gamma^+)\right)$ to be
\begin{equation*}\Zinvufrfneg\left(T(\gamma^-)(\hcylc,\tanghcyll)T(\gamma^+)\right)=\Zinvuf(\gamma^-)\Zinvufrfneg(\hcylc,\tanghcyll)\Zinvuf(\gamma^+),\end{equation*}
with respect to Definition~\ref{deffstconsttangframed} of $\Zinvufrfneg(\hcylc,\tanghcyll)$.\footnote{Alternatively, we can forget the orientations of $\tanghcyll$, $\gamma^-$, and $\gamma^+$ since $\Zinvufrf$ depends on the component orientations as in Proposition~\ref{proporcomp} for $\tanghcyll$, $\gamma^-$, and $\gamma^+$.}
\end{definition}

\begin{remark}
\label{deflimzqtang}
 Let $(\hcylc,\tanghcyll_q)$ be a $q$-tangle from a bottom limit configuration $\confc^- \in \ccompuptd{\finsetb^-}{\CC}$ to a top configuration $\confc^+ \in \ccompuptd{\finsetb^+}{\CC}$.
 These limit configurations are initial points $\gamma^{\pm}(0)=\confc^{\pm}$ of polynomial paths
 $\gamma^{\pm}$ of $\ccompuptd{\finsetb^{\pm}}{\CC}$,
 such that $\gamma^{\pm}\left(\left]0,1\right]\right) \subset \cinjuptd{\finsetb^{\pm}}{\CC}$ (as in the proof of Proposition~\ref{propfuncqbraid}). 
 This allows us to regard $\confc^{\pm}$ as a limit 
 \begin{equation*}\confc^{\pm}=\lim_{t \to 0} \gamma^{\pm}(t)\end{equation*}
 of injective configurations and view $\tanghcyll_q$ as a limit of
 framed tangles between injective configurations.
We can indeed view the $q$-tangle $\tanghcyll_q$ as the limit of 
 \begin{equation*}\tanghcyll_{q,\varepsilon}=T(\overline{\gamma^-\vert_{\left[0,\varepsilon\right]}})\tanghcyll_qT(\gamma^+\vert_{\left[0,\varepsilon\right]})\end{equation*} when $\varepsilon$ tends to $0$.
 Then $\Zinvufrfneg(\hcylc,\tanghcyll_q)$ can be defined alternatively to be
 \begin{equation*}\Zinvufrfneg(\hcylc,\tanghcyll_q)=\lim_{\varepsilon \to 0}\Zinvufrfneg(\tanghcyll_{q,\varepsilon}).\end{equation*}
 Indeed, the above consistent definition implies
 \begin{equation*}\Zinvufrfneg(\hcylc,\tanghcyll_q)=\Zinvuf(\gamma^-\vert_{\left[0,\varepsilon\right]})\Zinvufrfneg(\tanghcyll_{q,\varepsilon})\Zinvuf(\overline{\gamma^+\vert_{\left[0,\varepsilon\right]}}),\end{equation*}
 while Theorem~\ref{thmpoirbraid} and Proposition~\ref{propfuncqbraid} imply \begin{equation*}\lim_{\varepsilon \to 0}\Zinvuf(\gamma^-\vert_{\left[0,\varepsilon\right]})=\id\mbox{ and } \lim_{\varepsilon \to 0}\Zinvuf(\overline{\gamma^+\vert_{\left[0,\varepsilon\right]}})=\id.\end{equation*}
\end{remark}

We will construct variants of $\Zinvufrf$ in the spirit of Theorem~\ref{thmconststraight} in Section~\ref{secvarzinvf}.
These variants will allow us to prove Theorem~\ref{thmmainfunc} in Chapter~\ref{chappropzinvffunc}.

Proving Theorem~\ref{thmmainfunc} will require lengthy studies of compactifications
of configuration spaces, which are not manifolds with boundaries. These studies will help to get all
the nice and natural properties of $\Zinvufrf$, stated in Theorem~\ref{thmmainfunc}.
Robin Koytcheff, Brian Munson, and Ismar Voli\'{c} proposed another approach to obtain invariants of tangles and avoid our complicated configuration spaces in \cite{KoytcheffMunsonVolic}.

In the following three chapters, we will show the construction of $\Zinvufrf$ for $q$-tangles in detail, following the outline of this section.

\section{Properties of the functor\texorpdfstring{$\Zinvufrf$}{Zf}}
\label{secexppty}

In this section, we state the main natural properties of the functor $\Zinvufrf=(\Zinvufrf_n)_{n \in \NN}$, whose construction is outlined in the previous section.

\begin{theorem}
\label{thmmainfunc}
 The invariant $\Zinvufrf$ of $q$-tangles satisfies the following properties.
\begin{itemize}
\item $\Zinvufrf$ coincides with the invariant $\Zinvufrf$ of Definition~\ref{defZinvufrflink} for framed links in $\QQ$-spheres.
\item $\Zinvufrf$ coincides with the Poirier functor $Z^l$ of \cite{poirierv2} for combinatorial $q$-tangles in $\RR^3$.
\item \emph{Naturality:} If $\tanghcyll$ is a $q$-tangle with domain $\sourcetl$, then $\Zinvufrfneg(\tanghcyll)=\bigl(\Zinvufrf_k(\tanghcyll)\bigr)_{k\in \NN}$ is valued in $\Aavis(\sourcetl)=\Aavis(\sourcetl;\RR)$, and $\Zinvufrf_0(\tanghcyll)$ is the class of the empty diagram. If $\tanghcyll^{\prime}$ is a subtangle of $\tanghcyll$ with domain  $\sourcetl^{\prime}$, then $\Zinvufrfneg(\tanghcyll^{\prime})$ is obtained from $\Zinvufrfneg(\tanghcyll)$ by mapping all the diagrams with univalent vertices on $\sourcetl \setminus \sourcetl^{\prime}$ to zero and by forgetting $\sourcetl \setminus \sourcetl^{\prime}$.
\item \emph{Dependence on component orientations:}
If $\sourcetl$ and the components of $\tanghcyll$ are unoriented, then $\Zinvufrfneg(\tanghcyll)$ is valued in the space $\Aavis(\sourcetl)$ of Definition~\ref{defrkoruniv}, as in Proposition~\ref{proporcomp}. Otherwise, component orientation reversals affect $\Zinvufrf$ as in Proposition~\ref{proporcomp}.

\item \emph{Framing dependence:} For a $q$-tangle $\tanghcyll=(\hcylc,\sqcup_{j=1}^k K_j,\sqcup_{j=1}^k K_{j\parallel})$,
\begin{equation*}\prod_{j=1}^k\Bigl(\exp\bigl(-lk(K_j,K_{j\parallel})\alpha\bigr)\#_j\Bigr)\Zinvufrfneg(\tanghcyll)\end{equation*} is independent of the framing of $\tanghcyll$. It is denoted by $\Zinvuf(\hcylc,\sqcup_{j=1}^k K_j)$.\footnote{In particular, when a kink $\pkink$ is added to the $j$th component of a combinatorial $q$-tangle $\tanghcyll_c$, $\Zinvufrfneg(\tanghcyll_c)$ is changed to 
\begin{equation*}\Zinvufrfneg(\pkink \#_j\tanghcyll_c)  = \exp(\alpha)\#_j \Zinvufrfneg(\tanghcyll_c).\end{equation*} We similarly have
$\Zinvufrfneg(\nkink \#_j\tanghcyll_c)  = \exp(-\alpha)\#_j \Zinvufrfneg(\tanghcyll_c)$.}
\item \emph{Functoriality:} For two $q$-tangles $\tanghcyll_1$ and $\tanghcyll_2$ such that the bottom configuration of $\tanghcyll_2$ coincides with the top configuration of $\tanghcyll_1$, we have \begin{equation*}\Zinvufrfneg(\tanghcyll_1\tanghcyll_2)= \Zinvufrfneg(\tanghcyll_1) \Zinvufrfneg(\tanghcyll_2),\end{equation*}
with products obtained by stacking above in natural ways on both sides.
\item \emph{First duplication property:} Let $K$ be a component of a $q$-tangle $\tanghcyll$, then we have
\begin{equation*}\Zinvufrfneg\bigl(\tanghcyll( 2 \times K)\bigr)=\pi(2 \times K)^{\ast}\Zinvufrfneg(\tanghcyll)\end{equation*}
with respect to Notation~\ref{notationduplication}.
\item \emph{Second duplication property:} Let $\finsetb$ be a finite set, let $\confy $ be an element of $\ccompuptd{\finsetb}{\CC}$. Let $K$ be a component going from bottom to top in a $q$-tangle $\tanghcyll$. Recall $\tanghcyll(\confy \times K)= \tanghcyll\left(\left(\confy \times \left[0,1\right]\right)/K\right)$ from Notation~\ref{notcabling}.
 Then we have
\begin{equation*}\Zinvufrfneg\bigl(\tanghcyll( \confy \times K)\bigr)=\pi(\finsetb \times K)^{\ast}\Zinvufrfneg(\tanghcyll).\end{equation*}
\item \emph{Dependence on ambient orientation:} 
Let $s_{\frac12}$ be the orthogonal reflection with respect to the horizontal plane at height $\frac12$. Extend $s_{\frac12}$ from $\partial \hcylc$ to an orientation-reversing diffeomorphism $s_{\frac12}$ of $\hcylc$. Define the parallels $s_{\frac12}(K)_{\parallel}$ of interval components $s_{\frac12}(K)$ of $s_{\frac12}(\Link)$ so that they satisfy
$lk_{s_{\frac12}(\hcylc)}\bigl(s_{\frac12}(K),s_{\frac12}(K)_{\parallel}\bigr)=-lk_{\hcylc}\bigl(K,K_{\parallel}\bigr)$. Then
we have \begin{equation*}\Zinvufrf_n\Bigl(s_{\frac12}(\hcylc),s_{\frac12} \circ \tanghcyll\Bigr)=(-1)^n\Zinvufrf_n(\hcylc,\tanghcyll)\end{equation*} 
for all $n \in \NN$.

\item \emph{Symmetry:}  Let $\rho$ be a rotation of $\RR^3$ that preserves the standard homology cylinder $\drad{1} \times \left[0,1\right]$ (setwise). Let $\tanghcyll$ be a $q$-tangle of a rational homology cylinder $\hcylc$. Extend $\rho$ from $\partial \hcylc$ to an orientation-preserving diffeomorphism $\rho$ of $\hcylc$.
If the angle of $\rho$ is different from $0$ and $\pi$, assume that the interval components of $\tanghcyll$ go from bottom to top or from top to bottom. Then we have \begin{equation*}\Zinvufrfneg\bigl(\rho(\hcylc),\rho \circ \tanghcyll\bigr)=\Zinvufrfneg(\hcylc,\tanghcyll),\end{equation*}
where the parallels of interval components $\rho(K)$ of $\rho(\tanghcyll)$ are defined so that they satisfy $lk_{\rho(\hcylc)}\bigl(\rho(K),\rho(K)_{\parallel}\bigr)=lk_{\hcylc}\bigl(K,K_{\parallel}\bigr)$.\footnote{This statement applies to a diffeomorphism $\rho$ of $\hcylc$ restricting to $\partial \hcylc$ as the identity map, in particular.}
\item \emph{Cabling property:} Let $\finsetb$ be a finite set with cardinality greater than $1$. Let $\confy \in \ccompuptd{\finsetb}{\CC}$, let $\confy \times \left[0,1\right]$ denote the corresponding $q$-braid, and let $K$ be a strand of $\confy \times \left[0,1\right]$. Let $\tanghcyll$ be a $q$-tangle with domain $\sourcetl$.
Then $\Zinvufrfneg\left(\left(\confy \times \left[0,1\right]\right)\left(\tanghcyll/K\right)\right)$ is obtained from $\Zinvufrfneg(\tanghcyll)$ by the natural injection from $\Aavis(\sourcetl)$ to $\Aavis\left(\sqcup_{b \in \finsetb}\RR^{\{b\}}\left(\frac{\sourcetl}{K}\right)\right)$.
\item \emph{Full twist in degree one:} The expansion $\Zinvufrf_{\leq 1}$ up to degree $1$ of $\Zinvufrf$ satisfies 
\begin{equation*}\Zinvufrf_{\leq 1}\left(\fulltwistnonor\right)=1+\left[\onechordtwos\right],\end{equation*}
where the endpoints of the tangle are assumed to lie on $\RR \times \{0,1\}$.
\item \emph{Group-like behavior:} For any $q$-tangle $\tanghcyll$ and any integer $n$, we have \begin{equation*}\Delta_n\left(\Zinvufrf_n(\tanghcyll)\right)=\sum_{i=0}^n\Zinvufrf_i(\tanghcyll) \otimes \Zinvufrf_{n-i}(\tanghcyll)\end{equation*} with respect to the coproduct maps $\Delta_n$ of Section~\ref{seccoprod}.
\end{itemize}
\end{theorem}
The definition of $\Zinvufrf$ in Section~\ref{secdefzinvf} obviously extends the definition of $\Zinvufrf$ for tangles with empty boundary.
Note that the naturality property in Theorem~\ref{thmmainfunc} easily follows from the definition (which will be justified later). The behavior of $\Zinv(\hcylc,\tanghcyll,\tau)$ with respect to the coproduct can be observed from the definition, as in the proof of Proposition~\ref{propgrouplike}. According to Lemmas~\ref{lemAHopf} and \ref{lemexpprimgrouplike}, Definition~\ref{defxin}, and Proposition~\ref{propanom}, the correction factors are group-like. So the compatibility between the various products and the coproduct
ensures that $\Zinvufrf$ behaves as stated in Theorem~\ref{thmmainfunc} with respect to the coproduct,
for framed tangles between injective configurations. Then Remark~\ref{deflimzqtang} ensures that this also holds for general $q$-tangles.

\begin{corollary}
\label{corlogZinvfunc}
 If $\tanghcyll$ has at most one component,
let $p^c$ be the projection given by Corollary~\ref{corprojprim} from $\Aavis(\source)$ to the space $\Aavis^c(\source)$ of its primitive elements. 
Set \begin{equation*}\zinvufrfneg(\hcylc,\tanghcyll)=p^c\left(\Zinvufrfneg(\hcylc,\tanghcyll)\right).\end{equation*}
Then we have \begin{equation*}\Zinvufrfneg(\hcylc,\tanghcyll)=\exp\left(\zinvufrfneg(\hcylc,\tanghcyll)\right).\end{equation*}
\end{corollary}
\bp In these cases, Lemma~\ref{lemAHopf} guarantees that $\Aavis(\source)$ is a graded bialgebra, and Theorem~\ref{thmmainfunc} implies that $\Zinvufrfneg(\hcylc,\tanghcyll)$ is group-like. Conclude with Theorem~\ref{thmgrouplikexp}.\eop

The proof of Theorem~\ref{thmmainfunc} will be finished at the end of Section~\ref{secduprop}.
The multiplicativity of $\Zinvuf$ under connected sum of Theorem~\ref{thmconnsum}
is a direct corollary of the functoriality of $\Zinvufrf$ in the above statement.
The functoriality also implies that $\Zinvuf$ and $\Zinvufrf$ map tangles consisting of vertical segments in the standard cylinder to $1$. Consider such a trivial braid consisting of the two vertical segments 
$\{-\frac12\} \times \left[0,1\right]$ and $\{\frac12\} \times \left[0,1\right]$. Cable $\{-\frac12\} \times \left[0,1\right]$ by a $q$-tangle $(\hcylc_1,\tanghcyll_1)$, and cable $\{\frac12\} \times \left[0,1\right]$ by a $q$-tangle $(\hcylc_2,\tanghcyll_2)$. Call the resulting $q$-tangle $(\hcylc_1 \otimes \hcylc_2,\tanghcyll_1 \otimes \tanghcyll_2 )$. Formally, this tangle may be expressed as $\left(\left(\{-\frac12,\frac12\} \times \left[0,1\right]\right)\left(\frac{(\hcylc_1,\tanghcyll_1)}{\{-\frac12\} \times \left[0,1\right]} \right)\right)\left(\frac{(\hcylc_2,\tanghcyll_2)}{\{\frac12\} \times \left[0,1\right]}\right)$.

\begin{corollary}\label{cormonoidal} The functor $\Zinvufrf$ satisfies the following \emph{monoidality property} with respect to the above structure:
\begin{equation*}\Zinvufrfneg(\hcylc_1 \otimes \hcylc_2,\tanghcyll_1 \otimes \tanghcyll_2 )=\Zinvufrfneg(\hcylc_1,\tanghcyll_1) \otimes \Zinvufrfneg(\hcylc_2,\tanghcyll_2),\end{equation*}
where the product $\otimes$ of the right-hand side is simply induced by the disjoint union of diagrams.
\end{corollary}
\bp This is a consequence of the cabling property and the functoriality in Theorem~\ref{thmmainfunc}.
\eop

 More generally, Theorem~\ref{thmmainfunc} implies that the Poirier functor $Z^l$ is a good monoidal functor.
The multiplicativity of $\Zinvuf$ under connected sum of Theorem~\ref{thmconnsum} is also a consequence of Corollary~\ref{cormonoidal}.

\begin{remark}\label{rkfirstdup}
The first duplication property may be iterated. Note that $\pi(r \times K)^{\ast}$ is nothing but the composition of $(r-1)$ $\pi(2 \times K)^{\ast}$. Also note that iterating duplications $\frac{(\bulletpetit \bulletpetit)}{\bulletpetit}$ for configurations produces combinatorial configurations as in Example~\ref{exaqtang}, i.e., elements in the $0$-dimensional strata of some $\ccompuptd{\finsetb}{\RR}$ discussed in Example~\ref{exaStasheff}. For example, we have
\newcommand{\cupduptwo}{\begin{tikzpicture} [baseline=-.25cm]
\useasboundingbox (-.65,-.55) rectangle (.65,.2);
\draw (-.4,0) arc (-180:0:.4);
\draw (-.2,0) arc (-180:0:.2); 
\draw (-.5,-.15) node{\tiny 1} (-.05,-.1) node{\tiny 2};
\draw (-.5,.15) node{\scriptsize \textrm{(}} (-.1,.15) node{\scriptsize \textrm{)}} (.1,.15) node{\scriptsize \textrm{(}} (.5,.15) node{\scriptsize \textrm{)}};
\fill (-.4,.15) circle (1pt) (-.2,.15) circle (1pt) (.4,.15) circle (1pt) (.2,.15) circle (1pt);
\end{tikzpicture}}
\newcommand{\cupduptwoone}{\begin{tikzpicture} [baseline=-.25cm]
\useasboundingbox (-.65,-.55) rectangle (.65,.2);
\draw (-.4,0) arc (-180:0:.4);
\draw (-.5,-.15) node{\tiny 1};
\end{tikzpicture}}
\newcommand{\cupdupthree}{\begin{tikzpicture} [baseline=-.25cm]
\useasboundingbox (-.75,-.55) rectangle (.75,.2);
\draw (-.52,0) arc (-180:0:.52);
\draw (-.4,0) arc (-180:0:.4);
\draw (-.2,0) arc (-180:0:.2); 
\draw (-.65,.15) node{\scriptsize ((} (-.1,.15) node{\scriptsize \textrm{)}} (-.3,.15) node{\scriptsize \textrm{)}} (.1,.15) node{\scriptsize \textrm{(}} (.3,.15) node{\scriptsize \textrm{(}} (.65,.15) node{\scriptsize ))};
\fill (-.4,.15) circle (1pt) (-.2,.15) circle (1pt) (.4,.15) circle (1pt) (.2,.15) circle (1pt) (-.52,.15) circle (1pt) (.52,.15) circle (1pt);
\end{tikzpicture}}
\begin{equation*} \cup \left(2 \times \cup \right)=\cupduptwo\;\; \mbox{and} \;\; \cupduptwo \left(2 \times \cupduptwoone \right)=\cupdupthree.\end{equation*}
So, iterating twice the first duplication property implies
\begin{equation*}\Zinvufrfneg\left(\cupdupthree\right)=\pi\left(3 \times \cup \right)^{\ast}\Zinvufrfneg\left(\cup\right).\end{equation*}
\end{remark}
\begin{remark}\label{rksecondtdup}
The second duplication property, together with the behavior of $\Zinvufrf$ under component orientation reversal, yields a similar duplication property for a component $K$ going from top to bottom.
The behavior of $\Zinvufrf$ under orientation change, the functoriality, and the duplication properties allow us to generalize the cabling property to cablings of components $K$ going from bottom to top or from top to bottom in arbitrary $q$-tangles, by arbitrary $q$-tangles.
We may similarly generalize the cabling property to components $K$ going from top to top or from bottom to bottom, cabled by $q$-tangles in a rational homology cylinder whose bottom or top configurations are combinatorial configurations (as in Example~\ref{exaStasheff}).
In both cases, we can perform the insertion of the nontrivial part $T_i$ near an end of $K$ so that the result is a vertical composition of a tangle obtained by cabling a strand in a trivial vertical braid with $T_i$ and a possibly iterated duplication of the (main) tangle $T_m$ in the (main) rational homology cylinder $\hcylc_m$, as in Figure~\ref{figcablingbottom}.

\bfig
\centering
\begin{tikzpicture}[scale=1.4]
\draw (1,.95) node{\scriptsize{$T$}} (.4,.7) rectangle (1.6,1.2) (1,0) -- (1,.7) (.6,1.2) -- (.6,1.9) (1,1.2) -- (.6,1.9) (1.4,1.2) -- (1.4,1.9);
 \draw  [thin] (0,0) rectangle (2,1.9);
  \draw (0,.95) node[left]{\scriptsize $T_i=$};
 \begin{scope}[xshift=4.75 cm]
 \draw  [thin] (-1,0) rectangle (1,1.9);
\draw (.4,.5) arc (-90:80:.5 and .15);
\draw (.4,.5) arc (-90:-260:.5 and .15);
\draw [->] (-.4,0) -- (-.4,.9) arc (180:0:.4);
\draw (.4,0) -- (.4,.43) (.4,.57) -- (.4,.9);
\draw (.35,.95) node[right]{\scriptsize $K$};
\draw (-.95,.95) node[left]{\scriptsize $T_m=$};
 \end{scope}
 
  \begin{scope}[xshift=3.5 cm, yshift=-2.4 cm]
 \draw [thin] (-1,.25) rectangle (1,1.5);
\draw (.4,.5) arc (-90:65:.5 and .15);
\draw (.4,.5) arc (-90:-245:.5 and .15);
\draw (-.5,.25) -- (-.5,.9) arc (180:0:.5);
\draw (-.45,.25) -- (-.45,.9) arc (180:0:.45);
\draw (.45,.25) -- (.45,.44) (.45,.56) -- (.45,.9);
\draw (.5,.25) -- (.5,.44) (.5,.56) -- (.5,.9);
\draw [->] (-.35,.25) -- (-.35,.9) arc (180:0:.35);
\draw (.35,.25) -- (.35,.44) (.35,.56) -- (.35,.9);
\draw (-.95,.55) node[left]{\scriptsize $T_m\left(\frac{T_i}{K}\right)=$};
 \end{scope}
\begin{scope}[xshift=3.5 cm,  yshift=-2.8 cm]
 \draw  [thin] (-1,0) rectangle (1,.65);
 \draw (-.5,.65) -- (-.5,.42) (-.45,.65) -- (-.425,.42) (-.35,.65) -- (-.35,.42)  (-.65,.15) rectangle (-.2,.42) (-.425,.15) -- (-.425,0);
 \draw (-.425,.26) node{\scriptsize $T$};
 \draw (0,.3) node{\scriptsize $\otimes$} (.58,-.18) node{\scriptsize $)$} (.62,-.18) node{\scriptsize $)$} (.27,-.18) node{\scriptsize $($} (.39,-.18) node{\scriptsize $($};
 \draw  (.5,.65) -- (.5,0) (.45,.65) -- (.45,0) (.35,.65) -- (.35,0);
 \fill (.52,-.18) circle (1pt) (.43,-.18) circle (1pt) (.32,-.18) circle (1pt) (-.425,-.18) circle (1pt);
 \end{scope}
 
   \begin{scope}[xshift=5.93 cm, yshift=-2.4 cm]
\draw (.4,.5) arc (-90:80:.5 and .15);
\draw (.4,.5) arc (-90:-260:.5 and .15);
\draw [->] (-.4,-.4) -- (-.4,.9) arc (180:0:.4);
\draw (-.95,.55) node[left]{\scriptsize $=$};
\end{scope}
\begin{scope}[xshift=5.93 cm, yshift=-2.8 cm]
 \draw  [thin] (-1,0) rectangle  (1,1.9) (-1,.65) -- (1,.65) (.65,.2) rectangle (.15,.5);
 \draw (.4,.5) -- (.4,1.3) (.22,.2) -- (.39,0)  (.4,.2) --  (.41,0) (.58,.2) -- (.41,0);
 \node[rotate=180] at (.4,.32) {\scriptsize $T$};
 \draw (-.15,.3) node{\scriptsize $\otimes$} (.58,-.18) node{\scriptsize $)$} (.62,-.18) node{\scriptsize $)$} (.27,-.18) node{\scriptsize $($} (.39,-.18) node{\scriptsize $($};
 \fill (.52,-.18) circle (1pt) (.43,-.18) circle (1pt) (.32,-.18) circle (1pt) (-.425,-.18) circle (1pt);
 \end{scope}
 \end{tikzpicture}
\caption{Cabling a component going from bottom to bottom in two different ways}
\label{figcablingbottom}

\end{figure}
\end{remark}

\begin{remark}\label{rkcircledup}
The behavior of $\Zinvufrf$ when a component $K$ of a link $\Link \colon \source \longrightarrow \crats(\hcylc_m)$ is cabled by a semi-pure $q$-tangle $(\hcylc,\tanghcylt_i)$ can be described as follows.
\begin{itemize}
\item cut the domain of $K$ to replace it by a copy of $\RR$ using Proposition~\ref{propdiagrs},
\item duplicate the corresponding strand and $\Zinvufrfneg(\rats(\hcylc_m),\Link)$, accordingly, as in the duplication property above,
\item multiply the obtained element by $\Zinvufrfneg(\hcylc,\tanghcylt_i)$ i.e., concatenate the diagrams, naturally,
\item finally, close the domain of $\Link(\tanghcylt/K)$.
\end{itemize}
This follows easily from Theorem~\ref{thmmainfunc}, by viewing $(\rats,\Link)$ as a vertical composition of two tangles, where the bottom one is just a cup $\cup$. This cup is a trivial strand going from top to top in a standard cylinder. It is a part of $K$. We illustrate the process in Figure~\ref{figcablingcomplink}. In order to apply the iterated first duplication property, we change $\tanghcylt_i$ to a conjugate of $\tanghcylt_i^c$ whose bottom (or top) configuration is combinatorial. The result follows with $\Zinvufrfneg(\hcylc,\tanghcylt_i^c)$ instead of $\Zinvufrfneg(\hcylc,\tanghcylt_i)$. Lemma~\ref{lemdupcom} guarantees that changing $\Zinvufrfneg(\hcylc,\tanghcylt_i^c)$ to $\Zinvufrfneg(\hcylc,\tanghcylt_i)$ in the above recipe does not change the result. (This is consistent with the fact that the tangles $\Link(\tanghcylt_i^c/K)$ and $\Link(\tanghcylt_i/K)$ coincide.)

\bfig
\centering
\begin{tikzpicture}[scale=1.4]
\begin{scope}[xshift=-4 cm, yshift= -.6 cm]
\draw  [thin] (0,1.2) rectangle (.7,1.9);
\draw  [thin] (.08,1.35) rectangle (.62,1.75);
\draw (.35,1.55) node{\scriptsize{$T_i$}} 
(.15,1.2) -- (.15,1.35) (.15,1.9) -- (.15,1.75) 
(.35,1.2) -- (.35,1.35) (.35,1.9) -- (.35,1.75)
(.55,1.2) -- (.55,1.35) (.55,1.9) -- (.55,1.75);
\draw (0,1.555) node[left]{\scriptsize $T_i=$};
\end{scope}
\begin{scope}[xshift=-2 cm, yshift= .2 cm]
\draw  [thin] (0,0) rectangle (1,1.5);
\draw (.25,1.25) node{\scriptsize $\hcylc_m$};
\draw [dashed, rounded corners] (.2,.4) -- (.2,.9) -- (.7,.9) -- (.7,.4);
\draw [rounded corners,<-] (.2,.4) -- (.2,.1) -- (.7,.1) -- (.7,.4);
\draw (.35,.4) node{\scriptsize $K$};
\draw[dashed]
(.7,.65) ellipse (.2 and .15);
\draw (0,.75) node[left]{\scriptsize $L=$};
 \end{scope}
 
\begin{scope}[xshift=-.5 cm, yshift = .4 cm]
\draw  [thin] (0,-.4) rectangle (1,1.5) (0,0) -- (1,0);
\draw (.25,1.25) node{\scriptsize $\hcylc_m$};
\draw [dashed, rounded corners] (.2,.4) -- (.2,.9) -- (.7,.9) -- (.7,.4);
\draw [rounded corners,<-] (.2,.4) -- (.2,-.3) -- (.7,-.3) -- (.7,.4);
\draw (.35,.4) node{\scriptsize $K$};
\draw[dashed](.7,.65) ellipse (.2 and .15);
\draw (0,.55) node[left]{\scriptsize $=$};
 \end{scope}
 
 \begin{scope}[xshift=-3 cm, yshift =-2.2 cm]
\draw  [thin] (0,-.4) rectangle (1,-.1) (0,0) rectangle (1,1.5);
\draw (.25,1.25) node{\scriptsize $\hcylc_m$};
\draw [dashed, rounded corners] (.2,.4) -- (.2,.9) -- (.7,.9) -- (.7,.4);
\draw [rounded corners]  (.2,-.1) -- (.2,-.3) -- (.7,-.3) -- (.7,.4);
\draw[->] (.2,0) -- (.2,.4);
\draw (.35,.4) node{\scriptsize $\RR$};
\draw[dashed](.7,.65) ellipse (.2 and .15);
\draw (0,-.6) node[right]{\scriptsize Cut $K$};
\draw [thin,->] (0,-.6) arc (-90:-270:.275);
 \end{scope}
 
 \begin{scope}[xshift=-.5 cm,yshift=-2.5cm]
\draw [rounded corners] (-.5,.4) -- (-.5,0) -- (.5,0) -- (.5,1.1);
\draw [dashed, rounded corners] (-.5,1.1) -- (-.5,1.8) -- (.5,1.8) -- (.5,1.1);
\draw [rounded corners] (-.45,.4) -- (-.45,.05) -- (.45,.05) -- (.45,1.1);
\draw [dashed, rounded corners] (-.45,1.1) -- (-.45,1.75) -- (.45,1.75) -- (.45,1.1);
\draw [rounded corners] (-.3,.4) -- (-.3,.2) -- (.3,.2) -- (.3,1.1);
\draw [dashed, rounded corners] (-.3,1.1) -- (-.3,1.6) -- (.3,1.6) -- (.3,1.1);
\draw [dashed] (.4,1.4) ellipse (.3 and .15);
\draw (0,-.3) node{\scriptsize Duplicate};
 \end{scope}
 
 \begin{scope}[xshift=1.5 cm,yshift=-2.5cm]
\draw [rounded corners] (-.5,.3) -- (-.5,0) -- (.5,0) -- (.5,1.1);
\draw [dashed, rounded corners] (-.5,1.1) -- (-.5,1.8) -- (.5,1.8) -- (.5,1.1);
\draw [rounded corners] (-.45,.3) -- (-.45,.05) -- (.45,.05) -- (.45,1.1);
\draw [dashed, rounded corners] (-.45,1.1) -- (-.45,1.75) -- (.45,1.75) -- (.45,1.1);
\draw [rounded corners] (-.3,.3) -- (-.3,.2) -- (.3,.2) -- (.3,1.1);
\draw [dashed, rounded corners] (-.3,1.1) -- (-.3,1.6) -- (.3,1.6) -- (.3,1.1);
\draw [dashed] (.4,1.4) ellipse (.3 and .15);
\draw  [thin] (-.7,.55) rectangle (-.1,.95);
\draw (-.4,.75) node{\scriptsize{$T_i$}} 
(-.5,.55) -- (-.5,.4) (-.5,.95) -- (-.5,1.1) 
(-.4,.55) -- (-.45,.4) (-.4,.95) -- (-.45,1.1) 
(-.3,.55) -- (-.3,.4) (-.3,.95) -- (-.3,1.1);
\draw (0,-.3) node{\scriptsize Multiply};
 \end{scope}
 
 \begin{scope}[xshift=3.5 cm,yshift=-2.5cm]
\draw [rounded corners] (-.5,.4) -- (-.5,0) -- (.5,0) -- (.5,1.1);
\draw [dashed, rounded corners] (-.5,1.1) -- (-.5,1.8) -- (.5,1.8) -- (.5,1.1);
\draw [rounded corners] (-.45,.4) -- (-.45,.05) -- (.45,.05) -- (.45,1.1);
\draw [dashed, rounded corners] (-.45,1.1) -- (-.45,1.75) -- (.45,1.75) -- (.45,1.1);
\draw [rounded corners] (-.3,.4) -- (-.3,.2) -- (.3,.2) -- (.3,1.1);
\draw [dashed, rounded corners] (-.3,1.1) -- (-.3,1.6) -- (.3,1.6) -- (.3,1.1);
\draw [dashed] (.4,1.4) ellipse (.3 and .15);
\draw  [thin] (-.7,.55) rectangle (-.1,.95);
\draw (-.4,.75) node{\scriptsize{$T_i$}} 
(-.5,.55) -- (-.5,.4) (-.5,.95) -- (-.5,1.1) 
(-.4,.55) -- (-.45,.4) (-.4,.95) -- (-.45,1.1) 
(-.3,.55) -- (-.3,.4) (-.3,.95) -- (-.3,1.1);
\draw (0,-.3) node{\scriptsize Close};
 \end{scope}
 \end{tikzpicture}
\caption{Cabling a link component with a semi-pure tangle, step by step}
\label{figcablingcomplink}

\end{figure}

The first cutting step is not legitimate if $\Link$ has interval components. (Recall Exercise~\ref{exocounterexcut}). So the above recipe does not generalize to this case.\footnote{For a circle component $K$ in a general $q$-tangle $T_m$, and a semi-pure $q$-tangle $T_i$ as above, $T_m(T_i/K)$ is well-defined and it coincides with $T_m(T_i^c/K)$ for any conjugate $T_i^c$ of $T_i$. We can still write $T_m$ as a product $\left(\confc^-_m \times [0,1] \otimes \cup\right) T^{\prime}_m$, where $\cup$ is a part of $K$, pick a conjugate $T_i^c$ of $T_i$ whose bottom (or top) configuration is combinatorial, and compute $\Zinvufrfneg(T_m(T_i^c/K))$ from $\Zinvufrfneg(T^{\prime}_m)$ and $\Zinvufrfneg(T_i^c)$. Unfortunately, $\Zinvufrfneg(T^{\prime}_m)$ is not determined by $\Zinvufrfneg(T_m)$ anymore.}
\end{remark}

For a $q$-braid $\gamma \colon \left[0,1\right] \to \ccompuptd{\finsetb}{\CC}$, $\Zinvufrfneg(\gamma)=\Zinvuf(\gamma)$ stands for 
\begin{equation*}\Zinvufrfneg\bigl(T(\gamma)\bigr)=\Zinvufrfneg\bigl(\drad{1} \times \left[0,1\right],T(\gamma)\bigr).\end{equation*}                             
The following proposition, which leads to interesting cablings, is a corollary of Theorem~\ref{thmmainfunc}.

\begin{proposition}
\label{proptangtorus}
Let $q$ be a positive integer.
Let $\gamma_{1,q}$ be a braid represented by the map
\begin{equation*}\begin{array}{lll} \left[0,1\right] \times \underline{q}& \to& \CC \\
  (t,k) & \mapsto & \frac12 \exp \left(\frac{2i\pi (k +\chi(t))}{q}\right)
 \end{array}\end{equation*}
for a surjective map $\chi \colon \left[0,1\right] \to \left[0,1\right]$ with nonnegative derivative, which is constant in neighborhoods of $0$ and $1$, as in Figure~\ref{figbraidthreestrands} of $\gamma_{1,3}$. Recall Notation~\ref{notationduplication} and the anomaly $\alpha$ of Proposition~\ref{propanom}.
Then we have \begin{equation*}\Zinvufrfneg(\gamma_{1,q})=\exp\left( \pi\left(q \times \nochordonesorshort \right)^{\ast}(\frac{1}{q}\alpha) \right)\left(\exp(-\frac{1}{q}\alpha) \otimes \dots \otimes \exp(-\frac{1}{q}\alpha)\right).\end{equation*}
\end{proposition}
\bp
Let $\tilde{\gamma}_{1,1}$ be the trivial one-strand braid $K$ in the standard cylinder equipped with its parallel $K_{\parallel}$ such that $lk(K,K_{\parallel})=1$.
The framing dependence property in Theorem~\ref{thmmainfunc} implies $\Zinvufrfneg(\tilde{\gamma}_{1,1})=\exp(\alpha)$.
Let $\tilde{\gamma}_{q,q}$ be the $q$-tangle obtained by cabling $\tilde{\gamma}_{1,1}$ 
as in the second duplication property, by replacing the one-point configuration by the planar configuration of $\CC$ consisting of the $q$ points
$\frac12 \exp \left(\frac{2i\pi k}{q}\right)$, for $k \in \underline{q}$.
This duplication operation equips each strand $K_k$ of $\tilde{\gamma}_{q,q}$ with a parallel $K_{k\parallel,1}$ such that
$lk(K_k,K_{k\parallel,1})=1$. The $q$-tangle $\tilde{\gamma}_{q,q}$ coincides with $\gamma_{1,q}^q$ except for the framing since the standard framing of $\gamma_{1,q}^q$ equips $K_k$  with a parallel $K_{k\parallel}$ such that
$lk(K_k,K_{k\parallel})=0$.
According to the second duplication property, we have
\begin{equation*}\Zinvufrfneg(\tilde{\gamma}_{q,q})=\exp\left( \pi\left(q \times \nochordonesorshort\right)^{\ast}(\alpha) \right),\end{equation*}
whereas the framing dependence property implies
\begin{equation*}\Zinvufrfneg(\gamma_{1,q}^q)=\Zinvufrfneg(\tilde{\gamma}_{q,q})\bigl(\exp(-\alpha) \otimes \dots \otimes \exp(-\alpha)\bigr).\end{equation*}
By the invariance of $\Zinvufrf$ under rotation, $\Zinvufrfneg(\gamma_{1,q})$ is invariant under cyclic permutation of the strands. So the functoriality implies $\Zinvufrfneg(\gamma_{1,q}^q)=\Zinvufrfneg(\gamma_{1,q})^q$.
The result follows by unicity of a $q$th root of $\Zinvufrfneg(\gamma_{1,q}^q)$ with $1$ as degree $0$ part.
\eop

Dror Bar-Natan, Thang L\^e, and Dylan Thurston computed the Kontsevich integral of the trivial knot $O$ in \cite{bltwheels}.
Thus, Corollary~\ref{corlesunikontwoleg} allows one to express $\Zinvufrf$ for the unknot and for the torus knots as a function of the anomaly $\alpha$.
Note that the symmetry properties imply that $\Zinvufrfneg(\cap)$ vanishes in odd degree and that we have $\Zinvufrfneg(\cap)=\Zinvufrfneg(\cup)=\sqrt{\Zinvufrfneg(O)}$, where we implicitly use the natural isomorphism of Proposition~\ref{propdiagcom}.

\begin{lemma}
\label{lemalphatwoleg} We have
\begin{equation*}\Zinvufrfneg\left(\fulltwist \right)= \exp\left(\Psi(2\alpha)\left(\onechordtwosor\right)\right).\end{equation*}
\end{lemma}
\bp This lemma can be deduced from Corollary~\ref{corlesunikontwoleg}.
Below, as an exercise, we alternatively deduce it from Proposition~\ref{proptangtorus} and Theorem~\ref{thmmainfunc}, assuming that $\alpha$ is a two-leg element of $\Assis\left(\nochordonesorshort\right)$, but without assuming Corollary~\ref{corlesunikontwoleg}.
Since $\alpha$ is a two-leg element, we picture it as 
\begin{equation*}\alpha=\onechordwithboxonverticalinterval.\end{equation*}
So, using Lemma~\ref{lemsymvo}, we get
\begin{equation*}\pi\left(2 \times \nochordonesorshort\right)^{\ast}(\alpha)= 2 \onechordwithboxtwosor + \onechordwithboxonverticalinterval  \nochordonesor + \nochordonesor \onechordwithboxonverticalinterval.\end{equation*}
Since $\onechordwithboxonverticalinterval$ can be slid along its interval, we obtain
\begin{equation*}\exp\Bigl( \pi\bigl(2 \times \nochordonesorshort\bigr)^{\ast}(\alpha) \Bigr)=
\exp\left(\Psi(2\alpha)\left(\onechordtwosor\right)\right) \bigl(\exp( \alpha) \otimes \exp(\alpha)\bigr).\end{equation*}
\eop

\chapter{Invariance of \texorpdfstring{$\Zinvufrf$}{Zf} for long tangles}
\label{chapconszinvf}

In this chapter, we study appropriate compactifications of the configuration space $\check{C}(\crats(\hcylc),\tanghcyll;\Gamma)$ associated to a long tangle representative $\tanghcyll \colon \sourcetl \hookrightarrow \crats(\hcylc)$ and to
a Jacobi diagram $\Gamma$ with support $\sourcetl$.
These compactifications $C_{\tanghcyll}={C}(\rats(\hcylc),\tanghcyll;\Gamma)$
and $C^f_{\tanghcyll}={C}^f\!(\rats(\hcylc),\tanghcyll;\Gamma)$ are introduced in Definition~\ref{deffinercomptang}.
They allow us 
to prove Theorem~\ref{thmconvint}, which ensures that the integrals involved in the extension of $\Zinvufrf$ to long tangles converge, in Section~\ref{secstrucbyhand}. They also allow us to prove Theorem~\ref{thmfstconsttang}, which ensures the topological invariance of this extension of $\Zinvufrf$, in Section~\ref{secvarzinf}.
Our compactifications are locally diffeomorphic to products of smooth manifolds by singular subspaces of $\RR^n$ associated to trees.
We study these singular subspaces in Section~\ref{secsingmodtree} and show how Stokes' theorem applies to them in this preliminary section, which is independent of the rest of the book.
We describe the local structure of our compactifications and their codimension-one faces in Theorem~\ref{thmcomptang}.

\section{Singular models associated to trees}
\label{secsingmodtree}

\begin{definitions} \label{deffirsttree}
In this book, an \emph{oriented tree} is a tree $\CT$ as in Figure~\ref{figTgen} whose edges are oriented so that $\CT$ satisfies the following properties.
\begin{itemize}
\item There is exactly one vertex without outgoing edges. This vertex $\topt(\CT)$ is called the \emph{top} of $\CT$. It is also simply denoted by $\topt$ when $\CT$ is fixed.
\item The edges of $\CT$ are oriented towards $\topt(\CT)$. In other words, for any vertex $V$ of $\CT$ different from $\topt(\CT)$, the orientation of the edges in the injective path $\left[V,\topt(\CT)\right]$ from $V$ to $\topt(\CT)$ is induced by the orientation of $\left[V,\topt(\CT)\right]$.
\end{itemize}

Let $\CT$ be such an oriented tree.
A univalent vertex of $\CT$ with one outgoing edge is called a \emph{leaf} of $\CT$.
The set of leaves of $\CT$ is denoted by $L(\CT)$\index[N]{LT@$L(\CT)$ set of leaves}.
A \emph{node} of $\CT$ is a vertex with at least two ingoing edges.
A \emph{branch} of $\CT$ is an oriented injective path of oriented edges going from a leaf $\leafl$ to a node $\nodemaj$ or to the top $\topt$. Such a branch is denoted by $\left[\leafl,\nodemaj\right]$ or by $\left[\leafl,\topt\right]$. It is viewed as the subset of the set $E(\CT)$ of edges of $\CT$
between $\leafl$ and $\nodemaj$, or between $\leafl$ and $\topt$.
We denote the edge adjacent to a leaf $\leafl$ by $\mape(\leafl)$.
For any two vertices $\nodemaj_1$, $\nodemaj_2$ on the same branch $\left[\leafl,\topt\right]$, such that $\nodemaj_2$ is closer to $\topt$ than $\nodemaj_1$, $\left[\nodemaj_1,\nodemaj_2\right]$ is the set of edges between $\nodemaj_1$ and $\nodemaj_2$. These edges may contain $\nodemaj_1$ or $\nodemaj_2$ as an endpoint.
$\left]\nodemaj_1,\nodemaj_2\right]$ (resp. $\left[\nodemaj_1,\nodemaj_2\right[$)  denotes the set of edges of $\left[\nodemaj_1,\nodemaj_2\right]$ that do no contain $\nodemaj_1$ (resp. $\nodemaj_2$) as an endpoint.
For example, the set $\left[\leafl,\nodemaj\right] \setminus \{\mape(\leafl)\}$ is denoted by $\left]\leafl,\nodemaj\right]$. 
Similarly, for any two edges $\edgee_1$, $\edgee_2$ on the same branch $\left[\leafl,\topt\right]$, such that $\edgee_1$ is closer to $\leafl$ than $\topt$, $\left[\edgee_1,\edgee_2\right]$ is the set of edges between $\edgee_1$ and $\edgee_2$, including $\edgee_1$ and $\edgee_2$, and
$\left]\edgee_1,\edgee_2\right]$ (resp. $\left[\edgee_1,\edgee_2\right[$)  denotes the set $\left[\edgee_1,\edgee_2\right] \setminus \{\edgee_1\}$ (resp. $\left[\edgee_1,\edgee_2\right] \setminus \{\edgee_2\}$).
We also mix edges and vertices in this notation. For example, $\left]\edgee_1,\topt\right]$
is the set of edges between $\edgee_1$ and $\topt$ different from $\edgee_1$ that may contain $\topt$. The edges are ordered naturally on such an interval of edges.
The first is the first encountered when following the orientation of the interval induced by the orientation of the tree.

For two leaves $\leafl_1$ and $\leafl_2$ of $\CT$, $\nodemaj(\leafl_1,\leafl_2)$
denotes the node of $\CT$ such that $\left[\leafl_1,\topt\right]\cap\left[\leafl_2,\topt\right]=\left[\nodemaj(\leafl_1,\leafl_2),\topt\right] $.
For a subset $\CE$ of $E(\CT)$, $L(\CE)=L(\CE,\CT)$\index[N]{LE@$L(\CE)$ set of leaves} denotes the set of leaves $\leafl$ of 
$\CT$ such that $\left[\leafl,\topt\right]$ contains at least one edge of $\CE$. For a leaf
$\leafl$ in $L(\CE)$,  the closest edge to $\leafl$ in $\left[\leafl,\topt\right] \cap \CE$ is denoted by $\mape(\CE,\leafl)$\index[N]{eL@$\mape(\CE,\leafl)$ edge in a tree}.

\end{definitions}

\begin{examples}
 In Figure~\ref{figTgen}, the set of leaves is $L(\CT)=\{\leafl_0,\leafl_1, \dots, \leafl_5\}$.
 For $\CE=\{e_6,e_{13}\}$, we have $L(\CE)=L(\CT)$, $\mape(\CE,\leafl_0)=\mape(\CE,\leafl_3)=\mape(\CE,\leafl_4)=e_6$, and $\mape(\CE,\leafl_5)=e_{13}$.
 
\end{examples}
 
\bfig
\centering
\begin{equation*}
\begin{tikzpicture}
\draw[gray]  (2,4) -- (6.2,4) (-.6,1) -- (-.6,2) -- (.9,3) (.9,1) -- (.9,2) -- (.9,3) (2.4,1) -- (2.4,2) -- (.9,3) -- (3.4,4) (3.4,1) -- (3.4,4) (4.8,1) -- (4.8,4);
\fill (2,4) circle (1.5pt) (3.4,4) circle (1.5pt) (4.8,4) circle (1.5pt) (6.2,4) circle (1.5pt)  (-.6,1)  circle (1.5pt) (-.6,2)  circle (1.5pt) (.9,3)  circle (1.5pt) (.9,1)  circle (1.5pt) (.9,2)  circle (1.5pt) (2.4,1)  circle (1.5pt) (2.4,2)  circle (1.5pt) (3.4,1)  circle (1.5pt) (3.4,2.5) circle (1.5pt) (4.8,1)  circle (1.5pt) (4.8,2) circle (1.5pt) (4.8,3) circle (1.5pt);
\draw [gray,->] (2,4) -- (2.7,4);
\draw  (2,4.3) node{\scriptsize  $\leafl_0$} (2.7,4.25) node{\scriptsize $e_0$};
\draw [gray,->] (3.4,4) -- (4.1,4);
\draw (4.1,4.25) node{\scriptsize  $e_6$} (6.2,4.25) node{\scriptsize $\topt(\CT)$};
\draw [gray,->] (4.8,4) -- (5.5,4);
\draw  (5.4,4.25) node{\scriptsize  $e_7$};
\draw[gray,->] (-.6,1) -- (-.6,1.5);
\draw (-.95,1.5) node{\scriptsize $e_1$} (-.3,.9) node{\scriptsize$\leafl_1$};
\draw[gray,->] (.9,1) -- (.9,1.5);
\draw  (.55,1.5) node{\scriptsize $e_2$} (1.2,.9) node{\scriptsize$\leafl_2$};
\draw[gray,->] (2.4,1) -- (2.4,1.5);
\draw  (2.05,1.5) node{ \scriptsize$e_3$} (2.7,.9) node{\scriptsize$\leafl_3$};
\draw[gray,->] (3.4,1) -- (3.4,1.75);
\draw (3.7,1.75) node{\scriptsize $e_4$}
(3.7,.9) node{\scriptsize$\leafl_4$};
\draw[gray,->] (4.8,1) -- (4.8,1.5);
\draw (4.8,1.5) node[right]{\scriptsize $e_5=\mape(\leafl_5)$} 
(5.1,.9) node{\scriptsize$\leafl_5$};
\draw[gray,->] (4.8,2) -- (4.8,2.5);
\draw  (5.2,2.5) node{\scriptsize $e_{13}$};
\draw[gray,->]  (4.8,3) -- (4.8,3.5);
\draw (5.2,3.5) node{\scriptsize $e_{14}$};
\draw[gray,->]  (3.4,2.5) -- (3.4,3.25);
\draw (3.8,3.25) node{\scriptsize $e_{12}$};
\draw[gray,->]  (.9,3) -- (2.15,3.5);
\draw (1.5,3.6) node{\scriptsize $e_{9}$}
(1.05,2.95) node[right]{\scriptsize $\nodemaj(\leafl_1,\leafl_3)$};
\draw[gray,->] (2.4,2) -- (1.65,2.5);
\draw (2.2,2.5) node{\scriptsize $e_{11}$} ;
\draw [gray,->]  (.9,2) -- (.9,2.5);
\draw (1.2,2.35) node{\scriptsize $e_{10}$};
\draw[gray,->] (-.6,2) -- (.15,2.5);
\draw  (-.2,2.5) node{\scriptsize $e_8$};
\draw[very thick,->] (3.4,4) -- (4.8,4) (3.4,4) -- (4.1,4);
\draw[very thick,->] (4.8,2) -- (4.8,3) (4.8,2) -- (4.8,2.5);
\end{tikzpicture}\end{equation*}
\caption{A tree $\CT$ with a bold codimension-one system of edges}
\label{figTgen}

\end{figure}

This section is devoted to the study of the following space ${\setparamx}(\CT)$ associated to an oriented tree $\CT$.

\begin{definition} \label{defXT}
 For $\bigl((\varedge_e)_{e \in E(\CT)}\bigr) \in \left[0, \infty\right[^{E(\CT)}$
and for a branch $\left[\leafl,\nodemaj\right]$ of $\CT$, define \begin{equation*}\prodedge\bigl(\left[\leafl,\nodemaj\right]\bigr)=\prod_{e \in \left[\leafl,\nodemaj\right]}\varedge_e.\end{equation*}
Define ${\setparamx}(\CT)$\index[N]{X@${\setparamx}(\CT)$ singular space associated to a tree $\CT$} to be the set of the elements $\bigl((\varedge_e)_{e \in E(\CT)}\bigr)$ of $\left[0, \infty\right[^{E(\CT)}$ such that the equality 
\begin{equation*}\ast(\leafl_1,\leafl_2) \;\;\colon \;\;
\prodedge\bigl(\left[\leafl_1,\nodemaj(\leafl_1,\leafl_2)\right]\bigr)=\prodedge\bigl(\left[\leafl_2,\nodemaj(\leafl_1,\leafl_2)\right]\bigr)\end{equation*}
holds for any two leaves $\leafl_1$ and $\leafl_2$ of $\CT$. Set $\mathring{\setparamx}(\CT) ={\setparamx}(\CT) \cap \left]0,\infty\right[^{E(\CT)}$.
\end{definition}

\begin{example} \label{exanotcorner}
 The space ${\setparamx}(\CT)$ associated to the tree
 \begin{center}
\begin{tikzpicture} 
\draw (0,-.2) -- (2,.2) (4,-.2) -- (2,.2);
\draw [->] (3.5,.1) node{\scriptsize $e_{3}$} (4,-.2) -- (3.5,-.1);
\draw [->] (2.5,.3) node{\scriptsize $e_{4}$} (3,0) -- (2.5,.1);
\draw [->]  (.5,.1) node{\scriptsize $e_{1}$} (0,-.2) -- (.5,-.1);
\draw [->]  (1.5,.3) node{\scriptsize $e_{2}$} (1,0) -- (1.5,.1);
\fill (4,-.2) circle (1.5pt) (0,-.2) circle (1.5pt) (1,0) circle (1.5pt) (3,0) circle (1.5pt) (2,.2) circle (1.5pt);
\end{tikzpicture}
\end{center}
is
 $\{(u_1,u_2,u_3,u_4) \in \left[0,\infty\right[^4 \suchthat u_1u_2=u_3u_4\}$.  At its point $(0,0,0,0)$, the four half-lines $ \RR^+(1,0,0,0)$, $ \RR^+(0,1,0,0)$, $ \RR^+(0,0,1,0)$, and $\RR^+(0,0,0,1)$ embed in ${\setparamx}(\CT)$. Their four independent unit tangent vectors at $(0,0,0,0)$ generate $\RR^4$, but the complement of ${\setparamx}(\CT)$ is dense in $\left[0, \infty\right[^4$. So ${\setparamx}(\CT)$ is not a submanifold with ridges of $\left[0, \infty\right[^4$.
\end{example}

\begin{remarks}

 Let $\left[\leafl_1,\nodemaj\right]$ and $\left[\leafl_2,\nodemaj\right]$ be two branches of $\CT$ ending at the same vertex $\nodemaj$. If $\bigl((\varedge_e)_{e \in E(\CT)}\bigr) \in {\setparamx}(\CT)$, then we have $\prodedge\bigl(\left[\leafl_1,\nodemaj\right]\bigr)=\prodedge\bigl(\left[\leafl_2,\nodemaj\right]\bigr)$.
 
Let $\leafl_0$ be a leaf of $\CT$.
Then $\mathring{\setparamx}(\CT)$ is the set of the elements $\bigl((\varedge_e)_{e \in E(\CT)}\bigr)$ of $\left]0,\infty\right[^{E(\CT)}$
such that $\prodedge\bigl(\left[\leafl,\topt\right]\bigr)=\prodedge\bigl(\left[\leafl_0,\topt\right]\bigr)$ for any leaf $\leafl$ of $L(\CT) \setminus \{\leafl_0\}$.
Indeed, the equations $\ast(\leafl_1,\leafl_2)$ of Definition~\ref{defXT} are equivalent to $\prodedge\bigl(\left[\leafl_1,\topt\right]\bigr)=\prodedge\bigl(\left[\leafl_2,\topt\right]\bigr)$ when no variable is zero.
\end{remarks}

 \begin{definition} A \indexT{reducing system of edges} in an oriented tree $\CT$ is a set $\CE_r$ of edges such that $L(\CE_r) = L(\CT) \setminus \{\leafl_0\}$ for some leaf $\leafl_0$ and $\mape(\CE_r,.)$ is a bijection from $L(\CE_r)$ to $\CE_r$. For such a reducing system $\CE_r$, the leaf $\leafl_0$  of $L(\CT) \setminus L(\CE_r)$ is denoted by $\leafl_0(\CE_r)$. 
A \indexT{maximal free system of edges} of $\CT$ is the complement $E(\CT) \setminus \CE_r$ of a reducing system $\CE_r$ of edges of $\CT$ in $E(\CT)$.
\end{definition}

\begin{examples}
For example, for every leaf $\leafl_0$ of $\CT$, the set of edges adjacent to the leaves of $L(\CT) \setminus \{\leafl_0\}$ is a reducing system of edges.
In Figure~\ref{figTgen}, $\{e_1,e_{10},e_9,e_{12},e_{13}\}$ is also a reducing system of edges.
\end{examples}

\begin{lemma} \label{lemsysredcoord}
Let $\CE_r$ be a reducing system of edges of $\CT$.
Let $e_1$ be an edge of $\CE_r$.
Let $\leafl=\leafl(\CE_r,e_1)$ be the leaf of $\CT$ such that $e_1 = \mape(\CE_r,\leafl)$. There exists a unique leaf $\leafl^{\prime}=\leafl^{\prime}(\CE_r,e_1)$ of $\CT$ such that $e_1 \in \left[\leafl,\nodemaj(\leafl,\leafl^{\prime})\right]$ and $e_1$ is the only element of $\CE_r$ in $\left[\leafl,\nodemaj(\leafl,\leafl^{\prime})\right] \cup \left[\leafl^{\prime},\nodemaj(\leafl,\leafl^{\prime})\right]$. 
Then the equation $\ast(\leafl,\leafl^{\prime})$ may be written as 
 \begin{equation*}\varedge_{e_1} = \frac{\prodedge\bigl(\left[\leafl^{\prime},\nodemaj(\leafl,\leafl^{\prime})\right]\bigr)}{\prod_{e \in \left[\leafl,\nodemaj(\leafl,\leafl^{\prime})\right] \setminus \{ e_1\}} \varedge_e}=\varedge_{e_1}\left(E(\CT) \setminus \CE_r; (\varedge_e)_{e \in E(\CT) \setminus \CE_r}\right)\end{equation*} when $\eltx =\bigl((\varedge_e)_{e \in E(\CT)}\bigr) \in \mathring{\setparamx}(\CT)$.
 So it determines $\varedge_{e_1}$ in terms of the variables associated to $E(\CT) \setminus \CE_r$ in $\mathring{\setparamx}(\CT)$.
 Moreover, we have \begin{equation*}\frac{d \varedge_{e_1}}{\varedge_{e_1}} =\sum_{e \in \left[\leafl^{\prime},\nodemaj(\leafl,\leafl^{\prime})\right]}\frac{d \varedge_e}{\varedge_{e}} -\sum_{e \in \left[\leafl,\nodemaj(\leafl,\leafl^{\prime})\right] \setminus \{ e_1\}}\frac{d \varedge_e}{\varedge_{e}} \end{equation*} in $\mathring{\setparamx}(\CT)$.
\end{lemma}
\bp If there is an edge of $\CE_r$ after $e_1$ on $\left[\leafl,\topt\right]$, let $e_2= \mape(\CE_r,\leafl^{\prime})$ be the second edge of $\CE_r$ on $\left[\leafl,\topt\right]$.
This defines $\leafl^{\prime}$. Otherwise, set $\leafl^{\prime} = \leafl_0(\CE_r)$. See Figure~\ref{figlemsysredcoord}.  This proves the existence of $\leafl^{\prime}$ such that $e_1$ is the only element of $\CE_r$ in $\left[\leafl,\nodemaj(\leafl,\leafl^{\prime})\right] \cup \left[\leafl^{\prime},\nodemaj(\leafl,\leafl^{\prime})\right]$. For such an $\leafl^{\prime} \neq \leafl_0(\CE_r)$, $\mape(\CE_r,\leafl^{\prime})$ must be the second edge of $\CE_r$ on $\left[\leafl,\topt\right]$. So the condition of the statement determines $\leafl^{\prime}$.
\eop

\bfig
\centering
\begin{equation*}
\begin{tikzpicture}
\draw [dashed] (0,2) -- (2.6,2);
\draw [dotted] (3.6,2) -- (4,2);
\draw (0,2.25) node{\scriptsize $\leafl$} (.3,.8) node{\scriptsize $\leafl^{\prime}$} (.6,2) -- (1.6,2) (2.6,2) -- (3.6,2);
\draw [->] (.9,2.25) node[right]{\scriptsize $e_1$} (.6,2) -- (1.1,2);
\draw [->] (3.1,2.25) node{\scriptsize $e_2=\mape(\CE_r,\leafl^{\prime})$} (2.6,2) -- (3.1,2);
\fill (0,2) circle (1.5pt) (.6,2) circle (1.5pt) (1.6,2) circle (1.5pt) (2.1,2) circle (1.5pt) (2.6,2) circle (1.5pt) (3.6,2) circle (1.5pt) (0,.8) circle (1.5pt);
\draw [dashed] (0,.8) -- (2.1,2);
\draw (2,1.75) node[right]{\scriptsize $\nodemaj(\leafl,\leafl^{\prime})$};

\begin{scope}[xshift=7cm,yshift=-1.2cm]
\draw [dashed] (-1.5,3.2) -- (1.2,2);
\draw [dashed] (-1.5,2) -- (2,2);
\draw (-1.5,3.2) node[left]{\scriptsize $\leafl_0(\CE_r)$} (-1.5,2) node[left]{\scriptsize $\leafl$} (-.6,2) -- (0,2);
\draw [->] (-.3,1.8) node{\scriptsize $\edgee_1$} (-.6,2) -- (-.3,2);
\fill (-1.5,2) circle (1.5pt) (-.6,2) circle (1.5pt) (0,2) circle (1.5pt) (1.2,2) circle (1.5pt) (2,2) circle (1.5pt) (-1.5,3.2) circle (1.5pt);
\draw [dashed] (0,2) -- (2,2);
\draw (2,2.25) node{\scriptsize $\topt(\CT)$};
\draw (1.2,1.75) node{\scriptsize $\nodemaj(\leafl,\leafl_0(\CE_r))$};
\end{scope}
\end{tikzpicture}\end{equation*}
\caption{The two cases in Lemma~\ref{lemsysredcoord}}
\label{figlemsysredcoord}

\end{figure}

\begin{lemma} \label{lemsysred}
The set $\mathring{\setparamx}(\CT)$ is a smooth manifold of dimension
\begin{equation*}\dimd(\CT)=\cardlef{E(\CT)} -\cardlef{L(\CT)} +1.\end{equation*}  For every maximal free system $\CE_b$ of edges of $\CT$, $\mathring{\setparamx}(\CT)$ is freely parametrized by the (variables of) the edges of $\CE_b$.
For a subset $\CE$ of cardinality $\dimd(\CT)$, the form $\wedge_{e \in \CE}d\varedge_e$ is a nonvanishing volume form on $\mathring{\setparamx}(\CT)$ if and only if $ 
\CE$ is a maximal free system of edges of $\CT$.
\end{lemma}
\bp Let $\CE_b$ be a maximal free system of edges of $\CT$, and let $\CE_r$ be its complement. 
An edge $e$ of $\CE_r$ may be expressed as $\mape(\CE_r,\leafl(\CE_r,e))$ for a unique $\leafl(\CE_r,e)$ of $L(\CE_r)$, and any associated variable $\varedge_{e}$ can be expressed in terms of the variables associated to $\CE_b$ as in Lemma~\ref{lemsysredcoord}.
Furthermore, any $\bigl((\varedge_e)_{e \in E(\CT)}\bigr) \in \left]0,\infty\right[^{E(\CT)}$ such that 
\begin{equation*}\prodedge\bigl(\left[\leafl(\CE_r,e),\topt \right]\bigr)=\prodedge\bigl(\left[\leafl_0(\CE_r),\topt\right]\bigr),\end{equation*} for any $\edgee$ in $\CE_r$,
is in $\mathring{\setparamx}(\CT)$.
The equation defining $\varedge_{\edgee_1}$ in Lemma~\ref{lemsysredcoord}, for $\edgee_1 \in \CE_r$, implies
$\ast(\leafl(\CE_r,\edgee_1), \leafl^{\prime})$, where
$\leafl^{\prime}=\leafl_0(\CE_r)$ if $\edgee_1$ is the only edge of $\CE_r$ in 
$\left[\leafl(\CE_r,\edgee_1),\topt\right]$.
Otherwise, set $\leafl_{\edgee_1}=\leafl(\CE_r,\edgee_1)$ and let $\edgee_1$, $\edgee_2$, \dots, $\edgee_k$ denote the edges of $\CE_r$ on $\left[\leafl_{\edgee_1},\topt\right]$, where $\edgee_i=\mape(\CE_r,\leafl_{\edgee_i})$ denotes the $i^{th}$ edge that we meet from $\leafl_{\edgee_1}$ to $\topt$. Applying Lemma~\ref{lemsysredcoord} to $\edgee_i$ defines $\varedge_{e_i}$ by an equation, which is equivalent to 
$\ast(\leafl_{\edgee_i},\leafl_{\edgee_{i+1}})$ if $i<k$, and to $\ast(\leafl_{\edgee_k},\leafl_0(\CE_r))$ if $i=k$. These equations together imply $\ast(\leafl_{\edgee_1}, \leafl_0(\CE_r))$ for any $\leafl_{\edgee_1}$ of $L(\CT)$.

Therefore, any $\eltx=\bigl((\varedge_e)_{e \in E(\CT)}\bigr) \in \left]0,\infty\right[^{E(\CT)}$ satisfying the equations of  Lemma~\ref{lemsysredcoord} defining the $\varedge_{e}$ for $e \in \CE_r$ is in $\mathring{\setparamx}(\CT)$. Thus the space $\mathring{\setparamx}(\CT)$ is freely parametrized by the (variables of) the edges of $\CE_b$. So $\mathring{\setparamx}(\CT)$ is a smooth submanifold of $\RR^{E(\CT)}$ of dimension
\begin{equation*}\dimd(\CT)=\cardlef{E(\CT)} -\cardlef{L(\CT)} +1.\end{equation*}

Let $\CE$ be a subset of cardinality $\dimd(\CT)$, such that the form $\wedge_{e \in \CE}d\varedge_e$ is a nonvanishing volume form on $\mathring{\setparamx}(\CT)$. Set $\CE^c =E(\CT) \setminus \CE$. 
Let us show that the map $\mape(\CE^c,.) \colon L(\CE^c) \to \CE^c$ is injective.
If $\mape(\CE^c,\leafl_1)=\mape(\CE^c,\leafl_2)$, for two leaves $\leafl_1$ and $\leafl_2$, then $\left[\leafl_1,\nodemaj(\leafl_1,\leafl_2)\right] \cup \left[\leafl_2,\nodemaj(\leafl_1,\leafl_2)\right] \subset \CE$ and the relation $\prodedge\bigl(\left[\leafl_1,\nodemaj(\leafl_1,\leafl_2)\right]\bigr) = \prodedge\bigl(\left[\leafl_2,\nodemaj(\leafl_1,\leafl_2)\right]\bigr)$ gives rise to a nontrivial relation
between the forms $d\varedge_e$ for $e \in \left[\leafl_1,\nodemaj(\leafl_1,\leafl_2)\right] \cup \left[\leafl_2,\nodemaj(\leafl_1,\leafl_2)\right]$.

If $L\setminus L(\CE^c)$ has two distinct elements $\leafl_1$ and $\leafl_2$, then we similarly have a nontrivial linear relation between coordinate forms $d\varedge_e$ associated to edges of $\left[\leafl_1,\nodemaj(\leafl_1,\leafl_2)\right] \cup \left[\leafl_2,\nodemaj(\leafl_1,\leafl_2)\right]$. 
Therefore, the cardinality of $L(\CE^c)$ is at least $\cardlef{L(\CT)}-1$, which is the cardinality of $\CE^c$. So $\mape(\CE^c,.)$ is a bijection.
\eop

\begin{definition}
A \indexT{codimension-one system of edges} in an  oriented tree $\CT$ is a set $\CE$ of edges such that there is exactly one edge of $\CE$ in any path from a leaf to the top of $\CT$. 
\end{definition}

\begin{examples}
An example of a codimension-one system of edges of $\CT$ is the set of edges that start at leaves. In Figure~\ref{figTgen}, $\CE=\{e_6,e_{13}\}$, $\{e_7\}$, and $\{e_0,e_9,e_{12},e_5\}$ are other codimension-one systems of edges of $\CT$.
\end{examples}

\begin{lemma}
\label{lemcodimoneone}
For any codimension-one system $\CE_1$ of ${\setparamx}(\CT)$ and any edge $e_0 \in \CE_1$, $\CE_1 \setminus \{e_0\}$ can be completed to a reducing system that does not contain $e_0$.
\end{lemma}
\bp For each element $e$ of $\CE_1$, choose a leaf $\leafl_e$ such that $\mape(\CE_1,\leafl_e)=e$. Let $L_1=\{\leafl_e, e \in \CE_1\}$ be the set of these leaves, and set $L_2= L(\CT)  \setminus L_1$. Let $\CE_2=\mape(L_2)$ be the set of edges adjacent to the leaves of $L_2$. Then $\CE_1 \cup \CE_2 \setminus \{e_0\}$ is a reducing system of $\CT$.
\eop

\begin{lemma} \label{lemclosuresetparamx}
The closure of $\mathring{\setparamx}(\CT)$ in $\left[0, \infty\right[^{E(\CT)}$ is ${\setparamx}(\CT)$.

Let $\eltx = \bigl((\varedge_e)_{e \in E(\CT)}\bigr)$ be an element of 
${\setparamx}(\CT) \setminus \mathring{\setparamx}(\CT)$. Let $\CE(\eltx)$ denote the set of edges of $\CT$ such that $\varedge_e=0$. 
Then the image $\CE_1(x)$ of $\bigl(\mape(\CE(\eltx),.) \colon L(\CE(\eltx)) \to \CE(\eltx)\bigr)$ is a codimension-one system of edges of $\CT$. 

Let ${\setparamx}_{\CE(\eltx)}(\CT)$ be the set of
elements $\eltx^{\prime}$ of ${\setparamx}(\CT)$ such that $\CE(\eltx^{\prime})=\CE(\eltx)$.
Then ${\setparamx}_{\CE(\eltx)}(\CT)$ is a smooth manifold of dimension $\dimd(\CT)-1 -(\cardlef{\CE(\eltx)} - \cardlef{\CE_1(x)})$.
\end{lemma}
\bp Let $\eltx=\bigl((\varedge_e)_{e \in E(\CT)}\bigr) \in {\setparamx}(\CT)\setminus \mathring{\setparamx}(\CT)$. 
Let us first prove that the image $\CE_1(x)$ of $\mape(\CE(\eltx),.) \colon L(\CE(\eltx)) \to \CE(\eltx)$ is a codimension-one system of edges of $\CT$.
Since $\CE(\eltx) \neq \emptyset$, there is a leaf $\leafl$ of $\CT$ such that
$U\left(\left[\leafl,\topt\right] \right)=0$. This implies $U\left(\left[\leafl,\topt\right] \right)=0$ for all leaves. Thus, $L(\CE(\eltx))=L(\CT)$. Furthermore no branch $\left[\leafl,\topt\right]$ can contain more than one edge of $\CE_1(x)$.
Otherwise, the first two edges of $\CE_1(x)$ on such a branch $\left[\leafl_1,\topt\right]$ would be $\mape(\CE(\eltx),\leafl_1)$ and $\mape(\CE(\eltx),\leafl_2)$, respectively, and $\mape(\CE(\eltx),\leafl_1)$ would be the only edge of $\CE(\eltx)$ on $\left[\leafl_1,\nodemaj(\leafl_1,\leafl_2)\right] \cup \left[\leafl_2,\nodemaj(\leafl_1,\leafl_2)\right]$. Then we would have $U\left(\left[\leafl_1,\nodemaj(\leafl_1,\leafl_2)\right]\right)=0$ and $U\left(\left[\leafl_2,\nodemaj(\leafl_1,\leafl_2)\right]\right)\neq 0$. So $\ast(\leafl_1,\leafl_2)$ would not be satisfied. So $\CE_1(x)$ is a codimension-one system of edges of $\CT$.

Now, it suffices to prove the following two assertions. 
\begin{itemize}
\item $\eltx$ is in the closure of $\mathring{\setparamx}(\CT)$ in $\left[0, \infty\right[^{E(\CT)}$, and 
\item ${\setparamx}_{\CE(\eltx)}(\CT)$ is a smooth manifold of dimension $\dimd(\CT)-1 -(\cardlef{\CE(\eltx)} - \cardlef{\CE_1(x)})$.
\end{itemize}

Let us first prove them when $\mape(\CE(\eltx),.) \colon L(\CT) \to \CE_1(x)$ is injective.  
Define $x(t)=\bigl(\left(\varedge_e(t)\right)_{e \in E(\CT)}\bigr) \in \mathring{\setparamx}(\CT)$ from $x=\bigl((\varedge_e)_{e \in E(\CT)}\bigr)$ for $t \in \left]0,\infty\right[$, as follows. Pick $\leafl_0 \in L(\CT)$, and set $e_0=\mape(\CE(\eltx),\leafl_0)$.
Replace all the variables $\varedge_e$ for $e \in \CE(\eltx) \setminus \CE_1(x)$ by $t$, replace $\varedge_{e_0}$ by $t^k$ for some positive integer $k$, and leave the variables associated to the edges of $E(\CT) \setminus \CE(\eltx)$ (which are not zero) unchanged. For an edge $\edgef=\mape(\CE_1(x),\leafl_{\edgef})$ of $\CE_1(x)\setminus \{e_0\}$, set \begin{equation*}\varedge_{\edgef}(t) =\frac{\prodedge\bigl(\left[\leafl_0,\topt\right]\bigr)(t)}{\prod_{e \in \left[\leafl_{\edgef},\topt\right] \setminus \{\edgef\}}\varedge_e(t)}.\end{equation*}
Recall that $\edgef$ is the only edge of $\CE_1(x)$ on $\left[\leafl_{\edgef},\topt\right]$ and that $e_0$ is the only edge of $\CE_1(x)$ on $\left[\leafl_0,\topt\right]$.
Then $\varedge_{\edgef}(t)$ is equal to $\alpha t^{k+r(\edgef)}$ for some $\alpha>0$ and some $r(\edgef) \in \ZZ$. Choose $k$ so that $k+r(\edgef) \geq 1$ for any $\edgef \in \CE_1(x)$. Then $x(t)$ tends to $x$ when $t$ tends to zero. Furthermore, since 
all $\varedge_{\edgee}(t)$ are nonzero, the defining equations for the $\varedge_{\edgef}(t)$ are equivalent to the equations $\prodedge\bigl(\left[\leafl,\topt \right]\bigr)(t)=\prodedge\bigl(\left[\leafl_0,\topt \right]\bigr)(t)$.
They are satisfied for any $\leafl \in L(\CT)$, so $x(t) \in \mathring{\setparamx}(\CT)$.
This proves that $x$ is in the closure of $\mathring{\setparamx}(\CT)$.

The defining equations of ${\setparamx}(\CT)$ are satisfied as soon as the $\varedge_e$, for $e \in \CE_1(x)$,
vanish (still under the assumption that $\mape(\CE(\eltx),.) \colon L(\CE(\eltx)) \to \CE(\eltx)$ is injective). Therefore, ${\setparamx}_{\CE(\eltx)}(\CT)$ is a manifold freely parametrized by the variables corresponding to the edges of $E(\CT) \setminus \CE(x)$. Its dimension is $\cardlef{E(\CT)}  -\cardlef{\CE(\eltx)}- \cardlef{L(\CT)} +\cardlef{\CE_1(x)}= \dimd(\CT) -1 - (\cardlef{\CE(\eltx)} - \cardlef{\CE_1(x)})$.
The two assertions are proved when \begin{equation*}\mape\bigl(\CE(\eltx),.\bigr) \colon L\bigl(\CE(\eltx)\bigr) \to \CE(\eltx)\end{equation*} is injective.

In general, choose a leaf $\leafl_e$ such that $\mape(\CE(\eltx),\leafl_e)=e$
for each element $e$ of $\CE_1(x)$. Let $L_1=\{\leafl_e, e \in \CE_1(x)\} $ be the set of these leaves. Let $\CT_1$ be the subtree of $\CT$ such that $E(\CT_1)=\cup_{e \in \CE_1(x)} \left[\leafl_e,\topt\right]$ (so $L(\CT_1)=L_1$). Let $x_1$ be the natural projection of $x$ in $\left[0, \infty\right[^{E(\CT_1)}$. Note $\CE(x) \subseteq E(\CT_1)$ and
$\CE(x) = \CE(x_1)$.
Also note that the restriction of $\mape(\CE(x), .)$ to $L(\CT_1)$ is the map $\mape_{\CT_1}(\CE(x_1), .)$ associated to $\CT_1$, which is injective.

In particular, the first part of the proof expresses $x_1$ as a limit at $0$ of some continuous function $x_1(.) \colon \left]0,\infty\right[ \to \mathring{\setparamx}(\CT_1)$, such that $\varedge_{\edgee}(t)$ is constant for any edge $\edgee$ of $\CE(\CT_1) \setminus \CE(x)$. Define $x(.) \colon \left]0,\infty\right[ \to \left]0,\infty\right[^{E(\CT)}$ so that the variables $\varedge_{\edgee}(t)$ for $\edgee \notin \CE(x)$ are constant (and different from zero) and the variables $\varedge_{\edgee}(t)$ for $\edgee \in \CE(\CT_1)$ are the same for $x(t)$ and $x_1(t)$.
Let $L_2= L(\CT) \setminus L_1$. For $\leafl_2 \in L_2$, we have $\mape(\CE(\eltx),\leafl_2)=\mape(\CE(\eltx),\leafl_1(\leafl_2))$ for a unique $\leafl_1(\leafl_2)$ of $L_1$.
The equation $\ast(\leafl_2,\leafl_1(\leafl_2))$ between nonvanishing constant products holds for $x(t)$ for any $t \in \left]0,\infty\right[$, and it implies that $\prodedge\bigl(\left[\leafl_2,\topt\right]\bigr)(t) =\prodedge\bigl(\left[\leafl_1(\leafl_2),\topt\right]\bigr)(t)=\prodedge\bigl(\left[\leafl_0,\topt\right]\bigr)(t)$. So $x(.)$ is valued in $\mathring{\setparamx}(\CT)$. Its limit at $0$ is $x$. So $x$ is in the closure of $\mathring{\setparamx}(\CT)$.

Let $\CE_2=\mape(L_2)$ be the set of edges adjacent to the leaves of $L_2$. Set $\CE_3=E(\CT) \setminus (E(\CT_1) \cup \CE_2)$. Then any element $x^{\prime}$ of ${\setparamx}_{\CE(\eltx)}(\CT)$ is determined by its projection $x^{\prime}_1 \in {\setparamx}_{\CE(\eltx)}(\CT_1)$ and by the free nonzero variables associated to the edges of $\CE_3$. More precisely, for an edge $\edgee = \mape(\leafl_2 \in L_2)$, 
the equation $\ast(\leafl_2,\leafl_1(\leafl_2) )$ 
between nonvanishing products 
determines $\varedge^{\prime}_{\edgee}$ as a function of $x^{\prime}_1$ and the free nonzero variables associated to the edges of $\CE_3$.
For elements $\bigl((\varedge_e)_{e \in E(\CT)}\bigr) \in \{0\}^{\CE(\eltx)} \times \left]0,\infty\right[^{ E(\CT)\setminus \CE(\eltx)}$, if the equations $\ast(\leafl_2,\leafl_1(\leafl_2) )$  are satisfied for all $\leafl_2 \in L_2$, then all the equations $\ast(\leafl_2,\leafl^{\prime})$ for $\leafl_2 \in L_2$ and $\leafl^{\prime} \in L(\CT)$ are satisfied, as we prove below. Let $\leafl_2 \in L_2$, and let $\nodemaj( \leafl_1(\leafl_2))$ denote the closest node to $\mape(\CE(\eltx),\leafl_2)=\mape(\CE(\eltx),\leafl_1(\leafl_2))$ in $\left[\leafl_1(\leafl_2),\mape\bigl(\CE(\eltx),\leafl_2\bigr)\right]$. See Figure~\ref{figleaves}. Then $\ast(\leafl_2,\leafl_1(\leafl_2) )$ is equivalent to
$\prodedge\bigl(\left[\leafl_2,\nodemaj(\leafl_1(\leafl_2))\right]\bigr) =\prodedge\bigl(\left[\leafl_1(\leafl_2),\nodemaj(\leafl_1(\leafl_2)\right]\bigr)$.
In particular, if $\mape(\CE(\eltx),\leafl_2)=\mape(\CE(\eltx),\leafl^{\prime})$ for $\leafl^{\prime} \in L(\CT)$ as in Figure~\ref{figleaves}, then $\ast(\leafl^{\prime},\leafl_1(\leafl_2) )$ is equivalent to  
$\prodedge\bigl(\left[\leafl^{\prime},\nodemaj(\leafl_1(\leafl_2))\right]\bigr) =\prodedge\bigl(\left[\leafl_1(\leafl_2),\nodemaj(\leafl_1(\leafl_2)\right]\bigr)$.
So $\ast(\leafl_2,\leafl_1(\leafl_2) )$ and $\ast(\leafl^{\prime},\leafl_1(\leafl_2) )$ imply $\ast(\leafl_2,\leafl^{\prime})$.
If $\mape(\CE(\eltx),\leafl_2) \neq \mape(\CE(\eltx),\leafl^{\prime})$, then
$\mape(\CE(\eltx),\leafl_2) \in \left[\leafl_2,\nodemaj(\leafl_2,\leafl^{\prime})\right]$
and $\mape(\CE(\eltx),\leafl^{\prime}) \in \left[\leafl^{\prime},\nodemaj(\leafl_2,\leafl^{\prime})\right]$. So
$\ast(\leafl_2,\leafl^{\prime})$ is equivalent to $0=0$ and is satisfied.
Thus ${\setparamx}_{\CE(\eltx)}(\CT)$ is a smooth manifold whose dimension is 
$\cardlef{\CE_3} +\cardlef{E(\CT_1)}  -\cardlef{\CE(\eltx)}- \cardlef{L(\CT_1)} +\cardlef{\CE_1(x)}
=\cardlef{E(\CT)} - \cardlef{L_2}  -\cardlef{\CE(\eltx)}- \cardlef{L(\CT_1)} +\cardlef{\CE_1(x)}$.

\eop
\bfig
\centering

\begin{equation*}\begin{tikzpicture}
\begin{scope}[yscale=.4]
\draw[gray]  (-.8,4) -- (7.6,4)  (-.6,1) -- (-.6,2) -- (.9,3) (.9,1) -- (.9,2) -- (.9,3) (2.4,1) -- (2.4,2) -- (.9,3) -- (2,4) (3.4,1) -- (3.4,4) (4.8,1) -- (4.8,4);
\fill (-.8,4) ellipse (1.5pt and 4pt) (.6,4) ellipse (1.5pt and 4pt) (2,4) ellipse (1.5pt and 4pt) (3.4,4) ellipse (1.5pt and 4pt) (4.8,4) ellipse (1.5pt and 4pt) (6.2,4) ellipse (1.5pt and 4pt)  (-.6,1)  ellipse (1.5pt and 4pt) (-.6,2)  ellipse (1.5pt and 4pt) (.9,3)  ellipse (1.5pt and 4pt) (.9,1)  ellipse (1.5pt and 4pt) (.9,2)  ellipse (1.5pt and 4pt) (2.4,1)  ellipse (1.5pt and 4pt) (2.4,2)  ellipse (1.5pt and 4pt) (3.4,1)  ellipse (1.5pt and 4pt) (3.4,2.5) ellipse (1.5pt and 4pt) (4.8,1)  ellipse (1.5pt and 4pt) (4.8,2) ellipse (1.5pt and 4pt) (4.8,3) ellipse (1.5pt and 4pt) (7.6,4) circle (2.5pt);
\draw[dotted, gray] (6.9,4) -- (7.5,4);
\draw [gray,->] (2,4) -- (2.7,4);
\draw [gray,->] (-.8,4) -- (-.1,4);
\draw [gray,->] (.6,4) -- (1.3,4);
\draw [gray,->]
 (3.4,4) -- (4.1,4);
\draw [gray,->] (4.8,4) -- (5.5,4);
\draw [very thick,->] (6.2,4) -- (6.9,4);
\draw (-.8,4) node[left]{\scriptsize  $\leafl_1(\leafl_2)$} 
(6.9,4.6) node{\scriptsize  $\mape(\CE(x),\leafl_1(\leafl_2))$}  
(4.8,4.6)  node{\scriptsize  $\nodemaj(\leafl_1(\leafl_2))$} 
(-.6,1.5) node[right]{\scriptsize  $\mape(\leafl_2)$}
(6.9,3.2) node{\scriptsize  $=\mape(\CE(x),\leafl_2)$}
(6.9,2.3) node{\scriptsize  $=\mape(\CE(x),\leafl^{\prime}_2)$}
(2,4.6) node{\scriptsize  $\nodemaj(\leafl_2,\leafl_1(\leafl_2))$}
(1,3.2) node[left]{\scriptsize  $\nodemaj(\leafl_2,\leafl^{\prime}_1)$}
(2.4,.9) node[right]{\scriptsize $\leafl^{\prime}_1$} 
(3.4,.9) node[right]{\scriptsize $\leafl^{\prime}_2$} 
(-.6,.9) node[left]{\scriptsize$\leafl_2$}
(4.8,.9) node[right]{\scriptsize$\leafl^{\prime}_3$};
\draw[->]  (-.6,1) -- (-.6,1.5);
\draw[gray,->] (.9,1) -- (.9,1.5);
\draw[gray,->] (2.4,1) -- (2.4,1.5);
\draw[gray,->] (3.4,1) -- (3.4,1.75);
\draw[gray,->] (4.8,1) -- (4.8,1.5);
\draw[gray,->] (4.8,2) -- (4.8,2.5);
\draw[gray,->]  (4.8,3) -- (4.8,3.5);
\draw[gray,->]  (3.4,2.5) -- (3.4,3.25);
\draw[gray,->] (.9,3) -- (1.45,3.5);
\draw[gray,->] (2.4,2) -- (1.65,2.5);
\draw[gray,->]  (.9,2) -- (.9,2.5);
\draw[gray,->]  (-.6,2) -- (.15,2.5);
\draw[gray,->] (3.4,4) -- (4.8,4) (3.4,4) -- (4.1,4);
\draw[gray,->] (4.8,2) -- (4.8,3) (4.8,2) -- (4.8,2.5);
\end{scope}
\end{tikzpicture}\end{equation*}
\caption{Example of leaves $\leafl^{\prime}_i$ such that $\mape(\CE(x),\leafl^{\prime}_i)=\mape(\CE(x),\leafl_2)$ for the proof of Lemma~\ref{lemclosuresetparamx}}
\label{figleaves}

\end{figure}

\begin{lemma}
\label{lemcodimonex}
The codimension-one faces of ${\setparamx}(\CT)$ are in one-to-one correspondence with the codimension-one systems of edges of $\CT$. In a neighborhood of such an open face, ${\setparamx}(\CT)$ has the structure of a smooth manifold with boundary.
\end{lemma}
\bp According to Lemma~\ref{lemclosuresetparamx} above, if ${\setparamx}_{\CE}(\CT)$ is a nonempty manifold of dimension $\dimd(\CT)-1$, then $\CE$ is a codimension-one system. Let $\CE_c$ be a codimension-one system of edges of $\CT$.
For any edge of $e_0 \in \CE_c$, $\CE_c \setminus \{e_0\}$ can be completed to a reducing system $\CE_r$ that does not contain $e_0$ as in Lemma~\ref{lemcodimoneone}.
In particular, $\mathring{\setparamx}(\CT)$ is freely parametrized by the variables associated to the edges of $E(\CT) \setminus \CE_r$. When all these variables are nonzero except maybe $(\varedge_{e_0} \in \left[0,\varepsilon\right[)$, we get a local parametrization by $\left[0,\varepsilon\right[ \times \left]0,\varepsilon\right[^{E(\CT) \setminus (\CE_r \cup \{e_0\})}$ near the locus $\varedge_{e_0}=0$.
If $e \in \CE_c \setminus \{e_0\}$ and $e=\mape(\CE_r,\leafl_e)$, then
the variable $u_e$ may be expressed as \begin{equation*}u_e=\frac{\prod_{{\edgef} \in \left[\leafl_0(\CE_r),\nodemaj(\leafl_0(\CE_r),\leafl_e)\right] \setminus \{ e_0\}} \varedge_{\edgef}}
{\prod_{{\edgef} \in \left[\leafl_e,\nodemaj(\leafl_0(\CE_r),\leafl_e)\right] \setminus \{ e\}} \varedge_{\edgef}}\varedge_{e_0}\end{equation*}
with respect to this local parametrization.
\eop

We finish this section by proving a version of Stokes' theorem for spaces modelled by products of ${\setparamx}(\CT)$ by a manifold. 
Let us first introduce its statement given in Theorem~\ref{thmStokesonsetparamx}.

\begin{lemma}\label{lemconvonsetparamx} 
Let $\varepsilon \in \left]0,\infty\right[$. Let ${\setparamx}^{\varepsilon}(\CT)= {\setparamx}(\CT) \cap \left[0,\varepsilon\right]^{E(\CT)}$. Assume that ${\setparamx}(\CT)$ is oriented. Let $n\in \NN$. 
For any smooth form $\Omega$ on $\left[0,\varepsilon\right]^{E(\CT)} \times \left[0,1\right]^n$
of degree $(\dimd(\CT)+n)$, the integral $\int_{{\setparamx}^{\varepsilon}(\CT)\times \left[0,1\right]^n} \Omega$ of $\Omega$ along the interior of ${\setparamx}^{\varepsilon}(\CT)\times \left[0,1\right]^n$ is absolutely convergent.
\end{lemma}
For an ordered subset $\CE$  of $E(\CT)$, we set $\Omega_{\CE}=\wedge_{e \in \CE}d\varedge_e$, where the factors are ordered with respect to the order of $\CE$.

\bpo{Proof of lemma~\ref{lemconvonsetparamx}}
Any smooth form $\Omega$ on $\left[0,\varepsilon\right]^{E(\CT)} \times \left[0,1\right]^n$
of degree $(\dimd(\CT)+n)$ is a sum of forms $g_{\CE} \Omega_{\CE} \wedge (\wedge_{i=1}^n dx_i)$, for ordered subsets $\CE$ of $E(\CT)$ of cardinality $\dimd(\CT)$ and for smooth maps $g_{\CE}\colon \left[0,\varepsilon\right]^{E(\CT)} \times \left[0,1\right]^n \to \RR$,
up to forms that vanish identically on the interior of ${\setparamx}^{\varepsilon}(\CT)\times \left[0,1\right]^n$.
These forms are bounded on $\left[0,\varepsilon\right]^{E(\CT)} \times \left[0,1\right]^n$. They are zero on the interior of ${\setparamx}^{\varepsilon}(\CT)\times \left[0,1\right]^n$ unless $\CE$ is a maximal free system, according to Lemma~\ref{lemsysred}. When $\CE$ is a maximal free system, Lemma~\ref{lemsysred} ensures that $\mathring{\setparamx}^{\varepsilon}(\CT)$ is freely parametrized by the variables associated to the edges of $\CE$.
\eop

We restrict to the compact subspace ${\setparamx}^{\varepsilon}(\CT)$ of ${\setparamx}(\CT)$. 
The subspace ${\setparamx}^{\varepsilon}(\CT)$ is the closure of the open    manifold $\mathring{\setparamx}(\CT) \cap \left]0,\varepsilon\right[^{E(\CT)}$.
For a maximal free system $\CE_b$ of edges, $\mathring{\setparamx}(\CT)$ is freely parametrized by the variables $u_{\edgef}$ for $\edgef \in \CE_b$ as in Lemma~\ref{lemsysred}. For $e_1 \in
E(\CT) \setminus \CE_b$, the variable $u_{e_1}$ is the function
$u_{e_1}(\CE_b; V \in \left]0,\infty\right[^{\CE_b})$ of Lemma~\ref{lemsysredcoord}.
This function is a smooth quotient of monomials in the $u_{\edgef}$ for ${\edgef} \in \CE_b$.
Then $\mathring{\setparamx}(\CT) \cap  \left]0,\varepsilon\right]^{E(\CT)}$ is parametrized by the set of elements $V \in \left]0,\varepsilon\right]^{\CE_b}$ such that $u_e(\CE_b;U)\leq \varepsilon$ for any $e \in
E(\CT) \setminus \CE_b$. It is a subspace of $\left]0,\varepsilon\right]^{\CE_b}$ whose boundary can be stratified so that the open strata of the boundary are locally defined by an equation 
$u_{\edgee}(\CE_b; V \in \left]0,\infty\right[^{\CE_b})=\varepsilon$.
Let ${\partial}_{\varepsilon}{\setparamx}^{\varepsilon}(\CT)$ be the union of the corresponding open codimension-one faces of ${\setparamx}^{\varepsilon}(\CT)$.
Then the union $\check{\partial}{\setparamx}^{\varepsilon}(\CT)$ of the oriented open codimension-one faces of ${\setparamx}^{\varepsilon}(\CT)$ is \begin{equation*}\check{\partial}{\setparamx}^{\varepsilon}(\CT) = {\partial}_o{\setparamx}^{\varepsilon}(\CT) \cup {\partial}_{\varepsilon}{\setparamx}^{\varepsilon}(\CT), \end{equation*}
where ${\partial}_o{\setparamx}^{\varepsilon}(\CT)$ is the union of the intersection of the codimension-one faces of ${\setparamx}(\CT)$ of Lemma~\ref{lemcodimonex} with ${\setparamx}^{\varepsilon}(\CT)$.

\begin{theorem}\label{thmStokesonsetparamx} 
Let $\varepsilon \in \left]0,\infty\right[$. Recall ${\setparamx}^{\varepsilon}(\CT)= {\setparamx}(\CT) \cap \left[0,\varepsilon\right]^{E(\CT)}$. Assume that ${\setparamx}(\CT)$ is oriented. 
Let $n\in \NN$. Let $\omega$ be a smooth form on $\left[0,\varepsilon\right]^{E(\CT)} \times \left[0,1\right]^n$
of degree $(\dimd(\CT) -1 +n)$.
With the notation above, consider the dense part 
\begin{equation*}\check{\partial}\Bigl(\setparamx^{\varepsilon}(\CT)\times \left[0,1\right]^n\Bigr)=\left(\check{\partial}{\setparamx}^{\varepsilon}(\CT)\right)\times \left[0,1\right]^n \cup \left( (-1)^{\dimd(\CT)}\mathring{\setparamx}^{\varepsilon}(\CT) \times \partial \left( \left[0,1\right]^n\right)\right)\end{equation*}
of $\partial\left(\setparamx^{\varepsilon}(\CT)\times \left[0,1\right]^n\right)$.

Then the integral $\int_{\partial(\setparamx^{\varepsilon}(\CT)\times \left[0,1\right]^n)} \omega$ of $\omega$ along $\check{\partial}\left(\setparamx^{\varepsilon}(\CT)\times \left[0,1\right]^n\right)$ is absolutely convergent. Furthermore, Stokes' theorem
applies to this setting. So we have
\begin{equation*}\int_{\partial(\setparamx^{\varepsilon}(\CT)\times \left[0,1\right]^n)} \omega= \int_{{\setparamx}^{\varepsilon}(\CT)\times \left[0,1\right]^n} d\omega.\end{equation*}
\end{theorem}

Let us prepare the proof of Theorem~\ref{thmStokesonsetparamx} with a few lemmas.

\begin{lemma}
Let $\CE_1$ be an ordered subset of $E(\CT)$ of cardinality $\dimd(\CT)$. Let $g$ be a smooth function  on $\left[0,\varepsilon\right]^{E(\CT)} \times \left[0,1\right]^n$.
Then Theorem~\ref{thmStokesonsetparamx} holds for 
$\omega =g \Omega_{\CE_1} \wedge (\wedge_{i=2}^n dx_i)$.
\end{lemma}
\bp
We have
\begin{equation*}\begin{array}{ll}\int_{{\setparamx}^{\varepsilon}(\CT)\times \left[0,1\right]^n} d\omega &=(-1)^{\dimd(\CT)}\int_{\mathring{\setparamx}^{\varepsilon}(\CT)\times \left[0,1\right]^n} \frac{\partial g}{\partial x_1} \Omega_{\CE_1} \wedge (\wedge_{i=1}^n dx_i)\\
&= (-1)^{\dimd(\CT)} \int_{\mathring{\setparamx}^{\varepsilon}(\CT)\times \partial \left[0,1\right] \times  \left[0,1\right]^{n-1}} \omega\\
&= \int_{\partial\left({\setparamx}^{\varepsilon}(\CT)\times \left[0,1\right]^n\right)} \omega \end{array}\end{equation*}
since $\mathring{\setparamx}^{\varepsilon}(\CT)\times \partial \left[0,1\right] \times  \left[0,1\right]^{n-1}$
is the only part of $\partial\left({\setparamx}^{\varepsilon}(\CT)\times \left[0,1\right]^n\right)$ where the integral of $\omega$ does not vanish.
\eop

So it suffices to prove Theorem~\ref{thmStokesonsetparamx} for forms
$\omega =g \Omega_{\CE_2} \wedge (\wedge_{i=1}^n dx_i)$ for ordered subsets $\CE_2$ of $E(\CT)$ of cardinality $\dimd(\CT)-1$ such that $\Omega_{\CE_2}$ is not zero on $\mathring{\setparamx}^{\varepsilon}(\CT)$. (If $\Omega_{\CE_2}$ is zero on $\mathring{\setparamx}^{\varepsilon}(\CT)$, then it is also zero on ${\setparamx}^{\varepsilon}(\CT)$, and both sides of the equality to be proved are zero.)

\begin{lemma}\label{lemconvonsetparamxtwo}
Let $\CE_2$ be a subspace of $E(\CT)$ of cardinality $\dimd(\CT)-1$ such that $\Omega_{\CE_2}$ is not zero on $\mathring{\setparamx}^{\varepsilon}(\CT)$.
Let $\CE_2^c = E(\CT) \setminus \CE_2$. 
Then either $L(\CE_2^c)=L(\CT)$ or $L(\CE_2^c)=L(\CT) \setminus \{\leafl_0\}$ for a unique $\leafl_0 \in L(\CT)$.

If $L(\CE_2^c)=L(\CT)$, then $\mape(\CE_2^c,.)$ is a bijection from  $L(\CE_2^c)$ to $\CE_2^c$.

If $L(\CE_2^c)=L(\CT) \setminus \{\leafl_0\}$, then $\mape(\CE_2^c,.)$ is a bijection from  $L(\CE_2^c)$ to $\CE_2^c \setminus \{\edgef\}$ for a unique element $\edgef$ of $\CE_2^c$.
\end{lemma}
\bp
The set $L(\CT) \setminus L(\CE_2^c)$ cannot contain two distinct leaves $\ell_3$ and $\ell_4$, because Equation $\ast(\leafl_3,\leafl_4)$ would imply that $\Omega_{\CE_2}=0$.
Similarly, $\mape(\CE_2^c,.)$ cannot map two distinct leaves $\ell_3$ and $\ell_4$ of $L(\CE_2^c)$ to the same element. So $\mape(\CE_2^c,.)$ is injective.
Use that the cardinality of $\CE_2^c$ is $\cardlef{L(\CT)}$ to conclude.
\eop

For a subset $\CE^{\prime}$ of $\CE$, let $p_{\CE^{\prime}} \colon {\setparamx}^{\varepsilon}(\CT) \to \left[0,\varepsilon\right]^{\CE^{\prime}}$ denote the composition of the inclusion ${\setparamx}^{\varepsilon}(\CT) \hookrightarrow \left[0,\varepsilon\right]^{E(\CT)}$ with the natural projection.

\begin{lemma}\label{lemconvonsetparamxthree} Recall the assumptions of Lemma~\ref{lemconvonsetparamxtwo}.
Let $\CE_3$ be the set of edges ${\edgef}$ of $\CT$ such that 
$\CE_2 \cup \{\edgef\}$ is a maximal free system of edges of $\CT$.
For an edge $f$ of $\CE_3$ and an element $(u_f,U) \in \left]0,\varepsilon\right]^{\{\edgef\} \cup\CE_2}$, let $x_f(u_f,U)$ denote the element of $\mathring{\setparamx}(\CT)$ such that $p_{\{\edgef\} \cup\CE_2}(x_f(u_f,U))= (u_f,U)$.

Then there exist ${\edgef_2} \in \CE_3$ and piecewise smooth functions
$\funca$ and $\funcb$ from $\left]0,\varepsilon\right]^{\CE_2}$ to $\left[0,\varepsilon\right]$, such that ${\setparamx}^{\varepsilon}(\CT) \cap p_{\CE_2}^{-1}\left(]0,\varepsilon]^{\CE_2} \right)$ is the set \begin{equation*}\Bigl\{x_{\edgef_2}(u_{\edgef_2},U) \suchthat U \in \left]0,\varepsilon\right]^{\CE_2} \cap p_{\CE_2}\bigl({\setparamx}^{\varepsilon}(\CT)\bigr), u_{\edgef_2} \in \left[\funca(U), \funcb(U)\right]\Bigr\}.\end{equation*}
Furthermore, the space $\left]0,\varepsilon\right]^{\CE_2} \cap p_{\CE_2}\left({\setparamx}^{\varepsilon}(\CT)\right)$ is a subspace of
$\left]0,\varepsilon\right]^{\CE_2}$ determined by the condition $b(U)-a(U)\geq 0$ and conditions
$u_e(U) \leq \varepsilon$ for $e \notin \CE_3$, for functions $u_e$ of $U$ 
which are quotients of monomials in the variables $u_g$ for $g\in \CE_2$.
Moreover, the boundary of ${\setparamx}^{\varepsilon}(\CT)$ consists of
\begin{itemize}
 \item $(d(\CT)-1)$-dimensional strata where $u_{\edgez}$ vanishes for some ${\edgez} \in \CE_2$, 
 \item (negligible) strata of dimension less than $(d(\CT)-1)$, and
 \item $(d(\CT)-1)$-dimensional strata of 
\begin{equation*}\partial_{\funcb}\bigl({\setparamx}^{\varepsilon}(\CT)\bigr)
=\Bigl\{x_{\edgef_2}\bigl(\funcb(U),U\bigr) \suchthat U \in \left]0,\varepsilon\right]^{\CE_2} \cap p_{\CE_2}\bigl({\setparamx}^{\varepsilon}(\CT)\bigr)\Bigr\},\end{equation*}
and
\begin{equation*}-\partial_{\funca}\bigl({\setparamx}^{\varepsilon}(\CT)\bigr)
=-\Bigl\{x_{\edgef_2}\bigl(\funca(U),U\bigr) \suchthat U \in \left]0,\varepsilon\right]^{\CE_2} \cap p_{\CE_2}\bigl({\setparamx}^{\varepsilon}(\CT)\bigr) \Bigr\},\end{equation*}
which behave as standard codimension-one faces of ${\setparamx}^{\varepsilon}(\CT)$,\footnote{The signs before $\partial_{\funcb}$ and $\partial_{\funcb}$ correspond to the orientations provided that $du_{\edgef_2} \wedge \Omega_{\CE_2}$ orients ${\setparamx}^{\varepsilon}(\CT)$, which we assume without loss.}
\end{itemize}
with respect to some natural stratification.
\end{lemma}

\bpo{Proof of Theorem~\ref{thmStokesonsetparamx} assuming Lemma~\ref{lemconvonsetparamxthree}}
Consider a form $\omega =g \Omega_{\CE_2} \wedge (\wedge_{i=1}^n dx_i)$, for an ordered subset $\CE_2$ of $E(\CT)$  of cardinality $\dimd(\CT)-1$
and a smooth function $g \colon \left[0,\varepsilon\right]^{E(\CT)} \times \left[0,1\right]^n \to \RR$.
Assume that the $d\varedge_{\edgez}$, for ${\edgez} \in \CE_2$, are linearly independent in the $C^{\infty}(\mathring{\setparamx}^{\varepsilon}(\CT);\RR)$-module $\Omega^1(\mathring{\setparamx}^{\varepsilon}(\CT))$.
The integral of \begin{equation*}d\omega =\sum_{h \in E(\CT)} \frac{\partial g}{\partial \varedge_{\edgeh}} d\varedge_{\edgeh} \wedge\Omega_{\CE_2}\wedge (\wedge_{i=1}^n dx_i)\end{equation*} is absolutely convergent over the interior of
${\setparamx}^{\varepsilon} \times \left[0,1\right]^n$. 

With the notation of Lemma~\ref{lemconvonsetparamxthree}, we have
\begin{equation*}\int_{\partial(\setparamx^{\varepsilon}(\CT)\times \left[0,1\right]^n)} \omega= \int_{\partial_{\funcb}({\setparamx}^{\varepsilon}(\CT))\times \left[0,1\right]^n)}\omega - \int_{\partial_{\funca}({\setparamx}^{\varepsilon}(\CT))\times \left[0,1\right]^n)}\omega\end{equation*}
since $\omega$ vanishes along the $(d(\CT)-1)$-dimensional strata where $u_{\edgez}$ vanishes for some ${\edgez} \in \CE_2$.
On the other hand, we have
\begin{equation*}\int_{\setparamx^{\varepsilon}(\CT)\times \left[0,1\right]^n} d\omega= \int_{\left[0,1\right]^n}\left(\int_{\mathring{\setparamx}^{\varepsilon}(\CT)}\frac{\partial g_X}{\partial u_{\edgef_2} }du_{\edgef_2} \wedge \Omega_{\CE_2} \right) \wedge_{i=1}^n dx_i,\end{equation*}
where $g_{\setparamx}$ is the restriction of $g$ to ${\setparamx}^{\varepsilon}(\CT)\times \left[0,1\right]^n$, where $g_{\setparamx}$ and the $\varedge_{\edgeh}$ are functions of $\varedge_{\edgef_2}$ and of the $\varedge_{\edgez}$, for ${\edgez} \in \CE_2$. For a fixed implicit $(x_1, \dots,x_n)$, we compute
\begin{equation*}\begin{array}{ll}\int_{\mathring{\setparamx}^{\varepsilon}(\CT)}\frac{\partial g_X}{\partial u_{\edgef_2} }du_{\edgef_2} \wedge \Omega_{\CE_2}&=
\int_{ U \in \left]0,\varepsilon\right]^{\CE_2} \cap p_{\CE_2}\left({\setparamx}^{\varepsilon}(\CT)\right)} \left(\int_{u_{\edgef_2} \in \left[\funca(U), \funcb(U)\right]} \frac{\partial g_X}{\partial u_{\edgef_2} }du_{\edgef_2}\right) \Omega_{\CE_2}\\
&= \int_{ U \in \left]0,\varepsilon\right]^{\CE_2} \cap p_{\CE_2}\left({\setparamx}^{\varepsilon}(\CT)\right)} \bigl(g_X(b(U))-g_X(a(U))\bigr) \Omega_{\CE_2}.\end{array}\end{equation*}
\eop

Lemma~\ref{lemconvonsetparamxtwo} reduces the proof of Lemma~\ref{lemconvonsetparamxthree} to the proofs of the following two lemmas.

\begin{lemma}\label{lemconvonsetparamxthreeone}Lemma~\ref{lemconvonsetparamxthree} holds when $L(\CE_2^c)=L(\CT) \setminus \{\leafl_0\}$.
\end{lemma}

 \begin{lemma}\label{lemconvonsetparamxthreetwo}
 Lemma~\ref{lemconvonsetparamxthree} holds when  $L(\CE_2^c)=L(\CT).$
\end{lemma}

\bpo{Proof of Lemma~\ref{lemconvonsetparamxthreeone}}
Recall that $L(\CE_2^c)=L(\CT) \setminus \{\leafl_0\}$ and $\mape(\CE_2^c,.)$ is a bijection from  $L(\CE_2^c)$ to $\CE_2^c \setminus \{\edgef\}$ for a unique element $\edgef$ of $\CE_2^c$ from Lemma~\ref{lemconvonsetparamxtwo}.
There is a leaf $\leafl_1$ such that $\edgef \in \left[\leafl_1,\topt\right]$.
Let $\edgef_2$ be the last edge (the closest to $\edgef$) of $\CE_2^c$ in $\left[\leafl_1,\edgef\right[$. Then we have $\edgef_2=\mape(\CE_2^c,\leafl_2)$ for some leaf $\leafl_2$ of $L(\CE_2^c)$. If there is an edge of $\CE_2^c$ in $\left]\edgef,\topt\right]$, define the leaf $\leafl_3$
such that $\mape(\CE_2^c,\leafl_3)$ is the first edge of $\CE_2^c$ in $\left]\edgef,\topt\right]$ as in Figure~\ref{figStokes1}.
Otherwise, set $\leafl_3=\leafl_0$. Note that $\edgef$ and $\edgef_2$ are in $\CE_3$.

\bfig
\centering
\begin{equation*}
\begin{tikzpicture}
\draw [dashed] (.6,2) -- (4.6,2) (4.6,2.6) node[left]{\scriptsize $\leafl_0$} (4.6,2.6) -- (5.8,2) (.7,.8) node[left]{\scriptsize $\leafl_2$} (.7,.8) -- (1.6,2) (2.1,.8) -- (3,2) (3.7,.8) node[left]{\scriptsize $\leafl_3$} (3.7,.8) -- (4.6,2);
\draw (.6,2) -- (1.2,2) (2,2) -- (2.6,2) (3.4,2) -- (4,2) (4.6,2) -- (5.2,2) (2.4,1.2) -- (2.7,1.6);
\draw [dotted] (0,2) -- (.6,2) (5.8,2) node[right]{\scriptsize $\topt$} (5.2,2) -- (5.8,2);
\draw [->] (0,2) node[left]{\scriptsize $\leafl_1$} (.6,2) -- (.9,2);
\draw [->] (2.3,2.25) node{\scriptsize $\edgef_2$} (2,2) -- (2.3,2);
\draw [->] (3.7,2.25) node{\scriptsize $\edgef$} (3.4,2) -- (3.7,2);
\draw [->] (5.15,1.75) node{\scriptsize $\mape(\CE_2^c,\leafl_3)$}  (4.6,2) -- (4.9,2);
\draw [->] (2.4,1.2) -- (2.55,1.4);
\fill  (.6,2) circle (1.5pt) (1.2,2) circle (1.5pt) (2,2) circle (1.5pt) (2.6,2) circle (1.5pt) (3.4,2) circle (1.5pt) (4,2) circle (1.5pt) (4.6,2) circle (1.5pt) (5.2,2) circle (1.5pt) (2.4,1.2) circle (1.5pt) (2.7,1.6) circle (1.5pt) (4.6,2.6) circle (1.5pt) (5.8,2) circle (1.5pt) (0,2) circle (1.5pt) (.7,.8) circle (1.5pt) (2.1,.8) circle (1.5pt) (3.7,.8) circle (1.5pt) (1.6,2) circle (1.5pt) (3,2) circle (1.5pt);
\end{tikzpicture}\end{equation*}
\caption{The edge $f$, when $L(\CE_2^c) \neq L(\CT)$, in the proof of Theorem~\ref{thmStokesonsetparamx}}
\label{figStokes1}

\end{figure}
The product $\varedge_{\edgef_2}\varedge_{\edgef}$ is given by the expression 
\begin{equation*}\varedge_{\edgef_2}\varedge_{\edgef}=\frac{\prodedge\bigl(\left[\leafl_3,\nodemaj(\ell_2,\ell_3)\right]\bigr)}{\prod_{{\edgez} \in \left[\leafl_2,\nodemaj(\ell_2,\ell_3)\right] \setminus \{\edgef,\edgef_2\}}\varedge_{\edgez}}=\eta_2(U),\end{equation*}
in $U \in \left]0,\varepsilon\right]^{\CE_2} \cap p_{\CE_2}\left({\setparamx}^{\varepsilon}(\CT)\right)$. 
In particular, we have $\varedge_{\edgef}={\eta_2(U)}/{\varedge_{\edgef_2}}$, and the conditions $\varedge_{\edgef} \leq \varepsilon$ and $\varedge_{\edgef_2} \leq \varepsilon$ imply \begin{equation*}\frac{\eta_2(U)}{\varepsilon} \leq \varedge_{\edgef_2} \leq \varepsilon.\end{equation*}
Let $\edgee$ be an edge of $\CE_2^c$ different from $\edgef$. Then $\edgee$
may be expressed as $\edgee =\mape(\CE_2^c\setminus \{\edgef\},\leafl_e)$ for a unique $\leafl_e=\leafl(\CE_2^c\setminus \{\edgef\},\edgee)$.
As in Lemma~\ref{lemsysredcoord}, there is a leaf
 $\leafl^{\prime}=\leafl^{\prime}(\CE_2^c\setminus \{\edgef\},\edgee)$ of $\CT$ such that $\edgee$ is the only element of $\CE_2^c\setminus \{\edgef\}$ in $\left[\leafl_e,\nodemaj(\leafl_e,\leafl^{\prime})\right] \cup \left[\leafl^{\prime},\nodemaj(\leafl_e,\leafl^{\prime})\right]$.
 In particular, if $\edgee\neq\edgef_2 $, then $\edgef_2$ does not belong to $\left[\leafl_e,\nodemaj(\leafl_e,\leafl^{\prime})\right] \cup \left[\leafl^{\prime},\nodemaj(\leafl_e,\leafl^{\prime})\right]$. Note that $\edgef$ cannot be on $\left[\leafl^{\prime},\nodemaj(\leafl_e,\leafl^{\prime})\right]$.
 If $\edgef \notin \left[\leafl_e,\nodemaj(\leafl_e,\leafl^{\prime})\right]$, then $\varedge_e$ is a function of $U \in \left]0,\varepsilon\right]^{\CE_2}$, $\edgee \notin \CE_3$, and $\varedge_e$ is different from zero.
\bfig
\centering
\begin{equation*}
\begin{tikzpicture}
\draw [dashed] (-1.5,3.2) -- (.45,2);
\draw (-.85,2.8) -- (-.2,2.4);
\draw [dashed] (-1.5,2) -- (2,2);
\draw [->] (-.4,2.7) node{\scriptsize $\edgef_2$} (-.85,2.8) -- (-.525,2.6);
\draw (-1.7,3.2) node{\scriptsize $\leafl_2$} (-1.7,2) node{\scriptsize $\leafl_e$} (-1.7,.8) node{\scriptsize $\leafl^{\prime}$} (-.6,2) -- (0,2) (.9,2) -- (1.5,2) ;
\draw [->] (-.3,2.25) node{\scriptsize $\edgee$} (-.6,2) -- (-.3,2);
\draw [->] (1.2,2.25) node{\scriptsize $\edgef$} (.9,2) -- (1.2,2);
\fill (-1.5,2) circle (1.5pt) (-.6,2) circle (1.5pt) (0,2) circle (1.5pt) (.9,2) circle (1.5pt) (1.5,2) circle (1.5pt) (-1.5,.8) circle (1.5pt) (2,2) circle (1.5pt) (-1.5,3.2) circle (1.5pt) (.45,2) circle (1.5pt) (-.85,2.8) circle (1.5pt) (-.2,2.4) circle (1.5pt);
\draw [dashed] (-1.5,.8) -- (2,2);
\draw [dotted] (2,2) -- (2.4,2);
\fill (-1.5,.8) circle (1.5pt);
\draw (1.8,1.75) node[right]{\scriptsize $\nodemaj(\leafl_e,\leafl^{\prime})$};
\draw (.45,1.8) node{\scriptsize $\nodemaj(\leafl_e,\leafl_2)$};
\draw (0,.3) node{\scriptsize When $L(\CE_2^c) \neq L(\CT)$};
\draw (0,-.2) node{\scriptsize and $f \in \left[\leafl_e,\nodemaj(\leafl_e,\leafl^{\prime})\right]$};

\begin{scope}[xshift=4cm]
\draw [dashed] (0,2) -- (2.6,2);
\draw [dotted] (3.6,2) -- (4,2);
\draw (0,2.25) node{\scriptsize $\leafl_e$} (.3,.8) node{\scriptsize $\leafl$} (.6,2) -- (1.6,2) (2.6,2) -- (3.6,2);
\draw [->] (.9,2.25) node[right]{\scriptsize $e$} (.6,2) -- (1.1,2);
\draw [->] (3.1,1.75) node{\scriptsize $\mape(\CE_2^c,\leafl)$} (2.6,2) -- (3.1,2);
\fill (0,2) circle (1.5pt) (.6,2) circle (1.5pt) (1.6,2) circle (1.5pt) (2.1,2) circle (1.5pt) (2.6,2) circle (1.5pt) (3.6,2) circle (1.5pt) (0,.8) circle (1.5pt);
\draw [dashed] (0,.8) -- (2.1,2);
\draw (2,2.25) node{\scriptsize $\nodemaj(\leafl_e,\leafl)$};
\draw (2,.3) node{\scriptsize When $L(\CE_2^c) = L(\CT)$};
\draw (2,-.2) node{\scriptsize and $e \in \CE_2^c \setminus \CE_3$};
\end{scope}

\end{tikzpicture}\end{equation*}
\caption{More figures for the proof of Theorem~\ref{thmStokesonsetparamx}}
\label{fighf}

\end{figure}
If $\edgef \in \left[\leafl_e,\nodemaj(\leafl_e,\leafl^{\prime})\right]$, then $\nodemaj(\leafl_e,\leafl_2)$ is on $\left]\edgef_2,\edgef\right[$ (equivalently, $\edgef_2 \in \left[\leafl_2,\nodemaj(\leafl_e,\leafl_2)\right]$ and  $f \notin \left[\leafl_2,\nodemaj(\leafl_e,\leafl_2)\right]$), as in the left part of Figure~\ref{fighf}, $\edgee \in \CE_3$,
$\edgef_2$ is the unique element of  $\CE_2^c$ in $\left[\leafl_2,\nodemaj(\leafl_e,\leafl_2)\right]$,
and we have \begin{equation*}\varedge_e=\frac{\prodedge\bigl(\left[\leafl_2,\nodemaj(\leafl_e,\leafl_2)\right]\bigr)}{\prod_{\edgez \in \left[\leafl_e,\nodemaj(\leafl_e,\leafl_2)\right] \setminus \{e\}}\varedge_{\edgez}}.\end{equation*} So,
when $U \in \left]0,\varepsilon\right]^{\CE_2}\cap p_{\CE_2}\left({\setparamx}^{\varepsilon}(\CT)\right)$ is fixed, $\varedge_e$ is a linear function  $\varedge_{\edgee}=\eta_{\edgee,\edgef_2}(U) \varedge_{\edgef_2}$ of $\varedge_{\edgef_2}$.
In particular, we have $\varedge_{\edgef_2} \leq {\varepsilon}/{\eta_{\edgee,\edgef_2}(U)}$.
Set $\funca(U)= {\eta_2(U)}/{\varepsilon}$ and $\funcb(U)=\min_{\edgee \in \CE_3 \setminus \{\edgef\}}({\varepsilon}/{\eta_{\edgee,\edgef_2}(U)})$, where $\eta_{\edgef_2,\edgef_2}=1$. Then $U \in \left]0,\varepsilon\right]^{\CE_2}$ is in $p_{\CE_2}\left({\setparamx}^{\varepsilon}(\CT)\right)$ if and only if $\funca(U) \leq \funcb(U)$ and $\varedge_e(U) \leq \varepsilon$ for any $\edgee \in \CE_2^c \setminus \CE_3$. Furthermore, ${\setparamx}^{\varepsilon}(\CT) \cap p_{\CE_2}^{-1}\bigl(\left]0,\varepsilon\right]^{\CE_2} \bigr)$ is the set \begin{equation*}\Bigl\{x_{\edgef_2}(u_{\edgef_2},U) \suchthat U \in \left]0,\varepsilon\right]^{\CE_2} \cap p_{\CE_2}\bigl({\setparamx}^{\varepsilon}(\CT)\bigr), u_{\edgef_2} \in \left[\funca(U), \funcb(U)\right]\Bigr\}.\end{equation*} 
The boundary part $\partial_{\funca}({\setparamx}^{\varepsilon}(\CT))$, corresponds to $\varedge_{\edgef}=\varepsilon$. The boundary part $\partial_{\funcb}({\setparamx}^{\varepsilon}(\CT))$, along which $\varedge_{\edgee}=\varepsilon$ for some $\edgee\in \CE_3$, lies in the intersection of the nonsingular $\mathring{\setparamx}$ with loci $(\varedge_{\edgee}=\varepsilon)$ for some $\edgee \in \CE_3$. The locus $(\varedge_{\edgee}=\varepsilon)$ is transverse to $\mathring{\setparamx}$, for any $\edgee \in \CE_3$ since $\CE_2 \cup \{\edgee\}$ is a maximal free system for such an $\edgee$. So 
$\mathring{\setparamx} \cap (\varedge_{\edgee}=\varepsilon)$ is a manifold of dimension $(d(\CT)-1)$ for all
$\edgee \in \CE_3$. 
Let $\edgee$ and $\edgee^{\prime}$ be two distinct edges of $\CE_3$. Then $\varedge_{\edgee}=\varedge_{\edgee^{\prime}}=\varepsilon$ implies a nontrivial equation among variables $\varedge_{\edgeh}$ associated to edges $\edgeh$ of $\CE_2$ unless $\edgee$ and $\edgee^{\prime}$ are both adjacent to leaves and meet at a node, which implies that $\varedge_{\edgee}=\varedge_{\edgee^{\prime}}$ on $\mathring{\setparamx}$. We may assume that this configuration never occurs, without loss of generality, because the space ${\setparamx}(\CT)$ is canonically diffeomorphic to the space ${\setparamx}(\CT^{\prime})$ obtained from $\CT$ by erasing $\edgee^{\prime}$. Except for this special configuration, the loci $\mathring{\setparamx} \cap (\varedge_{\edgee}=\varedge_{\edgee^{\prime}}=\varepsilon)$ are of dimension less than $(d(\CT)-1)$, and we have a stratification of $\partial_{\funca}({\setparamx}^{\varepsilon}(\CT))$ and $\partial_{\funcb}({\setparamx}^{\varepsilon}(\CT))$, where the $(d(\CT)-1)$-dimensional strata are the loci of $\partial_{\funca}({\setparamx}^{\varepsilon}(\CT))$ and $\partial_{\funcb}({\setparamx}^{\varepsilon}(\CT))$ where $\varedge_{\edgee}=\varepsilon$ for (exactly) one $\edgee$ of $\CE_3$. These strata are smooth open $(d(\CT)-1)$-manifolds.
\eop

\bpo{Proof of Lemma~\ref{lemconvonsetparamxthreetwo}}
Let us now assume that $L(\CE_2^c)= L(\CT)$, and recall that $\mape(\CE_2^c,.)$ is a bijection from $L(\CT)$ to $\CE_2^c$ in this case.
The elements of $\CE_3$ are the edges $\edgee$ of $\CE_2^c$ such that there is no edge of $\CE_2^c$ on $\left]\edgee,\topt\right]$. In particular, $\CE_3$ is a codimension-one system of edges of $\CT$.

Let $\edgee=\mape(\CE_2^c,\leafl_e)$ be  an element of $\CE_2^c \setminus \CE_3$, let $\ell$ be a leaf of $\CT$ such that the first edge of $\CE_2^c$ on $\left]\edgee,\topt\right]$ is $\mape(\CE_2^c,\leafl)$. Then $\edgee$ is the only element of $\CE_2^c$ in $\left[\leafl_e,\nodemaj(\leafl,\leafl_e)\right] \cup \left[\leafl,\nodemaj(\leafl,\leafl_e)\right]$, as in Figure~\ref{fighf}. So $\varedge_{\edgee}$ depends only on the fixed variables of $\CE_2$, and it is not zero.

If $\edgef \in \CE_3$, let $\leafl(\edgef)$ denote the unique leaf such that $L(\CE_2^c \setminus \{\edgef\}) = L(\CT)\setminus \{\leafl(\edgef)\}$.
The variable $\varedge_{\edgef}$ is a linear function $\varedge_{\edgef}=\eta_{\edgef}(U)\varedge_{\edgef_2}$ of $\varedge_{\edgef_2}$ for one (arbitrary) $\edgef_2 \in \CE_3$, \begin{equation*}\eta_{\edgef}(U)=\frac{\prod_{{\edgez} \in \left[\leafl(\edgef_2),\topt\right] \setminus \{\edgef_2\}}\varedge_{\edgez}}{\prod_{{\edgez} \in \left[\leafl(f),\topt\right] \setminus \{\edgef\}}\varedge_{\edgez}}.\end{equation*}
Here, $\funca(U)=0$ and $\funcb(U)=\min_{\edgef \in \CE_3}({\varepsilon}/{\eta_{\edgef}(U)})$. So $\partial_{\funca}({\setparamx}^{\varepsilon}(\CT)) $ is the codimension-one face associated to $\CE_3$, as in Lemma~\ref{lemcodimonex}, while $\partial_{\funcb}({\setparamx}^{\varepsilon}(\CT))$ may be stratified as in the previous proof. 
\eop

The proof of Theorem~\ref{thmStokesonsetparamx} is finished. So we have the announced Stokes formula for all forms as in its statement.

\section{Configuration spaces of graphs on long tangles}
\label{secstrucbyhand}

An \emph{$\infty$-component} of an LTR $\tanghcyll$ is a connected component of the intersection of the image of $\tanghcyll$ with $\rats(\hcylc) \setminus \hcylc$. 
A univalent vertex $\eltv$ of a Jacobi diagram on the domain $\sourcetl$
of $\tanghcyll$ approaching $\infty$ moves on such an $\infty$-component.
Let $\confy_{\pm}(\eltv)$ denote the orthogonal projection on $\CC$ of the corresponding $\infty$-component of $\eltv$.
Then this $\infty$-component may be written as $\{\confy_{\pm}(\eltv)\} \times \left]1,\infty\right[$ or $\{\confy_{\pm}(\eltv)\} \times\left]-\infty,0\right[$. The projection $\confy_{\pm}(\eltv)$ depends on the considered $\infty$-component of $\eltv$, in general. When such a component is fixed, we simply denote the projection by $\confy(\eltv)$, and we speak of the $\infty$-component of a univalent vertex mapped to $\infty$.

\begin{definition}
\label{deffinercomptang}
Define the two-point compactification of $\RR=\left]-\infty, +\infty\right[$ to be $\left[-\infty, +\infty\right]$. 
For a Jacobi diagram $\Gamma$ on the domain $\sourcetl$
of $\tanghcyll$,
let $U_{ttbb}(\Gamma)$ be the set of vertices $\eltv$ of $U(\Gamma)$ such that the component $\tanghcyll(\sourcetl(\eltv))$ of $\eltv$ goes from top to top or from bottom to bottom.
Let $\overline{\sourcetl(\eltv)}$ denote the two-point compactification of the component $\sourcetl(\eltv)$ of $\sourcetl$.

The open manifold $\check{C}(\rats(\hcylc),\tanghcyll;\Gamma)$ embeds naturally in the product \begin{equation*} C_{V(\Gamma)}(\rats(\hcylc)) \times \prod_{\eltv \in U_{ttbb}(\Gamma)} \overline{\sourcetl(\eltv)}.\end{equation*}
Let $C^f_{\tanghcyll}=C^f\!(\rats(\hcylc),\tanghcyll;\Gamma)$ be the closure of $\check{C}(\rats(\hcylc),\tanghcyll;\Gamma)$ in this product. This closure maps naturally onto the closure $C_{\tanghcyll}={C}(\rats(\hcylc),\tanghcyll;\Gamma)$ of $\check{C}(\crats(\hcylc),\tanghcyll;\Gamma)$ in $C_{V(\Gamma)}(\rats(\hcylc))$.
An element of $C^f\!(\rats(\hcylc),\tanghcyll;\Gamma)$ is a configuration $\confc_{V(\Gamma)}$ of ${C}(\rats(\hcylc),\tanghcyll;\Gamma)$ equipped with the additional data of an $\infty$-component for each vertex of $U_{ttbb}(\Gamma) \cap \pbl(\confc_{V(\Gamma)})^{-1}(\infty)$.\footnote{This piece of data is automatically determined by $\confc_{V(\Gamma)}$ in most cases. However, it may happen that it is not. For example, when $\Gamma$ has a unique univalent vertex on a strand going from top to top, and when this vertex is mapped to $\infty$, the strand of the vertex may not be determined by the configuration in ${C}_{V(\Gamma)}(\rats(\hcylc))$.}
\end{definition}

In this section, we 
prove
the following theorem.
\begin{theorem}
\label{thmcomptang}
Let $\tanghcyll \colon \sourcetl \hookrightarrow \crats(\hcylc)$ be a long tangle representative, and let $\Gamma$
be a numbered degree $n$ Jacobi diagram with support $\sourcetl$ without looped edges. For any $\confc^0_{V(\Gamma)}$ in $C^f\!(\rats(\hcylc),\tanghcyll;\Gamma)$, there exist
\begin{itemize}
 \item a manifold $W$ with boundary and ridges, 
 \item a small $\varepsilon >0$,
 \item an oriented tree $\CT^0=\CT(\confc^0_{V(\Gamma)})$ (described in Notation~\ref{nottreeconftang}),
 with its associated singular space $\setparamx(\CT^0)$ defined in Definition~\ref{defXT}, and
 \item a smooth map \begin{equation*}\varphi \colon \left[0,\varepsilon\right[^{E(\CT^0)} \times W \to C_{V(\Gamma)}(\rats(\hcylc)),\end{equation*} whose restriction to $\bigl(\setparamx(\CT^0) \cap \left[0,\varepsilon\right[^{E(\CT^0)}\bigr) \times W$ is injective and induces a map
 \begin{equation*}\varphi^f \colon \bigl(\setparamx(\CT^0) \cap \left[0,\varepsilon\right[^{E(\CT^0)}\bigr) \times W \to C^f\!(\rats(\hcylc),\tanghcyll;\Gamma),\end{equation*}
\end{itemize}
such that $\Image(\varphi^f)$ is an open neighborhood of $\confc^0_{V(\Gamma)}$ in ${C}^f\!(\rats(\hcylc),\tanghcyll;\Gamma)$ where the $\infty$-component for each univalent vertex of the $\pbl(\confc^0_{V(\Gamma)})^{-1}(\infty)$ is the same as its $\infty$-component for $\confc^0_{V(\Gamma)}$. 

Let $\upvec$ denote the upward unit vertical vector.
The codimension-one open faces of ${C}^f\!(\rats(\hcylc),\tanghcyll;\Gamma)$ are 
\begin{itemize}
\item the faces corresponding to the collapse of a subgraph at one point in $\crats(\hcylc)$ as before,
\item the faces corresponding to a set of vertices mapped to $\infty$, for which the configuration up to dilation at $\infty$ is injective and does not map a point to $0$, as before,
\item additional faces called \emph{$T$-faces}\index[T]{Tfaces@$T$-faces} (for which ${C}(\rats(\hcylc),\tanghcyll;\Gamma)$ is not transverse to the ridges of $C_{V(\Gamma)}(\rats(\hcylc))$), where 
\begin{itemize}
\item a set of vertices $\finsetb \sqcup \sqcup_{j\in I}\finsetb_j$ is mapped to $\infty$, for a nonempty set I,
\item the corresponding configuration up to dilation from $\finsetb \sqcup \sqcup_{j\in I}\finsetb_j$ to $T_{\infty}\rats(\hcylc)$ maps each $\finsetb_j$ to a nonzero point of the vertical line, and it injects $\finsetb$ outside zero and the images of the $\finsetb_j$, which are distinct, 
\item each subset $\finsetb_j$ contains univalent vertices of at least $2$ distinct $\infty$-components,
\item for each $\finsetb_j$, the infinitesimal configuration of $\finsetb_j$ is an injective configuration of a Jacobi diagram on the lines that extend the half-lines above (resp. below) $\hcylc$, if $\finsetb_j$ is mapped to $\lambda \upvec$ for some $\lambda >0$ (resp. for some $\lambda <0$), up to global vertical translation. (No inversion is involved here.)\end{itemize}
\end{itemize}
\end{theorem}

Together with Lemma~\ref{lemconvonsetparamx} and Theorem~\ref{thmStokesonsetparamx}, Theorem~\ref{thmcomptang} implies the following lemma.

\begin{lemma}
\label{lemStokesconftang}
Theorem~\ref{thmconvint} is true.

Recall $C_{\tanghcyll}={C}(\rats(\hcylc),\tanghcyll;\Gamma)$
and $C^f_{\tanghcyll}={C}^f\!(\rats(\hcylc),\tanghcyll;\Gamma)$.
Let $\eta$ be a form of degree $(\dim(C_{\tanghcyll})-1)$ of ${C}_{V(\Gamma)}(\rats(\hcylc))$.
Then $\int_{C_{\tanghcyll}}d \eta$ is the sum $\sum_F\int_F \eta$, which runs over the codimension-one faces $F$ of $C^f_{\tanghcyll}$, oriented as such, and listed in Theorem~\ref{thmcomptang}.
\end{lemma}
\eopwobp

\begin{example}
\label{exathetatangtwo}
Let $K \colon \left]0,1\right[ \hookrightarrow \crats(\hcylc)$ be a (long) component of $\tanghcyll$. Assume that $K \colon \left]0,1\right[ \hookrightarrow \crats(\hcylc)$ goes from top to top. 
Let $d_1=-\{z_1\} \times \left[1,\infty\right]$ and $d_2=\{z_2\} \times \left[1,\infty\right]$ denote the vertical half-lines of $K$ above $\hcylc$, where $d_1$ is before $d_2$. 
Let $G=\{(h,k)\in \left]0,1\right]\times \RR \suchthat k+\frac{1}{h} \geq 1\}$.
Define the diffeomorphism \begin{equation*}\begin{array}{llll}g \colon &G &\to &d_1\times d_2\\
                           & (h,k) & \mapsto &\bigl((z_1,\frac{1}{h}), (z_2,k+\frac{1}{h})\bigr).
                           \end{array}\end{equation*}
This diffeomorphism $g$ extends as a continuous map \begin{equation*}g \colon G \cup \bigl(\{0\} \times \RR\bigr) \to {C}\bigl(\rats(\hcylc),\tanghcyll;\onechordsmalltK\bigr).\end{equation*} It maps $(0,k)$ to the limit $g(0,k)$ at $0$ in ${C}_2(\rats(\hcylc))$ of 
the $g(\left]0,\varepsilon\right]\times{k})$. The image of $g(0,k)$ under the canonical map from ${C}_2(\rats(\hcylc))$ to $\rats(\hcylc)^2$ is $(\infty,\infty)$.
The configuration in $T_{\infty}\rats(\hcylc)$ up to dilation  corresponding to $g(0,k)$ maps the two points to the same upward vertical vector.
We have $p_{\tau}(g(0,k))=\frac{(z_2-z_1,k)}{\norm{(z_2-z_1,k)}}$.
The image of
$\bigl(-g(\{0\}\times \RR)\bigr) $ under $p_{\tau}$ is the open half-circle from $\upvec$ to $-\upvec$ through
the direction of $(z_2-z_1)$.
This codimension-one face of ${C}^f\!(\rats(\hcylc),\tanghcyll;\onechordsmalltK)$ is an example of a $T$-face, for which $\finsetb \sqcup \sqcup_{j\in I}\finsetb_j=\finsetb_1$ is the pair of vertices of the graph.
Since this codimension-one face sits in a codimension-two face of $C_2(\rats(\hcylc))$, ${C}(\rats(\hcylc),\tanghcyll;\onechordsmalltK)$ is not transverse to the ridges of $C_2(\rats(\hcylc))$.
\end{example}

Recall that the elements of $C^f_{\tanghcyll}=C^f\!(\rats(\hcylc),\tanghcyll;\Gamma)$ are elements of the closure $C_{\tanghcyll}$ of $\check{C}(\crats(\hcylc),\tanghcyll;\Gamma)$ in $C_{V(\Gamma)}(\rats(\hcylc))$ equipped with the additional data of the $\infty$-components of the univalent vertices sent to $\infty$.

First note that the configuration space $C_{\tanghcyll}$ intersects $\pbl^{-1}(\crats(\hcylc)^{V(\Gamma)})$ as a smooth submanifold as in the case of links.
The only difference with the case of links occurs when some univalent vertices approach $\infty$.
Our configuration space is a local product of the space of the restrictions of the configurations to the points near $\infty$ and 
the space of the restrictions of the configurations to the other points, which is a smooth manifold with boundary whose structure has been studied in detail in Chapter~\ref{chapcompconf}.

Recall the orientation-reversing embedding $\phi_{\infty}$
\begin{equation*}\begin{array}{llll}\phi_{\infty}\colon &\RR^3 &\longrightarrow &S^3\\
& \mu (x \in S^2) & \mapsto & \left\{\begin{array}{ll} \infty \;&\;\mbox{if}\; \mu=0\\
\frac{1}{\mu}x \;&\;\mbox{otherwise.} \end{array}\right.\end{array}\end{equation*}
According to Corollary~\ref{cordescconftaucomp}, with the notation of Chapter~\ref{chapcompconf}, and especially those of Section~\ref{secbloinf}, an element $\confc_{V(\Gamma)}$ of $C_{V(\Gamma)}(\rats(\hcylc))$ consists of 
\begin{itemize}
\item a subset $\finsetv=\pbl(\confc_{V(\Gamma)})^{-1}(\infty)$ of $V(\Gamma)$,
\item an element $\confc_{V(\Gamma) \setminus \finsetv}=\confc_{V(\Gamma)}\vert_{V(\Gamma) \setminus \finsetv}$ of $C_{V(\Gamma) \setminus \finsetv}[\crats(\hcylc)]$,
\item an element $\confc=\confc_{V(\Gamma)}\vert_{\finsetv}$ of $C_{V}(\rats(\hcylc)) \cap \pbl^{-1}(\infty^{\finsetv})$,
\item a $\Delta$-parenthesization $\parentp^+=\parentp^+(\confc)$ of $V^+ = V \sqcup\{v_{\infty}\}$,
\item for each $\finseta \in \parentp^+$, an injective configuration \begin{equation*}T_{0} \phi_{\infty \ast} \circ w_{\finseta} \in \check{S}_{\kids(\finseta)}\Bigl(T_{\infty}\rats(\hcylc)\Bigr),\end{equation*} up to dilation and translation   (see Proposition~\ref{propstratinfty}). 
\end{itemize}

Proposition~\ref{proprest} describes the restriction maps, naturally.
As reminded above, the configuration space $C_{V(\Gamma)}(\rats(\hcylc))$ has a natural stratification induced by $\finsetv=\pbl(\confc)^{-1}(\infty)$, the parenthesization
associated to $\confc_{V(\Gamma) \setminus \finsetv}$ (as before Proposition~\ref{propdescconftautwo}) and the above $\Delta$-parenthesization $\parentp^+$ of $\finsetv^+$.
Each stratum has a well-defined dimension.

Below, we refine this partition induced on ${C}^f\!(\rats(\hcylc),\tanghcyll;\Gamma)$ by the stratification of $C_{V(\Gamma)}(\rats(\hcylc))$.

\begin{notation}
 \label{notationseparating}
As in Proposition~\ref{propstratinfty}, we define the totally ordered subset $\parentp^+_s$ to be the set of elements of $\parentp^+$ containing $v_{\infty}$ and
$\parentp_s= \{\finseta \setminus \{v_{\infty}\} \suchthat \finseta \in \parentp^+_s\}$.
We have \begin{equation*}\parentp_s=\{\finsetv=\finsetv(1), \finsetv(2), \dots, \finsetv(\sigma)\}\end{equation*} with $\finsetv(i+1)\subset\finsetv(i)$.
We let $\kids^s_d(\finsetv(i))$ denote the set of kids of $\finsetv(i)^+$ that do not contain $v_{\infty}$. We impose $w_{\finsetv(i)^+}(v_{\infty})=0$ and define $f_i$ to be the restriction of $w_{\finsetv(i)^+}$ to $\finsetv(i)$. So $f_i$ is an injective map from $\kids^s_d(\finsetv(i))$ to $\RR^3 \setminus \{0\}$ up to dilation. It belongs to $\sinjupdtcs\bigl(\RR^3,{\kids^s_d\bigl(\finsetv(i)\bigr)}\bigr)$.

The sets of $\parentp^+$ that contain a univalent vertex are called \emph{univalent}.
As in  Proposition~\ref{propstratinfty}, we set $\parentp_d=\parentp^+ \setminus \parentp^+_s$.
A \indexT{possibly separating set} associated to the above parenthesization $\parentp^+$ and to the data of the $\infty$-components of the elements of $\finsetv$ is a set $\finseta \in \parentp_d$ such that
\begin{itemize}
\item  $\finseta$ has at least two univalent vertices on different $\infty$-components and
\item each kid of $\finseta$ has all its univalent vertices on the same $\infty$-component.
\end{itemize}
Let $\parentp_{\cutx}$\index[N]{Parenthesizations!pacutX@$\parentp_{\cutx}$} denote the set of possibly separating sets associated to the above parenthesization $\parentp^+$.
A set $A$ of $\parentp_{\cutx}$ is \emph{separating} (with respect to $\confc_{V(\Gamma)}$) if it has at least two univalent kids $A_1$ and $A_2$
such that $w_{\finseta}(A_1)-w_{\finseta}(A_2)$ is not vertical.
The set of separating sets of $\confc_{V(\Gamma)}$ is denoted by $\parentp_x$\index[N]{Parenthesizations!pacutx@$\parentp_x$}.

Recall that $\upvec$ denotes the upward unit vertical vector.
Let $p_{\CC} \colon (\RR^3=\CC \times \RR) \to \CC$ denote the orthogonal projection onto the horizontal plane $\CC$, and let $p_{\RR} \colon \CC \times \RR \to \RR$  denote the orthogonal projection onto the vertical line $\RR$.
\end{notation}

We are going to prove the following two propositions.

\begin{proposition} \label{proppredefsepone}
The space $C^f_{\tanghcyll}=C^f\!(\rats(\hcylc),\tanghcyll;\Gamma)$ of Definition~\ref{deffinercomptang}
 is the space of configurations $\confc_{V(\Gamma)}$ of $C_{V(\Gamma)}(\rats(\hcylc))$ as above, equipped with $\infty$-components for the univalent vertices of $V=\pbl(\confc_{V(\Gamma)})^{-1}(\infty)$, such that the following conditions are satisfied.
\begin{enumerate}
 \item If the configuration $\confc_{V(\Gamma)}\vert_{U(\Gamma) \setminus (\finsetv \cap U(\Gamma))}$ is injective, then it factors through  the restriction to $U(\Gamma) \setminus (\finsetv \cap U(\Gamma))$ of a representative of $i_{\Gamma}$ that maps the univalent vertices of $\finsetv \cap U(\Gamma)$ to their $\infty$-components further than the elements of $U(\Gamma) \setminus (\finsetv \cap U(\Gamma))$. If the configuration $\confc_{V(\Gamma)}\vert_{U(\Gamma) \setminus (\finsetv \cap U(\Gamma))}$ is not injective, then it factors through a limit of such restrictions. In any case, the possible infinitesimal configurations of vertices of $U(\Gamma) \setminus (\finsetv \cap U(\Gamma))$ are locally ordered on the tangent space to their component, as in the case of links (see Sections~\ref{secfacecodimone} and \ref{secproofblowup}).
\item The $f_i$ map the elements
of $\kids^s_d\bigl(\finsetv(i)\bigr)$ that contain a univalent vertex on an $\infty$-component $y \times \left]1,\infty\right[$ (resp. $y \times \left]-\infty,0\right[$) to the half-line $\RR^+ \upvec$ (resp. $\RR^+ (-\upvec)$). 
\item If $\eltv_1 \in \finseta_1$ and $\eltv_2 \in \finseta_2$ are two univalent vertices of distinct kids $\finseta_1$ and $\finseta_2$ of an element $\finseta \in \parentp^+$, fix a  normalization of $w_{\finseta}$, and let $\vec{y}=\confy(\eltv_2)-\confy(\eltv_1)$.
\begin{itemize} 
 \item If $\vec{y}=0$ (that is if $\eltv_1$ and $\eltv_2$ are on the same $\infty$-component), and if  $\eltv_1$ is closer to $\infty$ than $\eltv_2$, then 
$w_{\finseta}(\finseta_2)-w_{\finseta}(\finseta_1)$ is a nonzero vertical vector, which may be expressed as $\alpha \upvec$, where  $\alpha$ is positive when the $\infty$-component
is above $\hcylc$, and $\alpha$ is negative 
when the $\infty$-component
is under $\hcylc$.
\item If $A \notin \parentp_{\cutx}$, then 
$w_{\finseta}(\finseta_2)-w_{\finseta}(\finseta_1)$ is also a nonzero vertical vector.
\item If $A \in \parentp_{\cutx}$ and if $\vec{y} \neq 0$, then $w_{\finseta}(\finseta_2)-w_{\finseta}(\finseta_1)$ may be expressed as $(\alpha \upvec + \beta \vec{y})$ for some nonzero pair $(\alpha,\beta)$ of $\RR \times \RR^+$. Furthermore,
if $\eltv_3 \in \finseta_3$ is a univalent vertex of another kid $\finseta_3$ of $\finseta$, then there exists $\alpha_3 \in \RR$ such that $w_{\finseta}(\finseta_3)-w_{\finseta}(\finseta_1)$ is equal to $(\alpha_3 \upvec + \beta \left(\confy(\eltv_3)-\confy(\eltv_1))\right)$.
\end{itemize} 
\end{enumerate}
\end{proposition}

\begin{proposition} \label{proppredefseptwo}
The space $C^f_{\tanghcyll}$ of Proposition~\ref{proppredefsepone}
 is partitioned 
by the data for a configuration $\confc_{V(\Gamma)}$ of
\begin{itemize}
\item the set $\finsetv=\pbl(\confc_{V(\Gamma)})^{-1}(\infty)$, 
\item the parenthesization $\parentp(\confc_{V(\Gamma)}\vert_{\vertsetv(\Gamma) \setminus\finsetv})$ of $\vertsetv(\Gamma) \setminus\finsetv$ associated to $\confc_{V(\Gamma)}\vert_{\vertsetv(\Gamma) \setminus\finsetv}$ (as before Proposition~\ref{propdescconftautwo}),
\end{itemize}
and, if $\vertsetv \neq \emptyset$,\footnote{When $\vertsetv =\emptyset$, the structure of $C^f_{\tanghcyll}$ near $\confc_{V(\Gamma)}$ is already known from Chapter~\ref{chapcompconf}.}
\begin{itemize}
\item the $\Delta$-parenthesization $\parentp^+=\parentp^+(\confc_{V(\Gamma)}\vert_{\finsetv})$ of $\vertsetv^+ =\vertsetv \sqcup \{v_{\infty}\}$ 
\item the data of the $\infty$-components of the univalent vertices that are mapped to $\infty$, and 
\item the set $\parentp_x$ of separating sets of $\parentp^+$.
\end{itemize}
The part associated to the above data is a smooth submanifold of $C_{V(\Gamma)}(\rats(\hcylc))$ of dimension
\begin{equation*}\cardbig{U(\Gamma)} +3 \cardbig{T(\Gamma)} - \cardbig{\parentp(\confc_{V(\Gamma)}\vert_{\vertsetv(\Gamma) \setminus\finsetv})}
-\cardbig{\parentp^+} + \cardbig{\parentp_x}.\end{equation*}
This partition is a stratification of $C^f_{\tanghcyll}$.
\end{proposition}

\begin{remark} \label{rkstrataTf}
 Proposition~\ref{proppredefseptwo} implies that the only codimension-one new parts --which  come necessarily from strata for which $\vertsetv \neq \emptyset$- come from the parts such that
 $\parentp^+=\{\vertsetv^+\} \sqcup \parentp_x$
and $\parentp_d=\parentp_x$.
They are the $T$-faces of Theorem~\ref{thmcomptang}, for which $\parentp_s=\{\vertsetv\}=\{\finsetb \sqcup \sqcup_{j\in I}\finsetb_j\}$ and $\parentp_d=\parentp_{\cutx}=\parentp_x=\{\finsetb_j \suchthat j\in I\}$.
\end{remark}

The rest of this section is devoted to the proofs of Theorem~\ref{thmcomptang} and the above two propositions.

Let $\confc^0_{V(\Gamma)}$ be a configuration of $C^f_{\tanghcyll}$ in a stratum as in the statement of Proposition~\ref{proppredefseptwo}. Let $\confc^0=\confc^0_{\vertsetv(\Gamma) }\vert_{\finsetv}$  
denote the restriction of $\confc^0_{V(\Gamma)}$ to $V(=V(\confc^0))$. 
 
It is easy to see that $\confc^0_{V(\Gamma)}$ has a neighborhood $N_{\Gamma}(\confc^0_{V(\Gamma)})$ in $C_{\tanghcyll}^f$ consisting of configurations mapping
\begin{itemize}
 \item $\finsetv$ to a fixed open neighborhood $N_{\infty}$ of $\infty$ in $\rats(\hcylc)$,
 \item the univalent vertices of $V$ to their $\infty$-component with respect to $\confc^0$,
 \item $\vertsetv(\Gamma) \setminus\finsetv$ to an open subspace $U$ of $\rats(\hcylc)\setminus N_{\infty}$.
\end{itemize}

Let $C^f_{\vertsetv}(N_{\infty},\tanghcyll,\confc^0)$ denote the space of restrictions to $V$ of configurations of ${C}(\rats(\hcylc),\tanghcyll;\Gamma)$ mapping univalent vertices of $V$ to $N_{\infty}$ and to their $\infty$-components (determined by $\confc^0$), and vertices of $\vertsetv(\Gamma) \setminus\finsetv$ to $U$.

The configuration $\confc^0_{V(\Gamma)}$ has a neighborhood $N_{\Gamma}(\confc^0_{V(\Gamma)})$ in $C_{\tanghcyll}^f$ that is a product of $C^f_{\vertsetv}(N_{\infty},\tanghcyll,\confc^0)$ by
a smooth submanifold $N_2$ of $C_{\vertsetv(\Gamma) \setminus\finsetv}[U]$.
In this product decomposition, the $N_2$-part contains the restriction of the configurations $\confc_{V(\Gamma)}$ to $\vertsetv(\Gamma) \setminus\finsetv$. 
This space has been studied before (see Proposition~\ref{propcompconfL}, Section~\ref{secblodiag}, and Theorem~\ref{thmcompuptd}), and the manifold $W$ of the statement of Theorem~\ref{thmcomptang} will be a product $W_{\finsetv} \times N_2$. 
This allows us to forget about the $N_2$-part.
We focus only on $\confc^0$ and on a neighborhood $N(\confc^0)$ of $\confc^0$ in $C_{V}(\rats(\hcylc))$.

The configuration $\confc^0$ is described by 
\begin{itemize}
\item a $\Delta$-parenthesization $\parentp^+=\parentp^+(\confc^0)=\parentp_s^+ \sqcup \parentp_d$ of $\finsetv^+$, where
$\parentp_s=\{\finsetv={\finsetv}(1),{\finsetv}(2), \dots, {\finsetv}(\sigma)\}
$ and ${\finsetv}(i+1)\subset {\finsetv}(i)$,
\item for any element $A$ of $\parentp^+$, an element $w^0_{\finseta} \colon K(A) \to \RR^3$ of the manifold $\tilde{W}_{\finseta}$ consisting of the injective maps $w_{\finseta} \colon K(A) \to \RR^3$ up to dilation and translation.
\end{itemize}

\begin{notation}
\label{notnormalizeta}
Let $i \in \overline{\sigma}$. Recall $\kids^s_d(V(i))=\kids\bigl(V(i)^+ \bigr) \cap \parentp_d$. Set $\kids^s(V(\sigma))=\kids^s_d(V(\sigma))$, and $\kids^s(V(i))=\kids^s_d(V(i)) \sqcup \{V(i+1)\}$ when $i\neq \sigma$.
We normalize $w^0_{V(i)^+}$ so that $w^0_{V(i)^+}(v_{\infty})=0$,  
and the restriction $f^0_i$ of $w^0_{V(i)^+}$ to $V(i)$ is in 
the manifold consisting of the maps $f_i \colon \kids^s(V(i)) \to \RR^3$ such that
\begin{itemize}
\item $f_i\bigl(V(i+1)\bigr)=0$, if $i\neq \sigma$,
\item $\sum_{A \in \kids^s_d(V(i))}\norm{f_i(A)}^2=1$,
\item $\norm{f_i(A)} > \eta$ for any $i$ and for any $A \in \kids^s_d\bigl(V(i)\bigr)$,
\item $\norm{f_i(A_2)-f_i(A_1)} > \eta$ for any two distinct elements $A_1$ and $A_2$ of $\kids^s\bigl(V(i)\bigr)$,
\end{itemize}
for some real number $\eta>0$,\footnote{Recall Notation~\ref{notationseparating}. The map $f_i$ represents an injective configuration $T_{0}\phi_{\infty} \circ f_i$ up to dilation of $\sinjupdtcs(T_{\infty}\rats(\hcylc),{\kids^s_d\bigl(\finsetv(i)\bigr)})$.}
We fix an open neighborhood  $W_i^s$ of  $f^0_i$ in the above manifold.
\end{notation}

\begin{notation}
\label{nothornnormaliz}
We choose univalent basepoints for univalent sets of $\parentp_d$. 
As usual, our basepoints also satisfy the conditions that for two elements $\finseta$ and $\finsetb$ of $\parentp_d$, such that $\finsetb \subset \finseta$, if $b(\finseta) \in \finsetb$, then $b(\finsetb)=b(\finseta)$.

For $\finseta \in \parentp_d$, we normalize the configurations of $\tilde{W}_{\finseta}$ as follows in a neighborhood $\tilde{N}_{\finseta}$ of a given $w^0_{\finseta} \in\tilde{W}_{\finseta}$. Choose a kid $k_n(\finseta)$ such that 
$|p_{\RR}(w^0_{\finseta}(k_n(\finseta))) -p_{\RR}(w^0_{\finseta}(b(\finseta)))|$ or $|p_{\CC}(w^0_{\finseta}(k_n(\finseta))) -p_{\CC}(w^0_{\finseta}(b(\finseta)))|$ is maximal in the set 
\begin{equation*}\biggl\{\Bigl|p_{\RR}\Bigl(w^0_{\finseta}\bigl(k\bigr)\Bigr) -p_{\RR}\Bigl(w^0_{\finseta}\bigl(b(\finseta)\bigr)\Bigr)\Bigr|, \Bigl|p_{\CC}\Bigl(w^0_{\finseta}\bigl(k\bigr)\Bigr) -p_{\CC}\Bigl(w^0_{\finseta}\bigl(b(\finseta)\bigr)\Bigr)\Bigr| \suchthat k \in \kids(\finseta) \biggr\},\end{equation*}
and call it the \emph{normalizing kid} of $\finseta$.
If 
\begin{equation*}\Bigl|p_{\RR}\Bigl(w^0_{\finseta}\bigl(k_n(\finseta)\bigr)\Bigr) -p_{\RR}\Bigl(w^0_{\finseta}\bigl(b(\finseta)\bigr)\Bigr)\Bigr| \geq \Bigl|p_{\CC}\Bigl(w^0_{\finseta}\bigl(k_n(\finseta)\bigr)\Bigr) -p_{\CC}\Bigl(w^0_{\finseta}\bigl(b(\finseta)\bigr)\Bigr)\Bigr|,\end{equation*}
then we say that $k_n(\finseta)$ is \indexT{vertically normalizing} or \emph{v-normalizing} and 
normalize the configurations $\confw_{\finseta}$ in a neighborhood of $w^0_{\finseta}$ by imposing \begin{equation*}\confw_{\finseta}(b(\finseta))=0\mbox{ and }|p_{\RR}(\confw_{\finseta}(k_n(\finseta)))|=1.\end{equation*}
Otherwise, we say that $k_n(\finseta)$ is \indexT{horizontally normalizing} or \emph{h-normalizing}, and normalize the configurations $\confw_{\finseta}$ in a neighborhood of $w^0_{\finseta}$ by imposing $\confw_{\finseta}(b(\finseta))=0$ and $|p_{\CC}(\confw_{\finseta}(k_n(\finseta)))|=1$. (These normalizations are compatible with the smooth structure of $C_{\vertsetv(\Gamma)}(\rats(\hcylc))$.)

In our neighborhood $\tilde{N}_{\finseta}$, we also impose that $\norm{\confw_{\finseta}(k)-\confw^0_{\finseta}(k)}<\varepsilon$, for a small $\varepsilon \in \left]0,1\right[$. So $\tilde{N}_{\finseta}$ is diffeomorphic to the product $W_{\finseta}$ of 
\begin{itemize}
 \item the product, over the nonnormalizing kids $k$ of $\finseta$ that do not contain $b(\finseta)$, of the open balls $\mathring{\ballb}(\confw^0_{\finseta}(k),\varepsilon)$ of radius $\varepsilon$ centered at $\confw^0_{\finseta}(k)$, by
 \item the set of $w_{\finseta}(k_n(\finseta))$
in $\mathring{\ballb}(w^0_{\finseta}(k_n(\finseta)),\varepsilon)$ such that $\left|p_{\RR}\bigl(w_{\finseta}(k_n(\finseta))\bigr)\right|=1$ (resp. such that $\left|p_{\CC}\bigl(w_{\finseta}(k_n(\finseta))\bigr)\right|=1$) if $k_n(\finseta)$ is v-normalizing (resp. if $k_n(\finseta)$ is h-normalizing).
\end{itemize}
We reduce the $\eta$ of Notation~\ref{notnormalizeta} if necessary, and we choose $\varepsilon$ so that $\norm{w_{\finseta}(\finsetb_2)-w_{\finseta}(\finsetb_1)} > \eta$ for any two distinct kids $\finsetb_1$ and $\finsetb_2$ of $\finseta$ in our normalized neighborhood $W_{\finseta}$.
Note that $\norm{\confw_{\finseta}(k)}<3$ for any $k \in  \kids(\finseta)$ in this neighborhood.
\end{notation}
All the considered maps $f_i=w_{V(i)^+}\vert_{V(i)}$, $w_{\finseta}$ are also considered as maps from $\finsetv^+$ to $\RR^3$, which are constant on the elements of $\kids^s\bigl(\finsetv(i)\bigr)$ and $\kids(\finseta)$, respectively, and which respectively map $\finsetv^+ \setminus \finsetv(i)$ 
and $\finsetv^+ \setminus \finseta$ to $0$. 

We use a chart $\psi$ of ${C}_{V}(\rats(\hcylc))$ of a neighborhood $N(\confc^0)$ of $\confc^0$ in  ${C}_{V}(\rats(\hcylc))$ mapping
\begin{equation*}\bigl((\mu_{\finseta})_{{\finseta} \in \parentp^+},(w_{\finseta})_{{\finseta} \in \parentp^+}\bigr) \in \left[0,\varepsilon\right[^{\parentp^+} \times \prod_{i \in \underline{\sigma}} W^s_i  \times \prod_{{\finseta} \in \parentp_d} W_{\finseta}\end{equation*}
to a configuration $\confc=\psi\left((\mu_{\finseta})_{{\finseta} \in \parentp^+},(w_{\finseta})_{{\finseta} \in \parentp^+}\right) \in {C}_{\finsetv}(\rats(\hcylc))$,
such that, when the $\mu_{\finseta}$ do not vanish, $\confc$ is the injective configuration
\begin{equation*}\confc = \phi_{\infty} \circ\left(\sum_{\finseta \in \parentp^+}\left(\prod_{\finsetd \in\parentp^+ \suchthat \finseta \subseteq \finsetd} \mu_{\finsetd} \right)w_{\finseta}\right).\end{equation*}
With this chart induced by Theorems~\ref{thmcompuptd}, \ref{thmcompdiaggen}, and \ref{thmcompconfbis}, we have \begin{equation*}\confc^0=\psi\bigl((\mu_{\finseta}^0=0),(w_{\finseta}^0)\bigr).\end{equation*}

\begin{notation}\label{notbigu}
 For $k \in \underline{\sigma}$, set $u_k=\mu_{V(k)^+}$, and set $U_k=\prod_{i=1}^k u_i$.
 For  $A \in \parentp_d$, let $k(A)$ be the maximal integer among the integers $k$ such that $A \subseteq V(k)$.
So $\psi$ maps
\begin{equation*}\bigl((u_i)_{i \in \underline{\sigma}}, (\mu_{\finseta})_{{\finseta} \in \parentp_d},(f_i)_{i \in \underline{\sigma}},(w_{\finseta})_{{\finseta} \in \parentp_d}\bigr) \in \left[0,\varepsilon\right[^{\sigma} \times \left[0,\varepsilon\right[^{\parentp_d}\times \prod_{i \in \underline{\sigma}} W^s_i  \times \prod_{{\finseta} \in \parentp_d} W_{\finseta}\end{equation*}
to a configuration $\confc=\psi\left((u_i), (\mu_{\finseta}),(f_i),(w_{\finseta})\right) \in {C}_{\finsetv}(\rats(\hcylc))$,
such that, when the $u_i$ and the $\mu_{\finseta}$ do not vanish, $\confc$ is the injective configuration
\begin{equation*}\confc = \phi_{\infty} \circ\left(\sum_{{\finsetv}(k) \in \parentp_s}U_k \left(f_{k} + \sum_{\finsetc \in \parentp_d \suchthat k(\finsetc)=k}\left(\prod_{\finsetd \in \parentp_d \suchthat  \finsetc \subseteq \finsetd} \mu_{\finsetd} \right) w_{\finsetc}\right)\right).\end{equation*} 
We also write \begin{equation*}\confc^0=\psi\bigl((u_i^0=0), (\mu_{\finseta}^0=0),(f_i^0),(w_{\finseta}^0)\bigr).\end{equation*}
\end{notation}

\begin{example}
\label{exaproofpropenferone}
In the special case of Example~\ref{exathetatangtwo}, with the graph
\begin{equation*}\Gamma=\begin{tikzpicture} \useasboundingbox (-.3,-.3) rectangle (.7,.3);
\draw [->,dash pattern=on 2pt off 2pt]  (0,-.4) -- (0,.6);
\draw (0,.5) node[left]{\scriptsize K};
\draw (0,-.2) node[right]{\scriptsize $v_1$};
\draw (0,.2) node[right]{\scriptsize $v_2$};
\draw (0,-.2) .. controls (-.3,-.2) and (-.3,.2) .. (0,.2);
\fill (0,-.2) circle (1.5pt) (0,.2) circle (1.5pt);
\end{tikzpicture},\end{equation*} consider configurations mapping $v_1$ to $-\{z_1\} \times \left[1,\infty\right]$ and $v_2$  to $\{z_2\} \times \left[1,\infty\right]$. When $\finsetv=\{v_1,v_2\}$ and $\parentp^+=\{\finsetv^+,\finsetv\}$, we have $\parentp_s=\parentp_d=\{\finsetv\}$. Set $f=f_1$, $u=u_1$, $\mu=\mu_{\finsetv}$, $b(\finsetv)=v_1$, $w=w_{\finsetv}$. We have $w(v_1)=0$, $\norm{f(v_1)}=1$, $f^0(v_1)=\upvec$, $\confc = \phi_{\infty} \circ\left(u \left(f +  \mu w \right)\right)$,
$\confc(v_1)=\phi_{\infty} (u f(v_1))=\frac1{u}f(v_1)$, and $\confc(v_2)=\phi_{\infty} (u (f(v_1) + \mu w(v_2)))$.
\end{example}

Back to the general proof of Theorem~\ref{thmcomptang} and Propositions \ref{proppredefsepone} and \ref{proppredefseptwo}, we will often reduce $\varepsilon$ and reduce the spaces $W_k^s$ and $W_{\finseta}$ to smaller manifolds, which are neighborhoods of $f_k^0$ and $w_{\finseta}^0$ in the initial manifolds $W_k^s$ and $W_{\finseta}$. In particular, we assume that the image $N(\confc^0)$ of $\psi$ is in ${C}_{V}(N_{\infty})$, and we set 
\begin{equation*}N^f_{\tanghcyll}(\confc^0) = N(\confc^0) \cap C^f_{\vertsetv}(N_{\infty},\tanghcyll,\confc^0).\end{equation*}

The intersection of $N^f_{\tanghcyll}(\confc^0)$ with $\check{C}_{V}(\rats(\hcylc))$ 
is determined by the conditions that univalent vertices belong to their $\infty$-components and that their order on the $\infty$-components is prescribed by the isotopy class of injections from $U(\Gamma)$ to $\sourcetl$.
These conditions are closed. So they still hold in $N^f_{\tanghcyll}(\confc^0)$.

In particular, the basepoints $b(\finseta)$ of the univalent elements $\finseta$ of $\kids^s_d({\finsetv}(k))$ must go to their $\infty$-components. We call this condition the \emph{first condition}.
We examine what this \say{first condition} imposes on the $f_k$ and prove the following two lemmas.

\begin{lemma}
\label{lemfkone}
For any $k\in \underline{\sigma}$, for any univalent element $\finseta$ of $\kids^s_d({\finsetv}(k))$, we have
 $f^0_k(A)=\pm \norm{f^0_k(A)}\upvec$, where $\norm{f^0_k(A)} \geq \eta$, and where
 the $\pm$ sign is $+$ if $A$ has a univalent vertex above $\hcylc$ and $-$ otherwise.
\end{lemma}
\bp
For an element $\finseta$ of $\kids^s_d({\finsetv}(k))$ and a configuration \begin{equation*}\confc=\psi\bigl((u_i), (\mu_{\finseta}),(f_i),(w_{\finseta})\bigr) \in {C}_{\finsetv}(\rats(\hcylc))\end{equation*} such that $U_k=\prod_{i=1}^k u_i \neq 0$, we have
\begin{equation*}\confc\bigl(b(A)\bigr)=\frac{1}{U_k}\frac{f_k(\finseta)}{\bignorm{f_k(\finseta)}^2}.\end{equation*}
So the condition $p_{\CC}(\confc\bigl(b(A)\bigr))=y\bigl(b(A)\bigr)$ is equivalent to the closed condition
\begin{equation}
\label{eqone}
 p_{\CC}\bigl(f_k(\finseta)\bigr)=U_k \bignorm{f_k(\finseta)}^2y\bigl(b(A)\bigr),
\end{equation}
where \begin{equation*}\Bignorm{f_k(\finseta)}^2=\Bignorm{p_{\CC}\bigl(f_k(\finseta)\bigr)}^2 + \Bignorm{p_{\RR}\bigl(f_k(\finseta)\bigr)}^2 \leq 1,\end{equation*}
with our normalization of Notation~\ref{notnormalizeta}.

So Equation~\ref{eqone} implies $\norm{p_{\CC}(f_k(A))}=O(U_k)$ (meaning that there exists $C \in \RR^{\ast+}$ such that $\norm{p_{\CC}(f_k(A))}\leq C U_k$).
In particular, since $U_k=0$ for $\confc^0$, we have $p_{\CC}(f^0_k(A))=0$ and $f^0_k(A)=\pm \norm{f^0_k(A)}\upvec$, where $\norm{f^0_k(A)} \geq \eta$.
Lemma~\ref{lemfkone} follows easily since the sign of $p_{\RR}(f_k(A))$ is constant on $N^f_{\tanghcyll}(\confc^0)$.
\eop

In particular, if $V(\sigma) \in \parentp_d$, then $f^0_{\sigma}\bigl(V(\sigma)\bigr)$ equals $\upvec$ if $V(\sigma)$ has a univalent vertex above $\hcylc$, and $f^0_{\sigma}\bigl(V(\sigma)\bigr)$ equals $(- \upvec)$ if $V(\sigma)$ has a univalent vertex under $\hcylc$.

\begin{lemma}
\label{lemfkonebis}
Let $v_1$ and $v_2$ be two vertices on some $\infty$-component of $\tanghcyll$. Assume that $v_1$ is closer to $\infty$ than $v_2$ and that $v_1 \in A_1$ and $v_2\in A_2$, for two different kids $A_1$ and $A_2$ of $V(k)$, where $k \in \underline{\sigma}$. Then we have
\begin{equation}\label{eqtwo}
 \norm{f^0_k(A_1)} < \norm{f^0_k(A_2)}.
\end{equation}
\end{lemma}
\bp 
The configuration $\confc^0$ is a limit at $0$ of a family $\confc(t)$ indexed by $t \in \left]0,\varepsilon\right[$ of injective configurations for which 
$\norm{\confc(t)(v_1)}>\norm{\confc(t)(v_2)}$, with
\begin{equation*}\confc(t)=\confc=\psi\bigl((u_i), (\mu_{\finseta}),(f_i),(w_{\finseta})\bigr).\end{equation*}
Therefore,  we have
$\norm{f^0_k(A_1)} \leq \norm{f^0_k(A_2)}$. Since $f^0_k(A_1) \neq f^0_k(A_2)$, the result follows.
\eop

In particular, we have $|p_{\RR}(f^0_k(A_1))| \leq |p_{\RR}(f^0_k(A_2))| -\eta$.
We possibly reduce $W^s_k$ by imposing $|p_{\RR}(f_k(\finseta))-p_{\RR}(f^0_k(\finseta))| < \varepsilon$ for some positive $\varepsilon$ such that $\varepsilon <\frac{\eta}{2}$. 
This condition ensures that the univalent vertices $b(\finseta)$, for the elements $\finseta$ of $\kids^s_d({\finsetv}(k))$, are well-ordered on any $\infty$-component.

\begin{lemma}
\label{lemfktwo}
Let $k \in \underline{\sigma}$ be such that $V(k) \notin \parentp_d$. Recall $U_k=\prod_{i=1}^ku_i$.
Let $A_{k,1}$, $A_{k,2}$, \dots, $A_{k,\ell(k)}$ denote the univalent elements of $\kids^s_{d}(V(k))$, where $\ell(k) \in \NN$.
Let $\kids^s_{td}(V(k))$ be the set of nonunivalent elements of $\kids^s_{d}(V(k))$. For $\finsetd \in \kids^s_{td}(V(k))$, let $\mathring{\ballb}_{\finsetd}=\mathring{\ballb}(f^0_k(\finsetd),\varepsilon)$ be the open ball of center $f^0_k(\finsetd)$ and of radius $\varepsilon$ in $\RR^3$.
\begin{itemize}
 \item 

If $\ell(k) \geq 1$, set 
\begin{equation*}W_k^L=\prod_{i=1}^{\ell(k)-1}\Bigl]p_{\RR}\bigl(f^0_k(\finseta_{k,i})\bigr)-\varepsilon,p_{\RR}\bigl(f^0_k(\finseta_{k,i})\bigr)+\varepsilon\Bigr[ \times \prod_{\finsetd \in \kids^s_{td}(V(k))}\mathring{\ballb}_{\finsetd}.\end{equation*}
Up to reducing $\varepsilon$, there is a smooth injective map 
\begin{equation*}\phi_k \colon [0,\varepsilon^k[ \times W_k^L \to W_k^s\end{equation*}
such that, up to reducing $N^f_{\tanghcyll}(\confc^0)$ (to a smaller open neighborhood of $\confc^0$ in $C^f_{\vertsetv}(N_{\infty},\tanghcyll,\confc^0)$), all elements \begin{equation*}\confc=\psi\left((u_i)_{i \in \underline{\sigma}}, (\mu_{\finseta})_{{\finseta} \in \parentp_d},(f_i)_{i \in \underline{\sigma}},(w_{\finseta})_{{\finseta} \in \parentp_d}\right)\end{equation*} of $N^f_{\tanghcyll}(\confc^0)$ satisfy the condition
\begin{equation*}f_k=\phi_k\biggl(\prod_{i=1}^ku_i,\Bigl(p_{\RR}\bigl(f_k(A_{k,i})\bigr)\Bigr)_{i\in \underline{\ell(k)-1}},\Bigl(f_k(\finsetd)\Bigr)_{\finsetd \in \kids^s_{td}(V(k))} \biggr).\end{equation*}
This condition is equivalent to 
\begin{equation*}p_{\CC}\Bigl(\confc\bigl(b(A_{k,i})\bigr)\Bigr) = y\bigl(b(A_{k,i})\bigr)\end{equation*} 
for any $i \in \underline{\ell(k)}$ when $\confc \in \check{C}_{V}(\rats(\hcylc))$, and it implies Equation~\ref{eqone} for any $\confc$ in $N(\confc^0)$.

\item If $\ell(k)=0$, set $W_k^L = W_k^s$.\end{itemize}

If $V(\sigma) \in \parentp_d$ and if $V(\sigma)$ is not univalent, also set $W_{\sigma}^L = W_{\sigma}^s$.

If $V(\sigma) \in \parentp_d$ and if $V(\sigma)$ is univalent, then all elements $\confc$ of $N^f_{\tanghcyll}(\confc^0)$ satisfy the condition \begin{equation*}p_{\CC}\Bigl(f_{\sigma}\bigl(V(\sigma)\bigr)\Bigr)=U_{\sigma}y\Bigl(b\bigl(V(\sigma)\bigr)\Bigr).\end{equation*}
This condition is equivalent to 
$p_{\CC}\bigl(c\bigl(b(V(\sigma))\bigr)\bigr) = y\bigl(b(V(\sigma))\bigr)$ when $\confc \in \check{C}_{V}(\rats(\hcylc))$.
These elements $\confc$ also satisfy
\begin{equation*}p_{\RR}\Bigl(f_{\sigma}\bigl(V(\sigma)\bigr)\Bigr)=\sqrt{1-U_{\sigma}^2\bigl|y\bigl(b\bigl(V(\sigma)\bigr)\bigr)\bigr|^2}p_{\RR}\Bigl(f^0_{\sigma}\bigl(V(\sigma)\bigr)\Bigr).\end{equation*}
In this case, set $W_{\sigma}^L=\{\ast_{\sigma}\}$.

Let $N_1$ denote the subspace of $N(\confc^0)$, where the first condition (stated before Lemma~\ref{lemfkone}) is satisfied. Then $N_1$ is a smooth manifold parametrized by \begin{equation*}\left[0,\varepsilon\right[^{\sigma} \times \left[0,\varepsilon\right[^{\parentp_d}\times \prod_{k \in \underline{\sigma}} W^L_k  \times \prod_{{\finseta} \in \parentp_d} W_{\finseta}.\end{equation*}

\end{lemma}
\bp The proof of Lemma~\ref{lemfkone} shows that
for $i  \in \underline{\ell(k) -1}$, $p_{\CC}(f_k(A_{k,i}))$ is an implicit function of $U_k=\prod_{i=1}^k u_i$ and $p_{\RR}(f_k(A_{k,i}))$, which is close to $\pm \norm{f^0_k(A_{k,i})} \neq 0$ on $N_1$ and $N^f_{\tanghcyll}(\confc^0)$. This implicit function is determined by Equation~\ref{eqone}.

Then the condition that $\sum_{A \in \kids^s_d(V(k))}\norm{f_k(A)}^2=1$ in $W_k^s$ 
determines $\norm{f_k(A_{k,\ell(k)})} \neq 0$ as a function of $U_k$, $\left(p_{\RR}(f_k(A_{k,i}))\right)_{i  \in \underline{\ell(k) -1}}$ and of the $f_k(\finsetd)$ for $\finsetd \in \kids^s_{td}(V(k))$. Now, Equation~\ref{eqone} determines $p_{\CC}(f_k(A_{k,\ell(k)}))$, which in turn determines $p_{\RR}(f_k(A_{k,\ell(k)}))$. This is how the map $\phi_k$ of the statement is constructed. It is easy to check that $\phi_k$ has the desired properties and that $N_1$ is parametrized naturally, as announced, using the maps $\phi_k$. \eop

In Example~\ref{exaproofpropenferone}, we have $p_{\CC}\bigl(f(v_1)\bigr)=u y(v_1)$ and \begin{equation*}p_{\RR}\bigl(f(v_1)\bigr)=\sqrt{1-u^2\norm{y(v_1)}^2}.\end{equation*}
So $f(v_1)$ is just a smooth function of the small parameter $u$.

We now restrict to the submanifold $N_1$ of $N(\confc^0)$ of Lemma~\ref{lemfktwo} and take care of the univalent basepoints of the kids of elements of $\parentp_d$ in the following lemmas.

\begin{lemma}
\label{lemptauconfczero}
For $\vecx \in S^2$, let $\syms(\vecx)$ denote the orthogonal reflection of $\RR^3$ with respect to the plane orthogonal to $\vecx$. 
Let $A$ be an element of $\parentp_d$. Recall that $k(A)$ is the maximal integer $k$ such that $A \subseteq V(k)$.
The restriction of $\confc^0$ to $\finseta$ maps $\finseta$ to \begin{equation*}\vecx^0_{\finseta}=\frac{f^0_{k(A)}(A)}{\norm{f^0_{k(A)}(A)}}
\in \partial \blowup{\rats(\hcylc)}{\infty}.\end{equation*} It is represented by $\syms(\vecx^0_{\finseta}) \circ w^0_{\finseta}$ up to translation and dilation, as a configuration of the ambient $\RR^3$ outside $\hcylc$.
If $\eltp$ and $\eltq$ are univalent vertices in two different kids of $A$, if they belong to an $\infty$-component $K^+$, and if $\eltp$ is closer to $\infty$ than $\eltq$,
then there exists $\alpha^0 \in \RR$ such that $|\alpha^0|>\eta$ and $w^0_{\finseta}(\eltq)-w^0_{\finseta}(\eltp)=\alpha^0 \upvec$, where $\alpha^0 >0$ if $K^+$ is above $\hcylc$, and $\alpha^0<0$ otherwise.
\end{lemma}
We introduce some notation before the proof.

\begin{notation}
\label{notpreldeflambda}
For an element $\finseta$ of $\parentp_d$ such that 
$k(A)=k$
and a configuration $\confc=\psi\left((u_i), (\mu_{\finseta}),(f_i),(w_{\finsetb})\right) \in N(\confc^0)$,
set \begin{equation*}M_{\finseta}=\prod_{\finsetd \in \parentp_d \suchthat  \finseta \subseteq \finsetd} \mu_{\finsetd},\end{equation*} and, for any element $\eltq$ of $\finseta$, define
 \begin{equation*}\tilde{f}_k(\eltq)= f_k(\finseta) + \sum_{\finsetc \in \parentp_d \suchthat \eltq \in \finsetc}M_{\finsetc} w_{\finsetc}(\eltq),\end{equation*}
so that we have $\confc(\eltq)=\frac{1}{U_k}\frac{\tilde{f}_k(\eltq)}{\norm{\tilde{f}_k(\eltq)}^2}$ when $U_k \neq 0$. Both $M_{\finseta}$ and $\tilde{f}_k$ depend on the configuration $\confc$.
\end{notation}

\bpo{Proof of Lemma~\ref{lemptauconfczero}}
The configuration $\confc^0_{V(\Gamma)}$ is a limit at $0$ of a family $\confc(t)$ indexed by $t \in \left]0,\varepsilon\right[$ of injective configurations of $N^f_{\tanghcyll}(\confc^0)$.
We have
\begin{equation*}\confc=\confc(t)=\psi\bigl((u_i), (\mu_{\finsetb}),(f_i),(w_{\finsetb})\bigr),\end{equation*} where the $u_i$ and the $\mu_{\finsetb}$ are positive. Set $k=k(\finseta)$.
Let $\eltp=b(\finseta)$ be the basepoint of $\finseta$, and let $\eltq$ be the basepoint of a kid of $\finseta$ that does not contain $\eltp$.
Since $\tilde{f}_k(\eltq)-\tilde{f}_k(\eltp) = M_{\finseta}w_{\finseta}(\eltq)$, we have
\begin{equation*}\norm{\tilde{f}_k(\eltq)}^2=\norm{\tilde{f}_k(\eltp)}^2 +2 M_{\finseta} \bigl\langle w_{\finseta}(\eltq) ,\tilde{f}_k(\eltp)\bigr\rangle +M_{\finseta}^2\norm{w_{\finseta}(\eltq)}^2.\end{equation*}
We get \begin{multline*}\confc(\eltq)-\confc(\eltp)=\frac{\norm{\tilde{f}_k(\eltp)}^2\tilde{f}_k(\eltq)-\norm{\tilde{f}_k(\eltq)}^2\tilde{f}_k(\eltp)}{U_k\norm{\tilde{f}_k(\eltq)}^2\norm{\tilde{f}_k(\eltp)}^2}\\
        =\frac{M_{\finseta}\norm{\tilde{f}_k(\eltp)}^2w_{\finseta}(\eltq)}{U_k\norm{\tilde{f}_k(\eltq)}^2\norm{\tilde{f}_k(\eltp)}^2}
  -     
\frac{2 M_{\finseta} \langle w_{\finseta}(\eltq),\tilde{f}_k(\eltp)\rangle +M_{\finseta}^2\norm{w_{\finseta}(\eltq)}^2}{U_k\norm{\tilde{f}_k(\eltq)}^2\norm{\tilde{f}_k(\eltp)}^2} \tilde{f}_k(\eltp).
\end{multline*}
When the $\mu_{\finsetb}$ tend to $0$ and when $\tilde{f}_k$ tends to ${f}_k^0$, $\tilde{f}_k(\eltp)$ and $\tilde{f}_k(\eltq)$ tend to ${f}^0_k(\finseta)$. Thus, \begin{equation*}\frac{U_k\norm{\tilde{f}_k(\eltq)}^2}{M_{\finseta}}\bigl(\confc(\eltq)-\confc(\eltp)\bigr)\end{equation*} tends to \begin{equation*}w_{\finseta}(\eltq) -2\bigl\langle w_{\finseta}(\eltq),\vecx^0_{\finseta}\bigr\rangle \vecx^0_{\finseta}=\syms\bigl(\vecx^0_{\finseta}\bigr)\bigl(w_{\finseta}(\eltq)\bigr).\end{equation*}
So $w_{\finseta}^0$ is the limit of the $\syms(\vecx^0_{\finseta}) \circ \confc\vert_{A}$ up to dilation and translation.

If $\finseta$ contains a univalent vertex of an $\infty$-component above $\hcylc$, then $\vecx^0_{\finseta}$ equals $\upvec$ according to Lemma~\ref{lemfkone}. In this case, if $\elta$ and $\eltq$ are univalent vertices in two different kids of $A$, if they belong to an $\infty$-component $K^+$, and if $\elta$ is closer to $\infty$ than $\eltq$, then we have $\confc(\eltq)-\confc(\elta)=-\alpha(t) \upvec$ for some positive $\alpha(t)$ for any $t>0$. So $(w_{\finseta}^0(\eltq)-w_{\finseta}^0(\elta))$, which is defined up to dilation, is a positive multiple of $\upvec$.
\eop

\begin{notation}
\label{notparentlambda}
Let $\parntpwidehat(=\parntpwidehat(\confc^0))$ denote the set of elements of $\parentp_d (=\parentp_d(\confc^0))$ that contain or are equal to an element of $\parentp_{\cutx}(=\parentp_{\cutx}(\confc^0))$. 
\end{notation}

\begin{lemma}
\label{lemdeflambda} Let $A \in \parentp_d$ be such that $k(A)=k$. 
Let $\confc$ be as in Notation~\ref{notpreldeflambda}. 
For any univalent kid $\finsetb$ of $A$ such that $y(b(\finsetb))-y(b(\finseta))=0$, if $\confc \in  N^f_{\tanghcyll}(\confc^0)$, we have
\begin{equation} \label{eqfour}
p_{\CC}\bigl(w_{\finseta}(\finsetb)\bigr)= U_k\Bigl(2\bigl\langle \tilde{f}_{k}\bigl(b(\finseta)\bigr),w_{\finseta}(\finsetb) \bigr\rangle
+M_{\finseta} \bignorm{w_{\finseta}(\finsetb)}^2 \Bigr)y\bigl(b(\finsetb)\bigr).
\end{equation}

Furthermore, as soon as $p_{\CC}(\confc(b(A)))=y(b(A))$ and $\confc \in \check{C}_{V}(\rats(\hcylc))$, Equation~\ref{eqfour} implies $p_{\CC}(\confc(b(\finsetb)))=y(b(\finsetb))$ for such a $\finsetb$.

If $\confc \in N^f_{\tanghcyll}(\confc^0)$, then there exists $\bigl(\lambda_{\finseta}=\lambda_{\finseta}(\confc)\bigr)_{A \in \parntplambda} \in (\RR^+)^{\parntplambda}$ such that the following properties are satisfied.
\begin{itemize}
\item For any univalent kid $\finsetb$ of $A \in \parntplambda$, we have
\begin{multline} \label{eqfive}
p_{\CC}\bigl(w_{\finseta}(\finsetb)\bigr)= \lambda_{A} \bignorm{\tilde{f}_{k(A)}\bigl(b(\finseta)\bigr)}^2 \Bigl(y\bigl(b(\finsetb)\bigr)-y\bigl(b(A)\bigr)\Bigr)\\
 + U_{k(A)} \Bigl(2\bigl\langle \tilde{f}_{k(A)}\bigl(b(\finseta)\bigr),w_{\finseta}(\finsetb) \bigr\rangle
+M_{\finseta} \bignorm{w_{\finseta}(\finsetb)}^2 \Bigr)y\bigl(b(\finsetb)\bigr).
\end{multline}
\item For $M \in \parentp_{\cutx}$, set $\tilde{\lambda}_M=\lambda_M\mu_M$.
For any $A \in \parntplambda \setminus \parentp_{\cutx}$, there exists $M \in \parentp_{\cutx}$ such that $M \subset A$ and we have
\begin{equation*}\lambda_{\finseta}=\tilde{\lambda}_M\prod_{\finsetd \in \parentp_d \suchthat  M \subset \finsetd \subset A} \mu_{\finsetd}.\end{equation*}
\item The map $\lambda_{A}$ is continuous on $N_{\Gamma}(c^0)$.
\item We have $\lambda_{A}M_{\finseta}=U_{k(A)}$.
\item For any three elements $A$, $\finsetb$, and $D$ in $\parentp_d$ such that $A$ and $\finsetb$ are in $\parntplambda$ and $A \cup \finsetb \subseteq D$, we have
\begin{equation*}\lambda_{\finseta} \prod_{C\in \parentp_d \suchthat A \subseteq C \subset D}\mu_C = \lambda_{\finsetb} \prod_{C\in \parentp_d \suchthat \finsetb \subseteq C \subset D}\mu_C .\end{equation*}
\item For a univalent kid $\finsetb$ of $A \in \parntplambda$, as soon as $p_{\CC}(\confc(b(A)))=y(b(A))$ holds and $\confc \in \check{C}_{V}(\rats(\hcylc))$, Equation~\ref{eqfive} implies $p_{\CC}(\confc(b(\finsetb)))=y(b(\finsetb))$.\end{itemize}
\end{lemma}
\bp Let $A \in \parentp_d$ be such that $k(A)=k$.
Consider a univalent kid $\finsetb$ of $\finseta$
and assume $p_{\CC}(\confc(b(A)))=y(b(A))$. This equality is equivalent to
\begin{equation*}p_{\CC}\Bigl(\tilde{f}_k\bigl(b(\finseta)\bigr)\Bigr)=U_k\bignorm{\tilde{f}_k\bigl(b(\finseta)\bigr)}^2 y\bigl(b(A)\bigr),\end{equation*}
with $\tilde{f}_k(b(\finseta))= f_k(A) + \sum_{\finsetc \in \parentp_d \suchthat \finseta \subset \finsetc}M_{\finsetc}w_{\finsetc}(A)$
and $\tilde{f}_k(b(\finsetb))=\tilde{f}_k(b(\finseta))+ M_{\finseta} w_{\finseta}(\finsetb)$.
So the condition $p_{\CC}(\confc(b(\finsetb)))=y(b(\finsetb))$ may be written as
\begin{equation*}p_{\CC}\Bigl(\tilde{f}_k\bigl(b(\finseta)\bigr)+ M_{\finseta} w_{\finseta}(\finsetb)\Bigr)=U_k\norm{\tilde{f}_k(b(\finsetb))}^2 y\bigl(b(\finsetb)\bigr),\end{equation*}
which is equivalent to
\begin{multline} \label{eqthree}
 M_{\finseta} p_{\CC}\bigl(w_{\finseta}(\finsetb)\bigr)
 = U_k\bignorm{\tilde{f}_k\bigl(b(\finseta)\bigr)}^2 \Bigl(y\bigl(b(\finsetb)\bigr)-y\bigl(b(A)\bigr)\Bigr)
\\ + U_k \Bigl(\bignorm{\tilde{f}_k(b(\finsetb))}^2-\bignorm{\tilde{f}_k\bigl(b(\finseta)\bigr)}^2\Bigr)y\bigl(b(\finsetb)\bigr),
\end{multline}
with \begin{equation*}\bignorm{\tilde{f}_k(b(\finsetb))}^2-\bignorm{\tilde{f}_k\bigl(b(\finseta)\bigr)}^2
= M_{\finseta} \Bigl(2\bigl\langle \tilde{f}_k\bigl(b(\finseta)\bigr),w_{\finseta}(\finsetb) \bigr\rangle
+M_{\finseta}\bignorm{w_{\finseta}(\finsetb)}^2 \Bigr).\end{equation*}

When $y(b(\finsetb))-y(b(\finseta))=0$ and $M_{\finseta} \neq 0$,
Equation~\ref{eqthree} simplifies to Equation~\ref{eqfour}, which also holds in the closure $C_{\tanghcyll}$.
In particular, we have $p_{\CC}(w_{\finseta}^0(b(\finsetb)))=0$. So we also get $|p_{\RR}(w_{\finseta}^0(b(\finsetb)))|\geq \eta$.

When 
$M_{\finseta} \neq 0$,
Equation~\ref{eqthree} is equivalent to
\begin{multline} \label{eqdeflambda}
p_{\CC}\bigl(w_{\finseta}(\finsetb)\bigr)-U_k \Bigl(2\bigl\langle \tilde{f}_k\bigl(b(\finseta)\bigr),w_{\finseta}(\finsetb) \bigr\rangle
+M_{\finseta} \bignorm{w_{\finseta}(\finsetb)}^2 \Bigr)y\bigl(b(\finsetb)\bigr) 
\\ = \frac{U_k}{M_{\finseta}} \bignorm{\tilde{f}_k\bigl(b(\finseta)\bigr)}^2 \Bigl(y\bigl(b(\finsetb)\bigr)-y\bigl(b(A)\bigr)\Bigr).  
\end{multline}
It tells that the left-hand side is colinear to $ \norm{\tilde{f}_k(b(\finseta))}^2 (y(b(\finsetb))-y(b(A)))$, and that the scalar product of these two vectors is nonnegative. This remains true in the closure $C_{\tanghcyll}$. When $y(b(\finsetb))-y(b(\finseta))\neq 0$, this uniquely defines $\lambda_{\finseta}=\lambda_{\finseta}(\finsetb,c)$ such that Equation~\ref{eqfive} holds for $\finsetb$. 
Furthermore, $\lambda_{\finseta}(\finsetb,c)$ is continuous on $C_{\tanghcyll}$,
we have \begin{equation*} \lim_{U_k \to 0}\lambda_{\finseta}(\finsetb,c) =\frac{\norm{p_{\CC}(w_{\finseta}(\finsetb))}}{\norm{\tilde{f}_k(b(\finseta))}^2\norm{y(b(\finsetb))-y(b(A))}},\end{equation*}
and $\lambda_{\finseta}(\finsetb,c)={U_k}/{M_{\finseta}}$ when $M_{\finseta} \neq 0$.
In particular, if $\finsetb^{\prime}$ is another univalent kid of $\finseta$ such that  $y(b(\finsetb^{\prime}))-y(b(\finseta))\neq 0$ and if $M_{\finseta} \neq 0$, then $\lambda_{\finseta}(\finsetb,c)=\lambda_{\finseta}(\finsetb^{\prime},c)$. This remains true in the closure $C_{\tanghcyll}$ when $M_{\finseta}=0$. This allows us to define the continuous $\lambda_{\finseta}=\lambda_{\finseta}(\finsetb,c)$ for the subset $\parntpplambda$ of $\parntplambda$ made of the sets $\finseta$ that have kids with basepoints on different components. These $\lambda_{\finseta}$ satisfy $\lambda_{\finseta}M_{\finseta} =U_{k(A)}$ when the parameters $\mu_D$ do no vanish. This also remains true in the closure $C_{\tanghcyll}$ when $M_{\finseta}=0$.
Observe $\parentp_{\cutx} \subseteq \parntpplambda$. For $\finseta \in \parntplambda \setminus \parntpplambda$, there exists $M \in \parentp_{\cutx}$ such that $M \subset A$. Define \begin{equation*}\lambda(\finseta,M)=\tilde{\lambda}_M\prod_{\finsetd \in \parentp_d \suchthat  M \subset \finsetd \subset A} \mu_{\finsetd}.\end{equation*}
So we have $\lambda(A,M)M_{\finseta} =U_{k(A)}$ and $\lambda(A,M)=\lambda(A,M^{\prime})$ for any subset $M^{\prime}$ of $A$ in $\parentp_{\cutx}$, when $M_{\finseta}\neq 0$. So $\lambda(A,M)=\lambda(A,M^{\prime})$ on $C_{\tanghcyll}$, and we can set $\lambda_{\finseta}=\lambda(A,M)$.
The other properties of the parameters $\lambda_{\finseta}$
are obvious when the parameters $\mu_D$ do no vanish. So they hold in $C_{\tanghcyll}$, and Lemma~\ref{lemdeflambda} is proved. 
(Note that the set $\parentp_{\cutx}$ of possibly separating sets is the subset of $\parntplambda$ consisting of its minimal sets with respect to the inclusion.)
\eop

\begin{lemma}
\label{lemfromczerocoor}
The configuration \begin{equation*}\confc^0=\psi\left((0)_{k \in \underline{\sigma}}, (0)_{{\finseta} \in \parentp_d},(f^0_k)_{k \in \underline{\sigma}},(w^0_{\finseta})_{{\finseta} \in \parentp_d} \right)\end{equation*} of $C^f_{\vertsetv}(N_{\infty},\tanghcyll,\confc^0)$ is such that 
\item \begin{itemize}
       \item we have $p_{\CC} \circ w^0_{\finseta}(\finsetb)=0$ for any univalent kid $\finsetb$ of $\finseta$ if $\finseta \in \parentp_d \setminus \parentp_{\cutx}$,
\item for any $\finseta \in \parentp_{\cutx}$, there exists $\lambda_{\finseta}^0 \geq 0$ such that
\begin{equation*}p_{\CC}\circ w^0_{\finseta}(\finsetb)=\lambda_{A}^0 \bignorm{{f}^0_{k(A)}(A)}^2 \Bigl(y\bigl(b(\finsetb)\bigr)-y\bigl(b(A)\bigr)\Bigr)\end{equation*} for any univalent kid $\finsetb$ of $A$.
      \end{itemize}
\end{lemma}
\bp
Lemma~\ref{lemdeflambda} implies that $\lambda_{\finseta}^0=0$ when $A \in \parntplambda \setminus \parentp_{\cutx}$.
\eop

\begin{example} \label{exaproofpropenfertwo}
Let us go back to Example~\ref{exaproofpropenferone}, where $\norm{y(v_2)-y(v_1)}$ is nonzero. When $\finseta=\finsetv$ and $\finsetb=\{v_2\}$, we have $\tilde{f}(v_2)=f(v_1) + \mu w(v_2)$. When $\confc$ is injective, Equation~\ref{eqthree} is equivalent to
\begin{equation*}p_{\CC}\bigl(w(v_2)\bigr)= \lambda \bignorm{\tilde{f}(v_1)}^2 \bigl(y(v_2)-y(v_1)\bigr) + u \Bigl(2\bigl\langle \tilde{f}(v_1),w(v_2) \bigr\rangle
+\mu \norm{w(v_2)}^2 \Bigr)y(v_2), \end{equation*}
with $\lambda=\lambda_{\finsetv}=\frac{u}{\mu}$. We get $p_{\CC}(w^0(v_2))= \lambda^0 (y(v_2)-y(v_1))$.
\end{example}

Let us now define the oriented tree $\CT^0$ (as in Definition~\ref{deffirsttree}) of the statement of Theorem~\ref{thmcomptang}.
\begin{notation} \label{nottreeconftang}
The set $E(\CT^0)$ of edges of $\CT^0$ is in one-to-one correspondence with $\{u_i\}_{i \in \underline{\sigma}} \cup \{ \mu_{\finseta}\}_{\finseta \in \parntpwidehat \setminus \parentp_x} \cup \{\lambda_{\finseta}\}_{\finseta \in \parntplambda \setminus \parentp_x}\cup \{\tilde{\lambda}_{\finseta}=\lambda_{\finseta}\mu_{\finseta}\}_{\finseta \in \parentp_{x}}$, and its edges are labeled by these variables.
So $E(\CT^0)$ is in one-to-one correspondence with the disjoint union of $\underline{\sigma}$, $\parentp_x$, and two disjoint copies of $\parntplambda \setminus \parentp_x$. (Recall that $\parntpwidehat$ is the set of elements of $\parentp_d$ that contain or are equal to an element of $\parentp_{\cutx}$.) 

The set of vertices of $\CT(\parentp)$ is in one-to-one correspondence with the disjoint union 
\begin{equation*}\parntpwidehat \sqcup \parentp_s^+ \sqcup  \{\rroot_1\} \sqcup \bigl\{\rroot_{\finseta} \suchthat \finseta \in \parntplambda \setminus \parentp_x\bigr\}.\end{equation*}
Its elements label the vertices. 

The edge labeled by $u_i$ ends at the vertex labeled by $V(i)^+$. It starts at the vertex labeled by $V(i-1)^+$ when $i>1$, and at the univalent vertex labeled by $\rroot_1$ when $i=1$. 

For $\finseta \in \parntplambda \setminus \parentp_x$, the edge labeled by $\lambda_{\finseta}$ starts from the univalent vertex labeled by $\rroot_{\finseta}$, and it goes to the vertex labeled by $\finseta$.

For $\finseta \in \parntpwidehat \setminus \parentp_x$ (resp. for $\finseta \in \parentp_x$), the edge labeled by $\mu_{\finseta}$ (resp. by $\tilde{\lambda}_{\finseta}$) starts at the vertex $\finseta$,  
It ends at the mother $\mother(\finseta)$ of $\finseta$. 
See Figure~\ref{figTconfzero} for an example of a tree $\CT^0$.
\end{notation}
\bfig
\centering
\begin{equation*}
\begin{tikzpicture}
\useasboundingbox (-1,0) rectangle (10.6,5);
\draw  (2,4) -- (6.2,4) (7.6,4) -- (9,4) (-.6,1) -- (-.6,2) -- (2,3) (.4,1) -- (2,3) -- (-.6,3.4) (1.4,1) -- (2,3) -- (3.4,4) (3.4,1) -- (3.4,4) (5.4,1) -- (4.8,3) (4.2,1) -- (4.8,3) -- (4.8,4) (8,0) --(8,1) -- (9,2) (9,0) -- (9,4) (9,3) -- (10.4,3) (9,2) -- (10.4,2);
\fill (2,4) circle (1.5pt) (3.4,4) circle (1.5pt) (4.8,4) circle (1.5pt) (6.2,4) circle (1.5pt) (7.6,4) circle (1.5pt) (9,4) circle (1.5pt)  (-.6,1)  circle (1.5pt) (-.6,2)  circle (1.5pt) (2,3)  circle (1.5pt) (.4,1)  circle (1.5pt) (1.4,1)  circle (1.5pt) (-.6,3.4)  circle (1.5pt) (3.4,1)  circle (1.5pt) (3.4,2.5) circle (1.5pt) (4.2,1)  circle (1.5pt) (5.4,1)  circle (1.5pt) (5.1,2) circle (1.5pt) (4.8,3) circle (1.5pt) (8,0)  circle (1.5pt) (8,1)  circle (1.5pt) (9,2)  circle (1.5pt) (9,0)  circle (1.5pt) (9,1) circle (1.5pt) (9,3) circle (1.5pt) (10.4,3)  circle (1.5pt) (10.4,2)  circle (1.5pt);
\draw [->] (.7,3.2) node[above]{\scriptsize  $\lambda_{m(A_1)}$} (-.6,3.4) node[above]{\scriptsize  $\rroot_{m(A_1)}$} (-.6,3.4) -- (.7,3.2);
\draw [->] (4.6,2.2) node[left]{\scriptsize  $\lambda_{m(A_5)}$} (4.2,1) -- (4.5,2);
\draw [->] (9.7,1.9) node[above]{\scriptsize  $\lambda_{m(A_9)}$} (10.4,2) -- (9.7,2);
\draw [->] (9.85,2.9) node[above]{\scriptsize  $\lambda_{m(m(A_9))}$} (10.4,3) -- (9.7,3);
\draw [->] (2,4.25) node{\scriptsize  $\rroot_1$} (2.7,4.25) node{\scriptsize $u_1$} (3.4,4.25) node{\scriptsize $V(1)^+$} (2,4) -- (2.7,4);
\draw [->] (4.1,4.25) node{\scriptsize  $u_2$} (4.8,4.25) node{\scriptsize $V(2)^+$} (3.4,4) -- (4.1,4);
\draw [->] (5.5,4.25) node{\scriptsize  $u_3$} (4.8,4) -- (5.5,4);
\draw [->] (8.3,4.25) node{\scriptsize  $u_{\sigma}$}  (9.2,4.25) node{\scriptsize $V(\sigma)^+$} (7.6,4) -- (8.3,4);
\draw (6.9,4) node{ $\dots$};
\draw[->] (-1,1.5) node{\scriptsize $\lambda_{A_1}$} (-.9,2.1) node{\scriptsize  $A_1$} (-.6,1) -- (-.6,1.5);
\draw[->] (1.2,2) node[left]{\scriptsize$\tilde{\lambda}_{A_2}$} (.4,1) -- (1.2,2);
\draw[->] (1.65,1.7) node[right]{\scriptsize $\tilde{\lambda}_{A_3}$} (1.4,1) -- (1.7,2);
\draw[->] (3.5,1.75) node[left]{\scriptsize $\lambda_{A_4}$} (3.4,1) -- (3.4,1.75);
\draw[->] (5.65,1.5) node{\scriptsize $\lambda_{A_5}$} (5.4,2) node{\scriptsize  $A_5$} (5.85,.9) node{\scriptsize$\rroot_{A_5}$} (5.4,1) -- (5.25,1.5);
\draw[->] (7.6,.5) node{\scriptsize $\lambda_{A_8}$} (8,0) -- (8,.5);
\draw[->] (8.65,.5) node{\scriptsize $\lambda_{A_9}$} (9,0) -- (9,.5);
\draw[->] (9.4,1.5) node{\scriptsize $\mu_{A_9}$} (9,1) -- (9,1.5);
\draw[->] (8.35,2.5) node{\scriptsize  $\mu_{\mother(A_9)}$} (9,2) -- (9,2.5);
\draw[->] (8.2,3.5) node{\scriptsize  $\mu_{\mother(\mother(A_9))}$} (9,3) -- (9,3.5);
\draw[->] (8.1,1.6) node{\scriptsize  $\mu_{A_8}$} (8,1) -- (8.5,1.5);
\draw[->] (5.4,2.5) node{\scriptsize  $\mu_{A_5}$} (5.4,3) node{\scriptsize  $\mother(A_5)$} (5.1,2) -- (4.95,2.5);
\draw[->] (5.5,3.5) node{\scriptsize  $\mu_{\mother(A_5)}$} (4.8,3) -- (4.8,3.5);
\draw[->] (3.8,3.35) node{\scriptsize  $\mu_{A_4}$} (3.4,2.5) -- (3.4,3.25);
\draw[->] (2.7,3.5) node[left]{\scriptsize  $\mu_{\mother(A_1)}$} (2,3) -- (2.7,3.5);
\draw[->] (.7,2.45) node[above]{\scriptsize  $\mu_{A_1}$} (-.6,2) -- (.7,2.5);
\draw (-.9,.9) node{\scriptsize$\rroot_{A_1}$} (.2,.9) node{\scriptsize $A_2$} (1.2,.9) node{\scriptsize $A_3$} (3.1,.9) node{\scriptsize$\rroot_{A_4}$} (3.1,2.5) node{\scriptsize $A_4$};
\end{tikzpicture}\end{equation*}
\caption{A tree $\CT^0=\CT(\confc^0)$}
\label{figTconfzero}

\end{figure}

According to Lemma~\ref{lemdeflambda}, the variables of \begin{equation*}\{u_i\}_{i \in \underline{\sigma}} \cup \{ \mu_{\finseta}\}_{\finseta \in \parntpwidehat \setminus \parentp_x} \cup \{\lambda_{\finseta}\}_{A \in \parntplambda \setminus \parentp_x}\cup \{\tilde{\lambda}_{\finseta}=\lambda_{\finseta}\mu_{\finseta}\}_{A \in \parentp_{x}}\end{equation*} satisfy the equations that define
$\setparamx(\CT^0)$.

\begin{lemma}\label{lemdimsetparamtzero}
 The dimension of $\setparamx(\CT^0)$ is $\dimd(\CT^0)=\cardbig{\parentp_s} +\cardbig{\parntpwidehat} -\cardbig{\parentp_x}$.
\end{lemma}
\bp
Recall $\dimd(\CT^0)=\cardbig{E(\CT^0)} - \cardbig{L(\CT^0)} + 1$ from Lemma~\ref{lemsysred}. 
The set of leaves different from $\rroot_1$ is in one-to-one correspondence with $\parntplambda$.
\eop

\begin{lemma}
\label{lemcoordparentpd}
Let $\confc^0=\psi\left((0), (0),(f^0_i),(w^0_{\finsetb})\right)$ be a configuration of ${C}_{V}(\rats(\hcylc))$ that satisfies the equations of Lemmas~\ref{lemfktwo} and \ref{lemdeflambda}. Then $\confc^0$ belongs to $ C^f_{\vertsetv}(N_{\infty},\tanghcyll,\confc^0)$.

Let $A$ be a univalent element of $\parentp_d$.
Let $\kids_{u}(A)$ denote the set of univalent kids of $A$ that do not contain $b(A)$, and let $\kids_{t}(A)$ denote the set of nonunivalent kids of $A$ that do not contain $b(A)$. Recall Notation~\ref{nothornnormaliz}.
For any element $\finsetd$ of $\kids_t(\finseta)$, let $B^{\prime}_{\finsetd}$ denote the open ball $\mathring{\ballb}(w^0_{\finseta}(\finsetd),\varepsilon)$ of radius $\varepsilon$ with center $w^0_{\finseta}(\finsetd)$ in $\RR^3$ if $\finsetd$ is not normalizing, and let $B^{\prime}_{\finsetd}$ denote the set of elements $w_{\finseta}(\finsetd)$
of $\mathring{\ballb}(w^0_{\finseta}(\finsetd),\varepsilon)$ such that $|p_{\RR}(w_{\finseta}(\finsetd))|=1$ (resp. such that $|p_{\CC}(w_{\finseta}(\finsetd))|=1$) if $\finsetd$ is v-normalizing (resp. if $\finsetd$ is h-normalizing).
For any element $\finsetd$ of $\kids_u(\finseta)$, let $J^{\prime}_{\finsetd}$ denote the interval $\left]p_{\RR}(w^0_{\finseta}(\finsetd))-\varepsilon,p_{\RR}(w^0_{\finseta}(\finsetd))+\varepsilon\right[$ if $\finsetd$ is not v-normalizing, and set $J^{\prime}_{\finsetd}=\{p_{\RR}(w^0_{\finseta}(\finsetd))\}$ if $\finsetd$ is v-normalizing.
Set \begin{equation*}W_{\finseta}^L=\prod_{\finsetd \in \kids_{u}(A)}J^{\prime}_{\finsetd} \times \prod_{\finsetd \in \kids_{t}(\finseta)}B^{\prime}_{\finsetd}.\end{equation*}
When $A$ is a nonunivalent element of $\parentp_d$, set $W^{\tanghcyll}_{\finseta}=W_{\finseta}$.
Recall that $W^{\tanghcyll}_k$ has been introduced in Lemma~\ref{lemfktwo}.

Assume that $2\varepsilon < \lambda^0_{\finseta}$ for all $\finseta \in \parentp_x$. Let $\parentp_{x,hn}$ be the set of elements $\finseta$ of $\parentp_x$ such that $k_n(\finseta)$ is a univalent horizontally normalizing kid. Set 
\begin{equation*}W= \left[0,\varepsilon\right[^{\parentp_d \setminus \parntpwidehat} \times \left( \prod_{k \in \underline{\sigma}} W^{\tanghcyll}_k  \right) \times \left(\prod_{{\finseta} \in \parentp_d} W^{\tanghcyll}_{\finseta} \right) \times \prod_{\finseta \in \parentp_x \setminus \parentp_{x,hn}}\left]\lambda^0_{\finseta}-\varepsilon,\lambda^0_{\finseta}+\varepsilon \right[.\end{equation*}
Recall Notation~\ref{nottreeconftang}. There exists a smooth map from $\left[0,\varepsilon\right[^{E(\CT^0)} \times W$ to ${C}_{V}(\rats(\hcylc))$ which restricts to $\bigl(\left[0,\varepsilon\right[^{E(\CT^0)} \cap \setparamx(\CT^0)\bigr) \times  W $ as a continuous injective map $\varphi$, whose image is an open neighborhood $N^f_{\tanghcyll}(\confc^0)$ of $\confc^0$ in $C^f_{\vertsetv}(N_{\infty},\tanghcyll,\confc^0)$.
\end{lemma}
\bp Here and in Notation~\ref{nottreeconftang}, the parameters $u_i$ and $\mu_{\finsetb}$ of $E(\CT^0)$ are the initial parameters of $\confc$ in the chart $\psi$ of ${C}_{V}(\rats(\hcylc))$. The parameters $\lambda_{\finseta}$ are defined in Lemma~\ref{lemdeflambda}. Lemma~\ref{lemdeflambda} implies that the parameters $\lambda_{\finseta}$ satisfy the equations of $\setparamx(\CT^0)$ for configurations in $C^f_{\vertsetv}(N_{\infty},\tanghcyll,\confc^0)$.  The factor of $\left[0,\varepsilon\right[^{\parentp_d \setminus \parntpwidehat}$ of $W$ contains the parameters $\mu_{\finsetb}$ for 
$\finsetb \in \parentp_d \setminus \parntpwidehat$. 

Lemma~\ref{lemfktwo} shows how to express the parameter $f_k$ of $\psi^{-1}(\confc) $, for configurations $\confc$ in $C^f_{\vertsetv}(N_{\infty},\tanghcyll,\confc^0)$, as a smooth function of $W \times \left[0,\varepsilon\right[^{\underline{\sigma}}$, where the factor $\left[0,\varepsilon\right[^{\underline{\sigma}}$ contains the $u_i$.

We now construct the $w_{\finseta}$ as smooth maps from $\left[0,\varepsilon\right[^{E(\CT^0)} \times W$ to $\cinjuptd{\kids(A)}{\RR^3}$,
in order to finish constructing a smooth map 
\begin{equation*}\varphi \colon \left[0,\varepsilon\right[^{E(\CT^0)} \times W \to C_{V(\Gamma)}\bigl(\rats(\hcylc)\bigr).\end{equation*}
The coordinates of the nonunivalent kids of $A$ and the vertical coordinates of the univalent kids of $A$ are part of $W$, and we do not change them.
We only need to determine the horizontal coordinate $p_{\CC}(w_{\finseta}(\finsetb))$ for any univalent kid $\finsetb$ of $A$ as a smooth map. 

Note that $\lambda_{\finseta}$ is a free parameter in $\left[0,\varepsilon\right[^{E(\CT^0)} \times W$ when $\finseta \in \parntplambda \setminus \parentp_{x,hn}$. 

Let $\finseta$ be a univalent element of $\parentp_d$. Set $k=k(A)$.
By induction, assume that $w_{\finsetc}$ and $\mu_{\finsetc}$ have already been constructed as smooth functions of $\left[0,\varepsilon\right[^{E(\CT^0)} \times W$ for any $\finsetc \in \parentp_d$ that contains $\finseta$, so that $\tilde{f}_k(b(\finseta))$
is a smooth function of $\left[0,\varepsilon\right[^{E(\CT^0)} \times W$ defined by
 the expression from Notation~\ref{notpreldeflambda}:
\begin{equation}
\label{eqftilde}
\tilde{f}_k\bigl(b(\finseta)\bigr)= f_k(\finseta) + \sum_{\finsetc \in \parentp_d \suchthat \finseta \subset \finsetc}M_{\finsetc}w_{\finsetc}(\finseta).
\end{equation}

If $\finseta \notin \parntplambda$, then 
Equation~\ref{eqfour} determines $p_{\CC}(w_{\finseta}(\finsetb))$ as a smooth implicit function of ($p_{\RR}(w_{\finseta}(\finsetb))$, which is determined by $W$, $\tilde{f}_k(b(\finseta))$, the $u_i$, and the $\mu_{\finsetc}$ for the $\finsetc \in \parentp_d$ such that $\finseta \subseteq \finsetc$.
If $\finseta \in \parntplambda$, then Equation~\ref{eqfive} determines $p_{\CC}(w_{\finseta}(\finsetb))$ as a smooth implicit function of the same parameters ($W$, $\tilde{f}_k(b(\finseta))$, the $u_i$, and the $\mu_{\finsetc}$ for the $\finsetc \in \parentp_d$ such that $\finseta \subseteq \finsetc$), and $\lambda_{\finseta}$.
The parameter $\mu_{\finseta}$ is among the parameters unless $\finseta \in \parentp_x$. If $\finseta \in \parentp_x \setminus  \parentp_{x,hn}$, then $\mu_{\finseta}$ is determined as
${\tilde{\lambda}_{\finseta}}/{\lambda_{\finseta}}$, where $\lambda_{\finseta}$ is a parameter of $W$. So $\mu_{\finseta}$ is this smooth function of our parameters. 
If $\finseta \in \parentp_{x,hn}$, then $p_{\CC}(w_{\finsetc}(k_n(\finseta)))$ is still a smooth implicit function of $W$, $\tilde{f}_k(b(\finseta))$, the $u_i$, the $\mu_{\finsetc}$ for the $\finsetc \in \parentp_d$ such that $\finseta \subseteq \finsetc$, 
and $\lambda_{\finseta}$, and the normalizing condition $|p_{\CC}(w_{\finseta}(k_n(\finseta)))|=1$ determines $\lambda_{\finseta}$ as a smooth implicit fonction of the given parameters.

We have constructed a smooth function $\varphi \colon \left[0,\varepsilon\right[^{E(\CT^0)} \times W \to C_{V(\Gamma)}(\rats(\hcylc)).$
Set \begin{equation*}\mathring{W}= \left]0,\varepsilon\right[^{\parentp_d \setminus \parntpwidehat} \times \times \left( \prod_{k \in \underline{\sigma}} W^{\tanghcyll}_k  \right) \times \left(\prod_{{\finseta} \in \parentp_d} W^{\tanghcyll}_{\finseta} \right) \times \prod_{\finseta \in \parentp_x \setminus \parentp_{x,hn}}\left]\lambda^0_{\finseta}-\varepsilon,\lambda^0_{\finseta}+\varepsilon \right[.\end{equation*}
Lemma~\ref{lemdeflambda} implies that $\varphi$ maps $\bigl(\left]0,\varepsilon\right[^{E(\CT^0)} \cap \setparamx(\CT^0) \bigr) \times \mathring{W}$ to $C^f_{\vertsetv}(N_{\infty},\tanghcyll,\confc^0)$ (up to reducing $\varepsilon$).
Since the closure of $\mathring{\setparamx}(\CT^0)$ in $\left[0, \infty\right[^{E(\CT^0)}$ is ${\setparamx}(\CT^0)$, according to Lemma~\ref{lemclosuresetparamx}, we have \begin{equation*}\varphi \Bigl(\bigl( \left[0,\varepsilon\right[^{E(\CT^0)} \cap \setparamx(\CT^0)\bigr) \times W\Bigr) \subseteq C^f_{\vertsetv}(N_{\infty},\tanghcyll,\confc^0),\end{equation*} too.
In particular, any $\confc^0$ that satisfies the equations of Lemmas~\ref{lemfktwo} and \ref{lemdeflambda} is in $C^f_{\vertsetv}(N_{\infty},\tanghcyll,\confc^0)$.

Lemma~\ref{lemdeflambda} also implies $N^f_{\tanghcyll}(\confc^0) \subseteq \varphi \bigl(\bigl( \left[0,\varepsilon\right[^{E(\CT^0)} \cap \setparamx(\CT^0)\bigr) \times W\bigr)$ for an open neighborhood $N^f_{\tanghcyll}(\confc^0)$ of $\confc^0$ in $C^f_{\vertsetv}(N_{\infty},\tanghcyll,\confc^0)$.
The injectivity of the restriction of $\varphi$ to $\bigl(\bigl( \left[0,\varepsilon\right[^{E(\CT^0)} \cap \setparamx(\CT^0)\bigr) \times W\bigr)$ comes from the fact that Equation~\ref{eqfive} determine the parameters $\lambda_{\finseta}$ for $A \in \parentp_x$ and the equations of Lemma~\ref{lemdeflambda} determine the others.\footnote{We can remove these other parameters $\lambda_{\finsetc}$, for $\finsetc \in \parntplambda\setminus \parentp_x$, from our parametrization for the statement. However, it was more convenient to keep them for the proof since the definition of the $\lambda_{\finseta}$, for $\finseta \in \parentp_{x,hn}$, involves the $\lambda_{\finsetc}$ for the $\finsetc$ that contain $\finseta$. These $\lambda_{\finsetc}$ are functions of the $\lambda_{\finseta}$ for which $\finseta \in \parentp_x$.}
\eop

\bpo{Proof of Proposition~\ref{proppredefsepone}}
Lemmas~\ref{lemfkone}, \ref{lemfkonebis}, \ref{lemptauconfczero}, and~\ref{lemfromczerocoor} show that a configuration $\confc^0_{V(\Gamma)}$ of $C_{\tanghcyll}$ must satisfy the conditions of the statement of Proposition~\ref{proppredefsepone}. According to Lemma~\ref{lemcoordparentpd}, if a configuration $\confc^0_{V(\Gamma)}$ satisfies the conditions of the statement of Proposition~\ref{proppredefsepone}, its restriction $\confc^0$ to $\vertsetv=\pbl(\confc_{V(\Gamma)})^{-1}(\infty)$ is in $C^f_{\vertsetv}(N_{\infty},\tanghcyll,\confc^0)$, so $\confc^0_{V(\Gamma)}$ is in $C_{\tanghcyll}$.
\eop

\bpo{Proof of Proposition~\ref{proppredefseptwo}}
The strata of Proposition~\ref{proppredefseptwo} are smooth submanifolds of $C_{V(\Gamma)}(\rats(\hcylc))$ since they correspond to the locus where all the variables associated to $\CT^0$ in Notation~\ref{nottreeconftang} are zero. (The parameters $\lambda_A$ are either zero or not zero on the whole stratum, in the charts of Lemma~\ref{lemcoordparentpd}.)
Let us assume $\vertsetv=\vertsetv(\Gamma)$. Then the codimension of the stratum of $\confc^0$ is $\dimd(\CT^0) +\cardbig{(\parentp_d \setminus \parntpwidehat)}$, which is $\cardlef{\parentp_s} + \cardlef{\parentp_d} -\cardlef{\parentp_x}$, according to Lemma~\ref{lemdimsetparamtzero}.
\eop

\begin{lemma}
\label{lemsubsecTface}
The codimension-one faces of ${C}^f\!(\rats(\hcylc),\tanghcyll;\Gamma)$ are those listed in Theorem~\ref{thmcomptang}.
In a neighborhood of these faces, ${C}^f\!(\rats(\hcylc),\tanghcyll;\Gamma)$ has the structure of a smooth manifold with boundary.

Let $\syms=\syms(\upvec)$ be the orthogonal reflection of $\RR^3$ with respect to the horizontal plane.
A configuration \begin{equation*}\confc^0_{V(\Gamma)}\vert_{V}=\confc^0=\left(T_{0}\phi_{\infty} \circ f_1^0,\left( T_0 \phi_{\infty \ast} \circ w^0_{\finseta}\right)_{A \in \parentp_{\cutx}}, \left(\lambda^0_{\finseta}\right)_{A \in \parentp_{\cutx}}\right)\end{equation*}
of a $T$-face is the limit at $t=0$ of a family of injective configurations $c(t)_{t \in \left]0,\varepsilon\right[}$ on the vertical parts of the tangle, far above or far below, such that $c(t)\vert_{A}=\syms \circ w^0_{\finseta}$ up to dilation and translation for any $\finseta \in \parentp_{\cutx}(=\parentp_d=\parentp_x)$. In particular, for an edge $e=(v_1,v_2)$ whose vertices are in $A$, we have \begin{equation*}p_{\tau}\circ p_e(\confc^0)=\frac{\syms \circ w^0_{\finseta}(v_2) -\syms \circ w^0_{\finseta}(v_1)}{\bignorm {\syms \circ w^0_{\finseta}(v_2) -\syms \circ w^0_{\finseta}(v_1) }}.\end{equation*}
For an edge $e=(v_1,v_2)$ whose vertices are in different kids of $V$, we have
\begin{multline*}p_{\tau}\circ p_e(\confc^0)=\frac{\phi_{\infty} \circ f_1^0(v_2) -\phi_{\infty} \circ f_1^0(v_1)}{\norm{ \phi_{\infty} \circ f_1^0(v_2) -\phi_{\infty} \circ f_1^0(v_1) }}\\=\frac{\norm{f_1^0(v_1)}^2f_1^0(v_2)- \norm{f_1^0(v_2)}^2f_1^0(v_1)}{\norm{\norm{f_1^0(v_1)}^2f_1^0(v_2)- \norm{f_1^0(v_2)}^2f_1^0(v_1)}}.\end{multline*}
For an edge $e=(v_1,v_2)$ such that $v_1 \in V$ and $v_2 \notin V$, we have \begin{equation*}p_{\tau}\circ p_e(\confc^0)=-\frac{f_1^0(v_1)}{\norm{f_1^0(v_1)}}.\end{equation*} For an edge $e=(v_1,v_2)$ such that $v_2\in V$ and $v_1 \notin V$, we have $p_{\tau}\circ p_e(\confc^0)=\frac{f_1^0(v_2)}{\norm{f_1^0(v_2)}}$.
\end{lemma}
\bp
Let $\partial_{\infty}({C}^f\!(\rats(\hcylc),\tanghcyll;\Gamma))$ be the subspace of ${C}^f\!(\rats(\hcylc),\tanghcyll;\Gamma)$ consisting of the configurations as above that map at least one univalent vertex to $\infty$.
Outside this subspace, the spaces ${C}^f\!(\rats(\hcylc),\tanghcyll;\Gamma)$ and ${C}(\rats(\hcylc),\tanghcyll;\Gamma)$ are the same, and ${C}^f\!(\rats(\hcylc),\tanghcyll;\Gamma)$ has the structure of a smooth manifold with ridges.
Recall from Remark~\ref{rkstrataTf} that the only codimension-one new parts are the $T$-faces of Theorem~\ref{thmcomptang}, where $\parentp_s=\{\finsetb \sqcup \sqcup_{j\in I}\finsetb_j\}$ and $\parentp_d=\parentp_{\cutx}=\parentp_x=\{\finsetb_j \suchthat j\in I\}$.
Lemmas~\ref{lemcoordparentpd} and~\ref{lemcodimonex} imply that these strata arise as codimension-one faces of ${C}^f\!(\rats(\hcylc),\tanghcyll;\Gamma)$, along which ${C}^f\!(\rats(\hcylc),\tanghcyll;\Gamma)$
is a smooth manifold with boundary. 

Lemma~\ref{lemfkone} implies ${f}^0_1(b)=\pm \norm{{f}^0_1(b)} \upvec$ for any univalent vertex $b$ of $\finsetb \sqcup \sqcup_{j\in I}\finsetb_j$.
Lemma~\ref{lemptauconfczero} implies that the restriction of $\confc^0$ to any $\finseta$ of $\parentp_{\cutx}=\{\finsetb_j \suchthat j\in I\}$ is represented by $\syms \circ w_{\finseta}^0$, up to translation and dilation, as a configuration of $\RR^3$.
Furthermore, Lemma~\ref{lemfromczerocoor} implies
\begin{equation*}p_{\CC}\bigl(w^0_{\finseta}(b)\bigr)
 = \lambda^0_{\finseta}\bignorm{{f}^0_1\bigl(b(A)\bigr)}^2 \Bigl(y(b)-y\bigl(b(A)\bigr)\Bigr)\end{equation*}
for any $A \in \parentp_{\cutx}$ and for any $b \in A$.
Therefore, the configuration $\confc^0$ is the limit of the following family $\confc(t)$ of configurations, indexed by $t \in \left]0,\varepsilon\right[$, where we have\\
$\confc(t)(b)=\frac{1}{t}\frac{{f}^0_1(b)}{\norm{{f}^0_1(b)}^2}$ for any trivalent vertex of $\finsetb$,\\
$\confc(t)(b)=(y(b),0)+\frac{1}{t}\frac{{f}^0_1(b)}{\norm{{f}^0_1(b)}^2}$ for any univalent vertex of $\finsetb$, and\\
$\confc(t)(b)=(y(b(\finseta)),0)+\frac{1}{t}\frac{{f}^0_1(b(\finseta))}{\norm{{f}^0_1(b(\finseta))}^2} +\frac{\syms \circ w_{\finseta}^0(b)}{\lambda^0_{\finseta}\norm{{f}^0_1(b(A))}^2}$ for any vertex $b$ of an element $\finseta$ of $\parentp_{\cutx}=\{\finsetb_j \suchthat j\in I\}$. 

So
$\confc^0$ is the limit at $t=0$ of the family of injective configurations $\confc(t)_{t \in \left]0,\varepsilon\right[}$, and $\confc(t)\vert_{A}=\syms \circ w^0_{\finseta}$ up to dilation and translation, for any $\finseta \in \parentp_{\cutx}$.
\eopwobp

Theorem~\ref{thmcomptang} is now proved. \eop

\section{Variations of integrals on configuration spaces of long tangles}
\label{secvarzinf}

In this section, we prove Theorem~\ref{thmfstconsttang}. 

\begin{lemma}
\label{lemetactwoprime}
 For any two propagating forms $\omega$ and $\omega^{\prime}$ of $C_2(\rats)$ (as in Definition~\ref{defpropagatortwo}) that coincide on $\partial C_2(\rats)$, there exists a one-form $\eta$ of $C_2(\rats)$ that vanishes on $\partial C_2(\rats)$ such that $\omega^{\prime}=\omega+d\eta$.
 In particular, for any two homogeneous propagating forms $\omega$ and $\omega^{\prime}$ of $C_2(\rats)$ as in Definition~\ref{defprophomogen} that coincide on $\ST \ballb_{\rats}$, there exists a one-form $\eta$ of $C_2(\rats)$ that vanishes on $\partial C_2(\rats)$ such that $\omega^{\prime}=\omega+d\eta$.
\end{lemma}
\bp Exercise. See the proof of Lemma~\ref{lemetactwo}.
\eop

\begin{lemma}
\label{lemvarzinfunc}
The element $\Zinv_n(\hcylc,\tanghcyll,(\omega(i)))$
of $\Aavis_n(\sourcetl)$ is independent of the chosen homogeneous propagating forms $\omega(i)$ of $(C_2(\rats(\hcylc)),\tau)$, under the assumptions of Theorem~\ref{thmfstconsttang}.

More generally, if the $\omega(i)$ are only assumed to be homogeneous propagating forms of $C_2(\rats(\hcylc))$, then
$\Zinv_n(\hcylc,\tanghcyll,(\omega(i)))$ depends only on $(\hcylc,\tanghcyll \cap \hcylc,\tau)$ and on the restrictions of the $\omega(i)$ to $\ST \hcylc$. 
\end{lemma}
\bp
By Lemma~\ref{lemetactwoprime}, it suffices to prove that $\Zinv$ does not vary when $\omega(i)$ is changed to $\omega(i) + d\eta$ for a one-form $\eta$ on $C_2(\rats(\hcylc))$ that vanishes on $\partial C_2(\rats(\hcylc))$.
Let $\Omega_{\Gamma}=\bigwedge_{e \in E(\Gamma)}p_e^{\ast}(\omega(j_E(e)))$ and let $\tilde{\Omega}_{\Gamma}$ be obtained from $\Omega_{\Gamma}$ by replacing $\omega(i)$ by $\eta$.
The variations of the integrals $\int_{({C}(\rats(\hcylc),\tanghcyll;\Gamma),o(\Gamma))}\Omega_{\Gamma}$ are computed with
 Stokes' theorem, as the sum over the codimension-one faces $F$ of ${C}(\rats(\hcylc),\tanghcyll;\Gamma)$ of $\int_{F}\tilde{\Omega}_{\Gamma}$, as allowed by Lemma~\ref{lemStokesconftang}.

These faces are the faces listed in Theorem~\ref{thmcomptang}.
The arguments of Lemmas~\ref{lemdiscon}, \ref{lemedge}, \ref{lemsym}, \ref{lemihx}, \ref{lemstu} allow us to get rid of all the faces, except the faces in which some vertices are at $\infty$, and the faces $\facee(\check{\Gamma}_{\finseta},\tanghcyll,\Gamma)$ in which $\check{\Gamma}_{\finseta}$ is a connected diagram on $\RR$ (these faces are components of $\facee(\finseta,\Link,\Gamma)$ as in the proof of Lemma~\ref{leminvtwo})
and $\Gamma$ is a diagram that contains $\check{\Gamma}_{\finseta}$ as a subdiagram on a component $\source_j$ of $\tanghcyll$.
The contribution of the faces $\facee(\check{\Gamma}_{\finseta},\tanghcyll,\Gamma)$ is zero when $i \notin j_E(E(\check{\Gamma}_{\finseta}))$ for dimension reasons. It is zero when $i \in j_E(E(\check{\Gamma}_{\finseta}))$ because $\eta$ vanishes on $\partial C_2(\rats(\hcylc))$. So we are left with the faces for which some vertices are at $\infty$.
Let $F$ be such a face.
Let $\finsetv$ be the set of vertices mapped to $\infty$ in $F$, let $E_{\infty}$ be the set of edges between elements of $\finsetv$, and let $E_m$ denote the set of edges with one end in $\finsetv$. When $i\in j_E(E_{\infty} \cup E_m)$, the contribution vanishes 
because $\eta$ vanishes on $\partial C_2(\rats(\hcylc))$.

Assume $i\notin j_E(E_{\infty} \cup E_m)$.
The face $F$ is diffeomorphic to a product by \begin{equation*}\check{C}_{\vertsetv(\Gamma) \setminus \finsetv}\bigl(\crats(\hcylc),\tanghcyll;\Gamma\bigr),\end{equation*} whose dimension is
\begin{equation*}3 \cardbig{T(\Gamma) \cap \bigl(\vertsetv(\Gamma) \setminus \finsetv\bigr)} + \cardbig{U(\Gamma) \cap \bigl(\vertsetv(\Gamma) \setminus \finsetv\bigr)},\end{equation*}
of a space $C_{\finsetv}$ of dimension 
$3 \cardlef{T(\Gamma) \cap \finsetv} + \cardlef{U(\Gamma) \cap \finsetv} -1,$
along which \begin{equation*}\bigwedge_{e \in E_{\infty} \cup E_m}p_e^{\ast}\Bigl(\omega\bigl(j_E(e)\bigr)\Bigr)\end{equation*} has to be integrated.
The degree $2\cardlef{E_{\infty} \cup E_m}$ of this form is bigger than the dimension of $C_{\finsetv}$ as a count of half-edges shows. So the faces for which some vertices are at $\infty$ (including the $T$-faces) do not contribute either.
\eop

\begin{proposition}
\label{propinvtwofonc}
Let $\tanghcyll \colon \sourcetl \hookrightarrow \crats(\hcylc)$ denote a long tangle representative in a rational homology cylinder. Let $\tau$ denote a parallelization of $\hcylc$ as in Definition~\ref{defparacyl}.
Let $\omega_0$ and $\omega_1$ be two homogeneous propagating forms of $C_2(\rats(\hcylc))$ (as in Definitions~\ref{defpropagatortwo} and \ref{defprophomogen}). Let $\tilde{\omega}$ be a closed $2$-form on $\left[0,1\right] \times \partial C_2(\rats(\hcylc))$ whose restriction $\tilde{\omega}(t)$ to $\{t\} \times \left(\partial C_2(\rats(\hcylc)) \setminus \ST  {\ballb}_{\rats(\hcylc)}\right)$ is $\projp_{\tau}^{\ast}(\omega_{S^2})$ for any $t \in \left[0,1\right]$, and such that the restriction of $\omega_i$ to $\partial C_2(\rats(\hcylc))$ is $\tilde{\omega}(i)$ for $i\in \{0,1\}$.
For any component $K_j$ of $\tanghcyll=\sqcup_{j=1}^k K_j$, define
$I_j=\sum_{\Gamma_{\finsetb} \in \CD^c(\RR)}\coefgambet_{\Gamma_{\finsetb}}I(\Gamma_{\finsetb},K_j,\tilde{\omega}),$
where \begin{equation*}I(\Gamma_{\finsetb},K_j,\tilde{\omega})=\int_{u\in \left[0,1\right]}\int_{w\in K_j \cap {\ballb}_{\rats(\hcylc)}}\int_{\cinjuptdanvec(T_{w}\crats(\hcylc),\vec{t}_w;{\Gamma}_{\finsetb})} \bigwedge_{e \in E(\Gamma_{\finsetb})}p_e^{\ast}\bigl(\tilde{\omega}(u)\bigr)\left[\Gamma_{\finsetb}\right]\end{equation*}
and $\vec{t}_w$ denotes the unit tangent vector to $K_j$ at $w$. \\(The notation $\cinjuptdanvec(T_{w}\crats(\hcylc),\vec{t}_w;{\Gamma}_{\finsetb})$ is introduced before Lemma~\ref{lemfacfacone}, and $\CD^c(\RR)$ is introduced at the beginning of Section~\ref{secdephom}.)
Define 
\begin{equation*}\zinv(\tilde{\omega})=\sum_{n\in \NN}\zinv_n\left(\left[0,1\right] \times \ST {\ballb}_{\rats(\hcylc)}; \tilde{\omega}\right)\end{equation*}
as in Corollary~\ref{corinvone}.
Then we have
\begin{equation*}\Zinv(\hcylc,\tanghcyll,\omega_1)=\left(\prod_{j=1}^k \exp\left(I_j\right) \#_j\right)\Zinv\left(\hcylc,\tanghcyll,\omega_0\right) \exp\bigl(\zinv(\tilde{\omega})\bigr).\end{equation*}
\end{proposition}
\bp According to Proposition~\ref{propinvtwomieux}, this statement holds when $\tanghcyll$ is a link, and when $\tanghcyll$ is the empty link in particular.
Using Notation~\ref{notationzZ}, it suffices to prove that
\begin{equation*}\Zinvlink(\hcylc,\tanghcyll,\omega_1)=\left(\prod_{j=1}^k \exp\left(I_j\right) \#_j\right)\Zinvlink\left(\hcylc,\tanghcyll,\omega_0\right)\end{equation*} since $\Zinv(\hcylc,\tanghcyll,\omega_1)=\Zinv(\crats(\hcylc),\emptyset,\omega_1)\Zinvlink(\hcylc,\tanghcyll,\omega_1)$ as in Lemma~\ref{lemmultcheck}.
As in the proof of Lemma~\ref{lemvarzinfunc} above, the only faces contributing to the variation of $\Zinvlink(\hcylc,\tanghcyll,\omega_t)$ are the faces $\facee(\check{\Gamma}_{\finseta},\tanghcyll,\Gamma)$ for which $\check{\Gamma}_{\finseta}$ is a connected diagram on $\RR$
and $\Gamma$ is a diagram that contains $\check{\Gamma}_{\finseta}$ as a subdiagram on a component $\source_j$ of $\tanghcyll$. Their contribution yields the result as in the proof of Lemma~\ref{lempropinvtwomieux}.
\eop

\begin{lemma}\label{lemvarithetatangtwo}
 Recall Definition~\ref{defIthetalong} of $I_{\theta}(K,\tau)$ for a long component $K \colon \RR \hookrightarrow \crats(\hcylc)$ of a tangle in a parallelized $\QQ$-cylinder $(\hcylc,\tau)$.
Let $\hcylc$ be a $\QQ$-cylinder. Let $(\tau_t)_{t \in \left[0,1\right]}$ be a smooth homotopy of parallelizations of $\hcylc$.
For any component $K$ of a tangle in $\hcylc$, we have
\begin{equation*}I_{\theta}(K,\tau_u) - I_{\theta}(K,\tau_0) = 2 \int_{\cup_{t \in \left[0,u\right]} p_{\tau_t}(\ST^+K)}\omega_{S^2}.\end{equation*}
\end{lemma}
\bp
When $K$ is closed, $I_{\theta}$ is defined in Lemma~\ref{lemdefItheta}, and the lemma follows from Proposition~\ref{propanom} and Lemma~\ref{lemvartauanom}.
Lemma~\ref{lemformprod} implies the existence of a closed $2$-form $\omega$ on $\left[0,1\right] \times C_2(\rats(\hcylc))$ that restricts to $\{t\} \times C_2(\rats(\hcylc))$ as a homogeneous propagating form of $(C_2(\rats(\hcylc)),\tau_t)$ for all $t \in \left[0,1\right]$. The integral of this form on $\partial \left(\left[0,u\right] \times C(\rats(\hcylc),K;\onechordsmallnonnum)\right)$ is zero, and it is half the difference between the two sides of the equality to be proved when $K$ is a long component. 
\eop

\bpo{Proof of Theorem~\ref{thmfstconsttang}}
Let $\omega$ be a homogeneous propagating form of $(C_2(\rats(\hcylc)),\tau)$. Let us study the variation of
$\Zinvlink(\hcylc,\tanghcyll,\tau)=\left(\Zinvlink_n(\hcylc,\tanghcyll,\omega)\right)_{n \in \NN}$ when $\tau$ varies inside its homotopy class.

Let $(\tau(t))_{t \in \left[0,1\right]}$ be a smooth homotopy of parallelizations of $\hcylc$. Set $\Zinvlink(t)=\Zinvlink(\hcylc,\tanghcyll,\tau(t)).$ 
Thanks to Proposition~\ref{propinvtwofonc}, we have
\begin{equation*}\frac{\partial}{\partial t}\Zinvlink(t)=\sum_{j=1}^k\left(\frac{\partial}{\partial t}\left(2 \int_{\cup_{u \in \left[0,t\right]} p_{\tau(u)}(\ST^+K_j)}\omega\right) \alpha \#_j \right)\Zinvlink(t)\end{equation*}
as in Lemma~\ref{lemvartauanom}. Lemma~\ref{lemvarithetatangtwo} implies
\begin{equation*}I_{\theta}\bigl(K_j,\tau(t)\bigr) - I_{\theta}\bigl(K_j,\tau(0)\bigr) = 2 \int_{\cup_{u \in \left[0,t\right]} p_{\tau(u)}(\ST^+K_j)}\omega\end{equation*} for any $j$.
As in Corollary~\ref{corvartauanom}, conclude that \begin{equation*}\prod_{j=1}^k\biggl(\exp\Bigl(-I_{\theta}\bigl(K_j,\tau(t)\bigr)\alpha\Bigr) \#_j\biggr)\Zinvlink(t)\end{equation*} is constant, and note $\Zinvlink_1(\hcylc,K_j,\tau(t))=\frac12 I_{\theta}\bigl(K_j,\tau(t)\bigr) \left[ \onechordsmalltseul\right]$ for an interval component $K_j$.

Proposition~\ref{propinvtwofonc} and Lemma~\ref{lemformprod}  
imply that changing the trivialization $\tau$ in a ball $B_{\tau}$ that does not meet the tangle does not change $\Zinvlink$ (where the form $\omega^{\partial}$ of Lemma~\ref{lemformprod} is easily assumed to pull back through the projection of $\left[0,1\right] \times (\partial C_2(\rats(\hcylc)) \setminus \ST(\ballb_{\tau}))$ onto $\partial C_2(\rats(\hcylc)) \setminus \ST(\ballb_{\tau})$ on $\left[0,1\right] \times (\partial C_2(\rats(\hcylc)) \setminus \ST(\ballb_{\tau}))$).
Then the proof of Theorem~\ref{thmfstconsttang} can be concluded like the proof of Theorem~\ref{thmfstconst} at the end of Section~\ref{secdephom}, with the following additional argument for the strands going from bottom to bottom or from top to top. In the proof of Theorem~\ref{thmfstconst}, we assumed that $\projp_{\partau}(\ST^+K_j)=v$ for some $v \in S^2$, and that $g$ maps $K_j$ to rotations with axis $v$ for any $j \in \underline{k}$, in order to ensure that $^{\partau}\psi(g^{-1})$ induces a diffeomorphism of $\cup_{w\in K_j}{\cinjuptdanvec(T_{w}\crats,\vec{t}_w;\check{\Gamma}_{\finsetb})}$. Without loss of generality, we instead assume $v = \upvec$ and $\projp_{\partau}(\ST^+K_j)=\pm \upvec$ for all components $K_j$ of $\tanghcyll$, except possibly in a neighborhood of the boundary of $\hcylc$, which is mapped to $1$ by $g$ (so that $\projp_{\partau}(\ST^+K_j)$ can move from $\pm \upvec$ to $\mp \upvec$ in this neighborhood).
\eop

\chapter{The invariant \texorpdfstring{$\Zinv$}{Z} as a holonomy for braids}
\label{chapzinvfbraid}

In this chapter, we interpret the extension of $\Zinvufrf$ to long tangles of the previous chapter as a holonomy for long braids, and we study it as such.

Recall the compactification $\ccompuptd{\finsetv}{\vecspt}$ of the space $\cinjuptd{\finsetv}{\vecspt}$ of injective maps from a finite set $\finsetv$ to a vector space $\vecspt$ up to translation and dilation, from Theorem~\ref{thmcompfacanom}.
Let $\finsetb$ be a finite set of cardinality at least $2$.
Let $\Gamma$ be a Jacobi diagram on a disjoint union of lines $\RR_{\eltb}$ indexed by elements $\eltb$ of $\finsetb$. Let $p_{\finsetb} \colon U(\Gamma)\to \finsetb$ be the natural map induced by $i_{\Gamma}$. We assume that $p_{\finsetb}$ is onto.
Let $\univu_{\eltb}=\univu_{\eltb}(\Gamma)=p_{\finsetb}^{-1}(\eltb)$ be the set of univalent vertices of $\Gamma$ sent to $\RR_b$ by $i_{\Gamma}$.
Let $\check{\CV}(\Gamma) \subset \cinjuptd{\vertsetv(\Gamma)}{\RR^3}$ be the quotient by the translations and the dilations of the space of injective maps $\confc$ from $\vertsetv(\Gamma)$ to $\RR^3=\CC \times \RR$
that map $\univu_{\eltb}(\Gamma)$ to a vertical line $\confy(\confc,\eltb) \times \RR $ for each $\eltb \in \finsetb$, with respect to the order induced by $i_{\Gamma}$, so that 
the planar configuration $\confy(\confc,.) \colon \finsetb \to \CC$ is injective.
Let $\CV(\Gamma)$ \index[N]{Vcal@${\CV}(\Gamma)$ configuration space} denote the closure of the image of $\check{\CV}(\Gamma)$ in $\ccompuptd{\vertsetv(\Gamma)}{\RR^3} \times \ccompuptd{\finsetb}{\CC}$ under the map $\bigl(\confc \mapsto (\confc,\confy(\confc,.))\bigr)$.

\section{On the structure of \texorpdfstring{$\CV(\Gamma)$}{V(Gamma)}}
\label{secstrucCV} 

In this section, we investigate the structure of $\CV(\Gamma)$, as we did in Section~\ref{secstrucbyhand} for ${C}(\rats(\hcylc),\tanghcyll;\Gamma)$.

\begin{lemma}
\label{lemcnsonecybelongstonu}
An element $(\confc,\confy)$ of $\ccompuptd{\vertsetv(\Gamma)}{\RR^3} \times \ccompuptd{\finsetb}{\CC}$ is in $\CV(\Gamma)$ if and only if
\begin{itemize}
 \item for any $ \eltb \in \finsetb$ and for any $ (v_1,v_2) \in \univu_{\eltb}^2$, the restriction of $\confc$ to $\{v_1,v_2\}$
is vertical, and its direction is that prescribed by $i_{\Gamma}$,
\item for any pair $(b_1,b_2)$ of distinct elements of $\finsetb$ and for any $(v_1,v_2)$ in $\univu_{\eltb_1}\times \univu_{\eltb_2}$, there exists $\beta \in \RR^+$ such that the restriction $\confc\vert_{\{v_1,v_2\}}$ of $\confc$ to $\{v_1,v_2\}$ satisfies 
\begin{equation*}p_{\CC} \bigl(\confc\vert_{\{v_1,v_2\}}(v_2)-\confc\vert_{\{v_1,v_2\}}(v_1)\bigr)=\beta \bigl(\confy\vert_{\{b_1,b_2\}}(b_2)-\confy\vert_{\{b_1,b_2\}}(b_1)\bigr)\end{equation*} 
(where $\beta$ depends on chosen representatives of $\confc$ and $\confy$), and
\item for any triple $(b_1,b_2,b_3)$ of distinct elements of $\finsetb$ and for any  $(v_1,v_2,v_3) \in \univu_{\eltb_1}\times \univu_{\eltb_2}\times \univu_{\eltb_3}$, there exists $\beta \in \RR^+$ such that
\begin{multline*}\Bigl(p_{\CC} \bigl( \confc\vert_{\{v_1,v_2,v_3\}}(v_2)-\confc\vert_{\{v_1,v_2,v_3\}}(v_1)\bigr),p_{\CC} \bigl( \confc\vert_{\{v_1,v_2,v_3\}}(v_3)-\confc\vert_{\{v_1,v_2,v_3\}}(v_1)\bigr)\Bigr)\\
=\beta \Bigl(\confy\vert_{\{b_1,b_2,b_3\}}(b_2)-\confy\vert_{\{b_1,b_2,b_3\}}(b_1),\confy\vert_{\{b_1,b_2,b_3\}}(b_3)-\confy\vert_{\{b_1,b_2,b_3\}}(b_1)\Bigr)\end{multline*} in $\CC^2$.
\end{itemize}
\end{lemma}

It is easy to see that an element of $\CV(\Gamma)$ must satisfy these conditions since they are closed and satisfied on $\check{\CV}(\Gamma)$. We will prove the converse after Lemma~\ref{lemstrucCV}.

Let $(\confc^0,\confy^0)$ satisfy the conditions of Lemma~\ref{lemcnsonecybelongstonu}.
We are going to study how a neighborhood $N(\confc^0,\confy^0)$ of $(\confc^0,\confy^0)$ in $\ccompuptd{\vertsetv(\Gamma)}{\RR^3} \times \ccompuptd{\finsetb}{\CC}$ intersects $\check{\CV}(\Gamma)$.

Let us introduce notation in the following long notation paragraph, which ends just before Lemma~\ref{lemcnsonecybelongstonubis}.

\begin{notation}
 \label{notstrucCV}

As in Theorem~\ref{thmcompuptd}, the configuration $\confc^0$ in $\ccompuptd{\vertsetv(\Gamma)}{\RR^3}$ is described by a $\Delta$-parenthesization $\parentp$ of $\vertsetv=\vertsetv(\Gamma)$ (as in Definition~\ref{defparenthesizationtwo})
and \begin{equation*}\left(c^0_{\finsetz} \in \cinjuptd{\kids(\finsetz)}{\RR^3}\right)_{\finsetz \in \parentp}.\end{equation*} The configurations $\confc$ in a neighborhood of $c^0$ in $\ccompuptd{\vertsetv}{\RR^3}$ may be expressed as
\begin{equation*}\begin{array}{ll}\confc\bigl((\mu_{\finsetz}), (\confc_{\finsetz})\bigr)&=\sum_{\finsetz \in \parentp}\left(\prod_{\finsety \in \parentp \suchthat  \finsetz \subseteq \finsety \subset \finsetv} \mu_{\finsety} \right) \confc_{\finsetz}\\
   &=\confc_{\finsetv}+ \sum_{\finsetz \in D(\finsetv)}\mu_{\finsetz}\bigl(\confc_{\finsetz} + \sum_{\finsety \in D(\finsetz)}\mu_{\finsety}\left(\confc_{\finsety}+\dots\right)\bigr)
  \end{array}
\end{equation*}
for $\left((\mu_{\finsetz})_{{\finsetz} \in \parentp \setminus \{\finsetv\}}, (\confc_{\finsetz})_{{\finsetz} \in \parentp}\right) \in \left[0,\varepsilon\right[^{\parentp \setminus \{\finsetv\}} \times \prod_{{\finsetz} \in \parentp} W_{\finsetz}$, as in Lemma~\ref{lemneighchartcomf}.

Similarly, the configuration $\confy^0$ in $\ccompuptd{\finsetb}{\CC}$ may be written as
\begin{equation*}\bigl(\confy^0_{\finsetd} \in \cinjuptd{\kids(\finsetd)}{\CC}\bigr)_{\finsetd \in \parentp_{\finsetb}}\end{equation*}
for a $\Delta$-parenthesization $\parentp_{\finsetb}$ of $\finsetb$, and the configurations $\confy$ in a neighborhood of $\confy^0$ in $\ccompuptd{\finsetb}{\CC}$ may be expressed as
\begin{equation*}\begin{array}{ll}\confy\bigl((u_{\finsetd}), (\confy_{\finsetd})\bigr)&=\sum_{\finsetd \in \parentp_{\finsetb}}\left(\prod_{\finsete \in \parentp_{\finsetb} \suchthat  \finsetd \subseteq \finsete \subset \finsetb} u_{\finsete} \right) \confy_{\finsetd}\\
   &=\confy_{\finsetb}+ \sum_{\finsetd \in D(\finsetb)}u_{\finsetd}\bigl(\confy_{\finsetd} + \sum_{\finsete \in D(\finsetd)}u_{\finsete}(\confy_{\finsete}+\dots)\bigr)
  \end{array}
\end{equation*}
for $((u_{\finsetd})_{{\finsetd} \in \parentp_{\finsetb} \setminus \{\finsetb\}}, (\confy_{\finsetd})_{{\finsetd} \in \parentp_{\finsetb}}) \in \left[0,\varepsilon\right[^{\parentp_{\finsetb} \setminus \{\finsetb\}} \times \prod_{{\finsetd} \in \parentp_{\finsetb}} N_{\finsetd}$.

We normalize the $\confy_{\finsetd}$ by choosing basepoints $b(\finsetd)$ for the $\finsetd \in \parentp_{B}$ and imposing $\confy_{\finsetd}(b(\finsetd))=0$ and 
$\sum_{\finsete \in \kids(\finsetd)}\norm{\confy_{\finsetd}(\finsete)}^2=1$.

For a set $\finsety$ of $\parentp$,  we choose a basepoint $b(\finsety)$,
which is univalent if there is a univalent vertex in $\finsety$.
We also choose a kid $k_n(\finsety)$ such that 
$|p_{\RR}(\confc^0_{\finsety}(k_n(\finsety))) -p_{\RR}(\confc^0_{\finsety}(b(\finsety)))|$ or $|p_{\CC}(\confc^0_{\finsety}(k_n(\finsety))) -p_{\CC}(\confc^0_{\finsety}(b(\finsety)))|$ is maximal in the set 
\begin{equation*}\biggl\{\Bigl|p_{\RR}\Bigl(\confc^0_{\finsety}(k)\Bigr) -p_{\RR}\Bigl(\confc^0_{\finsety}\bigl(b(\finsety)\bigr)\Bigr)\Bigr|, \Bigl|p_{\CC}\Bigl(\confc^0_{\finsety}(k)\Bigr) -p_{\CC}\Bigl(\confc^0_{\finsety}\bigl(b(\finsety)\bigr)\Bigr)\Bigr| \suchthat k \in \kids(\finsety) \biggr\}\end{equation*}
and we call it the \emph{normalizing kid} of $\finsety$.
If 
\begin{equation*}\Bigl|p_{\RR}\Bigl(\confc^0_{\finsety}(k_n(\finsety))\Bigr) -p_{\RR}\Bigl(\confc^0_{\finsety}\bigl(b(\finsety)\bigr)\Bigr)\Bigr| \geq \Bigl|p_{\CC}\Bigl(\confc^0_{\finsety}(k_n(\finsety))\Bigr) -p_{\RR}\Bigl(\confc^0_{\finsety}\bigl(b(\finsety)\bigr)\Bigr)\Bigr|,\end{equation*} then we say that $k_n(\finsety)$ is \indexT{vertically normalizing} or \emph{v-normalizing}, and we normalize the configurations $\confc_{\finsety}$ in a neighborhood of $\confc^0_{\finsety}$ by imposing $\confc_{\finsety}(b(\finsety))=0$ and $|p_{\RR}(\confc_{\finsety}(k_n(\finsety)))|=1$.
Otherwise, we say that $k_n(\finsety)$ is \indexT{horizontally normalizing} or \emph{h-normalizing}, and we first normalize the configurations $\confc_{\finsety}$ in a neighborhood of $\confc^0_{\finsety}$ by imposing \begin{equation*}\confc_{\finsety}\bigl(b(\finsety)\bigr)=0\mbox{ and }\Bigl|p_{\CC}\Bigl(\confc_{\finsety}\bigl(k_n(\finsety)\bigr)\Bigr)\Bigr|=1. \end{equation*} (These normalizations are compatible with the smooth structure of $\ccompuptd{\vertsetv}{\RR^3}$. In the case of a horizontally normalizing kid, they will be changed in Notation~\ref{notationnewnormhn}.)

In our neighborhood, we also impose $\norm{\confc_{\finsety}(k)-\confc^0_{\finsety}(k)}<\varepsilon$ for any kid $k$ of $\finsety$, for a small $\varepsilon \in \left]0,1\right[$. So the 
manifold $W_{\finsety}$ is diffeomorphic to the product of the product, over the nonnormalizing kids $k$ of $\finsety$ that do not contain $b(\finsety)$, of balls $\mathring{\ballb}(\confc^0_{\finsety}(k),\varepsilon)$, by the set of elements $\confc_{\finsety}(k_n(\finsety))$
of $\mathring{\ballb}(\confc^0_{\finsety}(k_n(\finsety)),\varepsilon)$, such that $|p_{\RR}(\confc_{\finsety}(k_n(\finsety)))|=1$ (resp. such that $|p_{\CC}(\confc_{\finsety}(k_n(\finsety)))|=1$) if $k$ is v-normalizing (resp. if $k$ is h-normalizing). Note that
$|p_{\CC}(\confc_{\finsety}(k))|<2$ for any $k \in  \kids(\finsety)$.

For a set $\finsety$ of $\parentp$, let $\finsetb(\finsety) \subseteq \finsetb$ be the set of (labels of) the components of its univalent vertices, and let $\hat{\finsetb}(\finsety)$ denote the smallest element of $\parentp_{\finsetb}$ such that $\finsetb(\finsety) \subseteq \hat{\finsetb}(\finsety)$ if $\cardlef{\finsetb(\finsety)} \geq 2$. If $\cardlef{\finsetb(\finsety)}=1$, set $\hat{\finsetb}(\finsety)=\finsetb(\finsety)$.
 If $\finsetb(\finsety) \neq \emptyset$, the set $\finsety$ is called \emph{univalent}.
 
For $\finsetd \in \parentp_{\finsetb}$, define the set $\parentp_{\cutx,\finsetd}$\index[N]{Parenthesizations!papacutxX@$\parentp_{\cutx,\finsetd}$} of elements
$\finsety$ of $\parentp$ such that $\hat{\finsetb}(\finsety)=\finsetd$ and $\hat{\finsetb}(\finsetz)\neq\finsetd$ for every daughter $\finsetz$ of $\finsety$. Note that any element $\finsety$ of $\parentp_{\cutx,\finsetd}$ has at least two kids $\finsety_a$ and $\finsety_b$ such that $\finsetb(\finsety_a) \neq \finsetd$ and $\finsetb(\finsety_b) \neq \finsetd$.
Set
\begin{equation*}\parentp^{\prime}_{\cutx}\index[N]{Parenthesizations!papacutX@$\parentp^{\prime}_{\cutx}$}= \cup_{\finsetd \in \parentp_{\finsetb}} \parentp_{\cutx,\finsetd}.\end{equation*}

Let $\widehat{\parentp^{\prime}}_{\cutx}$\index[N]{Parenthesizations!papacutXp@$\widehat{\parentp^{\prime}}_{\cutx}$} be the subset of $\parentp$ consisting of the univalent sets $\finsetz$ of $\parentp$ such that $\cardlef{\finsetb(\finsetz)} \geq 2$.\end{notation}

\begin{lemma}
\label{lemcnsonecybelongstonubis}
Let $(\confc^0,\confy^0)$ be a configuration parametrized as above, which satisfies the conditions of Lemma~\ref{lemcnsonecybelongstonu}.

For any $\finsety \in \widehat{\parentp^{\prime}}_{\cutx}$, there exists a unique $\lambda^0(\finsety) \in \RR^+$ 
such that for any two univalent kids $\finsety_1$ and $\finsety_2$ of $\finsety$ such that $\finsetb(\finsety_1) \subset \finsetd_1$ and $\finsetb(\finsety_2) \subset \finsetd_2$, where $\finsetd_1$ and $\finsetd_2$ are kids of $\hat{\finsetb}(\finsety)$, we have
\begin{equation*}p_{\CC}\bigl(\confc^0_{\finsety}(\finsety_2)-\confc^0_{\finsety}(\finsety_1)\bigr)=\lambda^0(\finsety)\bigl(\confy^0_{\hat{\finsetb}(\finsety)}(\finsetd_2)-\confy^0_{\hat{\finsetb}(\finsety)}(\finsetd_1)\bigr).\footnote{$\lambda^0(\finsety)$ depends on the fixed normalizations of $\confc^0_{\finsety}$ and  $\confy^0_{\hat{\finsetb}(\finsety)}$.}\end{equation*}
For any $\finsety \in \widehat{\parentp^{\prime}}_{\cutx} \setminus \parentp^{\prime}_{\cutx}$, we have $\lambda^0(\finsety)=0$.
Conversely, if an element $(\confc^0,\confy^0)$ satisfies the above properties and the first condition of Lemma~\ref{lemcnsonecybelongstonu}, then it satisfies the conditions of Lemma~\ref{lemcnsonecybelongstonu}.
\end{lemma}
\bp Let $ \finsety \in \parentp$. Assume $\cardlef{\finsetb(\finsety)} \geq 2$.
The first condition of Lemma~\ref{lemcnsonecybelongstonu} ensures that when $\eltv$ and $\eltv^{\prime}$ belong to $\finsety \cap U_b$, we have \begin{equation*}p_{\CC}\circ\confc^0_{\finsety}(\eltv)=p_{\CC}\circ\confc^0_{\finsety}(\eltv^{\prime})=(p_{\CC}\circ\confc^0_{\finsety})(b \in \finsetb).\end{equation*}
So, we view $p_{\CC}\circ\confc^0_{\finsety}$ as a map from $\finsetb(\finsety)$ to $\CC$.

Assume $\finsety \in \widehat{\parentp^{\prime}}_{\cutx}$. Then $\confy^0_{\hat{\finsetb}(\finsety)}\vert_{{\finsetb}(\finsety)}$ is not constant. So, the second and third conditions ensure that there exists a unique $\lambda^0(\finsety) \in \RR^+$ 
such that
$p_{\CC}\circ\confc^0_{\finsety}$ coincides with $\lambda^0(\finsety)\confy^0_{\hat{\finsetb}(\finsety)|{\finsetb}(\finsety)}$ up to translation.
Since $p_{\CC}\circ\confc^0_{\finsety}$ is constant on ${\finsetb}(\finsetz)$ for any kid $\finsetz$ of $\finsety$, if there exists such a kid $\finsetz$ such that $\hat{\finsetb}(\finsetz) = \hat{\finsetb}(\finsety)$, then $\lambda^0(\finsety)$ must be equal to zero.

The last assertion is an easy exercise.
\eop

In order to finish the proof of Lemma~\ref{lemcnsonecybelongstonu}, we are going to prove that the configurations that satisfy the conditions of its statement are in ${\CV}(\Gamma)$. We take a closer look at the structure of ${\CV}(\Gamma)$.

\begin{notation} \label{notationnewnormhn}
For a univalent $\finsety \in \parentp$, define $\deltb(\finsety) \in \finsetb$ such that $b(\finsety) \in U_{\deltb(\finsety)}$.
Also assume that $p_{\finsetb}\left(b\left(\vertsetv\right)\right)$ is the basepoint $b_0$ of $\finsetb$. As always, our basepoints satisfy that if $\finsetz \subset \finsety$ and if $b(\finsety) \in \finsetz$, then $b(\finsetz)= b(\finsety)$.
For $\finsetd \in \parentp_{\finsetb}$,
set $U_{\finsetd}=\prod_{\finsete \in \parentp_{\finsetb} \setminus \{\finsetb\} \suchthat \finsetd \subseteq \finsete}u_{\finsete}$. For $\deltb \in \finsetd$, set
\begin{equation*}\tilde{\confy}_{\finsetd}(\deltb) = \confy_{\finsetd}(\deltb) + \sum_{\finsete \in \parentp_{\finsetb} \suchthat  \deltb \in \finsete \subset \finsetd} \left(\prod_{\finsetf \in \parentp_{\finsetb} \suchthat  \finsete \subseteq \finsetf \subset \finsetd} u_{\finsetf}\right)\confy_{\finsete}(\deltb).\end{equation*}
For $\finsety \in \parentp_{\cutx}$,
set $M_{\finsety}=\prod_{\finsetz \in \parentp \setminus \{\finsetv\} \suchthat \finsety \subseteq \finsetz}\mu_{\finsetz}$.

Let $\parentp^{\prime}_{\cutx,hn}$ be the set of elements $\finsetz$ of $\widehat{\parentp^{\prime}}_{\cutx}$ such that $k_n(\finsetz)$ is a univalent horizontally normalizing kid. For an element  $\finsetz$
of $\parentp^{\prime}_{\cutx,hn}$, we have $\lambda^0(\finsetz) \neq 0$, which implies $\parentp^{\prime}_{\cutx,hn} \subseteq \parentp^{\prime}_{\cutx}$.
We change the normalizations of the $c_{\finsetz}$ for the elements $\finsetz$ of $\parentp^{\prime}_{\cutx,hn}$ in our neighborhood $N(\confc^0,\confy^0)$ of $(\confc^0,\confy^0)$ in $\ccompuptd{\vertsetv(\Gamma)}{\RR^3} \times \ccompuptd{\finsetb}{\CC}$ (which we reduce if needed) in $\ccompuptd{\finsetb}{\CC}$ so that 
 $\confc_{\finsetz}(b(\finsetz))=0$ and
 \begin{equation}
  p_{\CC}\Bigl(\confc_{\finsetz}\bigl(k_n(\finsetz)\bigr)\Bigr)=\lambda^0(\finsetz)\biggl(\tilde{\confy}_{\hat{\finsetb}(\finsetz)}\Bigl(d\bigl(k_n(\finsetz)\bigr)\Bigr)-\tilde{\confy}_{\hat{\finsetb}(\finsetz)}\Bigl(d(\finsetz)\Bigr)\biggr).
 \end{equation}
With this normalization, Lemma~\ref{lemcnsonecybelongstonubis} is still valid, and we have $\lambda(\finsetz)=\lambda^0(\finsetz)$ in our neighborhood. 

\end{notation}

\begin{lemma}
\label{lemdeflambdabraids}
 With the above normalizations and notation, there exist continuous maps $\lambda$ and $\lambda(Y)$, for $\finsety \in \parentp^{\prime}_{\cutx}$, from $N(\confc^0,\confy^0)\cap {\CV}(\Gamma)$ to $\RR^+$ satisfying the following properties.
For any $\finsety \in \parentp^{\prime}_{\cutx}$, we have
 $\lambda(\finsety)(\confc^0)=\lambda^0(\finsety)$ with the notation of Lemma~\ref{lemcnsonecybelongstonubis} and,
 for any configuration $(\confc,\confy)$ of $N(\confc^0,\confy^0)\cap {\CV}(\Gamma)$,
 \begin{itemize}
  \item we have $p_{\CC} \circ \confc(v)= \lambda \confy(\deltb)$ for any $\deltb \in \finsetb$ and for any $v \in U_{\deltb}$,
  \item for $\finsetd \in  \parentp_{\finsetb}$ and $\finsety \in \parentp_{\cutx,\finsetd}$,
  \begin{itemize}
  \item if $\finsety_a$ and $\finsety_b$ are two univalent kids of $\finsety$, then we have
  \begin{equation*}p_{\CC}\bigl( \confc_{\finsety}(\finsety_b)- \confc_{\finsety}(\finsety_a)\bigr)
=\lambda(\finsety) \Bigl(\tilde{\confy}_{\finsetd}\bigl(\deltb(\finsety_b)\bigr)-\tilde{\confy}_{\finsetd}\bigl(\deltb(\finsety_a)\bigr)\Bigr),\end{equation*}
 \item and the equality \begin{equation*}\ast(\finsety): \;\;\lambda(\finsety)M_{\finsety}=\lambda U_{\hat{\finsetb}(\finsety)}\end{equation*}
  holds.
 \end{itemize}
 \end{itemize}
 When $\finsetv \in \parentp^{\prime}_{\cutx}$, $\ast(\finsetv)$ is equivalent to $\lambda(\finsetv)=\lambda$.
 When $\finsetv \notin \parentp^{\prime}_{\cutx}$,
 we also let $\lambda(\finsetv)$ and $\lambda$ both denote $\lambda$, depending on the context.
\end{lemma}
\bp
In $N(\confc^0,\confy^0)\cap\check{\CV}(\Gamma)$, where $\prod_{\finsetd \in \parentp_{\finsetb} \setminus \{\finsetb\}}u_{\finsetd} \times \prod_{\finsetw \in \parentp \setminus \{\finsetv\}}\mu_{\finsetw } \neq 0$,
there exists some $\lambda >0$ such that $p_{\CC} \circ \confc(v)= \lambda \confy(b)$ for any $b \in U_b$ and for any $v \in U_b$.
This map $\lambda$, starting from $N(\confc^0,\confy^0)\cap\check{\CV}(\Gamma)$, can be extended continuously on $N(\confc^0,\confy^0)\cap {\CV}(\Gamma)$, as follows.
Let $b \in \finsetb$ be such that $b$ and $b_0$ belong to different kids of $\finsetb$, and let $v \in U_{b
}$. Then we have $p_{\CC} \left( \confc(v)\right)= \lambda\confy(b)$
on $N(\confc^0,\confy^0)\cap\check{\CV}(\Gamma)$, and $\confy(b)$ does not vanish on 
 $N(\confc^0,\confy^0)$. The closed condition that $p_{\CC} \left( \confc(v)\right)$ and $\confy(b
)$ are colinear and that their scalar product is nonnegative is satisfied on $N(\confc^0,\confy^0)\cap {\CV}(\Gamma)$. It allows us to define $\lambda(v)$ such that 
 $p_{\CC}\left(\confc(v)\right)= \lambda(v)\confy(b)$ on $N(\confc^0,\confy^0)\cap {\CV}(\Gamma)$
 and  $\lambda(v)$ is continuous. Now, since $\lambda(v)=\lambda$ is independent of $v$ as above on $N(\confc^0,\confy^0)\cap\check{\CV}(\Gamma)$, it is also on $N(\confc^0,\confy^0)\cap {\CV}(\Gamma)$. 
 Set $\lambda=\lambda(v)$.  Then we have $p_{\CC} \circ \confc(v)= \lambda \confy(b)$ for any $b \in U_b$ and for any $v \in U_b$ on $N(\confc^0,\confy^0)\cap\check{\CV}(\Gamma)$. This is also true on $N(\confc^0,\confy^0)\cap {\CV}(\Gamma)$. 
 
 Let $\finsety \in \parentp$. Let $\finsety_a$ and $\finsety_b$ be two univalent kids of $\finsety$. If $\confc \in N(\confc^0,\confy^0)\cap\check{\CV}(\Gamma)$, then we have
 \begin{equation*}M_{\finsety}\Bigl(p_{\CC}\bigl( \confc_{\finsety}(\finsety_b)- \confc_{\finsety}(\finsety_a)\bigr)\Bigr)
=\lambda U_{\hat{\finsetb}(\finsety)} \Bigl(\tilde{\confy}_{\hat{\finsetb}(\finsety)}\bigl(\deltb(\finsety_b)\bigr)-\tilde{\confy}_{\hat{\finsetb}(\finsety)}\bigl(\deltb(\finsety_a)\bigr)\Bigr).\end{equation*}

In particular, $p_{\CC}\left( \confc_{\finsety}(\finsety_b)- \confc_{\finsety}(\finsety_a)\right)$ and $\bigl(\tilde{\confy}_{\hat{\finsetb}(\finsety)}(\deltb(\finsety_b))-\tilde{\confy}_{\hat{\finsetb}(\finsety)}(\deltb(\finsety_a))\bigr)$ are colinear, and their scalar product is nonnegative on $N(\confc^0,\confy^0)\cap\check{\CV}(\Gamma)$. As above, as soon as there exist two kids $\finsety_a$ and $\finsety_b$ as above, such that $\deltb(\finsety_a)$ and $\deltb(\finsety_b)$ are in two distinct kids of $\hat{\finsetb}(\finsety)$, we can define the continuous function $\lambda(\finsety)$ such that
\begin{equation*}p_{\CC}\bigl( \confc_{\finsety}(\finsety_b)- \confc_{\finsety}(\finsety_a)\bigr)
=\lambda(\finsety) \Bigl(\tilde{\confy}_{\hat{\finsetb}(\finsety)}\bigl(\deltb(\finsety_b)\bigr)-\tilde{\confy}_{\hat{\finsetb}(\finsety)}\bigl(\deltb(\finsety_a)\bigr)\Bigr)\end{equation*}
for any two univalent kids
$\finsety_a$ and $\finsety_b$
of $\finsety$, and we have $\lambda(\finsety)M_{\finsety}=\lambda U_{\hat{\finsetb}(\finsety)}$.
\eop

Let $\parentp^{\prime}_{\finsetb}$ be the set of elements $\finsetd$ of $\parentp_{\finsetb}$ such that $\parentp_{\cutx,\finsetd} \neq \emptyset$. Let $\widehat{\parentp^{\prime}}_{\finsetb}$ be the set of elements $\finsetd$ of $\parentp_{\finsetb}$ that contain, or are equal to, an element of $\parentp^{\prime}_{\finsetb}$.

For any collections $(\finsetb_i)_{i \in \ZZ/n\ZZ}$, $(\finsetb^{+}_i)_{i \in \ZZ/n\ZZ}$, $(\finsetb^{\prime}_i)_{i \in \ZZ/n\ZZ}$ of (possibly equal) sets of $\widehat{\parentp^{\prime}}_{\finsetb}$ such that \begin{equation*}\left(\finsetb^{+}_{i-1} \cup \finsetb_i\right) \subseteq \finsetb^{\prime}_i,\end{equation*} for any collections $(\finsety_i)_{i \in \ZZ/n\ZZ}$, $(\finsety^{+}_i)_{i \in \ZZ/n\ZZ}$, $(\finsety^{\prime}_i)_{i \in \ZZ/n\ZZ}$ such that
\begin{equation*}\finsety_i \in \parentp_{\cutx,\finsetb_i}\;\mbox{,}\;\finsety^+_i \in \parentp_{\cutx,\finsetb^+_i}\;\mbox{,}\;\finsety^{\prime}_i \in \widehat{\parentp^{\prime}}_{\cutx}\;\mbox{, and}\; \left(\finsety^{+}_{i} \cup \finsety_i\right) \subseteq \finsety^{\prime}_i,\end{equation*}
we have
\begin{equation*}\prod_{i=1}^n \left(\lambda(\finsety_i)\frac{M_{\finsety_i}}{M_{\finsety^{\prime}_i}}\right) \prod_{i=1}^n \frac{U_{\finsetb^+_i}}{U_{\finsetb^{\prime}_{i+1}}} =\prod_{i=1}^n \frac{U_{\finsetb_i}}{U_{\finsetb^{\prime}_i}} \prod_{i=1}^n\left(\lambda(\finsety^+_i)\frac{M_{\finsety^+_i}}{M_{\finsety^{\prime}_i}}\right)\end{equation*}
for configurations of $N(\confc^0,\confy^0)\cap\check{\CV}(\Gamma)$. This equation is equivalent to the following equation, which also holds in $N(\confc^0,\confy^0)\cap {\CV}(\Gamma)$:
\begin{multline*}
\prod_{i=1}^n \left(\lambda(\finsety_i)\prod_{\finsetz \in \parentp \suchthat \finsety_i\subseteq \finsetz \subset \finsety^{\prime}_i}\mu_{\finsetz}\right)
\prod_{i=1}^n \left(\prod_{\finsete \in \parentp_{\finsetb} \suchthat \finsetb^+_i\subseteq \finsete \subset \finsetb^{\prime}_{i+1}}u_{\finsete}\right)\\
=
\prod_{i=1}^n \left(\lambda(\finsety^+_i)\prod_{\finsetz \in \parentp \suchthat \finsety^+_i\subseteq \finsetz \subset \finsety^{\prime}_i}\mu_{\finsetz}\right)
\prod_{i=1}^n \left(\prod_{\finsete \in \parentp_{\finsetb} \suchthat \finsetb_i\subseteq \finsete \subset \finsetb^{\prime}_i}u_{\finsete}\right).
\end{multline*}

 For $\finsety \in \parentp$, let $\kids_{u}(\finsety)$ denote the set of univalent kids of $\finsety$ that do not contain $b(\finsety)$, and let $\kids_{t}(\finsety)$ denote the set of nonunivalent kids of $\finsety$ that do not contain $b(\finsety)$.

We are going to prove the following lemma.

\begin{lemma}
\label{lemstrucCV} 
Let $(\confc^0,\confy^0)$ satisfy the conditions of Lemma~\ref{lemcnsonecybelongstonu} (and hence of Lemma~\ref{lemcnsonecybelongstonubis}).
There is a neighborhood of $(\confc^0,\confy^0)$ in ${\CV}(\Gamma)$ parametrized by the following variables:
\begin{enumerate}
 \item $\bigl((u_{\finsetd})_{{\finsetd} \in \parentp_{\finsetb} \setminus \{\finsetb\}}, (\confy_{\finsetd})_{{\finsetd} \in \parentp_{\finsetb}}\bigr) \in \left[0,\varepsilon\right[^{\parentp_{\finsetb} \setminus \{\finsetb\}} \times \prod_{{\finsetd} \in \parentp_{\finsetb}} N_{\finsetd}$,
 \item $(\mu_{\finsetz})_{{\finsetz} \in \parentp \setminus \{\finsetv\}}\in \left[0,\varepsilon\right[^{\parentp \setminus \{\finsetv\}}$,
 \item the $\confc_{\finsetz}$ for nonunivalent sets $\finsetz$ of $\parentp$,
 \item for univalent sets $\finsetz$ of $\parentp$, the $\confc_{\finsetz}(\finsety)$ for the kids $\finsety \in \kids_t(\finsetz)$ of $\finsetz$, and the $p_{\RR}\circ\confc_{\finsetz}(\finsety)$ for the kids $\finsety \in \kids_u(\finsetz)$ of $\finsetz$, 
 \item the parameter $\lambda$, and 
 \item for every element $\finsety$ of $\parentp^{\prime}_{\cutx}$, the parameter $\lambda(\finsety)$ defined in Lemma~\ref{lemdeflambdabraids}.
\end{enumerate}
These variables satisfy the following constraints and determine the $p_{\CC}\circ\confc_{\finsetz}(\finsety)$ for the univalent kids $\finsety \in \kids_u(\finsetz)$ of a set $\finsetz$ of $\parentp$ as follows.
\begin{enumerate}
\item For a vertically normalizing kid $\finsety$ of an element $\finsetz$ of $\parentp$, we have $p_{\RR}\circ\confc_{\finsetz}(\finsety)=p_{\RR}\circ\confc^0_{\finsetz}(\finsety)=\pm 1$.
\item 
For a horizontally normalizing nonunivalent kid $\finsety$ of an element $\finsetz$ of $\parentp$, we have $|p_{\CC}\circ\confc_{\finsetz}(\finsety)|=1$.
\item For a horizontally normalizing univalent kid $\finsety$ of an element $\finsetz$ of $\parentp$, we have \begin{equation*}
  p_{\CC}\Bigl(\confc_{\finsetz}\bigl(\finsety\bigr)\Bigr)=\lambda^0(\finsetz)\Bigl(\tilde{\confy}_{\hat{\finsetb}(\finsetz)}\bigl(d(\finsety)\bigr)-\tilde{\confy}_{\hat{\finsetb}(\finsetz)}\bigl(d(\finsetz)\bigr)\Bigr).
 \end{equation*}
\item  For any element $\finsety$ of $\parentp_{\cutx,\finsetd}$, we have 
\begin{equation*}\ast(\finsety): \;\;M_{\finsety}\lambda(\finsety)=\lambda U_{\finsetd}\end{equation*}
 where
 $M_{\finsety}=\prod_{\finsetz \in \parentp \setminus \{\finsetv\} \suchthat \finsety \subseteq \finsetz}\mu_{\finsetz}$
and $U_{\finsetd}=\prod_{\finsete \in \parentp_{\finsetb} \setminus \{\finsetb\} \suchthat \finsetd \subseteq \finsete}u_{\finsete}$.
\item For any collections $(\finsetb_i)_{i \in \ZZ/n\ZZ}$, $(\finsetb^{+}_i)_{i \in \ZZ/n\ZZ}$, $(\finsetb^{\prime}_i)_{i \in \ZZ/n\ZZ}$ of (possibly equal) sets of $\widehat{\parentp^{\prime}}_{\finsetb}$ such that \begin{equation*}\left(\finsetb^{+}_{i-1} \cup \finsetb_i\right) \subseteq \finsetb^{\prime}_i,\end{equation*} for any collections $(\finsety_i)_{i \in \ZZ/n\ZZ}$, $(\finsety^{+}_i)_{i \in \ZZ/n\ZZ}$, $(\finsety^{\prime}_i)_{i \in \ZZ/n\ZZ}$ such that
\begin{equation*}\finsety_i \in \parentp_{\cutx,\finsetb_i}\;\mbox{,}\;\finsety^+_i \in \parentp_{\cutx,\finsetb^+_i}\;\mbox{,}\;\finsety^{\prime}_i \in \widehat{\parentp^{\prime}}_{\cutx}\;  \mbox{, and}\; \left(\finsety^{+}_{i} \cup \finsety_i\right) \subseteq \finsety^{\prime}_i,\end{equation*}
we have
\begin{multline}
\label{eqcycconst}
\prod_{i=1}^n \left(\lambda(\finsety_i)\prod_{\finsetz \in \parentp \suchthat \finsety_i\subseteq \finsetz \subset \finsety^{\prime}_i}\mu_{\finsetz}\right)
\prod_{i=1}^n \left(\prod_{\finsete \in \parentp_{\finsetb} \suchthat \finsetb^+_i\subseteq \finsete \subset \finsetb^{\prime}_{i+1}}u_{\finsete}\right)\\
=
\prod_{i=1}^n \left(\lambda(\finsety^+_i)\prod_{\finsetz \in \parentp \suchthat \finsety^+_i\subseteq \finsetz \subset \finsety^{\prime}_i}\mu_{\finsetz}\right)
\prod_{i=1}^n \left(\prod_{\finsete \in \parentp_{\finsetb} \suchthat \finsetb_i\subseteq \finsete \subset \finsetb^{\prime}_i}u_{\finsete}\right).
\end{multline}
\item For any univalent kid $\finsetz_b$ of an element $\finsetz$ of $\parentp$ such that $\deltb(\finsetz_b)=\deltb(\finsetz)$, we have
 \begin{equation*}p_{\CC}\circ \confc_{\finsetz}(\finsetz_b)=p_{\CC}\circ \confc_{\finsetz}\bigl(b(\finsetz)\bigr)=0.\end{equation*} 
 \item For any univalent element $\finsetz$ of $\parentp$ such that $\cardlef{\finsetb(\finsetz)} \geq 2$, for any  maximal element $\finsety$ of $\parentp^{\prime}_{\cutx}$ such that $\finsety \subseteq \finsetz$ and $\hat{\finsetb}(\finsetz)=\hat{\finsetb}(\finsety)$,\footnote{Note that the minimal univalent elements $\finsetz$ of $\parentp$ such that $\cardlef{\finsetb(\finsetz)} \geq 2$ are in $\parentp^{\prime}_{\cutx}$.} for any univalent kid $\finsetz_b$ of $\finsetz$, $p_{\CC}\circ \confc_{\finsetz}(\finsetz_b)$ is equal to
\begin{equation}
\label{eqlambdafinsety}
p_{\CC}\circ \confc_{\finsetz}(\finsetz_b)=\frac{M_{\finsety}}{M_{\finsetz}}\lambda(\finsety) \Bigl(\tilde{\confy}_{\hat{\finsetb}(\finsetz)}\bigl(\deltb(\finsetz_b)\bigr)-\tilde{\confy}_{\hat{\finsetb}(\finsetz)}\bigl(\deltb(\finsetz)\bigr)\Bigr).
\end{equation}
\end{enumerate}
\end{lemma}

The analysis performed before Lemma~\ref{lemstrucCV} ensures that all the elements in our neighborhood
$N(\confc^0,\confy^0)$ of $(\confc^0,\confy^0)$ in $\ccompuptd{\vertsetv}{\RR^3} \times \ccompuptd{\finsetb}{\CC}$ that are in ${\CV}(\Gamma)$ are described by the parameters described in the statement and that they satisfy all the equations of the statement. (See the proof of Lemma~\ref{lemdeflambdabraids} for the last one.)
In order to prove Lemma~\ref{lemstrucCV}, we are going to prove that, conversely, the elements described in its statement are in  ${\CV}(\Gamma)$. Note that the elements described in this statement such that neither $\lambda$, nor the $\mu_Z$ nor the $u_D$ vanish, correspond to elements of $\cinjuptd{\vertsetv}{\RR^3} \times \cinjuptd{\finsetb}{\CC}$ and are in $\check{\CV}(\Gamma)$.

We define a tree $\CT(\parentp,\parentp_{\finsetb})$ with oriented edges from the tree whose vertices are the elements of $\widehat{\parentp^{\prime}}_{\cutx}$ and whose edges start at an element $\finsetz$ of $\widehat{\parentp^{\prime}}_{\cutx} \setminus \{\finsetv\}$, end at its mother $\mother(\finsetz)$, and are labeled by $\mu_{\finsetz}$,
\begin{itemize}
\item by gluing an edge, which arrives at the vertex $\finsety$, starts at a vertex labeled by $O_{\finsety}$, and which is labeled by $\lambda(\finsety)$, at the vertex $\finsety$, for any $\finsety$ of $\parentp^{\prime}_{\cutx}$,
\item by gluing an edge, which arrives at the vertex $\finsetv$, starts at a vertex labeled by $\finsetb$, and is labeled by $\lambda$, at the vertex
$\finsetv$,
\item by gluing to the root $\finsetb$ of this edge the subtree of $\parentp_{\finsetb}$ whose vertices are the elements of $\widehat{\parentp^{\prime}}_{\finsetb}$ and whose edges are labeled by the $u_{\finsetd}$ for $\finsetd \in \widehat{\parentp^{\prime}}_{\finsetb} \setminus \{\finsetb\}$, where the edge labeled by $u_{\finsetd}$ starts at $\finsetd$ and ends at its mother $\mother(\finsetd)$.
\end{itemize}

An example is drawn in Figure~\ref{figTbraidsgen}.
\bfig
\centering
\begin{equation*}
\begin{tikzpicture}[xscale=.9]
\draw (1.4,7) -- (4.2,4.2) -- (11.4,4.2) 
(1.4,4.2) -- (4.2,4.2) -- (2.8,2.8) -- (.1,2.8) 
(2.8,2.8) -- (2.8,0) -- (1.4,0) 
(1.4,1.4) -- (2.8,1.4) 
(6,4.2) -- (6,2.8) -- (5.6,1.4) -- (4.2,1.4)
(7.8,1.4) -- (6.4,1.4) -- (6,2.8)
(7.8,4.2) -- (9.6,2.8) -- (11.4,2.8)
(9.6,2.8) -- (9.6,1.4);
\fill (1.4,7) circle (1.5pt) (2.8,5.6) circle (1.5pt) (1.4,4.2) circle (1.5pt) (2.8,4.2) circle (1.5pt) (4.2,4.2) circle (1.5pt) (6,4.2) circle (1.5pt) (7.8,4.2) circle (1.5pt) (9.6,4.2) circle (1.5pt) (11.4,4.2) circle (1.5pt) (.1,2.8) circle (1.5pt) (1.4,2.8) circle (1.5pt) (2.8,2.8) circle (1.5pt) (6,2.8) circle (1.5pt) (9.6,2.8) circle (1.5pt) (11.4,2.8) circle (1.5pt) (1.4,1.4) circle (1.5pt) (2.8,1.4) circle (1.5pt) (9.6,1.4) circle (1.5pt) (1.4,0) circle (1.5pt) (2.8,0) circle (1.5pt) (4.2,1.4) circle (1.5pt) (5.6,1.4) circle (1.5pt) (7.8,1.4) circle (1.5pt) (6.4,1.4) circle (1.5pt); (11.4,2.8) circle (1.5pt)
\draw (7.8,4.2) circle (3pt) (9.6,4.2) circle (3pt) (11.4,4.2) circle (3pt) (11.4,2.8) circle (3pt) (9.6,1.4) circle (3pt);
\draw [->] (2.35,6.6) node{\scriptsize $\lambda(\finsetz_2)$} (2.9,5.8) node{\scriptsize $\finsetz_2$} (1.4,7) -- (2.1,6.3);
\draw [->] (3.5,5.2) node{\scriptsize $\mu_2$} (2.8,5.6) -- (3.5,4.9);
\draw [->] (2.1,4.5) node{\scriptsize $\lambda(\finsetz_3)$} (2.8,4.45) node{\scriptsize $\finsetz_3$} (1.4,4.2) -- (2.1,4.2);
\draw [->] (3.5,4.5) node{\scriptsize $\mu_3$} (4.4,4.45) node{\scriptsize $\finsetz_1$} (2.8,4.2) -- (3.5,4.2);
\draw [->] (5.1,4.5) node{\scriptsize $\mu_1$} (6,4.5) node{\scriptsize $\finsetv$} (4.2,4.2) -- (5.1,4.2);
\draw [->] (.7,3.1) node{\scriptsize $\lambda(\finsetz_5)$} (1.4,3.05) node{\scriptsize $\finsetz_5$} (.1,2.8) -- (.7,2.8);
\draw [->] (2.1,3.1) node{\scriptsize $\mu_5$} (2.6,3.05) node{\scriptsize $\finsetz_4$} (1.4,2.8) -- (2.1,2.8);
\draw [->] (1,1.4) node{\scriptsize $O_{\finsetz_6}$} (2.1,1.7) node{\scriptsize $\lambda(\finsetz_6)$} (3.1,1.4) node{\scriptsize $\finsetz_6$} (1.4,1.4) -- (2.1,1.4);
\draw [->] (1,0) node{\scriptsize $O_{\finsetz_7}$} (2.1,.3) node{\scriptsize $\lambda(\finsetz_7)$} (3.1,0) node{\scriptsize $\finsetz_7$} (1.4,0) -- (2.1,0);
\draw [->] (3.1,.7) node{\scriptsize $\mu_7$} (2.8,0) -- (2.8,.7);
\draw [->] (3.1,2.1) node{\scriptsize $\mu_6$} (2.8,1.4) -- (2.8,2.1);
\draw [->] (3.5,3.8) node{\scriptsize $\mu_4$} (2.8,2.8) -- (3.5,3.5);
\draw [->] (6.3,3.5) node{\scriptsize $\mu_8$} (6.3,2.8) node{\scriptsize $\finsetz_8$} (6,2.8) -- (6,3.5);
\draw [->] (4.9,1.6) node{\scriptsize $\lambda(\finsetz_{9})$} (5.6,1.2) node{\scriptsize $\finsetz_{9}$} (4.2,1.4) -- (4.9,1.4);
\draw [->] (7.1,1.6) node{\scriptsize $\lambda(\finsetz_{10})$} (6.4,1.2) node{\scriptsize $\finsetz_{10}$} (7.8,1.4) -- (7.1,1.4);
\draw [->] (5.5,2.1) node{\scriptsize $\mu_{9}$} (5.6,1.4) -- (5.8,2.1);
\draw [->] (6.6,2.1) node{\scriptsize $\mu_{10}$} (6.4,1.4) -- (6.2,2.1);
\draw [->] (6.9,4.4) node{\scriptsize $\lambda$} (7,4.6) node[right]{\scriptsize $\begin{array}{ll}\finsetb &=\hat{\finsetb}(\finsetz_2)\\&=\hat{\finsetb}(\finsetz_5)\end{array}$} (7.8,4.2) -- (6.9,4.2);
\draw [->] (8.7,4) node{\scriptsize $u_a$} (9.3,3.9) node[right]{\scriptsize $\finsetb_a =\hat{\finsetb}(\finsetz_6)$} (9.6,4.2) -- (8.7,4.2);
\draw [->] (10.5,4.4) node{\scriptsize $u_b$} (11.45,4.2) node[right]{\scriptsize $\finsetb_b =\hat{\finsetb}(\finsetz_7)$} (11.4,4.2) -- (10.5,4.2);
\draw [->] (9,3.5) node{\scriptsize $u_c$} (9.5,2.7) node[left]{\scriptsize $\finsetb_c$} (9.6,2.8) -- (8.7,3.5);
\draw [->] (10.5,3) node{\scriptsize $u_d$} (11.45,2.8) node[right]{\scriptsize $\finsetb_d =\hat{\finsetb}(\finsetz_{9})$} (11.4,2.8) -- (10.5,2.8);
\draw [->] (9.8,2.1) node{\scriptsize $u_e$} (9.65,1.4) node[right]{\scriptsize $\finsetb_e =\hat{\finsetb}(\finsetz_{3}) =\hat{\finsetb}(\finsetz_{10})$} (9.6,1.4) -- (9.6,2.1);

\end{tikzpicture}\end{equation*}
\caption{The tree $\CT(\parentp,\parentp_{\finsetb})$ associated to a configuration $\confc$}
\label{figTbraidsgen}

\end{figure}

\bpo{Proof of Lemma~\ref{lemcnsonecybelongstonu}} We prove that  our configuration $(\confc^0,\confy^0)$ of Lemma~\ref{lemcnsonecybelongstonubis} and \ref{lemstrucCV} is in ${\CV}(\Gamma)$ by exhibiting a family $(\confc(t),\confy(t))$ of configurations of $\check{\CV}(\Gamma)$, for $t \in \left]0,\varepsilon\right[$ tending to $0$, which tends to $(\confc^0,\confy^0)$.
Our family will be described by nonvanishing $u_{\finsetd}(t)$, $\mu_{\finsetz}(t)$ tending to $0$, nonvanishing $\lambda(\finsety)(t)$ tending to our given $\lambda^0(\finsety)$, and nonvanishing $\lambda(t)$ tending to our given $\lambda^0$, such that the equations \begin{equation*}\ast(\finsety)(t) \colon \;\;M_{\finsety}(t) \lambda(\finsety)(t) =\lambda(t)  U_{\finsetd}(t) \end{equation*} are satisfied for any element $\finsety$ of $\parentp_{\cutx,\finsetd}$. 
The parameters $(\confy_{\finsetd})_{{\finsetd} \in \parentp_{\finsetb}}$, 
$\confc_{\finsetz}(\finsety)$ for nonunivalent kids $\finsety$ of elements $\finsetz$ of $\parentp$,
$p_{\RR}(\confc_{\finsetz}(\finsety))$ for univalent kids $\finsety$ of elements $\finsetz$ of $\parentp$,
of the statement of Lemma~\ref{lemstrucCV} will be defined to be the same as those of $(\confc^0,\confy^0)$.

Set $\mu_{\finsety}(t)=t$ for any $\finsety \in \parentp \setminus \left(\parentp^{\prime}_{\cutx} \cup \{\finsetv\} \right)$,  
and $u_{\finsetd}(t)=t$ for any $ \finsetd \in \parentp_{\finsetb} \setminus \left(\parentp^{\prime}_{\finsetb}\cup \{\finsetb\} \right)$.

Denote $\lambda(t)$ by $u_{\finsetb}(t)$ to make notation homogeneous. See Figure~\ref{figTbraidsgen}. Recall that we have $\lambda(\finsetv)=\lambda$ when $\finsetv \in \parentp_{\cutx,\finsetb}$.

For $\finsetd \in \parentp^{\prime}_{\finsetb}$, we are going to define some integer $g(\finsetd) \geq 1$ and 
set $u_{\finsetd}(t)=t^{g(\finsetd)}$ when $\finsetd\neq \finsetb$.
We are going to set $u_{\finsetb}(t)=\lambda(t)=t^{g(\finsetb)}$ if $\lambda^0=0$, and $u_{\finsetb}(t)=\lambda(t)$, where $\lambda(t)=\lambda^0$ if $\lambda^0\neq 0$.

Let $\finsety \in \parentp^{\prime}_{\cutx}$. If $\finsety$ is maximal in $\parentp^{\prime}_{\cutx}$, 
set $\lambda^{\prime}(\finsety)=1$. 
If $\finsety$ is not maximal in $\parentp^{\prime}_{\cutx}$, let $s(\finsety)$ be the minimal element of $\parentp^{\prime}_{\cutx}$ that contains $\finsety$ strictly. If $\lambda^0(s(\finsety))=0$, set $\lambda^{\prime}(\finsety)(t)=1$. If $\lambda^0(s(\finsety))\neq 0$, 
set $\lambda^{\prime}(\finsety)(t)=\lambda(s(\finsety))(t)$.
For $\finsety \in \parentp^{\prime}_{\cutx}$, we are going to define some integer $g(\finsety) \geq 1$ and 
set 
\begin{itemize}
 \item $\lambda(\finsety)(t)=t$ and $\mu_{\finsety}(t)=t^{g(\finsety)}\lambda^{\prime}(\finsety)$ if $\lambda^0(\finsety)=0$,
 \item $\lambda(\finsety)(t)=\lambda^0(\finsety)$
and
$\mu_{\finsety}(t)=\frac{t^{g(\finsety)+1}}{\lambda^0(\finsety)}\lambda^{\prime}(\finsety)(t)$ if $\lambda^0(\finsety)\neq 0.$
\end{itemize}
(If $\finsety=\finsetv$, just forget about $\mu_{\finsetv}$, which is useless.)

Order the elements of $\parentp^{\prime}_{\finsetb}$, by calling them $\finsetd_1$, $\finsetd_2$, \dots $\finsetd_{\cardlef{\parentp^{\prime}_{\finsetb}}}$ so that $\finsetd_1$ is maximal in $\parentp^{\prime}_{\finsetb}$ (with respect to the inclusion) and $\finsetd_{i+1}$ is maximal in $\parentp^{\prime}_{\finsetb} \setminus \{ \finsetd_1,\finsetd_2, \dots, \finsetd_i\}$. Note that $\finsetb \in \parentp^{\prime}_{\finsetb}$, so $\finsetd_1=\finsetb$.

Recall $\lambda(t)=t^{g(\finsetb)}$ if $\lambda^0=0$ and $\lambda(t)\neq 0$ if $\lambda^0\neq 0$. Set $g(\finsetd)=1$ for any $ \finsetd \in \widehat{\parentp^{\prime}}_{\finsetb} \setminus \parentp^{\prime}_{\finsetb}$.

We are going to define the integers $g(\finsetd_i)$ and the integers $g(\finsety)$ for every $\finsety \in \parentp_{\cutx,\finsetd_i}$, for $i \in \underline{\cardlef{\parentp^{\prime}_{\finsetb}}}$, inductively, so that the following assertion $(\ast(i))$ holds for every $i$.

$(\ast(i))$: for every $\finsety \in \parentp_{\cutx,\finsetd_i}$, we have $M_{\finsety}(t)\lambda(\finsety)(t)=\lambda(t) U_{\finsetd_i}(t)=\lambda(t)t^{f(i)}$
 for all $t\in \left]0,\varepsilon\right[$, where $f(i)=\sum_{\finsetd \in \widehat{\parentp^{\prime}}_{\finsetb} \suchthat \finsetd_i \subseteq \finsetd \subset \finsetb}g(\finsetd)$.

Recall $\finsetd_1=\finsetb$ and note $f(1)=0$.
If $\finsetv \in \parentp_{\cutx,\finsetb}$, set $g(\finsetd_1)=1$. So $(\ast(1))$ is satisfied. (If $\lambda^0\neq 0$, then $g(\finsetd_1)$ is not used, and $(\ast(1))$ is satisfied.)
If $\finsetv \notin \parentp_{\cutx,\finsetb}$, then $\lambda^0=0$, according to Lemma~\ref{lemcnsonecybelongstonubis}. Define $g(\finsetd_1)$ to be the maximum over the elements $\finsety$ of $\parentp_{\cutx,\finsetb}$ of the number of elements of $\parentp$ that  contain $\finsety$ strictly, plus one.
Then define the $g(\finsety)$ for the elements $\finsety$ of $\parentp_{\cutx,\finsetb}$ so that $(\ast(1))$ holds.

Let $k \in \underline{\cardlef{\parentp^{\prime}_{\finsetb}}}$. 

Assume that we have defined $g(\finsetd_i)$ and the $g(\finsety)$ for all $\finsety \in \parentp_{\cutx,\finsetd_i}$, for all $i<k$, so that the assertions $(\ast(i))$ hold for any $i <k$. Let us define the $g(\finsety)$ for all $\finsety \in \parentp_{\cutx,\finsetd_k}$ and $g(\finsetd_k)$ so that $(\ast(k))$ holds.
For any $\finsety \in \parentp_{\cutx,\finsetd_k}$, $M_{\finsety}(t)\lambda(\finsety)(t)$ is equal to $\lambda(t)t^{g(\finsety)+h(\finsety)}$ for some integer $h(\finsety)$, and $\lambda(t) U_{\finsetd_k}(t)$ is equal to $\lambda(t)t^{g(\finsetd_k)+h(\finsetd_k)}$ for some integer $h(\finsetd_k)$.
Let $H(k)$ be the maximum over $h(\finsetd_k)$ and all the 
integers $h(\finsety)$ for $\finsety \in \parentp_{\cutx,\finsetd_k}$. Set $g(\finsetd_k)=H(k)+1-h(\finsetd_k)$ and $g(\finsety)=H(k)+1-h(\finsety)$, so that $(\ast(k))$ holds with $f(k)=H(k)+1$.

This process is illustrated in Figure~\ref{figTbraidsgenexlim}. 
Its result does not depend on the arbitrary order that respects our condition on the $\finsetd_i$.

\bfig
\centering
\begin{equation*}
\begin{tikzpicture}[xscale=.9]
\draw (1.4,7) -- (4.2,4.2) -- (11.4,4.2) 
(1.4,4.2) -- (4.2,4.2) -- (2.8,2.8) -- (.1,2.8) 
(2.8,2.8) -- (2.8,0) -- (1.4,0) 
(1.4,1.4) -- (2.8,1.4) 
(6,4.2) -- (6,2.8) -- (5.6,1.4) -- (4.2,1.4)
(7.8,1.4) -- (6.4,1.4) -- (6,2.8)
(7.8,4.2) -- (9.6,2.8) -- (11.4,2.8)
(9.6,2.8) -- (9.6,1.4);
\fill (1.4,7) circle (1.5pt) (2.8,5.6) circle (1.5pt) (1.4,4.2) circle (1.5pt) (2.8,4.2) circle (1.5pt) (4.2,4.2) circle (1.5pt) (6,4.2) circle (1.5pt) (7.8,4.2) circle (1.5pt) (9.6,4.2) circle (1.5pt) (11.4,4.2) circle (1.5pt) (.1,2.8) circle (1.5pt) (1.4,2.8) circle (1.5pt) (2.8,2.8) circle (1.5pt) (6,2.8) circle (1.5pt) (9.6,2.8) circle (1.5pt) (11.4,2.8) circle (1.5pt) (1.4,1.4) circle (1.5pt) (2.8,1.4) circle (1.5pt) (9.6,1.4) circle (1.5pt) (1.4,0) circle (1.5pt) (2.8,0) circle (1.5pt) (4.2,1.4) circle (1.5pt) (5.6,1.4) circle (1.5pt) (7.8,1.4) circle (1.5pt) (6.4,1.4) circle (1.5pt); (11.4,2.8) circle (1.5pt)
\draw (7.8,4.2) circle (3pt) (9.6,4.2) circle (3pt) (11.4,4.2) circle (3pt) (11.4,2.8) circle (3pt) (9.6,1.4) circle (3pt);
\draw [->] (2.35,6.6) node{\scriptsize $t$} (2.9,5.8) node{\scriptsize $\finsetz_2$} (1.4,7) -- (2.1,6.3);
\draw [->] (3.5,5.2) node{\scriptsize $t^2$} (2.8,5.6) -- (3.5,4.9);
\draw [->] (2.1,4.5) node{\scriptsize $t$} (2.8,4.45) node{\scriptsize $\finsetz_3$} (1.4,4.2) -- (2.1,4.2);
\draw [->] (3.5,4.5) node{\scriptsize $t^4$} (4.4,4.45) node{\scriptsize $\finsetz_1$} (2.8,4.2) -- (3.5,4.2);
\draw [->] (5.1,4.5) node{\scriptsize $t$} (6,4.5) node{\scriptsize $\finsetv$} (4.2,4.2) -- (5.1,4.2);
\draw [->] (.7,3.1) node{\scriptsize $\lambda^0(\finsetz_5)$} (1.4,2.55) node{\scriptsize $\finsetz_5$} (.1,2.8) -- (.7,2.8);
\draw [->] (2.1,3.1) node{\scriptsize $\frac{t^2}{\lambda^0(\finsetz_5)}$} (3.1,2.7) node{\scriptsize $\finsetz_4$} (1.4,2.8) -- (2.1,2.8);
\draw [->] (2.1,1.7) node{\scriptsize $\lambda^0(\finsetz_6)$} (3.1,1.4) node{\scriptsize $\finsetz_6$} (1.4,1.4) -- (2.1,1.4);
\draw [->] (2.1,.3) node{\scriptsize $t$} (3.1,0) node{\scriptsize $\finsetz_7$} (1.4,0) -- (2.1,0);
\draw [->] (2.8,.7) node[right]{\scriptsize $t\lambda^0(\finsetz_6)$} (2.8,0) -- (2.8,.7);
\draw [->] (2.8,2.1) node[right]{\scriptsize $\frac{t^3}{\lambda^0(\finsetz_6)}$} (2.8,1.4) -- (2.8,2.1);
\draw [->] (3.5,3.8) node{\scriptsize $t$} (2.8,2.8) -- (3.5,3.5);
\draw [->] (6.3,3.5) node{\scriptsize $t$} (6.3,2.8) node{\scriptsize $\finsetz_8$} (6,2.8) -- (6,3.5);
\draw [->] (4.9,1.6) node{\scriptsize $t$} (5.6,1.2) node{\scriptsize $\finsetz_{9}$} (4.2,1.4) -- (4.9,1.4);
\draw [->] (7.1,1.6) node{\scriptsize $\lambda^0(\finsetz_{10})$} (6.4,1.2) node{\scriptsize $\finsetz_{10}$} (7.8,1.4) -- (7.1,1.4);
\draw [->] (5.5,2.1) node{\scriptsize $t^4$} (5.6,1.4) -- (5.8,2.1);
\draw [->] (6.3,2.2) node[right]{\scriptsize $\frac{t^5}{\lambda^0(\finsetz_{10})}$} (6.4,1.4) -- (6.2,2.1);
\draw [->] (6.9,4.4) node{\scriptsize $t^4$} (7,4.6) node[right]{\scriptsize $\begin{array}{ll}\finsetb &=\hat{\finsetb}(\finsetz_2)\\&=\hat{\finsetb}(\finsetz_5)\end{array}$} (7.8,4.2) -- (6.9,4.2);
\draw [->] (8.7,4) node{\scriptsize $t$} (9.3,3.9) node[right]{\scriptsize $\finsetb_a =\hat{\finsetb}(\finsetz_6)$} (9.6,4.2) -- (8.7,4.2);
\draw [->] (10.5,4.4) node{\scriptsize $t^2$} (11.45,4.2) node[right]{\scriptsize $\finsetb_b =\hat{\finsetb}(\finsetz_7)$} (11.4,4.2) -- (10.5,4.2);
\draw [->] (9,3.5) node{\scriptsize $t$} (9.5,2.7) node[left]{\scriptsize $\finsetb_c$} (9.6,2.8) -- (8.7,3.5);
\draw [->] (10.5,3) node{\scriptsize $t$} (11.45,2.8) node[right]{\scriptsize $\finsetb_d =\hat{\finsetb}(\finsetz_{9})$} (11.4,2.8) -- (10.5,2.8);
\draw [->] (9.8,2.1) node{\scriptsize $t$} (9.65,1.4) node[right]{\scriptsize $\finsetb_e =\hat{\finsetb}(\finsetz_{3}) =\hat{\finsetb}(\finsetz_{10})$} (9.6,1.4) -- (9.6,2.1);
\end{tikzpicture}\end{equation*}
\caption{Interior configurations tending to our limit configuration $\confc^0$, whose associated tree $\CT(\parentp,\parentp_{\finsetb})$ is as in Figure~\ref{figTbraidsgen},
where $\lambda^0(\finsetz_2)=\lambda^0(\finsetz_3)=\lambda^0(\finsetz_7)= \lambda^0(\finsetz_9)=0$ and the other $\lambda^0(\finsetz_i)$ are not zero.}
\label{figTbraidsgenexlim}

\end{figure}

\eop

Let us now begin the longer proof of Lemma~\ref{lemstrucCV}, which ends after Sublemma~\ref{sublempartialorder}.
The elements of ${\CV}(\Gamma)$ obviously satisfy the equations of the statement of Lemma~\ref{lemstrucCV}.
Let us conversely prove that an element $\confc$ parametrized as in this statement belongs to ${\CV}(\Gamma)$ by exhibiting a sequence of configurations $\confc(t)$ of $\check{\CV}(\Gamma)$ that approaches it. 
Again, we index our family of configurations $\confc(t)$ tending to the given one $\confc$ by a variable $t$ tending to $0$.

For $\finsety \in \parentp \setminus \bigl(\widehat{\parentp^{\prime}}_{\cutx} \cup \{\finsetv\} \bigr)$, set $\mu_{\finsety}=\mu_{\finsety}(\confc)$ and
\begin{equation*}
 \mu_{\finsety}(t)= \left\{\begin{array}{ll} t & \mbox{if }\mu_{\finsety}=0\\
 
                            \mu_{\finsety} &\mbox{if }\mu_{\finsety}\neq 0.
                           \end{array}
\right.\end{equation*}
Similarly, for $ \finsetd \in \parentp_{\finsetb} \setminus \bigl(\widehat{\parentp^{\prime}}_{\finsetb}\cup \{\finsetb\} \bigr)$, set
$u_{\finsetd}=u_{\finsetd}(\confc)$ and
\begin{equation*}
 u_{\finsetd}(t)= \left\{\begin{array}{ll} t & \mbox{if }u_{\finsetd}=0\\
 
                            u_{\finsetd} &\mbox{if }u_{\finsetd}\neq 0.
                           \end{array}
\right.\end{equation*}

We will again focus on the parameters $\mu_{\finsetz}$ for $\finsetz \in \widehat{\parentp^{\prime}}_{\cutx} \setminus \{\finsetv\}$,  $u_{\finsetd}$ for $ \finsetd \in \widehat{\parentp^{\prime}}_{\finsetb} \setminus \{\finsetb\}$, $\lambda(\finsety)$ for $\finsety \in \parentp^{\prime}_{\cutx}$, and $\lambda$.
All of these parameters correspond to edges of our tree $\CT(\parentp,\parentp_{\finsetb})$. We denote the parameter corresponding to an edge $e$ by $\lambda(e)$ and the corresponding parameter for the family by $\lambda(e)(t)$. The other parameters of Lemma~\ref{lemstrucCV}
are fixed as in the proof of Lemma~\ref{lemcnsonecybelongstonu}.

When $\lambda(e) \neq 0$, set $\lambda(e)(t)=\lambda(e)$. 

Let $\parentp^{\prime}_{\cutx,c}$ be the set of elements $\finsety$ of $\parentp^{\prime}_{\cutx}$
such that $\lambda(\finsety)M_{\finsety}=0$.
If $\finsety \in \parentp^{\prime}_{\cutx} \setminus \parentp^{\prime}_{\cutx,c}$, then
$\bigl(\lambda(\finsety)M_{\finsety}=\lambda U_{\hat{\finsetb}(\finsety)}\bigr)$ is not zero.
Let $\parentp^{\prime}_{\finsetb,c}$ be the set of elements $\finsetd$ of $\parentp^{\prime}_{\finsetb}$
such that $\lambda U_{\finsetd}=0$.

The main equations \begin{equation*}\ast(\finsety)(t) : \lambda(\finsety)(t)M_{\finsety}(t)=\lambda(t) U_{\hat{\finsetb}(\finsety)}(t)\end{equation*} associated to elements of $\parentp^{\prime}_{\cutx} \setminus \parentp^{\prime}_{\cutx,c}$ are obviously satisfied.

For $\finsety \in \parentp^{\prime}_{\cutx}$, recall that $O_{\finsety}$ denotes the initial vertex of the edge of $\lambda(\finsety)$ in the tree $\CT(\parentp,\parentp_{\finsetb})$.
For $\finsety \in \parentp^{\prime}_{\cutx,c}$, define $e(\finsety)$ to be the closest edge to $O_{\finsety}$ such that $\lambda(e(\finsety))=0$ between $O_{\finsety}$ and the vertex $\finsetv$ in $\CT(\parentp,\parentp_{\finsetb})$.
For $\finsetd \in \parentp^{\prime}_{\finsetb,c}$, define $e(\finsetd)$ to be the closest edge to $\finsetd$ such that $\lambda(e(\finsetd))=0$ between $\finsetd$ and the vertex $\finsetv$.

For edges $e$ of $\CT(\parentp,\parentp_{\finsetb})$ such that $\lambda(e)=0$ that are not in $e(\parentp^{\prime}_{\cutx,c} \cup \parentp^{\prime}_{\finsetb,c})$, set $\lambda(e)(t)=t$.
For edges $e \in e(\parentp^{\prime}_{\cutx,c} \cup \parentp^{\prime}_{\finsetb,c})$, set $\lambda(e)(t)=r(e)t^{k(e)}$, where $(r(e),k(e)) \in \RR^+ \times \NN$, $r(e) \neq 0$, and $k(e) \neq 0$. We are going to show how to define pairs $(r(e),k(e))$, so that the equations 
$\ast(\finsety)(t)$ are satisfied, for all $t$ and for all $\finsety$ such that $\finsety \in \parentp^{\prime}_{\cutx}$. Since they imply the equations~\ref{eqcycconst}, because all the coefficients $\lambda(e)(t)$ are different from zero, this will conclude the proof.

If $\parentp^{\prime}_{\finsetb,c} =\emptyset$, then no parameter $\lambda(e)$ vanishes. So there is nothing to prove.

For a set $I$ of edges of $\CT(\parentp,\parentp_{\finsetb})$, $\lambda(I)$ denotes the product over the edges $e$ of $I$ of the $\lambda(e)$, while $\lambda_t(I)$ denotes the product over the edges $e$ of $I$ of the $\lambda(e)(t)$. With the notation of Definition~\ref{deffirsttree}, the main equations $\ast(\finsety)(t)$ may be written as
\begin{equation*}\ast(\finsety)(t) : \lambda_t\bigl([O_{\finsety},V]\bigr)=\lambda_t\bigl([\hat{\finsetb}(\finsety),V]\bigr).\end{equation*}

We also have the obvious sublemma.

\begin{sublemma}
\label{sublemcyccond} Condition~\ref{eqcycconst} of Lemma~\ref{lemstrucCV} may be rewritten as follows.
For any collections $(\finsety_i)_{i \in \ZZ/n\ZZ}$, $(\finsety^{+}_i)_{i \in \ZZ/n\ZZ}$, $(\finsety^{\prime}_i)_{i \in \ZZ/n\ZZ}$ such that $\finsety_i \in \parentp^{\prime}_{\cutx}$, $\finsety^+_i \in \parentp^{\prime}_{\cutx}$, $\finsety^{\prime}_i \in \widehat{\parentp^{\prime}}_{\cutx}$, and 
$\left(\finsety^{+}_{i} \cup \finsety_i\right) \subseteq \finsety^{\prime}_i$, and for any collection $(\finsetb^{\prime}_i)_{i \in \ZZ/n\ZZ}$ of sets of $\widehat{\parentp^{\prime}}_{\finsetb}$ such that
\begin{equation*}\bigl(\hat{\finsetb}(\finsety^{+}_{i-1}) \cup \hat{\finsetb}(\finsety_i)\bigr) \subseteq \finsetb^{\prime}_i,\end{equation*}
the products of the numbers over the arrows in one direction equals the product of the numbers over the arrows in the opposite direction in the cycle of Figure~\ref{figcyccond}.

\bfig
\centering

\begin{equation*} 
\begin{tikzpicture}
\useasboundingbox (.5,-.3) rectangle (13,5.4);
\draw
(2,4) node{\scriptsize $\finsetb^{\prime}_1$}
(2.9,4.4) node{\scriptsize $\lambda([\hat{\finsetb}(\finsety_1),\finsetb^{\prime}_1])$}
(3.8,4) node{\scriptsize $\finsety_1$}
(4.7,4.4) node{\scriptsize $\lambda([O_{\finsety_1},\finsety^{\prime}_1])$}
(5.6,4) node{\scriptsize $\finsety^{\prime}_1$} 
(6.5,4.3) node{\scriptsize $\lambda([O_{\finsety_1^+},\finsety^{\prime}_1])$}
(7.6,4) node{\scriptsize $\finsety^+_1$}
(8.7,4.4) node{\scriptsize $\lambda([\hat{\finsetb}(\finsety_1^+),\finsetb^{\prime}_2])$}
(9.8,4) node{\scriptsize $\finsetb^{\prime}_2$} 
(10.8,4.4) node{\scriptsize $\lambda([\hat{\finsetb}(\finsety_2),\finsetb^{\prime}_2])$}
(11.8,4) node{\scriptsize $\finsety_2$};
\draw [<-] (2.3,4) -- (3.5,4) ;
\draw [->] (4.1,4) -- (5.3,4) ;
\draw [<-] (5.9,4) -- (7.3,4) ;
\draw [->] (7.9,4) -- (9.5,4) ;
\draw [<-] (10.1,4) -- (11.5,4) ;
\draw [->] (12.1,4) -- (12.5,4) ;
\draw [dotted] (12.6,4) .. controls (13.2,4) and (12.3,3.2) .. (6,3.2) ..  controls (0,3.2) and (.6,2) .. (1.2,2);
\draw
(2,2) node{\scriptsize $\finsetb^{\prime}_i$}
(2.9,2.4) node{\scriptsize $\lambda([\hat{\finsetb}(\finsety_i),\finsetb^{\prime}_i])$}
(3.8,2) node{\scriptsize $\finsety_i$}
(4.7,2.4) node{\scriptsize $\lambda([O_{\finsety_i},\finsety^{\prime}_i])$}
(5.6,2) node{\scriptsize $\finsety^{\prime}_i$}
(6.5,2.3) node{\scriptsize $\lambda([O_{\finsety_i^+},\finsety^{\prime}_i])$} 
(7.6,2) node{\scriptsize $\finsety^+_i$}
(8.7,2.4) node{\scriptsize $\lambda([\hat{\finsetb}(\finsety_i^+),\finsetb^{\prime}_{i+1}])$}
(9.8,2) node{\scriptsize $\finsetb^{\prime}_{i+1}$} 
(11.8,2) node{\scriptsize $\finsety_{i+1}$} ;
\draw [->] (1.3,2) -- (1.7,2) ;
\draw [<-] (2.3,2) -- (3.5,2) ;
\draw [->] (4.1,2) -- (5.3,2) ;
\draw [<-] (5.9,2) -- (7.3,2) ;
\draw [->] (7.9,2) -- (9.4,2) ;
\draw [<-] (10.2,2) -- (11.4,2) ;
\draw [->] (12.2,2) -- (12.5,2) ;
\draw [dotted] (12.6,2) .. controls (13.2,2) and (12.3,1.2) .. (6,1.2) ..  controls (0,1.2) and (2,.4) .. (5.3,.4);
\draw
(5.6,.4) node{\scriptsize $\finsety^{\prime}_n$}
(7.6,.4) node{\scriptsize $\finsety^+_n$};
\draw [<-] (5.9,.4) -- (7.3,.4);
\draw [->] (7.9,.4) .. controls (8.5,.4) and (8.5,-.2) .. (2,-.2) ..  controls (0,-.2) and (0,4) .. (1.7,4);
 \end{tikzpicture}\end{equation*}
\caption{The cycle of Condition~\ref{eqcycconst} of Lemma~\ref{lemstrucCV}.}
\label{figcyccond}

\end{figure}

\end{sublemma}

We can redraw a typical part of the cycle of Figure~\ref{figcyccond} as in Figure~\ref{figcyccondbis} to emphasize the different natures of its segments and show how the sets at their ends determine the coefficients over the arrows.
\bfig
\centering
\begin{tikzpicture}[xscale=.9]
\draw
(1,1.6) node{\scriptsize $\finsetb^{\prime}_i$}
(1.22,.7) node{\scriptsize $\lambda([\hat{\finsetb}(\finsety_i),\finsetb^{\prime}_i])$}
(2.6,.2) node{\scriptsize $\hat{\finsetb}(\finsety_i)$}
(3.5,.6) node{\scriptsize $\finsety_i$}
(4.4,.2) node{\scriptsize $O_{\finsety_i}$}
(4.2,1.15) node{\scriptsize $\lambda([O_{\finsety_i},\finsety^{\prime}_i])$}
(6,1.6) node{\scriptsize $\finsety^{\prime}_i$}
(6.15,.7) node{\scriptsize $\lambda([O_{\finsety_i^+},\finsety^{\prime}_i])$} 
(7.6,.2) node{\scriptsize $O_{\finsety_i^+}$}
(8.5,.6) node{\scriptsize $\finsety_i^+$}
(9.4,.2) node{\scriptsize $\hat{\finsetb}(\finsety_i^+)$}
(8.9,1.15) node{\scriptsize $\lambda([\hat{\finsetb}(\finsety_i^+),\finsetb^{\prime}_{i+1}])$}
(11,1.6) node{\scriptsize $\finsetb^{\prime}_{i+1}$} 
(12.6,.2) node{\scriptsize $\hat{\finsetb}(\finsety_{i+1})$};
\draw [->] (2.6,.6) -- (1.2,1.4) ;
\draw [->] (4.4,.6) -- (5.8,1.4) ;
\draw [->] (7.6,.6) -- (6.2,1.4) ;
\draw [->] (9.4,.6) -- (10.8,1.4) ;
\draw [->] (12.6,.6) -- (11.2,1.4) ;
\end{tikzpicture}
\caption{Another representation of a part of the cycle of Figure~\ref{figcyccond}.}
\label{figcyccondbis}

\end{figure}

For $g \in e(\parentp^{\prime}_{\finsetb,c})$ and for $f\in e(\parentp^{\prime}_{\cutx,c})$, let $\parentp(f,g)$ be the set of elements $\finsety \in \parentp^{\prime}_{\cutx}$ such that $e(\finsety)=f$ and $e(\hat{\finsetb}(\finsety))=g$.
Note the following sublemma.

\begin{sublemma}
\label{sublemeqcvone} Recall the notation of Definition~\ref{deffirsttree} for sets of edges of a tree.
Let $g \in e(\parentp^{\prime}_{\finsetb,c})$, let $f\in e(\parentp^{\prime}_{\cutx,c})$. Assume that there exists an element $\finsety \in \parentp(f,g)$.
Set \begin{equation*}c(f,g)=\frac{\lambda([\hat{\finsetb}(\finsety),g[)}{\lambda([O_{\finsety},f[)}.\end{equation*}
Then $c(f,g)$ is a positive coefficient independent of the chosen element  $\finsety$ of $\parentp(f,g)$.
The equation 
\begin{equation*}\ast_t(f,g): \lambda_t\bigl(\left[f,\finsetv\right]\bigr)=c(f,g) \lambda_t\bigl(\left[g,\finsetv\right]\bigr)\end{equation*}
must be satisfied for our sequence of configurations $\confc(t)$ to be in $\check{\CV}(\Gamma)$.
Furthermore, if Equation~$\ast_t(f,g)$ is satisfied for all $t$, then for any $\finsety^+\in \parentp(f,g)$,
Equation~$\ast(\finsety^+)(t)$ is satisfied for all $t$.
\end{sublemma}
\bp With the given coefficient $c(f,g)$, Equation~$\ast_t(f,g)$ is equivalent to Equation~$\ast(\finsety)(t)$.
For any other element $\finsety^+$ of $\parentp(f,g)$, according to Condition~\ref{eqcycconst} of Lemma~\ref{lemstrucCV}---where $n=1$, $Y^{\prime}_1$ is the origin of $f$, and $\finsetb^{\prime}_1$ is the origin of $g$,
we have \begin{equation*}\lambda\bigl([O_{\finsety},f[\bigr)\lambda\bigl([\hat{\finsetb}(\finsety^+),g[\bigr)=\lambda\bigl([O_{\finsety^+},f[\bigr)\lambda\bigl([\hat{\finsetb}(\finsety),g[\bigr).\end{equation*} \eop

Define the equivalence relation $\sim$ on $e(\parentp^{\prime}_{\cutx,c} \cup \parentp^{\prime}_{\finsetb,c})$ to be the relation generated
by the equivalences: Whenever $g \in e(\parentp^{\prime}_{\finsetb,c})$ and $f \in e(\parentp^{\prime}_{\cutx,c})$, if
$\parentp(f,g)\neq \emptyset$, then we have $g \sim f$. 
When $\finsety \in \parentp(f,g)$, such a generating elementary equivalence will also be denoted by $g \sim_{\finsety} f$, and its inverse will be denoted by $f \sim_{\finsety} g$. In this case, set $c(g,f)=c(f,g)^{-1}$. 

\begin{sublemma}
\label{sublemeqcyc}
For any cycle
\begin{equation*}g=g_1  \sim_{\finsety_1} f_1 \sim_{\finsety^+_1} g_2  \sim_{\finsety_2} f_2 \sim_{\finsety^+_2} g_3 \dots f_n \sim_{\finsety^+_n} g_{n+1}=g\end{equation*} 
of elementary equivalences
for sequences $(g_i)_{i \in \ZZ/n\ZZ}$ of edges of $e(\parentp^{\prime}_{\finsetb,c})$ and $(f_i)_{i \in \ZZ/n\ZZ}$ of edges of $e(\parentp^{\prime}_{\cutx,c})$, we have
\begin{equation*}\prod_{i=1}^n c(g_i,f_i)c(f_i,g_{i+1}) = 1.\end{equation*}
\end{sublemma}
\bp
Apply Condition~\ref{eqcycconst} of Lemma~\ref{lemstrucCV}, as stated in Sublemma~\ref{sublemcyccond}, to the above sequences $(\finsety_i)$ and $(\finsety^+_i)$, where $Y^{\prime}_i$ is the origin $o(f_i)$ of $f_i$ and $\finsetb^{\prime}_i$ is the origin $o(g_i)$ of $g_i$. Note that Sublemma~\ref{sublemeqcvone} implies that, in each triangle of Figure~\ref{figcyccondtriangles}, the product of the coefficients over the edges of the boundary of the triangle that are oriented as part of that boundary is equal to the product of the coefficients over the edges with the opposite orientation.
\bfig
\centering
\begin{tikzpicture}[xscale=.9]
\draw
(1,1.6) node{\scriptsize $o(g_i)$}
(1.22,.7) node{\scriptsize $\lambda([\hat{\finsetb}(\finsety_i),g_i[)$}
(2.6,.2) node{\scriptsize $\hat{\finsetb}(\finsety_i)$}
(3.5,.6) node{\scriptsize $\finsety_i$}
(4.4,.2) node{\scriptsize $O_{\finsety_i}$}
(4.2,1.15) node{\scriptsize $\lambda([O_{\finsety_i},f_i[)$}
(6,1.6) node{\scriptsize $o(f_i)$}
(6.15,.7) node{\scriptsize $\lambda([O_{\finsety_i^+},f_i[)$} 
(7.6,.2) node{\scriptsize $O_{\finsety_i^+}$}
(8.5,.6) node{\scriptsize $\finsety_i^+$}
(9.4,.2) node{\scriptsize $\hat{\finsetb}(\finsety_i^+)$}
(8.9,1.15) node{\scriptsize $\lambda([\hat{\finsetb}(\finsety_i^+),g_{i+1}[)$}
(11,1.6) node{\scriptsize $o(g_{i+1})$} 
(12.6,.2) node{\scriptsize $\hat{\finsetb}(\finsety_{i+1})$}
(3.5,1.8) node{\scriptsize $c(g_i,f_i)$}
(8.5,1.8) node{\scriptsize $c(f_i,g_{i+1})$};
\draw [->] (2.6,.6) -- (1.2,1.4) ;
\draw [->] (4.4,.6) -- (5.8,1.4) ;
\draw [->] (7.6,.6) -- (6.2,1.4) ;
\draw [->] (9.4,.6) -- (10.8,1.4) ;
\draw [->] (12.6,.6) -- (11.2,1.4) ;
\draw [->] (2,1.6) -- (5,1.6);
\draw [->] (7,1.6) -- (10,1.6);
\end{tikzpicture}
\caption{Definition of $c(g_i,f_i)$ and $c(f_i,g_{i+1})$.}
\label{figcyccondtriangles}

\end{figure}
\eop

The following sublemma is an easy corollary of the previous one.
\begin{sublemma}
\label{sublemeqcvtwo}
Let  $e$ and $e^{\prime}$ be two elements of $e(\parentp^{\prime}_{\cutx,c} \cup \parentp^{\prime}_{\finsetb,c})$ such that $e \sim e^{\prime}$.
There exists a sequence $(e_i)_{i \in \underline{m}}$ of edges of $e(\parentp^{\prime}_{\cutx,c} \cup \parentp^{\prime}_{\finsetb,c})$ such that $e_1=e$, $e_m=e^{\prime}$, and $e_i \sim_{\finsetz_i} e_{i+1}$ for any $i \in \underline{m-1}$.
For such a sequence, set \begin{equation*}c(e,e^{\prime})=\prod_{i=1}^{m-1}c(e_i,e_{i+1}).\end{equation*}
Then $c(e,e^{\prime})$ is a positive coefficient independent of the chosen sequences as above.

The equation
\begin{equation*}\ast_t(e,e^{\prime}): \lambda_t\bigl(\left[e,\finsetv\right]\bigr)=c(e,e^{\prime}) \lambda_t\bigl(\left[e^{\prime},\finsetv\right]\bigr)\end{equation*}
must be satisfied for our sequence of configurations $\confc(t)$ to be in $\check{\CV}(\Gamma)$.
Furthermore, for any three elements $e$, $e_0$, and $e^{\prime}$ of $e(\parentp^{\prime}_{\finsetb,c})$ in the same equivalence class under $\sim$, the equations $\ast_t(e,e^{\prime})$ and $\ast_t(e^{\prime},e)$ are equivalent, and the equations $\ast_t(e,e_0)$ and $\ast_t(e^{\prime},e_0)$ imply $\ast_t(e,e^{\prime})$.
\end{sublemma}
\eopwobp

The following sublemma will allow us to define a partial order on the set of equivalence classes under $\sim$.
\begin{sublemma}
\label{sublempartialorder}
Let $k \in \NN$. It is not possible to find edges $e_1, \dots, e_k$, $e^{\prime}_1,\dots e^{\prime}_k$ of $e(\parentp^{\prime}_{\cutx,c} \cup \parentp^{\prime}_{\finsetb,c})$ such that
$e_j \sim e^{\prime}_j$ and $e_{j+1} \in \bigl]e^{\prime}_j,\finsetv\bigr]$ for all $j \in \ZZ/k\ZZ$.
\end{sublemma}
\bp
Assume that there exist edges $e_1, \dots, e_k$, $e^{\prime}_1,\dots e^{\prime}_k$  of $e(\parentp^{\prime}_{\cutx,c} \cup \parentp^{\prime}_{\finsetb,c})$ such that
$e_j \sim e^{\prime}_j$ and $e_{j+1} \in \bigl]e^{\prime}_j,\finsetv\bigr]$ for all $j \in \ZZ/k\ZZ$.
Let us picture this hypothesis as follows.

\begin{center}
\begin{tikzpicture}
\draw (2,.5) -- (3,1) (5,1) -- (6,1.5);
\draw[dotted] (-.1,.25) -- (.4,.5) (8,1.5) -- (8.5,1.75);
\draw (1.2,.5) node{\scriptsize $e_j \sim e^{\prime}_j$}
(4,1) node{\scriptsize $e_{j+1} \sim e^{\prime}_{j+1}$}
(7,1.5) node{\scriptsize $e_{j+2} \sim e^{\prime}_{j+2}$};
\end{tikzpicture}
\end{center}
Let us now construct a cycle as in Figure~\ref{figcyccond} under this hypothesis.
An elementary equivalence $g \sim_{\finsety} f$ where $\finsety \in \parentp(f,g)$ is represented by
a path $E(g,f)$ 
\begin{equation*}
\begin{tikzpicture}
\draw
(1.9,4) node{\scriptsize $o(g)$}
(2.9,4.4) node{\scriptsize $\lambda([\hat{\finsetb}(\finsety),g[)$}
(3.8,4) node{\scriptsize $\finsety$}
(4.7,4.4) node{\scriptsize $\lambda([O_{\finsety},f[)$}
(5.7,4) node{\scriptsize $o(f)$} ;
\draw [<-] (2.3,4) -- (3.5,4) ;
\draw [->] (4.1,4) -- (5.3,4) ;
\end{tikzpicture}\end{equation*}
where $o(g)$ is the set associated to the origin of $g$ and $o(f)$ is associated to the origin of $f$. (It is equal to $O_{\finsety}$ if $\lambda(\finsety)=0$, and it is the set associated to the initial point of $f$ otherwise.)
The inverse equivalence $f \sim_{\finsety} g$ is similarly represented by the following path $E(f,g)$
\begin{equation*}
\begin{tikzpicture}
\draw
(1.9,4) node{\scriptsize $o(f)$}
(2.9,4.4) node{\scriptsize $\lambda([O_{\finsety},f[)$}
(3.8,4) node{\scriptsize $\finsety$}
(4.7,4.4) node{\scriptsize $\lambda([\hat{\finsetb}(\finsety),g[)$}
(5.7,4) node{\scriptsize $o(g)$} ;
\draw [<-] (2.3,4) -- (3.5,4) ;
\draw [->] (4.1,4) -- (5.3,4) ;
\end{tikzpicture}.\end{equation*}
Again, the coefficients over the arrows are determined by the labels of the ends. In these cases, they do not vanish, and we picture the paths $E(g,f)$ and $E(f,g)$ as
$o(g) \hookleftarrow \finsety \hookrightarrow o(f)$ and $o(f) \hookleftarrow \finsety \hookrightarrow o(g)$, respectively.
Note that according to our assumptions, no $e_j$ can start at some $O_Y$.

If $e_j \neq e^{\prime}_j$ for all $j \in \ZZ/k\ZZ$, then our cycle of arrows as above is obtained by
assembling the following paths of arrows $A(e_j,e_{j+1})$ from $o(e_j)$ to $o(e_{j+1})$.
The path $A(e_j,e_{j+1})$ is obtained from a sequence $E(e_j,e^{\prime}_j)$ of paths $E(e,e^{\prime})$ of arrows associated to elementary equivalences 
by replacing the last arrow $X_j \hookrightarrow o(e^{\prime}_j)$, which ends at $o(e^{\prime}_j)$, by an arrow $X_j \rightarrow o(e_{j+1})$ ending at $o(e_{j+1})$. 
The coefficient of this arrow is obtained from the coefficient of $X_j \hookrightarrow o(e^{\prime}_j)$
by multiplication by $\lambda([e^{\prime}_j,e_{j+1} [)$, which vanishes in our case because it has a factor $\lambda(e^{\prime}_j)$.
According to the recalled criterion, the product of the coefficients over the arrows in the direction of our cycle
must be equal to the product of the coefficients over the arrows in the opposite direction.
The first product is zero because of its factors $\lambda(e^{\prime}_j)$. The second one is nonzero because it only contains nonzero factors associated to equivalences. Therefore, the lemma is proved when $e_j \neq e^{\prime}_j$ for any $j \in \ZZ/k\ZZ$. This case is ruled out.

We cannot have $e_j=e^{\prime}_j$ for all $j \in \ZZ/k\ZZ$. 

Up to permuting our indices cyclically, it suffices to rule out the case in which $k\geq 2$ and there exists $r \geq 1$ such that $e_r \neq e^{\prime}_r$, $e_j=e^{\prime}_j$ for all $j$ such that $r+1 \leq j \leq k$, and
$e_1 \neq e^{\prime}_1$. In this case, we define a path $A(e_r,e_1)$ by replacing the last arrow $X_r \hookrightarrow o(e^{\prime}_r)$ in $E(e_r,e^{\prime}_r)$ with $X_r \rightarrow o(e_1)$, and multiplying the corresponding coefficient by $\lambda([e^{\prime}_r,e_1[)=\prod_{j=r}^{k}\lambda([e^{\prime}_j,e_{j+1}[)$.
Similarly define paths $A(e_s,e_t)$, for all pair $(s,t)$ of integers such that $s<t<k$, $e_s \neq e^{\prime}_s$, $e_t \neq e^{\prime}_t$, and $e_j=e^{\prime}_j$ for all $j$ such that $s+1 \leq j \leq t-1$. Define the cycle yielding the contradiction by composing these paths, which include our former paths $A(e_j,e_{j+1})$ for which $e_j\neq e^{\prime}_j$ and $e_{j+1}\neq e^{\prime}_{j+1}$.
\eop

\bpo{End of the proof of Lemma~\ref{lemstrucCV}}
Sublemma~\ref{sublempartialorder} allows us to define the partial order $\succeq$ on the set ${E(\sim)}$ 
of equivalence classes of the relation $\sim$ on $e(\parentp^{\prime}_{\cutx,c} \cup \parentp^{\prime}_{\finsetb,c})$ such that
two equivalence classes $\overline{e}$ and $\overline{f}$ of $E(\sim)$ satisfy
$\overline{f} \succeq \overline{e}$ if and only if there exist a positive integer $k \geq 2$ and two sequences $(e_i)_{i\in \underline{k}\setminus \{1\}}$ and $(e^{\prime}_i)_{i\in \underline{k-1}}$ such that $e^{\prime}_1 \in \overline{e}$,
$e_k \in \overline{f}$, $e_j \sim e^{\prime}_j$ for all $j \in \underline{k-1}\setminus \{1\}$, and $e_{j+1} \in \bigl[e^{\prime}_j,\finsetv\bigr]$
 (i.e., $e_{j+1}$ is between the initial point of $e^{\prime}_j$ and $\finsetv$) for all $j \in \underline{k-1}$. For example, if $e_2 \in \bigl[e^{\prime}_1,\finsetv\bigr]$, then $\overline{e_2} \succeq \overline{e^{\prime}_1}$.

Fix an arbitrary total order on $E(\sim)$ compatible with the above partial order by writing
\begin{equation*}E(\sim)=\{\overline{g_i} \suchthat i \in \underline{m}\}\end{equation*}
so that for any $(i,j) \in \underline{m}^2$ such that $\overline{g_j} \succeq \overline{g_i}$, we have $j \leq i$.

We will pick one representative $g_i \in e(\parentp^{\prime}_{\finsetb,c})$ in each equivalence class $\overline{g_i}$ of $E(\sim)$
and define
$\lambda(g_i)(t)=r(g_i)t^{k(g_i)}$, for $i=1,\dots,m$, inductively.

For $e \in\overline{g_i}$, the equation 
\begin{equation*}\ast_t(e,g_i): \lambda_t\bigl(\left[e,\finsetv\right]\bigr)=c(e,g_i) \lambda_t\bigl(\left[g_i,\finsetv\right]\bigr),\end{equation*}
which must be satisfied, determines $\lambda(e)(t)=r(e)t^{k(e)}$ with $(r(e),k(e)) \in \RR^{+\ast} \times (\NN \setminus \{0\})$ as a function of $\lambda(g_i)(t)$ and of the $\lambda(e^{\prime})(t)=r(e^{\prime})t^{k(e^{\prime})}$
for edges $e^{\prime}$ that belong to $\cup_{j=1}^{i-1}\overline{g_j}$. (Recall that we fixed the coefficients of the other edges.)
More precisely, $k(e)-k(g_i)$ is a degree one polynomial in the variables $k(e^{\prime})$ for edges $e^{\prime}$ that belong to $\cup_{j=1}^{i-1}\overline{g_j}$,
which are already defined by induction.
In particular, if we choose $k(g_i)$ to be a sufficiently large integer and fix $r(g_i)=1$, then the $\lambda(e)(t)$ are uniquely determined 
so that the equations $\ast_t(e,g_i)$ are satisfied with $k(e) >0$ for any $e \in \overline{g_i}$.

Thus, once the induction is achieved, according to Sublemmas~\ref{sublemeqcvone} and \ref{sublemeqcvtwo},
Equation~$\ast(\finsety)(t)$ is satisfied for all $t$ for any $\finsety\in \parentp$.
This shows that the elements described in the statement of Lemma~\ref{lemstrucCV} are in ${\CV}(\Gamma)$ and finishes the proof of Lemma~\ref{lemstrucCV}. \eop

Let $\parentp_u$ denote the set of univalents elements of $\parentp$, which are the elements of $\parentp$ that contain at least one univalent vertex.

For a configuration $\confy \in \ccompuptd{\finsetb}{\CC}$ and a Jacobi diagram $\Gamma$ on a disjoint union of lines $\RR_{\eltb}$ indexed by elements $\eltb$ of $\finsetb$, let ${\CV}(\confy,\Gamma)$ \index[N]{Vcalp@${\CV}(\confy,\Gamma)$ configuration space} denote the preimage of $\confy$ under
\begin{equation*}p_{\ccompuptdempty{\finsetb}} \colon \CV(\Gamma) \to \ccompuptd{\finsetb}{\CC}.\end{equation*}

\begin{lemma}
\label{lemdimstratCV} With respect to the notation introduced before  Lemma~\ref{lemcnsonecybelongstonubis} and in Lemma~\ref{lemcnsonecybelongstonubis}, 
let $\parntppx$\index[N]{Parenthesizations!papcutx@$\parntppx$} be the subset of $\parentp^{\prime}_{\cutx}$ consisting of its elements $\finsety$ such that $\lambda^0(\finsety) \neq 0$.
The dimension of the stratum of $\confc^0$ in the fiber ${\CV}(\confy^0,\Gamma)$ over the configuration $\confy^0 \in \ccompuptd{\finsetb}{\CC}$ is $\cardlef{U(\Gamma)} + 3 \cardlef{T(\Gamma)}-\cardlef{\parentp} + \cardlef{\parntppx} -1$.
\end{lemma}
\bp
For any $\finsety \in \parentp \setminus \parentp_u$, the configuration $c_{\finsety}$ is defined up to global translation and dilation. (With our normalizations, the quotient by translations is replaced by the fact that we send basepoints to zero, and our horizontal or vertical normalization condition replaces the quotient by dilation.)
For any $\finsety \in  \parentp_u \setminus \parntppx$, the configuration $c_{\finsety}$, whose restriction to the set of univalent vertices is vertical, is defined up to global vertical translation and dilation.
For any $\finsety \in \parntppx$, we still have these two codimension-one normalization conditions on the configuration $c_{\finsety}$ when $\lambda(\finsety)$ is fixed. But varying the parameter $\lambda(\finsety)$ adds one to the dimension.

Say that a subset of $\finsetv(\Gamma)$ is \emph{trivalent} if it contains no univalent vertex. Write the set $\kids(\finseta)$ of kids of an element $\finseta$ of $\parentp$ 
as the union of its set $\kids^0_{u}(\finseta)$ of univalent kids and its set $\kids^0_{t}(\finseta)$ of trivalent kids (including the kid of the basepoint). \begin{equation*}\kids(\finseta)=\kids^0_{u}(\finseta) \sqcup \kids^0_{t}(\finseta).\end{equation*}
Then the dimension of the involved space of maps up to dilation and (possibly vertical) translation from $\kids(\finseta)$ to $\RR^3$ is 
\begin{itemize}
\item $3 \cardlef{\kids^0_{t}(\finseta)} -4 $ if $\finseta \in \parentp \setminus \parentp_u $ ,
\item $\cardlef{\kids^0_{u}(\finseta)} + 3 \cardlef{\kids^0_{t}(\finseta)} -2 $ if $\finseta \in \parentp_u \setminus  \parntppx$,
\item $\cardlef{\kids^0_{u}(\finseta)} + 3 \cardlef{\kids^0_{t}(\finseta)} -1 $ if $\finseta \in \parntppx$.
\end{itemize}
Let $\parentp_{ext}$ denote the union of $\parentp$ with the set of singletons of elements of $\finsetv$. So $\parentp_{ext}$ contains all the kids of elements of $\parentp$. The only element of $\parentp_{ext}$ that is not a kid is $\finsetv$. Let $\parentp_{ext,u}$ denote the set of univalent sets of $\parentp_{ext}$,
and let $\parentp_{ext,t}$ denote the set of trivalent sets of $\parentp_{ext}$.
The sum over $\finseta \in \parentp$ of the above dimensions is equal to
\begin{equation*}3 \cardbig{\parentp_{ext,t}} +  \cardbig{\parentp_{ext,u}} - 1 - 3 \cardbig{\parentp \setminus \parentp_u} - \cardbig{\parentp_u} - \cardbig{\parentp \setminus \parntppx},\end{equation*}
where  $\cardbig{\parentp_{ext,t}}=\cardbig{T(\Gamma)}+\cardbig{\parentp \setminus \parentp_u}$,
and $\cardbig{\parentp_{ext,u}}=\cardbig{U(\Gamma)}+\cardbig{\parentp_u}$.
\eop

Lemma~\ref{lemstrucCV} simplifies when $\confy^0 \in \cinjuptd{\finsetb}{\CC}$. It allows us to describe the structure of ${\CV}(\Gamma)$ over $\cinjuptd{\finsetb}{\CC}$ in the following lemma.

\begin{lemma}
\label{lemstrucCVinjtree} 
Let $\confy^0 \in \cinjuptd{\finsetb}{\CC}$. For any $(\confc^0,\confy^0) \in {\CV}(\Gamma)$,
there exist
\begin{itemize}
 \item a manifold $W$ with boundary and ridges, 
 \item an open neighborhood $N(\confy^0)$ of $\confy^0$ in $\cinjuptd{\finsetb}{\CC}$, 
 \item a small $\varepsilon >0$,
 \item an oriented tree $\CT_{\CV}(\confc^0)$, and
 \item a smooth map $\varphi_1\colon \left[0,\varepsilon\right[^{E(\CT_{\CV}(\confc^0))} \times W \times N(\confy^0) \to \ccompuptd{\vertsetv(\Gamma)}{\RR^3}$, 
\end{itemize}
such that the product of $ \varphi_1$ by the natural projection 
$p_{N(\confy^0)} \colon \left[0,\varepsilon\right[^{E(\CT_{\CV}(\confc^0))} \times W \times N(\confy^0) \to N(\confy^0)$
 has the following property. The map
\begin{equation*} \varphi_1 \times p_{N(\confy^0)} \colon \left[0,\varepsilon\right[^{E(\CT_{\CV}(\confc^0))} \times W \times N(\confy^0) \to \ccompuptd{\vertsetv(\Gamma)}{\RR^3} \times N(\confy^0)\end{equation*}
restricts to 
$\bigl(\setparamx(\CT_{\CV}(\confc^0)) \cap \left[0,\varepsilon\right[^{E(\CT_{\CV}(\confc^0))}\bigr) \times W \times N(\confy^0)$ as a bijection onto
an open neighborhood of $(\confc^0,\confy^0)$ in ${\CV}(\Gamma)$, where $\setparamx(\CT_{\CV}(\confc^0))$ was introduced in Definition~\ref{defXT}.
\end{lemma}
\bp According to Lemma~\ref{lemcnsonecybelongstonu}, $(\confc^0,\confy^0)$ has a neighborhood parametrized as in Lemma~\ref{lemstrucCV}. In such a parametrization, no $u_{\finsetd}$ is involved since $\confy^0=\confy_{\finsetv}^0 \in \cinjuptd{\finsetb}{\CC}$. 
Also note $\parentp^{\prime}_{\cutx} = \parentp_{\cutx,\finsetb}$. In particular, the elements of $\parentp^{\prime}_{\cutx}$ are pairwise disjoint sets.
Define the tree $\CT_{\CV}(\confc^0)$ associated to $(\confc^0,\confy^0)$ to be the tree obtained from $\CT(\parentp,\parentp_{\finsetb})$ by replacing the pair of edges respectively labeled by $\lambda(\finsety)$ and $\mu_{\finsety}$ by a single edge labeled by $\tilde{\lambda}(\finsety)=\lambda(\finsety)\mu_{\finsety}$ for each $\finsety \in \parntppx$. (In this tree, the right $\finsetb$-part in Figure~\ref{figTbraidsgen} is reduced to one edge labeled by $\lambda$, and the edges labeled by $\tilde{\lambda}(\finsety)$ for $\finsety \in \parntppx$  
start at univalent vertices.) The lemma is an easy consequence of Lemma~\ref{lemstrucCV}.
\eop

\begin{lemma} \label{lemfacesfiber}
 Let $\confy  \in \cinjuptd{\finsetb}{\CC}$. A codimension-one open face of ${\CV}(\confy,\Gamma)$ is a stratum $(\parentp,\parntppx)$ as in Lemma~\ref{lemdimstratCV} 
 such that
 \begin{itemize}
  \item either we have $\finsetv(\Gamma) \in \parntppx$, $\lambda \neq 0$, $\parntppx=\{\finsetv(\Gamma)\}$, and
$\parentp=\{\finsetv(\Gamma),\finsetc\}$ for some $\finsetc \subset \finsetv(\Gamma)$,
  \item or we have $\finsetv(\Gamma) \notin \parntppx$, $\lambda = 0$, and $\parentp=\{\finsetv(\Gamma)\} \cup \parntppx$ (where $\parntppx$ can be empty).
 \end{itemize}
 \end{lemma}
\bp According to Lemma~\ref{lemdimstratCV}, the faces of the $\partial{\CV}(\confy,\Gamma)$ with maximal dimension are such that $\cardbig{\parentp} =\cardbig{\parntppx} +1$, where $\finsetv(\Gamma) \in \parentp$.
\eop

Lemmas~\ref{lemcodimonex} and \ref{lemstrucCVinjtree} guarantee that ${\CV}(\confy,\Gamma)$ behaves as a codimension-one face of a manifold with boundary along a face as in Lemma~\ref{lemfacesfiber}.

\section{A one-form on \texorpdfstring{$\cinjuptd{\finsetb}{\CC}$}{SB(C)}}
\label{seconeform}

For a finite set $\finsetb$ of cardinality at least $2$, a configuration $\confy \in \ccompuptd{\finsetb}{\CC}$, and a Jacobi diagram $\Gamma$ on a disjoint union of lines $\RR_{\eltb}$ indexed by elements $\eltb$ of $\finsetb$ such that $p_{\finsetb} \colon U(\Gamma)\to \finsetb$ is onto, recall that ${\CV}(\confy,\Gamma)$ is the preimage of $\confy$ under
$p_{\ccompuptdempty{\finsetb}} \colon \CV(\Gamma) \to \ccompuptd{\finsetb}{\CC}.$
When $\confy \in \cinjuptd{\finsetb}{\CC}$, let $\check{\CV}(\confy,\Gamma)$ denote the quotient by vertical translations of $\check{C}(\RR^3,\confy^0 \times \RR;\Gamma)$ for a representative $\confy^0 \in  \check{C}_{\finsetb}\left[\drad{1}\right]$ of $\confy$.
Note that $\check{\CV}(\confy,\Gamma)$ is an open $T$-face of ${C}(\rats(\hcylc),\tanghcyll;\Gamma)$ for any tangle $\tanghcyll$ whose top configuration is $\confy$, as in Theorem~\ref{thmcomptang}, where the set $\finsetb$ of Theorem~\ref{thmcomptang} is empty, $I=\{j\}$, and $\parentp_s=\parentp_d=\parentp_x=\{\finsetv(\Gamma)\}$ with Notation~\ref{notationseparating}. 
Assume that $\Gamma$ is equipped with a vertex-orientation $o(\Gamma)$ as in Definition~\ref{defdia} and with an edge-orientation $o_E(\Gamma)$ of $H(\Gamma)$ as before Lemma~\ref{lemorc}. The space $\check{\CV}(\confy,\Gamma)$ is a smooth manifold of dimension $\cardlef{U(\Gamma)} + 3\cardlef{T(\Gamma)} -1$.
It is oriented as the part of the boundary of ${C}(\rats(\hcylc),\tanghcyll;\Gamma)$ that is the $T$-face in which all univalent vertices of $\Gamma$ tend to $\infty$ above $\hcylc$, for a tangle $\tanghcyll$ whose top configuration is $\confy$. 
(Here, $\Gamma$ is also viewed as a diagram on the domain $\sourcetl$ of $\tanghcyll$ by representing the original $[i_{\Gamma}]$ by a map $i_{\Gamma} \colon U(\Gamma) \hookrightarrow \left( \sqcup_{\eltb \in \finsetb} \left[1, \infty\right[_b \subset \sourcetl\right)$.)
Note that $\check{\CV}(\confy,\Gamma)$ is therefore oriented as the part of the boundary of $(-{C}(\rats(\hcylc),\tanghcyll;\Gamma))$ that is (minus) the $T$-face in which all univalent vertices of $\Gamma$ tend to $\infty$ below $\hcylc$, for a tangle $\tanghcyll$ whose bottom configuration is $\confy$.   The orientation of $\check{\CV}(\confy,\Gamma)$ depends on $o(\Gamma)$ and on $o_E(\Gamma)$, but it does not depend on the global orientations of the lines $\RR_{\eltb}$, which are only locally oriented by $o(\Gamma)$ near the images of the univalent vertices of $\Gamma$ as in Definition~\ref{defdia}.\footnote{The reader who prefers working with oriented strands can assume that the lines $\RR_{\eltb}$ are oriented from bottom to top and consider braids $L$ instead of tangles $L$, above, for the moment. However, since we will later need to allow various orientations for our strands $\RR_{\eltb}$, it is better to work with unoriented strands as much as possible.} 

Below, we define a one-form $\eta_{\Gamma}=\eta_{\Gamma,o(\Gamma)}$ on $\cinjuptd{\finsetb}{\CC}$ to be the integral of $\bigwedge_{e \in E(\Gamma)}p_{e,S^2}^{\ast}({\omega}_{S^2})$ along the interiors $\check{\CV}(\confy,\Gamma)$ of the compact fibers ${\CV}(\confy,\Gamma)$. 
We agree that the integral along the fiber of $dx \wedge \omega$ for a volume form $\omega$ of the fiber is $\bigl(\int_{\mbox{\scriptsize fiber}}\omega\bigr)dx$.
\begin{proposition}
 \label{propintfib}
 The integral of $\bigwedge_{e \in E(\Gamma)}p_{e,S^2}^{\ast}({\omega}_{S^2})$ along the interior $\check{\CV}(\confy,\Gamma)$ of the fiber ${\CV}(\confy,\Gamma)$ is absolutely convergent. It defines a smooth one-form $\eta_{\Gamma}$ on $\cinjuptd{\finsetb}{\CC}$.\footnote{Again, $\eta_{\Gamma}$ depends on the arbitrary vertex-orientation $o(\Gamma)$ of $\Gamma$, but the product $\eta_{\Gamma}\left[\Gamma\right]$ is independent of $o(\Gamma)$.}
 The definition of $\eta_{\Gamma}$ extends naturally to diagrams $\Gamma$
 on $\sqcup_{\eltb \in \finsetb}\RR_{\eltb}$ such that $p_{\finsetb} \colon U(\Gamma)\to \finsetb$ is not onto.\footnote{In general, $\eta_{\Gamma}$ pulls back through $\cinjuptd{p_{\finsetb}(U(\Gamma))}{\CC}$. So $\eta_{\Gamma}=0$ if $\cardlef{p_{\finsetb}(U(\Gamma))} <2$.} 

 Let $\gamma \colon \left[0,1\right] \to \cinjuptd{\finsetb}{\CC}$ be a smooth map. Orient  $p_{\ccompuptdempty{\finsetb}}^{-1}\left(\gamma(\left[0,1\right]\right))$ as the local product $\left[0,1\right] \times \mbox{fiber}$.\footnote{Note that it amounts to say that the $\left[0,1\right]$ factor replaces the upward vertical translation parameter of the quotient $\check{\CV}(\confy,\Gamma)$, as far as orientations are concerned.} 
 Then the integral 
 \begin{equation*}\int_{\left[0,1\right]} \gamma^{\ast}(\eta_{\Gamma})=\int_{p_{\ccompuptdempty{\finsetb}}^{-1}\left(\gamma(\left[0,1\right]\right))}\bigwedge_{e \in E(\Gamma)}p_{e,S^2}^{\ast}({\omega}_{S^2})\end{equation*}
 is absolutely convergent. The map
 $\left( t \mapsto \int_{\left[0,t\right]} \gamma^{\ast}(\eta_{\Gamma})\right)$ is differentiable, and we have
\begin{equation*}\frac{\partial}{\partial t} \left( \int_{\left[0,t\right]} \gamma^{\ast}(\eta_{\Gamma})\right)(u)=\eta_{\Gamma}\Bigl(\gamma(u),\frac{\partial}{\partial t}\gamma_u\Bigr).\end{equation*}
\end{proposition}
\bp In the proof, we assume that $p_{\finsetb}$ is onto, without loss of generality.
Lemma~\ref{lemstrucCVinjtree} implies that the integral $\int_{p_{\ccompuptdempty{\finsetb}}^{-1}\left(\gamma\left(\left[0,1\right]\right)\right)}\bigwedge_{e \in E(\Gamma)}p_{e,S^2}^{\ast}({\omega}_{S^2})$ is absolutely convergent. See the proof of Theorem~\ref{thmStokesonsetparamx}. Let us prove that the integral of $\bigwedge_{e \in E(\Gamma)}p_{e,S^2}^{\ast}({\omega}_{S^2})$ along the interior $\check{\CV}(\confy,\Gamma)$ of the fiber ${\CV}(\confy,\Gamma)$ is absolutely convergent and that it defines a smooth a one-form $\eta_{\Gamma}$ on $\cinjuptd{\finsetb}{\CC}$.

Let $\confy^0 \in \cinjuptd{\finsetb}{\CC}$, let $N(\confy^0)$ be small neighborhood of $\confy^0$ in $\cinjuptd{\finsetb}{\CC}$, and let \begin{equation*}\left(\zeta_r \colon N(\confy^0) \to \RR\right)_{r\in \underline{2 \cardlef{\finsetb} -3}}\end{equation*} be a system of coordinates on $N(\confy^0)$.
These coordinates give rise to associated smooth one-forms $d\zeta_r =d p_{N(\confy^0)} \circ \zeta_r$ on $\check{\CV}(N({\confy}^0),\Gamma)=\check{\CV}(\Gamma) \cap p_{\ccompuptdempty{\finsetb}}^{-1}(N(\confy^0))$.
For a local system $(f_1, \dots, f_k)$ of coordinates of the interior of a fiber $\check{\CV}({\confy},\Gamma)$, and a local product structure with the base, we also have associated forms $df_i$, which depend on the product structure. (Changing this product structure adds some combination of the $d\zeta_r$ and $df_j$ to $df_i$). We also have an associated volume form of the fiber
$\omega_F =df_1 \wedge \dots \wedge df_k$, which depends on the product structure too.
A $(k+1)$-form $\Omega$ on $\check{\CV}({N}({\confy}^0),\Gamma)$
may be expressed as 
$\Omega =\sum_{r=1}^{2 \cardlef{\finsetb} -3}  d \zeta_r\wedge (g_r\omega_F) + \sum_{i=1}^k\omega_{i,v}$,
where $\omega_{i,v}$ vanishes at the tangent vector $\xi_i$ to a curve of the fiber whose coordinates $f_j$ for $j \in \underline{k} \setminus \{i\}$ are constant ($\omega_{i,v}$ is expressed as a wedge product of coordinates forms that does not involve $df_i$, this decomposition is not canonical).
In order to check the convergence of the integral of the pull-back $\Omega$ of the form $\bigwedge_{e \in E(\Gamma)}p_{e,S^2}^{\ast}({\omega}_{S^2})$ on $\check{\CV}({N}({\confy}^0),\Gamma)$ along the fiber, it suffices to cover the fiber by finitely many neighborhoods as above, to express $\Omega$ as above with respect to the corresponding charts, and to check that the $g_r\omega_F$ (which are well-defined, up to forms whose integrals along the fiber vanish) and their derivatives with respect to the coordinates $\zeta_j$ are bounded in each of these neighborhoods.
Lemma~\ref{lemstrucCVinjtree} implies that $\Omega$ is the restriction to \begin{equation*}\left(\setparamx\bigl(\CT_{\CV}(\confc^0)\bigr) \cap \left[0,\varepsilon\right[^{E(\CT_{\CV}(\confc^0))}\right) \times W \times N(\confy^0)\end{equation*}
of a smooth form on $\left[0,\varepsilon\right[^{E(\CT_{\CV}(\confc^0))} \times W \times N(\confy^0)$, locally. The same holds for the $g_r\omega_F$ and all their iterated partial derivatives with respect to the $\zeta_{j}$.
These forms are bounded in a neighborhood of an arbitrary element $(\confc^0,\confy^0)$ of ${\CV}({\confy}^0,\Gamma)$ as in Lemma~\ref{lemstrucCVinjtree}.
Since the fiber ${\CV}({\confy}^0,\Gamma)$ is compact, the integral of $g_r\omega_F$ along the fiber is absolutely convergent, for any $r \in \underline{2 \cardlef{\finsetb} -3}$. So are the integrals of its iterated partial derivatives with respect to the $\zeta_{j}$ on $\check{\CV}(N({\confy}^0),\Gamma)$ for a small neighborhood $N({\confy}^0)$ of ${\confy}^0$ in $\cinjuptd{\finsetb}{\CC}$. This proves that $\eta_{\Gamma}$ is a well-defined smooth one-form on $\cinjuptd{\finsetb}{\CC}$. \eop

\begin{definition} \label{defetaholonomy}
For $k\in\NN \setminus \{0\}$, set \begin{equation*}\eta_{k,\finsetb}=\sum_{\Gamma \in \Davis^e_{k}(\sqcup_{\eltb \in \finsetb}\RR_{\eltb})}
\frac{(3k-\cardlef{E(\Gamma)})!}{(3k)!2^{\cardlef{E(\Gamma)}}}\eta_{\Gamma}\left[\Gamma\right] \in \Omega^1\bigl(\ccompuptd{\finsetb}{\CC}; \Aavis_k(\sqcup_{\eltb \in \finsetb}\RR_{\eltb})\bigr),\end{equation*}
where $\eta_{\Gamma}=0$ if $\cardlef{p_{\finsetb}(U(\Gamma))} <2$ or if $\Gamma$ is not connected,
and \begin{equation*}\eta_{\finsetb}=\sum_{k\in \NN \setminus \{0\}}\eta_{k,\finsetb}.\end{equation*} 
The form $\eta_{\finsetb}$ is a one-form on $\cinjuptd{\finsetb}{\CC}$ with coefficients in the space $\Aavis(\sqcup_{\eltb \in \finsetb}\RR_{\eltb})$ of Jacobi diagrams on $\sqcup_{\eltb \in \finsetb}\RR_{\eltb}$, which is treated as an unoriented manifold as in Definition~\ref{defrkoruniv}. The form $\eta_{\finsetb}$ will be regarded as a \emph{connection}. Let $p_{CS} \colon \check{C}_{\finsetb}\left[\drad{1}\right] \to \ccompuptd{\finsetb}{\CC}$ denote the natural projection.
For a path $\gamma \colon [a,b] \to \check{C}_{\finsetb}\left[\drad{1}\right]$, define the \indexT{holonomy} $\hol{\gamma}{\eta_{\finsetb}}$ of $\eta_{\finsetb}$ along $\gamma$ to be 
\begin{equation*}\hol{\gamma}{\eta_{\finsetb}}=\sum_{r=0}^{\infty} \int_{(t_1,\dots,t_r) \in [a,b]^r \suchthat  t_1 \leq t_2 \leq \dots \leq t_r}\bigwedge_{i=1}^r (p_{CS} \circ \gamma \circ p_i)^{\ast}(\eta_{\finsetb}),\end{equation*}
where $p_i(t_1,\dots,t_r)=t_i$, the wedge product of forms is performed as usual, and the diagrams are multiplied from bottom to top (from left to right) with respect to their order of appearance.
\end{definition}
 The degree $0$ part of $\hol{\gamma}{\eta_{\finsetb}}$ is the unit $\left[\emptyset\right]$ of $\Aavis_k(\sqcup_{\eltb \in \finsetb}\RR_{\eltb})$, and we have
\begin{equation*}\hol{\gamma}{\eta_{\finsetb}}=\left[\emptyset\right] + \sum_{r=1}^{\infty} \int_{(t_1,\dots,t_r) \in [a,b]^r \suchthat  t_1 \leq t_2 \leq \dots \leq t_r}\bigwedge_{i=1}^r (p_{CS} \circ \gamma \circ p_i)^{\ast}(\eta_{\finsetb}).\end{equation*}
This holonomy is valued in a space of diagrams on an unoriented domain as in  Definition~\ref{defrkoruniv}, Proposition~\ref{proporcomp}, and Remark~\ref{rkorunivtwo}. It satisfies the following properties.
\begin{itemize}
 \item For an orientation-preserving diffeomorphism $\psi \colon \left[c,d\right] \to [a,b]$, we have
 \begin{equation*}\hol{\gamma \circ \psi}{\eta_{\finsetb}}=\hol{\gamma}{\eta_{\finsetb}}.\end{equation*}
 \item When $\gamma_1\gamma_2$ is the path composition of $\gamma_1$ and $\gamma_2$, we have
\begin{equation*}\hol{\gamma_1\gamma_2}{\eta_{\finsetb}}=\hol{\gamma_1}{\eta_{\finsetb}}\hol{\gamma_2}{\eta_{\finsetb}}.\end{equation*}
\item We have $\frac{\partial}{\partial t}\hol{\gamma\vert_{\left[a,t\right]}}{\eta_{\finsetb}}=\hol{\gamma\vert_{\left[a,t\right]}}{\eta_{\finsetb}}\eta_{\finsetb}(\gamma^{\prime}(t))$.
\item  We have $\frac{\partial}{\partial t}\hol{\gamma\vert_{\left[t,b\right]}}{\eta_{\finsetb}}=-\eta_{\finsetb}(\gamma^{\prime}(t))\hol{\gamma\vert_{\left[t,b\right]}}{\eta_{\finsetb}}$.
\end{itemize}

The following lemma expresses the variation of the invariant $\Zinvuf$ of long tangles under isotopies that do not fix the bottom and top configurations. It uses the above holonomy and the anomaly $\alpha$ of Section~\ref{secanomalpha}.

\begin{lemma}
\label{lemholzinf}
Let $(h_t)_{t \in \left[0,1\right]}$ be an isotopy of $\crats(\hcylc)$ such that $h_t$ is the identity map on $\left(\CC \setminus \drad{1} \right)\times \RR$ for any $t$, $h_t$ may be expressed as $h^-_t \times \id_{\left]-\infty,0\right]}$ on $\CC \times \left]-\infty,0\right]$ for a planar isotopy $(h^-_t)_{t \in \left[0,1\right]}$, and $h_t$ may be expressed as $h^+_t \times \id_{\left[1,+\infty\right[}$ on $\CC \times \left[1,+\infty\right[$ for a planar isotopy $(h^+_t)_{t \in \left[0,1\right]}$. Assume that $h_0=\id$ and note that $h_t$ preserves $\hcylc$ setwise. 
Let $\tanghcyll$ be a long tangle representative of $\crats(\hcylc)$ whose bottom (resp. top) configuration is represented by a map $\confy^- \colon \finsetb^- \to \drad{1}$ (resp. $\confy^+ \colon \finsetb^+ \to \drad{1}$).
Let $J_{bb}$ denote the set of components of $\tanghcyll$ going from bottom to bottom, and let $J_{tt}$ denote the set of components of $\tanghcyll$ going from top to top. Set $\varepsilon(K_j)=-$ for $K_j \in J_{bb}$, and $\varepsilon(K_j)=+$ for $K_j \in J_{tt}$. 
For a component $K_j$ of $J_{bb} \cup J_{tt}$, let $\theta_j \colon \left[0,1\right] \to \RR$ be a path such that the difference $\bigl(h^{\varepsilon(K_j)}_t(\confy^{\varepsilon(K_j)}(K_j(1)))-h^{\varepsilon(K_j)}_t(\confy^{\varepsilon(K_j)}(K_j(0)))\bigr)$ is a positive multiple of the complex direction $\exp(i 2\pi\theta_j(t))$.
With the notation of Theorem~\ref{thmfstconsttang}, set $\Zinvuf(t)=\Zinvuf(\hcylc,h_t(\tanghcyll))$. Then we have
\begin{multline*}\left(\prod_{K_j \in J_{bb} \cup J_{tt}}\biggl(\exp\Bigl(-2\varepsilon(K_j)\bigl(\theta_j(t)-\theta_j(0)\bigr)\alpha\Bigr) \#_j\biggr)\right)\Zinvuf(t)\\=
\hol{h^-\vert_{\left[t,0\right]}\circ \confy^-}{\eta_{\finsetb^-}}\Zinvuf(0)\hol{h^+\vert_{\left[0,t\right]}\circ \confy^+}{\eta_{\finsetb^+}} .\end{multline*}
\end{lemma}
\bp Let $\tau$ be a parallelization of $\hcylc$. Set $\tanghcyll=(K_j)_
{j \in \underline{k}}$, and recall \begin{equation*}\Zinvuf(t)=\exp\Bigl(-\frac14 p_1(\tau)\ansothree\Bigr)\left(\prod_{j=1}^k\left(\exp\Bigl(-I_{\theta}\bigl(K_j(t),\tau\bigr)\alpha\Bigr) \#_j\right)\right)\Zinv\bigl(\hcylc,h_t(\tanghcyll),\tau\bigr).\end{equation*}
The algebraic boundary of the chain $\cup_{t\in \left[t_0,t_1\right]}{C}(\rats(\hcylc),h_t(\tanghcyll);\Gamma)$ is
\begin{equation*}{C}\bigl(\rats(\hcylc),h_{t_1}(\tanghcyll);\Gamma\bigr)-{C}\bigl(\rats(\hcylc),h_{t_0}(\tanghcyll);\Gamma\bigr)- \sum\left(\cup_{t\in \left[t_0,t_1\right]}F_t\right),\end{equation*} where the sum runs over the codimension-one faces $F_t$ of ${C}(\rats(\hcylc),h_{t}(\tanghcyll);\Gamma)$.\footnote{This chain is locally modelled by open subsets of ${C}(\rats(\hcylc),h_t(\tanghcyll);\Gamma) \times \left]t-\varepsilon,t+\varepsilon\right[$ unless the isotopies $h_t^{\pm}$ are degenerate.
See Theorem~\ref{thmcomptang} and Lemma~\ref{lemcodimonex}. Stokes' theorem applies thanks to Theorem~\ref{thmStokesonsetparamx}.}
Faces cancel as in Section~\ref{secvarzinf} except for the anomaly $\alpha$ faces and the 
faces for which some vertices are at $\infty$.

The variations due to the anomaly $\alpha$ faces contribute to $\frac{\partial}{\partial t}\Zinvlink(\hcylc,h_t(\tanghcyll),\tau)$
as \begin{equation*}\left(\sum_{j=1}^k\frac{\partial}{\partial t}\left( 2 \int_{\left[0,t\right]\times \ST^+K_j}p_{\tau}^{\ast}(\omega_{S^2})\right)\alpha \#_j\right)\Zinvlink\bigl(\crats,h_t(\tanghcyll),\tau\bigr)\end{equation*}
as in Lemma~\ref{lemvartauanom}.

When the bottom and top configurations are not fixed and when $K_j \in J_{bb} \cup J_{tt}$, we have 
\begin{equation*}I_{\theta}\bigl(K_j(u),\tau\bigr) - I_{\theta}\bigl(K_j(0),\tau\bigr) = 2 \int_{\cup_{t \in \left[0,u\right]} p_{\tau}(\ST^+K_j(t))\cup S(K_j(t))}\omega_{S^2},\end{equation*}
as in the proof of Lemma~\ref{lemvarithetatangtwo},
where $S(K_j(t))$ denotes the half-circle from $\varepsilon(K_j)\upvec$ to $-\varepsilon(K_j)\upvec$ through $\exp(2i \pi\theta_j(t))$ as in Lemma~\ref{lemvarithetatang}. So we get
\begin{equation*}I_{\theta}\bigl(K_j(t),\tau\bigr) - I_{\theta}\bigl(K_j(0),\tau\bigr) = 2 \int_{\cup_{u \in \left[0,t\right]} p_{\tau}(\ST^+K_j(u))}\omega_{S^2} - 2 \varepsilon(K_j) \bigl(\theta_j(t)-\theta_j(0)\bigr).\end{equation*}
Thus
\begin{equation*}\tilde{\Zinvuf}(t)=\left(\prod_{K_j \in J_{bb} \cup J_{tt}}\biggl(\exp\Bigl(-2\varepsilon(K_j)\bigl(\theta_j(t)-\theta_j(0)\bigr)\alpha \Bigr) \#_j\biggr)\right)\Zinvuf(t)\end{equation*}
gets no variation from the anomaly $\alpha$ faces as in Corollary~\ref{corvartauanom}.

Let
$F_t$ be a codimension-one face of ${C}(\rats(\hcylc),h_{t}(\tanghcyll);\Gamma)$ for which a subset $\finsetv$ of $\finsetv(\Gamma)$ is mapped to $\infty$. Such a face is either a $T$-face as in Theorem~\ref{thmcomptang} and Lemma~\ref{lemsubsecTface},
or a face $\facee_{\infty}(\finsetv,\Link,\Gamma)$ as around Notation~\ref{notsinjupdtcs}. Set $F=\cup_{t\in \left[t_0,t_1\right]}F_t$.

In both cases, an element $\confc^0$ of the face involves an injective configuration $T_{0}\phi_{\infty} \circ f_1^0$ from the kids of $\finsetv$ to $(T_{\infty}\rats(\hcylc) \setminus 0)$ up to dilation. 
Let $e=(v_1,v_2)$ be an edge.
If the vertices $v_1$ and $v_2$ are in different kids of $\finsetv$, then we have
\begin{equation*}\begin{array}{ll}p_{\tau}\circ p_e(\confc^0)&=\frac{\phi_{\infty} \circ f_1^0(v_2) -\phi_{\infty} \circ f_1^0(v_1)}{\norm{ \phi_{\infty} \circ f_1^0(v_2) -\phi_{\infty} \circ f_1^0(v_1) }}\\&=\frac{\norm{f_1^0(v_1)}^2f_1^0(v_2)- \norm{f_1^0(v_2)}^2f_1^0(v_1)}{\bignorm{\norm{f_1^0(v_1)}^2f_1^0(v_2)- \norm{f_1^0(v_2)}^2f_1^0(v_1)}}.\end{array}\end{equation*}
When $v_1 \in V$ and $v_2 \notin V$, then we have \begin{equation*}p_{\tau}\circ p_e(\confc^0)=-\frac{f_1^0(v_1)}{\norm{f_1^0(v_1)}}.\end{equation*} If  $v_2\in V$ and $v_1 \notin V$, we similarly have \begin{equation*}p_{\tau}\circ p_e(\confc^0)=\frac{f_1^0(v_2)}{\norm{f_1^0(v_2)}}.\end{equation*}

Let $E_{\infty}$ be the set of edges between elements of the set $\finsetv$ of vertices mapped to $\infty$ in $F$, and let $E_m$ denote the set of edges with one end in $\finsetv$.
The face $F_t$ is diffeomorphic to a product by $\check{C}_{\vertsetv(\Gamma) \setminus \finsetv}(\crats(\hcylc),h_t(\tanghcyll);\Gamma)$, whose dimension is
\begin{equation*}3 \cardbig{T(\Gamma) \cap (\vertsetv(\Gamma) \setminus \finsetv)} + \cardbig{U(\Gamma) \cap (\vertsetv(\Gamma) \setminus \finsetv)},\end{equation*}
of a space $C_{\finsetv,t}$ of dimension 
\begin{equation*}3 \cardbig{T(\Gamma) \cap \finsetv} + \cardbig{U(\Gamma) \cap \finsetv} -1 = 2\cardbig{E_{\infty}}+\cardbig{E_{m}}-1,\end{equation*}
and $\bigwedge_{e \in E_{\infty} \cup E_m}p_e^{\ast}(\omega(j_E(e)))$ has to be integrated along $\cup_{t \in \left[0,1\right]}C_{\finsetv,t}$, according to the above expression of $p_{\tau}\circ p_e$ for edges of $E_{\infty} \cup E_m$.
The degree of this form is $2\cardlef{E_{\infty} \cup E_m}$. So the face $F$ cannot contribute unless $E_m=\emptyset$.

Now the expression $p_{\tau}\circ p_e$ for edges of $E_{\infty}$ makes also clear that, if $f_1^0$ is changed to $\phi_{\infty}\circ T\circ \phi_{\infty} \circ f_1^0$, for a vertical translation $T$ such that $0$ is not in the image of $\phi_{\infty}\circ T\circ \phi_{\infty} \circ f_1^0$, then the image under $\prod_{e \in E_{\infty}}p_e$ is unchanged. So we have a one-parameter group acting on our face $F$ such that $\prod_{e \in E_{\infty}}p_e$  factors through this action. Unless $\finsetv^+$ has only one kid, this action is not trivial, and the quotient of the face by this action is of dimension strictly less than the face dimension.

Therefore, for the contributing faces, we have $E_m=\emptyset$ and $\finsetv^+$ has only one kid.
Thus, according to Lemma~\ref{lemsubsecTface}, we are left with the $T$-faces of Theorem~\ref{thmcomptang} (for which $\finsetb=\emptyset$ and $I$ has one element) for which $\parentp_x=\{\finsetv\}$. These faces yield the derivative $\frac{\partial}{\partial t}\tilde{\Zinvuf}=d\tilde{\Zinvuf} \left(\frac{\partial}{\partial t}\right)$ with
\begin{equation*}d \tilde{\Zinvuf}=
-\bigl((t \mapsto h^-_t\circ\confy^-)^{\ast}({\eta_{\finsetb^-}})\bigr)\tilde{\Zinvuf} + \tilde{\Zinvuf}\bigl((t \mapsto h^+_t\circ\confy^+)^{\ast}({\eta_{\finsetb^+}})\bigr).\end{equation*}
This proves the equality 
\begin{equation*}\tilde{\Zinvuf}(t)=\hol{h^-\vert_{\left[t,0\right]}\circ \confy^-}{\eta_{\finsetb^-}}\tilde{\Zinvuf}(0)\hol{h^+\vert_{\left[0,t\right]}\circ \confy^+}{\eta_{\finsetb^+}}\end{equation*}
and leads to the formula for $\Zinvuf$.\eop

\begin{corollary} \label{corholZinvufrf}
Under the hypotheses of Lemma~\ref{lemholzinf},
let $\tanghcyll_{\parallel}$ be a parallel of $\tanghcyll$.
Let $h_t(\tanghcyll)_{\parallel}=(h_t(K_j)_{\parallel})_j$ be the parallel of $h_t(\tanghcyll)$ such that
\begin{equation*}lk\bigl(h_t(K_j),h_t(K_j)_{\parallel}\bigr)-lk\bigl(K_j,K_{j\parallel}\bigr)=
-\varepsilon(K_j)2\bigl(\theta_j(t)-\theta_j(0)\bigr)\end{equation*} for any component $K_j$ of $J_{bb} \cup J_{tt}$,
and $lk(h_t(K),h_t(K)_{\parallel})=lk(K,K_{\parallel})$ for any other component $K$ of $\tanghcyll$.
Use Definition~\ref{deffstconsttangframed} to set \begin{equation*}\Zinvufrfneg(t)=\Zinvufrfneg\bigl(\hcylc,h_t(\tanghcyll),h_t(\tanghcyll)_{\parallel}\bigr).\end{equation*} 
Then we have \begin{equation*}\Zinvufrfneg(t)=\hol{h^-\vert_{\left[t,0\right]}\circ \confy^-}{\eta_{\finsetb^-}}\Zinvufrfneg(0)\hol{h^+\vert_{\left[0,t\right]}\circ \confy^+}{\eta_{\finsetb^+}} .\footnote{Again, the holonomies are considered as valued in spaces of diagrams on unoriented domains, where the vertex-orientation of Jacobi diagrams includes local orientations of strands, which can be made consistent with a global orientation induced by $\tanghcyll$, as in Definition~\ref{defrkoruniv}.}\end{equation*}
\end{corollary}
\bp See Definition~\ref{defselflkqtanggen}. \eop

A connection is \emph{flat} if its holonomy along a null-homotopic loop is trivial.

\begin{proposition} \label{propfuncbraid}
The connection $\eta_{\finsetb}$ is flat on $\cinjuptd{\finsetb}{\CC}$.
When $\gamma \colon \left[0,1\right] \to \check{C}_{\finsetb}\left[\drad{1}\right]$ is smooth with vanishing derivatives at $0$ and $1$, the image $T(\gamma)$ of the graph $\{(\gamma(t),t)\}$ of $\gamma$ in $\drad{1}\times\left[0,1\right]$
 is a tangle in $\drad{1}\times\left[0,1\right]$, and we have
\begin{equation*}\Zinvufrfneg\bigl(T(\gamma)\bigr)=\Zinvuf\bigl(T(\gamma)\bigr)=\hol{p_{CS} \circ \gamma}{\eta_{\finsetb}},\end{equation*}
where $p_{CS} \circ \gamma$ is the composition of $\gamma$ by the natural projection $p_{CS} \colon \check{C}_{\finsetb}\left[\drad{1}\right] \to \ccompuptd{\finsetb}{\CC}$.
For two framed tangles $(\hcylc_1,\tanghcyll_1)$ and $(\hcylc_2,\tanghcyll_2)$ such that the bottom of $\tanghcyll_2$ coincides with the top of $\tanghcyll_1$, if one of them is a braid $T(\gamma)$ as above, we have \begin{equation*}\Zinvufrfneg(\hcylc_1\hcylc_2,\tanghcyll_1\tanghcyll_2) = \Zinvufrfneg(\hcylc_1,\tanghcyll_1) \Zinvufrfneg(\hcylc_2,\tanghcyll_2),\end{equation*} with products obtained by stacking above in natural ways on both sides, reading from left to right.
\end{proposition}
\bp Applying Lemma~\ref{lemholzinf} when $\tanghcyll$ is a trivial braid, $h^-_t=h^-_0$ is constant, and $h^+_t\circ\confy^+=\gamma(t)$ shows $\Zinvuf\left(T(\gamma)\right)=\hol{p_{CS} \circ \gamma}{\eta_{\finsetb}}$.
Then the isotopy invariance of $\Zinvuf$ shows that $\eta_{\finsetb}$ is flat on $\cinjuptd{\finsetb}{\CC}$.
Applying Corollary~\ref{corholZinvufrf} when $h^-_t=h^-_0$ is constant and $\gamma(t)=h^+_t\circ\confy^+$ proves
\begin{equation*}\Zinvufrfneg\bigl(\hcylc,\tanghcyll T(\gamma)\bigr) = \Zinvufrfneg\bigl(\hcylc,\tanghcyll\bigr) \Zinvufrfneg\bigl(T(\gamma)\bigr).\end{equation*}
So we have $\Zinvufrfneg\left(T(\gamma_1)T(\gamma_2)\right) = \Zinvufrfneg\left(T(\gamma_1)\right) \Zinvufrfneg\left(T(\gamma_2)\right)$ for braids.

Applying Corollary~\ref{corholZinvufrf} when $h^+_t=h^+_0$ is constant and $\gamma(t)=h^-_{1-t}\circ\confy^-$ proves $\Zinvufrfneg(\hcylc_1\hcylc_2,\tanghcyll_1\tanghcyll_2) = \Zinvufrfneg(\hcylc_1,\tanghcyll_1) \Zinvufrfneg(\hcylc_2,\tanghcyll_2)$ when $(\hcylc_1,\tanghcyll_1)$ is a braid and $(\hcylc_2,\tanghcyll_2)$ is a framed tangle, too.
\eop

\begin{note} \label{noteKZ}
Recall the space $\check{C}_{\finsetb}\!\left[\CC\right]$ of injective planar configurations $\confy \colon \finsetb \hookrightarrow \CC$.
Let $\gamma \colon [0,1] \to \check{C}_{\finsetb}\!\left[\CC\right]$ be a braid.
Proposition~\ref{propfuncbraid} expresses $\Zinvuf(\gamma)=\Zinvuf\left(T(\gamma)\right)$ as the holonomy of the flat connection $\eta_{\finsetb}$ with coefficients in $\Aavis(\sqcup_{\eltb \in \finsetb}\RR_{\eltb})$ along $p_{CS} \circ \gamma$. 
Kontsevich defined his integral $Z^K$ as the holonomy of another flat connection $\eta_{KZ}$ along $\gamma$ for such a braid \cite[Section 4]{barnatan}. 
The involved \emph{Knizhnik–Zamolodchikov connection} $\eta_{KZ}$ is defined on the space $\check{C}_{\finsetb}\!\left[\CC\right]$. The coefficients of $\eta_{KZ}$ belong to the tensor product of the same diagram space $\Aavis(\sqcup_{\eltb \in \finsetb}\RR_{\eltb})$ by $\CC$.

Assume $\finsetb =\{b_1,b_2\}$.
Let $\gamma^c \colon \left]0,1\right] \to \check{C}_{\finsetb}\!\left[\CC\right]$ map
$t$ to $\gamma^c(t)$ with $\gamma^c(t)(b_1)=0$ and $\gamma^c(t)(b_2)=t$.
For $\varepsilon \in ]0,1[$, set $\gamma^c_{\varepsilon}=\gamma^c\vert_{[\varepsilon,1]}$.
Then the holonomy of $\eta_{\finsetb}$ along $\gamma^c_{\varepsilon}$ is trivial. However, the holonomy of $\eta_{KZ}$ along $\gamma^c_{\varepsilon}$ is not. It does not even converge when $\varepsilon$ tends to zero.
In particular, the connection $\eta_{KZ}$ does not factor through $p_{CS}$.

Nevertheless, the holonomy of $\eta_{KZ}$ can be regularized. It lead Thang L\^e and Jun Murakami to a definition of the Kontsevich integral for combinatorial $q$-braids $\gamma$ in \cite{LeMur}.  This definition is related to the limit $\Zinvufrfneg(T(\gamma))$ of the holonomy of $\eta_{\finsetb}$ as in Note~\ref{notelesunikonfunc}.
\end{note}

In the above proposition, we proved that $\eta_{\finsetb}$ is flat by proving that its holonomy is $\id=\left[\emptyset\right]$ on null-homotopic loops.
The flatness of a differentiable connection $\eta$ is often established by proving that its \emph{curvature} $(d\eta + \eta \wedge \eta)$ vanishes, instead. The following lemma shows how the curvature vanishing implies the homotopy invariance of the holonomy.

\begin{lemma} \label{lemrkflat} Set $\simplexr=\{(t_1,\dots,t_r) \in \left[0,1\right]^r \suchthat  t_1 \leq t_2 \leq \dots \leq t_r\}$.
Let $\eta$ be a one-form on $\cinjuptd{\finsetb}{\CC}$ with coefficients in the space $\Aavis(\sqcup_{\eltb \in \finsetb}\RR_{\eltb})$.
Let $\gamma \colon \left[0,1\right] \times \left[0,1\right] \rightarrow \cinjuptd{\finsetb}{\CC}$ be a homotopy mapping $(\left[0,1\right] \times \{0,1\}) \cup \{0\} \times \left[0,1\right]$ to a point. 
Set $\gamma_u(t)=\gamma(u,t)$. Then the holonomy of $\eta$ along $\gamma_1$ is \begin{equation*}\id +\sum_{r=1}^{\infty} \int_{\left[0,1\right]\times\simplexr}\left(\sum_{j=1}^{r}(-1)^{j-1} \bigwedge_{i=1}^r (\gamma \circ p_i)^{\ast}(\eta)\left(\frac{(\gamma \circ p_j)^{\ast}(d\eta + \eta \wedge \eta)}{(\gamma \circ p_j)^{\ast}(\eta)} \right)\right),\end{equation*}
where the fraction means that the denominator is replaced by the numerator in the preceding expression.

In particular, if $d\eta + \eta \wedge \eta$ vanishes, then the holonomy of $\eta$ is trivial along any null-homotopic loop.
\end{lemma}
\bp
Stokes' theorem allows us to compute
\begin{equation*}\hol{\gamma_1}{\eta}=\id +\sum_{r=1}^{\infty} \int_{\simplexr}\bigwedge_{i=1}^r (\gamma_1 \circ p_i)^{\ast}(\eta)\end{equation*}
by integrating $d\bigwedge_{i=1}^r (\gamma_u \circ p_i)^{\ast}(\eta)$ over $\left[0,1\right] \times \simplexr$ as follows. 
Set $F_0(\simplexr)=\{(0,t_2,\dots,t_r) \in \Delta^r \}$, $F_r(\Delta^r)=\{(t_1,t_2,\dots,t_{r-1},1) \in \Delta^r \}$, and
\begin{equation*}F_j(\Delta^r)=\bigl\{(t_1,\dots ,t_{r}) \in \Delta^r  \suchthat t_j=t_{j+1}\bigr\}.\end{equation*} for $j \in \underline{r-1}$. Then we have
\begin{equation*}\partial \simplexr=\sum_{j=0}^r (-1)^{j+1} F_j(\simplexr),\end{equation*}
and
\begin{equation*}\partial\left(\left[0,1\right] \times \simplexr\right)= \Bigl(\bigl(\partial \left[0,1\right]\bigr) \times \simplexr \Bigr) \cup \Bigl(\left[0,1\right] \times \partial \bigl(-\simplexr\bigr)\Bigr).\end{equation*}
So we get
\begin{multline*}\int_{\simplexr}\bigwedge_{i=1}^r (\gamma_1 \circ p_i)^{\ast}(\eta)-
\int_{\simplexr}\bigwedge_{i=1}^r (\gamma_0 \circ p_i)^{\ast}(\eta) \\+ \sum_{j=0}^r(-1)^{j}\int_{\left[0,1\right]\times F_j(\simplexr)}\bigwedge_{i=1}^r (\gamma \circ p_i)^{\ast}(\eta) = \int_{\left[0,1\right]\times \simplexr }d\left(\bigwedge_{i=1}^r (\gamma \circ p_i)^{\ast}(\eta)\right),\end{multline*}
where the faces $F_0$ and $F_r$ do not contribute since $\gamma$ maps $(\left[0,1\right] \times \{0,1\})$ to a point, and $\int_{\simplexr}\bigwedge_{i=1}^r (\gamma_0 \circ p_i)^{\ast}(\eta)$ similarly vanishes.
We obtain \begin{multline*}\int_{\simplexr}\bigwedge_{i=1}^r (\gamma_1 \circ p_i)^{\ast}(\eta)
\\= \sum_{j=1}^{r-1}(-1)^{j-1}\int_{\left[0,1\right]\times \simplexrm}\bigwedge_{i=1}^{r-1} (\gamma \circ p_i)^{\ast}(\eta)\left(\frac{(\gamma \circ p_j)^{\ast}(\eta \wedge \eta)}{(\gamma \circ p_j)^{\ast}(\eta)} \right) \\+\sum_{j=1}^{r}(-1)^{j-1} \int_{\left[0,1\right]\times \simplexr }\left(\bigwedge_{i=1}^r (\gamma \circ p_i)^{\ast}(\eta)\right)\left(\frac
{(\gamma \circ p_j)^{\ast}(d\eta)}{(\gamma \circ p_j)^{\ast}(\eta)} \right),\end{multline*}
which yields the result.
\eop

\begin{corollary}
\label{coretaBflat} We have
 $d\eta_{\finsetb} + \eta_{\finsetb} \wedge \eta_{\finsetb}=0$.
\end{corollary}
\bp Let us prove that the degree $k$ part $(d\eta_{\finsetb} + \eta_{\finsetb} \wedge \eta_{\finsetb})_k$ of $(d\eta_{\finsetb} + \eta_{\finsetb} \wedge \eta_{\finsetb})$ vanishes for any $k\in \NN$, by induction on the degree $k$. This is obviously true for $k=0$. Let us assume that $k>0$ and $(d\eta_{\finsetb} + \eta_{\finsetb} \wedge \eta_{\finsetb})_i$ vanishes for $i<k$.
For any disk $D=\gamma(\left[0,1\right]\times \Delta^{(1)})$ as in Lemma~\ref{lemrkflat}, the degree $k$ part of the holonomy of $\eta_{\finsetb}$ along $\partial D$, which vanishes, is the integral of 
$(d\eta_{\finsetb} + \eta_{\finsetb} \wedge \eta_{\finsetb})_k$ along $D$ according to Lemma~\ref{lemrkflat}. Therefore, the integral of 
$(d\eta_{\finsetb} + \eta_{\finsetb} \wedge \eta_{\finsetb})_k$ vanishes along any disk $D$, and $(d\eta_{\finsetb} + \eta_{\finsetb} \wedge \eta_{\finsetb})_k$ is zero.
\eop

Below, we compute $d\eta_{\finsetb}$ and sketch an alternative proof for the equality of 
Corollary~\ref{coretaBflat}.

\begin{lemma}\label{lemdetab}Along the open codimension-one faces of $\,{\CV}(\confy,\Gamma)$,
the integral of $\bigl(-\bigwedge_{e \in E(\Gamma)}p_{e,S^2}^{\ast}({\omega}_{S^2})\bigr)$ is absolutely convergent. This allows us to define the smooth two-form $(\confy \mapsto d\eta_{\Gamma}(\confy))$ on $\cinjuptd{\finsetb}{\CC}$ so that $d\eta_{\Gamma}(\confy)$ is the sum of these integrals along the open codimension-one faces of $\,{\CV}(\confy,\Gamma)$.
\end{lemma}
\bp The integral of $\bigwedge_{e \in E(\Gamma)}p_{e,S^2}^{\ast}({\omega}_{S^2})$ along these codimension-one faces of ${\CV}(\confy,\Gamma)$ is absolutely convergent and it defines a smooth two-form on $\cinjuptd{\finsetb}{\CC}$ as in the proof of Proposition~\ref{propintfib}.

To see that this two-form is $d\eta_{\Gamma}$, use a chart $\psi \colon \RR^s \to N(\confy)$ of a neighborhood of $\confy$ in $\cinjuptd{\finsetb}{\CC}$, with $s =2\cardlef{\finsetb} -3$. Let $\zeta_i$ denote the composition $p_i \circ \psi^{-1}$. So we have $\eta_{\Gamma}=\sum_{i=1}^s \eta_i d\zeta_i$ and \begin{equation*}d\eta_{\Gamma}=\sum_{(i,j) \in \underline{s}^2 \suchthat i<j}\left(\frac{\partial}{\partial \zeta_i} \eta_j - \frac{\partial}{\partial \zeta_j} \eta_i\right)d\zeta_i \wedge d\zeta_j,\end{equation*} 
where \begin{equation*}\frac{\partial}{\partial \zeta_2} \eta_1 = \lim_{\stackrel{t_2 \to 0}{t_2\in \left]0,\infty\right[}} \frac{1}{t_2} \left(\lim_{\stackrel{t_1 \to 0}{t_1\in \left]0,\infty\right[}}
\frac{1}{t_1} \int_{\psi(\left[0,t_1\right] \times \{(t_2,0,\dots,0)\})-\psi(\left[0,t_1\right] \times \{(0,0,\dots,0)\}) }\eta_{\Gamma}\right).\end{equation*}
We thus have
\begin{equation*}\frac{\partial}{\partial \zeta_1} \eta_2 - \frac{\partial}{\partial \zeta_2} \eta_1=\lim_{\stackrel{(t_1, t_2) \to 0}{(t_1,t_2) \in \left]0,\infty\right[^2} } \frac{1}{t_1t_2}\int_{\partial N(t_1,t_2)}\eta_{\Gamma},\end{equation*}
where $N(t_1,t_2)=\psi(\left[0,t_1\right]\times \left[0,t_2\right] \times \{0\})$. Furthermore, we have
\begin{equation*}\int_{\partial N(t_1,t_2)}\eta_{\Gamma}=
\int_{\partial p_{\ccompuptdempty{\finsetb}}^{-1}( N(t_1,t_2)) \setminus \left(\cup_{\confy \in N(t_1,t_2)} \partial {\CV}(\confy,\Gamma)\right)}\bigwedge_{e \in E(\Gamma)}p_{e,S^2}^{\ast}({\omega}_{S^2}).\end{equation*}
Lemma~\ref{lemstrucCVinjtree} and Theorem~\ref{thmStokesonsetparamx} imply that Stokes' theorem applies to the closed form $\bigwedge_{e \in E(\Gamma)}p_{e,S^2}^{\ast}({\omega}_{S^2})$. So we have $\int_{\partial p_{\ccompuptdempty{\finsetb}}^{-1}( N(t_1,t_2))}\bigwedge_{e \in E(\Gamma)}p_{e,S^2}^{\ast}({\omega}_{S^2})=0$ and \begin{equation*}\int_{ N(t_1,t_2)}d\eta_{\Gamma}= -\int_{y \in  N(t_1,t_2)}\int_{\partial {\CV}(\confy,\Gamma)}\bigwedge_{e \in E(\Gamma)}p_{e,S^2}^{\ast}({\omega}_{S^2}).\end{equation*}
\eop

\begin{remark}
\label{rkflatdirect}
The flatness condition $d\eta_{\finsetb}+\eta_{\finsetb}\wedge\eta_{\finsetb}=0$ can be proved from the definition of $d\eta_{\finsetb}$, which can be extracted from Definition~\ref{defetaholonomy} and Lemma~\ref{lemdetab}, as an exercise along the following lines. 

Let us use Lemma~\ref{lemdetab} to compute $d\eta_{\finsetb}$. Let $\confy \in \cinjuptd{\finsetb}{\CC}$.
According to Lemma~\ref{lemfacesfiber}, for a connected Jacobi diagram $\Gamma \in \Davis^e_{k}(\sqcup_{\eltb \in \finsetb}\RR_{\eltb})$, the faces of the $\partial{\CV}(\confy,\Gamma)$ with maximal dimension are 
of two types. Either $\finsetv(\Gamma) \notin \parntppx$, or $\finsetv(\Gamma) \in \parntppx$ and
$\parentp=\{\finsetv(\Gamma),\finsetc\}$. 

If $\finsetv(\Gamma) \notin \parntppx$, then we have $\parentp=\{\finsetv(\Gamma)\} \sqcup \parntppx$, and the configuration of the kids of $\finsetv(\Gamma)$ maps the univalent kids of $\finsetv(\Gamma)$ to the vertical line through the origin. It is defined up to vertical translation and dilation. Since $\Gamma$ is connected, if a univalent daughter $\finsety$ of $\finsetv(\Gamma)$ contains a trivalent vertex, then it contains a trivalent vertex that is bivalent, univalent, or $0$-valent in $\Gamma_{\finsety}$. This type of face does not contribute to $d\eta_{\finsetb}$ (as in Sections~\ref{seccancelfacessk} and~\ref{seccancelfaces}). (Note that the univalent daughters of $\finsetv(\Gamma)$ are in $\parntppx$. So they must have vertices on at least two strands, and $\Gamma_{\finsety}$ cannot be an edge between a univalent vertex and a trivalent one.)
For the faces in which all the daughters $\finsety$ of $\finsetv(\Gamma)$ contain only univalent vertices, the integrated form is determined by the configuration space of the kids of $\finsetv(\Gamma)$ up to dilation and (conjugates of) vertical translations since $\Gamma_{\finsety}$ does not contain chords. The dimension of this space is smaller than $3\cardlef{T(\Gamma)} + \cardlef{U(\Gamma)} -2$ if $\parntppx$ is not empty.
If $\parntppx$ is empty, then the face is independent of the planar configuration $\confy$. So it does not contribute to 
$d\eta_{\finsetb}$ either.

Thus, the only faces of the $\partial{\CV}(\confy,\Gamma)$ contributing to $d\eta_{\finsetb}$ are 
such that $\parentp=\{\finsetv(\Gamma),\finsetc\}$ and the univalent vertices of $\finsetc$ are on one strand.
Most of these faces cancel as an analysis similar to that performed in Section~\ref{seccancelfaces} shows. The only contributing faces are
the STU-faces that involve a collapse of an edge containing one univalent vertex such that the other two diagrams involved in the corresponding STU relation are not connected. Then these other two diagrams are disjoint unions of two components $\Gamma_1$ and $\Gamma_2$ on $\sqcup_{\eltb \in \finsetb}\RR_{\eltb}$. (If the involved diagrams were connected, then the corresponding faces would cancel by STU.) Consider a configuration $\confc_1$ of $\Gamma_1$ on $\sqcup_{\eltb \in \finsetb} \confy(\eltb) \times \RR_{\eltb}$
and a configuration $\confc_2$ of $\Gamma_2$ on $\sqcup_{\eltb \in \finsetb} \confy(\eltb) \times  \RR_{\eltb}$ for a planar configuration $\confy$. View $\Gamma_1$ far below $\Gamma_2$, and slide it vertically until it is far above.
During this sliding, there will be heights at which one univalent vertex of $\Gamma_1$ coincides with one univalent vertex of $\Gamma_2$.
(For a generic pair $(\confc_1,\confc_2)$, there are no heights at which more than one univalent vertex of $\Gamma_1$ coincides with one univalent vertex of $\Gamma_2$.)
Each collision corresponds to a configuration contributing to an incomplete STU relation in $d\eta_{\finsetb}$. Furthermore the sum of the corresponding graph classes is $\left(\left[\Gamma_1\right]\left[\Gamma_2\right]-\left[\Gamma_2\right]\left[\Gamma_1\right]\right)$. This roughly shows how $d\eta_{\finsetb}=-\eta_{\finsetb} \wedge \eta_{\finsetb}$.
\end{remark}

\chapter{Discretizable variants of \texorpdfstring{$\Zinvufrf$}{Zf} and extensions to \texorpdfstring{$q$-tangles}{q-tangles}}
\label{chapdiscrex}

We introduce and study variants of $\Zinvufrf$ involving nonhomogeneous propagating forms in Sections~\ref{secdischol}, \ref{secvarzinvf}, and \ref{secstraighttang}. These variants allow discrete computations from algebraic intersections as in Chapter~\ref{chaprat}. We will use them in the proofs of important properties of $\Zinvufrf$ in Chapter~\ref{chappropzinvffunc}, where we will finish the proof of Theorem~\ref{thmmainfunc}. 
Section~\ref{zinvfqbraid} is devoted to the extension of $\Zinvufrf$ and its variants to $q$-tangles. This extension relies on the theory of semi-algebraic sets. We review known facts about semi-algebraic structures and extract useful lemmas for our purposes in Section~\ref{secsemialg}.

Throughout this chapter, $N$ denotes some (large) fixed integer $\largen$, $\largen \geq 2$, and, for $i\in \underline{3\largen}$, $\tilde{\omega}(i,S^2)=(\tilde{\omega}(i,t,S^2))_{t\in \left[0,1\right]}$ is a closed $2$-form on $\left[0,1\right] \times S^2$ such that 
$\tilde{\omega}(i,0,S^2)$ is a volume-one form of $S^2$.
As in Definition~\ref{defnumdia}, for a finite set $\finseta$, an \emph{$\finseta$-numbered}\index[T]{Anumbered@$\finseta$-numbered Jacobi diagram} Jacobi diagram is a Jacobi diagram $\Gamma$ with oriented edges equipped with an injection $j_E \colon E(\Gamma) \hookrightarrow \finseta$, which numbers the edges.
Let $\Davis^e_{n,\finseta}(\sourcetl)$\index[N]{Diag@Diagram sets!DeA@$\Davis^e_{n,\finseta}(\sourcetl)$ $\finseta$-numbered} denote the set of $\finseta$-numbered degree $n$ Jacobi diagrams with support $\sourcetl$ without looped edges.

\section{Discretizable holonomies}
\label{secdischol} 

Let $\finsetb$ be a finite set. Let $\Gamma$ be
a $\underline{3\largen}$-numbered Jacobi diagram on $\sqcup_{\eltb \in \finsetb}\RR_{\eltb}$.
For an edge $e$ of $\Gamma$ from a vertex $v(e,1)$ to a vertex $v(e,2)$,
let
\begin{equation*}p_{e,S^2} \colon \left[0,1\right] \times \ccompuptd{\vertsetv(\Gamma)}{\RR^3} \to \left[0,1\right] \times S^2\end{equation*} send $(t,\confc)$ to
$\left(t,p_{S^2}\bigl((\confc(v(e,1)),\confc(v(e,2)))\bigr)\right)$.
The maps $p_{e,S^2}$ provide the form \begin{equation*}\bigwedge_{e \in E(\Gamma)}p_{e,S^2}^{\ast}\Bigl(\tilde{\omega}\bigl(j_E(e),S^2\bigr)\Bigr)\end{equation*} over $\left[0,1\right] \times \ccompuptd{\vertsetv(\Gamma)}{\RR^3}$. This form pulls back to provide smooth forms on the smooth strata of $\left[0,1\right] \times {\CV(\Gamma)}$.
Define the smooth form $\eta_{\Gamma}=\eta_{\Gamma}((\tilde{\omega}(i,S^2))_{i \in \underline{3\largen}})$ on $\left[0,1\right] \times \cinjuptd{\finsetb}{\CC}$ so that $\eta_{\Gamma}(t,\confy)$ is the integral of $\bigwedge_{e \in E(\Gamma)}p_{e,S^2}^{\ast}(\tilde{\omega}(j_E(e),S^2))$ along the interiors $\{t\} \times \check{\CV}(\confy,\Gamma)$ of the fibers $\{t\} \times {\CV}(\confy,\Gamma)$ of $\left[0,1\right] \times {\CV}(\Gamma)$, as in Proposition~\ref{propintfib}.

\begin{definitions}
 \label{defholtilde}
When $\finseta$ is a subset of $\underline{3\largen}$ with cardinality $3n$, with $n>0$, set
\begin{equation*}\eta_{\finsetb,\finseta}=\sum_{\Gamma \in \Davis^e_{n,\finseta}(\sqcup_{\eltb \in \finsetb}\RR_{\eltb})}
\coefgambet_{\Gamma}\eta_{\Gamma}\left[\Gamma\right] \in \Omega^1\Bigl(\left[0,1\right] \times \cinjuptd{\finsetb}{\CC}; \Aavis_{n}(\sqcup_{\eltb \in \finsetb}\RR_{\eltb})\Bigr),\end{equation*}
where $\coefgambet_{\Gamma}= \frac{(\cardlef{\finseta}-\cardlef{E(\Gamma)})!}{\cardlef{\finseta}!2^{\cardlef{E(\Gamma)}}}$.
The form \begin{equation*}\eta_{\finsetb,\finseta}=\eta_{\finsetb,\finseta}\Bigl(\bigl(\tilde{\omega}(i,S^2)\bigr)_{i \in \finseta}\Bigr)\end{equation*} pulls back to a one-form on $\left[0,1\right] \times \check{C}_{\finsetb}\left[\drad{1}\right]$ still denoted by $\eta_{\finsetb,\finseta}$ with coefficients
in $\Aavis_n(\sqcup_{\eltb \in \finsetb}\RR_{\eltb})$, which is again viewed as a space of diagrams on an unoriented domain as in Definition~\ref{defrkoruniv}, Proposition~\ref{proporcomp}, and Remark~\ref{rkorunivtwo}.
Set $\eta_{\finsetb,\emptyset}=0$.

Let $\eta_{\finsetb,\largen}$ denote $\eta_{\finsetb,\underline{3\largen}}$. If $\tilde{\omega}(i,0,S^2)=\omega_{S^2}$, then the restriction of $\eta_{\finsetb,\largen}$ to $\{0\} \times \cinjuptd{\finsetb}{\CC}$ is the form $\eta_{\largen,\finsetb}$ of Definition~\ref{defetaholonomy}.

For an integer $r$, and for a set $\finseta$ of cardinality $3n$, let $P_r(\finseta)$ \index[N]{PzqrA@$P_r(\finseta)$ set of $r$-tuples}
denote the set of $r$-tuples $(\finseta_1, \finseta_2, \dots,\finseta_r)$, where $\finseta_i$ is a subset
of $\finseta$ with a cardinality multiple of $3$, the $\finseta_i$ are pairwise disjoint, and their union is $\finseta$.

Let $n$ be a positive integer. Let $\finseta$ be a subset of cardinality $3n$ of $\underline{3\largen}
$. For a path $\gamma \colon \left[0,1\right] \to \left[0,1\right] \times \cinjuptd{\finsetb}{\CC}$, define the \emph{$\finseta$-holonomy}\index[T]{Aholonomy@$\finseta$-holonomy} $\lol{\gamma}{\eta_{\finsetb,\finseta}}$\index[N]{hollg@$\lol{\gamma}{\eta_{\finsetb,\finseta}}$ $\finseta$-holonomy} of $\eta_{\finsetb,\finseta}$ along $\gamma$ to be 
\begin{multline*}\lol{\gamma}{\eta_{\finsetb,\finseta}}\\
=\sum_{r=0}^{\infty} \sum_{(\finseta_1, \dots,\finseta_r) \in P_{r}(\finseta)} \frac{\prod_{i=1}^r\cardlef{\finseta_i}!}{\cardlef{\finseta}!}\int_{(t_1,\dots,t_r) \in \left[0,1\right]^r \suchthat  t_1 \leq t_2 \leq \dots \leq t_r}\bigwedge_{i=1}^r (\gamma \circ p_i)^{\ast}(\eta_{\finsetb,\finseta_i}),\end{multline*}
with $p_i(t_1,\dots,t_r)=t_i$ and $\lol{\gamma}{\eta_{\finsetb,\emptyset}}=\left[\emptyset\right]$.
Then we have
\begin{equation*}\lol{\gamma_1\gamma_2}{\eta_{\finsetb,\finseta}}=\sum_{(\finseta_1,\finseta_2) \in P_2(\finseta)}\frac{\cardlef{\finseta_1}!\cardlef{\finseta_2}!}{\cardlef{\finseta}!}\lol{\gamma_1}{\eta_{\finsetb,\finseta_1}}\lol{\gamma_2}{\eta_{\finsetb,\finseta_2}},\end{equation*}
\begin{equation*}\frac{\partial}{\partial t}\lol{\gamma\vert_{\left[0,t\right]}}{\eta_{\finsetb,\finseta}}=\sum_{(\finseta_1,\finseta_2) \in P_2(\finseta)}\frac{\cardlef{\finseta_1}!\cardlef{\finseta_2}!}{\cardlef{\finseta}!}\lol{\gamma\vert_{\left[0,t\right]}}{\eta_{\finsetb,\finseta_1}}\eta_{\finsetb,\finseta_2}\bigl(\gamma(t),\gamma^{\prime}(t)\bigr),\end{equation*}
and 
\begin{equation*}\frac{\partial}{\partial t}\lol{\gamma\vert_{\left[t,1\right]}}{\eta_{\finsetb,\finseta}}=-\sum_{(\finseta_1,\finseta_2) \in P_2(\finseta)}\frac{\cardlef{\finseta_1}!\cardlef{\finseta_2}!}{\cardlef{\finseta}!}\eta_{\finsetb,\finseta_1}\bigl(\gamma(t),\gamma^{\prime}(t)\bigr)\lol{\gamma\vert_{\left[t,1\right]}}{\eta_{\finsetb,\finseta_2}}.\end{equation*}
\end{definitions}

Note the following lemma about the behavior of the coefficients.

\begin{lemma}
 \label{lemcoefcoef} Let $\Gamma=\Gamma_1 \sqcup \Gamma_2$ be a Jacobi diagram equipped with an injection \begin{equation*}j_E \colon E(\Gamma) \hookrightarrow \finseta,\end{equation*} where $\card{\finseta}=3\deg(\Gamma)$. We have
\begin{equation*}\coefgambet_{\Gamma_1 \cup \Gamma_2}=\sum_{\begin{array}{l}{\scriptstyle (\finseta_1,\finseta_2) \in P_2(\finseta),}\\ {\scriptstyle j_E(E(\Gamma_1)) \subseteq \finseta_1,j_E(E(\Gamma_2)) \subseteq \finseta_2,}\\{\scriptstyle  \cardlef{\finseta_1} =3\deg(\Gamma_1),\cardlef{\finseta_2} =3\deg(\Gamma_2)}\end{array}} \frac{\cardlef{\finseta_1}!\cardlef{\finseta_2}!}{\cardlef{\finseta}!}\coefgambet_{\Gamma_1}\coefgambet_{\Gamma_2}\end{equation*}
\end{lemma}
\bp We have \begin{equation*}\coefgambet_{\Gamma_1 \cup \Gamma_2}=\frac{(\cardlef{\finseta} - \cardlef{E(\Gamma_1)} -\cardlef{E(\Gamma_2)} )!}{\cardlef{\finseta}!2^{\cardlef{E(\Gamma_1 \cup \Gamma_2)}}},\end{equation*} and the number of pairs $(\finseta_1,\finseta_2)$ in the sum is \begin{equation*}\frac{(\cardlef{\finseta} - \cardlef{E(\Gamma_1)} -\cardlef{E(\Gamma_2)} )!}{(\cardlef{\finseta_1} - \cardlef{E(\Gamma_1)})!(\cardlef{\finseta_2} - \cardlef{E(\Gamma_2)})!}.\end{equation*}
\eop

\begin{remark} We fixed the cardinality of the sets $A_i$ of edge indices to be $f_0(n_i)=3n_i$ for degree $n_i$ graphs.
We could have replaced $f_0$ by
any $f \colon \NN \to \NN$ such that $f(n) \geq 3n-2$ (so that $f(n) \geq \cardlef{j_E(E(\Gamma))}$ for any degree $n$ Jacobi diagram with at least two univalent vertices) 
and such that $f(n) + f(m) = f(n+m)$, and we would have obtained the same equalities as above. However, we had to fix such a map $f$ to get these equalities.
\end{remark}

If $\tilde{\omega}(i,0,S^2)=\omega_{S^2}$, then for a subset $\finseta$ of $\underline{3\largen}$ with cardinality $3n$ and for a path $\gamma \colon \left[0,1\right] \to \{0\} \times \cinjuptd{\finsetb}{\CC}$, $\lol{\gamma}{\eta_{\finsetb,\finseta}}$ is the degree $n$ part of $\hol{\gamma}{\eta_{\finsetb}}$.

\begin{notation}
\label{notomegarev}
For a finite set $\finseta$, let $\CD_{k,\finseta}^c(\RR)$ denote the set of connected
$\finseta$-numbered degree $k$ Jacobi diagrams with support $\RR$ without looped edges.
For a vertex-oriented graph $\check{\Gamma} \in \CD_{k,\finseta}^c(\RR)$, define the two-form $\omega(\check{\Gamma})$ on $\left[0,1\right]\times S^2$ to be \begin{equation*}\omega(\check{\Gamma})(t,v)=\int_{\{t\} \times {\anomq}(v;\check{\Gamma})} \bigwedge_{e \in E(\check{\Gamma})}p_{e,S^2}^{\ast}\bigl(\tilde{\omega}(j_E(e),S^2)\bigr),\end{equation*} where $p_{e,S^2}$ again abusively denotes $\id_{\left[0,1\right]} \times p_{e,S^2}$ and ${\anomq}(v;\check{\Gamma})$ is defined in Section~\ref{secanomalpha}.
For a subset $\finseta$ of $\underline{3\largen}$ of cardinality $3k$,
define the two-form $\omega(\finseta)$ on $\left[0,1\right]\times S^2$ with coefficients in $\Assis_k(\RR)$ to be
\begin{equation*}\omega(\finseta)=\omega\Bigl(\finseta,\bigl(\tilde{\omega}(i,S^2)\bigr)_{i \in \underline{3\largen}}\Bigr)=\sum_{\check{\Gamma} \in \CD_{k,\finseta}^c(\RR)}\coefgambet_{\check{\Gamma}} \omega(\check{\Gamma})\left[\check{\Gamma}\right].\end{equation*}

Here, we view $\Assis(\RR)$ as a space of Jacobi diagrams on the oriented $\RR$. So the sets $U(\Gamma)$ of univalent vertices in involved Jacobi diagrams $\Gamma$ are ordered by $i_{\Gamma}$, and the univalent vertices are oriented by the orientation of $\RR$.
Let $\revmaps_{\ast} \colon \Assis(\RR) \to \Assis(\RR)$ be the map that sends (the class of) a Jacobi diagram $\Gamma$ on $\RR$ to (the class of) the Jacobi diagram $\revmaps(\Gamma)$ obtained from $\Gamma$ by multiplying it by $(-1)^{\cardlef{U(\Gamma)}}$ and by reversing the order of $U(\Gamma)$.\footnote{To my knowledge, the map $\revmaps_{\ast}$ might be the identity map. If it is, real-valued Vassiliev invariants cannot distinguish an oriented knot from that obtained by reversing its orientation, as Greg Kuperberg proved in \cite{Kupinvertib}. In other words, there would be no odd Vassiliev invariants. Recall the first two questions in Section~\ref{secopq}.}\end{notation}

\begin{lemma}
\label{lemrevomega}
Let $\iota\colon \left[0,1\right]\times S^2 \to  \left[0,1\right]\times S^2$
map $(t,v)$ to $(t,\iota_{S^2}(v))$, where $\iota_{S^2}$ is the antipodal map of $S^2$.
With the above notation, the form \begin{equation*}\omega(\finseta)=\omega\Bigl(\finseta,\bigl(\tilde{\omega}(i,S^2)\bigr)_{i \in \underline{3\largen}}\Bigr)\end{equation*} is a closed form of $\left[0,1\right]\times S^2$ with coefficients in $\Assis_k(\RR)$, and we have $\iota^{\ast}(\omega(\finseta))=-\revmaps_{\ast}(\omega(\finseta))$ and $\iota^{\ast}(\omega(\finseta))=(-1)^k(\omega(\finseta))$.
\end{lemma}
\bp In order to prove that $\omega(\finseta)$ is closed, it suffices to prove that its integral vanishes on the boundary of any $3$-ball $B$ of $\left[0,1\right] \times S^2$. 
View \begin{equation*}{\anomq}_k(t,v)=\sum_{\check{\Gamma} \in \CD_{k,\finseta}^c(\RR)}\coefgambet_{\check{\Gamma}} \left[\check{\Gamma}\right] \left(\{t\} \times {\anomq}(v;\check{\Gamma}) \times (S^2)^{\finseta \setminus j_E(E(\check{\Gamma}))}\right)\end{equation*} as a $(6k-2)$-chain with coefficients in $\Assis_k(\RR)$. So
the integral $\int_{\partial B} \omega(\finseta)$ is the integral of the  pull-back of the closed form 
$\bigwedge_{a \in \finseta}\tilde{\omega}(a,S^2)$ under a natural map $P$ from $\cup_{(t,v) \in \partial B}{\anomq}_k(t,v)$ to $\left[0,1\right] \times (S^2)^{\finseta}$.

The analysis of the boundary of ${\anomq}(v;\check{\Gamma})$ in the proof of Proposition~\ref{propanom} shows that the codimension-one faces of the boundary of $P({\anomq}_k(t,v))$ can be glued. So the boundary of $P({\anomq}_k(t,v))$ vanishes algebraically.

Therefore, since the cycle $\cup_{(t,v) \in \partial B}P({\anomq}_k(t,v))$ bounds $\cup_{(t,v) \in B}P({\anomq}_k(t,v))$ in $\left[0,1\right] \times (S^2)^{\finseta}$, the integral $\int_{\partial B} \omega(\finseta)$ vanishes, and $\omega(\finseta)$ is closed.

For $\check{\Gamma} \in \CD_{k,\finseta}^c(\RR)$, recall that the class of $\revmaps(\check{\Gamma})$ is obtained from the class of $\check{\Gamma}$ by multiplying it by $(-1)^{\cardlef{U(\check{\Gamma})}}$ and by reversing the order of $U(\check{\Gamma})$. 
A configuration of ${\anomq}(v;\check{\Gamma})$ is naturally a configuration of ${\anomq}(-v;\revmaps(\check{\Gamma}))$.
Below, we consider $\check{\Gamma}$ and $\revmaps(\check{\Gamma})$ as diagrams on unoriented domains whose univalent vertices are equipped with matching local orientations of the domain at univalent vertices. 
The natural identification from ${\anomq}(v;\check{\Gamma})$ to ${\anomq}(-v;\revmaps(\check{\Gamma}))$ reverses the orientation since the local orientations at univalent vertices coincide, the quotients by dilation coincide, but the translations act in opposite directions.
Therefore, for any $2$-chain $\Delta$ of $\left[0,1\right] \times S^2$, we have $\int_{\Delta}\omega(\check{\Gamma})=- \int_{\iota(\Delta)}\omega(\revmaps(\check{\Gamma})),$ where $\iota(\Delta)$ is equipped with the orientation of $\Delta$. We get
\begin{equation*}\begin{array}{ll}\sum_{\check{\Gamma} \in \CD_{k,\finseta}^c(\RR)}\coefgambet_{\check{\Gamma}} \int_{\Delta} \omega(\check{\Gamma})\left[\check{\Gamma}\right] & = -\sum_{\check{\Gamma} \in \CD_{k,\finseta}^c(\RR)}\coefgambet_{\revmaps(\check{\Gamma})} \int_{\iota(\Delta)}\omega\bigl(\revmaps(\check{\Gamma})\bigr)\left[\check{\Gamma}\right]\\& = -\sum_{\check{\Gamma} \in \CD_{k,\finseta}^c(\RR)}\coefgambet_{\check{\Gamma}} \int_{\iota(\Delta)}\omega(\check{\Gamma})[\revmaps(\check{\Gamma})]\end{array}\end{equation*}
and
$\omega(\finseta)=-\iota^{\ast}(\revmaps_{\ast}(\omega(\finseta)))$.

In order to prove $\iota^{\ast}(\omega(\finseta))=(-1)^k(\omega(\finseta))$, use the notation and the arguments of the proof that $\alpha_{2n}=0$ in the proof of Proposition~\ref{propanom} to prove that for any $2$-chain $\Delta$ of $\left[0,1\right] \times S^2$ and any $\check{\Gamma} \in \CD_{k,\finseta}^c(\RR)$, we have
\begin{equation*}\int_{\Delta}\omega(\check{\Gamma})\left[\check{\Gamma}\right]=(-1)^k \int_{\iota(\Delta)}\omega(\check{\Gamma}^{eo})[\check{\Gamma}^{eo}].\end{equation*}
\eop

\begin{definition} \label{defAholonom}
For any oriented tangle component $K$, let $\ST^+K$ denote the fiber space over $K$ consisting of the tangent vectors to the knot $K$ of $\crats(\hcylc)$ that orient $K$ up to dilation, as in Section~\ref{secthetaknot}. When the ambient manifold is equipped with a parallelization $\tau$, define the one-form
$\eta(\finseta,p_{\tau}(\ST^+K))=\eta(\finseta,(\tilde{\omega}(i,S^2))_{i \in \underline{3\largen}},p_{\tau}(\ST^+K))$  on $\left[0,1\right]$ valued in $\Assis(\RR)$ to be
\begin{equation*}\eta\bigl(\finseta,p_{\tau}(\ST^+K)\bigr)(t)=\eta\bigl(\finseta,p_{\tau}(\ST^+K),\tau\bigr)(t)=\int_{\{t\} \times p_{\tau}(\ST^+K)}\omega\Bigl(\finseta,\bigl(\tilde{\omega}(i,S^2)\bigr)_{i \in \underline{3\largen}}\Bigr).\end{equation*}
More precisely, we have
\begin{equation*}\eta\bigl(\finseta,p_{\tau}(\ST^+K)\bigr)\Bigl(u,\frac{\partial}{\partial t}\Bigr)=\frac{\partial}{\partial t} \left( \int_{\left[0,t\right]\times p_{\tau}(\ST^+K)}\omega(\finseta)\right)(u)\end{equation*} as in Proposition~\ref{propintfib}.
Define the \emph{$\finseta$-holonomy} of $\eta(.,p_{\tau}(\ST^+K))$ along $[a,b] \subseteq \left[0,1\right]$ to be \begin{multline*}\Biglol{[a,b]}{\eta\bigl(\finseta,p_{\tau}(\ST^+K)\bigr)}=\\
\sum_{r=0}^{\infty} \sum_{(\finseta_1, \dots,\finseta_r) \in P_r(\finseta) } \frac{\prod_{i=1}^r\cardlef{\finseta_i}!}{\cardlef{\finseta}!}\int_{(t_1,\dots,t_r) \in [a,b]^r \suchthat  t_1 \leq t_2 \leq \dots \leq t_r}\bigwedge_{i=1}^r  p_i^{\ast}\Bigl(\eta\bigl(\finseta_i,p_{\tau}(\ST^+K)\bigr)\Bigr)\end{multline*}
\index[T]{Aholonomy@$\finseta$-holonomy}
\index[N]{holli@$\lol{[a,b]}{.}$ $\finseta$-holonomy}
with $p_i(t_1,\dots,t_r)=t_i$ and $\lol{[a,b]}{\eta(\emptyset,p_{\tau}(\ST^+K))}=\id=\left[\emptyset\right]$.
\end{definition}

\begin{lemma}
\label{lemrevholomega} For any oriented tangle component $K$, set $\ST^-K = \ST^+(-K)$. We have
\begin{multline*}\Biglol{[a,b]}{\eta\bigl(\finseta,-p_{\tau}(\ST^-K)\bigr)=-\eta\bigl(\finseta,p_{\tau}(\ST^-K)\bigr)}\\=\revmaps_{\ast}\left(\Biglol{[a,b]}{\eta\bigl(\finseta,p_{\tau}(\ST^+K)\bigr)}\right).\end{multline*}
\end{lemma}
\bp Lemma~\ref{lemrevomega} implies $\eta\bigl(\finseta,-p_{\tau}(\ST^-K)\bigr)=\revmaps_{\ast}\Bigl(\eta\bigl(\finseta,p_{\tau}(\ST^+K)\bigr)\Bigr)$. \eop

Also note the following lemma for later use.
\begin{lemma} 
\label{lemaddhol} Let $K(\source)$ be a connected oriented one-dimensional submanifold of a $3$-manifold $\rats(\hcylc)$ equipped with a parallelization $\tau$. Assume that $K=K(\source)$ is a finite union $K=\cup_{j=1}^kK_j$ of intervals $K_j$ with disjoint interiors. Let $\finseta$ be a subset of $\underline{3\largen}$ of cardinality multiple of $3$. Then we have
$\eta (\finseta,p_{\tau}(\ST^+K))=\sum_{j=1}^k\eta(\finseta,p_{\tau}(\ST^+K_j))
$
and \begin{multline*}\Biglol{\left[0,1\right]}{\eta\bigl(\finseta,p_{\tau}(\ST^+K)\bigr)}
\\= \sum_{(\finsetb_1, \dots,\finsetb_k) \in P_k(\finseta) } \frac{\prod_{j=1}^k\cardlef{\finsetb_j}!}{\cardlef{\finseta}!}\prod_{j=1}^k\Biglol{\left[0,1\right]}{\eta\bigl(\finsetb_j,p_{\tau}(\ST^+K_j)\bigr)}\end{multline*}
in $\Assis(\source)$.
\end{lemma}
\bp It suffices to prove the lemma when $k$ equals $2$.
Recall
\begin{equation*}\Biglol{\left[0,1\right]}{\sum_{j=1}^2\eta\bigl(.,p_{\tau}(\ST^+K_j)\bigr)}=\sum_{r=0}^{\infty} \sum_{(\finseta_1, \dots,\finseta_r) \in P_r(\finseta) } \frac{\prod_{i=1}^r\cardlef{\finseta_i}!}{\cardlef{\finseta}!}I(\finseta_1, \dots,\finseta_r),\end{equation*}
with \begin{multline*}I(\finseta_1, \dots,\finseta_r)\\=\int_{(t_1,\dots,t_r) \in \left[0,1\right]^r \suchthat  t_1 \leq t_2 \leq \dots \leq t_r}\bigwedge_{i=1}^r p_i^{\ast}\Bigl(\eta\bigl(\finseta_i,p_{\tau}(\ST^+K_1)\bigr)+\eta\bigl(\finseta_i,p_{\tau}(\ST^+K_2)\bigr)\Bigr).
     \end{multline*}
     Write 
     \begin{multline*}\bigwedge_{i=1}^r p_i^{\ast}\Bigl(\eta\bigl(\finseta_i,p_{\tau}(\ST^+K_1)\bigr)+\eta\bigl(\finseta_i,p_{\tau}(\ST^+K_2)\bigr)\Bigr)\\= \sum_{f \colon \{1,\dots,r\}\to\{1,2\}} \bigwedge_{i=1}^r p_i^{\ast}\Bigl(\eta\bigl(\finseta_i,p_{\tau}(\ST^+K_{f(i)})\bigr)\Bigr).
\end{multline*}
For $f$ as above, also write
\begin{multline*}
 I_j\bigl((\finseta_i)_{i \in f^{-1}(j)}\bigr)\\=\int_{(t_i)_{i \in f^{-1}(j)} \in \left[0,1\right]^{f^{-1}(j)} \suchthat  t_i \leq t_s \mbox{ \scriptsize when } i \leq s}\bigwedge_{i\in f^{-1}(j)}  p_i^{\ast}\Bigl(\eta\bigl(\finseta_i,p_{\tau}(\ST^+K_j)\bigr)\Bigr).
\end{multline*}
Decompose the terms $I_1((\finseta_i)_{i \in f^{-1}(1)})I_2((\finseta_i)_{i \in f^{-1}(2)})$ of the right-hand side of the equality to be proved with respect to the relative orders of the involved $t_i$ in $I_1$ and $I_2$.

The contribution to $\lol{\left[0,1\right]}{\sum_{j=1}^2\eta(.,p_{\tau}(\ST^+K_j))}$ of the terms such that
\begin{itemize}
\item $\{\finseta_i \suchthat  i \in \underline{r}\}$ is fixed as an unordered set,
\item $f$ is fixed as a map from this unordered set to $\{1,2\}$, and 
\item the partial order induced by the numbering in $\underline{r}$ is fixed over the sets $\{\finseta_i \suchthat  i \in f^{-1}(j)\}$ for $j=1, 2$, 
\end{itemize}
is also a sum over the possible total orders
on the $\{\finseta_i \suchthat  i \in \underline{r}\}$ that induce the given orders on the two subsets.

We can easily identify the involved coefficients of the above similar contributions to prove the result when $s$ equals $2$. 
 \eop

\section{Variants of \texorpdfstring{$\Zinvufrf$}{Zf} for tangles}
\label{secvarzinvf}

We now present alternative definitions of $\Zinvufrf$ involving nonhomogeneous propagating forms
associated to volume forms $\tilde{\omega}(i,1,S^2)=\omega_{S^2} +d\eta(i,S^2)$,
where $\omega_{S^2}$ is the homogeneous volume-one form of $S^2$. 

Let $\hcylc$ be a rational homology cylinder equipped with a parallelization $\tau$.
Let $\tanghcyll \colon \sourcetl \hookrightarrow \rats(\hcylc)$ denote a long tangle representative.
Recall that we fixed some (large) integer $\largen$, $\largen \geq 2$, and that
 $\bigl(\tilde{\omega}(i,S^2)=(\tilde{\omega}(i,t,S^2))_{t\in \left[0,1\right]}\bigr)_{i\in \underline{3\largen}}$ is a family of a closed $2$-forms on $\left[0,1\right] \times S^2$ such that 
$\tilde{\omega}(i,0,S^2)$ is a volume-one form of $S^2$.
For $i\in \underline{3\largen}$, let $\tilde{\omega}(i)=(\tilde{\omega}(i,t))_{t\in \left[0,1\right]}$ be a closed $2$-form on $\left[0,1\right] \times C_2(\rats(\hcylc))$ such that $\tilde{\omega}(i)=(\id_{\left[0,1\right]} \times p_{\tau})^{\ast}(\tilde{\omega}(i,S^2))$ on
$\left[0,1\right]\times \partial C_2(\rats(\hcylc))$. 

For a diagram $\Gamma \in \Davis^e_{k,\underline{3\largen}}(\sourcetl)$, define \begin{equation*}I\Bigl(\hcylc,\tanghcyll,\Gamma,o(\Gamma),\bigl(\tilde{\omega}(i,t)\bigr)_{i \in \underline{3\largen}}\Bigr)=\int_{(\check{C}(\crats(\hcylc),\tanghcyll;\Gamma),o(\Gamma))} \bigwedge_{e \in E(\Gamma)}p_e^{\ast}\Bigl(\tilde{\omega}\bigl(j_E(e),t\bigr)\Bigr),\end{equation*}
which converges, according to Theorem~\ref{thmconvint}.

\begin{theorem}
\label{thmtangconstcompar} Let $\tanghcyll \colon \sourcetl \hookrightarrow \rats(\hcylc)$ denote the long tangle associated to a tangle in a rational homology cylinder equipped with a parallelization $\tau$. Let $\{K_j\}_{j \in I}$ be the set of components of $\tanghcyll$. 
Assume that the bottom and top configurations of $\tanghcyll$ are respectively represented by maps $\confy^- \colon \finsetb^- \to \drad{1}$ and $\confy^+ \colon \finsetb^+ \to \drad{1}$. 

Let $\largen \in \NN$.
For $i\in \underline{3\largen}$, let $\tilde{\omega}(i,S^2)=(\tilde{\omega}(i,t,S^2))_{t\in \left[0,1\right]}$ be a closed $2$-form on $\left[0,1\right] \times S^2$ such that $\tilde{\omega}(i,0,S^2)$  is 
a volume-one form of $S^2$,\footnote{In this chapter, we apply the theorem only when $\tilde{\omega}(i,0,S^2)$ is the standard homogeneous volume-one form $\omega_{S^2}$ on $S^2$. However, we use this general statement in the next chapter.} 
and let $\tilde{\omega}(i)=(\tilde{\omega}(i,t))_{t\in \left[0,1\right]}$ be a closed $2$-form on $\left[0,1\right] \times C_2(\rats(\hcylc))$ such that $\tilde{\omega}(i)=(\id_{\left[0,1\right]} \times p_{\tau})^{\ast}(\tilde{\omega}(i,S^2))$ on
$\left[0,1\right]\times \partial C_2(\rats(\hcylc))$.
For a subset $\finseta$ of $\underline{3\largen}$ with cardinality $3k$, set
\begin{equation*}\Zinv\Bigl(\hcylc,\tanghcyll,\tau,\finseta,\bigl(\tilde{\omega}(i,t)\bigr)_{i \in \finseta}\Bigr)=\sum_{\Gamma \in \Davis^e_{k,\finseta}(\sourcetl)}\coefgambet_{\Gamma}I\Bigl(\hcylc,\tanghcyll,\Gamma,\bigl(\tilde{\omega}(i,t)\bigr)_{i \in \finseta}\Bigr)\left[\Gamma\right] \in \Aavis_{k}(\sourcetl)\end{equation*}
and \begin{equation*}\Zinv(\hcylc,\tanghcyll,\tau,\finseta)(t)=\Zinv\Bigl(\hcylc,\tanghcyll,\tau,\finseta,\bigl(\tilde{\omega}(i,t)\bigr)_{i \in \finseta}\Bigr).\end{equation*}

Then we have \begin{equation*}\Zinv(\hcylc,\tanghcyll,\tau,\finseta)(t)=\sum_{\aleph=(\finseta_1,\finseta_2,\finseta_3,(\finseta_{K_j})_{j \in I}) \in P_{3+\cardlef{I}}(\finseta)}\zeta_d(\aleph) \Zinv(\aleph,t),\end{equation*}
with \begin{equation*}\zeta_d(\aleph)=\frac{ \cardbig{\finseta_1}!\cardbig{\finseta_2}!\cardbig{\finseta_3}! \prod_{j \in I}\cardbig{\finseta_{K_j}}! }{\cardlef{\finseta}!}\end{equation*}
and 
\begin{multline*}\Zinv(\aleph,t)= 
\left(\prod_{j\in I}\Biglol{\left[0,t\right]}{\eta\bigl(\finseta_{K_j},p_{\tau}(\ST^+K_j)\bigr)} \#_j\right)\\
\lol{\left[t,0\right] \times \confy^- }{\eta_{\finsetb^-,\finseta_1}}\Zinv(\hcylc,\tanghcyll,\tau,\finseta_2)(0)\lol{\left[0,t\right]\times \confy^+ }{\eta_{\finsetb^+,\finseta_3}}.\end{multline*}
The terms of this formula belong to spaces of diagrams on unoriented one-manifolds as in Definition~\ref{defrkoruniv}, except for the term $\lol{\left[0,t\right]}{\eta(\finseta_{K_j},p_{\tau}(\ST^+K_j))}$ and its action $\#_j$ for which we first pick an orientation of the $K_j$, which we may forget afterwards.\footnote{Both sides are independent of the domain orientations in the sense of the last sentence of Proposition~\ref{proporcomp} and Remark~\ref{rkorunivtwo}, thanks to Lemma~\ref{lemrevholomega}.}
The formula implies that $\Zinv(\hcylc,\tanghcyll,\tau,\finseta)(t)$ depends only on $(\tilde{\omega}(i,t,S^2))_{i \in \finseta}$ for any $t$. (It also depends on $\tau$ and on the specific embedding $\tanghcyll$.) It will be denoted by $\Zinv(\hcylc,\tanghcyll,\tau,\finseta,(\tilde{\omega}(i,t,S^2))_{i \in \finseta})$.
\end{theorem}
\bp
Compute $\frac{\partial}{\partial t}\Zinv(\finseta,t)=d\Zinv(\finseta,.) \left(\frac{\partial}{\partial t}\right)$ as 
 in Lemma~\ref{leminvtwo} with the help of Proposition~\ref{propinvone}, using the same analysis of faces as in the proof of Lemma~\ref{lemholzinf}, to find
\begin{equation*}\begin{array}{ll}d\Zinv(\finseta,.)=&
\sum_{(\finseta_1,\finseta_2) \in P_2(\finseta)}\frac{\cardlef{\finseta_1}!\cardlef{\finseta_2}!}{\cardlef{\finseta}!}\sum_{j\in I}\left(\eta(\finseta_1,p_{\tau}(\ST^+K_j))\#_j\right) \Zinv(\finseta_2,t)\\
&+\sum_{(\finseta_1,\finseta_2) \in P_2(\finseta)}\frac{\cardlef{\finseta_1}!\cardlef{\finseta_2}!}{\cardlef{\finseta}!}\Zinv(\finseta_1,t)(\{t\} \times\confy^+ )^{\ast}({\eta_{\finsetb^+,\finseta_2}})
\\
&-\sum_{(\finseta_1,\finseta_2) \in P_2(\finseta)}\frac{\cardlef{\finseta_1}!\cardlef{\finseta_2}!}{\cardlef{\finseta}!}(\{t\} \times \confy^-)^{\ast}({\eta_{\finsetb^-,\finseta_1}}) \Zinv(\finseta_2,t),
\end{array}\end{equation*}
with $\eta(\emptyset,p_{\tau}(\ST^+K_j))=0$ and $\eta_{\finsetb^+,\emptyset}=0$. (Recall Lemma~\ref{lemcoefcoef} for the behavior of the coefficients.) This
shows that both sides of the equality to be proved vary in the same way when $t$ varies. Since they take the
same value at $t=0$, the formula is proved.
Apply the formula when $\tilde{\omega}(i,0,S^2)$ is the standard form ${\omega}_{S^2}$, and use Lemma~\ref{lemformprod} together with the isotopy invariance of $\Zinv$ of Theorem~\ref{thmfstconsttang} to see that $\Zinv(\hcylc,\tanghcyll,\tau,\finseta)(t)$ depends only on the $\tilde{\omega}(i,S^2)$ for $i \in \finseta$. So it depends only on the $\tilde{\omega}(i,t,S^2)$.
\eop

Let us introduce some notation in order to rephrase Theorem~\ref{thmtangconstcompar}.
View $\Zinv(\hcylc,\tanghcyll,\tau,.)(t)$ as a map from the set $\parentp_{(3)}(\underline{3\largen})$ of subsets of $\underline{3\largen}$ with cardinality multiple of $3$ to
$\Aavis_{\leq \largen}(\sourcetl)=
\oplus_{k=0}^{\largen}\Aavis_{k}(\sourcetl)$, which maps $\emptyset$ to the class of the empty diagram.
Similarly, consider $\lol{\left[0,t\right]}{\eta(.,p_{\tau}(\ST^+K_j))}$,
$\lol{\left[t,0\right] \times \confy^-}{\eta_{\finsetb^-,.}}$ and $\lol{\left[0,t\right] \times\confy^+ }{\eta_{\finsetb^+,.}}$ as maps from $\parentp_{(3)}(\underline{3\largen})$ to spaces of diagrams, which map the empty set to the class of the empty diagram.
The values of these maps can be multiplied as in the statement of the theorem using the structures of the space of diagrams.
\begin{definition} \label{defprodbiz}
For such maps $z_1$ and $z_2$ from $P_{(3)}(\underline{3\largen})$ to spaces of diagrams, define their product $(z_1 z_2)_{\sqcup}$ to be the map with domain $\parentp_{(3)}(\underline{3\largen})$ such that
\begin{equation*}(z_1 z_2)_{\sqcup}(\finseta)=\sum_{(\finseta_1,\finseta_2) \in P_2(\finseta)}\frac{\cardlef{\finseta_1}!\cardlef{\finseta_2}!}{\cardlef{\finseta}!}z_1(\finseta_1)z_2(\finseta_2)\end{equation*}
whenever the products $z_1(\finseta_1)z_2(\finseta_2)$ make sense.
\end{definition}
This product is associative, and $(z_1 z_2 z_3)_{\sqcup}$ denotes
\begin{equation*}\bigl((z_1 z_2)_{\sqcup}z_3\bigr)_{\sqcup}=\bigl(z_1 (z_2z_3)_{\sqcup}\bigr)_{\sqcup}.\end{equation*}
The maps that send all nonempty elements of $P_{(3)}(\underline{3\largen})$ to $0$ and the empty element to the class of the empty diagram are neutral
for this product. They are denoted by $\id$.
With the above definition, we may write the equality of Theorem~\ref{thmtangconstcompar} as
\begin{multline*}\Zinv(.)(t)=\\
\left(\left(\prod_{j\in I}\biglol{\left[0,t\right]}{\eta(.,\ST^+_j)} \#_j\right)
\lol{\left[t,0\right] \times \confy^-}{\eta_{\finsetb^-,.}}\Zinv(.)(0)\lol{\left[0,t\right] \times \confy^+}{\eta_{\finsetb^+,.}} \right)_{\sqcup} \end{multline*}
with $\Zinv(.)(t)=\Zinv(\hcylc,\tanghcyll,\tau,.)(t)$ and $\ST^+_j=p_{\tau}(\ST^+K_j)$.

\section{Straight tangles}
\label{secstraighttang}

\begin{definition}
 \label{defstraighttang}
Recall that $\upvec$ denotes the vertical unit vector.
A tangle $\tanghcyll \colon \sourcetl \hookrightarrow \hcylc$ is \emph{straight}\index[T]{straight!tangle} with respect to $\tau$ if
\begin{itemize}
 \item $p_{\tau}(\ST^+K) \subset \{-\upvec,\upvec\}$ for closed components $K$ and components $K$ going from bottom to top or from top to bottom,
\item for any interval component $\source_j$, the map $p_{\tau}$ sends the unit tangent vectors to $\tanghcyll(\source_j)$ to the vertical half great circle $S_{WE}$ from $-\upvec$ to $\upvec$ that contains the west-east direction (to the right), or to the vertical half great circle $S_{EW}$ from $-\upvec$ to $\upvec$ that contains the east-west direction (to the left). See Figure~\ref{figSWE}.
\end{itemize}

\bfig
\centering
\begin{equation*}
\begin{tikzpicture}
\draw
(1.15,0) node{\scriptsize $S_{WE}$}
(-1.15,0) node{\scriptsize $S_{EW}$}
(-.2,-1) node{\scriptsize $-\upvec$}
(.2,1) node{\scriptsize $\upvec$};
\draw [->] (0,-.7) arc (-90:0:.7);
\draw [->] (0,-.7) arc (270:180:.7);
\draw (.7,0)  arc (0:180:.7);
\fill (0,-.7) circle (1.5pt) (0,.7) circle (1.5pt);
\end{tikzpicture}
\end{equation*}
\caption{The half-circles $S_{WE}$ and $S_{EW}$.}
\label{figSWE}
\end{figure}

Orient $S_{WE}$ and $S_{EW}$ from $-\upvec$ to $\upvec$.
Straight tangles with respect to $\tau$ get the following framing induced by $\tau$.
For any $k \in K$, $p_{\tau}(\ST_k^+K)$ is an element $v_k$ of the vertical circle $S_{WE} \cup (-S_{EW})$, which intersects the horizontal plane $\CC$ in the real line.
Let $\rho_{i,\pi/2}(v_k)$ be the image of $v_k$ under the rotation of angle $\pi/2$ about the axis $i$ ($(i\in\CC)$ points toward the sheet). Then $K_{\parallel}$ is the parallel of $K$ obtained by pushing $K$ in the direction of the section $\left(k \mapsto \tau(\rho_{i,\pi/2}(v_k))\right)$ of the unit normal bundle of $K$.
(This is consistent with the conventions of Definition~\ref{defparalleltang}.)
\end{definition}

\begin{remark}
A boundaryless straight tangle is a straight link in the sense of Definition~\ref{defstraightlink}. However, the converse is not true. The present definition is more restrictive.
\end{remark}

The following proposition generalizes Lemma~\ref{lemdegloctheta} to interval components.

\begin{proposition}
\label{proplktangles}
Let $K$ be a component of a straight $q$-tangle in a parallelized homology cylinder $(\hcylc,\tau)$. 
Then we have \begin{equation*}I_{\theta}(K,\tau)=I_{\theta}(-K,\tau)=lk(K,K_{\parallel})\end{equation*}
with the notation of Definitions~\ref{defselflkqtangone} and \ref{defselflkqtanggen} for $lk(K,K_{\parallel})$, and of Definition~\ref{defIthetalong} and Lemma~\ref{lemdefItheta} for $I_{\theta}$. 
\end{proposition} 

In order to prove Proposition~\ref{proplktangles}, we describe some propagating forms, which will also be useful in the next chapter.

\begin{notation}
\label{notratsrest}
For an interval $I$ of $\RR$ that contains $\left[0,1\right]$
and for a real number $x \in \left[1,+\infty\right[$,
$\rats_{x,I}(\hcylc)$ denotes the part that replaces $\drad{x} \times I$ in $\rats(\hcylc)$, and $\rats^c_{x,I}(\hcylc)=S^3 \setminus (\dorad{x} \times \mathring{I})$ denotes the closure of its complement.

Let $\chi_{\hcylc}$ be a smooth map from $\crats(\hcylc)$ to $\left[0,1\right]$ that sends $\hcylc=\rats_{1,\left[0,1\right]}(\hcylc)$ to $1$ and $\rats^c_{2,\left[-1,2\right]}(\hcylc)$ to $0$.
Define \begin{equation*}\begin{array}{llll}\pi_{\hcylc} \colon &\crats(\hcylc) &\to &\RR^3\\
         &x & \mapsto &(1-\chi_{\hcylc}(x))x,
        \end{array}
\end{equation*} with $0x=0$.
The map 
\begin{equation*}\begin{array}{llll}p \colon &\biggl(\crats(\hcylc)\biggr)^2 \setminus \biggl(\crats_{2,\left[-1,2\right]}(\hcylc)^2 \cup \diagon\Bigl(\crats^c_{2,\left[-1,2\right]}(\hcylc)^2\Bigr)\biggr) &\to &S^2\\
         &(x,y) & \mapsto & \frac{\pi_{\hcylc}(y)-\pi_{\hcylc}(x)}{\norm{\pi_{\hcylc}(y)-\pi_{\hcylc}(x)}}
        \end{array}\end{equation*}
extends to \begin{equation*}D(p)=C_2\bigl(\rats(\hcylc)\bigr) \setminus \mathring{C}_2\bigl(\rats_{2,\left[-1,2\right]}(\hcylc)\bigr).\end{equation*}
When a parallelization $\tau$ of $\hcylc$ is given, the corresponding extension of $p$ to \begin{equation*}D(p_{\tau})=D(p) \cup \ST\rats_{2,\left[-1,2\right]}(\hcylc)\end{equation*} is denoted by $p_{\tau}$.
\end{notation}

\begin{lemma}
\label{lemextformfunct}
For any volume-one $2$-form $\omega(S^2)$ on $S^2$ and any parallelization $\tau$ of $\hcylc$ (as in Definition~\ref{defparacyl}), there exists a propagating form $\omega$ of $(C_2(\crats(\hcylc)),\tau)$ that restricts to $D(p_{\tau})$ as $p_{\tau}^{\ast}(\omega(S^2))$.
For any $\vecx \in S^2$, there exists a propagating chain $F$ of $(C_2(\crats(\hcylc)),\tau)$ that restricts to $D(p_{\tau})$ as $p_{\tau}^{-1}(\vecx)$.
\end{lemma}
\bp
We again need an extension to the interior of
${C}_2(\rats_{2,\left[-1,2\right]}(\hcylc))$ of a closed $2$-form defined on the boundary. Since this space is a $6$-manifold with ridges with the same homology as $S^2$, the form extends as a closed form.
\eop

\bpo{Proof of Proposition~\ref{proplktangles}}
Lemma~\ref{lemdegloctheta} leaves us with the case of interval components $K$. In this case, we have $I_{\theta}(K,\tau)=2 I(\onechordsmallnonnum,\omega_{S^2})$ and
$I_{\theta}(-K,\tau)= I_{\theta}(K,\tau)$.

First assume that $K$ goes from bottom to top.
In this case, Lemma~\ref{lemvarithetatang} guarantees that $I(\onechordsmallnonnum,\omega(S^2))$ is independent of the chosen volume-one form $\omega(S^2)$ of $S^2$. Choose a volume-one form $\omega_{0,0}$ of $S^2$ $\varepsilon$-dual to $p_{\tau}^{-1}(i)$ for the complex horizontal direction $i$ and for a small positive number $\varepsilon$. (Recall Definition~\ref{defformdual}.)  Set $\omega_0(S^2) =\frac1{2}({\omega_{0,0} -\iota^{\ast}(\omega_{0,0})})$, and choose a propagating form $\omega_0$ that restricts to $D(p_{\tau})$ as $p_{\tau}^{\ast}(\omega_0(S^2))$ as in Lemma~\ref{lemextformfunct}.
Assume that $\varepsilon$ is sufficiently small so that $p_{\tau}$ maps $(K \times K_{\parallel}) \cap D(p_{\tau})$ outside the support of $\omega_0(S^2)$.
Define 
knots $C$ and $C_{\parallel}$ such that 
\begin{itemize}
\item the knots $C$ and $C_{\parallel}$ respectively coincide with $K$ and $K_{\parallel}$ on $\rats_{9,[-9,9]}(\hcylc)$, 
\item 
$(C,C_{\parallel})$ is isotopic to the pair $(\hat{K},\hat{K}_{\parallel})$ of Definition~\ref{defselflkqtangone} of $lk(K,K_{\parallel})$ so that we have $lk(K,K_{\parallel})=lk(C,C_{\parallel})$,
and
\item $p_{\tau}$ maps $(C \times C_{\parallel}) \cap D(p_{\tau})$ 
outside the support of $\omega_0(S^2)$ so that we have $\int_{C \times C_{\parallel}}\omega_0=\int_{K \times K_{\parallel}}\omega_0$.
\end{itemize}
Our choice of $\omega_0$ also allows us to let $K_{\parallel}$ approach and replace $K$ without changing the rational integral $\int_{K \times K_{\parallel}}\omega_0$. So we have $\int_{K \times K_{\parallel}}\omega_0
=\int_{K\times K \setminus \diag}\omega_0.$
Then we get \begin{equation*}\begin{array}{ll}lk(K,K_{\parallel})&=\int_{C \times C_{\parallel}}\omega_0=\int_{K \times K_{\parallel}}\omega_0
=\int_{K\times K \setminus \diag}\omega_0\\&=I\left(\onechordsmallnonnum,\omega_0(S^2)\right)+I\left(\onechordsmallnonnumorrev,\omega_0(S^2)\right)\\&=I\left(\onechordsmallnonnum,\omega_0(S^2)\right)+I\left(\onechordsmallnonnum,-\iota^{\ast}(\omega_0(S^2))\right)
=2I\left(\onechordsmallnonnum,\omega_0(S^2)\right).\end{array}\end{equation*}
This proves
Proposition~\ref{proplktangles} when $K$ goes from bottom to top.
Since we have $I_{\theta}(-K,\tau)=I_{\theta}(K,\tau)$, Definition~\ref{defselflkqtangone} allows us to deduce it
when $K$ goes from top to bottom.

Let us now assume that $K$ goes from bottom to bottom (resp. from top to top). 
Lemma~\ref{lemlkindepor} reduces the proof to the case $p_{\tau}(\ST^+K) \subseteq S_{WE}$.

Let $S(K)$ be the half-circle from $-\upvec$ to $\upvec$ (resp. from $\upvec$ to $-\upvec$) through the direction of $\left[K(0),K(1)\right]$.
Lemma~\ref{lemvarithetatang} implies that $I(\onechordsmallnonnum,\omega(S^2))$ depends only on the integrals of the restriction of $\omega(S^2)$ to the components of $S^2 \setminus (p_{\tau}(\ST^+K)\cup S(K))$, with the notation of Lemma~\ref{lemvarithetatang}. In particular, when $S(K)=\pm S_{WE}$ (i.e., when the segment $\left[K(0),K(1)\right]$ is directed and oriented as the real line, as in Figure~\ref{figgammaparthetazero} 
), it does not depend on $\omega(S^2)$, and we conclude as above.

Otherwise, $S^2 \setminus (p_{\tau}(\ST^+K)\cup S(K))$ has two connected components $A_1$ and $A_2$ such that
$\partial \overline{A}_1=p_{\tau}(\ST^+K)\cup S(K)=-\partial \overline{A}_2$.
For $j \in \underline{2}$, let $I_j=I(\onechordsmallnonnum,\omega_j(S^2))$, where $\omega_j(S^2)$ is a volume-one form supported on $A_j$. Then $I_2$ is a rational number since it is the intersection of a propagating chain with boundary $p_{\tau}^{-1}(\vecx)$ for $\vecx \in A_2$ with 
$C(\rats(\hcylc,K;\onechordsmallnonnum))$ in ${C}_2(\rats(\hcylc))$.\footnote{When $\hcylc=\drad{1} \times \left[0,1\right]$ and $\tau = \taust$, the map $p_{\tau}$ extends to ${C}_2(\rats(\hcylc))$ and $I_2$ is the integral local degree of this extended $p_{\tau}$ at a point of $A_2$.}
According to Lemma~\ref{lemvarithetatang}, for any volume-one form $\omega(S^2)$ of $S^2$, we have
\begin{equation*}\begin{array}{ll} I(\onechordsmallnonnum,\omega(S^2))&=I_1 +\int_{A_1}\bigl(\omega(S^2)-\omega_1(S^2)\bigr)=I_1 +\int_{A_1}\omega(S^2)-1\\&
  =I_2 + \int_{A_1}\bigl(\omega(S^2)-\omega_2(S^2)\bigr)=I_2 +\int_{A_1}\omega(S^2). 
  \end{array}
\end{equation*}
The rational numbers $I_1$ and $I_2$ are not changed when $(K,\tau)$ moves continuously so that the angle
from the real positive half-line to $\bigvec{K(0)K(1)}$ varies in $\left]0,2\pi\right[$ and the trivialization varies accordingly so that $K$ remains straight. Therefore, $I(\onechordsmallnonnum,\omega_{S^2})$ varies like the area of $A_1$,
which we compute now.
Assume that the direction of $\left[K(0),K(1)\right]$ coincides with the direction of $\exp(2i\pi\theta)$ for $\theta \in \left]0,1\right[$.

If $K$ goes from bottom to bottom,
then $A_1$ is the set of vectors $v$ of $S^2$ such that $p_{\CC}(v)=\lambda \exp(2i\pi\alpha)$ for some $(\lambda,\alpha) \in \left]0,+\infty\right[\times \left]0,\theta\right[$, and the area of $A_1$ is $\theta$.
So we have $I_{\theta}(K,\tau)=I_1+I_2 + 2\theta-1$.

If $K$ goes from top to top,
then $A_1$ is the set of vectors $v$ of $S^2$ such that $p_{\CC}(v)=\lambda \exp(2i\pi\alpha)$ for some $(\lambda,\alpha) \in \left]0,+\infty\right[\times \left]\theta, 1\right[$, and the area of $A_1$ is $1-\theta$.
So we have $I_{\theta}(K,\tau)=I_1+I_2 + 1 - 2\theta$.

Let us treat the case $\theta=\frac12$ when $K$ goes from bottom to bottom. With the notation of Definition~\ref{defselflkqtanggen}, set $\theta_1=-\theta_0=\frac12$.  
So $\hat{K}_{\parallel,\theta_0,\theta_1}$ looks as in Figure~\ref{figgammapar}, and it is isotopic to the parallel $-\widehat{(-K)}_{\parallel}$, which looks as in Figure~\ref{figgammaparthetazero}.
In this case, the interior of $A_1$ contains the direction $i$. We can assume that $\omega_1(S^2)=\omega_{0,0}$ and compute
$lk(K, K_{\parallel})=lk(C, C_{\parallel}) $ for a pair  $(C,C_{\parallel})$ of parallel closed curves that coincides with $(K,-(-K)_{\parallel})$ in a big neighborhood of $\hcylc$ and that lies in a vertical plane orthogonal to $i$ outside that big neighborhood. We get
\begin{equation*}\begin{array}{ll}lk(K, K_{\parallel})&=lk(C, C_{\parallel})=\int_{K\times K \setminus \diag}\omega_{0,0}= I\left(\onechordsmallnonnum,\omega_{0,0}\right)+I\left(\onechordsmallnonnumorrev,\omega_{0,0}\right) \\&=I\left(\onechordsmallnonnum,\omega_{0,0}\right)+I\left(\onechordsmallnonnum,-\iota^{\ast}(\omega_{0,0})\right)=I_1+I_2=I_{\theta}(K,\tau).\end{array}\end{equation*}
So Proposition~\ref{proplktangles} holds in the case $\theta=\frac12$ when $K$ goes from bottom to bottom. It similarly holds in the case $\theta=\frac12$ when $K$ goes from top to top.

We now deduce the case in which $\theta \in \left]0,1\right[$ from the case $\theta=\frac12$.
Recall that the rational numbers $I_1$ and $I_2$ are unchanged under isotopies that make $\theta$ vary continuously in $\left]0,1\right[$. Define $\theta_1=\theta$ and $\theta_0=\theta-1$,
so that the arcs $\alpha_0$ and $\alpha_1$ of Definition~\ref{defselflkqtanggen} vary continuously as $\theta$ varies from $0$ to $1$.
Then the isotopy class of the pair $(\widehat{K(\theta)}, \widehat{K(\theta)}_{\parallel,\theta_0,\theta_1})$ is unchanged. So when $K$ goes from bottom to bottom, we have $lk(K(\theta), K(\theta)_{\parallel})=lk(K(\frac12), K(\frac12)_{\parallel})+2\theta-1$ and $I_{\theta}(K({\theta}),\tau)=I_{\theta}(K({\frac12}),\tau)+2\theta-1$.
When $K$ goes from top to top, $lk(K(\theta), K(\theta)_{\parallel})-I_{\theta}(K,\tau)$ is similarly fixed when $\theta$ varies. So Proposition~\ref{proplktangles} holds in any case.
\eop

View the anomaly $\ansothree$ of Section~\ref{secansothree} as the map from $P_{(3)}(\underline{3\largen})$ to $\Assis(\RR)$ that sends any subset of $\underline{3\largen}$ with cardinality $3k$ to $\beta_k$.

With the notation of Definition~\ref{defprodbiz}, we get the following corollary of Theorem~\ref{thmtangconstcompar}, Theorem~\ref{thmfstconsttang}, and Proposition~\ref{proplktangles}.

\begin{theorem}
 \label{thmdefsanstauvarzinvf}
Let $\tanghcyll$ be a straight tangle with respect to a parallelization $\tau$. Let $J_{bb}=J_{bb}(\tanghcyll)$ (resp. $J_{tt}=J_{tt}(\tanghcyll)$) denote the set of components of $\tanghcyll$ going from bottom to bottom (resp. from top to top). For $K \in J_{bb} \cup J_{tt}$, the \emph{orientation of $K$ induced by $\tau$} is the orientation of $K$ such that $p_{\tau}(\ST^+K) \subseteq S_{WE}$. 
Under the assumptions of Theorem~\ref{thmtangconstcompar}, \begin{equation*}\biggl(\Zinv\Bigl(\hcylc,\tanghcyll,\tau,.,\bigl(\tilde{\omega}(i,1)\bigr)_{i \in \underline{3\largen}}\Bigr)\exp\Bigl(-\frac14 p_1(\tau)\ansothree(.)\Bigr)\biggr)_{\sqcup}\end{equation*}
depends only on the $\tilde{\omega}(i,1,S^2)$, on the boundary-fixing diffeomorphism class of $(\hcylc,\tanghcyll)$, on the orientations of the components of $J_{bb} \cup J_{tt}$ induced by $\tau$, and on the parallel $\tanghcyll_{\parallel}$ of $\tanghcyll$ induced by $\tau$.

It is denoted by $\Zinvufrfneg(\hcylc,\tanghcyll,\tanghcyll_{\parallel},.,(\tilde{\omega}(i,1,S^2))_{i \in \underline{3\largen}})$.
Set
\begin{equation*}\Zinvufrf_{\leq \largen}(\hcylc,\tanghcyll,\tanghcyll_{\parallel},.)=\Zinvufrfneg\bigl(\hcylc,\tanghcyll,\tanghcyll_{\parallel},.,(\omega_{S^2})_{i \in \underline{3\largen}}\bigr).\end{equation*} Then $\Zinvufrf_{\leq \largen}(\hcylc,\tanghcyll,\tanghcyll_{\parallel},.)$
maps any subset of $\underline{3\largen}$ of cardinality $3k$ to the degree $k$ part $\Zinvufrf_{k}(\hcylc,\tanghcyll,\tanghcyll_{\parallel})$ of the invariant $\Zinvufrfneg(\hcylc,\tanghcyll,\tanghcyll_{\parallel})$ of Definition~\ref{deffstconsttangframed}.  We will drop the subscript \say{$\leq \largen$} from $\Zinvufrf_{\leq \largen}(\hcylc,\tanghcyll,\tanghcyll_{\parallel},.)$.
\end{theorem}
\bp 
First note that $\Zinvufrfneg(\hcylc,\tanghcyll,\tanghcyll_{\parallel},.,(\omega_{S^2})_{i \in \underline{3\largen}})$ maps any subset of $\underline{3\largen}$ of cardinality $3k$ to  $\Zinvufrf_{k}(\hcylc,\tanghcyll,\tanghcyll_{\parallel})$ as stated.

Let $K$ be a component of $\tanghcyll$. If $K$ is a knot, or a component going from bottom to top or from top to bottom, then the form $\eta(\finseta_{K},p_{\tau}(\ST^+K))$ of Definition~\ref{defAholonom} vanishes.
Since the components of $J_{bb} \cup J_{tt}$ are equipped with the orientation induced by $\tau$, $\eta(\finseta_{K},p_{\tau}(\ST^+K))$ is the same for all components $K$ of $J_{tt}$, and it is independent of $\tanghcyll$. It is denoted by $\eta(\finseta_{K},S_{WE})$. Similarly, we have $\eta(\finseta_{K},p_{\tau}(\ST^+K))=-\eta(\finseta_{K},S_{WE})$ for all components $K$ of $J_{bb}$. Therefore, the factor $\left(\prod_{j\in I}\lol{\left[0,1\right]}{\eta(\finseta_{K_j},p_{\tau}(\ST^+K_j))} \#_j\right)$ in Theorem~\ref{thmtangconstcompar} is equal to
\begin{equation*}\left(\prod_{K_j \in J_{bb}}\biglol{\left[0,1\right]}{-\eta(\finseta_{K_j},S_{WE})} \#_j\right) \left(\prod_{K_j \in J_{tt}}\biglol{\left[0,1\right]}{\eta(\finseta_{K_j},S_{WE})} \#_j\right),\end{equation*}
and $\Zinvufrfneg(\hcylc,\tanghcyll,\tanghcyll_{\parallel},.,(\tilde{\omega}(i,1,S^2))_{i \in \underline{3\largen}})$ is determined by $\Zinvufrfneg(\hcylc,\tanghcyll,\tanghcyll_{\parallel})$ and by the holonomies, which depend only on the $\tilde{\omega}(i,.,S^2)$.
\eop

\begin{definition}
\label{defnotzfvariant}
With the notation of Theorem~\ref{thmdefsanstauvarzinvf} and under its assumptions, Theorems~\ref{thmtangconstcompar} and \ref{thmfstconsttang} together with Proposition~\ref{proplktangles} imply that
\begin{equation*}\left(\prod_{j=1}^k\Bigl(\exp\bigl(-lk(K_j,K_{j\parallel})\alpha\bigr)\#_j\Bigr)  \Zinvufrfneg\bigl(\hcylc,\tanghcyll,\tanghcyll_{\parallel},.,(\tilde{\omega}(i,1,S^2))_{i \in \underline{3\largen}}\bigr)\right)_{\sqcup}\end{equation*} is independent of the framing of $\tanghcyll$ with $\tanghcyll=(K_j)_{j \in \underline{k}}$ and $\tanghcyll_{\parallel}=(K_{j\parallel})_{j \in \underline{k}}$.  It is denoted by $\Zinvuf(\hcylc,\tanghcyll,.,(\tilde{\omega}(i,1,S^2))_{i \in \underline{3\largen}})$. It a priori depends on the orientations of the components of $J_{bb} \cup J_{tt}$.

The data of an orientation for the components of $J_{bb} \cup J_{tt}$ is called a \emph{$J_{bb,tt}$-orientation}, and $\tanghcyll$ is said to be \emph{$J_{bb,tt}$-oriented} when it is equipped with such an orientation.

All the involved products are as in Definition~\ref{defprodbiz}, and
$lk(K_j,K_{j\parallel})\alpha$ is considered as a function of subsets of $\underline{3\largen}$
with cardinality multiple of $3$, which depends only on the degree.
View the invariant of Theorem~\ref{thmfstconsttang} as such a function $\Zinvuf(\hcylc,\tanghcyll)(.)$. Then, according to Theorem~\ref{thmtangconstcompar} ---applied when when $\tilde{\omega}(i,0,S^2)=\omega_{S^2}$, we get
\begin{multline}\label{eqchangeformzinv}
\Zinvuf\Bigl(\hcylc,\tanghcyll,.,\bigl(\tilde{\omega}(i,1,S^2)\bigr)_{i \in \underline{3\largen}}\Bigr)=\\
\left(\begin{array}{r} \left(\prod_{K_j \in J_{bb}}\biglol{\left[0,1\right]}{-\eta(.,S_{WE})} \#_j\right) \left(\prod_{K_j \in J_{tt}}\biglol{\left[0,1\right]}{\eta(.,S_{WE})} \#_j\right)\\
\lol{\left[1,0\right] \times \confy^- }{\eta_{\finsetb^-,.}}\Zinvuf(\hcylc,\tanghcyll)(.)\lol{\left[0,1\right]\times \confy^+ }{\eta_{\finsetb^+,.}}\end{array}\right)_{\sqcup}.
\end{multline}

This allows us to extend
the definition of $\Zinvufrfneg(.,.,(\tilde{\omega}(i,1,S^2))_{i \in \underline{3\largen}})$ for $J_{bb,tt}$-oriented framed tangles that are not represented by straight tangles so that the equality 
\begin{multline*}\Zinvufrfneg\Bigl(\hcylc,\tanghcyll,\tanghcyll_{\parallel},.,\bigl(\tilde{\omega}(i,1,S^2)\bigr)_{i \in \underline{3\largen}}\Bigr)\\=
\left(\prod_{j=1}^k\Bigl(\exp\bigl(lk(K_j,K_{j\parallel})\alpha\bigr)\#_j\Bigr)  \Zinvuf\Bigl(\hcylc,\tanghcyll,.,\bigl(\tilde{\omega}(i,1,S^2)\bigr)_{i \in \underline{3\largen}}\Bigr)\right)_{\sqcup}\end{multline*}
holds for all these $J_{bb,tt}$-oriented framed tangles, and we have
\begin{multline} \label{eqchangeformzinvrf}
\Zinvufrfneg\Bigl(\hcylc,\tanghcyll,\tanghcyll_{\parallel},.,\bigl(\tilde{\omega}(i,1,S^2)\bigr)_{i \in \underline{3\largen}}\Bigr)\\=
\left(\begin{array}{r} \left(\prod_{K_j \in J_{bb}}\biglol{\left[0,1\right]}{-\eta(.,S_{WE})} \#_j\right) \left(\prod_{K_j \in J_{tt}}\biglol{\left[0,1\right]}{\eta(.,S_{WE})} \#_j\right)\\
\lol{\left[1,0\right] \times \confy^- }{\eta_{\finsetb^-,.}}\Zinvufrfneg(\hcylc,\tanghcyll,\tanghcyll_{\parallel})(.)\lol{\left[0,1\right]\times \confy^+ }{\eta_{\finsetb^+,.}}\end{array}\right)_{\sqcup}.
\end{multline}
\end{definition}

In order to compute $\Zinvufrfneg(\hcylc,\tanghcyll,\tanghcyll_{\parallel},.,.)$ from the discretizable definition of $\Zinv(\hcylc,\tanghcyll,\tau,.,.)$ in Theorem~\ref{thmtangconstcompar}, we first represent $\tanghcyll$ as a straight tangle with another induced parallel $\tanghcyll^{\prime}=(K_j^{\prime})_{j \in \underline{k}}$ (but with the same $J_{bb,tt}$-orientation induced by the parallelization), and we correct by setting
\begin{multline*} \Zinvufrfneg\Bigl(\hcylc,\tanghcyll,\tanghcyll_{\parallel},.,\bigl(\tilde{\omega}(i,1,S^2)\bigr)_{i \in \underline{3\largen}}\Bigr)\\=
\left(\prod_{j=1}^k\Bigl(\exp\bigl(lk(K_j,K_{j\parallel}-K_j^{\prime})\alpha\bigr)\#_j\Bigr)  \Zinvufrfneg\Bigl(\hcylc,\tanghcyll,\tanghcyll^{\prime},.,\bigl(\tilde{\omega}(i,1,S^2)\bigr)_{i \in \underline{3\largen}}\Bigr)\right)_{\sqcup}.\end{multline*}

\begin{remark}
\label{rknotzfvariant}
Definition~\ref{defnotzfvariant} is not canonical because of the arbitrary choice of $S_{WE}$. The defined invariant may not have the usual natural dependence on the component orientations (as in Proposition~\ref{proporcomp}). 
Indeed, Definition~\ref{defnotzfvariant} involves the $J_{bb,tt}$-orientation. So it is not symmetric under the reversal of a component orientation.
See Remark~\ref{rknotzfvarianttwo} for further explanations.
\end{remark}

\begin{lemma}\label{lemholzinftwo}
Let $(\omega(i))_{i\in \underline{3\largen}}$ denote a fixed family of propagating forms of $(C_2(\rats(\hcylc)),\tau)$. 
These propagators may be expressed as $\tilde{\omega}(i,1)$ for forms $\tilde{\omega}(i)$ as in Theorem~\ref{thmtangconstcompar} (thanks to Lemma~\ref{lemformprod}).
Set
\begin{equation*}\Zinv(\hcylc,\tanghcyll,\tau,\finseta)=\Zinv(\hcylc,\tanghcyll,\tau,\finseta)(1)\end{equation*} 
with the notation of Theorem~\ref{thmtangconstcompar}.
Let $h_t$ be an isotopy of $\crats(\hcylc)$ that is the identity on $\left(\CC \setminus \drad{1} \right)\times \RR$ for any $t$ and that restricts to $\CC \times \left]-\infty,0\right]$ and $\CC \times \left[1,+\infty\right[$ as isotopies $h^-_t \times \id_{\left]-\infty,0\right]}$ and $h^+_t \times \id_{\left[1,+\infty\right[}$, for planar isotopies $h^-_t$ and $h^+_t$. Assume $h_0=\id$.
Let $\tanghcyll$ be a long tangle of $\crats(\hcylc)$ whose bottom (resp. top) configuration is represented by a map $\confy^- \colon \finsetb^- \to \drad{1}$ (resp. $\confy^+ \colon \finsetb^+ \to \drad{1}$). Let $(\tau_t)_{t\in \left[0,1\right]}$ be a smooth homotopy of parallelizations of $\hcylc$ such that $p_{\tau_t}\vert_{\ST^+(\tanghcyll)}$ is constant with respect to $t$.
With the notation of Definition~\ref{defprodbiz},
for $\finseta \in P_{(3)}(\underline{3\largen})$, set
 $\Zinv(t,\finseta)=\Zinv(\hcylc,h_t(\tanghcyll),\tau_t,\finseta)$. Then we have
 \begin{equation*}\Zinv(t,.)=\left(\lol{h^-_{\left[t,0\right]}\circ \confy^-}{\eta_{\finsetb^-,.}}\Zinv(0,.)\lol{h^+_{\left[0,t\right]}\circ \confy^+}{\eta_{\finsetb^+,.}}\right)_{\sqcup}. \end{equation*}
\end{lemma}
\bp The proof is similar to that of Lemma~\ref{lemholzinf}. \eop

Lemma~\ref{lemholzinftwo} implies that $\bigl(\lol{h^-_{\left[t,0\right]}\circ \confy^-}{\eta_{\finsetb^-,.}}\lol{h^-_{\left[0,t\right]}\circ \confy^-}{\eta_{\finsetb^-,.}}\bigr)_{\sqcup}$ is neutral for the product of Definition~\ref{defprodbiz}.

The following proposition can be proved as Proposition~\ref{propfuncbraid}.

\begin{proposition} \label{propfuncbraidtwo}
With the notation and assumptions of Theorem~\ref{thmtangconstcompar},
when \begin{equation*}\gamma \colon \left[0,1\right] \to \check{C}_{\finsetb}\left[\drad{1}\right]\end{equation*} is smooth with vanishing derivatives at $0$ and $1$, deform the standard parallelization of $\RR^3$ to a homotopic parallelization $\tau$ such that $T(\gamma)$
 is straight with respect to $\tau$ at any time of the homotopy. For any subset $\finseta$ of $\underline{3\largen}$ with cardinality $3k$, we have
\begin{equation*}\Zinv\Bigl(\hcylc_0=\drad{1}\times\left[0,1\right],T(\gamma),\tau,\finseta,\bigl(\tilde{\omega}(i,1,S^2)\bigr)_{i\in \underline{3\largen}}\Bigr)=\lol{\{1\}\times p_{CS} \circ \gamma}{\eta_{\finsetb,\finseta}},\end{equation*}
where $p_{CS}$ is the natural projection $\check{C}_{\finsetb}\left[\drad{1}\right] \to \ccompuptd{\finsetb}{\CC}$.
Thus we also have \begin{equation*}\begin{array}{ll}\Zinvuf\Bigl(\gamma,\finseta,\bigl(\tilde{\omega}(i,1,S^2)\bigr)_{i\in \underline{3\largen}}\Bigr)&=\Zinvufrfneg\Bigl(\hcylc_0,T(\gamma),T(\gamma)_{\parallel},\finseta,\bigl(\tilde{\omega}(i,1,S^2)\bigr)_{i\in \underline{3\largen}}\Bigr)\\&=\lol{\{1\} \times p_{CS} \circ \gamma}{\eta_{\finsetb,\finseta}}.\end{array}\end{equation*}
\end{proposition}

Proposition~\ref{propfuncbraid} could be fully generalized to this setting, too, but we will prove more general functoriality properties in Section~\ref{secfunc}.

\begin{lemma} \label{lemthmconnecgenzero}
Let $u \in \left[0,1\right]$. Let
 $\gamma \colon \left[0,1\right] \to \{u\} \times \cinjuptd{\finsetb}{\CC}$ be a smooth path. We have the following properties.
 \begin{itemize}
  \item $\lol{\gamma}{\eta_{\finsetb,\finseta}}$ depends only on $\gamma(0)$, $\gamma(1)$, the $\tilde{\omega}(i,S^2)$ for $i \in \finseta$, and the homotopy class of $\gamma$ relatively to $\partial \gamma$ in $\{u\} \times \cinjuptd{\finsetb}{\CC}$.
  \item If $(\overline{\gamma} \colon t \mapsto \gamma(1-t))$ denotes the inverse of $\gamma$ with respect to the path composition, then 
$\bigl(\lol{\gamma}{\eta_{\finsetb,.}}\lol{\overline{\gamma}}{\eta_{\finsetb,.}}\bigr)_{\sqcup}$ is neutral with respect to the product of Definition~\ref{defprodbiz}.
\item Let $w$, $t$, and $\varepsilon$ be three elements of $\left[0,1\right[$, such that 
$\left[w,w+ \varepsilon\right]$ and $\left[t,t+\varepsilon\right]$ are subsets of $\left[0,1\right]$. Let $\ell$ be the boundary of the square $\left[w,w+ \varepsilon\right] \times p_{CS} \circ\gamma\left(\left[t,t+\varepsilon\right]\right)$ of $\left[0,1\right] \times \cinjuptd{\finsetb}{\CC}$, then $\lol{\ell}{\eta_{\finsetb,.}}$ is trivial.
 \end{itemize}
\end{lemma}
\bp 
The first assertion is a direct consequence of Proposition~\ref{propfuncbraidtwo}.
It implies that for any $u \in \left[0,1\right]$, for any homotopically trivial loop $\ell$ of $\{u\} \times \cinjuptd{\finsetb}{\CC}$,  $\lol{\ell}{\eta_{\finsetb,.}}$ is neutral with respect to the product of Definition~\ref{defprodbiz}. Since $\gamma \overline{\gamma}$ is such a loop, it implies the second assertion.
The third assertion is a direct consequence of Lemma~\ref{lemholzinftwo} and Proposition~\ref{propfuncbraidtwo}.
\eop

\section{Semi-algebraic structures on some configuration spaces}
\label{secsemialg}

We would like to extend the definitions of our connections $\eta$ of Sections~\ref{seconeform} and \ref{secdischol} on $\ccompuptd{\finsetb}{\CC}$, in order to extend the definition of $\Zinvufrf$ to $q$-tangles.
Unfortunately, I do not know whether the connections $\eta$ extend as differentiable forms on $\ccompuptd{\finsetb}{\CC}$.
However, we will be able to extend the definitions of their holonomies and prove that these holonomies
along paths make sense (as $\sqrt{t^{-1}}$ may be integrated on $\left[0,1\right]$ though $\sqrt{t^{-1}}$ is not defined at $0$).
In order to do that, we will need to prove that integrals over singular spaces converge absolutely. Our proofs rely on the theory of semi-algebraic sets. We review
the results of this theory that we will use below. Our primary reference is 
\cite[Section 1.4 and Chapter 2]{BCR} by Jacek Bochnak, Michel Coste, and Marie-Fran{\c{c}}oise Roy.

\begin{definition}\cite[Definition 2.1.4]{BCR}
\label{defsemialgset}
 A \emph{semi-algebraic subset} of $\RR^n$ is a subset of the form 
\begin{equation*}\bigcup_{i=1}^s\left(\bigcap_{j=1}^{r_i}\bigl\{x \in \RR^n \suchthat  f_{i,j}(x) <0\bigr\} \cap \bigcap_{j=r_i+1}^{s_i}\bigl\{x \in \RR^n \suchthat  f_{i,j}(x)=0\bigr\}\right) \end{equation*}
for an integer $s$, $2s$ integers $r_1$, \dots , $r_s$, $s_1$, \dots, $s_s$, such that $s_i\geq r_i$ for any $i \in \underline{s}$, and $\sum_{i=1}^ss_i$ real polynomials $f_{i,j}$ in the natural coordinates of $x$.
A \indexT{semi-algebraic set} is a semi-algebraic subset of $\RR^n$ for some $n \in \NN$.
\end{definition}

The set of semi-algebraic subsets of $\RR^n$ is obviously stable under finite union, finite intersection, and taking complements. 
The set of semi-algebraic sets is stable under finite products.

\begin{theorem}
\label{thmbcrone} Semi-algebraic sets also satisfy the following deeper properties, proved in \cite{BCR}.
\begin{description} 
\item \cite[Theorem 2.2.1]{BCR} Let $S$ be a semi-algebraic subset of $\RR^{n+1}$. Let $\Pi \colon \RR^{n+1} \to \RR^n$ be the projection onto the space of the first $n$ coordinates. Then $\Pi(S)$ is a semi-algebraic subset of $\RR^n$.
\item \cite[Proposition 2.2.2]{BCR} The closure and the interior of a semi-algebraic set are semi-algebraic sets.
\end{description}
\end{theorem}
\eopwobp

\begin{definition} \cite[Definition 2.2.5]{BCR}
 A map from a semi-algebraic subset of $\RR^n$ to a semi-algebraic subset of $\RR^m$ is \emph{semi-algebraic} if its graph is semi-algebraic in $\RR^{n+m}$.
\end{definition}

The following proposition \cite[Proposition 2.2.7]{BCR} can be deduced from Theorem~\ref{thmbcrone} above as an exercise.
\begin{proposition}
Let $f$ be a semi-algebraic map from a semi-algebraic set $A$ to a semi-algebraic subset of $\RR^n$. For any semi-algebraic subset $S$ of $A$,  $f(S)$ is semi-algebraic. For any semi-algebraic subset $S$ of $\RR^n$, $f^{-1}(S)$ is semi-algebraic. The composition of two composable semi-algebraic maps is semi-algebraic.
\end{proposition}
\eopwobp

As an example, which will be useful very soon, we prove the following proposition.

\begin{proposition}
\label{propsemialgstrucconfS}
Let $\finsetv$ denote a finite set of cardinality at least $2$.
Let ${\vecspt}$ be a vector space of dimension $\dimdel$. The manifold $\ccompuptd{\finsetv}{\vecspt}$ of Theorem~\ref{thmcompfacanom} has a canonical structure of a semi-algebraic set. For any subset $\finseta$ of $\finsetv$ with $\cardlef{\finseta} \geq 2$, the restriction map $\ccompuptd{\finsetv}{\vecspt} \to \ccompuptd{\finseta}{\vecspt}$ is semi-algebraic with respect to the canonical structures. 
\end{proposition}
\bp The charts of Lemma~\ref{lemnormfibdiag} provide canonical semi-algebraic structures on $\cuptd{\finsetv}{\vecspt}$ and $\cinjuptd{\finsetv}{\vecspt}$. The restriction maps from $\cinjuptd{\finsetv}{\vecspt}$ to $\cinjuptd{\finseta}{\vecspt}$ are semi-algebraic with respect to these structures.
The description of $\ccompuptd{\finsetv}{\vecspt}$ as the closure of the image of $\cinjuptd{\finsetv}{\vecspt}$ in $\prod_{\finseta \in \CP_{\geq 2}}\cuptd{\finseta}{\vecspt}$ of Lemma~\ref{lemdescclosccompuptd} makes clear that $\ccompuptd{\finseta}{\vecspt}$ has a natural semi-algebraic structure, thanks to Theorem~\ref{thmbcrone}.
\eop

We also have the following easy lemma.

\begin{lemma}
\label{lemcompcalV}
The space $\check{\CV}(\Gamma)$ of Chapter~\ref{chapzinvfbraid} and its compactification $\CV(\Gamma)$ carry natural structures of semi-algebraic sets. The projection $p_{\ccompuptdempty{\finsetb}}$ from $\CV(\Gamma)$ to $\ccompuptd{\finsetb}{\CC}$ and its projections to the $\ccompuptd{e}{\RR^3}$ for ordered pairs $e$ of $V(\Gamma)$ are semi-algebraic maps with respect to these structures.
For any configuration $\confy \in \ccompuptd{\finsetb}{\CC}$, the spaces $\check{\CV}(\confy,\Gamma)$ and ${\CV}(\confy,\Gamma)$ are semi-algebraic.
\end{lemma}
\eopwobp

\begin{lemma} \label{lemdifsemialg} 
 Let $f \colon \left]a,b\right[^d \to \RR$ be a $C^1$ semi-algebraic map. Then its partial derivatives $\frac{\partial f}{\partial x_i}$ are semi-algebraic functions.
\end{lemma}
\bp Let $(e_1, \dots, e_d)$ be a basis of $\RR^d=\{\sum_{i=1}^d x_ie_i\}$.
As in \cite[Proposition 2.9.1]{BCR}, note that the set 
\begin{equation*}\biggl\{\Bigl(t,x,f(x),\bigl(f(x+te_i)-f(x)\bigr)/t\Bigr) \suchthat  t \in \left]0,1\right],x \in \left]a,b\right[^d, x+te_i \in \left]a,b\right[^d\biggr\}\end{equation*}
is semi-algebraic. So are its closure, the locus $(t=0)$ of this closure, and its projection to the graph of the partial derivative of $f$ with respect to $x_i$.
\eop

\begin{lemma}
\label{lemcritsetsemialg} 
 Let $f$ be a semi-algebraic smooth map from an open hypercube $\left]0,1\right[^d$ to
$\RR^n$.
Then the \emph{critical set} of $f$, which is the subset of $\left]0,1\right[^d$ for which $f$ is not a submersion, is semi-algebraic.
\end{lemma}
\bp 
According to Lemma~\ref{lemdifsemialg}, the partial derivatives $\frac{\partial p_j \circ f}{\partial x_i}$ with respect to the factors of $\RR^d$ of the $p_j \circ f$ for the projections $p_j$ on the factors of $\RR^n$ are semi-algebraic. It is easy to see that the product and the sum of two real-valued semi-algebraic maps are semi-algebraic. Being in the critical set may be written as \say{For any subset $I$ of $\underline{d}$ of cardinality $n$, the determinant $\det\bigl[\frac{\partial p_j \circ f}{\partial x_i}(x)\bigr]_{i \in I, j \in \underline{n} }$ is equal to zero.}
\eop

An essential property of semi-algebraic sets, which we are going to use, is the following 
decomposition theorem \cite[Proposition 2.9.10]{BCR}.

\begin{theorem}
\label{thmdecsemialg}
Let $S$ be a semi-algebraic subset of $\RR^n$. Then, as a set, $S$ is the disjoint union of finitely many smooth semi-algebraic submanifolds, each semi-algebraically diffeomorphic to an open hypercube $\left]0,1\right[^d$.
\end{theorem}

 The \emph{dimension} of a semi-algebraic set is the maximal dimension of a hypercube in a decomposition as above. It is proved in \cite[Section 2.8]{BCR} that it does not depend on the decomposition. According to \cite[Theorem 2.8.8]{BCR}, the dimension of the image of a semi-algebraic set of dimension $d$ under a semi-algebraic map is smaller than or equal to $d$. According to \cite[Proposition 2.8.13]{BCR}, if $A$ is a semi-algebraic set of dimension $\dim(A)$, then we have $\dim(\overline{A}\setminus A) < \dim(A)$.

The following lemma is a corollary of Theorem~\ref{thmdecsemialg}.

\begin{lemma}
\label{lemconvintsemalg}
 Let $f$ be a continuous semi-algebraic map from a compact semi-algebraic set $A$ of dimension $d$ to a semi-algebraic smooth manifold $B$ with boundary.
 Let $\omega$ be a smooth differential form of degree $d$ on $B$.
Assume that the restriction of $f$ to each piece of a decomposition as in Theorem~\ref{thmdecsemialg} is smooth.
Then the integrals $\int_{\Delta}f^{\ast}(\omega)$ of $f^{\ast}(\omega)$ over the open pieces $\Delta$ of dimension $d$ of such a decomposition converge absolutely, and $\int_{A}f^{\ast}(\omega)$ is well defined to be the sum of these $\int_{\Delta}f^{\ast}(\omega)$.
\end{lemma}
\bp 
It suffices to prove the lemma when $\omega$ is supported on a subset
$[-1,1]^n$ of an open subset of $B$ semi-algebraically diffeomorphic to $B_{k,n}=\left]-2,1\right]^k \times\left]-2,2\right[^{n-k}$. Indeed, using
a partition of unity allows us to write $\omega$ as a finite sum of such forms around
the compact $f(A)$. This allows us to reduce the proof to the case
$B=B_{k,n}$. Now a degree $d$ differential form on $B_{k,n}$ is a sum over the
parts $J$ of cardinality $d$ of $\underline{n}$
 of pull-backs of degree $d$ forms on $B_J=B_{J,k,n}=\left]-2,1\right]^{\underline{k} \cap J} \times\left]-2,2\right[^{(\underline{n} \setminus  \underline{k}) \cap J}$ multiplied by smooth functions on $B_{k,n}$, which are bounded on their compact supports. This
allows us to reduce the proof to the case in which $\omega$ is such a pull-back of a form $\omega_J$ on $B_J$, under the projection $p_J \colon B_{k,n} \to B_J$, multiplied by a bounded function $g_J$ on $B_{k,n}$.

Decompose $A$ as in Theorem~\ref{thmdecsemialg}. It suffices to prove that the integral of $\omega$ over each hypercube $H$ of dimension $d$ converges absolutely. Let $f_J$ denote $p_J \circ f$.
Consider the closure $\overline{H} \subset A$ of the hypercube $H$ in $A$, set $\partial H = \overline{H} \setminus H$. Then $\partial H$ and its image $ f_J(\partial H)$ in $B_J$ are algebraic subsets of $B_J$ of dimension less than $d$, according to the dimension properties recalled before the lemma. 
Therefore, the form $f_J^{\ast}(\omega_J)$ vanishes on the dimension $d$ pieces of the intersection of $H$ with the semi-algebraic compact set $f_J^{-1}(f_J(\partial H))$.
Let $\Sigma(f_J)$ be the set of critical points of $f_{J}\vert_{H}$. 
According to Lemma~\ref{lemcritsetsemialg}, $\Sigma(f_J)$ is semi-algebraic.
According to the Morse--Sard theorem~\ref{thmMorseSard}, $f_J(\Sigma(f_J))$, which is semi-algebraic, is of zero measure. Therefore, its dimension is less than $d$.
Now, $B_J \setminus (f_J(\Sigma(f_J) \cup \partial H) \cup \partial B_J)$ is an open semi-algebraic subset of $B_J$, which therefore has a finite number of connected components according to Theorem~\ref{thmdecsemialg}. On each of these connected components, the local degree of $f_J$ is finite because $\overline{H}$ is compact and the points in the preimage of a regular value are isolated. Our assumptions make this local degree locally constant.
Indeed, for a point $y$ in such a component, there exists a small $d$-dimensional disk $D(y)$ around $y$ whose preimage contains disk neighborhoods of the points of the preimage, each mapped diffeomorphically to $D(y)$. The image of $\overline{H}$ minus these open disks is a compact that does not meet $y$. Therefore, there is a smaller disk around $y$ that is not met by this compact.

Then $\int_{\overline{H}}f^{\ast}(\omega)$ is the integral of $\omega_J$ weighted by this bounded local degree and by a multiplication by $g_J \circ f$. So it is absolutely convergent.
\eop

Recall that an \emph{open simplex} in $\RR^n$ is a subset of the form $v_1\dots v_k=\{\sum_{i=1}^k t_iv_i \suchthat  t_i \in \left]0,1\right], \sum_{i=1}^k t_i=1\}$, where $v_1,\dots, v_k$ are affinely independent points in $\RR^n$. The faces of $v_1\dots v_k$ are the simplices $v_{i_1}\dots v_{i_j}$ for subsets $\{i_1,\dots,i_j\} \subset \underline{k}$.
A \emph{locally finite simplicial complex} in $\RR^n$ is a locally finite 
collection $K$ of disjoint open simplices such that each face of a simplex of $K$ belongs to $K$. For such a complex $K$, $|K|$ denotes the union of the simplices of $K$.

The following Lojasiewicz triangulation theorem \cite[Theorem 1, p. 463, \S 3]{Lojasiewicz} ensures that a compact semi-algebraic set may be viewed as a topological chain (as in Subsection~\ref{subalgint}).

\begin{theorem}
\label{thmloja} 
For any locally finite collection $\{B_i\}$ of semi-algebraic subsets of $\RR^n$, there exist a locally finite simplicial complex $K$ of $\RR^n$ such that $|K|=\RR^n$ and a homeomorphism $\tau$ from $\RR^n$ to $\RR^n$ such that
\begin{itemize}
 \item for any open simplex $\sigma$ of $K$, $\tau(\sigma)$ is an analytic submanifold of $\RR^n$, and $\tau\vert_{\sigma}$ is an analytic isomorphism from $\sigma$ to $\tau(\sigma)$,
\item for any open simplex $\sigma$ of $K$ and any $B_i$ of the collection $\{B_i\}$, we have $\tau(\sigma) \subset B_i$ or $\tau(\sigma) \subset \RR^n \setminus B_i$.
\end{itemize}
\end{theorem}

\section{Extending \texorpdfstring{$\Zinvufrf$}{Zf} to \texorpdfstring{$q$-tangles}{q-tangles}}
\label{zinvfqbraid}

In Chapter~\ref{chapzinvfbraid}, the behavior of $\Zinvuf$ on braids, which are paths in $\cinjuptd{\finsetb}{\CC}$ for some finite set $\finsetb$ was discussed. Recall that $\Zinvuf$ and $\Zinvufrf$ coincide for braids. In this section,
we extend $\Zinvufrf$ and its variants of Section~\ref{secvarzinvf} to paths of $\ccompuptd{\finsetb}{\CC}$, where $\finsetb$ is a finite set.
This is already mostly done in \cite{poirierv2}, where the main ideas come from. However, our presentation
is different, and it provides additional statements and explanations.

Our extension to paths of $\ccompuptd{\finsetb}{\CC}$ will allow us to define the extension of $\Zinvufrf$ to $q$-tangles in rational homology cylinders so that Proposition~\ref{propfuncbraid} is still valid in the setting of $q$-tangles.

Recall the semi-algebraic subsets ${\CV}(\Gamma)$ and ${\CV}(\confy,\Gamma)$ of $\ccompuptd{\vertsetv(\Gamma)}{\RR^3} \times \ccompuptd{\finsetb}{\CC}$, introduced in Chapter~\ref{chapzinvfbraid}, for
a $\underline{3\largen}$-numbered Jacobi diagram $\Gamma$ on $\sqcup_{\eltb \in \finsetb}\RR_{\eltb}$.
Both $\ccompuptd{\vertsetv(\Gamma)}{\RR^3}$ and $\ccompuptd{\finsetb}{\CC}$ are stratified by $\Delta$-parenthesizations according to Theorem~\ref{thmcompuptd}.
Let $\parentp_{\finsetb}$ be a $\Delta$-parenthesization of $\finsetb$. Let $\parentp$
be a $\Delta$-parenthesization of $\vertsetv(\Gamma)$. Set \begin{equation*}{\CV}_{\parentp_{\finsetb},\parentp}(\Gamma)={\CV}(\Gamma) \cap \left(\ccompuptd{\vertsetv(\Gamma),\parentp}{\RR^3} \times \ccompuptd{\finsetb,\parentp_{\finsetb}}{\CC}\right).\end{equation*}
An element of $\ccompuptd{\vertsetv(\Gamma),\parentp}{\RR^3}$ is denoted by $(c_{\finsety})_{{\finsety} \in \parentp}$, where $c_{\finsety} \in \cinjuptd{\kids(\finsety)}{\RR^3}$.
An element of $\ccompuptd{\finsetb,\parentp_{\finsetb}}{\CC}$ is denoted by $(y_{\finsetd})_{{\finsetd} \in \parentp_{\finsetb}}$, where $y_{\finsetd} \in \cinjuptd{\kids(\finsetd)}{\CC}$.
Fix $\Gamma$ and $\parentp_{\finsetb}$.
Recall the natural map $p_{\finsetb} \colon U(\Gamma)\to \finsetb$ induced by $i_{\Gamma}$.
Let $\confc=((c_{\finsety})_{{\finsety} \in \parentp}, (y_{\finsetd})_{{\finsetd} \in \parentp_{\finsetb}}) \in {\CV}_{\parentp_{\finsetb},\parentp}(\Gamma)$, and let ${\finsety}$ be in the set $\widehat{\parentp^{\prime}}_{\cutx}$ introduced in  Notation~\ref{notstrucCV}.
Then we have \begin{equation*}p_{\CC} \circ c_{\finsety}\vert_{U(\Gamma) \cap \finsety}=\lambda(\finsety) \left(y_{\hat{\finsetb}(\finsety)} \circ p_{\finsetb}- y_{\hat{\finsetb}(\finsety)} \circ p_{\finsetb}\bigl(b(\finsety)\bigr)\right)\end{equation*}
for some $\lambda(\finsety) \geq 0$, with the normalizations of Notation~\ref{notstrucCV}.
Recall the set $\parntppx(\confc)$ of separating sets $\finsety$ for $\confc$ from Lemma~\ref{lemdimstratCV}. Such a separating set is an element $\finsety$ of $\widehat{\parentp^{\prime}}_{\cutx}$ with $\lambda(\finsety) \neq 0$.
For a subset $\parntppx$ of $\widehat{\parentp^{\prime}}_{\cutx}$, set 
${\CV}_{\parentp_{\finsetb},\parentp,\parntppx}(\Gamma)=\{\confc \in {\CV}_{\parentp_{\finsetb},\parentp}(\Gamma) \suchthat \parntppx(\confc)=\parntppx\}$.
We use the data $(\parentp_{\finsetb},\parentp,\parntppx)$ to stratify ${\CV}(\Gamma)$ (or $\left[0,1\right] \times {\CV}(\Gamma)$, whose strata will be the products by $\left[0,1\right]$ of the strata of ${\CV}(\Gamma)$ by definition).
Recall that for any $\finsetd \in \parentp_{\finsetb}$,
the elements $\finsety$ of $\parntppx$ with $\hat{\finsetb}(\finsety)=\finsetd$ are minimal with respect to the inclusion among the elements of $\parentp$ with $\hat{\finsetb}(\finsety)=\finsetd$.

For $\confy \in \ccompuptd{\finsetb,\parentp_{\finsetb}}{\CC}$, for a $\Delta$-parenthesization $\parentp$ of $\vertsetv(\Gamma)$, and for a subset $\parntppx$ of $\parentp$, set 
\begin{equation*}{\CV}(\confy,\Gamma,\parentp,\parntppx)={\CV}(\confy,\Gamma)\cap {\CV}_{\parentp_{\finsetb},\parentp,\parntppx}(\Gamma).\end{equation*}
Recall from Lemma~\ref{lemdimstratCV} that
when ${\CV}(\confy,\Gamma,\parentp,\parntppx)$ is not empty, its dimension is
$\cardlef{U(\Gamma)} +3 \cardlef{T(\Gamma)} - 1 -\cardlef{\parentp \setminus \parntppx}.$
In particular, the dimension of ${\CV}(\confy,\Gamma)$ is at most $2\cardlef{E(\Gamma)}-1$.

Fix the family $(\tilde{\omega}(i,S^2)=(\tilde{\omega}(i,t,S^2))_{t\in \left[0,1\right]})_{i \in  \underline{3\largen}}$ of closed $2$-forms on $\left[0,1\right] \times S^2$ once for all in this section, and assume that $\tilde{\omega}(i,0,S^2)=\omega_{S^2}$ for all $i$.
For an edge $e$ of $\Gamma$,
recall the map
\begin{equation*}p_{e,S^2} \colon \left[0,1\right] \times \ccompuptd{\vertsetv(\Gamma)}{\RR^3} \to \left[0,1\right] \times S^2,\end{equation*} sending $(t,\confc \in \cinjuptd{\vertsetv(\Gamma)}{\RR^3})$ to
$\left(t,p_{S^2}\bigl((\confc(v(e,1)),\confc(v(e,2)))\bigr)\right)$.
Also recall the $(2\cardlef{E(\Gamma)})$-form 
\begin{equation*}\Omega_{\Gamma}=\bigwedge_{e \in E(\Gamma)}p_{e,S^2}^{\ast}\Bigl(\tilde{\omega}\bigl(j_E(e),S^2\bigr)\Bigr)\end{equation*} over $\left[0,1\right] \times \ccompuptd{\vertsetv(\Gamma)}{\RR^3}$. This form pulls back to provide smooth forms on the smooth strata of $\left[0,1\right] \times {\CV(\Gamma)}$.

Let $\finseta$ denote a subset of $\underline{3\largen}$ with cardinality $3n$.
An \emph{ordered $r$-component $\finseta$-numbered Jacobi diagram} $\Gamma^{(r)}$ on $\sqcup_{\eltb \in \finsetb}\RR_{\eltb}$ is a degree $n$ $\finseta$-numbered Jacobi diagram $\Gamma^{(r)}$ on $\sqcup_{\eltb \in \finsetb}\RR_{\eltb}$ that has $r$ connected components $\Gamma_1$, \dots
$\Gamma_r$ and such that $i_{\Gamma}$ is represented by an injection of $V(\Gamma)$ that maps all
univalent vertices of $\Gamma_i$ before (or below) the univalent vertices of  $\Gamma_{i+1}$ for any $i \in \underline{r-1}$.
The data of such an ordered $r$-component $\finseta$-numbered Jacobi diagram $\Gamma^{(r)}$ is equivalent to the data of an $r$-tuple $(\Gamma_1,\dots,\Gamma_r)$ of $\finseta$-numbered connected Jacobi diagrams with pairwise disjoint $j_E(E(\Gamma_i))$ such that the sum of the degrees of the $\Gamma_i$ is $n$.
Let $\Davis^{e,r}_{n,\finseta}(\sqcup_{\eltb \in \finsetb}\RR_{\eltb})$ denote the set of these ordered $r$-component $\finseta$-numbered Jacobi diagrams $\Gamma^{(r)}$ on $\sqcup_{\eltb \in \finsetb}\RR_{\eltb}$.\footnote{The notation $\Davis^{e,r}_{n,\finseta}$ has a redundancy since the cardinality of $\finseta$ is $3n$. We keep the redundancy for consistency here because we will use other spaces of numbered Jacobi diagrams, where the degree is not determined by the cardinality of the set of indices, in Chapter~\ref{chappropzinvffunc}.} 

Such a diagram provides the $(2\cardlef{E(\Gamma^{(r)})})$-form
\begin{equation*}\Omega_{\Gamma^{(r)}}=\bigwedge_{i=1}^r P_{i}^{\ast}(\Omega_{\Gamma_i})\end{equation*}  on (the smooth strata of) $\prod_{i=1}^r\left(\left[0,1\right] \times {\CV(\Gamma_i)}\right)$,
where \begin{equation*}P_i \colon \prod_{i=1}^r\bigl(\left[0,1\right] \times {\CV(\Gamma_i)}\bigr) \to \left[0,1\right] \times {\CV(\Gamma_i)}\end{equation*} is the projection onto the $i$th factor.
This form is also the pullback of a smooth form on $\left[0,1\right]^r \times (S^2)^{E(\Gamma^{(r)})}$, by a semi-algebraic map.
Recall $\simplexr=\{(t_1,\dots,t_r) \in \left[0,1\right]^r  \suchthat 0\leq t_1 \leq t_2 \dots \leq t_r \}$.
A semi-algebraic path $\gamma \colon \left[0,1\right] \to \left[0,1\right] \times \ccompuptd{\finsetb}{\CC}$ induces
the semi-algebraic map 
\begin{equation*}\begin{array}{llll}\gamma^{(r)} \colon &\simplexr &\to &\bigl(\left[0,1\right] \times \ccompuptd{\finsetb}{\CC}\bigr)^r\\
   &(t_1,\dots,t_r) & \mapsto & \bigl(\gamma(t_1),\dots,\gamma(t_r)\bigr).
  \end{array}\end{equation*}
Consider the product $P_{\Gamma^{(r)}}\colon \prod_{i=1}^r\left(\left[0,1\right] \times {\CV(\Gamma_i)}\right) \to (\left[0,1\right] \times \ccompuptd{\finsetb}{\CC})^r$ of natural projections. Assume that $\gamma$ is injective.
Set $C(\Gamma^{(r)},\gamma)=P_{\Gamma^{(r)}}^{-1}(\gamma^{(r)}(\simplexr))$.
Then $C(\Gamma^{(r)},\gamma)$ is a semi-algebraic subset of $\prod_{i=1}^r\left(\left[0,1\right] \times {\CV(\Gamma_i)}\right)$ of dimension at most $2 \cardlef{E(\Gamma^{(r)})}$ whose $2 \cardlef{E(\Gamma^{(r)})}$-dimensional strata are oriented canonically, as soon as the Jacobi diagrams $\Gamma_i$ are: Fix an arbitrary vertex-orientation for the $\Gamma_i$. The set $C(\Gamma^{(r)},\gamma)$ is locally oriented as the product of the $C(\Gamma_i^{(1)},\gamma)$ for $i \in \underline{r}$.
The parameter $t_i$ replaces the translation parameter in $\CV(\gamma(t_i),\Gamma_i)$.

Define the \emph{$\finseta$-holonomy} $\lol{.}{\eta_{\finsetb,\finseta}}$ along injective semi-algebraic paths $\gamma$ of $\left[0,1\right] \times \ccompuptd{\finsetb}{\CC}$, with respect to our family $\left(\tilde{\omega}(i,S^2)\right)_{i \in  \underline{3\largen}}$, to be
\begin{equation*}\lol{\gamma}{\eta_{\finsetb,\finseta}}=\left[\emptyset\right] +\sum_{r=1}^{\infty} \sum_{\Gamma^{(r)} \in \Davis^{e,r}_{n,\finseta}(\sqcup_{\eltb \in \finsetb}\RR_{\eltb})} \coefgambet_{\Gamma^{(r)}}\int_{C(\Gamma^{(r)},\gamma)}\Omega_{\Gamma^{(r)}}[\Gamma^{(r)}],\end{equation*}
with
\begin{equation*}\coefgambet_{\Gamma^{(r)}}=\frac{\bigl(\cardbig{\finseta} -\cardbig{E(\Gamma^{(r)})}\bigr)!}{\card{\finseta}!2^{\card{E(\Gamma^{(r)})}}}.\end{equation*}
(Again, we fix an arbitrary vertex-orientation for the components $\Gamma_i$ of each $\Gamma^{(r)}$, and $\lol{\gamma}{\eta_{\finsetb,\finseta}}$ is independent of our choices.)

The involved integrals make sense as soon as $\gamma$ is semi-algebraic, thanks to Lemma~\ref{lemconvintsemalg}, which also justifies the following lemma.

\begin{lemma}
\label{lemconvthmpoirbraid}
For any injective semi-algebraic path $\gamma$ of $\left[0,1\right] \times \ccompuptd{\finsetb}{\CC}$, \begin{equation*}\lim_{\varepsilon \to 0} \biglol{\gamma\vert_{\left[\varepsilon, 1-\varepsilon\right]}}{\eta_{\finsetb,\finseta}}\end{equation*} makes sense, and it is equal to $\lol{\gamma}{\eta_{\finsetb,\finseta}}$.
\end{lemma}
\eopwobp

Together with the identification $\Zinvuf(T(\gamma))=\hol{p_{CS} \circ \gamma}{\eta_{\finsetb}}$ for braids provided by Proposition~\ref{propfuncbraid}, Lemma~\ref{lemconvthmpoirbraid} implies the
convergence part of Theorem~\ref{thmpoirbraid}.
The above convergent integrals extend Definition~\ref{defholtilde} of $\lol{\gamma}{\eta_{\finsetb,\finseta}}$ for injective semi-algebraic paths in $\left[0,1\right] \times \cinjuptd{\finsetb}{\CC}$. Note the following easy lemma.

\begin{lemma}
\label{lemmultlol}
The $\finseta$-holonomy $\lol{\gamma}{\eta_{\finsetb,.}}$, which is valued in
$\Aavis_{n}(\sqcup_{\eltb \in \finsetb}\RR_{\eltb})$, extends naturally to noninjective semi-algebraic paths. This holonomy is multiplicative under path composition with respect to the product of Definition~\ref{defprodbiz}.
\end{lemma}
\eopwobp

Recall that $\tilde{\omega}(i,0,S^2)$ is the standard homogeneous volume-one form on $S^2$. When $\gamma$ is valued in $\{0\} \times \ccompuptd{\finsetb}{\CC}$, there is no need to number the diagram edges since $\tilde{\omega}(i,0,S^2)=\omega_{S^2}$ for all $i$, and we simply have
\begin{equation*}\lol{\gamma}{\eta_{\finsetb}}=\sum_{r=0}^{\infty} \sum_{k=0}^{\infty}\sum_{\Gamma^{(r)} \in \Davis^{e,r}_{k,\underline{3k}}(\sqcup_{\eltb \in \finsetb}\RR_{\eltb})} \coefgambet_{\Gamma^{(r)}}\int_{C(\Gamma^{(r)},\gamma)}\Omega_{\Gamma^{(r)}} [\Gamma^{(r)}]\in \Aavis(\sqcup_{\eltb \in \finsetb}\RR_{\eltb}).\end{equation*} In this case, $\lol{\gamma}{\eta_{\finsetb}}$ is nothing but the Poirier functor $Z^l$ of \cite[Section 1.4]{poirierv2} applied to $\gamma$.
The projection in 
$\Aavis_{n}(\sqcup_{\eltb \in \finsetb}\RR_{\eltb})$
of $\lol{\gamma}{\eta_{\finsetb}}$, which coincides with the holonomy $\hol{\gamma}{\eta_{\finsetb}}$ defined in Section~\ref{seconeform} when $\gamma$ is valued in $\{0\} \times \cinjuptd{\finsetb}{\CC}$, coincides with $\lol{\gamma}{\eta_{\finsetb,\finseta}}$, when $\cardlef{\finseta}=3n$, in this case of homogeneous forms.

Recall that $\ccompuptd{\finsetb}{\CC}$ is a smooth manifold with ridges, which can also be equipped with a semi-algebraic structure for which the local charts provided in Theorem~\ref{thmcompuptd} are semi-algebraic maps. 
In such a trivialized open simply connected subspace, any two points can be connected by a semi-algebraic path. In particular, any two points $a$ and $b$ of $\ccompuptd{\finsetb}{\CC}$ are connected by a semi-algebraic 
path $\gamma \colon \left[0,1\right] \to \ccompuptd{\finsetb}{\CC}$ such that $\gamma(0)=a$ and $\gamma(1)=b$.
Furthermore, any path from $a$ to $b$ of $\ccompuptd{\finsetb}{\CC}$ can be $C^0$-approximated by a homotopic semi-algebraic path. So, any homotopy class of paths from $a$ to $b$ has a semi-algebraic representative.

Now Theorem~\ref{thmpoirbraid} is a direct corollary of the following one, which will be proved after Lemma~\ref{lemthmconnecgenbis}. This theorem is a mild generalization of \cite[Proposition 9.2]{poirierv2}, thanks to Lemma~\ref{lemconvthmpoirbraid}.

\begin{theorem}
\label{thmconnecgen}
Let  $\gamma \colon \left[0,1\right] \to \left[0,1\right] \times \ccompuptd{\finsetb}{\CC}$ be a semi-algebraic path. Then
$\lol{\gamma}{\eta_{\finsetb,\finseta}}$ depends only on $\gamma(0)$, $\gamma(1)$, the $\tilde{\omega}(i,S^2)$ for $i \in \finseta$, and the homotopy class of $\gamma$ relatively to $\partial \gamma$.
\end{theorem}

According to Lemmas~\ref{lemthmconnecgenzero} and \ref{lemmultlol}, Theorem~\ref{thmconnecgen} holds for smooth paths $\gamma$ of $\{u\} \times \cinjuptd{\finsetb}{\CC}$ and their piecewise smooth compositions.
We now prove the following other particular case of Theorem~\ref{thmconnecgen}.

\begin{lemma} \label{lemthmconnecgen}
Let $\gamma \colon \left[0,1\right]^2\to \left[0,1\right] \times \ccompuptd{\finsetb}{\CC}$
be a semi-algebraic homotopy such that 
\begin{itemize}
\item $\gamma$ is injective on $\left[0,1\right] \times \left]0,1\right[$,
\item $\gamma_u(t)=\gamma(u,t)$ belongs to $\left[0,1\right] \times \cinjuptd{\finsetb}{\CC}$ for any $(u,t) \in \left]0,1\right]\times \left]0,1\right[$,
 \item we have $\gamma_u(0)=\gamma_0(0)$
and $\gamma_u(1)=\gamma_0(1)$ for all $u \in \left[0,1\right]$,

\item $\gamma_0(t)$ is in a fixed stratum of $\left[0,1\right] \times \ccompuptd{\finsetb}{\CC}$ for $t \in \left]0,1\right[$, where a \emph{stratum} of $\left[0,1\right] \times \ccompuptd{\finsetb}{\CC}$ is the product by $\left[0,1\right]$ of a stratum of $\ccompuptd{\finsetb}{\CC}$ associated to a $\Delta$-parenthesization.
\end{itemize}
Then we have
\begin{equation*}\lol{\gamma_0}{\eta_{\finsetb,\finseta}}=\lol{\gamma_1}{\eta_{\finsetb,\finseta}}.\end{equation*}
\end{lemma}

Lemma~\ref{lemthmconnecgen} is the direct consequence of Lemmas~\ref{lemcurvetaab} and \ref{lemthmconnecgenzerobis} below. The proof of Lemma~\ref{lemcurvetaab} uses Lemma~\ref{lemdetabbis} and the following sublemma.

\begin{sublemma}
Under the assumptions of Lemma~\ref{lemthmconnecgen},
let $\Gamma^{(r)}$ be an element of $\Davis^{e,r}_{n,\finseta}(\sqcup_{\eltb \in \finsetb}\RR_{\eltb})$, recall $C(\Gamma^{(r)},\gamma_u)=P_{\Gamma^{(r)}}^{-1}(\gamma_u^{(r)}(\simplexr))$,
and let
\begin{equation*}C\bigl(\Gamma^{(r)},(\gamma_u)\bigr)=\cup_{u\in \left[0,1\right]}C\bigl(\Gamma^{(r)},\gamma_u\bigr)\end{equation*} be the associated semi-algebraic set of dimension $2 \cardlef{E(\Gamma^{(r)})}+1$. Then
 the codimension-one boundary of $C(\Gamma^{(r)},(\gamma_u))$ is \begin{equation*}C(\Gamma^{(r)},\gamma_1) - C(\Gamma^{(r)},\gamma_0)
-\cup_{u\in \left[0,1\right]} \partial C(\Gamma^{(r)},\gamma_u),\end{equation*}
with \begin{equation*}\partial C\bigl(\Gamma^{(r)},\gamma_u\bigr)=\partial_C C\bigl(\Gamma^{(r)},\gamma_u\bigr) + \partial_{\Delta} C\bigl(\Gamma^{(r)},\gamma_u\bigr),\end{equation*}
\begin{equation*}\partial_C C\bigl(\Gamma^{(r)},\gamma_u\bigr)=
\pm \cup_{(t_1,\dots,t_r) \in \simplexr}\partial P_{\Gamma^{(r)}}^{-1}\Bigl(\gamma_u^{(r)}\bigl((t_1,\dots,t_r)\bigr)\Bigr),\end{equation*}
and \begin{equation*}\partial_{\Delta} C\bigl(\Gamma^{(r)},\gamma_u\bigr)=\pm P_{\Gamma^{(r)}}^{-1}\bigl(\gamma_u^{(r)}(\partial \simplexr)\bigr).\end{equation*}
\end{sublemma}
\bp 
When the image of $(\gamma_u)$ is in $\left[0,1\right] \times \cinjuptd{\finsetb}{\CC}$, it follows from Lemmas~\ref{lemstrucCVinjtree} and \ref{lemcodimonex}.
Let us prove that it is still true for our homotopies $(\gamma_u)$.
The part $C(\Gamma^{(r)},\gamma_1)$ is in the boundary as before. We can ignore the contributions of the extremities of $\gamma_1$ since they belong to parts of dimension at most $2 \cardlef{E(\Gamma^{(r)})}-1$, thanks to Lemma~\ref{lemdimstratCV}.
For the part coming from $\cup_u\partial C(\Gamma^{(r)},\gamma_u)$ in the $2 \cardlef{E(\Gamma)}$-dimensional boundary, we may restrict to $u \in \left]0,1\right[$ for dimension reasons, which we do. So this part is in the boundary as before, too

The part over $\gamma_0$ of the codimension-one boundary of $C(\Gamma^{(r)},(\gamma_u))$ is included in $C(\Gamma^{(r)},\gamma_0)$. Let us prove that the corresponding algebraic boundary is indeed $-C(\Gamma^{(r)},\gamma_0)$ when $\gamma_0\left(\left]0,1\right[\right)$ is in some stratum of $\left[0,1\right] \times \partial \ccompuptd{\finsetb}{\CC}$, associated to a parenthesization $\parentp_B$ of $\finsetb$. 

Let $t_i \in \left]0,1\right[$. In a neighborhood $[0,\eta[ \times N(t_i)$ of $(0,t_i)$ in $\left[0,1\right]^2$, $\gamma_u(t)=\gamma(u,t)$ is expressed as
\begin{equation*}\gamma_u(t)=\Bigl(\bigl(\confy_{\finsetd}(u,t)\bigr)_{{\finsetd} \in \parentp_{\finsetb}} ,\bigl(u_{\finsetd}(u,t)\bigr)_{{\finsetd} \in \parentp_{\finsetb} \setminus \{\finsetb\}}\Bigr),\end{equation*}
where we have $u_{\finsetd}(0,t)=0$ for all $t \in N(t_i)$.

Let $\confc^0$ be a point of the $(2 \cardlef{E(\Gamma_i)}-1)$-dimensional open part of $\CV(\gamma_0(t_i),\Gamma_i)$ with $\gamma^{\prime}_0(t_i)\neq 0$. According to Lemma~\ref{lemdimstratCV}, $\parentp(\confc^0)$ equals $\parntppx(\confc^0)$.
We express a neighborhood of $\confc^0$ in $\cup_{(u,t) \in [0,\eta[ \times N(t_i)}\CV(\gamma_u(t),\Gamma_i)$, as a product by $[0,\eta[ \times N(t_i)$, as follows.
We use the above coordinates $\confy_{\finsetd}(u,t)$, $u_{\finsetd}(u,t)$ of the base, and the parameters $\confc_{\finsetz}$,
$p_{\RR}\circ\confc_{\finsetz}(\finsety)$, $\lambda(\finsety)$ listed in the fourth, fifth and sixth sets of variables
of Lemma~\ref{lemstrucCV}. We have $\lambda(V(\Gamma))(\confc^0)=\lambda(\confc^0)$ and $\lambda(\confc^0) \neq 0$ since we have $\parentp=\parntppx$.
Let $\finsety$ be an element of $\parentp \setminus \{V(\Gamma)\}$. We similarly have $\lambda(\finsety)(\confc^0)\neq 0$. Let $\finsety^+$ be the smallest set of $\parentp$ that strictly contains $\finsety$. Set $\finsety^{\prime}=\finsety^+$. Set $B_1=\hat{B}(\finsety)$ and $B^{\prime}_1=B^+_1=\hat{B}(\finsety^+)$. Equation~\ref{eqcycconst} in Lemma~\ref{lemstrucCV} applied to these sets when $n=1$ may be written as
$\lambda(\finsety)\mu_{\finsety}=\lambda(\finsety^+)u_{\hat{B}(\finsety)}$. It implies that $\mu_{\finsety}=\bigl({\lambda(\finsety^+)}/{\lambda(\finsety)}\bigr)u_{\hat{B}(\finsety)}$ is determined by the listed parameters, and that all $\mu_{\finsety}$ are zero over $\gamma_0(N(t_i))$ in our neighborhood of $\confc^0$. So we have $\parentp=\parntppx=\parentp(\confc^0)$ over $\gamma_0(N(t_i))$ in our neighborhood of $\confc^0$. In particular, if $\finsety \subsetneq \finsety^{\prime}$, then $\hat{B}(\finsety) \subsetneq \hat{B}(\finsety^{\prime})$, and all equations \ref{eqcycconst} are satisfied, there. Over $\gamma(\left]0,\eta\right[ \times N(t_i))$, no $u_D$ vanishes so that no $\mu_{\finsety}$ vanishes either in our neighborhood of $\confc^0$, and the equations \ref{eqcycconst} are implied by the equations $\ast(\finsety)$, which are implied by the equations $\mu_{\finsety}=({\lambda(\finsety^+)}/{\lambda(\finsety)})u_{\hat{B}(\finsety)}$ and $\lambda(V(\Gamma))=\lambda$.
So we get a neighborhood of $\confc^0$ in $\cup_{(u,t) \in [0,\eta[ \times N(t_i)}\CV(\gamma_u(t),\Gamma_i)$ parametrized by $(u,t)$ and by the parameters $\confc_{\finsetz}$,
$p_{\RR}\circ\confc_{\finsetz}(\finsety)$, $\lambda(\finsety)$ in the fourth, fifth and sixth lists of variables
of Lemma~\ref{lemstrucCV}. For any horizontally normalizing kid $\finsety$ of $\parentp$,  remove $\lambda(\finsety)=\lambda^0(\finsety)$ from the parameters. Use the first and second constraints of Lemma~\ref{lemstrucCV} to remove other superfluous parameters $p_{\RR}\circ\confc_{\finsetz}(\finsety)$, $|p_{\CC}\circ\confc_{\finsetz}(\finsety)|$ and get a free system of parameters.
Thus, we obtain an open $(2 \cardlef{E(\Gamma_i)})$-dimensional neighborhood $\CO$ of $\confc^0$ in $\cup_{t \in N(t_i)}\CV(\gamma_0(t),\Gamma_i)$
and a local open embedding of the product of $[0,\eta[ \times \CO$ into $\cup_{(u,t) \in [0,\eta[ \times N(t_i)}\CV(\gamma_u(t),\Gamma_i)$. So
$\left(-C(\Gamma^{(r)},\gamma_0)\right)$ is the algebraic boundary of $C(\Gamma^{(r)},(\gamma_u))$ over $\gamma_0$.
\eop

\begin{lemma} \label{lemdetabbis}
For a Jacobi diagram in $\Davis^{e}_{n,\finseta}(\sqcup_{\eltb \in \finsetb}\RR_{\eltb})$, for an element $(t,\confy)$ of $\left[0,1\right] \times \cinjuptd{\finsetb}{\CC}$, $d\eta_{\Gamma}(t,\confy)$ is the 
 integral of $\bigl(-\bigwedge_{e \in E(\Gamma)}p_{e,S^2}^{\ast}(\tilde{\omega}(j_E(e),S^2))\bigr)$ along the interiors of the codimension-one faces of $\{t\} \times {\CV}(\confy,\Gamma)$.
\end{lemma}
\bp See the proof of Lemma~\ref{lemdetab}.
\eop 

Recall from the beginning of Section~\ref{seconeform} that the fiber ${\CV}(\confy,\Gamma)$ is oriented so that the orientation of ${\CV}(\confy,\Gamma)$ preceded
by the upward translation parameter---which replaces the parametrization of the paths along which we integrate---matches the usual orientation of our configuration spaces, induced as in Lemma~\ref{lemorc}.

For a subset $\finseta_i$ of $\finseta$, we have \begin{equation*}d\eta_{\finsetb,\finseta_i}=\sum_{\Gamma_i \in \Davis^e_{n_i,\finseta_i}(\sqcup_{\eltb \in \finsetb}\RR_{\eltb})}
\coefgambet_{\Gamma_i}d\eta_{\Gamma_i}\left[\Gamma_i\right],\end{equation*}
and we set
\begin{equation*}(\eta\wedge\eta)_{\finsetb,\finseta_i}=\sum_{(\Gamma,\Gamma^{\prime}) \in \Davis^{e,2}_{n_i,\finseta_i}(\sqcup_{\eltb \in \finsetb}\RR_{\eltb})} \coefgambet_{\Gamma \sqcup \Gamma^{\prime}}\eta_{\Gamma} \wedge\eta_{\Gamma^{\prime}} \left[\Gamma\right]\left[\Gamma^{\prime}\right].\end{equation*}

\begin{lemma} \label{lemcurvetaab}
Let $(\gamma_u(t))_{u \in \left[0,1\right]}$
be a semi-algebraic homotopy satisfying the assumptions of Lemma~\ref{lemthmconnecgen}
or a smooth homotopy valued in  $\left[0,1\right] \times \cinjuptd{\finsetb}{\CC}$. Assume that
$\gamma_u(0)$ and $\gamma_u(1)$ do not depend on $u$.

Then we have
\begin{multline*}\lol{\gamma_1}{\eta_{\finsetb,\finseta}}-\lol{\gamma_0}{\eta_{\finsetb,\finseta}}\\
 =\sum_{r=0}^{\infty} \sum_{(\finseta_1, \dots,\finseta_r) \in P_{r}(\finseta)} \frac{\prod_{i=1}^r\cardlef{\finseta_i}!}{\cardlef{\finseta}!}\int_{\left[0,1\right] \times \simplexr}\sum_{i=1}^r\delta(i,\finseta_1, \dots,\finseta_{r}),
\end{multline*}
with \begin{equation*}\delta(i,\finseta_1, \dots,\finseta_{r})=(-1)^{i-1}\bigwedge_{j=1}^r (\gamma \circ p_j)^{\ast}(\eta_{\finsetb,\finseta_j})\left(\frac{(\gamma \circ p_i)^{\ast}\bigl(d\eta_{\finsetb,\finseta_i} + (\eta\wedge\eta)_{\finsetb,\finseta_i}\bigr)}{(\gamma \circ p_i)^{\ast}(\eta_{\finsetb,\finseta_i})}\right),\end{equation*}
where the fraction means that $(\gamma \circ p_i)^{\ast}(\eta_{\finsetb,\finseta_i})$ is replaced by $(\gamma \circ p_i)^{\ast}(d\eta_{\finsetb,\finseta_i}+(\eta\wedge\eta)_{\finsetb,\finseta_i})$.
\end{lemma}
\bp
Set $\partial_C C(\Gamma^{(r)},(\gamma_u))=-\cup_{u \in \left]0,1\right[}\partial_C C(\Gamma^{(r)},\gamma_u)$
and \begin{equation*}\partial_{\Delta} C\bigl(\Gamma^{(r)},(\gamma_u)\bigr)=-\cup_{u \in \left]0,1\right[}\partial_{\Delta} C\bigl(\Gamma^{(r)},\gamma_u\bigr).\end{equation*}
Since \begin{equation*}\partial C\bigl(\Gamma^{(r)},(\gamma_u)\bigr)=C\bigl(\Gamma^{(r)},\gamma_1\bigr) - C\bigl(\Gamma^{(r)},\gamma_0\bigr)+\partial_C C\bigl(\Gamma^{(r)},(\gamma_u)\bigr) +\partial_{\Delta} C\bigl(\Gamma^{(r)},(\gamma_u)\bigr)\end{equation*} is a null-homologous cycle,
we get
\begin{equation*}\lol{\gamma_1}{\eta_{\finsetb,\finseta}}-\lol{\gamma_0}{\eta_{\finsetb,\finseta}}=\bigdelol{(\gamma_u)}{d\eta_{\finsetb,\finseta}} +\bigdelol{(\gamma_u)}{(\eta\wedge \eta)_{\finsetb,\finseta}},
\end{equation*}
with
\begin{equation*}\bigdelol{(\gamma_u)}{d\eta_{\finsetb,\finseta}}\stackrel{\mbox{\scriptsize def}}=-
\sum_{r=0}^{\infty} \sum_{\Gamma^{(r)} \in \Davis^{e,r}_{n,\finseta}(\sqcup_{\eltb \in \finsetb}\RR_{\eltb})} \coefgambet_{\Gamma^{(r)}}\int_{\partial_C C(\Gamma^{(r)},(\gamma_u))}\Omega_{\Gamma^{(r)}}\left[\Gamma^{(r)}\right]\end{equation*}
and 
\begin{equation*}\bigdelol{(\gamma_u)}{(\eta\wedge \eta)_{\finsetb,\finseta}}\stackrel{\mbox{\scriptsize def}}=-
\sum_{r=0}^{\infty} \sum_{\Gamma^{(r)} \in \Davis^{e,r}_{n,\finseta}(\sqcup_{\eltb \in \finsetb}\RR_{\eltb})} \coefgambet_{\Gamma^{(r)}}\int_{\partial_{\Delta} C(\Gamma^{(r)},(\gamma_u))}\Omega_{\Gamma^{(r)}}\left[\Gamma^{(r)}\right].\end{equation*}
Let us study these terms.
Since 
\begin{equation*}P_{\Gamma^{(r)}}^{-1}\Bigl(\gamma_u^{(r)}\bigl((t_1,\dots,t_r)\bigr)\Bigr)=\pm \prod_{i=1}^r \biggl( \Bigl\{p_{\left[0,1\right]}\bigl(\gamma_u(t_i)\bigr) \Bigr\}\times \CV\Bigl(p_{\cinjuptd{\finsetb}{\CC}}\bigl(\gamma_u(t_i)\bigr),\Gamma_i\Bigr)\biggr)\end{equation*} is diffeomorphic to $\pm \prod_{i=1}^r \CV(\gamma_u(t_i),\Gamma_i)$
(forgetting the natural $p_{\cinjuptd{\finsetb}{\CC}}$), we have
\begin{equation*}\partial P_{\Gamma^{(r)}}^{-1}\Bigl(\gamma_u^{(r)}\bigl((t_1,\dots,t_r)\bigr)\Bigr)=\pm \sum_{i=1}^r \partial \CV\bigl(\gamma_u(t_i),\Gamma_i\bigr) \times \prod_{j\in \underline{n}\setminus \{i\}} \CV\bigl(\gamma_u(t_j),\Gamma_j\bigr).\end{equation*}
According to Lemma~\ref{lemdetabbis}, we have
\begin{equation*}\delol{(\gamma_u)}{d\eta_{\finsetb,\finseta}}=\sum_{r=0}^{\infty} \sum_{(\finseta_1, \dots,\finseta_r) \in P_{r}(\finseta)} \frac{\prod_{i=1}^r\cardlef{\finseta_i}!}{\cardlef{\finseta}!}\int_{\left[0,1\right] \times \simplexr}\sum_{i=1}^r\alpha(i,\finseta_1, \dots,\finseta_{r}),\end{equation*}
with \begin{equation*}\alpha(i,\finseta_1, \dots,\finseta_{r})=(-1)^{i-1}\bigwedge_{j=1}^r (\gamma \circ p_j)^{\ast}(\eta_{\finsetb,\finseta_j})\left(\frac{(\gamma \circ p_i)^{\ast}(d\eta_{\finsetb,\finseta_i})}{(\gamma \circ p_i)^{\ast}(\eta_{\finsetb,\finseta_i})}\right)\end{equation*} and $d\eta_{\finsetb,\emptyset}=0$.

Let us give a few details about the involved signs. Recall that the space $C(\Gamma^{(r)},\gamma_u)$ is oriented locally as the product of the $C(\Gamma_i^{(1)},\gamma_u)$ for $i \in \underline{r}$, which are oriented so that the parameter $t$ replaces the translation parameter in $\CV(\gamma_u(t),\Gamma_i)$. Thus, the boundary $\partial_C$, along which we integrate, is
locally diffeomorphic to \begin{equation*}\left(-\cup_{u \in \left]0,1\right[} \left( \prod_{j=1}^{i-1}C\bigl(\Gamma_j^{(1)},\gamma_u\bigr) \right)\times \partial_C C\bigl(\Gamma_i^{(1)},\gamma_u\bigr) \times \left(\prod_{j=i+1}^{r}C\bigl(\Gamma_j^{(1)},\gamma_u\bigr)\right)\right),\end{equation*} where the dimension of $C(\Gamma_j^{(1)},\gamma_u)$ is even. When rewriting such an integral as an integral over $\left[0,1\right] \times \Delta^{r}$ of the two-form $(\gamma \circ p_i)^{\ast}(d\eta_{\finsetb,\finseta_i})$ and one-forms $ (\gamma \circ p_j)^{\ast}(\eta_{\finsetb,\finseta_j})$, one must take into account the fact that
the two-form will be integrated along the product by $\left[0,1\right]$ of the interval parametrized by $t_i$. This gives rise to the factor $(-1)^{i-1}$.

Recall $\partial \simplexr=\sum_{i=0}^r (-1)^{i+1} F_i(\simplexr)$,
with $F_0(\simplexr)=\{(0,t_2,\dots,t_r) \in \Delta^r \}$, $F_r(\simplexr)=\{(t_1,t_2,\dots,t_{r-1},1) \in \simplexr \}$, and
\begin{equation*}F_i(\simplexr)=\bigl\{(t_1,\dots,t_i,t_i,t_{i+1},\dots ,t_{r-1}) \in \simplexr \bigr\}\end{equation*}
for $i \in \underline{r-1}$.
Observe that the faces $F_0$ and $F_r$ do not contribute to $\delol{(\gamma_u)}{(\eta\wedge \eta)_{\finsetb,\finseta}}$. Indeed for $F_0$, the directions of the edges of $\Gamma_1$
are in the image of $\CV(\gamma_0(0),\Gamma_1)$, which is $(2\cardlef{E(\Gamma_1)} -1)$-dimensional.
The contribution of the faces $F_i$ yields
\begin{multline*}\bigdelol{(\gamma_u)}{(\eta\wedge \eta)_{\finsetb,\finseta}}\\
=\sum_{r=1}^{\infty} \sum_{(\finseta_1, \dots,\finseta_{r-1}) \in P_{r-1}(\finseta)} \frac{\prod_{i=1}^{r-1}\cardlef{\finseta_i}!}{\cardlef{\finseta}!}\int_{\left[0,1\right] \times \Delta^{r-1}}\sum_{i=1}^{r-1}\tilde{\delta}_d(i,\finseta_1, \dots,\finseta_{r-1}),\end{multline*}
with \begin{equation*}\tilde{\delta}_d(i,\finseta_1, \dots,\finseta_{r-1})= (-1)^{i-1}\bigwedge_{j=1}^{r-1} (\gamma \circ p_j)^{\ast}(\eta_{\finsetb,\finseta_j})\left(\frac{(\gamma \circ p_i)^{\ast}\bigl((\eta\wedge\eta)_{\finsetb,\finseta_i}\bigr)}{(\gamma \circ p_i)^{\ast}(\eta_{\finsetb,\finseta_i})}\right).\end{equation*}
\eop

Note that
$(\eta\wedge\eta)_{\finsetb,\finseta_i}=\sum_{(\Gamma,\Gamma^{\prime}) \in \Davis^{e,2}_{n_i,\finseta_i}(\sqcup_{\eltb \in \finsetb}\RR_{\eltb})} \coefgambet_{\Gamma \sqcup \Gamma^{\prime}}\eta_{\Gamma} \wedge\eta_{\Gamma^{\prime}} \left[\Gamma\right]\left[\Gamma^{\prime}\right]$ is valued in the space of primitive elements of $\Aavis_{n_i}(\sqcup_{\eltb \in \finsetb}\RR_{\eltb})$. (Use $\eta_{\Gamma} \wedge\eta_{\Gamma^{\prime}} \left[\Gamma\right]\left[\Gamma^{\prime}\right] +\eta_{\Gamma^{\prime}} \wedge\eta_{\Gamma} \left[\Gamma^{\prime}\right]\left[\Gamma\right]=\eta_{\Gamma} \wedge\eta_{\Gamma^{\prime}}\left(\left[\Gamma\right]\left[\Gamma^{\prime}\right]-\left[\Gamma^{\prime}\right]\left[\Gamma\right]\right)$.)
 
 \begin{lemma} \label{lemthmconnecgenzerobis}
The form $d\eta_{\finsetb,\finseta}+(\eta\wedge\eta)_{\finsetb,\finseta}$ on $\left[0,1\right] \times \cinjuptd{\finsetb}{\CC}$ vanishes identically for any subset $\finseta$ of $\underline{3N}$ with cardinality $3n$.
\end{lemma}
\bp We proceed by induction on $n$ as in the proof of Corollary~\ref{coretaBflat},
using Lemma~\ref{lemcurvetaab} and Lemma~\ref{lemthmconnecgenzero}, which guarantees that $\lol{.}{\eta_{\finsetb,\finseta}}$ vanishes along homotopically trivial loops in $\{u\} \times \cinjuptd{\finsetb}{\CC}$ and along boundaries of squares $\left[w,w+ \varepsilon\right] \times p_{CS} \circ\gamma\left(\left[t,t+\varepsilon\right]\right)$. (This flatness condition can also be proved directly as an exercise, as in Remark~\ref{rkflatdirect}.)
\eop

Lemma~\ref{lemthmconnecgen} is proved.
\eopwobp

We now generalize Lemma~\ref{lemthmconnecgen} as follows.

\begin{lemma}
 \label{lemthmconnecgenbis}
Let $\gamma$ and $\delta$ be two semi-algebraic paths of $\left[0,1\right]\times \ccompuptd{\finsetb}{\CC}$ homotopic relatively to $\{0,1\}$. Assume $\gamma\left(\left]0,1\right[\right) \subset \cinjuptd{\finsetb}{\CC}$ and $\delta\left(\left]0,1\right[\right) \subset \cinjuptd{\finsetb}{\CC}$. Then we have
\begin{equation*}\lol{.}{\eta_{\finsetb,\finseta}}(\gamma)=\lol{.}{\eta_{\finsetb,\finseta}}(\delta).\end{equation*}
\end{lemma}
\bp According to Lemma~\ref{lemthmconnecgenzero}, it suffices to take care of homotopies near the endpoints. Thanks to Lemma~\ref{lemmultlol}, it suffices to prove that there exist $t, t^{\prime}>0$, and a path $\varepsilon$ from $\gamma(t)$ to 
$\delta(t^{\prime})$ in $\left[0,1\right]\times \cinjuptd{\finsetb}{\CC}$ such that $\gamma\vert_{\left[0,t\right]}\varepsilon \overline{\delta\vert_{\left[0,t^{\prime}\right]}}$ is null-homotopic and the holonomy along
$\gamma\vert_{\left[0,t\right]}\varepsilon \overline{\delta\vert_{\left[0,t^{\prime}\right]}}$ is one. When $t$ and $t^{\prime}$ are small enough, the images of $\gamma\vert_{\left[0,t\right]}$ and $\delta\vert_{\left[0,t^{\prime}\right]}$ lie in a subset equipped with a local semi-algebraic chart as in Theorem~\ref{thmcompuptd}, from which it is easy to construct semi-algebraic interpolations in products of sphere pieces and intervals.
Furthermore, there is no loss of generality in assuming that $\gamma\vert_{\left[0,t\right]}\varepsilon$ and $\delta\vert_{\left[0,t^{\prime}\right]}$ meet only
at $\gamma(0)$ and $\delta(t^{\prime})$, and that straight interpolation provides a boundary-fixing semi-algebraic homotopy from $\gamma\vert_{\left[0,t\right]}\varepsilon$ to $\delta\vert_{\left[0,t^{\prime}\right]}$, which satisfies the injectivity hypotheses of Lemma~\ref{lemthmconnecgen}. (Otherwise, we could use an intermediate $\gamma^{\prime}\vert_{\left[0,t\right]}.$) Thus, Lemma~\ref{lemthmconnecgen} allows us to prove 
$\lol{\gamma\vert_{\left[0,t\right]}\varepsilon}{\eta_{\finsetb,\finseta}}=\lol{\delta\vert_{\left[0,t^{\prime}\right]}}{\eta_{\finsetb,\finseta}}$.
\eop

\bpo{Proof of Theorem~\ref{thmconnecgen}} 
Lemma~\ref{lemthmconnecgenbis} allows us to define a map $\tilde{h}_{\finseta}$ induced by $\lol{.}{\eta_{\finsetb,\finseta}}$ from homotopy classes of paths with fixed boundaries of $\left[0,1\right] \times \ccompuptd{\finsetb}{\CC}$ to $\Aavis_{ {\cardlef{\finseta}}/{3}}(\sqcup_{\eltb \in \finsetb}\RR_{\eltb})$ as follows.
For a path $\gamma \colon \left[0,1\right] \to \ccompuptd{\finsetb}{\CC}$, set
$\tilde{h}_{\finseta}(\gamma)=\lol{.}{\eta_{\finsetb,\finseta}}(\delta)$ for any semi-algebraic path $\delta$ of $\left[0,1\right]\times \ccompuptd{\finsetb}{\CC}$ that is homotopic to $\gamma$ relatively to $\{0,1\}$ and such that $\delta\left(\left]0,1\right[\right) \subset \cinjuptd{\finsetb}{\CC}$.

Now, it suffices to prove that $\lol{\gamma}{\eta_{\finsetb,\finseta}}$ coincides with $\tilde{h}_{\finseta}(\gamma)$ for any semi-algebraic path $\gamma \colon \left[0,1\right] \to \left[0,1\right] \times \ccompuptd{\finsetb}{\CC}$.
Recall that a stratum of $\left[0,1\right] \times \ccompuptd{\finsetb}{\CC}$ is the product by $\left[0,1\right]$ of a stratum of $\ccompuptd{\finsetb}{\CC}$ associated to a $\Delta$-parenthesization.
The preimage under $\gamma$ of such a stratum of $\left[0,1\right] \times \ccompuptd{\finsetb}{\CC}$ is semi-algebraic. So $\gamma$ is a finite composition of paths whose interiors lie in a fixed stratum of $\ccompuptd{\finsetb}{\CC}$ (according to the Lojasiewicz theorem \ref{thmloja}). Thus, thanks to Lemma~\ref{lemmultlol}, it suffices to prove that $\lol{\gamma}{\eta_{\finsetb,\finseta}}$ coincides with $\tilde{h}_{\finseta}(\gamma)$
for any injective semi-algebraic path $\gamma$ whose interior lies in a fixed stratum of $\ccompuptd{\finsetb}{\CC}$ and in a subset equipped with a local semi-algebraic chart as in Theorem~\ref{thmcompuptd}.
Such a path can be deformed by sending the vanishing coordinates of $\gamma(t)$ in the $\left[0,\varepsilon\right[$ factors in such a chart to $\varepsilon (\frac12 - |\frac12 -t|)u$ for $u \in\left[0,1\right]$ giving rise to a semi-algebraic homotopy $(\gamma_u(t))_{u \in \left[0,1\right]}$ such that $\tilde{h}_{\finseta}(\gamma)=\lol{\gamma_1}{\eta_{\finsetb,\finseta}}$. So Lemma~\ref{lemthmconnecgen} implies $\lol{\gamma_0}{\eta_{\finsetb,\finseta}}=\lol{\gamma_1}{\eta_{\finsetb,\finseta}}$. \eop

Theorem~\ref{thmconnecgen} allows us to set the following definition. 

\begin{definition} \label{defgenlol}
 For any continuous path $\gamma \colon \left[0,1\right] \to \left[0,1\right] \times \ccompuptd{\finsetb}{\CC}$, define $\lol{.}{\eta_{\finsetb,\finseta}}(\gamma)$ to be \begin{equation*}\lol{.}{\eta_{\finsetb,\finseta}}(\gamma)=\lol{.}{\eta_{\finsetb,\finseta}}(\delta)\end{equation*} for any semi-algebraic path $\delta$ of $\left[0,1\right]\times \ccompuptd{\finsetb}{\CC}$ homotopic to $\gamma$ relatively to $\{0,1\}$.
\end{definition}

We can now generalize \cite[Proposition 1.18]{poirierv2} for braids.
\begin{proposition}
\label{propdubraid}
Let $\finsetb$ and $\finsetc$ be two finite sets. Let $\eltb_0 \in \finsetb$.

 Let $\gamma_{\finsetb}$ be a path of $\left[0,1\right] \times \ccompuptd{\finsetb}{\CC}$
and let $\gamma_{\finsetc}$ be a path of $\left[0,1\right] \times \ccompuptd{\finsetc}{\CC}$.  
Let $\gamma_{\finsetb}(\gamma_{\finsetc}/b_0)=\gamma_{\finsetb}(\gamma_{\finsetc}/K_{b_0})$ be the $q$-braid obtained by cabling the strand $K_{b_0}$ of $b_0$ in $T(\gamma_{\finsetb})$ by $T(\gamma_{\finsetc})$, as in Notation~\ref{notcabling}.
Then we have \begin{equation*}\biglol{\gamma_{\finsetb}(\gamma_{\finsetc}/b_0)}{\eta_{\finsetb\left(\frac{\finsetc}{\eltb_0}\right),.}}=\left(\pi(\finsetc \times \eltb_0)^{\ast}\bigl(\lol{\gamma_{\finsetb}}{\eta_{\finsetb,.}}\bigr)\lol{\gamma_{\finsetc}}{\eta_{\finsetc,.}}\right)_{\sqcup},\end{equation*}
where $\pi(\finsetc \times \eltb_0)^{\ast}$ denotes the duplication of the strand $\RR_{b_0}$ for diagrams as in Notation~\ref{notationduplication} and we use the product of Definition~\ref{defprodbiz}.
\end{proposition}
\bp 
Thanks to Definition~\ref{defgenlol}, we assume that $\gamma_{\finsetc}$ and $\gamma_{\finsetb}$ are semi-algebraic, without loss of generality. 
Since a semi-algebraic path is a path composition of finitely many semi-algebraic paths whose interiors lie in a fixed stratum of $\left[0,1\right] \times \ccompuptd{.}{\CC}$ (according to the Lojasiewicz theorem \ref{thmloja}), 
the functoriality of Lemma~\ref{lemmultlol} allows us to assume furthermore that
the image of the interior of $\gamma_{\finsetc}$ lies in a fixed stratum of $\left[0,1\right] \times \ccompuptd{\finsetc}{\CC}$
and that the image of the interior of $\gamma_{\finsetb}$ lies in a fixed stratum of $\left[0,1\right] \times \ccompuptd{\finsetb}{\CC}$, for the proof. (Recall the commutation lemma~\ref{lemdupcom}.)

Let $\finsetb({\finsetc}/{\eltb_0})$ be the set obtained from $\finsetb$ by replacing $\eltb_0$ by $\finsetc$. 
We refer to Lemma~\ref{lemdimstratCV} for the description of the stratification of $\CV(\Gamma)$ and the dimensions of the fibers for connected diagrams on $\finsetb({\finsetc}/{\eltb_0}) \times \RR$. When computing the \say{holonomy} of $\eta_{\finsetb({\finsetc}/{\eltb_0}),.}$ we integrate over products of one-parameter families
$\CV(\confy,\Gamma)$ with $\confy \in \gamma_{\finsetb}(\gamma_{\finsetc}/b_0)\left(\left]t_i-\varepsilon,t_i+\varepsilon\right[\right)$, locally.
We may restrict to the strata of $\CV(\confy,\Gamma)$ of dimension $(2\cardlef{E(\Gamma)}-1)$, which are described in Lemma~\ref{lemdimstratCV}.

Consider a connected diagram $\Gamma$ on $\finsetb({\finsetc}/{\eltb_0}) \times \RR$ together with a $\Delta$-parenthesization $\parentp$ of its vertices corresponding to such a stratum of configurations. 
Since $\parentp$ is equal to $\parntppx$, all elements of $\parentp$ are univalent. Assume $\parentp \neq \{ V(\Gamma) \}$.
Let $\Gamma_{\parentp}$ be obtained from $\Gamma$ by identifying all the vertices in a daughter $\finseta$ of $V(\Gamma)$ to a single vertex $\eltv_{\finseta}$ and by erasing the edges between two elements in $\finseta$, for each $\finseta$ of $D(V(\Gamma))$. Then $\Gamma_{\parentp}$ is connected and its vertices $\eltv_{\finseta}$ move along vertical lines. Let $U(\Gamma_{\parentp})$ and  $T(\Gamma_{\parentp})$ respectively denote the set of univalent vertices of $\Gamma_{\parentp}$ distinct from the $\eltv_{\finseta}$ and the set of trivalent vertices of $\Gamma_{\parentp}$ distinct from the $\eltv_{\finseta}$. The dimension of the one-parameter family of configurations of the vertices of $\Gamma_{\parentp}$ up to vertical translation is $\cardlef{D(V(\Gamma))} +\cardlef{U(\Gamma_{\parentp})} + 3 \cardlef{T(\Gamma_{\parentp})}$. The form $\bigwedge_{e \in E(\Gamma_{\parentp})}p_{e,S^2}^{\ast}(\tilde{\omega}(j_E(e),S^2))$ factors through this one-parameter family. Thus, a count of half-edges shows that the stratum cannot contribute unless the vertices $\eltv_{\finseta}$ are univalent in $\Gamma_{\parentp}$. 

Now assume that all the vertices $\eltv_{\finseta}$ are univalent in $\Gamma_{\parentp}$. Let $e_{\finseta}$ be the edge of $\Gamma_{\parentp}$ with maximal label such that $e_{\finseta}$ is adjacent to a vertex $\eltv_{\finseta}$ for some $\finseta \in D(V(\Gamma))$. The subgraph $\Gamma_{\finseta}$ of $\Gamma$ consisting of the vertices of $\finseta$ and the edges of $\Gamma$ between two such vertices is connected. It has one bivalent vertex $b$ in $\Gamma_{\finseta}$ (which is the end of $e_{\finseta}$ in $\Gamma \cap \finseta$). Its configurations are considered up to vertical translations. Their contribution is
opposite to that of the configurations of the graph $\Gamma^{\prime}_{\finseta}$ obtained from $\Gamma_{\finseta}$ by exchanging the labels and possibly the orientations of the two edges of $\Gamma_{\finseta}$ that contain $b$, as in Lemma~\ref{lemsym}.

Thus, the $\Delta$-parenthesization of $\Gamma$ is $\{V(\Gamma)\}$ in the strata that may contribute. If $p_{\finsetb({\finsetc}/{\eltb_0})}(U(\Gamma)) \subset \finsetc$, then $\Gamma$ is a diagram on $\finsetc \times \RR$, which contributes as in $\lol{\gamma_{\finsetc}}{\eta_{\finsetc,.}}$. Otherwise, the projection to the horizontal plane of the vertices of $U(\Gamma) \cap p_{\finsetb({\finsetc}/{\eltb_0})}^{-1}(\finsetc)$ is reduced to a point. So all diagrams obtained from these diagrams $\Gamma$ by changing the map from $U(\Gamma) \cap p_{\finsetb({\finsetc}/{\eltb_0})}^{-1}(\finsetc)$ to $\finsetc$ arbitrarily contribute together to $\pi(\finsetc \times \eltb_0)^{\ast}(\lol{\gamma_{\finsetb}}{\eta_{\finsetb,.}})$ as desired, locally. 
We get the proposition since the two kinds of diagrams commute thanks to Lemma~\ref{lemdupcom}.
\eop

\begin{definition}
\label{defvarzfqtang} Recall Proposition~\ref{propfuncbraidtwo} and Definition~\ref{defgenlol}.
 For a $q$-braid (representative) $\gamma \colon \left[0,1\right] \to \ccompuptd{\finsetb}{\CC}$, set 
\begin{equation*}\Zinvufrfneg\Bigl(\gamma,.,\bigl(\tilde{\omega}(i,1,S^2)\bigr)_{i\in \underline{3\largen}}\Bigr)=\lol{\{1\} \times \gamma }{\eta_{\finsetb,.}}.\end{equation*}
For a $q$-tangle \begin{equation*}T=T(\gamma^-)(\hcylc,\tanghcyll,\tanghcyll_{\parallel})T(\gamma^+)\end{equation*} 
such that $\gamma^-$ and $\gamma^+$ are $q$-braids, and $(\hcylc,\tanghcyll)$ is a $J_{bb,tt}$-oriented framed tangle whose bottom configuration is $\gamma^-(1)$ and whose top configuration is $\gamma^+(0)$, as in Definition~\ref{deftanglelong},
for $\largen \in \NN$, and for a family $\left({\omega}(i,S^2)\right)_{i\in \underline{3\largen}}$ of volume-one forms of $S^2$, set \begin{multline*}\Zinvufrfneg\Bigl(T,.,\bigl({\omega}(i,S^2)\bigr)_{i \in \underline{3\largen}}\Bigr)\\                                                                                           =\biggl(\Zinvufrfneg\Bigl(\gamma^-,.,\bigl({\omega}(i,S^2)\bigr)\Bigr)\Zinvufrfneg\Bigl(\hcylc,\tanghcyll,\tanghcyll_{\parallel},.,\bigl({\omega}(i,S^2)\bigr)\Bigr)\Zinvufrfneg\Bigl(\gamma^+,.,\bigl({\omega}(i,S^2)\bigr)\Bigr)\biggr)_{\sqcup}                                                                                                        \end{multline*}
with the notation of Definition~\ref{defprodbiz} and Definition~\ref{defnotzfvariant}.
\end{definition}
Lemma~\ref{lemmultlol} and the isotopy invariance of Theorems~\ref{thmdefsanstauvarzinvf} and \ref{thmconnecgen} ensure that the definition above is consistent.

Theorem~\ref{thmtangconstcompar} 
allows us to express the variation of \begin{equation*}\Zinvufrfneg\Bigl(T, .,\bigl(\omega(i,S^2)\bigr)_{i \in \underline{3\largen}}\Bigr)\end{equation*} when $(\omega(i,S^2))_{i \in \underline{3\largen}}$ varies for framed straight tangles with injective top and bottom configurations.\footnote{See also the proof of Theorem~\ref{thmdefsanstauvarzinvf}.}
As a corollary of Theorem~\ref{thmconnecgen}, this expression generalizes to $q$-tangles. We get the following theorem.

\begin{theorem}
\label{thmtangconstcomparbis} Let $\largen \in \NN$.
For $i\in \underline{3\largen}$, let $\tilde{\omega}(i,S^2)=(\tilde{\omega}(i,t,S^2))_{t\in \left[0,1\right]}$ be a closed $2$-form on $\left[0,1\right] \times S^2$ with
$\int_{\{0\} \times S^2}\tilde{\omega}(i,0,S^2)=1$.
Let $T$ denote a $J_{bb,tt}$-oriented $q$-tangle.
Assume that the bottom and top configurations of $T$ are elements $\confy^-$ of $\ccompuptd{\finsetb^-}{\CC}$ and $\confy^+$ of $\ccompuptd{\finsetb^+}{\CC}$, respectively.
Then we have
\begin{multline*}\Zinvufrfneg\Bigl(T,.,\bigl(\tilde{\omega}(i,1,S^2)\bigr)_{i \in \underline{3\largen}}\Bigr)\\
 =\left(\begin{array}{r} \left(\prod_{K_j \in J_{bb}}\biglol{\left[0,1\right]}{-\eta(.,S_{WE})} \#_j\right) \left(\prod_{K_j \in J_{tt}}\biglol{\left[0,1\right]}{\eta(.,S_{WE})} \#_j\right)\\
\lol{\left[1,0\right] \times \confy^- }{\eta_{\finsetb^-,.}}\Zinvufrfneg\Bigl(T,.,\bigl(\tilde{\omega}(i,0,S^2)\bigr)\Bigr)\lol{\left[0,1\right]\times \confy^+ }{\eta_{\finsetb^+,.}}\end{array}\right)_{\sqcup}.
\end{multline*}
\end{theorem}
\eopwobp

Recall that the forms $\eta$ and their holonomies introduced in Definitions~\ref{defholtilde} and \ref{defAholonom} depend on 
the forms $\tilde{\omega}(i,S^2)$. In the end of this chapter, we will apply Theorem~\ref{thmtangconstcomparbis}, only
when $\tilde{\omega}(i,0,S^2)$ is the standard form $\omega_{S^2}$. So $\Zinvufrfneg(T,.,(\tilde{\omega}(i,0,S^2))$ is simply $\Zinvufrfneg(T,.)$.

The following lemma is also easy to prove.
\begin{lemma}
\label{lemrevcompgen}
The behavior of $\Zinvufrf$ under reversing the orientation of a closed component is the same as that described in Proposition~\ref{proporcomp}.
Its behavior under reversing the orientation of component going from bottom to top or from top to bottom in a $q$-tangle $T$ is similar:
Let $K$ be an oriented component going from bottom to top or from top to bottom in a $q$-tangle $T$.
For a Jacobi diagram $\Gamma$ on the domain $\sourcetl$ of $T$, let $U_K(\Gamma)$ denote the set of univalent vertices of $\Gamma$ mapped to the domain $\RR_K$ of $K$. This set is ordered by the orientation of $\RR_K$.
When the orientation of $K$ is changed, $\Zinvufrfneg(T)$ is modified by reversing the orientation of $\RR_K$ (that is reversing the order of $U_K(\Gamma)$) in classes $\left[\Gamma\right]$ of diagrams $\Gamma$ on $\sourcetl$, and multiplying these classes by $(-1)^{\cardlef{U_K(\Gamma)}}$ in $\Aavis(\sourcetl)$, simultaneously.
\end{lemma}
\eopwobp

In other words, we can forget the orientation of closed components and components going from bottom to top or from top to bottom, and regard $\Zinvufrfneg(\Link)$ as valued in spaces of diagrams where the domains of these components are not oriented, as in Definitions~\ref{defdia} and \ref{defrkoruniv}. However, we a priori need a $J_{bb,tt}$-orientation as Lemma~\ref{lemetaWEEW} below indicates.

\begin{lemma} \label{lemetaWEEW}
Let $\finseta$ be a set of cardinality $3k$. Recall Definition~\ref{defAholonom}.
The form $\eta(\finseta,S_{WE} \cup (- S_{EW}))$ is zero when $k$ is odd.
If $k$ is even, then
we have \begin{equation*}\eta\bigl(\finseta,S_{WE} \cup (- S_{EW})\bigr) = -\revmaps_{\ast}\Bigl(\eta\bigl(\finseta,S_{WE} \cup (- S_{EW})\bigr)\Bigr).\footnote{In particular, $\eta(\finseta,S_{WE} \cup (- S_{EW}))$ is also zero when $k$ is even, if $\revmaps_{\ast}$ is the identity map. But this is unknown to me.}\end{equation*}
If $\eta(\finseta,S_{WE} \cup (- S_{EW}))$ is zero, then $\Zinvufrfneg(.,.,(\tilde{\omega}(i,1,S^2))_{i \in \underline{3\largen}})$ is independent of the
$J_{bb,tt}$-orientation of the tangles (as in Proposition~\ref{proporcomp}).
\end{lemma}
\bp The circle $S_{WE} \cup (- S_{EW})$ is the great circle $\partial D$ of $S^2$ that is the boundary of the hemisphere $D$ of $S^2$ centered at $(-i)$. By Definition~\ref{defAholonom}, we have
$\eta(\finseta,\partial D)(t)=\int_{\{t\} \times \partial D}\omega(\finseta)$
with
\begin{equation*}\int_{\left[0,t\right] \times \partial D}\omega(\finseta) =\int_{\left(\partial \left[0,t\right]\right) \times D} \omega(\finseta) =\int_{\{t\} \times D} \omega(\finseta) -\int_{\{0\} \times D} \omega(\finseta)\end{equation*}
and
\begin{equation*}\int_{\{t\} \times S^2} \omega(\finseta) =  \int_{\{t\} \times D}\Bigl( \omega(\finseta) - \iota^{\ast}\bigl(\omega(\finseta)\bigr)\Bigr)=\bigl(1+(-1)^{k+1}\bigr)\int_{\{t\} \times D}\omega(\finseta)\end{equation*}
according to 
Lemma~\ref{lemrevomega}. In particular, when $k$ is odd, the integral \begin{equation*}\int_{\{t\} \times D}\omega(\finseta)=\frac12 \int_{\{t\} \times S^2}\omega(\finseta)\end{equation*} is independent of $t$ and $\eta(\finseta,\partial D)$ is zero. When $k$ is even, Lemma~\ref{lemrevomega} implies $\omega(\finseta)=-\revmaps_{\ast}(\omega(\finseta))$.

Changing the orientation of a component $K$ of $J_{bb}$ amounts to construct the invariant by imposing the condition $p_{\tau}(\ST^+K) \subseteq S_{EW}$ rather than $p_{\tau}(\ST^+K) \subseteq S_{WE}$. So this replaces the factor $\lol{\left[0,1\right]}{-\eta(.,S_{WE})}$ associated to $K$ with $\lol{\left[0,1\right]}{-\eta(.,S_{EW})}$ in the formula of Theorem~\ref{thmtangconstcomparbis}. See Theorem~\ref{thmdefsanstauvarzinvf} and Definition~\ref{defnotzfvariant}. This amounts to multiply by $\lol{\left[0,1\right]}{\eta(.,S_{WE} \cup - S_{EW})}$ on the component of $K$.
\eop

\begin{remarks}
\label{rknotzfvarianttwo}
Similarly, if we had imposed that $p_{\tau}$ maps the unit tangent vectors to components of $J_{bb}\cup J_{tt}$ to the vertical half great circle $S(\theta)$ from $-\upvec$ to $\upvec$ that contains the complex direction $\exp(2i\pi \theta)$, for $\theta \in \left]0,1\right[$, in our definition of straight tangles in Section~\ref{secstraighttang}, then $S(\theta)$ would replace $S_{WE}$ in the formula of Theorem~\ref{thmtangconstcomparbis}, and 
$\Zinvufrfneg(T,.,(\tilde{\omega}(i,1,S^2))_{i \in \underline{3\largen}})$ would have been multiplied by \begin{multline*}\left(\prod_{K_j \in J_{bb}}\lolsp{\left[0,1\right]}{\eta\Bigl(.,S_{WE}\cup \bigl(- S(\theta)\bigr)\Bigr)} \#_j\right) \\ \left(\prod_{K_j \in J_{tt}}\lolsp{\left[0,1\right]}{\eta\bigl(.,S(\theta) \cup (-S_{WE})\bigr)} \#_j\right).\end{multline*} 

In particular, when $T$ is a tangle with only one component going from bottom to bottom, $\Zinvufrfneg(T,.,(\tilde{\omega}(i,1,S^2))_{i \in \underline{3\largen}})$ would have been multiplied by $\lolsp{\left[0,1\right]}{\eta(.,S_{WE} \cup (- S(\theta)))}$. 

With the notation of Definition~\ref{defAholonom} and Lemma~\ref{lemrevomega}, \begin{multline*}\lolsp{\left[0,1\right]}{\eta\Bigl(\underline{3},S_{WE} \cup \bigl(- S(\theta)\bigr)\Bigr)}=\int_{\left[0,1\right]\times (S_{WE} \cup (- S(\theta)))} \omega(\underline{3})\\
=\frac16\sum_{i=1}^3\int_{\left[0,1\right]\times (S_{WE} \cup (- S(\theta)))}\bigl(\tilde{\omega}(i,S^2)-\iota^{\ast}\tilde{\omega}(i,S^2)\bigr)\left[ \onechordsmalltseul \right]
\end{multline*}
would have been added to $\Zinvufrfneg(T,\underline{3},(\tilde{\omega}(i,1,S^2))_{i \in \underline{3}})$.
For a chain $D$ of $S^2$ bounded by $(S_{WE} \cup (- S(\theta)))$, we have
\begin{equation*}\int_{\left[0,1\right]\times (S_{WE} \cup (- S(\theta)))}\tilde{\omega}(i,S^2)= \int_{(\partial \left[0,1\right] )\times D}\tilde{\omega}(i,S^2).\end{equation*}
In particular, if $\theta$ is in $\left]1/2,1\right[$, then we can choose $D$ and $\tilde{\omega}(i,1,S^2)$ so that $\left(\tilde{\omega}(i,1,S^2)-\iota^{\ast}\tilde{\omega}(i,1,S^2)\right)$ is supported outside $D \cup \iota_{S^2}(D)$ for any $i$. In this case, we get \begin{equation*}\lolsp{\left[0,1\right]}{\eta\Bigl(\underline{3},S_{WE} \cup \bigl(- S(\theta)\bigr)\Bigr)}=-\int_{\{0\}\times D} \omega_{S^2}=\theta-1.\end{equation*}
Thus, as claimed in Remark~\ref{rknotzfvariant}, Definition~\ref{defnotzfvariant} is not canonical. 

The above calculation does not rule out the alternative choice of the vertical half great circle $S(\frac12)=S_{EW}$ for our definition. This choice would multiply $\Zinvufrfneg(T,.,(\tilde{\omega}(i,1,S^2))_{i \in \underline{3\largen}})$ by $\lol{\left[0,1\right]}{\eta(.,S_{WE} \cup (- S_{EW}))}$, which is zero in degree $1$, according to Lemma~\ref{lemetaWEEW}.

It might be tempting to modify the definition of $\Zinvufrfneg(T,.,(\tilde{\omega}(i,1,S^2))_{i \in \underline{3\largen}})$, by multiplying it by \begin{equation*}\left(\prod_{K_j \in J_{bb}}\biglol{\left[0,1\right]}{\eta(.,S_{WE})} \#_j\right) \left(\prod_{K_j \in J_{tt}}\biglol{\left[0,1\right]}{-\eta(.,S_{WE})} \#_j\right).\end{equation*} Unfortunately, $\lol{\left[0,1\right]}{\eta(.,S_{WE})}$ depends on the closed $2$-forms $\tilde{\omega}(i,S^2)$ of $\left[0,1\right] \times S^2$, and not only on the $\tilde{\omega}(i,1,S^2)$. Indeed, assume that all $\tilde{\omega}(i,S^2)$ coincide with each other and change all of them by adding $d \eta_S$, for a one-form $\eta_S$ of $\left[0,1\right] \times S^2$ supported on the product of $\left[1/4,3/4\right]$ by a small neighborhood of $\upvec$.
Then the variation of the degree one part of $2\lol{\left[0,1\right]}{\eta(.,S_{WE})}$ maps $\underline{3}$ to $\int_{\left[0,1\right]\times S_{WE}}\left(d \eta_S-\iota^{\ast}(d \eta_S)\right)$,
with 
\begin{equation*}\int_{\left[0,1\right]\times S_{WE}}d \eta_S = \int_{\partial (\left[0,1\right]\times S_{WE})} \eta_S=-\int_{\left[0,1\right]\times \{\upvec\}}\eta_S\end{equation*} and 
\begin{equation*}\begin{array}{ll}\int_{\left[0,1\right]\times S_{WE}}-\iota^{\ast}(d \eta_S)= -\int_{\left[0,1\right]\times \iota_{S^2}(S_{WE})}d \eta_S&
=\int_{\partial (\left[0,1\right]\times -\iota_{S^2}(S_{WE}))} \eta_S\\
  &=  -\int_{\left[0,1\right]\times \{\upvec\}}\eta_S      . \end{array}                                                                                                                                      
\end{equation*} (In Theorem~\ref{thmtangconstcompar}, the factors $\lol{\left[t,0\right] \times \confy^{-} }{\eta_{\finsetb^-,.}}$ and $\lol{\left[0,t\right] \times \confy^{+} }{\eta_{\finsetb^+,.}}$ also depend on $\tilde{\omega}(i,S^2)$. However, both types of dependences cancel each other.) 
\end{remarks}

\chapter{Justifying the properties of \texorpdfstring{$\Zinvufrf$}{Zf}}
\label{chappropzinvffunc}

Recall Definition~\ref{defzqtang} of the invariant $\Zinvufrf$.
So far we have succeeded in constructing this 
 invariant $\Zinvufrf$ of $q$-tangles, invariant under boundary-fixing diffeomorphisms, generalizing both 
the invariant $\Zinvufrf$ for framed links in $\QQ$-spheres
and the Poirier functor $Z^l$ for $q$-tangles in $\RR^3$. The framing dependence of Theorem~\ref{thmmainfunc} comes from Definition~\ref{deffstconsttangframed}.

The behavior of $\Zinvuf$ and $\Zinvufrf$ under orientation changes of the components described in the statement of Theorem~\ref{thmmainfunc} can be justified as in the case of links in rational homology spheres treated in Section~\ref{secfstpropzinv}.

\begin{lemma}
 The invariant $\Zinvufrf$ behaves as prescribed by Theorem~\ref{thmmainfunc} under the diffeomorphisms $s_{\frac12}$ and $\rho$.
\end{lemma}
\bp 
Let $\tanghcyll$ be a tangle representative as in Theorem~\ref{thmfstconsttang}.
We proceed as in the proof of Proposition~\ref{propinvdiffnat}, except that we need to take care of the facts that, for $\psi =s_{\frac12}$ or $\psi =\rho$, $\psi_{\ast}(\partau)$ is not asymptotically standard, and that $s_{\frac12}$ reverses the orientation.
Therefore, we use $\partau^{\prime}=\psi_{\ast}(\partau) \circ (\id_{\crats(\hcylc)} \times \psi^{-1}_{\RR^3})=T\psi \circ \partau \circ (\psi^{-1} \times \psi^{-1}_{\RR^3})$, as an asymptotically standard
parallelization of $\rats=\rats(\hcylc)$. So we have $(\rho_{\ast}^{-1})^{\ast}\projp_{\partau}^{\ast}(\omega_{S^2})=\projp_{\partau^{\prime}}^{\ast}(\omega_{S^2})$ and $(\sigma_{\frac12\ast}^{-1})^{\ast}\projp_{\partau}^{\ast}(\omega_{S^2})=-\projp_{\partau^{\prime}}^{\ast}(\omega_{S^2})$. If $\omega$ is a homogeneous propagating form of $\bigl(C_2(\rats),\partau\bigr)$, then $(\rho_{\ast}^{-1})^{\ast}(\omega)$ is a homogeneous propagating form of $(C_2(\rho(\rats)),\partau^{\prime})$ and $(-\sigma_{\frac12\ast}^{-1})^{\ast}(\omega)$ is a homogeneous propagating form of $(C_2(s_{\frac12}(\rats)),\partau^{\prime})$. Let us now focus on the case $\psi=s_{\frac12}$, 
since the case $\psi =\rho$ is similar, but simpler.
For any Jacobi diagram $\Gamma$ on the domain of $\tanghcyll$, equipped with an implicit orientation $o(\Gamma)$, compute $I=I\bigl(s_{\frac12}(\rats),s_{\frac12}(\tanghcyll),\Gamma,(-\sigma_{\frac12\ast}^{-1})^{\ast}(\omega)\bigr)$. We have
\begin{equation*}\begin{array}{ll}I&=\int_{\check{C}(s_{\frac12}(\rats),s_{\frac12}(\tanghcyll);\Gamma)} \bigwedge_{e \in E(\Gamma)}p_e^{\ast}\bigl(-(\sigma_{\frac12\ast}^{-1})^{\ast}(\omega)\bigr)\\
  &=(-1)^{\cardlef{E(\Gamma)}}\int_{\check{C}(s_{\frac12}(\rats),s_{\frac12}(\tanghcyll);\Gamma)} (\sigma_{\frac12\ast}^{-1})^{\ast}\bigl(\bigwedge_{e \in E(\Gamma)}p_e^{\ast}(\omega)\bigr)\\
&=(-1)^{\cardlef{E(\Gamma)} + \cardlef{T(\Gamma)}}I\left(\rats,\tanghcyll,\Gamma,\omega\right).
  \end{array}
 \end{equation*}
Thus, we get $\Zinv_n(s_{\frac12}(\hcylc),s_{\frac12}(\tanghcyll),
\partau^{\prime})=(-1)^n\Zinv_n(\hcylc,\tanghcyll,\partau)$ for all $n \in \NN$.
In particular, Corollary~\ref{corThetap} and Propositions~\ref{propThetaorrev} and \ref{propthetazoneone} imply $p_1(\partau^{\prime})=-p_1(\partau)$ and 
$I_{\theta}(K_j,\partau^{\prime})=-I_{\theta}(K_j,\partau)$ for any component $K_j$ of $\tanghcyll$.
So we have $\Zinv_n(s_{\frac12}(\hcylc),s_{\frac12}(\tanghcyll))=(-1)^n\Zinv_n(\hcylc,\tanghcyll)$ for all $n \in \NN$ since the anomalies $\alpha$ and $\beta$ vanish in even degrees, thanks to Propositions~\ref{propdeftwoanom} and \ref{propanom}.
Now, recall Proposition~\ref{proplktangles}, and note that if a component $K$ is straight with respect to $\partau$, then $s_{\frac12}(K)$ is straight with respect to $\partau^{\prime}$.
In particular, the condition 
\begin{equation*}lk_{s_{\frac12}(\hcylc)}\left(s_{\frac12}(K),s_{\frac12}(K)_{\parallel}\right)=-lk_{\hcylc}\left(K,K_{\parallel}\right)\end{equation*}
is realizable and natural, and
we get the desired equality 
for a framed tangle from an injective bottom configuration to an injective top configuration.
Thanks to Remark~\ref{deflimzqtang}, it is still true for a $q$-tangle.  
\eop

We are thus left with the proofs of the functoriality, the duplication properties, and the cabling property  
to finish the proof of Theorem~\ref{thmmainfunc}. These proofs will occupy four sections of this chapter, which will end with a section describing other properties of $\Zinvufrf$.
The corresponding properties of variants of $\Zinvuf$ and $\Zinvufrf$ involving nonhomogeneous propagating forms will be treated simultaneously since they are often easier to prove, and since we are going to use them to prove some of the results for homogeneous propagating forms.

\section{Transversality and rationality}
\label{sectransrat}

In this section, we generalize the rationality results of Chapter~\ref{chaprat} to the tangle case. The generalization will be useful in the proofs of the properties later.

Let $S^2_H$ \index[N]{SatwoH@$S^2_H$ subset of $S^2$} denote the subset of $S^2$ consisting of the vectors whose vertical coordinate is in $]-\frac12,\frac12[$.

\begin{proposition}
\label{proptranstang}
Let $(\hcylc,\tau)$ be a parallelized rational homology cylinder.
Let $\tanghcyll \colon \sourcetl \hookrightarrow \hcylc$ be a long tangle of $\crats(\hcylc)$.
Let $\largen$ be an integer greater than $1$.
Then there exist $(X_1,X_2,\dots,X_{3\largen}) \in (S^2_H)^{3\largen}$, $\largem \in \left[1,+\infty\right[$, and propagating chains $\propP(i)$ of $(C_2(\crats(\hcylc)),\tau)$ for $i \in \underline{3\largen}$ such that
\begin{itemize}
 \item $\propP(i)$ intersects the domain
$D(p_{\tau})$ of Notation~\ref{notratsrest} as $p_{\tau}^{-1}(X_i)$,
\item the $\propP(i) \cap {C}_2(\rats_{\largem,\left[-\largem,\largem\right]}(\hcylc))$ are in general $3\largen$-position with respect to $\tanghcyll$, with the natural generalization of the notion of Definition~\ref{defgenthreenpos} (where $\rats_{\largem,\left[-\largem,\largem\right]}(\hcylc)$ replaces $\rats$, with Notation~\ref{notratsrest}),
\item the intersections \begin{equation*}I_S\Bigl(\Gamma,\bigl(\propP(i)\bigr)_{i \in \underline{3\largen}}\Bigr)=\bigcap_{e \in E(\Gamma)}p_e^{-1}\Bigl(\propP\bigl(j_E(e)\bigr)\Bigr)\end{equation*} in 
$C(\rats(\hcylc),\tanghcyll;\Gamma)$ are transverse and located in 
${C}(\rats_{\largem,\left[-\largem,\largem\right]}(\hcylc),\tanghcyll;\Gamma)$
for any $\Gamma \in \Davis^e_{\underline{3\largen}}(\source)=\cup_{k \in \NN}\Davis^e_{k,\underline{3\largen}}(\source)=\cup_{k=0}^{3\largen}\Davis^e_{k,\underline{3\largen}}(\source)$,
\item for any $\alpha >0$, there exists $\beta >0$ such that
\begin{equation*}\bigcap_{e \in E(\Gamma)}p_e^{-1}\biggl(N_{\beta}\Bigl(\propP\bigl(j_E(e)\bigr)\Bigr)\biggr) \subset N_{\alpha}\biggl(I_S\Bigl(\Gamma,\bigl(\propP(i)\bigr)_{i \in \underline{3\largen}}\Bigr)\biggr)\end{equation*} for any $\Gamma \in \Davis^e_{\underline{3\largen}}(\source)$, with the notation of Definition~\ref{defformdual},
\item there exists an open ball $B_X$ around $(X_1,X_2,\dots,X_{3\largen})$ in $(S^2)^{\underline{3\largen}}$ such that $B_X \subset (S^2_H)^{\underline{3\largen}}$ and for any 
$(Y_1,Y_2,\dots,Y_{3\largen}) \in B_X$, there exist propagating chains $\propP(i)(Y_i)$ of $(C_2(\crats(\hcylc)),\tau)$ satisfying all the above conditions with respect to $Y_i$ with the same $\largem$.
\end{itemize}
The set of $(X_1,X_2,\dots,X_{3\largen}) \in (S^2_H)^{3\largen}$ such that there exist \begin{equation*}\largem=\largem(X_1,X_2,\dots,X_{3\largen}) \in \left[1,+\infty\right[\end{equation*} and propagating chains $\propP(i)$ of $(C_2(\crats(\hcylc)),\tau)$ satisfying the above conditions is dense in $(S^2_H)^{3\largen}$.
\end{proposition}

In order to prove the proposition, we begin by producing some \begin{equation*}(\vecw_1,\vecw_2,\dots,\vecw_{3\largen}) \in (S^2_H)^{3\largen}\end{equation*} (in a given neighborhood of some $(\vecw^0_1,\vecw^0_2,\dots,\vecw^0_{3\largen})$ in $(S^2_H)^{3\largen}$) and \begin{equation*}\largem(\vecw_1,\vecw_2,\dots,\vecw_{3\largen}) \in \left[1,+\infty\right[.\end{equation*}

For a one-manifold $\sourcetl$ and a finite set $\finseta$, the set of connected
$\finseta$-numbered degree $k$ Jacobi diagrams with support $\sourcetl$ without looped edges is denoted by $\CD_{k,\finseta}^c(\sourcetl)$.

\begin{lemma}
\label{lemexisbonform} Let $\largen$ be an integer greater than $1$. Let $\finsetb$ be a finite set.
Let $\confy \colon \finsetb \hookrightarrow \drad{1}$ be a planar configuration.
For a $\underline{3\largen}$-numbered Jacobi diagram $\Gamma$ on $\sqcup_{b \in \finsetb}\RR_b$, define the semi-algebraic map 
\begin{equation*}g(\Gamma) \colon {C}\bigl(S^3,\confy(\finsetb) \times \RR;\Gamma\bigr) \times (S^2)^{\underline{3\largen} \setminus j_E(E(\Gamma))} \to (S^2)^{\underline{3\largen}}\end{equation*}
to be the product $\bigl(\prod_{e \in E(\Gamma)} p_{e,S^2}\bigr) \times \id\left((S^2)^{\underline{3\largen} \setminus j_E(E(\Gamma))}\right)$.

The subset $\CO(\largen,\confy)$ of $(S^2_H)^{3\largen}$ of points that are in the complement of the images of the maps $g(\Gamma)$ for all $\Gamma \in \cup_{k=1}^{3\largen}\CD_{k,\underline{3\largen}}^c\left(\sqcup_{b \in \finsetb}\RR_b\right)$ is dense and open.\footnote{ Though we are going to study $\Zinvufrf$ up to degree $\largen$, higher degree diagrams will occur in our proofs. See the proof of Lemma~\ref{lemconftangboun}.}
\end{lemma}
\bp
It suffices to prove that the complement of the image of the map $g(\Gamma)$ is open and dense for any of the finitely many graphs $\Gamma \in \cup_{k=1}^{3\largen}\CD_{k,\underline{3\largen}}^c(\RR)$.
The dimension of ${C}(S^3,\confy(\finsetb) \times \RR;\Gamma)$ is the same as the dimension of $(S^2)^{j_E(E(\Gamma))}$.
The quotient of $\check{C}(S^3,\confy(\finsetb) \times \RR;\Gamma)$ by global vertical translations
is a semi-algebraic set with dimension one less. Thus, the image of $g(\Gamma)$ is a compact semi-algebraic subset of $(S^2_H)^{3\largen}$ of codimension at least one. Its complement is thus an open dense semi-algebraic subset of $(S^2_H)^{3\largen}$.
\eop

Note the following easy lemma.
\begin{lemma}
\label{lemineqnorm}
 Let $a$ and $h$ denote two vectors of $\RR^n$ such that $a$ and $a+h$ are different from $0$.
Then we have \begin{equation*}\lefnorm{\frac{1}{\norm{a+h}}(a+h) - \frac{1}{\norm{a}}a} \leq 2\frac{\norm{h}}{\norm{a}}.\end{equation*}
\end{lemma}
\bp The left-hand side can be written as \begin{equation*}\lefnorm{\frac{1}{\norm{a}}h +\left(\frac{1}{\norm{a+h}} - \frac{1}{\norm{a}} \right)(a+h) }.\end{equation*}
So it is less than or equal to \begin{equation*}\left(\frac{\norm{h}}{\norm{a}} + \frac{\left|\norm{a}-\norm{a+h}\right|}{\norm{a}}\right).\end{equation*}

\eop

\begin{lemma}
\label{lemconftangboun} 
Under the assumptions of Lemma~\ref{lemexisbonform}, assume that the image of  the configuration $\confy \colon \finsetb \hookrightarrow \drad{1}$ contains $0$. 
Equip $(S^2_H)^{3\largen}$ with the distance coming from the Euclidean norm of $(\RR^3)^{\underline{3\largen}}$. Let $\varepsilon \in \bigl]0,\frac1{20^{12\largen}}\bigr[$ be such that the ball
$B(\vecw_1,\vecw_2,\dots,\vecw_{3\largen})$ centered at $(\vecw_1,\vecw_2,\dots,\vecw_{3\largen})$ of radius
$24 \largen \varepsilon^{\frac1{12 \largen}}$ of $(S^2_H)^{3\largen}$ sits in the subset $\CO(\largen,\confy)$ of Lemma~\ref{lemexisbonform}.

Let $\tanghcyll \colon \sourcetl \hookrightarrow \hcylc$ be a long tangle of $\crats(\hcylc)$ whose bottom and top planar configurations are subconfigurations (i.e., restrictions) of $\confy$.

Let $\diskstwo(\vecw_i,\varepsilon)$ be the disk of radius $\varepsilon$ centered at $\vecw_i$ in $S^2$.
For $i \in \underline{3\largen}$ and $Y_i \in \diskstwo(\vecw_i,\varepsilon)$, let $\propP(Y_i)$ be a propagating chain of $(C_2(\rats(\hcylc)),\tau)$ that coincides with
$p_{\tau}^{-1}(Y_i)$ on the domain $D(p_{\tau})$ of Notation~\ref{notratsrest}.
Then for any $\Gamma \in \cup_{k=1}^{3\largen}\CD_{k,\underline{3\largen}}^c(\sourcetl)$, we have
\begin{equation*} \bigcap_{e \in E(\Gamma)}p_e^{-1}\Bigl(\cup_{Y_{j_E(e)}\in B(\vecw_{j_E(e)},\varepsilon)}\propP\bigl(Y_{j_E(e)}\bigr)\Bigr)
\subset {C}\Bigl(\rats_{\frac{1}{\sqrt{\varepsilon}},\left[-\frac{1}{\sqrt{\varepsilon}},\frac{1}{\sqrt{\varepsilon}}\right]}(\hcylc),\tanghcyll;\Gamma\Bigr).\end{equation*}

For $i \in \underline{3\largen}$, let $\omega_{i,S^2}$ be a two-form of $S^2$
supported in 
$\diskstwo(\vecw_i,\varepsilon)$. Then for any family $(\omega_i)_{i \in \underline{3 \largen}}$ of closed propagating forms $\omega_i$ of $(C_2(\rats(\hcylc)),\tau)$ restricting to $D(p_{\tau})$ as $p_{\tau}^{\ast}(\omega_{i,S^2})$,
the support of $\bigwedge_{e \in E(\Gamma)}p_e^{\ast}(\omega(j_E(e))) $ is included in ${C}(\rats_{\frac{1}{\sqrt{\varepsilon}},\left[-\frac{1}{\sqrt{\varepsilon}},\frac{1}{\sqrt{\varepsilon}}\right]}(\hcylc),\tanghcyll;\Gamma)$.
\end{lemma}
\bp
Fix a connected $\underline{3 \largen}$-numbered Jacobi diagram $\Gamma$ on the domain $\sourcetl$ of $\tanghcyll$. Note that $\Gamma$ has at most $6 \largen$ vertices.
Define a sequence $\alpha_1, \dots,\alpha_{6\largen}$ 
by \begin{equation*}\alpha_k=\varepsilon^{\frac{1-k}{12\largen}}.\end{equation*}

Since $\varepsilon$ is smaller than $\frac1{20^{12\largen}}$, we have $\varepsilon^{-\frac1{12\largen}}>20$ and $\alpha_{2}> 20$.

Define an open covering of $\check{C}(\crats(\hcylc),\tanghcyll;\Gamma)$ associated to colorings of the vertices of $\Gamma$ by colors blue and $k$, for $k \in \underline{6 \largen}$, such that
\begin{itemize}
\item vertices of color $1$ go to (or belong to)
$\mathring{\rats}_{3,\left[-2,3\right]}(\hcylc)$,
 \item blue vertices and vertices of color $k\geq 2$ do not go to  
 $\rats_{2,\left[-1,2\right]}(\hcylc)$,
 with Notation~\ref{notratsrest},
\item any vertex of color $2$ is connected to a vertex of color $1$ by an edge of $\Gamma$, 
\item any vertex of color $2$ is at a distance smaller than 
$5\alpha_2$
from $(0,0) \in \RR^3=\CC \times \RR$ (with respect to the Euclidean norm of $\RR^3$),
\item for $k$ such that $2 \leq k \leq  6 \largen -1$, any vertex of color $(k+1)$ is connected to a vertex of color $k$ by an edge of length smaller than 
$5\alpha_{k+1}$,\footnote{This edge length makes sense since vertices of color $k\geq 2$ belong to  
 $\rats^c_{2,\left[-1,2\right]}(\hcylc) \subset \RR^3$.}
\item when there is an edge of $\Gamma$ between a blue vertex and a vertex colored by $1$, the distance between the blue vertex and $(0,0)$ is greater than 
$3\alpha_2$ 
(with respect to the Euclidean distance of $\RR^3$),
\item when there is an edge of $\Gamma$ between a blue vertex and a vertex colored by $k$ for  $2 \leq k \leq  6 \largen -1,$ the distance between the two vertices is greater than 
$3\alpha_{k+1}$. (Since $\Gamma$ has at most $6\largen$ vertices, if there is a blue vertex, then no vertex can be colored by $6\largen$.)
\end{itemize}
The subset $U(\colorc)$ of $\check{C}(\crats(\hcylc),\tanghcyll;\Gamma)$, consisting of the configurations that satisfy the above conditions with respect to a coloring $\colorc$ of the vertices is open.
Let us prove that $\check{C}(\crats(\hcylc),\tanghcyll;\Gamma)$ is covered by these sets.
For a configuration $\confc$, color by $1$ its vertices 
that are in  
$\mathring{\rats}_{3,\left[-2,3\right]}(\hcylc)$.
Then color by $2$ all the still uncolored vertices $v$ connected to a vertex of color $1$ (by an edge of $\Gamma$) such that $d(v,(0,0)) <5\alpha_2$. 
Continue by coloring all the possible uncolored vertices connected to a vertex of color $2$ by an edge of length smaller than $5\alpha_{3}$ by $3$, and so on, in order to end up with a coloring, which obviously satisfies the above conditions, by coloring the uncolored vertices blue.

Note that the distance between a vertex colored by $k \geq 2$ and the point $(0,0) \in \RR^3$ is less than
\begin{equation*}\left(5 \sum_{i=2}^k \alpha_i= 5\frac{\alpha_{k+1}-\alpha_2}{\varepsilon^{-\frac{1}{12\largen}}-1} \right)< 5\frac{\alpha_{k+1}}{\frac56\varepsilon^{-\frac{1}{12\largen}}} \leq 6 \alpha_k\leq 6 \alpha_{6\largen},\end{equation*}
where since $\varepsilon^{-\frac1{12\largen}}$ is greater than $20$, we have
\begin{equation*}\left(6 \alpha_{6\largen} = 6  (\varepsilon^{-\frac1{12\largen}})^{6\largen-1} \right)\leq \left (\bigl(\varepsilon^{-\frac1{12\largen}}\bigr )^{6\largen} = \frac{1}{\sqrt{\varepsilon}}\right).\end{equation*}
Thus, the vertices colored by some $k$ are in $\rats_{\frac{1}{\sqrt{\varepsilon}},\left[-\frac{1}{\sqrt{\varepsilon}},\frac{1}{\sqrt{\varepsilon}}\right]}(\hcylc)$.

Let us prove that an open set $U(\colorc)$ associated to a coloring $\colorc$ for which the color blue appears
cannot intersect $\bigcap_{e \in E(\Gamma)}p_e^{-1}(\cup_{Y_{j_E(e)}\in B(\vecw_{j_E(e)},\varepsilon)}\propP(Y_{j_E(e)}))$. Fix such a coloring $\colorc$. Remove from $\Gamma$ all the edges that do not contain a blue vertex (without removing their ends). Let $\Gamma_b$ be a connected component with at least one blue vertex of the obtained graph. It has blue vertices, which are trivalent or univalent (in $\Gamma$ and $\Gamma_b$). The blue univalent vertices go to $\confy \times \RR$. Color its other vertices colored by some $k$ red.
Red vertices may have $1$, $2$, or $3$ adjacent edges in $\Gamma_b$.  Let $\Gamma^{\prime}_b$ be the uni-trivalent graph obtained by blowing up $\Gamma_b$ at its red vertices by replacing such a vertex by a red univalent vertex for each adjacent edge. Color the edges between blue vertices blue, and the edges between a blue vertex and a red one purple.
To a configuration of $U(\colorc)$ in 
$ \bigcap_{e \in E(\Gamma)}p_e^{-1}(\cup_{Y_{j_E(e)}\in B(\vecw_{j_E(e)},\varepsilon)}\propP(Y_{j_E(e)}))$
 associate the configuration of $\Gamma^{\prime}_b$ obtained by sending all the red vertices to $\origino=(0,0)$, leaving the positions of the blue vertices unchanged.
Thus, 
\begin{itemize}
 \item the direction of a blue edge numbered by $i$ is at a distance less than $\varepsilon$ from $\vecw_i$,
\item the direction of a purple edge numbered by $p$ is at a distance less than $(\varepsilon+4\varepsilon^{\frac1{12\largen}})$ from $\vecw_p$.
\end{itemize}
Let us justify the second assertion. Let $\vertexb$ denote the blue vertex of the purple edge numbered by $p$, and let $\vertexr$ be its red vertex in the configured graph $\Gamma$.
Assume that $\vertexr$ is colored by $k$ with $k \geq 2$. Then we have
\begin{equation*}d(\vertexb,\vertexr)>3\alpha_{k+1},\; (\origino,\vertexb) \in D(p_{\tau}),\;(\vertexr,\vertexb) \in D(p_{\tau}),\; \mbox{and} \; d(\origino,\vertexr)<6\alpha_{k}.\end{equation*} 
Since the configuration consisting of $\vertexb$ and $\vertexr$ in $\Gamma$ is in some $\propP(Y_p)$, we have
\begin{equation*}\left\lvert\pm \frac{1}{\norm{\bigvec{\vertexr\vertexb}}}\bigvec{\vertexr\vertexb} - \vecw_p\right\rvert \leq \varepsilon.\end{equation*}
Apply Lemma~\ref{lemineqnorm} with $a=\bigvec{\vertexr\vertexb}$ and $h=\bigvec{\origino\vertexr}$ to obtain
\begin{equation*}\lefnorm{\frac{1}{\norm{\bigvec{\origino\vertexb}}}\bigvec{\origino\vertexb} - \frac{1}{\norm{\bigvec{\vertexr\vertexb}}}\bigvec{\vertexr\vertexb}} \leq \frac{2\norm{\bigvec{\origino\vertexr}}}{\norm{\bigvec{\vertexr\vertexb}}} < \frac{4\alpha_{k}}{\alpha_{k+1}} \leq 4\varepsilon^{\frac1{12\largen}}.\end{equation*}
This proves the second assertion when $\vertexr$ is colored by $k \geq 2$.

Assume that $\vertexr$ is colored by $1$.
Then we have $d(\vertexb,\origino)>3\alpha_{2}$. 
When the configuration consisting of $\vertexb$ and $\vertexr$ in $\Gamma$ is in some $\propP(Y_p)$, there exists an $\vertexs$ in $\drad{3} \times \left[-2,3\right]$ such that
the direction of $\pm \bigvec{\vertexs\vertexb}$ is at a distance less than $\varepsilon$ from $\vecw_p$. Here, we have $d(\origino,\vertexs)<5\alpha_{1}$. So the direction of $\pm \bigvec{\origino\vertexb}$ is still
at a distance less than $(\varepsilon+4\varepsilon^{\frac1{12\largen}})$ from $\vecw_p$.
Indeed Lemma~\ref{lemineqnorm} applied with $a=\bigvec{\origino\vertexb}$ and $h=\bigvec{\vertexs\origino}$ yields
\begin{equation*}\lefnorm{\frac{1}{\norm{\bigvec{\origino\vertexb}}}\bigvec{\origino\vertexb} - \frac{1}{\norm{\bigvec{\vertexs\vertexb}}}\bigvec{\vertexs\vertexb}} 
 \leq \frac{2\norm{\bigvec{\origino\vertexs}}}{\norm{\bigvec{\origino\vertexb}}} \leq \frac{10\alpha_{1}}{3\alpha_{2}} < 4\varepsilon^{\frac1{12\largen}}.\end{equation*}

Therefore, the directions of the edges numbered by $i$ 
of the configured graph $\Gamma^{\prime}_b$ are at a distance less than 
$(\varepsilon + 4\varepsilon^{\frac1{12\largen}} \leq 8\varepsilon^{\frac1{12\largen}})$ 
from the $\vecw_i$.
But the directions of these edges cannot be in the image of 
$\prod_{e \in E(\Gamma^{\prime}_b)} p_{e,S^2}$ 
according to our conditions. 
Indeed, together with $(\vecw_i)_{i \in \underline{3N} \setminus j_E(E(\Gamma^{\prime}_b))}$,
they form a $3N$-tuple that is at a distance less than 
$3N\times (8\varepsilon^{\frac1{12\largen}})$ from $(\vecw_i)_{i \in \underline{3N}}$.

Thus 
$ \bigcap_{e \in E(\Gamma)}p_e^{-1}(\cup_{Y_{j_E(e)}\in B(\vecw_{j_E(e)},\varepsilon)}\propP(Y_{j_E(e)}))$
does not intersect the open subsets of the coverings that use the blue color.

It is now easy to conclude.
\eop

\bpo{Proof of Proposition~\ref{proptranstang}}
Fix an $\varepsilon$ as in Lemma~\ref{lemconftangboun}.
For a diagram $\Gamma$ of $\Davis^e_{\underline{3\largen}}(\sourcetl)$ and a subset $E$ of $E(\Gamma)$, the map \begin{equation*}q({\Gamma},E)=\prod_{e \in E}p_{\tau} \circ p_e \colon {C}\Bigl(\rats_{\frac{1}{\varepsilon},[-\frac{1}{\varepsilon},\frac{1}{\varepsilon}]}(\hcylc),\tanghcyll;\Gamma\Bigr) \cap \bigcap_{e \in E}p_e^{-1}\bigl(D(p_{\tau})\bigr) \to (S^2)^{E}\end{equation*}  
has an open dense set of regular values. The product of this set by $(S^2)^{\underline{3\largen}\setminus E}$ is also open and dense. So is the intersection $\CI_q$ over all such pairs $(\Gamma,E)$. Thus, there exist $(\vecx_1,\vecx_2,\dots,\vecx_{3\largen})$ in this intersection and $\alpha \in \left]0,\varepsilon\right]$ such that $\prod_{i=1}^{3\largen}B(\vecx_{i},\alpha) \subset  \CI_q \cap \prod_{i=1}^{3\largen}B(\vecw_{i},\varepsilon)$ for the $B(\vecw_{i},\varepsilon)$ of Lemma~\ref{lemconftangboun}.
For $(\vecy_1,\vecy_2,\dots,\vecy_{3\largen}) \in \prod_{i=1}^{3\largen}B(\vecx_{i},\alpha)$, let $\propP(\vecy_i)$ be a propagator of $(C_2(\crats(\hcylc)),\tau)$ restricting to $D(p_{\tau})$ as $p_{\tau}^{-1}(\vecy_i)$, for each $i \in \underline{3\largen}$, as in Lemma~\ref{lemextformfunct}.
Then the $\propP(\vecy_i)$ can be put in general $3\largen$ position as in Section~\ref{secexistransv}, by changing
them only on $\mathring{C}_2(\rats_{2,\left[-1,2\right]}(\hcylc))$ since they satisfy the general position conditions on the boundaries.
Thus, Proposition~\ref{proptranstang} holds with $\largem=\frac{1}{\sqrt{\varepsilon}}$.
For a given $\alpha >0$, the existence of a $\beta >0$ 
such that
\begin{equation*}\bigcap_{e \in E(\Gamma)}p_e^{-1}\biggl(N_{\beta}\Bigl(\propP\bigl(j_E(e)\bigr)\Bigr)\biggr) \subset N_{\alpha}\biggl(I_S\Bigl(\Gamma,\bigl(\propP(i)\bigr)_{i \in \underline{3\largen}}\Bigr)\biggr)\end{equation*}  for any $\Gamma \in \Davis^e_{\underline{3\largen}}(\source)$ can be proved as in the end of the proof of Lemma~\ref{lemepsneigint}.
\eop

\begin{corollary}
\label{corratvariant}
 For forms $\omega_i$ $\beta$-dual (as in Definition~\ref{defformdual}) to the $\propP(i)$ of Proposition~\ref{proptranstang}, for any subset $\finseta$ of $\underline{3\largen}$ with cardinality $3k$, 
\begin{equation*}\Zinv\bigl(\hcylc,\tanghcyll,\tau,\finseta,(\omega_i)_{i \in \underline{3\largen}}\bigr),\end{equation*} which is defined in Theorem~\ref{thmtangconstcompar}, is rational.
\end{corollary}
\bp As in Lemma~\ref{lemexistformtransv}, the involved configuration space integrals can be computed as algebraic intersections of rational preimages of the $\propP(i)$.
\eop

\section{Functoriality}
\label{secfunc}

In this chapter, we prove the functoriality of $\Zinvufrf$, which implies the multiplicativity of $\Zinvuf$ under connected sum. A reader only interested in the latter proof can read the proof by replacing the set $\finsetb$ of strands by $\{0\}$ and by viewing $\crats(\hcylc_j)$ as an asymptotically standard $\RR^3=\CC \times \RR$, identified with $\CC \times \RR$ outside $\drad{1} \times \left[0,1\right]$, where $\drad{1}$ is the unit disk of $\CC$, for $j\in \underline{2}$.

\begin{proposition}
\label{propfunc}
Let $\largen \in \NN$. Let $(\hcylc_1,\tanghcyll_1,\tau_1)$ and $(\hcylc_2,\tanghcyll_2,\tau_2)$ be two composable tangles. There exist volume-one forms $\omega(i,S^2)$ of $S^2$ for $i \in \underline{3\largen}$ such that 
\begin{multline*}\Zinv\Bigl(\hcylc_1\hcylc_2,\tanghcyll_1\tanghcyll_2,\tau_1\tau_2, \finseta,\bigl({\omega}(i,S^2)\bigr)_{i \in \finseta}\Bigr) =\\\sum_{(\finseta_1,\finseta_2) \in P_2(\finseta)}\frac{\cardlef{\finseta_1}!\cardlef{\finseta_2}!}{\cardlef{\finseta}!} \Zinv\Bigl(\hcylc_1,\tanghcyll_1,\tau_1,\finseta_1,\bigl({\omega}(i,S^2)\bigr)\Bigr) \Zinv\Bigl(\hcylc_2,\tanghcyll_2,\tau_2,\finseta_2,\bigl({\omega}(i,S^2)\bigr)\Bigr)\end{multline*}
for any subset $\finseta$ of $\underline{3\largen}$ with cardinality $3n$, with the notation of Theorem~\ref{thmtangconstcompar}.
\end{proposition}
\bp Let $\confy \colon \finsetb \hookrightarrow \drad{1}$ be a planar configuration whose image contains $0$ and the images of the bottom and top configurations of $\tanghcyll_1$ and $\tanghcyll_2$.
Let $B(X_1,X_2,\dots,X_{3\largen})$ be a ball centered at $(X_1,X_2,\dots,X_{3\largen})$ of radius
$24 \largen \varepsilon^{\frac1{12 \largen}}$ of $(S^2_H)^{3\largen}$ sitting in the subset $\CO(\largen,\confy)$ of Lemma~\ref{lemexisbonform}, with $\varepsilon \in \left]0,\frac1{20^{12\largen}}\right[$ as in Lemma~\ref{lemconftangboun}.
For $i \in \underline{3\largen}$, let $\omega(i,S^2)$ be a volume-one form on $S^2$ supported on a disk $\diskstwo(X_i,\varepsilon)$ of $S^2_H$ of radius $\varepsilon$ centered at $X_i$.
 Define an extension $\omega_1(i)$ of $p_{\tau_1}^{\ast}(\omega(i,S^2))$ to $C_2(\rats(\hcylc_1))$ as in Lemma~\ref{lemextformfunct}.
Let $m_{\varepsilon}$ denote the multiplication by $\varepsilon$ in $\RR^3$.
Let $m_{1,\varepsilon}$ be a diffeomorphism from $\crats(\hcylc_1)$ to a manifold denoted by $\crats(\varepsilon\hcylc_1)$ such that $m_{1,\varepsilon}$ coincides with $m_{\varepsilon}$ on $\RR^3 \setminus \Int(\drad{1})\times\left]0,1\right[$.
Call $(\varepsilon \hcylc_1,\varepsilon \tanghcyll_1)$ the intersection of $m_{1,\varepsilon}((\hcylc_1,\tanghcyll_1))$ with the part that replaces $\drad{1} \times \left[0,1\right]$ in $\crats(\varepsilon\hcylc_1)$. Note that $\crats(\varepsilon\hcylc_1)$ is standard outside $\drad{\varepsilon}\times\left[0,\varepsilon\right]$.
Use forms $(m_{1,\varepsilon}^{-1})^{\ast}(\omega_1(i))$ for $(\varepsilon \hcylc_1,\varepsilon \tanghcyll_1)$. So we have
\begin{equation*}\Zinv\Bigl(\varepsilon\hcylc_1,\varepsilon\tanghcyll_1,m_{1,\varepsilon \ast}(\tau_1),\finseta_1,\bigl({\omega}(i,S^2)\bigr)\Bigr)=\Zinv\Bigl(\hcylc_1,\tanghcyll_1,\tau_1,\finseta_1,\bigl({\omega}(i,S^2)\bigr)\Bigr) \end{equation*} for any subset $\finseta_1$ of $\underline{3\largen}$ with cardinality multiple of $3$, with
$m_{1,\varepsilon \ast}(\tau_1)=\tang m_{1,\varepsilon} \circ \tau_1 \circ \left(m_{1,\varepsilon}^{-1} \times m_{\varepsilon}^{-1}\right)$.
Define $\omega_2(i)$ on $C_2(\crats(\hcylc_2))$, so that $\omega_2(i)$ extends $p_{\tau_2}^{\ast}(\omega(i,S^2))$, as in Lemma~\ref{lemextformfunct}.
Let $T_{2,\varepsilon}$ be the vertical translation of $\RR^3$ by $(0,0,1-\varepsilon)$.
Let $m_{2,\varepsilon}$ be a diffeomorphism from $\crats(\hcylc_2)$ to a manifold denoted by $\crats(\varepsilon\hcylc_2)$ such that $m_{2,\varepsilon}$ coincides with the composition $T_{2,\varepsilon} \circ m_{\varepsilon}$ on $\RR^3 \setminus \Int(\drad{1})\times\left]0,1\right[$.
 Call $(\varepsilon \hcylc_2,\varepsilon \tanghcyll_2)$ the intersection of  $m_{2,\varepsilon}((\hcylc_2,\tanghcyll_2))$ with the part that replaces $\drad{1} \times \left[0,1\right]$ in $\crats(\varepsilon\hcylc_2)$. Note that $\crats(\varepsilon\hcylc_2)$ is standard outside $\drad{\varepsilon}\times\left[1-\varepsilon,1\right]$.
Use forms $(m_{2,\varepsilon}^{-1})^{\ast}(\omega_2(i))$ for $(\varepsilon \hcylc_2,\varepsilon \tanghcyll_2)$. So we have
\begin{equation*}\Zinv\Bigl(\varepsilon\hcylc_2,\varepsilon\tanghcyll_2,m_{2,\varepsilon \ast}(\tau_2),\finseta_2,\bigl({\omega}(i,S^2)\bigr)\Bigr)=\Zinv\Bigl(\hcylc_2,\tanghcyll_2,\tau_2,\finseta_2,\bigl({\omega}(i,S^2)\bigr)\Bigr)\end{equation*}
for any subset $\finseta_2$ of $\underline{3\largen}$ with cardinality multiple of $3$, where
$m_{2,\varepsilon \ast}(\tau_2)=\tang m_{2,\varepsilon} \circ \tau_2 \circ \left(m_{2,\varepsilon}^{-1} \times m_{\varepsilon}^{-1}\right)$.

Now, let $(\varepsilon\hcylc_1\hcylc_2,\varepsilon\tanghcyll_1\tanghcyll_2)$ be obtained from 
$(\varepsilon\hcylc_1,\varepsilon\tanghcyll_1)$ by inserting \begin{equation*}(\varepsilon\hcylc_2,\varepsilon\tanghcyll_2) \cap \rats_{\varepsilon,\left[1-\varepsilon,1\right]}(\varepsilon\hcylc_2)\end{equation*} instead of $\drad{\varepsilon}\times\left[1-\varepsilon,1\right]$. Here and below, we use a natural extension of the $\rats_{.,.}$ notation introduced in \ref{notratsrest}.
Define the propagator $\omega(i)$ of $(\rats(\varepsilon\hcylc_1\hcylc_2),\tau_1\tau_2)$,
\begin{itemize}
\item to coincide with $(m_{1,\varepsilon}^{-1})^{\ast}(\omega_1(i))$ on $C_2(\rats(\varepsilon\hcylc_1\hcylc_2)\setminus \rats_{2\varepsilon,\left[1-2\varepsilon,1+\varepsilon\right]}(\varepsilon\hcylc_1\hcylc_2))$,
\item to coincide with $(m_{2,\varepsilon}^{-1})^{\ast}(\omega_2(i))$ on $C_2(\rats(\varepsilon\hcylc_1\hcylc_2)\setminus \rats_{2\varepsilon,\left[-\varepsilon,2\varepsilon\right]}(\varepsilon\hcylc_1\hcylc_2))$, and
\item to be zero on $
\rats_{.1,\left[-.1,.1\right]}(\varepsilon\hcylc_1\hcylc_2) \times \rats_{.1,\left[.9,1.1\right]}(\varepsilon\hcylc_1\hcylc_2)$ and \begin{equation*}\rats_{.1,\left[.9,1.1\right]}(\varepsilon\hcylc_1\hcylc_2) \times \rats_{.1,\left[-.1,.1\right]}(\varepsilon\hcylc_1\hcylc_2).\end{equation*}
 \end{itemize}
This definition is consistent because the form ${\omega}(i,S^2)$ is supported in $S^2_H$. 
Then compute $\Zinv(\hcylc_1\hcylc_2,\tanghcyll_1\tanghcyll_2,\tau_1\tau_2, \finseta,({\omega}(i,S^2))_{i \in \finseta})$
as \begin{equation*}\Zinv\Bigl(\varepsilon\hcylc_1\hcylc_2,\varepsilon\tanghcyll_1\tanghcyll_2,\tau_1\tau_2, \finseta,\bigl({\omega}(i,S^2)\bigr)_{i \in \finseta}\Bigr)\end{equation*} with these propagators $\omega(i)$.
We are going to prove the following lemma under the above hypotheses.
\begin{lemma}
\label{lemindephor}
For any $\underline{3 \largen}$-numbered Jacobi diagram $\Gamma_1$
on the domain of $\tanghcyll_1$ of degree at most $\largen$, 
the form on ${C}(\rats(\varepsilon\hcylc_1),\varepsilon\tanghcyll_1;\Gamma_1)$ \begin{equation*}\bigwedge_{e \in E(\Gamma_1)}
\bigl(p_e m_{1,\varepsilon}^{-1}\bigr)^{\ast}\Bigl(\omega_1\bigl(j_E(e)\bigr)\Bigr)\end{equation*}
is supported on $C_{\vertsetv(\Gamma_1)}\left(\rats_{.1,\left[-.1,.1\right]}(\varepsilon\hcylc_1)\right) \cap {C}(\rats(\varepsilon\hcylc_1),\varepsilon\tanghcyll_1;\Gamma_1)$.

For any $\underline{3 \largen}$-numbered Jacobi diagram $\Gamma_2$
on the domain of $\tanghcyll_2$ of degree at most $\largen$, 
the form $\bigwedge_{e \in E(\Gamma_2)}
(p_e m_{2,\varepsilon}^{-1})^{\ast}(\omega_2(j_E(e)))$ on ${C}(\rats(\varepsilon\hcylc_2),\varepsilon\tanghcyll_2;\Gamma_2)$
is supported
on $C_{\vertsetv(\Gamma_2)}\left(\rats_{.1,\left[.9,1.1\right]}(\varepsilon\hcylc_2)\right) \cap {C}(\rats(\varepsilon\hcylc_2),\varepsilon\tanghcyll_2;\Gamma_2)$.

For any $\underline{3 \largen}$-numbered Jacobi diagram $\Gamma$
on the domain of $\tanghcyll_1\tanghcyll_2$ of degree at most $\largen$, the form $\bigwedge_{e \in E(\Gamma)}
p_e^{\ast}(\omega(j_E(e)))$ on ${C}(\rats(\varepsilon\hcylc_1\hcylc_2),\varepsilon\tanghcyll_1\tanghcyll_2;\Gamma)$
is supported on \begin{equation*}\cup_{\{\finsetv_1,\finsetv_2\} \in P_2(\Gamma)}C_{\finsetv_1}\left(\rats_{.1,\left[-.1,.1\right]}(\varepsilon\hcylc_1\hcylc_2)\right) \times C_{\finsetv_2}\left( \rats_{.1,\left[.9,1.1\right]}(\varepsilon\hcylc_1\hcylc_2)\right),\end{equation*}
where $P_2(\Gamma)$ denotes the set of partitions $\{\finsetv_1,\finsetv_2\}$ of $V(\Gamma)$ into two disjoint subsets $\finsetv_1$ and $\finsetv_2$ such that no edge of $\Gamma$ has one vertex in $\finsetv_1$ and the other in $\finsetv_2$.
\end{lemma}
Assuming Lemma~\ref{lemindephor}, we can conclude the proof of Proposition~\ref{propfunc} as follows.
Lemma~\ref{lemindephor} implies that both sides of the equality to be proved are sums over pairs $(\Gamma_1,\Gamma_2)$
of $\finseta$-numbered diagrams such that $\Gamma_1$ is a diagram on the domain of $\tanghcyll_1$, $\Gamma_2$ is a diagram on the domain of $\tanghcyll_2$, and $j_E(E(\Gamma_1)) \cap j_E(E(\Gamma_2))=\emptyset$, of terms
\begin{equation*}I\Bigl(\hcylc_1,\tanghcyll_1,\Gamma_1,\bigl({\omega}_1(i)\bigr)_{i \in \finseta}\Bigr)I\Bigl(\hcylc_2,\tanghcyll_2,\Gamma_2,\bigl({\omega}_2(i)\bigr)_{i \in \finseta}\Bigr)\left[\Gamma_1\right]\left[\Gamma_2\right].\end{equation*}
Using Lemma~\ref{lemcoefcoef} to identify the coefficients finishes the proof of Proposition~\ref{propfunc} up to the proof of Lemma~\ref{lemindephor}, which follows.
\eop

\bpo{Proof of Lemma~\ref{lemindephor}}
The first two assertions follow from Lemma~\ref{lemconftangboun}. Let us focus on the third one.
Fix a $\underline{3 \largen}$-numbered Jacobi diagram $\Gamma$ of degree at most $\largen$ on the domain of $\tanghcyll_1\tanghcyll_2$.
For $i \in \underline{6\largen}$, set $\beta_i=\varepsilon \alpha_i$ with the sequence $\alpha_i=\varepsilon^{\frac{1-i}{12\largen}}$ of the proof of Lemma~\ref{lemconftangboun}.
Define an open covering of $\check{C}(\crats(\varepsilon\hcylc_1\hcylc_2),\varepsilon\hcylc_1\hcylc_2;\Gamma)$ associated to colorings of the vertices by colors blue, $(1,k)$, and $(2,k)$, with $k \in \underline{6 \largen}$, such that
\begin{itemize}
 \item blue vertices and vertices of color $(j,k)$ with $j \in \underline{2}$ and $k\geq 2$ do not go to $\rats_{2\varepsilon,\left[-\varepsilon,2\varepsilon\right]}(\varepsilon\hcylc_1\hcylc_2) \cup \rats_{2\varepsilon,\left[1-2\varepsilon,1+\varepsilon\right]}(\varepsilon\hcylc_1\hcylc_2)$,
\item vertices of color $(1,1)$ go to $\mathring{\rats}_{3\varepsilon,\left[-2\varepsilon,3\varepsilon\right]}(\varepsilon\hcylc_1\hcylc_2)$,
\item vertices of color $(2,1)$ go to $\mathring{\rats}_{3\varepsilon,\left[1-3\varepsilon,1+2\varepsilon\right]}(\varepsilon\hcylc_1\hcylc_2)$,
\item for $j \in \underline{2}$, any vertex of color $(j,2)$ is connected by an edge of $\Gamma$ to a vertex of color $(j,1)$ and is at a distance smaller than $5\beta_2$ from $(0,j-1)$ (with respect to the Euclidean norm of $\RR^3$),
\item for $j \in \underline{2}$ and $k$ such that $2 \leq k\leq 6 \largen -1$, any vertex of color $(j,k+1)$ is connected to a vertex of color $(j,k)$ by an edge of length smaller than $5\beta_{k+1}$,
\item when there is an edge of $\Gamma$ between a blue vertex and a vertex colored by $(j,1)$ for $j \in \underline{2}$, the distance between the blue vertex and $(0,j-1)$ is greater than $3\beta_2$, and
\item when there is an edge of $\Gamma$ between a blue vertex and a vertex colored by $(j,k)$ with $j \in \underline{2}$ and $2 \leq k\leq 6 \largen -1$, the distance between the two vertices is greater than $3\beta_{k+1}$.
\end{itemize}
The subset $U(\colorc)$ of $\check{C}(\crats(\varepsilon\hcylc_1\hcylc_2),\varepsilon\hcylc_1\hcylc_2;\Gamma)$ consisting of the configurations that satisfy the above conditions with respect to a coloring $\colorc$ of the vertices is open,
and $\check{C}(\crats(\varepsilon\hcylc_1\hcylc_2),\varepsilon\hcylc_1\hcylc_2;\Gamma)$ is covered by these sets as in the proof of Lemma~\ref{lemconftangboun}. The only additional thing to notice is that a vertex could not be simultaneously colored by $(1,k)$ and by $(2,k^{\prime})$ since a vertex colored by
$(j,k)$ is at a distance less than 
\begin{equation*}6 \beta_k \leq \sqrt{\varepsilon} \leq \frac1{20^6}\end{equation*}
from $(0,j-1)$.
In particular, the vertices colored by $(1,k)$ are in \begin{equation*}\rats_{.1,\left[-.1,.1\right]}(\varepsilon\hcylc_1\hcylc_2),\end{equation*}
and the vertices colored by $(2,k)$ are in $\rats_{.1,\left[.9,1.1\right]}(\varepsilon\hcylc_1\hcylc_2)$.

So the form $\bigwedge_{e \in E(\Gamma)}
p_e^{\ast}(\omega(j_E(e)))$ vanishes on open sets corresponding to colorings for which a vertex $(1,k)$ is connected to a vertex $(2,k^{\prime})$ (by some edge of $\Gamma$), according to the conditions before Lemma~\ref{lemindephor}.

As in the proof of Lemma~\ref{lemconftangboun}, we prove that our form vanishes on open sets $U(\colorc)$ associated to colorings for which the color blue appears. Fix such a coloring $\colorc$ and remove from $\Gamma$ all the edges that do not contain a blue vertex. Let $\Gamma_b$ be a connected component of this graph with at least one blue vertex. It has blue vertices, which are trivalent or univalent in $\Gamma$ and $\Gamma_b$. The blue univalent vertices lie on $\varepsilon \confy \times \RR$.
The other vertices of $\Gamma_b$ are either colored by some $(1,k)$, in which case we color them yellow, or by some $(2,k)$, in which case we color them red.
Red and yellow vertices may have $1$, $2$, or $3$ adjacent edges in $\Gamma_b$.  Let $\Gamma^{\prime}_b$ be the uni-trivalent graph obtained by blowing up $\Gamma_b$ at its yellow and red vertices by replacing such a vertex with a univalent vertex of the same color for each adjacent edge. Color the edges between blue vertices blue, the edges between a blue vertex and a yellow one green,
and the edges between a blue vertex and a red one purple.
To a configuration of $U(\colorc)$ in the support of $\bigwedge_{e \in E(\Gamma)}
p_e^{\ast}(\omega(j_E(e)))$, associate the configuration of $\Gamma^{\prime}_b$ obtained by sending all the yellow vertices to $\origino$ and all the red ones to $(0,1)$, leaving the positions of the blue vertices unchanged.
Thus, 
\begin{itemize}
 \item the direction of a blue edge numbered by $i$ is in the support of $\omega(i,S^2)$ at a distance less than $\varepsilon$ from $X_i$,
\item the direction of a green edge numbered by $g$ is at a distance less than $(\varepsilon+4\varepsilon^{\frac1{12\largen}})$ from $X_g$, (as in
 the proof of Lemma~\ref{lemconftangboun}),
\item the direction of a purple edge numbered by $p$ is at a distance less than $(\varepsilon+4\varepsilon^{\frac1{12\largen}})$ from $X_p$.
\end{itemize}

However, the directions of the edges of $\Gamma^{\prime}_b$ cannot be in the image of $\prod_{e \in E(\Gamma)} p_{e,S^2}$ according to our conditions in the beginning of the proof of Proposition~\ref{propfunc} (the $\varepsilon$ rescaling of $\confy$ does not change the image). Therefore, the support of $\bigwedge_{e \in E(\Gamma)}
p_e^{\ast}(\omega(j_E(e)))$ does not intersect the open subsets of the covering that use the blue color.
It is now easy to conclude.

\eop

For an integer $\largen$, $\Zinv_{\leq \largen}$ denotes the truncation of $\Zinv$ valued in
$\Aavis_{\leq \largen}(\source)=\prod_{j=0}^{\largen} \Aavis_j(\source)$.

\begin{theorem}
\label{thmfuncgen} Let $\largen$ be a natural number, and let $(\omega_{i,S^2})_{i \in \underline{3\largen}}$ be a family of volume-one forms of $S^2$.
For any two composable tangles $(\hcylc_1,\tanghcyll_1,\tau_1)$ and $(\hcylc_2,\tanghcyll_2,\tau_2)$ in parallelized rational homology cylinders,
 we have
\begin{multline*}
 \Zinv\Bigl(\hcylc_1\hcylc_2,\tanghcyll_1\tanghcyll_2,\tau_1\tau_2,.,\bigl(\omega_{i,S^2}\bigr)_{i \in \underline{3\largen}}\Bigr)\\
 =\biggl( \Zinv\Bigl(\hcylc_1,\tanghcyll_1,\tau_1,.,\bigl(\omega_{i,S^2}\bigr)\Bigr) \Zinv\Bigl(\hcylc_2,\tanghcyll_2,\tau_2,.,\bigl(\omega_{i,S^2}\bigr)\Bigr)\biggr)_{\sqcup}
\end{multline*}
with the notation of Theorem~\ref{thmtangconstcompar} and Definition~\ref{defprodbiz}.

For any two composable $J_{bb,tt}$-oriented $q$-tangles $T_1$ and $T_2$, we also have
\begin{equation*}\Zinvufrfneg\bigl(T_1T_2,.,(\omega_{i,S^2})_{i \in \underline{3\largen}}\bigr) =\Bigl( \Zinvufrfneg\bigl(T_1,.,(\omega_{i,S^2})\bigr) \Zinvufrfneg\bigl(T_2,.,(\omega_{i,S^2})\bigr)\Bigr)_{\sqcup},\end{equation*}
with the notation of Definition~\ref{defvarzfqtang}, and
\begin{equation*}\Zinvufrfneg(T_1)\Zinvufrfneg(T_2)=\Zinvufrfneg(T_1T_2).\end{equation*}
\end{theorem}
\bp 
Let us prove the first assertion.
Apply Theorem~\ref{thmtangconstcompar}, with \begin{equation*}\tilde{\omega}(i,1,S^2)=\omega_{i,S^2}\;\;\; \mbox{and}\;\;\;  \tilde{\omega}(i,0,S^2)={\omega}(i,S^2)\end{equation*} with the form ${\omega}(i,S^2)$ of Proposition~\ref{propfunc}.
We get \begin{multline*}\Zinv_{\leq \largen}\bigl(\hcylc_1\hcylc_2,\tanghcyll_1\tanghcyll_2,\tau_1\tau_2,.,(\omega_{i,S^2})\bigr)
  \\=\left(\begin{array}{l}
\left(\prod_{j\in I}\Biglol{\left[0,1\right]}{\eta\bigl(.,p_{\tau}(\ST^+K_j)\bigr)} \#_j\right)\\
\lol{\left[1,0\right] \times \confy_1^-}{\eta_{\finsetb^-,.}}\Zinv\Bigl(\hcylc_1\hcylc_2,\tanghcyll_1\tanghcyll_2,\tau_1\tau_2,.,\bigl({\omega}(i,S^2)\bigr)\Bigr)\lol{\left[0,1\right] \times \confy_2^+}{\eta_{\finsetb^+,.}}\end{array}\right)_{\sqcup},    
       \end{multline*}
where $\confy_i^-$ (resp. $\confy_i^+$) represents the bottom (resp. top) configuration of $\tanghcyll_i$, with $\confy_2^-=\confy_1^+$ and
\begin{multline*}\Zinv\Bigl(\hcylc_1\hcylc_2,\tanghcyll_1\tanghcyll_2,\tau_1\tau_2,.,\bigl({\omega}(i,S^2)\bigr)\Bigr)
\\=\biggl( \Zinv\Bigl(\hcylc_1,\tanghcyll_1,\tau_1,.,\bigl({\omega}(i,S^2)\bigr)\Bigr) \Zinv\Bigl(\hcylc_2,\tanghcyll_2,\tau_2,.,\bigl({\omega}(i,S^2)\bigr)\Bigr)\biggr)_{\sqcup}.\end{multline*}
A neutral factor $\bigl(\lol{\left[0,1\right] \times \confy_1^+}{\eta_{\finsetb_1^+,.}}\lol{\left[1,0\right] \times \confy_2^-}{\eta_{\finsetb_2^-,.}}\bigr)_{\sqcup}$ can be inserted in the middle. So the first equality of the statement becomes clear, up to the behavior of the factors $\lol{\left[0,1\right]}{\eta(.,p_{\tau}(\ST^+K_j))}$ of Definition~\ref{defAholonom}. For these factors, note that a component $K$ of $\tanghcyll_1\tanghcyll_2$ consists of a bunch of components $K_k$ from $\tanghcyll_1$ and $\tanghcyll_2$, for $k$ in a finite set $E$, and that $\eta(\finseta,p_{\tau}(\ST^+K))$ is the sum of the corresponding $\eta(\finseta,p_{\tau}(\ST^+K_k))$.
Lemma~\ref{lemaddhol} ensures that we have \begin{equation*}\Biglol{\left[0,1\right]}{\sum_{k \in E}\eta\bigl(.,p_{\tau}(\ST^+K_k)\bigr)}
=\left(\prod_{k \in E}\lolsp{\left[0,1\right]}{\eta\bigl(.,p_{\tau}(\ST^+K_k)\bigr)}\right)_{\sqcup}\end{equation*} in the commutative algebra
$\Assis(\RR)$.
This finishes the proof of the first assertion of Theorem~\ref{thmfuncgen}. 

Orient $\tanghcyll_1$ and $\tanghcyll_2$ in a compatible way.
Recall $p_1(\tau_1\tau_2)=p_1(\tau_1) + p_1(\tau_2)$. Also recall the associativity of the product $()_{\sqcup}$. According to the first assertion applied to straight tangles with the induced parallelization,
if $\tanghcyll_1$ and $\tanghcyll_2$ are framed by parallels $\tanghcyll_{1\parallel}$ and $\tanghcyll_{2\parallel}$ induced by parallelizations $\tau_k$ 
such that $p_{\tau_k}(\ST^+\tanghcyll_k) \subset S_{WE}$ for $k \in \underline{2}$, then we have  \begin{multline*}\Zinvufrfneg\bigl(\hcylc_1\hcylc_2,\tanghcyll_1\tanghcyll_2,(\tanghcyll_1\tanghcyll_2)_{\parallel},.,(\omega_{i,S^2})_{i \in \underline{3\largen}}\bigr) \\=\Bigl( \Zinvufrfneg\bigl(\hcylc_1,\tanghcyll_1,\tanghcyll_{1\parallel},.,(\omega_{i,S^2})\bigr) \Zinvufrfneg\bigl(\hcylc_2,\tanghcyll_2,\tanghcyll_{2\parallel},.,(\omega_{i,S^2})\bigr)\Bigr)_{\sqcup}.\end{multline*}
This generalizes to any pair $((\tanghcyll_1,\tanghcyll_{1\parallel}),(\tanghcyll_2,\tanghcyll_{2\parallel}))$ of parallelized $J_{bb,tt}$-oriented tangles  
with the invariant $\Zinvufrf$ of framed tangles of Definition~\ref{defnotzfvariant}, as follows.
When $(\tanghcyll_1,\tanghcyll_{1\parallel})$ is not representable as a straight tangle with respect to a parallelization, then $(\tanghcyll_1,\tanghcyll_{1\parallel+1})$ is, where $\tanghcyll_{1\parallel+1}$ is the parallel of $\tanghcyll_1$ such that $(\tanghcyll_{1\parallel+1}-\tanghcyll_{1\parallel})$ is homologous to a positive meridian of $\tanghcyll_1$ in a tubular neighborhood of $\tanghcyll_1$ deprived of $\tanghcyll_1$. Thus, the known behavior of $\Zinvufrf$ under such a framing change yields the second equality of the statement when $T_1$ and $T_2$ have injective bottom and top configurations. According to Definitions~\ref{defvarzfqtang} and \ref{defgenlol}, the second equality is also true when $T_1$ and $T_2$ are $q$-braids, thanks to the multiplicativity of $\lol{.}{.}$ with respect to the product of Definition~\ref{defprodbiz} in Lemma~\ref{lemmultlol}.
Thus, Definition \ref{defvarzfqtang}
of $\Zinvufrf$ implies the second equality for general $q$-tangles. The third equality is a direct consequence of the second one when $\omega_{i,S^2}=\omega_{S^2}$ for all $i$.
 \eop

\section{Insertion of a tangle in a trivial \texorpdfstring{$q$-braid}{q--braid}.}

In this section, we prove the following result, which is the cabling property of Theorem~\ref{thmmainfunc} generalized to all variants of the invariant $\Zinvufrf$ of Definition~\ref{defvarzfqtang}.

\begin{proposition}
\label{propcabtrivbraid}
Let $\finsetb$ be a finite set with cardinality greater than $1$. Let $\confy \in \ccompuptd{\finsetb}{\CC}$, let $\confy \times \left[0,1\right]$ denote the corresponding $q$-braid, and let $K$ be a strand of $\confy \times \left[0,1\right]$. Let $\tanghcyll$ be a $q$-tangle with domain $\sourcetl$.
Then \begin{equation*}\Zinvufrfneg\bigl(\left(\confy \times \left[0,1\right]\right)\left(\tanghcyll/K\right)\bigr)\end{equation*} is obtained from $\Zinvufrfneg(\tanghcyll)$ by the natural injection from $\Aavis(\sourcetl)$ to $\Aavis\bigl(({\finsetb} \times \RR)(\sourcetl/K)\bigr)$.

Furthermore, if $\tanghcyll$ is $J_{bb,tt}$-oriented, for any $\largen \in \NN$, for any subset $\finseta$ of $\underline{3\largen}$ whose cardinality is a multiple of $3$, and for any family of volume-one forms $(\omega_{i,S^2})_{i \in \underline{3\largen}}$,
$\Zinvufrfneg\bigl(\left(\confy \times \left[0,1\right]\right)\left(\tanghcyll/K\right), \finseta,(\omega_{i,S^2})_{i \in \underline{3\largen}}\bigr)$ is obtained from $\Zinvufrfneg\bigl(\tanghcyll,\finseta,(\omega_{i,S^2})_{i \in \underline{3\largen}}\bigr)$ by the natural injection from $\Aavis(\sourcetl)$ to $\Aavis\bigl(({\finsetb} \times \RR)(\sourcetl/K)\bigr)$.
\end{proposition}

We first prove the following particular case of Proposition~\ref{propcabtrivbraid}.

\begin{lemma}
Under the hypotheses of Proposition~\ref{propcabtrivbraid}, furthermore assume that $\confy \in \cinjuptd{\finsetb}{\CC}$ and $\tanghcyll$ is a $J_{bb,tt}$-oriented framed tangle represented by a tangle embedding
\begin{equation*}\tanghcyll \colon \sourcetl \hookrightarrow \hcylc\end{equation*}
with injective bottom and top configurations. Let $\largen$ be a positive integer.
Then there exists a family $(\omega(i,S^2))_{i \in \underline{3\largen}}$ of volume-one forms of $S^2$ such that
\begin{equation*}\Zinvufrfneg\Bigl(\left(\confy \times \left[0,1\right]\right)\left(\tanghcyll/K\right), \finseta,\bigl(\omega(i,S^2)\bigr)_{i \in \underline{3\largen}}\Bigr)\end{equation*} is obtained from $\Zinvufrfneg(\tanghcyll,\finseta,(\omega(i,S^2))_{i \in \underline{3\largen}})$ by the natural injection from $\Aavis(\sourcetl)$ to $\Aavis\bigl(({\finsetb} \times \RR)(\sourcetl/K)\bigr)$ for any subset $\finseta$ of $\underline{3\largen}$ whose cardinality is a multiple of $3$.
\end{lemma}
\bp Without loss of generality, translate and rescale $\confy$ so that $K=\{0\} \times \left[0,1\right]$.
Let $\eta \in \left]0,1\right[$ be the distance between $K$ and the other strands of $\confy \times \left[0,1\right]$.
Because of the known variation of $\Zinvufrf$ under framing changes, there is no loss of generality in assuming that $\tanghcyll$ is straight with respect to a parallelization $\tau$, which we do. 
Let $\confy_1 \colon \finsetb \hookrightarrow \drad{1}$ be a planar configuration whose image contains $\confy$ and the images of the bottom and top configurations of $\tanghcyll$.
Let $B(X_1,X_2,\dots,X_{3\largen})$ be a ball centered at $(X_1,X_2,\dots,X_{3\largen})$ of radius
$24 \largen \varepsilon^{\frac1{12 \largen}}$ of $(S^2_H)^{3\largen}$ sitting in the subset $\CO(\largen,\confy_1)$ of Lemma~\ref{lemexisbonform}, with $\varepsilon \in \left]0,\frac1{20^{12\largen}}\right[$ as in Lemma~\ref{lemconftangboun}.
For $i \in \underline{3\largen}$, let $\omega(i,S^2)$ be a volume-one form on $S^2$ supported on a disk $\diskstwo(X_i,\varepsilon)$ of $S^2_H$.
 Define $\omega_1(i)$ on $C_2(\rats(\hcylc))$ so that $\omega_1(i)$ coincides with $p_{\tau}^{\ast}(\omega(i,S^2))$ on $D(p_{\tau})$, as in Lemma~\ref{lemextformfunct}.
Perform a global homothety $m_{1,\eta\varepsilon}$ of $\crats(\hcylc)$, of ratio $\eta\varepsilon$, where $\drad{1}\times\left[0,1\right]$ is consequently changed to $\drad{\eta\varepsilon}\times\left[0,\eta\varepsilon\right]$.
Call $(\eta\varepsilon \hcylc,\eta\varepsilon \tanghcyll)$ the intersection of the image of the long tangle $(\hcylc,\tanghcyll)$ by this homothety with the part that replaces $\drad{1} \times \left[0,1\right]$, which is now standard outside $\drad{\eta\varepsilon}\times\left[0,1\right]$.
Use forms $(m_{1,\eta\varepsilon}^{-1})^{\ast}(\omega_1(i))$ for $(\eta\varepsilon \hcylc,\eta\varepsilon \tanghcyll)$. So we have
\begin{equation*}\Zinvufrfneg\Bigl(\eta\varepsilon\hcylc,\eta\varepsilon\tanghcyll,m_{1,\eta\varepsilon\ast}(\tau),.,\bigl({\omega}(i,S^2)\bigr)\Bigr)=\Zinvufrfneg\Bigl(\hcylc,\tanghcyll,\tau,.,\bigl({\omega}(i,S^2)\bigr)\Bigr).\end{equation*}

Let $\left(\confy \times \left[0,1\right]\right)\left(\eta\varepsilon \tanghcyll/K\right)$ be the tangle obtained from $\left(\confy \times \left[0,1\right]\right)$ by letting 
$\rats_{\eta\varepsilon,\left[0,1\right]}(\eta\varepsilon\hcylc)$ replace $\drad{\eta\varepsilon}\times\left[0,1\right]$.
Graphs that do not involve vertices of $\rats_{2\eta\varepsilon,\left[-{\eta\varepsilon},2\eta\varepsilon\right]}(\eta\varepsilon\hcylc)$ cannot contribute to \begin{equation*}\Zinvufrfneg\Bigl(\eta\varepsilon\hcylc,\left(\confy \times \left[0,1\right]\right)\left(\eta\varepsilon \tanghcyll/K\right), \finseta,\bigl(\omega(i,S^2)\bigr)_{i \in \underline{3\largen}}\Bigr).\end{equation*}
As in the proof of Lemma~\ref{lemconftangboun}, the only contributing graphs are located in $\crats_{\eta\sqrt{\varepsilon},[-\eta\sqrt{\varepsilon},\eta\sqrt{\varepsilon}]}(\eta\varepsilon\hcylc)$.

We conclude that $\Zinvufrfneg(\eta\varepsilon\hcylc,\left(\confy \times \left[0,1\right]\right)\left(\eta\varepsilon \tanghcyll/K\right), \finseta,(\omega(i,S^2))_{i \in \underline{3\largen}})$ is obtained from $\Zinvufrfneg(\hcylc,\tanghcyll,\finseta,(\omega(i,S^2))_{i \in \underline{3\largen}})$ by the natural injection from $\Aavis(\sourcetl)$ to $\Aavis\bigl(({\finsetb} \times \RR)(\sourcetl/K)\bigr)$.
Since this is also true when $\eta$ is replaced by a smaller $\eta^{\prime}$, this is also true
when $\eta\varepsilon \tanghcyll$ is replaced by a legal composition $\gamma^{-}(\eta\varepsilon \tanghcyll)\gamma^{+}$ for braids $\gamma^{-}$ and $\gamma^{+}$ with constant projections in $\ccompuptd{\finsetb^{\pm}}{\
CC}$, which respectively go from $\eta^{\prime} \confy^-$ to $\eta \confy^-$, and from $\eta \confy^+$ to $\eta^{\prime} \confy^+$ (up to adjusting the parallelizations), thanks to the isotopy invariance of $\Zinvufrf$.
Therefore, this is also true at the limit, when $\eta^{\prime}$ tends to zero, thanks to Lemma~\ref{lemconvthmpoirbraid}. See Definition~\ref{defvarzfqtang}.
\eop

\begin{corollary}
\label{corquatdouze}
Proposition~\ref{propcabtrivbraid} is true when $\confy\in \cinjuptd{\finsetb}{\CC}$.
\end{corollary}
\bp
Recall
Theorem~\ref{thmtangconstcomparbis}, which expresses the variation of \begin{equation*}\Zinvufrfneg\Bigl(\bigl(\confy \times [0,1]\bigr)(\tanghcyll/K), \finseta,(\omega_{i,S^2})_{i \in \underline{3\largen}}\Bigr)\end{equation*} when $(\omega_{i,S^2})_{i \in \underline{3\largen}}$ varies, for $q$-tangles. This variation is given by the insertion of factors on components going from bottom to bottom or from top to top, which are identical in both sides of the implicit equality to be proved, and $\finsetd$-holonomies for the bottom and top configurations, for $\finsetd \subseteq \finseta$. The $\finsetd$-holonomies satisfy the duplication property of Proposition~\ref{propdubraid}. The $\finsetd$-holonomies of the bottom and top configurations of $\tanghcyll$ contribute in the same way to both sides of the equality. The $\finsetd$-holonomies of the bottom and top configurations of $\confy \times \left[0,1\right]$ are inverse to each other.
After the insertion, they are duplicated both at the top and at the bottom on possibly different numbers of strands.
Let $\finsetb^+$ (resp. $\finsetb^-$) be the set of upper (resp. lower) $\infty$-components of $\tanghcyll$. 
Lemma~\ref{lemcom} ensures that for any diagram $\Gamma$ on $\sourcetl$ and for any duplication $\pi(\finsetb^+ \times K^+)^{\ast}$ (resp. $\pi(\finsetb^- \times K^-)^{\ast}$) of the upper part $K^+$ (resp. lower part $K^-$) of the long strand of $K$ by $\finsetb^+ \times \left[1,+\infty\right[$ (resp. $\finsetb^- \times \left]-\infty,0\right] $) of a diagram $\Gamma^{\prime}$ on $\finsetb \times \RR$,
we have \begin{equation*}\Gamma \pi(\finsetb^+ \times K^+)^{\ast}(\Gamma^{\prime})=\pi(\finsetb^- \times K^-)^{\ast}(\Gamma^{\prime})\Gamma\end{equation*}
in $\Aavis\bigl(({\finsetb} \times \RR)(\sourcetl/K)\bigr)$.
Thus, the \say{holonomies} $\pi(\finsetb^- \times K^-)^{\ast}\lol{\left[1,0\right] \times \confy}{\eta_{\finsetb,.}}$ and $\pi(\finsetb^+ \times K^+)^{\ast}\lol{\left[0,1\right] \times \confy}{\eta_{\finsetb,.}}$ cancel, and Proposition~\ref{propcabtrivbraid} is true when $\confy\in \cinjuptd{\finsetb}{\CC}$ as soon as the bottom and top configurations of $\tanghcyll$ may be represented by injective configurations.
When $\tanghcyll=T(\gamma^-)(\hcylc,\tanghcyll,\tanghcyll_{\parallel})T(\gamma^+)$ is a general $q$-tangle and $\gamma^-$ and $\gamma^+$ are paths of configurations, $\left(\confy \times \left[0,1\right]\right)\left(\frac{\tanghcyll}{K}\right)$ is equal to
\begin{equation*}\left(\confy \times \left[0,1\right]\right)\left(\frac{T(\gamma^-)}{K}\right)\left(\confy \times \left[0,1\right]\right)\left(\frac{(\hcylc,\tanghcyll,\tanghcyll_{\parallel})}{K}\right)\left(\confy \times \left[0,1\right]\right)\left(\frac{T(\gamma^+)}{K}\right).\end{equation*}
So the result follows using the Functoriality theorem~\ref{thmfuncgen} and the cabling theorem for $q$-braids (Proposition~\ref{propdubraid}).
\eop

\bpo{Proof of Proposition~\ref{propcabtrivbraid}}
Corollary~\ref{corquatdouze} leaves us with the case in which $\confy$ is a limit configuration.
To treat this case, pick a path
$\gamma \colon \left[0,1\right] \to \ccompuptd{\finsetb}{\CC}$ such that $\gamma(1)=\confy$ and 
$\gamma\left(\left[0,1\right[\right) \subset \cinjuptd{\finsetb}{\CC}$, view $\confy \times \left[0,1\right]$ as the path composition $\overline{\gamma} \gamma(0) \gamma$ where $\gamma(0)$ is thought of as a constant map. If $\confy^{-}$ and $\confy^{+}$ respectively denote the bottom and top configurations of $\tanghcyll$ and if $K^{-}$ and $K^{+}$ denote the strand of $K$ in $\overline{\gamma}$ and in $\gamma$, respectively, then
$\left(\confy \times \left[0,1\right]\right)\left(\frac{\tanghcyll}{K}\right)$ may be expressed as \begin{equation*}T(\overline{\gamma})\left(\frac{\confy^{-}\times \left[0,1\right]}{K^{-}}\right)\left(\gamma(0) \times \left[0,1\right]\right)\left(\frac{\tanghcyll}{K}\right)T({\gamma})\left(\frac{\confy^{+}\times \left[0,1\right]}{K^{+}}\right).\end{equation*}

Use the functoriality theorem~\ref{thmfuncgen}, Corollary~\ref{corquatdouze}, the cabling theorem for $q$-braids (Proposition~\ref{propdubraid}), and the commutation argument in the  above proof to conclude.
\eop

\section{Duplication property}
\label{secduprop}

We are about to show how  $\Zinvufrf$ and all its variants behave under a general parallel duplication of a component going from bottom to top in a tangle.

\begin{proposition}
\label{propduptang}
Let $K$ be a component going from bottom to top or from top to bottom in a $q$-tangle $\tanghcyll$ in a rational homology cylinder $\hcylc$.
Let $\confy $ be an element of $\ccompuptd{\finsetb}{\CC}$ for a finite set $\finsetb$.
Let $\tanghcyll( \confy \times K)$ be the tangle obtained by duplicating $K$ as in Section~\ref{secintqtangle}.
Then we have
\begin{equation*}\Zinvufrfneg\bigl(\tanghcyll( \confy \times K)\bigr)=\pi(\finsetb \times K)^{\ast}\Zinvufrfneg(\tanghcyll)\end{equation*} with a natural extension of Notation~\ref{notationduplication}.
Furthermore, for any $\largen \in \NN$, for any family $(\omega_{i,S^2})_{i \in \underline{3\largen}}$ of volume-one forms of $S^2$, and for any subset $\finseta$ of $\underline{3\largen}$ whose cardinality is a multiple of $3$, we have
\begin{equation*}\Zinvufrfneg\bigl(\hcylc,\tanghcyll( \confy \times K), \finseta,(\omega_{i,S^2})_{i \in \underline{3\largen}}\bigr)=\pi(\finsetb \times K)^{\ast}\Zinvufrfneg\bigl(\hcylc,\tanghcyll,\finseta,(\omega_{i,S^2})_{i \in \underline{3\largen}}\bigr)\end{equation*}
for any $J_{bb,tt}$-orientation of $\tanghcyll$.

\end{proposition}

In order to prove this proposition, we are going to prove the following lemmas.
\begin{lemma}
\label{lempropduptangred}
Let $\tanghcyll \colon \sourcetl \hookrightarrow \hcylc$ be a straight tangle in a parallelized rational homology cylinder $(\hcylc,\tau)$. Let $K$ be a component of $\tanghcyll$ going from bottom to top. Let $\confy$ be an element of $\cinjuptd{\finsetb}{\CC}$ for a finite set $\finsetb$.  Let $\largen \in \NN$.
There exists a family of volume-one forms $(\omega(i,S^2))_{i \in \underline{3\largen}}$ such that
\begin{equation*}\Zinv\Bigl(\hcylc,\tanghcyll(\confy \times K),\tau, \finseta,\bigl(\omega(i,S^2)\bigr)_{i \in \underline{3\largen}}\Bigr)=\pi(\finsetb \times K)^{\ast}\Zinv\Bigl(\hcylc,\tanghcyll,\tau,\finseta,\bigl(\omega(i,S^2)\bigr)_{i \in \underline{3\largen}}\Bigr)\end{equation*} for any subset $\finseta$ of $\underline{3\largen}$ whose cardinality is a multiple of $3$.
\end{lemma}

\begin{lemma}
\label{lemlemimplpropduptangred}
Lemma~\ref{lempropduptangred} implies Proposition~\ref{propduptang}.
\end{lemma}
\bp The known behavior of $\Zinv$ under strand orientation changes for components going from bottom to top of Lemma~\ref{lemrevcompgen} allows us to reduce the proof to the case in which $K$ goes from bottom to top. Lemma~\ref{lempropduptangred}, Theorem~\ref{thmtangconstcomparbis}, and Proposition~\ref{propdubraid} imply that Proposition~\ref{propduptang} holds when $\tanghcyll$ is a straight tangle (with injective bottom and top configurations) and when $\confy \in \cinjuptd{\finsetb}{\CC}$. Then the duplication property for braids of Proposition~\ref{propdubraid} and the functoriality
imply that Proposition~\ref{propduptang} holds if $\confy \in \cinjuptd{\finsetb}{\CC}$ for any $J_{bb,tt}$-oriented $q$-tangle $\tanghcyll$ for which $(K,K_{\parallel})$ can be represented by a straight knot with respect to a parallelization $\tau$ of $\hcylc$ and its associated parallel.
Therefore, Proposition~\ref{propduptang} is also true if $\confy \in \ccompuptd{\finsetb}{\CC}$ by iterating the duplication process as soon as $(K,K_{\parallel})$ is representable by a straight knot.
In particular, it is true when $K$ is a strand of a trivial braid whose framing has been changed so that $lk(K,K_{\parallel})=2$. (Recall Lemma~\ref{lemstraightlk} and Proposition~\ref{proplktangles}.) Thanks to the functoriality of $\Zinvufrf$, since an element whose degree $0$ part is $1$ is determined by its square, Proposition~\ref{propduptang} is also true when $K$ is a strand of a trivial braid whose framing has been changed so that $lk(K,K_{\parallel})=1$.
If our general $(K,K_{\parallel})$ is not representable, then Proposition~\ref{propduptang} is true when $\tanghcyll$ is composed by a trivial braid such that the framing of the strand $I$ that extends $K$ is changed so that $lk(I,I_{\parallel})=-1$.
So it is also true for $\tanghcyll$. 
\eop

Let us introduce some notation for the proof of Lemma~\ref{lempropduptangred}. Choose a tubular neighborhood
\begin{equation*}N_{\eta_0}(K)=\drad{\eta_0} \times \RR_K\end{equation*} of $K=\{0\} \times \RR_K$ in $\rats(\hcylc) \setminus(\tanghcyll(\sourcetl) \setminus K)$ for some $\eta_0$ such that $0<10\eta_0 <1$, where $\drad{\eta_0}$ denotes the disk of radius $\eta_0$ centered at $0$ in $\CC$.
Assume that
the trivialization $\partau$ maps $(d \in\drad{\eta_0} ,k \in \RR_K,e_1=\upvec)$ to an oriented tangent vector to $d \times K$, and that $\partau$ maps $(d,k,(e_2,e_3))$ to the standard frame $(1,i)$ of $\drad{\eta_0} (\times k) \subset \CC$.
Pick a representative $\confy$ of $\confy$ in $\check{C}_{\finsetb}[\drad{1/2}]$.
For $\eta \in \left]0,\eta_0\right]$, let $\tanghcyll(\eta^{2} \confy \times K)$ denote
the tangle obtained from $\tanghcyll$ by replacing $\{0\} \times \RR_K$ by $\eta^{2}\confy \times \RR_K$ in $\drad{\eta_0} \times \RR_K$.

Let us now reduce the proof of Lemma~\ref{lempropduptangred} to the proof of the following lemma.
\begin{lemma}
\label{lempropduptangredbaby}
 There exist $\eta_1 \in \left]0,\eta_0\right]$ and volume-one forms $(\omega(i,S^2))$ of $S^2$ for $i \in \underline{3\largen}$ such that
\begin{equation*}\Zinv\Bigl(\hcylc,\tanghcyll(\eta^{2} \confy \times K),\tau, \finseta,\bigl(\omega(i,S^2)\bigr)_{i \in \underline{3\largen}}\Bigr)=\pi(\finsetb \times K)^{\ast}\Zinv\Bigl(\hcylc,\tanghcyll,\tau,\finseta,\bigl(\omega(i,S^2)\bigr)_{i \in \underline{3\largen}}\Bigr)\end{equation*} for any $\eta \in \left]0,\eta_1\right]$ and any subset $\finseta$ of $\underline{3\largen}$ whose cardinality is a multiple of three.
\end{lemma}

\begin{lemma}
Lemma~\ref{lempropduptangredbaby} implies Lemma~\ref{lempropduptangred}.
\end{lemma}
\bp Definition~\ref{defvarzfqtang} and Lemma~\ref{lemconvthmpoirbraid} allow us to write 
\begin{multline*}
\Zinvufrfneg\Bigl(\hcylc,\tanghcyll(\confy \times K), \finseta,\bigl(\omega(i,S^2)\bigr)_{i \in \underline{3\largen}}\Bigr)\\
\begin{array}{l}
=\lim_{\eta \to 0}\Zinvufrfneg\bigl(\hcylc,\tanghcyll(\eta^{2}\confy \times K), \finseta,(\omega(i,S^2))_{i \in \underline{3\largen}}\bigr)\\
=\pi(\finsetb \times K)^{\ast}\Zinvufrfneg\bigl(\hcylc,\tanghcyll,\finseta,(\omega(i,S^2))_{i \in \underline{3\largen}}\bigl).\end{array}
\end{multline*}
\eop

To prove Lemma~\ref{lempropduptangredbaby}, we need some preliminary lemmas, which involve the following new type of Jacobi diagram.
A \emph{special Jacobi $\underline{3\largen}$-diagram} on $\finsetb \times \RR$ is a connected graph $\Gamma_s$ without looped edges with univalent vertices, trivalent vertices, and one bivalent vertex, equipped with an injection $j_E$ from its set $E(\Gamma_s)$ of edges into $\underline{3\largen}$ and with an isotopy class of injections $j_{\Gamma_s}$ from its set $U(\Gamma_s)$ of univalent vertices into $\finsetb \times \RR$. The space of these diagrams is denoted by 
$\Davis^{e,{\textbf{special}}}_{\underline{3\largen}}(\finsetb \times \RR)$.
For a special Jacobi $\underline{3\largen}$-diagram $\Gamma_s$ with univalent vertices on at least two strands, the space of configurations of $\Gamma_s$ with respect to $\confy$ is the space $\check{\CV}(\confy,\Gamma_s)$ of injections of the set $V(\Gamma_s)$ of vertices of $\Gamma_s$ into $\CC \times \RR$ whose restriction to $U(\Gamma_s)$ is the composition of $\confy \times 1_{\RR}$ with an injection from $U(\Gamma_s)$ into $\finsetb \times \RR$ in the isotopy class $\left[j_{\Gamma_s}\right]$, up to vertical translation, for our representative $\confy \in \check{C}_{\finsetb}[\drad{\frac12}]$.
This space $\check{\CV}(\confy,\Gamma_s)$ is similar to former spaces $\check{\CV}(\confy,\Gamma)$ of Section~\ref{seconeform} and is compactified as in Chapter~\ref{chapzinvfbraid}. See also Lemma~\ref{lemcompcalV}. Its compactification is its closure in $\ccompuptd{\vertsetv(\Gamma)}{\RR^3}$. If $\Gamma_s$ has no univalent vertices or univalent vertices on one strand, then configurations are also considered up to dilation, and the compactification is again the closure in $\ccompuptd{\vertsetv(\Gamma)}{\RR^3}$. The configurations of $\check{\CV}(\confy,\Gamma_s)$ are normalized so that a vertex of $\Gamma_s$ is sent to $\drad{1/2} \times \{0\}$. 

\begin{lemma}
\label{lemregxifunc}

Let $\finsetb_{\infty}$ be a finite set. Let $\confy_{\infty} \colon \finsetb_{\infty} \hookrightarrow \drad{1}$ be a planar configuration.
Let $\largen$ be a natural number.
The set $\CO_s(\largen,\confy_{\infty})$ of points $(X_i)_{i \in \underline{3\largen}}$ of $(S^2_H)^{3\largen}$ that are regular values of
\begin{itemize}
 \item the maps $g(\Gamma)$ of Lemma~\ref{lemexisbonform} associated to $\underline{3\largen}$-numbered Jacobi diagrams $\Gamma$ on $\sqcup_{b \in \finsetb_{\infty}}\RR_b$ and to the configuration $\confy_{\infty}$ and
\item similar maps $g(\Gamma_s)$ associated to special $\underline{3\largen}$-numbered Jacobi diagrams $\Gamma_s$ (as above with one bivalent vertex) on $\sqcup_{b \in \finsetb_{\infty}}\RR_b$ and to the configuration $\confy_{\infty}$
\end{itemize}
is a dense open subset of $(S^2_H)^{3\largen}$.
\end{lemma}
\bp The arguments of Lemma~\ref{lemexisbonform} allow us to prove that the images
of the above maps $g(\Gamma)$ and of the maps $g(\Gamma_s)$, when the univalent vertices of $\Gamma_s$ are on one strand or when $\Gamma_s$ has no univalent vertices, are compact semi-algebraic subsets of $(S^2_H)^{3\largen}$ of codimension at least one. The complement of the union of these images is therefore open and dense.

For a special $\Gamma_s$ with univalent vertices on at least $2$ strands, the images under $g(\Gamma_s)$ of the boundary of the configuration space and of the parts where $g(\Gamma_s)$ is not a submersion are also compact semi-algebraic subsets of $(S^2_H)^{3\largen}$ of codimension at least one.
\eop

\begin{lemma}
\label{lempropforminsidetube}
Assume that the configuration  $\confy_{\infty}$ of Lemma~\ref{lemregxifunc} contains the configuration $\confy \colon \finsetb \hookrightarrow \drad{1}$, the bottom configuration $\confy_-$ of $\tanghcyll$, the top configuration  $\confy_+$ of $\tanghcyll$, and $0$ as subconfigurations.

Let $(X_i)_{i \in \underline{3\largen}}$ be the center of a tiny ball of radius $\varepsilon_0 >0$ in the set $\CO_s(\largen,\confy_{\infty})$ defined in Lemma~\ref{lemregxifunc}.
There exist $\varepsilon_1 \in \bigl]0,\frac{\varepsilon_0}{\sqrt{3\largen}}\bigr[$ and $\largem_1 \in \left]1,+\infty\right[$ such that for any family $(\omega(i,S^2))_{i \in \underline{3\largen}}$ of volume-one forms of $S^2$ supported inside disks $\diskstwo(X_i,\varepsilon)$ of $S^2$ centered at $X_i$ of radius $\varepsilon \in \left]0,\varepsilon_1\right]$, 
\begin{itemize}
 \item for any Jacobi diagram $\Gamma \in \Davis^e_{\underline{3\largen}}(\finsetb_{\infty} \times \RR)$ and for any special Jacobi diagram $\Gamma \in \Davis^{e,{\textbf{special}}}_{\underline{3\largen}}(\finsetb \times \RR)$ with univalent vertices on at most one strand, the support of $\bigwedge_{e \in E(\Gamma)}(p_{S^2}\circ p_e)^{\ast}\bigl(\omega(j_E(e),S^2)\bigr)$ in
${C}(\RR^3, \confy_{\infty} \times K;\Gamma)$ is empty,
\item for any special Jacobi diagram $\Gamma \in \Davis^{e,{\textbf{special}}}_{\underline{3\largen}}(\finsetb \times \RR)$,\begin{itemize}
\item the support of $\wedge_{e \in E(\Gamma)}(p_{S^2}\circ p_e)^{\ast}(\omega(j_E(e),S^2))$ in the space
$\check{\CV}(\confy,\Gamma)$ introduced before Lemma~\ref{lemregxifunc} is contained in disjoint open subsets where $\prod_{e \in E(\Gamma)}(p_{S^2}\circ p_e)$ is a diffeomorphism onto \begin{equation*}\prod_{e \in E(\Gamma)}\mathring{\diskstwo}(X_{j_E(e)},\varepsilon) \end{equation*} and the distance of the images of a vertex under two configurations is at most $\largem_1 \varepsilon$,\footnote{In $\check{\CV}(\confy,\Gamma)$, configurations are normalized so that univalent vertices go to $\confy(\finsetb) \times \RR$ for our fixed representative of $\confy$ and one vertex goes to  $\confy(\finsetb) \times \{0\}$. The distance is the Euclidean distance of $\RR^3$. }
and 
\item the images of the vertices under the configurations of the support are contained in $\drad{\largem_1} \times \left[-\largem_1,\largem_1\right]$.
\end{itemize}
\end{itemize}
\end{lemma}
\bp The definition of $\CO_s(\largen,\confy_{\infty})$ in Lemma~\ref{lemregxifunc} and the hypotheses on $(X_i)_{i \in \underline{3\largen}}$ in Lemma~\ref{lemexispropfunc} guarantee that the first assertion is satisfied for all
$\varepsilon_1 \in \bigl]0,\frac{\varepsilon_0}{\sqrt{3\largen}}\bigr[$. Proceed as in Section~\ref{secexistformtransv},
to reduce the support to disjoint subsets, where $\prod_{e \in E(\Gamma)}(p_{S^2}\circ p_e)$ is a diffeomorphism onto $\prod_{e \in E(\Gamma)}\mathring{\diskstwo}(X_{j_E(e)},\varepsilon_1) $. (It is simpler here.) Then there exists an $\largem$ such that the distance of the images of a vertex under two configurations is at most $\largem\varepsilon$ on the preimage of $\prod_{e \in E(\Gamma)}\mathring{\diskstwo}(X_{j_E(e)},\varepsilon) $ under such a diffeomorphism, for any $\varepsilon \in \left]0,\varepsilon_1\right]$. \eop

Choose a Riemannian metric on $\crats(\hcylc)$, which coincides with the standard metric of $\RR^3$ outside $\crats_{1,\left[0,1\right]}(\hcylc)$, and 
assume that this Riemannian metric restricts as the natural product metric on $N_{\eta_0}(K)$, locally. (Reduce $\eta_0$ if necessary.)

Let $C_{2,\leq 10\eta_0}(N_{\eta_0}(K))$ denote the closure in $C_2(\rats(\hcylc))$ of the space of pairs of points $(x_1,x_2) \in N_{\eta_0}(K)^2$ at a distance less than $10 \eta_0$ from each other.
Naturally extend $p_{\tau}$ to $C_{2,\leq 10\eta_0}(N_{\eta_0}(K))$ by viewing $\drad{\eta_0} \times \RR_K$ as a subspace of $\RR^3$, locally, with the usual formula $p_{\tau}=\frac1{\norm{x_2-x_1}}(x_2-x_1)$.

Let $\Gamma \in \Davis^e_{\underline{3\largen}}(\sourcetl)$.
Set $U_K(\Gamma)=j_{\Gamma}^{-1}(\RR_K)$, where $\RR_K$ is viewed as the domain of $K$.
For $\eta \in \left]0,\eta_0\right]$, let ${C}(\rats(\hcylc),\tanghcyll,\eta;\Gamma)$ be the configuration space obtained from ${C}(\rats(\hcylc),\tanghcyll;\Gamma)$ by replacing the condition that the restriction of the configurations to the set $U_K(\Gamma)$ of univalent vertices on $\RR_K$ factors through $K \circ j_{\Gamma}$,
for a $\Gamma$-compatible injection $j_{\Gamma}$ from $U(\Gamma)$ to $\sourcetl$,
with the condition that $U_K(\Gamma)$ is mapped to the interior of $N_{\eta}(K)$. 
(In other words, the conditions on the restriction of a configuration $\confc$ to $U(\Gamma)$ now only impose that $\confc(U_K(\Gamma)) \subset \mathring{N}_{\eta}(K)$, and that $\confc\vert_{U(\Gamma)\setminus U_K(\Gamma)}$ may be expressed as $L \circ j_{\Gamma}\vert_{U(\Gamma)\setminus U_K(\Gamma)}$ for some $\Gamma$-compatible $j_{\Gamma}$.) There is a natural projection $p_{K}$ from this configuration space ${C}(\rats(\hcylc),\tanghcyll,\eta;\Gamma)$ to $\dorad{\eta}^{U_K(\Gamma)}$,
and ${C}(\rats(\hcylc),\tanghcyll;\Gamma)$ is contained in the preimage $p_{K}^{-1}((0)^{U_K(\Gamma)})$ in
${C}(\rats(\hcylc),\tanghcyll,\eta;\Gamma)$. (This preimage also contains ${C}(\rats(\hcylc),\tanghcyll;\tilde{\Gamma})$ for Jacobi diagrams $\tilde{\Gamma}$ different from $\Gamma$ because the linear order of the vertices of $U_K(\Gamma)$ is not induced by $j_{\Gamma}$.)

\begin{lemma}
\label{lemexispropfunc}
Let $(X_i)_{i \in \underline{3\largen}}$ and $\varepsilon_1$ be as in Lemma~\ref{lempropforminsidetube}.
There exist propagating chains $\propP(i)$ of $(C_2(\rats(\hcylc)),\tau)$ and $\varepsilon_1$-dual propagating forms $\omega(i)$ of $(C_2(\rats(\hcylc)),\tau)$, for $i \in \underline{3\largen}$, and $\eta_2 \in \left]0,\eta_0\right]$ such that
\begin{itemize}
\item we have $\propP(i) \cap D(p_{\tau})=p_{\tau}^{-1}(X_i) \cap D(p_{\tau})$ for any $i \in \underline{3\largen}$,
\item we have $\propP(i) \cap C_{2,\leq 10\eta_2}(N_{\eta_2}(K))=p_{\tau}^{-1}(X_i) \cap C_{2,\leq 10\eta_2}(N_{\eta_2}(K))$, for the above natural extension of $p_{\tau}$ on $C_{2,\leq 10\eta_2}(N_{\eta_2}(K))$,
\item the $\propP(i)$ are in general $3\largen$-position with respect to $(\crats(\hcylc),\tanghcyll,\tau)$
(again as in Definition~\ref{defgenthreenpos}), and
\item $\omega(i)$ restricts to $D(p_{\tau}) \cup C_{2,\leq 10\eta_2}(N_{\eta_2}(K))$ as $p_{\tau}^{\ast}(\omega(i,S^2))$ for some form $\omega(i,S^2)$ as in Lemma~\ref{lempropforminsidetube}.
\end{itemize}

Let  $\Gamma \in \Davis^e_{\underline{3\largen}}(\sourcetl)$.

For $(\omega(i))_{i \in \underline{3\largen}}$ as above and $\eta\in \left]0,\eta_2\right]$, let $\Supp(\Gamma,\eta;(\omega(i)))$ denote the support of $\bigwedge_{e\in E(\Gamma)}p_e^{\ast}(\omega(j_E(e)))$ in ${C}(\rats(\hcylc),\tanghcyll,\eta;\Gamma)$.

For $(\propP(i))_{i \in \underline{3\largen}}$ as above, for any configuration $\confc$ in the discrete set
\begin{equation*}I_S\Bigl(\tanghcyll,\Gamma,\bigl(\propP(i)\bigr)\Bigr)=C\bigl(\rats(\hcylc),\tanghcyll;\Gamma\bigr) \bigcap \cap_{e \in E(\Gamma)}p_e^{-1}\Bigl(\propP\bigl(j_E(e)\bigr)\Bigr),\end{equation*} and for any edge $e$ of $E(\Gamma)$, $p_e(c)$ is in the interior of a $4$-cell $\Delta_e(c)$ of $\propP(j_E(e))$, and ${C}_2(\rats(\hcylc))$ is diffeomorphic to $\drad{\varepsilon_1} \times \Delta_e(c)$ near $p_e(c)$, where $\drad{\varepsilon_1}$ is the local fiber of a tubular neighborhood of $\Delta_e(c)$.  Let $p_{c,e}$ be the associated projection onto $\drad{\varepsilon_1}$ in a neighborhood of $c$ in ${C}(\rats(\hcylc),\tanghcyll,\eta;\Gamma)$. Without loss of generality, assume that $\omega(j_E(e))$ may be expressed as $p_{c,e}^{\ast}(\omega_{\varepsilon}(c,e))$ locally, for a $2$-form $\omega_{\varepsilon}(c,e)$ supported on $\drad{\varepsilon}$, such that $\int_{\drad{\varepsilon}} \omega_{\varepsilon}(c,e)$ is the rational coefficient of the above $4$-cell $\Delta_e(c)$ in $\propP(j_E(e))$. 

With this notation, there exist $\varepsilon \in \left]0,\varepsilon_1\right]$, propagating chains $\propP(i)$, $\varepsilon$-dual propagating forms $\omega(i)$ as above, $\eta_3 \in \left]0,\eta_2\right]$, and $\largem_2>1$ such that for any $\eta \in \left]0,\eta_3\right]$ and for any Jacobi diagram $\Gamma \in \Davis^e_{\underline{3\largen}}(\sourcetl)$,
\begin{itemize}
 \item $\Supp(\Gamma,\eta;(\omega(i)))$ is contained in disjoint submanifolds, indexed by the configurations $\confc$ of $I_S(\Gamma,(\propP(i)))$, and which contain those, where $p_K \times \bigl(p_E=\prod_{e \in E(\Gamma)} p_{c,e}\bigr)$ is a diffeomorphism onto $\dorad{\eta}^{U_K(\Gamma)} \times \dorad{\varepsilon}^{E(\Gamma)}$, and where the distance of the images of a vertex under two configurations of $(p_K\times p_E)^{-1}\bigl(\dorad{\eta}^{U_K(\Gamma)} \times \{W\}\bigr)$ is at most $\largem_2 \eta$ for any $W \in \dorad{\varepsilon}^{E(\Gamma)}$,
\item the involved configurations map vertices of $V(\Gamma) \setminus U_K(\Gamma)$ at a distance greater than $9 \largem_2 \eta$ from $K$, and
\item they map two distinct vertices of $U_K(\Gamma)$ at a distance greater than $9 \largem_2 \eta$ from each other.
\end{itemize}
\end{lemma}
\bp 
The existence of the $\propP(i)$ in general $3\largen$-position with prescribed behavior near the boundary can be proved as in Section~\ref{secexistransv}. 
Fix a graph $\Gamma \in \Davis^e_{\underline{3\largen}}(\sourcetl)$.
Once the $\propP(i)$ are in general $3\largen$-position, for such a given graph $\Gamma$, 
$I_S(\Gamma,(\propP(i)))$
consists of a finite number of isolated intersection points at which trivalent vertices cannot be on $K$. Indeed, this would correspond to a degenerate configuration for a graph in which $3$ univalent vertices replace the trivalent vertex on $K$.

The existence of a family $(\omega(i))_{i \in \underline{3\largen}}$ of propagating forms of $\rats(\hcylc)$, $\varepsilon_1$-dual to $\propP(i)$, which may be expressed as $\projp_{\tau}^{\ast}(\omega(i,S^2))$ on $D(p_{\tau}) \cup  C_{2,\leq 10\eta_2}(N_{\eta_2}(K))$, with respect to a family $(\omega(i,S^2))_{i \in \underline{3\largen}}$ of $2$-forms of volume $1$ supported inside a disk $\diskstwo(X_i,\varepsilon)$, as in Lemma~\ref{lempropforminsidetube}, can be proved as in Section~\ref{secexistformtransv}.

Let $\confc \in I_S(\Gamma,(\propP(i)))$.
Transversality implies that the restriction to a neighborhood of $\confc$ in $C(\rats(\hcylc),\tanghcyll;\Gamma)$ of $p_E$ is a submersion, with the notation of the statement. 
This implies that the restriction of $p_K \times p_E$ to a neighborhood of $\confc$ in ${C}(\rats(\hcylc),\tanghcyll,\eta;\Gamma)$ is a submersion, too, so that it is a local diffeomorphism in a neighborhood of $\confc$ in ${C}(\rats(\hcylc),\tanghcyll,\eta;\Gamma)$. After reducing $\eta$ and $\varepsilon$, we obtain a neighborhood $N(\confc)$ of $\confc$ in ${C}(\rats(\hcylc),\tanghcyll,\eta;\Gamma)$ such that $p_K \times p_E$ is a diffeomorphism from $N(\confc)$ onto $\dorad{\eta}^{U_K(\Gamma^)} \times \dorad{\varepsilon}^{E(\Gamma)}$ for all $\confc\in I_S(\Gamma,(\propP(i)))$. 
The compact intersection of $C(\rats(\hcylc),\tanghcyll;\Gamma)$ with the complement of $\cup_{\confc \in I_S(\Gamma,(\propP(i)))} N(\confc)$ is mapped outside $\prod_{e \in E(\Gamma)} \propP(j_E(e))$ by $\prod_{e \in E(\Gamma)}p_e$. Therefore, its image avoids a neighborhood of $\prod \propP(j_E(e))$. Reducing $\varepsilon$ allows us to assume that it avoids the compact closure of the support of $\bigwedge_{e\in E(\Gamma)}p_e^{\ast}(\omega(j_E(e)))$. We may now reduce $\eta$ so that ${C}(\rats(\hcylc),\tanghcyll,\eta;\Gamma)$ does not meet this compact closure outside $\cup_{\confc \in I_S(\Gamma,(\propP(i)))} N(\confc)$.
This can be achieved simultaneously for all the finitely many considered $\Gamma$.

Now, the Jacobians of the corresponding inverse local diffeomorphisms (viewed as maps from $\drad{\eta}^{U_K(\Gamma)} \times \prod_{e\in E(\Gamma)} \drad{\varepsilon}$ to 
$\rats(\hcylc)^{V(\Gamma)}\setminus \diagon$) are bounded (after reducing $\varepsilon$ and $\eta$ if necessary). 

So we get a $\largem_2$ 
for which the distance of the images of a vertex under two configurations of any $\left(p_K\times p_E\right)^{-1}\bigl(\dorad{\eta}^{U_K(\Gamma)} \times \{W\}\bigr)$ in a connected component of $\Supp(\Gamma,\eta;(\omega(i)))$ is at most $\largem_2 \eta$, for all the finitely many considered $\Gamma$.
It is easy to reduce $\eta_3$ so that the last two conditions are satisfied with our given $\largem_2$.
\eop

Set $\Davis^e_{\leq \largen,\underline{3\largen}}(.)=\cup_{k=0}^{\largen}\Davis^e_{k,\underline{3\largen}}(.)$. 
For $\Gamma \in \Davis^e_{\leq \largen,\underline{3\largen}}(\sourcetl(\finsetb \times \RR_K))$, let \begin{equation*}I_S\bigl(\tanghcyll(\eta^{2} \confy \times K),\Gamma,(\propP(i))_{i \in \underline{3\largen}}\bigr)\end{equation*} denote the set of configurations $\confc$ of \begin{equation*}{C}\bigl(\rats(\hcylc),\tanghcyll(\eta^{2} \confy \times K);\Gamma\bigr) \bigcap \bigcap_{e \in E(\Gamma)}p_e^{-1}\bigl(\propP(j_E(e))\bigr)\end{equation*} with respect to the propagating chains $\propP(i)$.

The following crucial lemma justifies the introduction of special Jacobi diagrams.

\begin{lemma}
\label{lemcharconf}
Let $\eta_1$ be the minimum in the set $\{\eta_3,\frac1{(2 \largen -1)8\largen \largem_1}, \frac1{100 \largen \largem_2}\}$ of positive numbers introduced in Lemmas~\ref{lempropforminsidetube} and \ref{lemexispropfunc}.
For any family $(\propP(i))_{i \in \underline{3\largen}}$ of propagating chains as in Lemma~\ref{lemexispropfunc}, 
for any \begin{equation*}\Gamma_{\finsetb} \in \Davis^e_{\leq \largen,\underline{3\largen}}\bigl(\sourcetl(\finsetb \times \RR_K)\bigr),\end{equation*} for any configuration 
$\confc_{\eta_1}$ of the set $I_S(\tanghcyll(\eta_1^{2} \confy \times K),\Gamma_{\finsetb},(\propP(i))_{i \in \underline{3\largen}})$, there exists a continuous map 
\begin{equation*}\begin{array}{lll}\left]0,\eta_1\right] & \to & {C}_V(\Gamma_{\finsetb})\bigl(\rats(\hcylc)\bigr)\\
   \eta & \mapsto & \confc_{\eta}
  \end{array}\end{equation*}
such that, for any $\eta \in \left]0,\eta_1\right]$, $\confc_{\eta}$ belongs to $I_S(\tanghcyll(\eta^{2} \confy \times K),\Gamma_{\finsetb},(\propP(i))_{i \in \underline{3\largen}})$ and
the graph $\Gamma_{\finsetb}$ configured by $\confc_{\eta}$ is the union of 
\begin{itemize}
 \item (small red) special Jacobi diagrams $\Gamma_s$ on $\finsetb \times \RR_K$ of diameter less than $10 \eta$ configured on $\eta^{2} \confy \times K$ in  $N_{\eta}(K)$ (with univalent vertices on at least two strands of $\finsetb \times \RR_K$), and
 \item a uni-trivalent (blue and purple) graph $\plainGamma$ on the domain $\sourcetl$ of $\tanghcyll$ configured so that its
univalent vertices are  
\begin{itemize}
\item either univalent vertices of $\Gamma_{\finsetb}$ on 
$(\tanghcyll \setminus K) \cup (\eta^{2} \confy \times K)$
\item or trivalent vertices of $\Gamma_{\finsetb}$ attached to a bivalent vertex of a (small red) special graph $\Gamma_s$, as in Figure~\ref{figdupJacspec}.
\end{itemize}
\end{itemize}

\bfig
\centering
\begin{tikzpicture}
\useasboundingbox (0,0) rectangle (6.2,3);
\draw [dashed] (0,0) -- (0,.5) .. controls (0,.7) and (.3,1) .. (.5,1)  .. controls (.7,1) and (.85,.7) .. (1,.7) .. controls (1.15,.7) and (1.5,.8).. (1.5,1) -- (1.5,3)
 (.2,0) -- (.2,.5) .. controls (.2,.65) and (.35,.8) .. (.5,.8)  .. controls (.65,.8) and (.8,.5) .. (1,.5) .. controls (1.2,.5) and (1.7,.7) .. (1.7,1) -- (1.7,3);
\draw [dashed] (5,1.1) arc (-180:180:.6);
 \draw [violet, thick] (.5,1) -- (1.7,2.8) (.2,.3) -- (3.6,.4) (1.9,2.1) -- (4,2.3);
\draw [red, thin] (1.5,2) -- (1.95,1.85) -- (1.9,2.1) -- (1.7,2.2)  (1.5,2) --  (1.95,1.85) -- (1.7,1.8);
\draw [blue, very thick] (5,1.1) -- (3.6,.4) --  (4,2.3) -- (5.6,1.7);
\fill [red] (.5,1) circle (1.5pt) (1.7,2.8) circle (1.5pt) (1.5,2) circle (1.5pt) (1.7,2.2) circle (1.5pt) (1.9,2.1) circle (1.5pt) (1.7,1.8) circle (1.5pt) (1.95,1.85) circle (1.5pt) (.2,.3) circle (1.5pt);
\fill [blue] (5,1.1) circle (2.5pt) (5.6,1.7) circle (2.5pt) (4,2.3) circle (2.5pt) (3.6,.4) circle (2.5pt);
\draw (1.9,2.1) node[above]{\scriptsize $\alpha_1$};
\end{tikzpicture}
\caption{A configured diagram $\Gamma_{\finsetb}$ with its small red vertices, its small thin red edges, its big blue vertices (the rightmost four), and its thick blue edges (the rightmost three)}
\label{figdupJacspec}
\end{figure}

Furthermore, the configuration $\confc_{\eta}$ arises as a transverse intersection point, and the intersections that involve at least one (red) special Jacobi diagram cancel algebraically. For a fixed $\eta$, the remaining configurations are in natural one-to-one correspondence---independent of $\eta$---with triples $(\plainGamma, f, \confc^{\prime})$,
where $\plainGamma \in \Davis^e_{\leq \largen,\underline{3\largen}}(\sourcetl)$, 
$U_K(\plainGamma)$ is the set of univalent vertices of $\Gamma_{\finsetb}$ on $\finsetb \times \RR_K$,
$f \in \finsetb^{U_K(\plainGamma)}$, and $\confc^{\prime} \in I_S(\tanghcyll,\plainGamma,(\propP(i))_{i \in \underline{3\largen}})$.
The inverse of this natural one-to-one correspondence maps $(\plainGamma, f, \confc^{\prime})$ to a pair $(\plainGamma_f, \confc_{\eta})$, where \begin{itemize}
\item $\plainGamma_f$ is a Jacobi diagram on $\sourcetl(\finsetb \times \RR_K)$ obtained from $\plainGamma$ by changing the (isotopy class of the) injection from $U_K(\plainGamma)$ to $\RR_K$ to the injection from $U_K(\plainGamma)$ to $\finsetb \times \RR_K$ that maps a vertex $u$ to $f(u) \times \RR_K$ so that the order of vertices on each strand of $\finsetb \times \RR_K$ is induced by their former order on $\RR_K$,
\item $\confc_{\eta}$ belongs to $I_S(\tanghcyll(\eta^{2} \confy \times K),\plainGamma_f,(\propP(i))_{i \in \underline{3\largen}})$, 
\item  $d(\confc_{\eta}(v), \confc^{\prime}(v))$ is smaller than $\eta$ for any vertex $v$ of $\plainGamma$, 
\end{itemize}
and the sign of the algebraic intersection at $\confc_{\eta}$ is the same as the sign of the algebraic intersection at $\confc^{\prime}$, with respect to consistent vertex-orientations of $\Gamma$ and $\Gamma_f$. 

Furthermore, for any $\eta \in \left]0,\eta_1\right]$, the support $\Supp(\tanghcyll(\eta^{2} \confy \times K),\Gamma_{\finsetb};(\omega(i)))$ of $\bigwedge_{e\in E(\Gamma_{\finsetb})}p_e^{\ast}(\omega(j_E(e)))$ in
${C}(\rats(\hcylc),\tanghcyll(\eta^{2} \confy \times K);\Gamma_{\finsetb})$ consists of disjoint neighborhoods $N_{\tanghcyll(\eta^{2} \confy \times K)}(\confc_{\eta})$ of configurations $\confc_{\eta}$ as above, where projections $p_{\confc_{\eta},e}$ as in the statement of Lemma~\ref{lemexispropfunc} make sense, and such that the restriction of $\prod_{e\in E(\Gamma_{\finsetb})}p_{\confc_{\eta},e}$ to $N_{\tanghcyll(\eta^{2} \confy \times K)}(\confc_{\eta})$ is a diffeomorphism onto $\dorad{\varepsilon}^{E(\Gamma_{\finsetb})}$.
\end{lemma}
\bp Let $\Gamma_{\finsetb} \in \Davis^e_{\leq \largen,\underline{3\largen}}(\sourcetl(\finsetb \times \RR_K))$. 
Instead of starting with a configuration $\confc_{\eta_1} \in I_S(\tanghcyll(\eta_1^{2} \confy \times K),\Gamma_{\finsetb},(\propP(i))_{i \in \underline{3\largen}})$, as in the statement,
we consider a configuration $\confc$ in  
$\Supp(\tanghcyll(\eta_4^{2} \confy \times K),\Gamma_{\finsetb};(\omega(i)))$, 
for some $\eta_4 \in \left]0,\eta_1\right]$.
View $\Gamma_{\finsetb}$ as a graph configured by $\confc$. So its vertices become elements of $\rats(\hcylc)$.
Color the vertices of $\Gamma_{\finsetb}$ in $N_{{\eta_1}^2}(K)$ with $(\mbox{red},1)$.
Next color by $(\mbox{red},2)$ its uncolored vertices at a distance less than $4{\eta_1}^2\largem_1$ from a vertex colored by $(\mbox{red},1)$.
For $k \geq 2$, inductively color the still uncolored vertices that are at a distance less than $4{\eta_1}^2\largem_1$ from a vertex $\eltv$ colored by $(\mbox{red},k)$, by $(\mbox{red},k+1)$.

Define the map
\begin{equation*}\begin{array}{llll}
   \rhcylr \colon & \crats(\hcylc) & \rightarrow & \left[0,\eta_1\right]\\
& (z_D,t) \in \drad{{\eta_1}} \times \RR_K &\mapsto &|z_D|\\ 
& x \in \rats(\hcylc) \setminus (\dorad{{\eta_1}} \times \RR_K) &\mapsto &{\eta_1}.
\end{array}
\end{equation*}

Note that a vertex $\eltv$ colored by $(\mbox{red},k)$ with $k\geq 2$ satisfies ${\eta_1}^2 \leq \rhcylr(\eltv) \leq 4k{\eta_1}^2 \largem_1$, by induction. So we have $\rhcylr(\eltv)\leq 8\largen{\eta_1}^2\largem_1 \leq {\eta_1}$ for such a vertex.

Color the vertices that are still uncolored after this algorithm blue.
Color the edges between two blue vertices blue.
Color the edges between a red vertex and a blue one purple.
Also color the edges between red vertices at a distance greater than $8\largen{\eta_1}^2\largem_1$ purple.
Color the remaining edges between two red vertices red.

Remove the red edges between red vertices, and the red vertices that do not belong to a purple edge from $\Gamma_{\finsetb}$.
Blow up the obtained graph $\tilde{\plainGamma}$ at red vertices that belong to at least two purple edges, so that these red vertices become univalent. A red vertex that belongs to $r$ purple edges is transformed into $r$ red vertices during this process.
Let $\plainGamma^{\prime}$ be the obtained configured uni-trivalent graph with blue and purple edges. Its red vertices are in $N_{{\eta_1}}(K)$.
Let $U_K(\plainGamma^{\prime})$ denote the set of red vertices of $\plainGamma^{\prime}$.
The restriction of $\confc$ to $V(\plainGamma^{\prime})$ is in 
${C}(\rats(\hcylc),\tanghcyll,{\eta_1};\plainGamma^{\prime}) \bigcap \Supp\bigl( \bigwedge_{e\in E(\plainGamma^{\prime})}p_e^{\ast}(\omega(j_E(e)))\bigr)$.
So, according to Lemma~\ref{lemexispropfunc}, it is in one of the
 disjoint submanifolds
 where $p_K \times p_E$ restricts as a diffeomorphism 
 onto $\dorad{{\eta_1}}^{U_K(\plainGamma^{\prime})} \times \dorad{\varepsilon}^{E(\plainGamma^{\prime})}$, with $p_E=\prod_{e \in E(\plainGamma^{\prime})} p_{c,e}$, and where the distance of the images of a vertex under two configurations of $\left(p_K\times p_E \right)^{-1}\bigl(\dorad{\eta}^{U_K(\Gamma^{\prime})} \times \{W\}\bigr)$ is at most $\largem_2 \eta$ for any $W \in \dorad{\varepsilon}^{E(\Gamma^{\prime})}$. Set $W_0=p_E(\confc)$. (For $\confc \in I_S(\tanghcyll(\eta_4^{2} \confy \times K),\Gamma_{\finsetb},(\propP(i)))$, we have $W_0=0=(0)_{e \in E(\plainGamma^{\prime})}$.)  
 Then $\confc^{\prime}=(p_K \times p_E)^{-1}((0)_{U_K(\plainGamma^{\prime})},W_0)$ is a configuration
 of a graph $\plainGamma$ on $\sourcetl$, obtained from $\plainGamma^{\prime}$ by adding the data of an isotopy of injections of $U_K(\plainGamma^{\prime})$ into $\RR_K$, where $U_K(\plainGamma)=U_K(\plainGamma^{\prime})$.
According to Lemma~\ref{lemexispropfunc}, no collision of vertices of $\plainGamma$ can occur. So the red vertices of $\plainGamma^{\prime}$ were univalent in $\tilde{\plainGamma}$  (there was no need to blow them up)
and $\tilde{\plainGamma}=\plainGamma^{\prime}$.

Furthermore, $\confc^{\prime}$ maps two red vertices at a distance at least $9 \largem_2 {\eta_1}$ from each other. In particular, two red vertices of $\tilde{\plainGamma}$ are at a distance at least $7 \largem_2 {\eta_1}$ from each other, with respect to $\confc$.

The univalent vertices of $\tilde{\plainGamma}$ are either univalent vertices of $\Gamma_{\finsetb}$, sent to $\eta_4^{2} \confy \times K$ or to $\tanghcyll \setminus K$ by $\confc$, or trivalent vertices of $\Gamma_{\finsetb}$, which belong to a bivalent vertex of a red subgraph of $\Gamma_{\finsetb}$.
Let $\Gamma_R$ be the subgraph of $\Gamma_{\finsetb}$ consisting of its red vertices and of its red edges.

Let $\Gamma_{R,1}$ be a connected component of $\Gamma_R$ such that $\Gamma_{R,1}$ is not reduced to a univalent vertex.
Since two vertices of $\Gamma_{R,1}$ are at a distance at most $(2 \largen -1)8\largen{\eta_1}^2\largem_1 \leq {\eta_1}$ from each other, there is a most one red vertex of $\tilde{\plainGamma}$ in $\Gamma_{R,1}$, and $\confc$ sends $\Gamma_{R,1}$ to a part of $N_{{\eta_1}}(K)$ identified with a part of $\RR^3$ of diameter less than $10{\eta_1}$. So such a $\Gamma_{R,1}$ configured by $\confc$ may be viewed as a graph with straight edges, directed by the $X_i$ if $\confc \in I_S(\tanghcyll(\eta_4^{2} \confy \times K),\Gamma_{\finsetb},(\propP(i)))$, and by the  $\tilde{W}_i \in \diskstwo(X_i,\varepsilon)$ in general. In particular, Lemma~\ref{lempropforminsidetube} implies that $\Gamma_{R,1}$ must be a configured special Jacobi diagram in $\drad{{\eta_1}^2\largem_1} \times \left[x-{\eta_1}^2\largem_1,x+{\eta_1}^2\largem_1\right]$. Its configuration $\eta_4^2\confc_{R,1}$ is determined up to translation along $\RR_K$.
The projection of the bivalent vertex $\alpha_1$ of $\Gamma_{R,1}$ to $\drad{{\eta_1}^2\largem_1}$ is $p_{\CC}({\eta_4}^2\confc_{R,1}(\alpha_1))$. 

Let $\CA$ denote the set of bivalent vertices of $\Gamma_R$.
Write the corresponding configured special Jacobi diagrams  $(\Gamma_{R,\alpha},\confc_{R,\alpha})_{\alpha \in \CA}$.
Let $U_K(\Gamma_{\finsetb})=j_{\Gamma_{\finsetb}}^{-1}(\finsetb \times \RR_K)$. Note the natural inclusions $U_K(\plainGamma)\subseteq U_K(\Gamma_{\finsetb}) \sqcup \CA$ and $\CA \subseteq U_K(\plainGamma)$.
Let $f_{\Gamma_{\finsetb}} \colon U_K(\Gamma_{\finsetb}) \to \finsetb$ be the map that sends $u \in j_{\Gamma_{\finsetb}}^{-1}(\{b\} \times \RR_K)$ to $b$. Let $f$ denote the restriction of $f_{\Gamma_{\finsetb}}$ to $U_K(\plainGamma) \setminus\CA$.

So far, our analysis allows us to associate \begin{equation*}\Phi(\Gamma_{\finsetb},\confc) =\Bigl(\plainGamma,\confc^{\prime},\CA \subset U_K(\plainGamma), f \colon U_K(\plainGamma) \setminus\CA \to \finsetb, (\Gamma_{R,\alpha},\confc_{R,\alpha})_{\alpha \in \CA}\Bigr) \end{equation*}
to our configured graph $(\Gamma_{\finsetb},\confc)$ as above,
where \begin{itemize}
       \item $\plainGamma \in \Davis^e_{\leq \largen,\underline{3\largen}}(\sourcetl)$,
       \item $\confc^{\prime} \in \Supp\bigl(\tanghcyll,\plainGamma;(\omega(i))\bigr) = {C}\bigl(\rats(\hcylc),\tanghcyll;\plainGamma\bigr) \bigcap \Supp\bigl( \bigwedge_{e\in E(\plainGamma)}p_e^{\ast}(\omega(j_E(e)))\bigr)$, 
       \item $p_E(\confc^{\prime})=p_E(\confc)=W_0$,
       \item $\confc^{\prime} \in I_S\bigl(\tanghcyll,\plainGamma,(\propP(i))\bigr)$
    when $\confc \in I_S\bigl(\tanghcyll(\eta_4^{2} \confy \times K),\Gamma_{\finsetb},(\propP(i))\bigr)$,
       \item the $\Gamma_{R,\alpha}$ are special Jacobi diagrams on $\finsetb \times \RR_K$ whose disjoint union is numbered in $\underline{3\largen} \setminus j_E(E(\plainGamma))$,
       \item the $\confc_{R,\alpha}$ are configurations of these diagrams with respect to $\confy \times \RR$ in \begin{equation*}\bigcap_{e\in E(\Gamma_{R,\alpha})}(p_{S^2}\circ p_e^{-1})(\tilde{W}_{j_E(e)}),\end{equation*} (or in $\bigcap_{e\in E(\Gamma_{R,\alpha})}(p_{S^2}\circ p_e^{-1})(X_{j_E(e)})$ when $\confc \in I_S(\tanghcyll(\eta_4^{2} \confy \times K),\Gamma_{\finsetb},(\propP(i)))$).
      \end{itemize}

Let us now show how to reconstruct $(\Gamma_{\finsetb},\confc)$ from the data
\begin{equation*}\bigl(\plainGamma,\confc^{\prime},\CA \subset U_K(\plainGamma), f \colon U_K(\plainGamma) \setminus\CA \to \finsetb, (\Gamma_{R,\alpha},\confc_{R,\alpha})_{\alpha \in \CA}\bigr).\end{equation*}

The graph $\Gamma_{\finsetb}$ is obtained from $\Gamma$ and the $\Gamma_{R,\alpha}$
by gluing $\plainGamma$ and the $\Gamma_{R,\alpha}$ at the vertices $\alpha$. The components of its univalent vertices are determined by those of the univalent vertices of $\Gamma_{R,\alpha}$ and of $\plainGamma$ outside $\CA$, and by $f$. Their order on a strand of $\finsetb \times \RR_K$ is the restriction to the vertices of such a strand of the order induced on $U_K(\Gamma_{\finsetb})$ by letting the ordered set of the univalent vertices of $\Gamma_{R,\alpha}$ replace $\alpha$, for every $\alpha \in \CA$, in the ordered set of the univalent vertices of $\plainGamma$.

Below, we construct a smooth embedding
\begin{equation*}{\textbf{c}} \colon \left]0,\eta_1\right] \times \dorad{\varepsilon}^{E(\Gamma_{\finsetb})} \hookrightarrow {C}_{V(\Gamma_{\finsetb})}\bigl(\rats(\hcylc)\bigr)\end{equation*}
such that $\confc(\eta,{W}_{\finsetb})$ belongs to ${C}(\rats(\hcylc),\tanghcyll(\eta^{2} \confy \times K);\Gamma_{\finsetb})$ for any $(\eta,{W}_{\finsetb}) \in \left]0,\eta_1\right] \times \dorad{\varepsilon}^{E(\Gamma_{\finsetb})}$ and our configuration $\confc$ is equal to ${\textbf{c}}(\eta_4,W_{\finsetb,0})$.

Lemma~\ref{lempropforminsidetube} guarantees that $p_{E,\alpha}=\prod_{e \in E(\Gamma_{R,\alpha})}(p_{S^2}\circ p_e)$ is a diffeomorphism from a neighborhood of $\confc_{R,\alpha}$ in the configuration space of configurations of $\Gamma_{R,\alpha}$ on $\confy \times \RR_K$ up to translation along $\RR_K$ onto $\prod_{e \in E(\Gamma_{R,\alpha})}\mathring{\diskstwo}(X_{j_E(e)},\varepsilon)$, for every $\alpha  \in \CA$. Use an implicit diffeomorphism
$\phi_e\colon {\diskstwo}(X_{j_E(e)},\varepsilon) \to \drad{\varepsilon}$ to identify ${\diskstwo}(X_{j_E(e)},\varepsilon)$ with 
$\drad{\varepsilon}$. Set $W_{\alpha,0}(=(\phi_e \circ p_{S^2}\circ p_e(\confc_{R,\alpha}))_{e \in E(\Gamma_{R,\alpha})})=p_{E,\alpha}(\confc_{R,\alpha})$.
Lemma~\ref{lemexispropfunc} provides a neighborhood of $\confc^{\prime}$ in ${C}(\rats(\hcylc),\tanghcyll,\eta_1;\plainGamma)$, where $p_K \times p_E$ is a diffeomorphism onto $\dorad{\eta_1}^{U_K(\plainGamma)} \times \dorad{\varepsilon}^{E(\plainGamma)}$. Recall $p_E(\confc)=W_0$.
These diffeomorphisms assemble to form a diffeomorphism $\Psi_{\eta}$ from a neighborhood of $\confc(\eta,W)$ in ${C}(\rats(\hcylc),\tanghcyll(\eta^{2} \confy \times K);\Gamma_{\finsetb})$
to $\dorad{\varepsilon}^{E(\Gamma_{\finsetb})}$, whose inverse $\Psi_{\eta}^{-1}={\textbf{c}}(\eta,.)$ is described below.
Write $W_{\finsetb}=(W_e)_{e \in E(\Gamma_{\finsetb})}=((W_{\alpha})_{\alpha \in \CA},W)$, with
$W_{\alpha}=(W_e)_{e \in E(\Gamma_{R,\alpha})}$ and $W=(W_e)_{e \in E(\plainGamma)}$.
$\Psi_{\eta}^{-1}(W)$ is constructed from representatives $\eta^2\confc_{R,\alpha}(\alpha)$ (where the vertical translation parameter is fixed) of the $\eta^2 p_{E,\alpha}^{-1}(W_{\alpha})$ by 
assembling them with
\begin{equation*}\confc_1(\eta,W)=(p_K \times p_E)^{-1}\Biggl(\biggl(\Bigl(\eta^2\confy\bigl(f(u)\bigr)\Bigr)_{u \in U_K(\plainGamma) \setminus \CA}, \Bigl(p_{\CC}\bigl(\eta^2\confc_{R,\alpha}(\alpha)\bigr)\Bigr)_{\alpha \in \CA}\biggr),W\Biggr)\end{equation*}
so that the height (projection onto $\RR_K$) $p_{\RR}(\eta^2\confc_{R,\alpha}(\alpha))$
of $\alpha$ in $\eta^2 \confc_{R,\alpha}(\alpha)$
 coincides with $p_{\RR}(\confc_1(\eta,W)(\alpha))$.

 This ensures that ${\textbf{c}}(\eta,W_{\finsetb})=\Psi_{\eta}^{-1}(W_{\finsetb})$ arises as a transverse intersection, for any $W_{\finsetb}$. This is true when $W_{\finsetb}=0$, in particular. So it is true for the configuration $\confc_{\eta_1}$ of the statement, which may be expressed as ${\textbf{c}}(\eta_1,0)$. The family $\confc_{\eta}$ is the continuous family $\Psi_{\eta}^{-1}(0)$, in this case. The sign of the corresponding intersection is the sign of the Jacobian of $\Psi_{\eta}$.
 Since we started with an arbitrary configuration in $\Supp(\tanghcyll(\eta_4^{2} \confy \times K),\Gamma_{\finsetb};(\omega(i)))$ for some $\eta_4 \in \left]0,\eta_1\right]$, the above arguments also prove the final assertion of the lemma. Let us finally focus on the claimed algebraic cancellation.
From now on, we set $W=0$.

Assume that $\CA \neq \emptyset$. Every $\alpha \in \CA$ is contained in one edge $e(\alpha)$ of $\plainGamma$. Choose $\alpha_0$ in $\CA$ such that $j_E(e(\alpha_0))$ is minimal.
Let $s(\Gamma_{R,\alpha_0})$ be obtained from $\Gamma_{R,\alpha_0}$ by exchanging the labels of the two edges $e_1$ and $e_2$ of $\Gamma_{R,\alpha_0}$ that contain $\alpha_0$, and by reversing their orientations if they both come from $\alpha_0$ or go to $\alpha_0$ as in Lemma~\ref{lemsym}. Let $s_{\alpha_0}(\Gamma_{\finsetb})$ be obtained from 
$\Gamma_{\finsetb}$ by performing the same changes.
Let $s(\confc_{R,\alpha_0})$ be obtained from $\confc_{R,\alpha_0}$ by changing the position of $\confc_{R,\alpha_0}(\alpha_0)$ by a central symmetry with respect to the middle of the two other ends of $e_1$ and $e_2$. 
The intersection point associated to the configured graph \begin{equation*}\bigl(s_{\alpha_0}(\Gamma_{\finsetb}), \confc_2\bigr)=\Phi^{-1}\Bigl(\plainGamma,\confc^{\prime},\CA, f, (\Gamma_{R,\alpha},\confc_{R,\alpha})_{\alpha \in \CA \setminus \{\alpha_0\}}, s(\Gamma_{R,\alpha_0}),s(\confc_{R,\alpha_0}),0\Bigr)\end{equation*}
and the intersection point associated to $(\Gamma_{\finsetb},\confc)$ cancel algebraically, as in Lemma~\ref{lemsym}. (Our process defines an involution on the configured graphs $(\Gamma_{\finsetb},\confc)$ with $\CA \neq \emptyset$ such that a configured graph
and its image cancel.)

Therefore, the configured graphs $(\Gamma_{\finsetb},\confc)$ that contribute to the intersection are the graphs for which $\CA = \emptyset$. They are obtained from some $\plainGamma \in \Davis^e_{\leq \largen,\underline{3\largen}}(\sourcetl)$, some $\confc^{\prime}$, and some $f\colon U_K(\plainGamma) \to \finsetb$ as in the statement.

\eop

\begin{samepage}
\bpo{Proof of Lemma~\ref{lempropduptangredbaby}} 
Lemma~\ref{lemcharconf} implies that 
\begin{itemize}\item 
for any $\plainGamma \in \Davis^e_{\leq \largen,\underline{3\largen}}(\sourcetl)$,
$I(\hcylc,\tanghcyll,\plainGamma,o(\plainGamma),({\omega}(i))_{i \in \underline{3\largen}})$ is the algebraic intersection $I(\hcylc,\tanghcyll,\plainGamma,o(\plainGamma),(\propP(i))_{i \in \underline{3\largen}})$ of the preimages of the propagating chains  $\propP(i)$ in ${C}(\rats(\hcylc),\tanghcyll;\plainGamma)$ with respect to $\plainGamma$,
\item for any $\eta \in \left]0,\eta_1\right]$ and for any
 $\Gamma_{\finsetb} \in \Davis^e_{\leq \largen,\underline{3\largen}}(\sourcetl(\finsetb \times \RR_K))$, we have
\begin{multline*}I\Bigl(\hcylc,\tanghcyll(\eta^{2} \confy \times K),\Gamma_{\finsetb},o(\Gamma_{\finsetb}),\bigl({\omega}(i)\bigr)\Bigr)=\\I\Bigl(\hcylc,\tanghcyll(\eta^{2} \confy \times K),\Gamma_{\finsetb},o(\Gamma_{\finsetb}),\bigl(\propP(i)\bigr)\Bigr), \end{multline*} and
 \item for any $\eta \in \left]0,\eta_1\right]$ and for any subset $\finseta$ of $\underline{3\largen}$ with cardinality $3k$, we have
\begin{multline*}\sum_{\plainGamma \in \Davis^e_{k,\finseta}(\sourcetl)}\coefgambet_{\plainGamma}I\Bigl(\hcylc,\tanghcyll,\plainGamma,\bigl(\propP(i)\bigr)_{i \in \finseta}\Bigr)\pi(\finsetb \times K)^{\ast}\bigl([\plainGamma]\bigr)
\\= \sum_{\Gamma_{\finsetb} \in \Davis^e_{k,\finseta}(\sourcetl(\finsetb \times \RR_K))}\coefgambet_{\Gamma_{\finsetb}}I\Bigl(\hcylc,\tanghcyll(\eta^{2} \confy \times K),\Gamma_{\finsetb},\bigl(\propP(i)\bigr)_{i \in \finseta}\Bigr)\left[\Gamma_{\finsetb}\right].\end{multline*}
\end{itemize} 
 \eop
\end{samepage}

Now, both the second duplication property and the first duplication property  for components going from bottom to top or from top to bottom of Theorem~\ref{thmmainfunc} are proved.
Below, we prove the first duplication property, more generally, in the doubling case.

\begin{lemma}
\label{leminvcap}
Let $\invtrivcar$ be the element of $\Aavis(\RR)$ obtained from 
\begin{equation*}\Zinvufrfneg\left(\assocnum \right) \in \Aavis\biggl(\sourceass\biggr)\end{equation*} by inserting $ \sourceass$ in $\RR$ as indicated by the picture $\sourceassr$. Let $\invtrivcar^{-\frac12}$ be the unique element of $\Aavis(\RR)$ whose degree $0$ part is $1$ such that $(\invtrivcar^{-\frac12})^2\invtrivcar$ equals $1$.
Then we have \begin{equation*}\Zinvufrfneg\Bigl(\tancapred\Bigr)=\Zinvufrfneg\Bigl(\tancuplow\Bigr)=\invtrivcar^{-\frac12}.\end{equation*}
\end{lemma}
\bp Use Theorem~\ref{thmmainfunc} except for the first duplication property, which is about to be proved. The symmetry and the isotopy invariance respectively imply \begin{equation*}\Zinvufrfneg\Bigl(\tancapred\Bigr)=\Zinvufrfneg\Bigl(\tancuplow\Bigr) \mbox{ and }
\Zinvufrfneg\left(\snake\right)= \Zinvufrfneg\left(\snakestr\right)=1.\end{equation*}
The functoriality implies \begin{equation*}\Zinvufrfneg\left(\snake\right)=\Zinvufrfneg\left(\snakedashed\right)=\Zinvufrfneg\Bigl(\botsnake\Bigr)\Zinvufrfneg\left(\assocnum\right)\Zinvufrfneg\Bigl(\topsnake\Bigr).\end{equation*}
The cabling property implies that $\Zinvufrfneg\bigl(\botsnake\bigr)$ is obtained from $\Zinvufrfneg\bigl(\tancuplow\bigr)$ by the map induced by the natural injection from $\tancuplow$
to $\botsnake$. $\Zinvufrfneg\bigl(\topsnake\bigr)$ is obtained similarly from $\Zinvufrfneg\bigl(\tancapred\bigr)$. Since the insertions of $\Zinvufrfneg\bigl(\tancuplow\bigr)$ and
$\Zinvufrfneg\bigl(\tancapred\bigr)$ can be performed at arbitrary places according to Proposition~\ref{propdiagcom}, we have
\begin{equation*}\Zinvufrfneg\Bigl(\tancapred\Bigr)\Zinvufrfneg\Bigl(\tancuplow\Bigr)\invtrivcar =1\end{equation*}
in the algebra $\Aavis(\RR)$. In this algebra, an element whose degree $0$ part is the class of the empty diagram is determined by its inverse. Its square also determines it.
\eop
\begin{lemma}
\label{lemdupcupcap}
The first duplication property of Theorem~\ref{thmmainfunc} is true when $K$ is the unique component of the tangle \tancapred\, or the unique component of \tancuplow. In other words, we have
\begin{equation*}\Zinvufrfneg\left(2 \times \tancuplow\right)=\pi(2 \times \RR)^{\ast}\Bigl(\Zinvufrfneg\bigl(\tancuplow\bigr)\Bigr)\end{equation*} and $\Zinvufrfneg\left(2 \times \tancapred\right)=\pi(2 \times \RR)^{\ast}\left(\Zinvufrfneg(\tancapred)\right)$.
\end{lemma}
\bp Again we have $\Zinvufrfneg\bigl(2 \times\tancuplow\bigr)=\Zinvufrfneg\bigl(2 \times\tancapred\bigr)$ by symmetry.
Thus, $\Zinvufrfneg\bigl(2 \times \tancapred\bigr)$ can be computed from $\Zinvufrfneg\left(2 \times \snake \right)$ as $\Zinvufrfneg\left(\tancapred\right)$ is computed from $\Zinvufrfneg\left( \snake \right)$ in the proof of Lemma~\ref{leminvcap}.
Indeed the boxes $\Zinvufrfneg\bigl(2 \times \tancapred\bigr)$ and $\Zinvufrfneg\bigl(2 \times \tancuplow\bigr)$ can slide across the duplicated strands of \begin{multline*}\pi\Bigl(2 \times \sourceass\Bigr)^{\!\ast}\left(\Zinvufrfneg\left(\assocnum \right)\right)=\Zinvufrfneg\left(\assocnum\Bigl( 2 \times \sourceassone{1}\Bigr)\Bigl( 2 \times \sourceassone{2}\Bigr)\Bigl( 2 \times \sourceassone{3}\Bigr)\right)\\ 
=\pi\Bigl(2 \times \sourceassone{1}\Bigr)^{\!\!\ast} \pi\Bigl(2 \times \sourceassone{2}\Bigr)^{\!\!\ast}\pi\Bigl(2 \times \sourceassone{3}\Bigr)^{\!\!\ast}\left(\Zinvufrfneg\left(\assocnum \right)\right),\end{multline*} according to Lemma~\ref{lemdupcom}, so that we get
\begin{equation*}\Zinvufrfneg\left(2 \times \tancapred\right)^2\pi(2 \times \RR)^{\ast}(\invtrivcar)=1\end{equation*}
in the algebra $\Aavis(\RR \sqcup \RR)$.
Since $\pi(2 \times \RR)^{\ast}$ is an algebra morphism from $\Aavis(\RR)$ to $\Aavis(\RR \sqcup \RR)$, Lemma~\ref{leminvcap} implies
\begin{equation*}\pi(2 \times \RR)^{\ast}\left(\Zinvufrfneg\bigl(\tancapred\bigr)^2\right)\pi(2 \times \RR)^{\ast}(\invtrivcar)=1\end{equation*}
in the algebra $\Aavis(\RR \sqcup \RR)$.
Since the multiplication by
an element whose degree $0$ part is $1$ is injective, and since an element whose degree $0$ part is $1$ is determined by its square, we get $\Zinvufrfneg\bigl(2 \times \tancapred\bigr)=\pi(2 \times \RR)^{\ast}\bigl(\Zinvufrfneg(\tancapred)\bigr)$ as desired.
\eop

\begin{proposition}
\label{propduptangtwo}
Let $K$ be a component of a $q$-tangle $\tanghcyll$ in a rational homology cylinder $\hcylc$.
Let $\tanghcyll( 2 \times K)$ be the tangle obtained by duplicating $K$ as in Section~\ref{secintqtangle}.
Then we have
\begin{equation*}\Zinvufrfneg\bigl(\tanghcyll(2 \times K)\bigr)=\pi(2 \times K)^{\ast}\Zinvufrfneg(\tanghcyll).\end{equation*}
\end{proposition}
\bp A tangle $\tanghcyll$, with a strand $K$ going from bottom to bottom, can be written as a composition
\begin{equation*}\tanghcyll_1\tanghcyll_2= \begin{tikzpicture} [baseline=-.45cm]
\useasboundingbox (-.3,-.55) rectangle (1.9,.25);
\draw [->] (-.12,-.15) .. controls (-.12,.1) and (.12,.1) .. (.12,-.15);
\draw (.7,-.15) -- (.7,.25)  (1.2,-.15) -- (1.2,.25)  (1.7,-.15) -- (1.7,.25)
(-.3,-.55) rectangle (1.9,-.15) (.8,-.35) node{\scriptsize $\tanghcyll_1$};
\end{tikzpicture}\end{equation*}
of some tangle $\tanghcyll_1$, with a cabling $\tanghcyll_2$ of a trivial braid by the replacement of a strand by $\tancapred$, where $K$ is the concatenation of one strand of $\tanghcyll_1$ going from bottom to top, $\tancapred$, and another strand of $\tanghcyll_1$, which goes from top to bottom.
The statement for such a pair follows from Lemma~\ref{lemdupcupcap}, Proposition~\ref{propcabtrivbraid}, and Proposition~\ref{propduptang}, using functoriality.
The case in which $K$ goes from top to top can be treated similarly, by sending $\tancuplow$ below. So the proposition is proved.
\eop

Theorem~\ref{thmmainfunc} is now proved. 
\eopwobp

\begin{remark}
 Theorem~\ref{thmtangconstcomparbis} and Proposition~\ref{propdubraid} do not allow me to prove
\begin{equation*}\Zinvufrfneg\Bigl(\hcylc,\tanghcyll( 2 \times K), \finseta,\bigl(\omega(i,S^2)\bigr)_{i \in \underline{3\largen}}\Bigr)=\pi(2 \times K)^{\ast}\Zinvufrfneg\Bigl(\hcylc,\tanghcyll,\finseta,\bigl(\omega(i,S^2)\bigr)_{i \in \underline{3\largen}}\Bigr)\end{equation*}
for a $J_{bb,tt}$-oriented $q$-tangle $\tanghcyll$, an integer $\largen$, a subset $\finseta$ of $\underline{3\largen}$ whose cardinality is a multiple of $3$, and a family of volume-one forms $(\omega(i,S^2))_{i \in \underline{3\largen}}$.
Indeed, I can see no reason to believe that $\pi(2 \times \RR)^{\ast}\bigl(\lol{\left[0,t\right]}{\eta(.,S_{WE})}\bigr)$ is the product of twice $\lol{\left[0,t\right]}{\eta(.,S_{WE})}$ on the two strands of $2 \times \RR$.
Unfortunately, as noticed in Remark~\ref{rknotzfvarianttwo}, I do not know how to get rid of our noncanonical normalization of 
$\Zinvufrfneg(\hcylc,\tanghcyll,\finseta,(\omega(i,S^2))_{i \in \underline{3\largen}})$  and of the corresponding factors $\lol{\left[0,t\right]}{\eta(.,S_{WE})}$, which might not behave well under duplication.
\end{remark}

\section{Behavior of \texorpdfstring{$\Zinvufrf$}{Zf} with respect to the coproduct}

The behavior of $\Zinvufrf$ with respect to the coproduct described in Theorem~\ref{thmmainfunc} is justified after the statement of Theorem~\ref{thmmainfunc}. Proposition~\ref{propgrouplikefunc} below shows how this behavior generalizes to the variants of $\Zinvufrf$.

Before stating and proving it, let us prove that
the compatibility between product and coproduct implies the following preliminary lemma.
\begin{lemma}
\label{lemprodcoprodbiz}
Say that a map $F$ from the set $\parentp_{(3)}(\underline{3\largen})$ of subsets of $\underline{3\largen}$ whose cardinalities are multiple of $3$ to a space of Jacobi diagrams is \emph{group-like} if
 \begin{equation*}\Delta\bigl(F(\finsetb)\bigr)=\sum_{(\finsetb_1,\finsetb_2) \in P_2(\finsetb)}\frac{\cardlef{\finsetb_1}!\cardlef{\finsetb_2}!}{\cardlef{\finsetb}!}F(\finsetb_1) \otimes F(\finsetb_2)\end{equation*} for any element $\finsetb$ of $\parentp_{(3)}(\underline{3\largen})$.
 
Let $F$ and $G$ be two maps from $\parentp_{(3)}(\underline{3\largen})$ to spaces $\CA_F$ and $\CA_G$ of Jacobi diagrams, such that there is a product from $\CA_F \times \CA_G$ to a space of Jacobi diagrams $\CA_{FG}$.\footnote{  This assumption guarantees that the product $(FG)_{\sqcup}$ of Definition~\ref{defprodbiz} makes sense.}
If $F$ and $G$ are group-like, then $(FG)_{\sqcup}$ is group-like, too.
\end{lemma}
\bp
Let $\finseta \in \parentp_{(3)}(\underline{3\largen})$. We have
\begin{multline*}\Delta\bigl((FG)_{\sqcup}(\finseta)\bigr)=\sum_{(\finsetb,\finsetc) \in P_2(\finseta)}\frac{\cardlef{\finsetb}!\cardlef{\finsetc}!}{\cardlef{\finseta}!}\Delta\bigl(F(\finsetb)\bigr)\Delta\bigl(G(\finsetc)\bigr)
\\=\sum_{(\finsetb_1,\finsetb_2,\finsetc_1,\finsetc_2) \in P_4(\finseta)}\frac{ \cardlef{\finsetb_1}! \cardlef{\finsetb_2}!\cardlef{\finsetc_1}!\cardlef{\finsetc_2}!}{\cardlef{\finseta}!}(F(\finsetb_1)G(\finsetc_1) \otimes F(\finsetb_2)G(\finsetc_2)),
\end{multline*}
with \begin{equation*}(FG)_{\sqcup}(\finseta_1)= \sum_{(\finsetb_1,\finsetc_1) \in P_2(\finseta_1)}\frac{ \cardlef{\finsetb_1}!\cardlef{\finsetc_1}!}{\cardlef{\finseta_1}!}F(\finsetb_1) G(\finsetc_1).\end{equation*}
So we get
\begin{equation*}\Delta\bigl((FG)_{\sqcup}(\finseta)\bigr)=\sum_{(\finseta_1,\finseta_2) \in P_2(\finseta)}\frac{\cardlef{\finseta_1}!\cardlef{\finseta_2}!}{\cardlef{\finseta}!}(FG)_{\sqcup}(\finseta_1) \otimes (FG)_{\sqcup}(\finseta_2).\end{equation*}
\eop

\begin{proposition}
\label{propgrouplikefunc}
For any $q$-tangle $\tanghcyll$ in a rational homology cylinder $\hcylc$, 
for any $\largen \in \NN$, for any subset $\finseta$ of cardinality $3k$ of $\underline{3\largen}$, and for any family $(\omega(i,S^2))_{i \in \underline{3\largen}}$ of volume-one forms of $S^2$, we have
\begin{multline*}\Delta_k\biggl(\Zinvufrfneg\Bigl(\hcylc,\tanghcyll, \finseta,\bigl(\omega(i,S^2)\bigr)_{i \in \underline{3\largen}}\Bigr)\biggr)=\\\sum_{j=0}^k\sum_{\begin{array}{l}{\scriptstyle \finseta_1 \subset\finseta, }\\
{\scriptstyle  \cardlef{\finseta_1}=3j,}\\ 
{\scriptstyle \finseta_2 =\finseta \setminus \finseta_1}\end{array}}
\frac{\cardlef{\finseta_1} ! \cardlef{\finseta_2} !}{\cardlef{\finseta} !}\Zinvufrfneg\Bigl(\hcylc,\tanghcyll, \finseta_1,\bigl(\omega(i,S^2)\bigr)\Bigr) \otimes \Zinvufrfneg\Bigl(\hcylc,\tanghcyll, \finseta_2,\bigl(\omega(i,S^2)\bigr)\Bigr)\end{multline*}
with the coproduct maps $\Delta_n$ defined in Section~\ref{seccoprod}.
\end{proposition}
\bp
Observe that the statement of Proposition~\ref{propgrouplikefunc} is valid for $q$-tangles that can be represented as straight tangles with respect to a parallelization $\tau$ such that $p_1(\tau)=0$, by Theorem~\ref{thmdefsanstauvarzinvf}.
The coefficients are treated as in Lemma~\ref{lemcoefcoef}.

Say that a map $F$ from $\parentp_{(3)}(\underline{3\largen})$ to a space $\CA_F$ of Jacobi diagrams is \emph{cardinality-determined} if it maps any element $\finseta$ of $\parentp_{(3)}(\underline{3\largen})$ to a degree ${\cardlef{\finseta}}/{3}$ element $F_{{\cardlef{\finseta}}/{3}}$ which depends only on the cardinality of $\finseta$.
The truncation $F_{\leq \largen}=(F_k)_{k\in \NN \suchthat k\leq \largen}$ until degree $N$ of any element $(F_k)_{k\in \NN}$ of $\CA_F$ can be viewed as such a cardinality-determined map.
Note that such a truncation of a group-like element is group-like in the sense of Lemma~\ref{lemprodcoprodbiz}.

Let $T$ be a trivial q-braid (represented by a constant path) except for the framing of one of its strands $K$, which is 
$lk(K,K_{\parallel})=1$ instead of $0$. Definition~\ref{defnotzfvariant} implies $\Zinvufrfneg(\drad{1} \times \left[0,1\right],T, .,(\omega(i,S^2))) = \exp(\alpha)\#_K\left[\emptyset\right]$.
In particular, this expression does not depend on $(\omega(i,S^2))_{i \in \underline{3\largen}}$, and $\Zinvufrfneg(\drad{1} \times \left[0,1\right],T, .,(\omega(i,S^2)))$ is group-like. Similarly, $\exp_{\leq \largen}\left(-\frac14 p_1(\tau)\ansothree(.)\right)_{\sqcup}$ is group-like. Therefore,
Lemma~\ref{lemprodcoprodbiz} allows us to conclude the proof of Proposition~\ref{propgrouplikefunc} for framed tangles with injective bottom and top configurations. Use Lemma~\ref{lemconvthmpoirbraid} to conclude for general $q$-tangles.
\eop

\begin{proposition}
\label{propwithZinvlink}
Let $(\omega(i,S^2))_{i \in \underline{3\largen}}$ be a family of volume-one forms of $S^2$. Let $\hcylc$ be a rational homology cylinder.
Let $\tanghcyll$ be a $J_{bb,tt}$-oriented $q$-tangle of $\hcylc$.
Let $\Zinvlinkufrfneg\left(\hcylc,\tanghcyll,.,(\omega(i,S^2))\right)$ 
\index[N]{ZZ@$\Zinvuf$ and some variants (see also the summary in the next pages)!ZuTT@$\Zinvlinkuf$}
denote the projection of
$\Zinvufrfneg\left(\hcylc,\tanghcyll,.,(\omega(i,S^2))\right)$ on $\Assis(\source)$.
Then we have \begin{equation*}\Zinvufrfneg\Bigl(\hcylc,\tanghcyll,.,\bigl(\omega(i,S^2)\bigr)_{i \in \underline{3\largen}}\Bigr)=\biggl(\Zinvlinkufrfneg\Bigl(\hcylc,\tanghcyll,.,\bigl(\omega(i,S^2)\bigr)_{i \in \underline{3\largen}}\Bigr)\Zinvufrf_{\leq N}\Bigl(\hcylc,\emptyset\Bigr)\biggr)_{\sqcup}.\end{equation*}
\end{proposition}
\bp
$\Zinvufrfneg\bigl(\hcylc,\emptyset, .,(\omega(i,S^2))_{i \in \underline{3\largen}}\bigr)=\Zinvuf(\rats(\hcylc))$ is independent of $(\omega(i,S^2))_{i \in \underline{3\largen}}$. See Theorems~\ref{thmconststraight} and~\ref{thmdefsanstauvarzinvf}.\eop

\section{A proof of universality}
\label{secproofthmbn}

In this section, we apply the properties of $\Zinvufrf$ to
study its variation under crossing changes and 
generalize Theorem~\ref{thmbn}.

Define a \emph{singular tangle representative with $n$ double points} to be an oriented $1$-dimensional 
manifold $\tanghcyll(\sourcetl)$ immersed in $\hcylc$ such that \begin{itemize}
\item the boundary of $\tanghcyll(\sourcetl)$  sits in the interior of $\drad{1} \times \{0,1\}$,
\item $\tanghcyll(\sourcetl)$ meets a neighborhood $N\left(\partial (\drad{1} \times \left[0,1\right]\right))$ as vertical segments, and
\item the only singular points of the immersion are $n$ double points \doublep, for which the directions of the two meeting branches generate a \emph{tangent plane}. 
\end{itemize}
Define a \emph{singular tangle with $n$ double points} to be an equivalence class of such representatives under the equivalence relation defined as in the nonsingular case in Definition~\ref{deftangleone}, by adding the adjective singular.

Extend the invariant $\Zinvuf$ to unframed singular tangles by the local rule
\begin{equation*}\Zinvuf\Bigl(\doublepcirc\Bigr)=\Zinvuf\Bigl(\pccirc\Bigr)-\Zinvuf\Bigl(\nccirc\Bigr).\end{equation*}
This local rule relates the invariants $\Zinvuf$ of three singular tangles that coincide outside
the represented ball and are as in the pictures in this ball.

Define the \indexT{chord diagram} $\Gamma_C(\tanghcyll)$ associated to an unframed singular tangle $\tanghcyll$ with $n$ double points to be the following diagram on the domain $\sourcetl$ of $\tanghcyll$.
The diagram $\Gamma_C(\tanghcyll)$ has $2n$ vertices. Its vertices are univalent and located at the preimages of the double points. It has $n$ edges, one between each pair of preimages of a double point.  These edges are called \emph{chords}. (The chords are attached on the left-hand side of the oriented domain $\sourcetl$, when orientations of univalent vertices are needed, as in Definition~\ref{defrkoruniv}.) 
\begin{notation}
\label{notzbar}
Respectively denote the images of $\Zinvuf$ and of $\projassis \circ \Zinvuf={\Zinvlinkuf}$ under the quotient by the $1T$-relation by $\Zinvufmodis$\index[N]{ZZ@$\Zinvuf$ and some variants (see also the summary in the next pages)!ZZZb@$\Zinvufmodis$ mod 1T}  and $\Zinvlinkufmodis$\index[N]{ZZ@$\Zinvuf$ and some variants (see also the summary in the next pages)!ZZZcb@$\Zinvlinkufmodis$}. (So $\Zinvufmodis(\tanghcyll)$ and  $\Zinvlinkufmodis(\tanghcyll)$ respectively belong to $\Aavis(\sourcetl)/(1T)$ and $ \Assis(\sourcetl)/(1T)$.)
\end{notation}

In this section, we prove the following theorem, which is a generalization of Theorem~\ref{thmbn} from knots to links and tangles.

\begin{theorem}
\label{thmunivsingtang} Let $n$ be a natural number.
For any unframed singular tangle $\tanghcyll$ with $n$ double points in a rational homology cylinder $\hcylc$, the expansion $\Zinvufmodis_{\leq n-1}(\tanghcyll)$ up to degree $n-1$ of $\Zinvufmodis(\tanghcyll)$ vanishes, and its expansion $\Zinvufmodis_{\leq n}(\tanghcyll)$ up to degree $n$ is equal to
\begin{equation*}\Zinvufmodis_{\leq n}(\tanghcyll)=\left[\Gamma_C(\tanghcyll)\right]\end{equation*}
in $\Aavis(\sourcetl)/(1T)$. So we have $\Zinvlinkufmodis_{\leq n}(\tanghcyll)=[\Gamma_C(\tanghcyll)]$
in 
$\Assis(\sourcetl)/(1T)$.
\end{theorem}

Theorem~\ref{thmbn} and its proof, whose easiest part is presented in Section~\ref{secchordd}, generalize to $k$-component oriented links, with numbered components, to produce the following corollary to Theorem~\ref{thmunivsingtang}. The isomorphism of the corollary was first shown by Dror Bar-Natan and Maxim Kontsevich \cite{barnatan}.

\begin{corollary}
With the notation of Section~\ref{secdeffintype}, for any $(k,n) \in \NN^2$,
$\Zinvlinkufmodis_{\leq n}$ induces an isomorphism from $\frac{\CF_n(\CK_k;\QQ)}{\CF_{n+1}(\CK_k;\QQ)}$ to $\Assis(\sqcup_{i=1}^k (S^1)_i)/(1T)$,
where $(S^1)_i$ is the copy of $S^1$ associated to the $i^{th}$ component of a link.
\end{corollary}
\eopwobp

In order to prove Theorem~\ref{thmunivsingtang}, we first define framed singular tangles and extend $\Zinvufrf$ to these tangles.
A \emph{parallelization} of a singular tangle is an isotopy class of parallels as in the nonsingular case, with the same restrictions near the boundary, where the parallel of a neighborhood of a double point is on one side of the tangent plane of the double point \doubleppetittexte. Recall that there are two ways of desingularizing \doubleppetittexte, the \emph{positive} one for which \doubleppetittexte is replaced by \pcpetittexte, and the \emph{negative} one for which \doubleppetittexte is replaced by \ncpetittexte.
In particular, every desingularization of such a singular tangle gets a natural parallelization from the parallelization of the singular tangle. Locally, the parallel of each branch is well-defined.

In general, we define the self-linking number of a component of a singular framed tangle as before, where the \emph{components} of a singular framed tangle are in natural one-to-one correspondence with the components of its domain. Let $\doublepointtang$ be a double point for which both branches belong to a component $K_j$. Let $\tanghcyll(\doublepointtang,+)$ and $\tanghcyll(\doublepointtang,-)$ respectively denote the positive and negative desingularizations of $\tanghcyll$ at $\doublepointtang$. Then the self-linking numbers of $K_j$ in $\tanghcyll$ and in these two desingularizations are related by
\begin{equation*}lk((K_j,K_{j\parallel}) \subset \tanghcyll(\doublepointtang,+))=lk((K_j,K_{j\parallel}) \subset \tanghcyll) +1\end{equation*}
and $lk((K_j,K_{j\parallel}) \subset \tanghcyll(\doublepointtang,-))=lk((K_j,K_{j\parallel}) \subset \tanghcyll) -1$.

We formally extend $\Zinvufrf$ to (framed) singular $q$-tangles by the formula
\begin{equation*}\Zinvufrfneg\Bigl(\doublepcirc\Bigr)=\Zinvufrfneg\Bigl(\pccirc\Bigr)-\Zinvufrfneg\Bigl(\nccirc\Bigr),\end{equation*}
where the parallels of the three tangles are supposed to be behind and to match on the boundary of the ball.

As an example, we have
\begin{equation*}\Zinvufrfneg\left(\fulltwistsing\right)=\Zinvufrfneg\left(\fulltwist\right)-\Zinvufrfneg\left(\fulltwistfake\right)=\Zinvufrfneg\left(\fulltwist\right)-\Zinvufrfneg\left(\twostraightstrands\right),\end{equation*}
where the endpoints of the tangles are assumed to lie on the real line.
So Lemmas~\ref{lemtwostraightstrands} and \ref{lemfulltwist} yield \begin{equation*}\Zinvufrf_{\leq 1}\left(\fulltwistsing\right)=\left[\onechordtwosor\right].\end{equation*} 

Note that $\Zinvufrf$ is now a functor on the category of singular $q$-tangles, which satisfies the cabling property and the duplication properties of Theorem~\ref{thmmainfunc} provided that the components involved in a double point are not duplicated.                                                                                                                                                                                                                                                                                                                                                                                                                                                      

\begin{proposition}
\label{propfundsingtang} Let $n$ be a natural number. For any singular $q$-tangle $\tanghcyll$ with $n$ double points, we have
\begin{equation*}\Zinvufrf_{\leq n}\left(\tanghcyll\right)=\left[\Gamma_C\left(\tanghcyll\right)\right].\end{equation*}
\end{proposition}
\bp In the proof below, we evaluate the lowest degree part that does not vanish in 
$\Zinvufrfneg\left(\tanghcyll\right)$ for various singular $q$-tangles. Note that this part
is unchanged when such a singular $q$-tangle is multiplied by a nonsingular $q$-tangle (except for the modification of the domains) since the lowest degree part that does not vanish for a nonsingular $q$-tangle is the class of the empty diagram.
In particular, the lowest degree nonvanishing part is independent of the bottom and top configurations of our $q$-tangles, which will not be specified.

Applying the cabling property of $\Zinvufrf$ to the following cable of a trivial braid with three strands yields
\begin{equation*}\Zinvufrf_{\leq 1}\left(\twostraightstrands \fulltwistsing\right)=\nochordtwosor \onechordtwosor.\end{equation*}
We thus get
\begin{equation*}\Zinvufrf_{\leq 1}\left(\doublecuptwistsing\right)=\doublecupnotwistsource\end{equation*}
by functoriality as desired for such a $q$-tangle.

Starting with a trivial braid and successively cabling some of its strands by replacing $n$ of them by $\doublecuptwistsing$, we find 
\begin{equation*}\Zinvufrf_{\leq n}\left(\doublecuptwistsing \dots \doublecuptwistsing \twostraightstrandsnonorlong \dots \twostraightstrandsnonorlong \right)=
\doublecupnotwistsource \dots \doublecupnotwistsource \twostraightstrandsnonorlong \dots \twostraightstrandsnonorlong.\end{equation*}
(Here, we first apply the cabling property to the first strand of a trivial braid $\id_{n+k}$ and regard the resulting tangle $T_1$ as a product $\id_{n+k-1}T_1 $, where $\id_{n+k-1}$ is a trivial braid. We next apply the cabling property to the first strand of $\id_{n+k-1}$. We keep going to find the obtained formula with the help of the functoriality property.)

Since every singular $q$-tangle with $n$ double points can be written as a product of a $q$-tangle as above and a nonsingular $q$-tangle by moving the double points below, the proposition follows.
\eop

We are ready to deduce Theorem~\ref{thmunivsingtang} from Proposition~\ref{propfundsingtang}.

\bpo{Proof of Theorem~\ref{thmunivsingtang}}
Assume that the singular $q$-tangle $\tanghcyll$ with $n$ double points has $k$ components $K_i$ for $i=1,\dots, k$.

The case $n=0$ is obvious.
Assume that $n=1$. Let $s_i$ denote the self-linking number
of $K_i$ in $\tanghcyll$. Let $\tanghcyll^+$ be the positive desingularization of $\tanghcyll$, and let $\tanghcyll^-$ be the negative desingularization of $\tanghcyll$.
Let $s_i^+$ (resp. $s_i^-$) denote the self-linking number
of the $i^{th}$ component in $\tanghcyll^+$ (resp. in $\tanghcyll^-$).
Recall that $\Zinvufmodis$ is the image of $\Zinvuf$ under the quotient by the $1T$-relation. Also recall \begin{equation*}\Zinvuf(\tanghcyll^+)=\prod_{j=1}^k\left(\exp(-s_j^+\alpha)\#_j\right)\Zinvufrfneg(\tanghcyll^+).\end{equation*}
If the two strands involved in the double point belong to two distinct components, then we have $s_j^+=s_j=s_j^-$ for any $j \in \underline{k}$ and \begin{equation*}\Zinvuf\bigl(\pc\bigr) - \Zinvuf\bigl(\nc\bigr)=\prod_{j=1}^k\bigl(\exp(-s_j\alpha)\#_j\bigr)\left(\Zinvufrfneg\bigl(\pc\bigr) - \Zinvufrfneg\bigl(\nc\bigr)\right).\end{equation*}
So the result follows since the lowest degree part of any $\exp(-s_j\alpha)$ is the class of the empty diagram.

If the two strands involved in the double point belong to the same component $K_i$, then we have $s_i^+=s_i^- +2$ and
\begin{multline*}
 \Zinvuf\bigl(\pc\bigr) - \Zinvuf\bigl(\nc\bigr)=\prod_{j=1}^k\left(\exp(-s_j^+\alpha)\#_j\right)\left(\Zinvufrfneg\bigl(\pc\bigr) - \Zinvufrfneg\bigl(\nc\bigr)\right)\\+ \left(\prod_{j=1}^k\left(\exp(-s_j^+\alpha)\#_j\right)-\prod_{j=1}^k\left(\exp(-s_j^-\alpha)\#_j\right)\right)\Zinvufrfneg\bigl(\nc\bigr),
\end{multline*}
where the $\bigl(\prod_{j=1}^k\left(\exp(-s_j^+\alpha)\#_j\right)-\prod_{j=1}^k\left(\exp(-s_j^-\alpha)\#_j\right)\bigr)$ \say{factor} begins with its degree one part, which is $-2\alpha_1=-\left[ \onechordsmalltseul \right]$. We get
\begin{equation*}\Zinvufmodis_{\leq 1}\bigl(\pc\bigr) - \Zinvufmodis_{\leq 1}\bigl(\nc\bigr)=\left[ \doublepdashchord \right]\end{equation*}
as desired.
Note that this equality would be wrong if $\Zinvufmodis$ were replaced by $\Zinvuf$, that is, without moding out by $1T$.

Let us now conclude the proof by induction on $n$. Assume that the result is known for singular $q$-tangles with less than $n$ double points.
Let $\doublepointset$ denote the set of double points of $\tanghcyll$. For $i \in \underline{k}$, let $\doublepointset_i$ denote the set of double points for which both branches belong to $K_i$.
For a subset $I$ of $\doublepointset$, the $q$-tangle obtained from $\tanghcyll$ by performing negative desingularizations on double points of $I$ and positive ones on double points of $\doublepointset \setminus I$ is denoted by $\tanghcyll_I$, and $s_{j,I}$ denotes the self-linking number of the component $K_{j,I}$ in $\tanghcyll_I$.

Then $\Zinvufmodis(\tanghcyll)$ is equal to 
\begin{multline*}\sum_{I \subseteq \doublepointset}(-1)^{\cardlef{I}}\Zinvufmodis(\tanghcyll_I)
=
\sum_{I \subseteq \doublepointset}(-1)^{\cardlef{I}}\prod_{j=1}^k\left(\exp(-s_{j,I}\alpha)\#_j\right)\Zinvufrfneg(\tanghcyll_I) \;\mbox{mod 1T}\\
=
\sum_{I \subseteq \doublepointset}(-1)^{\cardlef{I}}
\left(\prod_{j=1}^k\left(\exp(-s_{j,I}\alpha)\#_j\right)-\prod_{j=1}^k\left(\exp(-s_{j,\emptyset}\alpha)\#_j\right)\right)\Zinvufrfneg(\tanghcyll_I)
\\
+\prod_{j=1}^k\left(\exp(-s_{j,\emptyset}\alpha)\#_j\right) \sum_{I \subseteq \doublepointset}(-1)^{\cardlef{I}}\Zinvufrfneg(\tanghcyll_I)\;\mbox{mod 1T}.
\end{multline*}

The lowest degree term of the last line is $\left[\Gamma_C(\tanghcyll) \right]$. So it suffices to prove that the previous line does not contain terms of degree less than $n+1$ (mod 1T).
This previous line can be rewritten as
\begin{equation*}T_2=\sum_{I \subseteq \doublepointset}(-1)^{\cardlef{I}}\biggl(1-\prod_{j=1}^k\Bigl(\exp\bigl( (s_{j,I}-s_{j,\emptyset})\alpha\bigr)\#_j\Bigr) \biggr)\Zinvuf(\tanghcyll_I)\;\mbox{mod 1T},\end{equation*}
with $s_{j,I}-s_{j,\emptyset} =-2 \cardlef{I \cap \doublepointset_j}= \sum_{\doublepointtang \in I \cap \doublepointset_j}(-2)$.
So $T_2$ can be rewritten as 
\begin{equation*}T_2=\sum_{I \subseteq \doublepointset}(-1)^{\cardlef{I}}\biggl(1-\prod_{j=1}^k\Bigl( \prod_{\doublepointtang \in I \cap\doublepointset_j} \exp( -2 \alpha_{\doublepointtang})\#_j\Bigr) \biggr)\Zinvufmodis(\tanghcyll_I),\end{equation*} where $\alpha_{\doublepointtang}$ is a copy of $\alpha$.

Let $F$ be the set of maps $f\colon \cup_{j=1}^k \doublepointset_j \to \NN$ such that $f\left(\cup_{j=1}^k \doublepointset_j\right) \neq \{0\}$. Let $F(I)$ be the set of maps $f$ of $F$ such that
$f(\doublepointtang)=0$ for any $\doublepointtang \notin I$.
Then we have 
\begin{equation*}1-\prod_{j=1}^k\Bigl( \prod_{\doublepointtang \in I \cap\doublepointset_j} \exp( -2 \alpha_{\doublepointtang})\#_j\Bigr) = - \sum_{f \in F(I)}D_f\end{equation*}
with \begin{equation*}D_f=\prod_{j=1}^k\biggl(\Bigl( \prod_{\doublepointtang \in \doublepointset_j}\frac{(-2 \alpha)^{f(\doublepointtang)}}{f(\doublepointtang)!}\Bigr)\#_j\biggr).\end{equation*} 
Now, $T_2$ can be rewritten as 
\begin{equation*} \sum_{f \in F} D_f\sum_{I \suchthat I \subseteq \doublepointset, f\in F(I)}(-1)^{\cardlef{I}+1}\Zinvufmodis(\tanghcyll_I).\end{equation*}
Set $N(f)=\{\doublepointtang \suchthat \doublepointtang \in \cup_{j=1}^k \doublepointset_j, f(\doublepointtang) \neq 0\}$.
 Let $\tanghcyll(N(f),-)$ denote the singular $q$-tangle with $(n- \cardlef{N(f))})$ double points obtained from $\tanghcyll$ by desingularizing the double points of $N(f)$ in a negative way.
 Since the condition $f\in F(I)$ is equivalent to the condition 
$N(f) \subseteq I$, we get
\begin{equation*}\sum_{I \suchthat I \subseteq \doublepointset, f\in F(I)}(-1)^{\cardlef{I}+1}\Zinvufmodis(\tanghcyll_I)=(-1)^{\cardlef{N(f)}+1}\Zinvufmodis\bigl(\tanghcyll(N(f),-)\bigr).\end{equation*}
Note that $(-2 \alpha)^{f(\doublepointtang)}$ is of degree at least $3$  as soon as $f(\doublepointtang) \neq 0$ (when working modulo 1T since $\alpha_1=\alpha_2=0$).
So the degree of $D_f$ is at least $3\cardlef{N(f)}$.
By induction, the degree of $\Zinvufmodis(\tanghcyll(N(f),-))$ is $n-\cardlef{N(f)}$.
Therefore, the parts of $T_2$ of degree at most $n$ vanish. \eop

\part{Universality}

\chapter{The main universality statements and their corollaries}
\label{chapuniv} 

\section{Universality with respect to Lagrangian-preserving surgeries}
\label{secstateunivlag}

Let us recall some definitions quickly surveyed in the book introduction.

\begin{definition}
\label{defhhh}
 An \emph{integer (resp. rational) homology handlebody} of genus $g$ is a compact oriented $3$-manifold $A$ with the same integral (resp. rational) homology as the usual solid handlebody $\handlebody_g$ of Figure~\ref{fighandlebodyg}.
The {\em Lagrangian\/} $\CL_A$ of a compact $3$-manifold $A$ is the kernel of the map induced by the inclusion from $H_1(\partial A;\QQ)$ to $H_1(A;\QQ)$.
\end{definition}

\begin{exo}
 Show that if $A$ is a rational homology handlebody of genus $g$, then $\partial A$ is a connected genus $g$ surface. (See Appendix~\ref{sechomology}, where some basic properties of homology are recalled.)
\end{exo}

In Figure~\ref{fighandlebodyg}, the Lagrangian of $\handlebody_g$ is freely generated by the classes of the curves $a_i$.

\begin{definition}
\label{deflagsur}
An {\em integral (resp. rational) Lagrangian-Preserving (or LP) surgery\/} $(A^{\prime}/A)$ is the replacement of an integer (resp. rational) homology handlebody $A$, embedded in the interior of a $3$-manifold $M$, by another such $A^{\prime}$, whose boundary $\partial A^{\prime}$ is identified with $\partial A$
by an orientation-preserving diffeomorphism sending $\CL_A$ to $\CL_{A^{\prime}}$.
The manifold $M(A^{\prime}/A)$ obtained by such an LP-surgery is
\begin{equation*}M(A^{\prime}/A) = \bigl(M \setminus \Int(A)\bigr) \cup_{\partial A} A^{\prime}.\end{equation*}
(This only defines the topological structure of $M(A^{\prime}/A)$, but $M(A^{\prime}/A)$ is equipped with its unique smooth structure.)
\end{definition}

We present an interesting example of an integral Lagrangian-preserving surgery in Subsection~\ref{sublagboun}.
The Matveev Borromean surgery of Section~\ref{secfintypZ} is another example of integral LP-surgery.

\begin{lemma}
If $(A^{\prime}/A)$ is an integral (resp. rational) LP-surgery, then the homology of $M(A^{\prime}/A)$ with $\ZZ$-coefficients (resp. with $\QQ$-coefficients) is canonically isomorphic to
$H_{\ast}(M;\ZZ)$ (resp. to $H_{\ast}(M;\QQ)$). If $M$ is a $\QQ$-sphere, if $(A^{\prime}/A)$ is a rational LP-surgery, and if $(J,K)$ is a two-component link of $M \setminus A$, then the linking number of $J$ and $K$ in $M$ and the linking number of $J$ and $K$ in $M(A^{\prime}/A)$ coincide.
\end{lemma}
\bp Exercise. \eop

Let $(A^{\prime}/A)$ be a rational LP-surgery in a punctured rational homology sphere $\crats$. Let \begin{equation*}\partial_{MV} \colon H_2(A \cup_{\partial A}
-A^{\prime};\QQ) \rightarrow {\mathcal{L}}_A\end{equation*} be the morphism that maps the class of a closed surface in the closed 3--manifold $(A \cup_{\partial A}
-A^{\prime})$ to the boundary of its intersection with $A$.
The Mayer--Vietoris long exact sequence (see Theorem~\ref{thmMV}) shows that the above canonical morphism $\partial_{MV}$ is an isomorphism.
This isomorphism carries the algebraic triple intersection of surfaces to a trilinear antisymmetric form $\CI_{AA^{\prime}}$ on $\CL_A$. Explicitly, we have
\begin{equation*}\CI_{AA^{\prime}}(a_{i},a_{j},a_{k})=\bigl\langle \partial_{MV}^{-1}(a_i), \partial_{MV}^{-1}(a_j), \partial_{MV}^{-1}(a_k)\bigr\rangle_{\!A \cup
-A^{\prime}\,}.\end{equation*}

Let $(a_1,a_2,\dots,a_g)$ be a basis of $\CL_A$ and let $z_1,\dots,z_g$ be (curves representing) homology classes of $\partial A$ such that the system $(z_1,z_2,\dots,z_g)$ is dual to $(a_1,a_2,\dots,a_g)$ with respect to $\langle,\rangle_{\!\partial A\,}$:
\begin{equation*}\langle a_i,z_j \rangle_{\!\partial A\,}=\delta_{ij}=\left\{\begin{array}{ll} 1 & \mbox{if}\; i=j\\ 0 & \mbox{if}\; i\neq j.\end{array}\right.\end{equation*} 
Note that $(z_1,\dots,z_g)$
is a basis of $H_1(A;\QQ)$.

Represent $\CI_{AA^{\prime}}$
by the following combination $T\left(\CI_{AA^{\prime}}\right)$ of tripods 
whose three univalent vertices form an ordered set

\begin{equation*}T\left(\CI_{AA^{\prime}}\right)=\sum_{(i,j,k)  \in \underline{g}^3  \suchthat  i <j <k}\CI_{AA^{\prime}}(a_{i},a_{j},a_{k})
\smalltripod{${z}_i$}{${z}_j$}{${z}_k$},
\end{equation*}
where the tripods are considered up to the relations \begin{equation*}\smalltripod{$x$}{$y$}{$z$}=\smalltripod{$y$}{$z$}{$x$}=-\smalltripod{$y$}{$x$}{$z$}=-\smalltripod{$x$}{$-y$}{$z$}\;\;.\end{equation*}

\begin{notation}
\label{notcontract}
Let $G$ be a graph with $2k$ oriented trivalent vertices and with univalent vertices. Assume that the univalent vertices of $G$ are decorated with disjoint
curves of a punctured $\QQ$-sphere $\crats$.
Let $P(G)$ be the set of partitions of the set of univalent vertices of $G$ into disjoint pairs.

For $p \in P(G)$, identifying the two vertices of each pair provides a vertex-oriented trivalent Jacobi diagram $\Gamma_p$. Let $\ell(p)$ be the product, over the disjoint pairs of $p$, of the linking numbers of the curves corresponding to the two vertices in a pair.
We get an element $[\ell(p)\Gamma_p]$ of $\Aavis_k(\emptyset)$.
Define
\begin{equation*}\langle \langle G \rangle \rangle=\sum_{p \in P(G)} \bigl[\ell(p)\Gamma_p\bigr].\end{equation*} 
The contraction $\langle \langle . \rangle \rangle$ is linearly extended to linear combinations
of graphs. The disjoint union of combinations of graphs is bilinear.
\end{notation}

The universality theorem with respect to Lagrangian-preserving surgeries is the following one. It was proved in \cite{lessumgen} for the invariant $\Zinvuf$ of rational homology spheres. 
The statement below is more general since it applies to the invariant $\Zinvuf$ of Theorem~\ref{thmfstconsttang}, which satisfies the properties stated in Theorem~\ref{thmmainfunc}. Nevertheless, its proof reproduced in this book, is identical to the proof of the preprint \cite{lessumgen}, except for some editorial improvements.
(\cite{lessumgen} has never been submitted for publication.)

\begin{theorem} 
\label{thmmainunivlag} 
Let $\tanghcyll$ be a $q$-tangle representative in a rational homology cylinder $\hcylc$. Let $x$ be a positive integer. Let $\sqcup_{i=1}^{x}A^{(i)}$ be a disjoint union of rational homology handlebodies embedded in $\hcylc \setminus \tanghcyll$.
Let $(A^{(i)\prime}/A^{(i)})$ be rational LP-surgeries in $\hcylc$. For a subset $I$ of $\underline{x}$, let 
$\hcylc_I=\hcylc\left((A^{(i)\prime}/A^{(i)})_{i \in I}\right)$ be the rational homology cylinder obtained from $\hcylc$ by performing the LP-surgeries that replace $A^{(i)}$ by $A^{(i)\prime}$ for $i \in I$.
Set
$X=\left[\hcylc,\tanghcyll;(A^{(i)\prime}/A^{(i)})_{i \in\underline{x}}\right]$ and \begin{equation*}\Zinvuf_{n}(X)=  \sum_{I\subseteq \underline{x}}(-1)^{x+\cardlef{I}}\Zinvuf_n\left(\hcylc_I,\tanghcyll \right).\end{equation*}

If $2n <x $, then $\Zinvuf_{n}(X)$ vanishes. If $2n =x $, then we have
\begin{equation*}\Zinvuf_{n}(X)=\left[\Bigl\langle \Bigl\langle \bigsqcup_{i \in \underline{x}} T\bigl(\CI_{A^{(i)}A^{(i)\prime}}\bigr) \Bigr\rangle \Bigr\rangle \right].\end{equation*}
\end{theorem}

Before proving Theorem~\ref{thmmainunivlag}, we discuss some of its consequences and variants.
Section~\ref{secproofsurcas} shows that Theorem~\ref{thmmainunivlag} yields a direct proof of a surgery formula for the Theta invariant, as in \cite[Section 9]{lesformagt}.
The article \cite{lesformagt} presents many other surgery formulae implied by Theorem~\ref{thmmainunivlag}, which are not reproduced in this book.

Section~\ref{secfintypZ} shows how Theorem~\ref{thmmainunivlag} implies that $\Zinvuf$ restricts to a universal finite type invariant of integer homology $3$-spheres.
In Section~\ref{secfintypQ}, we review the Moussard classification of finite type invariants of rational homology $3$-spheres \cite{moussardAGT}, and we show how $\Zinvuf$ can be augmented to provide a universal finite type invariant of rational homology $3$-spheres, too, following an idea of Gw\'ena\"el Massuyeau.
Section~\ref{seccasson} shows how Theorem~\ref{thmmainunivlag} also implies that the invariant $\frac{1}{6}\Theta$ is the Casson--Walker invariant. Assuming this identification, the surgery formula of Section~\ref{secproofsurcas} is nothing but a consequence of the Casson--Walker surgery formula of \cite{wal}.
So a reader who does not need examples can skip Section~\ref{secproofsurcas}, at first.

We sketch the proof of Theorem~\ref{thmmainunivlag} in Section~\ref{secskproofthmmainunivlag}. We complete the details of the proof in the following two chapters.

Theorem~\ref{thmmainunivlag} and the universality theorems~\ref{thmbn} and~\ref{thmunivsingtang} for knots or tangles are put together in Theorem~\ref{thmunivmix}, which generalizes all of them.

Section~\ref{secmixuniv} shows how the main ingredients of the proof of Theorem~\ref{thmmainunivlag} also lead to Theorem~\ref{thmmainunivlagtwolegs}, which allows us to compute the degree $2$ part of $\Zinvlinkuf$, for any null-homologous knot, in Theorem~\ref{thmcompztwo}, with the help of the contents of Section~\ref{secproofsurcas}. 

Section~\ref{secDehnsur} below gives some background about Dehn surgeries. We use this background in Sections~\ref{secproofsurcas} and \ref{secfintypZ}.

\section{On Dehn surgeries}
\label{secDehnsur}

In this section, we define the manifold $\rats_{(K;p/q)}$ obtained by $p/q$-surgery on a $\QQ$-sphere $\rats$ along a knot $K$. We also introduce the lens spaces $L(p,q)$ and give some examples of surgeries on links used in Section~\ref{secfintypZ}.

Let $K$ be a knot in a $3$-manifold $M$, and let $N(K)$ be a tubular neighborhood of $K$.
The \emph{exterior} $E(K)$ of $K$ is the closure of the complement of $N(K)$ in $M$. \index[T]{exterior of a knot}
Let $\mu$ be a nonseparating simple closed curve of the boundary $\partial N(K)$ of $N(K)$.
The manifold obtained from $M$ by \emph{Dehn surgery}\index[T]{Dehn surgery} on $K$ with respect to $\mu$ is the union $E(K) \cup_{\partial  N(K)} T$ of $E(K)$ and a solid torus $T$, where $E(K)$ and $T$ are glued along $\partial N(K)$ by an orientation-reversing homeomorphism from $\partial T$ to $\partial N(K)$, which maps a meridian of $T$ to $\mu$. The result is then smoothed in a standard way. \index[T]{Dehn surgery} (Since the gluing of $T$ can be achieved by gluing a meridian disk of $T$ along $\mu$, thickening it, and gluing a $3$-dimensional ball to the resulting boundary, this surgery operation is well-defined.\footnote{The operation of gluing a $3$-ball $B^3$ along $S^2$ is well-defined, because any homeomorphism $f$ from $S^2$ to $S^2$ extends to $B^3$ as the homeomorphism that maps $tx$ to $tf(x)$ for $t \in \left[0,1\right]$ and $x\in S^2$.})

\begin{example}
 \label{exaStwotimesSone}
As the reader can check, the manifold obtained by Dehn surgery on the unknot $U$ of $S^3$ with respect to its meridian $m(U)$ is $S^2 \times S^1$.
\end{example}

Let $K$ be a knot in a rational homology sphere. If $K$ is null-homologous, then $K$ has a unique parallel $\ell(K)$ such that $lk(K,\ell(K))=0$. This parallel is called the \emph{preferred longitude}\index[T]{preferred longitude of a knot} of $K$. (A parallel of $K$ is also called a \indexT{longitude} of $K$.)
Let $\mu$ be a simple closed curve in the boundary $\partial N(K)$ of a tubular neighborhood of $K$ such that $\mu$ does not separate $\partial N(K)$.
 The class of the curve $\mu$ in $H_1(\partial N(K))$ may be expressed as $pm(K) +q\ell(K)$, where $m(K)$ is the meridian of $K$. The coefficient of the Dehn surgery along $K$ with respect to $\mu$ is ${p}/{q}$. We refer to this Dehn surgery as the $p/q$-surgery on $K$.
This coefficient ${p}/{q}$ may be expressed as ${lk(K,\mu)}/{\langle m(K),\mu\rangle_{\!\partial N(K)\,}}$. The $p/q$-surgery along a non-necessarily null-homologous knot $K$ in a rational homology $3$-sphere $\rats$ is the Dehn surgery with respect to a nonseparating simple closed curve $\mu$ of $\partial N(K)$ such that \begin{equation*}\frac{lk(K,\mu)}{\langle m(K),\mu\rangle_{\!\partial N(K)\,}}=\frac{p}{q}.\end{equation*}
Let $\rats_{(K;p/q)}$ \index[N]{Rmanif@$3$-manifolds!Rpq@$\rats_{(K;p/q)}$ Dehn surgery} denote the result of a $p/q$-Dehn surgery on $\rats$ along $K$.
As shown in Example~\ref{exaStwotimesSone}, we have $S^3_{(U;0)}=S^2 \times S^1$.

According to a theorem independently proved by Raymond Lickorish \cite{Lickorishsurg} and Andrew Wallace \cite{Wallace} in 1960, every closed oriented $3$-manifold can be obtained from $S^3$ by surgery along a link of $S^3$ whose components are equipped with integers. (Surgeries are  performed simultaneously along all the components of the link.) In \cite{rourke}, Colin Rourke gave a quick and elementary proof of this result.

\begin{examples}
 \label{exaStwotimesSoneconn}
As the reader can check, the manifold obtained by Dehn surgery on the trivial link of $S^3$ with $g$ components, all equipped with the coefficient $0$, is the connected sum of $g$ copies of $S^2 \times S^1$. Furthermore, this connected sum is homeomorphic to the manifold $\handlebody_g \cup_{\id_{\partial \handlebody_g}} (-\handlebody_g)$.

\bfig 
\centering
\begin{tikzpicture} [scale=0.15]
\draw (0,-1.5) circle (2.5);
\draw [white,line width=6pt,rotate=120] (0,-1.5) circle (2.5);
\draw[rotate=120] (0,-1.5) circle (2.5);
\draw [rotate=-120,white,line width=6pt] (-1.77,0.27) arc (135:190:2.5);
\draw [rotate=-120,white,line width=6pt] (1.77,0.27) arc (45:90:2.5);
\draw[rotate=-120] (0,-1.5) circle (2.5);
\draw [white,line width=6pt] (-1.77,0.27) arc (135:190:2.5);
\draw [white,line width=6pt] (1.77,0.27) arc (45:90:2.5);
\draw (-1.77,0.27) arc (135:190:2.5);
\draw (1.77,0.27) arc (45:90:2.5);
\end{tikzpicture}
\caption{Borromean link}\label{figborro}
\end{figure}

As a more challenging exercise, the reader can prove the following fact. The manifold obtained by Dehn surgery on the \emph{Borromean link} of $S^3$, represented in Figure~\ref{figborro}, whose components are equipped with the coefficient $0$, is diffeomorphic to $(S^1)^3$. A hint can be found in \cite[Example 13.1.5]{thurstonwngt}.
\end{examples}

Let $p$ and $q$ be two coprime integers, $p>0$. View $S^3$ as the unit sphere of $\CC^2$. The \emph{lens space} $L(p,q)$ is the quotient of $S^3$ by the action of $\ZZ/p\ZZ$ on $S^3$, where the generator $[1]$ of $\ZZ/p\ZZ$ acts on a unit vector $(z_1,z_2)$ of $\CC^2$ by mapping it to $\bigl(\exp({2i\pi}/{p})z_1,\exp({2i\pi q}/{p})z_2\bigr)$.

Let us study more surgeries along \emph{unknots} or \emph{trivial knots}, which are knots that bound an embedded disk, and prove the following well-known lemma.

\begin{lemma}
 \label{lemsurgunknot}
Let $k$ be an integer. Let $U$ be a trivial knot.
Then $S^3_{(U;1/k)}\cong S^3$. More generally, for any pair $(a,b)$ of coprime integers such that $a>0$, $S^3_{(U;a/(b+ka))}$ is diffeomorphic to $S^3_{(U;a/b)}$, and we have \begin{equation*}S^3_{(U;a/b)}=L(a,-b).\end{equation*}
If $U$ is a trivial knot of a $3$-manifold $M$, then $M_{(U;a/b)}$ is the connected sum $M \# L(a,-b)$.
\end{lemma}
\bp
The exterior $E$ of the unknot $U$ in $S^3$ is a solid torus whose meridian $m(E)$ is the preferred longitude $\ell(U)$ of $U$. The meridian $m(U)$ is a longitude $\ell(E)$ of $E$. Performing $({a}/{b})$-surgery along $U$ on $S^3$ amounts to gluing a solid torus with meridian 
\begin{equation*}\mu=a m(U) +b \ell(U)=b m(E) +a \ell(E)\end{equation*} to $E$, where $\langle m(E),\mu  \rangle_{\!\partial E\,}=a$ and $\langle \mu,\ell(E)  \rangle_{\!\partial E\,}=b$.
The manifold $S^3_{(U;a/b)}$ is the union of two solid tori $E$ and $T$ glued by a homeomorphism from $(-\partial T)$ to $\partial E$ mapping the meridian of $T$ to a curve $\mu$ as above.

Let $k$ be an integer. Then $\left(\ell(E) - k m(E)\right)$ is another longitude of $E$. This shows that for any coprime integers $a$ and $b$, the manifold $S^3_{(U;a/(b+ka))}$ is diffeomorphic to $S^3_{(U;a/b)}$.
In particular, we have $S^3_{(U;1/k)}\cong S^3$.

For a trivial knot $U$ in a $3$-manifold $M$, the manifold $M_{(U;a/b)}$ is the connected sum of $S^3_{(U;a/b)}$ and $M$. (The connected sum replaces a ball in the interior of the above solid torus $E$ by the exterior of a ball containing $U$ in $M$.)

Below, $S^3$ is viewed as the sphere of $\CC^2$ with radius $\sqrt{2}$.
The action of $\ZZ/p\ZZ$ on $S^3$ that defines the lens space $L(p,q)$ preserves the solid torus $|z_1| \leq |z_2|$ and the solid torus $|z_1| \geq |z_2|$. Let $E$ be the quotient of the second torus. We have \begin{equation*}E=\frac{\left\{\Bigl({\exp\bigl(\frac{2i\pi t}{p} \bigr)}{\sqrt{2 - |z_2|^2}}, z_2\Bigr) \suchthat  t\in \left[0,1\right], z_2 \in \CC, |z_2| \leq 1 \right\}}{ \Bigl({\sqrt{2 - |z_2|^2}}, z_2\Bigr) \sim \Bigl({\exp\bigl(\frac{2i\pi}{p} \bigr)}{\sqrt{2 - |z_2|^2}}, z_2\exp\bigl(\frac{2i\pi q}{p} \bigr)\Bigr)}.\end{equation*}
The meridian of the solid torus $E$ is 
\begin{equation*}m(E)=\Bigl\{\bigl(1, \exp(2i\pi u)\bigr) \suchthat  u\in \left[0,1\right]\Bigr\}\end{equation*}
and the possible (homology classes of) longitudes of $E$ are all the $\ell(E) + k m(E)$ for $k\in \ZZ$,
where \begin{equation*}\ell(E)=\left\{\Bigl(\exp\bigl(\frac{2i\pi t}{p} \bigr), \exp\bigl(\frac{2i\pi tq}{p} \bigr)\Bigr) \suchthat  t\in \left[0,1\right] \right\}.\end{equation*}
The boundary of $E$ is oriented as $(-S^1) \times S^1$.
The quotient of the solid torus $|z_1| \leq |z_2|$ is also a solid torus whose meridian is 
\begin{equation*}\begin{array}{ll}\mu&=\Bigl\{\bigl(\exp(2i\pi s), 1\bigr) \suchthat  s\in [0,1
]\Bigr\}\\&=\cup_{j=0}^{p-1} \left\{\left(\exp\bigl(\frac{2i\pi (t+j)}{p} \bigr),1\right) \suchthat  t\in \left[0,1\right]\right\}
   \\&=\cup_{j=0}^{p-1} \left\{\left(\exp\bigl(\frac{2i\pi t}{p} \bigr), \exp\left(-\frac{2i\pi qj}{p} \right)\right) \suchthat  t\in \left[0,1\right]\right\}.
  \end{array}
 \end{equation*}
So we have $\langle m(E),\mu  \rangle_{\!\partial E\,}=p$ and $\langle \mu,\ell(E)  \rangle_{\!\partial E\,}=-q$.
(There are $|q|$ pairs $(t,j)$ with $t\in \left[0,1\right[$ and $j \in \{0,1,\dots, p-1\}$ such that $(t+j)\in \frac{p}{q} \ZZ$.) 
\eop

\begin{remark} A homeomorphism from $M$ to $M_{(U;1/k)}$ can also be directly described as follows.
Let $D$ be a disk bounded by $U$, and let $d$ be a smaller disk inside $D$. The disk $D$ is parametrized by the disk of radius $2$ in $\CC$, and the unit disk parametrizes $d$.
The exterior $E$ of $U$ is homeomorphic to \begin{equation*}\Bigl(M \setminus \bigl(\mathring{D} \times \left]0,2\pi\right[\bigr)\Bigr) \cup_{d \times \{0,2\pi\}} \bigl(d \times \left[0,2\pi\right]\bigr)\end{equation*} by a homeomorphism mapping the meridian of $U$ to 
\begin{equation*}\bigl(\{1\} \times \left[0,2\pi\right]\bigr) \cup \bigl(\left[1,2\right] \times \{2\pi\}\bigr) \cup \bigl(-\{2\} \times \left[0,2\pi\right]\bigr) \cup \bigl(-\left[1,2\right] \times \{0\} \bigr).\end{equation*}
See Figure~\ref{figdoublecylinder}.

\bfig
\centering
\begin{tikzpicture} \useasboundingbox (-3,-.8) rectangle (3,1.8);
\draw (-.6,1.7) node{\scriptsize $\partial D \times \{2\pi\}$}
(0,-.25) node{\tiny $\partial d \times \{0\}$};
\begin{scope}[yshift=1cm]
\draw [thick,->] (-1.2,-1) -- (-1.2,0) (90:1.2 and .4) arc (90:450:1.2 and .4);
\draw (90:.3 and .1) arc (90:450:.3 and .1);
\end{scope}
\draw [thick,dashed] (0:1.2 and .4) arc (0:180:1.2 and .4);
\draw [thick,dashed] (180:1.2 and .4) arc (180:360:1.2 and .4);
\draw [dashed] (-.3,1) -- (-.3,0) (-90:.3 and .1) arc (-90:270:.3 and .1);
\draw [very thick,gray,dashed,->] (.3,.5) -- (.3,1) (1.2,0) -- (.3,0) -- (.3,.5);
\draw [very thick,gray,->]  (.5,.5) node{\scriptsize $m$} (1,.4) node{\scriptsize $m$} (1.2,.5) -- (1.2,0) (.3,1) -- (1.2,1)  -- (1.2,.5);
\end{tikzpicture}
\caption{The gray image $m$ of $m(U)$ in ${D} \times \left[0,2\pi\right]$}
\label{figdoublecylinder}
\end{figure}
The homeomorphism of $E$ restricting to $\bigl(M \setminus \bigl(\mathring{D} \times \left[0,2\pi\right]\bigr)\bigr)$ as the identity map and sending $(z,\theta) \in d \times \left[0,2\pi\right]$ to $(z\exp(i k \theta),\theta)$ sends the above meridian $m(U)$ to a curve homologous to $m(U) +k\ell(U)$. So this homeomorphism extends to provide a homeomorphism from $M$ to $M_{(U;\frac1k)}$.
\end{remark}

As another example, we prove the following standard lemma.
\begin{lemma}
\label{lemDehnmeridianzero}
Let $\merid$ be a meridian of a knot $K$ in a $3$-manifold $A$.
Equip $\merid$ with its preferred longitude $\ell(\merid)$.\footnote{It is equivalent to equip $\merid$ with the coefficient $0$ since $\merid$ lies in a ball.} Equip $K$ with a curve $\mu$ parallel to $K$.
Then the Dehn surgery on $((K;\mu),(\merid;\ell(\merid))$ does not change the $3$-manifold $A$.
\end{lemma}
\bp
Before and after the surgery, the two involved tori can be glued together along an annulus, whose core is a meridian of one of the knots and a longitude of the other, to form a solid torus in which $\ell(\merid)$ bounds a disk.
\eop

We refer the reader to the book \cite[Chapter 9, G, H]{Rolfsen} by Dale Rolfsen for many other examples of surgeries.

\section{Direct proof of a surgery formula for \texorpdfstring{$\Theta$}{Theta}}
\label{secproofsurcas}

In this section, we apply Theorem~\ref{thmmainunivlag} to compute $\Theta(\rats_{(K;p/q)})-\Theta(\rats) +\Theta(L(p,q))$ for any null-homologous knot $K$ in Proposition~\ref{propcaswsurfor}. In order to prove Proposition~\ref{propcaswsurfor}, we describe a special LP-surgery introduced in \cite[Section 9]{lesformagt}. We will also use this special LP-surgery to compute the degree $2$ part of $\Zinvlinkuf$ for a null-homologous knot in Theorem~\ref{thmcompztwo}.

Let $K$ be a null-homologous knot $K$ in a $3$-manifold $M$.
A \indexT{Seifert surface} of $K$ is a compact connected oriented surface $\Sigma$ embedded in $M$ such that the boundary $\partial \Sigma$ of $\Sigma$ is $K$.
A \indexT{symplectic basis} for the $H_1$ of such a Seifert surface is a basis $(x_1,y_1, \dots,x_g,y_g)$, as in Figure~\ref{figsympbasSeif}, where $\langle x_i,y_i\rangle_{\!\Sigma\,}=1$, for $i \in \underline{g}$.

\bfig
\centering
\begin{tikzpicture} \useasboundingbox (-1.1,-1.4) rectangle (7.5,1.3); On baisse tout d'1.6 et on recule de 1.1
\draw [-] (0:1.1) arc (0:180:1.1) (0:.7) arc (0:180:.7);
\draw [->,gray] (0:.9) arc (0:-120:.9) (0:.9) arc (0:240:.9); 
\draw (240:.7) node{\scriptsize $y_1$};
\draw (-.7,0) .. controls (-.7,-.1) ..  (-.3,-.1) .. controls (.1,-.1) .. (.1,0) (.5,0) .. controls (.5,-.1) ..  (.6,-.1) .. controls (.7,-.1) .. (.7,0) (1.1,0) .. controls (1.1,-.1) ..  (1.5,-.1) .. controls (1.9,-.1) .. (1.9,0);
\begin{scope}[xshift=1.2cm]
\draw [draw=white,double=black,very thick]  (0:1.1) arc (0:180:1.1) (0:.7) arc (0:180:.7);
\draw [draw=white,double=gray,very thick] (0:.9) arc (0:180:.9);
\draw [->,gray] (-60:.9) arc (-60:0:.9) (180:.9) arc (180:300:.9);
\draw (300:.75) node{\scriptsize $x_1$};
\end{scope}
\begin{scope}[xshift=5.2cm]
\draw [-] (0:1.1) arc (0:180:1.1) (0:.7) arc (0:180:.7);
\draw [->,gray] (0:.9) arc (0:-120:.9) (0:.9) arc (0:240:.9);
\draw (240:.7) node{\scriptsize $y_g$};
\draw (-.7,0) .. controls (-.7,-.1) ..  (-.3,-.1) .. controls (.1,-.1) .. (.1,0) (.5,0) .. controls (.5,-.1) ..  (.6,-.1) .. controls (.7,-.1) .. (.7,0) (1.1,0) .. controls (1.1,-.1) ..  (1.5,-.1) .. controls (1.9,-.1) .. (1.9,0);
\end{scope}
\begin{scope}[xshift=6.4cm]
\draw [draw=white,double=black,very thick]  (0:1.1) arc (0:180:1.1) (0:.7) arc (0:180:.7);
\draw [draw=white,double=gray,very thick] (0:.9) arc (0:180:.9);
\draw [->,gray] (-60:.9) arc (-60:0:.9) (180:.9) arc (180:300:.9);
\draw (300:.7) node{\scriptsize $x_g$};
\end{scope}
\draw [->] (2.3,0) .. controls (2.3,-.1) ..  (2.9,-.1) (3.5,-.1) .. controls (4.1,-.1) .. (4.1,0) (3.2,-1.1) -- (6.4,-1.1) .. controls (7,-1.1) and (7.5,-.6) .. (7.5,0)  (3.2,-1.15) node[above]{\scriptsize $K=\partial \Sigma$} (3.8,-.3) node{\scriptsize $\Sigma$} (-1.1,0) .. controls (-1.1,-.6) and (-.6,-1.1) .. (0,-1.1) -- (3.2,-1.1)  ;
\draw [dotted] (2.9,-.1) -- (3.5,-.1); 
\end{tikzpicture}
\caption{Symplectic basis of a Seifert surface}\label{figsympbasSeif}
\end{figure}

\begin{proposition}
\label{propcaswsurfor}
Let $K$ be a null-homologous knot in a rational homology sphere $\rats$. Let $\Sigma$ be a Seifert surface of $K$ in $\rats$, and let
$(x_1,y_1, \dots,x_g,y_g)$ be a symplectic basis of $\Sigma$. For a curve $c$ of $\Sigma$, let $c^+$ denote its push-off in the direction of the positive normal to $\Sigma$. Set
\begin{equation*}a_2(\Sigma) = \sum_{(i,j) \in \underline{g}^2} \left(lk(x_i,x_j^+)lk(y_i,y_j^+)-lk(x_i,y_j^+)lk(y_i,x_j^+)\right).\end{equation*}
Then we have
\begin{equation*} \Theta\bigl(\rats_{(K;p/q)}\bigr) =\Theta\bigl(\rats\bigr) - \Theta\bigl(L(p,q)\bigr) + 6 \frac{q}{p} a_2(\Sigma).\end{equation*}
\end{proposition}

We will prove this proposition in Subsection~\ref{subproofcaswalsurg} in these words. Thus, this proposition implies that $a_2(\Sigma)$ is an invariant of $K$. This invariant will be denoted by $a_2(K)$. It is equal to $\frac12{\Delta_K^{\prime\prime}(1)}$, where $\Delta_K$ is the Alexander polynomial of $K$ normalized so that $\Delta_K(t)=\Delta_K(t^{-1})$ and $\Delta_K(1)=1$. 

\begin{definition}
\label{defAlex}
 Here is a possible quick definition of the \emph{Alexander polynomial} $\Delta_K$ of the null-homologous knot $K$. Rewrite the symplectic basis $(x_i,y_i)_{i \in \underline{g}}$ as the basis $(z_j)_{j \in \underline{2g}}$ such that $z_{2i-1}=x_i$ and $z_{2i}=y_i$ for $i \in \underline{g}$. Let $V=\left[lk(z_j,z_k^+)\right]_{(j,k) \in  \underline{2g}^2}$ denote the associated Seifert matrix, and let $^tV$ denote its transpose, then we have \begin{equation*}\Delta_K(t)=\det\left(t^{1/2}V- t^{-1/2}\, ^tV\right).\end{equation*}
\end{definition}

See \cite{Alexanderpoly} or \cite[Chapter 2]{lespup} for other definitions of the Alexander polynomial, which will be mentioned later but not used in this book anymore.

\begin{remark} Proposition~\ref{propcaswsurfor} is also a consequence of the identification of $\Theta$ with $6\lambda_{CW}$ in Theorem~\ref{thmThetaeqlambda}, which is proved independently in Section~\ref{seccasson}, and of the Walker surgery formula proved in \cite[Theorem~5.1]{wal}.
 \end{remark}

\subsection{A Lagrangian-preserving surgery associated to a Seifert surface}
\label{sublagboun}

\begin{definition}
 Let $c(S^1)$ be a curve embedded in the interior of an oriented surface $F$, and let $c(S^1) \times [-1,1]$ be a collar neighborhood of $c(S^1)$ in $F$.
A \emph{right-handed (resp. left-handed) Dehn twist} about the curve $c(S^1)$ is a 
homeomorphism
 of $F$ that coincides with the identity
map of $F$ outside $c(S^1) \times [-1,1]$ and that maps \begin{equation*}\bigl(c(z),t\bigr) \in c(S^1) \times [-1,1]\mbox{ to }
\biggl(c\Bigl(z\exp\bigl(if(t)\bigr)\Bigr),t\biggr)\end{equation*} for $f(t)=\pi(t+1)$ (resp. for $f(t)=-\pi(t+1)$).
\end{definition}

Let $\Sigma$ be a Seifert surface of a knot $K$ in a manifold $M$. 
Consider an annular neighborhood $[-3,0] \times K$ of $(\{0\} \times K=)K =\partial \Sigma$ in $\Sigma$,
a small disk $D$ inside $\left]-2,-1\right[ \times K$, and a small open disk $d$ whose closure is in the interior of $D$.
Let $F=\Sigma \setminus d$. Let $h_{F}$ be the composition of the two left-handed Dehn twists on $F$
about $c=\partial D$ and $K_2=\{-2\} \times K$ with the right-handed one about $K_1=\{-1\} \times K$. See Figure~\ref{fig4}.

\bfig
\centering
\begin{tikzpicture} \useasboundingbox (-.5,-.3) rectangle (6.3,3.3);
\draw [thick,->] (5.9,1.5) .. controls (5.9,2.2) .. (4.2,2.2) ..  controls (3.2,2.2) and (2.5,3) .. (1.5,3) ..  controls (.7,3) and (0,2.2) .. (0,1.5) .. controls (0,.8) and (.7,0) .. (1.5,0) .. controls (2.5,0) and (3.2,.8) .. (4.2,.8) .. controls (5.9,.8) .. (5.9,1.5);
\draw [->] (1.5,.2) .. controls (2.4,.2) and (3.3,1) .. (4.2,1) .. controls (5.7,1) .. (5.7,1.5) .. controls (5.7,2) .. (4.2,2) ..  controls (3.3,2) and (2.4,2.8) .. (1.5,2.8) ..  controls (.8,2.8) and (.2,2) .. (.2,1.5) .. controls (.2,1) and (.8,.2) .. (1.5,.2);

\draw [gray,->] (3.6,1.5) .. controls (3.6,2) and (3,2.3) .. (2,2.3) .. controls (1,2.3) and (.4,2) .. (.4,1.5) ..  controls (.4,1) and (1,.7) .. (2,.7) ..  controls (3,.7) and (3.6,1) .. (3.6,1.5);
\draw (1,1.2) .. controls (1.08,1.25) .. (1.08,1.35) .. controls (1.08,1.45) and (1,1.5) .. (1,1.6) .. controls (1,1.75) and (1.1,1.9) .. (1.2,1.9) .. controls (1.3,1.9) and (1.4,1.75) .. (1.4,1.6) ..  controls (1.4,1.5) and (1.32,1.45) .. (1.32,1.35) .. controls (1.32,1.25) .. (1.4,1.2);
\draw (1.2,1.7) .. controls (1.25,1.6) .. (1.2,1.5); 
\draw (1.25,1.78) .. controls (1.23,1.76) .. (1.2,1.7) .. controls (1.15,1.6) .. (1.2,1.5) .. controls (1.23,1.44) .. (1.25,1.42);
\begin{scope}[xshift=1.6cm]
\draw (1,1.2) .. controls (1.08,1.25) .. (1.08,1.35) .. controls (1.08,1.45) and (1,1.5) .. (1,1.6) .. controls (1,1.75) and (1.1,1.9) .. (1.2,1.9) .. controls (1.3,1.9) and (1.4,1.75) .. (1.4,1.6) ..  controls (1.4,1.5) and (1.32,1.45) .. (1.32,1.35) .. controls (1.32,1.25) .. (1.4,1.2);
\draw (1.2,1.7) .. controls (1.25,1.6) .. (1.2,1.5); 
\draw (1.25,1.78) .. controls (1.23,1.76) .. (1.2,1.7) .. controls (1.15,1.6) .. (1.2,1.5) .. controls (1.23,1.44) .. (1.25,1.42);
\end{scope}
\fill [gray] (4.9,1.5) circle (.2);
\draw [<-,gray] (5.3,1.5) arc (0:-60:.4);
\draw [gray] (5.3,1.5) arc (0:300:.4);
\draw (3.9,1.5) node{\scriptsize $K_2$} (5.4,1.5) node{\scriptsize $c$} (6.15,1.5) node{\scriptsize $K$} (1.5,.4) node{\scriptsize $K_1$} (1.55,2.5) node{\scriptsize $\Sigma$} (4.9,1.5) node{\scriptsize $d$} (2,1.5) node{\scriptsize \dots};
\end{tikzpicture}
\caption{$K$, $\Sigma$, $F$, $c$, $K_1$, and $K_2$}\label{fig4}
\end{figure}

View $F$ as $F \times \{0\}$ in the boundary of a handlebody 
$A_{F}=F \times \left[-1,0\right]$ of $M$.
Extend $h_F$
to a homeomorphism $h_A$ of $\partial A_{F}$ that is the identity map
outside $F \times \{0\}$.
Let $A_{F}^{\prime}$ be a copy of $A_{F}$.
Identify $\partial A^{\prime}_F$ with $\partial A_{F}$ with 
\begin{equation*}h_{A} \colon \partial A^{\prime}_F \rightarrow \partial A_{F}.\end{equation*}

Define the {\em surgery associated to $\Sigma$\/} to be the surgery $(A_{F}^{\prime}/A_{F})$
associated to $(A_{F}, A_{F}^{\prime}; h_{A})$.
If $j$ denotes the embedding from $\partial A_{F}$ to $M$. This surgery
replaces 
\begin{equation*}M =A_{F} \cup_{j} \left( M \setminus \mathring{A}_{F} \right) \end{equation*}
by
\begin{equation*}M_F=A_{F}^{\prime} \cup_{j \circ h_{A}} \left(M \setminus \mathring{A}_{F} \right) .\end{equation*}

\begin{proposition}
\label{prophomafaprimef}
With the above notation, the surgery $(A_{F}^{\prime}/A_{F})$ associated to $\Sigma$ is a Lagrangian-preserving surgery with the following properties.
There is a homeomorphism from $M_F$ to $M$, 
\begin{itemize}
\item which extends the identity map of
\begin{equation*} M \setminus \bigl(\left[-3,0\right] \times K \times \left[-1,0\right] \bigr),\end{equation*}
\item which transforms a curve passing through $d \times \left[-1,0\right]$ by a band sum with $K$,
\item which transforms a $0$-framed meridian $m$ of $K$ passing through $d \times \left[-1,0\right]$, viewed as a curve of $M \setminus \mathring{A}_{F}$ (which may be expressed as $h_A^{-1}(m)$ in $A_{F}^{\prime}$),
to a $0$-framed copy of $K$ isotopic to the framed curve $h_A^{-1}(m)$ of Figure~\ref{fig5} (with the framing induced by $\partial A_F$).
\end{itemize}

\bfig
\centering
\begin{tikzpicture} 
\draw [thick,->] (4,.8) .. controls (5.5,.8) .. (5.8,1.5) .. controls (5.95,1.85) .. (5.95,2.02) .. controls (5.95,2.05) and (6.1,2.2) .. (4.6,2.2) ..  controls (3.6,2.2) and (2.7,3) .. (1.7,3) ..  controls (.9,3) and (0,2.2) .. (0,1.5) .. controls (0,.8) and (.2,0) .. (1,0) .. controls (2,0) and (3,.8) .. (4,.8);
\draw [->] (1,.2) .. controls (1.9,.2) and (3.2,1) .. (4.07,1) .. controls (5.38,1) .. (5.6,1.5) .. controls (5.82,2) .. (4.33,2) ..  controls (3.4,2) and (2.4,2.8) .. (1.7,2.8) ..  controls (1,2.8) and (.2,2) .. (.2,1.5) .. controls (.2,1) and (.3,.2) .. (1,.2);
\draw [gray,->] (3.6,1.5) .. controls (3.8,2) and (3.3,2.2) .. (2.3,2.2) .. controls (1.3,2.2) and (.6,2) .. (.4,1.5) ..  controls (.2,1) and (.7,.8) .. (1.7,.8) ..  controls (2.7,.8) and (3.4,1) .. (3.6,1.5);
\draw (1,1.2) .. controls (1.08,1.25) .. (1.08,1.35) .. controls (1.08,1.45) and (1,1.5) .. (1,1.6) .. controls (1,1.75) and (1.1,1.9) .. (1.2,1.9) .. controls (1.3,1.9) and (1.4,1.75) .. (1.4,1.6) ..  controls (1.4,1.5) and (1.32,1.45) .. (1.32,1.35) .. controls (1.32,1.25) .. (1.4,1.2);
\draw (1.2,1.7) .. controls (1.25,1.6) .. (1.2,1.5); 
\draw (1.25,1.78) .. controls (1.23,1.76) .. (1.2,1.7) .. controls (1.15,1.6) .. (1.2,1.5) .. controls (1.23,1.44) .. (1.25,1.42);
\begin{scope}[xshift=1.6cm]
\draw (1,1.2) .. controls (1.08,1.25) .. (1.08,1.35) .. controls (1.08,1.45) and (1,1.5) .. (1,1.6) .. controls (1,1.75) and (1.1,1.9) .. (1.2,1.9) .. controls (1.3,1.9) and (1.4,1.75) .. (1.4,1.6) ..  controls (1.4,1.5) and (1.32,1.45) .. (1.32,1.35) .. controls (1.32,1.25) .. (1.4,1.2);
\draw (1.2,1.7) .. controls (1.25,1.6) .. (1.2,1.5); 
\draw (1.25,1.78) .. controls (1.23,1.76) .. (1.2,1.7) .. controls (1.15,1.6) .. (1.2,1.5) .. controls (1.23,1.44) .. (1.25,1.42);
\end{scope}
\draw (4.9,1.5) circle (.2);
\draw [->,gray] (4.5,1.5) .. controls (4.42,1.3) and (4.65,1.15) .. (4.8,1.15) .. controls (4.95,1.15) and (5.22,1.3) .. (5.3,1.5) .. controls (5.38,1.7) and (5.15,1.85) .. (5,1.85) .. controls (4.85,1.85) and (4.58,1.7) .. (4.5,1.5);
\draw (3.9,1.5) node{\scriptsize $K_2$} (4.35,1.45) node{\scriptsize $c$} (4,.6) node{\scriptsize $K$} (1,.4) node{\scriptsize $K_1$} (1.75,2.5) node{\scriptsize $\Sigma$} (4.9,1.5) node{\scriptsize $d$} (2,1.5) node{\scriptsize \dots}  (5.6,.7) node{\scriptsize $m$};
\draw [dashed] (4.7,1.5) -- (4.7,.5) (5.1,1.5) -- (5.1,.5);
\draw [thick] (0,1.5) -- (0,.5) (5.95,2.02) -- (5.95,1.02);
\begin{scope}[yshift=-1cm]
\draw [dashed] (4.9,1.5) circle (.2);
\draw [thick] (0,1.5) .. controls (0,.8) and (.2,0) .. (1,0) .. controls (2,0) and (3,.8) .. (4,.8) .. controls (5.5,.8) .. (5.8,1.5) .. controls (5.95,1.85) .. (5.95,2.02);
\end{scope}
\draw [very thick,gray,->] (5.8,.8)  -- (5.8,.5) -- (5.1,.5) (5.1,1.5) -- (5.8,1.5) -- (5.8,.8);
\draw [very thick,dashed,gray]  (5.1,1.5) -- (5.1,.5);
\begin{scope}[xshift=3cm,yshift=-3.5cm]
\draw [thick,->] (4,.8) .. controls (5.5,.8) .. (5.8,1.5) .. controls (5.95,1.85) .. (5.95,2.02) .. controls (5.95,2.05) and (6.1,2.2) .. (4.6,2.2) ..  controls (3.6,2.2) and (2.7,3) .. (1.7,3) ..  controls (.9,3) and (0,2.2) .. (0,1.5) .. controls (0,.8) and (.2,0) .. (1,0) .. controls (2,0) and (3,.8) .. (4,.8);
\draw [->] (3.6,1.5) .. controls (3.8,2) and (3.3,2.2) .. (2.3,2.2) .. co
ntrols (1.3,2.2) and (.6,2) .. (.4,1.5) ..  controls (.2,1) and (.7,.8) .. (1.7,.8) ..  controls (2.7,.8) and (3.4,1) .. (3.6,1.5);
\draw (1,1.2) .. controls (1.08,1.25) .. (1.08,1.35) .. controls (1.08,1.45) and (1,1.5) .. (1,1.6) .. controls (1,1.75) and (1.1,1.9) .. (1.2,1.9) .. controls (1.3,1.9) and (1.4,1.75) .. (1.4,1.6) ..  controls (1.4,1.5) and (1.32,1.45) .. (1.32,1.35) .. controls (1.32,1.25) .. (1.4,1.2);
\draw (1.2,1.7) .. controls (1.25,1.6) .. (1.2,1.5); 
\draw (1.25,1.78) .. controls (1.23,1.76) .. (1.2,1.7) .. controls (1.15,1.6) .. (1.2,1.5) .. controls (1.23,1.44) .. (1.25,1.42);
\end{scope}
\begin{scope}[xshift=4.6cm,yshift=-3.5cm]
\draw (1,1.2) .. controls (1.08,1.25) .. (1.08,1.35) .. controls (1.08,1.45) and (1,1.5) .. (1,1.6) .. controls (1,1.75) and (1.1,1.9) .. (1.2,1.9) .. controls (1.3,1.9) and (1.4,1.75) .. (1.4,1.6) ..  controls (1.4,1.5) and (1.32,1.45) .. (1.32,1.35) .. controls (1.32,1.25) .. (1.4,1.2) (-1,-.6) node{\scriptsize $A_F$};
\draw (1.2,1.7) .. controls (1.25,1.6) .. (1.2,1.5); 
\draw (1.25,1.78) .. controls (1.23,1.76) .. (1.2,1.7) .. controls (1.15,1.6) .. (1.2,1.5) .. controls (1.23,1.44) .. (1.25,1.42);
\end{scope}
\begin{scope}[xshift=3cm,yshift=-3.5cm]
\draw (4.9,1.5) circle (.2);
\draw (3.9,1.5) node{\scriptsize $K_2$} (4,.6) node{\scriptsize $K$} (1,.5) node{\scriptsize $h_A^{-1}(m)$} (1.75,2.5) node{\scriptsize $\Sigma$} (4.9,1.5) node{\scriptsize $d$} (2,1.5) node{\scriptsize \dots};
\draw [dashed] (4.7,1.5) -- (4.7,.5) (5.1,1.5) -- (5.1,.5);
\draw [thick] (0,1.5) -- (0,.5) (5.95,2.02) -- (5.95,1.02);
\draw [very thick,gray,->] (5.8,.8)  -- (5.8,.5) -- (5.1,.5) (5.8,1.5) -- (5.8,.8);
\draw [very thick,dashed,gray] (5.1,1.5) -- (5.1,.5);
\draw [very thick,gray,->] (5.1,1.5) .. controls (5.4,1.5) and (5.3,1.15) .. (4.8,1.15) .. controls (4.65,1.15) and (4.42,1.3) .. (4.5,1.5) .. controls (4.58,1.7) and (4.85,1.85) .. (5,1.85) .. controls (5.3,1.85) and (5.3,2) .. (4.33,2) ..  controls (3.4,2) and (2.4,2.8) .. (1.7,2.8) ..  controls (1,2.8) and (.2,2) .. (.2,1.5) .. controls (.2,1) and (.3,.2) .. (1,.2);
\draw [very thick,gray] (1,.2) .. controls (1.9,.2) and (3.2,1) .. (4.07,1) .. controls (5.38,1) ..(5.56,1.45) .. controls (5.6,1.5) .. (5.8,1.5);
\end{scope}
\begin{scope}[xshift=3cm,yshift=-4.5cm]
\draw [dashed] (4.9,1.5) circle (.2);
\draw [thick] (0,1.5) .. controls (0,.8) and (.2,0) .. (1,0) .. controls (2,0) and (3,.8) .. (4,.8) .. controls (5.5,.8) .. (5.8,1.5) .. controls (5.95,1.85) .. (5.95,2.02);
\end{scope}
\end{tikzpicture}
\caption{$m$ and $h_A^{-1}(m)$}\label{fig5}
\end{figure}
\end{proposition}

\bp
Observe that $h_A\vert_{(F \times \{-1,0\} )\cup \left(\partial \Sigma \times \left[-1,0\right]\right)}$ extends to $\Sigma \times \left[-1,0\right]$ as 
\begin{equation*}\begin{array}{llll} h \colon &\Sigma \times \left[-1,0\right] & \rightarrow & \Sigma \times \left[-1,0\right]\\
&(\sigma,t) & \mapsto & h(\sigma,t)=\bigl(h_t(\sigma),t\bigr),
\end{array}\end{equation*}
where $h_{-1}$ is the identity map of $\Sigma$,
$h_0$ is the extension of $h_F$ by the identity map on $d$,\footnote{This extension is isotopic to the identity map of $\Sigma$.}
$h_t$ coincides with the identity map outside $[-5/2,-1/2]\times K(S^1)$,
and $h_t$ is defined as follows on $[-5/2,-1/2]\times K(S^1)$.\\
\noindent $\bullet$ When $t \leq -1/2$, $h_t$ coincides with the identity map $h_{-1}$ outside the disk $D$, whose elements are written as $D(z\in \CC)$, with $|z|\leq 1$. The elements of $d$ are the $D(z)$ for $|z|<1/2$. On $D$, $h_t$ describes the isotopy between the identity map and the left-handed Dehn twist about $\partial D$ located on $\{D(z) \suchthat 1/2 \leq |z| \leq 1\}$. We have
\begin{equation*}\begin{array}{lll}
h_t(z \in D) & = z \exp\bigl(i\pi(2t+2)4(|z|-1)\bigr) & \mbox{if} \; |z| \geq 1/2 \\
h_t(z \in D) & = z \exp\bigl(-2i\pi(2t+2)\bigr) & \mbox{if} \; |z| \leq 1/2.\\
\end{array}\end{equation*}
\noindent $\bullet$ When $t \geq -1/2$,  $h_t$ describes the following isotopy between $h_{-1/2}$ and the composition $h_0$ of $h_{-1/2}$ with the left-handed Dehn twist about $K_2$ and the right-handed Dehn twist about $K_1$, where the first twist is supported on $[-5/2,-2]\times K(S^1)$ and the second one is supported on $[-1,-1/2]\times K(S^1)$:
\begin{equation*}\begin{array}{ll}
h_t\bigl(u,K(z)\bigr) & = \Biggl(u,K\biggl(z \exp\Bigl(i(2t+1)\bigl(4\pi (u+5/2)\bigr)\Bigr)\biggr)\Biggr) \\& \hspace{6cm} 
\mbox{if} \; -5/2 \leq u \leq -2, \\
h_t\bigl(u,K(z)\bigr) & =  h_{-1/2}\Biggl(u,K\biggl(z \exp\Bigl(i(2t+1)(2\pi)\Bigr)\biggr)\Biggr)\\& \hspace{6cm} 
\mbox{if} \; -2 \leq u \leq -1, \\
h_t\bigl(u,K(z)\bigr) & = \Biggl(u,K\biggl(z \exp\Bigl(i(2t+1)\bigl(4\pi (-u-1/2)\bigr)\Bigr)\biggr)\Biggr) \\& \hspace{6cm}
\mbox{if} \;-1 \leq u \leq -1/2.
\end{array}\end{equation*}
Now, $M_F$ is naturally homeomorphic to 
\begin{equation*}\biggl( A^{\prime}_F\cup_{h\vert_{\partial A^{\prime}_F \setminus \left(\partial d \times \left[-1,0\right]\right)}} \Bigl(M \setminus \Int\bigl(\Sigma \times \left[-1,0\right]\bigr) \Bigr) \biggr)  \cup_{\partial \left(d \times \left[-1,0\right]\right)} \left(d \times \left[-1,0\right]\right) ,\end{equation*}
and hence to
\begin{equation*}\bigl(\Sigma \times \left[-1,0\right]\bigr) \cup_{h\vert_{\partial\left(\Sigma \times \left[-1,0\right]\right)}} \Bigl(M \setminus \Int\bigl(\Sigma \times \left[-1,0\right]\bigr) \Bigr),   \end{equation*}
which is mapped homeomorphically to $M$ by the identity map outside $\Sigma \times \left[-1,0\right]$, and by $h$ on $\Sigma \times \left[-1,0\right]$. Therefore, we indeed have a homeomorphism from $M_F$ to $M$. This homeomorphism is the identity map outside $[-3,0]\times K \times \left[-1,0\right]$.
It maps $d \times \left[-1,0\right]$ to a cylinder running along $K$ after being twisted negatively.

View the meridian $m$ as a curve of $M \setminus \mathring{A}_{F}$  with its framing induced by the boundary of $A_F$. Thicken $m$ as a band $\left[0,1\right] \times m$ in $\partial A_F$. Assume that a part of this band lies in the vertical boundary $\partial d \times \left[-1,0\right]$ of $d \times \left[-1,0\right]$ and can be written as a rectangle $[\theta,\theta^{\prime}] \times \left[-1,0\right] \subset \partial d \times \left[-1,0\right]$.
Then $h$ sends $[\theta,\theta^{\prime}] \times \left[-1,-1/2\right]$ to some other rectangle in $\partial d \times \left[-1,-1/2\right]$. The image under $h$ of $[\theta,\theta^{\prime}] \times [-1/2,0]$ together with a small additional piece of the thickened meridian can be isotoped in a tubular neighborhood of $K$ to a band on $\partial A_F$, which is first vertical in $\partial d \times [-1/2,0]$, and which then runs along $K$.
So $h$ sends the framed meridian $m$ to a curve isotopic to 
$h_A^{-1}(m)$ in a tubular neighborhood of $K$ with the framing induced by the boundary of $A_F$. See Figure~\ref{fig5}.

Now, $H_1(\partial A_{F})$ is generated by the generators of $H_1(\Sigma)\times \{0\}$, the generators of $H_1(\Sigma)\times \{-1\}$,
and the homology classes of $c=\partial D$ and $m$. Among them, the only generator that could be affected by $h_A$ is the class of $m$, which is not.
Thus $h_A$ acts trivially on $H_1(\partial A_{F})$, and the defined surgery is an LP-surgery.
\eop

Let $\Sigma \times \left[-1,2\right]$ be an extension of the previous neighborhood of $\Sigma$. Set $B_F=F \times \left[1,2\right]$. 
Define the homeomorphism $h_B$ of $\partial B_F$ to be the identity map anywhere except 
on $F \times \{1\}$, where it coincides with the homeomorphism $h_F$ of $F$, with
the obvious identification.

Let $B_{F}^{\prime}$ be a copy of $B_{F}$.
Identify $\partial B^{\prime}_{F}$ with $\partial B_{F}$ with 
\begin{equation*}h_B \colon \partial B^{\prime}_{F} \rightarrow \partial B_{F}.\end{equation*}

Define the {\em inverse surgery associated to $\Sigma$\/} to be the surgery 
associated to $(B_{F}, B_{F}^{\prime})$ (or $(B_{F}, B_{F}^{\prime}; h_B)$).
We can apply the previous study to this surgery by using the central symmetry of $\left[-1,2\right]$.

The following obvious lemma, which we will not prove, justifies the terminology.

\begin{lemma}
With the above notation, performing both surgeries $(B_{F}^{\prime}/B_{F})$ and $(A_{F}^{\prime}/A_{F})$ affects neither $M$ nor the curves in the complement
of $F \times \left[-1,2\right]$ (up to isotopy). Performing $(A_{F}^{\prime}/A_{F})$ (resp. $(B_{F}^{\prime}/B_{F})$) changes a $0$-framed meridian of $K$ passing through $d \times \left[-1,2\right]$
into a $0$-framed copy of $K$ (resp. $(-K)$).
\end{lemma}

\begin{lemma}
\label{lemtriplag}
Let $(x_i, y_i)_{i=1, \dots, g}$ be a symplectic basis of $\Sigma$. Then the tripod combination $T(\CI_{A_FA_F^{\prime}})$ associated
to the surgery
$(A_{F}^{\prime}/A_{F})$
is \begin{equation*}T(\CI_{A_FA_F^{\prime}})=\sum_{i=1}^g \bigtripod{$c=\partial D$}{$y_i$}{$x_i$}.\end{equation*}

For a curve $\gamma$ of $F$, let $\gamma^+$ denote $\gamma \times \{1\}$.
The tripod combination $T(\CI_{B_FB_F^{\prime}})$ associated
to the surgery
$(B_{F}^{\prime}/B_{F})$
is \begin{equation*}T(\CI_{B_FB_F^{\prime}})=-\sum_{i=1}^g 
\bigtripod{$c^+$}{$y_i^+$}{$x_i^+$}.\end{equation*}
\end{lemma}
\bp For a curve $\gamma$ of $F$, $\gamma^-$ denotes $\gamma \times \{-1\}$. In order to compute the intersection form of $(A_F \cup -A^{\prime}_F)$,
use the basis $\bigl(m,(x_i - x_i^-)_{i\in \underline{g}},(y_i-y_i^-)_{i\in \underline{g}}\bigr)$
of the Lagrangian of $A_F$. The system $\bigl(c,(y_i)_{i\in \underline{g}}, (- x_i)_{i\in \underline{g}}\bigr)$ is dual to this basis.
The only curve of the Lagrangian basis modified by $h_A$ is $m$, and $h_A^{-1}(m)$ may be expressed as a path composition $mc^{-1}K_2$. 
Let $D_m$ be a disk of $A_F$ bounded by $m$.
The isomorphism $\partial_{MV}^{-1}$ from $\CL_{A_F}$ to $H_2(A_F \cup -A^{\prime}_F)$ satisfies:
\begin{equation*}\begin{array}{ll}
\partial_{MV}^{-1}(x_i - x_i^-) &= S(x_i)=-\bigl(x_i \times \left[-1,0\right]\bigr) \cup \bigl(x_i \times \left[-1,0\right] \subset A^{\prime}_F\bigr),\\
\partial_{MV}^{-1}(y_i - y_i^-) &=S(y_i)=-\bigl(y_i \times \left[-1,0\right]\bigr) \cup \bigl(y_i \times \left[-1,0\right] \subset A^{\prime}_F\bigr), \\
\partial_{MV}^{-1}(m) &=S_A(m)=D_m -\Bigl(\bigl(\Sigma \setminus (\left]-2,0\right] \times K)\bigr)\cup (D_m \subset A^{\prime}_F) \Bigr),
\end{array}\end{equation*}
where the given expression of $\partial_{MV}^{-1}(m)$
must be completed in \begin{equation*}\partial A_F \cap \bigl([-2,0] \times K \times \left[-1,0\right]\bigr)\end{equation*} so that the boundary of $\partial_{MV}^{-1}(m)$ actually vanishes, as it does algebraically.

Since $x_i$ intersects only $y_i$ among the curves $x_j$ and $y_j$ for $j \in \underline{g}$, the surface $S(x_i)$ intersects only $S(y_i)$ and $S_A(m)$ in our basis of $H_2(A_F \cup -A^{\prime}_F)$. The algebraic intersection of
$S(x_i)$, $S(y_i)$, and $S_A(m)$ is $-1$.

For the surgery $(B_{F}^{\prime}/B_{F})$, use the reflection $T \colon \left[-1,2\right] \to \left[-1,2\right]$ such that $T(x)=1-x$. The induced reflection
$T_F=1_F \times T$ of $F \times \left[-1,2\right]$ maps $A_F$ onto $-B_F$. Use the image by $T_F$ of the above basis of $\CL_{A_F}$ for $\CL_{B_F}$. The system \begin{equation*}\left(-c^+,(-y_i)_{i\in \underline{g}}, (x_i^+)_{i\in \underline{g}}\right)\end{equation*} is dual to the obtained basis of $\CL_{B_F}$. (Since $T_F$ reverses the orientation, the intersection numbers on $\partial A_F$ are multiplied by $-1$.) Use the images under $T_F$ of the former surfaces. Their triple intersection numbers are the same since their positive normals and the ambient orientation are reversed.
\eop

\subsection{A direct proof of the Casson surgery formula}
\label{subproofcaswalsurg}

In this subsection, we prove Proposition~\ref{propcaswsurfor}, assuming Theorem~\ref{thmmainunivlag}, which will be proved independently.

Note the following easy, well-known lemma.
\begin{lemma}
\label{lemvarlk}
The variation of the linking number of two knots $J$ and $K$ in a rational homology $3$-sphere $\rats$ after a $p/q$-surgery
on a knot $V$ disjoint from $J \sqcup K$ in $\rats$ is given by the following formula.
\begin{equation*}lk_{\rats_{(V;p/q)}}(J,K)=lk_{\rats}(J,K)-\frac{q}{p}lk_{\rats}(V,J)lk_{\rats}(V,K).\end{equation*}
\end{lemma}
\bp The $p/q$-surgery
on $V$ is the surgery with respect to a curve $\mu_V \subset \partial N(V)$. Set $q_V=\langle m(V),\mu_V\rangle_{\!\partial N(V)\,}$ and $p_V=lk(V,\mu_V)$. We have $\frac{p}{q}=\frac{p_V}{q_V}$.
In $H_1(\rats \setminus (V \cup K);\QQ)$, we have \begin{equation*}J=lk_{\rats}(J,K)m(K) + lk_{\rats}(V,J)m(V) \;\mbox{and}\;
\mu_V=p_Vm(V) + q_Vlk_{\rats}(V,K)m(K).\end{equation*} Since $\mu_V$ vanishes in $H_1(\rats_{(V;p/q)} \setminus  K;\QQ)$, we get
\begin{equation*}J=lk_{\rats}(J,K)m(K) -\frac{q}{p}lk_{\rats}(V,K) lk_{\rats}(V,J)m(K)\end{equation*} in $H_1(\rats_{(V;p/q)} \setminus  K;\QQ)$.
\eop

\bpo{Proof of Proposition~\ref{propcaswsurfor} assuming Theorem~\ref{thmmainunivlag}}
Recall that $K$ bounds a Seifert surface $\Sigma$ in a rational homology sphere $\rats$.
Let $\Sigma \times \left[-1,2\right]$ be a collar of $\Sigma$ in $\rats$, and let $(A^{\prime}/A)=(A_{F}^{\prime}/A_{F})$ and $(B^{\prime}/B)=( B_{F}^{\prime}/B_{F})$ be the LP-surgeries of Subsection~\ref{sublagboun}. Let $U$ be a meridian of $K$
passing through $d \times \left[-1,2\right]$, such that performing one of the two surgeries
transforms $U$ into $\pm K$ and performing both or none of them leaves $U$ unchanged.
We have \begin{equation*}\begin{array}{ll}\Zinvuf_1\Bigl(\bigl[\rats_{(U;p/q)},\emptyset;A^{\prime}/A,B^{\prime}/B\bigr]\Bigr)&=2\Zinvuf_1\left(\rats_{(U;p/q)}\right)-2\Zinvuf_1\left(\rats_{(K;p/q)}\right)\\
&=\Bigl[\bigl\langle \bigl\langle T\left(\CI_{AA^{\prime}}\right) \sqcup T\left(\CI_{BB^{\prime}}\right) \bigr\rangle \bigr\rangle_{\!\!\rats_{(U;p/q)\,}}\Bigr].\end{array}\end{equation*}
According to Lemma~\ref{lemtriplag},
the tripods associated
to the surgery
$(A, A^{\prime})$ and to the surgery
$(B, B^{\prime})$,
are \begin{equation*}\sum_{i=1}^{g} \bigtripod{$c$}{$y_i$}{$x_i$}\mbox{ and }\sum_{j=1}^{g} \bigtripodrev{$c^{+}$}{$y_j^{+}$}{$x_j^{+}$},\end{equation*} respectively.
Among the curves of the tripods in the right-hand side, the only curve linking $c$ algebraically in $\rats_{(U;p/q)_{i\in N}}$
is $c^{+}$ with a linking number $-q/p$. Therefore, the vertices labeled by $c$ and $c^{+}$ must be paired together with coefficient $-q/p$. We get
\begin{equation*}\left\langle\left\langle \bigtripod{$c$}{$y_i$}{$x_i$} \bigtripodrev{$c^{+}$}{$y_j^{+}$}{$x_j^{+}$} \right\rangle\right\rangle= -\frac{q}{p}\left( lk(x_i,x_j^{+})lk(y_i,y_j^{+}) - lk(x_i,y_j^{+})lk(y_i,x_j^{+})\right) \Bigl[\tata\Bigr].\end{equation*}
So we have \begin{equation*}\Bigl[\bigl\langle \bigl\langle T\left(\CI_{AA^{\prime}}\right) \sqcup T\left(\CI_{BB^{\prime}}\right) \bigr\rangle \bigr\rangle_{\!\!\rats_{(U;p/q)\,}}\Bigr]= -\frac{q}{p} a_2(\Sigma) \Bigl[\tata\Bigr]\end{equation*} and
\begin{equation*}\Zinvuf_1\left(\rats_{(K;p/q)}\right)=\Zinvuf_1\left(\rats_{(U;p/q)}\right) +\frac{q}{2p} a_2(\Sigma) \Bigl[\tata \Bigr].\end{equation*}
Corollary~\ref{corthetazone} implies $\Zinvuf_1(\rats)=\frac{1}{12}\Theta(\rats)\left[\tata\right]$, where $\left[\tata\right]\neq 0$ in $\Aavis(\emptyset)$.
Also recall that $\Theta$ is additive under connected sum according to Corollary~\ref{corconnsumTheta}, and that Proposition~\ref{propThetaorrev} implies $\Theta(L(p,-q))=-\Theta(L(p,q))$.
The result follows, thanks to Lemma~\ref{lemsurgunknot}.
\eop

\section{Finite type invariants of \texorpdfstring{$\ZZ$-}{integer homology }spheres}
\label{secfintypZ}

In this section, we state the fundamental theorem of finite type invariants for integer homology $3$-spheres due to Thang L\^e \cite{le}, and we show how we may use Theorem~\ref{thmmainunivlag} in its proof. This shows in what sense Theorem~\ref{thmmainunivlag} implies that $\Zinvuf$ restricts to a \emph{universal} finite type invariant of integer homology $3$-spheres.
In order to do this, we first follow Mikhail Goussarov \cite{ggp} and Kazuo Habiro \cite{habiro} and construct surjective maps from
$\Aavis_n(\emptyset)$ to $\CF_{2n}(\CM)/\CF_{2n+1}(\CM)$.

\paragraph{Mapping \texorpdfstring{$\Aavis_n(\emptyset)$}{trivalent diagrams} to \texorpdfstring{$\CF_{2n}(\CM)/\CF_{2n+1}(\CM)$}{the graded space}.}

Let $\Gamma$ be a degree $n$ trivalent Jacobi diagram whose vertices are numbered in $\underline{2n}$.
Let $\Sigma(\Gamma)$ be an oriented surface containing $\Gamma$ in its interior such that $\Sigma(\Gamma)$ is a regular neighborhood of $\Gamma$ in $\Sigma(\Gamma)$. Equip $\Gamma$ with its vertex-orientation induced by the orientation of $\Sigma(\Gamma)$.
Embed $\Sigma(\Gamma)$ in a ball inside $\RR^3$.
Replace neighborhoods
\begin{equation*}
\Nedget\mbox{ of the edges }\edget\mbox{ by neighborhoods }\Nedgehopf \mbox{ of }\edgehopf\mbox{.}\end{equation*}
Thus, $\Sigma(\Gamma)$ is transformed into a collection of disjoint oriented surfaces \begin{equation*}\Sigma(Y) = \normyep\mbox{,}\end{equation*} one for each trivalent vertex. 
The graph \smally equipped with its framing induced by $\Sigma(Y)$ is called a \emph{$Y$-graph}. Its looped edges are called \emph{leaves}.
Thickening the $\Sigma(Y)$ transforms each of them into a standard genus $3$ handlebody $\aborb$.
The handlebody $\aborb$ has three handles with meridians $m_j$ and longitudes $\ell_j$, such that $\langle m_i, \ell_j\rangle_{\!\partial \aborb\,} =\delta_{ij}$. Its longitudes $\ell_j$ are on $\Sigma(Y)$ as in Figure~\ref{fighandlebodygtrois}.

\bfig
\centering
 \begin{tikzpicture} \useasboundingbox (-2.5,-2.5) rectangle (2.5,2.2);
\draw [very thick] (30:.4) -- (0,0) -- (0,-.4) (150:.4) -- (0,0) (30:1.1) circle (.7) (150:1.1) circle (.7)  (-90:1.1) circle (.7);
\draw [->] (30:.6) arc (-150:210:.5);
\draw [->] (150:.6) arc (-30:330:.5);
\draw [->] (-90:.6) arc (-270:90:.5);
\draw (30:.8) node{\tiny $\ell_1$} (150:.8) node{\tiny $\ell_2$} (-90:.8) node{\tiny $\ell_3$} (30:2.2) node{\tiny $m_1$} (149:2.25) node{\tiny $m_2$} (-90:2.2) node{\tiny $m_3$};
\draw [->,draw=white,double=black,very thick] (30:1.6) .. controls (29:1.6) and (26:1.7) .. (26:1.8) .. controls (26:1.9) and (29:2) .. (30:2);
\draw (26:1.8) .. controls (26:1.9) and (29:2) .. (30:2);
\draw [>-] (30:2) .. controls (31:2) and (34:1.95) .. (34:1.9);
\draw (30:1.6) .. controls (31:1.6) and (34:1.65) .. (34:1.7);
\draw [->,draw=white,double=black,very thick] (150:1.6) .. controls (149:1.6) and (146:1.7) .. (146:1.8) .. controls (146:1.9) and (149:2) .. (150:2);
\draw (146:1.8) .. controls (146:1.9) and (149:2) .. (150:2);
\draw [>-] (150:2) .. controls (151:2) and (154:1.95) .. (154:1.9);
\draw (150:1.6) .. controls (151:1.6) and (154:1.65) .. (154:1.7);
\draw [->,draw=white,double=black,very thick] (-90:1.6) .. controls (-91:1.6) and (-94:1.7) .. (-94:1.8) .. controls (-94:1.9) and (-91:2) .. (-90:2);
\draw (-94:1.8) .. controls (-94:1.9) and (-91:2) .. (-90:2);
\draw [>-] (-90:2) .. controls (-89:2) and (-86:1.95) .. (-86:1.9);
\draw (-90:1.6) .. controls (-89:1.6) and (-86:1.65) .. (-86:1.7);
\end{tikzpicture}
\caption{Meridians and longitudes of the genus $3$ handlebody $\aborb$}
\label{fighandlebodygtrois}
\end{figure}

The Matveev \emph{Borromean surgery} on $\aborb$ is the Dehn surgery on the six-component link $L_6$ inside $\aborb$ with respect to the parallels of its components that are parallel in Figure~\ref{figdelborro}. Sergei Matveev studied it in \cite{matveev}.
\bfig 
\centering
\begin{tikzpicture} [scale=0.15]
\begin{scope}
\newcommand{\feuille}[1]{
\draw[rotate=#1] (0,0) -- (0,-8);
\draw[rotate=#1] (0,-11) circle (3);}
\feuille{0}
\feuille{120}
\feuille{-120}
\draw (3,-4) node{$\Gamma$};
\end{scope}
\draw[very thick,->] (21.5,-3) -- (23.5,-3);
\begin{scope}[xshift=1200]
\newcommand{\bras}[1]{
\draw[rotate=#1] (0,-1.5) circle (2.5);
\draw [rotate=#1,white,line width=8pt] (-0.95,-4) -- (0.95,-4);
\draw[rotate=#1] {(0,-11) circle (3) (1,-3.9) -- (1,-7.6)};
\draw[rotate=#1,white,line width=6pt] (-1,-5) -- (-1,-8.7);
\draw[rotate=#1] {(-1,-3.9) -- (-1,-8.7) (-1,-8.7) arc (-180:0:1)};}
\bras{0}
\draw [white,line width=6pt,rotate=120] (0,-1.5) circle (2.5);
\bras{120}
\draw [rotate=-120,white,line width=6pt] (-1.77,0.27) arc (135:190:2.5);
\draw [rotate=-120,white,line width=6pt] (1.77,0.27) arc (45:90:2.5);
\bras{-120}
\draw [white,line width=6pt] (-1.77,0.27) arc (135:190:2.5);
\draw [white,line width=6pt] (1.77,0.27) arc (45:90:2.5);
\draw (-1.77,0.27) arc (135:190:2.5);
\draw (1.77,0.27) arc (45:90:2.5);
\draw (3.5,-4.5) node{$L_6$};
\end{scope}
\end{tikzpicture}
\caption{$Y$-graph and associated LP-surgery}\label{figdelborro}
\end{figure}

\begin{lemma}
\label{lemMatveevtorus}
The Matveev Borromean surgery changes the handlebody $\aborb$ to an integer homology handlebody $\aborbp$ with the same boundary and Lagrangian as $\aborb$.
The manifold $\aborb \cup_{\partial} (-\aborbp)$ is diffeomorphic to $(S^1)^3$, and there exists $\varepsilon_S=\pm 1$ such that
\begin{equation*}{T}(\CI_{\aborb\aborbp})=\varepsilon_S 
\begin{tikzpicture} \useasboundingbox (-.2,-.4) rectangle (1.1,0);
\begin{scope}[yshift=-.3cm]
\draw [thick] (.4,0) -- (0,0) -- (.4,.3) (0,0) -- (.4,-.3);
\draw (.4,.3) node[right]{\footnotesize ${\ell}_3$} (.4,0) node[right]{\footnotesize ${\ell}_2$} (.4,-.3) node[right]{\footnotesize ${\ell}_1$};
\fill (0,0) circle (1.5pt);
\end{scope}
\end{tikzpicture}.\end{equation*}
(The sign $\varepsilon_S$ is well-determined by the data. We do not need its explicit value.)
\end{lemma}
\bp First observe that the Lagrangian of $\aborbp$ is the same as the Lagrangian of $\aborb$.
As in Example~\ref{exaStwotimesSoneconn}, $\aborb \cup_{\partial} (-\aborb)$ is obtained from $S^3$ 
by surgery on three $0$-framed meridians of three handles of $\aborb$, where $\aborb$ is embedded in $S^3$ in a standard way. So $\aborbp \cup_{\partial} (-\aborb)$ is obtained by surgery on the zero-framed nine-component link obtained from the six-component link $L_6$ of Figure~\ref{figdelborro} by adding a meridian for each outermost component of $L_6$. Lemma~\ref{lemDehnmeridianzero} implies that $\aborbp \cup_{\partial} (-\aborb)$ is obtained by surgery on the zero-framed Borromean link. Therefore, according to Example~\ref{exaStwotimesSoneconn}, $\aborbp \cup_{\partial} (-\aborb)$ is diffeomorphic to $(S^1)^3$. So is $\aborb \cup_{\partial} (-\aborbp)$. Easy homological computations imply that $\aborbp$ is an integer homology handlebody.
\eop

Denote the Borromean surgeries associated to the $Y$-graphs corresponding to the vertices of $\Gamma$ by $(A^{(i)\prime}/A^{(i)})$.
With the notation of Section~\ref{secdeffintype},
define $\psi_n(\Gamma)$ to be the class of 
\begin{equation*}\bigl[S^3;(A^{(i)\prime}/A^{(i)})_{i\in \underline{2n}}\bigr]=\sum_{I \subseteq \underline{2n}} (-1)^{\cardlef{I}} S^3\bigl((A^{(i)\prime}/A^{(i)})_{i\in I} \bigr)\end{equation*}
in ${\CF_{2n}(\CM)}/{\CF_{2n+1}(\CM)}$. The coefficient field $\KK$ of Section~\ref{secdeffintype} for $\CF_{2n}(\CM)$ is $\RR$, from now on.\footnote{For the statements involving only invariants $\Zinvuf$ valued in spaces of Jacobi diagrams with rational coefficients (when no interval components are involved),  we can fix the coefficient field to be $\QQ$, provided that we also restrict the coefficient field of our related spaces of Jacobi diagrams to be $\QQ$.}

In \cite[Theorem 4.13, Section 4]{ggp}, Stavros Garoufalidis, Mikhail Goussarov, and Michael Polyak proved the following theorem.

\begin{theorem}[Garoufalidis, Goussarov, Polyak] \label{thmggp}
Let $n \in \NN$. For a degree $n$ trivalent Jacobi diagram $\Gamma$, the element $\psi_n(\Gamma)$ of ${\CF_{2n}(\CM)}/{\CF_{2n+1}(\CM)}$ constructed above depends only on the class of $\Gamma$ in $\Aavis_n(\emptyset)$, and the map
\begin{equation*}\psi_n \colon \Aavis_n(\emptyset)  \rightarrow \frac{\CF_{2n}(\CM)}{\CF_{2n+1}(\CM)} \end{equation*}
is surjective. Furthermore, we have $\frac{\CF_{2n+1}(\CM)}{\CF_{2n+2}(\CM)}=\{0\}$.
\end{theorem}

Assuming the above theorem,
the following {\em L\^e fundamental theorem on finite type invariants of $\ZZ$-spheres\/} becomes a corollary of Theorem~\ref{thmmainunivlag}.

\begin{theorem}[L\^e]
\label{thmle}
There exists a family $\left(\Zgen_n \colon \CF_0(\CM) \rightarrow \Aavis_n(\emptyset) \right)_{n \in \NN}$ of linear maps such that
\begin{itemize}
\item $\Zgen_n(\CF_{2n+1}(\CM))=0$,
\item the restriction $\overline{\Zgen}_n$ to ${\CF_{2n}(\CM)}/{\CF_{2n+1}(\CM)}$ of the morphism induced by $\Zgen_n$ on ${\CF_{0}(\CM)}/{\CF_{2n+1}(\CM)}$ to $\Aavis_n(\emptyset)$ is a left inverse of $\psi_n$.
\end{itemize}
In particular, for any $n \in \NN$, we have \begin{equation*}\frac{\CF_{2n}(\CM)}{\CF_{2n+1}(\CM)} \cong \Aavis_n(\emptyset)\mbox{ and }\frac{\CI_{2n}(\CM)}{\CI_{2n-1}(\CM)}\cong \Aavis^{\ast}_n(\emptyset).\end{equation*}
\end{theorem}
An invariant $\Zgen$ satisfying the properties in the statement of Theorem~\ref{thmle} above
is called a \emph{universal finite type invariant of $\ZZ$-spheres}. \index[T]{universal!finite type invariant}
In order to prove Theorem~\ref{thmle}, Thang L\^e proved that the L\^e--Murakami--Ohtsuki invariant $Z^{LMO}=(Z^{LMO}_n)_{n\in \NN}$ of \cite{lmo} is a universal finite type invariant of $\ZZ$-spheres in \cite{le}.

As a corollary of Theorem~\ref{thmmainunivlag}, we get the following Kuperberg--Thurston theorem \cite{kt}.
\begin{theorem}[Kuperberg, Thurston]
\label{thmkt}
 The restriction of $\Zinvuf$ to $\ZZ$-spheres is a universal finite type invariant of $\ZZ$-spheres.
\end{theorem}
\bp
Theorem~\ref{thmmainunivlag} ensures $\Zinvuf_n(\CF_{2n+1}(\CM))=0$. This reduces the proof of Theorem~\ref{thmkt} to the proof of the following lemma. \eop

\begin{lemma}
\label{lemkt}
 For any trivalent Jacobi diagram $\Gamma$, we have ${\Zinvuf}_n \circ \psi_n\bigl(\left[\Gamma\right]\bigr)=\left[\Gamma\right]$.
\end{lemma}
\bp Let us show how this lemma follows from Theorem~\ref{thmmainunivlag}. Number the vertices of $\Gamma$ in $\underline{2n}$. Call $(A^{(i)\prime}/A^{(i)})$ the Borromean surgery associated to the vertex $i$.
The associated tripod is
\begin{equation*}T\left(\CI_{A^{(i)}A^{(i)\prime}}\right)= \varepsilon_S \bigtripod{$\ell_1^{(i)}$}{$\ell_2^{(i)}$}{$\ell_3^{(i)}$}\end{equation*} for the fixed $\varepsilon_S =\pm 1$ of Lemma~\ref{lemMatveevtorus}.
Embed the tripods $T\left(\CI_{A^{(i)}A^{(i)\prime}}\right)$ into the graph $\Gamma$, naturally, so that the half-edge of $\ell_k^{(i)}$ is on the half-edge that gave rise to the leaf of $\ell_k^{(i)}$ in the $Y$-graph associated to $i$. 
In order to contribute to \begin{equation*}\Bigl\langle \Bigl\langle \bigsqcup_{i \in \underline{2n}} T\left(\CI_{A^{(i)}A^{(i)\prime}}\right) \Bigr\rangle \Bigr\rangle_{\!\!n\,},\end{equation*} a partition must pair a half-edge associated to a leaf of some $\ell_k^{(i)}$ with the half-edge of the only leaf that links $\ell_k^{(i)}$, which is the other half-edge of the same edge.
We get
\begin{equation*}\biggl[\Bigl\langle \Bigl\langle \bigsqcup_{i \in \underline{2n}} T\bigl(\CI_{A^{(i)}A^{(i)\prime}}\bigr) \Bigr\rangle \Bigr\rangle_{\!\!n\,}\biggr] = \left[\Gamma\right].\end{equation*}
\eop

\begin{remark}
In the original work of Thang L\^e \cite{le} and in the article \cite{ggp}, the primary filtration used for the space of $\ZZ$-spheres is defined from Borromean surgeries rather than from integral LP-surgeries. In \cite{aucles}, Emmanuel Auclair and I proved that the two filtrations coincide. We also proved that a universal finite type invariant of $\ZZ$-spheres automatically satisfies the more general formula of Theorem~\ref{thmmainunivlag} for any $X=\bigl[\crats;(A^{(i)\prime}/A^{(i)})_{i \in\underline{x}}\bigr]$ such that $\rats$ is a $\ZZ$-sphere and the $(A^{(i)\prime}/A^{(i)})$ are integral LP-surgeries in $\crats$.
In \cite{ggp}, Stavros Garoufalidis, Mikhail Goussarov, and Michael Polyak compare other filtrations of the space of $\ZZ$-spheres, including the original filtration of Tomotada Ohtsuki using surgeries on algebraically split links. This original Ohtsuki filtration defined in his introduction of finite invariants of $\ZZ$-spheres \cite{ohtkno} gives rise to the same notion of real-valued finite-type invariants.
\end{remark}

\section{Finite type invariants of \texorpdfstring{$\QQ$-}{rational homology }spheres}
\label{secfintypQ}

For a $\QQ$-sphere $\rats$, the cardinality of $H_1(\rats;\ZZ)$ is the product over the prime numbers $p$ of $p^{\nu_p(\rats)}$, where $\nu_p(\rats)$ is called the \emph{$p$-valuation of the order of $H_1(\rats;\ZZ)$}. In \cite[Proposition 1.9]{moussardAGT}, Delphine Moussard proved that $\nu_p$ is a degree $1$ invariant of $\QQ$-spheres with respect to $\CO^{\QQ}_{\CL}$, which is defined in Section~\ref{secdeffintype}. She also proved \cite[Corollary 1.10]{moussardAGT} that the degree $1$ invariants of $\QQ$-spheres with respect to $\CO^{\QQ}_{\CL}$ are (possibly infinite) linear combinations of the invariants $\nu_p$ and of a constant map.

Define an \emph{augmented trivalent Jacobi diagram} to be the disjoint union of a trivalent Jacobi diagram and a finite number of isolated $0$-valent vertices equipped with prime numbers. The \emph{degree}
of such a diagram is half the number of its vertices. It is a half-integer.
For a half-integer $h$, let $\aaug_h$ denote the quotient of the $\QQ$-vector space generated by
degree $h$ augmented trivalent Jacobi diagrams, by the Jacobi relation and the antisymmetry relation.
The product induced by the disjoint union turns $\aaug=\prod_{h\in \frac12 \NN}\aaug_h$ to a graded algebra.
In \cite{moussardAGT}, Delphine Moussard proved that
\begin{equation*}\frac{\CF_{n}(\CM_{\QQ})}{\CF_{n+1}(\CM_{\QQ})} \cong \aaug_{n/2}\end{equation*} for any integer $n$.
Her proof used
the configuration space integral $Z_{KKT}$, described in \cite{kt} and \cite{lesconst}, and the splitting formulae of \cite{lessumgen} stated in Theorem~\ref{thmmainunivlag}. See \cite[Theorem 1.7]{moussardAGT}. The invariant $Z_{KKT}$ is the restriction to $\QQ$-spheres of the invariant $\Zinvuf$ described in this book.

The maps $\psi_n$ of Section~\ref{secfintypZ} can be generalized to canonical maps 
\begin{equation*}\psi_h \colon \aaug_{h} \to \frac{\CF_{2h}(\CM_{\QQ})}{\CF_{2h+1}(\CM_{\QQ})}\end{equation*} as follows.
For any prime number $p$, let $B_p$ be a rational homology ball such that $|H_1(B_p;\ZZ)|=p$.
Let $\Gamma^a$ be the disjoint union of a degree $k$ trivalent Jacobi diagram $\Gamma$ and $r$ isolated $0$-valent vertices $v_j$ equipped with prime numbers $p_j$ for $j \in \underline{r}$.
Embed $\Gamma^a$ in $\RR^3$. Thicken it, replace $\Gamma$ by $2k$ genus $3$ handlebodies $A^{(i)}$ associated to the vertices of $\Gamma$ as in Section~\ref{secfintypZ}, and replace each vertex $v_j$ by a small ball $B(v_j)$ around it so that the $B(v_j)$ and the $A^{(i)}$ form a family of $2k+r$ disjoint rational homology handlebodies.
Define $\psi_{k+r/2}(\Gamma^a)$ to be the class of
\begin{equation*}\Bigl[S^3;\bigl(A^{(i)\prime}/A^{(i)}\bigr)_{i\in \underline{2k}}, \bigl(B_{p_j}/B(v_j)\bigr)_{j\in \underline{r}}\Bigr]\mbox{ in }\frac{\CF_{2k+r}(\CM_{\QQ})}{\CF_{2k+r+1}(\CM_{\QQ})},\end{equation*} where the $\bigl(A^{(i)\prime}/A^{(i)}\bigr)$ are the Borromean surgeries associated to the $Y$-graphs corresponding to the vertices of $\Gamma$ as in Section~\ref{secfintypZ}.
\begin{lemma}
The map $\psi_h$ is well-defined.
\end{lemma}
\bp According to \cite[Lemma 6.11]{moussardAGT}, if $B^{\prime}_{p_j}$ is a $\QQ$-ball whose $H_1(.;\ZZ)$ has the same cardinality as $H_1(B_{p_j};\ZZ)$, then
\begin{equation*}\bigl(S^3(B^{\prime}_{p_j}/B^3)-S^3(B_{p_j}/B^3)\bigr)\end{equation*} belongs to $\CF_{2}(\CM_{\QQ})$.
This guarantees that $\psi_{r/2+k}(\Gamma^a)$ does not depend on the chosen balls $B_{p_j}$. Thus, Theorem~\ref{thmggp}
implies that $\psi_h$ is well-defined.
\eop

This map $\psi_h$ is canonical.
The generators of ${\CF_{2h}(\CM_{\QQ})}/{\CF_{2h+1}(\CM_{\QQ})}$ exhibited in \cite[Section 6.2 and Proposition 6.9]{moussardAGT} are in the image of $\psi_h$. So $\psi_h$ is surjective.

Let $\aaugc_h$ denote the subspace of $\aaug_h$ generated by connected degree $h$ diagrams. So, if $\aaugc_h \neq 0$, then $h \in \NN$ or $h=1/2$. Set $\aaugc =\prod_{h \in \frac12\NN}\aaugc_n$.
Let $\zaug$ denote the $\aaugc$-valued invariant $\zaug$ of $\QQ$-spheres such that, for any $\QQ$-sphere $\rats$, we have
\begin{itemize}
 \item $\zaug_0(\rats)=0$,
\item $\zaug_{1/2}(\rats)=\sum_{p\; \mbox{\scriptsize prime}} \nu_p(\rats) \bullet_p$, and
\item $\zaug_n(\rats)$ is the natural projection $\zaug_n(\rats,\emptyset)=p^c(\Zinvuf_n(\rats,\emptyset))$ of $\Zinvuf_n(\rats,\emptyset)$ to the subspace $\Aavisc_n(\emptyset)$ of $\Aavis_n(\emptyset)$ generated by connected diagrams. (Recall from Notation~\ref{notationzZ} that the projection $p^c$ maps disconnected diagrams to $0$.)
\end{itemize}

Define an $\aaug$-valued invariant $\bzaug=(\bzaug_n)_{n \in \frac12\NN}$ to be
$\bzaug =\exp(\zaug)$ for the $\aaugc$-valued invariant $\zaug$ (meaning $\bzaug(\rats) =\exp(\zaug(\rats))$ for any $\QQ$-sphere $\rats$).

As noticed by Gw\'ena\"el Massuyeau, the Moussard fundamental theorem for {\em finite type invariants of $\QQ$-spheres\/} can be stated as follows.

\begin{theorem}[Moussard]
\label{thmmoussmass}
The family $\bigl(\bzaug_h \colon \CF_0(\CM_{\QQ}) \rightarrow \aaug_{h} \bigr)_{h \in \frac12 \NN }$ of linear maps is such that, for any $h \in \frac12\NN$,
\begin{itemize}
\item we have $\bzaug_h\bigl(\CF_{2h+1}(\CM_{\QQ})\bigr)=0$, and
\item $\bzaug_h$ induces a left inverse to $\psi_h$ 
from $\frac{\CF_{2h}(\CM_{\QQ})}{\CF_{2h+1}(\CM_{\QQ})}$ to $\aaug_{h}$.
\end{itemize}
In particular, we have $\frac{\CF_{n}(\CM_{\QQ})}{\CF_{n+1}(\CM_{\QQ})} \cong \aaug_{n/2}$ and $\frac{\CI_{n}(\CM_{\QQ})}{\CI_{n-1}(\CM_{\QQ})}\cong (\aaug_{n/2})^{\ast}$ for any $n \in \NN$.
\end{theorem}
\bp The first assertion follows from Theorem~\ref{thmmainunivlag} and Lemma~\ref{lemmultinvfin}. Let $\Gamma^a$ be the disjoint union of a degree $k$ trivalent Jacobi diagram $\Gamma$ and $r$ isolated $0$-valent vertices $v_j$ equipped with prime numbers $p_j$ for $j \in \underline{r}$.
Let $\psi_{\Gamma}$
be a representative of $\psi_{k+r/2}(\Gamma^a)$ in $\CF_{2k+r}(\CM_{\QQ})$.
Let us prove 
\begin{equation*}\bzaug_{\leq k+r/2}(\psi_{\Gamma})=[\Gamma^a] 
\end{equation*}
Write $\psi_{\Gamma}$ as $\bigl[S^3;(A^{(i)\prime}/A^{(i)})_{i\in \underline{2k}}, (B_{p_j}/B(v_j))_{j\in \underline{r}}\bigr]$.
If $r=0$, then $\psi_{\Gamma}$ belongs to $\CF_{2k+r}(\CM)$ and $\bzaug(\psi_{\Gamma}) $ is equal to $\Zinvuf(\psi_{\Gamma})$. So we have $\bzaug_k(\psi_{\Gamma})=\left[\Gamma\right]$, thanks to Lemma~\ref{lemkt}. (All the involved manifolds are $\ZZ$-spheres.)
The general case follows by induction. If $r>0$, let $\Gamma^{\prime}$ be obtained from $\Gamma^a$ by forgetting the vertex $v_r$. Let $\psi_{\Gamma^{\prime}}$ be obtained from $\psi_{\Gamma}$ by forgetting the surgery $(B_{p_r}/B(v_r))$. This surgery is nothing but a connected sum with $S_{p_r}=B_{p_r}\cup_{S^2} B^3$.
Since $\bzaug$ is multiplicative under connected sum according to Theorem~\ref{thmconnsum}, we have
\begin{equation*}\bzaug(\psi_{\Gamma})=\bzaug(\psi_{\Gamma^{\prime}})\bigl(\bzaug(S_{p_r})-1\bigr).\end{equation*}
Then, identifying the nonvanishing terms with minimal degree yields
\begin{equation*}\bzaug_{\leq k+r/2}(\psi_{\Gamma})=\bzaug_{k+(r-1)/2}(\psi_{\Gamma^{\prime}})[\bullet_{p_r}],\end{equation*}
which allows us to conclude the proof. 
Indeed, any $\lambda \in \bigl({\CF_{n}(\CM_{\QQ})}/{\CF_{n+1}(\CM_{\QQ})}\bigr)^{\ast}$ extends to the linear form
$\lambda \circ \psi_{n/2} \circ \bzaug_{n/2}$ of $\CI_{n}(\CM_{\QQ})$. So the natural injection \begin{equation*}\frac{\CI_{n}(\CM_{\QQ})}{\CI_{n-1}(\CM_{\QQ})}\hookrightarrow \left(\frac{\CF_{n}(\CM_{\QQ})}{\CF_{n+1}(\CM_{\QQ})}\right)^{\ast}\end{equation*} is surjective.
\eop

According to Corollary~\ref{corthetazone}, we have
\begin{equation*}\Zinvuf_1(\rats,\emptyset)=\frac{1}{12}\Theta(\rats)\left[\tata\right].\end{equation*}
In particular, the invariant $\Theta$ is of degree at most $2$ with respect to $\CO_{\CL}^{\QQ}$ according to Theorem~\ref{thmmainunivlag}. Furthermore, Lemma~\ref{lemkt}
implies \begin{equation*}\Theta\Bigl(\psi_1\bigl(\tata\bigr)\Bigr)=12.\end{equation*}

The following easy corollary of Theorem~\ref{thmmoussmass} can be proved as Corollary~\ref{coruniv}.

\begin{corollary}
\label{corunivdegtwo}
For any real-valued invariant $\nu$ of $\QQ$-spheres of degree at most $2$ with respect to $\CO_{\CL}^{\QQ}$, there exist real numbers $a_{\theta}$, $a_0$, $a_p$ for any prime number $p$, and $a_{p,q}$ for any pair $(p,q)$ of prime numbers such that $p\leq q$, such that
\begin{equation*}\nu(\rats)=a_0 + \sum_{p\; \mbox{\scriptsize prime}}a_p \nu_p(\rats) + \sum_{p,q \;\mbox{\scriptsize prime} \suchthat  p\leq q}a_{p,q}\nu_p(\rats)\nu_q(\rats) + a_{\theta}\Theta(\rats)\end{equation*}
for any $\QQ$-sphere $\rats$.
\end{corollary}

Note that the above infinite sums of the statement do not cause problems since they are finite when applied to a $\QQ$-sphere $\rats$.

According to Proposition~\ref{propThetaorrev} (or to Theorem~\ref{thmorchangman}),  we have
$\Theta(-\rats)=-\Theta(\rats)$ for any $\QQ$-sphere $\rats$.

\begin{theorem}
\label{thmcharTheta}
Let $\nu$ be a real-valued invariant of $\QQ$-spheres such that
\begin{itemize}
\item the invariant $\nu$ is of degree at most $2$ with respect to $\CO_{\CL}^{\QQ}$ and
\item we have $\nu(-\rats)=-\nu(\rats)$ for any $\QQ$-sphere $\rats$.
\end{itemize}
Then there exists a real number $a_{\theta}$ such that 
$\nu=a_{\theta}\Theta$.
\end{theorem}
\bp
Apply Corollary~\ref{corunivdegtwo}. We have \begin{equation*}(\nu-a_{\theta}\Theta)(-\rats)=(\nu-a_{\theta}\Theta)(\rats)=-(\nu-a_{\theta}\Theta)(\rats)\end{equation*} for any $\QQ$-sphere $\rats$. We get $\nu=a_{\theta}\Theta$.
\eop

\begin{remark}
A similar result was proved in \cite[Proposition 6.2]{lessumgen} without using the Moussard theorem.
\end{remark}

\section{Identifying \texorpdfstring{$\Theta$}{Theta} with the Casson--Walker invariant}
\label{seccasson}

In 1984, Andrew Casson introduced an invariant of $\ZZ$-spheres, which counts the conjugacy classes of irreducible $SU(2)$-representations of their fundamental groups using Heegaard splittings. See \cite{akmc,gm,mar}.
This invariant lifts the Rohlin $\mu$-invariant of Definition~\ref{defmuRohlin} from ${\ZZ}/{2\ZZ}$ to $\ZZ$.
In 1988, Kevin Walker generalized the Casson invariant to $\QQ$-spheres in \cite{wal}.
Here, the Casson--Walker invariant $\lambda_{CW}$ is normalized as in \cite{akmc,gm,mar} for integer homology $3$-spheres, and as $\frac{1}{2}\lambda_W$ for rational homology $3$-spheres, where $\lambda_W$ is the Walker normalisation in \cite{wal}.
\cite[Lemma 3.1]{wal} implies
\begin{equation*}\lambda_{CW}(-\rats)=-\lambda_{CW}(\rats)\end{equation*} for any $\QQ$-sphere $\rats$.
\cite{lesinv} shows that the Casson--Walker generalization satisfies the same splitting formulae as $\frac{1}{6}\Theta$. So $\lambda_{CW}$ is of degree at most $2$ with respect to $\CO_{\CL}^{\QQ}$, and we have \begin{equation*}\lambda_{CW}\Bigl(\psi_1\bigl(\tata\bigr)\Bigr)=2.\end{equation*}
(This is a consequence of \cite[Theorem 1.3]{lesinv}.)

As a direct corollary of Theorem~\ref{thmcharTheta}, we obtain the following theorem first proved by Greg Kuperberg and Dylan Thurston in \cite{kt} for $\ZZ$-spheres, and generalized to $\QQ$-spheres in \cite[Section 6]{lessumgen}. See \cite[Theorem 2.6]{lessumgen}.

\begin{theorem}
\label{thmThetaeqlambda}
 We have $\Theta=6\lambda_{CW}$.
\end{theorem}
\bp Recall that Lemma~\ref{lemkt} and Corollary~\ref{corthetazone} imply $\Theta(\psi_1(\tata))=12.$
\eop

\section{Sketch of the proof of \texorpdfstring{Theorem~\ref{thmmainunivlag}}{the splitting formulae}}
\label{secskproofthmmainunivlag}

Fix $\sqcup_{i=1}^{x}A^{(i)}$ and a representative of $\tanghcyll$ in $\hcylc$ as in the statement of Theorem~\ref{thmmainunivlag}. We have $\sqcup_{i=1}^{x}A^{(i)}\subset \hcylc \setminus \tanghcyll$.

For $I\subseteq \underline{x}$, recall $\hcylc_I=\hcylc\left((A^{(i)\prime}/A^{(i)})_{i \in I}\right)$. Set $\crats_I=\crats\left(\hcylc_I \right)$, $\rats_I=\rats\left(\hcylc_I\right)$,  $\crats=\crats_{\emptyset}=\crats(\hcylc)$, and $\rats=\rats_{\emptyset}=\rats(\hcylc)$.

For any part $X$ of $\rats_I$, $C_2(X)$ denotes the preimage of $X^2$ under the blowdown map from $C_2(\rats_I)$ to $\rats_I^2$.

In order to prove Theorem~\ref{thmmainunivlag}, we will compute the $\Zinvuf_n(\hcylc_I,\tanghcyll)$ with antisymmetric homogeneous propagating forms $\omega_I$ on the $C_2(\rats_I)$ such that the $\omega_I$ coincide with each other as much as possible (with respect to Definition~\ref{defantisympropform} of antisymmetric propagating forms). More precisely, for any subsets $I$ and $J$ of $\underline{x}$, our forms will satisfy \begin{equation*}\omega_I=\omega_J \;\; \mbox{on}\;\; C_2\Bigl(\bigl(\rats \setminus \cup_{i \in I \cup J} \Int(A^{(i)})\bigr) \cup_{i \in I \cap J} A^{(i)\prime}\Bigr).\end{equation*}

When dealing with integral LP-surgeries, such forms will be associated with parallelizations $\tau_I$ of the $\hcylc_I$ (as in Definition~\ref{defparacyl}), which coincide as much as possible with each other, i.e., such that 
\begin{equation*}\tau_I=\tau_J\mbox{ on }\Bigl(\bigl(\rats \setminus \cup_{i \in I \cup J} \Int(A^{(i)})\bigr) \cup_{i \in I \cap J} A^{(i)\prime}\Bigr)
\times \RR^3.\end{equation*}

Unfortunately, such a consistent choice of parallelizations is not always possible for rational LP-surgeries. See Section~\ref{secwhypseudop} and Example~\ref{exabadlag}.
To remedy this problem, we will make the definition of $\Zinvuf$ more flexible by allowing more general propagating forms associated with generalizations of parallelizations, called \emph{pseudo-parallelizations}.

We define these pseudo-parallelizations in Chapter~\ref{chappseudopar}, and we show that they satisfy the following properties.
\begin{itemize}
\item They generalize parallelizations. They are genuine parallelizations outside a link tubular neighborhood, inside which they can be thought of as an average of genuine parallelizations.
\item A parallelization defined near the boundary of a rational homology handlebody always extends to this rational homology handlebody as a pseudo-parallelization. (See Lemma~\ref{lempseudoparextend}.)
\item When $\crats$ is an asymptotic rational homology $\RR^3$, a pseudo-parallelization $\tilde{\tau}$ of $\crats$ induces a homotopy class of special complex trivializations $\tilde{\tau}_{\CC}$ of $T\crats \otimes_{\RR} \CC$, which has a Pontrjagin number $p_1(\tilde{\tau}_{\CC})$. Outside the link tubular neighborhood considered above, $\tilde{\tau}_{\CC}$ is $\tilde{\tau}\otimes_{\RR} 1_{\CC}$. We set $p_1(\tilde{\tau})=p_1(\tilde{\tau}_{\CC})$. (See Definitions~\ref{defpseudotrivpone} and \ref{defponepseudotrivpone}.)
\item The notion of a \emph{homogeneous propagating form} of $(C_2(\rats),\tilde{\tau})$ is presented in Definition~\ref{defpropagatorpseudohom}. 
This definition allows us to extend the definition of $\Zinvuf$ of Theorem~\ref{thmfstconsttang} using pseudo-parallelizations instead of parallelizations as follows.
For any long tangle representative \begin{equation*}\tanghcyll \colon \sourcetl \hookrightarrow \crats(\hcylc)\end{equation*} in a rational homology cylinder equipped with a pseudo-parallelization $\tilde{\tau}$ restricting to a neighborhood of the image of $\tanghcyll$ as a genuine parallelization,
for any $n \in \NN$, and for any family $(\omega(i))_{i\in \underline{3n}}$ of homogeneous propagating forms of $(C_2(\rats(\hcylc)),\tilde{\tau})$, the sum
\begin{equation*}\Zinv_n\Bigl(\hcylc,\tanghcyll,\bigl(\omega(i)\bigr)\Bigr)=\sum_{\Gamma \in \Davis^e_n(\sourcetl)}\coefgambet_{\Gamma}I\Bigl(\hcylc,\tanghcyll,\Gamma,\bigl(\omega(i)\bigr)_{i \in \underline{3n}}\Bigr)\left[\Gamma\right] \in \Aavis_n(\sourcetl)\end{equation*}
depends only on $(\hcylc,\tanghcyll,p_1(\tilde{\tau}_{\CC}))$ and on the $I_{\theta}(K_j,\tilde{\tau})$, which are defined as in Lemma~\ref{lemdefItheta} and Definition~\ref{defIthetalong}, for the components $K_j$, $j \in \underline{k}$, of $\tanghcyll$. It is denoted by $\Zinv_n(\hcylc,\tanghcyll,\tilde{\tau})$.
Set \begin{equation*}\Zinv(\hcylc,\tanghcyll,\tilde{\tau})=\bigl(\Zinv_n(\hcylc,\tanghcyll,\tilde{\tau})\bigr)_{n\in \NN} \in \Aavis(\sourcetl).\end{equation*}
Then Theorem~\ref{thmfstconsttangpseudo} ensures that we have \begin{equation*}\Zinvuf(\hcylc,\tanghcyll)=\exp\Bigl(-\frac14 p_1(\tilde{\tau}_{\CC})\ansothree\Bigr)\prod_{j=1}^k\Bigl(\exp\bigl(-I_{\theta}(K_j,\tilde{\tau})\alpha\bigr)\#_j\Bigr) \Zinv(\hcylc,\tanghcyll,\tilde{\tau}).\end{equation*}
\end{itemize}

In the general case, we begin by choosing pseudo-parallelizations $\tau_I$ of the $\hcylc_I$ such that 
$\tau_I$ and $\tau_J$ coincide on \begin{equation*}\Bigl(\bigl(\rats \setminus \cup_{i \in I \cup J} \Int(A^{(i)})\bigr) \cup_{i \in I \cap J} A^{(i)\prime}\Bigr)
\times \RR^3.\end{equation*}

A reader only interested in the cases for which pseudo-parallelizations are not necessary, as in the applications of Sections~\ref{secproofsurcas} and \ref{secfintypZ}, can skip Chapter~\ref{chappseudopar}, and substitute the word pseudo-parallelization with parallelization in the rest of the proof below and in Chapter~\ref{chapsimnormprop}.

Set $\tau_{\emptyset}=\tau$. The following lemma relates the $p_1(\tau_I)$.

\begin{lemma}
\label{lempofi}
Set $p(i)=p_1(\tau_{\{i\}})-p_1(\tau)$.
 For any subset $I$ of $\underline{x}$, we have
\begin{equation*}p_1(\tau_I)=p_1(\tau) + \sum_{i \in I}p(i).\end{equation*}
\end{lemma}
\bp Proceed by induction on the cardinality of $I$. The lemma is obviously true if $\cardlef{I}$ is zero or one.
Assume that $\cardlef{I} \geq 2$. Let $j \in I$. It suffices to prove $p_1(\tau_I) - p_1(\tau_{I \setminus \{j\}})=p_1(\tau_{\{j\}})-p_1(\tau)$.
This follows by applying twice the second part of Proposition~\ref{proppont}, where ${\manifm}_0=A^{(j)}$, ${\manifm}_1=A^{(j)\prime}$, and $D=\hcylc$ or $D=\hcylc_{I \setminus \{j\}}$. The first application identifies $\left(p_1(\tau_{\{j\}})-p_1(\tau)\right)$ with $p_1(\tau\vert_{A^{(j)}},\tau_{\{j\}|A^{(j)\prime}})$. The second one yields the conclusion.
\eop

For any $i \in \underline{x}$, fix pairwise disjoint simple closed curves $(a^i_j)_{j=1,\dots,g_i}$ and pairwise disjoint simple closed curves  $(z^i_j)_{j=1,\dots,g_i}$ on $\partial A^{(i)}$ such that
\begin{equation*}\CL_{A^{(i)}}=\oplus_{j=1}^{g_i} [a^i_j]\end{equation*} 
and 
\begin{equation*}\bigl\langle a^i_j,z^i_k \bigr\rangle_{\!\!\partial A^{(i)}\,}=\delta_{jk}=\left\{\begin{array}{ll} 0 &  \mbox{if} \; j \neq k\\
1 & \mbox{if} \; j = k.\end{array}\right.\end{equation*}

Let $[-4,4] \times \left(\sqcup_{i \in N} \partial A^{(i)}\right)$ be a tubular neighborhood of $(\sqcup_{i \in N} \partial A^{(i)})$ 
in $\hcylc$. This neighborhood intersects $A^{(i)}$ as $[-4,0] \times \partial A^{(i)}$.
Let $[-4,0] \times \partial A^{(i)}$ be a neighborhood of $\partial A^{(i)\prime}=\partial A^{(i)}$ in $A^{(i)\prime}$.
The manifold $\hcylc_{\{i\}}=\hcylc_i$ is obtained from $\hcylc$ by removing $\left(A^{(i)} \setminus
\left(\left]-4,0\right] \times \partial A^{(i)}\right)\right)$ and by gluing back $A^{(i)\prime}$ along $\left]-4,0\right] \times \partial A^{(i)}$.

Let $\eta_{[-1,1]}$ \index[N]{etamin@$\eta_{[-1,1]}$} be a one-form with compact support in $\left]-1,1\right[$ such that  \begin{equation*}\int_{[-1,1]}\eta_{[-1,1]}=1. \end{equation*}
Let $a^i_j  \times [-1,1]$ be a tubular neighborhood of $a^i_j$ in $\partial A^{(i)}$. 
Let $\eta(a^i_j)$ \index[N]{etazaij@$\eta(a^i_j)$ one-form on $A_I^{(i)}$} be a closed one-form on $A^{(i)}$ such that
the support of $\eta(a^i_j)$ intersects $[-4,0] \times \partial A^{(i)}$ inside  $[-4,0] \times \left(a^i_j  \times [-1,1]\right)$, where $\eta(a^i_j)$ can be written as
\begin{equation*}\eta(a^i_j)=\projp_{[-1,1]}^{\ast}(\eta_{[-1,1]}),\end{equation*} 
with the projection $\projp_{[-1,1]}\colon [-4,0] \times \left(a^i_j  \times [-1,1]\right) \to [-1,1]$ to the $[-1,1]$ factor. Let $\eta(a^i_j)$ also denote a closed one-form on $A^{(i)\prime}$ which can be written in the same way on $[-4,0] \times \partial A^{(i)}$.
Note that the forms $\eta(a^i_j)$ on $A^{(i)}$ and $A^{(i)\prime}$ induce a closed one-form 
on $(A^{(i)} \cup_{\partial A^{(i)}} -A^{(i)\prime})$ that restrict to the previous ones on $A^{(i)}$ and $(-A^{(i)\prime})$. This one-form is also denoted by $\eta(a^i_j)$.
The form $\eta(a^i_j)$ on $(A^{(i)} \cup_{\partial A^{(i)}} -A^{(i)\prime})$ is Poincar\'e dual to the homology class $\partial_{MV}^{-1}(a^i_{j})$ in $(A^{(i)} \cup_{\partial A^{(i)}} -A^{(i)\prime})$, with the notation introduced before Theorem~\ref{thmmainunivlag}. Define the part $A_I^{(i)}$ of $\hcylc_I$ to be $A_I^{(i)}=A^{(i)}$ if $i \notin I$ and $A_I^{(i)}=A^{(i)\prime}$ if $i \in I$.\index[N]{AIi@$A_I^{(i)}$ $3$-handlebody}

The following proposition is the key to the proof of Theorem~\ref{thmmainunivlag}. Its proof is more complicated than I expected. We give it in Chapter~\ref{chapsimnormprop}. 

\begin{proposition}
\label{propnormasim}
There exist homogeneous antisymmetric propagating forms $\omega_I$ \index[N]{omegaI@$\omega_I$ propagating form} of $(C_2(\rats_I),\tau_I)$ with the following properties.
\begin{itemize}
 \item For any subsets $I$ and $J$ of $\underline{x}$, $\omega_I$ and $\omega_J$ coincide on \begin{equation*} C_2\Bigl(\bigl(\rats \setminus \cup_{i \in I \cup J} \Int(A^{(i)})\bigr) \cup_{i \in I \cap J} A^{(i)\prime}\Bigr),\end{equation*}
\item For any $(i,k) \in \underline{x}^2$ such that $i \neq k$, $\omega_I$ can be written as follows on $A_I^{(i)} \times A_I^{(k)}$:
\begin{equation*}\omega_I=\sum_{\substack {j=1,\dots,g_i\\ \ell=1,\dots,g_k }}lk(z^i_j,z^k_{\ell})\projp_{A_I^{(i)}}^{\ast}\bigl(\eta(a^i_j)\bigr) \wedge \projp_{A_I^{(k)}}^{\ast}\bigl(\eta(a^k_{\ell})\bigr),\end{equation*}
where $\projp_{A_I^{(i)}}\colon A_I^{(i)} \times A_I^{(k)} \to A_I^{(i)}$ and  $\projp_{A_I^{(k)}}\colon A_I^{(i)} \times A_I^{(k)} \to A_I^{(k)}$ again denote the natural projections onto the factor corresponding to the subscript.
\end{itemize}
\end{proposition}

Recall Notation~\ref{notcontract}. Let $\Gamma$ be an oriented Jacobi diagram without univalent vertices. 
When $x$ is even and when $G=\bigsqcup_{i \in \underline{x}} T\left(\CI_{A^{(i)}A^{(i)\prime}}\right)$,
define 
\begin{equation*}\langle \langle G  \rangle \rangle_{\!\Gamma\,}=\sum_{p \in P(G) \suchthat \Gamma_p\; \mbox{\scriptsize isomorphic to}\; \Gamma} \left[\ell(p)\Gamma_p\right],\end{equation*} 
where the sum runs over the $p$ such that $\Gamma_p$ is isomorphic to $\Gamma$ as a nonoriented trivalent graph.

Assuming Proposition~\ref{propnormasim}, one can prove the following lemma.
\begin{lemma}
\label{lemthindiv}
Let $\omega_I$ be forms as in Proposition~\ref{propnormasim}. Let $\Gamma$ be an oriented Jacobi diagram on $\sourcetl$.
If $\Gamma$ has less than $x$ trivalent vertices, then we have 
\begin{equation*}\sum_{I \subseteq \underline{x}}(-1)^{\cardlef{I}}I\bigl(\rats_I,\tanghcyll,\Gamma,o(\Gamma),(\omega_I)\bigr)=0.\end{equation*}
If $\Gamma$ is a trivalent Jacobi diagram with $x$ vertices, then we have
\begin{equation*}\sum_{I \subseteq \underline{x}}(-1)^{\cardlef{I}}I\bigl(\rats_I,\tanghcyll,\Gamma,(\omega_I)\bigr)\left[\Gamma\right]=\cardlef{\Aut(\Gamma)}\Bigl\langle \Bigl\langle  \bigsqcup_{i \in \underline{x}} T\left(\CI_{A^{(i)}A^{(i)\prime}}\right) \Bigr\rangle \Bigr\rangle_{\!\!\Gamma\,}.\end{equation*}
\end{lemma}
\bp Let $n \in \NN$. Let $\Gamma$ be an oriented degree $n$ Jacobi diagram on $\sourcetl$. Let us compute \begin{equation*}\Delta=\sum_{I \subseteq \underline{x}}(-1)^{\cardlef{I}}I\bigl(\rats_I,\tanghcyll,\Gamma,o(\Gamma),(\omega_I)\bigr).\end{equation*}
Number the vertices of $\Gamma$ in $\underline{2n}$ arbitrarily.
So the open configuration space $\check{C}(\crats_I,\tanghcyll;\Gamma)$ becomes a submanifold of $\crats_I^{2n}$. The order of the vertices orders the oriented local factors (some of which are tangle components) of $\check{C}(\crats_I,\tanghcyll;\Gamma)$. Thus, it orients $\check{C}(\crats_I,\tanghcyll;\Gamma)$.
Orient the edges of $\Gamma$ so that the edge-orientation of $H(\Gamma)$ and the vertex-orientation of $\Gamma$ induce the above orientation of $\check{C}(\crats_I,\tanghcyll;\Gamma)$ as in Lemma~\ref{lemorc}.

Let $i \in \underline{x}$. The forms $\bigwedge_{e \in E(\Gamma)}p_e^{\ast}(\omega_I)$ over \begin{equation*}\check{C}(\crats_I,\tanghcyll;\Gamma) \cap (\crats_I \setminus A^{(i)}_I)^{2n}\end{equation*}
are identical for $I=K$ and $I=K\cup\{i\}$ for any $K \subseteq \underline{x} \setminus \{i\}$.
Since their integrals enter the sum $\Delta$ with opposite signs, they cancel each other.
This argument allows us to get rid of the contributions of the integrals over the
$\check{C}(\crats_I,\tanghcyll;\Gamma) \cap (\crats_I \setminus A^{(1)}_I)^{2n}$
for any $I \subseteq \underline{x}$. The contributions over the 
\begin{equation*}\Bigl(\check{C}(\crats_I,\tanghcyll;\Gamma) \setminus \bigl(\check{C}(\crats_I,\tanghcyll;\Gamma) \cap (\crats_I \setminus A^{(1)}_I)^{2n} \bigr) \cap (\crats_I \setminus A^{(2)}_I)^{2n}\end{equation*} cancel in the same way. Iterating,
we get rid of the contributions of the integrals over the
$\check{C}(\crats_I,\tanghcyll;\Gamma) \cap (\crats_I \setminus A^{(i)}_I)^{2n}$
for any $i \in\underline{x}$, and for any $I \subseteq \underline{x}$.
Thus, we are left with the contributions of the integrals over the subsets $Q_I$ of \begin{equation*}\check{C}(\crats_I,\tanghcyll;\Gamma) \subset \crats_I^{2n}\end{equation*} with the following property:\\
For any $i \in \underline{x}$, any element of $Q_I$ is sent to $A^{(i)}_I$ under at least one of the $(2n)$ projections onto $\crats_I$.
These subsets $Q_I$ are empty if $2n <x$. Thus, the lemma is proved
when $\Gamma$ has less than $x$ trivalent vertices. 

Assume $\Gamma$ is trivalent and $x=2n$. Then
 $Q_I$ is equal to
\begin{equation*}\cup_{\sigma \in \mathfrak{S}_{2n}}\prod_{i=1}^{2n}A_I^{(\sigma(i))},\end{equation*}
where $\mathfrak{S}_{2n}$ is the set of permutations of $\underline{2n}$.
We get
\begin{equation*}\Delta=\sum_{\sigma \in \mathfrak{S}_{2n}}\Delta_{\sigma}\end{equation*}
with 
\begin{equation*}\Delta_{\sigma}=\sum_{I \subseteq \underline{2n}}(-1)^{\cardlef{I}} \int_{\prod_{i=1}^{2n}A_I^{(\sigma(i))}}\bigwedge_{e \in E(\Gamma)}p_e^{\ast}(\omega_I).\end{equation*}
Let us compute $\Delta_{\sigma}$.
For any $i \in \underline{x}$, $p_i \colon \check{C}(\crats_I,\tanghcyll;\Gamma) \longrightarrow \crats_I$
denotes the projection onto the $i$th factor.
When $e$ is an oriented edge from the vertex $x(e) \in V(\Gamma)$ to $y(e) \in V(\Gamma)$, we have
\begin{multline*}{p_e^{\ast}(\omega_I)\vert_{\prod_{i=1}^{2n}A_I^{(\sigma(i))}}=\\
\sum_{\substack{j=1,\dots,g_{\sigma(x(e))}\\ \ell=1,\dots,g_{\sigma(y(e))} }}
lk\Bigl(z^{\sigma(x(e))}_j,z^{\sigma(y(e))}_{\ell}\Bigr)p_{x(e)}^{\ast}\Bigl(\eta\bigl(a^{\sigma(x(e))}_j\bigr)\Bigr) 
\wedge p_{y(e)}^{\ast}\Bigl(\eta\bigl(a^{\sigma(y(e))}_{\ell}\bigr)\Bigr),}\end{multline*}
where the vertices are regarded as elements of $\underline{2n}$ via the numbering.
Recall the sets $E(\Gamma)$ and $H(\Gamma)$ of edges and half-edges of $\Gamma$.
For a half-edge $c$, let $v(c)$ denote the label of the vertex contained in $c$.

Let $F_{\sigma}$ denote the set of maps $f$
from $H(\Gamma)$ to $\NN$ such that ${f(c)} \in \{1,2,\dots,g_{\sigma(v(c))}\}$ for any $c \in H(\Gamma)$.
For such a map $f$, $f(x(e))$ (resp. $f(y(e))$) denotes the value of $f$ at the half-edge of $e$ that contains $x(e)$ (resp. $y(e)$). We have
\begin{equation*}\Delta_{\sigma}=
\sum_{f \in F_{\sigma}}\biggl(\prod_{e \in E(\Gamma)}lk\Bigl(z^{\sigma(x(e))}_{f(x(e))},z^{\sigma(y(e))}_{f(y(e))}\Bigr)\biggr) I(f)\end{equation*}
with
\begin{equation*}I(f)= \int_{\prod_{i=1}^{2n}(A^{(\sigma(i))} \cup -A^{(\sigma(i))\prime})} \bigwedge_{e \in E(\Gamma)} p_{x(e)}^{\ast}\Bigl(\eta\bigl(a^{\sigma(x(e))}_{f(x(e))}\bigr)\Bigr) \wedge p_{y(e)}^{\ast}\Bigl(\eta\bigl(a^{\sigma(y(e))}_{f(y(e))}\bigr)\Bigr).\end{equation*}
We also have
\begin{equation*}I(f)= \prod_{i=1}^{2n}\int_{(A^{(\sigma(i))} \cup -A^{(\sigma(i))\prime})} \bigwedge_{c \in H(\Gamma) \suchthat v(c)=i} p_{i}^{\ast}\Bigl(\eta\bigl(a^{\sigma(i)}_{{f(c)}}\bigr)\Bigr),\end{equation*}
where we order the exterior product's factors over the half-edges of $v^{-1}(i)$ according to the vertex-orientation of $i$.
Observe
\begin{equation*} \int_{A^{(\sigma(i))} \cup (-A^{(\sigma(i))\prime})} \bigwedge_{c \in  v^{-1}(i)}\eta\bigl(a^{\sigma(i)}_{f(c)}\bigr)= \CI_{{A^{(\sigma(i))}}A^{(\sigma(i))\prime}}\Bigl(\bigotimes_{c \in v^{-1}(i)} a^{\sigma(i)}_{f(c)}\Bigr),\end{equation*}
where the factors of the tensor product are ordered according to the vertex-orientation of $i$, again. Indeed, the closed form $\eta\bigl(a^{\sigma(i)}_{f(c)}\bigr)$ is dual to the homology class $\partial_{MV}^{-1}\bigl(a^{\sigma(i)}_{f(c)}\bigr)$ in $A^{(\sigma(i))} \cup (-A^{(\sigma(i))\prime})$, with the notation introduced before Theorem~\ref{thmmainunivlag}.

Summarizing, we get
\begin{equation*}\Delta_{\sigma}=\sum_{f \in F_{\sigma}}\biggl(\biggl(\prod_{e \in E(\Gamma)}lk\Bigl(z^{\sigma(x(e))}_{f(x(e))},z^{\sigma(y(e))}_{f(y(e))}\Bigr)\biggr)\biggl(\prod_{i \in \underline{2n}}\CI_{{A^{(\sigma(i))}}A^{(\sigma(i))\prime}}\Bigl(\bigotimes_{c \in v^{-1}(i)} a^{\sigma(i)}_{f(c)}\Bigr)\biggr)\biggr).\end{equation*}

We may restrict the sum to the subset $\tilde{F}_{\sigma}$ of $F_{\sigma}$ consisting of the maps $f$ of $F_{\sigma}$ that restrict to $v^{-1}(i)$ as injections for any $i$.

Finally, $\Delta$ is a sum running over all the ways of renumbering the vertices of $\Gamma$ by elements of $\underline{x}$ (via $\sigma$) and of coloring the half-edges $c$ of $v^{-1}(i)$ by three distinct curves $z^{\sigma(i)}_{f(c)}$ via $f$.
In particular, a pair $(\sigma,f)$ provides  
a tripod \begin{equation*}\bigtripod{${z}^y_j$}{${z}^y_k$}{${z}^y_{\ell}$}\end{equation*}
 such that $1 \leq j < k< \ell \leq g_y$ for any $y \in \underline{x}$,
and it provides a pairing of the ends of the univalent vertices of the tripods, giving rise
to the graph $\Gamma$ with a possibly different vertex-orientation.
The vertices of the obtained graph are furthermore numbered by the numbering of the vertices of $\Gamma$, and its edges are identified with the original edges of $\Gamma$.

Fix a set of tripods associated to the elements of $\underline{x}$ as above and a pairing of their univalent vertices giving rise to $\Gamma$ as a nonoriented graph.
Then there are $\cardlef{\Aut(\Gamma)}$ ways of numbering its vertices and edges to get a graph isomorphic to $\Gamma$. So the pairing occurs $\cardlef{\Aut(\Gamma)}$ times.
\eop

\bpo{Proof of Theorem~\ref{thmmainunivlag} up to the unproved assertions of this section\\... restated precisely at the end of the proof}
Lemma~\ref{lemthindiv} and Proposition~\ref{propdefhomog} easily imply
\begin{equation*}\begin{array}{llll}
\sum_{I \subseteq \underline{x}}(-1)^{\cardlef{I}}\Zinv_n\bigl(\hcylc_I,\tanghcyll,(\omega_I)\bigr)&=& 0 &\mbox{if}\; 2n < x,\\
&=&\Bigl[\bigl\langle \bigl\langle \bigsqcup_{i \in \underline{x}} T\bigl(\CI_{A^{(i)}A^{(i)\prime}}\bigr) \bigr\rangle \bigr\rangle \Bigr]&\mbox{if}\; 2n =x .\end{array}\end{equation*}
Theorem~\ref{thmfstconsttangpseudo} (or Theorem~\ref{thmfstconsttang} if pseudo-parallelizations are not required) implies
\begin{equation*}\Zinvuf\left(\hcylc,\tanghcyll \right)=\exp\left(-\frac14 p_1(\tau)\ansothree\right)\prod_{j=1}^k\Bigl(\exp\bigl(-I_{\theta}(K_j,\tau)\alpha\bigr)\#_j\Bigr) \Zinv(\hcylc,\tanghcyll,\tau).\end{equation*}
Set \begin{equation*}Y_n=\sum_{I \subseteq \underline{x}}(-1)^{\cardlef{I}}\left(\exp\left(-\frac14 p_1(\tau_I)\ansothree\right)\Zinv\right)_{\!n}\Bigl(\hcylc_I,\tanghcyll,(\omega_I)\Bigr),\end{equation*}
where $(.)_n$ stands for the degree $n$ part.
Note that the framing corrections $\prod_{j=1}^k\left(\exp(-I_{\theta}(K_j,\tau_I)\alpha)\#_j\right)$ do not depend on $I$. Therefore,
it suffices to prove
\begin{equation*}Y_n=\sum_{I \subseteq \underline{x}}(-1)^{\cardlef{I}}\Zinv_n\Bigl(\hcylc_I,\tanghcyll,(\omega_I)\Bigr)\end{equation*} when $2n \leq x$. We have
\begin{multline*}\left(\exp\left(-\frac14 p_1(\tau_I)\ansothree\right)\Zinv\right)_{\!n}\bigl(\hcylc_I,\tanghcyll,(\omega_I)\bigr)=\Zinv_n\bigl(\hcylc_I,\tanghcyll,(\omega_I)\bigr) \\+ \sum_{j<n}\Zinv_j\bigl(\hcylc_I,\tanghcyll,(\omega_I)\bigr) P_{n-j}(I),\end{multline*}
with  \begin{equation*}P_{n-j}(I) = \biggl( \exp\Bigl(-\frac14 \bigl(p_1(\tau) +\sum_{i \in I}p(i)\bigr)\bigl(\ansothree_1 + \ansothree_3 + \dots\bigr) \Bigr)\biggr)_{\!n-j},\end{equation*}
thanks to Lemma~\ref{lempofi}.

This element $P_{n-j}(I)$ of $\Aavis_{n-j}(\emptyset) $
can be expanded as a combination $\sum m_{\Gamma,g,K}\left[\Gamma\right]$,
where 
\begin{itemize}
\item the $\Gamma$ are degree $(n-j)$ Jacobi diagrams,
\item the $K$ are subsets of $I$ with cardinality at most $n-j$,
\item the $m_{\Gamma,g,K}$ are monomials of degree
at most $(n-j)$ in $p_1(\tau)$ and in the $p(i)$ for $i \in K$,
\item the $m_{\Gamma,g,K}$ do not depend on $I$,
\item for every $i \in K$, the variable $p(i)$ actually occurs in the monomial $m_{\Gamma,g,K}$, and 
\item for a subset $J$ of $\underline{x}$, $m_{\Gamma,g,K}\left[\Gamma\right]$ appears in $P_{n-j}(J)$ if and only if $K$ is a subset of $J$.
\end{itemize}

Thus, we can write the sum of the undesired terms in $Y_n$ by factoring out the $m_{\Gamma,g,K}\left[\Gamma\right]$. The factor of $m_{\Gamma,g,K}\left[\Gamma\right]$ is
\begin{equation*}\sum_{I  \suchthat  K \subseteq I \subseteq \underline{x}}(-1)^{\cardlef{I}}\Zinv_j\bigl(\hcylc_I,\tanghcyll,(\omega_I)\bigr).\end{equation*}
This sum actually runs over the subsets of $\underline{x} \setminus K$. The cardinality of $\underline{x} \setminus K$ is at least $x +j - n$. The inequalities $2n \leq x$ and $j<n$ imply $j<n \leq x-n$, and hence $2j<x +j - n$. Therefore, the beginning of the proof ensures that the above factor of $m_{\Gamma,g,K}\left[\Gamma\right]$ is zero. Hence we have $Y_n=\sum_{I \subseteq \underline{x}}(-1)^{\cardlef{I}}\Zinv_n\left(\hcylc_I,\tanghcyll,(\omega_I)\right).$
This concludes the reduction of the proof of Theorem~\ref{thmmainunivlag} to the proof of Proposition~\ref{propnormasim} given in Chapter~\ref{chapsimnormprop}, and to the proofs that pseudo-parallelizations and associated
propagating forms exist and satisfy the announced properties, given in Chapter~\ref{chappseudopar}.
\eop

\section{Mixed universality statements}
\label{secmixuniv}

We can mix the statements of Theorem~\ref{thmunivsingtang} and \ref{thmmainunivlag}
to get the following statement, which covers both of them.

\begin{theorem}
\label{thmunivmix} Let $y $, $z \in \NN$. Recall $\underline{y}=\{1,2,\dots,y\}$. Set $(\underline{z}+y)=\{y+1,y+2,\dots,y+z\}$.
Let $\tanghcyll$ be a singular $q$-tangle representative in a rational homology cylinder $\hcylc$, whose double points are numbered by $\underline{y}$ and sitting in balls $B_b$ of desingularizations for $b \in \underline{y}$.
For a subset $I$ of $\underline{y}$, let $\tanghcyll_I$ denote the $q$-tangle obtained from $\tanghcyll$ by performing negative desingularizations on double points of $I$ and positive ones on double points of $\underline{y} \setminus I$ in the balls $B_b$. Let $\sqcup_{i=y+1}^{y+z}A^{(i)}$ be a disjoint union of rational homology handlebodies embedded in $\hcylc \setminus (\tanghcyll \cup_{b=1}^{y} B_b)$.
Let $(A^{(i)\prime}/A^{(i)})$ be rational LP-surgeries in $\hcylc$. 

Set
$X=\bigl[\hcylc,\tanghcyll;(A^{(i)\prime}/A^{(i)})_{i \in\underline{z}+y}\bigr]$ and, using Notation~\ref{notzbar}, \begin{equation*}\Zinvufmodis_{n}(X)=  \sum_{I\subseteq \underline{y+z}}(-1)^{\cardlef{I}}\Zinvufmodis_n\Bigl(\hcylc\bigl((A^{(i)\prime}/A^{(i)})_{i \in I \cap  (\underline{z}+y)}\bigr),\tanghcyll_{I \cap  \underline{y}} \Bigr).\end{equation*}
If $2n <2y+z $, then $\Zinvufmodis_{n}(X)$ vanishes. If $2n =2y+z $, then we have
\begin{equation*}\Zinvufmodis_{n}(X)=\left[\Bigl\langle \Bigl\langle \bigsqcup_{i \in \underline{z} +y} T\left(\CI_{A^{(i)}A^{(i)\prime}}\right) \Bigr\rangle \Bigr\rangle_{\!\!\crats(\hcylc)\;}
\sqcup \Gamma_C(\tanghcyll)\right].\end{equation*}
\end{theorem}
\bpo{Proof assuming Theorem~\ref{thmmainunivlag}}
Write $(\hcylc,\tanghcyll)$ as a product of a tangle $\tanghcyll_1$ of the form \begin{equation*}\left(\doublecuptwistsing \dots \doublecuptwistsing \twostraightstrandsnonorlong \dots \twostraightstrandsnonorlong \right)\end{equation*} in the standard rational homology cylinder $\drad{1} \times \left[0,1\right]$, by a nonsingular tangle $\tanghcyll_2$ in $\hcylc$, so that $\sqcup_{i=y+1}^{y+z}A^{(i)}$ is in the latter factor. This can be achieved by moving the double points below. Then Proposition~\ref{propfundsingtang} and Theorem~\ref{thmmainunivlag} imply that
\begin{equation*}\sum_{I\subseteq \underline{y+z}}(-1)^{\cardlef{I}}\Zinvufrf_n\Bigl(\hcylc\bigl((A^{(i)\prime}/A^{(i)})_{i \in I \cap  (\underline{z}+y)}\bigr),\tanghcyll_{I \cap  \underline{y}} \Bigr)\end{equation*}
satisfies the conclusions of the statement with $\Zinvufrf$ instead of $\Zinvufmodis$, by functoriality.
Now, to prove Theorem~\ref{thmunivmix} it suffices to \say{multiply the proof of Theorem~\ref{thmunivsingtang} (at the end of Section~\ref{secproofthmbn}) by $\bigl[\bigl\langle \bigl\langle \bigsqcup_{i \in \underline{z} +y} T(\CI_{A^{(i)}A^{(i)\prime}}) \bigr\rangle \bigr\rangle_{\!\!\crats(\hcylc)\,}\bigr]$} .
\eop

In order to prove more interesting mixed universality properties,
we say more about the normalization of the propagating forms of Proposition~\ref{propnormasim}.

Recall that $[-4,4] \times \partial A^{(i)}$ denotes a regular neighborhood of $\partial A^{(i)}$
embedded in $\hcylc$, that intersects $A^{(i)}$ as $ [-4,0] \times \partial A^{(i)}$.
All the neighborhoods $[-4,4] \times \partial A^{(i)}$ are disjoint from each other and from $\tanghcyll$.
Throughout this paragraph, we use the corresponding  coordinates on the image of this implicit embedding.

For $t \in [-4,4]$, set \index[N]{Ait@$A^{(i)}_t$ $3$-handlebody}
\begin{equation*}A^{(i)}_t=\left\{\begin{array}{ll} A^{(i)} \cup  \bigl(\left[0,t\right] \times \partial A^{(i)}\bigr) & \mbox{if}\;\; t \geq 0\\
 A^{(i)} \setminus \bigl( \left]t,0\right] \times \partial A^{(i)}\bigr) & \mbox{if} \;\;t \leq 0.\end{array} \right.\end{equation*}
We have $\partial A^{(i)}_t =\{t\} \times \partial A^{(i)}.$

For $i \in \underline{x}$, choose a basepoint $p^i$ in $\partial A^{(i)}$
outside the neighborhoods $a^i_j \times [-1,1]$ of the $a^i_j$ and outside neighborhoods $z^i_j \times [-1,1]$ the $z^i_j$.
Fix a path $\left[p^i,q^i\right]$ from $p^i$ to a point $q^i$ of $\partial \hcylc$ in
\begin{equation*}\hcylc \setminus \left(\tanghcyll \cup \Int(A^{(i)}) \cup \cup_{k \suchthat k\neq i} A^{(k)}_4 \right)\end{equation*} 
so that the paths $\left[p^i,q^i\right]$ are pairwise disjoint.
Choose a closed $2$-form $\omega(p^i)$ on $\left(C_1(\rats_I) \setminus \Int(A^{(i)})\right)$ such that
\begin{itemize}
\item the integral of $\omega(p^i)$ \index[N]{omegapi@$\omega(p^i)$ $2$-form} along a closed surface of 
$(\hcylc_I \setminus \Int(A^{(i)}))$ is its algebraic intersection with $\left[p^i,q^i\right]$,
\item the support of $\omega(p^i)$ intersects $\left(\hcylc_I \setminus \Int(A^{(i)})\right)$
inside a tubular neighborhood of $\left[p^i,q^i\right]$ disjoint from \begin{equation*}\left( \cup_{k \suchthat k\neq i} A^k_4\right)\cup \tanghcyll \cup \Biggl(\left[0,4\right] \times \biggl( \cup_{j=1}^{g_i} \Bigl( \left(a^i_j \times [-1,1]\right) \cup \left(z^i_j \times [-1,1]\right) \Bigr)\biggr)  \Biggr).\end{equation*}
\item $\omega(p^i)$ restricts as the usual volume form $\omega_{S^2}$ on $\partial C_1(\rats_I)=S^2$.
\end{itemize}

For $i \in \underline{x}$, for $j=1, \dots, g_i$, the curve $\{4\} \times a^i_j$ bounds a rational chain $\Sigma(a^i_j)$ in 
$A^{(i)}_4$ and a rational chain $\Sigma^{\prime}(a^i_j)$ in 
$A^{(i)\prime}_4$. When viewed as a chain in $\hcylc_I$, such a chain is denoted by $\Sigma_I(a^i_j)$. We have $\Sigma_I(a^i_j)=\Sigma(a^i_j)$ if $i \notin I$, and $\Sigma_I(a^i_j)=\Sigma^{\prime}(a^i_j)$ if $i \in I$.
The form $\eta(a^i_j)$ supported on $[-4,4] \times  a^i_j \times [-1,1]$ in $A^{(i)}_{I,4} \setminus A^{(i)}_{I,-4}$ may be expressed as $\eta(a^i_j)=\projp_{[-1,1]}^{\ast}\left(\eta_{[-1,1]}\right)$ there. Thus, it extends naturally to $A^{(i)}_{I,4}=\bigl(A^{(i)}_{I}\bigr)_4$, as a closed form dual to the chain $\Sigma_I(a^i_j)$. 

For $i \in \underline{x}$, for $j=1, \dots, g_i$, $z^i_j$ bounds a rational chain in $\hcylc_I$. Therefore, it cobounds a rational chain $\Sigma_I(\check{z}^i_j)$ in $\bigl(\hcylc_I \setminus \mathring{A}_I^{(i)}\bigr) \setminus \bigl(\cup_{i=1}^x\left[p^i,q^i\right] \bigr) $ with a combination of $a^i_{\ell}$ with rational coefficients. We have
\begin{equation*}\partial \Sigma_I(\check{z}^i_j)= z^i_j - \sum_{j=1}^{g_i} lk_e\bigl(z^i_j,\{-1\} \times z^i_{\ell}\bigr)a^i_{\ell}=\check{z}^i_j.\end{equation*}
Furthermore, $\Sigma_I(\check{z}^i_j)$ may be assumed to intersect $\mathring{A}_I^{(k)}$ as 
\begin{equation*}\sum_{m=1}^{g_k}lk\bigl({z}^i_j,{z}^k_m\bigr) \Sigma_I(a^k_m),\end{equation*} for $k \neq i$.
There is a closed one-form $\eta_I(z^i_j)$ dual to $\Sigma_I(\check{z}^i_j)$ in $\bigl(\hcylc_I \setminus \mathring{A}_I^{(i)}\bigr)$, such that
$\eta_I(z^i_j)$ is supported near $\Sigma_I(\check{z}^i_j)$ and outside the supports of $\omega(p^i)$ and of the other $\omega(p^k)$,
and we have
\begin{equation*}\eta_I(z^i_j)=\sum_{m=1}^{g_k}lk\bigl({z}^i_j,{z}^k_m\bigr) \eta(a^k_m)\end{equation*} on $\mathring{A}_I^{(k)}$ for $k \neq i$. The integral of $\eta_I(z^i_j)$ along a closed curve of $\bigl(\hcylc_I \setminus \mathring{A}_I^{(i)}\bigr)$ is its linking number with $z^i_j$ in $\hcylc_I$.

We will prove the following proposition in Chapter~\ref{chapsimnormprop}.
\begin{proposition}
\label{propnorma}
The antisymmetric propagating forms $\omega_I$ \index[N]{omegaI@$\omega_I$ propagating form} of $\left(C_2(\rats_I),\tau_I\right)$ of Proposition~\ref{propnormasim} can be chosen so that
\begin{enumerate} 
\item for every $i \in \underline{x}$, the restriction of $\omega_I$ to 
\begin{equation*}\biggl(A_I^{(i)} \times \left(C_1(\rats_I) \setminus A^{(i)}_{I,3}\right) \biggr) \subset C_2(\rats_I)\end{equation*}
is equal to
\begin{equation*}\sum_{j \in \underline{g_i}} p_1^{\ast}\left(\eta_I(a^i_j)\right) \wedge p_2^{\ast}\left(\eta_I(z^i_j)\right) + p_2^{\ast}\left(\omega(p^i)\right),\end{equation*}
where $p_1$ and $p_2$ respectively denote the first and second projection of $A_I^{(i)} \times (C_1(\rats_I) \setminus A^{(i)}_{I,3})$ to $C_1(\rats_I)$, and
\item for every $i$, for any $j \in \{1,2, \dots, g_i\}$, we have  \begin{equation*}\int_{\Sigma_I(a^i_j) \times p^i} \omega_I=0,\end{equation*}
where $p^i \in \partial A_I^{(i)}$ and $\partial \Sigma_I(a^i_j) \subset \{4\} \times \partial A_I^{(i)}$.
\end{enumerate}
\end{proposition}

A \emph{two-leg Jacobi diagram} is a uni-trivalent Jacobi diagram with two univalent vertices, called \emph{legs}\index[T]{leg}. When these legs are colored by possibly noncompact connected components $K_j$ of a tangle $\tanghcyll$, a two-leg diagram gives rise to a diagram on the domain $\sourcetl$ of the tangle by attaching the legs to the corresponding components. The class in $\Aavis(\sourcetl)$ of this diagram is well-defined.
Indeed, according to Lemma~\ref{lemoneleg}, Jacobi diagrams with one univalent vertex vanish in 
$\Aavis(\sourcetl)$. Therefore, the STU relation guarantees that if the same noncompact component colors the two legs, changing their order with respect to the orientation component does not change the diagram class. (See also Lemma~\ref{lemsymvo}.)

Generalize the contraction of trivalent graphs associated to LP-surgeries of Notation~\ref{notcontract} 
to graphs with legs as follows.

\begin{notation} \label{notGamma2}
Let $\tanghcyll \colon \sourcetl \hookrightarrow \hcylc$ be a long tangle representative in a rational homology cylinder $\hcylc$.

Let $G$ be a graph with oriented trivalent vertices and with two kinds of univalent vertices, the \emph{decorated ones} and the \emph{legs}, where the components of the legs of $G$ are edges from a leg to a decorated univalent vertex.
The \emph{legs} are univalent vertices on $\sourcetl$.
The decorated univalent vertex in a leg segment is decorated with the leg component.
The other \emph{decorated univalent vertices} of $G$ are decorated with disjoint
curves of $\crats=\crats(\hcylc)$ disjoint from the image of $\tanghcyll$. Such a curve $c$ bounds a compact oriented surface $\Sigma(c)$ in $\hcylc$, and its linking number with a component $K_u$ of $\tanghcyll$ is the algebraic intersection $\langle K_u,\Sigma(c)\rangle$.
Let $\check{P}(G)$ be the set of partitions of the set of decorated univalent vertices of $G$ in disjoint pairs.

For $p \in \check{P}(G)$, identifying the two decorated vertices of each pair provides a vertex-oriented Jacobi diagram $\Gamma_p$ on $\sourcetl$. Multiplying it by the product $\ell(p)$ over the pairs of $p$ of the linking numbers of the curves that decorate the two vertices yields an element $\left[\ell(p)\Gamma_p\right]$ of $\Aavis(\sourcetl)$.

Define \begin{equation*}\bigl\langle \bigl\langle G \check{\rangle} \bigr\rangle= \projassis\Bigl( \sum_{p \in \check{P}(G)} \bigl[\ell(p)\Gamma_p\bigr] \Bigr) \in \Assis(\sourcetl)\end{equation*}
Extend this contraction to linear combinations
of graphs, linearly.
Assume that the components of $\tanghcyll$ are numbered in $\underline{k}$. To a pair $(u,v)$ of 
$\underline{k}$, associate the univalent graph $G_{uv}=\guv$ consisting of two distinct segments, the first one with its leg on $K_u$, and the second one with its leg on $K_v$. The legs are considered as distinct even if there is an automorphism of $G_{uu}=\guu$ that maps a leg to the other when $K_u$ is a circle. We think of the leg of the first segment as the first leg of $G_{uv}$, and the other leg is the second leg of $G_{uv}$.

For a finite collection $(A^{(i)\prime}/A^{(i)})_{i \in\underline{2x}}$ of disjoint rational LP-surgeries in $\crats(\hcylc) \setminus \tanghcyll$,
define the element 
\begin{equation*}\left[\Gamma_{(2)}\left(\hcylc,\tanghcyll;(A^{(i)\prime}/A^{(i)})_{i \in\underline{2x}}\right)\right]=\frac12\sum_{(u,v) \in \underline{k}^2}\biggl[\Bigl\langle \Bigl\langle \bigsqcup_{i \in \underline{2x}} T\left(\CI_{A^{(i)}A^{(i)\prime}}\right) \sqcup {G}_{uv}\check{\bigr\rangle} \Bigr\rangle\biggr]\end{equation*}
of $\Assis(\sourcetl)$.
\end{notation}

\begin{examples}
Assume \begin{equation*}T\left(\CI_{A^{(1)}A^{(1)\prime}}\right)=\bigtripod{${z}^1_1$}{${z}^1_2$}{${z}^1_{3}$} \;\mbox{and}\; T\left(\CI_{A^{(2)}A^{(2)\prime}}\right)=\bigtripod{${z}^2_1$}{${z}^2_2$}{${z}^2_{3}$},\end{equation*}
where $lk({z}^1_i,{z}^2_j)=0$ as soon as $i \neq j$, and $lk({z}^1_i,K_u)=lk({z}^2_i,K_u)=0$ if $i\neq 3$.
If $K_u$ is an interval, then we get
\begin{equation*}\Bigl\langle \Bigl\langle \bigsqcup_{i \in \underline{2}} T\left(\CI_{A^{(i)}A^{(i)\prime}}\right) \sqcup {G}_{uu}\check{\bigr\rangle} \Bigr\rangle
= -2 lk\left({z}^1_1,{z}^2_1\right)lk\left({z}^1_2,{z}^2_2\right)lk\left({z}^1_3,K_u\right)lk\left({z}^2_3,K_u\right) \left[\vertbubblu\right].\end{equation*}
If $K_u$ is a circle, then we get
\begin{equation*}\Bigl\langle \Bigl\langle \bigsqcup_{i \in \underline{2}} T\left(\CI_{A^{(i)}A^{(i)\prime}}\right) \sqcup {G}_{uu}\check{\bigr\rangle} \Bigr\rangle
= -2 lk({z}^1_1,{z}^2_1)lk({z}^1_2,{z}^2_2)lk({z}^1_3,K_u)lk({z}^2_3,K_u) \left[\twochordbubblu\right].\end{equation*}
\end{examples}

\begin{theorem} 
\label{thmmainunivlagtwolegs}
Let $\tanghcyll \colon \sourcetl \hookrightarrow \hcylc$ be a long tangle representative in a rational homology cylinder $\hcylc$. Let $\sqcup_{i=1}^{2x}A^{(i)}$ be a disjoint union of rational homology handlebodies embedded in $\hcylc \setminus \tanghcyll$.
Let $(A^{(i)\prime}/A^{(i)})$ be rational LP-surgeries in $\hcylc$. Set
$X=\left[\hcylc,\tanghcyll;(A^{(i)\prime}/A^{(i)})_{i \in\underline{2x}}\right]$ and \begin{equation*}\Zinvlinkuf_{n}(X)=  \sum_{I\subseteq \underline{2x}}(-1)^{\cardlef{I}}\Zinvlinkuf_n\left(\hcylc_I,\tanghcyll \right),\end{equation*}
where $\hcylc_I=\hcylc\left((A^{(i)\prime}/A^{(i)})_{i \in I}\right)$ is the rational homology cylinder obtained from $\hcylc$ by performing the LP-surgeries that replace $A^{(i)}$ by $A^{(i)\prime}$ for $i \in I$. Assume that $x \neq 0$.
If $n <x+1 $, then $\Zinvlinkuf_{n}(X)$ vanishes. If $n=x+1$, then we have
\begin{equation*}\Zinvlinkuf_{n}(X)=\Bigl[\Gamma_{(2)}\bigl(\hcylc,\tanghcyll;(A^{(i)\prime}/A^{(i)})_{i \in\underline{2x}}\bigr)\Bigr],\end{equation*}
where $\Zinvlinkuf$ is defined in Proposition~\ref{propwithZinvlink}. See also Notation~\ref{notAcheck} and Notation~\ref{notGamma2}.
\end{theorem}
\begin{remark}
 This result holds modulo 1T when $x=0$.
\end{remark}

\bpo{Proof of Theorem~\ref{thmmainunivlagtwolegs} assuming Proposition~\ref{propnorma}}
Recall that $\Zinvlinkuf$ takes its values in $\Assis(\sourcetl)$, where the diagrams with connected trivalent components vanish. Therefore, Theorem~\ref{thmmainunivlag} implies the result when $n<x+1$. Assume that $n=x+1$. Since the framing correction terms involving $\ansothree$ vanish in $\Assis(\sourcetl)$, and since the other correction terms
are the same for all the $\Zinvlinkuf\left(\hcylc_I,\tanghcyll \right)$, we get
\begin{equation*}\Zinvlinkuf_{x+1}(X)=\sum_{I\subseteq \underline{2x}}(-1)^{\cardlef{I}}\Zinvlink_n\bigl(\hcylc_I,\tanghcyll,(\omega_I) \bigr).\end{equation*}
Let $\Gamma$ be a Jacobi diagram of degree $x+1$ on $\sourcetl$ that contributes to $\Zinvlinkuf_{x+1}(X)$. Since its class does not vanish in $\Assis_{x+1}(\sourcetl)$, each component of $\Gamma$ must contain at least two univalent vertices because one-leg diagrams vanish in $\Assis(\sourcetl)$. Furthermore, as in the proof of Lemma~\ref{lemthindiv}, the configurations must involve at least one point in some $A_I^{(i)}$, for each $i$ in $\underline{2x}$. Such a point must be a trivalent vertex position. Therefore, the graph $\Gamma$ has at least $2x$ trivalent vertices.
Finally, $\Gamma$ must be a connected two-leg Jacobi diagram of degree $x+1$ on $\sourcetl$.
Assume that the univalent vertices of $\Gamma$ are on components $K_u$ and $K_v$.
Order the set of univalent vertices of $\Gamma$ so that the first is on $K_u$ and the second is on $K_v$. Let $\Gamma^U$ denote the graph $\Gamma$ equipped with this order.
Number the trivalent vertices of $\Gamma$ in $\underline{2x}$. Orient the open configuration space $\check{C}(\crats_I,\tanghcyll;\Gamma)$ as an open oriented submanifold of $K_u \times K_v \times \crats_I^{2x}$ with respect to the order of $V(\Gamma)$ induced by the numbering of its elements.
Orient the edges and the trivalent vertices of $\Gamma$ so that the induced orientation of $\check{C}(\crats_I,\tanghcyll;\Gamma)$ matches the previous one.

As in the proof of Lemma~\ref{lemthindiv}, we have
\begin{equation*}\Delta(\Gamma)=\sum_{I \subseteq \underline{2x}}(-1)^{\cardlef{I}}I(\rats_I,\tanghcyll,\Gamma,o(\Gamma),(\omega_I))\left[\Gamma\right]=\sum_{\sigma \in \mathfrak{S}_{2x}}\Delta_{\sigma}\left[\Gamma\right],\end{equation*} with
\begin{equation*}\Delta_{\sigma}=\sum_{I \subseteq \underline{2x}}(-1)^{\cardlef{I}} \int_{D(I,\sigma)}\bigwedge_{e \in E(\Gamma)}p_e^{\ast}\left(\omega_I\right),\end{equation*}
and
\begin{equation*}D(I,\sigma)=K_u \times K_v \times \prod_{i=1}^{2x}A_I^{(\sigma(i))},\end{equation*}
when $K_u \neq K_v$, or when $K_u$ is a closed component. When $K_u=K_v$ and $K_u$ is a long component, the factor $K_u \times K_v$ must be replaced by one of its two subsets of pairs with a fixed order on $K_u$.
When $e$ is an oriented edge between two trivalent vertices, recall the expression of $p_e^{\ast}(\omega_I)$ from 
Proposition~\ref{propnormasim}.
Without loss of generality, assume that the legs are the first half-edges of their edges. With the projections $p_i:\check{C}(\crats_I,\tanghcyll;\Gamma) \longrightarrow \crats_I$, Proposition~\ref{propnorma} implies
\begin{equation*}p_e^{\ast}\left(\omega_I\right)\vert_{D(I,\sigma)}=\sum_{j \in \underline{g_{\sigma(y(e))}}} p_{x(e)}^{\ast}\left(\eta(z^{\sigma(y(e))}_j)\right) \wedge p_{y(e)}^{\ast}\left(\eta(a^{\sigma(y(e))}_{j})\right)
\end{equation*} for an edge from a leg $x(e)$ to a trivalent vertex $y(e)$.
In this case, the edge $e$ is associated to $x(e)$ and denoted by $e(x(e))$.
Let $H_T(\Gamma)$ denote the subset of $H(\Gamma)$ consisting of the half-edges that contain a trivalent vertex, and let $E_T(\Gamma)$ denote the set of edges of $\Gamma$ between trivalent vertices.
Let $F_{\sigma}$ denote the set of maps $f$
from $H_T(\Gamma)$ to $\NN$ such that ${f(c)}$ belongs to $\{1,2,\dots,g_{\sigma(v(c))}\}$ for any $c \in H_T(\Gamma)$. We have
\begin{equation*}\Delta_{\sigma}=
\sum_{f \in F_{\sigma}}\left(\prod_{e \in E_T(\Gamma)}lk\left(z^{\sigma(x(e))}_{f(x(e))},z^{\sigma(y(e))}_{f(y(e))}\right)\right) I(f),\end{equation*}
where $I(f)$ is equal to
\begin{equation*}\int_{K_u} \eta\left(z^{\sigma(y(u))}_{f(y(u))}\right) \int_{K_v} \eta\left(z^{\sigma(y(v))}_{f(y(v))}\right)\times \prod_{i=1}^{2x}\int_{(A^{(\sigma(i))} \cup -A^{(\sigma(i))\prime})} \bigwedge_{c \in H(\Gamma) \suchthat v(c)=i} p_{i}^{\ast}\left(\eta(a^{\sigma(i)}_{{f(c)}})\right),\end{equation*}
when $K_u \neq K_v$, or when $K_u$ is a closed component.
Recall \begin{equation*}\int_{K_u} \eta\left(z^{\sigma(y(u))}_{f(y(u))}\right)=lk\left(K_u,z^{\sigma(y(u))}_{f(y(u))}\right).\end{equation*}
Summarizing, when $K_u \neq K_v$ or when $K_u$ is a closed component, we have
\begin{equation*}\Delta_{\sigma}=\sum_{f \in F_{\sigma}}\left(\left(\prod_{e \in E(\Gamma)}lk(e;f)\right)\left(\prod_{i \in \underline{2x}}\CI_{{A^{(\sigma(i))}}A^{(\sigma(i))\prime}}\left(\bigotimes_{c \in v^{-1}(i)} a^{\sigma(i)}_{f(c)}\right)\right)\right),\end{equation*}
where $lk(e;f)=
lk\left(z^{\sigma(x(e))}_{f(x(e))},z^{\sigma(y(e))}_{f(y(e))}\right)$ when $e \in E_T(\Gamma)$,
and \begin{equation*}lk(e;f)=
lk\left(K_{x(e)},z^{\sigma(y(e))}_{f(y(e))}\right)\end{equation*} when $x(e)$ is univalent.

Finally, $\Delta(\Gamma)$ is a sum, running over all the ways of renumbering the trivalent vertices of $\Gamma$ by elements of $\underline{2x}$ (via $\sigma$), and of coloring the half-edges $c$ of $v^{-1}(i)$ by three distinct curves $z^{\sigma(i)}_{f(c)}$ via $f$.

In particular, a pair $(\sigma,f)$ provides 
a tripod \begin{equation*}\bigtripod{${z}^y_j$}{${z}^y_k$}{${z}^y_{\ell}$}\end{equation*}
for any $y \in \underline{2x}$ such that $1 \leq j < k< \ell \leq g_y$.
It also provides a pairing of the ends of the univalent vertices of the tripods and of the legs on $K_u$ and $K_v$ (the first one and the second one when $v=u$), which gives rise
to the graph $\Gamma$ with a possibly different vertex-orientation.
The vertices of the obtained graph are furthermore numbered by the numbering of the vertices of $\Gamma$, and its edges are identified with the original edges of $\Gamma$. The order of $U(\Gamma)$ is induced by the order on the legs of $G_{uv}$.

Let $\Aut(\Gamma^U)$ denote the set of automorphisms of $\Gamma$ that fix the univalent vertices of $\Gamma$ (pointwise). (The set $\Aut(\Gamma^U) $ is distinct from $\Aut(\Gamma)$ if and only if $u$=$v$, $K_u$ is a closed component, and there exists an automorphism of $\Gamma$ that exchanges its two univalent vertices.)

Fix a set of tripods associated to the elements of $\underline{2x}$ as above and a pairing of their univalent vertices and the legs on $K_u$ and $K_v$ (which are distinguished as the first one and the second one when $v=u$).
Then there are exactly $\cardlef{\Aut(\Gamma^U)}$ ways of numbering its vertices and edges to get a graph isomorphic to $\Gamma^U$ by an isomorphism that fixes the univalent vertices. So the pairing occurs $\cardlef{\Aut(\Gamma^U)}$ times.

Let $G$ be a family of $2x$ tripods 
\begin{equation*}\bigtripod{${z}^y_j$}{${z}^y_k$}{${z}^y_{\ell}$},\end{equation*}
one for each $y \in \underline{2x}$, such that $1 \leq j < k< \ell \leq g_y$.
Define 
\begin{equation*}\bigl[\langle \langle  G \sqcup {G}_{uv} \rangle \rangle_{\Gamma^U}\bigr]=\sum_{p \in \check{P}(G \sqcup {G}_{uv}) \suchthat \Gamma_p \; \mbox{\scriptsize isomorphic to}\; \Gamma^U} \left[\ell(p)\Gamma_p\right],\end{equation*}
where the sum runs over the $p$ such that $\Gamma_p$ is isomorphic to $\Gamma^U$, as a nonoriented uni-trivalent graph on $\sourcetl$, equipped with a fixed order on $U(\Gamma)$, with the notation introduced before the statement of Theorem~\ref{thmmainunivlagtwolegs}.
Similarly define
\begin{equation*}\bigl[\langle \langle  G \sqcup {G}_{uv} \rangle \rangle_{\Gamma}\bigr]=\sum_{p \in \check{P}(G \sqcup {G}_{uv}) \suchthat \Gamma_p \; \mbox{\scriptsize isomorphic to}\; \Gamma} \left[\ell(p)\Gamma_p\right],\end{equation*}
where the sum runs over the $p$ such that $\Gamma_p$ is isomorphic to $\Gamma$, as a nonoriented uni-trivalent graph on $\sourcetl$.
(We forget the order on $U(\Gamma)$.)

We get \begin{equation*}\Delta(\Gamma)=\cardlef{\Aut(\Gamma^U)}\left[\left\langle \left\langle \bigsqcup_{i \in \underline{2x}} T\left(\CI_{A^{(i)}A^{(i)\prime}}\right) \sqcup {G}_{uv} \right\rangle \right\rangle_{\!\!\Gamma^U\,}\right].\end{equation*}

If $u\neq v$, then $\Gamma^U$ and $\Gamma$ coincide, and we have
\begin{equation*}\frac{\Delta(\Gamma)}{\cardlef{\Aut(\Gamma)}}=\frac12 \Bigl(
\Bigl[\bigl\langle \bigl\langle \bigsqcup_{i \in \underline{2x}} T(\CI_{A^{(i)}A^{(i)\prime}}) \sqcup {G}_{uv}\bigr\rangle \bigr\rangle_{\!\Gamma\,}\Bigr] + \Bigl[\bigl\langle \bigl\langle \bigsqcup_{i \in \underline{2x}} T(\CI_{A^{(i)}A^{(i)\prime}}) \sqcup {G}_{vu}\bigr\rangle \bigr\rangle_{\!\Gamma\,}\Bigr] \Bigr).\end{equation*}

Assume that $u=v$ and that $K_u$ is a circle.
If there exists an automorphism of $\Gamma$ that exchanges its univalent vertices, then we have $\cardlef{\Aut(\Gamma)}=2\cardlef{\Aut(\Gamma^U)}$ and $\langle \langle . \rangle \rangle_{\Gamma^U}=\langle \langle . \rangle \rangle_{\Gamma}$. Otherwise, we have $\cardlef{\Aut(\Gamma)}=\cardlef{\Aut(\Gamma^U)}$ and $\left[\langle \langle . \rangle \rangle_{\Gamma}\right]=2\left[\langle \langle. \rangle \rangle_{\Gamma^U}\right]$. So we obtain
\begin{equation*}\Delta(\Gamma)=\frac12 \cardlef{\Aut(\Gamma)}\left[\left\langle \left\langle \bigsqcup_{i \in \underline{2x}} T\left(\CI_{A^{(i)}A^{(i)\prime}}\right) \sqcup {G}_{uu} \right\rangle \right\rangle_{\!\Gamma\,}\right]\end{equation*}
in any case.

It remains to study the case in which the univalent vertices of $\Gamma$ belong to the same noncompact component $K_u$. In this case, we compute the sum $\Delta(\Gamma)+\Delta(\Gamma^s)$,
where $\Gamma^s$ is obtained from $\Gamma$ by exchanging the order of its univalent vertices on
$K_u$. 

Again, we find \begin{equation*}\Delta(\Gamma)+\Delta(\Gamma^s)=\cardlef{\Aut(\Gamma)}\left[\left\langle\left\langle  \bigsqcup_{i \in \underline{2x}} T\left(\CI_{A^{(i)}A^{(i)\prime}}\right) \sqcup {G}_{uu} \right\rangle \right\rangle_{\!\!\Gamma\,}\right],\end{equation*}
where the contraction $\langle\langle.\rangle \rangle_{\Gamma}$ keeps only the graphs that are isomorphic to $\Gamma$ (as a graph with an ordered pair of free legs). 
(Recall $\left[\Gamma^s\right]=\left[\Gamma\right]$ in $\Assis(\sourcetl)$). 
If $\Gamma$ and $\Gamma^s$ are isomorphic, then we have
$\Delta(\Gamma)=\Delta(\Gamma^s)$ and
 \begin{equation*}\Delta(\Gamma)=\frac12 \cardlef{\Aut(\Gamma)}\left[\left\langle\left\langle  \bigsqcup_{i \in \underline{2x}} T\left(\CI_{A^{(i)}A^{(i)\prime}}\right) \sqcup {G}_{uu} \right\rangle \right\rangle_{\!\!\Gamma\,}\right].\end{equation*}
Otherwise, we have $\cardlef{\Aut(\Gamma)}=\cardlef{\Aut(\Gamma^s)}$
and $\left[\langle\langle.\sqcup {G}_{uu}\rangle \rangle_{\Gamma}\right]=\left[\langle\langle.\sqcup {G}_{uu}\rangle \rangle_{\Gamma^s}\right]$. So we get
\begin{equation*}\sum_{\substack{\Gamma \;\mbox{\scriptsize as above} \\ \mbox{\scriptsize with $2$ univalent vertices on}\; K_u}} \frac{\Delta(\Gamma)}{\cardlef{\Aut(\Gamma)}}=\frac12\left[\left\langle\left\langle  \bigsqcup_{i \in \underline{2x}} T\left(\CI_{A^{(i)}A^{(i)\prime}}\right) \sqcup {G}_{uu} \right\rangle \right\rangle\right]\end{equation*} in any case.
\eop

As in Section~\ref{secfintypZ}, we can associate an alternate sum of tangles to a framed embedding of a Jacobi diagram on a tangle $\tanghcyll$ into a rational homology cylinder.

Let $\Gamma$ be such a Jacobi diagram, whose trivalent vertices are numbered in $\underline{x}$.
Let $\Sigma(\Gamma)$ be an oriented surface containing $\Gamma$ in its interior such that $\Sigma(\Gamma)$ is a regular neighborhood of $\Gamma$ in $\Sigma(\Gamma)$. Equip $\Gamma$ with its vertex-orientation induced by the orientation of $\Sigma(\Gamma)$.
Embed $\Sigma(\Gamma)$ in $\hcylc$, so that $\tanghcyll$ intersects $\Sigma(\Gamma)$ near univalent vertices as in \begin{equation*}\Nchordrotate\mbox{ or }\Nunivv,\end{equation*} $\tanghcyll$ is tangent to $\Sigma(\Gamma)$ at univalent vertices, and $\tanghcyll$ does not meet $\Sigma(\Gamma)$ outside such neighborhoods of the univalent vertices. Note that the embedding of $\Sigma(\Gamma)$ induces a local orientation of $\tanghcyll$ as in Definition~\ref{defrkoruniv}. Replace a chord between two univalent vertices by a crossing change so that
\begin{equation*}\Nchord\mbox{ encodes the singular point }\Nchordsing \end{equation*} associated to the positive crossing change from \begin{equation*}\Nchordnegwithoutsurf\mbox{ to }\Nchordpos.\end{equation*}  Replace neighborhoods \begin{equation*}\Nunivv\mbox{ of the edges }\univv \end{equation*} between a trivalent vertex and a univalent vertex by neighborhoods \begin{equation*}\Nunivvhopf\mbox{ of }\univvhopf.\end{equation*}
Finally, replace neighborhoods of the edges between trivalent vertices as in Section~\ref{secfintypZ}.

Thus, $\Sigma(\Gamma)$ transforms $\tanghcyll$ into a singular tangle $\tanghcyll^s$ whose double points are associated to the chords
of $\Gamma$ as above, equipped with a collection of disjoint oriented surfaces \begin{equation*}\Sigma(Y) = \normyep,\end{equation*} associated to trivalent vertices of $\Gamma$. The oriented surfaces $\Sigma(Y)$ are next thickened to become framed genus $3$ handlebodies.

Define $\psi(\Sigma(\Gamma))$ to be $\left[\hcylc,\tanghcyll^s;(A^{(i)\prime}/A^{(i)})_{i\in \underline{x}}\right]$, where the surgeries $(A^{(i)\prime}/A^{(i)})$ associated to the trivalent vertices of $\Gamma$ are defined as in Section~\ref{secfintypZ}. Set $\Zinvlinkuf\bigl(\psi(\Sigma(\Gamma))\bigr)=
\sum_{I\subseteq \underline{x}}(-1)^{x+\cardlef{I}}\Zinvlinkuf_n\left(\hcylc_I,\tanghcyll^s \right)$,
where $\hcylc_I=\hcylc\left((A^{(i)\prime}/A^{(i)})_{i \in I}\right)$.

\begin{corollary} Let $n$ be a positive integer.
 Let $\Gamma$ be a degree $n$ Jacobi diagram with two univalent vertices. Let $\Sigma(\Gamma)$ be a regular neighborhood of $\Gamma$  embedded in $\hcylc$ as above.
 Then we have $\Zinvlinkuf_{\leq n}\bigl(\psi(\Sigma(\Gamma))\bigr)=\left[\Gamma\right]$.
\end{corollary}
\bp This follows from Theorem~\ref{thmmainunivlagtwolegs} as in the proof of Lemma~\ref{lemkt}.
\eop

\begin{remark}
This corollary could be true for more general Jacobi diagrams. 
\end{remark}

To finish this section, we apply Theorem~\ref{thmmainunivlagtwolegs} with the LP-surgery of Subsection~\ref{sublagboun} to compute the degree $2$ part of $\Zinvlinkuf$ for a null-homologous knot and prove Theorem~\ref{thmcompztwo} below.
Recall the definition of $a_2(K)$ from the statement of Proposition~\ref{propcaswsurfor} and the lines that follow.
\begin{theorem}
\label{thmcompztwo}
Let $K$ be a null-homologous knot in a rational homology sphere $\rats$. 
Then we have
\begin{equation*} \Zinvlinkuf_2(\rats,K) =\left(\frac1{24}- a_2(K)\right)\left[\tripodnonor\right].\end{equation*}
\end{theorem}
\bp
As in the proof of Proposition~\ref{propcaswsurfor} in Subsection~\ref{subproofcaswalsurg}, let $K$ bound a Seifert surface $\Sigma$ equipped with a symplectic basis $(x_i,y_i)_{i \in \underline{g}}$ in $\rats$.
Let $\Sigma \times \left[-1,2\right]$ be a collar of $\Sigma$ in $\rats$, and let $(A^{\prime}/A)=(A_{F}^{\prime}/A_{F})$ and $(B^{\prime}/B)=( B_{F}^{\prime}/B_{F})$ be the LP-surgeries of Subsection~\ref{sublagboun}. Let $U$ be a meridian of $K$
passing through $d \times \left[-1,2\right]$. According to Proposition~\ref{prophomafaprimef}, $(\rats(A^{\prime}/A),U)$ is diffeomorphic to $(\rats,K)$. Similarly, $(\rats(B^{\prime}/B),U)$ is diffeomorphic to $(\rats,-K)$, while $(\rats(A^{\prime}/A,B^{\prime}/B),U)$ is diffeomorphic to $(\rats,U)$.
So we have \begin{equation*}\begin{array}{ll}\Zinvlinkuf_2\bigl(\left[\rats,U;A^{\prime}/A,B^{\prime}/B\right]\bigr)&=2\Zinvlinkuf_2\left(\rats,U\right)-\Zinvlinkuf_2\left(\rats,K\right)-\Zinvlinkuf_2\left(\rats,-K\right)\\
&=\left[\Gamma_{(2)}\left(\rats,U;A^{\prime}/A,B^{\prime}/B\right)\right],\end{array}\end{equation*}
with
\begin{equation*}\left[\Gamma_{(2)}\left(\rats,U;A^{\prime}/A,B^{\prime}/B\right)\right]=\frac12 \left[\left\langle \left\langle  T\left(\CI_{AA^{\prime}}\right) \sqcup T\left(\CI_{BB^{\prime}}\right) \sqcup \twochordfreeleg \right\rangle \right\rangle\right] \in \Aavis_{2}(\sourcetl),\end{equation*}
\begin{equation*}T\left(\CI_{AA^{\prime}}\right)=\sum_{i=1}^{g} \bigtripod{$c$}{$y_i$}{$x_i$},\end{equation*}
and \begin{equation*}T\left(\CI_{BB^{\prime}}\right)=\sum_{j=1}^{g} \bigtripodrev{$c^{+}$}{$y_j^{+}$}{$x_j^{+}$},\end{equation*}
according to Lemma~\ref{lemtriplag}. We get
\begin{equation*}\left[\Gamma_{(2)}\left(\rats,U;A^{\prime}/A,B^{\prime}/B\right)\right] =  a_2(K)\left[\twochordbubble \right].\end{equation*}
Since $\Assis_2(S^1)$ is generated by the chord diagrams $\twoischord$ and $\twoxchord$, which are symmetric with respect to the orientation change on $S^1$, we have $\Zinvlinkuf_2(\rats,K)=\Zinvlinkuf_2(\rats,-K)$. So we obtain
\begin{equation*}2\Zinvlinkuf_2\left(\rats,U\right)-2\Zinvlinkuf_2\left(\rats,K\right)=a_2(K)\left[\twochordbubble \right] = 2a_2(K)\left[\tripodnonor\right].\end{equation*}
According to Example~\ref{examplecomconfinttwobis} and the multiplicativity of $\Zinvuf$ under connected sum of Theorem~\ref{thmconnsum}, we get
$\Zinvlinkuf_2\left(\rats,U\right)=\frac1{24}\left[\tripodnonor\right]$.
\eop

\begin{remark} This theorem  generalizes a result of Enore Guadagnini, Maurizio Martellini, and Mihail Mintchev in \cite{gmm} for the case $\rats=S^3$, to any rational homology sphere $\rats$. In the case of $S^3$, the known proof relies on the facts that $\Zinvlinkuf_2$ is of degree $2$ and that the space of real-valued  knot invariants of degree at most $2$ is generated by $a_2$ and a constant nonzero invariant. This uses the fact that any knot can be unknotted by crossing changes. This is no longer true in general rational homology spheres since crossing changes do not change the homotopy class. Our proof is more direct, and our result holds for any null-homologous knot in a rational homology sphere.\end{remark}

Recall the Conway weight system $w_C$ from Example~\ref{exaweightConway}, the Alexander polynomial of Definition~\ref{defAlex}, and Proposition~\ref{propchordtriv}. A \emph{long knot} $\check{K}$ of a $\QQ$-sphere $\rats(\hcylc)$ is an embedding of $\RR$ into $\crats(\hcylc)$, whose image intersects the complement of $\hcylc$, as the vertical embedding $(j_{\RR} \colon t \mapsto (0,0,t))$ does. Replace $j_{\RR} \left(\RR \setminus \left]0,1\right[\right)$ by an arc of $\crats \setminus \mathring{\hcylc}$ from $(0,0,1)$ to $(0,0,0)$, which cobounds an embedded topological disk in $\crats \setminus \mathring{\hcylc}$ with an arc of $\partial\hcylc$ with the same ends. This provides a knot $\overline{\check{K}}$, whose isotopy class is well-defined. In \cite{Leturcq3}, David Leturcq proves the following theorem. Together with the properties of the functor $\Zinvufrf$ of the third part of the book, this Leturcq theorem implies Theorem~\ref{thmcompztwo}.

\begin{theorem}[Leturcq]
\label{thmLeturcqAlex} 
For any long knot $\check{K}$ in a rational homology sphere $\crats$, such that $\overline{\check{K}}$ is null-homologous, we have the following equality in $\RR[[h]]$:
\begin{equation*}\sum_{n \in \NN} w_C \bigl(\Zinvlinkuf_n(\rats,\check{K})\bigr) h^n = \Delta_{\overline{\check{K}}}\bigl(\exp(h)\bigr).\end{equation*}
\end{theorem}

Leturcq's proof of this theorem relies on a direct computation with appropriate propagating forms. In \cite{Leturcq2}, David Leturcq obtains a similar theorem for the Bott--Cattaneo--Rossi invariants of higher dimensional knots \cite{cattaneorossi,Leturcq1}. This Leturcq theorem in higher dimensions generalizes a theorem of Tadayuki Watanabe \cite{watanabeAlexander}, who proved it for ribbon knots.

\chapter{More flexible definitions of \texorpdfstring{$\Zinvuf$}{Z} using pseudo-parallelizations}
\label{chappseudopar}

This chapter presents the \emph{pseudo-parallelizations} introduced in Section~\ref{secskproofthmmainunivlag} and defined in Section~\ref{defpseudopar}.
Our proof of the universality theorem~\ref{thmmainunivlag} uses
these generalizations of parallelizations first introduced in \cite[Section 4.3 and 4.2]{lessumgen} and furthermore studied in \cite[Section~10]{lesbetaone} and \cite[Sections 7 to~10]{lesuniveq}. 

Pseudo-parallelizations allow us to give more flexible definitions for our invariants $\Zinvuf$ and $\Zinvlinkuf$. 
In Theorem~\ref{thmfstconsttangpseudo}, we generalize the definition
of the invariant $\Zinvuf$ of Theorem~\ref{thmfstconsttang}, and thus the definition of the $q$-tangle invariant $\Zinvufrf$
of Definition~\ref{defzqtang}, by allowing 
propagating forms associated with pseudo-parallelizations.
We conclude the proof of Theorem~\ref{thmfstconsttangpseudo} in Section~\ref{secproofthmfstconsttangpseudo}. This proof involves all the previous sections.

In Chapter~\ref{chappseudoparmuchmore}, we present variants of the definition of $\Zinvufrf$ involving nonhomogeneous propagating forms or propagating chains associated with pseudo-parallelizations.
We will not use these variants in the proof of Theorem~\ref{thmmainunivlag}.

\section{Why we need pseudo-parallelizations}
\label{secwhypseudop}

This section explains why a parallelization of the exterior of a $\QQ$-handlebody $A$ does not necessarily extend to $A^{\prime}$ after a rational LP-surgery $(A^{\prime}/A)$. It justifies why I could not avoid this chapter and some of its difficulties.

\begin{notation}
\label{nottwistacrosscurve} Let $\alpha$ be a smooth map from $[-1,1]$ to $[0,2 \pi]$ such that $\alpha$ maps $[-1,-1+\varepsilon]$ to $0$ for some $\varepsilon  \in \left]0,\frac18\right[$,
it increases from $0$ to $2\pi$ on $[-1+\varepsilon, 1-\varepsilon]$, and $\alpha(-u)+\alpha(u)=2\pi$ for any $u \in [-1,1]$.

Let $\Sigma$ be a surface, orientable or not.  Let $\gamma$ be a two-sided curve properly embedded in $\Sigma$ and equipped with a collar $\gamma \times [-1,1]$ in $\Sigma$. Let $\rho_{\alpha(u)}=\rho(\alpha(u),(0,0,1))$ denote the rotation of $\RR^3$ of angle $\alpha(u)$ 
whose axis is directed by $(0,0,1)$. 
Then $\CR_{\gamma}$ denotes the map from $\Sigma$ to $SO(3)$ sending 
$\Sigma \setminus \bigl(\gamma \times [-1,1]\bigr)$ to $\id$
and satisfying
\begin{equation*}\CR_{\gamma}\bigl(c \in \gamma,u \in [-1,1]\bigr)=\rho_{\alpha(u)}.\end{equation*}
The homotopy class of $\CR_{\gamma}$ is well-defined.
\end{notation}

\begin{lemma}
 \label{lemextsurfSOthree} Let $\RR P^2_{\pi}$ denote the nonorientable submanifold of $SO(3)$ of the rotations of angle $\pi$.
 Let $\Sigma$ be a connected surface, orientable or not. Let $f\colon \partial \Sigma \to SO(3)$ be a map transverse to $\RR P^2_{\pi}$. Then $f$ extends to $\Sigma$ if and only if the $\ZZ/2\ZZ$-valued algebraic intersection
 $\langle f(\partial \Sigma),\RR P^2_{\pi}\rangle_{\!SO(3)\,}$ of $f(\partial \Sigma)$ and $\RR P^2_{\pi}$ in $SO(3)$ is zero.
\end{lemma}

\begin{remark}
 \label{rqextsurfSOthree} Recall $\pi_1(SO(3))=\ZZ/2\ZZ$ from Section~\ref{secannpont}.
The homotopy class of a loop of $SO(3)$ is determined by its $\ZZ/2\ZZ$-valued algebraic intersection with $\RR P^2_{\pi}$. This proves Lemma~\ref{lemextsurfSOthree} when $\Sigma$ is a disk.
\end{remark}

\bpo{Proof of Lemma~\ref{lemextsurfSOthree}} When $f$ extends to $\Sigma$, we may choose an extension $f_{\Sigma}$ transverse to $\RR P^2_{\pi}$. Then 
$f(\partial \Sigma) \cap \RR P^2_{\pi}$ bounds the one-manifold $f_{\Sigma}(\Sigma) \cap \RR P^2_{\pi}$. So the $\ZZ/2\ZZ$-valued intersection $\langle f(\partial \Sigma),\RR P^2_{\pi}\rangle_{\!SO(3)\,}$ is zero.

Conversely, assume that $\langle f(\partial \Sigma),\RR P^2_{\pi}\rangle_{\!SO(3)\,}$ is zero. 
Then there is a disjoint union $\gamma$ of intervals embedded in $\Sigma$  and transverse to $\partial \Sigma$ with the following properties.
\begin{itemize}
 \item The boundary $\partial \gamma$ of $\gamma$ is in $\partial \Sigma$.
 \item The connected components $K$ of $\partial \Sigma$ such that $\langle f(K),\RR P^2_{\pi}\rangle=0$ do not meet $\gamma$
 \item The connected components $K$ of $\partial \Sigma$ such that $\langle f(K),\RR P^2_{\pi}\rangle=1$ meet $\gamma$ at one point (of $\partial \gamma$).
\end{itemize}
Then the restriction to $\partial \Sigma$  of the map $\CR_{\gamma}$ associated to $(\Sigma, \gamma)$ as in Notation~\ref{nottwistacrosscurve} is homotopic to $f$. Therefore, $f$ extends to $\Sigma$.
\eop

Recall that a \emph{framed knot} is a knot equipped with a parallel (up to homotopy). Equivalently, it is a knot equipped with a normal nonzero vector field $\vec{n}_K$.
Let $K$ be framed knot in an oriented $3$-manifold $M$. Let $\vec{t}_K$ be a tangent vector of $K$ that induces the orientation of $K$.
These data induce the direct trivialization $\tau_K$ of $TM\vert_{K}$ (up to homotopy) such that $\tau_K(e_1)=\vec{t}_K$ and $\tau_K(e_2)=\vec{n}_K$.
The homotopy class of the trivialization $\tau_K$ is well-defined and does not 
depend on the orientation of $K$.

\begin{lemma}
\label{lemtaubounding} Let $K$ be a framed knot
bounding a possibly nonorientable compact surface $\Sigma$ in an oriented $3$-manifold $M$. Assume that $\Sigma$ (or, more precisely, the parallel of $K$ on $\Sigma$) induces the framing of $K$.
Let $\tau$ be a trivialization of the tangent space of $M$ over $\Sigma$. Then the restriction of $\tau$ to $K$ is not homotopic to $\tau_K$.
\end{lemma}
\bp 
Let us prove that the homotopy class of the restriction of $\tau$ to $K$ is independent of the trivialization $\tau$ of $TM$ over $\Sigma$.
Any other trivialization of $TM$ over $\Sigma$
may be written as $\tau \circ \psi_{\RR}(f)$ for a map $f$ from $\Sigma$ to $SO(3)$, with the notation of Section~\ref{secannpont}. Lemma~\ref{lemextsurfSOthree} and Remark~\ref{rqextsurfSOthree} imply that the restriction $f\vert_{K}$ of such a map $f$ is homotopically trivial. 
Thus, the homotopy class of the restriction of $\tau$ to $K$ is independent of the trivialization $\tau$ of $TM$ over $\Sigma$. It is also independent of the oriented $3$-manifold $M$ that contains $\Sigma$.
Since the tangent bundle of an oriented $3$-manifold over a possibly nonorientable closed surface is trivializable, the homotopy class of the restriction of $\tau$ to $K$ is also independent of $\Sigma$. Therefore, it is enough to prove the lemma when $\Sigma$ is a disk, and it is obvious in this case.
\eop

\begin{definition}\label{deftaubounding}
If $(K,\tau_K)$ is a framed knot in an oriented $3$-manifold $M$ and if $\tau$ is a trivialization of the restriction of $TM$ to $K$, we say that
$K$ is \emph{$\tau$-bounding}\index[T]{taubounding@$\tau$-bounding} if $\tau$ is not homotopic to $\tau_K$. (This notion is independent of the whole manifold $M$, depending only on what happens in a tubular neighborhood of $K$.)
\end{definition}

\begin{definition}
\label{deftwistacrosscurve} Let $\Sigma$ be a surface, orientable or not. Let $\gamma$ be a two-sided curve properly embedded in $\Sigma$. 
Recall the map $\CR_{\gamma}$ from $\Sigma$ to $SO(3)$ from Notation~\ref{nottwistacrosscurve}.
Define the \emph{twist map} $\CT_{\gamma}$ \emph{across} $\gamma$ to be the map from $\Sigma \times \RR^3$ to itself such that
\begin{equation*}\CT_{\gamma}(y \in \Sigma;\cvarM\in \RR^3)=\bigl(y;\CR_{\gamma}(y)(\cvarM)\bigr).\end{equation*}
Assume that $\Sigma$ is embedded in an oriented $3$-manifold $M$.
For a trivialization $\tau$ of $TM\vert_{\Sigma}$, the \emph{twist of $\tau\vert_{\Sigma}$ across $\gamma$} is the trivialization $\tau\vert_{\Sigma} \circ \CT_{\gamma}$ (defined up to homotopy).
\end{definition}

Let $A$ be a compact oriented connected $3$-manifold 
with boundary $\partial A$.
Define the {\em $\ZZ/2\ZZ$-Lagrangian\/} $\CL_A^{\ZZ/2\ZZ}$ of $A$ \index[N]{Lagrangians!LAGZ2@$\CL_A^{\ZZ/2\ZZ}$} to be the kernel
\begin{equation*}\CL_A^{\ZZ/2\ZZ}=\Ker\bigl(H_1(\partial A;\ZZ/2\ZZ) \longrightarrow
H_1( A;\ZZ/2\ZZ)\bigr)\end{equation*}
of the map induced by the inclusion map.
This is a Lagrangian subspace of $(H_1(\partial A;\ZZ/2\ZZ);\langle ., . \rangle )$.

 The curves of an oriented compact surface $F$ embedded in an oriented $3$-manifold are naturally framed by the surface $F$: They are framed 
 by a nonzero normal vector field tangent to $F$.
\begin{proposition}
\label{propuseless}
Let $\partial A$ be a connected oriented compact surface.
Let $\tau$ be a trivialization of \/$T(\partial A \times [-2,2])$.
Then there exists a unique map 
\begin{equation*}\phi_{\tau} \colon H_1(\partial A;\ZZ/2\ZZ) \longrightarrow \frac{\ZZ}{2\ZZ}\end{equation*}
such that 
\begin{enumerate}
\item for a connected curve $x$ of $\partial A=\partial A \times \{0\}$, we have
$\phi_{\tau}(x)=0$ if and only if $x$ (framed by $\partial A$) is $\tau$-bounding
and,
\item for any pair $(x,y)$ of elements of $H_1(\partial A;\ZZ/2\ZZ)$, we have
\begin{equation*}\phi_{\tau}(x+y)=\phi_{\tau}(x) + \phi_{\tau}(y) + \langle x, y \rangle_{\!\partial A\,} .\end{equation*}
\end{enumerate}
The map $\phi_{\tau}$ satisfies the following properties.
\begin{itemize}
\item 
Let $x$ be a disjoint union of curves in $\partial A$. Assume that $\partial A \times [-2,2]$ is embedded in an oriented $3$-manifold $M$, where $x$ bounds a connected surface $\Sigma$, orientable or not. 
If $\Sigma$ and $\partial A$ induce the same framing of $x$, then
$\tau\vert_{x}$ extends to $\Sigma$ as a trivialization of $TM\vert_{\Sigma}$ if and only if $\phi_{\tau}(x)$ is zero.
\item Let $c$ be curve of $\partial A$. Let $\twistriv_c$ denote the twist across $c$ of Definition~\ref{deftwistacrosscurve}. For any $x \in H_1(\partial A;\ZZ/2\ZZ)$, we have
\begin{equation*}\phi_{ \tau \circ \twistriv_c}(x)=\phi_{\tau}(x) + \langle x, c \rangle_{\!\partial A\,}.\end{equation*}
\item When $A$ is a compact oriented connected $3$-manifold 
with boundary $\partial A$, the trivialization $\tau$ extends as a trivialization over $A$ if and only if
$\phi_{\tau}(\CL_A^{\ZZ/2\ZZ})=\{0\}$.
\end{itemize}
\end{proposition}
\bp For a disjoint union $x=\sqcup_{i=1}^nx_i$ of connected curves $x_i$ on $\partial A$, define $\phi_{\tau}\left(x\right)=\sum_{i=1}^n \phi_{\tau}(x_i)$ in $\ZZ/2\ZZ$, with $\phi_{\tau}(x_i)=0$ if $x_i$ (framed by $\partial A$) is $\tau$-bounding, and $\phi_{\tau}(x_i)=1$ otherwise. 

Assume that such a framed disjoint union $x$ bounds a connected surface $\Sigma$ in an oriented $3$-manifold $M$. Also assume that $\Sigma$ and $\partial A$ induce the same framing of $x$.
Let us prove that $\tau\vert_{x}$ extends to $\Sigma$ as a trivialization of $TM\vert_{\Sigma}$ if and only if $\phi_{\tau}(x)$ is zero.

When $\phi_{\tau}(x)=0$, group all the curves $x_i$ such that $\phi_{\tau}(x_i)=1$ by pairs. Make each such pair bound an annulus that induces the framing. Make each curve $x_i$ such that $\phi_{\tau}(x_i)=0$ bound a disk which induces the framing.
Let $\hat{\Sigma}$ denote the union of $\Sigma$ with the above disks and annuli. 
The restriction to $\hat{\Sigma}$ of the tangent bundle of an oriented $3$-manifold $M$ in which $\hat{\Sigma}$ embeds is independent of $M$. It admits a trivialization $\hat{\tau}$. Extend $\tau$ to $\hat{\Sigma} \setminus \mathring{\Sigma}$. Write the restriction of $\hat{\tau}$ to $(\hat{\Sigma} \setminus \mathring{\Sigma}) \times \RR^3$ as $\hat{\tau}=\tau \circ \psi_{\RR}(f)$ for a map $f$ from $\hat{\Sigma} \setminus \mathring{\Sigma}$ to $SO(3)$. According to Lemma~\ref{lemextsurfSOthree}, the intersection of $f(x)$ and $\RR P^2_{\pi}$ is zero, and 
the map $f$ extends to $\hat{\Sigma}$. Thus $\tau = \hat{\tau}\circ \psi_{\RR}(f)^{-1}$ also extends to $\Sigma$.

When $\phi_{\tau}(x)=1$, assume that $\phi_{\tau}(x_1)=1$ without loss of generality. Group the other curves $x_i$ such that $\phi_{\tau}(x_i)=1$ pairwise, make them cobound a disjoint union of annuli, and make the curves $x_i$ such that $\phi_{\tau}(x_i)=0$ bound disks, in a framed way as above. Let $\hat{\Sigma}$ denote the union of $\Sigma$ with these disks and these annuli.
The boundary of $\hat{\Sigma}$ is $x_1$.
The trivialization $\tau$ still extends to $\hat{\Sigma} \setminus \mathring{\Sigma}$. If $\tau$ extends to $\Sigma$, then it extends to $\hat{\Sigma}$, and $x_1$ is $\tau$-bounding. So $\phi_{\tau}(x_1)=0$, which is absurd. Therefore, $\tau$ does not extend to $\Sigma$.

Let us prove that our above definition of $\phi_{\tau}(x)$ depends only on the class of $x$ in $H_1(\partial A;\ZZ/2\ZZ)$.
Let $x$ be an embedded (possibly nonconnected) curve in $\partial A$.
Let $y$ be another such curve in $\partial A \times \{-1\}$ homologous to $x$ modulo $2H_1(\partial A;\ZZ)$. Then there exists a framed (possibly nonorientable) cobordism between $x$ and $y$ in $\partial A \times [-1,1]$, and it is easy to see that 
$\phi_{\tau}(x)=0$ if and only if $\phi_{\tau}(y)=0$.

Let us check that $\phi_{\tau}$ behaves as predicted under addition.
Because we are dealing with elements of $H_1(\partial A;\ZZ/2\ZZ)$, we can
consider representatives $x$ and $y$ of $x$ and $y$ 
such that $x$ is connected and intersects $y$ at most once. Next, the known additivity under disjoint union reduces the proof to the case in which $x$ and $y$ are connected and $x$ and $y$ intersect once.
Note that both sides of the equality to be proved vary in the same way under trivialization changes. 
Consider the punctured torus neighborhood of $x \cup y$ and
a trivialization $\tau$ that restricts to the punctured torus as the direct sum of a parallelization of the torus and the normal vector to $\partial A$. Then 
we have $\phi_{\tau}(x+y)=\phi_{\tau}(x)= \phi_{\tau}(y)=1$.
We leave the last two assertions to the reader.
\eop

\begin{example}
\label{exabadlag}
For any $\QQ$-handlebody $A$, there exists a Lagrangian subspace $\CL^{\ZZ}$ of $(H_1(\partial A;\ZZ);\langle ., . \rangle )$, such that
$\CL_A=\CL^{\ZZ} \otimes \QQ$. The following example shows that
$\CL_A^{\ZZ/2\ZZ}$ \index[N]{Lagrangians!LAGZ2@$\CL_A^{\ZZ/2\ZZ}$} may differ from $\CL^{\ZZ} \otimes {\ZZ/2\ZZ}$.

Let $\CM$ be a M\"obius band embedded in the interior of a solid torus $D^2 \times S^1$ so that
the core of the solid torus is the core of $\CM$. 
Embed $D^2 \times S^1$ into $S^2 \times S^1=D^2 \times S^1 \cup_{\partial D^2 \times S^1} (-D^2 \times S^1)$ as the first copy. Orient the knot $\partial \CM$ so that $\partial \CM$ pierces twice $S^2 \times 1$ positively.
Let $m$ be the meridian of $\partial \CM$. Let $\ell$ be the parallel of $\partial \CM$ induced by $\CM$.
Let $A$ be the exterior of the knot $\partial \CM$ in 
$S^2 \times S^1$.
The reader can check that $A$ is a $\QQ$-handlebody, as an exercise. Observe
$\CL_A^{\ZZ}=\ZZ[2m]$, $\CL_A=\QQ[m]$, and $\CL_A^{\ZZ/2\ZZ}=\ZZ/2\ZZ[\ell]$.

Let $A_{0}$ be the solid torus with boundary $\partial A = \partial A_0$, where $m$ bounds a disk. We have $\CL_{A_0}=\QQ[m]$. Equip $A_{0}$ with a parallelization $\tau$ such that $\phi_{\tau}(\ell)=1$.
According to Proposition~\ref{propuseless}, the restriction to $\partial A$ of $\tau$ does not extend to $A$.
\end{example}

\section{Definition of pseudo-parallelizations}
\label{defpseudopar}

\begin{definition}
\label{defpseudotriv} Recall $\varepsilon$ and the maps $\alpha$ and $\rho_{\alpha(u)}$ from Notation~\ref{nottwistacrosscurve}. Set $N\left(\partial [-1,1]\right)=[-1,1]\setminus \left]-1+\varepsilon, 1-\varepsilon\right[$.
A {\em pseudo-parallelization\/} $\tilde{\tau}=(N(\gamma);\tau_e,\tau_b)$ of an oriented $3$-manifold $\Aman$ with possible boundary consists of 
\begin{itemize}
\item a framed link $\gamma$ of the interior of $\Aman$, which will be called \emph{the link of the pseudo-parallelization} $\tilde{\tau}$, equipped with a neighborhood $N(\gamma)=[a,b] \times \gamma \times [-1,1]$, for real numbers $a$ and $b$ such that $a<b$ and $\varepsilon < \frac{b-a}{4}$,
\item a parallelization $\tau_e$ of 
$\Aman$ outside $\left]a+\varepsilon,b-\varepsilon\right[ \times \gamma \times \left]-1+\varepsilon,1-\varepsilon\right[$,
\item a parallelization $\tau_b \colon N(\gamma) \times \RR^3 \rightarrow TN(\gamma)$ of $N(\gamma)$ such that
\begin{equation*}\tau_b=\left\{\begin{array}{ll} \tau_e & \mbox{over } [b-\varepsilon,b] \times \gamma \times [-1,1]\mbox{ and }[a,b] \times \gamma \times N\left(\partial [-1,1]\right) \\
 \tau_e \circ  \CT_{\gamma}& \mbox{over }  [a,a+\varepsilon] \times \gamma \times [-1,1],
\end{array} \right.\end{equation*}
\end{itemize}
where 
\begin{equation*}\CT_{\gamma}\bigl(t,c \in \gamma,u \in [-1,1];\cvarM\in \RR^3\bigr)=\bigl(t,c,u;\rho_{\alpha(u)}(\cvarM)\bigr).\end{equation*}
\end{definition}

\begin{lemma}
\label{lempseudoparextend}
 Let $A$ be a compact oriented $3$-manifold and let $\tau$ be a trivialization of $TA$ defined on a collar $[-4,0]\times \partial A$ of $\partial A(=\{0\} \times \partial A)$.
Then there is a pseudo-parallelization of $A$ that extends the restriction of $\tau$ to $\left[-1,0\right]\times \partial A$.
\end{lemma}
\bp
There exists a trivialization $\tau^{\prime}$ of $TA$ on $A$.
After a homotopy of $\tau$ around $\{-2\} \times \partial A$, there exists a union $\gamma \times [-1,1]$ of annuli of $\{-2\} \times \partial A$ such that $\tau = \tau^{\prime} \circ \CT_{\gamma}$ on $\{-2\} \times \partial A$. (When $\partial A$ is connected, $\gamma$ can be assumed to be connected, too.)
Consider the neighborhood $N(\gamma)=[-2,-1]\times \gamma \times[-1,1]$ of $\gamma$. Define $\tau_e$ to coincide with $\tau$ on $([-2,0]\times \partial A) \setminus \Int(N(\gamma))$, and with $\tau^{\prime}$ on $A \setminus (\left]-2,0\right]\times \partial A)$.
Define $\tau_b$ to coincide with $\tau$ on $N(\gamma)$.
\eop

\begin{definition}
\label{defpseudotrivpone}[Trivialization $\tilde{\tau}_{\CC}$ of $T\Aman \otimes_{\RR} \CC$]
Define a smooth map
\begin{equation*}F_U\colon  [a,b] \times [-1,1]  \longrightarrow SU(3)\end{equation*} such that
\begin{equation*}F_U(t,u) = \left\{\begin{array}{ll}\id&
 \mbox{if}\; |u|>1-\varepsilon\\
\rho_{\alpha(u)} &  \mbox{if}\; t<a+\varepsilon\\
\id &  \mbox{if}\; t>b-\varepsilon.\end{array}\right.\end{equation*}
Since $\pi_1(SU(3))$ is trivial, it is possible to define such a smooth map.
Define a trivialization $\tilde{\tau}_{\CC}$ of $T\Aman \otimes_{\RR} \CC$, associated with the pseudo-parallelization $\tilde{\tau}$ of Definition~\ref{defpseudotriv}, as follows.
\begin{itemize}
\item On $\bigl(\Aman\setminus N(\gamma)\bigr) \times \CC^3$, we have $\tilde{\tau}_{\CC} =\tau_e \otimes 1_{\CC}$,
\item Over $[a,b] \times \gamma \times [-1,1]$,  we have
$\tilde{\tau}_{\CC}(t,c,u;\cvarM) =\tau_b\bigl(t,c,u;F_U(t,u)^{-1}(\cvarM)\bigr)$.
\end{itemize}
Since $\pi_2(SU(3))$ is trivial, the homotopy class of $\tilde{\tau}_{\CC}$
is well-defined.
\end{definition}

\begin{definition}
\label{defponepseudotrivpone}Let ${\manifm}_0$ and ${\manifm}_1$ 
be two compact connected oriented $3$-manifolds
whose boundaries $\partial {\manifm}_0$ and $\partial {\manifm}_1$ have collars identified by a diffeomorphism.
Let $\tau_0$ be a pseudo-parallelization of ${\manifm}_0$, which restricts to a collar neighborhood of $\partial {\manifm}_0$ as a genuine trivialization.
Let $\tau_1$ be a pseudo-parallelization of ${\manifm}_1$ that coincides with $\tau_0$ on this collar neighborhood.
We use the definition of Proposition~\ref{proppontdef} of relative Pontrjagin numbers and define $p_1(\tau_0,\tau_1)$ to be $p_1(\tau_{0,\CC},\tau_{1,\CC})$.
Let $\crats$ be a rational homology $\RR^3$. A pseudo-parallelization $\tilde{\tau}$ of $\crats$ is \emph{asymptotically standard} if it coincides with the standard parallelization $\taust$ of $\RR^3$ outside $\ballb_{\rats}$. (Recall Definition~\ref{defparasyst}.) For such an asymptotically standard pseudo-parallelization, set 
$p_1(\tilde{\tau})=p_1\left((\taust)\vert_{B^3},\tilde{\tau}\vert_{\ballb_{\rats}}\right)$.
\end{definition}

\begin{definition}
\label{defhombounformpseudo}[Homogeneous boundary form associated with $\tilde{\tau}$]
Let $\tilde{\tau}=(N(\gamma);\tau_e,\tau_b)$ be a pseudo-parallelization of a  $3$-manifold $\Aman$.
 Recall $\varepsilon$ and the map $\alpha$ from Notation~\ref{nottwistacrosscurve} and Definition~\ref{defpseudotriv}. Define a smooth map 
\begin{equation*}F \colon [a,b] \times [-1,1]  \longrightarrow  SO(3)\end{equation*}
such that
\begin{equation*}F(t,u)  =  \left\{\begin{array}{ll}\id&
 \mbox{if}\; |u|>1-\varepsilon\\
\rho_{\alpha(u)} &  \mbox{if}\; t<a+\varepsilon\\
\rho_{-\alpha(u)} &  \mbox{if}\; t>b-\varepsilon .\end{array}\right.\end{equation*}
Since the restriction of $F$ to the boundary of $[a,b] \times [-1,1]$ is trivial in $\pi_1(SO(3))$, it is possible to define such a smooth map $F$.

Let $p_{\tau_b}=p(\tau_b)$ denote the projection from $\ST N(\gamma)$ to $S^2$ induced by $\tau_b$. We have $p_{\tau_b}(\tau_b(t,c,u;\cvarM \in S^2))=\cvarM$.
Define
\begin{equation*}\begin{array}{llll} {F}(\gamma,\tau_b) \colon &  [a,b] \times   \gamma \times [-1,1] \times S^2 & \longrightarrow & [a,b] \times   \gamma \times [-1,1] \times S^2\\
&(t,c,u; \svarM) & \mapsto & \bigl(t,c,u; F(t,u)(\svarM)\bigr).\end{array}\end{equation*}

Define the closed two-form \index[N]{omegazctaub@$\omega(\gamma,\tau_b)$ $2$-form} $\omega(\gamma,\tau_b)$ on $\ST\left([a,b] \times \gamma \times [-1,1]\right)$ to be
\begin{equation*}\omega(\gamma,\tau_b)=\frac{p\bigl(\tau_b \circ \CT_{\gamma}^{-1}\bigr)^{\ast}(\omega_{S^2}) +p\bigl(\tau_b \circ {F}(\gamma,\tau_b)^{-1}\bigr)^{\ast}(\omega_{S^2})}2.\end{equation*}

The \emph{homogeneous boundary form} associated to $(\tilde{\tau},F)$ is the following closed $2$-form $\omega(\tilde{\tau},F)$ on $\ST \Aman$.
\begin{equation*}\omega(\tilde{\tau},F)=\left\{ \begin{array}{ll}p_{\tau_e}^{\ast}\left(\omega_{S^2}\right)& \mbox{on}\; \ST\bigl(\Aman \setminus N(\gamma)\bigr)\\
\omega(\gamma,\tau_b)
&\mbox{on}\;\ST\bigl(N(\gamma)\bigr).
\end{array}\right.\end{equation*}
A \emph{homogeneous boundary form} of $(\ST \Aman,\tilde{\tau})$ is a homogeneous boundary form associated to $(\tilde{\tau},F)$ for some $F$ as above.
\end{definition}

We will justify the consistency of Definition~\ref{defhombounformpseudo} by applying the following lemma with the constant map $\kappa$ with value one. We will use the general lemma in Lemma~\ref{lemappprimtwist}.

\begin{lemma}
\label{lemprimtwist}
Let $(e_1=(1,0,0),e_2=(0,1,0),e_3=(0,0,1))$ denote the standard basis 
of $\RR^3$. 
Let $v_i:\RR^3 \longrightarrow \RR$ denote the 
$i$th coordinate with respect to this basis.
Let $\rho_{\theta}=\rho_{\theta,e_3}$ \index[N]{Rtheta@$\rho_{\theta}$} denote the rotation of $\RR^3$ of angle $\theta$
whose axis is directed by $e_3$.
For $k \in \ZZ$, define \index[N]{Tcal@$\CT_k$ twist} \begin{equation*}\begin{array}{llll}\CT_k:&\RR \times S^2 &\longrightarrow &S^2\\
& (\theta,\vecx) & \mapsto & \rho_{k\theta}(\vecx). \end{array}\end{equation*}
Let $\kappa \colon [-1,1] \to \RR$ be a smooth map. Consider the associated volume form $(\kappa \circ v_3)\omega_{S^2}$ on $S^2$. It is invariant under the rotations $\rho_{\theta}$.
Then we have
\begin{equation*}\CT_k^{\ast}\bigl((\kappa \circ v_3)\omega_{S^2}\bigr)=\CT_0^{\ast}\bigl((\kappa \circ v_3)\omega_{S^2}\bigr) 
+ \frac{k(\kappa \circ v_3)}{4\pi} d\theta \wedge dv_3\end{equation*}
\end{lemma}
\bp Recall the homogeneous two-form $\omega_{S^2}$ on $S^2$
with total area $1$. 
When $\vecx \in S^2$, 
and when $v$ and $w$ are two tangent vectors of $S^2$ at $\vecx$, we have
\begin{equation*}\omega_{S^2}(v \wedge w)=\frac{1}{4\pi}\det(\vecx,v,w),\end{equation*}
where $\vecx \wedge v \wedge w = 
\det(\vecx,v,w) e_1 \wedge e_2 \wedge e_3$ in $\bigwedge^3\RR^3$.

Since $\rho_{\theta}$ preserves the area in $S^2$ and leaves $v_3$ invariant, the restrictions of 
$\CT_k^{\ast}\left((\kappa \circ v_3)\omega_{S^2}\right)$ and $\CT_0^{\ast}\left((\kappa \circ v_3)\omega_{S^2}\right)$ coincide on 
$\bigwedge^2T_{(\theta,\vecx)}(\{\theta \}\times S^2)$.
Therefore, we are left with the computation of 
\begin{equation*}\Bigl(\CT_k^{\ast}\bigl((\kappa \circ v_3)\omega_{S^2}\bigr)-\CT_0^{\ast}\bigl((\kappa \circ v_3)\omega_{S^2}\bigr)\Bigr)(u \wedge v)\end{equation*}
when $u \in T_{(\theta,\vecx)}(\RR \times \{\vecx \})$
and $v \in T_{(\theta,\vecx)}(\{\theta \}\times S^2)$. We have
\begin{equation*}\CT_0^{\ast}\bigl((\kappa \circ v_3)\omega_{S^2}\bigr)(u \wedge v)=0,\end{equation*}
\begin{equation*}\CT_k^{\ast}\bigl((\kappa \circ v_3)\omega_{S^2}\bigr)_{(\theta,\vecx)}(u \wedge v)=\frac{\kappa \circ v_3(\vecx)}{4\pi}\det\bigl(\rho_{k\theta}(\vecx),
T_{(\theta,\vecx)}\CT_k(u),T_{(\theta,\vecx)}\CT_k(v)\bigr),\end{equation*}
and $T_{(\theta,\vecx)}\CT_k(v)=\rho_{k\theta}(v)$.
Since $\rho_{k\theta}$ preserves the volume in $\RR^3$, we get
\begin{equation*}\CT_k^{\ast}\bigl((\kappa \circ v_3)\omega_{S^2}\bigr)_{(\theta,\vecx)}(u \wedge v)=\frac{\kappa \circ v_3(\vecx)}{4\pi}\det\bigl(\vecx,
\rho_{-k\theta}(T_{(\theta,\vecx)}\CT_k(u)),v\bigr).\end{equation*}
Let $\vecx_i$ stand for $v_i(\vecx)$. We have
\begin{equation*}T_{(\theta,\vecx)}\CT_k(u)=kd\theta(u)\rho_{k\theta+\pi/2}(\vecx_1 e_1 +\vecx_2 e_2).\end{equation*}
We obtain
\begin{equation*}\CT_k^{\ast}\left(\omega_{S^2}\right)_{(\theta,\vecx)}(u \wedge v)=\frac{kd\theta(u)}{4\pi}\det(\vecx,
-\vecx_2 e_1 +\vecx_1 e_2,v),\end{equation*} and therefore
\begin{equation*}\begin{array}{ll}\CT_k^{\ast}(\omega_{S^2})(u \wedge .)
&=\frac{kd\theta(u)}{4\pi}
\det\left(\begin{array}{ccc}\vecx_1 & -\vecx_2 & dv_1\\
\vecx_2 & \vecx_1 & dv_2\\
\vecx_3 & 0 & dv_3\end{array}\right)\\
&=\frac{kd\theta(u)}{4\pi}
\bigl(-\vecx_3\vecx_1dv_1 -\vecx_3\vecx_2 dv_2+
 (1-\vecx_3^2)dv_3\bigr)\\
&=\frac{kd\theta(u)}{4\pi} dv_3.\end{array}\end{equation*}
\eop

\bpo{Proof of the consistency of Definition~\ref{defhombounformpseudo}} It suffices to prove 
\begin{equation*}p\left(\tau_b \circ \CT_{\gamma}^{-1}\right)^{\ast}\left(\omega_{S^2}\right) +p\left(\tau_b \circ \CT_{\gamma}\right)^{\ast}(\omega_{S^2})=2p\left(\tau_b\right)^{\ast}\left(\omega_{S^2}\right)\end{equation*}
on $\ST\left([b-\varepsilon,b] \times \gamma \times [-1,1]\right)$, where we have
 \begin{equation*}p_{\tau_b\circ \CT_{\gamma}^{\pm 1}}\bigl(\tau_b(t,c,u;\cvarM)\bigr)= p_{\tau_b\circ \CT_{\gamma}^{\pm1}}\Bigl(\tau_b\circ \CT_{\gamma}^{\pm1}\bigl(t,c,u;\rho_{\mp \alpha(u)}(\cvarM)\bigr)\Bigr)=\rho_{\mp \alpha(u)}(\cvarM) .\end{equation*}
 Set $\tilde{p}_{\tau_b}=p_{[-1,1]} \times {p}_{\tau_b} \colon \ST\bigl([b-\varepsilon,b] \times \gamma \times [-1,1]\bigr) \to [-1,1] \times S^2$. We have
 \begin{equation*}p\left(\tau_b\right) =\CT_{0} \circ (\alpha \times \id_{S^2}) \circ \tilde{p}_{\tau_b}\mbox{ and }p\left(\tau_b \circ \CT_{\gamma}^{\pm 1}\right) =\CT_{\mp 1} \circ (\alpha \times \id_{S^2}) \circ \tilde{p}_{\tau_b}.\end{equation*}
We get $p\left(\tau_b \circ \CT_{\gamma}^{\pm 1}\right)^{\ast}\left(\omega_{S^2}\right) =\left( (\alpha \times \id_{S^2}) \circ \tilde{p}_{\tau_b}\right)^{\ast} \left( \CT_{\mp 1}^{\ast}\left(\omega_{S^2}\right)\right) $. \\
Thus, Lemma~\ref{lemprimtwist}
implies that Definition~\ref{defhombounformpseudo} is consistent. \eop

\begin{definition}
\label{defpropagatorpseudohom}
Let $\crats$ be a rational homology $\RR^3$, equipped with an asymptotically standard pseudo-parallelization $\tilde{\tau}$.
A \emph{homogeneous propagating form}\index[T]{homogeneous!propagating form} of $(C_2(\rats),\tilde{\tau})$ is a propagating form of $C_2(\rats)$ (as in Definition~\ref{defpropagatortwo}) that coincides with a homogeneous boundary form of $(\ST \crats,\tilde{\tau})$ as in Definition~\ref{defhombounformpseudo} on $\ST \crats$.
\end{definition}

\begin{lemma}
\label{lemexistpropagatorpseudohom}
Such homogeneous propagating forms exist for any $(\rats,\tilde{\tau})$.
\end{lemma}
\bp See Section~\ref{secprop}.
\eop

As in Definition~\ref{defparacyl},
a \emph{pseudo-parallelization} of a rational homology cylinder $\hcylc$ is a pseudo-parallelization of $\crats(\hcylc)$ that agrees with the standard parallelization of $\RR^3$ outside $\hcylc$.

The main result of this chapter is the following theorem.

\begin{theorem}

\label{thmfstconsttangpseudo}
Let $\hcylc$ be a rational homology cylinder. Let $\tau=(N(\gamma);\tau_e,\tau_b)$ be a pseudo-parallelization of $\hcylc$.
Let $\tanghcyll \colon \sourcetl \hookrightarrow \crats(\hcylc) \setminus N(\gamma)$ be a long tangle representative in $\crats(\hcylc) \setminus N(\gamma)$.

Definition~\ref{defIthetalong} of $I_{\theta}(K,\partau)$ naturally extends for such a pseudo-parallelization $\tau$ when $K$ is a component of $\tanghcyll$. 

With this extended definition of $I_{\theta}$ and with Definition~\ref{defponepseudotrivpone} of $p_1(\tau)$,
Theorem~\ref{thmfstconst} and Theorem~\ref{thmfstconsttang} also hold when $\tau$ is a pseudo-parallelization $\tau=(N(\gamma);\tau_e,\tau_b)$ of $\hcylc$ such that 
$N(\gamma)$ does not meet the image of the long tangle representative (or the link) $\tanghcyll \colon \sourcetl \hookrightarrow \crats(\hcylc)$.
\end{theorem}

In order to prove Theorem~\ref{thmfstconsttangpseudo} in Section~\ref{secproofthmfstconsttangpseudo}, we prove some preliminary lemmas in the next sections.

\section{Integration of homogeneous propagating forms along surfaces}

\begin{definition}
\label{defeulernumb}
Let $\Sigma$ be a compact oriented surface with boundary. Let $\Sigma_0$ denote the image of the zero section in the tangent bundle $T\Sigma$ of $\Sigma$. Let $\vecx$ be a nowhere vanishing section of $T\Sigma$ along the boundary of $\Sigma$. 
Let $\tilde{\vecx}$ be an extension of $\vecx$ over $\Sigma$ whose image $\tilde{\vecx}(\Sigma)$ in $T \Sigma$ is transverse to $\Sigma_0$.
The sections $\Sigma_0$ and $\tilde{\vecx}(\Sigma)$ are naturally oriented by $\Sigma$.
The \emph{relative Euler number}\index[T]{relative!Euler number} $\chi(\vecx;\Sigma)$ is their algebraic intersection $\langle \tilde{\vecx}(\Sigma), \Sigma_0 \rangle_{\!T \Sigma\,}$ in $T \Sigma$.
\end{definition}

Note that this definition makes sense since all the extensions of $\vecx$ are homotopic relatively to $\partial \Sigma$.
This Euler number
is an obstruction to extending $\vecx$ over $\Sigma$ as a nowhere vanishing section of $T \Sigma$.
Here are some other well-known properties of this number.

\begin{lemma}
\label{lemptyeulernumb} Let $\Sigma$ be a compact oriented surface with boundary, and let $\vecx$ be a section of $\ST \Sigma$ along the boundary of $\Sigma$.
\begin{itemize}
\item If $\Sigma$ is connected and if $\chi(\vecx;\Sigma)=0$, then $\vecx$ extends as a nowhere vanishing section of  $T \Sigma$.
\item If $\vecx$ is tangent to the boundary of $\Sigma$, then $\chi(\vecx;\Sigma)$ is the Euler characteristic $\chi(\Sigma)$ of $\Sigma$.
\item More generally,
let \/$a^{(1)}$, \dots, $a^{(k)}$ denote the \/$k$ connected components of the boundary \/$\partial \Sigma$ of $\Sigma$. 
For \/$i=1,\dots,k$, the unit bundle $\ST\Sigma\vert_{a^{(i)}}$ of $T\Sigma\vert_{a^{(i)}}$
is an \/$S^1$-bundle over \/$a^{(i)}$ with a canonical trivialization induced by \/$Ta^{(i)}$. Let \/$d(\vecx,a^{(i)})$ be the degree of the projection on the fiber \/$S^1$ of this bundle of the section \/$\vecx$, with respect to this canonical trivialization. Then we have \begin{equation*}\chi(\vecx;\Sigma)=\sum_{i=1}^k d(\vecx,a^{(i)}) + \chi(\Sigma).\end{equation*}
\end{itemize}
\end{lemma}
\bp First observe all these properties when $\Sigma$ is a disk.
When $\Sigma$ is connected, there is a disk $D$ that contains all the zeros of an arbitrary generic extension of $\tilde{\vecx}$ of $\vecx$. If $\chi(\vecx;\Sigma)=0$, then $\chi(\tilde{\vecx}\vert_{\partial D};D)=0$, and $\tilde{\vecx}\vert_{\partial D}$ extends to $D$ as a nowhere vanishing section. So $\vecx$ extends to $\Sigma$ as a nowhere vanishing section.

Let $\ST^+\partial \Sigma$ denote the unit vector field of $\partial \Sigma$ that is tangent to $\partial \Sigma$ and induces its orientation. Let us prove that the following equality $(\ast(\Sigma))$
\begin{equation*} \chi(\ST^+\partial \Sigma;\Sigma)=\chi(\Sigma)\end{equation*}
holds for a general $\Sigma$.

For $i=1,2$, let $\Sigma_i$ be a compact oriented surface, and let $c_i$ be a connected component of $\partial \Sigma_i$. Set $\Sigma=\Sigma_1 \cup_{c_1 \sim -c_2}\Sigma_2$. 
Since the section $\ST^+ c_1$ is homotopic to $(-\ST^+ c_1)$ as a section of $\ST \Sigma$, we have
\begin{equation*}\chi(\ST^+\partial \Sigma;\Sigma)=\chi(\ST^+\partial \Sigma_1;\Sigma_1)+\chi(\ST^+\partial \Sigma_2;\Sigma_2).\end{equation*}
Assume that $(\ast(\Sigma_1))$ holds.
Since $\chi(\Sigma)=\chi(\Sigma_1)+\chi(\Sigma_2)$,
the equalities $(\ast(\Sigma_2))$ and $(\ast(\Sigma))$ are equivalent.
Since $S^1 \times S^1$ is parallelizable, the equality $(\ast(S^1 \times S^1))$ holds.
So $(\ast(S^1 \times S^1\setminus \mathring{D}^2))$ holds. The general case follows easily.

The third property of $\chi(\vecx;\Sigma)$ is an easy consequence of the previous one.
\eop

\begin{lemma}
\label{lemsplusstau}
Recall the vectors $e_2=(0,1,0)$ and $e_3=(0,0,1)$ of $\RR^3$.
Let $\Sigma$ be a compact oriented surface immersed in a \/$3$-manifold $M$ equipped with a parallelization $\tau$. Assume that $\tau(. \times e_3)$ is a positive normal to $\Sigma$ along $\partial \Sigma$.
Let \/$s_+(\Sigma) \subset \ST M$ (resp. $s_-(\Sigma) \subset \ST M$)  be the graph of the section of \/$\ST M\vert_{\Sigma}$ in \/$\ST M$ associated to the positive (resp. negative) normal to \/$\Sigma$. Let $s_{\tau}(\Sigma; e_3)$ be the graph of the section $\tau(\Sigma \times \{e_3\})$.
Then the cycles \begin{equation*}2\bigl(s_+(\Sigma)-s_{\tau}(\Sigma; e_3)\bigr) - \chi\bigl(\tau(. \times e_2)\vert_{\partial \Sigma};\Sigma\bigr)\ST M\vert_{\ast} \end{equation*}
and \begin{equation*}2\bigl(s_-(\Sigma)-s_{\tau}(\Sigma; -e_3)\bigr) + \chi\bigl(\tau(. \times e_2)\vert_{\partial \Sigma};\Sigma\bigr)\ST M\vert_{\ast} \end{equation*} of \/$\ST M\vert_{\Sigma}$ are null-homologous in \/$\ST M\vert_{\Sigma}$.
\end{lemma}
\bp Consider the involution $\iota_{\Sigma}$ of $\ST M\vert_{\Sigma}$ that sends a vector to its opposite. This involution reverses the orientation of a fiber $\ST M\vert_{\ast}$. Therefore, it sends 
the second cycle (with $s_-(\Sigma)$) to the first one (with $s_+(\Sigma)$).
Thus, it suffices to prove that \begin{equation*}2\bigl(s_+(\Sigma)-s_{\tau}(\Sigma; e_3)\bigr) - \chi\bigl(\tau(. \times e_2)\vert_{\partial \Sigma};\Sigma\bigr)\ST M\vert_{\ast} \end{equation*} is a null-homologous cycle in \/$\ST M\vert_{\Sigma}$. 
The trivialization $\tau$ can be homotoped so that $\tau(. \times e_3)$ is a positive normal to $\Sigma$ over a one-skeleton of $\Sigma$. Therefore, we can assume that $\tau(. \times e_3)$ is a positive normal to $\Sigma$ over the complement of a disjoint union of disks embedded in the interior of $\Sigma$.
So it suffices to prove that
$2\left(s_+(D)-s_{\tau}(D; e_3)\right) - \chi\left(\tau(. \times e_2)\vert_{\partial D};D\right)\ST M\vert_{\ast} $
is null-homologous, for such a disk $D$.
Let $\tau_0$ be a trivialization of $\ST M\vert_{D}$ such that $\tau_0(.,e_3)$ is the positive normal to $D$ and 
$\tau_0(.,e_2)$ is tangent to $\Sigma$. The trivialization $\tau_0$ identifies $\ST M\vert_{D}$ with $D \times S^2$, and we have \begin{equation*}\bigl[s_+(D)\bigr]-\bigl[s_{\tau}(D;e_3)\bigr]=-\bigl\langle s_{\tau}(D;e_3),s_{\tau_0}(D;-e_3) \bigr\rangle_{\!\ST M\vert_{D}\,} [\ST M\vert_{\ast}].\end{equation*}
Let $[e_2,e_3]$ denote the shortest arc of great circle from $e_2$ to $e_3$ on $S^2$.
We have $\langle s_{\tau}(D;e_3),s_{\tau_0}(D;-e_3) \rangle_{\!\ST M\vert_{D}\,}=\langle s_{\tau}(D;e_2),s_{\tau_0}(D;-e_3) \rangle_{\!\ST M\vert_{D}\,}$ because $s_{\tau}(D;e_3)\cap s_{\tau_0}(D;-e_3)$ and $s_{\tau}(D;e_2)\cap s_{\tau_0}(D;-e_3)$ cobound $\tau(D \times [e_2,e_3]) \cap s_{\tau_0}(D;-e_3)$ in the interior of $\ST M\vert_{D}$.
We similarly get \begin{equation*}\langle s_{\tau}(D;e_3),s_{\tau_0}(D;-e_3) \rangle_{\!\ST M\vert_{D}\,}=\langle s_{\tau}(D;-e_2),s_{\tau_0}(D;-e_3) \rangle_{\!\ST M\vert_{D}\,}.\end{equation*}
Applying the involution $\iota_{\Sigma}$, we obtain
\begin{equation*}2\Bigl(\bigl[s_+(D)\bigr]-\bigl[s_{\tau}(D;e_3)\bigr]\Bigr)=\bigl\langle s_{\tau}(D;e_2),s_{\tau_0}(D;e_3)-s_{\tau_0}(D;-e_3) \bigr\rangle_{\!\ST M\vert_{D}\,}\bigl[\ST M\vert_{\ast}\bigr].\end{equation*}
The projection to $\RR^2 \times \{0\}$ of $\tau(. \times e_2)$ is an extension of the section $\tau(. \times e_2)$ of $T D$ to $D$. Its intersection with the zero section of $T D$ is the above intersection number.\eop

\begin{definition}
\label{defhomotoppseudo}
 A \emph{homotopy} from a pseudo-parallelization $(N(\gamma);\tau_e,\tau_b)$ to another such $(N(\gamma);\tau^{\prime}_e,\tau^{\prime}_b)$ is a homotopy from the pair $(\tau_e,\tau_b)$ to the pair $(\tau^{\prime}_e,\tau^{\prime}_b)$ such that \begin{equation*}\tau_b=\left\{\begin{array}{ll} \tau_e & \mbox{over}\; \partial\left([a,b] \times \gamma \times [-1,1]\right)\setminus \left(\{a\} \times \gamma \times [-1,1]\right) \\
 \tau_e \circ  \CT_{\gamma}& \mbox{over}\; \{a\} \times \gamma \times [-1,1]
\end{array} \right.\end{equation*} at any time. 
\end{definition}

\begin{proposition}
\label{propintsplusstau}
Let $\Sigma$, $e_3$, $s_+(\Sigma)$, and $s_-(\Sigma)$ be as in Lemma~\ref{lemsplusstau}. Let $\tilde{\tau}$ be a pseudo-parallelization restricting to a neighborhood of $\partial \Sigma$ as a genuine parallelization of $M$ such that $\tilde{\tau}(. \times e_3)$ is the positive normal to $\Sigma$ along $\partial \Sigma$. Let $\omega(\tilde{\tau})$ be a homogeneous boundary form associated with $\tilde{\tau}$ (as in Definition~\ref{defhombounformpseudo}). Then we have
\begin{equation*}\int_{s_+(\Sigma)}\omega(\tilde{\tau})=\frac12\chi\bigl(\tilde{\tau}(. \times e_2)\vert_{\partial \Sigma};\Sigma\bigr)\end{equation*}
and
\begin{equation*}\int_{s_-(\Sigma)}\omega(\tilde{\tau})=-\frac12\chi\bigl(\tilde{\tau}(. \times e_2)\vert_{\partial \Sigma};\Sigma\bigr).\end{equation*}
\end{proposition}
\bp Observe that $\omega(\tilde{\tau})$ is a closed form, which satisfies $\int_{\ST M\vert_{\ast}}\omega(\tilde{\tau})=1$. So Lemma~\ref{lemsplusstau} implies 
\begin{equation*}
\int_{s_{\pm}(\Sigma)}\omega(\tilde{\tau})=\mp \frac12\chi\bigl(\tilde{\tau}(. \times e_2)\vert_{\partial \Sigma};\Sigma\bigr) +\int_{s_{\tilde{\tau}}(\Sigma; \pm e_3)}\omega(\tilde{\tau}).
\end{equation*}
When $\tilde{\tau}$ is a genuine parallelization, the term $\int_{s_{\tilde{\tau}}(\Sigma; \pm e_3)}\omega(\tilde{\tau})$ is zero, and the proposition follows.

In general, for $\tilde{\tau}=(N(\gamma)=[a,b] \times \gamma \times [-1,1];\tau_e,\tau_b)$, the integral $\int_{s_+(\Sigma)}\omega(\tilde{\tau})$ is invariant under an isotopy of $\Sigma$ that fixes $\partial \Sigma$ since $\omega(\tilde{\tau})$ is closed.
It is also invariant under a homotopy of $(\tau_e,\tau_b)$ as in Definition~\ref{defhomotoppseudo} that is fixed in a neighborhood of $\partial \Sigma$.
(See Lemma~\ref{lemprimhom}.)

In particular, there is no loss of generality in assuming that $\Sigma$ meets $N(\gamma)$ along disks $D_c=[a,b] \times \{c\} \times [-1,1]$, and that ${\tau_b}(. \times e_3)$ is the positive normal to $D_c$ along $\partial D_c$ for these disks. Thus, thanks to the good behavior of the two sides of the equality to be proved under gluings along circles that satisfy the boundary conditions, it suffices to prove the proposition when $\Sigma$ is a meridian disk $D_c$ of $\gamma$ (with its corners smoothed) such that ${\tau_b}(. \times e_3)$ is the positive normal to $\Sigma$ along $\partial \Sigma$.
On $\ST M\vert_{D_c}$, the form $\omega(\tilde{\tau})$ is then equal to
\begin{equation*}\omega(\gamma,\tau_b)=\frac{p\bigl(\tau_b \circ \CT_{\gamma}^{-1}\bigr)^{\ast}\left(\omega_{S^2}\right) +p\bigl(\tau_b \circ {F}(\gamma,\tau_b)^{-1}\bigr)^{\ast}\left(\omega_{S^2}\right)}2,\end{equation*}
where $p(\tau_b \circ \CT_{\gamma}^{-1})^{\ast}\left(\omega_{S^2}\right)$ and $p\left(\tau_b \circ {F}(\gamma,\tau_b)^{-1}\right)^{\ast}\left(\omega_{S^2}\right)$ are propagating forms respectively associated with the parallelizations $\tau_b \circ \CT_{\gamma}^{-1}$ and $\left(\tau_b \circ {F}(\gamma,\tau_b)^{-1}\right)$, and $\chi\left(\tau_b \circ \CT_{\gamma}^{-1}(. \times e_2)\vert_{\partial \Sigma};\Sigma\right)=\chi\left(\tau_b(. \times e_2)\vert_{\partial \Sigma};\Sigma\right)$.

Therefore, we have
\begin{equation*}\int_{s_+(\Sigma)}\omega(\tilde{\tau})=\frac14\biggl(\chi\Bigl(\tau_b(. \times e_2)\vert_{\partial \Sigma};\Sigma\Bigr)+ \chi\Bigl(\bigl(\tau_b\circ {F}(\gamma,\tau_b)^{-1}\bigr)(. \times e_2)\vert_{\partial \Sigma};\Sigma\Bigr)\biggr).\end{equation*}
Thanks to Lemma~\ref{lemptyeulernumb}, this average is $\frac12 \left(\chi\left(\tau_e(. \times e_2)\vert_{\partial \Sigma};\Sigma\right)\right)$. This concludes the computation of $\int_{s_+(\Sigma)}\omega(\tilde{\tau})$. The computation of $\int_{s_-(\Sigma)}\omega(\tilde{\tau})$ is similar.
\eop

\section[Anomalous terms for pseudo-parallelizations]{Anomalous terms \\ for pseudo-parallelizations}

\begin{proposition}
\label{propextexistspseudo}
Let $A$ be a compact $3$-manifold equipped with two pseudo-parallelizations $\tau_0$ and $\tau_1$ that coincide with a common genuine parallelization along a regular neighborhood of $\partial A$. 
There exists a closed 
$2$-form $\omega$ on $\left[0,1\right] \times \ST A$ that restricts 
\begin{itemize}
 \item to $\{0\} \times \ST A$ as a homogeneous boundary form $\omega(\tau_0)$ of $(\ST A,\tau_0)$,
\item to $\{1\} \times \ST A$ as a homogeneous boundary form $\omega(\tau_1)$ of $(\ST A,\tau_1)$,
\item to $\left[0,1\right] \times \ST A\vert_{\partial A}$ as $p_{\ST A}^{\ast}\left(\omega(\tau_0)\right)$ with respect to the natural projection $p_{\ST A} \colon \left[0,1\right] \times \ST A \to \ST A$.
\end{itemize}
\end{proposition}
\bp Without loss of generality, assume that $A$ is connected.
Set $X=\left[0,1\right] \times \ST A$. Then $X$ is diffeomorphic to $\left[0,1\right] \times A \times S^2$ by a diffeomorphism induced by a parallelization $\tau$. The closed two-form $\omega$ is defined consistently on $\partial X$. It suffices to prove that the coboundary map $\partial$ of the long exact cohomology sequence associated to the pair $(X,\partial X)$ maps the class of $\omega\vert_{\partial X}$ to $0$ in
$H^3(X,\partial X)$. Since $H_3(X,\partial X)$
is Poincar\'e dual to \begin{equation*}H_3(X) \cong \bigl(H_1(A) \otimes H_2(S^2)\bigr) \oplus \bigl(H_3(A) \otimes H_0(S^2)\bigr),\end{equation*} it is generated \begin{itemize}
\item by classes of the form $\left[0,1\right] \times s_+(\Sigma)$ for surfaces $\Sigma$ of $A$ such that $\partial \Sigma \subset \partial A$ and for graphs $s_+(\Sigma)$ of sections in $\ST A$ associated to positive normals of the $\Sigma$, and 
\item by $\left[0,1\right] \times \{a\} \times S^2$ for some $a \in A$, when $\partial A=\emptyset$.
\end{itemize}
The evaluation of $\partial \left[\omega\vert_{\partial X}\right]$ on these classes is the evaluation of $\left[\omega\vert_{\partial X}\right]$ on their boundary. It is clearly zero for $\partial \left[0,1\right] \times \{a\} \times S^2$ since
\begin{equation*}\int_{\{(1,a)\} \times S^2}\omega=\int_{\{(0,a)\} \times S^2}\omega=
1.\end{equation*}
Let us conclude the proof by proving 
\begin{equation*}\int_{\partial (\left[0,1\right] \times  s_+(\Sigma))}\omega=0\end{equation*} for a surface $\Sigma$ as above transverse to $\{-1\}\times \partial A$.
This integral is invariant under the homotopies of $(\tau_0,\tau_1)$ that fix $(\tau_0,\tau_1)$
near $\partial A$. (See Lemma~\ref{lemprimhom}.)
Therefore, we can assume that $\tau_0=\tau_1$ in a neighborhood $[-2,0]\times \partial A$ of $\partial A$ in $A$, that $\tau_0$ is a genuine parallelization in this neighborhood, and that the positive normal to $\Sigma$ is $\tau_0(. \times e_3)$ along $\Sigma \cap \left(\{-1\}\times \partial A\right)$. Set $A_{-1}=A\setminus (\left]-1,0\right] \times \partial A)$ and $\Sigma_{-1}=\Sigma \cap A_{-1}$. Extend $\omega$ so that $\omega=p_{\ST A}^{\ast}\left(\omega(\tau_0)\right)$ over $\left[0,1\right] \times \ST A\vert_{[-2,0] \times \partial A}$. We have $\int_{\partial (\left[0,1\right] \times  s_+(\Sigma))}\omega=\int_{\partial (\left[0,1\right] \times  s_+(\Sigma_{-1}))}\omega$.
Proposition~\ref{propintsplusstau} ensures 
$\int_{\{0\} \times  s_+(\Sigma_{-1})}\omega=\int_{\{1\} \times  s_+(\Sigma_{-1})}\omega$.
Since $\int_{\left[0,1\right] \times \partial s_+(\Sigma_{-1})}\omega=0$, we get $\int_{\partial (\left[0,1\right] \times  s_+(\Sigma))}\omega=0.$
\eop

\begin{proposition}
\label{propextdefzn} Let $A$ be a compact oriented $3$-manifold equipped with three pseudo-parallelizations $\tau_0$, $\tau_1$, and $\tau_2$ that coincide with a common genuine parallelization along a regular neighborhood of $\partial A$. Let $n \in \NN$. Let $\omega(\tau_0)$ and $\omega(\tau_1)$ be homogeneous boundary forms respectively associated with $(\ST A,\tau_0)$ and $(\ST A,\tau_1)$.

Under the assumptions of Proposition~\ref{propextexistspseudo}, as in Corollary~\ref{corinvone}, set \begin{equation*}z_n(\left[0,1\right] \times \ST A;\omega)=\sum_{\Gamma \in \Davis^{c}_n}\coefgambet_{\Gamma}\int_{\left[0,1\right]\times \cinjuptd{\finsetv(\Gamma)}{T A}}\bigwedge_{e \in E(\Gamma)}p_e^{\ast}\left(\omega\right)\left[\Gamma\right].\end{equation*}
If $A$ embeds in a rational homology $3$-ball, then $z_n(\left[0,1\right] \times \ST A;\omega)$ 
depends only on the pseudo-parallelizations $\tau_0$ and $\tau_1$.
It is denoted by
$z_n(A;\tau_0,\tau_1)$, and the following properties are satisfied.
\begin{itemize} \item If $n$ is even, then $z_n(A;\tau_0,\tau_1)=0$.
\item If $B$ is a compact oriented $3$-manifold embedded in the interior of $A$, if $\tau_0$ and $\tau_1$ coincide on a neighborhood of $A\setminus B$, and if $\tau_0$ restricts to a neighborhood of $\partial B$ as a genuine parallelization, then $z_n(B;\tau_0,\tau_1)=z_n(A;\tau_0,\tau_1)$.
\item If $\tau_0$ and $\tau_1$ are actual parallelizations, then we have \begin{equation*}z_n(A;\tau_0,\tau_1)=\frac{p_1(\tau_0,\tau_1)}{4}\ansothree_n.\end{equation*}
\item We have $z_n(A;\tau_0,\tau_2)=z_n(A;\tau_0,\tau_1) + z_n(A;\tau_1,\tau_2)$.
(In particular, we have $z_n(A;\tau_0,\tau_1)=-z_n(A;\tau_1,\tau_0)$.)
\item For any orientation-preserving bundle isomorphism $\Psi$ of $\ST A$ over the identity map of $A$, we have
$z_n(A;\Psi \circ \tau_0,\Psi \circ \tau_1)=z_n(A;\tau_0,\tau_1)$.
\item For any orientation-preserving diffeomorphism $\psi$ from $A$ to another  compact oriented $3$-manifold $B$, we have
\begin{equation*}z_n\Bigl(B;T\psi \circ \tau_0 \circ \bigl(\psi^{-1} \times 1_{\RR^3}\bigr),T\psi \circ \tau_1 \circ \left(\psi^{-1} \times 1_{\RR^3}\right) \Bigr)=z_n(A;\tau_0,\tau_1).\end{equation*}
\item If $\tau^{\prime}_1$ is homotopic to $\tau_1$ relatively to $\partial A$ in the sense of Definition~\ref{defhomotoppseudo}, then we have \begin{equation*}z_n(A;\tau_0,\tau^{\prime}_1)=z_n(A;\tau_0,\tau_1).\end{equation*}
\item For any orientation-preserving diffeomorphism $\psi_1$ of $A$ isotopic to the identity map of $A$ relatively to $\partial A$, we have
\begin{equation*}z_n\Bigl(A;\tau_0,T\psi_1 \circ \tau_1\circ \bigl(\psi_1^{-1} \times 1_{\RR^3}\bigr) \Bigr)=z_n(A;\tau_0,\tau_1),\end{equation*}
where $\psi_1$ is used to carry the required parametrization of $N(\gamma)$.
\end{itemize}

\end{proposition}
\bp In this proof, the manifold $A$ embeds in a rational homology $3$-ball.
Lemma~\ref{lemsymeven} implies
that $z_n(A;\tau_0,\tau_1)=0$ when $n$ is even. Assume that $n$ is odd from now on.
Let us first prove that
$z_n(\left[0,1\right] \times \ST A;\omega)$ does not depend on the closed extension $\omega$ when $A$ is a rational homology $3$-ball and when $\tau_0$ is standard near $\partial A$. 
According to Lemma~\ref{lemexistpropagatorpseudohom}, $\omega(\tau_0)$ (resp. $\omega(\tau_1)$) extends to a homogeneous propagating form of $(C_2(S^3(A/\ballb_{S^3})),\tau_0)$ (resp.  of $(C_2(S^3(A/\ballb_{S^3})),\tau_1)$). 
Set $X=\left[0,1\right] \times C_2(S^3(A/\ballb_{S^3}))$. The above extensions of $\omega(\tau_0)$ and $\omega(\tau_1)$ together with $\omega$ (extended as $\projp_{\taust}^{\ast}\left(\omega_{S^2}\right)$ on $\left[0,1\right] \times \left(\partial C_2(S^3(A/\ballb_{S^3})) \setminus \ST A\right)$) determine a closed $2$-form of $\partial X$. This form extends as a closed form on $X$ by Lemma~\ref{lemformprod}.
Then Corollary~\ref{corinvone} implies that $z_n(\left[0,1\right] \times \ST A;\omega)$ does not depend on $\omega$ when $A$ is a rational homology $3$-ball and when $\tau_0$ is standard near $\partial A$. 
In general, embed $A$ in the interior of such a space $\ballb_{\rats}$. The pseudo-parallelization $\tau_0$ extends to $\ballb_{\rats}$ as a pseudo-parallelization standard near $\partial \ballb_{\rats}$ according to Lemma~\ref{lempseudoparextend}. The form $\omega$ of Proposition~\ref{propextexistspseudo} extends to $\left[0,1\right] \times \ST \ballb_{\rats}$ as $p_{\ST \ballb_{\rats}}^{\ast}\left(\omega(\tau_0)\right)$ on $\left[0,1\right] \times \ST (\ballb_{\rats} \setminus A)$.
Then $z_n(\left[0,1\right] \times \ST (\ballb_{\rats} \setminus A);\omega)=0$ since 
$\bigwedge_{e \in E(\Gamma)}p_e^{\ast}\left(\omega\right)$ pulls back through $\cinjuptd{\finsetv(\Gamma)}{T (\ballb_{\rats} \setminus A)}$ for any involved $\Gamma$.
So $z_n(\left[0,1\right] \times \ST A;\omega)=z_n(\left[0,1\right] \times \ST \ballb_{\rats};\omega)$ is independent of $\omega\vert_{\left[0,1\right] \times \ST A}$. 

Set $z_n(A;\omega(\tau_0),\omega(\tau_1) )=z_n(\left[0,1\right] \times \ST A;\omega)$.
We easily observe \begin{equation*}z_n\bigl(A;\omega(\tau_0),\omega(\tau_2)\bigr)=z_n\bigl(A;
\omega(\tau_0),\omega(\tau_1)\bigr) + z_n\bigl(A;\omega(\tau_1),\omega(\tau_2)\bigr).\end{equation*}

Assume $\tau_0=\tau_1=(N(\gamma);\tau_e,\tau_b)$. Assume that the forms $\omega(\tau_0)=\omega(\tau_0,F_0)$ and $\omega(\tau_1)=\omega(\tau_0,F_1)$ of Definition~\ref{defhombounformpseudo} are obtained from one another by changing the map $F=F_0\colon [a,b] \times [-1,1]\to SO(3)$ to another one $F_1$.
There exists a homotopy $F_t$ from $F_0$ to $F_1$. Such a homotopy induces a homotopy
\begin{equation*}\begin{array}{llll} {F}_t(\gamma,\tau_b) \colon &  [a,b] \times   \gamma \times [-1,1] \times S^2 & \longrightarrow & [a,b] \times   \gamma \times [-1,1] \times S^2\\
&(s,c,u; \svarM) & \mapsto & \bigl(s,c,u; F_t(s,u)(\svarM)\bigr).\end{array}\end{equation*}
Then we have $z_n(A;\omega(\tau_0),\omega(\tau_1) )=z_n(N(\gamma);\omega(\tau_0,F_0),\omega(\tau_0,F_1) )$.
Use $\tau_b$ to identify $\ST N(\gamma)$ with $[a,b] \times \gamma \times [-1,1] \times S^2$, and define $\omega(\gamma,\tau_b)$ on $\left[0,1\right] \times [a,b] \times \gamma \times [-1,1] \times S^2$
with respect to the formula for $\omega(\gamma,\tau_b)$ in Definition~\ref{defhombounformpseudo} by
\begin{equation*}\omega(\gamma,\tau_b)=\frac{p\bigl(\tau_b \circ \CT_{\gamma}^{-1}\bigr)^{\ast}\left(\omega_{S^2}\right) +p\bigl(\tau_b \circ {F}_.(\gamma,\tau_b)^{-1}\bigr)^{\ast}\left(\omega_{S^2}\right)}2.\end{equation*}
This formula does not depend on the coordinate along $\gamma$. So $\omega$ pulls back through a projection from $\left[0,1\right] \times \ST N(\gamma)$ to $\left[0,1\right] \times [a,b] \times [-1,1] \times S^2$. The 
$\bigwedge_{e \in E(\Gamma)}p_e^{\ast}(\omega)$ pull back through $\left[0,1\right] \times \cinjuptd{\finsetv(\Gamma)}{T N(\gamma)\vert_{[a,b] \times \{c\} \times [-1,1]}}$. Hence $z_n(\left[0,1\right] \times \ST N(\gamma);\omega)$ vanishes. 
This proves that $z_n(N(\gamma);\omega(\tau_0,F_0),\omega(\tau_0,F_1) )$ vanishes. We conclude that
$z_n(A;\omega(\tau_0),\omega(\tau_1) )$ depends only on $\tau_0$ and $\tau_1$.

For an orientation-preserving bundle isomorphism $\Psi$ of $\ST A$ over the identity map of $A$, the pseudo-parallelization 
\begin{equation*}\Psi\circ\left(\tau_0=(N(\gamma);\tau_e,\tau_b)\right)=\left(N(\gamma);\Psi\circ\tau_e,\Psi\circ\tau_b\right)\end{equation*}
makes unambiguous sense. The following commutative diagram
\begin{equation*}
\xymatrix{\ST A  \ar@/^1pc/[drr]^{p(\Psi \circ \tau)}\\
\ST A \ar[u]_{\Psi} \ar[r]^{\tau^{-1}} \ar@/_1pc/[rr] _{p(\tau)}& A\times S^2 \ar[ul]_{(\Psi \circ \tau)} \ar[r]^{p_{S^2}} & S^2}
\end{equation*}
shows that $p(\tau_b)=p(\Psi\circ\tau_b) \circ \Psi$. So we have
\begin{equation*}\omega(\Psi\circ\tau_0,F)=\left(\Psi^{-1}\right)^{\ast}\omega(\tau_0,F).\end{equation*}
The form $\omega$ on $\left[0,1\right] \times \ST A$ can be pulled back by the orientation-preserving $1_{\left[0,1\right]} \times \Psi^{-1}$, similarly. We get
\begin{equation*}z_n(A;\Psi \circ \tau_0,\Psi \circ \tau_1)=z_n(A;\tau_0,\tau_1).\end{equation*}

Proposition~\ref{propdeftwoanom} implies that $z_n(A;\tau_0,\tau_1)=\frac{p_1(\tau_0,\tau_1)}{4}\ansothree_n$ as soon as $\tau_0$ and $\tau_1$ are actual parallelizations and $A$ embeds in a rational homology ball to which $\tau_0$ extends as a genuine parallelization. Embed $A$ in a rational homology ball $B_{\rats}$. Let $\tau$ be a parallelization of $B_{\rats}$. Then $\tau\vert_{A}=\Psi \circ \tau_0$ for an orientation-preserving bundle isomorphism $\Psi$ of $\ST A$. We have $z_n(A;\Psi \circ\tau_0,\Psi \circ\tau_1)=\frac{p_1(\Psi \circ\tau_0,\Psi \circ\tau_1)}{4}\ansothree_n$. The above behavior of $z_n$ under an orientation-preserving bundle isomorphism of $\ST A$ allows us to
conclude that $z_n(A;\tau_0,\tau_1)=\frac{p_1(\tau_0,\tau_1)}{4}\ansothree_n$ as soon as $\tau_0$ and $\tau_1$ are actual parallelizations.

For any orientation-preserving diffeomorphism $\psi$ from $A$ to $B$, the reader can check 
\begin{equation*}z_n(B;T\psi \circ \tau_0\circ \left(\psi^{-1} \times 1_{\RR^3}\right),T\psi \circ \tau_1\circ \left(\psi^{-1} \times 1_{\RR^3}\right))=z_n(A;\tau_0,\tau_1)\end{equation*} as above.

If $\tau_1$ is homotopic to $\tau_0$ in the sense of Definition~\ref{defhomotoppseudo}, then there exists a map $\Psi \colon \left[0,1\right] \times \ST A \to \ST A$ such that $\left(t\mapsto \Psi(t,.) \circ \tau_0\right)$ is a homotopy of pseudo-parallelizations from $\tau_0$ to $\tau_1$.
So we have \begin{equation*}z_n(A;\tau_0,\tau_1)=z_n\Bigl(\left[0,1\right] \times \ST A;\omega=\bigl(\Psi^{-1}\bigr)^{\ast}\bigl(\omega(\tau_0,F)\bigr)\Bigr).\end{equation*}
The form $\omega$ pulls back through a map from $\left[0,1\right] \times \ST A$ to $\ST A$, and the forms
$\bigwedge_{e \in E(\Gamma)}p_e^{\ast}(\omega)$ pull back through $\cinjuptd{\finsetv(\Gamma)}{T A}$. So $z_n(\left[0,1\right] \times \ST A;\omega)$ vanishes.

Let $\psi \colon \left[0,1\right]\times A \to A$ be an isotopy, which maps $(t,u)$ to $\psi_t(u)$, such that $\psi_0=\id$. It induces the homotopy
\begin{equation*}\begin{array}{llll}\Psi \colon &\left[0,1\right]\times \ST A &\to &\ST A \\
  & (t,u) & \mapsto & \tau_1 \circ \left(\psi_t \times 1_{S^2}\right) \circ \tau_1^{-1} \circ T\psi_t^{-1}(u),
  \end{array}\end{equation*} 
which satisfies 
\begin{equation*}\begin{array}{ll}p\bigl(T\psi_t \circ \tau_1 \circ (\psi_t^{-1}\times 1_{S^2})\bigr)&=p_{S^2} \circ \tau_1^{-1} \circ \tau_1 \circ \bigl(T\psi_t \circ \tau_1 \circ (\psi_t^{-1}\times 1_{S^2})\bigr)^{-1}\\&=p(\tau_1)\circ \Psi(t,.).\end{array}\end{equation*}
We get
\begin{multline*}z_n\Bigl(A;\tau_1=T\psi_0 \circ \tau_1\circ \bigl(\psi_0^{-1} \times 1_{\RR^3}\bigr),T\psi_1 \circ \tau_1\circ\bigl(\psi_1^{-1} \times 1_{\RR^3}\bigr)\Bigr)\\\
=z_n\Bigl([0,1]\times \ST A;\bigl(\Psi\bigr)^{\ast}\bigl(\omega(\tau_1)\bigr)\Bigr)=0.\end{multline*}
\eop

The main result of this section is the following theorem.

\begin{theorem}
\label{thmznponepseudo}
Let $A$ be a compact $3$-manifold equipped with two pseudo-parallelizations $\tau_0$ and $\tau_1$ that coincide with a common genuine parallelization along a regular neighborhood of $\partial A$. Assume that $A$ embeds in a rational homology $3$-ball. 
With the notation of Proposition~\ref{propextdefzn}, we have
\begin{equation*}z_n(A;\tau_0,\tau_1)=\frac{p_1(\tau_0,\tau_1)}{4}\ansothree_n \end{equation*} for any natural integer $n$.
\end{theorem}

To prove this theorem, we will prove it in special cases and show that these special cases are sufficient to get a complete proof.

\begin{lemma}
\label{lempseudodoubletorus}
Assume $A=[2,9] \times \gamma \times [-2,2]$. Equip $A$ with a pseudo-parallelization
$\tau_0=(N(\tilde{\gamma});\tau_e,\tau_b)$ such that
\begin{equation*}N(\tilde{\gamma})=[3,5] \times \gamma \times [-1,1] \sqcup [6,8] \times \gamma \times [-1,1].\end{equation*}
Then there exists a parallelization $\tau_1$ of $A$ that coincides with $\tau_e$ in a neighborhood of $\partial A$.
Furthermore, we have
\begin{equation*}z_n(A;\tau_0,\tau_1)=\frac{p_1(\tau_0,\tau_1)}{4}\ansothree_n\end{equation*} for any such parallelization.
\end{lemma}
\bp We first prove the lemma for some chosen pseudo-parallelizations $\tilde{\tau}_0$ and $\tilde{\tau}_1$ satisfying the assumptions and behaving as \say{products by $\gamma$}.
For these pseudo-parallelizations, this product behavior will imply 
$p_1(\tilde{\tau}_0,\tilde{\tau}_1)=0$ and $z_n(A;\tilde{\tau}_0,\tilde{\tau}_1)=0$.
Define the parallelization
\begin{equation*}\begin{array}{llll}\tau_A \colon& A \times \RR^3 &\to& \ST A\\
   &(s_0,c_0,u_0;e_1)&\mapsto&\frac{d}{ds}(s,c_0,u_0)(s_0,c_0,u_0)\\
   &(s_0,c_0,u_0;e_2)&\mapsto&\frac{d}{dc}(s_0,c,u_0)(s_0,c_0,u_0)\\
   &(s_0,c_0,u_0;e_3)&\mapsto&\frac{d}{du}(s_0,c_0,u)(s_0,c_0,u_0)
  \end{array}\end{equation*} of $A$.
Define $\tilde{\tau}_e \colon \bigl( A \setminus \mathring{N}(\tilde{\gamma}) \bigr) \times \RR^3 \to  \ST \bigl( A \setminus \mathring{N}(\tilde{\gamma})\bigr) $ by
\begin{equation*}\tilde{\tau}_e=\left\{
\begin{array}{ll}\tau_A &\mbox{on}\; A \setminus \left(\left[2,8\right[ \times \gamma \times \left]-1,1\right[\right) \times \RR^3\\
 \tau_A \circ \CT_{\gamma}^{-1} &\mbox{on}\; [5,6] \times \gamma \times [-1,1]\\
 \tau_A \circ \CT_{\gamma}^{-2} &\mbox{on}\; [2,3] \times \gamma \times [-1,1],
  \end{array}\right.\end{equation*}
where $\CT_{\gamma}(t,c \in \gamma,u \in [-1,1];\cvarM\in \RR^3)=(t,c,u;\rho_{\alpha(u)}(\cvarM))$ as in Definition~\ref{defpseudotriv}.
Define $\tilde{\tau}_b  \colon {N}(\tilde{\gamma}) \times \RR^3 \to  \ST {N}(\tilde{\gamma})$ by
\begin{equation*}\tilde{\tau}_b=\left\{
\begin{array}{ll}
 \tau_A &\mbox{on}\; [6,8] \times \gamma \times [-1,1]\\
 \tau_A \circ \CT_{\gamma}^{-1} &\mbox{on}\; [3,5] \times \gamma \times [-1,1].
  \end{array}\right.\end{equation*}
Set $\tilde{\tau}_0=(N(\tilde{\gamma});\tilde{\tau}_e,\tilde{\tau}_b)$. Define a map
\begin{equation*}\begin{array}{llll}\tilde{F} \colon& [3,8] \times [-1,1] &\to& SO(3)\\
   &(3,u)&\mapsto&\rho_{-2\alpha(u)}\\
  & (8,u)&\mapsto&1_{SO(3)}\\
   &(t,\pm 1)&\mapsto&1_{SO(3)}.
  \end{array}\end{equation*}
Finally define $\tilde{\tau}_1$ such that
\begin{equation*}\tilde{\tau}_1=\left\{
\begin{array}{ll}\tau_A &\mbox{on}\; A \setminus \left( \left[2,8\right[ \times \gamma \times \left]-1,1\right[ \right) \times \RR^3\\
 \tau_A \circ \CT_{\gamma}^{-2} &\mbox{on}\; [2,3] \times \gamma \times [-1,1]
  \end{array}\right.\end{equation*}
and \begin{equation*}\tilde{\tau}_1(s,c,u;\fvarM)=\Bigl({\tau}_A\bigl(s,c,u;\tilde{F}(s,u)(\fvarM)\bigr)\Bigr)\end{equation*} when $(s,u) \in [3,8] \times [-1,1].$

Let us prove $p_1(\tilde{\tau}_0,\tilde{\tau}_1)=0$.
The involved trivializations of $T(\left[0,1\right] \times A) \otimes C$ on $\partial \left[0,1\right] \times A$ are obtained from the natural parallelization $T\left[0,1\right] \oplus \tau_A$ by composition by a map from $\partial (\left[0,1\right] \times A)= \gamma \times \partial(\left[0,1\right] \times [2,9] \times [-2,2])$ to $SU(4)$, which does not depend on the coordinate along $\gamma$. Since $\pi_2(SU(4))=\{0\}$, this map extends to $SU(4)$.

Let us similarly prove $z_n(A;\tilde{\tau}_0,\tilde{\tau}_1)=0$. 
Set $Y=\left[0,1\right] \times S^2 \times [2,9] \times [-2,2]$. The parallelization $\tau_A$ identifies $\left[0,1\right] \times \ST A$ with $\gamma \times Y$.
We have $z_n(A;\tilde{\tau}_0,\tilde{\tau}_1)=z_n(\left[0,1\right] \times \ST A;\omega)$ for a closed two-form $\omega$ whose restriction $\omega\vert_{\partial}$ to $\partial(\left[0,1\right] \times \ST A)$ factors through the projection of $\gamma \times \partial Y$ onto $\partial Y$. The involved closed $2$-form on $\partial Y$ extends to $Y$ as a closed form $\omega_{Y}$ since  $\omega\vert_{\partial}$ extends to the whole $\left[0,1\right] \times \ST A$, according to Proposition~\ref{propextexistspseudo}. Then $\omega$ can be chosen as the pull-back of $\omega_{X}$ under the projection of $\gamma \times X$ onto $X$. Thus, the forms
$\bigwedge_{e \in E(\Gamma)}p_e^{\ast}(\omega)$ again pull back through a projection onto a space of dimension smaller than the degree of the forms. So $z_n(\left[0,1\right] \times \ST A;\omega)$ vanishes.

There exists an orientation-preserving bundle isomorphism $\Psi$ of $\ST A$ over the identity map of $A$,
such that $\tau_0=\Psi \circ \tilde{\tau}_0$. The parallelization $\tau_1=\Psi \circ \tilde{\tau}_1$ satisfies the assumptions of the lemma. We have
$z_n(A;\Psi \circ \tilde{\tau}_0,\Psi \circ \tilde{\tau}_1)=z_n(A;\tilde{\tau}_0,\tilde{\tau}_1)=0$ and
$p_1(\Psi \circ \tilde{\tau}_0,\Psi \circ \tilde{\tau}_1)=p_1(\tilde{\tau}_0,\tilde{\tau}_1)=0$. So we have $z_n(A;\tau_0,\tau_1)=\frac{p_1(\tau_0,\tau_1)}{4}\ansothree_n$.
We conclude for any another parallelization $\tau^{\prime}_1$ that coincides with $\tau_e$ near $\partial A$, because $z_n(A;\tau_1,\tau^{\prime}_1)=\frac{p_1(\tau_1,\tau^{\prime}_1)}{4}\ansothree_n$.
\eop

\begin{lemma}
\label{lembounsurfpseudo}
Let $A$ be a compact oriented $3$-manifold that embeds in a rational homology $3$-ball.
Let $[-7,0] \times \partial A$ be a collar neighborhood of $A$.
Let $\gamma \times [-2,2]$ be a disjoint union of annuli in $\partial A$. Set $N(\gamma)=[-2,-1] \times \gamma \times \left[-1,1\right]$.
Let $\tau_0=(N(\gamma);\tau_e,\tau_b)$ be a pseudo-parallelization of $A$ that coincides with the restriction of a parallelization $\tau_1$ of $A$ in a neighborhood of $\partial A$.
Then we have
\begin{equation*}z_n(A;\tau_0,\tau_1)=\frac{p_1(\tau_0,\tau_1)}{4}\ansothree_n .\end{equation*}
\end{lemma}
\bp Recall $A_{-2}= A\setminus (\left]-2,0\right] \times \partial A )$. Figure~\ref{figproofpseudoone} shows the schema of the proof.

\bfig
\centering
\begin{tikzpicture} 
\draw (6,2) rectangle (10,5) (8,3.5) node{\scriptsize Box to be computed} (8,3) node{\scriptsize over $\left[0,1\right] \times A$};
\draw [dashed] (0,2) rectangle (6,4) (3,3.5) node{\scriptsize Box over $\left[-1,0\right] \times A_{-2}$} (3,3) node{\scriptsize isomorphic to the} (3,2.5) node{\scriptsize box to be computed};
\draw [dashed] (0,4) rectangle (6,5) (3,4.65) node{\scriptsize Trivial product by $\left[-1,0\right]$  over } (3,4.35) node{\scriptsize $\left[-1,0\right] \times  [-2,0] \times \partial A$};
\draw (-.1,1.7) node[right]  {\scriptsize $(A,\tau_{-1})$} (-.5,1.4) node[right] {\scriptsize as in Lemma~\ref{lempseudodoubletorus}} (6,1.7) node {\scriptsize $(A,\tau_0)$} (10,1.7) node  {\scriptsize $(A,\tau_1)$};
\draw (-.1,2.3) node[left]{\scriptsize $-7$} (-.1,4) node[left]{\scriptsize $-2$}(-.1,5) node[left]{\scriptsize $0$};
\fill [white] (-.15,4.3) rectangle (.15,4.75);
\draw (-.15,4.3) rectangle (.15,4.75) (0,4.5) node {\scriptsize $\gamma$};
\begin{scope}[xshift=6cm]
 \fill [white] (-.15,4.3) rectangle (.15,4.75);
\draw (-.15,4.3) rectangle (.15,4.75) (0,4.5) node {\scriptsize $\gamma$};
\end{scope}
\begin{scope}[yshift=-1cm]
 \fill [white] (-.15,4.3) rectangle (.15,4.75);
\draw (-.15,4.3) rectangle (.15,4.75) (0,4.5) node {\scriptsize $\gamma$};
\end{scope}
\end{tikzpicture}
\caption{Schema of proof for Lemma~\ref{lembounsurfpseudo}}
\label{figproofpseudoone}
\end{figure}
Let $f \colon [-7,0] \to [-7,-2]$ be a diffeomorphism
such that $f(t)=t-2$ when $t\geq -3$, and $f(t)=t$ when $t\leq -6$.
Let $\psi\colon  A \to A_{-2}$ be a diffeomorphism restricting to $A \setminus ( \left]-7,0\right] \times \partial A )$ as the identity map, and mapping $(t,x) \in [-7,0] \times \partial A$ to $(f(t),x)$.

There exists a bundle isomorphism $\Phi$ of $\ST A_{-2}$ over the identity map of $A_{-2}$ such that $\tau_e\vert_{A_{-2}}=\Phi \circ T\psi \circ \tau_1\circ \left(\psi^{-1} \times 1_{\RR^3}\right)$.
Let $\tau_{-1}$ be the pseudo-parallelization of $A$ that coincides with $\tau_0$ over 
$[-2,0] \times \Sigma$, and with $\Phi \circ T\psi \circ \tau_0\circ \left(\psi^{-1} \times 1_{\RR^3}\right)$
over $A_{-2}$.
Then $\tau_{-1}$ is a parallelization outside $[-7,0] \times \gamma \times [-2,2]$. Since $[-7,0] \times \gamma \times [-2,2]$ satisfies the hypotheses of Lemma~\ref{lempseudodoubletorus} up to reparametrization, Lemma~\ref{lempseudodoubletorus} and Proposition~\ref{propextdefzn}
ensure that 
\begin{equation*}z_n(A;\tau_{-1},\tau_1)=\frac{p_1(\tau_{-1},\tau_1)}{4}\ansothree_n.\end{equation*}

To conclude, we prove that $z_n(A;\tau_0,\tau_1)=\frac12 z_n(A;\tau_{-1},\tau_1)$ and 
$p_1(\tau_0,\tau_1)=\frac12 p_1(\tau_{-1},\tau_1)$.
The element $z_n(A;\tau_{-1},\tau_0)$ of $\Aavis_n(\emptyset)$ can be written as
\begin{equation*}z_n\bigl(A_{-2};\Phi \circ T\psi \circ \tau_0\circ (\psi^{-1} \times 1_{\RR^3}),\Phi \circ T\psi \circ \tau_1\circ (\psi^{-1} \times 1_{\RR^3})\bigr)
=z_n(A;\tau_0,\tau_1).\end{equation*}
We similarly have $p_1(\tau_{-1},\tau_0)=p_1(\tau_0,\tau_1)$.
\eop

\begin{lemma}
\label{lemSigmapseudoone}
Let $\Sigma$ be a compact oriented surface.
Let $\gamma_0$ and $\gamma_1$ be two disjoint unions of curves of $\Sigma$ with respective tubular neighborhoods $\gamma_0 \times [-1,1]$ and $\gamma_1 \times [-1,1]$. 
Set $A=[0,3]\times \Sigma$, $N(\gamma_0)=\left[1,2\right] \times \gamma_0 \times [-1,1]$ and $N(\gamma_1)=\left[1,2\right] \times \gamma_1 \times [-1,1]$.

There exist two pseudo-parallelizations $\tau_0=(N(\gamma_0);\tau_{e,0},\tau_{b,0})$ and $\tau_1=(N(\gamma_1);\tau_{e,1},\tau_{b,1})$ that coincide near $\partial A$ if and only if $\gamma_0$ and $\gamma_1$ 
have the same class in $H_1(\Sigma;\ZZ/2\ZZ)$.
In this case, we have \begin{equation*}z_n(A;\tau_0,\tau_1)=\frac{p_1(\tau_0,\tau_1)}{4}\ansothree_n .\end{equation*}
\end{lemma}
\bp
Without loss of generality, assume that $\Sigma$ is connected and that the complement of $\left(\gamma_0 \times [-1,1]\right) \cup \left(\gamma_1 \times [-1,1]\right)$ in $\Sigma$ is not empty.
Let $\tau$ be a parallelization of $A=[0,3]\times \Sigma$. 

Assume that $\tau_0$ and $\tau_1$ are two pseudo-parallelizations as in the statement, which coincide near $\partial A$. Let us prove that $\gamma_0$ and $\gamma_1$ are homologous modulo $2H_1(\Sigma;\ZZ)$.

If the boundary of $\Sigma$ is empty, choose a disk $D$ of $\Sigma$ outside $\gamma_0 \times [-1,1] \cup \gamma_1 \times [-1,1]$ and assume that $\tau_{e,0}$ and $\tau_{e,1}$ coincide on $[0,3]\times D$, without loss of generality.
This allows us to assume $\partial \Sigma \neq \emptyset$, without loss of generality, by possibly removing the interior of $D$ from $\Sigma$.
For $i \in \{0,1\}$ and for the pseudo-parallelization $\tau_i$ of the statement, write $\tau_{e,i}$ 
as $\tau \circ \psi_{\RR}(g_i)$ for some $g_i \colon A \setminus \mathring{N}(\gamma_i) \to SO(3)$, with the notation of Section~\ref{secannpont}. Then the restriction of $g_i$ to a meridian curve of $\gamma_i$ is not homotopic to a constant loop. The maps $g_0$ and $g_1$  coincide near $\partial A$. Let $c \colon \left[0,1\right] \to \Sigma$ be a path such that $c(0)$ and $c(1)$ are in $\partial \Sigma$. Then 
the restriction of $g_i$ to \begin{equation*}\Bigl(\{3\} \times c\bigl([0,1]\bigr)\Bigr) \cup \Bigl(-[0,3] \times c(1)\Bigr) \cup \Bigl(-\{0\} \times c\bigl([0,1]\bigr)\Bigr) \cup \Bigl([0,3] \times c(0)\Bigr)\end{equation*}
is null-homotopic if and only if the mod $2$ intersection of $c$ with $\gamma_i$ is trivial.
So $\gamma_0$ and $\gamma_1$ must be homologous modulo $2H_1(\partial A;\ZZ)$.

Conversely, assume that $\gamma_0$ and $\gamma_1$ are homologous modulo $2H_1(\Sigma;\ZZ)$. Define $g_0 \colon A \setminus \mathring{N}(\gamma_0)\to SO(3)$ to be the map sending $A \setminus \left(\left[0,2\right] \times \gamma_0 \times \left]-1,1\right[\right)$ to the identity of $SO(3)$, and mapping
$(t,c,u) \in \left(\left[0,1\right] \times \gamma_0 \times [-1,1]\right)$ to $\rho_{-\alpha(u)}$.
First define $g_1\colon A \setminus \mathring{N}(\gamma_1)\to SO(3)$ 
on $\bigl(A \setminus (\left[0,1\right[ \times \Sigma) \bigr) \setminus \mathring{N}(\gamma_1)$
so that $g_1$ sends $\bigl(A \setminus ([0,1] \times \Sigma) \bigr) \setminus \mathring{N}(\gamma_1)$ to the identity of $SO(3)$ and $g_1$ maps 
$(t,c,u) \in \left(\{1\} \times \gamma_1 \times [-1,1]\right)$ to $\rho_{-\alpha(u)}$. Since the classes of $\gamma_0$ and $\gamma_1$ 
coincide in $H_1(\Sigma;\ZZ/2\ZZ)$, the restrictions 
to $\{1\} \times \Sigma$ of $g_1$ and $g_0$ are homotopic on the one-skeleton of $\Sigma$, and hence on $\Sigma$. 
This allows us to extend $g_1$ to $\left[0,1\right] \times \Sigma$ so that 
$g_0$ and $g_1$ coincide in a neighborhood of $\partial A$.
Then, for $i\in \{0,1\}$, there exists $\tau_{b,i}$ such that $\tau_i=(N(\gamma_i);\tau_{e,i}=\tau \circ \psi_{\RR}(g_i),\tau_{b,i})$ is a pseudo-parallelization.
The pseudo-parallelizations $\tau_0$ and $\tau_1$ coincide near $\partial A$.

Set $B=\left[0,6\right]\times \Sigma$ and $N(\gamma^{\prime}_0)=[4,5] \times \gamma_0 \times [-1,1]$.
Extend $\tau_0$ to $B$ as a pseudo-parallelization $\tau_{0,B}=(N(\gamma_0) \sqcup N(\gamma^{\prime}_0);\tau_{e,B,0},\tau_{b,B,0})$. Extend $\tau_1$ to $B$ as a pseudo-parallelization $\tau_{1,B}$ that coincides with $\tau_{0,B}$ on $[3,6]\times \Sigma$.
We have $z_n(A;\tau_0,\tau_1)=z_n(B;\tau_{0,B},\tau_{1,B})$ and $p_1(\tau_0,\tau_1)=p_1(\tau_{0,B},\tau_{1,B})$.

According to Lemma~\ref{lempseudodoubletorus}, there is a parallelization $\tau_2$ of $B$ that coincides with $\tau_{0,B}$ near $\partial B$, and we have $z_n(B;\tau_{0,B},\tau_2)=\frac{p_1(\tau_{0,B},\tau_2)}{4}\ansothree_n$.

Apply Lemma~\ref{lembounsurfpseudo} 
to prove that $z_n(B;\tau_{1,B},\tau_2)=\frac{p_1(\tau_{1,B},\tau_2)}{4}\ansothree_n$.
(To apply Lemma~\ref{lembounsurfpseudo} as it is stated in 
$B$, first rotate $N(\gamma_0)=\left[1,2\right] \times \gamma_0 \times [-1,1]$ around $\gamma_0$ by an isotopy that sends
$1 \times [-1,1]$ to $2 \times [1,-1]$ and apply Proposition~\ref{propextdefzn}.) 
This implies that
$z_n(B;\tau_{0,B},\tau_{1,B})=\frac{p_1(\tau_{0,B},\tau_{1,B})}{4}$.
\eop

\begin{lemma}
\label{lemSigmapseudotwo}
Let $\Sigma$ be a compact connected oriented surface with boundary.
Let $\gamma$, $\gamma_0$, and $\gamma_1$ be three disjoint unions of curves of $\Sigma$ with respective tubular neighborhoods $\gamma \times [-1,1]$, $\gamma_0 \times [-1,1]$, and $\gamma_1 \times [-1,1]$.\footnote{The disjoint union of curves $\gamma$ may intersect $\gamma_0$ and $\gamma_1$, and $\gamma_0$ may intersect $\gamma_1$.} Assume that $[\gamma_1]=[\gamma_0]+\left[\gamma\right]$ in $H_1(\Sigma;\ZZ/2\ZZ)$.
Set $A=\left[0,6\right]\times \Sigma$, $N(\gamma)=[4,5] \times \gamma \times [-1,1]$, $N(\gamma_0)=\left[1,2\right] \times \gamma_0 \times [-1,1]$, and $N(\gamma_1)=\left[1,2\right] \times \gamma_1 \times [-1,1]$.
Let $\tau_0=(N(\gamma) \sqcup N(\gamma_0);\tau_{e,0},\tau_{b,0})$ and $\tau_1=(N(\gamma_1);\tau_{e,1},\tau_{b,1})$
be two pseudo-parallelizations which coincide near $\partial A$.
Then we have \begin{equation*}z_n(A;\tau_0,\tau_1)=\frac{p_1(\tau_0,\tau_1)}{4}\ansothree_n .\end{equation*}
\end{lemma}
\bp Lemma~\ref{lemSigmapseudoone} allows us to choose arbitrary representatives of $\left[\gamma\right]$ and $\left[\gamma_0\right]$ for the proof, without loss of generality.
In particular, there is no loss of generality in assuming that $\gamma_0$ is connected and that the intersection of $\gamma$ and $\gamma_0$ has no more than one point.

If $\gamma$ and $\gamma_0$ are disjoint, then we can perform an isotopy in $A$ to lower $\gamma$, and the result is a direct consequence of Lemma~\ref{lemSigmapseudoone}.

Assume that the intersection of $\gamma$ and $\gamma_0$ has one transverse point.
Attach two copies $\Sigma$ and $\Sigma^{\prime}$ of $\Sigma$ to a disk $D$ along intervals $I$ and $I^{\prime}$.
Let $\gamma^{\prime}$, $\gamma^{\prime}_0$, and $\gamma^{\prime}_1$ be the respective copies of $\gamma$, $\gamma_0$, and $\gamma_1$ in $\Sigma^{\prime}$.
Let $\tilde{\Sigma}=\Sigma \cup D \cup \Sigma^{\prime}$ as in Figure~\ref{figsigmatildepseudo}.

\bfig
\centering
\begin{tikzpicture}[scale=.5] \useasboundingbox (-1.1,-1.4) rectangle (16.1,1.3);
\draw [-] (0:1.1) arc (0:180:1.1) (0:.7) arc (0:180:.7);
\draw (-.7,0) .. controls (-.7,-.1) ..  (-.3,-.1) .. controls (.1,-.1) .. (.1,0) (.5,0) .. controls (.5,-.1) ..  (.6,-.1) .. controls (.7,-.1) .. (.7,0) (1.1,0) .. controls (1.1,-.1) ..  (1.5,-.1) .. controls (1.9,-.1) .. (1.9,0);
\begin{scope}[xshift=1.2cm]
\draw [draw=white,double=black,very thick]  (0:1.1) arc (0:180:1.1) (0:.7) arc (0:180:.7);
\end{scope}
\begin{scope}[xshift=4.2cm]
\draw [-] (0:1.1) arc (0:180:1.1) (0:.7) arc (0:180:.7);
\draw (-.7,0) .. controls (-.7,-.1) ..  (-.3,-.1) .. controls (.1,-.1) .. (.1,0) (.5,0) .. controls (.5,-.1) ..  (.6,-.1) .. controls (.7,-.1) .. (.7,0) (1.1,0) .. controls (1.1,-.1) ..  (1.5,-.1) .. controls (1.9,-.1) .. (1.9,0);
\end{scope}
\begin{scope}[xshift=5.4cm]
\draw [draw=white,double=black,very thick]  (0:1.1) arc (0:180:1.1) (0:.7) arc (0:180:.7);
\end{scope}
\draw [-] (2.3,0) .. controls (2.3,-.1) ..  (2.5,-.1) .. controls (3.1,-.1) .. (3.1,0) (5.4,-1.1) .. controls (6,-1.1) and (6.5,-.6) .. (6.5,0)   (2.8,-.5) node{\scriptsize $\Sigma$} (-1.1,0) .. controls (-1.1,-.6) and (-.6,-1.1) .. (0,-1.1);
\draw [dashed] (0,-1.1) -- (5.4,-1.1);
\begin{scope}[xshift=9cm]
\draw [-] (0:1.1) arc (0:180:1.1) (0:.7) arc (0:180:.7);
\draw (-.7,0) .. controls (-.7,-.1) ..  (-.3,-.1) .. controls (.1,-.1) .. (.1,0) (.5,0) .. controls (.5,-.1) ..  (.6,-.1) .. controls (.7,-.1) .. (.7,0) (1.1,0) .. controls (1.1,-.1) ..  (1.5,-.1) .. controls (1.9,-.1) .. (1.9,0);
\begin{scope}[xshift=1.2cm]
\draw [draw=white,double=black,very thick]  (0:1.1) arc (0:180:1.1) (0:.7) arc (0:180:.7);
\end{scope}
\begin{scope}[xshift=4.2cm]
\draw [-] (0:1.1) arc (0:180:1.1) (0:.7) arc (0:180:.7);
\draw (-.7,0) .. controls (-.7,-.1) ..  (-.3,-.1) .. controls (.1,-.1) .. (.1,0) (.5,0) .. controls (.5,-.1) ..  (.6,-.1) .. controls (.7,-.1) .. (.7,0) (1.1,0) .. controls (1.1,-.1) ..  (1.5,-.1) .. controls (1.9,-.1) .. (1.9,0);
\end{scope}
\begin{scope}[xshift=5.4cm]
\draw [draw=white,double=black,very thick]  (0:1.1) arc (0:180:1.1) (0:.7) arc (0:180:.7);
\end{scope}
\draw [-] (2.3,0) .. controls (2.3,-.1) ..  (2.5,-.1) .. controls (3.1,-.1) .. (3.1,0) (5.4,-1.1) .. controls (6,-1.1) and (6.5,-.6) .. (6.5,0) (2.8,-.5) node{\scriptsize $\Sigma^{\prime}$} (-1.1,0) .. controls (-1.1,-.6) and (-.6,-1.1) .. (0,-1.1);
\draw [dashed] (0,-1.1) -- (5.4,-1.1);
\end{scope}
\draw (5.4,-1.1) -- (9,-1.1);
\draw (0,-1.1) .. controls (-.6,-1.1) .. (-.6,-1.4) .. controls (-.6,-1.7) .. (0,-1.7) -- (14.4,-1.7) .. controls (15,-1.7) .. (15,-1.4) .. controls (15,-1.1) .. (14.4,-1.1) (7.5,-1.4) node{\scriptsize $D$};
\end{tikzpicture}
\caption{$\tilde{\Sigma}$}\label{figsigmatildepseudo}
\end{figure}

Let $B =\left[0,6\right]\times \tilde{\Sigma}$. Let $\tau_{B,0}$ be a pseudo-parallelization of $B$ extending the pseudo-parallelization $\tau_0$ used both for $A$ and for $A^{\prime}=\left[0,6\right]\times \Sigma^{\prime}$. Let $\tau_{B,1}$ be a pseudo-parallelization that coincides with $\tau_1$ on $A$ and on $A^{\prime}$, and with $\tau_{B,0}$ on $\left[0,6\right]\times D$.
Then we have $z_n(A;\tau_0,\tau_1)=\frac12 z_n(B;\tau_{B,0},\tau_{B,1})$ and $p_1(\tau_0,\tau_1)=\frac12p_1(\tau_{B,0},\tau_{B,1})$.
Since the intersection of $(\gamma\cup\gamma^{\prime})$ and $(\gamma_0\cup\gamma^{\prime}_0)$ is zero modulo $2H_1(\partial A;\ZZ)$, the homology classes of these curves can be represented by curves that do not intersect. So we have $z_n(B;\tau_{B,0},\tau_{B,1})=\frac{p_1(\tau_{B,0},\tau_{B,1})}{4}\ansothree_n$.
\eop

\bpo{Proof of Theorem~\ref{thmznponepseudo}}
Let us first prove the theorem when $A$ is a rational homology ball, according to the schema of Figure~\ref{figproofpseudotwo}. Then there exists a parallelization $\tau_2$ of $A$ that coincides with $\tau_0=(N(\gamma);\tau_e,\tau_b)$ in a neighborhood of $\partial A$.
Thicken the neighborhood $N(\gamma)$ to $[a-7,b+7] \times \gamma \times [-2,2]$.
Add bands to $\gamma \times [-2,2]$ so that the disjoint union $\gamma \times [-2,2]$ is embedded in a connected oriented surface $\Sigma$ of $A$ with one boundary component. Let $[a-7,b+7] \times \Sigma$ be embedded in $A$ so that this parametrization matches the previous one.

\bfig
\centering
\begin{tikzpicture} 
\draw (-6,-.5) -- (-6,0) rectangle (-2.8,5) -- (-2.8,5.5)
(-6,5.5) -- (-6,5) 
(-2.8,0) -- (-2.8,-.5) 
(-6,2) -- (-2.8,2) 
(-6,3) -- (-2.8,3) 
(-2.8,1) -- (3,1) 
(3,-.5) -- (3,5.5) 
(-2.8,2) -- (.6,2)  
(-2.8,4) -- (3,4) 
(.6,4) -- (.6,1) 
(-5.9,5) node [left]{\scriptsize $\{b+7\} \times \Sigma$}
(-5.9,4) node [left]{\scriptsize $\{b+4\} \times \Sigma$}
(-5.9,3) node [left]{\scriptsize $\{b\} \times \Sigma$}
(-5.9,2) node [left]{\scriptsize $\{a\} \times \Sigma$}
(-5.9,1) node [left]{\scriptsize $\{a-4\} \times \Sigma$}
(-5.9,0) node [left]{\scriptsize $\{a-7\} \times \Sigma$} 
(-4.4,2.65)node {\scriptsize Trivial product over}
(-4.4,2.35) node {\scriptsize $\left[0,1\right] \times  [a,b] \times \Sigma$}
(-6,-.65) node {\scriptsize $(A,\tau_0)$}
(-2.8,-.65) node {\scriptsize $(A,\tau_3)$}
(2.7,-.65) node {\scriptsize $(A,\tau_2)$}
(-4.4,1.15) node{\scriptsize Box treated} (-4.4,.85) node{\scriptsize by Lemma~\ref{lempseudodoubletorus}}
(-4.4,4.15) node{\scriptsize Box treated} (-4.4,3.85) node{\scriptsize by Lemma~\ref{lempseudodoubletorus}}
(.6,.4) node{\scriptsize Box over $C$} (.6,.1) node{\scriptsize treated by Lemma~\ref{lembounsurfpseudo}}
(.6,4.9) node{\scriptsize Box over $C$} (.6,4.6) node{\scriptsize treated by Lemma~\ref{lembounsurfpseudo}}
(-1.1,3.15) node{\scriptsize Box treated} (-1.1,2.85) node{\scriptsize by Lemma~\ref{lemSigmapseudotwo}}
(1.8,2.65) node{\scriptsize Box treated} (1.8,2.35) node{\scriptsize by Lemma~\ref{lempseudodoubletorus}}
(-1.1,1.65) node{\scriptsize Trivial product over}
(-1.1,1.35) node{\scriptsize $\left[0,1\right] \times  [a-4,a] \times \Sigma$}
(-4.4,5.2) node{\scriptsize Trivial product}
(-4.4,-.2) node{\scriptsize Trivial product};
\begin{scope}[xshift=-6cm,yshift=2cm]
 \fill [white] (-.15,.3) rectangle (.15,.75);
\draw (-.15,.3) rectangle (.15,.75) (0,.5) node {\scriptsize $\gamma$};
\end{scope}
\begin{scope}[xshift=-2.8cm,yshift=2cm]
 \fill [white] (-.15,.3) rectangle (.15,.75);
\draw (-.15,.3) rectangle (.15,.75) (0,.5) node {\scriptsize $\gamma$};
\end{scope}
\begin{scope}[xshift=-2.8cm]
 \fill [white] (-.15,.3) rectangle (.15,.75);
\draw (-.15,.3) rectangle (.15,.75) (0,.5) node {\scriptsize $\gamma_a$};
\end{scope}
\begin{scope}[xshift=-2.8cm,yshift=1cm]
 \fill [white] (-.15,.3) rectangle (.15,.75);
\draw (-.15,.3) rectangle (.15,.75) (0,.5) node {\scriptsize $\gamma_a$};
\end{scope}
\begin{scope}[xshift=-2.8cm,yshift=3cm]
 \fill [white] (-.15,.3) rectangle (.15,.75);
\draw (-.15,.3) rectangle (.15,.75) (0,.5) node {\scriptsize $\gamma_b$};
\end{scope}
\begin{scope}[xshift=-2.8cm,yshift=4cm]
 \fill [white] (-.15,.3) rectangle (.15,.75);
\draw (-.15,.3) rectangle (.15,.75) (0,.5) node {\scriptsize $\gamma_b$};
\end{scope}
\begin{scope}[xshift=-.8cm]
\fill [white] (-2.15,-.15) rectangle (-1.85,.15);
\draw (-2,.0) node {\scriptsize $\tau_0$};
\end{scope}
\begin{scope}[xshift=-.8cm,yshift=2cm]
\fill [white] (-2.15,-.15) rectangle (-1.85,.15);
\draw (-2,.0) node {\scriptsize $\tau_0$};
\end{scope}
\begin{scope}[xshift=-.8cm,yshift=3cm]
\fill [white] (-2.15,-.15) rectangle (-1.85,.15);
\draw (-2,.0) node {\scriptsize $\tau_0$};
\end{scope}
\begin{scope}[xshift=-.8cm,yshift=5cm]
\fill [white] (-2.15,-.15) rectangle (-1.85,.15);
\draw (-2,.0) node {\scriptsize $\tau_0$};
\end{scope}
\begin{scope}[xshift=-.8cm,yshift=1cm]
\fill [white] (-2.15,-.15) rectangle (-1.85,.15);
\draw (-2,.0) node {\scriptsize $\tau_2$};
\end{scope}
\begin{scope}[xshift=-.8cm,yshift=4cm]
\fill [white] (-2.15,-.15) rectangle (-1.85,.15);
\draw (-2,.0) node {\scriptsize $\tau_2$};
\end{scope}
\begin{scope}[xshift=.6cm,yshift=1cm]
 \fill [white] (-.15,.3) rectangle (.15,.75);
\draw (-.15,.3) rectangle (.15,.75) (0,.5) node {\scriptsize $\gamma_a$};
\end{scope}
\begin{scope}[xshift=.6cm,yshift=3cm]
 \fill [white] (-.15,.3) rectangle (.15,.75);
\draw (-.15,.3) rectangle (.15,.75) (0,.5) node {\scriptsize $\gamma_a$};
\end{scope}
\begin{scope}[xshift=2.6cm,yshift=2.5cm]
\fill [white] (-2.15,-.15) rectangle (-1.85,.15);
\draw (-2,.0) node {\scriptsize $\tau_4$};
\end{scope}
\end{tikzpicture}
\caption{Schema of proof for Theorem~\ref{thmznponepseudo}}
\label{figproofpseudotwo}
\end{figure}

After a possible homotopy of $\tau_2$,
there exist annuli $\gamma_a \times [-1,1]$ and $\gamma_b \times [-1,1]$ in $\Sigma$ such that 
\begin{itemize}
\item $\tau_0$ coincides with $\tau_2$ in a neighborhood of $[a-7,b+7] \times \partial \Sigma$,
\item $\tau_{0}\vert_{\{a-4\} \times \Sigma}=\tau_2 \circ \CT_{\gamma_a}$, and
\item $\tau_{0}\vert_{\{b+4\} \times \Sigma}=\tau_2 \circ \CT_{\gamma_b}^{-1}$.
\end{itemize}

Set \begin{multline*}N(\gamma_3)= N(\gamma) \sqcup \bigl([b+2,b+3] \sqcup [b+5,b+6]\bigr) \times \gamma_b \times [-1,1]\\ 
\sqcup \bigl([a-6,a-5] \sqcup [a-3,a-2]\bigr) \times \gamma_a \times [-1,1].\end{multline*}
Let $\tau_3=(N(\gamma_3);\tau_{3,e},\tau_{3,b})$ be a pseudo-parallelization, which coincides with 
$\tau_0$ outside $([b+1,b+7]) \times \gamma_b \times [-1,1]
\sqcup ([a-7,a-1]) \times \gamma_a \times [-1,1]$, and which coincides with $\tau_2$ on $\{a-4\} \times \Sigma$ and on $\{b+4\} \times \Sigma$.
According to Lemma~\ref{lempseudodoubletorus} and to Proposition~\ref{propextdefzn}, we have
\begin{equation*}z_n(A;\tau_0,\tau_3)=\frac{p_1(\tau_0,\tau_3)}{4}\ansothree_n, \end{equation*}
as the left part of Figure~\ref{figproofpseudotwo} shows.

Set $B=[a-4,b+4] \times \Sigma$ and $C=A \setminus \left(\left]a-4,b+4\right[ \times \mathring{\Sigma}\right)$. We have
\begin{equation*}z_n(A;\tau_2,\tau_3)=z_n(B;\tau_2,\tau_3) + z_n\left(C;\tau_2,\tau_3\right),\end{equation*}
and $p_1$ decomposes similarly.

Let us now prove \begin{equation*}z_n\left(C;\tau_2,\tau_3\right)=\frac14 p_1(\tau_{2}\vert_{C},\tau_{3}\vert_{C})\ansothree_n.\end{equation*}
To do this, we apply Lemma~\ref{lembounsurfpseudo}, after an isotopy of $[b+4,b+7] \times \Sigma$ which sends $[b+5,b+6] \times \gamma_b \times [-1,1]$ to itself (at the end) so that $\{b+5\} \times \gamma_b \times [-1,1]$ is sent to $\{b+6\} \times \gamma_b \times \left(-[-1,1]\right)$, and $\{b+6\} \times \gamma_b \times [-1,1]$ is sent to $\{b+5\} \times \gamma_b \times \left(-[-1,1]\right)$.

Let $\tau_4=(([a-3,a-2] \sqcup [b+2,b+3]) \times \gamma_a \times [-1,1];\tau_{4,e},\tau_{4,b})$ be a pseudo-parallelization of $B$ that coincides with $\tau_2$ in a neighborhood of $\partial B$. Lemma~\ref{lempseudodoubletorus} implies \begin{equation*}z_n\left(B;\tau_4,\tau_2\right)=\frac14 p_1(\tau_4,\tau_{2}\vert_{B})\ansothree_n.\end{equation*}
Since $\left[\gamma_a\right] + \left[\gamma\right] + \left[\gamma_b\right]=0$ in $H_1(\Sigma;\ZZ/2\ZZ)$, Lemma~\ref{lemSigmapseudotwo} and Proposition~\ref{propextdefzn} imply $z_n\left(B;\tau_3,\tau_4\right)=\frac14 p_1(\tau_{3}\vert_{B},\tau_4)\ansothree_n$.
We get $z_n\left(B;\tau_2,\tau_3\right)=\frac14 p_1(\tau_{2}\vert_{B},\tau_3\vert_{B})\ansothree_n$, and hence
\begin{equation*}z_n(A;\tau_2,\tau_3)=\frac14 p_1(\tau_{2},\tau_{3})\ansothree_n.\end{equation*}
So we have $z_n(A;\tau_0,\tau_2)=\frac14 p_1(\tau_{0},\tau_{2})\ansothree_n$.
For the same reasons, we have $z_n(A;\tau_1,\tau_2)=\frac14 p_1(\tau_{1},\tau_{2})\ansothree_n$. Hence the lemma is proved when $A$ is a rational homology ball.

In general, the manifold $A$ is assumed to embed into a rational homology ball $B$, the pseudo-parallelization $\tau_{0}$ on $A$ extends to a pseudo-parallelization $\tilde{\tau}_0$ of $B$, and the pseudo-parallelization $\tau_{1}$ over $A$ extends to a pseudo-parallelization $\tilde{\tau}_1$ of $B$, which coincides with $\tilde{\tau}_0$ over $B\setminus \mathring{A}$.
We have
\begin{equation*}z_n(A;\tau_0,\tau_1)=z_n(B;\tilde{\tau}_0,\tilde{\tau}_1)=\frac{p_1(\tilde{\tau}_0,\tilde{\tau}_1)}{4}\ansothree_n= \frac{p_1(\tau_0,\tau_1)}{4}\ansothree_n.\end{equation*}
\eop

\section{Proof of \texorpdfstring{Theorem~\ref{thmfstconsttangpseudo}}{the main theorem}}
\label{secproofthmfstconsttangpseudo}

\begin{proposition}
\label{propfstconsttangpseudo}
Let $\crats$ be a rational homology $\RR^3$ equipped with an asymptotically standard pseudo-parallelization $\tau$. Let $\omega({\tau})$ be a homogeneous propagating form of $(C_2(\rats),{\tau})$. Let $n \in \NN$. With the notation of Corollary~\ref{cordefzinvuf} and Notation~\ref{notationzZ}, we have
\begin{equation*}\zinv_n\bigl(\crats,\omega({\tau})\bigr)=\zinvuf_n(\rats) +\frac14 p_1(\tau) \ansothree_n \end{equation*}
and
\begin{equation*}\Zinvuf(\rats)=\Zinv\bigl(\crats,\emptyset,\omega({\tau})\bigr)\exp\left(- \frac{p_1(\tau)}{4}\ansothree\right).\end{equation*}
\end{proposition}
\bp Let $\tau_0$ and $\tau_1$ be two pseudo-parallelizations of $\crats$ standard outside ${\ballb}_{\rats}$. Let $\omega$ be a $2$-form on $\left[0,1\right] \times \ST {\ballb}_{\rats}$ as in Lemma~\ref{lemformprod}.
Corollary~\ref{corinvone} implies
\begin{equation*}\zinv_n\bigl(\crats,\omega(\tau_1)\bigr)-\zinv_n\bigl(\crats,\omega(\tau_0)\bigr)=\zinv_n\bigl([0,1] \times \ST {\ballb}_{\rats}; \omega \bigr),\end{equation*}
while Theorem~\ref{thmznponepseudo} implies
\begin{equation*}\zinv_n\bigl([0,1] \times \ST {\ballb}_{\rats}; \omega \bigr)=\frac14\bigl(p_1(\tau_1)-p_1(\tau_0)\bigr) \ansothree_n.\end{equation*}
Conclude with Corollary~\ref{cordefzinvuf}.
\eop

\bpo{Proof of Theorem~\ref{thmfstconsttangpseudo}}
Theorem~\ref{thmfstconsttang} implies \begin{equation*}\Zinv(\hcylc,\tanghcyll,\tau)= \exp\left(\frac14 p_1(\tau)\ansothree\right)\prod_{j=1}^k\Bigl(\exp\bigl(I_{\theta}(K_j,\tau)\alpha\bigr)\#_j\Bigr) \Zinvuf(\hcylc,\tanghcyll)\end{equation*} for any actual parallelization $\tau$ of $\hcylc$. We want to prove the same equality when
$\tau$ is a pseudo-parallelization that is an actual parallelization over a tubular neighborhood $N(\tanghcyll)$ of $\tanghcyll$. Recall Lemma~\ref{lemmultcheck}. 
Proposition~\ref{propfstconsttangpseudo} leaves us with the proof that
\begin{equation*}\Zinvlink(\hcylc,\tanghcyll,\tau^{\prime})= \prod_{j=1}^k\Bigl(\exp\bigl(I_{\theta}(K_j,\tau^{\prime})\alpha\bigr)\#_j\Bigr) \Zinvlinkuf(\hcylc,\tanghcyll)\end{equation*} for any pseudo-parallelization $\tau^{\prime}$ of $\hcylc$ that is an actual parallelization over a tubular neighborhood $N(\tanghcyll)$ of $\tanghcyll$.

First assume that the restriction of $\tau^{\prime}$ to the tubular neighborhood $N(\tanghcyll)$ of $\tanghcyll$ extends to an actual parallelization $\tau$ of $\hcylc$. In this case, we apply Proposition~\ref{propinvtwofonc} with a closed $2$-form $\tilde{\omega}$ on $\left[0,1\right] \times \ST \hcylc$ as in its statement.
Thus, we get
\begin{equation*}\Zinvlink(\hcylc,\tanghcyll,\tau^{\prime}) = \left(\prod_{j=1}^k \exp\left(I_j\right) \#_j\right) \Zinvlink(\hcylc,\tanghcyll,\tau),\end{equation*}
where $I_j$ is defined from $\tilde{\omega}$ for any component $K_j$ of $\tanghcyll=\sqcup_{j=1}^k K_j$ in Proposition~\ref{propinvtwofonc}. Let $\projp_{\ST N(\tanghcyll)} \colon \left[0,1\right] \times \ST N(\tanghcyll) \to \ST N(\tanghcyll)$ denote the projection to the second factor.
According to Proposition Proposition~\ref{propextexistspseudo}, the form  $\tilde{\omega}$ may be expressed as $\projp_{\ST N(\tanghcyll)}^{\ast} \projp_{\tau}^{\ast}\left(\omega_{S^2}\right)$ over $\left[0,1\right] \times \ST N(\tanghcyll)$.
This factorization implies that the $I_j$ vanish. So, we have
$\Zinvlink(\hcylc,\tanghcyll,\tau^{\prime})=\Zinvlink(\hcylc,\tanghcyll,\tau)$. The degree one part of this equality implies $I_{\theta}(K_j,\tau^{\prime})=I_{\theta}(K_j,\tau)$. Hence the theorem is proved when the restriction of $\tau^{\prime}$ to the tubular neighborhood $N(\tanghcyll)$ of $\tanghcyll$ extends to an actual parallelization $\tau$ of $\hcylc$.

Recall that $\drad{r}$ is the disk of the complex numbers of module less than or equal to $r$.
Let $\taust$ denote both the standard parallelization of $\RR^3$ and its restriction to $\drad{4} \times [-2,2]$. Consider the neighborhood \begin{equation*}N(\gamma_2)=(D_3 \setminus \mathring{D}_1) \times [-1,1]\end{equation*} of $(\gamma_2=\partial D_2 \times \{0\})$. Let $\tau_2=(N(\gamma_2);\tau_e,\tau_b)$ be a pseudo-parallelization of $\drad{4} \times [-2,2]$ that coincides with $\taust$ in a neighborhood of $\partial (\drad{4} \times [-2,2])$ and that maps $e_3$ to the vertical direction of $\{0\} \times [-2,2]$ along $\{0\} \times [-2,2]$.

Let $\taust^{\hcylc}$ and $\tau_2^{\hcylc}$ be two pseudo-parallelizations of $\hcylc$ that satisfy the following set $(\ast)(\taust,\tau_2,\hcylc,\tanghcyll,K_j)$ of assumptions: They are actual parallelizations over a tubular neighborhood $N(\tanghcyll)$ of $\tanghcyll$.
There is an embedding of $\drad{4} \times [-2,2]$ in the rational homology cylinder $\hcylc$ equipped with the long tangle representative $\tanghcyll$ so that (the image of) $\drad{4} \times [-2,2]$ intersects (the image of) $\tanghcyll$ along $\{0\} \times [-2,2]$ and the orientations of $\{0\} \times [-2,2]$ and $\tanghcyll$ match. 
With respect to this embedding, $\taust^{\hcylc}$ and $\tau_2^{\hcylc}$ respectively coincide with $\taust$ and $\tau_2$ on the image of $\drad{4} \times [-2,2]$, and they coincide with each other elsewhere. The component of $\tanghcyll$ that intersects $\drad{4} \times [-2,2]$ is denoted by $K_j$.

According to Proposition~\ref{propextexistspseudo}, there exists a closed $2$-form $\tilde{\omega}$ on $\left[0,1\right] \times \ST \left(\drad{4} \times [-2,2]\right)$ that restricts to \begin{equation*}\Bigl(\{0\} \times \ST \bigl(\drad{4} \times [-2,2]\bigr)\Bigr) \cup \Bigl(\left[0,1\right] \times \ST\bigl(\drad{4} \times [-2,2]\bigr)\vert_{\partial \left(\drad{4} \times [-2,2]\right) }\Bigr)\end{equation*} as $\projp_{\taust}^{\ast}\left(\omega_{S^2}\right)$, and to $\{1\} \times \ST \left(\drad{4} \times [-2,2]\right)$ as a homogeneous propagating form of $C_2(\rats(\hcylc),\tau_2^{\hcylc})$ does. This closed $2$-form is actually independent of $(\hcylc,\tanghcyll)$, so is the induced quantity $I(\tau_2)$ 
of Proposition~\ref{propinvtwofonc} such that \begin{equation*}\Zinvlink(\hcylc,\tanghcyll,\tau_2^{\hcylc})=\exp\bigl(I(\tau_2)\bigr) \#_j \Zinvlink(\hcylc,\tanghcyll,\taust^{\hcylc}).\end{equation*}
In particular, the difference $(I_{\theta}(K_j,\tau_2^{\hcylc})-I_{\theta}(K_j,\taust^{\hcylc}))$ is a constant $\ell(\tau_2)$. It can be obtained from the degree one part of $I(\tau_2)$.

Apply this computation when $\rats(\hcylc)$ is $SO(3)$, $(\tanghcyll=K)$ is a knot, the homology class of $K$ represents the generator of $H_1(SO(3);\ZZ)=\ZZ/2\ZZ$, and  $\taust^{\hcylc}$ is an actual parallelization.
The two homotopy classes of parallelizations of $N(K)$ are obtained from one another by composition by the restriction to $N(K) \times \RR^3$ of the map \begin{equation*}\begin{array}{lll}SO(3) \times \RR^3 & \to& SO(3) \times \RR^3\\                                                                                                                                                                                                                                                         (\rho,x) & \mapsto & (\rho,\rho(x)).                                                                                                                                                                                                                                            \end{array}\end{equation*}
So, all parallelizations of $N(K)$ extend to $\hcylc$.
In particular, $\tau^{\hcylc}_{2}\vert_{N(K)}$ extends to $\hcylc$ as a parallelization standard near $\partial \hcylc$.

Then the first studied case implies $I(\tau_2)=\left(I_{\theta}(K,\tau_2^{\hcylc})-I_{\theta}(K,\taust^{\hcylc})\right)\alpha=\ell(\tau_2)\alpha$. So we get
$I(\tau_2)=\left(I_{\theta}(K_j,\tau_2^{\hcylc})-I_{\theta}(K_j,\taust^{\hcylc})\right)\alpha$, and hence
\begin{multline*}\prod_{\ell=1}^k\Bigl(\exp\bigl(-I_{\theta}(K_{\ell},\tau_2^{\hcylc})\alpha\bigr)\#_{\ell}\Bigr)\Zinvlink(\hcylc,\tanghcyll,\tau_2^{\hcylc})\\
=\prod_{\ell=1}^k\Bigl(\exp\bigl(-I_{\theta}(K_{\ell},\taust^{\hcylc})\alpha\bigr)\#_{\ell}\Bigr)\Zinvlink(\hcylc,\tanghcyll,\taust^{\hcylc})\end{multline*}
for any two pseudo-parallelizations  $\taust^{\hcylc}$ and $\tau_2^{\hcylc}$ that satisfy
$(\ast)(\taust,\tau_2,\hcylc,\tanghcyll,K_j)$.

Let $\tau^{\prime}$ be a pseudo-parallelization of $\hcylc$ that coincides with an actual parallelization $\tau_N$ of $\hcylc$ on $N(\tanghcyll)$. Let $\tau$ be a parallelization of $\hcylc$. The restrictions of $\tau_N$ and $\tau$ to $\partial N(\tanghcyll)$ are homotopic along the meridians of $\tanghcyll$. They differ by the action of the generator of $\pi_1(SO(2))$ along parallels on components $K_j$ for $K_j$ in some finite set $\finseta$. If $\finseta=\emptyset$, the first studied case implies 
\begin{equation*} \prod_{j=1}^k\Bigl(\exp\bigl(-I_{\theta}(K_j,\tau^{\prime})\alpha\bigr)\#_j\Bigr) \Zinvlink(\hcylc,\tanghcyll,\tau^{\prime})=\Zinvlinkuf(\hcylc,\tanghcyll).\end{equation*}

Otherwise, equip each component $K_j$ of $\finseta$ with one embedding of $\drad{4} \times [-2,2]$ whose image meets $\tanghcyll$ in $K_j$ along $\{0\} \times [-2,2]$, so that the orientations of $\{0\} \times [-2,2]$ and $\tanghcyll$ match.
Perform a homotopy of $\tau^{\prime}$ as in Definition~\ref{defhomotoppseudo} to transform $\tau^{\prime}$ to a pseudo-parallelization $\tau^{\prime \prime}$ such that $\tau^{\prime \prime}$ is induced by $\taust$ on the image of the above embeddings.
For $\finsetb \subseteq \finseta$, let $\tau^{\prime \prime}_{\finsetb}$ be obtained from $\tau^{\prime \prime}$
by changing $\taust$ to $\tau_2$ on the images of the embeddings of $\drad{4} \times [-2,2]$ that meet an element of $\finsetb$.
Then for $\finsetb \subseteq \finseta$ and for $K_{\ell} \in \finseta \setminus \finsetb$, the pseudo-parallelizations $\tau^{\prime \prime}_{\finsetb}$ and $\tau^{\prime \prime}_{\finsetb \cup \{K_{\ell}\}}$ satisfy  $(\ast)(\tau^{\prime \prime}_{\finsetb},\tau^{\prime \prime}_{\finsetb \cup \{K_{\ell}\}},\hcylc,\tanghcyll,K_{\ell})$. Therefore, we have
\begin{multline*}\prod_{j=1}^k\Bigl(\exp\bigl(-I_{\theta}(K_j,\tau^{\prime \prime}_{\finsetb \cup \{K_{\ell}\}})\alpha\bigr)\#_j\Bigr)\Zinvlink(\hcylc,\tanghcyll,\tau^{\prime \prime}_{\finsetb \cup \{K_{\ell}\}})\\=\prod_{j=1}^k\Bigl(\exp\bigl(-I_{\theta}(K_j,\tau^{\prime \prime}_{\finsetb})\alpha\bigr)\#_j\Bigr)\Zinvlink(\hcylc,\tanghcyll,\tau^{\prime \prime}_{\finsetb}).\end{multline*}
We get
\begin{equation*}\prod_{j=1}^k\Bigl(\exp\bigl(-I_{\theta}(K_j,\tau^{\prime \prime}_{\finseta})\alpha\bigr)\#_j\Bigr)\Zinvlink(\hcylc,\tanghcyll,\tau^{\prime \prime}_{\finseta})=\prod_{j=1}^k\Bigl(\exp\bigl(-I_{\theta}(K_j,\tau^{\prime \prime})\alpha\bigr)\#_j\Bigr)\Zinvlink(\hcylc,\tanghcyll,\tau^{\prime \prime}) \end{equation*}
by induction on $\cardlef{\finseta}$.
Since $\tau^{\prime \prime}_{\finseta}$ and $\tau$ are homotopic on $N(\tanghcyll)$, the first studied case
implies \begin{equation*}\prod_{j=1}^k\Bigl(\exp\bigl(-I_{\theta}(K_j,\tau^{\prime \prime}_{\finseta})\alpha\bigr)\#_j\Bigr)\Zinvlink(\hcylc,\tanghcyll,\tau^{\prime \prime}_{\finseta})=\Zinvlinkuf(\hcylc,\tanghcyll).\end{equation*}
So it suffices to prove that 
\begin{equation*}\prod_{j=1}^k\Bigl(\exp\bigl(-I_{\theta}(K_j,\tau^{\prime \prime})\alpha\bigr)\#_j\Bigr)\Zinvlink(\hcylc,\tanghcyll,\tau^{\prime \prime})=\prod_{j=1}^k\Bigl(\exp\bigl(-I_{\theta}(K_j,\tau^{\prime })\alpha\bigr)\#_j\Bigr)\Zinvlink(\hcylc,\tanghcyll,\tau^{\prime})\end{equation*} 
or that 
$\prod_{j=1}^k\left(\exp(-I_{\theta}(K_j,\tau^{\prime })\alpha)\#_j\right)\Zinvlink(\hcylc,\tanghcyll,\tau^{\prime})$ is invariant by a homotopy of $\tau^{\prime}$, as in Definition~\ref{defhomotoppseudo}, supported in a ball
where $\tau^{\prime}$ is a genuine parallelization (namely, around an image of $\drad{4} \times [-2,2]$).

Again, the effect on $\Zinvlink(\hcylc,\tanghcyll,\tau^{\prime})$ of such a homotopy depends only on the homotopy
inside the ball, according to Propositions~\ref{propinvtwofonc} and \ref{propextexistspseudo}. Since such a ball equipped with the homotopy may be inserted in a tangle equipped with a genuine trivialization $\tau$, we conclude that \begin{equation*}\prod_{j=1}^k\left(\exp(-I_{\theta}(K_j,\tau^{\prime })\alpha)\#_j\right)\Zinvlink(\hcylc,\tanghcyll,\tau^{\prime})\end{equation*} is indeed invariant under the above homotopies. This is sufficient to conclude the proof of Theorem~\ref{thmfstconsttangpseudo}.
\eop

We give
more definitions of $\Zinv$ involving non-necessarily homogeneous propagating forms and pseudo-parallelizations in Chapter~\ref{chappseudoparmuchmore}.

\chapter{Simultaneous normalization of propagating forms}
\label{chapsimnormprop}

This chapter is devoted to the proof of Propositions~\ref{propnormasim} and \ref{propnorma}.
As shown in Section~\ref{secskproofthmmainunivlag}, this is sufficient to prove Theorem~\ref{thmmainunivlag}.
In this chapter, we use real coefficients for homology and cohomology unless otherwise mentioned.

\section{Sketch}

First note that the homogeneous boundary form of Definition~\ref{defhombounformpseudo} defined on $\partial C_2(\rats)$
(or the form $\projp_{\partau}^{\ast}\left(\omega_{S^2}\right)$ as in 
Definition~\ref{defpropagatorone} when pseudo-parallelizations are not involved)
is antisymmetric on $\partial C_2(\rats)$ as in Definition~\ref{defantisympropform}. So it extends as a closed antisymmetric 2-form $\omega=\omega_{\emptyset}$ on $C_2(\rats)$ as in 
Lemma~\ref{lemiotaform}.
Also note that if the restriction of $\omega_I$ to 
\begin{equation*}A_I^{(i)} \times \bigl(C_1(\rats_I) \setminus A^{(i)}_{I,3}\bigr) \subset C_2(\rats_I)\end{equation*}
equals
\begin{equation*}\sum_{j \in \underline{g_i}} p_1^{\ast}\left(\eta_I(a^i_j)\right) \wedge p_2^{\ast}\left(\eta_I(z^i_j)\right) + p_2^{\ast}\left(\omega(p^i)\right)\end{equation*} 
as stated in Proposition~\ref{propnorma}, then
the restriction of $\omega_I$ to $A_I^{(i)} \times A_I^{(k)}$ equals
\begin{equation*}\sum_{(j,\ell) \in \underline{g_i} \times \underline{g_k}}lk\left(z^i_j,z^k_{\ell}\right)\projp_{A_I^{(i)}}^{\ast}\left(\eta(a^i_j)\right) \wedge \projp_{A_I^{(k)}}^{\ast}\left(\eta(a^k_{\ell})\right),\end{equation*}
for $k\neq i$, as required in Proposition~\ref{propnormasim}.

To arrange the propagating forms $\omega_I$ as in Propositions~\ref{propnormasim} and \ref{propnorma}, we will first show how to make $\omega$ satisfy the conditions of Proposition~\ref{propnorma}, with respect to the notation before Proposition~\ref{propnorma}.
More precisely, we will prove the following proposition in Subsection~\ref{proofnorma}.

\begin{proposition}
\label{propnormababy}
Let $\tilde{\omega}$ be a propagating form of $C_2(\rats)$ as in Definition~\ref{defpropagatortwo}. Its restriction to $\partial C_2(\rats)\setminus \ST \ballb_{\rats}$ may be expressed as $\projp_{\tau}^{\ast}\left(\tilde{\omega}_{S^2}\right)$, for some volume-one form $\tilde{\omega}_{S^2}$ of $S^2$. Let $\tilde{\omega}(p^i)$ (resp. $\tilde{\omega}(p^i)_{\iota}$) be a degree-two form on $\left(C_1(\rats) \setminus \Int(A^{(i)})\right)$ that satisfies the same properties as the form $\omega(p^i)$ (introduced before Proposition~\ref{propnorma}) except that it restricts to $\partial C_1(\rats)=S^2$ as $\tilde{\omega}_{S^2}$ (resp. as $-\iota_{S^2}^{\ast}\tilde{\omega}_{S^2}$) instead of the usual volume form $\omega_{S^2}$.\footnote{In our proof of Propositions~\ref{propnormasim} and \ref{propnorma}, we use Proposition~\ref{propnormababy} only when $\tilde{\omega}_{S^2}=\omega_{S^2}$ and $\tilde{\omega}(p^i)=\tilde{\omega}(p^i)_{\iota}=\omega(p^i)$, but the general statement is useful in other related work.} If $\tilde{\omega}$ is antisymmetric, then assume that $\tilde{\omega}(p^i)_{\iota}=\tilde{\omega}(p^i)$.

Then there exists a propagating form $\omega$ of $C_2(\rats)$ such that
\begin{enumerate} 
\item the form $\omega$ coincides with $\tilde{\omega}$ on $\partial C_2(\rats)$,
\item for every $i \in \underline{x}$, the restriction of $\omega$ to 
\begin{equation*}A^{(i)} \times \bigl(C_1(\rats) \setminus A^{(i)}_{3}\bigr) \subset C_2(\rats)\end{equation*}
is equal to \begin{equation*}\sum_{j=1}^{g_i} p_1^{\ast}\left(\eta(a^i_j)\right) \wedge p_2^{\ast}\left(\eta(z^i_j)\right) + p_2^{\ast}\bigl(\tilde{\omega}(p^i)\bigr),\end{equation*}
where $p_1$ and $p_2$ respectively denote the first and the second projection of $A^{(i)} \times (C_1(\rats) \setminus A^{(i)}_{3})$ to $C_1(\rats)$,
and the restriction of $\omega$ to 
\begin{equation*}\bigl(C_1(\rats) \setminus A^{(i)}_{3}\bigr) \times A^{(i)} \subset C_2(\rats)\end{equation*}
is equal to \begin{equation*}\sum_{j \in \underline{g_i}} p_1^{\ast}\left(\eta(z^i_j)\right) \wedge p_2^{\ast}\left(\eta(a^i_j)\right) - p_1^{\ast}\bigl(\tilde{\omega}(p^i)_{\iota}\bigr),\end{equation*}
\item for every $i \in \underline{x}$, for any $j \in \{1,2, \dots, g_i\}$, we have \begin{equation*}\int_{\Sigma(a^i_j) \times p^i} \omega=0\;\;\;\; \mbox{and} \;\;\;\; \int_{p^i \times \Sigma(a^i_j)} \omega=0,\end{equation*}
where $p^i \in \partial A^{(i)}$ and $\partial \Sigma(a^i_j) \subset \{4\} \times \partial A^{(i)}$, and
\item $\omega$ is antisymmetric if $\tilde{\omega}$ is.
\end{enumerate}

\end{proposition}

Assume that Proposition~\ref{propnormababy} is proved. This is the goal of Subsection~\ref{proofnorma}. 
Recall that we will use Proposition~\ref{propnormababy} only when $\tilde{\omega}_{S^2}=\omega_{S^2}$ and $\tilde{\omega}(p^i)=\tilde{\omega}(p^i)_{\iota}=\omega(p^i)$. Also recall that the $\eta(a^i_j)$ are defined both in $A^{(i)}$ and $A^{(i)\prime}$ and they are identical near $\partial A^{(i)}$ and $\partial A^{(i)\prime}$.
Finally recall that $\omega(p^i)$ is supported outside $\cup_{k \in \underline{x}} \Int(A^{(k)})$ and the $\eta(z^i_j)$ restrict to the $A^{(k)}$ as a combination of $\eta(a^k_{\ell})$ (fixed by the linking numbers).
When changing some $A^{(i)}$ into some $A^{(i)\prime}$ with the same Lagrangian, it is easy to change the restrictions of $\omega$ inside the parts ($A^{(r)} \times (C_1(\rats) \setminus A^{(r)}_{3})$ or $(C_1(\rats) \setminus A^{(r)}_{3}) \times A^{(r)} \subset C_2(\rats)$) mentioned in the statement of Proposition~\ref{propnormababy}.
 Indeed, all the forms $\eta(a^i_j)$, $\eta(z^i_j)$, and $\omega(p^i)$ can be defined on the parts of the $\rats_I$ where they are needed so that these forms coincide with each other whenever it makes sense, and so that they have the properties that were required for $\rats$.
Define $\omega_0(\rats_I)$ on
\begin{multline*}
D\bigl(\omega_0(\rats_I)\bigr)= \index[N]{Dfomegazero@$D(\omega_0(\rats_I))$ domain in $C_2$} \\
\biggl(C_2(\rats_I) \setminus \Bigl( \cup_{i \in I} \pbl^{-1}\bigl((A^{(i)\prime}_{-1} \times A^{(i)\prime}_3) \cup (A^{(i)\prime}_3 \times A^{(i)\prime}_{-1}) \bigr) \Bigr)  \biggr) \cup \pbl^{-1}\bigl(\diagon(\crats_I)\bigr) \end{multline*}
so that
\begin{enumerate}
\item we have $\omega_0(\rats_I)=\omega$ on $C_2\bigl(\rats \setminus (\cup_{i \in I}A^{(i)\prime}_{-1})\bigr)$,
\item when $i \in I$, we have \begin{equation*}\omega_0(\rats_I)=\sum_{i=1}^{g_i} p_1^{\ast}\left(\eta(a^i_j)\right) \wedge p_2^{\ast}(\eta(z^i_j)) + p_2^{\ast}\left(\omega(p^i)\right)\end{equation*} on $A^{(i)\prime} \times (\crats_I \setminus A^{(i)\prime}_3)$,
\item  when $i \in I$, we have $\omega_0(\rats_I)=-\iota^{\ast}\bigl(\omega_0(\rats_I)\bigr)$ on $(\crats_I \setminus A^{(i)\prime}_3) \times A^{(i)\prime} $, and
\item  on $\partial C_2(\rats_I)$, the form $\omega_0(\rats_I)$ coincides with the homogeneous boundary form $\omega({\tau}_I,F)$ of Definition~\ref{defhombounformpseudo}, for a map $F$ that is the same for all $I \subseteq \underline{x}$.
\end{enumerate}
Note that this definition is consistent. 

Set $\rats_i=\rats_{\{i\}}$ and $D_A(\omega_0(\rats_i))=C_2(A^{(i)\prime}_4) \cap D(\omega_0(\rats_i))$.
\begin{lemma}
\label{lemker}
With the above notation, for any $i \in \underline{x}$, the cohomology class of
$\omega_0(\rats_i)$ vanishes on the kernel
 of the map
\begin{equation*}H_2\Bigl(D_A\bigl(\omega_0(\rats_i)\bigr)\Bigr) \longrightarrow H_2\Bigl(C_2\bigl(A^{(i)\prime}_4\bigr)\Bigr)\end{equation*} induced by the inclusion.
\end{lemma}
 
This lemma was surprisingly difficult to prove for me. It will be proved in Subsection~\ref{secprooflemker}. Assume it for the moment. Then (the cohomology class of) $\omega_0(\rats_i)$  is in the image of the natural map
\begin{equation*}H^2\Bigl(C_2\bigl(A^{(i)\prime}_4\bigr)\Bigr) \longrightarrow H^2\Bigl(D_A\bigl(\omega_0(\rats_i)\bigr)\Bigr).\end{equation*}
So  $\omega_0(\rats_i)$ extends to a closed form $\omega_1(\rats_i)$ on $C_2(A^{(i)\prime}_4)$.
Change this form to $\omega_{\{i\}}=\frac12\bigl(\omega_1(\rats_i) - \iota^{\ast}\left(\omega_1(\rats_i)\right)\bigr)$ to get an antisymmetric homogeneous propagating form of $\bigl(C_2(\rats_i),\tau_{\{i\}}\bigr)$.
For any $I \subseteq \underline{x}$, define
\begin{equation*}\omega_I= \left\{ \begin{array}{ll} \omega_{0}(\rats_I) \; &\mbox{on}\; C_2(\rats_I) \setminus \Bigl( \cup_{i \in I} \pbl^{-1}\bigl((A^{(i)\prime}_{-1} \times A^{(i)\prime}_4) \cup (A^{(i)\prime}_4 \times A^{(i)\prime}_{-1}) \bigr) \Bigr)\\
\omega_{\{i\}} &\mbox{on}\;C_2(A^{(i)\prime}_4) \; \mbox{for}\; i \in I.
\end{array}\right.\end{equation*} 
This definition is consistent since the $C_2(A^{(i)\prime}_4)$ do not intersect.
The forms $\omega_I$ satisfy the properties of
Proposition~\ref{propnormasim}. In order to finish the proof of Proposition~\ref{propnorma}, up to Proposition~\ref{propnormababy}
and Lemma~\ref{lemker},
let us prove that $\int_{\Sigma^{\prime}(a^i_j) \times p^i} \omega_{\{i\}}=0$.
Note
\begin{equation*}\int_{\Sigma^{\prime}(a^i_j) \times p^i} \omega_{\{i\}}=\int_{\Sigma^{\prime}(a^i_j) \times (p^i\times \{4\})} \omega_{\{i\}} +  \int_{\partial \Sigma^{\prime}(a^i_j) \times \left(p^i\times \left[0,4\right]\right)} \omega_{\{i\}}.\end{equation*} 
The same formula applied to $\Sigma(a^i_j)$ instead of $\Sigma^{\prime}(a^i_j)$ yields
\begin{equation*}0=\int_{\Sigma(a^i_j)  \times (p^i\times \{4\})} \omega +  \int_{\partial \Sigma(a^i_j) \times (p^i\times \left[0,4\right])} \omega. \end{equation*}
The prescribed behavior of the forms on $A_{I}^{(i)} \times \bigl(\crats_I \setminus A_{I3}^{(i)}\bigr)$ implies
\begin{equation*}\int_{(\Sigma^{\prime}(a^i_j)\cap A^{(i)\prime}_0) \times (p^i\times \{4\})} \omega_{\{i\}}=\int_{(\Sigma(a^i_j)\cap A^{(i)}_0) \times (p^i\times \{4\})} \omega=0.\end{equation*} 
Since $\omega$ and $\omega_{\{i\}}$ coincide on $C_2\left(A_4^{(i)}\setminus A^{(i)}_{-1}\right)$, we also have 
\begin{equation*}\int_{\left(\left(-\partial \Sigma(a^i_j)\times \left[0,4\right]\right) \times (p^i\times \{4\})\right) \cup \left(\partial \Sigma(a^i_j) \times \left(p^i\times \left[0,4\right]\right)\right)} (\omega_{\{i\}} - \omega)=0.\end{equation*}
This shows $\int_{\Sigma^{\prime}(a^i_j) \times p^i} \omega_{\{i\}}=0$.

Thus, Proposition~\ref{propnormababy}
and Lemma~\ref{lemker} imply Propositions~\ref{propnormasim} and \ref{propnorma}. Their proofs occupy the next two subsections.

\section{Proof of \texorpdfstring{Proposition~\ref{propnormababy}}{the normalization proposition}}
\label{proofnorma}
The homology classes of the $(z_j^i \times (4 \times a_k^i))_{(j,k) \in \{1,\dots,g_i\}^2}$ and $(p^i \times \partial C_1(\rats))$ form a basis of
\begin{equation*}H_2\Bigl(A^{(i)} \times \bigl(C_1(\rats) \setminus A^{(i)}_3\bigr)\Bigr)=\Bigl(H_1(A^{(i)}) \otimes H_1(\rats \setminus A^{(i)})\Bigr) \oplus H_2\Bigl(C_1(\rats) \setminus A^{(i)}\Bigr).\end{equation*}
According to Lemma~\ref{lemlkprop}, the evaluation of the cohomology class of any propagating form of $C_2(\rats)$
at these classes is $lk(z_j^i, (4 \times a_k^i))=\delta_{kj}$ for the first ones
and $1$ for the last one. 
In particular, the form of the statement integrates correctly on this basis.

Let us first prove Proposition~\ref{propnormababy} when $\underline{x}=\{1\}$. Set $A^1=A$, and forget about the 
superfluous superscripts $1$. Let $\omega_0$ be a propagating two-form of $C_2(\rats)$ that restricts to $\partial C_2(\rats)\setminus \ST \ballb_{\rats}$ as $\projp_{\tau}^{\ast}\left(\tilde{\omega}_{S^2}\right)$,
and let $\omega_b$ be the closed 2-form defined on 
$\left(A_1 \times (C_1(\rats) \setminus \Int(A_2))\right)$ by the statement (extended naturally). Since this form $\omega_b$ integrates correctly on \begin{equation*}H_2\Bigl(A_1 \times \bigl(C_1(\rats) \setminus \Int(A_2)\bigr)\Bigr),\end{equation*} there exists a one-form $\eta$ on
$\left(A_1 \times (C_1(\rats) \setminus \Int(A_2))\right)$ such that
$\omega_b=\omega_0 +d \eta$.

This form $\eta$ is closed on $A_1 \times \partial C_1(\rats)$.
Since $H^1\left(A_1 \times (C_1(\rats) \setminus \Int(A_2))\right)$ maps surjectively to 
$H^1(A_1 \times \partial C_1(\rats))$, we may extend $\eta$ as a closed one-form $\tilde{\eta}$
on $\left(A_1 \times (C_1(\rats) \setminus \Int(A_2))\right)$. Changing $\eta$ into 
$(\eta-\tilde{\eta})$ turns $\eta$ to a primitive of $(\omega_b-\omega_0)$
that vanishes on $A_1 \times \partial C_1(\rats)$.
Let $\chi\colon C_2(\rats) \to \left[0,1\right]$ be a smooth function supported in 
$\left(A_1 \times \left(C_1(\rats) \setminus \Int(A_2)\right)\right)$ and constant with the value 
$1$ on $\left(A \times \left(C_1(\rats) \setminus A_3\right)\right)$.
Set \begin{equation*}\omega_a=\omega_0 +d \chi \eta.\end{equation*}
Then $\omega_a$ is a closed form that has the required form on $(A \times (C_1(\rats) \setminus A_3))$. Furthermore, the restrictions of 
$\omega_a$ and $\omega_0$ agree on $\partial C_2(\rats)$ since $d \chi \eta$
vanishes there (because $\eta$ vanishes on $A_1 \times \partial C_1(\rats)$).

Adding to $\eta$ a combination $\eta_c$ of the closed forms $p_2^{\ast}\left(\eta(z_j)\right)$,
which vanish on $A_1 \times \partial C_1(\rats)$, does not change the above properties,
but adds 
\begin{equation*}\int_{p \times ( [2,4] \times a_j)} d (\chi \eta_c)=\int_{p \times (4 \times a_j)} \eta_c\end{equation*}
to
$\int_{p \times \Sigma(a_j)} \omega_a$.
Therefore, since the $p_2^{\ast}\left(\eta(z_j)\right)$ generate the dual of $\CL_A$, we may choose $\eta_c$ so that all the $\int_{p \times \Sigma(a_j)} \omega_a$ vanish.
After this step, $\omega_a$ is a closed form, which takes the prescribed values on 
\begin{equation*}PS_a=\partial C_2(\rats) \cup \bigl(A \times (C_1(\rats) \setminus A_3)\bigr),\end{equation*}
and such that all the $\int_{p \times \Sigma(a_j)} \omega_a$ vanish.
To make $\omega_a$ take the prescribed values on $\iota(PS_a)$ and integrate as required on $\Sigma(a_j) \times p$,
we apply 
similar modifications to $\omega_a$ on the symmetric part $(C_1(\rats) \setminus \Int(A_2)) \times A_1$. The support of these modifications is disjoint
from the support of the previous ones. Therefore, they do not interfere and transform 
$\omega_a$ into a closed form $\omega_b$ with the following additional properties:
\begin{itemize}
\item the form $\omega_b$ has the prescribed form on $(C_1(\rats) \setminus A_3) \times A$, 
\item we have $\int_{\Sigma(a_j) \times p } \omega_b =0$, for all $j=1, \dots g_1$.
\end{itemize}

Thus, the form $\omega=\omega_b$ (resp. $\omega=\frac{\omega_b - \iota^{\ast}\left(\omega_b\right)}{2}$ if antisymmetry is desired)
has all the required properties, and Proposition~\ref{propnormababy} is proved for $\underline{x}=\{1\}$.

We now proceed by induction on $x$.
We start with a $2$-form $\omega_0$ that satisfies all the hypotheses with $\underline{x-1}$ instead of $\underline{x}$. By the first
step, we also assume that we have a $2$-form $\omega_b$ satisfying all the hypotheses with $\{x\}$ instead of $\underline{x}$, with the enlarged  $A_1^{(x)}$ replacing $A^{(x)}$. 

Now, we proceed similarly.
There exists a one-form $\eta$ on
$C_2(\rats)$ such that
$\omega_b=\omega_0 +d \eta$.
The exact sequence
\begin{equation*} 0=H^1\bigl(C_2(\rats)\bigr) \longrightarrow H^1\bigl(\partial C_2(\rats)\bigr) \longrightarrow H^2\bigl(C_2(\rats), \partial C_2(\rats)\bigr) \cong H_4\bigl(C_2(\rats)\bigr)=0 \end{equation*}
implies that $H^1(\partial C_2(\rats))$ is trivial. Therefore, the form $\eta$ is exact
on $\partial C_2(\rats)$. Thus, we may assume that $\eta$ vanishes on $\partial C_2(\rats)$, which we do. Let  \begin{equation*}\chi\colon C_2(\rats) \to \left[0,1\right]\end{equation*} be a smooth function supported in 
$\bigl(A^{(x)}_1 \times (C_1(\rats) \setminus \Int(A^{(x)}_2))\bigr)$ and constant with the value 
$1$ on $\bigl(A^{(x)} \times (C_1(\rats) \setminus A^{(x)}_3)\bigr)$.
Again, we are going to modify $\eta$ by some closed forms
so that \begin{equation*}\omega_a=\omega_0 +d \chi \eta\end{equation*}
has the prescribed value on 
\begin{multline*}PS_a=\partial C_2(\rats) \cup \biggl(\cup_{k=1}^x \Bigl(A^{(k)} \times \bigl(C_1(\rats) \setminus A^{(k)}_3\bigr)\Bigr) \biggr) \\ 
\cup \biggl(\cup_{k=1}^{x-1}\Bigl( \bigl(C_1(\rats) \setminus A^{(k)}_3\bigr) \times A^{(k)}\Bigr)\biggr).\end{multline*}
Our form $\omega_a$ is as required everywhere except possibly on
\begin{equation*}\Bigl(A^{(x)}_1 \times \bigl(C_1(\rats) \setminus \Int(A^{(x)}_2)\bigr)\Bigr) \setminus
\Bigl(A^{(x)} \times \bigl(C_1(\rats) \setminus A^{(x)}_3\bigr)\Bigr)\end{equation*}
and in particular on the intersection of this domain with the domains 
where it was normalized previously, which is included
in \begin{equation*}A^{(x)}_1 \times \bigl(\partial C_1(\rats) \cup (\cup_{k=1}^{x-1}A^{(k)})\bigr).\end{equation*}

Recall that $\eta$ vanishes on $A^{(x)}_1 \times \partial C_1(\rats)$.
Our assumptions also imply that $\eta$ is closed on $A^{(x)}_1 \times A^{(k)}$, for any $k <x$.
Let us prove that they imply that $\eta$ is exact on $A^{(x)}_1 \times A^{(k)}$, for any $k <x$.
To do that, it suffices to check the following two assertions.
\begin{enumerate}
\item For any $j =1,\dots,g_x$, we have
$\int_{z^x_j \times p^k} \eta=0.$
\item For any $j =1,\dots,g_k$, we have
$\int_{p^x \times z^k_j} \eta=0.$
\end{enumerate}
Let us prove the first assertion. Let $\infty(v) \in \partial C_1(\rats)$, and let $\left[p^k,\infty(v)\right]$ be a path from $p^k$ to $\infty(v)$ in $C_1(\rats)$ that intersects $\hcylc$ like the path $\left[p^k,q^k\right]$ introduced before Proposition~\ref{propnorma}.
Since $\int_{z^x_j \times \infty(v)} \eta=0$, we have
\begin{equation*}\int_{z^x_j \times p^k} \eta= \int_{\partial (z^x_j \times \left[p^k,\infty(v)\right])}
\eta=\int_{ z^x_j \times \left[p^k,\infty(v)\right]}
(\omega_b-\omega_0),\end{equation*}
where $\int_{ z^x_j \times \left[p^k,\infty(v)\right]}
\omega_b=0$ because the supports of the $\eta(z^x_{\ell})$ do not intersect $\left[p^k,\infty(v)\right]$. 
Let us compute \begin{equation*}\int_{ z^x_j \times \left[p^k,\infty(v)\right]}
\omega_0= - \int_{ \Sigma(z^x_j) \times \partial \left[p^k,\infty(v)\right]} \omega_0
=\int_{ \Sigma(z^x_j) \times \{p^k\}} \omega_0.\end{equation*}
The last integral vanishes because
\begin{enumerate}
\item the surface $\Sigma(z^x_j)$ intersects $A^{(k)}_4$ as copies of $\Sigma(a^k_{\ell})$,
\item we have $\int_{ \Sigma(a^k_{\ell}) \times p^k}\omega_0=0$, thanks to the third condition of Proposition~\ref{propnormababy}, and
\item the integral of $\omega_0$ also vanishes on the remaining part of $\Sigma(z^x_j) \times p^k$ since $\omega_0$ is determined on $((C_1(\rats) \setminus A^{(k)}_4) \times A^{(k)})$ and since the support of $\omega(p^k)$ is disjoint from $\Sigma(z^x_j)$.
\end{enumerate}
Let us prove the second assertion, namely $\int_{p^x \times z^k_j} \eta=0$ for $j \in \underline{g_k}$.
Again, we have $\int_{\infty(v) \times z^k_j} \eta=0$ since $\eta$ vanishes on $\partial C_2(\rats)$. We get
\begin{equation*}\int_{p^x \times z^k_j} \eta=- \int_{ \left[p^x,\infty(v)\right] \times z^k_j }
(\omega_b-\omega_0).\end{equation*}
The integral $\int_{ \left[p^x,\infty(v)\right] \times z^k_j }
\omega_0$ is zero because of the behavior of $\omega_0$ on $(C_1(\rats) \setminus A^{(k)}_4) \times A^{(k)}$. We are left with the computation of 
\begin{equation*}\int_{ \left[p^x,\infty(v)\right] \times z^k_j }
\omega_b= \int_{\left( \partial \left[p^x,\infty(v)\right] \right) \times \Sigma(z^k_j) }
\omega_b= - \int_{ \{p^x\} \times \Sigma(z^k_j) } \omega_b.\end{equation*}
Again, we know that this integral is zero along the intersection 
of $\{p^x\} \times \Sigma(z^k_j) $ with $A^{(x)} \times (C_1(\rats) \setminus A^{(x)}_4)$ because 
 $\Sigma(z^k_j)$ does not meet the support of $\omega(p^x)$. We conclude
because $\int_{ \{p^x\} \times \Sigma(a^x_{\ell}) } \omega_b=0$ and because $\Sigma(z^k_j)$ intersects $A^{(x)}_4$ along copies of $ \Sigma(a^x_{\ell})$.

Since $\eta$ is exact on $A^{(x)}_1 \times \bigl(\cup_{k=1}^{x-1}A^{(k)}\bigr)$, we can assume that it vanishes identically there.

Thus, the form $\omega_a$ takes the prescribed values on $A^{(x)} \times (C_1(\rats) \setminus A^{(x)}_4)$, it coincides with $\omega_0$ where $\omega_0$ was prescribed, and it integrates
correctly along the $\Sigma(a^k_{\ell}) \times p^k$ and their symmetric with respect to $\iota$, for $k\neq x$.
Let us now modify $\eta$ so that 
the integrals of $\omega_a$ along the $\{p^x\} \times \Sigma(a^x_{\ell})$ vanish for $\ell=1, \dots, g_x$, too.
We do this by adding to $\eta$  a linear combination $(-\eta_c)$ of $p_2^{\ast}\left(\eta(z^x_j)\right)$ that vanishes on the $A^{(x)}_1 \times A^{(k)}$, for $k<x$, so that we do not change the properties obtained above.
Let $f\colon H_1(\rats \setminus \Int(A^{(x)})) \longrightarrow \RR$ be the linear map 
defined by
\begin{equation*}f(a^x_{\ell})=\int_{\{p^x\} \times \Sigma(a^x_{\ell})}\omega_a.\end{equation*}
Set $\eta_c=\sum_{\ell=1}^{g_x} f(a^x_{\ell})p_2^{\ast}\left(\eta(z^x_{\ell})\right)$. We have $f(y)=\int_{p^x \times y}\eta_c$ for any $y \in \CL_{A^{(x)}}$.

Fix $k<x$. For any $j\in \underline{g_k}$, we have
\begin{multline*}f(z^k_j)=\sum_{\ell=1}^{g_x}lk\left(z^k_j, z^x_{\ell}\right)f(a^x_{\ell})=\int_{\{p^x\} \times \Sigma(z^k_j)} \omega_a\\=\int_{\{\infty(v)\} \times \Sigma(z^k_j)} \omega_a - \int_{\left[p^x,\infty(v)\right] \times z^k_j} \omega_a=0.\end{multline*}
The restriction of  $\eta_c$ to $A^{(x)}_1 \times A_k$ may be expressed as \begin{equation*}\sum_{\ell=1}^{g_x} f(a^x_{\ell}) \sum_{j=1}^{g_k} lk\left(z^k_j, z^x_{\ell}\right) p_2^{\ast}\left(\eta(a^k_j)\right)=\sum_{j=1}^{g_k}f(z^k_j)p_2^{\ast}\left(\eta(a^k_j)\right)=0.\end{equation*}
So $\eta_c$ vanishes on the $A^{(x)}_1 \times A^{(k)}$ for any $k<x$.
Changing $\eta$ into $(\eta-\eta_c)$ does not change $\omega_a$ on the prescribed set but removes $\int_{\{p^x\} \times  \Sigma(a^x_{\ell})}d\chi \eta_c=\int_{\{p^x\} \times (4 \times  a^x_{\ell})}\eta_c=f(a^x_{\ell})$
from $\int_{\{p^x\} \times  \Sigma(a^x_{\ell})} \omega_a$, which becomes $0$. 

After this step, $\omega_a$ is a closed form taking the prescribed values on $PS_a$ such that the integrals of $\omega_a$ along the $(\{p^x\} \times \Sigma(a^x_{\ell}))$ vanish for $\ell=1, \dots, g_i$.
To make $\omega_a$ take the prescribed values on $\iota(PS_a)$, 
we apply 
similar modifications to $\omega_a$ on the symmetric part $(C_1(\rats) \setminus \mathring{A}^{(x)}_2) \times A^{(x)}_1$. Again, the support of these modifications is disjoint
from the support of the previous ones. Thus, they do not interfere. They transform 
$\omega_a$ to a closed form $\omega_c$ with the following additional properties.
\begin{itemize}
\item The form $\omega_c$ has the prescribed form on $(C_1(\rats) \setminus A^{(x)}_3) \times A^{(x)}$.
\item We have $\int_{\Sigma(a^x_j) \times p^x } \omega_c =0$ for all $j=1, \dots g_x$.
\end{itemize}

Now, the form $\omega=\omega_c$ (resp. $\omega=\frac{\omega_c - \iota^{\ast}(\omega_c)}{2}$ if antisymmetry is desired)
has all the required properties, and Proposition~\ref{propnormababy} is proved.

\eopwobp

\section{Proof of \texorpdfstring{Lemma~\ref{lemker}}{the crucial lemma}}
\label{secprooflemker}
In this section, we conclude the proofs of Propositions~\ref{propnormasim} and \ref{propnorma} by proving Lemma~\ref{lemker}.
We first state some homological lemmas.

\begin{lemma}
\label{lemhomdiagS}
Let $S$ be a closed (oriented) surface.
Let $S$ and $S^+$ be two copies of $S$. Let $(c_i)_{i\in \underline{2g}}$ and $(c^{\ast}_i)_{i\in \underline{2g}}$ be two dual bases
of $H_1(S;\ZZ)$ such that 
$\langle c_i, c^{\ast}_j\rangle_{\!S\,}=\delta_{ij}$. Let $\ast$ be a point of $S$.
Set $\diagon(S \times S^+)=\{(x,x^+) \suchthat  x \in S\}$.
We have the following equality in $H_2(S \times S^+;\ZZ)$:
\begin{equation*}\left[\diagon(S \times S^+)\right]=\left[\ast \times S^+\right] +\left [S \times \ast^+\right]+ \sum_{i=1}^{2g} \left[c_i \times c^{\ast+}_i\right].\end{equation*}
\end{lemma}
\bp We have
\begin{equation*}H_2(S \times S^+;\ZZ)=\ZZ\left[\ast \times S^+\right] \oplus \ZZ\left[S \times \ast^+\right] \oplus \bigoplus_{(i,j) \in \underline{2g}^2}\ZZ\left[c_i \times c^{\ast +}_j\right].\end{equation*}
The dual basis of the above basis with respect to the intersection form is \begin{equation*}\Bigl(\left[S \times \ast^+\right],\left[\ast \times S^+\right], \bigl([c^{\ast}_i \times c^{+}_j]\bigr)_{(i,j) \in \underline{2g}^2}\Bigr).\end{equation*}
To get the coordinates of $[\diagon(S \times S^+)]$ with respect to the first basis,
we compute its intersection numbers with the second one. We have
\begin{equation*}\left\langle\left[\diagon(S \times S^+)\right],\left[c^{\ast}_i \times c^{+}_i\right]\right\rangle=\pm 1.\end{equation*}
The tangent space to $\diagon(S \times S^+)$ is naturally parametrized by
$(u_i, v^{\ast}_i, u_i, v^{\ast}_i)$, and the tangent space to $[c^{\ast}_i \times c^{+}_i]$ is naturally parametrized by $(0,w^{\ast}_i,x_i,0)$. So the intersection sign is the sign of the permutation 
\begin{equation*}(u,v,w,x) \mapsto (u,w,x,v),\end{equation*} which is $+1$.
\eop

\begin{lemma} 
\label{lemhomdiagSbry} Let $\Sigma$ be a connected compact oriented surface with one boundary component $J(S^1)$ equipped with a basepoint $\ast=J(1)$.
Let $(c_i)_{i\in \underline{2g}}$ and $(c^{\ast}_i)_{i\in \underline{2g}}$ be two dual bases
of $H_1(\Sigma;\ZZ)$ such that 
$\langle c_i, c^{\ast}_j\rangle=\delta_{ij}$.
Let $\Sigma$ and $\Sigma^+$ be two copies of $\Sigma$. Set $J^+=J^+(S^1)=\partial \Sigma^+$.
Define the subspaces $J\times_{\ast,\leq}J^+$ and $J\times_{\ast,\geq}J^+$ of $J \times J^+$
to be \begin{equation*}J\times_{\ast,\leq}J^+=\Bigl\{\bigl(J(\exp(2i\pi t)),J(\exp(2i\pi u))\bigr) \suchthat  (t,u)\in \left[0,1\right]^2,t\leq u\Bigr\}\end{equation*} and
\begin{equation*}J\times_{\ast,\geq}J^+=\Bigl\{\bigl(J(\exp(2i\pi t)),J(\exp(2i\pi u))\bigr) \suchthat  (t,u)\in \left[0,1\right]^2,t\geq u\Bigr\}.\end{equation*}
Let $\diagon(\Sigma \times \Sigma^+)$ be the subspace $\{(x,x) \suchthat  x \in \Sigma\}$ of $\Sigma \times \Sigma^+$.
Then the chains \begin{equation*}C_{\ast,\leq}(\Sigma,\Sigma^+)=\diagon(\Sigma \times \Sigma^+)-\ast \times \Sigma^+ -\Sigma \times \ast^+ -J\times_{\ast,\leq}J^+\end{equation*} and
\begin{equation*}C_{\ast,\geq}(\Sigma,\Sigma^+)=\diagon(\Sigma \times \Sigma^+)-\ast \times \Sigma^+ -\Sigma \times \ast^+ +J\times_{\ast,\geq}J^+\end{equation*}
are cycles, and
we have
\begin{equation*}\bigl[C_{\ast,\leq}(\Sigma,\Sigma^+)\bigr]=\bigl[C_{\ast,\geq}(\Sigma,\Sigma^+)\bigr] =\sum_{i=1}^{2g} \left[c_i \times c^{\ast+}_i\right]\end{equation*} in $H_2(\Sigma \times \Sigma^+;\ZZ)$.
\end{lemma}
\bp Since $\partial (J\times_{\ast,\leq}J^+) =\diagon(J \times J^+) -\ast \times J^+ -J \times \ast^+$, the chain
 $C_{\ast,\leq}(\Sigma,\Sigma^+)$ is a cycle. 
Consider the closed surface $S$ obtained from $\Sigma$ by gluing a disk $D$ along $J$. According to Lemma~\ref{lemhomdiagS}, we have
\begin{equation*}\left[\diagon(S \times S^+)\right]=\left[\ast \times S^+\right] + \left[S \times \ast^+\right]+ \sum_{i=1}^{2g} \left[c_i \times c^{\ast +}_i\right]\end{equation*} in $H_2(S\times S^+;\ZZ)$.
This implies \begin{equation*}[C_{\ast,\leq}(\Sigma,\Sigma^+)-C_{\ast,\leq}(-D,(-D)^+)]=\sum_{i=1}^{2g} \left[c_i \times c^{\ast +}_i\right]\end{equation*} in $H_2(S \times S^+;\ZZ)$.
Since the cycle $C_{\ast,\leq}(-D,(-D)^+)$ lies in $D\times D^+$, it is null-homologous there.
Since $H_2(\Sigma \times \Sigma^+;\ZZ)$ injects naturally into 
$H_2(S \times S^+;\ZZ)$, we deduce that $\left[C_{\ast,\leq}(\Sigma,\Sigma^+)\right] =\sum_{i=1}^{2g} \left[c_i \times c^{\ast +}_i\right]$ in $H_2(\Sigma \times \Sigma^+;\ZZ)$. The proof for $C_{\ast,\geq}(\Sigma,\Sigma^+)$ is similar.
\eop

Consider a rational homology handlebody $A$ with a collar
$[-4,0] \times \partial A$ of its boundary.
For $s \in [-4,0]$, recall
$A_s= A \setminus  ( \left]s,0\right] \times \partial A )$ and $ \partial A_s =\{s\} \times \partial A$.
Fix pairwise disjoint simple closed curves $(a_i)_{i=1,\dots,g_A}$ and pairwise disjoint simple closed curves  $(z_i)_{i=1,\dots,g_A}$ on $\partial A$ such that
$\CL_{A}=\oplus_{i=1}^{g_A} [a_i]$
and 
$\langle a_i,z_j\rangle_{\!\partial A\,}=\delta_{ij}$.

Consider a curve $a$ of $\partial A$ disjoint from the curves $a_i$. The class of $a$ is in $\CL_A$. Let $k \in \NN \setminus \{0\}$ be its order $H_1(A;\ZZ)$. 
Let $\Sigma=k\Sigma(a)$ be a surface of $A$ immersed in $A$ bounded by $ka$. Assume that $\Sigma$
intersects $ \left[-1,0\right] \times \partial A $ as $k$ copies of $\left[-1,0\right] \times a$, the interior $\Int(A_{-1})$ of $A_{-1}$ as an embedded surface, and $\Int(A_{-1})\setminus A_{-2}$ as $k$ disjoint annuli. See the thick part of Figure~\ref{figlemfondunivone}. For $s \in [-2,0]$, set $\Sigma_{s}=\Sigma \cap A_s$.

\begin{lemma}
\label{homctimec}
With the above notation,
let $(c_i)_{i=1, \dots 2g}$ and $(c^{\ast}_i)_{i=1, \dots, 2g}$ be two dual bases of ${H_1(\Sigma_{-2};\ZZ)}/{H_1( \partial \Sigma_{-2};\ZZ)}$ such that $\langle c_i, c^{\ast}_j\rangle=\delta_{ij}$. Represent $(c_i)_{i=1, \dots 2g}$ and $(c^{\ast}_i)_{i=1, \dots, 2g}$ by curves $(c_i)_{i=1, \dots 2g}$ and $(c^{\ast}_i)_{i=1, \dots, 2g}$ of $\Sigma_{-2}$. Let $\Sigma_{-2} \times [-1,1]$ denote a tubular neighborhood of $\Sigma_{-2}=\Sigma_{-2} \times \{0\}$ in $A_{-2}$. For a curve $\sigma$ of $\Sigma_{-2}$, the curve $\sigma \times \{1\}$ is denoted by $\sigma^+$.

Then $\sum_{i=1}^{2g} c_i \times c^{\ast}_i$ is homologous to $\sum_{(j,{\ell}) \in \underline{g_A}^2\setminus \diagon}\langle \Sigma,\Sigma(a_j),\Sigma(a_{\ell}) \rangle_{\!A\,} z_j \times z_{\ell}$ in $A^2$. 
Furthermore, the sum $\sum_{i=1}^{2g} c_i \times c^{\ast +}_i$ is homologous to \begin{equation*}\sum_{(j,{\ell}) \in \underline{g_A}^2 \setminus \diagon}\langle \Sigma,\Sigma(a_j),\Sigma(a_{\ell}) \rangle_{\!A\,} z_j \times z_{\ell} -g\ST A\vert_{\ast}\end{equation*} in $C_2(A)$.
\end{lemma}
\bp Assume that $\Sigma$ and the $\Sigma(a_j)$ are transverse to each other. For $(j,{\ell}) \in \{1, \dots, g_A\}^2$, set $\gamma_{\Sigma j}=\Sigma\cap\Sigma(a_j)$ and $\gamma_{\Sigma \ell}=\Sigma\cap\Sigma(a_{\ell})$. If $j \neq \ell$, also set $\gamma_{j\ell}=\Sigma(a_j)\cap\Sigma(a_{\ell})$.
We have
\begin{equation*}c_i=\sum_{j=1}^{g_A}\langle c_i, \Sigma(a_j) \rangle_{\!A\,} z_j=\sum_{j=1}^{g_A}\langle c_i , \gamma_{\Sigma j} \rangle_{\!\Sigma\,} z_j \mbox{ and }  
c^{\ast}_i=\sum_{\ell=1}^{g_A}\langle c^{\ast}_i, \gamma_{\Sigma \ell} \rangle_{\!\Sigma\,} z_{\ell}\end{equation*} in $H_1(A)$.
This implies
\begin{equation*}c_i \times c^{\ast}_i
=\sum_{(j,\ell) \in \underline{g_A}^2}\langle c_i, \gamma_{\Sigma j} \rangle_{\!\Sigma\,}\langle c^{\ast}_i, \gamma_{\Sigma \ell} \rangle_{\!\Sigma\,} z_j \times z_{\ell} \end{equation*} in $H_2(A^2)$.
On the other hand, we have
\begin{equation*}\gamma_{\Sigma j}=\sum_{i=1}^{2g}\langle c_i, \gamma_{\Sigma j} \rangle_{\!\Sigma\,}c^{\ast}_i \mbox{ and } 
\gamma_{\Sigma \ell}=-\sum_{i=1}^{2g}\langle c^{\ast}_i, \gamma_{\Sigma \ell}\rangle_{\!\Sigma\,}c_i\end{equation*} in $H_1(\Sigma_{-2})/H_1(\partial \Sigma_{-2})$.
This implies \begin{equation*}\langle \gamma_{\Sigma j},\gamma_{\Sigma \ell} \rangle_{\!\Sigma\,}=\sum_{i=1}^{2g}\langle c_i, \gamma_{\Sigma j} \rangle_{\!\Sigma\,}\langle c^{\ast}_i,  \gamma_{\Sigma \ell} \rangle_{\!\Sigma\,}.\end{equation*}
In particular,
we have
\begin{equation*}\sum_{i=1}^{2g}\langle c_i, \gamma_{\Sigma j} \rangle_{\!\Sigma\,}\langle c^{\ast}_i,  \gamma_{\Sigma j} \rangle_{\!\Sigma\,}=0\end{equation*} for any $j \in \{1, \dots, g_A\}$.
If $j \neq \ell$, then
$\langle \Sigma,\Sigma(a_j),\Sigma(a_{\ell}) \rangle_{\!A\,}=\langle \gamma_{\Sigma j},\gamma_{\Sigma \ell} \rangle_{\!\Sigma\,}$. So $\alpha=\sum_{i=1}^{2g} c_i \times c^{\ast+}_i$ is homologous to $\sum_{(j,{\ell}) \in \underline{g_A}^2 \setminus \diagon}\langle \Sigma,\Sigma(a_j),\Sigma(a_{\ell}) \rangle_{\!A\,} z_j \times z_{\ell}$  
as announced.

Let us now compute the homology class $[\alpha]$ of $\alpha$ in $H_2(C_2(A))$.
For a curve $\sigma$ of $\Sigma_{-2}$, both $a \times \sigma^+$ and $\sigma \times a^+$ are null-homologous in $C_2(A)$. Hence, the class $[\alpha]$ depends only on the class of $\sum_{i=1}^{2g} c_i \otimes c^{\ast+}_i$ in $H_1(\Sigma_{-2})/H_1(\partial \Sigma_{-2}) \otimes H_1(\Sigma_{-2}^+)/H_1(\partial \Sigma_{-2}^+)$. This class is determined by the property that for any two closed curves $e$ and $f$ of $\Sigma$, we have $\langle e \times f^+,\alpha \rangle_{\!\Sigma \times \Sigma^+\,}=-\langle e,f \rangle_{\!\Sigma\,}$. So the class $[\alpha] \in H_2(C_2(A))$ is independent of the dual bases $(c_i)$ and $(c^{\ast}_i)$.
In particular, we have $[\alpha]=\bigl[\sum_{i=1}^{2g}  c^{\ast}_i \times (-c^{+}_i)\bigr]$. Set \begin{equation*}\beta=\sum_{(j,{\ell}) \in \underline{g_A}^2 \setminus \diagon}\bigl\langle \Sigma,\Sigma(a_j),\Sigma(a_{\ell}) \bigr\rangle_{\!\!A\,} z_j \times z_{\ell}-g \ST A\vert_{\ast}.\end{equation*}
The previous computation tells us that the difference $[\alpha-\beta]$
is a rational multiple of $[\ST A\vert_{\ast}]$ in $H_2(C_2(A))$. To evaluate this multiple, we embed $A$ in a rational homology ball
obtained from $A$ by adding thickened disks along neighborhoods of the $z_i$. Embed this rational homology ball in a rational homology sphere $\rats$.
We get 
\begin{multline*}[\alpha-\beta]=\frac{1}{2}\sum_{i=1}^{2g}\bigl(lk_{\rats}(c_i,c^{\ast+}_i)-lk_{\rats}(c^{\ast-}_i, c_i)\bigr)[\ST A\vert_{\ast}]
\\-\sum_{(j,{\ell}) \in \{1, \dots, g_A\}^2\setminus \diagon}\bigl\langle \Sigma,\Sigma(a_j),\Sigma(a_{\ell})\bigr\rangle_{\!\!A\,} lk_{\rats}(z_j,z_{\ell})[\ST A\vert_{\ast}]\end{multline*} in $H_2(C_2(\rats);\RR)$.
Since $lk_{\rats}(z_j,z_{\ell})=0$,
the second row vanishes and we obtain $[\alpha-\beta]=-g[\ST A\vert_{\ast}]$ in $H_2(C_2(\rats);\RR)$.
Since $[\ST A\vert_{\ast}]\neq 0$ in $H_2(C_2(\rats);\RR)$, this equality also holds in $H_2(C_2(A);\RR)$.
\eop

We now define a cycle $F^2(\Sigma(a))$ of $\partial C_2(A)$ associated to the surface $\Sigma=k\Sigma(a)$ introduced before Lemma~\ref{homctimec}.
Let $\left(a \times [-1,1]\right)$ be a tubular neighborhood of $a$ in $\partial A$.
Let $p(a) \in a$. View $a$ as the image of a map $a\colon \left[0,1\right] \rightarrow a$ such that $a(0)=a(1)=p(a)$.
Set
\begin{equation*}\Sigma^+=\Sigma_{-1} \cup k\Bigl\{\bigl(t-1, a(\alpha),t\bigr) \suchthat (t,\alpha)\in\left[0,1\right]^2\Bigr\}.\end{equation*}
So $\partial \Sigma^+= k a^+$, where $a^+= a \times \{1\}$.
Set
\begin{equation*}p(a)^+=(p(a),1)=(0,p(a),1) \in a \times [-1,1] \subset \bigl(\partial A =\{0\} \times \partial A\bigr).\end{equation*}
Recall $a\times_{p(a),\geq}a^+=\{((a(v),0),(a(w),+1)) \suchthat (v,w)\in \left[0,1\right]^2,v\geq w \}$.
Let $T(a)$ be the closure of $\{((a(v),0),(a(v),t)) \suchthat (t,v)\in \left]0,1\right] \times \left[0,1\right]\}$ (oriented by $(t,v)$) in $\partial C_2(A)$.
Let $s_+(\Sigma)$ be the positive normal section of $\ST(A)\vert_{\Sigma}$. Let $g$ be the genus of $\Sigma$. Set $e(\Sigma(a)=\frac1{k}{\Sigma})= \frac1{k}(g+k-1)$. We have
\begin{equation*}e(\Sigma(a))=\frac{-\chi(\Sigma)}{2k} + \frac{1}{2}.\end{equation*}

\begin{lemma}
\label{lemf(a)}
With the above notation, the chain
\begin{multline*}F^2(\Sigma(a))= \frac{1}{k}s_+(\Sigma) + T(a) - \Bigl(p(a) \times \frac{1}{k}\Sigma^+\Bigr) -\Bigl(\frac{1}{k}\Sigma \times p(a)^+\Bigr)  + \Bigl(a\times_{p(a),\geq}a^+\Bigr)\\ 
\\+ e\bigl(\Sigma(a)\bigr)[\ST A\vert_{\ast}] -\sum_{(j,{\ell}) \in \{1, \dots, g_A\}^2\setminus \diagon}\bigl\langle \Sigma(a),\Sigma(a_j),\Sigma(a_{\ell}) \bigr\rangle z_j \times z_{\ell}\end{multline*}
is a cycle, and it is null-homologous in $C_2(A)$.
\end{lemma}
\bp
For $k=1$ (when we are dealing with integral homology handlebodies, for example), it is a direct consequence of Lemma~\ref{lemhomdiagSbry} and Lemma~\ref{homctimec} above.
Let us now focus on the case $k>1$. Observe that $F^2(\Sigma(a))$ is a cycle.
Without loss of generality, assume
\begin{equation*}\Sigma \cap \bigl([-2,-1] \times \partial A\bigr)=\Bigl\{\bigl(t-2,a(\alpha),(1-t)\frac{j-1}{k}\bigr) \suchthat (t,\alpha)\in \left[0,1\right]^2,j\in \underline{k}\Bigr\}.\end{equation*}
For the proof, we change $\Sigma^+$ to a surface (still denoted by $\Sigma^+$)  with the same boundary as follows.

Recall $\Sigma_{-2}=\Sigma \cap A_{-2}$. We have $\partial \Sigma_{-2} =\cup_{j=1}^k \left(\{-2\}\times a \times\{\frac{j-1}{k}\}\right)$. Let $\Sigma_{-2}^+$ be a parallel copy of $\Sigma$ on its positive side with boundary
$\partial \Sigma_{-2}^+ =\cup_{j=1}^k \left(\{-2\}\times a \times\{\frac{j-1/2}{k}\}\right)$.
Redefine $\Sigma^+$ so that we have $\Sigma^+\cap A_{-2}=\Sigma_{-2}^+$,
\begin{equation*}\Sigma^+ \cap \bigl([-1,0] \times \partial A\bigr)=k\bigl([-1,0] \times a \times \{1\}\bigr),\end{equation*}
and
\begin{equation*}\Sigma^+ \cap \bigl([-2,-1]\times\partial A\bigr)=\Bigl\{\bigl(t-2,a(\alpha),(1-t)\frac{j-1/2}{k}+t\bigr) \suchthat (t,\alpha)\in \left[0,1\right]^2,j\in\underline{k}\Bigr\}\end{equation*}
as in Figure~\ref{figlemfondunivone}, which represents $\Sigma \cap \left([-2,0]\times p(a) \times[-1,1]\right) $ as the thick lines and $\Sigma^+ \cap \left([-2,0] \times p(a) \times[-1,1]\right)$ as the thin lines when $k=3$.
Observe that this modification of $\Sigma^+$ changes neither the boundary of $\Sigma^+$ nor the
class of the resulting cycle $F^2(\Sigma(a))$ in $H_2(C_2(A))$.

\bfig
\centering
\begin{tikzpicture} \useasboundingbox (0,0) rectangle (2.4,2.3);
\draw [line width=2pt] (0,0) -- (0,2) (0.8,0) -- (0,1) (1.6,0) -- (0,1) (0,2) node[left]{\scriptsize $p(a)$} (0,1) node[left]{\scriptsize $(-1,p(a),0)$} (0,0) node[left]{\scriptsize $(-2,p(a),0)$};
\fill (0,0) circle (2pt) (0.8,0) circle (2pt) (0,1) circle (2pt) (1.6,0) circle (2pt) (0,2) circle (2pt);
\draw (0.4,0) -- (2.4,1) -- (2.4,2) (1.2,0) -- (2.4,1) (2,0) -- (2.4,1) (2.4,2) node[right]{\scriptsize $p(a)^+$} (2.4,1) node[right]{\scriptsize $(-1,p(a),1)$};
\fill (0.4,0) circle (1.5pt) (1.2,0) circle (1.5pt) (2.4,1) circle (1.5pt) (2,0) circle (1.5pt) (2.4,2) circle (1.5pt);
\end{tikzpicture}
\caption{How $\Sigma$ and $\Sigma^+$ intersect $\left([-2,0]\times p(a) \times[-1,1]\right)$ in the proof of Lemma~\ref{lemf(a)}.}\label{figlemfondunivone}
\end{figure}

Let $S$ be the closed surface obtained from $\Sigma_{-2}$ by
gluing abstract disks $D_j$ with respective boundaries $\{-2\}\times (-a) \times\{\frac{j-1}{k}\}$ on 
$\partial \Sigma_{-2}$.
Let $S^+$ be similarly obtained from $\Sigma_{-2}^+$ by
gluing abstract disks $D_j^+$ with boundaries $\{-2\}\times (-a) \times\{\frac{j-1/2}{k}\}$ on 
$\partial \Sigma_{-2}^+$.
For $j=1,\dots, k$, set $p_j=(-2,p(a),\frac{j-1}{k}) \in \partial A_{-2}$ and $p^+_j=(-2,p(a),\frac{j-1/2}{k})$.
Recall the dual bases $(c_i)$ and $(c^{\ast}_i)$ of Lemma~\ref{homctimec}.
Lemma~\ref{lemhomdiagS} implies that
\begin{equation*}C(S)=\diagon(S\times S^+)- p_1 \times S^+ -S \times p^+_k - \sum_{i=1}^{2g} c_i \times c^{\ast +}_i\end{equation*}
is null-homologous in $H_2(S\times S^+)$. Choose closed representatives $c_i$ of the classes $c_i$ 
in the interior of $\Sigma_{-2}$ so that  $\left(\Sigma_{-2} \setminus \cup_{i=1}^{2g}c_i\right)$ is connected.
Let $[p_1,p_j]$ denote a path from $p_1$ to $p_j$ in $\left(\Sigma_{-2} \setminus \cup_{i=1}^{2g}c_i\right)$. Let $[p_j^+,p_k^+]$ be a path  from $p_j^+$ to $p_k^+$ in $\left(\Sigma^+_{-2} \setminus \cup_{i=1}^{2g}c^+_i\right)$.
Add to $C(S)$ the null-homologous cycles
\begin{equation*}\partial \bigl(-[p_1,p_j] \times D_j^+\bigr)=p_1\times D_j^+ -p_j\times D_j^+
+[p_1,p_j]\times \partial D_j^+,\end{equation*} 
\begin{equation*}\partial \bigl(D_j \times [p_j^+,p_k^+]\bigr)=D_j \times p^+_k -D_j \times p^+_j +\partial D_j \times [p_j^+,p_k^+],\end{equation*}
and the null-homologous cycles $(-C_{\ast,\leq}(D_j,D_j^+))$ of Lemma~\ref{lemhomdiagSbry}, for $j=1,\dots,k$.
This addition transforms $C(S)$ to the null-homologous cycle
\begin{equation*}\begin{array}{ll}C(\Sigma_{-2})=&\diagon(\Sigma_{-2}\times \Sigma_{-2}^+)- p_1 \times \Sigma_{-2}^+ -\Sigma_{-2} \times p^+_k - \sum_{i=1}^{2g} c_i \times c^{\ast +}_i\\& +\sum_{j=1}^k\bigl(\partial D_j \times [p_j^+,p_k^+]+[p_1,p_j]\times \partial D_j^+ +\partial D_j \times_{p(a),\leq} \partial D_j^+\bigr).\end{array}\end{equation*}
Let us deform the cycle $C(\Sigma_{-2})$ continuously in $\Sigma \times \Sigma^+$ to move the level $\{-2\} \times \partial A$ to the level $\{0\} \times \partial A$.
During such a deformation, when $s$ tends to $(-1)$,
the path $\left[\{s\} \times p_j^+,\{s\} \times p_k^+\right]$ becomes a loop $[p_j^+,p_k^+]_{-1}$ on $\Sigma^+_{-1}=\Sigma^+ \cap A_{-1}$, and the path  $\left[\{s\} \times p_1,\{s\} \times p_j\right]$ becomes a loop $[p_1,p_j]_{-1}$ on $\Sigma_{-1}$.
Thus the cycle $C(\Sigma_{-2})$ can be naturally deformed in $\Sigma \times \Sigma^+$ to the following still null-homologous cycle:
\begin{multline*}C(\Sigma)=\diagon(\Sigma \times \Sigma^+)- p(a) \times \Sigma^+ -\Sigma \times p(a)^+ - \sum_{i=1}^{2g} c_i \times c^{\ast +}_i\\ +\sum_{j=1}^k\bigl((-a) \times [p_j^+,p_k^+]_{-1}+[p_1,p_j]_{-1}\times (-a^+) \bigr) + k \left(a \times_{p(a),\geq} a^+\right).\end{multline*}
Since $a$ bounds $\frac{1}{k}\Sigma$, the cycle $(-a) \times [p_j^+,p_k^+]_{-1}$ is homologous to the cycle $\langle \frac{-1}{k}\Sigma,[p_j^+,p_k^+]_{-1} \rangle_{\!A\,} \ST A\vert_{\ast}$, in $C_2(A)$.
Similarly, $[p_1,p_j]_{-1}\times (-a^+)$ is homologous to $\langle \frac{-1}{k}\Sigma^+,[p_1,p_j]_{-1} \rangle_{\!A\,} \ST A\vert_{\ast}$, in $C_2(A)$.
Intersections occur where $\Sigma$ and $\Sigma^+$ intersect, in $([-2,-1] \times \partial A)$, as shown in Figure~\ref{figlemfondunivtwo}, where the positive normal to $\Sigma$ goes from left to right.

\bfig
\centering
\begin{tikzpicture} 
\draw [line width=2pt] (0,0) -- (0,1) (0.8,0) -- (0,1) (1.6,0) -- (0,1) (0,1) node[left]{\scriptsize $(-1,p(a),0)$} (0,0) node[left]{\scriptsize $(-2,p(a),0)$};
\fill (0,0) circle (2pt) (0.8,0) circle (2pt) (0,1) circle (2pt) (1.6,0) circle (2pt);
\draw (2.4,1) node[right]{\scriptsize $(-1,p(a),1)$};
\draw [->] (2.4,1) -- (0.4,0);
\draw [->] (2,0) -- (2.4,1);
\draw [dashed] (1.2,0) -- (2.4,1);
\draw [dotted] (0.4,0) .. controls (-.2,-.3) and (1.8,-.5) .. (2,0);
\end{tikzpicture}
\caption{The intersection $\Sigma \cap \Sigma^+$ and the loop $[p_j^+,p_k^+]_{-1}$ ($j=1$, $k=3$)}\label{figlemfondunivtwo}
\end{figure}
For any $j \in \underline{k}$, we have
\begin{equation*}\Bigl\langle \frac{-1}{k}\Sigma,[p_j^+,p_k^+]_{-1} \Bigr\rangle_{\!\!A\,}=\frac{k-j}{k}\mbox{
and }\Bigl\langle \frac{-1}{k}\Sigma^+,[p_1,p_j]_{-1} \Bigr\rangle_{\!\!A\,}=\frac{j-1}{k}.\end{equation*}
Therefore, Lemma~\ref{homctimec} implies that the null-homologous cycle $C(\Sigma)$ is homologous to 
\begin{multline*}\diagon(\Sigma \times \Sigma^+)- p(a) \times \Sigma^+ -\Sigma \times p(a)^+ - \sum_{(j,{\ell}) \in \{1, \dots, g_A\}^2\setminus \diagon}\langle \Sigma,\Sigma(a_j),\Sigma(a_{\ell}) \rangle z_j \times z_{\ell} \\+g\ST A\vert_{\ast} + (k-1)\ST A\vert_{\ast} + k \left(a \times_{p(a),\geq} a^+\right).\end{multline*}
This cycle is naturally homologous to $kF^2(\Sigma(a))$, which is therefore homologous to zero.
\eop

\begin{lemma}
\label{lemhomctwoa}
If $A$ is a rational homology handlebody such that $H_1(A)=\oplus_{j=1}^{g(A)}\RR[z_j]$, then we have
\begin{equation*}H_3(C_2(A))=\oplus_{j=1}^{g(A)} \RR\left[\ST A\vert_{z_j}\right].\end{equation*}
\end{lemma}
\bp
The configuration spaces $C_2(A)$ and $C_2(\mathring{A})$ have the same homotopy type, which is the homotopy type of $\mathring{A}^{2} \setminus \diagon$. We have
$H_3(\mathring{A}^2)=H_4(\mathring{A}^2)=0$ and
$H_4(\mathring{A}^{2},\mathring{A}^{2} \setminus \diagon)\cong H_4(\mathring{A}^2 \times B^3,\mathring{A}^2 \times S^2 )=\oplus_{j=1}^{g(A)} \RR[z_j \times B^3]$.
\eop

\begin{lemma}
\label{lemgenker}
Let $i \in \underline{x}$. For any $j\in \underline{g_i}$, assume that the chains $\Sigma(a^i_j)$ and $\Sigma^{\prime}(a^i_j)$ defined before Proposition~\ref{propnorma} intersect $[-1,4]\times \partial A^{(i)}$
as $[-1,4] \times a^i_j$. Also assume that they may be respectively expressed as $\frac1{k^{\prime}}\Sigma^{\prime}$ and $\frac1{k}\Sigma$, for immersed surfaces $\Sigma$ and $\Sigma^{\prime}$, which respectively intersect $\Int(A^{(i)}_{-1})$ and $\Int(A^{(i)\prime}_{-1})$ as embedded surfaces (as before Lemma~\ref{homctimec}).
Fix $p(a^i_j)$ on $\{4\} \times a^i_j$.
Then the classes of the cycles $F^2(\Sigma^{\prime}(a^i_j))$ of $\partial C_2(A^{(i)\prime}_4)$ defined in Lemma~\ref{lemf(a)}, for $j\in \underline{g_i}$ generate the kernel
 of the map induced by the inclusion
\begin{equation*}H_2\Bigl(D_A\bigl(\omega_0(\rats_i)\bigr)\Bigr) \longrightarrow H_2\Bigl(C_2\bigl(A^{(i)\prime}_4\bigr)\Bigr)\end{equation*} 
for the domain $D_A(\omega_0(\rats_i))$ defined before Lemma~\ref{lemker}.
\end{lemma}
\bp First note that $D_A(\omega_0(\rats_i))$ is homotopically equivalent by retraction to $\partial C_2(A^{(i)\prime}_4)$.
The cycles $F^2(\Sigma^{\prime}(a^i_j))$ lie in $\partial C_2(A^{(i)\prime}_4))$ and bound chains $G^3(a^i_j)$ in $C_2(A^{(i)\prime}_4)$ according to Lemma~\ref{lemf(a)}. These chains can be assumed to be transverse to the boundary. So they satisfy 
\begin{equation*}\Bigl\langle \bigl[G^3(a^i_j)\bigr], \bigl[\ST A^{(i)\prime}\vert_{z^i_k}\bigr]\Bigr\rangle_{\!\!C_2(A^{(i)\prime}_4)\,}=\pm \Bigl\langle F^2\bigl(\Sigma^{\prime}(a^i_j)\bigr), \ST A^{(i)\prime}\vert_{\{0\} \times z^i_k}\Bigr\rangle_{\!\!\partial C_2(A^{(i)\prime}_4)\,}= \pm \delta_{jk}.\end{equation*}
Therefore, Poincar\'e duality and Lemma~\ref{lemhomctwoa} imply
\begin{equation*}H_3\bigl(C_2(A^{(i)\prime}_4),\partial C_2(A^{(i)\prime}_4)\bigr)=\oplus _{j=1}^{g_i} \RR \left[G^3(a^i_j)\right].\end{equation*}

 The boundary map of the long exact sequence associated to the pair $\bigl(C_2(A^{(i)\prime}_4), \partial C_2(A^{(i)\prime}_4)\bigr)$ sends the cycle $\bigl[G^3(a^i_j)\bigr]$ of $\bigl(C_2(A^{(i)\prime}_4), \partial C_2(A^{(i)\prime}_4)\bigr)$ to $\bigl[F^2(\Sigma^{\prime}(a^i_j))\bigr]$.
\eop

\bpo{Proof of Lemma~\ref{lemker}}
According to Lemma~\ref{lemgenker}, it suffices to prove \begin{equation*}\int_{F^2(\Sigma^{\prime}(a^i_j))}\omega_0(\rats_i)=0\end{equation*}
for any $i \in \underline{x}$, and for any $j \in \underline{g_i}$.
Fix $i \in \underline{x}$ and $j \in \underline{g_i}$. Set $a=\{4\} \times a^i_j$. Let $F^{\prime}$ denote the cycle $F^2(\Sigma^{\prime}(a))$ of $\partial C_2(A^{(i)\prime}_4)$ associated to $\Sigma^{\prime}=k^{\prime}\Sigma^{\prime}(a)$ and to $p(a)=p(a^i_j)$, and let $F$ denote the cycle $F^2(\Sigma(a))$ of $\partial C_2(A^{(i)}_4)$ similarly associated to $\Sigma(a)$ and to $p(a)$.
We have
\begin{multline*}F^{\prime}= T(a) + a\times_{p(a),\geq}a^+ - p(a) \times \frac{1}{k^{\prime}}\Sigma^{\prime+} -\frac{1}{k^{\prime}}\Sigma^{\prime} \times p(a)^+ +\frac{1}{k^{\prime}}s_+(\Sigma^{\prime}) \\+ e(\Sigma^{\prime}(a))[\ST A\vert_{\ast}]
-\sum_{(p,q) \in \{1, \dots g_i\}^2\setminus \diagon}\langle \Sigma^{\prime}(a),\Sigma^{\prime}(a^i_p),\Sigma^{\prime}(a^i_q) \rangle(\{4\} \times  z^i_p) \times (\{4\} \times   z^i_q)\end{multline*}
Set $\Sigma^{\prime+}_{-1}=\Sigma^{\prime+} \cap A^{(i)\prime}_{-1}$.
The integral of $\omega_0(\rats_i)$ along 
\begin{equation*} - p(a) \times  \frac{1}{k^{\prime}}\Sigma^{\prime+}_{-1} -\frac{1}{k^{\prime}}\Sigma^{\prime}_{-1} \times p(a)^+\end{equation*}
is zero 
because of the prescribed form of $\omega_0(\rats_i)$ on 
\begin{equation*}\bigl( (\crats \setminus A^{(i)\prime}_3) \times A^{(i)\prime}\bigr) \cup \bigl(A^{(i)\prime} \times (\crats \setminus A^{(i)\prime}_3) \bigr) .\end{equation*}
For the form $\omega$ of Proposition~\ref{propnormababy}, the same argument implies
\begin{equation*}\int_{- p(a) \times \Sigma^+_{-1}(a^i_j) -\Sigma_{-1}(a^i_j) \times p(a)^+}\omega=0.\end{equation*} The part
\begin{multline*}C=T(a) + a\times_{p(a),\geq}a^+ -\bigl( p(a) \times \frac{1}{k^{\prime}}(\Sigma^{\prime+} \setminus \Sigma^{\prime+}_{-1})\bigr) -\bigl(\frac{1}{k^{\prime}}(\Sigma^{\prime}\setminus \Sigma^{\prime}_{-1}) \times p(a)^+\bigr) \\+\frac{1}{k^{\prime}}s_+(\Sigma^{\prime} \setminus \Sigma^{\prime}_{0})\end{multline*} of 
$F^{\prime}$ (or $F$)
lies in the intersection of $C_2(A^{(i)}_4)$ and $C_2(A^{(i)\prime}_4)$, inside which $\omega=\omega_0(\rats_i)$.
So we have \begin{equation*}\int_{C}\omega=\int_{C}\omega_0(\rats_i).\end{equation*}
Recall \begin{equation*}\int_{(\{4\} \times  z^i_p) \times (\{4\} \times   z^i_q)}\omega =  lk(z^i_p,z^i_q)\end{equation*} for any $(p,q) \in \{1, \dots g_i\}^2\setminus \diagon$.
Since $lk(z^i_p,z^i_q)=lk(z^i_q,z^i_p)$, and since \begin{equation*}\bigl\langle \Sigma(a),\Sigma(a^i_p),\Sigma(a^i_q) \bigr\rangle=-\bigl\langle \Sigma(a),\Sigma(a^i_q) ,\Sigma(a^i_p)\bigr\rangle,\end{equation*}
the integral of $\omega$ along 
\begin{equation*}\sum_{(p,q) \in \{1, \dots g_i\}^2\setminus \diagon}\bigl\langle \Sigma(a),\Sigma(a^i_p),\Sigma(a^i_q) \bigr\rangle\bigl(\{4\} \times  z^i_p\bigr) \times \bigl(\{4\} \times   z^i_q\bigr)\end{equation*} is equal to zero.
The integral of $\omega_0(\rats_i)$ along 
\begin{equation*}\sum_{(p,q) \in \{1, \dots g_i\}^2\setminus \diagon}\bigl\langle \Sigma^{\prime}(a),\Sigma^{\prime}(a^i_p),\Sigma^{\prime}(a^i_q) \bigr\rangle\bigl(\{4\} \times  z^i_p\bigr) \times \bigl(\{4\} \times   z^i_q\bigr)\end{equation*}
 vanishes similarly.

Since $\int_{F}\omega=0$, we have \begin{multline*}\int_{F^{\prime}}\omega_0(\rats_i)=\int_{F^{\prime}}\omega_0(\rats_i)-\int_{F}\omega
\\=\int_{\frac{1}{k^{\prime}}s_+(\Sigma^{\prime}_{0})}\omega_0(\rats_i) - \int_{s_+(\Sigma_{0}(a^i_j))}\omega +e\bigl(\Sigma^{\prime}(a)\bigr)-e\bigl(\Sigma(a)\bigr).\end{multline*}
Recall that $\tau$ coincides with $\tau_i$ on $\left[0,4\right] \times \partial A^{(i)}$. When $\tau$ maps $e_3$ to the positive normal to $\Sigma^{\prime+}_{0}$ along $\partial \Sigma^{\prime+}_{0}$, Proposition~\ref{propintsplusstau} and Lemma~\ref{lemptyeulernumb} imply 
\begin{equation*}\int_{\frac{1}{k^{\prime}}s_+(\Sigma^{\prime}_{0})}\omega_0(\rats_i)= \frac12 d\bigl(\tau(. \times e_2),\{0\} \times a\bigr) + \frac{1}{2k^{\prime}}\chi(\Sigma^{\prime}_{-2}).\end{equation*}
We get \begin{multline*}\int_{\frac{1}{k^{\prime}}s_+(\Sigma^{\prime}_{0})}\omega_0(\rats_i)+e\bigl(\Sigma^{\prime}(a)\bigr)= \frac12 d\bigl(\tau(. \times e_2),\{0\} \times a\bigr) + \frac{1}{2}\\=\int_{s_+(\Sigma_{0}(a^i_j))}\omega +e\bigl(\Sigma(a)\bigr)\end{multline*}
and $\int_{F^{\prime}}\omega_0(\rats_i)=0$.
When $\tau$ does not map $e_3$ to the positive normal to $\Sigma^{\prime+}_{0}$ along $\partial \Sigma^{\prime+}_{0}$, perform
a simultaneous homotopy on $\tau$ and $\tau_i$ to make this happen without changing \begin{equation*}\int_{\frac{1}{k^{\prime}}s_+(\Sigma^{\prime}_{0})}\omega_0(\rats_i) - \int_{s_+(\Sigma_{0}(a^i_j))}\omega.\end{equation*} Thus, the above proof still implies $\int_{F^{\prime}}\omega_0(\rats_i)=0$.
\eop

\chapter{Much more flexible definitions of \texorpdfstring{$\Zinvuf$}{Z}}
\label{chappseudoparmuchmore}

\section{More propagating forms associated to pseudo-parallelizations}
\label{secpseudopropform}

In this section, we define nonhomogeneous propagating forms associated with pseudo-parallelizations,
and we give more flexible definitions of $\Zinv$
involving these forms. 
In Section~\ref{secpseudosec}, we define propagating chains associated with pseudo-parallelizations, which allow discrete computations of $\Zinv$ associated with pseudo-parallelizations.
As in Chapter~\ref{chaprat} and Section~\ref{sectransrat}, the corresponding discrete definition of $\Zinv$ are justified by using nonhomogeneous propagating forms $\varepsilon$-dual (as in Definition~\ref{defformdual}) to these propagating chains.

Let $\Aman$ be an oriented $3$-manifold with possible boundary, equipped with a pseudo-parallelization $\tilde{\tau}=(N(\gamma);\tau_e,\tau_b)$ as in Definition~\ref{defpseudotriv}.
Let $\omega_s$ be a $2$-form of $S^2$ invariant under the rotations around the vertical axis such that $\int_{S^2}\omega_s=1$.
Let $\omega_i$ be a $2$-form of $S^2$ such that $\int_{S^2}\omega_i=1$.
Let $\eta_{i,s,1}$ be a $1$-form of $S^2$ such that $\omega_i = \omega_s +d\eta_{i,s,1}$. Let $\varepsilon \in \left]0,1/2\right[$ be the small positive number of Definitions~\ref{defpseudotriv} and \ref{defhombounformpseudo}. Let $\varepsilon_i \in \left]0,\varepsilon/2\right[$, and let $k$ be a large integer greater than $3$. Let $p(\tau_b)$ denote the projection from $UN(\gamma)$ to $S^2$ induced by $\tau_b$. We have \begin{equation*}p(\tau_b)\bigl(\tau_b(t,c,u;\cvarM \in S^2)\bigr)=\cvarM .\end{equation*}

\begin{lemma}\label{lemprelpseudoeta}
 Under the hypotheses above, there exists a one-form 
$\eta_{i,s}$ on $\ST(\Aman)$ such that
\begin{itemize}
 \item $\eta_{i,s}= p(\tau_e)^{\ast}\left(\eta_{i,s,1}\right)$ on $\ST\Bigl(\Aman \setminus \bigl([a,b-\varepsilon_i+\varepsilon_i^k] \times \gamma  \times [-1,1]\bigr)\Bigr)$,
 \item $\eta_{i,s}= p(\tau_b)^{\ast}\left(\eta_{i,s,1}\right)$ on 
$[a,b]\times \gamma\times N\left(\partial [-1,1]\right) \times S^2$, 

\item $\eta_{i,s}=\frac12 \Bigl(p\bigl(\tau_b \circ \CT_{\gamma}^{-1}\bigr)^{\ast}\left(\eta_{i,s,1}\right) +p\bigl(\tau_b \circ F(\gamma,\tau_b)^{-1}\bigr)^{\ast}\left(\eta_{i,s,1}\right)\Bigr)$ on \\
$\left[a,b-\varepsilon_i-\varepsilon_i^k\right]\times \gamma  \times [-1,1] \times S^2$, and
\item $\eta_{i,s}$ pulls back through $p_{[a,b]} \times p_{[-1,1]} \times p(\tau_b)$ on
$\ST\bigl(N(\gamma)=[a,b]\times \gamma  \times [-1,1] \bigr)$.
\end{itemize}
The following addition is useless in this book but useful for studying equivariant invariants as in \cite{lesuniveq}.
Let $c$ and $d$ be two elements of $\left]-1,1\right[$ such that $c <d$.
Let $v_3$ denote the projection on the third coordinate in $\RR^3$, and let $S^2_{\left]c,d\right[}=\{ X \in S^2 \suchthat v_3(X)\in \left]c,d\right[\}$.
If $\omega_s$ and $\omega_i$ are compactly supported in $S^2_{\left]c,d\right[}$, then $\eta_{i,s,1}$ can be chosen so that it is also compactly supported in $S^2_{\left]c,d\right[}$, and $\eta_{i,s}$ can be chosen so that it is compactly supported in $[a,b]\times \gamma  \times [-1,1] \times S^2_{\left]c,d\right[}$. 
\end{lemma}
\bp The existence of $\eta_{i,s,1}$ compactly supported in $S^2_{\left]c,d\right[}$, when $\omega_s$ and $\omega_i$ are compactly supported in $S^2_{\left]c,d\right[}$, comes from the fact that $H^2_c(S^2_{\left]c,d\right[})=\RR$, according to Theorem~\ref{thmDeRhamcompact}. 

Extend $\eta_{i,s}$ to \begin{equation*}C (\times \gamma) \times S^2=\left[b-\varepsilon_i-\varepsilon_i^k,b-\varepsilon_i+\varepsilon_i^k\right] (\times \gamma)  \times [-1+\frac{\varepsilon}{2}, 1-\frac{\varepsilon}{2}] \times S^2\end{equation*} as follows.
Cover $S^2$ by three open spaces $S^2_{\left]c,d\right[}$, $N(N)=\{ X \in S^2 \suchthat v_3(X) >\frac{c+d}{2}\}$, and $N(S)=\{ X \in S^2 \suchthat v_3(X) <\frac{c+d}{2}\}$.
Pick a corresponding partition of unity $(\chi_{S^2_{\left]c,d\right[}}, \chi_{N(S)}, \chi_{N(N)})$ of functions respectively
compactly supported on $S^2_{\left]c,d\right[}$, $N(N)$, and $N(S)$ whose sum is one. Let $S$ stand for $S^2_{\left]c,d\right[}$, $N(N)$, or $N(S)$.
On the product $C (\times \gamma) \times S$   
the form $\eta_{i,s}$ is a combination of basic standard one-forms. Smoothly extend the $\RR$-valued coordinate functions of the forms from $\left(\partial (C (\times \gamma))\right) \times S$ to $C (\times \gamma) \times S$ in order to obtain $\eta_{i,s}$ with the required properties, so that $\eta_{i,s}$ is compactly supported in $[a,b]\times \gamma  \times [-1,1]  \times S^2_{\left]c,d\right[}$ when $\eta_{i,s,1}$ is compactly supported in $S^2_{\left]c,d\right[}$.
\eop

Recall the notation of Definition~\ref{defhombounformpseudo} and extend
$\tau_e$ to \begin{equation*}U\left([b-\varepsilon_i+\varepsilon_i^k,b]  \times \gamma  \times [-1,1]\right),\end{equation*} so that it coincides with $\tau_b$, there. 

\begin{definition}
\label{defpseudoformprel} Under the hypotheses of the beginning of this section,
define $\omega=\omega(\tilde{\tau},\omega_i,k,\varepsilon_i,\eta_{i,s})$ on $\ST A$ to satisfy
\begin{equation*}
\omega=p(\tau_e)^{\ast}\left(\omega_i\right)  \mbox{ on } \ST\Bigl(\Aman \setminus \bigl([a,b-\varepsilon_i+\varepsilon_i^k] \times \gamma  \times [-1,1]\bigr)\Bigr),\end{equation*}
\begin{multline*}\omega=\frac12 \Bigl(p(\tau_b \circ \CT_{\gamma}^{-1})^{\ast}\left(\omega_i\right) + p(\tau_b \circ F(\gamma,\tau_b)^{-1})^{\ast}\left(\omega_i\right)\Bigr)\\ \mbox{on }\ST\bigl([a,b-\varepsilon_i-\varepsilon_i^k] \times \gamma  \times [-1,1]\bigr),\end{multline*}
and 
\begin{equation*}
\omega=
p(\tau_b)^{\ast}\left(\omega_s\right) + d \eta_{i,s}\mbox{ on }\ST\bigl([b-\varepsilon_i-2\varepsilon_i^k,b-\varepsilon_i+2\varepsilon_i^k]  \times \gamma  \times \left]-1,1\right[\bigr).
\end{equation*} 
\end{definition}

\begin{lemma}
\label{lemappprimtwist}
Definition~\ref{defpseudoformprel} of  $\omega(\tilde{\tau},\omega_i,k,\varepsilon_i,\eta_{i,s})$ is consistent.
Furthermore, we have $\omega(\tilde{\tau},\omega_i,k,\varepsilon_i,\eta_{i,s})=\omega(\tilde{\tau},\omega_s,k,\varepsilon_i,0)+ d\eta_{i,s}$ on $\ST A$.
\end{lemma}
\bp
Any $2$-form $\omega_s$ invariant under the rotations around the vertical axis may be expressed as $(\kappa \circ v_3)\omega_{S^2}$. Therefore, when $\eta_{i,s}$ is zero, Lemma~\ref{lemprimtwist} ensures that the definitions of $\omega(\tilde{\tau},\omega_s,k,\varepsilon_i,0)$  match on $[b-\varepsilon_i-2\varepsilon_i^k,b-\varepsilon_i-\varepsilon_i^k]  \times \gamma  \times \left]-1,1\right[$. Conclude by observing that $\omega(\tilde{\tau},\omega_i,k,\varepsilon_i,\eta_{i,s})=\omega(\tilde{\tau},\omega_s,k,\varepsilon_i,0)+ d\eta_{i,s}$ on $\ST N(\gamma)$.
\eop

\begin{definition}[General boundary form associated with a pseudo-parallelization $\tilde{\tau}$] 
\label{defpseudoform}
Let $\Aman$ be an oriented $3$-manifold equipped with a pseudo-parallelization $\tilde{\tau}$.
A \emph{boundary form} of $(\Aman,\tilde{\tau})$ is a $2$-form on $\ST \Aman$ that may be expressed as $\omega(\tilde{\tau},\omega_i,k,\varepsilon_i,\eta_{i,s})$ for some $(\omega_i,k,\varepsilon_i,\eta_{i,s})$ as in Definition~\ref{defpseudoformprel}.
When $\eta_{i,s,1}$, $\omega_s$, and $\omega_i$ are compactly supported in $S^2_{\left]c,d\right[}$, and when the intersection of the support of $\eta_{i,s}$ with $[b-\varepsilon_i-2\varepsilon_i^k,b-\varepsilon_i+2\varepsilon_i^k]\times \gamma  \times [-1,1] \times S^2$ is a compact subspace of $[b-\varepsilon_i-2\varepsilon_i^k,b-\varepsilon_i+2\varepsilon_i^k]\times \gamma  \times [-1,1] \times S^2_{\left]c,d\right[}$, the boundary form is said to be \emph{adapted to} $S^2_{\left]c,d\right[}$.
\end{definition}

Note that Definition~\ref{defpseudoform} coincides with Definition~\ref{defhombounformpseudo} of a homogeneous boundary form when $\omega_i=\omega_s=\omega_{S^2}$ and $\eta_{i,s}=0$.

\begin{lemma}
 \label{lempropintsplusstau}
Under the assumptions of Proposition~\ref{propintsplusstau}, for any boundary form $\omega$ of $(M,\tilde{\tau})$ as in Definition~\ref{defpseudoform}, we have
\begin{equation*}\int_{s_+(\Sigma)}\omega=\frac12\chi\bigl(\tilde{\tau}(. \times e_2)\vert_{\partial \Sigma};\Sigma\bigr)\end{equation*}
and \begin{equation*}\int_{s_-(\Sigma)}\omega=-\frac12\chi\bigl(\tilde{\tau}(. \times e_2)\vert_{\partial \Sigma};\Sigma\bigr).\end{equation*}
\end{lemma}
\bp Thanks to Lemma~\ref{lemsplusstau}, Lemma~\ref{lempropintsplusstau} holds when $\Sigma$ does not meet $N(\gamma)$.
Therefore, as in the proof of Proposition~\ref{propintsplusstau}, it suffices to prove Lemma~\ref{lempropintsplusstau} when $\Sigma$ is a meridian of the link $\gamma$ of the pseudo-parallelization $\tilde{\tau}$. Let us treat this case. 
When $\omega_i=\omega_s$ and $\eta_{i,s}=0$, the proof of Proposition~\ref{propintsplusstau} applies. In general, let $\tilde{\omega}$ be the form obtained from $\omega$ by changing $\omega_i$ to $\omega_s$ and $\eta_{i,s}$ to $0$.
The form $\omega$ may be written as $\tilde{\omega}+ d\eta_{i,s}$ on $\ST N(\gamma)$, where $\eta_{i,s}$ is expressed as $p(\tau_e)^{\ast}\left(\eta_{i,s,1}\right)$ along $s_+(\Sigma)\vert_{\partial \Sigma}$. We have $\int_{s_+(\Sigma)}\omega-\int_{s_+(\Sigma)}\tilde{\omega}=\int_{s_+(\Sigma)\vert_{\partial \Sigma}}p(\tau_e)^{\ast}\left(\eta_{i,s,1}\right)$. This is zero 
since $p(\tau_e)$ maps $s_+(\Sigma)\vert_{\partial \Sigma}$ to a point.
We similarly get $\int_{s_-(\Sigma)}\omega=\int_{s_-(\Sigma)}\tilde{\omega}$.
\eop

Theorem~\ref{thmznponepseudo} generalizes as follows to these boundary forms.

\begin{theorem}
\label{thmznponepseudogen}
Let $\Aman$ be a compact $3$-manifold that embeds in a rational homology $3$-ball.
Assume that $\Aman$ is equipped with two pseudo-parallelizations $\tau_0$ and $\tau_1$ that coincide with a common genuine parallelization along a regular neighborhood of $\partial \Aman$. 
For $i\in \underline{3n}$, let $\omega_{0,i}(\tau_0)$ be a boundary form of $(\Aman,\tau_0)$ and let $\omega_{1,i}(\tau_1)$ be a boundary form of $(\Aman,\tau_1)$.
There exists a closed 
$2$-form $\omega(i)$ on $\left[0,1\right] \times \ST \Aman$ that restricts 
\begin{itemize}
 \item to $\{0\} \times \ST \Aman$ as $\omega_{0,i}(\tau_0)$
\item to $\{1\} \times \ST \Aman$ as $\omega_{1,i}(\tau_1)$,
\item to $\left[0,1\right] \times \ST \Aman\vert_{\partial \Aman}$ as $(\id_{\left[0,1\right]} \times p_{\tau_0})^{\ast}\left(\omega_{S,i}\right)$ with respect to a closed $2$-form $\omega_{S,i}$ on $\left[0,1\right] \times S^2$ and to the projection $\id_{\left[0,1\right]} \times p_{\tau_0} \colon \left[0,1\right] \times \ST \Aman\vert_{\partial \Aman} \to \left[0,1\right] \times S^2$.
\end{itemize}
Let $n$ be a natural integer.
As in Corollary~\ref{corinvone}, set \begin{equation*}z_n\Bigl(\left[0,1\right] \times \ST A;\bigl(\omega(i)\bigr)_{i\in \underline{3n}}\Bigr)=\sum_{\Gamma \in \Davis^{c}_n}\coefgambet_{\Gamma}\int_{\left[0,1\right]\times \cinjuptd{\finsetv(\Gamma)}{T A}}\bigwedge_{e \in E(\Gamma)}p_e^{\ast}\Bigl(\omega\bigl(j_E(e)\bigr)\Bigr)\left[\Gamma\right].\end{equation*}
Then we have
\begin{equation*}z_n\Bigl(\left[0,1\right] \times \ST A;\bigl(\omega(i)\bigr)_{i\in \underline{3n}}\Bigr)=\frac{p_1(\tau_0,\tau_1)}{4}\ansothree_n .\end{equation*}
\end{theorem}
\bp The existence of $\omega_{S,i}$ is proved in Lemma~\ref{lemformprod}.
The proof of the existence of the form $\omega(i)$ with its prescribed properties is obtained
from the proof of Proposition~\ref{propextexistspseudo} by replacing Proposition~\ref{propintsplusstau} with Lemma~\ref{lempropintsplusstau}.

Let $t$ stand for the coordinate in $\left[0,1\right]$.
Assume that $\omega_{0,i}(\tau_0)$ is a homogeneous boundary form and that we have $\tau_0=\tau_1$. Choose $\omega_{S,i}=\omega_{S^2} + d (t\eta_{i,s,1})$
with the notation before Lemma~\ref{lemprelpseudoeta}.
Then the form $\omega(i)$ can be chosen so that it is equal to
\begin{equation*}
p(\tau_e)^{\ast}\left(\omega_{S,i}\right)\mbox{ on }\left[0,1\right] \times \ST\Bigl(\Aman \setminus \bigl([a,b-\varepsilon_i+\varepsilon_i^k] \times \gamma  \times [-1,1]\bigr)\Bigr),
\end{equation*}
\begin{equation*}
p(\tau_b)^{\ast}\left(\omega_{S^2}\right) + d (t\eta_{i,s}) \mbox{ on }\left[0,1\right] \times \ST \bigl([b-\varepsilon_i-\varepsilon_i^k,b-\varepsilon_i+\varepsilon_i^k]  \times \gamma  \times \left]-1,1\right[\bigr) \mbox{, and}
\end{equation*}
\begin{equation*}\frac12\Bigl({p\bigl(\tau_b \circ \CT_{\gamma}^{-1}\bigr)^{\ast}\left(\omega_{S,i}\right)} + {p\bigl(\tau_b \circ F(\gamma,\tau_b)^{-1}\bigr)^{\ast}\left(\omega_{S,i}\right)}\Bigr)\end{equation*}
on $\left[0,1\right] \times \ST \bigl([a,b-\varepsilon_i-\varepsilon_i^k] \times \gamma  \times [-1,1]\bigr)$.
In this case, the parts over $\left(\Aman \setminus \left([a,b] \times \gamma  \times [-1,1]\right)\right)$ in $z_n(\left[0,1\right] \times \ST \Aman;(\omega(i))_{i\in \underline{3n}})$ cancel because the forms do not depend on the factor $\Aman$ (which is factored out via the parallelization $\tau_e$), and the parts
over $\left([a,b] \times \gamma  \times [-1,1]\right)$ also cancel because the forms do not depend on
the factor $\gamma$.
So we simply find $z_n(\left[0,1\right] \times \ST A;(\omega(i))_{i\in \underline{3n}})=0$ as announced.

In general, as in the beginning of the proof of Proposition~\ref{propextdefzn}, and as in Corollary~\ref{corinvone}, the sum
$z_n(\left[0,1\right] \times \ST A;(\omega(i)))$ depends only on the restriction of the $\omega(i)$ to $\partial\left({\left[0,1\right] \times \ST \Aman}\right)$. 
The above arguments can also be used to prove that $z_n(\left[0,1\right] \times \ST A;(\omega(i))_{i\in \underline{3n}})$ is independent of the forms $\omega_{S,i}$. So $z_n(\left[0,1\right] \times \ST A;(\omega(i))_{i\in \underline{3n}})$ depends only on the pairs $(\omega_{0,i}(\tau_0),\omega_{1,i}(\tau_1))$. Denote it by $z_n\left((\omega_{0,i}(\tau_0),\omega_{1,i}(\tau_1))_{i\in \underline{3n}}\right)$. Observe 
\begin{multline*}z_n\Bigl(\bigl(\omega_{0,i}(\tau_0),\omega_{2,i}(\tau_2)\bigr)_{i\in \underline{3n}}\Bigr)=z_n\Bigl(\bigl(\omega_{0,i}(\tau_0),\omega_{1,i}(\tau_1)\bigr)_{i\in \underline{3n}}\Bigr) \\+ z_n\Bigl(\bigl(\omega_{1,i}(\tau_1),\omega_{2,i}(\tau_2)\bigr)_{i\in \underline{3n}}\Bigr)\end{multline*}
again.
Conclude with the study of the special case and with Theorem~\ref{thmznponepseudo}.
\eop

\begin{definition}
 \label{defpseudopropform}
Let $\crats$ be a rational homology $\RR^3$ equipped with an asymptotically standard pseudo-parallelization $\tilde{\tau}$.
A \emph{propagating form}\index[T]{propagating!form} of $(C_2(\rats),\tilde{\tau})$ is a propagating form of $C_2(\rats)$ (as in Definition~\ref{defpropagatortwo}) that coincides with a boundary form (of Definition~\ref{defpseudoform}) of $(\crats,\tilde{\tau})$ on $\ST \crats$.
\end{definition}

\begin{lemma} \label{lemIthetapseudo}
 Let $\crats$ be a rational homology $\RR^3$ equipped with an asymptotically standard pseudo-parallelization $\tau=(N(\gamma);\tau_e,\tau_b)$.
Let $K \colon \source \hookrightarrow \crats \setminus N(\gamma)$ be a knot embedding, which is straight with respect to $\tau\vert_{\crats \setminus N(\gamma)}$, as in Definition~\ref{defstraightlink}. Let $K_{\parallel,\tau}$ denote the parallel of $K$ induced by $\tau$, as in Lemma~\ref{lemdegloctheta}.
For any propagating form $\omega_p$ of $\bigl(C_2(\rats),\partau\bigr)$, we have
\begin{equation*}\int_{C(K;\onechordsmall{$S^1$})}\omega_p=I_{\theta}(K,\partau)=lk(K,K_{\parallel,\partau,Y}).\end{equation*}
\end{lemma}
\bp Lemma~\ref{lemdegloctheta} extends to this case. See also Lemmas~\ref{lemvaritheta} and \ref{lemappprimtwist}.
\eop

Theorem~\ref{thmconststraight} generalizes as follows to pseudo-parallelizations:

\begin{theorem}
\label{thmconststraightpseudo}
Let $\crats$ be a rational homology $\RR^3$ equipped with an asymptotically standard pseudo-parallelization $\tau=(N(\gamma);\tau_e,\tau_b)$.
Let $\Link \colon \source \hookrightarrow \crats \setminus N(\gamma)$ be a link embedding, which is straight with respect to $\tau\vert_{\crats \setminus N(\gamma)}$. Let $\Link_{\parallel,\tau}$ denote the parallel of $\Link$ induced by $\tau$.
For any $i\in \underline{3n}$, let $\omega(i)$ be a propagating form of $(C_2(\rats),\tau)$.
Then \begin{equation*}\Zinv^s_n\Bigl(\crats,\Link,\bigl(\omega(i)\bigr)_{i \in \underline{3n}}\Bigr)=\sum_{\Gamma \in \Davis^e_n(\source)}\coefgambet_{\Gamma}I\Bigl(\rats,\Link,\Gamma,\bigl(\omega(i)\bigr)_{i \in \underline{3n}}\Bigr)\left[\Gamma\right] \in \Aavis_n(\source).\end{equation*}
is independent of the chosen
$\omega(i)$. Set \begin{equation*}\Zinv^s_n(\crats,\Link,\tau)=\Zinv^s_n\Bigl(\crats,\Link,\bigl(\omega(i)\bigr)_{i \in \underline{3n}}\Bigr)\end{equation*} and $\Zinv^s(\crats,\Link,\tau)=\left(\Zinv^s_n(\crats,\Link,\tau)\right)_{n \in \NN}$. 
Then we have \begin{equation*}\Zinv^s(\crats,\Link,\tau)=\Zinv(\crats,\Link,\tau) =\Zinvufrfneg(\crats,\Link,\Link_{\parallel,\tau}) \exp\left(\frac{p_1(\tau)}{4}\ansothree\right),\end{equation*}
where $\Zinv(\crats,\Link,\tau)$ is defined in Theorem~\ref{thmfstconsttangpseudo} (as in Theorem~\ref{thmfstconst}),
with Definition~\ref{defponepseudotrivpone} for $p_1(\tau)$, and Definition~\ref{defZinvufrflink} for $\Zinvufrfneg(\crats,\Link,\Link_{\parallel,\tau})$.
\end{theorem}

Theorem~\ref{thmconststraightpseudo} is a consequence of the following generalization
of Theorem~\ref{thmtangconstcompar} to pseudo-parallelizations:

\begin{theorem}
\label{thmtangconstcomparpseudo}
Let $\hcylc$ be a rational homology cylinder equipped with a pseudo-parallelization $\tau=(N(\gamma);\tau_e,\tau_b)$.
Let $\tanghcyll \colon \sourcetl \hookrightarrow \rats(\hcylc) \setminus N(\gamma)$ be a long tangle representative associated to a tangle in $\hcylc$. Let $\{K_j\}_{j \in I}$ be the set of components of $\tanghcyll$.
Assume that the bottom (resp. top) configuration of $\tanghcyll$ is represented by a map $\confy^- \colon \finsetb^- \to \drad{1}$ (resp. $\confy^+ \colon \finsetb^+ \to \drad{1}$). 

Let $\largen \in \NN$.
For $i\in \underline{3\largen}$, let $\tilde{\omega}(i,S^2)=(\tilde{\omega}(i,t,S^2))_{t\in \left[0,1\right]}$ be a closed $2$-form on $\left[0,1\right] \times S^2$ such that $\int_{S^2}\tilde{\omega}(i,0,S^2)=1$.
There exists a closed $2$-form $\tilde{\omega}(i)=(\tilde{\omega}(i,t))_{t\in \left[0,1\right]}$ on $\left[0,1\right] \times C_2(\rats(\hcylc))$ such that \begin{equation*}\tilde{\omega}(i)=(\id_{\left[0,1\right]} \times p_{\tau_e})^{\ast}\left(\tilde{\omega}(i,S^2)\right)\end{equation*} on
$\left[0,1\right]\times \left(\partial C_2(\rats(\hcylc)) \setminus \ST N(\gamma) \right)$ and
the restriction of $\tilde{\omega}(i,t)$ to $\ST \crats(\hcylc)$ is a boundary form of $(\crats(\hcylc),{\tau})$ as in Definition~\ref{defpseudoform} for all $t \in \left[0,1\right]$.

For such a family $(\tilde{\omega}(i))_{i\in \underline{3\largen}}$, and for a subset $\finseta$ of $\underline{3\largen}$ with cardinality $3k$, set
\begin{equation*}\Zinv\Bigl(\hcylc,\tanghcyll,\tau,\finseta,\bigl(\tilde{\omega}(i,t)\bigr)_{i \in \finseta}\Bigr)=\sum_{\Gamma \in \Davis^e_{k,\finseta}(\sourcetl)}\coefgambet_{\Gamma}I\Bigl(\hcylc,\tanghcyll,\Gamma,\bigl(\tilde{\omega}(i,t)\bigr)_{i \in \finseta}\Bigr)\left[\Gamma\right] \in \Aavis_{k}(\sourcetl).\end{equation*}
Then  \begin{equation*}\Zinv(\hcylc,\tanghcyll,\tau,\finseta)(t)=\Zinv\Bigl(\hcylc,\tanghcyll,\tau,\finseta,\bigl(\tilde{\omega}(i,t)\bigr)_{i \in \finseta}\Bigr)\end{equation*} depends only on $(\tilde{\omega}(i,t,S^2))_{i \in \finseta}$ for any $t$ (and on $(\hcylc,\tanghcyll,\tau)$). 

It will be denoted by $\Zinv(\hcylc,\tanghcyll,\tau,\finseta,(\tilde{\omega}(i,t,S^2))_{i \in \finseta})$.
When $\tilde{\omega}(i,t,S^2)$ is the standard homogeneous form $\omega_{S^2}$ on $S^2$ for any $i$ and for all $t \in \left[0,1\right]$, the map
$\Zinv(\hcylc,\tanghcyll,\tau,.)(t)$ sends any subset of $\underline{3\largen}$ with cardinality $3k$ to $\Zinv_k(\hcylc,\tanghcyll,\tau)$.

Furthermore, with the notation of Definition~\ref{defprodbiz}, for any given orientation of $\tanghcyll$, we have \begin{multline*}\Zinv(\hcylc,\tanghcyll,\tau,.)(t)=\\
\biggl(\Bigl(\prod_{j\in I}\biglol{\left[0,t\right]}{\eta(.,\ST^+_j)} \#_j\Bigr)
\lol{\left[t,0\right] \times \confy^-}{\eta_{\finsetb^-,.}}\Zinv(\hcylc,\tanghcyll,\tau,.)(0)\lol{\left[0,t\right] \times \confy^+}{\eta_{\finsetb^+,.}} \biggr)_{\sqcup},\end{multline*}
where $\ST^+_j=p_{\tau}(\ST^+K_j)$.
\end{theorem}
\bp Fix $i \in \underline{3\largen}$. Pick a one-form $\tilde{\eta}_{i,s,1}$  on $\left[0,1\right] \times S^2$ such that $\tilde{\omega}(i,S^2)=p_{S^2}^{\ast}\left(\omega_{S^2}\right)+d\tilde{\eta}_{i,s,1}$. Use this form to construct a one-form
$\tilde{\eta}_{i,s}(t)$ on $\left[0,1\right] \times \ST\left([a,b]\times \gamma \times [-1,1]\right)$ that pulls back through $p_{\left[0,1\right]} \times p_{[a,b]} \times p_{[-1,1]} \times p(\tau_b)$ and such that its restriction to $\{t\} \times \ST(N(\gamma ))$ satisfies the conditions of Lemma~\ref{lemprelpseudoeta}, with respect to $\tilde{\eta}_{i,s,1}\vert_{\{t\} \times \ST(N(\gamma))}$. 

Next use $\tilde{\eta}_{i,s}$ to construct the restriction of $\tilde{\omega}(i)$ to $\left[0,1\right]\times \partial C_2(\rats(\hcylc))$, such that $\tilde{\omega}(i,t)\vert_{\{t\} \times \ST \crats(\hcylc)}= \omega({\tau},\tilde{\omega}(i,t,S^2),k,\varepsilon_i,\tilde{\eta}_{i,s}(t))$. The existence of $\tilde{\omega}(i)$ follows, as in Lemma~\ref{lemformprod}.

Use the proof of Theorem~\ref{thmtangconstcompar} to prove the variation formula that expresses $\Zinv(\hcylc,\tanghcyll,\tau,\finseta,(\tilde{\omega}(i,t))_{i \in \finseta})$ as a function of $\Zinv(\hcylc,\tanghcyll,\tau,\finseta,(\tilde{\omega}(i,0))_{i \in \finseta})$ for forms $\tilde{\omega}(i,t)$ as in the statement.

Then note that if $\tilde{\omega}(i,0,S^2)$ is the standard homogeneous form $\omega_{S^2}$ on $S^2$ for any $i$, then
$\Zinv(\hcylc,\tanghcyll,\tau,.)(0)$ maps any subset of $\underline{3\largen}$ with cardinality $3k$ to $\Zinv_k(\hcylc,\tanghcyll,\tau)$, by Theorem~\ref{thmfstconsttangpseudo}.
Thanks to Lemma~\ref{lemformprod}, this shows that $\Zinv(\hcylc,\tanghcyll,\tau,\finseta)(t)$ depends only on the forms $\tilde{\omega}(i,S^2)$, and therefore only on the forms $\tilde{\omega}(i,t,S^2)$.
\eop

\bpo{Proof of Theorem~\ref{thmconststraightpseudo}} 
Apply the formula of Theorem~\ref{thmtangconstcomparpseudo}.
There are no factors $\lol{\left[t,0\right] \times \confy^-}{\eta_{\finsetb^-,.}}$ and $\lol{\left[0,t\right] \times \confy^+}{\eta_{\finsetb^+,.}}$ because $\tanghcyll$ is a link.
Since $\tanghcyll$ is a straight link with respect to $\tau$, the factors $\lol{\left[0,t\right]}{\eta(.,p_{\tau}(\ST^+K_j)}$ also vanish, even though the notion of straight link of Definition~\ref{defstraightlink} is less restrictive than the notion of straight tangle of Definition~\ref{defstraighttang}. Theorem~\ref{thmconststraightpseudo} follows from Theorems~\ref{thmtangconstcomparpseudo} and~\ref{thmfstconsttangpseudo}, and Lemma~\ref{lemIthetapseudo}.\eop

\begin{remark}
When $\tilde{\omega}(i,1,S^2)$ is the standard homogeneous form $\omega_{S^2}$ on $S^2$ for any $i$,
the variation formula of Theorem~\ref{thmtangconstcomparpseudo} yields alternative expressions of $\Zinv$.
\end{remark}

\section[Pseudo-sections associated with pseudo-parallelizations]{Pseudo-sections associated with \texorpdfstring{\\}{} pseudo-parallelizations}
\label{secpseudosec}

For an asymptotically standard parallelization $\tau$ of a punctured rational homology $3$-sphere $\crats$, a propagating chain of $(C_2(\rats),\tau)$ is defined in Definition~\ref{defpropagatorone}
to be a $4$-chain $\propP$ of $C_2(\rats)$
such that $\partial \propP =p_{\tau}^{-1}(\cvarM)$ for some $\cvarM \in S^2$, 
where $p_{\tau}^{-1}(\cvarM) \subset\partial C_2(\rats)$, the intersection $p_{\tau}^{-1}(\cvarM) \cap \left(\partial C_2(\rats) \setminus \ST\crats \right)$ is independent of $\tau$, and \begin{equation*}p_{\tau}^{-1}(\cvarM) \cap  \ST \crats = \tau(\crats \times \{\cvarM\}).\end{equation*}

The image in $\ST \Aman$ of the restriction of a section $\tau(\crats \times \cvarM)=\tau(\crats \times \{\cvarM\})$ to a part $\Aman$ of $\crats$ is denoted by $s_{\tau}(\Aman;\cvarM)$.

In this section, we define pseudo-sections $s_{\tilde{\tau}}(\crats;\cvarM)$ associated with pseudo-parallelizations $\tilde{\tau}$. A \emph{propagating chain} of $(C_2(\rats),\tilde{\tau})$ is a $4$-chain $\propP$ of $C_2(\rats)$
such that \begin{equation*}\partial \propP =\bigl(p_{\tau}\vert_{\partial C_2(\rats) \setminus \ST\crats}\bigr)^{-1}(\cvarM) \cup s_{\tilde{\tau}}(\crats;\cvarM)\end{equation*} for some $\cvarM \in S^2$. Thus, the pseudo-sections $s_{\tilde{\tau}}(\crats;\cvarM)$ will play the same role as the sections $s_{\tau}(\crats;\cvarM)$ in the more flexible definition of $\Zinvuf$ presented in Section~\ref{secdefZpseudosec}.

\begin{definition}[Pseudo-sections $s_{\tilde{\tau}}(.;\cvarM)$]
\label{defpseudosec} 
Recall the map ${F}(\gamma,\tau_b)$ of Definition~\ref{defhombounformpseudo} and the notation of Definition~\ref{defpseudotriv}.

Let $\cvarM \in S^2$. Let $S^1(\cvarM)$ be the intersection of $S^2$ with the plane orthogonal to the axis generated by $(0,0,1)$ that contains $\cvarM$. So $S^1(\cvarM)$ is a circle or a point.
Let $\chainproppseudo{\cvarM}$ be a $2$-dimensional chain in $[-1,1] \times S^1(\cvarM)$ whose
boundary is 
$\{\left(u,\rho_{-\alpha(u)}(\cvarM)\right), u \in[-1,1]\}
+\{\left(u,\rho_{\alpha(u)}(\cvarM)\right), u \in[-1,1]\} - 2[-1,1]\times \{\cvarM\}$,
as in Figure~\ref{figchainproppseudo}.

\bfig
\centering
\begin{tikzpicture} \useasboundingbox (0,0) rectangle (1.2,1.2);
\fill[gray] (0,1.2) -- (1.2,1.2) -- (.6,.6) --  (0,1.2);
\draw[pattern= {horizontal lines}] (0,0) -- (1.2,0) -- (.6,.6) --  (0,0);
\draw [->]  (0,0) -- (0,1.2) (1.2,.6) -- (1.2,1.2) (1.2,0) -- (1.2,.6);
\draw [very thick,->] (.6,0) -- (0,0) (1.2,0) -- (.6,0);
\draw [very thick,->] (.6,1.2) -- (0,1.2) (1.2,1.2) -- (.6,1.2);
\draw (0,1.2) -- (1.2,0) (1.2,1.2) -- (0,0);
\draw [->] (0,0) -- (.3,.3);
\draw [->] (.6,.6) -- (.9,.9);
\draw [->] (0,1.2) -- (.3,.9);
\draw [->] (.6,.6) -- (.9,.3);
\draw (0,0) node[left]{\scriptsize $-1$} (1.2,0) node[right]{\scriptsize $1$} (1.2,.6) node[right]{\scriptsize $S^1(\cvarM)$};
\end{tikzpicture}
\caption{The chain $\chainproppseudo{\cvarM}$ in $[-1,1] \times S^1(\cvarM)$}\label{figchainproppseudo}
\end{figure}
Then $s_{\tilde{\tau}}(\Aman;\cvarM)$ is the following $3$-cycle of $(\ST \Aman, \ST \Aman\vert_{\partial \Aman})$
\begin{multline*}s_{\tilde{\tau}}(\Aman;\cvarM)=  s_{\tau_e}\bigl(\Aman\setminus \mathring{N}(\gamma);\cvarM\bigr) \\ +
\frac12\Bigl(s_{\tau_b \circ \CT_{\gamma}^{-1}}\bigl(N(\gamma);\cvarM\bigr) + s_{\tau_b \circ F(\gamma,\tau_b)^{-1}}\bigl(N(\gamma);\cvarM\bigr) + \{b\} \times \gamma \times \chainproppseudo{\cvarM}\Bigr),
\end{multline*}
where $\ST \Aman\vert_{\{b\} \times \gamma \times [-1,1]}$ is identified with $\{b\} \times \gamma \times [-1,1] \times S^2$ by $\tau_b$ (or in the same way by $\tau_e$).

We also introduce small deformations of these sections, associated with $\varepsilon_i$ such that $0\leq\varepsilon_i < \varepsilon$ (with respect to the $\varepsilon$ of Definition~\ref{defpseudotriv}) as follows.
Let $N(\gamma,\varepsilon_i)=\left[a,b-\varepsilon_i\right] \times \gamma \times [-1,1] \times S^2$. Extend $\tau_e$ over $[b-\varepsilon,b] \times \gamma \times [-1,1]$ so that it coincides with $\tau_b$, there.
Then  
$s_{\tilde{\tau}}(\Aman;\cvarM,\varepsilon_i)$ is the following $3$-cycle of $(\ST \Aman, \ST \Aman\vert_{\partial \Aman})$
\begin{equation*}\begin{array}{ll}s_{\tilde{\tau}}(\Aman;\cvarM,\varepsilon_i)=&s_{\tau_e}\bigl(\Aman\setminus \mathring{N}(\gamma,\varepsilon_i);\cvarM\bigr)\\& + 
\frac12\Bigl(s_{\tau_b \circ \CT_{\gamma}^{-1}}\bigl(N(\gamma,\varepsilon_i);\cvarM\bigr) +  s_{\tau_b \circ F(\gamma,\tau_b)^{-1}}\bigl(N(\gamma,\varepsilon_i);\cvarM\bigr) \Bigr)
\\&+ 
\frac12  \{b-\varepsilon_i\} \times \gamma \times \chainproppseudo{\cvarM}.\end{array}
\end{equation*}

\end{definition}

When $\Sigma$ is a $2$-chain that intersects $N(\gamma)$ along sections $N_c(\gamma)=[a,b]\times\{c\}\times [-1,1]$ (which are oriented as meridian disks of $(-\gamma)$), set \begin{equation*}s_{\tilde{\tau}}(\Sigma;\cvarM)=s_{\tilde{\tau}}(\Aman;\cvarM)\cap \ST \Aman\vert_{\Sigma}.\end{equation*} So $s_{\tilde{\tau}}(N_c(\gamma);\cvarM)$ equals
\begin{equation*}\frac12\Bigl(s_{\tau_b \circ \CT_{\gamma}^{-1}}\bigl(N_c(\gamma);\cvarM\bigr) + s_{\tau_b \circ F(\gamma,\tau_b)^{-1}}\bigl(N_c(\gamma);\cvarM\bigr) - \{b\} \times \{c\} \times \chainproppseudo{\cvarM}\Bigr).
\end{equation*}
Note that $\chainproppseudo{e_3}$ lies in $[-1,1]\times \{e_3\}$. So 
$s_{\tau}(\Aman;\pm e_3)$ is simply given by
\begin{multline*}s_{\tau}(\Aman;\pm e_3)=s_{\tau_e}\bigl(\Aman\setminus \mathring{N}(\gamma);\pm e_3\bigr) \\+ 
\frac12\Bigl({s_{\tau_b \circ \CT_{\gamma}^{-1}}\bigl(N(\gamma);\pm e_3\bigr) + s_{\tau_b \circ F(\gamma,\tau_b)^{-1}}\bigl(N(\gamma);\pm e_3\bigr)}\Bigr).
\end{multline*}

Below, we discuss common properties of homology classes of sections and pseudo-sections.

\begin{lemma}
\label{lemdegstwoso}
Let $\qvarM \in S^2$.
 Let $\Phi$ be a map from the unit disk $D^2$ of $\CC$ to $SO(3)$ such that $\Phi(\exp(i\beta))$ is the rotation $\rho_{2\beta,\qvarM}$ of angle $2\beta$ whose axis is directed by $\qvarM$ for $\beta \in \left[0,2\pi\right]$.
Then the map \begin{equation*}\begin{array}{llll}\Phi_{\qvarM}=\Phi(\cdot)(\qvarM) \colon &D^2  &\rightarrow & S^2\\
  &z&\mapsto &\Phi(z)(\qvarM)
 \end{array}\end{equation*}
sends $\partial D^2$ to $\qvarM$, and the degree of the induced map
from $D^2/\partial D^2$ to $S^2$ is $(-1)$.
\end{lemma}
\bp
First note that the above degree does not depend on $\Phi$ on the interior of $D^2$ since $\pi_2(SO(3))=0$.
View the restriction of $\Phi_{\qvarM}$ to $\partial D^2$, as the path composition of the maps $\left(\beta \in \left[0,2\pi\right] \mapsto  \rho_{\beta,\qvarM} \right)$ and the inverse of $\left(\beta \in \left[0,2\pi\right] \mapsto \rho_{\beta,-\qvarM} \right)$ (which is homotopic to the first map).
Consider an arc $\xi$ of a great circle of $S^2$ from $-\qvarM$ to $\qvarM$. Then $\Phi$ can be regarded as the map from $\xi \times \left[0,2\pi\right]$ to $SO(3)$ that maps $(\cvarM,\beta)$ to $\rho_{\beta,\cvarM}$. So
the only preimage of $-\qvarM$ under $\Phi_{\qvarM}$ is $(\cvarM_0,\pi)$, where $\cvarM_0 \in \xi$ and $\cvarM_0 \perp \qvarM$. The local degree is easily seen to be $(-1)$.
\eop

Lemma~\ref{lemsplusstau} generalizes as follows to pseudo-parallelizations, with the notation of Definition~\ref{defeulernumb}.

\begin{proposition}
\label{propsplusplustau} 
Recall $e_3=(0,0,1) \in \RR^3$.
Let $\Sigma$ be a compact oriented surface immersed in a \/$3$-manifold $M$ equipped with a pseudo-parallelization $\tau=(N(\gamma);\tau_e,\tau_b)$. Assume that $\tau$ is an actual parallelization around the boundary \/$\partial \Sigma$ of $\Sigma$ and that $\tau(. \times e_3)$ is a positive normal to $\Sigma$ along $\partial \Sigma$.
Let \/$s_+(\Sigma) \subset \ST M$ be the section of \/$\ST M\vert_{\Sigma}$ in \/$\ST M$ associated to the positive normal field \/$n$ to \/$\Sigma$, which
coincides with \/$\tau(\partial \Sigma \times e_3)$ on $\partial \Sigma$.
Then \begin{equation*}2\bigl(s_+(\Sigma)-s_{\tau}(\Sigma; e_3)\bigr) - \chi\bigl(\tau(. \times e_2)\vert_{\partial \Sigma};\Sigma\bigr)\ST M\vert_{\ast} \end{equation*}
and 
\begin{equation*}2\bigl(s_-(\Sigma)-s_{\tau}(\Sigma; - e_3)\bigr) + \chi\bigl(\tau(. \times e_2)\vert_{\partial \Sigma};\Sigma\bigr)\ST M\vert_{\ast} \end{equation*}
are null-homologous cycles in \/$\ST M\vert_{\Sigma}$.
\end{proposition}
\bp
Lemma~\ref{lemsplusstau} gives the result for actual parallelizations.
Let us prove it for pseudo-parallelizations.
Since $s_{\tau}(\Sigma;\pm e_3)=s_{\tau}(M;\pm e_3) \cap \ST M\vert_{\Sigma}$ is well-defined for any surface $\Sigma$, as soon as the above intersection is transverse, there is no loss of generality in assuming that $\Sigma$ meets $N(\gamma)$ along sections $N_c(\gamma)$, for some $c \in \gamma$.

Since the cycles of the statement behave well under gluing (or cutting) surfaces along curves that satisfy the boundary assumptions, and since these assumptions are easily satisfied after a homotopy as in Definition~\ref{defhomotoppseudo} (when the cutting process is concerned), it suffices to prove the lemma for a meridian disk $\Sigma$ of
$\gamma$. 
For such a disk, $s_{\tau}(\Sigma; e_3)$
is the average of two genuine sections corresponding to trivializations $\tau_1$ and $\tau_2$. So, with the notation of Lemma~\ref{lemptyeulernumb}, the chain
\begin{equation*}2\bigl(s_+(\Sigma)-s_{\tau}(\Sigma; e_3)\bigr) - \biggl(\frac{d\bigl(\tau_1(. \times e_2),\partial \Sigma\bigr)+d\bigl(\tau_2(. \times e_2),\partial \Sigma\bigr)}{2}  +\chi(\Sigma)\biggr)\ST M\vert_{\ast}\end{equation*}
is a null-homologous cycle.
We conclude that the first cycle is null-homologous because the exterior trivialization $\tau_e$ satisfies
\begin{equation*}d\bigl(\tau_e(. \times e_2),\partial \Sigma\bigr) =\frac{d\bigl(\tau_1(. \times e_2),\partial \Sigma\bigr)+d\bigl(\tau_2(. \times e_2),\partial \Sigma\bigr)}{2}.\end{equation*} The second cycle can be treated similarly.
\eop

The obvious property that {\em genuine sections \/$s_{\tau}(\Sigma;\fvarM)$ and $s_{\tau}(\Sigma;\svarM)$ corresponding to distinct \/$\fvarM$ and \/$\svarM$ of \/$S^2$ are disjoint\/} generalizes as follows for pseudo-sections.

\begin{lemma}
\label{lemintpseudosec}
Let $\Sigma$ be an oriented surface embedded in a \/$3$-manifold $A$ equipped with a pseudo-parallelization $\tilde{\tau}=(N(\gamma);\tau_e,\tau_b)$ such that $\Sigma$ intersects $N(\gamma)$ only along sections 
$N_{c_i}(\gamma)=[a,b] \times \{c_i\} \times [-1,1]$ in the interior of $\Sigma$.
If $\svarM \in S^2$ and if $\tvarM \in S^2 \setminus S^1(\svarM)$, then
the algebraic intersection $\langle s_{\tilde{\tau}}(\Sigma;\svarM),s_{\tilde{\tau}}(\Sigma;\tvarM)\rangle_{\!\ST A\vert_{\Sigma}\,}$ of $s_{\tilde{\tau}}(\Sigma;\svarM)$ and $s_{\tilde{\tau}}(\Sigma;\tvarM)$ in $\ST A\vert_{\Sigma}$ is zero.
\end{lemma}
\bp Recall that $S^1(\cvarM)$ is the intersection circle of $S^2$ with the plane orthogonal to the axis generated by $e_3=(0,0,1)$ that contains $\cvarM$.
Consider the contribution to $\langle s_{\tilde{\tau}}(\Sigma;\svarM),s_{\tilde{\tau}}(\Sigma;\tvarM)\rangle_{\!\ST A\vert_{\Sigma}\,}$ of an intersection point $c$ of $\gamma$ with $\Sigma$.
Such a contribution may be expressed as $\pm \frac{1}{4} (d_Y(Z)+d_Z(Y))$ 
with 
\begin{equation*}d_Y(Z)=\Bigl\langle s_{\tau_b \circ F(\gamma,\tau_b)^{-1}}\bigl(N_c(\gamma);\tvarM\bigr),s_{\tau_b \circ \CT_{\gamma}^{-1}}\bigl(N_c(\gamma);\svarM\bigr) \Bigr\rangle \end{equation*} 
and $d_Z(Y)=\bigl\langle s_{\tau_b \circ F(\gamma,\tau_b)^{-1}}\bigl(N_c(\gamma);\svarM\bigr),s_{\tau_b \circ \CT_{\gamma}^{-1}}\bigl(N_c(\gamma);\tvarM\bigr) \bigr\rangle$.
We prove that both intersection numbers $d_Y(Z)$ and $d_Z(Y)$ are $\pm 1$ and that their signs are opposite.
Consider an arc $\xi$ of great circle from $e_3$ to $-e_3$.
Identify $\xi$ with $[a,b]$ via an orientation-preserving diffeomorphism, so that the map $F$ of Definition~\ref{defhombounformpseudo} can be regarded as the map that maps $\left(V \in \xi,u \in [-1,1]\right)$ to the rotation $\rho_{\alpha(u),V}$ with axis $V$ and angle $\alpha(u)$.
Then the intersection
\begin{equation*}d_Y(Z)=\Bigl\langle s_{\tau_b \circ F(\gamma,\tau_b)^{-1}}\bigl(N_c(\gamma);\tvarM\bigr),s_{\tau_b \circ \CT_{\gamma}^{-1}}\bigl(N_c(\gamma);\svarM\bigr) \Bigr\rangle \end{equation*} 
in \begin{equation*}\ST A\vert_{N_c(\gamma)} =_{\tau_b \circ F(\gamma,\tau_b)^{-1}} N_c(\gamma) \times S^2\end{equation*}
is the algebraic intersection 
\begin{equation*}\Bigl\langle N_c(\gamma)\times \{\tvarM\},F(\gamma,\tau_b) \circ \CT_{\gamma}^{-1}\bigl(N_c(\gamma) \times \{\svarM\}\bigr) \Bigr\rangle_{\!N_c(\gamma) \times S^2\,}.\end{equation*} 
It is the degree of the map 
\begin{equation*}\begin{array}{llll}f_{\svarM}\colon &\xi \times [-1,1]&\rightarrow &S^2\\&(V,u)&\mapsto &\rho_{\alpha(u),V}\circ \rho_{\alpha(-u),e_3} (\svarM)\end{array}\end{equation*} at $\tvarM$. The other algebraic intersection $d_Z(Y)$
is the degree of the map $f_{\tvarM}$ at $\svarM$.
The boundary of the image of $f_{\svarM}$ is $- 2 S^1(\svarM)$, where $S^1(\svarM)$ is oriented as the boundary of the closure of the connected component of $S^2 \setminus S^1(\svarM)$ that contains $e_3$. Therefore, the degree increases by $2$ from the component of $S^2 \setminus S^1(\svarM)$ that contains $e_3$ to the component of $(-e_3)$. Furthermore, the degree of $f_{\svarM}$ on the component of $e_3$ is independent of $\svarM\neq e_3$. Lemma~\ref{lemdegstwoso} implies $\deg(f_{- e_3})=\deg(\Phi_{-e_3})=-1$. So the degree of $f_{\svarM}$ at $e_3$ is $-1$ when $\svarM\neq e_3$. Therefore, the degree of $f_{\svarM}$ on the component of $-e_3$, which is independent of $\svarM \neq (-e_3)$, is $1$.
Thus, the degree of $f_{\svarM}$ at $\tvarM$ and the degree of $f_{\tvarM}$ at $\svarM$ are opposite.
\eop

\begin{lemma}
\label{lemdiffz}
Let $\Aman$ be a rational homology handlebody equipped with two pseudo-parallelizations $\tilde{\tau}_0$ and $\tilde{\tau}_1$ that coincide with the same genuine parallelization near the boundary of $\Aman$. Let $\cvarM \in S^2$. Let $\eta \in \left[0,\varepsilon\right[$. 
There exists a rational $4$-chain $H(\tilde{\tau}_0,\tilde{\tau}_1,\cvarM,\eta)$ in $U\Aman$ such that
\begin{equation*}\partial H(\tilde{\tau}_0,\tilde{\tau}_1,\cvarM,\eta)=s_{\tilde{\tau}_1}(\Aman;\cvarM,\eta) - s_{\tilde{\tau}_0}(\Aman;\cvarM,\eta).\end{equation*}
\end{lemma}
\bp Let $\svarM \in S^2 \setminus S^1(\cvarM)$.
Let us show that $C=s_{\tilde{\tau}_1}(\Aman;\cvarM,\eta) - s_{\tilde{\tau}_0}(\Aman;\cvarM,\eta)$ vanishes in 
\begin{equation*}H_3(U\Aman;\QQ)=H_1(\Aman;\QQ) \otimes H_2(S^2;\QQ).\end{equation*}
To do so, we prove that the algebraic intersections of $C$ with $s_{\tilde{\tau}_0}(S;\svarM)$ vanish for $2$-cycles $S$ of $(\Aman,\partial \Aman)$ that generate $H_2(\Aman,\partial \Aman)$. We assume $\eta=0$ without loss.
Of course $s_{\tilde{\tau}_0}(\partial \Aman;\cvarM)$ does not intersect $s_{\tilde{\tau}_0}(S;\svarM)$. According to Lemma~\ref{lemintpseudosec}, the chain $s_{\tilde{\tau}_0}(S;\cvarM)$ does not intersect $s_{\tilde{\tau}_0}(S;\svarM)$ in $U\Aman\vert_{S}$ algebraically. So $s_{\tilde{\tau}_0}(\Aman;\cvarM)$ does not intersect $s_{\tilde{\tau}_0}(S;\svarM)$ algebraically, either. Since Proposition~\ref{propsplusplustau} guarantees that $(s_{\tilde{\tau}_0}(S;\cvarM)-s_{\tilde{\tau}_1}(S;\cvarM))$ bounds in $U\Aman\vert_{S}$, the chain $s_{\tilde{\tau}_1}(S;\cvarM)$ does not intersect $s_{\tilde{\tau}_0}(S;\svarM)$ algebraically in $U\Aman\vert_{S}$. So $s_{\tilde{\tau}_1}(\Aman;\cvarM)$ does not intersect $s_{\tilde{\tau}_0}(S;\svarM)$ algebraically. This proves the existence of the desired $4$-chain $H(\tilde{\tau}_0,\tilde{\tau}_1,\cvarM,\eta)$.
\eop

\section{Definition of \texorpdfstring{$\Zinvuf$}{Z} with respect to pseudo-sections}
\label{secdefZpseudosec}

\begin{definition}
 \label{defpseudopropchain}
Let $\crats$ be a rational homology $\RR^3$ equipped with an asymptotically standard pseudo-parallelization $\tilde{\tau}$.
A \emph{propagating chain}\index[T]{propagating!chain} of $(C_2(\rats),\tilde{\tau})$ is a propagating chain of $C_2(\rats)$ (as in Definition~\ref{defpropagatortwo}) whose boundary intersects $\ST \crats$ as
a chain $s_{\tilde{\tau}}(\crats;\cvarM,\varepsilon_i)$ as in Definition~\ref{defpseudosec}.
\end{definition}

Again, for a given $s_{\tilde{\tau}}(\crats;\cvarM,\varepsilon_i)$, a propagating chain whose boundary intersects $\ST \crats$ as 
 $s_{\tilde{\tau}}(\crats;\cvarM,\varepsilon_i)$ exists because $H_3(C_2(\rats);\QQ)=0$.

\begin{lemma}
\label{lempropchaindualpseudopropsec}
Let $\crats$ be a rational homology $\RR^3$ equipped with an asymptotically standard pseudo-parallelization $\tilde{\tau}$. For any positive number $\alpha$, and for
any propagating chain $\propP$ of $(C_2(\rats),\tilde{\tau})$, transverse to $\partial C_2(\rats)$, there exists a propagating form
 of $(C_2(\rats),\tilde{\tau})$ (as in Definition~\ref{defpseudopropform}) that is $\alpha$-dual to $\propP$ (as in Definition~\ref{defformdual}).
\end{lemma}
\bp There is no loss of generality in assuming that $\propP$ intersects a collar $\left[-1,0\right] \times \partial  C_2(\rats)$ as $\left[-1,0\right] \times \partial \propP$.

Assume that the boundary of $\propP$ intersects $\ST \crats$ as
a chain $s_{\tilde{\tau}}(\crats;\cvarM,\varepsilon_i)$. Let $\omega_i$ be volume-one form of $S^2$, which is $\alpha_{\cvarM}$-dual to $\cvarM$ for some small $\alpha_{\cvarM} \in \left]0,\alpha\right[$.
Set \begin{equation*}I(\varepsilon_i)=\left[b-\varepsilon_i-\varepsilon_i^k,b-\varepsilon_i+\varepsilon_i^k\right]\end{equation*} and $I_2(\varepsilon_i)=[b-\varepsilon_i-2\varepsilon_i^k,b-\varepsilon_i+2\varepsilon_i^k]  \times [-1,1]  \times S^2$.

Forms $\alpha$-dual to a given chain can be constructed as in Lemma~\ref{lemconstrformalphadual}.
We use this process to construct a form $\omega_1$ $\alpha$-dual to $\propP$, which pulls back through $\partial  C_2(\rats)$ on $\left[-1,0\right] \times \partial  C_2(\rats)$, which may be expressed as in Definition~\ref{defpseudoformprel} outside $\left[-1,0\right] \times \ST \left(I(\varepsilon_i)  \times \gamma  \times [-1,1]\right)$, and which factors through  $I(\varepsilon_i)  \times [-1,1]  \times S^2$ on $\left[-1,0\right] \times \ST \bigl(I(\varepsilon_i)  \times \gamma  \times [-1,1]\bigr)$.

It remains to see that $\omega_1$ can be assumed to have the form prescribed by 
Definition~\ref{defpseudopropform}  
on $\left[-1,0\right] \times \ST \left(I(\varepsilon_i)  \times \gamma  \times [-1,1]\right)$.

We also have a form $\omega_2$ as in Definition~\ref{defpseudopropform}, which also pulls back through $\partial  C_2(\rats)$ on $\left[-1,0\right] \times \partial  C_2(\rats)$ and which coincides with $\omega_1$ on $\left[-1,0\right] \times \ST \left(\crats \setminus \left(I(\varepsilon_i)  \times \gamma  \times [-1+\varepsilon,1-\varepsilon]\right)\right)$.

Both forms factor through $I_2(\varepsilon_i)  \times [-1,1]  \times S^2$, and they are cohomologous there. So there exists a one-form $\eta$ on $I_2(\varepsilon_i)  \times [-1,1]  \times S^2$ such that
$\omega_2-\omega_1=p^{\ast}(d \eta)$, where $p$ is the natural projection that forgets the $\left[-1,0\right]$ factor and the $\gamma$ factor, and 
$d\eta$ vanishes outside $I(\varepsilon_i)  \times [-1+\varepsilon,1-\varepsilon]  \times S^2$.

To conclude the proof, it suffices to prove that such $\eta$ can be chosen so that $\eta$ is zero outside $I(\varepsilon_i)  \times [-1+\varepsilon,1-\varepsilon]  \times S^2$, too.
Therefore, it suffices to prove that the cohomology class of $\eta$ is zero on \begin{equation*}\Bigl(\bigl(I_2(\varepsilon_i)  \times [-1,1]\bigr) \setminus \bigl( I(\varepsilon_i)  \times [-1+\varepsilon,1-\varepsilon]\bigr) \Bigr)\times S^2.\end{equation*} Thus, it suffices to prove that the integral of $\eta$ along \begin{equation*} \partial \bigl(I_2(\varepsilon_i)  \times [-1,1] \times \{\svarM\}\bigr)\end{equation*} is zero for some $\svarM \in S^2$.
This integral is also the integral of $(\omega_2-\omega_1)$ along $I_2(\varepsilon_i)  \times [-1,1] \times \{\svarM\}$ or along $\left(N_c(\gamma)=[a,b]  \times \{c\}  \times [-1,1]\right) \times \{\svarM\}$, where the integral of $\omega_1$ is the algebraic intersection of $s_{\tau_b}([a,b]  \times \{c\} \times [-1,1];\svarM)$ with $\propP$.

Note that nothing depends on the trivialization $\tau_b$, which identifies 
\begin{equation*}\ST \left([a,b]  \times \gamma  \times [-1,1]\right)\end{equation*} with $\left([a,b]  \times \gamma  \times [-1,1]\right) \times S^2$. Therefore, there is no loss of generality in assuming that $\tau_b$ maps $e_3$ to the direction of $(-\gamma)$, which is the positive normal to the meridian disks, so that
$(N_c(\gamma) \times \{e_3\})=s_+(N_c(\gamma))$.

Assume $\cvarM \neq e_3$ and choose $\svarM=e_3$.
Then it suffices to prove $\int_{s_+(N_c(\gamma))}\omega_1=\int_{s_+(N_c(\gamma))}\omega_2$.
Lemma~\ref{lempropintsplusstau} implies
\begin{equation*}\int_{s_+(N_c(\gamma))}\omega_2=\frac12\chi\left({\tau}_e(. \times e_2)\vert_{\partial N_c(\gamma)};N_c(\gamma)\right).\end{equation*}
Let us compute
$\int_{s_+(N_c(\gamma))}\omega_1= \langle s_+(N_c(\gamma)) ,P \rangle$.
Thanks to Proposition~\ref{propsplusplustau},
\begin{equation*}2\Bigl(s_+\bigl(N_c(\gamma)\bigr)-s_{\tilde{\tau}}\bigl(N_c(\gamma); e_3\bigr)\Bigr) - \chi\left({\tau}_e(. \times e_2)\vert_{\partial N_c(\gamma)};N_c(\gamma)\right)\ST \crats\vert_{\ast} \end{equation*}
is null-homologous in $\ST \crats\vert_{\ast}$. So its algebraic intersection with $\propP$ is zero.
The algebraic intersection of $s_{\tilde{\tau}}(N_c(\gamma); e_3)$ with $\propP$ is also zero according to Lemma~\ref{lemintpseudosec}.
We thus get $\int_{s_+(N_c(\gamma))}\omega_1=\int_{s_+(N_c(\gamma))}\omega_2$, and the lemma is proved when $e_3 \neq \cvarM$. 

When $e_3 = \cvarM$, choose $\svarM=-e_3$ and conclude by computing the integrals along $s_-(N_c(\gamma))$, similarly.
\eop

For a family of propagating chains $\propP_i$ of $(C_2(\rats),\tilde{\tau})$ in general $3n$-position with respect to $\Link$ as in Definition~\ref{defgenthreenpos}, and for propagating forms $\omega(i)$, which are $\alpha$-dual to them for a sufficiently small $\alpha$ as in Lemma~\ref{lempropchaindualpseudopropsec}, 
the integrals involved in $\Zinv_n^s(\crats,\Link,\tilde{\tau})$ in Theorem~\ref{thmconststraightpseudo} can be computed as algebraic intersections of preimages of 
transverse propagating chains of $(C_2(\rats),\tilde{\tau})$ as in Sections~\ref{secfromintegtointer} and \ref{sectransrat}.

This provides the announced discrete definition of $\Zinv_n^s(\crats,\Link,\tilde{\tau})$ with respect to pseudo-parallelizations.

\bookmarksetup{startatroot}
\appendix



\chapter{Some basic algebraic topology}
\label{chaphomology}

This appendix reviews the main well-known results in algebraic topology used throughout the book.
The proofs of these results can be found in several books about basic algebraic topology, such as \cite{greenberg,HatcherAT,Spanier}, by Marvin Greenberg, Allen Hatcher, and Edwin Spanier.
Here we state only the weak versions used in the book, with some sketches of proofs and some hints to show how things work.

\section{Homology}
\label{sechomology}

We first review some properties of homology.

A \emph{topological pair} $(X,Y)$ consists of a topological space $X$ and a subspace $Y$ of $X$.
A \emph{map} $f\colon (X,Y) \to (A,B)$ between such pairs is a continuous map $f$ from $X$ to $A$ sending $Y$ to $B$.
A topological space $X$ is identified with the pair $(X,\emptyset)$.

The coefficient ring $\Lambda$ of the homology theories $H(.)=H(.;\Lambda)$ considered in this book is $\ZZ/2\ZZ$, $\ZZ$, $\QQ$, or $\RR$.

A \emph{covariant functor} $H$ from the category of topological pairs 
$(X,Y)$ to the category of graded $\Lambda$-modules and homomorphisms of degree $0$ maps any topological pair $(X,Y)$ 
to a sequence $H(X,Y;\Lambda)=H(X,Y)=(H_k(X,Y))_{k\in \ZZ}$ of $\Lambda$-modules. It associates to any map $f\colon (X,Y) \to (A,B)$ between pairs a sequence of $\Lambda$-linear morphisms $H(f)=(H_k(f) \colon H_k(X,Y)\to H_k(A,B))_{k\in \ZZ}$
so that if $g\colon (A,B) \to (V,W)$ is another map between pairs $H_k(g \circ f)=H_k(g) \circ H_k(f)$, and 
$H_k$ maps the identity map of $(X,Y)$ to the identity map of $H_k(X,Y)$.

The homology theories $(H,\partial)$ that we consider consist of 
\begin{itemize}
\item a \emph{covariant functor} $H$ from the category of topological pairs 
$(X,Y)$ to the category of graded $\Lambda$-modules and homomorphisms of degree $0$, and
\item 
sequences $\partial_k(X,Y)$ of linear maps $\partial_k(X,Y) \colon H_k(X,Y)\to H_{k-1}(Y)$ associated to topological pairs 
$(X,Y)$ such that for any map $f\colon (X,Y) \to (A,B)$, we have
\begin{equation*}\partial_k(A,B) \circ H_k(f) = H_{k-1}(f\vert_{Y}) \circ \partial_k(X,Y).\end{equation*}
\end{itemize}

The considered homology theories satisfy the following Eilenberg--Steenrod axioms \cite[Chap. 4, Sec. 8, p.199]{Spanier}. 

\begin{itemize}
\item {\textbf{Homotopy axiom:}} For two topological pairs $(X,Y)$ and $(A,B)$, if 
\begin{equation*}\begin{array}{llll}f \colon& \left[0,1\right] \times X &\to &A \\
&(t,x)& \mapsto & f_t(x)=f(t,x)\end{array}\end{equation*} is a continuous map sending
$\left[0,1\right] \times Y$ to $B$, then $H(f_0)$ is equal to $H(f_1)$.
\item {\textbf{Exactness axiom:}} For any topological pair $(X,Y)$ with inclusion maps $i\colon Y \hookrightarrow X$ and $j\colon (X,\emptyset) \hookrightarrow (X,Y)$,
there is a long exact sequence: 
\begin{equation*}\dots \hflv{\partial_{k+1}(X,Y)} H_k(Y) \hflv{H_k(i)} H_k(X) \hflv{H_k(j)} H_k(X,Y) \hflv{\partial_{k}(X,Y)} H_{k-1}(Y) \hflv{H_{k-1}(i)} \dots\footnote{Such a sequence of homomorphisms is \emph{exact}\index[T]{exact!sequence} if and only if the image of a homomorphism is equal to the kernel of the following one.}\end{equation*}
\item {\textbf{Excision axiom:}} For any topological pair $(X,Y)$, if $U$ is an open subset of $X$ whose closure lies in the interior of $Y$, then the inclusion map 
$e \colon (X \setminus U,Y \setminus U) \to (X,Y)$ induces isomorphisms 
$H_k(e) \colon H_k(X \setminus U,Y \setminus U) \to H_k(X,Y)$ for any $k \in \ZZ$.
\item {\textbf{Dimension axiom:}} If $X$ has one element $x$, then $H_k(X)\cong \{0\}$ for any $k \neq 0$. Furthermore, $H_0(X)$ is isomorphic to $\Lambda$ and is equipped with a canonical generator denoted by $[x]$. So we write $H_0\left(\{x\}\right)=\Lambda\left[x\right]$.
\end{itemize}

An example of such a homology theory is the \emph{singular homology} described in \cite[Chap. 10 and 13]{greenberg} and in \cite[Section 2.1]{HatcherAT}.

The spaces considered in this book are homotopy equivalent to finite CW-complexes. The above Eilenberg--Steenrod axioms characterize the homology of these spaces.
Let us show a few examples of properties of homology and computations from these axioms.
First note that the functoriality implies that $H$ maps a homeomorphism $f\colon (X,Y) \to (A,B)$ between topological pairs to an isomorphism $H(f)$.
Also note that the homotopy axiom implies that the homologies of $\RR^n$ and its unit ball $B^n$ are isomorphic to the homology of a point, which is determined
by the dimension axiom.

Here are other easy consequences of the Eilenberg--Steenrod axioms.

\begin{proposition}
\label{proprdis}
Let $X$ and $Y$ be two topological spaces. Let $i_X$ and $i_Y$ denote their respective inclusion maps into their disjoint union
$X \sqcup Y$. Then for any $k \in \ZZ$,
$H_k(i_X) \colon H_k(X) \to H_k(X \sqcup Y)$ and $H_k(i_Y) \colon H_k(Y) \to H_k(X \sqcup Y)$ are injective,
and we have \begin{equation*}H_k(X \sqcup Y)=H_k(i_X)\bigl(H_k(X)\bigr) \oplus H_k(i_Y)\bigl(H_k(Y)\bigr).\end{equation*}
\end{proposition}
\bp
The excision isomorphism $H_k(e_X) \colon H_{k}(X)=H_{k}(X ,\emptyset) \to H_{k}(X \sqcup Y,Y)$ may be expressed as
\begin{equation*}H_k(e_X)= H_k\bigl(j\colon (X \sqcup Y,\emptyset) \hookrightarrow (X \sqcup Y,Y)\bigr) \circ H_k(i_X \colon X \hookrightarrow X \sqcup Y).\end{equation*}
So $H_k(j)$ is surjective and $H_k(i_X)$
is injective. Similarly, $H_k(i_Y)$ is injective. Thus, the long exact sequence of
$(X \sqcup Y,Y)$ yields short exact sequences
\begin{equation*}0 \rightarrow H_k(Y) \hflv{H_k(i_Y)}
H_k(X \sqcup Y) 
\hflv{H_k(e_X)^{-1} \circ H_k(j)} H_k(X)
\rightarrow 0\end{equation*}
for any integer $k$.
Since $\left(H_k(e_X)^{-1} \circ H_k(j)\right) \circ H_k(i_X)$ is the identity map, the exact sequence splits, and we get the result. 
\eop

\begin{proposition}
\label{proplonghomseqtriple}
Let $X$ be a topological space. Let $Y$ and $Z$ be two subspaces of $X$ such that $Z\subseteq Y$. Let
\begin{equation*}\bigl(\partial_{k}=\partial_{k}(X,Y,Z)\bigr) \colon H_k(X,Y) \to H_{k-1}(Y,Z)\end{equation*}
be the composition of $\partial_{k}(X,Y) \colon H_k(X,Y) \to H_{k-1}(Y)$ and the map from $H_{k-1}(Y)$ to $H_{k-1}(Y,Z)$ induced by the inclusion.
Then the long sequence 
\begin{equation*}\xymatrix{H_{k+1}(X,Y) \ar[r]^{\partial_{k+1}} & H_k(Y,Z) \ar[r]^{H_k(i)} & H_k(X,Z) \ar[r]^{H_k(j)} & H_k(X,Y) \ar[lld]_{\partial_{k}} &\\
& H_{k-1}(Y,Z) \ar[r]_{H_{k-1}(i)} & H_{k-1}(X,Z) \ar[r]^{H_{k-1}(j)} & H_{k-1}(X,Y) \ar[r]^{\partial_{k-1}}&  H_{k-2}(),}
\end{equation*}
where the $H_k(i)$ and the $H_k(j)$ are induced by inclusion, 
is exact. It is called the \emph{long exact sequence of homology of the triple} $(X,Y,Z)$.
\end{proposition}
\bsp To prove that $H_k(j) \circ H_k(i)=0$, note that the inclusion map from $(Y,Z)$ to $(X,Y)$ factors through $(Y,Y)$. Next
chase in commutative diagrams using functoriality and the exactness axiom, which is this exact sequence when $Z=\emptyset$.
\eopwobp

Let $X$ be a topological space with $k$ connected components. Choose a basepoint $x_i$ in each connected component $X_i$. Let $Y$ denote the set of these basepoints equipped with the natural discrete topology. Let $i\colon Y \hookrightarrow X$ be the inclusion map, and let $p$ be the map sending each connected component to its basepoint. By functoriality $H(p)\circ H(i)$ is the identity map. So $H(i)$ is injective. In particular $H_0(i)$ injects $\oplus_{x_i \in Y}\Lambda[x_i]$ into $H_0(X)$.  The element $H_0(i)([x_i])$ is also denoted by $[x_i]$. The homotopy axiom implies that the element $[x_i]$ of $H_0(X_i)$ is independent of the basepoint $x_i$ of $X_i$ if $X_i$ is path-connected.

\begin{lemma}
\label{lemrrmin} We have 
\begin{equation*}H_k\bigl(\RR, \RR \setminus \{0\}\bigr)=\left\{\begin{array}{ll}
\Lambda \partial_1^{-1}\bigl(\RR,\RR \setminus \{0\}\bigr)\bigl([1]-[-1]\bigr)\cong \Lambda& \textrm{if}\;k=1\\
\{0\} &  \textrm{otherwise.}
\end{array} \right. \end{equation*}
\end{lemma}
\bp
Proposition~\ref{proprdis} and the observations before imply
\begin{equation*}H_0\bigl(\RR \setminus \{0\}\bigr)= \Lambda[-1] \oplus \Lambda[+1].\end{equation*}
The morphism $H_0(i)$ induced by inclusion from $H_0(\RR \setminus \{0\})$ to $H_0(\RR) =\Lambda[1]$
maps $[-1]$ and $[+1]$ to the preferred generator $[-1]=[0]=[1]$ of $H_0(\RR)$.
The long exact sequence associated to $(\RR,\RR \setminus \{0\},\emptyset)$ implies that $\partial_1(\RR,\RR \setminus \{0\})\colon H_1(\RR,\RR \setminus \{0\}) \longrightarrow H_0(\RR \setminus \{0\})$
is an isomorphism onto its image, which is the kernel $\Lambda\left([+1]-[-1]\right)$ of $H_0(i)$. 
The long exact sequence also implies $H_q(\RR, \RR \setminus \{0\})=\{0\} $ when $q \neq 1$.
\eop

The reader can also apply the Eilenberg--Steenrod axioms and the above observations to prove 
the following proposition.

\begin{proposition}
\label{prophnrnsetm}
Let $n \in \NN \setminus \{0\}$ and let $k\in \ZZ$.
Recall \begin{equation*}\ballb^n=\Bigl\{x=(x_1,\dots,x_n) \in \RR^n \suchthat  \lVert x \rVert^2 =\sum_{i=1}^n x_i^2 \leq 1\Bigr\}\end{equation*} and $S^{n-1}=\partial \ballb^n$.  

Assume that $n\geq 2$. Set
$\ballb_+^{n-1}=\{x \in S^{n-1}  \suchthat  x_1 \geq 0\}$, $\ballb_-^{n-1}=\{x \in S^{n-1}  \suchthat  x_1 \leq 0\}$, and $S^{n-2}_+=\{x \in S^{n-1}  \suchthat  x_1 = 0\}$.
The morphisms of the sequence
\begin{multline*}H_k(\RR^n,\RR^n \setminus \{0\})  \stackrel{H_k(i)}{\longleftarrow} H_k(\ballb^n,S^{n-1})
\stackrel{\partial}{\longrightarrow} H_{k-1}(S^{n-1},\ballb_-^{n-1}) \\
\longleftarrow H_{k-1}(\ballb_+^{n-1},S^{n-2}_+) \stackrel{H_{k-1}(p)}{\longrightarrow} H_{k-1}(\RR^{n-1},\RR^{n-1} \setminus \{0\}),\end{multline*}
where \begin{itemize}
\item the unlabeled morphism and $H_k(i)$ are induced by inclusions, 
\item $\partial$ is the morphism $\partial_k(\ballb^n,S^{n-1},\ballb_-^{n-1})$ of Proposition~\ref{proplonghomseqtriple}, and
\item $p:\ballb_+^{n-1} \longrightarrow \RR^{n-1}$ forgets the first coordinate and shifts the numbering of the remaining ones by $(-1)$,\end{itemize}
are isomorphisms.
In particular, when $k=n$, their composition   
\begin{equation*}D_n \colon H_n(\RR^n,\RR^n \setminus \{0\}) \to H_{n-1}(\RR^{n-1},\RR^{n-1} \setminus \{0\})\end{equation*}
is an isomorphism. 
Let $[\RR,\RR \setminus \{0\}]=\partial_1^{-1}(\RR,\RR \setminus \{0\})([1]-[-1])$ be the preferred generator of 
$H_1(\RR, \RR \setminus \{0\})$ of Lemma~\ref{lemrrmin}. Inductively define \begin{equation*}\bigl[\RR^n,\RR^n \setminus \{0\}\bigr]=D_n^{-1}\bigl[\RR^{n-1},\RR^{n-1} \setminus \{0\}\bigr].\end{equation*} Then we have
\begin{equation*}H_k\bigl(\RR^n,\RR^n \setminus \{0\}\bigr)=\left\{\begin{array}{ll}
\Lambda \bigl[\RR^n,\RR^n \setminus \{0\}\bigr] \cong \Lambda& \textrm{if}\;k=n\\
\{0\} &  \textrm{otherwise.}
\end{array} \right.\end{equation*} for any positive integer $n$.
For $n \geq 1$, set \begin{equation*}[\ballb^n,S^{n-1}]=H_n(i)^{-1}\left([\RR^n,\RR^n \setminus \{0\}]\right).\end{equation*} Then we have
\begin{equation*}H_k(\ballb^n,S^{n-1})=\left\{\begin{array}{ll}
\Lambda [\ballb^n,S^{n-1}] \cong \Lambda& \textrm{if}\;k=n\\
\{0\} &  \textrm{otherwise.}
\end{array} \right. \end{equation*}
\end{proposition}
\eopwobp

As a corollary, we get the following proposition:

\begin{proposition}
\label{prophomologyspheres}
 Let $n \in \NN$. Assume that $n\geq 1$. Set
 \begin{equation*}[S^n]=\partial_{n+1}(\ballb^{n+1},S^{n})\left([\ballb^{n+1},S^{n}]\right).\end{equation*} We have
 \begin{equation*}H_k(S^n)=\left\{\begin{array}{ll}
\Lambda [S^n] \cong \Lambda& \textrm{if}\;k=n\\
\Lambda [(1,0,\dots,0)] \cong \Lambda& \textrm{if}\;k=0\\
\{0\} &  \textrm{otherwise.}
\end{array} \right. \end{equation*}
\end{proposition}
\bp Use the exact sequence associated to $(\ballb^{n+1},S^{n})$ and the previous proposition.
\eop

\begin{proposition}
\label{prophomologyori}
Let $n \in \NN$. Let $\phi$ be a diffeomorphism from $\RR^n$ to $\RR^n$ sending $0$ to $\phi(0)=0$.
If $\phi$ preserves the orientation, then $\phi$ induces the identity map on $H_n(\RR^n,\RR^n \setminus \{0\})$.
Otherwise, we have \begin{equation*}H_n(\phi)\Bigl(\bigl[\RR^n,\RR^n \setminus \{0\}\bigr]\Bigr)=-\bigl[\RR^n,\RR^n \setminus \{0\}\bigr].\end{equation*}
\end{proposition}
\bp When $n=1$, this is easy to see
with the generator \begin{equation*}\bigl[\RR,\RR \setminus \{0\}\bigr]=\partial_1\bigl(\RR,\RR \setminus \{0\}\bigr)^{-1}\bigl([+1]-[-1]\bigr).\end{equation*}
Let $\iota_n \colon \RR^n \to \RR^n$ map $(x_1, x_2, \dots , x_{n-1}, x_n)$ to $(x_1, x_2, \dots , x_{n-1}, -x_n)$. With the notation of Proposition~\ref{prophnrnsetm}, we have
$H_{n-1}(\iota_{n-1}) \circ D_n=D_n \circ H_{n}(\iota_{n})$. We inductively deduce $H_{n}(\iota_{n})([\RR^n,\RR^n \setminus \{0\}])=-[\RR^n,\RR^n \setminus \{0\}]$. Any orientation-reversing linear isomorphism of $\RR^n$ is 
homotopic to $\iota_n$ through linear isomorphisms. Therefore, we have $H_{n}(\phi)([\RR^n,\RR^n \setminus \{0\}])=-[\RR^n,\RR^n \setminus \{0\}]$ for any such isomorphism $\phi$. We similarly observe that $H_{n}(\phi)$ sends $[\RR^n,\RR^n \setminus \{0\}]$ to itself for an orientation-preserving linear isomorphism $\phi$. Diffeomorphisms preserving $0$ are homotopic to their linear derivative at $0$
through maps preserving $\RR^n \setminus \{0\}$, near $0$. So the result follows for general diffeomorphisms thanks to the excision axiom.
\eop

Let $X$ be a topological space equipped with a triangulation $T$ as in Subsection~\ref{submorelowdif}. The Eilenberg--Steenrod axioms also imply that the homology $H(X)$ of $X$ can be computed as follows.\footnote{The reader can prove it, but it takes more space.}
Equip each simplex of the triangulation with an arbitrary orientation, and let $C_k(T)=C_k(T;\Lambda)$ denote the $\Lambda$-module freely generated by the simplices of dimension $k$ of $T$.
The \indexT{algebraic boundary} of such a $k$-dimensional simplex $\Delta$ 
is the sum \begin{equation*}\partial_k (\Delta)=\sum_{\delta}\varepsilon(\Delta,\delta)\delta,\end{equation*} running over all the $(k-1)$-dimensional simplices $\delta$ in the boundary of $\Delta$, where
$\varepsilon(\Delta,\delta)$ is $1$ if $\delta$ is oriented as (part of) the boundary of $\Delta$ (with the outward normal first convention as usual), and $(-1)$ otherwise.
Define the \indexT{boundary map} $\partial_k \colon C_k(T) \to C_{k-1}(T)$ to be the linear map sending a $k$-dimensional simplex $\Delta$ to its algebraic boundary $\partial_k (\Delta)$.
Then we have $\partial_k \circ \partial_{k+1}=0$. So $(C(T),\partial)=(C_k(T),\partial_k)_{k\in \ZZ}$ is a \emph{chain complex}. Its homology
\begin{equation*}H\bigl(C(T),\partial\bigr)=\left(H_k(C(T),\partial)=\frac{\Ker \partial_k}{\Image \partial_{k+1}}\right)_{k \in \ZZ}\end{equation*}
is canonically isomorphic to the homology of $X$.

The elements of $\Ker \partial_k$ are the \emph{$k$-dimensional cycles}\index[T]{cycle} of $C(T)$. The elements of $\Image \partial_{k+1}$ are the \emph{$k$-dimensional boundaries} of $C(T)$.
The elements of $C_k(T)=C_k(T;\Lambda)$ are called the \emph{simplicial chains} of $T$ of dimension $k$ with coefficients in $\Lambda$

Thus, the homology of an $n$-dimensional manifold $M$ that can be equipped with a triangulation $T$ vanishes in degrees higher than $n$ and in negative degrees. (The existence of a homology theory satisfying the axioms and Proposition~\ref{prophnrnsetm} imply that the notion of dimension is well-defined for topological manifolds.)
Let $C$ be an $n$-dimensional cycle in such a connected manifold $M$, and let $\Delta_1$ and $\Delta_2$ be two $n$-simplices of $T$ intersecting along an $(n-1)$-simplex.
Then the coefficients of $\Delta_1$ and $\Delta_2$ in $C$ must coincide if the orientations of $\Delta_1$ and $\Delta_2$ are consistent along $\Delta_1 \cap \Delta_2$. They must be opposite to each other if the orientations are not consistent.

In particular, when $\Lambda={\ZZ}/{2\ZZ}$, if $M$ is connected, then the existence of a nonzero $n$-dimensional cycle implies that the boundary of $M$ is empty and that $M$ is compact.
When $\Lambda=\ZZ$, $\QQ$ or $\RR$, if $M$ is connected, then the existence of a nonzero $n$-dimensional cycle furthermore implies that $M$ is orientable. 

Assume that $M$ is a compact, $n$-dimensional, connected, triangulable, oriented manifold with empty boundary, and assume that the $n$-simplices of $T$ are equipped with the induced orientation. Then the sum of all these simplices is a cycle. Its homology class is called the \emph{fundamental class} of $M$. It is denoted by $[M]$. (When $\Lambda={\ZZ}/{2\ZZ}$, this definition does not require an orientation of $M$.) Let $\ballb^n$ be an n-dimensional closed ball embedded in $M$ by an orientation-preserving embedding. Let $0$ be the center of this ball.
Thanks to the excision axion, the inclusions induce isomorphisms from $H_n(\ballb^n,\ballb^n \setminus \{0\})$ to $H_n(\RR^n,\RR^n \setminus \{0\})$, and from $H_n(\ballb^n,\ballb^n \setminus \{0\})$ to $H_n(M,M \setminus \{0\})$. Let $[\ballb^n,\ballb^n \setminus \{0\}]$ denote the generator of $H_n(\ballb^n,\ballb^n \setminus \{0\})$ mapped to $[\RR^n,\RR^n \setminus \{0\}]$ by the first isomorphism and let $[M,M \setminus \{0\}]$ denote the image of $[\ballb^n,\ballb^n \setminus \{0\}]$ under the second isomorphism. 
Then the inclusion induces an isomorphism from $H_n(M)$ to $H_n(M,M \setminus \{0\})$, which maps the generator $[M]$ of $H_n(M)$ to the generator $[M,M \setminus \{0\}]$ of $H_n(M,M \setminus \{0\})$.
This defines the generator $[M]$ of $H_n(M)$ independently of a triangulation.

Let $X$ be a topological space equipped with a triangulation $T$ as above. Let $Y$ be a closed subspace of $X$ that is a union of simplices of $T$. Let $T_Y$ be the corresponding triangulation of $Y$. Set $C_k(T,T_Y)={C_k(T)}/{C_k(T_Y)}$ and define $\partial_k(T,T_Y) \colon C_k(T,T_Y) \to C_{k-1}(T,T_Y)$ to be the map induced by the previous boundary map $\partial_k$.
Then the homology $H(X,Y)$ is canonically isomorphic to the homology of the chain complex $(C(T,T_Y),\partial(T,T_Y))$.

When $M$ is a connected, compact, oriented $n$-dimensional manifold with boundary, equipped with a triangulation $T$ whose $n$-dimensional simplices are oriented by the orientation of $M$, the sum of all $n$-dimensional simplices of $T$ is a cycle of $C_n(M,\partial M)$. Its homology class is called the \emph{fundamental class} of $(M,\partial M)$. It is denoted by $\left[M,\partial M\right]$. Again, if $\ballb^n$ is an n-dimensional closed ball embedded in $M$ by an orientation-preserving embedding, then the inclusions induce isomorphisms
from $H_n(\ballb^n,\ballb^n \setminus \{0\})$ to $H_n(M,M \setminus \{0\})$ and from $H_n(M,\partial M)$ to $H_n(M,M \setminus \{0\})$, and the image of $[\ballb^n,\ballb^n \setminus \{0\}]$ under the first isomorphism coincides with the image of $\left[M,\partial M\right]$ under the second one.
Is is denoted by $[M,M \setminus \{0\}]$.

These considerations allow us to talk about \emph{the homology class of a compact oriented $p$-dimensional submanifold $P$ of a manifold $M$}. It is the image of $\left([P] \in H_p(P)\right)$ in $H_p(M)$ under the map induced by the inclusion, and it is often still denoted by $[P]$.
When $P$ is a compact oriented $p$-dimensional manifold with boundary embedded in a topological space $X$ so that $\partial P$ is embedded in a subspace $Y$ of $X$, we define the class $[P,\partial P]$ of $(P,\partial P)$ in $H_p(X,Y)$, similarly.

With these conventions, when $M$ is a connected, compact, oriented, $n$-dimensional manifold with boundary, such that $M$ can be equipped with a triangulation, the boundary map $\partial_n$ in the homology sequence of the pair $(M,\partial M)$ maps $\left[M,\partial M\right]$
to the class $[\partial M]$, where $[\partial M]$ is the sum of the classes of the connected components of $\partial M$ equipped with the orientation induced by the orientation of $M$ with respect to the outward normal first convention.

All the manifolds considered in this book can be equipped with triangulations.
If $M$ is a manifold equipped with a triangulation $T$, and if $P$ is a $p$-dimensional closed oriented manifold embedded in $M$ that is a union of simplices of $T$, then the homology class of $P$ in $M$ vanishes if and only if the cycle that is the sum of the simplices of dimension $p$ of $P$ (equipped with the orientation of $P$) is the (algebraic) boundary of a simplicial chain of $T$ of dimension $p+1$.

Homologies with various coefficients are related by the universal coefficient theorem. See \cite[29.12]{greenberg}, \cite[Thm. 3A.3]{HatcherAT}, or \cite[Chap. 5, Sec.2]{Spanier}, for example. 
Here is an excerpt of this theorem.

\begin{theorem}
\label{thmcoefunivhom}
When $\Lambda=\QQ$ or $\RR$, we have
\begin{equation*}H_k(X;\Lambda)=H_k(X;\ZZ)\otimes_{\ZZ}\Lambda\end{equation*} for any topological space $X$.
\end{theorem}

In this book, we mostly use cohomology with coefficients in $\QQ$, ${\ZZ}/{2\ZZ}$ or $\RR$. In these cases, it can be defined by the following excerpt of the universal coefficient theorem for cohomology, which can be found in \cite[Chap.5, Sec.5, Thm. 3, page 243]{Spanier} and in \cite[Thm. 3.2]{HatcherAT}, for example.

\begin{theorem}
\label{thmcoefunivcohom}
 When $\Lambda$ is a field, we have
 \begin{equation*}H^k(X,Y;\Lambda)=\Hom_{\Lambda}\bigl(H_k(X,Y;\Lambda),\Lambda\bigr)\end{equation*} for any $k \in \ZZ$ and for any topological space $X$.
\end{theorem}

Note the sign $=$ in the above theorems, meaning that the identifications are canonical. For a continuous map $f\colon (X,Y) \to (A,B)$  and an integer $k \in \ZZ$, the morphism
\begin{equation*}H^k(f;\Lambda) \colon H^k(A,B;\Lambda) \to H^k(X,Y;\Lambda)\end{equation*}
maps a linear map $g$ of $H^k(A,B;\Lambda)$ to $g \circ H_k(f;\Lambda)$.

For a general $\Lambda$, and for a pair $(X,Y)$ of topological spaces such that $X$ is equipped with a triangulation $T$ as above, and $Y$ is equipped with a subtriangulation $T_Y$ of $T$ as before, the cohomology $H^{\ast}(X,Y;\Lambda)$ of $(X,Y)$ is the cohomology of the complex $\left(C^{\ast}\left(T,T_Y;\Lambda\right),\partial^{\ast}\left(T,T_Y;\Lambda\right)\right)$, where $C^{k}(T,T_Y;\Lambda)$ is equal to $\Hom(C_{k}(T,T_Y;\Lambda);\Lambda)$ and $\partial^k \colon C^{k}(T,T_Y;\Lambda) \to C^{k+1}(T,T_Y;\Lambda)$ maps a linear form $g$ to $g \circ \partial_{k+1}$. We have
\begin{equation*}H^{k}(X,Y;\Lambda)=\frac{\Ker(\partial^k)}{\Image(\partial^{k-1})}.\end{equation*}
The elements of $\Ker(\partial^k)$ are the \emph{$k$-dimensional cocycles}\index[T]{cocycle}, and the elements of $\Image(\partial^{k-1})$ are the \emph{$k$-dimensional coboundaries}\index[T]{coboundary}.

Here is a weak version of the Poincar\'e duality theorem.  See \cite[Chap. 26 to 28, in particular (28.18)]{greenberg} or \cite[Thm. 3.43]{HatcherAT}.

\begin{theorem}
\label{thmpoincdual}
Let $M$ be a compact, $n$-dimensional manifold with a possible boundary. Let $H$ denote the singular homology. Then there are canonical Poincar\'e duality isomorphisms
 from $H^k(M,\partial M;{\ZZ}/{2\ZZ})$ to $H_{n-k}(M;{\ZZ}/{2\ZZ})$, and from $H^k(M;{\ZZ}/{2\ZZ})$ to $H_{n-k}(M,\partial M;{\ZZ}/{2\ZZ})$, for any $k\in \ZZ$.

If $M$ is oriented, then for any $\Lambda \in {{\ZZ}/{2\ZZ}, \ZZ, \QQ, \RR}$, there are canonical Poincar\'e duality isomorphisms
from $H^k(M,\partial M;\Lambda)$ to $H_{n-k}(M;\Lambda)$, and from $H^k(M;\Lambda)$ to $H_{n-k}(M,\partial M;\Lambda)$. 

Let $A$ be an oriented smooth submanifold of $M$ of dimension $n-k$ such that $\partial A \subset \partial M$. Let $P^{-1}([A,\partial A]) \in H^k(M;\Lambda)$ denote the image of its class under the inverse of such a  Poincar\'e duality isomorphism. Let $B$ be the class of a closed oriented $k$-dimensional submanifold or a simplicial $k$-dimensional cycle transverse to $A$. Then
the evaluation $P^{-1}([A,\partial A])(B)$ of $P^{-1}([A,\partial A])$ at $B$ is the algebraic intersection $\langle B,A \rangle_{\!M\,}$.
\end{theorem}

Here is a weak version of the K\"unneth theorem. See \cite[Chap.5, Sec.3, Thm. 10 (and Chap.5, Sec.2, Lem. 5)]{Spanier} or \cite[Section 3.B]{HatcherAT}. 
\begin{theorem}
 \label{thmKun}
Let $H$ denote the singular homology with coefficients in a commutative field $\Lambda$.
Then for any two topological spaces $X$ and $Y$, for any $k \in \NN$, we have
\begin{equation*}H_k(X\times Y)=\oplus_{i=0}^kH_i(X) \otimes_{\Lambda}H_{k-i}(Y).\end{equation*}
\end{theorem}
Again, the sign $=$ means that the identification is canonical. For embeddings of oriented closed manifolds $P$ into $X$ and $Q$ into $Y$, the tensor product $[P] \otimes [Q]$ of the homology classes of their images, represents the homology class of $P \times Q$.

We end this section by stating a weak version of the following Mayer--Vietoris exact sequence, which can be recovered from the Eilenberg--Steenrod axioms. See \cite[Chap. 17, Thm. 17.6]{greenberg}, for example.

\begin{theorem}
\label{thmMV} Let $X$ be a topological space. Let $A$ and $B$ be subspaces of $X$ whose union $A\cup B$ equals $X$. Let $i_A \colon A\cap B \hookrightarrow A$, $i_B \colon A\cap B \hookrightarrow B$, $j_A \colon A \hookrightarrow X$, $j_B \colon B \hookrightarrow X$ denote the inclusion maps.
 Assume that $j_A$ and $j_B$ induce isomorphisms from $H(A,A\cap B)$ to $H(X,B)$ and from $H(B,A\cap B)$ to $H(X,A)$.
Let $\partial_{MV,k+1} \colon H_{k+1}(X) \to H_k(A\cap B)$ be the composition of the map from $H_{k+1}(X)$ to $H_{k+1}(X,B)$ induced by inclusion, the inverse of the isomorphism from $H_{k+1}(A,A\cap B)$
to $H_{k+1}(X,B)$, and the boundary map $\partial_{k+1}(A,A\cap B)\colon H_{k+1}(A,A\cap B) \to H_k(A\cap B)$ of the long exact sequence of $(A,A\cap B)$.
Then there is a long exact sequence
\begin{equation*}\dots \to H_{k+1}(X) \hflv{\partial_{MV,k+1}}  H_k(A\cap B) \hflv{i_{MV,k}} H_k(A) \oplus H_k(B) \hflv{j_{MV,k}}  H_k(X) \hflv{\partial_{MV,k}} \dots\end{equation*}
such that $j_{MV,k}(\alpha \in H_k(A) ,\beta \in H_k(B))=H_k(j_A)(\alpha)+H_k(j_B)(\beta)$ and
$i_{MV,k}(\gamma)=(H_k(i_A)(\gamma),-H_k(i_B)(\gamma))$.
\end{theorem}

\section{Homotopy groups}
\label{sechomotopygroups}

Let $X$ be a topological space equipped with a basepoint $x$. Let $n$ be a positive integer.
The set of homotopy classes of maps from $\left[0,1\right]^n$ to $X$ mapping $\partial (\left[0,1\right]^n)$ to $x$ is denoted by $\pi_n(X,x)$. It is a \indexT{homotopy group}.
Its product sends a pair $\left([f],\left[g\right]\right)$ of homotopy classes of maps $f$ and $g$ such that $f$ sends $\left[1/2,1\right]\times \left[0,1\right]^{n-1}$ to $x$ and $g$ sends $\left[0,1/2\right]\times \left[0,1\right]^{n-1}$ to $x$ to
the class $[f]\left[g\right]$ of the map that coincides with $f$ on $\left[0,1/2\right]\times \left[0,1\right]^{n-1}$ and with $g$ on $\left[1/2,1\right]\times \left[0,1\right]^{n-1}$. This product is commutative when $n\geq 2$. The set of path-connected components of $X$
is denoted by $\pi_0(X)$.

\begin{remark}
 Classically, the set $\pi_n(X,x)$ is defined as the set of homotopy classes of maps from $(S^n,(1,0,\dots,0))$ to $(X,x)$. The two definitions coincide since $S^n$ is homeomorphic to the quotient ${\left[0,1\right]^{n}}/(\partial \left[0,1\right]^{n})$ for $n\geq 1$.
\end{remark}

\begin{examples}
 Let $k$ and $n$ be positive integers, such that $1\leq k \leq n$. A standard approximation theorem \cite[Chapter 2, Theorem 2.6, p. 49]{hirsch} implies that any continuous map from $S^k$ to $S^n$ is homotopic to a smooth map.
 The Morse--Sard theorem~\ref{thmMorseSard} ensures that if $k<n$, any smooth map is valued in the complement of a point in $S^n$, which is contractible. Therefore, $\pi_k(S^n,\ast)=\{1\}$ if $1\leq k <n$.
 The reader can develop these arguments to prove that there is a canonical isomorphism from $\pi_n(S^n,\ast)$ to $\ZZ$, which maps the homotopy class of a smooth map from $(S^n,\ast)$ to itself to its degree (see Definition~\ref{defdegr}).
 
A weak version of the \emph{Hurewicz theorem} 
relating homotopy groups to homology groups, 
ensures that for any path-connected topological space $X$ equipped with a basepoint $x$, $H_1(X;\ZZ)$ is the abelianization of $\pi_1(X,x)$. See \cite[(12.1)]{greenberg}, for example.
\end{examples}

A map $p \colon E \to B$ is called a \emph{weak fibration} if it satisfies the following \emph{homotopy lifting property with respect to cubes}:

For any integer $n \in \NN$, for any pair $(h_0 \colon \left[0,1\right]^n \times \{0\} \to E,H\colon \left[0,1\right]^{n+1} \to B)$ of continuous maps such that $H\vert_{\left[0,1\right]^n \times \{0\}}=p\circ h_0$, there
exists a continuous extension $h$ of $h_0$ to $\left[0,1\right]^{n+1}$ such that $H=p\circ h$.

To such a weak fibration, we associate the following long exact sequence in homotopy \cite[Chap. 7, Sec. 2, Thm. 10]{Spanier}.

\begin{theorem}
\label{thmlongseqhomotopy}
Let  $p \colon E \to B$ be a weak fibration. Let $e \in E$ be a basepoint of $E$. Let $b=p(e)$ denote its image under $p$, and let $F=p^{-1}(b)$ denote the fiber over $b$. Then we have the long exact sequence
\begin{multline*}\dots \pi_{n+1}(B,b) \to \pi_n(F,e) \to \pi_n(E,e) \to \pi_n(B,b) \to \pi_{n-1}(F,e) \dots\\
\dots  \to \pi_1(B,b) \to \pi_0(F) \to \pi_0(E) \to \pi_0(B),\end{multline*}
where the maps between the $\pi_n$ are respectively induced by the inclusion $F \hookrightarrow E$ and by $p$,
and the map from $\pi_n(B,b) \to \pi_{n-1}(F,e)$ is constructed as follows.
An element of $\pi_n(B,b)$ is represented by a map $H \colon \left[0,1\right]^n \to B$, which has a lift $h \colon \bigl[0,1\bigr]^n \to E$ that maps $[0,1]^{n-1} \times \{0\} \cup \bigl(\partial [0,1]^{n-1} \times [0,1]\bigr)$ to $e$.
Such an element is mapped to the homotopy class of the restriction of $h$ to $\left[0,1\right]^{n-1} \times \{1\}$.
The last three maps of the exact sequence are just maps between sets. Exactness means that the preimage of the component of the basepoint is the image of the previous map.
\end{theorem}

\begin{remark}
To define the map from $\pi_n(B,b)$ to $\pi_{n-1}(F,e)$ in the above theorem, we implicitly used the fact that the pair \begin{equation*}\Bigl(\left[0,1\right]^n,\left[0,1\right]^{n-1} \times \{0\} \cup \bigl(\partial \left[0,1\right]^{n-1} \times \left[0,1\right]\bigr)\Bigr)\end{equation*} is homeomorphic to $\bigl(\left[0,1\right]^n,\left[0,1\right]^{n-1} \times \{0\}\bigr)$ for all $n \in \NN \setminus \{0\}$.
\end{remark}

A path from a point $x$ to another point $y$ of $X$ induces an isomorphism from $\pi_n(X,x)$ to $\pi_n(X,y)$ for any integer $n$. So the basepoint is frequently omitted from the notation $\pi_n(X,x)$ when $X$ is path-connected.

A map $p \colon E \to B$ is called a \indexT{covering map} if every $b \in B$ has a neighborhood $U$ such that $p^{-1}(U)$
is a disjoint union of subsets of $E$, each of which is mapped homeomorphically onto $U$ by $p$. Such a covering map is an example of a weak fibration, for which $p$ induces isomorphisms from $\pi_n(E,e)$ to $\pi_n(B,p(e))$ for any $e \in E$ and any $n\geq 2$.

\chapter{Differential forms and de Rham cohomology}
\label{chapDeRham}

Here are a few well-known results about differential forms used throughout the book.
More detailed accounts can be found in the books \cite[Chapters XI and XII]{Godbillon} by Claude Godbillon and \cite{BT} by Raoul Bott and Loring Tu.

\section{Differential forms}
Let $M$ be a smooth manifold with possible boundary and ridges.
Degree $0$ forms of $M$ are smooth functions on $M$.
The differential $df=Tf \colon TM \rightarrow \RR$ of a smooth map $f$ from $M$ to $\RR$ is an example of a degree-one form of $M$.
In general, a \emph{degree $k$ differential form} on $M$ is a smooth section of the bundle $\bigwedge^k (TM)^{\ast}=\Hom(\bigwedge^k TM ;\RR)$ over $M$. So, such a form $\omega$ maps any $m \in M$ to an element $\omega(m)$ of $\Hom(\bigwedge^k T_mM ;\RR)$, smoothly, in the sense below.
On an open part $U$ of a manifold identified with an open subspace of $\RR^n$ by a chart $\phi \colon U \rightarrow \RR^n$ mapping $u \in U$ to $\phi(u)=(\phi_1(u),\dots,\phi_n(u))$, the degree $k$ forms are
uniquely expressed as \begin{equation*}\sum_{(i_1, \dots, i_k) \in \NN^k \suchthat 1 \leq i_1<  \dots < i_k\leq n}f_{i_1 \dots i_k} d \phi_{i_1} \wedge \dots \wedge d \phi_{i_k}.\end{equation*}
for smooth maps $f_{i_1 \dots i_k} \colon U \rightarrow  \RR$.
The vector space of degree $k$ differential forms on $M$ is denoted by $\Omega^k(M)$.
The differential $T\psi \colon TM \rightarrow TN$ of a smooth map $\psi$ from a manifold $M$ to a manifold $N$ induces
the \emph{pull-back} of differential forms \begin{equation*}\begin{array}{llll}\psi^{\ast} \colon & \Omega^k(N) & \rightarrow & \Omega^k(M)\\
&\omega &\mapsto &\bigl(x \mapsto \omega({\psi(x)}) \circ \bigwedge^k T_x\psi\bigr). \end{array}\end{equation*}
The antisymmetric bilinear exterior product $\wedge$ equips $\oplus_{k \in \NN}\Omega^k(M)$ with a structure of graded algebra, such that $\psi^{\ast}(\omega \wedge \omega^{\prime})=\psi^{\ast}(\omega) \wedge \psi^{\ast}(\omega^{\prime})$ for any two forms $\omega$ and $\omega^{\prime}$ on $N$.
This graded algebra is equipped with a unique operator $d \colon \Omega^k(M) \rightarrow \Omega^{k+1}(M)$ such that
\begin{itemize}
 \item $(df=Tf)$ is the differential of $f$ for any $f \in \Omega^0(M)$, 
 \item we have $d\circ d(f)=0$ for any $f \in \Omega^0(M)$,
\item for any $\alpha \in \Omega^{|\alpha|}(M)$ and $\beta \in \Omega^{|\beta|}(M)$, we have \begin{equation*}d(\alpha \wedge \beta)=d\alpha \wedge \beta + (-1)^{|\alpha|}\alpha \wedge d\beta ,\end{equation*} where $|\alpha|$ and $|\beta|$ denote the respective integral degrees of $\alpha$ and $\beta$, and
\item the derivation operator $d$ commutes with the above pull-backs, we have $d  \psi^{\ast}(\omega)=\psi^{\ast}(d\omega)$.
\end{itemize}
Then we have $d\circ d=0$.

The \emph{support}\index[T]{support!of a form} of a differentiable form is the closure of the set where it does not vanish.

Let $x_i \colon \RR^k \rightarrow \RR$ be the usual coordinate functions.
The \emph{integral} of a degree $k$ differential form $\omega=f dx_1 \wedge \dots \wedge dx_k$ over a $k$-dimensional compact submanifold $C$ of $\RR^k$ with boundary and ridges (like $\left[0,1\right]^k$) is $\int_C \omega= \int_C f dx_1  \dots dx_k$.
For any smooth map $\psi$ of $\RR^k$ that restricts to an orientation-preserving diffeomorphism from $C$ to its image, we have
\begin{equation*}\int_C\psi^{\ast}(\omega)=\int_{\psi(C)} \omega .\end{equation*}
Thanks to this property, we unambiguously define the integral of a $k$-form $\omega$ over any $k$-dimensional compact submanifold $C$ of a manifold $M$, identified to a subspace of $\RR^k$ by a diffeomorphism $\phi \colon C \rightarrow \RR^k$ onto its image, to be \begin{equation*}\int_C\omega=\int_{\phi(C)} \phi^{-1 \ast}(\omega).\end{equation*} This definition extends naturally to general compact manifolds with boundaries and ridges. 

One of the most useful theorems in this book is the following Stokes theorem. See \cite[Theorem 3.5, Page 31]{BT}.

\begin{theorem}[Stokes' theorem]
 Let $\omega$ be a degree $d$ form on an oriented  smooth compact $(d+1)$-manifold $M$ with boundary $\partial M$, then we have
\begin{equation*}\int_{\partial M} \omega = \int_M d\omega.\end{equation*}
\end{theorem}
This theorem applies to manifolds with ridges, and $\int_{\partial M} \omega$ is the sum of the $\int_{C} \omega$, over the codimension zero faces $C$ of $\partial M$. (The closures of these faces are $d$-manifolds with boundaries.)

\section{De Rham cohomology}
\label{secderhamcohom}

A degree $k$ differential form $\omega$ on $M$ is \emph{closed}\index[T]{closed!form} if $d\omega=0$. It is \emph{exact}\index[T]{exact!form} if $\omega \in d \Omega^{k-1}(M)$.
Define the degree $k$ de Rham cohomology\index[T]{de Rham cohomology} module of $M$ to be
\begin{equation*}H_{dR}^k(M)=\frac{\Ker\left(d \colon \Omega^{k}(M) \rightarrow \Omega^{k+1}(M)\right)}{d\bigl(\Omega^{k-1}(M)\bigr)}.\end{equation*}

For a compact submanifold $N$ (without boundary) of the interior of $M$, or for a connected component $N$ of $\partial M$, as in \cite[Chapter XII]{Godbillon}, let ${\Omega}^{k}(M,N)$ be the kernel of the restriction map from 
$\Omega^{k}(M)$ to $\Omega^{k}(N)$. The relative degree $k$ de Rham cohomology module of $(M,N)$ is
\begin{equation*}H_{dR}^k(M,N)=\frac{\Ker\left(d \colon {\Omega}^{k}(M,N) \rightarrow {\Omega}^{k+1}(M,N)\right)}{d\bigl({\Omega}^{k-1}(M,N)\bigr)}.\end{equation*} 
We have a natural short exact sequence of chain complexes
\begin{equation*} 0 \to {\Omega}^{k}(M,N) \to \Omega^{k}(M) \to \Omega^{k}(N) \to 0.\end{equation*}
This sequence
 induces a natural long exact cohomology sequence \begin{equation*}\to H_{dR}^{k-1}(N) \to H_{dR}^{k}(M,N) \to H_{dR}^{k}(M) \to H_{dR}^{k}(N)  \to \end{equation*}
where the maps from $H_{dR}^{k}(M,N)$ to $H_{dR}^{k}(M)$ and from $H_{dR}^{k}(M)$ to $H_{dR}^{k}(N)$ are induced by the restrictions and the cohomology class $[\eta]$ in $H_{dR}^{k-1}(N)$ of a closed form $\eta$ of $N$ is mapped to $\bigl([d\tilde{\eta}] \in H_{dR}^{k}(M,N)\bigr)$ for an extension $\tilde{\eta}$ of $\eta$ to $M$.

The degree $k$ forms with compact support in $M$ also form a subspace $\Omega_c^{k}(M)$ of $\Omega^{k}(M)$, 
and the degree $k$ de Rham cohomology module with compact support of $M$ is
\begin{equation*}H_{dR,c}^k(M)=\frac{\Ker\left(d \colon \Omega_c^{k}(M) \rightarrow \Omega_c^{k+1}(M)\right)}{d\bigl(\Omega_c^{k-1}(M)\bigr)}.\end{equation*}

For any smooth map $\psi$ from $M$ to another manifold $M^{\prime}$ that maps $N$ to a submanifold $N^{\prime}$, the pull-back $\psi^{\ast}$ induces maps still denoted by $\psi^{\ast}$ from $ H_{dR}^k(N^{\prime})$ to $H_{dR}^k(N)$, and from $H_{dR}^k(M^{\prime},N^{\prime})$ to $ H_{dR}^k(M,N)$. If $\phi$ is another such smooth map from the pair $(M^{\prime},N^{\prime})$ to another such $(M^{\prime \prime},N^{\prime \prime})$, then we have
\begin{equation*}\psi^{\ast} \circ \phi^{\ast}=(\phi \circ \psi)^{\ast}.\end{equation*}

When such a map $\psi$ is \emph{proper} (i.e., when the preimage of a compact is compact), the map $\psi$ also induces $\psi^{\ast} \colon H_{dR,c}^k(M^{\prime}) \rightarrow H_{dR,c}^k(M)$.

The following standard lemma implies that $\psi^{\ast} \colon H_{dR}^k(M^{\prime}) \rightarrow H_{dR}^k(M)$ depends only on the homotopy class of $\psi \colon M \to M^{\prime}$.
\begin{lemma}
\label{lemprimhom} Let $V$ and $W$ be two smooth manifolds, and
 let \begin{equation*}\begin{array}{llll}h\colon & \left[0,1\right] \times V &\rightarrow & W\\
       & (t,v) & \mapsto & h_t(v)
       \end{array}
\end{equation*} be a smooth homotopy. Let $\omega$ be a degree $d$ closed form on $W$.
Then we have \begin{equation*}h_t^{\ast}(\omega)-h_0^{\ast}(\omega)=d \eta_t(h,\omega)\end{equation*}
for any $t\in \left[0,1\right]$,
where $\eta_t(h,\omega)$ is the following degree $(d-1)$ form on $V$.
For $u\in \left[0,1\right]$, let $i_u \colon V \rightarrow \left[0,1\right] \times V$ map $v \in V$ to $i_u(v)=(u,v)$.
Let $h^{\ast}(\omega)\bigl((u,v)\bigr)(\frac{\partial}{\partial t}\wedge.)$ be obtained from $h^{\ast}(\omega)$ by evaluating it at the tangent vector $\frac{\partial}{\partial t}$ to $\left[0,1\right] \times \{v\}$ at $(u,v)$. (Thus, $h^{\ast}(\omega)\bigl((u,v)\bigr)(\frac{\partial}{\partial t}\wedge.)$ is the value at $(u,v) \in \left[0,1\right] \times V$ of a degree $(d-1)$ form of $\left[0,1\right] \times V$.) Then we define
\begin{equation*}\eta_t(h,\omega)(v)=\int_0^t i_u^{\ast}\left(h^{\ast}(\omega)\bigl((u,v)\bigr)\left(\frac{\partial}{\partial t}\wedge.\right)\right)du.\end{equation*}
\end{lemma}
\bp Observe $h_u^{\ast}(\omega)=i_u^{\ast}\left(h^{\ast}(\omega)\right)$ and write
\begin{equation*}h^{\ast}(\omega)=\omega_1 + dt \wedge h^{\ast}(\omega)\left(\frac{\partial}{\partial t}\wedge .\right),\end{equation*}
where we have $h_u^{\ast}(\omega)=i_u^{\ast}(\omega_1)$.
Since $\omega$ is closed, $d h^{\ast}(\omega)$ vanishes. Therefore, we obtain
\begin{equation*}0=d\omega_1 -dt \wedge d\Bigl(h^{\ast}(\omega)\bigl(\frac{\partial}{\partial t}\wedge .\bigr)\Bigr).\end{equation*}
On the other hand, with the natural projection $p_V \colon \left[0,1\right] \times V \rightarrow V$, we can write \begin{equation*}(d\omega_1)\bigl((u,v)\bigr)=\Bigl(dt \wedge \frac{\partial}{\partial t} (\omega_1) + p_V^{\ast}\left(d h_u^{\ast}(\omega)\right)\Bigr)\bigl((u,v)\bigr).\end{equation*}
We get \begin{equation*}i_u^{\ast}\left(\frac{\partial}{\partial t} (\omega_1)\right)=i_u^{\ast}\left(d\left(h^{\ast}(\omega)\left(\frac{\partial}{\partial t}\wedge .\right)\right)\right)\end{equation*}
and
\begin{equation*}\frac{\partial}{\partial u} h_u^{\ast}(\omega) =\frac{\partial}{\partial u}i_u^{\ast}\left(\omega_1\right)= i_u^{\ast}\left(\frac{\partial}{\partial t} (\omega_1)\right)=i_u^{\ast}\left(d\left(h^{\ast}(\omega)\left(\frac{\partial}{\partial t}\wedge.\right)\right)\right).\end{equation*}
Since we have $h_t^{\ast}(\omega)(v)-h_0^{\ast}(\omega)(v)=\int_0^t \frac{\partial}{\partial u} h_u^{\ast}(\omega)(v) du$, the lemma follows.
\eop

In particular, if there exist smooth maps $f \colon M \rightarrow N$ and $g \colon N \rightarrow M$ such that $f \circ g$ is smoothly homotopic to the identity map of $N$, and $g \circ f$ is smoothly homotopic to the identity map of $M$, then $f^{\ast}$ is an isomorphism from 
$H_{dR}^k(N)$ to $H_{dR}^k(M)$ for any $k$. In particular, all the smoothly contractible manifolds,
such as $\RR^n$, have the same de Rham cohomology as the point. So, for such a manifold $C$,
$H_{dR}^k(C)=\{0\}$ for any $k\neq 0$, and $H_{dR}^0(C)=\RR$.

More generally, in 1931, Georges de Rham identified the de Rham cohomology of a smooth $n$-dimensional manifold $M$ to its singular cohomology with coefficients in $\RR$.
See \cite[Pages 205-207]{Warnerseconded}, or \cite[Theorems 8.9 Page 98, 15.8 page 191]{BT}.

The de Rham isomorphism sends the cohomology class $[\omega]$ of a closed degree $k$ form to a linear map $[\omega]_{dR}$ of $\Hom(H_k(M;\RR);\RR)=H^k(M;\RR)$. The linear map
$[\omega]_{dR}$ sends the homology class of a closed oriented $k$-dimensional submanifold $N$ of $M$ to $[\omega]_{dR}([N])=\int_N\omega$. 

If $M$ has a smooth triangulation $T$, any $k$-form $\omega$ similarly defines a simplicial cochain, i.e., an element of $C^k(T;\RR)=\Hom(C_k(T;\RR);\RR))$,
by integration over the simplices. Stokes' theorem guarantees that the induced map from $\Omega^k(M)$ to $C^k(T;\RR)$ commutes with differentials, and induces a morphism from $H_{dR}^k(M)$ to $H^k(M;\RR)=\Hom(H_k(M;\RR);\RR)$.

We also have the following theorem \cite[page 197]{GreubHalperinVanston} or \cite[(5.4) and Remark 5.7]{BT}.

\begin{theorem}
\label{thmDeRhamcompact}
For any oriented manifold $M$ without boundary of dimension $n$ (whose cohomology is not necessarily finite-dimensional), the morphism from 
$\Omega^k(M;\RR)$ to $\Hom(\Omega^{n-k}_c(M);\RR)$\begin{equation*}\left(\omega \mapsto \Bigl(\omega_2 \mapsto \int_{M}\omega \wedge \omega_2\Bigr)\right)\end{equation*} induces an isomorphism from 
$H^k(M;\RR)$ to $\Hom(H^{n-k}_c(M);\RR)$,
for any integer $k$.
\end{theorem}

In particular, when the real cohomology of $M$ is finite-dimensional, so is its homology, and $H^{n-k}_c(M)$ is isomorphic to $H_k(M;\RR)$.
Below, we exhibit the image $[\omega_A]$ of the homology class of an oriented compact submanifold $A$ with trivial normal bundle in $H_k(M;\RR)=\Hom(H^{k}(M);\RR)$ under a canonical isomorphism from $H_k(M;\RR)$ to $H^{n-k}_c(M)$.

Let $A$ be such a $k$-dimensional submanifold without boundary of the manifold $M$, and let $N(A)=B^{n-k} \times A $ be a tubular neighborhood of $A$ in $M$. Let $\omega_B$ be an $(n-k)$-form of $B^{n-k}$ supported in the interior of $B^{n-k}$ such that $\int_{B^{n-k}} \omega_B=1$. Use the natural projection $p_B \colon B^{n-k} \times A \to B^{n-k}$ to pull back $\omega_B$ on $N(A)$ and define a form $\omega_A$ to coincide with $p_B^{\ast}(\omega_B)$ on $N(A)$ and to be zero outside $N(A)$.
The cohomology class $[\omega_A]$ of $\omega_A$ is independent of the choice of $\omega_B$.
The integral of the closed form $\omega_A$ over a compact $(n-k)$ submanifold $B$ of $M$ transverse to $A$ is the algebraic intersection $\langle B,A\rangle_{\!M\,}$ of $B$ and $A$.
Note that the support of $\omega_A$ may be chosen in an arbitrarily small neighborhood of $A$.
In the words of Definition~\ref{defformdual}, the form $\omega_A$ may be chosen to be $\alpha$-dual to $A$ for an arbitrarily small positive number $\alpha$.

\begin{lemma}
\label{lemconstrformalphadual} Let $M$ be a compact manifold of dimension $n$ with possible boundary.
Assume that $M$ is equipped with a smooth triangulation.
The above process can be extended to produce canonical \emph{Poincar\'e duality isomorphisms} from $H_{a}(M,\partial M;\RR)$ to $H_{dR}^{n-a}(M)$, where 
such a Poincar\'e duality isomorphism maps the class of an $a$-dimensional cycle $A$ of $(M,\partial M)$ to the class of a closed $(n-a)$-form $\omega_A$ $\alpha$-dual to $A$ for an arbitrarily small positive number $\alpha$, as follows.
\end{lemma}
\bp
Let $A$ be a linear combination of $a$-dimensional simplices of the triangulation $T$ of $M$ such that the algebraic boundary of $A$ is in $C(T_{\partial M})$. Such an $A$ is 
a \emph{simplicial $a$-cycle} of $(M,\partial M)$.
Its \emph{support} $\overline{A}$ is the union in $M$ of the closed simplices with a nonzero coefficient of $A$.
We construct a form $\omega_A$ $\alpha$-dual to $A$ as follows.
Let $A^{(k)}$ be the intersection of $\overline{A}$ with the $k$-skeleton of $T$. ($A^{(k)} = \emptyset$ if $k< 0$.)
Let $N(A^{(a-1)})$ be a small neighborhood of $A^{(a-1)}$.

First construct a closed form $\omega_A$ $\alpha$-dual to $A$ outside
$N(A^{(a-1)})$ so that its support is in an arbitrarily small neighborhood of $A$, as explained above in the case of closed manifolds with trivial normal bundles.

Then extend $\omega_A$ around each $(a-1)$-simplex $\Delta$ of $A$, outside a small neighborhood $N(A^{(a-2)})$ of $A^{(a-2)}$ as follows.

Without loss of generality, assume that the intersection of a neighborhood of $\Delta$ with the complement of $N(A^{(a-2)})$ is diffeomorphic to $\Delta^{\prime} \times D^{n-a+1}$, for some interior $\Delta^{\prime}$ of an $(a-1)$-simplex in the interior of $\Delta$ (some \emph{truncation} of $\Delta$). See Figure~\ref{figneighformdual}. The $(n-a)$-form $\omega_A$ is defined and closed on a neighborhood of $\Delta^{\prime} \times \partial D^{n-a+1}$. Its integral along 
$\{x\} \times \partial D^{n-a+1}$ is the algebraic intersection $\langle \partial D^{n-a+1},A\rangle_{\!M\,}$. It is the coefficient of $\Delta$ in $\partial A$, up to sign. Since this coefficient is zero, the form $\omega_A$ is exact in a neighborhood of $\Delta^{\prime} \times \partial D^{n-a+1}$. So $\omega_A$ may be written as
$d \eta$ on this neighborhood. Choose a map $\chi$ on $\RR^{n-a+1}$ that takes the value $1$ outside a small neighborhood of $0$ and that vanishes in a smaller neighborhood of $0$ so that $d \chi \eta$ makes sense on $\Delta^{\prime} \times D^{n-a+1}$ and extends $\omega_A$ as a closed form.

\bfig
\centering
\begin{tikzpicture}[scale=.5] \useasboundingbox (-5,-1.2) rectangle (5,2);
\begin{scope}[xshift=-2.5cm]
\fill [lightgray] (-1.8,-1) rectangle (1.8,1.9);

\draw [thin] (-30:1.4) -- (90:1.4) -- (-150:1.4) -- (-30:1.4) (-1.8,-.7) -- (1.8,-.7) 
(90:1.4) -- (90:1.9);
\draw [line width=2pt] (-150:1.4) -- (-30:1.4);
\fill [white, rounded corners] (-1.9,2) -- (-.3,2) -- (-.3,1.4) -- (-1.32,-.4) -- (-1.9,-.4) -- (-1.9,2);
\fill [white, rounded corners] (1.9,2) -- (.3,2) -- (.3,1.4) -- (1.32,-.4) -- (1.9,-.4) -- (1.9,2);
\fill [white, rounded corners] (-30:1.1) -- (90:1.1) -- (-150:1.1) -- (-30:1.1);
\fill [white, rounded corners] (90:1.1) -- (-150:1.1) -- (-30:1.1) -- (90:1.1);
\draw [very thin, <-] (.25,1.6) -- (.6,1.6) node[right]{\scriptsize $N(A^{(1)})$};
\draw [very thin, <-] (0,-.63) -- (0,-.3) node[above]{\scriptsize $\Delta$};
\draw [very thin, <-] (-1.65,-.65) -- (-1.65,-.3) node[above]{\scriptsize $A$};
\draw [very thin, <-] (-.97,-.1) -- (-1.45,-.1);
\end{scope}
\begin{scope}[xshift=2.5cm]
\fill [lightgray] (-1.8,-1) rectangle (1.8,1.9);
\draw [line width=2pt] (-.9,-.7) -- (.9,-.7);
\fill [white, rounded corners] (-1.9,2) -- (-.3,2) -- (-.3,1.4) -- (-1.32,-.4) -- (-1.9,-.4) -- (-1.9,2);
\fill [white, rounded corners] (1.9,2) -- (.3,2) -- (.3,1.4) -- (1.32,-.4) -- (1.9,-.4) -- (1.9,2);
\fill [white, rounded corners] (-30:1.1) -- (90:1.1) -- (-150:1.1) -- (-30:1.1);
\fill [gray] (90:1.4) circle (.45);
\fill [gray] (-30:1.4) circle (.45);
\fill [gray] (-150:1.4) circle (.45);
\draw [thin] (-30:1.4) -- (90:1.4) -- (-150:1.4) -- (-30:1.4) (-1.8,-.7) -- (1.8,-.7) 
(90:1.4) -- (90:1.9);
\draw [very thin, <-] (.25,1.45) -- (.6,1.6) node[right]{\scriptsize $N(A^{(0)})$};
\draw [very thin, <-] (1.2,-.35) -- (1.2,1.3);
\draw [very thin, <-] (0,-.63) -- (0,-.3) node[above]{\scriptsize $\Delta^{\prime}$};
\end{scope}
\end{tikzpicture}
\caption{The neighborhoods $N(A^{(1)})$ and $N(A^{(0)})$, a simplex $\Delta$ and its truncation $\Delta^{\prime}$, when $A$ is one-dimensional}\label{figneighformdual}
\end{figure}

This process allows us to define a closed form $\omega_A$ $\alpha$-dual to $A$ outside a small neighborhood $N(A^{(a-2)})$ of $A^{(a-2)}$. 
Iterate the process in order to extend such a form outside $N(A^{(a-3)})$, outside $N(A^{(a-4)})$, \dots,
and on the whole $M$. Note that the forms will be automatically exact around the truncated smaller cells since we have $H^{n-a}(\partial D^{n-a+k})=0$ when $k>1$.
\eop

In this setting, the correspondence between chain boundaries and the differentiation operator $d$ can be roughly seen as follows. Let $A$ be a compact oriented $a$-dimensional submanifold  with boundary of a manifold $M$. Assume that the normal bundle of $A$ is trivial. Let $[-1,1] \times \partial A$ be embedded in the interior of $M$ so that $\{0\} \times \partial A=\partial A$ and $\left[-1,0\right] \times \partial A$ is a neighborhood of $\partial A$ in $A$. Set $A^+=A \cup_{\left[-1,0\right] \times \partial A} [-1,1] \times \partial A$. Choose a map $\chi \colon A^+ \to [0,1]$ sending $A$ to $1$ and a neighborhood of $\partial A^+$ to $0$ such that $\chi$ factors through the projection onto $\left[0,1\right]$ on $\left[0,1\right] \times \partial A$. Let $N(A^+)=B^{n-a} \times A^+$ be the total space of the normal bundle to $A^+ =A \cup_{\partial A} \left(\left[0,1\right] \times \partial A\right)$
embedded in $M$. Let $p_{A^+}\colon N(A^+) \to A^+$ and $p_B\colon N(A^+) \to B^{n-a}$ denote the natural projections. Let $\omega_B$ be an $(n-a)$-form of $B^{n-a}$ supported in the interior of $B^{n-a}$ such that $\int_{B^{n-a}} \omega_B=1$. The form \begin{equation*}\omega_A=(\chi \circ p_{A^+}) p_B^{\ast}(\omega_B)\end{equation*} extends as a smooth form that vanishes outside $N(A^+)$. The form $d\omega_A$ is supported in $B^{n-a} \times [0,1] \times \partial A$. Its integral along a chain $C$ whose boundary does not meet $B^{n-a} \times [0,1] \times \partial A$ is $\pm \langle C , \partial A\rangle_{M}$.

\def\cprime{$'$}
\providecommand{\bysame}{\leavevmode ---\ }
\providecommand{\og}{``}
\providecommand{\fg}{''}
\providecommand{\smfandname}{\&}
\providecommand{\smfedsname}{\'eds.}
\providecommand{\smfedname}{\'ed.}
\providecommand{\smfmastersthesisname}{M\'emoire}
\providecommand{\smfphdthesisname}{Th\`ese}

\newpage

\printindex[T]
\printindex[N]

\chapter*{Summarizing the main definitions of $\Zinv$}
\label{SummaryZ}
\addcontentsline{toc}{chapter}{Summarizing the definitions of  \texorpdfstring{$\Zinv$}{Z}}

In the informal summary below, we first recall the various definitions of $\Zinv$ and its variants ($\Zinvlink$, $\zinv$, $\Zinvufrf$) for a given link embedding $\Link = \sqcup_{j=1}^k K_j$ into a $\QQ$-sphere $\rats$. 
The punctured $\crats =\rats \setminus \{\infty\}$ is equipped with a parallelization $\tau$ that makes $(\crats,\tau)$ an asymptotic rational homology $\RR^3$ as in Definition~\ref{defasyrathommRthree}.

For $n \in \NN$, for a family $(\omega(i))_{i \in \underline{3n}})$ of  propagating forms of $C_2(\rats)$ as in Section~\ref{secprop}, we have
\begin{equation*}\Zinv_n\Bigl(\crats,\Link,\bigl(\omega(i)\bigr)_{i \in \underline{3n}}\Bigr)\stackrel{\mbox{\scriptsize Not. \ref{notationzZ}}}{=}\sum_{\Gamma \in \Davis^e_n(\source)}\coefgambet_{\Gamma}I\Bigl(\rats,\Link,\Gamma,\bigl(\omega(i)\bigr)_{i \in \underline{3n}}\Bigr)\left[\Gamma\right] \in \Aavis_n(\source).\footnote{When $\omega(i)=\omega$ for all $i$, we set $\Zinv_n(\crats,\Link,\omega)=\Zinv_n(\crats,\Link,(\omega(i))_{i \in \underline{3n}})$ as in Proposition~\ref{propdefhomog}. The use of a family of distinct $\omega(i)$ allows us to discretize the definition of $\Zinvuf$ in Chapter~\ref{chaprat}, thanks to Version~\ref{thmconststraight} of Theorem~\ref{thmfstconst}.}\end{equation*}
If the $\omega(i)$ are propagating forms of $(C_2(\rats),\tau)$ as in Definition~\ref{defpropagatorone}, and if they are homogeneous as in Definition~\ref{defprophomogen}, then we have
\begin{equation*}\Zinv_n(\crats,\Link,\tau)\stackrel{\mbox{\scriptsize Thm. \ref{thmfstconst}}}{=} \Zinv_n\Bigl(\crats,\Link,\bigl(\omega(i)\bigr)_{i \in \underline{3n}}\Bigr)\end{equation*} and
\begin{equation*}\Zinv(\crats,\Link,\tau)=\bigl(\Zinv_n(\crats,\Link,\tau) \bigr)_{n\in \NN}=\sum_{n\in \NN}
\Zinv_n(\crats,\Link,\tau).\end{equation*}
With the anomalies $\alpha$ of Section~\ref{secanomalpha} and $\beta$ of Section~\ref{secansothree}, the Pontragin number $p_1(\partau)$ of Definition~\ref{deffirstpontcrats}, and the definition of $I_{\theta}$ in Lemma~\ref{lemdefItheta}, we have 
\begin{equation*}\Zinvuf(\rats,\Link)\stackrel{\mbox{\scriptsize Thm. \ref{thmfstconst}}}{=}\exp\Bigl(-\frac14 p_1(\partau)\ansothree\Bigr)\prod_{j=1}^k\Bigl(\exp\bigl(-I_{\theta}(K_j,\partau)\alpha\bigr)\#_j\Bigr) \Zinv(\crats,\Link,\partau),\footnote{The similar letters $\Zinv$ and $\Zinvuf$ respectively denote a function involving auxiliary data and the induced invariant.}\end{equation*}
and \begin{equation*}\Zinvuf(\rats)=\Zinvuf(\rats,\emptyset).\end{equation*}
With the projection $p^c \colon \Aavis(\emptyset) \to \Aavis^c(\emptyset)$ on the connected part of Notation~\ref{notationzZ}, we have \begin{equation*}\zinv(\rats) = p^c\bigl(\Zinvuf(\rats)\bigr)\;\;\mbox{and}\;\;\Zinvuf(\rats)=\exp\bigl(\zinv(\rats)\bigr).\end{equation*}
With the projection $\projassis \colon \Aavis(\source) \to \Assis(\source)$ of Notation~\ref{notationzZ}, which forgets diagrams without univalent vertices,
we have \begin{equation*}\Zinvlinkuf = \projassis \circ \Zinvuf\;\;\mbox{and}\;\;
\Zinvuf(\rats,\Link)=\Zinvuf(\rats)\Zinvlinkuf(\rats,\Link).\end{equation*}
When $\Link = \sqcup_{j=1}^k K_j$ has a parallel $\Link_{\parallel} = \sqcup_{j=1}^k K_{j\parallel}$, the \emph{framed} version $\Zinvufrf$ of $\Zinvuf$ satisfies
\begin{equation*}\Zinvufrf(\crats,\Link,\Link_{\parallel})
\stackrel{\mbox{\scriptsize Def. \ref{defZinvufrflink}}}{=}
\prod_{j=1}^k\Bigl(\exp\bigl(lk(K_j,K_{j\parallel})\alpha\bigr)\#_j\Bigr) \Zinvuf(\rats,\Link).\end{equation*}

Now, we give a similar summary for the extension of $\Zinvuf$ to tangles in $\QQ$-cylinders of the book's third part. We fix such a $\QQ$-cylinder $\hcylc$, the associated $\QQ$-sphere $\rats(\hcylc)$, a tangle $\tanghcyll$ in $\hcylc$, and the associated long tangle also denoted by $\tanghcyll$ in $\crats(\hcylc)$ as in Definition~\ref{defLTR}.
Theorem~\ref{thmfstconsttang} defines $\Zinv(\hcylc,\tanghcyll,\tau)$ (resp. $\Zinv_n(\hcylc,\tanghcyll,(\omega(i))_{i \in \underline{3n}})$)  to be a natural extension of $\Zinv(\rats(\hcylc),\tanghcyll,\tau)$ (resp. $\Zinv_n(\rats(\hcylc),\tanghcyll,(\omega(i))_{i \in \underline{3n}})$) from the case where $\tanghcyll$ is a link in $\hcylc$ to the case where $\tanghcyll$ is any tangle of 
$\hcylc$ canonically extended to a long tangle as in Definition~\ref{defLTR}. We have
\begin{equation*}\Zinvuf(\hcylc,\tanghcyll)=\exp\Bigl(-\frac14 p_1(\tau)\ansothree\Bigr)\prod_{j=1}^k\Bigl(\exp\bigl(-I_{\theta}(K_j,\tau)\alpha\bigr)\#_j\Bigr) \Zinv(\hcylc,\tanghcyll,\tau),\end{equation*}
with respect to Definition~\ref{defIthetalong} for $I_{\theta}(K_j,\tau)$ when $K_j$ is a long component.  Definition~\ref{deffstconsttangframed} of
\begin{equation*}\Zinvufrf\bigl(\hcylc,(\tanghcyll,\tanghcyll_{\parallel})\bigr)=\prod_{j=1}^k\Bigl(\exp\bigl(lk\left(K_j,K_{j\parallel}\right)\alpha\bigr)\#_j\Bigr)\Zinvuf(\hcylc,\tanghcyll).\end{equation*} extends from tangles to $q$-tangles, which are cobordisms between limit configurations, by a natural limit process as in Remark~\ref{deflimzqtang}.

Again, discrete computations and proofs of properties involve \emph{straight} tangles as in Section~\ref{secstraighttang}
 and distinct propagating forms \emph{dual} to propagating chains as in Chapter~\ref{chaprat}. Theorems~\ref{thmtangconstcompar} and \ref{thmdefsanstauvarzinvf}, and Definition~\ref{defvarzfqtang} describe variants of $\Zinvuf$ associated to such data. These variants depend on 2-forms over $S^2$ in a way described in Theorem~\ref{thmtangconstcomparbis}.
 \end{document}